\def\includeall{\renewcommand{\includeonly}[1]{}}
	   \includeall

\documentclass[reqno]{amsbook} 

\usepackage{Intro}

\begin{document}

%%% set up shortrefs for Exercises 
%%% need to do after \begin{document} for compatibility with hyperref
\makeatletter
 \global\let\precref@getlabel=\cref@getlabel
 \gdef\cref@getlabel#1#2{% 
\@ifundefined{r@#1}{\precref@getlabel{#1}{#2}}%
  {\ifthenelse{\equal{\thesection}{\sectionref{#1}}}%
   {\precref@getlabel{short@#1}{#2}}{\precref@getlabel{#1}{#2}}}}
 \newcommand{\oldref}{} \global\let\oldref\ref
 \renewcommand{\ref}[1]{%
 	\@ifundefined{section@#1}{\oldref{#1}}%
		{\ifthenelse{\equal{\thesection}{\sectionref{#1}}}%
 			{\oldref{short@#1}}{\oldref{#1}}}}
 \makeatother

\frontmatter

%!TEX root = IntroArithGrps.tex

\begingroup

\pagestyle{empty}

%\hrule \medskip 

\pdfbookmark[-1]{Frontmatter}{frontmatter}
\pdfbookmark[0]{Title page}{titlepage}

\ 
\vskip 0.75in

%\leftskip=0.575in
\thispagestyle{empty}

\hbox
{%\ifportable\hskip-0.4 in\fi
\vbox
{{\input titlefont.tex%% make the title really huge
\bfseries%\sffamily
%\global\setbox0\hbox{A}
%\leftskip = \wd0
\global\dimen0=\wd0
{\noindent \hbox{\hskip 0.365in} \hskip 7pt
\sffamily\small Introduction to}
\mathversion{bold}
 {\Large
\\[0pt] \hbox{\hskip 0.365in} 
\llap{A\kern0.7pt}\hbox{R\kern1pt IT\kern0pt H\kern0pt M\kern0pt E\kern-1pt T\kern0pt I\kern-3pt C}
\dimen0=0pt
\\[0pt] \hbox{\hskip 0.365in} 
\llap{G\kern-2pt}R\kern-5pt O\kern-3pt U\kern0.5pt P\kern-1.5pt S 
% I don't like the look of Lucida's sans-serif "G" , so use serif for the title
}
\vskip 3in \sffamily
\noindent  \hbox{\hskip 0.365in} \hskip 7pt
\small Dave Witte Morris
\par
} % end of large font
}}

\vfil

\endgroup
\newpage

\begingroup
\thispagestyle{empty}

\ 

\vskip 1in

\begin{center}
 
 Dave Witte Morris
% \\ Department of Mathematics and Computer Science
 \\ University of Lethbridge
% \\ Lethbridge, Alberta, T1K~3M4, Canada
 
 \vskip 0.125in
 
\textsf{\href{http://people.uleth.ca/~dave.morris/}{http://people.uleth.ca/\!{\mathversion{bold}\larger$\sim$}dave.morris/}}

\vfil

\normalfont\sffamily\noindent
 \hskip 5pt \ifportable Arxiv \fi Version 1.0\ of April 2015  % !!!
% \\ \hbox{\hskip 5pt}Please send comments and corrections to 
%\\  \hbox{\hskip 25pt} \href{mailto:Dave.Morris@uleth.ca}{\texttt{Dave.Morris@uleth.ca}}
	% update the date and version number !!!

\vfil

\ifportable\else
published by

\smallskip

{\larger Deductive Press}
\\ \url{http://deductivepress.ca/}

\bigskip

\setbox0\hbox{(paperback)}

ISBN: 978-0-9865716-0-2 \hbox to \wd0{(paperback)\hss}
\\[\smallskipamount] ISBN: 978-0-9865716-1-9 \hbox to \wd0{(hardcover)\hss}

\fi % \ifportable

\vfil

\centerline{%
	\boxit{\raise0.2pt\hbox{\href{http://creativecommons.org/publicdomain/zero/1.0/}%
	 		{\ifportable\hbox{\includegraphics[scale=0.225]{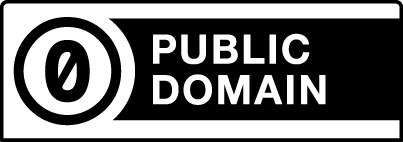}}\else
			 \hbox{\includegraphics[scale=1.2]{PDF/cc-zero-bw.pdf}}\fi}}%
	\vbox{\hbox{\ 
 	To the extent possible under law, Dave Witte Morris %
	}\hbox{\ 
	has waived all copyright and related or neighboring }%
	\hbox{\ 
	rights to this work.}%
	}\hskip-3pt}\hskip3.5pt} % get the box properly centered !!!

\medskip

{\it
You can copy, modify, and distribute this work, even for commercial
\\ purposes, all without asking permission. For more information, see}
	\\ \url{http://creativecommons.org/publicdomain/zero/1.0/}

\bigskip\bigskip

\ifportable \else

A PDF file of this book is available for download from
\\
	\url{http://arxiv.org/src/math/0106063/anc/}
\\[5pt]
The Latex source files are also available from
\\
\url{http://arxiv.org/abs/math/0106063}

\bigskip\bigskip
{\smaller
Cover photo \copyright\ \href{http://Depositphotos.com/portfolio-1987283.html}{Depositphotos.com/jbryson}. Used with permission.} % @@@

\fi % \ifportable

\end{center}
 
% \bigskip

\eject

\thispagestyle{empty}

\ \pdfbookmark[0]{Dedication}{dedication}
\vfil

\begin{center}
 \itshape
\larger

to my father
\\[5pt] E.\ Kendall Witte
\\ \smaller 1926 --- 2013

\end{center}

 \vfil

% \newpage
\endgroup

%!TEX root = IntroArithGrps.tex

%% TABLE OF CONTENTS

%%% manually choose better page breaks in the TOC
\newif\iftoceject \tocejectfalse
\newcommand{\toceject}{\ifportable\else\iftoceject\eject\fi\fi}

\cleardoublepage
\setcounter{page}{1}

\makeatletter

\newcommand{\tocchapterhead}[1]{\topskip 0pc\relax
  \begingroup
  \fontsize{\@xivpt}{18}\bfseries\sffamily\centering
  #1\par \endgroup
  \skip@34\p@ \advance\skip@-\normalbaselineskip
  \vskip\skip@ }
  
  \newcommand{\tocchapterstyle}{\bfseries\mathversion{bold}}

 \renewcommand{\l@part}[2]{\vskip0.25in \centerline{\bf \sffamily#1\hss}}

\newlength{\tocchaplabelwidth}

\renewcommand{\tocchapter}[3]{\sffamily\@ifnotempty{#2}{\tocchapterstyle}%
 \setbox0\hbox{#1 99.\enspace}\setlength{\tocchaplabelwidth}{\wd0}
  \indentlabel{\hbox to \wd0{\@ifnotempty{#2}{\ignorespaces#1 #2.}\hss}}#3}
 
\renewcommand{\tocappendix}[3]{\tocchapter{#1}{#2}{#3}}

\newlength{\intertoccorrection}

\makeatother

% list of chapters needs less vertical space, so it will fit on one page
\setlength{\intertoccorrection}{-5.275pt}

\begingroup %% Print a list of the numbered chapters (not frontmatter or indexes)
\pdfbookmark[0]{List of Chapters}{chapterlist}

\def\\{ \ignorespaces} % don't want line breaks in the chapter titles
\makeatletter \@fileswfalse \makeatother %% don't erase the TOC file
\renewcommand{\contentsname}{List of Chapters}
 \makeatletter
\renewcommand{\tocchapterstyle}{} % don't need chapter names to be bold
\renewcommand{\@afterheading}{\vskip-0.3in} % less space after "List of Chapters" heading
\renewcommand{\l@section}[2]{} % don't print the sections in this toc
\renewcommand{\l@schapter}[2]{} % don't print unnumbered chapters in this toc
\renewcommand{\@makeschapterhead}{\tocchapterhead}
 \makeatother
 \tableofcontents
\bigskip
 \endgroup % resume listing sections in the toc

 \makeatletter \@fileswtrue \makeatother %% now can erase the TOC file

\tocejecttrue % we do force page breaks in the full Table of Contents

% table of contents can use default vertical spacing
\setlength{\intertoccorrection}{0pt}

\begingroup
\def\\{ \ignorespaces} % don't want line breaks in the chapter titles

\makeatletter
\newlength{\tocsectlabelwidth}
\setbox0\hbox{\S99.9.\hskip0.3em}\global\setlength{\tocsectlabelwidth}{\wd0}%

\renewcommand{\tocchapter}[3]{%
\setbox0\hbox{#2M.\hskip0.3em}\global\setlength{\tocsectlabelwidth}{\wd0}%
\sffamily\@ifnotempty{#2}{\tocchapterstyle}%
  \indentlabel{\@ifnotempty{#2}{\ignorespaces#1 #2.\enspace}}#3}

\renewcommand{\tocsection}[3]{\normalfont
\setbox0\hbox{\S#2}
  \indentlabel{\hbox to \tocsectlabelwidth{\@ifnotempty{#2}{\hss\ignorespaces#1 
  \ifdim\wd0<15pt \hphantom{0}\fi\ifnotappendixtoc\S\fi#2.\ }}}#3}

\notappendixtoctrue % want \S on section lines, unless in appendix

\cleardoublepage
\thispagestyle{empty}
\refstepcounter{section}% to correct the PDF bookmarks
 \chaptermark{Contents}
\pdfbookmark[0]{Table of Contents}{toc}%
\@makeschapterhead{Contents}

\@input{\jobname.toc}
    \if@filesw
      \newwrite\tf@toc
      \immediate\openout\tf@toc \jobname.toc\relax
      \fi
\makeatother
\cleardoublepage
\endgroup

 %!TEX root = IntroArithGrps.tex

\refstepcounter{section}% to correct the PDF bookmarks

\chapter*{Preface}
\label{Preface}

This is a book about arithmetic subgroups of semisimple Lie groups, which means that we will discuss the group $\SL(n,\integer)$, and certain of its subgroups. By definition, the subject matter combines algebra (groups of matrices) with number theory (properties of the integers). However, it also has important applications in geometry. In particular, arithmetic groups arise in classical differential geometry as the fundamental groups of locally symmetric spaces. (See \cref{WhatisLocSymmChap,GeomIntroRank} for an elaboration of this line of motivation.) They also provide important examples and test cases in geometric group theory.

My intention in this text is to give a fairly gentle introduction to several of the main methods and theorems in the subject. There is no attempt to be encyclopedic, and proofs are usually only sketched, or only carried out for an illustrative special case. 
Readers with sufficient background will learn much more from \cite{MargulisBook} and \cite{PlatonovRapinchukBook} (written by the masters) than they can find here.

%The goal is to to acquaint graduate students and other non-experts with the background to understand seminars and papers that make use of some of the basic concepts, and perhaps even make use of some of these ideas themselves. Anyone needing a deeper understanding of the topics touched on here should consult the masters, but I hope this work will make those books less forbidding.
 
The book assumes knowledge of algebra, analysis, and topology that might be taught in a first-year graduate course, plus some acquaintance with Lie groups. (\cref{SSGrpsChap} quickly recounts the essential Lie theory, and \cref{BackChap} lists the required facts from graduate courses.) Some individual proofs and examples assume additional background (but may be skipped).

Generally speaking, the chapters are fairly independent of each other (and they all have their own bibliographies), so there is no need to read the book linearly.
To facilitate making a plan of study, the bottom of each chapter's first page states the main prerequisites that are not in \cref{SSGrpsChap,BackChap}.
Individual chapters (or, sometimes, sections) could be assigned for reading in a course or presented in a seminar. 
(The book has been released into the public domain, so feel free to make copies for such purposes.) Notes at the end of each chapter have suggestions for further reading. (Many of the subjects have been given book-length treatments.)
 Several topics (such as amenability and Kazhdan's property~$(T)$) are of interest well beyond the theory of arithmetic groups. 

Although this is a long book, some very important topics have been omitted. In particular, there is almost no discussion of the cohomology of arithmetic groups, even though it is a subject with a long history and continues to be a very active field. (See the lecture notes of Borel \cite{Borel-IntroCoho} for a recent survey.)
Also, there is no mention at all of automorphic forms. (Recent introductions to this subject include \cite{Deitmar-AutoForms} and \cite{PCMI-AutoForms}.)

Among the other books on arithmetic groups, the authoritative monographs of Margulis \cite{MargulisBook} and Platonov-Rapinchuk \cite{PlatonovRapinchukBook} have already been mentioned. They are essential references, but would be difficult reading for my intended audience. 
Some works at a level more comparable to this book include:
	\begin{itemize} \itemsep=\medskipamount
	\medskip

	\item[\cite{Borel-IntroGrpArith}] This classic gives an explanation of reduction theory (discussed here in \cref{ReductionChap}) and some of its important consequences.
	
	\item[\cite{Humphreys-ArithGrps}] This exposition covers reduction theory (at a more elementary level than \cite{Borel-IntroGrpArith}), adeles, ideles, and fundamentals of the Congruence Subgroup Property (mentioned here in \fullcref{MargNormSubgrpRems}{CSP}).
	
	\item[\cite{Ji-ArithGrpsWhatWhyHow}] This extensive survey touches on many more topics than are covered here (or even in \cite{MargulisBook} and \cite{PlatonovRapinchukBook}), with 60 pages of references.

	\item[\cite{MaclachlanReid-ArithHyp3Mflds}] This monograph thoroughly discusses arithmetic subgroups of the groups $\SL(2,\real)$ and $\SL(2,\complex)$.

	\item[\cite{RaghunathanBook}] This is an essential reference (along with \cite{MargulisBook} and \cite{PlatonovRapinchukBook}). It is the standard reference for basic properties of lattices in Lie groups (covered here in \cref{BasicLatticesChap}). It also has proofs of the Godement Criterion (discussed here in \cref{GodementSect}), the existence of both cocompact and noncocompact arithmetic subgroups (discussed here in \cref{CocptSect}), and reduction theory for arithmetic groups of $\rational$-rank one (discussed here in \cref{ReductionChap}). It also includes several topics not covered here, such as cohomology vanishing theorems, and lattices in non-semisimple Lie groups.

	\item[\cite{Sury-CSP}] This textbook provides an elementary introduction to the Congruence Subgroup Property (which has only a brief mention here in \fullcref{MargNormSubgrpRems}{CSP}).

	\item[\cite{ZimmerBook}] After developing the necessary prerequisites in ergodic theory and representation theory, this monograph provides proofs of three major theorems of Margulis: Superrigidity, Arithmeticity, and Normal Subgroups (discussed here in \cref{MargulisSuperChap,NormalSubgroupChap}). It also proves a generalization of the superrigidity theorem that applies to ``Borel cocycles\zz.''
	
	\end{itemize}

\medskip

\hfill \textit{Dave Morris}

 %\hfill \textit{Lethbridge, Alberta, Canada}

 \hfill \textit{April 2015} % !!!

 \newpage

\putintocfalse % do NOT want these references in the Table of Contents
\makemarksfalse % do not want these references to change the funning head

\vfil
 \eject

%!TEX root = IntroArithGrps.tex

%\cleardoublepage

%\thispagestyle{plain}

\refstepcounter{section}% to correct the PDF bookmarks

\begingroup
 \makeatletter \@openrightfalse %acknowledgments don't need cleardoublepage
 \makeatother
% \tolerance=5000
 
\putintocfalse % do NOT want this in the Table of Contents
\chapter*{Acknowledgments}
%\begin{ack}

 In writing this book, I have received major help from discussions with
G.\,A.\,Margulis, 
M.\,S.\,Raghunathan, 
S.\,G.\,Dani, 
T.\,N.\,Venkataramana, 
Gopal Prasad, 
Andrei Rapinchuk,
Robert\,J.\,Zimmer,
Scot Adams, 
and 
Benson Farb.
I have also benefited from the comments
and suggestions of many other colleagues, including
 Ian Agol,
 Marc Burger,
 Indira Chatterji,
 Yves Cornulier,
 Alessandra Iozzi,
 Anders Karlsson,
 Sean Keel,
 Nicolas Monod,
 Hee Oh,
 Alan Reid,
 Yehuda Shalom,
 Hugh Thompson,
 Shmuel Weinberger,
 Barak Weiss,
 and 
 Kevin Whyte.
 They may note that some of the
remarks they made to me have been reproduced
here almost verbatim.

I would like to thank \'E.\,Ghys, D.\,Gaboriau and the other
members of the \'Ecole Normale Sup\'erieure of Lyon for the
stimulating conversations and invitation to speak that were
the impetus for this work, and B.\,Farb, A.\,Eskin, and other
members of the University of Chicago mathematics
department for encouraging me to write this introductory text,
and for the opportunity to lecture from it to a stimulating audience on numerous occasions.

I am also grateful to 
the University of Chicago, 
the \'Ecole Normale Sup\'erieure of Lyon, 
the University of Bielefeld, 
the Isaac Newton Institute for Mathematical Sciences, 
the University of Michigan, 
the Tata Institute for Fundamental Research,
the University of Western Australia,
and
the Mathematical Sciences Research Institute,
for their warm
hospitality while I worked on various parts of the book. 
The preparation of this manuscript was partially
supported by research grants from the National Science
Foundation of the USA and the National Science and Engineering
Research Council of Canada.
% \end{ack}

I also thank my wife, Joy Morris, for her patience and encouragement during the (too many!)\ years this project was underway.

\endgroup

  	\CCtrue

\mainmatter

 %!TEX root = IntroArithGrps.tex

\part{Introduction} \label{IntroductionPart}

 %!TEX root = IntroArithGrps.tex

\mychapter{What is a\texorpdfstring{\\}{ }Locally Symmetric Space?}
\label{WhatisLocSymmChap}

\prereqs{understanding of geodesics, and other concepts of Differential Geometry.}

In this chapter, we give a geometric introduction to the notion of a
symmetric space or a locally symmetric space, and explain the
central role played by simple Lie groups and their lattice subgroups.
(Since geometers are the target audience here, we assume familiarity with differential geometry that will not be needed in other parts of the book.)
This material is not a prerequisite for reading any of the later
chapters, except \cref{GeomIntroRank}; it is intended to provide
a geometric motivation for the study of lattices in semisimple Lie
groups. Since arithmetic subgroups are the primary examples of lattices, this also motivates the main topic of the rest of the book.

\section{Symmetric spaces}

Recall that a \defit{Riemannian manifold} is
a smooth manifold~$M$, together with the choice of an inner
product $\langle \cdot \mid \cdot \rangle_x$ on the tangent
space~$T_x M$, for each $x \in M$, such that $\langle \cdot
\mid \cdot \rangle_x$ varies smoothly as $x$ varies.
The nicest Riemannian manifolds are homogeneous. This means that every
point looks exactly like every other point:

\begin{defn} \label{HomogDefn}
 A Riemannian manifold~$X$ is a \defit{homogeneous space} if its
isometry group $\Isom(X)$ acts transitively. That is, for every $x,y
\in X$, there is an isometry~$\phi$ of~$X$, such that $\phi(x) = y$.
 \end{defn}

\begin{notation} \label{GoNotation}
We use \nindex{$G^\circ$ = identity component of~$G$}$G^\circ$ to denote the \term{identity component} of the group~$G$.
\end{notation}

\begin{egs}
 Here are some elementary examples of (simply connected) homogeneous
spaces.
 \begin{enumerate} \itemsep=\smallskipamount
 
 \item The round sphere $S^n 
 = \{\, x \in \real^{n+1} \mid \|x\| = 1 \,\}$. Rotations are the only
orientation-preserving isometries of~$S^n$, so we have $\Isom(S^n)^\circ =
\SO(n+1)$. Any point on~$S^n$ can be rotated to any other point, so
$S^n$ is homogeneous.

 \item Euclidean space $\real^n$. Every
orientation-preserving isometry of~$\real^n$ is a combination of a
translation and a rotation, and this implies that $\Isom(\real^n)^\circ = \SO(n) \ltimes
\real^n$. Any point in~$\real^n$ can be
translated to any other point, so $\real^n$ is homogeneous. 

 \item The hyperbolic plane $\hyperbolic^2
 = \{\, z \in \complex \mid \Im z > 0 \,\}$,
 where the inner product on $T_z \hyperbolic^2$ is given by
 $$ \langle u \mid v \rangle_{\hyperbolic^2}
 = \frac{1}{4(\Im z)^2} \, \langle u \mid v \rangle_{\real^2} .$$
 It is not difficult to show that 
 	$$  \text{$\Isom(\hyperbolic^2)^\circ$ is isomorphic to 
	$\PSL(2,\real)^\circ = \SL(2,\real)/\{\pm 1\}$,} $$
by noting that
$\SL(2,\real)$ acts on $\hyperbolic^2$ by linear-fractional
transformations $z \mapsto (az+b)/(cz+d)$, and confirming, by
calculation, that these linear-fractional transformations preserve
the hyperbolic metric.

 \item Hyperbolic space
 $\hyperbolic^n = \{\, x \in \real^n \mid x_n > 0 \,\}$,
 where the inner product on $T_x \hyperbolic^n$ is given by
 $$ \langle u \mid v \rangle_{\hyperbolic^n}
 = \frac{1}{4x_n^2} \, \langle u \mid v \rangle_{\real^n} .$$
 It is not difficult to see that $\hyperbolic^n$ is homogeneous
\csee{HnisHomog}. One can also show that that
the group $\Isom(\hyperbolic^n)^\circ$ is isomorphic to
$\SO(1,n)^\circ$ \csee{HyperModel}.
 \item A cartesian product of any combination of the above
\csee{prod->homog}.
 \end{enumerate}
 \end{egs}

\begin{defns}
 Let $\phi \colon X \to X$.
 \begin{enumerate}
 \item We say that $\phi$ is \defit[involutive
isometry]{involutive} (or that $\phi$ is an \defit{involution}) if
$\phi^2 = \Id\mk$.
 \item A \defit{fixed point} of~$\phi$ is a point $p \in X$, such that
$\phi(p) = p$.
 \item A fixed point~$p$ of~$\phi$ is \defit[isolated fixed
point]{isolated} if there is a neighborhood~$U$ of~$p$, such that $p$
is the only fixed point of~$\phi$ that is contained in~$U$.
 \end{enumerate}
 \end{defns}

Besides an isometry taking $x$ to~$y$, each of the above spaces also
has a nice involutive isometry that fixes~$x$.
 \begin{enumerate} \itemsep=\smallskipamount
 \item Define $\phi_1 \colon S^n \to S^n$ by
 $$ \phi_1(x_1,\ldots,x_{n+1}) = (-x_1,\ldots,-x_n,x_{n+1}) .$$
 Then $\phi_1$ is an isometry of~$S^n$, such that $\phi_1$ has only
two fixed points: namely, $e_{n+1}$ and~$-e_{n+1}$, where $e_{n+1} =
(0,0,\ldots,0,1)$. Therefore, $e_{n+1}$ is an isolated fixed point
of~$\phi_1$.
 \item Define $\phi_2 \colon \real^n \to \real^n$ by $\phi_2(x) =
-x$. Then $\phi_2$ is an isometry of~$\real^n$, such that $0$ is the
only fixed point of~$\phi_2$.
 \item Define $\phi_3 \colon \hyperbolic^2 \to \hyperbolic^2$ by
$\phi_3(z) = -1/z$. Then $i$ is the only fixed point of~$\phi_3$.
 \item There are involutive isometries of~$\hyperbolic^n$ that have a
unique fixed point \csee{invol(BallModel)}, but they are somewhat
difficult to describe in the upper-half-space model that we are using.
 \end{enumerate}
The existence of such an isometry is the additional condition that is
required to be a symmetric space.

\begin{defn} \label{symmDefn}
 A Riemannian manifold~$X$ is a \defit[symmetric!space]{symmetric space} if
 \begin{enumerate}
 \item $X$ is connected,
 \item \label{symmDefn-homog}
 $X$ is homogeneous, and
 \item \label{symmDefn-fp}
 there is an involutive isometry~$\phi$ of~$X$, such that $\phi$ has
at least one isolated fixed point.
 \end{enumerate}
 \end{defn}

\begin{rem}
 If $X$ is a symmetric space, then all points of~$X$ are essentially
the same, so, for each $x \in X$ (not only for \emph{some} $x \in X$),
there is an isometry~$\phi$ of~$X$, such that $\phi^2 = \Id$ and $x$
is an isolated fixed point of~$\phi$ \csee{symm->xfixed}.
Conversely, if Condition~\pref{symmDefn-fp} is replaced with this
stronger assumption, then Condition~\pref{symmDefn-homog} can be
omitted \csee{xfixed->symm}.
 \end{rem}

We constructed examples of involutive isometries of $S^n$, $\real^n$,
and~$\hyperbolic^n$ that have an isolated fixed point. The following
proposition shows that no choice was involved: the involutive
isometry with a given isolated fixed point~$p$ is unique, if it
exists. Furthermore, in the exponential coordinates at~$p$, the
involution must simply be the map $x \mapsto -x$.

\begin{prop} \label{phiunique}
 Suppose $\phi$ is an involutive isometry of a Riemmanian
manifold~$X$, and suppose $p$ is an isolated fixed point of~$\phi$.
Then
 \begin{enumerate}
 \item \label{phiunique-dphi}
 $d\phi_p = - \Id$, and
 \item \label{phiunique-phi}
 for every geodesic $\gamma$ with~$\gamma(0) = p$, we have $\phi
\bigl( \gamma(t) \bigr) = \gamma(-t)$, for all $t \in \real$.
 \end{enumerate}
 \end{prop}

\begin{proof}
 \pref{phiunique-dphi} From the Chain Rule, and the fact that $\phi(p) =
p$, we have
 $$ d(\phi^2)_p = d\phi_{\phi(p)} \circ d\phi_p = (d\phi_p)^2 .$$
 Also, because $\phi^2 = \Id$, we know that $d(\phi^2)_p = d \Id_p =
\Id$. We conclude that $(d\phi_p)^2 = \Id$; hence, the linear
transformation $d\phi_p \colon T_p X \to T_p X$ satisfies the polynomial
equation $x^2 - 1 = 0$. 

Suppose $d\phi_p \neq -\Id$. (This will lead to a contradiction.)
Since the polynomial $x^2 - 1$ has no repeated roots, we know that
$d\phi_p$ is diagonalizable. Furthermore, because $1$ and~$-1$ are
the only roots of $x^2 - 1$, we know that $1$ and~$-1$ are the only
possible eigenvalues of $d\phi_p$. Therefore, because $d\phi_p \neq -\Id$,
we conclude that $1$ is an eigenvalue; so we may choose some
nonzero $v \in T_p X$, such that $d\phi_p(v) = v$. Let $\gamma$ be the
geodesic with $\gamma(0) = p$ and $\gamma'(0) = v$. Then, because
$\phi$ is an isometry, we know that $\phi \circ \gamma$ is also a
geodesic. We have 
 $$(\phi \circ \gamma)(0) = \phi \bigl( \gamma(0) \bigr)
 = \phi(p) = p = \gamma(0) $$
 and
 $$ (\phi \circ \gamma)'(0)
 = d \phi_{\gamma(0)} \bigl( \gamma'(0) \bigr)
 = d \phi_p (v)
 = v
 = \gamma'(0) .$$
 Since every geodesic is uniquely determined by prescribing its initial
position and its initial velocity, we conclude that $\phi \circ \gamma
= \gamma$. Therefore, $\phi \bigl( \gamma(t) \bigr) = \gamma(t)$, so
$\gamma(t)$ is a fixed point of~$\phi$, for every~$t$. This
contradicts the fact that the fixed point $p = \gamma(0)$ is isolated.

\pref{phiunique-phi} Define $\overline{\gamma}(t) = \gamma(-t)$, so
$\overline{\gamma}$ is a geodesic. Because $\phi$ is an isometry, we
know that $\phi \circ \gamma$ is also a geodesic. We have
 $$(\phi \circ \gamma)(0) = \phi \bigl( \gamma(0) \bigr)
 = \phi(p) = p = \overline{\gamma}(0) $$
 and, from~\pref{phiunique-dphi},
 $$ (\phi \circ \gamma)'(0)
 = d \phi_{\gamma(0)} \bigl( \gamma'(0) \bigr)
 = -\gamma'(0)
 = \overline{\gamma}'(0) .$$
 Since a geodesic is uniquely determined by prescribing its initial
position and its initial velocity, we conclude that $\phi \circ
\gamma = \overline{\gamma}$, as desired.
 \end{proof}

\begin{defn}
 Let $M$ be a Riemannian manifold, and let $p \in M$. 
It is a basic fact of differential geometry that there is a neighborhood~$V$ of~$0$ in $T_pM$, such that the exponential map $\exp_p$ maps~$V$ diffeomorphically onto a
neighborhood~$U$ of~$p$ in~$M$. By making $V$ smaller, we may assume it is:
	\begin{itemize}
	\item \defit[symmetric!neighborhood]{symmetric} (that is, $-V = V$),
	and
	\item \defit[star-shaped neighborhood]{star-shaped} (that is, $tV \subseteq V$, for $0\le t < 1$).
	\end{itemize}
The \defit{geodesic symmetry} at~$p$ is the diffeomorphism~$\tau$ of~$U$ that is defined by
 $$ \tau \bigl( \exp_p(v) \bigr) = \exp_p(-v) ,$$
 for all $v \in V$.

In other words, for each geodesic~$\gamma$ in~$M$, such that
$\gamma(0) = p$, and for all $t \in \real\mk$, such that $t \gamma'(0)
\in V$, we have $\tau \bigl( \gamma(t) \bigr) = \gamma(-t)$.
 \end{defn}

\begin{note}
 The geodesic symmetry~$\tau$ is a local diffeomorphism, but, for
most manifolds~$M$, it is \emph{not} a local isometry
\ccf{locsymm<>geodsymm}.
 \end{note}

In this terminology, the preceding proposition shows that if an
involutive isometry~$\phi$ has a certain point~$p$ as an isolated
fixed point, then, locally, $\phi$ must agree with the geodesic
symmetry at~$p$. This has the following easy consequence, which is
the motivation for the term \defit[symmetric!space]{symmetric space}.

\begin{cor} \label{symm<>geodsymm}
 A connected Riemannian manifold~$M$ is a symmetric space if and only
if, for each $p \in M$, the geodesic symmetry at~$p$ extends to an
isometry of~$M$.
 \end{cor}

\begin{exercises}

\item \label{HnisHomog}
 Show that $\hyperbolic^n$ is homogeneous. 
 \hint{For any $t
\in \real^+$, the dilation $x \mapsto tx$ is an isometry
of~$\hyperbolic^n$. Also, for any $v \in \real^{n-1}$, the
translation $x \mapsto x + v$ is an isometry of~$\hyperbolic^n$.}

\item \label{HnBallModel}
 Let $B^n = \{\, x \in \real^n \mid \|x\| < 1 \,\}$ be the open unit
ball in~$\real^n$, equip $T_x B^n$ with the inner product
 $$ \langle u \mid v \rangle_{B^n}
 = \frac{1}{ \bigl( 1 - \|x\|^2 \bigr)^2} \, \langle u \mid v
\rangle_{\real^n} ,$$
 and let $e_n = (0,0,\ldots,0,1) \in \real^n$.
 Show that the map $\phi \colon B^n \to \hyperbolic^n$ defined by
 $$ \phi(x) = \frac{x + e_n}{\|x + e_n\|^2} - \frac{1}{2} e_n$$
 is an isometry from~$B_n$ onto~$\hyperbolic^n$. (In geometric terms,
$\phi$ is obtained by composing a translation with the inversion
centered at the south pole of~$B^n$.)

\item \label{invol(BallModel)}
 Show that $x \mapsto -x$ is an isometry of~$B^n$ (with respect to the
Riemannian metric $\langle \cdot \mid \cdot \rangle_{B^n}$ defined in \cref{HnBallModel}).

\item \label{HyperModel}
 For $u,v \in \real^{n+1}$, define
 $$ \langle u \mid v \rangle_{1,n}
 = u_0 v_0 - \sum_{j=1}^n u_j v_j .$$
 (Note that, for convenience, we start our numbering of the
coordinates at~$0$, rather than at~$1$.) Let
 $$ X_{1,n}^+ = \{\, x \in \real^{n+1} \mid \langle x \mid x
\rangle_{1,n} = 1, \ x_0 > 0 \,\} ,$$
 so $X_{1,n}^+$ is one sheet of a 2-sheeted hyperboloid. Equip $T_x
X_{1,n}^+$ with the inner product obtained by restricting $\langle
\cdot \mid \cdot \rangle_{1,n}$ to this subspace. 
 \begin{enumerate}
 \item Show that the bijection $\psi \colon B^n \to X_{1,n}^+$
defined by 
 $$ \psi(x) = \frac{1}{1 - \|x\|^2} \, (1,x) $$
 is an isometry. (Note that this implies that the restriction of
$\langle \cdot \mid \cdot \rangle_{1,n}$ to $T_x X_{1,n}^+$ is
positive definite, even though $\langle \cdot \mid \cdot
\rangle_{1,n}$ is not positive definite on all of~$\real^{n+1}$.)
 \item Show $\SO(1,n)^\circ$ acts transitively on~$X_{1,n}^+$ by
isometries.
 \end{enumerate}

\item \label{Stab(Hn)} 
 For $G = \SO(1,n)^\circ = \Isom(\hyperbolic^n)^\circ$, show there is
some $p \in \hyperbolic^n$, such that $\Stab_G(p) = \SO(n)$.
 \hint{This is easy in the hyperboloid model~$X_{1,n}^+$.}

\item \label{prod->homog}
 Show that if $X_1,X_2,\ldots,X_n$ are homogeneous spaces, then the
cartesian product $X_1 \times X_2 \times \cdots \times X_n$ is also
homogeneous.

\item Show that every homogeneous space is \defit{geodesically
complete}. That is, for every geodesic segment $\gamma \colon
(-\epsilon, \epsilon) \to X$, there is a doubly-infinite geodesic
$\overline{\gamma} \colon \real \to X$, such that
$\overline{\gamma}(t) = \gamma(t)$ for all $t \in (-\epsilon,
\epsilon)$.

\item Show that if $X_1,\ldots,X_n$ are symmetric spaces, then the
cartesian product $X_1 \times X_2 \times \cdots \times X_n$ is also a
symmetric space.

\item \label{symm->xfixed}
 Show that if $X$ is a symmetric space, then, for each $x \in X$,
there is an isometry~$\phi$ of~$X$, such that $\phi^2 = \Id$ and $x$ is
an isolated fixed point of~$\phi$.

\item \label{xfixed->symm}
 Let $X$ be a connected Riemannian manifold, and assume, for each $x
\in X$, that there is an isometry~$\phi$ of~$X$, such that $\phi^2 =
\Id$ and $x$ is an isolated fixed point of~$\phi$. Show that $X$ is
homogenous, and conclude that $X$ is a symmetric space.

\item Show that the real projective space $\real P^n$ (with the metric
that makes its universal cover a round sphere) has an involutive
isometry~$\phi$, such that $\phi$ has both an isolated fixed point,
and a fixed point that is not isolated. Is $\real P^n$ a symmetric
space?

\end{exercises}

\section{How to construct a symmetric space} \label{ConstructSymm}

In this section, we describe how Lie groups are used to construct
symmetric spaces. Let us begin by recalling the well-known
group-theoretic structure of any homogeneous space.

Suppose $X$ is a connected homogeneous space, and let $G =
\Isom(X)^\circ$. Because $\Isom(X)$ is transitive on~$X$, and $X$ is
connected, we see that $G$ is transitive on~$X$
\csee{Xconn->idcomp}, so we may identify $X$ with the coset space
$G/K$, where $K$ is the stabilizer of some point in~$X$. Note that
$K$ is compact \csee{StabCpct}.

Conversely, if $K$ is any compact subgroup of any Lie group~$G$, then
there is a $G$-invariant Riemannian metric on $G/K$
\csee{G/K->Riem}, so $G/K$ (with this metric) is a homogeneous space.
 (For
any closed subgroup~$H$ of~$G$, the group~$G$ acts transitively on
the manifold $G/H$, by diffeomorphisms. However, when $H$ is not
compact, $G$ usually does not act by isometries of any Riemannian
metric on $G/H$, so there is no reason to expect $G/H$ to be a
homogeneous space in the sense of \cref{HomogDefn}.)

 \begin{eg}  \ 
 \noprelistbreak
 \begin{enumerate}
 
 \item For $X = S^n$, we have $G = \SO(n+1)$, and we may let 
 	$$ \text{$K = \Stab_G(e_{n+1}) = \SO(n)$, \ so \ $S^n = \SO(n+1)/\SO(n)$} .$$
Note that,
letting $\sigma$ be the diagonal matrix 
 $$ \sigma = \diag(-1,-1,\ldots,-1,1) ,$$ we
have $\sigma^2 = \Id$, and
  $K = \czer_G(\sigma)$ is the centralizer of~$\sigma$ in~$G$.
  
 \item For $X = \real^n$, we have $G = \SO(n) \ltimes \real^n$, and
we may let 
	$$ \text{$K = \Stab_G(0) = \SO(n)$, \ so \ $\real^n = \bigl( \SO(n)
\ltimes \real^n)/\SO(n)$} .$$
Note that the map $\sigma \colon (k,v)
\mapsto (k,-v)$ is an automorphism of~$G$, such that $\sigma^2 =
\Id$, and 
 $$ \czer_G(\sigma) = \{\, g \in G \mid \sigma(g) = g \,\} = K .$$
 
 \item For $X = \hyperbolic^2$, we have $G \approx \SL(2,\real)$, and
we may let 
	$$ \text{$K = \Stab_G(i) \approx \SO(2)$, \ so \ 
 $\hyperbolic^2 = \SL(2,\real)/\SO(2)$} .$$
  
 \item For $X = \hyperbolic^n$, we have $G = \SO(1,n)^\circ$, and we
may take $K = \SO(n)$ \csee{Stab(Hn)}.  Note that, for
$\sigma = \diag(1,-1,-1,\ldots,-1)$, we have
$\sigma^2 = \Id$, and
  $K = \czer_G(\sigma)$.
 \end{enumerate}
 Therefore, in each of these cases, there is an automorphism~$\sigma$
of~$G$, such that $K$ is the centralizer of~$\sigma$. (In other words,
$K = \{\, k \in G \mid \sigma(k) = k \,\}$ is the set of fixed points of~$\sigma$ in~$G$.) The following
proposition shows, in general, that a slightly weaker condition makes
$G/K$ symmetric, not just homogeneous.
 \end{eg}

\begin{prop} \label{G/K->symm}
 Let
 \begin{itemize}
 \item $G$ be a connected Lie group,
 \item $K$ be a compact subgroup of~$G$, and
 \item $\sigma$ be an involutive automorphism of~$G$, such that $K$
is an open subgroup of $\czer_G(\sigma)$.
 \end{itemize}
 Then $G/K$ can be given the structure of a symmetric space, such that
the map $\tau(gK) = \sigma(g) K$ is an involutive isometry of $G/K$
with $eK$ as an isolated fixed point.
 \end{prop}

\begin{proof}
 To simplify the proof slightly, let us assume that $K = \czer_G(\sigma)$
\csee{KinC(sigma)}. 

Because $K$ is compact, we know there is a $G$-invariant Riemmanian
metric on $G/K$ \csee{G/K->Riem}. Then, because $\langle \tau
\rangle$ is finite, and normalizes~$G$, it is not difficult to see
that we may assume this metric is also $\tau$-invariant
\csee{FiniteAvgMetric}. (This conclusion can also be reached by
letting $G^+ = \langle \sigma \rangle \ltimes G$ and $K^+ = \langle
\sigma \rangle \times K$, so $K^+$ is a compact subgroup of~$G^+$,
such that $G^+/K^+ = G/K$.) Therefore, $\tau$ is an involutive isometry of
$G/K$.

Suppose $gK$ is a fixed point of~$\tau$, with $g \approx e$. Then
$\sigma(g) \in g K$, so we may write $\sigma(g) = gk$, for some $k
\in K$. Since $\sigma$ centralizes~$k$ (and $\sigma$ is an
automorphism), we have
 $$\sigma^2(g) = \sigma \bigl( \sigma(g) \bigr)
 = \sigma(gk) = \sigma(g) \, \sigma(k) 
 = (gk) (k)= g k^2 .$$
 On the other hand, we know $\sigma^2(g) = g$ (because $\sigma$ is
involutive), so we conclude that $k^2 = e$.

Since $g \approx e$, and $\sigma(e) = e$, we have $\sigma(g) \approx
g$, so $k = g^{-1} \sigma(g) \approx e$. Since $k^2 = e$, we conclude
that $k = e$. (There is a neighborhood~$U$ of~$e$ in~$G$, such that,
for every $u \in U \smallsetminus \{e\}$, we have $u^2 \neq e$.) 
Therefore $\sigma(g) = gk = ge = g$, so $g \in \czer_G(\sigma) = K$;
hence, $gK = eK$.
 \end{proof}

Conversely, for any symmetric space~$X$, there exist $G$, $K$,
and~$\sigma$ as in \cref{G/K->symm}, such that $X$ is
isometric to $G/K$ \csee{symm->G/K}. 

\begin{eg} \label{SLnSymm}
 Let $G = \SL(n,\real)$, $K = \SO(n)$, and define $\sigma(g) =
(g^{-1})^T$ (the transpose-inverse). Then $\sigma^2 = 1$ and
$\czer_G(\sigma) = K$, so the theorem implies that $G/K$ is a symmetric
space. Let us describe this space somewhat more concretely.

Recall that any real symmetric matrix~$A$ can be diagonalized
over~$\real$. In particular, all of its eigenvalues are real. If all the
eigenvalues of~$A$ are strictly positive, then we say that $A$ is
\defit[positive!definite]{positive definite}.

Let 
 $$X = \{\, A \in \SL(n,\real) \mid 
 \mbox{$A$ is symmetric and positive definite} \,\} ,$$
 and define $\alpha \colon G \times X \to X$ by 
 $\alpha(g,x) = g x g^T$. 
 Then:
 \begin{enumerate} \renewcommand{\theenumi}{\alph{enumi}}
 \item \label{SLnSymm-act}
 $\alpha$ defines an action of $G$ on~$X$; i.e., we have
$\alpha(gh,x) = \alpha \bigl( g, \alpha(h,x) \bigr)$ for all $g,h \in G$
and $x \in X$.
\item \label{SLnSymm-G/K}
This action is transitive, and we have $K = \Stab_G(\Id)$, so $X$ may be
identified with $G/K$.
 \item \label{SLnSymm-TX}
 $T_{\Id} X = 
 \{\, u \in \Mat_{n\times n}(\real) \mid 
 \mbox{$u$ is symmetric and $\trace(u) = 0$} \,\} $. (By definition,
we have $X \subseteq \SL(n,\real)$. The condition $\trace(u) = 0$ is obtained by differentiating the restriction $\det(A) = 1$.)
 \item \label{SLnSymm-<>}
 The inner product $\langle u \mid v \rangle = \trace(uv)$ on
$T_{\Id} X$ is $K$-invariant, so it may be extended to a $G$-invariant
Riemannian metric on~$X$.
 \item \label{SLnSymm-tau}
 The map $\tau \colon X \to X$, defined by $\tau(A) = A^{-1}$, is an
involutive isometry of~$X$, such that 
 $\tau \bigl( \alpha(g,x) \bigr) = \sigma(g) \, \tau(x)$ for all $g
\in G$ and $x \in X$.
 \end{enumerate}
 \end{eg}

\begin{eg}
 Other examples of symmetric spaces are:
 \begin{enumerate}
 \item $\SL(n,\complex)/\SU(n)$, and
 \item $\SO(p,q)^\circ / \bigl( \SO(p) \times \SO(q) \bigr) $.
 \end{enumerate}
 \end{eg}

These are special cases of a consequence of
\cref{G/K->symm} that will be stated after we introduce some terminology.

\begin{defns} \ 
\noprelistbreak
	\begin{enumerate}
	\item A symmetric space~$X$ is \defit[irreducible!symmetric
space]{irreducible} if its universal cover is not isometric to any
nontrivial product $X_1 \times X_2$.
	\item A Riemannian manifold is \defit[flat!manifold]{flat} if its curvature tensor is identically zero, or, equivalently, if every point in~$X$ has a neighborhood that is isometric to an open subset of the Euclidean space~$\real^n$.
	\end{enumerate}
 \end{defns}

\begin{prop} \label{noncpctsymm<>simplegrp}
 Let $G$ be a connected, noncompact, simple Lie group with finite
center. Then $G$ has a maximal compact subgroup~$K$ \textup(which
is unique up to conjugacy\textup), and $G/K$ is a simply
connected, noncompact, irreducible symmetric space.
Furthermore, $G/K$ has non-positive sectional curvature and is not flat.

Conversely, any noncompact, non-flat, irreducible symmetric space is
of the form $G/K$, where $G$ is a connected, noncompact, simple Lie
group with trivial center, and $K$ is a maximal compact subgroup
of~$G$.
 \end{prop}

\begin{rem}
 Let $K$ be a compact subgroup of a connected, simple Lie group~$G$
with finite center, such that $G/K$ is a symmetric space
\ccf{G/K->symm}. \Cref{noncpctsymm<>simplegrp} shows that if $G$ is not compact,
then $K$ must be a maximal compact subgroup of~$G$, which is
essentially unique. 

On the other hand, if $G$ is compact, then the subgroup~$K$ may not
be unique, and may not be maximal. For example, both
$\SO(n)/\SO(n-1)$ and $\SO(n)/\{e\}$ are symmetric spaces. The former
is a round sphere, which has already been mentioned. The latter is a
special case of the fact that every connected, compact Lie group is a
symmetric space \csee{X=G}.

\'E.\,Cartan obtained a complete list of all the symmetric spaces
(both compact and noncompact) by finding all of the simple Lie
groups~$G$ \csee{RealSimpleGrps}, and determining, for each of them, which compact subgroups~$K$ can arise in
\cref{G/K->symm}.
 \end{rem}

\begin{exercises}

\item \label{Xconn->idcomp}
 Suppose a topological group~$G$ acts transitively (and continuously)
on a connected topological space~$M$. Show that the identity component
$G^\circ$ is transitive on~$M$.

\item \label{StabCpct}
 Let $\{g_n\}$ be a sequence of isometries of a connected, complete Riemannian manifold~$M$, and assume there exists $p \in M$, such that $g_n p = p$ for all~$n$.) Show there is a subsequence $\{g_{n_k}\}$ of~$\{g_n\}$ that \term[converge uniformly on compact sets]{converges uniformly on compact subsets} of~$M$. (That is, there is some isometry~$g$ of~$M$, such that, for every $\epsilon > 0$ and every compact subset~$C$ of~$M$, there exists $k_0$, such that $d(g_{n_k} c, gc) < \epsilon$ for all $c \in C$ and all $k > k_0$.)
\hint{This is a special case of the \thmindex{Arzelà-Ascoli}Arzelà-Ascoli Theorem. For each $c \in C$, the sequence $\{g_n c\}$ is bounded, and therefore has a convergent subsequence. By Cantor diagonalization, there is a subsequence that works for all~$c$ in a countable, dense subset of~$C$.}

\item \label{cpct->orthog}
 Let $K$ be a compact group, and let $\rho \colon K \to \GL(n,\real)$
be a continuous homomorphism. Show that there is a $K$-invariant inner
product $\langle \cdot \mid \cdot \rangle_K$ on~$\real^n$; that is,
such that
 $ \bigl\langle \rho(k) u \mid \rho(k) v \bigr\rangle_K =
\bigl\langle u \mid v \bigr\rangle_K $
 for all $k \in K$ and all $u,v \in \real^n$.
\hint{Define \ 
 $ \langle u \mid v \rangle_K
 = \int_K \bigl\langle \rho(k) u \mid \rho(k) v \bigr\rangle \,
d\mu(k) $,
 \ where $\mu$ is Haar measure on~$K$.}

\item \label{G/K->Riem}
 Let $K$ be a compact subgroup of a Lie group~$G$. Use
\cref{cpct->orthog} to show that there is a $G$-invariant
Riemannian metric on $G/K$.
\hint{A $G$-invariant Riemannian
metric on $G/K$ is determined by the inner product it assigns to the tangent space
$T_{eK}(G/K)$.}

\item \label{KinC(sigma)}
 Complete the proof of \cref{G/K->symm}, by removing the
simplifying assumption that $K = \czer_G(\sigma)$.

\item \label{FiniteAvgMetric}
 Let $F$ be a finite group of diffeomorphisms (not necessarily
isometries) of a Riemannian manifold $\bigl( M, \langle \cdot \mid
\cdot \rangle_x \bigr)$. Define a new inner product $\langle \cdot
\mid \cdot \rangle'_x$ on each tangent space $T_xM$ by
 $$ \langle u \mid v \rangle'_x
 = \sum_{f \in F} \langle df_x(u) \mid df_x(v) \rangle_{f(x)} .$$
 \begin{enumerate}
 \item Show that the Riemannian metric $\langle \cdot \mid \cdot
\rangle'$ on~$M$ is $F$-invariant.
 \item Show that if $G$ is a group of isometries of $\bigl( M,
\langle \cdot \mid \cdot \rangle_x \bigr)$, and $G$~is normalized
by~$F$, then $\langle \cdot \mid \cdot \rangle'$ is $G$-invariant.
 \end{enumerate}

\item \label{symm->G/K}
 For any symmetric space~$X$, show that there exist $G$, $K$,
and~$\sigma$ as in \cref{G/K->symm}, such that $X$ is
isometric to $G/K$. 
\hint{Suppose $\tau$ is an involutive isometry of~$X$ with an
isolated fixed point~$p$. Let $G = \Isom(X)^\circ$ and $K =
\Stab_G(p)$. Define $\sigma(g) = \tau g \tau$. Show $K \subset
\czer_G(\sigma)$ and, using the fact that $p$ is isolated, show that $K$
contains the identity component of $\czer_G(\sigma)$.}

\item Verify assertions
 \pref{SLnSymm-act}, \pref{SLnSymm-G/K}, \pref{SLnSymm-TX},
\pref{SLnSymm-<>}, and \pref{SLnSymm-tau}
 of \cref{SLnSymm}.
 \hint{To prove transitivity in \pref{SLnSymm-G/K}, you may assume that every symmetric matrix $A$ is diagonalizable by an \emph{orthogonal} matrix. That is, there exists $g$, such that $g A g^{-1}$ is diagonal and $g g^\transpose = \Id$. Note that every positive-definite diagonal matrix has a square root that is also a diagonal matrix.}

\item \label{IsomFinComps}
Show that if $X$ is a connected homogeneous space, then
$\Isom(X)$ has only finitely many connected components.
\hint{Every component of $G = \Isom(X)$ intersects the compact group $\Stab_g(x)$.}

\item \label{X=G}
 Show that if $G$ is compact, then there is a $G$-invariant
Riemannian metric on~$G$ that makes $G$ a symmetric space.
 \hint{The involutive isometry is $g \mapsto g^{-1}$.}

\end{exercises}

\section{Locally symmetric spaces}

The gist of the following definition is that a locally symmetric
space is a Riemannian manifold that is locally isometric to a
symmetric space; that is, every point has a neighborhood that is
isometric to an open subset of some symmetric space.

\begin{defn}
 A complete Riemannian manifold~$M$ is \defit[locally!symmetric]{locally symmetric} if
its universal cover is a symmetric space. In other words, there is a
symmetric space~$X$, and a group~$\Gamma$ of isometries of~$X$, such
that 
 \begin{enumerate}
 \item $\Gamma$ acts freely and properly discontinuously on~$X$, and
 \item $M$ is isometric to $\Gamma \backslash X$.
 \end{enumerate}
 \end{defn}

\begin{rem} \label{locsymm<>geodsymm}
 At every point of a symmetric space, the geodesic symmetry $\gamma(t) \mapsto
\gamma(-t)$ extends to an isometry of the entire manifold
\csee{symm<>geodsymm}. In a locally symmetric space, the geodesic
symmetry~$\tau$ at each point is an isometry on its domain, but it may
not be possible to extend~$\tau$ to an isometry that is well-defined
on the entire manifold; that is, the geodesic symmetry is only a
\emph{local} isometry. That is the origin of the term \defit[locally!symmetric]{locally
symmetric}.
 \end{rem}

\begin{eg}
 Define $g \colon \hyperbolic^2 \to \hyperbolic^2$ by $g(z) = z+1$,
let $\Gamma = \langle g \rangle$, and let $M = \Gamma \backslash
\hyperbolic^2$. Then (obviously) $M$ is locally symmetric.

However, $M$ is not symmetric. We provide several different
geometric proofs of this fact, in order to illustrate the important
distinction between symmetric spaces and locally symmetric spaces.
(It can also be proved group-theoretically \csee{GrpCuspNotSymm}.)
 The manifold~$M$ is a cusp:
 
 $$\includegraphics{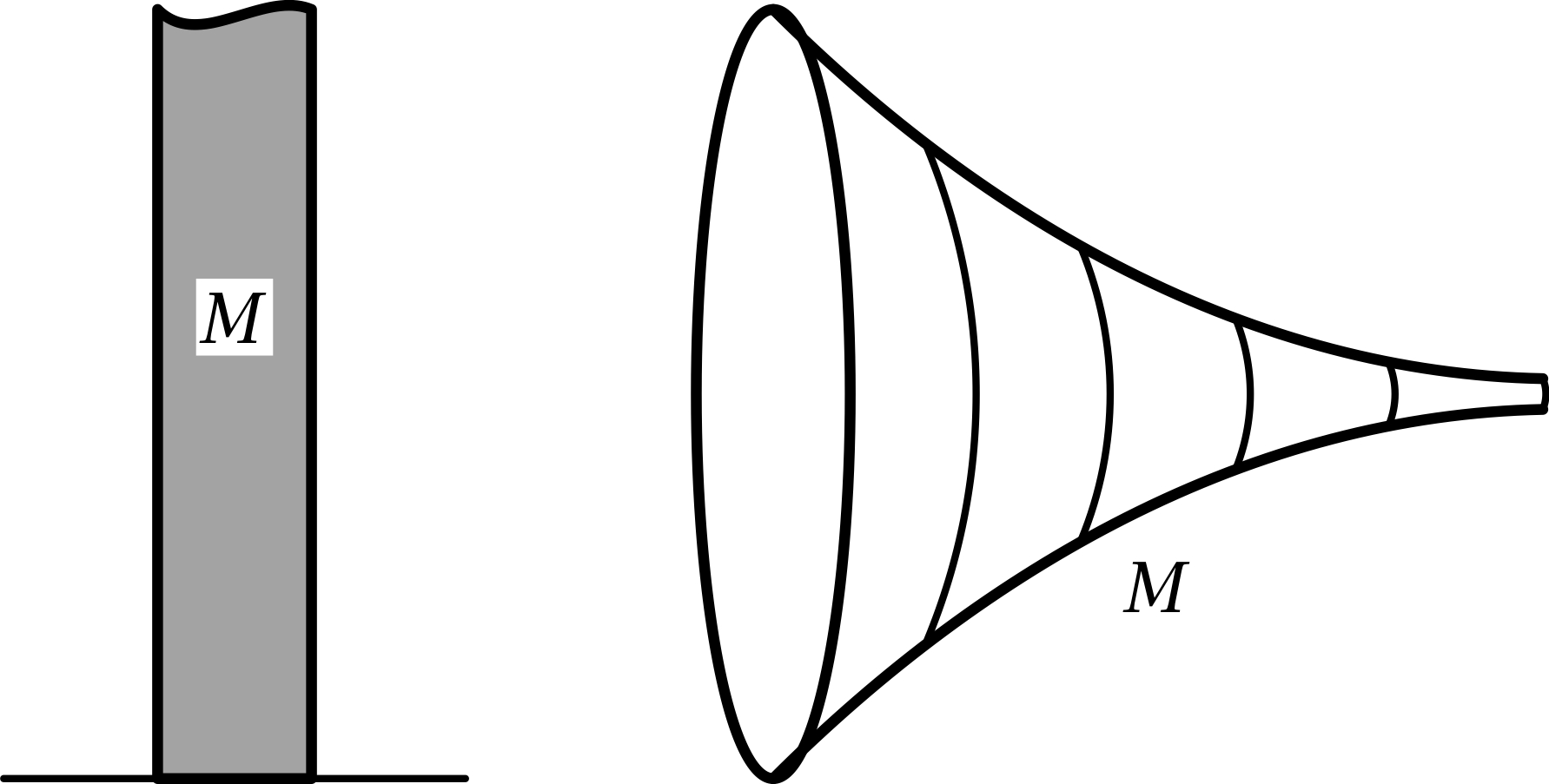}$$
%texpreamble
%("  \usepackage{amsmath}
% \usepackage[LY1]{fontenc}
% \usepackage[expert,LY1,mylucidascale]{mylucidabr}
% ");
%defaultpen(  fontcommand("\normalfont") + fontsize(10) ); 
%
%from graph access *;
%unitsize(0.75cm);
%
%real linethick = 1.5;
%
%real w = 1, h = 5, M = 3, axis = 1;
%
%draw( (-axis,0)--(w + axis, 0), linewidth(1) );
%fill( (0,0)--(w,0)--(w,h){WNW}..{NW}(0,h)--cycle, gray(0.7) );
%draw( (0,0)--(w,0)--(w,h){WNW}..{NW}(0,h)--cycle, linewidth(linethick) );
%
%real bb = 0.25;
%fill( shift(w/2,M)*( (-bb,-bb)--(-bb,bb)--(bb,bb)--(bb,-bb) -- cycle), white);
%label( "$M$", (w/2,M) );
%
%
%currentpen = linewidth(linethick);
%
%real gap = 2, a = 5, b = 0.1, cusp;
%
%cusp = w + axis + gap;
%
%draw( (cusp,0){NE}..{(1,0.02)}(cusp + a,(h/2)-b) );
%draw( (cusp,h){SE}..{(1,-0.02)}(cusp + a,(h/2)+b) );
%draw( ellipse( (cusp,h/2), 0.5, h/2) ) ;
%
%currentpen = linewidth(1);
%
%//draw( (cusp, 0){NNW}..{NNE}(cusp, h) );
%//draw( (cusp, 0){NNE}..{NNW}(cusp, h) );
%draw( (cusp+(a/5), 0.9){NNE}..{NNW}(cusp+(a/5), h-0.9) );
%draw( (cusp+2*(a/5), 1.55){NNE}..{NNW}(cusp+2*(a/5), h-1.55) );
%draw( (cusp+3*(a/5), 2){NNE}..{NNW}(cusp+3*(a/5), h-2) );
%draw( (cusp+4*(a/5), 2.3){NNE}..{NNW}(cusp+4*(a/5), h-2.3) );
%draw( (cusp+5*(a/5), 2.4){NNE}..{NNW}(cusp+5*(a/5), h-2.4) );
%
%label( "$M$", (6.5,1.25) );

\begin{enumerate}

\item Any point far out in the cusp lies on a short loop that is not
null-homotopic, but points at the other end do not lie on such a
loop. Therefore, $M$ is not homogeneous, so it cannot be symmetric.

\item The geodesic symmetry performs a $180^\circ$ rotation. Therefore, if
it is a well-defined diffeomorphism of~$M$, it must
interchange the two ends of the cusp. However, one end is thin, and
the other end is (very!) wide, so no isometry can interchange these
two ends. Hence, the geodesic symmetry (at any point) is not an
isometry, so $M$ is not symmetric.

\item Let us show, directly, that the geodesic symmetry at some
point $p \in \hyperbolic^2$ does not factor through to a well-defined
map on $\Gamma \backslash \hyperbolic^2 = M$.

 \begin{itemize}
 \item Let $x = -1 + i$ and $y = 1+i$, and let $p \in i\real$ be the
midpoint of the geodesic segment joining~$x$ and~$y$:
%\begin{figure}[h]
$$ \includegraphics{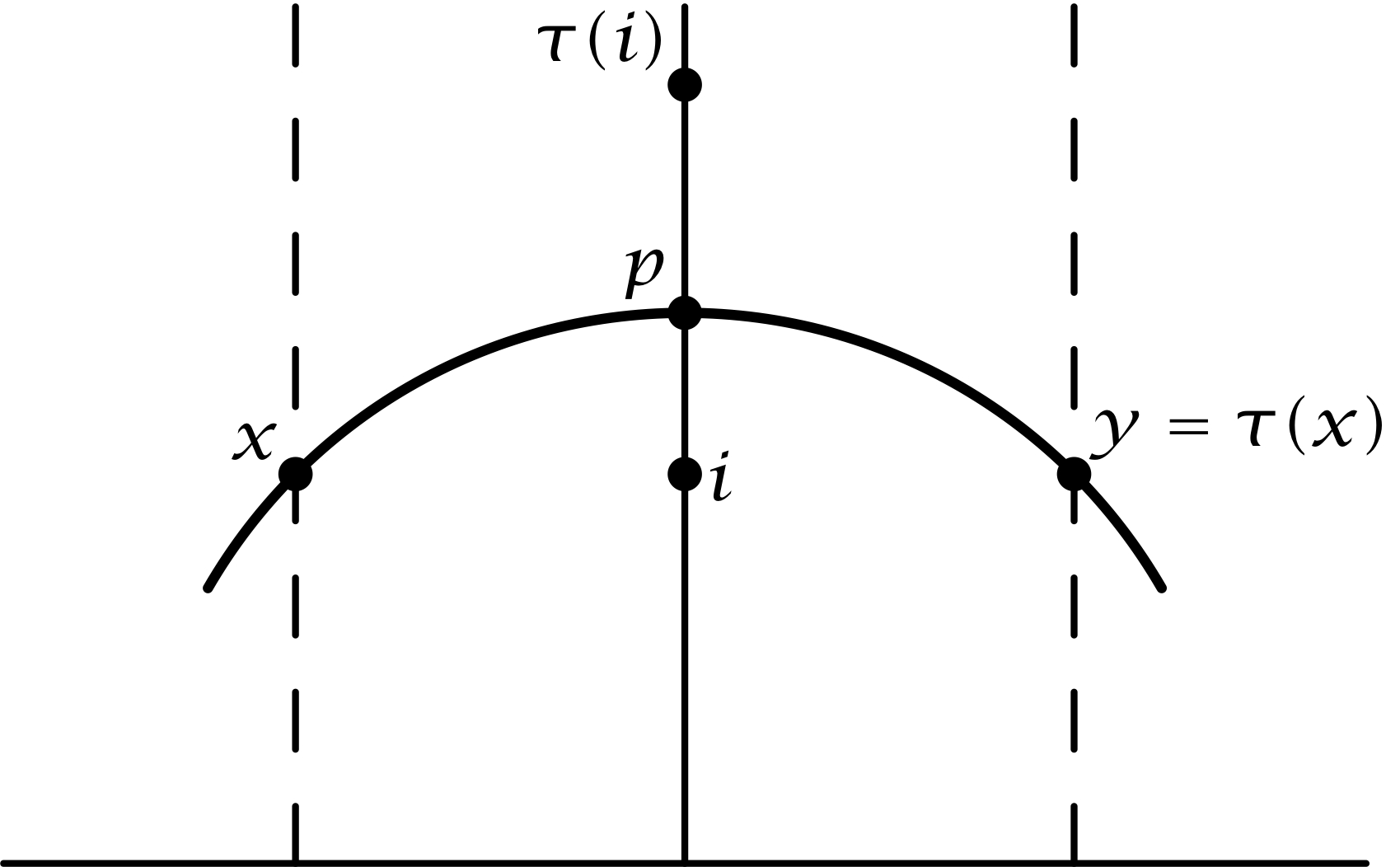} $$
% \caption{The geodesic symmetry~$\tau$ at~$p$.}
%\label{CuspNotSymm}
% \end{figure}
%texpreamble
%(" \usepackage[LY1]{fontenc}
% \usepackage[expert, LY1, mylucidascale]{mylucidabr} % I adjusted the scaling
% \usepackage{amsmath}
% ");
%defaultpen(  fontcommand("\normalfont") + fontsize(10) ); 
%
%from graph access *;
%unitsize(2cm);
%
%real r = sqrt(2), a = -1, b = 1;
%
%currentpen= defaultpen + linewidth(1.5);
%
%pair x = (a,b), y = (-a,b), p = (0,r), i = (0,1), ti = (0, 2) ;
%
%real w = 1.75*x.x, h = 1.1*ti.y;
%
%draw( (w,0)--(-w,0), linewidth(1) );
%draw( (0,0)--(0,h), linewidth(1) );
%
%draw( (x.x,0)--(x.x,h), dashed + linewidth(1) );
%draw( (y.x,0)--(y.x,h), dashed + linewidth(1) );
%
%label( "$x$" , x , NW);
%label( "$y = \tau(x)$" , y , NE);
%label( "$p$" , p , NW);
%label( "$i$" , i , E);
%label( "$\tau(i)$" , ti , NW);
%draw( x , linewidth(5) );
%draw( y , linewidth(5) );
%draw( p , linewidth(5) );
%draw( i , linewidth(5) );
%draw( ti , linewidth(5) );
%
%draw( arc( (0,0), r, 30, 150 ) );

 \item Let $\tau$ be the geodesic symmetry at~$p$. Then $\tau(x) = y
= 1 + i$.
 \item Because the imaginary axis is a geodesic, we have $\tau(i) =
ai$, for some $a > 1$.
 \item Now $i = x + 1 = g(x)$, so $x$ and~$i$ represent the same point
in~$M$. However, $\tau(i) - \tau(x) = -1 + (a-1)i$ is not an
integer (it is not even real), so $\tau(x)$ and~$\tau(i)$ do
\textbf{not} represent the same point in~$M$.  Therefore, $\tau$ does not
factor through to a well-defined map on~$M$.
 \end{itemize}

 \end{enumerate}
 \end{eg}

 \begin{rems} \ \label{locsymm<>curv}
 \noprelistbreak
 \begin{enumerate}
 \item Some authors do not require $M$ to be complete in their
definition of a locally symmetric space. This would allow the
universal cover of~$M$ to be an open subset of a symmetric space,
instead of the entire symmetric space.
 \item A more intrinsic (but more advanced) definition is that a
complete, connected Riemannian manifold~$M$ is \defit[locally!symmetric]{locally
symmetric} if and only if the curvature tensor of~$M$ is invariant
under all parallel translations, and $M$ is complete. 
 \end{enumerate}
 \end{rems}

Any complete, connected manifold of constant negative curvature is a
locally symmetric space, because the universal cover of such a manifold is $\hyperbolic^n$  (after normalizing the curvature to be~$-1$).
As a generalization of this, we are interested in locally symmetric
spaces~$M$ whose universal cover~$\widetilde{M}$ is of
\defit{noncompact type}, with no flat factors; that is, such that
each irreducible factor of $\widetilde{M}$ is noncompact (and not
flat).  From \cref{noncpctsymm<>simplegrp}, we see, in this
case, that $\widetilde{M}$ can be written in the form $\widetilde{M} = 
G/K$, where $G = G_1 \times \cdots \times G_n$ is a product of
noncompact simple Lie groups, and $K$~is a maximal compact subgroup
of~$G$. We have $M = \Gamma \backslash \widetilde{M}$, for some discrete subgroup~$\Gamma$ of $\Isom \bigl(  \widetilde{M}  \bigr)$. We know that $\Isom \bigl(  \widetilde{M} \bigr)$ has only finitely many connected components \csee{IsomFinComps}, so, if we replace $M$ with an appropriate finite cover, we can arrange that $\Gamma \subset \Isom \bigl( \, \widetilde{M} \, \bigr)^\circ = G$. Then
	$$ \text{$M = \Gamma \backslash G/K$, and $\Gamma$~is a discrete subgroup of~$G$.} $$

A topologist may like $M$ to be compact, but it turns out that a
very interesting theory is obtained by making the weaker assumption
that $M$ has finite volume. Hence, the subgroup~$\Gamma$ should be
chosen so that $\Gamma \backslash G/K$ has finite volume. Because
$\Gamma \backslash G$ is a principal $K$-bundle over $\Gamma
\backslash G/K$, and $K$~has finite measure, it is not difficult to
see, from \thmindex{Fubini's}{Fubini's Theorem}, that $\Gamma \backslash G$ has finite
volume \csee{lattice<>mu(X/Gamma)}. This leads to the following
definition.

\begin{defn}
 A subgroup~$\Gamma$ of~$G$ is a \defit[lattice!subgroup]{lattice} in~$G$ if $\Gamma$~is discrete and $\Gamma \backslash G$ has finite volume
\textup(which respect to the Haar measure on~$G$\textup). 
 \end{defn}

\begin{eg}
If $\Gamma$ is discrete and $\Gamma
\backslash G$ is compact, then $\Gamma$ is a lattice in~$G$, because
any compact Riemannian manifold obviously has finite volume.
\end{eg}

\begin{eg} \label{SL2Zlatt}
 $\SL(2,\integer)$ is a lattice in $\SL(2,\real)$. To see this, let
 \begin{equation} \label{SL2ZFundDom}
 \fund = \{\, z \in \hyperbolic^2 \mid \mbox{$|z|
\ge 1$ and $-1/2 \le \Re z \le 1/2$} \,\} 
 \end{equation}
 \csee{FundDomSL2R}. It is well known (though not obvious) that
$\fund$ is a fundamental domain for the action of $\SL(2,\integer)$
on~$\hyperbolic^2$ \csee{SL2ZinFundDom,SL2ZBdryFundDom}; it
therefore suffices to show that $\fund$ has finite volume, or, more
precisely, finite hyperbolic area. 

\begin{figure}[t]
 \includegraphics{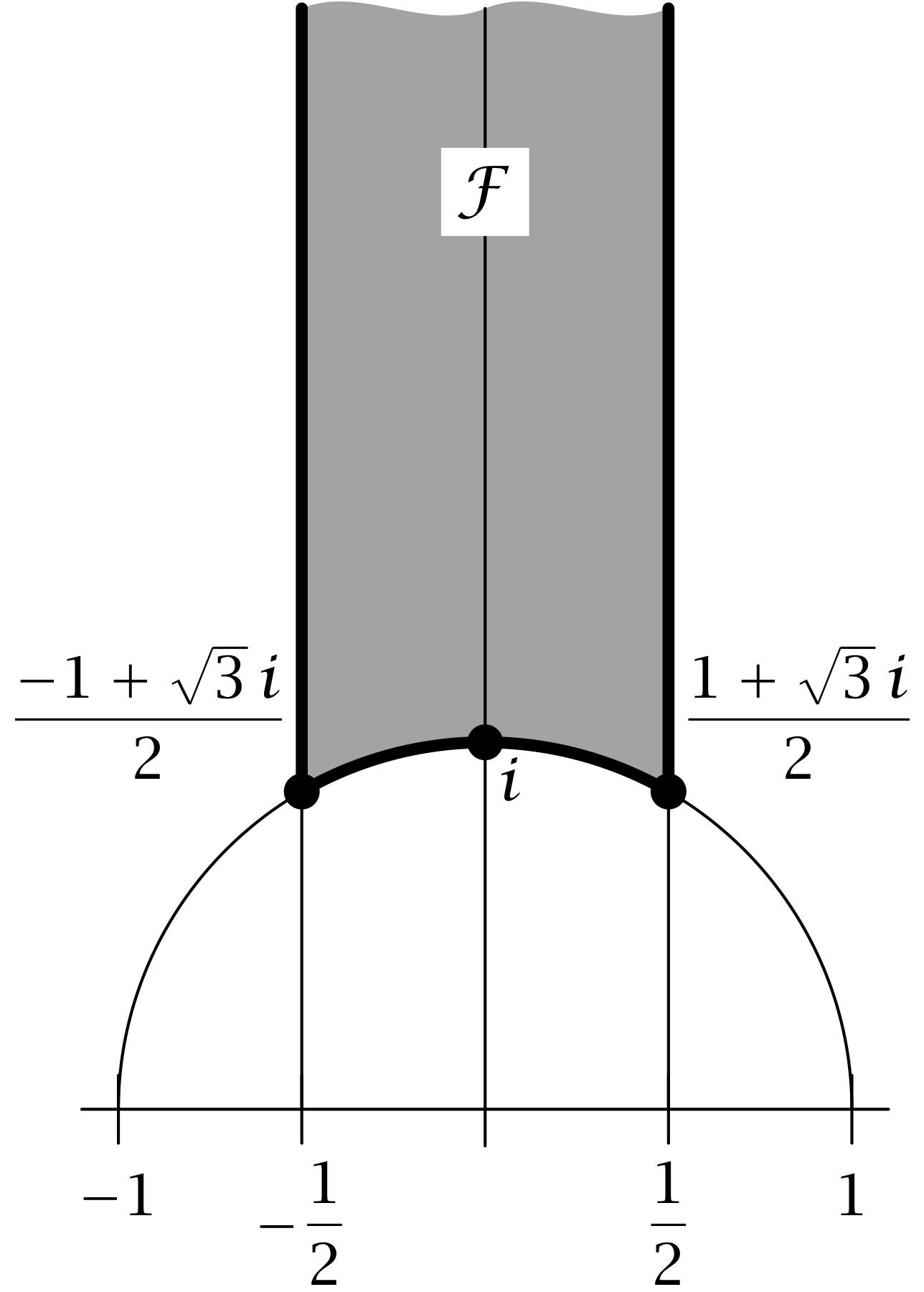}
 \caption{A fundamental domain~$\fund$ for $\SL(2,\integer)$ in
$\SL(2,\real)$.} \label{FundDomSL2R}
 \end{figure}
%texpreamble
%(" \usepackage[LY1]{fontenc}
% \usepackage[expert, LY1, mylucidascale]{mylucidabr} % I adjusted the scaling
% \usepackage{amsmath}
% \everymath{\displaystyle}
% ");
%defaultpen(  fontcommand("\normalfont") + fontsize(10) ); 
%
%from graph access *;
%
%size(0,3inch);
%pair i = (0,1);
%real h = sqrt(3)/2;
%pair a = (-1/2, h);
%pair b = (1/2, h);
%real top = 3;
%filldraw( (1/2,h) .. i .. (-1/2,h)--(-1/2,top){ENE}..{ENE}(0,top){ENE}..{ENE}(1/2,top)--cycle , gray(0.7) , invisible  );
%draw( (-1/2,h)--(-1/2,top) , linewidth(2) );
%draw( (1/2,h)--(1/2,top) , linewidth(2) );
%draw( (-1/2,h) .. i .. (1/2,h) , linewidth(2) );
%draw( (-1,0){N}..(-1/2,h) .. i .. (1/2,h) .. {S}(1,0));
%draw( (-1/2,0) -- (-1/2,h) );
%draw( (1/2,0) -- (1/2,h) );
%dotfactor = 12;
%dot( i ); label( "$i$", (0,1), SE );
%dot( a ); label( "$\frac{-1 + \sqrt{3} \, i}{2}$", a , NW );
%dot( b ); label( "$\frac{1 + \sqrt{3} \, i}{2}$", b , NE );
%real[] xticklist = {-1, -1/2, 1/2, 1};
%string labelfunc(real x){ 
%	if (x < -0.75){return "$-1$";} 
%	else if (x > 0.75) { return "$1$" ;}
%	else if (x < 0) { return "$-\frac{1}{2}$" ;}
%	else { return "$\frac{1}{2}$" ;}
%	}
%xaxis(-1.1, 1.1, Ticks(xticklist, ticklabel=labelfunc) );
%yaxis(-0.1, top, true);
%
%real fundw = 0.12, fundh = 2.5;
%fill( (-fundw, fundh - fundw) -- (fundw, fundh - fundw) -- (fundw, fundh + fundw) --  (-fundw, fundh + fundw) -- cycle, white );
%label( "$\mathcal{F}$",  (0,fundh) );

The hyperbolic area~$dA$ of an infinitesimal rectangle is the product
of its hyperbolic length and its hyperbolic width. If the Euclidean
length is $dx$ and the Euclidean width is $dy$, and the rectangle is
located at the point $x + iy$, then, by definition of the hyperbolic
metric, the hyperbolic length is $(dx)/(2y)$ and the hyperbolic width
is $(dy)/(2y)$. Therefore, 
 $$ dA = \frac{dx \, dy}{4 y^2} .$$
 Since $\Im z \ge \sqrt{3}/2$ for all $z \in \fund$,
we have
 $$ \vol(\fund)
 = \int_{x + iy \in \fund} \, dA
 \le \int_{\sqrt{3}/2}^\infty \int_{-1/2}^{1/2} \, \frac{dx \, dy}{4
y^2}
 =  \frac{1}{4} \int_{\sqrt{3}/2}^\infty \frac{1}{y^2} \, dy 
 < \infty .$$

Unfortunately, $\SL(2,\integer) \backslash \hyperbolic^2$ is not a
locally symmetric space, because $\SL(2,\integer)$ does not act
freely on~$\hyperbolic^2$ (so the quotient space is not a Riemannian
manifold). However, there are finite-index subgroups of
$\SL(2,\integer)$ that do act freely \ccf{torsionfree}, and these
provide interesting locally symmetric spaces.
 \end{eg}

Calculations similar to (but more complicated than)
\cref{SL2Zlatt} show:
 \begin{itemize}
 \item $\SL(n,\integer)$ is a lattice in $\SL(n,\real)$, and
 \item $\SO(p,q) \cap \SL(n,\integer)$ is a lattice in $\SO(p,q)$.
 \end{itemize}
 As in the example of $\SL(2,\integer) \backslash \hyperbolic^2$, the
hard part is to find a fundamental domain for $\Gamma \backslash G$
(or an appropriate approximation of a fundamental domain); then it is
not difficult to see that its volume is finite.
 These are special cases of the following general theorem, which
implies that every simple Lie group has a lattice.

\begin{thm}[{(\thmindex{arithmetic subgroups are lattices}{Arithmetic subgroups are lattices} \csee{arith->latt})}]
\label{Arith->Latt}
 Assume 
 \begin{itemize}
 \item $G = G_1 \times \cdots \times G_m$ is a product of simple Lie
groups,
 \item $G \subseteq \SL(\ell,\real)$, and 
 \item $G \cap \SL(\ell,\rational)$ is dense in~$G$.
 \end{itemize} Then $G_{\integer} = G \cap \SL(\ell,\integer)$ is a
lattice in~$G$.
 \end{thm}

Lattices constructed by taking the integer points of~$G$ in this way
are said to be \defit[arithmetic!subgroup]{arithmetic} \csee{ArithDefn}. (For most
simple Lie groups, these are the only lattices \csee{MargulisArith}.)
When $\ell$ is large, there is more than one way to embed $G$ in
$\SL(\ell,\real)$, and we will see that different embeddings can lead
to quite different intersections with $\SL(\ell,\integer)$. In particular, if $G$ is a noncompact, simple Lie group, then:
	\begin{itemize}
	\item By taking
an appropriate embedding of~$G$ in some $\SL(\ell,\real)$, we will
construct a lattice~$\Gamma$ in~$G$, such that $\Gamma \backslash G$ is \textbf{not} compact \csee{GHasNoncpctLatt}.
	\item By taking a different
embedding, we will construct a different lattice $\Gamma'$, such that
$\Gamma' \backslash G$ is compact \csee{GHasCpctLatt}.
	\end{itemize}

We will also see that algebraic properties of~$\Gamma$ influence the geometry of the corresponding locally symmetric space~$M$. 
In particular, the structure of~$\Gamma$ determines whether $M$ is compact or not.
(For example, the ``Godement Criterion'' \pref{GodementNoCpctFactor} implies that $M$~is compact if and only if every element of~$\Gamma$ is a diagonalizable matrix over~$\complex$.)
%More generally, we will see how group-theoretic
%properties of $\Gamma$ influence the large-scale structure of~$M$ \csee{LargeScaleSect}.
Much more generally, the following important theorem implies that
every geometric property of~$M$ is faithfully reflected in some
group-theoretic property of~$\Gamma$. 

\begin{thm}[(\thmindex{Mostow Rigidity}{Mostow Rigidity Theorem} \csee{MostowChap})] \label{MostowIrred}
 Let $M_1$ and $M_2$ be finite-volume locally symmetric spaces\/ \textup(not both $2$-dimensional\/\textup), such
that
 \begin{itemize}
 \item the universal covers of~$M_1$ and~$M_2$ are neither compact,
nor flat, nor reducible,
and
\item the volumes of $M_1$ and~$M_2$ are normalized\/ \textup(i.e., $\vol M_1 = \vol M_2 = 1$\textup).
 \end{itemize}
  If\/ $\pi_1(M_1) \iso \pi_1(M_2)$, then $M_1$ is isometric to~$M_2$.

 In fact, every homotopy equivalence is homotopic to an isometry.
 \end{thm}

The theorem implies that locally symmetric spaces
have no nontrivial deformations, which is why it is called a ``rigidity'' theorem:

\begin{cor}
Let $\{g_t\}$ be a continuous family of Riemannian metrics on a manifold~$M$ with\/ $\dim M > 2$, such that, for each~$t$:
	\begin{itemize}
	\item $(M,g_t)$ is a finite-volume locally symmetric space  whose universal cover is neither compact, nor flat, nor reducible,
	and
	\item $\vol(M,g_t) = 1$.
	\end{itemize}
Then $(M,g_t)$ is isometric to $(M,g_0)$, for every~$t$.
\end{cor}

\begin{defn}
 A locally symmetric space is \defit[irreducible!locally
symmetric space]{irreducible} if no \emph{finite} cover of~$M$ can be
written as a nontrivial cartesian product $M_1 \times M_2$.
 \end{defn}

 It is important to note that the universal cover of an irreducible
locally symmetric space need not be an irreducible symmetric space.
In other words, there can be lattices in $G_1 \times \cdots \times
G_n$ that are not of the form $\Gamma_1 \times \cdots \times
\Gamma_n$ \csee{SL(2Z[sqrt2])}. 

\begin{rem}
\Cref{MostowIrred} (and the corollary) can be generalized to the case where
only~$M_1$, rather than the universal cover of~$M_1$, is irreducible.
However, this requires the hypotheses to be strengthened: it
suffices to assume that no irreducible factor of~$M_1$ or~$M_2$ is
either compact or flat or $2$-dimensional. Furthermore, the conclusion needs to be weakened:
rather than simply multiplying by a single scalar to normalize
the volume, there can be a
different scalar on each irreducible factor of the universal cover.
\end{rem}

\begin{exercises}

\item
 Let
 \begin{itemize}
 \item $X$ be a simply connected symmetric space, 
 \item $\Gamma \backslash X$ be a locally symmetric space whose
universal cover is~$X$ (so $\Gamma$ is a discrete group of isometries
that acts freely and properly discontinuously on~$X$), and
 \item $\tau$ be an isometry of~$X$.
 \end{itemize}
 Show that if $\tau$ factors through to a well-defined map on $\Gamma
\backslash X$, then $\tau$ normalizes~$\Gamma$ (that is, $\tau
\gamma \tau^{-1} \in \Gamma$, for every $\gamma \in \Gamma\mk$).

\item \label{GrpCuspNotSymm}
 Define $g \colon \hyperbolic^2 \to \hyperbolic^2$ by $g(z) = z+1$. 
 \begin{enumerate}
 \item Show the geodesic symmetry~$\tau$ at~$i$ is given by
$\tau(z) = -1/z$.
 \item Show that $\tau$ does not normalize $\langle g \rangle$.
 \item Conclude that $\tau$ does not factor through to a well-defined
map on $\langle g \rangle \backslash \hyperbolic^2$.
 \end{enumerate}

\item
 Let
 \begin{itemize}
 \item $X$ be a simply connected symmetric space, and
 \item $\Gamma \backslash X$ be a locally symmetric space whose
universal cover is~$X$ (so $\Gamma$ is a discrete group of isometries
that acts freely and properly discontinuously on~$X$).
 \end{itemize}
 Show that $X$ is homogeneous if and only if the normalizer
$\nzer_G(\Gamma)$ is transitive on~$X$, where $G = \Isom(X)$.

\item Let $M = \Gamma \backslash G/K$ be a locally symmetric space,
and assume that $G$ has no compact factors. Show that if
$\nzer_G(\Gamma)/\Gamma$ is finite, then $\Isom(M)$ is finite.

\item Show that if $K$ is any compact subgroup of a Lie group~$G$,
then there is a unique (up to a scalar multiple) $G$-invariant Borel
measure~$\nu$ on $G/K$, such that $\nu(C) < \infty$, for every
compact subset~$C$ of $G/K$.

\item \label{lattice<>mu(X/Gamma)}
 Let
\begin{itemize}
 \item $K$ be a compact subgroup of a Lie group~$G$, and
 \item $\Gamma$ be a discrete subgroup of~$G$ that acts freely on
$G/K$.
 \end{itemize}
 Show that $\Gamma \backslash G$ has finite volume if and only if
$\Gamma \backslash G/K$ has finite volume.

\item \label{SL2ZinFundDom}
 Let $\Gamma = \SL(2,\integer)$, and define $\mathcal{F} \subset
\hyperbolic^2$ as in~\eqref{SL2ZFundDom}. Show, for each $p \in
\hyperbolic^2$, that there is some $\gamma \in \Gamma$ with
$\gamma(p) \in \mathcal{F}$.
 \hint{If $\Im \gamma(p) \le \Im p$ for all $\gamma \in \Gamma$, and
$-1/2 \le \Re p \le 1/2$, then $p \in \mathcal{F}$.}

\item \label{SL2ZBdryFundDom}
 Let $\Gamma = \SL(2,\integer)$, and define $\mathcal{F} \subset
\hyperbolic^2$ as in~\eqref{SL2ZFundDom}. Show, for $z,w \in
\mathcal{F}$, that if there exists $\gamma \in \Gamma$ with
$\gamma(z) = w$, then either $z = w$ or $z,w \in \bdry \mathcal{F}$.
 \hint{Assume $\Im w \le z$. Then $|\gamma_{2,1}z + \gamma_{2,2}| \le
1$. Hence $|\gamma_{2,1}| \in \{0,1\}$. If $|\gamma_{2,1}| = 1$ and
$\gamma_{2,2} \neq 0$, then $|\Re z| = 1/2$, so $z \in \bdry
\mathcal{F}$. If $|\gamma_{2,1}| = 1$ and
$\gamma_{2,2} = 0$, then $w = (az-1)/z$. Since $|\Re(1/z)|
\le |\Re z| \le 1/2$, and $|\Re w| \le 1/2$, we see that either $\Re
z = 1/2$ or $w = -1/z$.}

 \end{exercises}

\begin{notes}

 Either of Helgason's books \cite{HelgasonBookOld, HelgasonBook} is a
good reference for the geometric material on symmetric spaces and
locally symmetric spaces, the connection with simple Lie groups, and
much more. Lattices are the main topic of Raghunathan's book
\cite{RaghunathanBook}. 

\Cref{Arith->Latt} is a result of Borel and Harish-Chandra
\cite{BorelHarishChandra} that will be proved in \cref{SLnZLattChap,ReductionChap}.

\Cref{MostowIrred} combines work of Mostow
\cite{MostowRig}, Prasad \cite{PrasadRig}, and Margulis
\cite{Margulis-DiscGrpMot}. We will discuss it in  \cref{MostowChap}.

\Cref{SL2Zlatt} appears in many number theory texts, including \cite[\S7.1.2, pp.~77--79]{Serre-CourseArith}. Our hints
for \cref{SL2ZinFundDom,SL2ZBdryFundDom} are taken
from \cite[Prop.~4.4, pp.~181--182]{PlatonovRapinchukBook}.

\end{notes}

 %!TEX root = IntroArithGrps.tex

\mychapter{\texorpdfstring{Geometric Meaning of\\$\real$-rank and $\rational$-rank}%
	{Geometric Meaning of R-rank and Q-rank}}
\label{GeomIntroRank}

\prereqs{locally symmetric spaces (\cref{WhatisLocSymmChap}) and other differential geometry.}

This chapter, like the previous one, is motivational. It is not a
prerequisite for later chapters.

\section{Rank and real rank} \label{IntroRrankSect}

Let $X$ be a symmetric space \csee{symmDefn}. 
%(That is, for each point $p \in X$, there is an isometry $\phi$ of~$X$, such that the derivative $d \phi_p$ is $-\Id$ on the tangent space $T_p X$.) 
For example, $X$ could
be a Euclidean space $\real^n$, or a round sphere $S^n$, or a hyperbolic
space $\hyperbolic^n$, or a product of any combination of these.

As is well known, the rank of~$X$ is a natural number that
describes part of the geometry of~$X$, namely, the
dimension of a maximal flat.

\begin{defn}
 A \defit[flat!in a symmetric space]{flat} in~$X$ is a connected, totally geodesic,
flat submanifold of~$X$.
 \end{defn}

\begin{defn}
	\nindex{$\rank X$ = dimension of maximal flat}%
 $\rank X$ is the largest natural number~$r$, such that
$X$ contains an $r$-dimensional flat.
 \end{defn}

Let us assume that $X$ has no flat factors. (That is, the universal
cover of~$X$ is not isometric to a product of the form $Y \times
\real^n$. Mostly, we will be interested in the case where
$X$ also does not have any compact factors.)

Let $G = \Isom(X)^\circ$. Then $G$ acts transitively on~$X$, and there
is a compact subgroup~$K$ of~$G$, such that $X = G/K$. Because $X$ has no
flat factors, $G$ is a connected, semisimple, real Lie group with
trivial center (see~\S\ref{ConstructSymm}). (We remark that $G$ is
isomorphic to a closed subgroup of $\SL(\ell,\real)$, for some~$\ell$.)

The real rank can be understood similarly. It is an invariant of~$G$
that is defined algebraically \csee{RrankChap}, but it has the following geometric
interpretation.

\begin{thm} \label{Rrank-geometric}
	\nindex{$\Rrank G$ = maximal dimension of closed, simply connected flat}%
 $\Rrank G$ is the largest natural number~$r$, such that
$X$ contains a closed, \textbf{simply connected}, $r$-dimensional
flat.
 \end{thm}
 
\begin{warn}
By \emph{closed}, we simply mean that the flat contains all of its accumulation points, not that it is compact. (A closed, simply connected flat is homeomorphic to some Euclidean space~$\real^r$.)
\end{warn}

For example, if $X$ is compact, then every closed, totally
geodesic, flat subspace of~$X$ must be a torus,
not~$\real^n$, so $\Rrank G = 0$. On the other hand, if $X$
is not compact, then $X$ has unbounded geodesics (for
example, if $X$ is irreducible, then every geodesic goes to
infinity), so $\Rrank G \ge 1$. Hence:
 \centerline{$\Rrank G = 0
  \qquad\Leftrightarrow \qquad 
  \text{$X$~is compact}
 .$}
 Thus, there is a huge difference between $\Rrank G = 0$
and $\Rrank G > 0$, because no one would mistake a
compact space for a noncompact one.

\begin{rem}
 $\Rrank G = \rank X$ if and only if $X$ has no compact
factors.
 \end{rem}

There is also an important difference between $\Rrank G =
1$ and $\Rrank G > 1$. The following proposition is an important
example of this.

\begin{defn}
 $X$ is \defit{two-point homogeneous} if, whenever $(x_1,x_2)$ and
$(y_1,y_2)$ are two pairs of points in~$X$ with $d(x_1,x_2) =
d(y_1,y_2)$, there is an isometry~$g$ of~$X$ with $g(x_1) = y_1$ and
$g(x_2) = y_2$.
 \end{defn}

If $\Rrank G > 1$, then there exist maximal flats $H_1$ and~$H_2$ that
intersect nontrivially. On the other hand, there also exist some pairs
$x_1,x_2$, such that $\{x_1,x_2\}$ is not contained in the intersection
of any two (distinct) maximal flats. This establishes one direction of
the following result.

\begin{prop} \label{2ptHomog}
 Assume $X$ is noncompact and irreducible. The symmetric space~$X$ is
two-point homogeneous if and only if\/ $\Rrank G = 1$.
 \end{prop}

The following is an infinitesimal version of this result.

\begin{prop} \label{Rrank1->transitive}
 Assume $X$ is noncompact and irreducible. 
 The action of~$G$ on the set of unit tangent vectors
of~$X$ is transitive iff\/ $\Rrank G = 1$. 
 \end{prop}

\begin{cor}
 $\Rrank  \SO(1,n) = 1$.
 \end{cor}

\begin{proof}
 For $G = \SO(1,n)$, we have $X = \hyperbolic^n$. The stabilizer
$\SO(n)$ of a point in~$\hyperbolic^n$ acts transitively on the unit
tangent vectors at that point. So $G$ acts transitively on the unit
tangent vectors of~$X$.
 \end{proof}
 
More generally, it can be shown that $\Rrank \bigl( \SO(m,n) \bigr) = \min\{m,n\}$. Also,
$\Rrank \bigl( \SL(n,\real) \bigr) = n-1$. Although they may not be obvious geometrically, these real ranks are easy to calculate from the algebraic
definition that will be given in \cref{RrankChap}.

\begin{rem}
For every~$r$, there is a difference between $\Rrank G =
r$ and $\Rrank G > r$, but this difference is less
important as $r$~grows larger: the three main cases are
$\Rrank G = 0$, $\Rrank G = 1$, and $\Rrank G \ge 2$.
(This is analogous to the situation with smoothness
assumptions: countless theorems require a function to be
$C^0$ or $C^1$ or~$C^2$, but far fewer theorems require a
function to be, say,~$C^7$, rather than only~$C^6$.)
\end{rem}

%Real rank has the following implication for the geometry of
%$\Gamma \backslash X$:
%
%\begin{prop}
% Assume $X$ has no compact factors. Then:
% \begin{enumerate}
% \item There is a dense geodesic in $\Gamma \backslash X$.
% \item The geodesic flow on the unit tangent bundle $T^1(\Gamma
%\backslash X)$ has a dense orbit if and only if\/ $\Rrank G = 1$.
% \end{enumerate}
% \end{prop}

\begin{exercises}

\item Show $\Rrank(G_1 \times G_2) = \Rrank G_1 + \Rrank G_2$.

\item Assume $\Rrank G = 1$. Show $X$ is irreducible if and only if
$X$~has no compact factors.

\item Show that if $X$ is reducible, then $X$ is \emph{not} two-point
homogeneous. (Do not assume the fact about maximal flats that
was mentioned, without proof, before \cref{2ptHomog}.)

 \end{exercises}

\section{\texorpdfstring{$\rational$}{Q}-rank} \label{QrankIntroSect}

Now let $\Gamma \backslash X$ be a locally symmetric space
modeled on~$X$, and assume that
$\Gamma \backslash X$ has
finite volume. Hence, $\Gamma$ is a (torsion-free) discrete subgroup
of~$G$, such that $\Gamma \backslash G$ has finite volume;
in short, $\Gamma$ is a \defit[lattice!subgroup]{lattice} in~$G$.

The real rank depends only on~$X$, so it is not affected by the choice of
a particular lattice~$\Gamma$. We now describe an analogous
algebraically defined invariant, $\Qrank\Gamma$, that does depend
on~$\Gamma$, and therefore distinguishes between some of the various locally
homogeneous spaces that are modeled on~$X$. We will mention some of the
geometric implications of $\rational$-rank, leaving a more detailed
discussion to later chapters.

\begin{thm}[\csee{DivTorusSect,LargeScaleSect}] \label{QrankFlats} \ 
\noprelistbreak
\begin{enumerate}
 \item \label{QrankFlats-exists}
 $\Qrank\Gamma$ is the largest natural number~$r$, such
that some finite cover of\/ $\Gamma \backslash X$ contains a closed, simply connected,
$r$-dimensional flat.
\item \label{QrankFlats-BddDist}
$\Qrank\Gamma$ is the smallest natural number~$r$, for which there
exists collection of finitely many closed, $r$-dimensional flats, such that all of\/ $\Gamma \backslash X$ is within a bounded distance of the union of these flats. 
\end{enumerate}
 \end{thm}

\begin{rem} \label{QrankPossiblesMentioned}
 It is clear from \fullcref{QrankFlats}{exists}
that $\Qrank\Gamma$ always exists (and is finite). Furthermore,
$0 \le \Qrank\Gamma \le \Rrank G$. Although not so obvious, it can be shown that the extreme values are always
attained: there are lattices $\Gamma_c$ and~$\Gamma_s$ in~$G$ with
$\Qrank\Gamma_c = 0$ and $\Qrank\Gamma_s = \Rrank G$ \csee{GHasCpctLatt,Qrank=Rrank}. So it is
perhaps surprising that there may be gaps in between. 
(For example,  if $G \iso \SO(2,n)$, with $n \ge 5$, and
$n$~is odd, then $\Rrank G = 2$, but \cref{QrankGap} shows
there does not exist a lattice~$\Gamma$ in~$G$, such that $\Qrank\Gamma = 1$.)
 \end{rem}

\begin{eg}[\csee{QrankEg}] \label{QrankIntroEgs}
 From the algebraic definition, which will appear in
\cref{QrankChap}, it is easy to calculate 
 $$\Qrank \bigl( \SO(m,n)_{\integer} \bigr) = \min\{m,n\} = \Rrank \bigl(
\SO(m,n) \bigr)$$
 and 
 $$\Qrank \bigl( \SL(n,\integer) \bigr) = n-1 = \Rrank \bigl(
\SL(n,\real) \bigr) .$$
 \end{eg}

As for the real rank, the biggest difference is between spaces where the
invariant is zero and those where it is nonzero, because this is again
the distinction between a compact space and a noncompact one:

\begin{thm}[\csee{Qrank0Ex}] \label{Qrank0<>cocpct}
 $\Qrank\Gamma = 0$ iff\/ $\Gamma
\backslash X$~is compact.
 \end{thm}

\fullCref{QrankFlats}{BddDist} implies that the $\rational$-rank of~$\Gamma$ is
directly reflected in the large-scale geometry of~$\Gamma \backslash X$,
as described by the asymptotic cone of~$\Gamma \backslash X$.
 Intuitively, the asymptotic cone of a metric space is obtained
by looking at it from a large distance. For example, if $\Gamma
\backslash X$ is compact, then, as we move farther away, the manifold
appears smaller and smaller (see the illustration below). % !!! \csee{shrinkcpct}. 
In the limit, the manifold shrinks to a point.
\begin{align} \label{shrinkcpct} % cref calls this an equation, instead of a figure !!!
 \lower0.75\baselineskip\hbox{\includegraphics{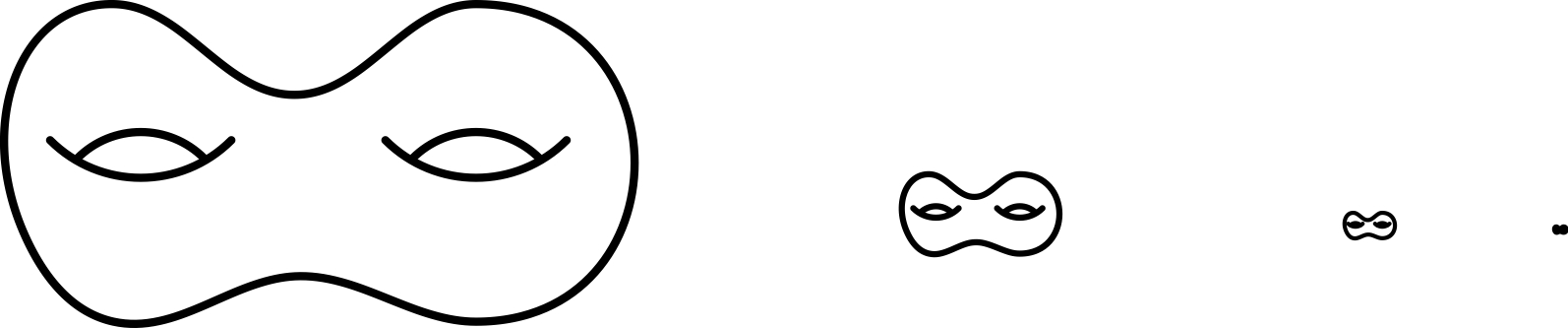}}
\end{align}
%%%\begin{figure}[h]
%%% \includegraphics{PDF/shrinkcpct.jpg}
%%% \caption{Looking at a compact manifold from farther and farther away.}
%%%% \label{shrinkcpct}
%%% \end{figure}
%texpreamble
%("  \usepackage{amsmath}
% \usepackage[LY1]{fontenc}
% \usepackage[expert,LY1,mylucidascale]{mylucidabr}
% ");
%defaultpen(  fontcommand("\normalfont") + fontsize(10) ); 
%
%from graph access *;
%unitsize(0.39cm);
%
%currentpen = linewidth(1.25);
%
%void hole(real x, real s) {
%	draw( (x+s,s){SW}..{NW}(x-s, s) );
%	draw( (x+s*0.7,s*0.8){NW}..{SW}(x-s*0.7, s*0.8) );
%	}	
%
%
%void M(real m, real s, real thickness){
%	currentpen = linewidth(thickness);
%	draw( shift(m,0)* (scale(s)*((0,-1){W}..(-2,-0.5)..(-4,-1)..(-5,0)..(-4,2.5)..(-2,1.5)..{E}(0,2.5)
%		..{W}(0,-1) )));
%	hole(m,s); hole(m-3.7*s, s);
%	}
%
%M(0, 1, 1); 
%M(6, 0.25, 0.75);
%M(10, 0.08, 0.5);
%M(12, 0.02, 0.5);

An intuitive understanding is entirely sufficient for our purposes here,
but, for the interested reader, we provide a more formal definition.

\begin{defn} \label{AsympConeDefn}
 The \defit{asymptotic cone} of a metric space $(M,d)$ is the
limit space 
 $$\lim_{\epsilon \to 0^+} \bigl( (M, \epsilon d), p \bigr) ,$$
if the limit exists. Here, $p$ is an arbitrary (but fixed!) point of~$M$,
and the limit is with respect to Gromov's Hausdorff distance. (Roughly
speaking, a large ball around~$p$ in $(M, \epsilon d)$ is $\delta$-close
to being isometric to a large ball around a certain (fixed) point~$p_0$
in the limit space $(M_0,d_0)$.)
 \end{defn}

\begin{egs} \  \label{TanConeEgs}
\noprelistbreak
 \begin{enumerate}
 \item If $\Gamma \backslash X$ is compact, then the
asymptotic cone of $\Gamma \backslash X$ is a
point, as is illustrated in \pref{shrinkcpct}. This point is a $0$-dimensional simplicial
complex, which is a geometric manifestation of the
fact that $\Qrank\Gamma = 0$.
 \item \label{TanConeEgs-rank1} 
If $\Rrank G = 1$, and $\Gamma \backslash X$
is not compact, then, as is well known, $\Gamma \backslash X$ has
finitely many cusps. The asymptotic cone of a cusp
is a ray, so the asymptotic cone of $\Gamma
\backslash X$ is a ``\term[star (of finitely many rays)]{star}'' of finitely many rays emanating
from a single vertex \csee{shrinkcusps}. Therefore, the asymptotic cone
of $\Gamma \backslash X$ is a $1$-dimensional simplicial
complex. This manifests the fact that $\Qrank\Gamma = 1$.
 \end{enumerate}
 \end{egs}

\begin{figure}[h]
 \centerline{\includegraphics{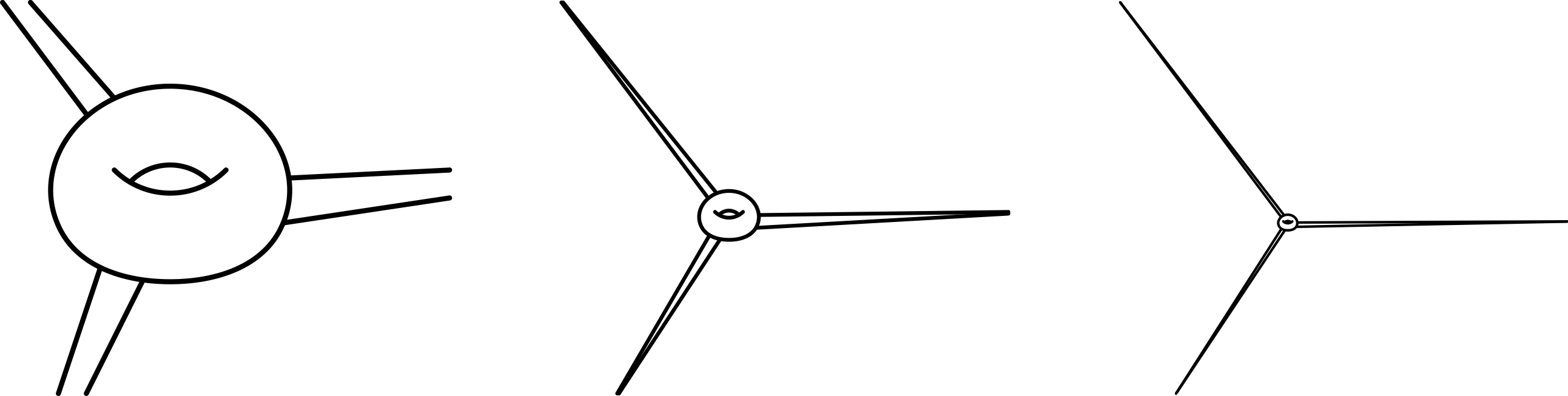}}
 \caption{Looking at a manifold with cusps from farther and farther
away.}
 \label{shrinkcusps}
 \end{figure}
%texpreamble
%("  \usepackage{amsmath}
% \usepackage[LY1]{fontenc}
% \usepackage[expert,LY1,mylucidascale]{mylucidabr}
% ");
%defaultpen(  fontcommand("\normalfont") + fontsize(10) ); 
%
%from graph access *;
%unitsize(0.39cm);
%
%currentpen = linewidth(1.25);
%
%void hole(real x, real s) {
%	draw( (x+s,s){SW}..{NW}(x-s, s) );
%	draw( (x+s*0.7,s*0.8){NW}..{SW}(x-s*0.7, s*0.8) );
%	}	
%
%
%void M(real m, real s, real thickness){
%	currentpen = linewidth(thickness);
%	draw( (m-3, 4)..(m,0) ); draw( (m-3+0.5*s^2, 4)..(m+s,0) );
%	draw( (m-2, -3)..(m-s,0) ); draw( (m-2+0.5*s^2, -3)..(m,0) );
%	draw( (m+5, s)..(m,0.75*s) ); draw( (m+5, s - 0.5*s^2)..(m,0.75*s-s) );
%
%	fill( shift(m,0)* (scale(s)*((0,-1){W}..(-2,0)..(0,2.5)..(2,0)..cycle )), white);
%	draw( shift(m,0)* (scale(s)*((0,-1){W}..(-2,0)..(0,2.5)..(2,0)..cycle )));
%	hole(m,s);
%	}
%
%M(0, 1, 1); 
%M(10, 0.25, 0.75);
%M(20, 0.08, 0.5);

\begin{thm}[\csee{LargeScaleRems}] \label{HattoriThm}
 The asymptotic cone of\/ $\Gamma \backslash X$ is a
simplicial complex whose dimension is\/ $\Qrank\Gamma$.
 \end{thm}

\begin{eg} \label{TanConeSL3Z}
 Let $G = \SL(3,\real)$ and $\Gamma = \SL(3,\integer)$. From
\cref{HattoriThm}, we see that the asymptotic cone of
$\Gamma \backslash G/K$ is a $2$-dimensional simplicial complex. In
fact, it turns out to be (isometric to) the sector
 $$ \bigset{ (x,y) \in \real^2 }{ 0 \le y \le \frac{\sqrt{3}}{2} x }
 .$$
 (It is not a coincidence that this sector is a Weyl chamber of
the Lie algebra $\LieSL(3,\real)$.)
 \end{eg}

\begin{rems} \ 
\noprelistbreak
 \begin{enumerate}
 \item If $\Qrank \Gamma = 1$, then the asymptotic cone of
$\Gamma \backslash X$ is a star of finitely many rays emanating from the
origin \fullccf{TanConeEgs}{rank1}. Note that this intersects the unit sphere
in finitely many points.
 \item In general, if $\Qrank \Gamma = k$, then  the unit sphere
contains a certain simplicial complex~$\mathcal{T}_\Gamma$ of dimension
$k - 1$, such that the asymptotic cone of $\Gamma \backslash X$
is the union of all the rays emanating from the origin that pass
through~$\mathcal{T}_\Gamma$.
 \item For $\Gamma = \SL(3,\integer)$, the simplicial
complex~$\mathcal{T}_\Gamma$ is a single edge \ccf{TanConeSL3Z}. In
general, the \term{Tits building}~$\mathcal{T}_G$ is a certain
simplicial complex defined from the parabolic $\rational$-subgroups
of~$G$, and $\mathcal{T}_\Gamma$ can be obtained from~$\mathcal{T}_G$ by
modding out the action of~$\Gamma$.
\item The asymptotic cone is also known as ``\index{tangent cone at infinity|indsee{asymptotic cone}}{tangent cone at infinity}\zz.''
 \end{enumerate}
 \end{rems}

\begin{rem} \label{BScohodim}
Although we will not prove this, the $\rational$-rank is directly reflected in the cohomology of $\Gamma \backslash X$. Namely, let $c$ be the cohomological dimension of $\Gamma \backslash X$. Because 
$\Gamma \backslash X$ is a manifold of dimension $\dim X$, we have 
 $c = \dim X$ if and only if $\Gamma \backslash X$ is compact. So the deficiency $\dim
X - c$ is, in some sense, a measure of how far $\Gamma \backslash X$ is
from being compact. This measure is precisely $\Qrank\Gamma$ (if $X$ has no compact factors).
\end{rem}

%\begin{thm} \label{}
% Assume $X$ has no compact factors. Then the cohomological dimension of
%$\Gamma \backslash X$ is $(\dim X) - \bigl(\Qrank\Gamma \bigr)$.
% \end{thm}

\begin{exercises}

\item \Cref{Qrank0<>cocpct} states that if $\Gamma
\backslash X$ is compact, then $\Qrank \Gamma = 0$.
	\begin{enumerate}
	\item Prove this directly from \fullcref{QrankFlats}{exists}.
	\item Prove this directly from \fullcref{QrankFlats}{BddDist}.
	\end{enumerate}

\end{exercises}

\begin{notes}

Helgason's book \cite{HelgasonBook} has a thorough treatment of
rank and $\real$-rank.

\fullCref{QrankFlats}{BddDist} was proved by B.\,Weiss \cite{Weiss-Qrank}.

\Cref{Qrank0<>cocpct} was proved  for arithmetic lattices by
Borel and Harish-Chandra \cite{BorelHarishChandra} and, independently, by
Mostow and Tamagawa \cite{MostowTamagawa}. For non-arithmetic lattices, this theorem is part of the \emph{definition} of $\rational$-rank.

A more precise version of \cref{HattoriThm} (providing a
description of the geometry of the simplicial complex) was proved by
Hattori \cite{Hattori}. 
Proofs also appear in \cite{JiMacpherson-GeomCpct} and~\cite{Leuzinger-TitsGeom}.

 \Cref{BScohodim} is due to Borel and Serre
\cite{BorelSerre-corners}.

\end{notes}

 %!TEX root = IntroArithGrps.tex

\mychapter{Brief Summary}
\label{SummaryChap}

\makemarksfalse % don't let the sections in this chapter change the running heads

\prereqs{none for most of the non-proof material.}

This book is about arithmetic subgroups, and other lattices, in semisimple Lie groups.
Given a lattice~$\Gamma$ in a semisimple Lie group~$G$, we will investigate both the algebraic structure of~$\Gamma$, and properties of the corresponding homogeneous space $G/\Gamma$. We will also study the close relationship between~$G$ and~$\Gamma$. For example, we will see that~$G$ is essentially the only semisimple group in which $\Gamma$ can be embedded as a lattice (``\thmindex{Mostow Rigidity}{Mostow Rigidity Theorem}''), and, conversely, we will usually be able to make a list of all the lattices in~$G$ (``\thmindex{Margulis!Arithmeticity}{Margulis Arithmeticity Theorem}''). 

%\begin{rem}
%Geometers interested in the locally symmetric space $\Gamma
%\backslash G/K$ usually place additional
%restrictions on $G$ and~$\Gamma$:
% \begin{itemize}
% \item Geometers assume that the center of $G$ is trivial,
%for otherwise $G$ does not act \emph{faithfully} as a group
%of isometries of the symmetric space $G/K$.
% \item Geometers assume that $\Gamma$ is torsion-free (that
%is, that $\Gamma$ has no nontrivial elements of finite
%order), for otherwise $\Gamma$ does not act freely on $G/K$.
% \end{itemize}
%\end{rem}

This chapter provides a very compressed outline of the material in this book. To help keep it brief, let us assume, for the remainder of the chapter, that
	$$ \text{\mathversion{bold}\bf$G$ is a noncompact, simple Lie group, and $\Gamma$ is a lattice in~$G$} .$$
This means \csee{LatticeDefn}:
	\begin{itemize}
	\item $\Gamma$ is a discrete subgroup of~$G$, 
	and 
	\item the homogeneous space $G/\Gamma$ has finite volume (with respect to the Haar measure on~$G$).
	\end{itemize} 
(If $G/\Gamma$ is compact, which is a very important special case, we say $\Gamma$ is \defit[cocompact subgroup]{cocompact}.)

\notocsection*{\cref{IntroductionPart}. Introduction}
All three chapters in this part of the book are entirely optional; none of the material will be needed later (although some examples and remarks do refer back to it). \Cref{WhatisLocSymmChap,GeomIntroRank} provide geometric motivation for the study of arithmetic groups, by explaining the connection with locally symmetric spaces. The present chapter (\cref{SummaryChap}) is a highly condensed version of the entire book.

\notocsection*{\cref{FundamentalsPart}. Fundamentals}
This part of the book presents definitions and other foundational material for the study of arithmetic groups.

\subsection*{\cref{BasicLatticesChap}. Basic Properties of Lattices} 
This chapter presents a few important definitions, including the notions of \emph{lattice subgroups}, \emph{commensurable subgroups}, and \emph{irreducible lattices}. It also proves a number of fundamental algebraic and geometric consequences of the assumption that $\Gamma$ is a lattice, including the following.

\smallbreak

\pref{GammaUnip->notcpct} Recall that an element~$u$ of $\SL(n,\real)$ is \emph{unipotent} if its characteristic polynomial is $(x-1)^n$ (or, in other words, its only eigenvalue is~$1$).
If $G/\Gamma$ is compact, then $\Gamma$~does not have any nontrivial unipotent elements.
This is proved by combining the Jacobson-Morosov Lemma \pref{JacobsonMorosov} with the observation that if a sequence $c_i \Gamma$ leaves all compact sets, and $U$~is a precompact set in~$G$, then, after passing to a subsequence, the sets $Uc_1 \Gamma, Uc_2 \Gamma, \ldots$ are all disjoint. However, $G/\Gamma$ has finite volume, so it cannot have infinitely many disjoint open sets that all have the same volume.

\smallbreak

\pref{GammaNotinConn} (\emph{Borel Density Theorem}) $\Gamma$ is not contained in any connected, proper, closed subgroup of~$G$. 
Assuming that $G/\Gamma$ is compact, the key to proving this is to note that if $\rho \colon G \to \GL(m,\real)$ is any continuous homomorphism, $u$~is any unipotent element of~$G$, and $v \in \real^m$, then the coordinates of the vector $\rho(u^k) v$ are polynomial functions of~$k$. However, if $G/\Gamma$ is compact, and $v$~happens to be $\rho(\Gamma)$-invariant, then the coordinates are all bounded. Since every bounded polynomial is constant, we conclude that every $\rho(\Gamma)$-invariant vector is $\rho(G)$-invariant. From this, the desired conclusion follows by looking at the action of~$G$ on exterior powers of its Lie algebra.

\smallbreak

\pref{GammaFinPres}  $\Gamma$ is finitely presented. When $G/\Gamma$ is compact, this follows from the fact that the fundamental group of any compact manifold is finitely presented. For the noncompact case, it follows from the existence of a nice fundamental domain for the action of~$\Gamma$ on~$G$ (which will be explained in \cref{ReductionChap}).

\smallbreak

\pref{torsionfree} (\emph{Selberg's Lemma}) $\Gamma$ has a torsion-free subgroup of finite index. For example, if $\Gamma = \SL(3,\integer)$, then the desired torsion-free subgroup can be obtained by choosing any prime $p \ge 3$, and taking the matrices in~$\Gamma$ that are congruent to the identity matrix, modulo~$p$.

\smallbreak

\pref{GammaResidFinite} $\Gamma$ is residually finite. For example, if $\Gamma = \SL(3,\integer)$, then no nontrivial element of~$\Gamma$ is in the intersection of the finite-index subgroups used in the preceding paragraph's % still preceding paragraph @@@
proof of Selberg's Lemma.

\smallbreak

\pref{FreeInGamma} (\emph{Tits Alternative}) $\Gamma$ contains a nonabelian free subgroup. This is proved by using the \emph{Ping-Pong Lemma} \pref{PingPong}, which, roughly speaking, states that if homeomorphism~$a$ contracts all of the space toward one point, and homeomorphism~$b$ contracts all of the space toward a different point, then the group generated by $a$ and~$b$ is free.

\smallbreak

\pref{MooreErgBasicThm} (\emph{Moore Ergodicity Theorem}) If $H$ is any noncompact, closed subgroup of~$G$, then every real-valued, $H$-invariant, measurable function on $G/\Gamma$ is constant (a.e.). 
The general case will be proved in \cref{MooreErgPfSect},
%The proof of this (and other theorems of this type) will be discussed in \cref{ErgodicChap} (Ergodic Theory).
but suppose, for example, that $G = \SL(2,\real)$, $H = \{a^s\}$ is the group of diagonal matrices, and $f$ is an $H$-invariant function that, for simplicity, we assume is uniformly continuous. If we let $\{u^t\}$ be the group of upper-triangular matrices with $1$'s on the diagonal, then we have
	$$ u^t \cdot f = u^t a^s \cdot f = a^s u^{e^{-s} t} \cdot f \stackrel{s \to \infty}{\ \longrightarrow\ } a^s u^0 \cdot f = f ,$$
so $f$~is invariant under~$\{u^t\}$. Similarly, it is also invariant under the group of lower-triangular matrices. So $f$ is $G$-invariant, and therefore constant. 
%A similar argument applies when $f$~is merely measurable, not uniformly continuous.

\subsection*{\cref{ArithGrpsChap}. What is an Arithmetic Group?} 
Roughly speaking, an \emph{arithmetic subgroup}~$G_\integer$ of~$G$ is obtained by embedding $G$ in some $\SL(\ell,\real)$, and taking the resulting set of integer points of~$G$. That is, $G_\integer$ is the intersection of~$G$ with $\SL(\ell,\integer)$. However, in order for~$G_\integer$ to be called an arithmetic subgroup, the embedding $G \hookrightarrow \SL(\ell,\integer)$ is required to satisfy a certain technical condition (``defined over~$\rational$'') \csee{DefdQDefn}.

\smallbreak

\pref{arith->latt} Every arithmetic subgroup of~$G$ is a lattice in~$G$. This fundamental fact will be proved in \cref{SLnZLattChap,ReductionChap}.

\smallbreak

\pref{MargulisArith} (\emph{Margulis Arithmeticity Theorem})
Conversely, if $G$ is neither $\SO(1,n)$ nor $\SU(1,n)$, then every lattice in~$G$ is an arithmetic subgroup. Therefore, in most cases, ``arithmetic subgroup'' is synonymous with ``lattice\zz.''
This amazing theorem will be proved in \cref{MargArithPf}.

It is a hugely important result. The definition of ``lattice'' is quite abstract, but a fairly explicit list of all the lattices in~$G$ can be obtained by combining this theorem with the classification of arithmetic subgroups that will be given in \cref{ArithClassicalChap}.

\smallbreak

\pref{GodementCriterion} (\emph{Godement Compactness Criterion}) 
$G/G_\integer$ is compact if and only if the identity element is the only unipotent element of~$G_\integer$. The direction ($\Rightarrow$) is very elementary and was proved in the previous chapter \see{GammaUnip->notcpct}. The converse uses the same main idea, combined with the simple observation that if a polynomial has integer coefficients, and all of its roots are close to~$1$, then all of its roots are exactly equal to~$1$.

\smallbreak

\pref{RestrictScalarsSect} The embedding of~$G$ in $\SL(\ell,\real)$ is not at all unique, and different embeddings can yield quite different arithmetic subgroups~$G_\integer$. One very important method of constructing non-obvious embeddings is called \emph{Restriction of Scalars}. It starts by choosing a field~$F$ that is a finite extension of~$\rational$. If we think of~$F$ as a vector space over~$\rational$, then it can be identified with some~$\rational^n$, in such a way that the ring~$\ints$ of algebraic integers of~$F$ is identified with~$\integer^n$. This implies that the group $G_{\ints}$ is isomorphic to $G'_\integer$, where $G'$ is a semisimple group that has~$G$ as one of its factors. Therefore, this method allows arithmetic subgroups to be constructed not only from ordinary integers, but also from algebraic integers.

\subsection*{\cref{EgArithGrpsChap}. Examples of Arithmetic Groups} 
This chapter explains how to construct many arithmetic subgroups of $\SL(2,\real)$, %\csee{ArithLattSL2}, 
$\SO(1,n)$, %\csee{ArithSO1nSect}, 
and $\SL(n,\real)$,
% \csee{NonLattinSL3Sect,CocpctLattSL3R,LattSlnRSect}, 
by using unitary groups and quaternion algebras (and other division algebras). (Restriction of scalars is also used for some of the cocompact ones.) It will be proved in \cref{ArithClassicalChap} that these fairly simple constructions actually produce all of the arithmetic subgroups of these groups.

\pref{NonArithSO1nSect} There exist non-arithmetic lattices in $\SO(1,n)$ for every~$n$. This was proved by M.\,Gromov and I.\,Piatetski-Shapiro. They ``glued together'' two arithmetic lattices to create a ``hybrid'' lattice that is not arithmetic.

\subsection*{\cref{SLnZLattChap}. $\SL(n,\integer)$ is a lattice in $\SL(n,\real)$} 
This chapter explains two different proofs of the fundamental fact (already mentioned in \cref{arith->latt}) that $G_\integer$ is a lattice in~$G$, in the illustrative special case where $G = \SL(n,\real)$ and $G_\integer = \SL(n,\integer)$.

\smallbreak

The first proof is quite short and elementary, and is presented fairly completely. It constructs a nice set that is (approximately) a fundamental domain for the action of~$\Gamma$ on~$G$.
%, by using the Iwasawa decomposition $G = KAN$. This proof 
The key notion is that of a \emph{Siegel set}.  We begin with the Iwasawa decomposition $G = KAN$.%
\noprelistbreak
	\begin{itemize}
	\item $K = \SO(n)$ is a maximal compact subgroup of~$G$.
	\item The group~$A$ of diagonal matrices in~$G$ is isomorphic to~$\real^{n-1}$, so we can think of it as a real vector space. Under this identification, the ``simple roots'' are linear functionals $\alpha_1, \ldots,\alpha_{n-1}$ on~$A$. Choose any $t \in \real$, and let
		$$ A_t = \{\, a \in A \mid \text{$\alpha_i(a) \ge t$ for all~$i$} \,\} ,$$
	so $A_t$ is a polyhedral cone in~$A$.
	
	\item $N$ is the group of upper-triangular matrices with $1$'s on the diagonal.
	\item Choose any compact subset~$N_0$ of~$N$.
	\end{itemize}
Then the product $\Siegel = N_0 \, A_t \,K$ is a Siegel set \csee{SiegelSLnZSect}. It depends on the choice of~$t$ and~$N_0$.

A straightforward calculation shows that every Siegel set has finite volume \csee{SiegelSLnRFinMeas}. It is also not terribly difficult to find a Siegel set~$\Siegel$ with the property that $G_\integer  \cdot \Siegel = G$ \csee{SiegelFundDomSLnZ}. This implies that $G/G_\integer$ has finite volume, so $G_\integer$ is a lattice in~$G$, or, in other words, $\SL(n,\integer)$ is a lattice in $\SL(n,\real)$.

Unfortunately, some difficulties arise when generalizing this method to other groups, because it is more difficult to use Siegel sets to construct an appropriate fundamental domain in the general case. The main ideas will be explained in \cref{ReductionChap}.

\smallbreak

So we also present a different proof that is much easier to generalize \csee{SLNZISLATTSlickSect}. Namely, the general case is quite easy to prove if one accepts the following key fact that was proved by Margulis: If 
	\begin{itemize}
	\item $u^t$ is any unipotent $1$-parameter subgroup of\/ $\SL(n,\real)$,
	and
	\item $x \in \SL(n,\real)/\SL(n,\integer)$,
	\end{itemize}
then there is a compact subset~$C$ of\/ $\SL(n,\real)/\SL(n,\integer)$, and some $\epsilon > 0$, such that at least $\epsilon \%$ of the orbit $\{u^t x\}_{t \in \real}$ is in the set~$C$ \csee{DaniMargulisUnipReturns}.

\notocsection*{\cref{ConceptsPart}. Important Concepts}
This part of the book explores several fundamental ideas that are important not only for their applications to arithmetic groups, but much more generally.

\smallbreak

\subsection*{\cref{RrankChap}.  Real rank}
This chapter defines the real rank of~$G$, which is an important invariant in the study of semisimple Lie groups. It also describes some consequences of assuming that the real rank is at least two, and presents the definition and basic structure of the minimal parabolic subgroups of~$G$.

\smallbreak

\subsection*{\cref{QrankChap}.  $\rational$-rank}
This chapter, unlike the others in this part of the book, discusses a topic that is primarily of interest in the theory of arithmetic groups (and related algebraic groups). Largely parallel to \cref{RrankChap}, it defines the $\rational$-rank of~$\Gamma$, describes some consequences of assuming that the $\rational$-rank is at least two, and presents the definition and basic structure of the minimal parabolic $\rational$-subgroups of~$G$.

\smallbreak

\subsection*{\cref{QuasiChap}. Quasi-isometries}
Any finite generating set~$S$ for~$\Gamma$ yields a metric~$d_S$ on~$\Gamma$: the distance from~$x$ to~$y$ is the minimal number of elements of~$S$ that need to be multiplied together to obtain $x^{-1} y$. Unfortunately, this ``word metric'' is not canonical, because it depends on the choice of the generating set~$S$. However, it is  well-defined up to a bounded factor, so, to get a geometric object that is uniquely determined by~$\Gamma$, we consider two metric spaces to be equivalent (or \emph{quasi-isometric}) if there is a map between them that only distorts distances by a bounded factor \csee{QuasiIsomDefn}. 

\smallbreak

\pref{CocpctActIsQI} Some quasi-isometries arise from cocompact actions: it is not difficult to see that if $\Gamma$ acts cocompactly, by isometries on a (nice) space~$X$, then there is a quasi-isometry from~$\Gamma$ to~$X$. Thus, for example, any cocompact lattice in $\SO(1,n)$ is quasi-isometric to the hyperbolic space~$\hyperbolic^n$.

\smallbreak

\pref{GromovHyperGrpsSect}
$\Gamma$ is \emph{Gromov hyperbolic} if and only if $\Rrank G = 1$ and $\Gamma$~is compact, except that all lattices in $\SL(2,\real)$ are hyperbolic, not only the cocompact ones. One direction is a consequence of the well-known fact that $\integer \times \integer$ is not contained in any hyperbolic group. The other direction (for the cocompact case) is a special case of the fact that the fundamental group of any closed manifold of strictly negative sectional curvature is hyperbolic.

\smallbreak

\subsection*{\cref{UnitaryRepChap}. Unitary representations}
This chapter presents some basic concepts in the theory of unitary representations, the study of group actions on Hilbert spaces. The Moore Ergodicity Theorem \pref{MooreErgBasicThm} is proved in \cref{MooreErgPfSect}, and the ``induced representations'' defined in \cref{InducedRepSect} will be used in \cref{LattInTSect} to prove that $\Gamma$ has Kazhdan's Property~$(T)$ if $\Rrank G \ge 2$.

\smallbreak

\pref{DecayMatCoeffSimple} (\emph{Decay of matrix coefficients}) If $\pi$ is a continuous homomorphism from~$G$ to the unitary group of a Hilbert space~$\Hilbert$, then 
	$$ \text{$\lim_{\|g\|\to \infty} \langle \pi(g) \phi \mid \psi \rangle = 0$, \ for all $\phi,\psi \in \Hilbert$} .$$
This yields the Moore Ergodicity Theorem \pref{MooreErgBasicThm} as an easy corollary, and the proof is based on the existence of $a \in G$ and (unipotent) subgroups $U^+$ and~$U^-$ of~$G$, such that $\langle U^+, U^-\rangle = G$ and $a^n u a^{-n} \to e$ as $n \to \infty$ (or $-\infty$), for all $u \in U^+$ (or $U^-$, respectively). 

\smallbreak

\pref{PeterWeyl}
Every unitary representation of any compact Lie group is a direct sum of finite-dimensional, irreducible unitary representations.

 \smallbreak

\pref{IntIrredAbel}
Every unitary representation of any abelian Lie group is a direct integral of one-dimensional unitary representations.

 \smallbreak

\pref{UnitaryGDirectInt}
Every unitary representation of~$G$ is a direct integral of irreducible unitary representations.

\subsection*{\cref{AmenableChap}. Amenable Groups}
Amenability is such a fundamental notion that it has very many quite different definitions, all of which determine exactly the same class of groups \csee{AmenDefn,AmenEquiv}. One useful choice is that a group~$\Lambda$ is \emph{amenable} if every continuous action of~$\Lambda$ on a compact, metric space has a finite, invariant measure. 

\smallbreak

\pref{FreeSubgrp->Nonamen} The fact that the lattice~$\Gamma$ contains a nonabelian free subgroup \csee{FreeInGamma} implies that it is not amenable. This is because subgroups of amenable groups are amenable \csee{SubgrpAmen}, and free groups do have actions (such as the actions described in the Ping-Pong Lemma \pref{PingPong}) that do not have a finite, invariant measure.

\smallbreak

Even so, amenability plays an important role in the study of~$\Gamma$, through the following observation:

\smallskip

\pref{G/amen->Meas(X)} (\emph{Furstenberg Lemma})
If $P$ is an amenable subgroup of~$G$, and we have a continuous action of~$\Gamma$ on some compact, metric space~$X$, then there exists a measurable, $\Gamma$-equivariant map from $G/P$ to the space $\Prob(X)$ of measures~$\mu$ on~$X$, such that $\mu(X) = 1$. To prove this, let $\mathcal{F}$ be the set of measurable, $\Gamma$-equivariant maps from $G$ to $\Prob(X)$. With an appropriate weak topology, this is a compact, metrizable space, and $P$~acts on it by translation on the right. Since $P$ is amenable, there is a $P$-invariant, finite measure~$\mu$ on~$\mathcal{F}$. The barycenter of this measure is a fixed point of~$P$ in~$\mathcal{F}$, and this fixed point is a function on~$G$ that factors through to a well-defined $\Gamma$-equivariant map from $G/P$ to $\Prob(X)$.

\smallbreak

\subsection*{\cref{KazhdanTChap}. Kazhdan's Property $(T)$}
To say $\Gamma$ has Kazhdan's property~$(T)$ means that if a unitary representation of~$\Gamma$ does not have any (nonzero) vectors that are fixed by~$\Gamma$, then it does not \emph{almost-invariant vectors}, that is, vectors that are moved only a small distance by the elements of any given finite subset of~$\Gamma$ \csee{KazhdanTDefn}.

\smallbreak

\pref{T+amen->finite} Kazhdan's property~$(T)$ is, in a certain sense, the antithesis of amenability: a discrete group cannot have both properties unless it is finite.
This is because the regular representation of any amenable group has almost-invariant vectors.

\smallbreak

\pref{KazhdanEasy} Every discrete group with Kazhdan's property~$(T)$ is finitely generated. To see this, let $\Hilbert = \bigoplus_F \LL2(\Lambda/F)$, where $F$ ranges over all the finitely generated subgroups of~$\Lambda$. Then, by construction, every finite subset of~$\Lambda$ fixes some nonzero vector in~$\Hilbert$.

\smallbreak

\pref{WhichGKazhdan} $G$ has Kazhdan's property~$(T)$, unless $G$ is either $\SO(1,n)$ or $\SU(1,n)$.
To prove this for $G = \SL(3,\real)$, first note that the semidirect product $\SL(2,\real) \ltimes \real^2$ can be embedded in~$G$. Also note that there are elements $a$ and~$b$ of $\SL(2,\real)$, such that, if $Q$ is any of the $4$~quadrants of~$\real^2$, then either $aQ$ or~$bQ$ is disjoint from~$Q$ (except for the $0$ vector). Applying this to the Pontryagin dual of~$\real^2$ implies that if a representation of the semidirect product $\SL(2,\real) \ltimes \real^2$ has almost-invariant vectors, then it must have a nonzero vector that is invariant under~$\real^2$. This vector must be invariant under all of $\SL(3,\real)$, by a generalization of the Moore Ergodicity Theorem that is called the \emph{Mautner phenomenon} \pref{MautnerPhenom}.

\smallbreak

\pref{Kazhdan:G->Gamma} $\Gamma$ has Kazhdan's property~$(T)$, unless $G$ is either $\SO(1,n)$ or $\SU(1,n)$. Any unitary representation~$\pi$ of~$\Gamma$ can be ``induced'' to a representation~$\pi_\Gamma^G$ of~$G$. If $\pi$ has almost-invariant vectors, then the induced representation has almost-invariant vectors, and, since $G$ has Kazhdan's property~$(T)$, this implies that $\pi_\Gamma^G$ has $G$-invariant vectors. Any such vector must come from a $\Gamma$-invariant vector in~$\pi$.

\smallbreak

\pref{T<>FH} A group has Kazhdan's property~$(T)$ if and only if every action of the group by (affine) isometries on any Hilbert space has a fixed point. 
This is not at all obvious, but here is the proof of one direction. 

Suppose $\Gamma$ does not have Kazhdan's property~$(T)$, so there exists a unitary representation of~$\Gamma$ on some Hilbert space~$\Hilbert$ that has almost-invariant vectors, but does not have invariant vectors. 
Choose an increasing chain $F_1 \subseteq F_2 \subseteq \cdots$ of finite subsets whose union is all of~$\Gamma$. Since $\Hilbert$ has almost-invariant vectors, there exists a unit vector $v_n \in \Hilbert$, such that $\| f v_n - v_n\| < 1/2^n$ for all $f \in F_n$.
Now, define $\alpha \colon \Gamma \to \Hilbert^\infty$ by 
	$$ \alpha(g)_n = n\bigl( gv_n - v_n \bigr) .$$
Then $\alpha$ is a $1$-cocycle, so defining
	$ g * v = gv + \alpha(g)$
yields an action of~$\Gamma$ on the Hilbert space~$\Hilbert^\infty$. Since $\Hilbert$ has no nonzero invariant vectors, it is not difficult to see that $\alpha$ is an unbounded function on~$\Gamma$, so $\alpha$ is not a coboundary. This implies that the corresponding action on~$\Hilbert^\infty$ has no fixed points.

\smallbreak

\subsection*{\cref{ErgodicChap}. Ergodic Theory}
\emph{Ergodic Theory} can be defined as the measure-theoretic study of group actions. In this category, the analogue of the transitive actions are the so-called \emph{ergodic} actions, for which every measurable, invariant function is constant (a.e.) \csee{ErgodicDefn}. 

%\smallbreak
%
%\pref{MooreErgodicity} (\emph{Moore Ergodicity Theorem}) 
%If $H$ is any noncompact, closed subgroup of~$G$, then the action of~$H$ on  $G/\Gamma$ is ergodic. 
%This was already stated in \cref{MooreErgBasicThm}, but this section provides the proof \csee{MooreErgPfSect}.

\smallbreak

\pref{PointwiseErgThm} (\emph{Pointwise Ergodic Theorem})
If $\integer$ acts ergodically on~$X$, with finite invariant measure, and $f$~is any $L^1$-function on~$X$, then the average of~$f$ on almost every $\integer$-orbit is equal to the average of~$f$ on the entire space~$X$.

\smallbreak

\pref{ErgodicDecomp}
Every measure-preserving action of~$G$ can be measurably decomposed into a union of ergodic actions.
%If $G$ acts continuously on a complete metric space~$X$, with a $\sigma$-finite, invariant measure~$\mu$ measure-preserving action of~$G$ can be decomposed (measurably) as a union of ergodic actions.

\smallbreak

\pref{GMixing}
If the action of~$G$ on a space~$X$ is ergodic, with a finite, invariant measure, then the action of~$G$ on $X \times X$ is also ergodic.

\notocsection*{\cref{MajorPart}. Major Results}

Here are some of the major theorems in the theory of arithmetic groups.

\subsection*{\cref{MostowChap}. Mostow Rigidity Theorem} \ 
% check to make sure there's not a page break here !!!

\smallskip

\pref{MostowRigidity} (\emph{Mostow Rigidity Theorem}) 
Suppose $\Gamma_i$ is a lattice in~$G_i$, for $i = 1,2$, and $\varphi \colon \Gamma_1 \to \Gamma_2$. If $G_i$ has trivial center and no compact factors, and is not $\PSL(2,\real)$, then $\varphi$ extends to an isomorphism $\overline\varphi \colon G_1 \to G_2$. 

In most cases, the desired conclusion is a consequence of the Margulis Superrigidity Theorem, which will be discussed in \cref{MargulisSuperChap}. 
However, a different proof is needed when $G_1 = G_2 = \SO(1,n)$ (and some other cases). Assuming that the lattices are cocompact, the proof uses the fact (mentioned in \cref{CocpctActIsQI}) that $\Gamma_1$ and~$\Gamma_2$ are quasi-isometric to~$\hyperbolic^n$. Comparing the two embeddings yields a quasi-isometry $\varphi$ from $\hyperbolic^n$ to itself. By proving that this quasi-isometry induces a map on the boundary that is conformal (i.e., preserves angles), it is shown that the two embeddings are conjugate by an isometry of~$\hyperbolic^n$.

\smallbreak

\pref{LotsOfLattsSL2} Mostow's theorem does not apply to $\PSL(2,\real)$: in this group, there are uncountably many lattices that are isomorphic to each other, but are not conjugate. This follows from the fact that there are uncountably many different right-angled hexagons in the hyperbolic plane~$\hyperbolic^2$. A compact surface of genus~$g$ can be constructed by gluing $4g-4$ of these hexagons together, in such a way that the fundamental group is a cocompact lattice in $\PSL(2,\real)$. The uncountably many different hexagons yield uncountably many non-conjugate lattices.

\smallbreak

\pref{QI->GIso} From Mostow's Theorem, we know that lattices in two different groups $G_1$ and~$G_2$ cannot be isomorphic. In fact, the lattices cannot even be quasi-isometric.
Some ideas in the proof of this fact are similar to the argument of Mostow's theorem, but we omit the details.

\subsection*{\cref{MargulisSuperChap}. Margulis Superrigidity Theorem} \ 

\smallskip

\pref{MargSuperStatementSect} (\emph{Margulis Superrigidity Theorem}) Suppose $\rho \colon \Gamma \to \GL(n,\real)$ is a homomorphism. If $G$ is neither $\SO(1,n)$ nor $\SU(1,n)$, and  mild hypotheses are satisfied, then $\rho$ extends to a homomorphism $\overline\rho \colon G \to \GL(n,\real)$. 

Assuming $\Rrank G \ge 2$, a proof is presented in \cref{SuperPfSect}. Start by letting $H$ be the Zariski closure of $\rho(\Gamma)$, and let $Q$ be a parabolic subgroup of~$H$. Furstenberg's Lemma \pref{G/amen->Meas(X)} provides a $\Gamma$-equivariant map $\psi \colon G/P \to \Prob(\real \projective^n)$. By using ``proximality\zz,'' $\psi$~can be promoted to a map $\widehat\psi \colon G/A \to \real^n$ (where $A$ is a maximal $\real$-split torus of~$G$). Thus, we have an $A$-invariant (measurable) section of the flat vector bundle over $G/\Gamma$ that is associated to~$\varphi$. Since $G$ is generated by the centralizers of nontrivial, connected subgroups of~$A$, this implies there is a finite-dimensional, $G$-invariant space of sections of the bundle, from which it follows that $\varphi$ has the desired extension to a homomorphism defined on all of~$G$.

\smallbreak

\pref{MostowRigidityIrred} This theorem of Margulis is a strengthening of the Mostow Rigidity Theorem \pref{MostowRigidity}, because the homomorphism $\rho$ is not required to be an isomorphism.  (On the other hand, Mostow's theorem applies to the groups $\SO(1,n)$ and $\SU(1,n)$, which are not allowed in the superrigidity theorem.)

\smallbreak

\pref{Super->VecBdlTrivial} In geometric terms, the superrigidity theorem implies (under mild hypotheses) that flat vector bundles over $G/\Gamma$ become trivial on a finite cover.

\smallbreak

\pref{MargArithPf} (\emph{Margulis Arithmeticity Theorem}) If $G$ is neither $\SO(1,n)$ nor $\SU(1,n)$, then the superrigidity theorem implies that every lattice in~$G$ is an arithmetic subgroup (as was stated without proof in \cref{MargulisArith}). 

The basic idea of the proof is that if there is some $\rho(\gamma)$ with a matrix entry that is transcendental, then composing $\rho$ with arbitrary elements of the Galois group $\Gal(\complex/\rational)$ would result in uncountably many different $n$-dimensional representations of~$\Gamma$. Since $G$ has only finitely many representations of each dimension, this would contradict superrigidity. Thus, we conclude that $\rho(\Gamma) \subseteq \GL(n,\overline\rational)$. By using a $p$-adic version of the superrigidity theorem, $\overline\rational$ can be replaced with~$\integer$.

\smallbreak

\pref{SuperRank1Sect} For groups of real rank one, the proof of superrigidity described in \cref{SuperPfSect} does not apply, because $A$ does not have any nontrivial, proper subgroups. Instead, a more geometric approach is used (but only a brief sketch will be provided). Let $X$ and~$Y$ be the symmetric spaces associated to~$G$ and~$H$, respectively, where $H$ is the Zariski closure of $\rho(\Gamma)$. By minimizing a certain energy functional, one can show there is a harmonic $\Gamma$-equivariant map $\psi \colon X \to Y$. Then, by using the geometry of $X$ and~$Y$, it can be shown that this harmonic map must be a totally geodesic embedding. This provides an embedding of the isometry group of~$X$ in the isometry group of~$Y$. In other words, an embedding of~$G$ in~$H$.

\subsection*{\cref{NormalSubgroupChap}. Normal Subgroups of~$\Gamma$} \ 

\smallskip

\pref{MargNormalSubgrpsThm} If $\Rrank G \ge 2$, then $\Gamma$ is almost simple. More precisely, every normal subgroup of~$\Gamma$ either is finite, or has finite index.
This is proved by showing that if $N$ is any infinite, normal subgroup of~$\Gamma$, then the quotient $\Gamma/N$ is amenable. Since $\Gamma/N$ has Kazhdan's property~$(T)$ (because we saw in \cref{Kazhdan:G->Gamma} that $\Gamma$ has this property), this implies $\Gamma / N$ is finite.

\smallbreak

\pref{Rrank1->GammaNotAlmSimple} On the other hand, if $\Rrank G = 1$, then $\Gamma$ is very far from being simple --- there are many, many infinite normal subgroups of~$\Gamma$. In fact, $\Gamma$~is ``SQ-universal\zz,'' which means that if $\Lambda$ is any finitely generated group, then there is a normal subgroup~$N$ of~$\Gamma$, such that $\Lambda$ is isomorphic to a subgroup of~$\Gamma/N$ \csee{Rank1SQUniv}.

\goodbreak

\subsection*{\cref{ArithClassicalChap}. Arithmetic Subgroups of Classical Groups} \ 
The main result of this chapter is the table on \cpageref{IrredInG} that provides a list of all of the arithmetic subgroups of~$G$ (unless $G$ is either an exceptional group or a group whose complexification~$G_\complex$ is isogenous to $\SO(8,\complex)$). Inspection of the list establishes several results that were stated without proof in previous chapters.

\smallskip

\pref{QFormsOfSLnSect} It was stated without proof in \cref{LattSlnRSect} that every arithmetic subgroup of $\SL(n,\real)$ is either a special linear group or a unitary group (if we allow division algebras in the construction). The proof of this fact is based on a calculation of the group cohomology of Galois groups (or \emph{Galois cohomology}, for short). To introduce this method in a simpler setting, it is first proved that the only $\real$-forms of the complex Lie group $\SL(n,\complex)$ are $\SL(n,\real)$, $\SL(n/2, \quaternion)$, and $\SU(k,\ell)$ \csee{GaloisCohoRealFormsSect}.

\smallbreak

\pref{QFormClassicalSect} The same methods show that all the $\rational$-forms of any classical group~$G$ are classical groups (except that there is a problem when $G$ is a real form of $\SO(8,\complex)$ \csee{D4weird}). However, we do not provide the calculations.

\smallbreak

\pref{Isotypic->irred} We say that a semisimple group $H = G_1 \times \cdots \times G_r$ is \emph{isotypic} if all the simple factors of~$H_\complex$ are isogenous to each other.
A theorem of Borel and Harder \pref{BorelHarderLocGlob} on Galois cohomology implies that if $H$ is isotypic, then it has an arithmetic subgroup that is \emph{irreducible}: it is not commensurable to a nontrivial direct product $\Gamma_1 \times \Gamma_2$. (The converse follows from the Margulis Arithmeticity Theorem unless $H$ is either $\SO(1,n) \times K$ or $\SU(1,n) \times K$.)

\subsection*{\cref{ReductionChap}. Construction of a Coarse Fundamental Domain}
This chapter presents some of the main ideas involved in the construction of a nice subset of~$G$ that approximates a fundamental domain for $G/\Gamma$ (when $\Gamma$ is an arithmetic subgroup). This generalizes the construction for $\Gamma = \SL(n,\integer)$ that was explained in \cref{SLnZLattChap}.

As in \cref{SLnZLattChap}, the key notion is that of a \emph{Siegel set}. 
The main difference is that, instead of the maximal $\real$-split torus~$A$, we must work with a subtorus~$T$ of~$A$ that is $\rational$-split, not merely $\real$-split:
	\begin{itemize}
	\item $K$ is a maximal compact subgroup of~$G$ (same as before),
	\item $S$ is a maximal $\rational$-split torus in~$G$,
	\item $S_t$ is a sector in~$S$,
	\item $P$ is a minimal parabolic $\rational$-subgroup of~$G$ that contains~$S$,
	and
	\item $P_0$ is a compact subset of~$P$.
	\end{itemize}
Then $K S_t P_0$ is a \emph{Siegel set} for~$\Gamma$ in~$G$ \csee{ArithSiegelDefn}.

It may not be possible to find a Siegel set~$\Siegel$, such that $\Siegel \cdot \Gamma = G$ \csee{SiegelNotFundEg}. (When $\dim S = 1$, this is because each Siegel set can only cover one cusp, and $G/\Gamma$ may have several cusps.) However, there is always a finite union of (translates of) Siegel sets that will suffice \csee{ReductThyArithGrps}.

\smallbreak

The existence of a nice set~$\fund$, such that $\fund \cdot \Gamma = G$, has  important consequences, such as the fact that $\Gamma$ is finitely presented \csee{FinPresSect}. This fact was stated in \cref{GammaFinPres}, but could only be proved for the cocompact case there.

\subsection*{\cref{RatnerChap}. Ratner's Theorems on Unipotent Flows}
If $\{a^t\}$ is any $1$-parameter group of diagonal matrices in~$G$, then there are $\{a^t\}$-orbits in $G/\Gamma$ that have bad closures: the closure is a fractal. M.\,Ratner proved that if a subgroup~$V$ is generated by $1$-parameter unipotent subgroups, then it is much better behaved: the closure of every $V$-orbit is a $C^\infty$~submanifold of $G/\Gamma$ \csee{Ratner-OrbitClosure}. 

\smallbreak

This theorem has important consequences in geometry and number theory. As a sample application in the theory of arithmetic groups, we mention that if $\Gamma_1$ and~$\Gamma_2$ are any two lattices in~$G$, then the subset $\Gamma_1 \, \Gamma_2$ of~$G$ is either discrete or dense \csee{ProdLattsDense}. This is proved by letting $\Gamma = \Gamma_1 \times \Gamma_2$ in $G \times G$, and letting $V$ be the diagonal embedding of~$G$ in the same group.

\smallbreak

Ratner proved that the actions of $1$-parameter unipotent subgroups on $G/\Gamma$ also have nice measurable properties: every finite, invariant probability measure is the Haar measure on a closed orbit of some subgroup of~$G\mk$ \csee{Ratner-MeasClass}, and every dense orbit is uniformly distributed \csee{Ratner-Equidistribution}.

\smallbreak

We will not prove Ratner's theorems, but some of the ideas in the proof will be described. One of the main ingredients is called ``shearing'' \csee{RatnerShearingSect}. For example, suppose $G = \SL(2,\real)$ and $V = \{u^t\}$ is a $1$-parameter unipotent subgroup. Then the key point is that if $x$ and~$y$ are two nearby points in $G/\Gamma$ (and are not on the same $\{u^t\}$-orbit), then the fastest relative motion between the two points is along the $V$-orbits. More precisely, there is some $t$, such that $u^t x$ is close to either $u^{t+1} y$ or $u^{t-1} y$.

\notocsection*{Appendices}
The main text is followed by three appendices. The first two (\cref{SSGrpsChap,BackChap}) recall some facts that are used in the main text. The third (\cref{SarithChap}) defines the notion of \emph{$S$-arithmetic group}, and quickly summarizes how the results on arithmetic groups extend to this more general setting. 

\makemarkstrue % resume updates of the running heads

%!TEX root = IntroArithGrps.tex

\part{Fundamentals} \label{FundamentalsPart}
  \addtocontents{toc}{\protect\toceject} % @@@
 %!TEX root = IntroArithGrps.tex

\mychapter{Basic Properties of Lattices}
\label{BasicLatticesChap}

\prereqs{none.}

This book is about lattices in semisimple Lie groups (with emphasis
on the ``arithmetic'' ones).

\bigskip

\setcounter{equation}{-1} % standing assumps will be numbered 3.0

%% emphasize this! (box with gray border is used here @@@)
\definecolor{grayborder}{gray}{0.6}
\definecolor{whitebkgnd}{gray}{1}
\newlength{\grayborder} \grayborder=10pt \relax 
\newlength{\graygap} \graygap=9pt \relax 
\newlength{\assumpw} \assumpw=\textwidth
	\advance\assumpw -2\grayborder\relax 
	\advance\assumpw -2\graygap\relax
\newbox\assumpbox
\setbox\assumpbox\hbox{\begin{minipage}{\assumpw}%
\begin{standassump} \label{standassump}
Throughout this book:
 \noprelistbreak \parindent=2.2em
 \begin{enumerate}
 \item \label{standassump-G}
 \textbf{\mathversion{bold}$G$ is a linear, semisimple Lie group}
 (see \cref{SSLieGrpDefnSect} for an explanation of these terms),
 with only finitely many connected components,
 and
 \item \label{standassump-Gamma}
 \textbf{\mathversion{bold}$\Gamma$ is a lattice in~$G$} (see \cref{LatticeDefn}).
 \end{enumerate}
 Similar restrictions apply to the symbols $G_1$, $G_2$,
$G'$, $\Gamma_1$, $\Gamma_2$, $\Gamma'$, etc.
\end{standassump}
\end{minipage}%
}
\newlength{\grayh} \grayh=\ht\assumpbox \advance\grayh\dp\assumpbox
	\advance\grayh2\grayborder\relax \advance\grayh2\graygap\relax
\vbox{
\vbox to 0pt{\hbox to \textwidth{\hss\color{grayborder}\vrule width \textwidth height \grayh \hss}\vss}%
\vskip-\baselineskip
\newlength{\whitew} \whitew=\assumpw \advance\whitew2\graygap\relax
\newlength{\whiteh} \whiteh=\grayh \advance\whiteh-2\grayborder\relax
\vbox to 0pt{\vskip\grayborder\hbox to \textwidth{\hss\color{whitebkgnd}\vrule width \whitew height \whiteh \hss}\vss}
\vskip \grayborder \vskip \graygap \vskip-1pt % @@@
\centerline{\copy\assumpbox}
\vskip \grayborder \vskip \graygap
}

\begin{rem}
Without losing any of the main ideas, it may be assumed,
throughout, that $G$ is either $\SL(n,\real)$ or
$\SO(m,n)$ (or a product of these groups), 
%. Much of the material is of interest even in the special case $G = \SO(1,n)$, 
but it is best if the reader is also acquainted with the other
``classical groups\zz,'' such as unitary groups and symplectic groups \csee{ClassicalDefn}. 
\end{rem}

Three definitions in this \lcnamecref{BasicLatticesChap} are very important: lattice subgroups
\pref{LatticeDefn}, commensurable subgroups
\pref{CommensDefn}, and irreducible lattices
(\ref{irreducibleLattice} and~\ref{prodirredlatt}). The
rest of the material in this chapter may not be essential
for a first reading, and can be referred back to when
necessary. However, if the reader has no prior experience
with lattices, then the basic properties discussed
in~\cref{LattDefnSect} will probably be helpful.

\section{Definition} \label{LattDefnSect}

\begin{lem} \label{ExistsFundDom}
  If $\Lambda$ is a discrete subgroup of~$G$, then 
there is a {\normalfont\defit[fundamental!domain!strict]{strict fundamental domain}} for $G/\Lambda$ in~$G$.
That is, there is a Borel subset~$\fund$ of~$G$, such that
the natural map $\fund \to G/\Lambda$, defined
by $g \mapsto g \Lambda$, is bijective.
 \end{lem}

\begin{proof}
 Since $\Lambda$ is discrete, there is a nonempty, open
subset~$U$ of~$G$, such that $(U^{-1} \, U) \cap \Lambda =
\{e\}$. Since $G$ is second countable (or, if you prefer,
since $G$ is $\sigma$-compact), there is a sequence
$\{g_n\}$ of elements of~$G$, such that $\bigcup_{n=1}^\infty g_n
U = G$. Let
 $$ \fund = \bigcup_{n = 1}^\infty \left( 
 \raise1pt\hbox{$\displaystyle % eliminate some blank space % @@@
 g_n U 
\smallsetminus \bigcup_{i< n} g_{i} U \Lambda
$}
  \right) .$$
 Then $\fund$ is obviously Borel, and it is a strict fundamental domain for $G/\Lambda$ 
 \csee{FisFundDom}.
 \end{proof}

\begin{rem} \ 
\noprelistbreak
\begin{enumerate}
\item The above \lcnamecref{ExistsFundDom} is stated for the space $G/\Lambda$ of left cosets of~$\Lambda$, but, in some situations, it is more natural to work with the space $\Lambda \backslash G$ of right cosets. In this book, we will feel free to use whichever is most convenient at a particular time, and leave it to the reader to translate between the two, by using the fact that the function $g \Lambda \mapsto \Lambda g^{-1}$ is a homeomorphism from $G/\Lambda$ to $\Lambda \backslash G$ \csee{LeftCosetsToRightCosets,FundDomRightCosets}. Our choice will usually be determined by the preference for most mathematicians to write their actions on the left. (Therefore, if $G$ is acting, then we will tend to use $G/\Lambda$, but if we are thinking of~$\Lambda$ as acting on~$G$, then we usually consider the quotient $\Lambda \backslash G$.)

\item Definitions in the literature vary somewhat, but saying that a subset~$\fund$ of~$G$ is a \defit[fundamental!domain]{fundamental domain} for $G/\Lambda$ typically means:
	\begin{enumerate}
	\item $\fund \Lambda = G$,
	\item $\fund$ is a closed set that is nice: its interior $\interior\fund$ is dense in~$\fund$, and its boundary $\fund \smallsetminus \interior\fund$ has measure~$0$,
	and
	\item $\fund \lambda \cap \interior\fund = \emptyset$, for all nonidentity $\lambda \in \Lambda$.
	\end{enumerate}
It is not difficult to see that if $\fund$ is a fundamental domain, then it has a Borel subset $\fund'$, such that $\fund'$ is a strict fundamental domain and $\fund \smallsetminus \fund'$ has measure~$0$. This means that, for many purposes (such as calculating integrals), it suffices to have a fundamental domain, rather than finding a set that is precisely a strict fundamental domain.
\end{enumerate}
\end{rem}

\begin{prop} \label{HaarOnHomog}
 Let $\Lambda$ be a discrete subgroup of~$G$, and let
$\mu$ be Haar measure on~$G$. There is a unique
\textup(up to a scalar multiple\textup)
$\sigma$-finite, $G$-invariant Borel measure~$\nu$ on
$G/\Lambda$.
More precisely:
 \begin{enumerate}
 \item \label{HaarOnHomog-nu}
 For any strict fundamental domain~$\fund$, the
measure~$\nu$ can be defined by 
  \begin{equation} \label{HomogFromHaar}
 \nu(A/\Lambda) = \mu( A \cap \fund) ,
 \end{equation}
 for every Borel set~$A$ in~$G$, such that $A \Lambda = A$.
 \item \label{HaarOnHomog-mu}
 Conversely, for $A \subseteq G$, we have
 \begin{equation} \label{HaarFromHomog}
 \mu(A) = \int_{G/\Lambda} \#(A \cap x \Lambda)
\, d \nu(x \Lambda) .
 \end{equation}
 \end{enumerate}
 \end{prop}

\begin{proof} 
 See \cref{HaarOnHomogIsGinv,HaarOnHomog->Haar} for \pref{HaarOnHomog-nu}
and~\pref{HaarOnHomog-mu}. The uniqueness of~$\nu$ follows
from~\pref{HaarOnHomog-mu} and the uniqueness of the Haar
measure~$\mu$.
 \end{proof}

\begin{rem}
 We always assume that the $G$-invariant measure~$\nu$ on
$G/\Lambda$ is normalized so that
\pref{HomogFromHaar} and~\pref{HaarFromHomog} hold.
 \end{rem}

\begin{cor} \label{inject->samemeas}
 Let $\Lambda$ be a discrete subgroup of~$G$, and let
$\phi \colon G \to G/\Lambda$ be the natural
quotient map $\phi(g) = g \Lambda$. If $A$ is a Borel
subset of~$G$, such that the restriction $\phi|_A$ is
injective, then $\nu \bigl( \phi(A) \bigr) = \mu(A)$.
 \end{cor}

\begin{rems} \ \label{VolOnG/Gamma}
 \noprelistbreak
 \begin{enumerate}
 \item \label{VolOnG/Gamma-smooth}
 The Haar measure~$\mu$ on~$G$ is given by a smooth
volume form, so the associated measure~$\nu$ on $G/\Lambda$ is also given by a volume form. Therefore, we say
that $G/\Lambda$ has \defit[finite volume
(homogeneous space)]{finite volume} if $\nu(G/\Lambda) < \infty$.
 \item The assumption that $\Lambda$ is discrete cannot be
eliminated from \cref{HaarOnHomog}. However, a $G$-invariant measure on $G/\Lambda$ can be constructed under the weaker assumption that $\Lambda$ is closed and
unimodular \csee{HaarOnHomog<>unimod}. 
%However, \pref{HomogFromHaar} and \pref{HaarFromHomog} are not valid in this generality.
 \end{enumerate}
 \end{rems}

\begin{defn} \label{LatticeDefn}
 A subgroup~$\Gamma$ of~$G$ is a
\defit[lattice!subgroup]{lattice} in~$G$ if
 \begin{itemize}
 \item $\Gamma$ is a discrete subgroup of~$G$, and
 \item $G/\Gamma$ has finite volume.
 \end{itemize}
 \end{defn}

\begin{rem}
The definition is not vacuous: we will explain in \cref{GHasLatt} that $G$~does have at least one lattice (in fact, infinitely many), although part of the proof will be postponed to \cref{SLnZLattChap}.
\end{rem}

\begin{prop} \label{Latt<>fund}
 Let $\Lambda$ be a discrete subgroup of~$G$, and let
$\mu$ be Haar measure on~$G$. The following are equivalent:
 \begin{enumerate}
 \item \label{Latt<>fund-latt}
 $\Lambda$ is a lattice in~$G$.
 \item \label{Latt<>fund-fund}
 There is a strict fundamental domain~$\fund$ for $G/\Lambda$, with $\mu(\fund) < \infty$.
 \item \label{Latt<>fund-otherside}
 There is a strict fundamental domain~$\fund'$ for $\Lambda \backslash G$, with $\mu(\fund') < \infty$.
 \item \label{Latt<>fund-coarse}
 There is a Borel subset $C$ of~$G$, such that $C \Lambda = G$ and $\mu(C) < \infty$.
 \end{enumerate}
 \end{prop}

\begin{proof}
 ($\ref{Latt<>fund-latt} \Leftrightarrow
\ref{Latt<>fund-fund}$) From \cref{HomogFromHaar}, we
have $\nu(G/\Lambda) = \mu(\fund)$. Therefore,
$G/\Lambda$ has finite volume if and only if
$\mu(\fund) < \infty$.

($\ref{Latt<>fund-fund} \Leftrightarrow \ref{Latt<>fund-otherside}$)
If $\fund$ is any strict fundamental domain for $G/\Lambda$, then $\fund^{-1}$ is a strict fundamental domain for $\Lambda \backslash G$ \csee{FundDomRightCosets}. 
Since $G$ is unimodular, we have $\mu(\fund^{-1}) = \mu(\fund)$ \csee{mu(inverse)}. 

 ($\ref{Latt<>fund-fund} \Rightarrow
\ref{Latt<>fund-coarse}$) Obvious.

 ($\ref{Latt<>fund-coarse} \Rightarrow
\ref{Latt<>fund-latt}$) We have $C \cap x \Lambda \neq
\emptyset$, for every $x \in G$, so, from
\pref{HaarFromHomog}, we see that
 \begin{align*} \nu(G/\Lambda) 
 = \int_{G/\Lambda} 1 \, d\nu(x \Lambda)
 \le \int_{G/\Lambda} \#(C \cap x \Lambda) \,
d\nu(\Lambda x)
 = \mu(C)
 < \infty 
 . & \qedhere \end{align*}
 \end{proof}

\begin{eg}
 As mentioned in \cref{SL2Zlatt}, $\SL(2,\integer)$ is
a lattice in $\SL(2,\real)$.
 \end{eg}

\begin{defn}
 A closed subgroup~$\Lambda$ of~$G$ is
\defit[cocompact subgroup]{cocompact} (or
\defit[uniform subgroup]{uniform}) if $G/\Lambda$ is compact.
 \end{defn}

\begin{cor} \label{Latt<>cpct/finind} \ 
 \noprelistbreak
 \begin{enumerate}
 \item \label{Latt<>cpct/finind-cpct}
 Every cocompact, discrete subgroup of~$G$ is a lattice.
 \item \label{Latt<>cpct/finind-finind}
 Every finite-index subgroup of a lattice is a lattice.
 \end{enumerate}
 \end{cor}

\begin{proof}
 \Cref{cocpct->latt,finindLatt->latt}.
 \end{proof}

\begin{rem}
 Lattices in~$G$ are our main interest, but we will
occasionally encounter lattices in Lie groups~$H$ that
are not semisimple. If $H$ is unimodular, then all of the
above results remain valid with $H$ in the place of~$G$.
In contrast, if $H$ is not unimodular, then
\cref{HaarOnHomog} may fail: there may exist a discrete subgroup~$\Lambda$, such that there is no
$H$-invariant Borel measure on $H/\Lambda$.
Instead, there is sometimes only a semi-invariant
measure~$\nu$:
 $$ \nu(hA) = \Delta(h) \, \nu(A) ,$$
 where $\Delta$~is the modular function of~$H$
\csee{SemiinvariantMeasure}. 
%This is sufficient to
%determine whether $H/\Lambda$ has finite volume
%or not, so \cref{LatticeDefn} applies.
 \end{rem}

For completeness, let us specifically state the following
concrete generalization of \cref{LatticeDefn}
\cf{Latt<>fund}.

\begin{defn} \label{LatticeGeneralDefn}
 A subgroup~$\Lambda$ of a Lie group~$H$ is a
\defit[lattice!subgroup]{lattice} in~$H$ if
 \begin{itemize}
 \item $\Lambda$ is a discrete subgroup of~$H$, and
 \item there is an $H$-invariant measure~$\nu$ on $H/\Lambda$, such that $\nu(H/\Lambda) < \infty$.
% a Borel subset $C$ of~$H$, such that
%$\Lambda C = H$ and $\mu(C) < \infty$, where $\mu$~is the
%left Haar measure on~$H$.
 \end{itemize}
 \end{defn}

\begin{eg}
 $\integer^n$ is a cocompact lattice in~$\real^n$. 
 \end{eg}

\begin{prop}
 If a Lie group~$H$ has a lattice, then $H$~is unimodular.
 \end{prop}

\begin{proof}
 Let $\fund$ be a strict fundamental domain for $H/\Lambda$.
 The proof of \cref{HaarOnHomog} shows $\nu(A/\Lambda) = \nu(A \cap \fund)$, for every Borel set~$A$ in~$G$, such that $A \Lambda = A$. Then \cref{SemiinvariantMeasure} implies $\nu(hA/\Lambda) = \Delta(h) \, \nu(A/\Lambda)$. In particular, we see that $\nu(H/\Lambda) = \Delta(h) \, \nu(H/\Lambda)$, by letting $A = H$  (and noting that $hH = H$).
Since $\nu(H/\Lambda) < \infty$, this implies $\Delta(h) = 1$, as desired.
 \end{proof}

\begin{exercises}[Recall that (in accordance with the Standing Assumptions \pref{standassump}), $\Gamma$~is a lattice in~$G$, and $G$~is a semisimple Lie group.]

\item Show that $\Gamma$ is finite if and only if $G$~is
compact.

\item \label{FisFundDom}
 Complete the proof of \cref{ExistsFundDom}; that
is, show that $\fund$ is a strict fundamental domain.

\item \label{LeftCosetsToRightCosets}
Define $f \colon G/\Lambda \to \Lambda \backslash G$ by $f(g \Lambda) = \Lambda g^{-1}$. Show that $f$ is a homeomorphism.

\item \label{FundDomRightCosets}
Show \cref{ExistsFundDom} easily implies an analogous statement that applies to right cosets. More precisely, show that if 
	\begin{itemize}
	\item $\Lambda$ is a discrete subgroup of~$G$, 
	\item $\fund$ is a strict fundamental domain for $G/\Lambda$, 
	and
	\item $\fund^{-1} = \{\, x^{-1} \mid x \in \fund \,\}$,
	\end{itemize}
then the natural map $\fund^{-1} \to \Lambda \backslash G$, defined
by $g \mapsto \Lambda g$, is bijective.

\item \label{mu(inverse)}
Show that $\mu(A^{-1}) = \mu(A)$ for every Borel subset~$A$ of~$G$.
\hint{Defining $\mu'(A) = \mu(A^{-1})$ yields a $G$-invariant measure on~$G$. The uniqueness of Haar measure implies $\mu' = \mu$. Where did you use the fact that $G$ is unimodular?}

\item \label{mu(fund)}
Let 
\noprelistbreak
	\begin{itemize}
	\item $\Lambda$ be a discrete subgroup of~$G$,
	\item $\fund$ and~$\fund'$ be strict fundamental domains for $G/\Lambda$,
	\item $\mu$ be Haar measure on~$G$,
	and
	\item $A$ be a Borel subset of~$G$.
	\end{itemize}
Show:
	 \begin{enumerate}
	 \item For each $g \in G$, there is a unique $\lambda \in \Lambda$, such that $g \lambda \in \fund$.
	 \item For each $\lambda \in \Lambda$, if we let
	 	 $A_\lambda = \{\, a \in A \mid a \lambda \in \fund \,\}$, then $A_\lambda$ is Borel, and $A$~is the disjoint union of the sets $\{\, A_\lambda \mid
	\lambda \in \Lambda \,\}$.
	 \item $\mu(\fund) = \mu(\fund')$.
	 \item \label{mu(fund)-F=F'}
	 If $A \Lambda = A$, then $\mu(A \cap \fund) = \mu(A \cap
	\fund')$.
	 \end{enumerate}

\item \label{HaarOnHomogIsGinv}
 Show, for every Haar measure~$\mu$ on~$G$, that the Borel
measure~$\nu$ defined in \fullcref{HaarOnHomog}{nu} is
$G$-invariant.
 \hint{For any $g \in G$, the set $g \fund$ is a strict
fundamental domain. From \fullcref{mu(fund)}{F=F'}, we
know that $\nu$~is independent of the choice of the strict
fundamental domain~$\fund$.}

\item \label{HaarOnHomog->Haar}
 If $\Lambda$ is a discrete subgroup
of~$G$, and $\nu$ is a $\sigma$-finite, $G$-invariant
Borel measure on $G/\Lambda$, show that the Borel
measure~$\mu$ defined in \fullcref{HaarOnHomog}{mu} is
$G$-invariant.

\item \label{HaarOnHomog<>unimod}
 Let $H$ be a closed subgroup of~$G$. Show that there is a
$\sigma$-finite, $G$-invariant Borel measure~$\nu$ on
$G/H$ if and only if $H$~is unimodular.
 \hint{($\Rightarrow$) For a left Haar measure~$\rho$
on~$H$, define a left Haar measure~$\mu$ on~$G$ by
 $$ \mu(A) = \int_{G/H} \rho(x^{-1} A \cap H) \,
d \nu(xH) .$$
 Then $\mu(A) = \Delta_H(h) \, \mu(Ah)$ for $h \in H$,
where $\Delta_H$ is the \term{modular function} of~$H$. Since
$G$ is unimodular, we must have $\Delta_H \equiv 1$.}

\item \label{finext->latt}
 Show that if $\Lambda$ is a discrete subgroup of~$G$
that contains~$\Gamma$, then $\Lambda$ is a lattice
in~$G$, and $\Gamma$ has finite index in~$\Lambda$.
\hint{Let $\fund$ be a strict fundamental domain for $G/\Lambda$, and let $F$~be a set of coset representatives for $\Gamma$ in~$\Lambda$. Then $\fund \cdot F$ is a strict fundamental domain for $G/\Gamma$, and therefore has finite measure.}

\item \label{CpctInHomog}
 Let $\Lambda$ be a discrete subgroup of~$G$. Show
that a subset~$A$ of $G/\Lambda$ is precompact
if and only if there is a compact subset~$C$ of~$G$, such
that $A \subseteq C\Lambda / \Lambda$.
 \hint{($\Leftarrow$) The continuous image of a compact
set is compact. ($\Rightarrow$) Let $\mathcal{U}$ be a
cover of~$G$ by precompact, open sets.} 

\item \label{cocpct->latt}
 Prove \fullcref{Latt<>cpct/finind}{cpct}.
 \hint{\cref{CpctInHomog} and
\fullcref{Latt<>fund}{coarse}.}

\item \label{finindLatt->latt}
 Prove \fullcref{Latt<>cpct/finind}{finind}.
 \hint{\cref{Latt<>fund}. A finite union of sets
of finite measure has finite measure.}

\item \label{SemiinvariantMeasure}
 Let 
 \begin{itemize}
 \item $H$ be a Lie group,
 \item $\Lambda$ be a discrete subgroup of~$H$,
 \item $\mu$ be the right Haar measure on~$H$,
 and
 \item $\fund$ be a strict fundamental domain for $H/\Lambda$. %such that $\mu(\fund) < \infty$.
 \end{itemize}
 Define a $\sigma$ Borel measure~$\nu$ on $H/\Lambda$ by 
 $ \nu(A/\Lambda) = \mu( A \cap \fund) $,
 for every Borel set~$A$ in~$H$, such that $A \Lambda =
A$. Show $\nu(h A / \Lambda) =
\Delta(h) \, \nu(A/\Lambda)$, where $\Delta$~is
the modular function of~$H$.
 \hint{Cf.\ \cref{HaarOnHomogIsGinv}.}

%\item Suppose $\Gamma_1$ and~$\Gamma_2$ are lattices
%in~$G$, such that $\Gamma_1 \subseteq \Gamma_2$. Show that
%$\Gamma_1$ has finite index in~$\Gamma_2$.

\item Show that every discrete, cocompact subgroup of every Lie group is a lattice.
\hint{Define $\nu$ as in \cref{SemiinvariantMeasure}. Since %$H/\Lambda$ is compact, we know 
$\nu(H/\Lambda)< \infty$ (why?), we must have $\Delta(h) = 1$.}

\end{exercises}

\section{Commensurability and isogeny} \label{CommSect}

We usually wish
to ignore the minor differences that come from passing to a finite-index
subgroup. The following
definition describes the resulting equivalence relation.

\begin{defn} \label{CommensDefn}
 We say that two subgroups~$\Lambda_1$ and~$\Lambda_2$ of a group~$H$ are
\defit{commensurable} if $\Lambda_1 \cap \Lambda_2$ is a
finite-index subgroup of both~$\Lambda_1$ and~$\Lambda_2$.
This is an equivalence relation on the collection of all
subgroups of~$H$ \csee{CommEquiv}.
 \end{defn}

\begin{egs} \ \label{CommensEg}
\noprelistbreak
 \begin{enumerate}
 \item Two cyclic subgroups $a \integer$
and~$b \integer$ of~$\real$ are commensurable if and
only if $a$ is a nonzero rational multiple of~$b$; therefore,
commensurability of subgroups generalizes the classical
notion of commensurability of real numbers.

 \item \label{CommensEg-latt}
 It is easy to show that every subgroup commensurable
to a lattice is itself a lattice. (For example, this follows from
\fullcref{Latt<>cpct/finind}{finind} and
\cref{finext->latt}.)
 \end{enumerate}
 \end{egs}

The analogous notion for Lie groups (with finite center and finitely many connected components) is called ``isogeny:''

\begin{defns} \  \label{IsogenyDefn}
\noprelistbreak
	\begin{enumerate}
	\item $G_1$ is \defit{isogenous} to~$G_2$ if some finite cover of~$(G_1)^\circ$ is isomorphic to some finite cover of~$(G_2)^\circ$. This is an equivalence relation.
	\item A (continuous) homomorphism $\varphi \colon G_1 \to G_2$ is an \defit{isogeny} if it is an isomorphism modulo finite groups. More precisely:
	\begin{itemize}
	\item the kernel of~$\varphi$ is finite,
	and
	\item the image of~$\varphi$ has finite index in~$G_2$.
	\end{itemize}
	\end{enumerate}
\end{defns}

\begin{rem}
The following are equivalent:
	\begin{enumerate}
	\item $G_1$ is isogenous to~$G_2$.
	\item $\Ad (G_1)^\circ \iso \Ad (G_2)^\circ$.
	\item $G_1$ and~$G_2$ are \defit[locally!isomorphic]{locally isomorphic}, that is, the Lie
algebras $\Lie G_1$ and~$\Lie G_2$ are isomorphic.
	\item There is an isogeny from some finite cover of~$(G_1)^\circ$ to~$G_2$.
	\end{enumerate}
\end{rem}

The normalizer of a subgroup is very important in group
theory. Because we are ignoring finite groups, the
following definition is natural in our context.

\begin{defn}
 An element~$g$ of~$G$ \emph{commensurates}~$\Gamma$ if
$g \Gamma g^{-1}$ is commensurable to~$\Gamma$. Let
 $$\Comm_G(\Gamma) = \{\, g \in G \mid \text{$g$
commensurates~$\Gamma$} \,\} .$$
 This is called the \defit{commensurator} of~$\Gamma$.
 %(or the \defit{commensurability subgroup} of~$\Gamma$).
 \end{defn}

\begin{rem}
 The commensurator of~$\Gamma$ is sometimes much larger than the
normalizer of~$\Gamma$. For example, let $G =
\SL(n,\real)$ and $\Gamma = \SL(n,\integer)$. Then
$\nzer_G(\Gamma)$ is commensurable to~$\Gamma$
\csee{latticenormalizer}, but $\Comm_G(\Gamma)$ contains $\SL(n,\rational)\mk$\csee{SLQCommSLZ}, so $\Comm_G(\Gamma)$
is dense in~$G$, even though $\nzer_G(\Gamma)$ is discrete.
Therefore, in this example (and, more generally, whenever
$\Gamma$ is ``arithmetic''), $\nzer_G(\Gamma)$ has infinite
index in $\Comm_G(\Gamma)$.

On the other hand, if $G = \SO(1,n)$, then it is known
that there are examples in which $\Gamma$, $\nzer_G(\Gamma)$,
and $\Comm_G(\Gamma)$ are commensurable to each other
(see \cref{Comm(Gamma)discrete}
and \cref{NonarithInSO1n}).
 \end{rem}

\begin{defn}
 We say that two groups $\Lambda_1$ and~$\Lambda_2$ are
\defit[commensurable!abstractly]{abstractly commensurable} if some finite-index subgroup of~$\Lambda_1$ is isomorphic to some finite-index subgroup
of~$\Lambda_2$.
 \end{defn}

Note that if $\Lambda_1$ and~$\Lambda_2$ are
commensurable, then they are abstractly commensurable, but not
conversely.

\begin{exercises}

\item \label{CommEquiv}
 Verify that commensurability is an equivalence relation.

\item If $\Gamma_1$ is commensurable
to~$\Gamma_2$, show $\Comm_G(\Gamma_1) =
\Comm_G(\Gamma_2)$.

\end{exercises}

\section{Irreducible lattices} \label{IrredLattSect}

Note that $\Gamma_1 \times \Gamma_2$ is a lattice in $G_1
\times G_2$. A lattice that can be decomposed as a product
of this type is said to be \defit[reducible!lattice]{reducible}.

\begin{defn} \label{irreducibleLattice}
 $\Gamma$ is \defit[irreducible!lattice]{irreducible} if $\Gamma N$ is dense
in~$G$, for every noncompact, closed, normal subgroup~$N$
of~$G^\circ$ (and $\Gamma$ is infinite, or, equivalently, $G$~is not compact).
 \end{defn}

\begin{eg}
If $G$ is simple (and not compact), then
every lattice in~$G$ is irreducible.  Conversely, if $G$
is not simple, then not every lattice
in~$G$ is irreducible.  To see this, assume, for
simplicity, that $G$ is connected and has trivial center (and is not compact). Then we may write
$G$ as a nontrivial direct product $G = G_1 \times G_2$,
where each of $G_1$ and~$G_2$ is semisimple. If we let
$\Gamma_i$ be any lattice in~$G_i$, for $i = 1,2$, then
$\Gamma_1 \times \Gamma_2$ is a reducible lattice in~$G$.
\end{eg}

The following \lcnamecref{prodirredlatt} shows (under mild assumptions)
that every lattice is commensurable to a product of irreducible
lattices. Therefore, the preceding example provides essentially
the only way to construct reducible lattices, so most
questions about lattices can be reduced to the irreducible
case. We postpone the proof, because
it relies on some results from later in this \lcnamecref{BasicLatticesChap}.

\begin{prop}[(see proof on \cpageref{prodirredlattPf})] \label{prodirredlatt}
 Assume 
 \noprelistbreak
 \begin{itemize}
 \item $G$ has trivial center,
  and 
 \item $\Gamma$ projects densely into the maximal compact
factor of~$G$. 
 \end{itemize}
 Then there is a direct-product decomposition $G = G_1 \times
\cdots \times G_r$, such that $\Gamma$ is commensurable
to\/ $\Gamma_1 \times \cdots \times \Gamma_r$, where\/
$\Gamma_i = \Gamma \cap G_i$, and\/ $\Gamma_i$ is an
irreducible lattice in~$G_i$, for each~$i$.
 \end{prop}

For readers familiar with locally symmetric spaces, these results can be restated in the following geometric
terms. 

\begin{defn}
 Recall that a locally symmetric space $\Gamma \backslash
X$ is \defit[irreducible!locally symmetric space]{irreducible} if there do not exist (nontrivial)
locally symmetric spaces $\Gamma_1 \backslash X_1$ and
$\Gamma_2 \backslash X_2$, such that the product
 $ (\Gamma_1 \backslash X_1) \times (\Gamma_2 \backslash
X_2)$
 finitely covers $\Gamma \backslash X$.
 \end{defn}

The following is obvious by induction on $\dim X$.

\begin{prop}
 There exist locally symmetric spaces\/
 $\Gamma_1 \backslash X_1, \ldots, \Gamma_r \backslash X_r$
 that are irreducible, such
that the product\/
 $ (\Gamma_1 \backslash X_1) \times \cdots \times (\Gamma_r
\backslash X_r)$
 finitely covers $\Gamma \backslash X$.
 \end{prop}

The following is a restatement of
\cref{prodirredlatt} (in the special case where $G$
has no compact factors).

\begin{prop} 
 Let $M$ be an irreducible locally symmetric space, such
that the universal cover~$X$ of~$M$ has no compact
factors, and no flat factors. For any nontrivial cartesian
product decomposition $X = X_1 \times X_2$ of~$X$, the
image of~$X_1$ is dense in~$M$.
 \end{prop}

We will see in \cref{SL(2Z[sqrt2])} that
$\SL(2,\real) \times \SL(2,\real)$ has an irreducible
lattice (for example, a lattice isomorphic to
$\SL\bigl( 2,\integer[\sqrt{2}] \bigr)$). 
More generally, \cref{Isotypic->irred} shows
that $G$ has an irreducible lattice if all the simple factors of the ``\term{complexification}'' of~$G$ are isogenous to each other. The converse is proved in \cref{irred->isotypic}, under the additional assumption that $G$ has no compact factors.

\begin{exercises}

\item \label{irred->CpctFactor}
 Show that if $\Gamma$ is irreducible, %and $G$~is not compact, 
then $\Gamma$ projects densely into the maximal
compact factor of~$G$.

\end{exercises}

\section{\texorpdfstring{Unbounded subsets of $\Gamma \backslash G$}{Unbounded subsets of G/Γ}} % should be Γ\G, but that gives a TeX error @@@
\label{DivergeSect}

 Geometrically, looking at the fundamental
domain described in \cref{SL2Zlatt} makes it clear that the sequence
$\{ni\}$ tends to~$\infty$ in $\SL(2,\integer) \backslash
\hyperbolic^2$. In this section, we give an algebraic
criterion that determines whether or not a sequence tends
to~$\infty$ in $G/\Gamma$, without any need for
a fundamental domain.

Recall that the \defit{injectivity radius} of a Riemannian
manifold~$X$ is the maximal $r \ge 0$, such that, for
every $x \in X$, the exponential map is a diffeomorphism
on the open ball of radius~$r$ around~$x$. If $X$~is
compact, then the injectivity radius is nonzero. The
following proposition shows that the converse holds in the
special case where $X = \Gamma \backslash G/K$ is locally
symmetric of finite volume.

\begin{prop} \label{Cpct<>InjRad}
 For $g \in G$, define $\phi_g \colon G \to G/\Gamma$ by 
 $\phi_g(x) = xg\Gamma$. The homogeneous
space $G/\Gamma$ is compact if and only if
there is a nonempty, open subset~$U$ of~$G$, such that,
for every $g \in G$, the restriction $\phi_g|_U$
of~$\phi_g$ to~$U$ is injective.
 \end{prop}

\begin{proof}
 ($\Rightarrow$) Define $\phi \colon G \to G/\Gamma$ 
 by $\phi(x) = x \Gamma$. Then $\phi$~is a
covering map, so, for each $p \in G/\Gamma$,
there is a connected neighborhood~$V_p$ of~$p$, such that
the restriction of~$\phi$ to each component of
$\phi^{-1}(V_p)$ is a diffeomorphism onto~$V_p$. Since
$\{\, V_p \mid p \in G/\Gamma\}$ is an open
cover of $G/\Gamma$, and $G/\Gamma$
is compact, there is a connected neighborhood~$U$ of~$e$
in~$G$, such that, for each $p \in G/\Gamma$,
there is some $p' \in G/\Gamma$, with $Up
\subseteq V_{p'}$ \csee{LebesgueNumberG}. Then $\phi_g|_U$
is injective, for each $g \in G$.

($\Leftarrow$) We prove the contrapositive. Let $U$ be any
nonempty, precompact, open subset of~$G$. (We wish to show,
for some $g \in G$, that $\phi_g|_U$ is not injective.)
 If $C$ is any compact subset of $G/\Gamma$,
then, because $G/\Gamma$ is not compact, we
have 
 $$ (G/\Gamma) \smallsetminus (U^{-1} C) \neq
\emptyset .$$
 Hence, we may inductively construct a sequence
$\{g_n\}$ of elements of~$G$, such that the open sets
$\phi_{g_1}(U), \phi_{g_2}(U), \ldots$ are pairwise
disjoint. Since $G/\Gamma$ has finite volume,
these sets cannot all have the same volume, so, for
some~$n$, the restriction $\phi_{g_n}|_U$ is not injective
\csee{inject->samemeas}. 
 \end{proof}

Let us restate this geometric result in algebraic terms.

\begin{notation}
 For elements $a$ and~$b$ of a group~$H$, and subsets~$A$
and~$B$ of~$H$, let%
	\nindex{${}^b \! a = b a b^{-1}$ for $a,b \in G$}%
	\nindex{${}^G\Gamma$ = $\{\, {}^g \gamma \mid \gamma \in \Gamma, \ g \in G \,\}$}
\begin{align*}
	{}^b \! a &= b a b^{-1} , 
	&
	{}^B \! a &= \{\, {}^b \!a \mid b \in B\,\}, \\
 	{}^b \! A &= \{\, {}^b \!a \mid a \in A \,\}, 
 	&
	{}^B \! A &= \{\, {}^b \!a \mid a\in A, b \in B \,\}
 	. \end{align*}
  \end{notation}

\begin{cor} \label{G/GammaCpct<>NoAccPt}
 $G/\Gamma$ is compact if and only if the
identity element~$e$ is \textbf{not} an accumulation point
of\/ ${}^G\Gamma$.
 \end{cor}

\begin{proof}
 We have
\begin{align*} 
\text{$\phi_g|_U$ is injective}
 &\quad \Leftrightarrow \quad 
 \nexists u_1,u_2 \in U, \ 
 \text{$u_1 g \Gamma = u_2 g \Gamma$ \  and \ $u_1 \neq u_2$}
 \\&\quad \Leftrightarrow \quad 
 {}^g \Gamma \cap (U^{-1} U) = \{e\} 
 . \qedhere \end{align*}
 \end{proof}

This has the following interesting consequence.

\begin{cor} \label{GammaUnip->notcpct}
 If\/ $\Gamma$ has a nontrivial, \index{unipotent!element}{unipotent element}, then
$G/\Gamma$ is \textbf{not} compact.
 \end{cor}

\begin{proof}
 If $u$~is a nontrivial, unipotent element
of~$\Gamma$, then, from the \thmindex{Jacobson-Morosov}{Jacobson-Morosov Lemma}
\pref{JacobsonMorosov}, we know there is a continuous
homomorphism $\phi \colon \SL(2,\real) \to G$, with 
 $\phi 
 \begin{Smallbmatrix}
 1 & 1 \\
 0 & 1
 \end{Smallbmatrix}
 = u$.  
 Let 
 $ a = \phi 
 \begin{Smallbmatrix}
  1/2 & 0 \\
 0 & 2
 \end{Smallbmatrix}
 \in G $.
 Then
 \goodbreak % @@@
 \begin{align*}
  a^{n} \, u \, a^{-n} 
 &= \phi 
 \left( 
 \begin{bmatrix}
 2^{-n} & 0 \\
 0 & 2^n
 \end{bmatrix}
 \begin{bmatrix}
 1 & 1 \\
 0 & 1
 \end{bmatrix}
 \begin{bmatrix}
 2^n & 0 \\
 0 & 2^{-n}
 \end{bmatrix}
 \right)
 \\&= 
 \phi
 \left( \begin{bmatrix}
 1 & 2^{-2n} \\
 0 & 1
 \end{bmatrix} \right)
 \to
 \phi
 \left( \begin{bmatrix}
 1 & 0 \\
 0 & 1
 \end{bmatrix} \right)
 = e .\end{align*}
 Therefore, $e$~is an accumulation point of ${}^G u$, so
\cref{G/GammaCpct<>NoAccPt} implies that $G/\Gamma$ is not compact.
 \end{proof}

\begin{rems} \  \label{GodementConverse}
\noprelistbreak
	\begin{enumerate}
	\item If $G$ has no compact factors, then the converse of
\cref{GammaUnip->notcpct} is true.
However, we will prove this only in the special case where $\Gamma$ is
``arithmetic'' \csee{GodementSect}.
	\item In general (without any assumption on compact factors), it can be shown that $G/\Gamma$ is compact if and only if every element of~$\Gamma$ is semisimple \csee{CocpctAllSS}.
	\end{enumerate}
\end{rems}

The proofs of \cref{Cpct<>InjRad} and
\cref{G/GammaCpct<>NoAccPt} establish the following
more general version of those results.

\begin{prop}[\csee{diverge<>contractEx}] \label{diverge<>contract}
 Let $\Lambda$ be a lattice in a Lie group~$H$, and let
$C$ be a subset of~$H$. The image of~$C$ in $H/\Lambda$ is 
precompact if and only if the identity
element~$e$ is \textbf{not} an accumulation point of\/
${}^C \! \Lambda$.
 \end{prop}

The following is a similar elementary result that applies
to the important special case where $G = \SL(\ell,\real)$
and $\Gamma = \SL(\ell,\integer)$, without relying on the
fact that $\SL(\ell,\integer)$ is a lattice. 
%For
%convenience, the result discusses $G/\Gamma$, rather than
%$\Gamma \backslash G$ (because we write $gv$, not~$vg$,
%for $g \in G$ and $v \in \real^\ell$). Since $\Gamma
%\backslash G$ is diffeomorphic to $G/\Gamma$ (via the map
%$\Gamma g \mapsto g^{-1} \Gamma$), this change of notation
%is of no real significance.

\begin{prop}[(\thmindex{Mahler Compactness Criterion}Mahler Compactness Criterion)]
\label{MahlerCpct}
 Let $C \subseteq \SL(\ell,\real)$. The image of~$C$ in\/
$\SL(\ell,\real)/\SL(\ell,\integer)$ is precompact if and
only if\/ $0$~is \textbf{not} an accumulation point of 
 $$ C \integer^\ell = \{\, cv \mid c \in C, v \in
\integer^\ell  \,\} .$$
 \end{prop}

\begin{proof}
 ($\Rightarrow$)  Since the image of~$C$ in
$\SL(\ell,\real)/\SL(\ell,\integer)$ is precompact, there
is a compact subset~$C_0$ of~$G$, such that $C \subseteq C_0
\SL(\ell,\integer)$ \csee{CpctInHomog}. There is no harm in
assuming that $C = C_0 \SL(\ell,\integer)$ (by enlarging~$C$). Then
 $ C \bigl( \integer^\ell \smallsetminus \{0\} \bigr) 
 = C_0 \bigl( \integer^\ell \smallsetminus \{0\} \bigr) 
 $
 is closed (since $\integer^\ell \smallsetminus \{0\}$,
being discrete, is closed and $C_0$~is compact), so $C
\bigl( \integer^\ell \smallsetminus \{0\} \bigr)$ contains
all of its accumulation points. In addition, since $0$~is
fixed by every element of~$C$, we know that $0 \notin C
\bigl( \integer^\ell \smallsetminus \{0\} \bigr)$.
Therefore, $0$~is not an accumulation point of $C \bigl(
\integer^\ell \smallsetminus \{0\} \bigr)$.

($\Leftarrow$) To simplify the notation (while retaining
the main ideas), let us assume $\ell = 2$
\csee{MahlerCpctPfExer}. Suppose $\{g_n\}$ is a sequence
of elements of~$\SL(2,\real)$, such that $0$~is
\textbf{not} an accumulation point of $\bigcup_{n=1}^\infty
g_n \integer^2$. We wish to show there is a
sequence $\{\gamma_n\}$ of elements of $\SL(2,\integer)$,
such that $\{g_n \gamma_n\}$ has a convergent subsequence.

For each~$n$, let 
\noprelistbreak
 \begin{itemize}
 \item $v_n \in \integer^2 \smallsetminus \{0\}$, such
that $\|g_n v_n\|$ is minimal,
 \item $\pi_n \colon \real^2 \to \real g_n v_n$ and
 $\pi_n^\perp \colon \real^2 \to (\real g_n v_n)^\perp$
 be the orthogonal projections, and
 \item $w_n \in \integer^2 \smallsetminus \real v_n$, such
that $\| \pi_n^\perp(g_n w_n)\|$ is minimal. 
 \end{itemize}
 By replacing $w_n$ with $w_n + k v_n$, for an appropriately chosen $k \in
\integer$, we may assume $\|\pi_n(g_n w_n)\| \le \|g_n
v_n\|/2$. Then, since the minimality of $\|g_n v_n\|$ implies
$\| g_n v_n \| \le \| g_n w_n \|$, we have 
 $$ \| g_n v_n \|
% \le \| g_n w_n \|
 \le \| \pi_n^\perp(g_n w_n)\| + \|\pi_n(g_n w_n)\|
 \le \| \pi_n^\perp(g_n w_n)\| + \frac{\| g_n v_n \|}{2}
,$$
 so
 \begin{equation} \label{LowerBdOnWn}
 \|\pi_n^\perp(g_n w_n)\| \ge \frac{\| g_n v_n \|}{2} .
 \end{equation}

Let $C$ be the convex hull of $\{0,v_n,w_n\}$ and
(thinking of $v_n$ and~$w_n$ as column vectors) let
$\gamma_n = \bigl[ v_n \  w_n \bigr] \in \Mat_{2\times 2}(\integer)$.
From the minimality of $\|g_n v_n\|$ and $\|
\pi_n^\perp(g_n w_n)\|$, we see that $C \cap \integer^2 =
\{0,v_n,w_n\}$ \csee{min->ConvHull}, so $\det \gamma_n =
\pm 1$ \csee{DetOne}.
Therefore, perhaps after replacing $w_n$ with~$-w_n$, we
have $\gamma_n \in \SL(2,\integer)$. Since 
	$\gamma_n \! \left[ \begin{smallmatrix} 1 \\ 0 \end{smallmatrix} \right] = v_n$ 
and $\gamma_n \! \left[ \begin{smallmatrix} 0 \\ 1 \end{smallmatrix} \right] = w_n$, 
we may assume, by
replacing $g_n$ with~$g_n \gamma_n$, that 
 $$ \mbox{$v_n = \left[ \begin{smallmatrix} 1 \\ 0 \end{smallmatrix} \right]$
 \qquad and\qquad 
 $w_n = \left[ \begin{smallmatrix} 0 \\ 1 \end{smallmatrix} \right]$.} $$

Note that 
 \begin{equation} \label{det=gvxperp}
 \|\pi_n^\perp (g_n w_n)\| \cdot \|g_n v_n \|
 = \det g_n
 = 1 .
 \end{equation}
 By combining this with \pref{LowerBdOnWn}, we see that
$\{ g_n v_n\}$ is a bounded sequence, so, by passing to a
subsequence, we may assume $g_n v_n$ converges to some
vector~$v$. By assumption, we have $v \neq 0$.

Now, from \pref{det=gvxperp}, and the fact that $\|g_n v_n
\| \to \|v\|$ is bounded away from~$0$, we see that
$\|\pi_n^\perp (g_n w_n)\|$ is bounded. Because
$\|\pi_n(g_n w_n)\|$ is also bounded, we conclude that
$\|g_n w_n\|$ is bounded. Hence, by passing to a
subsequence, we may assume $g_n w_n$ converges to some
vector~$w$. From \pref{LowerBdOnWn}, we know that
$\|\pi_n^\perp(g_n w_n)\| \not\to 0$, so $w \notin \real
v$.

Since $v \neq 0$ and $w \notin \real v$, there is some $g
\in \GL(\ell,\real)$ with $g \! \left[ \begin{smallmatrix} 1 \\ 0 \end{smallmatrix} \right] = v$ and $g \! \left[ \begin{smallmatrix} 0 \\ 1 \end{smallmatrix} \right] = w$.
We have 
 $$ g_n \! \left[ \begin{smallmatrix} 1 \\ 0 \end{smallmatrix} \right] = g_n v_n \to v = g \! \left[ \begin{smallmatrix} 1 \\ 0 \end{smallmatrix} \right] $$
 and, similarly, $g_n \!  \left[ \begin{smallmatrix} 0 \\ 1 \end{smallmatrix} \right] \to g \! \left[ \begin{smallmatrix} 0 \\ 1 \end{smallmatrix} \right]$, so $g_n x \to gx$
for all $x \in \real^2$. Therefore, $g_n \to g$, as desired.
 \end{proof}
 
\begin{exercises}

\item \label{LebesgueNumberG}
 Suppose a Lie group~$H$ acts continuously on a compact
topological space~$M$, and $\mathcal{V}$~is an open cover
of~$M$. Show that there is a neighborhood~$U$ of~$e$
in~$H$, such that, for each $m \in M$, there is some $V
\in \mathcal{V}$ with $Um \subseteq V$.
 \hint{This is analogous to the fact that every open cover
of a compact metric space has a ``\term[Lebesgue!number]{Lebesgue number}\zz.''
Each $m \in M$ is contained in some $V_m \in \mathcal{V}$. Choose $V_m'$ containing~$m$, and a neighborhood $U'_m$ of~$e$, such that $U'_m V_m' \subseteq V_m$. Cover $M$ with finitely many $V_m'$.}

\item \label{diverge<>contractEx}
Prove \cref{diverge<>contract}.
\hint{See the proofs of \cref{Cpct<>InjRad} and
\cref{G/GammaCpct<>NoAccPt}.}

\item \label{LattInH->ProperMap}
 Use \cref{diverge<>contract} to show that if
$H$ is a closed subgroup of~$G$, such that $H/\Gamma$
is a lattice in~$H$, then the natural inclusion map 
 $H/(\Gamma \cap H) \hookrightarrow G/\Gamma$ is proper.
 \hint{It suffices to show that if $C$ is a subset of~$H$,
such that the image of~$C$ in $G/\Gamma$ is
precompact, then the image of~$C$ in $H/(\Gamma \cap H)$ is also precompact.}

\item Let $G = \SL(2,\real)$, $\Gamma = \SL(2,\integer)$,
and $A$~be the subgroup consisting of all the diagonal
matrices in~$G$. Show that the natural inclusion map 
 $A/(\Gamma \cap A) \hookrightarrow G/\Gamma$ is proper, but $\Gamma \cap A$ is \emph{not}
a lattice in~$A$. (Therefore, the converse of
\cref{LattInH->ProperMap} does not hold.)

\item Let $G = \SL(3,\real)$, $\Lambda = \SL(3,\integer)$, and $a = \diag(1/2, \, 1, \, 2)
% \begin{bmatrix}
% 1/2 & 0 & 0 \\
% 0 & 1 & 0 \\
% 0 & 0 & 2
% \end{bmatrix}
 \in H$.
 Show that $a^n \Lambda \to \infty$ in $G/\Lambda$ 
 as $n \to \infty$. That is, show, for each compact
subset~$C$ of $G/\Lambda$, that, for all
sufficiently large~$n$, we have $a^n \Lambda \notin C$.
(For the purposes of this exercise, do \emph{not} assume
that $\Lambda$ is a lattice in~$G$.)

%\item Show that if $\Gamma$ contains a nontrivial
%unipotent element, then $\Gamma \backslash G$ is not
%compact.

%\item Prove ($\Rightarrow$) of \cref{MahlerCpct}.

\item \label{MahlerCpctPfExer}
 Prove \cref{MahlerCpct}($\Leftarrow$)  % is it still leftarrow !!!
  without assuming $\ell = 2$.
 \hint{Extend the definition of $v_n$ and~$w_n$ to an
inductive construction of vectors $u_{1,n},\ldots,u_{\ell,n} \in
\integer^\ell$.}

\item \label{min->ConvHull}
 Suppose that $v$~and~$w$ are linearly independent vectors
in~$\real^\ell$, and $x = av + bw$, with $a,b \ge 0$
and $a+b \le 1$. Show that either
 \begin{itemize}
 \item $x \in \{v,w\}$, or
 \item $\|x\| < \|v\|$, or
 \item $x \notin \real v$ and $d(x,\real v) < d(w,\real v)$.
 \end{itemize}
 \hint{Either $b = 1$, or $b = 0$, or $0< b < 1$.}
 
 \item \label{DetOne}
 Let $C$ be the convex hull of $\{0,v,w\}$, where $v$ and~$w$ are linearly independent vectors in $\integer^2$. Show that if $C \cap \integer^2 = \{0,v,w\}$, then $\det \bigl[ v \  w \bigr] = \pm 1$.
 \hint{Let $P$ be the the convex hull of $\{0,v,w,v+w\}$, so $\bigl| \det \bigl[ v \  w \bigr] \bigr|$ is the area of~$P$. If this area is $> 1$, then the translates of~$P$ by elements of~$\integer^2$ cannot be a tiling of~$\real^2$.}
 
\item \label{NGammaclosed}
 Let $H$ be a closed subgroup of~$G$. 
 	\begin{enumerate}
	\item \label{NGammaclosed-closed}
	Show that if $\Gamma \cap H$ is a lattice in~$H$, then $H \Gamma$ is
closed in~$G$.
	\item \label{NGammaclosed-converse}
	Show that the converse holds if $H$~is normal in~$G$.
	\end{enumerate}
\hint{\pref{NGammaclosed-closed}~\Cref{LattInH->ProperMap}.
\pref{NGammaclosed-converse}~Since $G/\Gamma$ is a bundle over $G/(H\Gamma)$ with fiber $H\Gamma/\Gamma$, Fubini's Theorem implies that $H\Gamma/\Gamma$ has finite volume. So the $H$-equivariantly homeomorphic space $H/(\Gamma \cap H)$ also has finite volume.}

\item \label{CpctOpenSubgrp->LattsCocpct}
Suppose 
	\begin{itemize}
	\item $\Lambda$ is a non-cocompact lattice in a topological group~$H$,
	and
	\item $H$ has a compact, open subgroup~$K$.
	\end{itemize}
Show that $\Lambda$ has a nontrivial element of finite order.
\hint{Since $K$ is compact and $\Lambda$ is discrete, it suffices to show that some conjugate of~$\Lambda$ intersects~$K$ nontrivially.}

\end{exercises}

\section{Borel Density Theorem and some consequences}
\label{BDTSect}

The results in this \lcnamecref{BDTSect} require the minor assumption that $\Gamma$ projects densely into the maximal compact factor of~$G$. This hypothesis is automatically satisfied (vacuously) if $G$ has no compact factors.
Recall that $G^\circ$ denotes the identity component of~$G$ \csee{GoNotation}.

\begin{thm}[(Borel)]
\label{BDT}
 Assume
 \begin{itemize}
 \item $\Gamma$ projects densely into the maximal compact factor of~$G$,
 \item $V$ is a finite-dimensional vector space
over\/~$\real$ or\/~$\complex$, and
 \item $\rho \colon G \to \GL(V)$ is a continuous
homomorphism.
 \end{itemize}
 Then:
 \begin{enumerate}
 \item \label{BDT-vector}
 Every $\rho(\Gamma)$-invariant vector in~$V$ is $\rho(G^\circ)$-invariant.
 \item \label{BDT-subspace}
 Every $\rho(\Gamma)$-invariant subspace of~$V$ is $\rho(G^\circ)$-invariant.
 \end{enumerate}
 \end{thm}

 \begin{proof}
 The proof is not difficult, but, in order to get to the applications more quickly, we will postpone it until \cref{BDTPfSect}. For now, to illustrate the main idea, let us just prove~\pref{BDT-vector}, in the
special case where $G/\Gamma$ is compact (and $G$ is connected).
Assume also that
$G$ has no compact factors \csee{BDTCpctFactor}; then $G$
is generated by its unipotent elements \csee{<unip>=G}, so
it suffices to show that $v$ is invariant under $\rho(u)$,
for every nontrivial unipotent element~$u$ of~$G$. Because
$\rho(u)$ is unipotent \csee{rho(unip)},
we know that $\rho(u^n) v$ is a polynomial function of~$n$ \csee{un=poly}.
However, because $G/\Gamma$ is compact and $\rho(\Gamma) v
= v$, we also know that $\rho(G) v$ is compact, so
 $\{\, \rho(u^n) v \mid n \in \natural \,\}$ is bounded. Every
bounded polynomial is constant, so we conclude that
$\rho(u^n) v = v$ for all~$n$; in particular, $\rho(u)v =
\rho(u^1) v = v$, as desired.
 \end{proof}

\begin{cor} \label{BDT-normalize}
 Assume\/ $\Gamma$ projects densely into the maximal compact
factor of~$G$.
 If $H$ is a connected, closed subgroup of~$G$ that is
normalized by\/~$\Gamma$, then $H$ is normal in~$G^\circ$.
 \end{cor}

\begin{proof}
%Since every closed subgroup of a Lie group is a Lie subgroup \csee{ClosedSubgrpIsLie}, we know that 
The Lie algebra $\Lie H$ of~$H$ is a vector subspace
of the Lie algebra~$\Lie G$ of~$G$. Also, because $\Gamma$ normalizes~$H$, 
we know that $\Lie H$ is invariant
under $\Ad_G \Gamma$. From \fullcref{BDT}{subspace}, we
conclude that $\Lie H$ is invariant under $\Ad G^\circ$. Since $H$ is connected,
this implies that $H$ is a normal subgroup of~$G^\circ$.
 \end{proof}

\begin{cor} \label{C(Gamma)=Z(G)}
 If\/ $\Gamma$ projects densely into the maximal compact factor of~$G$ \textup(and $G$ is connected\/\textup), then $\czer_G(\Gamma) = Z(G)$.
 \end{cor}

\begin{proof}
Recall that $G \subseteq \SL(\ell,\real)$, for some~$\ell$ \csee{standassump}.
 Let $V = \Mat_{\ell \times \ell}(\real)$ be the vector
space of all real $\ell \times \ell$ matrices, so $G
\subseteq V$. For $g \in G$ and $v \in V$, define $\rho(g) v
= g v g^{-1}$, so $\rho \colon G \to \GL(V)$ is a
continuous homomorphism. If $c \in \czer_G(\Gamma)$, then
 $ \rho(\gamma) c = \gamma c \gamma^{-1} = c$
 for every $\gamma \in \Gamma$, so
\fullcref{BDT}{vector} implies that $\rho(G) c = c$.
Therefore $c \in Z(G)$.
 \end{proof}

\begin{cor} \label{NinZ(G)}
 Assume\/ $\Gamma$ projects densely into the maximal compact
factor of~$G$ \textup(and $G$ is connected\/\textup).
 If $N$ is a finite, normal subgroup of\/~$\Gamma$, then $N
\subseteq Z(G)$.
 \end{cor}

\begin{proof}
 The quotient $\Gamma/ \czer_\Gamma(N)$ is finite,
because it embeds in the finite group $\Aut(N)$,
 so $\czer_\Gamma(N)$ is a lattice in~$G$
\fullcsee{Latt<>cpct/finind}{finind}.
 Then, because $N \subseteq \czer_G \bigl( \czer_\Gamma(N) \bigr)$,
\cref{C(Gamma)=Z(G)} implies $N \subseteq Z(G)$.
 \end{proof}

\begin{cor} \label{latticenormalizer}
 If\/ $\Gamma$ projects densely into the maximal compact
factor of~$G$, then\/ $\Gamma$ has finite index in
its normalizer $\nzer_G(\Gamma)$.
 \end{cor}

\begin{proof}
By passing to a subgroup of finite index, we may assume $G$ is connected.
 Because $\Gamma$ is discrete, the identity component
$\nzer_G(\Gamma)^\circ$ of $\nzer_G(\Gamma)$ must
centralize~$\Gamma$. So $\nzer_G(\Gamma)^\circ \subseteq
\czer_G(\Gamma) = Z(G)$ is finite. On the other hand,
$\nzer_G(\Gamma)^\circ$ is connected. Therefore,
$\nzer_G(\Gamma)^\circ$ is trivial, so $\nzer_G(\Gamma)$ is
discrete. Hence $\Gamma$ has finite index in
$\nzer_G(\Gamma)$ \csee{finext->latt}.
 \end{proof}

\begin{cor}[(\thmindex{Borel Density}{Borel Density Theorem})]
\label{BDT(Zardense)}
If\/ $\Gamma$ projects
densely into the maximal compact factor of~$G$ \textup(and $G$ is connected\/\textup), then\/
$\Gamma$ is Zariski dense in~$G$. That is, if $Q \in
\real[x_{1,1},\ldots,x_{\ell,\ell}]$ is a polynomial
function on $\Mat_{\ell \times \ell}(\real)$, such that
$Q(\Gamma) = 0$, then $Q(G) = 0$.
 \end{cor}

\begin{proof}
 Let 
 $$ \mathcal{Q} = \{\, Q \in
\real[x_{1,1},\ldots,x_{\ell,\ell}] \mid Q(\Gamma) = 0
\,\} .$$
From the definition of~$\mathcal{Q}$, it is obvious that $\Gamma$ is contained in the corresponding variety $\Var(\mathcal{Q})$ \csee{AlgicGrpDefn}. Since
$\Var(\mathcal{Q})$ has only finitely many connected
components \csee{Zar->AlmConn}, this implies that
$\Var(\mathcal{Q})^\circ$ is a connected subgroup of~$G$
that contains a finite-index subgroup of~$\Gamma$. Hence
\cref{BDT-normalize} implies that
$\Var(\mathcal{Q})^\circ = G$ \csee{GammaNotinConn}, so
$G \subseteq \Var(\mathcal{Q})$, as desired.
 \end{proof}

With the above results, we can now provide the proof that was postponed from \cpageref{prodirredlatt}:
 
 \begin{proof}[Proof of \cref{prodirredlatt}] \label{prodirredlattPf}
  We may assume $\Gamma$ is reducible (otherwise, let $r =
1$). Hence, there is some noncompact, connected, closed,
normal subgroup~$N$ of~$G$, such that $N \Gamma$ is
\emph{not} dense in~$G$; let $H$ be the closure of
$N\Gamma$, and let $H_1 = H^\circ$. Because $\Gamma \subset
H$, we know that $\Gamma$ normalizes $H_1$, so
$H_1$ is a normal subgroup of~$G$ \csee{BDT-normalize,irred->CpctFactor}).

Let $\Lambda_1 = H_1 \cap \Gamma$. By
definition, $H_1$ is open in~$H$, so the subgroup $H_1 \Gamma$ is also
open in~$H$. It is therefore closed,
%On the other hand, we have $N \subseteq H_1$,
%and $N\Gamma$ is dense in~$H$ (by the definition of~$H$), so $H_1 \Gamma$ is dense in~$H$. Therefore $H_1 \Gamma = H$. Therefore $H_1 \Gamma$ is
%closed in~$G$, 
so $\Lambda_1$ is a lattice in~$H_1$
\csee{NGammaclosed}. 

Because $H_1$ is normal in~$G$ and $G$ is semisimple (with
trivial center), there is a normal subgroup~$H_2$ of~$G$,
such that $G = H_1 \times H_2$ \csee{G=NxH}. Let $\Lambda
= H_1 \cap (H_2 \Gamma)$ be the projection of~$\Gamma$
to~$H_1$. Now $\Gamma$ normalizes~$\Lambda_1$, and $H_2$
centralizes~$\Lambda_1$, so $\Lambda$ must normalize
$\Lambda_1$. Therefore \cref{latticenormalizer}
implies that $\Lambda$ is discrete (hence closed), so $H_2
\Gamma = \Lambda \times H_2$ is closed, so $\Lambda_2 = H_2
\cap \Gamma$ is a lattice in~$H_2$ \csee{NGammaclosed}.

Because $\Lambda_1$ is a lattice in~$H_1$ and $\Lambda_2$
is a lattice in~$H_2$, we know that $\Lambda_1 \times
\Lambda_2$ is a lattice in $H_1 \times H_2 = G$. Hence,
$\Lambda_1 \times \Lambda_2$ has finite index in~$\Gamma$
\csee{finext->latt}.

By induction on $\dim G$, we may write 
	$$ \text{$H_1 = G_1 \times \cdots \times G_s$
	 \ and \ 
	 $H_2 = G_{s+1} \times \cdots \times G_r$} 
	 , $$
so that $\Gamma \cap G_i$ is an irreducible lattice
in~$G_i$, for each~$i$.
 \end{proof}

\begin{rem} \label{BDTNotProj}
For simplicity, the statement of \cref{BDT} assumes that $\Gamma$ projects densely into the maximal compact factor of~$G$. Without this assumption, the proof of \fullcref{BDT}{vector} establishes the weaker conclusion that $v$ is $\rho(S)$-invariant, for every noncompact, simple factor~$S$ of~$G$. This leads to alternate versions of the corollaries that make no assumption about the compact factor of~$G$. For example, the analogue of \cref{BDT-normalize} states that 
if $H$ is a connected, closed subgroup of~$G$ that is
normalized by~$\Gamma$, then $H$ is normalized by every noncompact, simple factor of~$G$.
\end{rem}

\begin{exercises}

\item \label{BDTCpctFactor}
 Prove \fullref{BDT}{vector}, under the assumption that $G/\Gamma$ is compact (but allowing $G$ to have compact factors).
 \hint{The above proof shows that $v$ is invariant under the the product~$N$ of all the noncompact factors of~$G$. So it is invariant under the closure of $N\Gamma$, which is~$G$.}

\item \label{<unip>=G}
Show that if $G$ is connected, and has no compact factors, then it is generated by its unipotent elements.
\hint{Consider each simple factor of~$G$ individually. The conjugates of a unipotent element are also unipotent.}

 \item \label{rho(unip)}
Suppose $\rho \colon G \to \SL(\ell,\real)$ is a continuous homomorphism. Show that if $u$ is unipotent, then $\rho(u)$ is unipotent.
\hint{The Jacobson-Morosov Lemma \pref{JacobsonMorosov} allows you to assume $G = \SL(2,\real)$ and $u = 
 \begin{Smallbmatrix}
 1 & 1 \\
 0 & 1
 \end{Smallbmatrix}
 $. Then a sequence of conjugates of~$u$ converges to~$\Id$, so the
characteristic polynomial of $\rho(u)$ is the same as the
characteristic polynomial of $\Id$.}

\item \label{un=poly}
Show that if $u$ is a unipotent element of $\SL(\ell,\real)$ and $v \in \real^\ell$, then each coordinate of the vector $u^n \, v$ is a polynomial function of~$n$.
\hint{Let $u = \Id + T$, where $T^{\ell+1} = 0$. Then
 $u^n \, v = (\Id + T)^n v = \sum_{k=0}^\ell \binom{n}{k} T^k v$.}

\item \label{GammaNotAbelian}
 Show that if $G$ is not compact, then $\Gamma$ is not
abelian.

\item \label{GammaNotSolvable}
 Generalizing \cref{GammaNotAbelian}, show that if
$G$ is not compact, then $\Gamma$ is not solvable.

\item \label{CommGammaInfinite}
 Strengthening \cref{GammaNotAbelian}, show that if
$G$ is not compact, then the commutator subgroup
$[\Gamma,\Gamma]$ is infinite. 

\item Assume the hypotheses of \cref{BDT}, and that $G$ is connected. For
definiteness, assume that $V$ is a real vector space. For
any subgroup~$H$ of~$G$, let $V[H]$ be the $\real$-span of 
 $\{\, \rho(h) \mid h
\in H \,\}$
 in $\End(V)$.
 Show that $V[\Gamma] = V[G]$.
 
 \item \label{IdCompN(Gamma)inK}
 Show the identity component of $\nzer_G(\Gamma)$ is contained in the maximal compact factor of~$G$.
 \hint{Apply \cref{latticenormalizer} to $G/K$, where $K$ is the maximal compact factor.}

\item \label{GammaHasNonUnip}
 Show that if $G$ is not compact, then $\Gamma$ has an
element that is not unipotent.
 \hint{Any unipotent element~$\gamma$ of $\SL(\ell,\real)$ satisfies the polynomial $(x - 1)^\ell = 0$.}

\item \label{GammaNotinConn}
 Assume $G$ has no compact factors. Show that if $H$ is a
connected, closed subgroup of~$G$ that contains a finite-index
subgroup of~$\Gamma$, then $H = G^\circ$.
\hint{$H$ is normalized by $\Gamma \cap H$, so $H \normal G^\circ$.}

\item Assume $G$ has trivial center and no compact factors. Show
that $\Gamma$ is reducible if and only if there is a
finite-index subgroup~$\Gamma'$ of~$\Gamma$ such that
$\Gamma'$ is isomorphic to $A \times B$, for some infinite
groups $A$~and~$B$.\par
\noindent {\smaller(Actually, although you do not need to prove it, there is no need to assume the center of~$G$ is trivial. This is because $\Gamma$~has a subgroup of finite index that is torsion free \csee{torsionfree}, and therefore does not intersect the center of~$G$.)\par}

\item \label{NcapGammaFinite}
 Show that if $\Gamma$ is irreducible, then $N \cap
\Gamma$ is finite, for every connected, closed,
normal subgroup~$N$ of~$G$, such that $G/N$ is not
compact. 
 \hint{See the proof of \cref{prodirredlatt}.}
 
\item Let $\rho_1$ and~$\rho_2$ be finite-dimensional,
real representations of~$G$. Assume $G$ is connected, and has no compact
factors. Show that if the restrictions $\rho_1|_\Gamma$
and $\rho_2|_\Gamma$ are isomorphic, then $\rho_1$
and~$\rho_2$ are isomorphic.
\hint{We are assuming $\rho_i \colon G \to \GL(n,\real)$, for some~$n$, and that there is some $A \in \GL(n,\real)$, such that $\rho_1(g) = A \,\rho_2(g) A^{-1}$, for all $\gamma \in \Gamma$. You wish to show there is some $A' \in \GL(n,\real)$, such that the same condition holds for all $g \in G$, with $A'$ in the place of~$A$. The Borel Density Theorem implies that you may take $A' = A$.}

\end{exercises}

\section{Proof of the Borel Density Theorem}
\label{BDTPfSect}

The proof of the Borel Density Theorem \pref{BDT} is based on the contrast between two behaviors. On the one hand, if $u$~is a unipotent matrix, and $v$~is a vector that is not fixed by~$u$, then some component of $u^n \, v$ is a nonconstant polynomial, and therefore tends to~$\pm\infty$ with~$n$. 
On the other hand, the following observation implies that if $v$~is fixed by~$\Gamma$, then some subsequence converges to a finite limit.

\begin{lem}[(\thmindex{Poincar\'e Recurrence}{Poincar\'e Recurrence Theorem})]
\label{PoincareRecurThm}
 Let
 \begin{itemize}
 \item $(X,d)$ be a metric space,
 \item $T \colon X \to X$ be a homeomorphism, and
 \item $\mu$ be a $T$-invariant measure on~$X$, such that $\mu(X) < \infty$.
 \end{itemize}
 Then, for almost every $x \in X$, there is a
sequence $n_k \to \infty$, such that $T^{n_k} x \to x$.
 \end{lem}

\begin{proof}
 Let
 $$ A_\epsilon = \{\, a \in X \mid \forall m>0, \
d(T^m x, x) > \epsilon \,\}.$$
 It suffices to show $\mu(A_\epsilon) = 0$ for every
$\epsilon$.

Suppose $\mu(A_{\epsilon}) > 0$. Then we may choose a
subset~$B$ of~$A_\epsilon$, such that $\mu(B) > 0$ and
$\operatorname{diam}(B) < \epsilon$. Because the sets $B,
T^{-1}B, T^{-2} B, \ldots$ all have the same measure, and
$\mu(X) < \infty$, they cannot all be disjoint:
there exists $m < n$, such that $T^{-m} B \cap T^{-n} B
\neq \emptyset$. By applying $T^{n}$, we may assume $n =
0$. For $x \in T^{-m} B \cap B$, we have $T^m x \in B$ and
$x \in B$, so 
 $$d(T^m x, x) \le \operatorname{diam}(B) < \epsilon .$$
 This contradicts the definition of~$A_\epsilon$.
 \end{proof}
 
 \begin{rem}
Part~\pref{BDT-vector} of \cref{BDT} is a corollary of Part~\pref{BDT-subspace}. Namely, if $v$~is $\rho(\Gamma)$-invariant, then the $1$-dimensional subspace $\real v$ (or $\complex v$) is also invariant, so \pref{BDT-subspace} implies that the subspace is $\rho(G)$-invariant. Since $G$ has no nontrivial homomorphism to the abelian group~$\real^\times$ (or~$\complex^\times$), this implies that the vector~$v$ is $\rho(G)$-invariant. 
\end{rem}

However, we will provide a direct proof of~\pref{BDT-vector}, since it is quite short (and a little more elementary than the proof of~\pref{BDT-subspace}).

\begin{proof}[Proof of \fullcref{BDT}{vector}]
Suppose $v$ is a vector in~$V$ that is fixed by $\rho(\Gamma)$. It suffices to show, for every unipotent $u \in G$, that $v$~is fixed by $\rho(u)$.

Since $u$ is $\rho(\Gamma)$-invariant, the map $\rho$ induces a well-defined map $\overline\rho \colon G/\Gamma \to V$, defined by $\overline\rho(g\Gamma) = \rho(g) v$. Since $\overline\rho$ is $G$-equivariant, it pushes the $G$-invariant, finite measure~$\nu$ on $G/\Gamma$ to a $\rho(G)$-invariant, finite measure $\overline\nu$ on~$V$. Therefore, \cref{PoincareRecurThm} implies, for a.e.\ $g \in G$, that $\bigl\{ \rho(u^n g) v \bigr\}$ has a convergent subsequence.

However, each component of the vector $\rho(u^n g) v$ is a polynomial function of~$n$ \csee{un=poly}. Therefore, the preceding paragraph implies, for a.e.\ $g \in G$, that $\rho(u^n g) v$ is constant (independent of~$n$). This means that $\rho(g) v$ is fixed by $\rho(u)$. Since this is true for a.e.~$g$, we conclude, by continuity, that it is true for all~$g$, including $g = e$. Hence, $v$ is fixed by~$u$.
\end{proof}

To prepare for the proof of \fullcref{BDT}{subspace}, we make a few observations about the action of~$G$ on the projectivization of~$V$.

\begin{prop} \label{mu(PV)}
 Assume
 \begin{itemize}
 \item $G$ has no compact factors,
 \item $V$ is a finite-dimensional vector space over
$\real$ or~$\complex$,
 \item $\rho \colon G \to \GL(V)$ is a continuous homomorphism, and
 \item $\mu$ is a $\rho(G)$-invariant measure
on the projective space $\projective(V)$.%
	\nindex{$\projective(V)$ = projectivization of vector space~$V$}
 \end{itemize}
 If $\mu \bigl( \projective(V) \bigr) < \infty$, then $\mu$ is supported on the set of fixed points of $\rho(G^\circ)$.
 \end{prop}

\begin{proof}
 We know that $G^\circ$ is generated by its unipotent elements
\csee{<unip>=G}, so it suffices to show that  $\mu$ is
supported on the set of fixed points of $\rho(u)$, for
every unipotent element~$u$ of~$G$.

Let 
\noprelistbreak
	\begin{itemize}
	\item $u$ be a unipotent element of~$G$, 
	\item $T = \rho(u) - \Id$,
	and
	\item $v \in V \smallsetminus \{0\}$. 
	\end{itemize}
Then $T$ is nilpotent (because $\rho(u)$ is unipotent
\csee{rho(unip)}), so there is some integer
$r \ge 0$, such that $T^r v \neq 0$, but $T^{r+1}v = 0$.
 We have 
 $$\rho(u) T^r v 
 = (\Id + T)(T^r v)
 = T^r v + T^{r+1}v
 = T^r v + 0
 = T^r v ,$$
 so $[T^r v] \in \projective(V)$ is a fixed point for
$\rho(u)$. Also, for each $n \in \natural$, we have
 $$ \rho(u^n) [v]
 = \left[ \sum_{k=0}^r \binom{n}{k} T^k v
\right]
 = \left[ \binom{n}{r}^{-1}
 \sum_{k=0}^r \binom{n}{k} T^k v
\right]
 \to [T^r v] $$
 (because, for $k < r$, we have $\binom{n}{k}/\binom{n}{r}
\to 0$ as $n \to \infty$).
 Therefore, $\rho(u^n) [v]$ converges to a fixed point of
$\rho(u)$, as $n \to \infty$.

The Poincar\'e Recurrence Theorem \pref{PoincareRecurThm} implies,
for $\mu$-almost every $[v] \in \projective(V)$, that
there is a sequence $n_k \to \infty$, such that
$\rho(u^{n_k}) [v] \to [v]$. On the other hand, the preceding paragraph tells us that $\rho(u^{n_k}) [v]$
converges to a fixed point of $\rho(u)$. Therefore,
$\mu$-almost every element of $\projective(V)$ is a fixed
point of $\rho(u)$. In other words, $\mu$ is supported on
the set of fixed points of $\rho(u)$, as desired.
 \end{proof}

The assumption that $G$ has no compact factors cannot be
omitted from \cref{mu(PV)}. For example, the usual
Lebesgue measure is an $\SO(n)$-invariant, finite
measure on~$S^{n-1}$, but $\SO(n)$ has no fixed points
on~$S^{n-1}$. We can, however, make the following weaker
statement. 

\begin{cor} \label{mu(PV)cpct}
 Assume
 \begin{itemize}
 \item $V$ is a finite-dimensional vector space over
$\real$ or~$\complex$,
 \item $\rho \colon G \to \GL(V)$ is a continuous homomorphism, and
 \item $\mu$ is a $\rho(G)$-invariant measure on the projective space $\projective(V)$.
 \end{itemize}
 If $\mu \bigl( \projective(V) \bigr) < \infty$, then there is a cocompact, closed, normal subgroup~$G'$ of~$G$, such that $\mu$ is supported on the set of fixed points of $\rho(G')$.
 \end{cor}

\begin{proof}
 Let $K$ be the maximal connected, compact, normal
subgroup of~$G$, and write $G \approx G' \times K$, for
some closed, normal subgroup~$G'$ of~$G$. Then $G'$ has no
compact factors, so we may apply \cref{mu(PV)} to the
restriction $\rho|_{G'}$.
 \end{proof}

It is now easy to prove the other part of \cref{BDT}:

\begin{proof}[Proof of \fullcref{BDT}{subspace}]
By passing to a subgroup of finite index, we may assume $G$ is connected. For simplicity, let us also assume that $G$ has no compact factors \csee{BDT(subspace)cpct}.

Suppose $W$ is a subspace in~$V$ that is fixed by $\rho(\Gamma)$, and let $d = \dim W$. Note that $\rho$ induces a continuous homomorphism $\widehat\rho \colon G \to \GL(\bigwedge^d V)$, and, since $W$ is $\rho(\Gamma)$-invariant, the $1$-dimensional subspace $\bigwedge^d W$ is $\widehat\rho(\Gamma)$-invariant. Hence, $\widehat\rho$ induces a well-defined map $\overline\rho \colon G/\Gamma \to \projective \bigl( \bigwedge^d V \bigr)$, with $\overline\rho(e\Gamma) = [\bigwedge^d W]$. Then, since $\overline\rho$ is $G$-equivariant, it pushes the $G$-invariant, finite measure~$\nu$ on $G/\Gamma$ to a $\widehat\rho(G)$-invariant, finite measure $\overline\nu$ on $\projective(\bigwedge^d V)$. Then \cref{mu(PV)} tells us that $\overline\rho(G/\Gamma)$ is contained in the set of fixed points of $\overline\rho(G)$. In  particular, $[\bigwedge^d W] = \overline\rho(e\Gamma)$ is fixed by $\widehat\rho(G)$. This means that $W$ is $\rho(G)$-invariant.
 \end{proof}

\begin{rem} \label{BDTnotDiscrete}
The proofs of the two parts of \cref{BDT} never use the fact that the lattice~$\Gamma$ is discrete. Therefore, $\Gamma$ can be replaced with any closed subgroup~$H$ of~$G$, such that there is a $G$-invariant, finite measure on $G/H$, and $H$~projects densely into the maximal compact factor of~$G$.
\end{rem}

\begin{exercises}

\item \label{BDT(subspace)cpct}
 Complete the proof of \fullcref{BDT}{subspace}, by
removing the assumption that $G$ has no compact factors.
 \hint{See the hint to \cref{BDTCpctFactor}.}

%\item \label{Isom(M)finite}
% Suppose $M$ is a locally symmetric space, such that
%no irreducible factor of the universal cover of~$M$ is
%either flat or compact. Show that $\Isom(M)$ is finite. 

\item \label{BDT-G/H}
 Let $H$ be a closed subgroup of~$G$ that projects densely
into the maximal compact factor of~$G$.
 Show that if there is a $G$-invariant, finite measure on
$G/H$, then the identity component $H^\circ$ is a normal subgroup
of~$G^\circ$.
\hint{\Cref{BDTnotDiscrete} and the proof of \cref{BDT-normalize}.}

%\item Let $H$ be a closed, proper subgroup of~$G$, and assume $G$ is connected, simple, and noncompact. Show that if $H$ is not discrete, then there does not exist a $G$-invariant, finite measure on $G/H$.
%(If $H$ is not unimodular, then there is not
%even a $\sigma$-finite $G$-invariant measure on $G/H$
%\csee{HaarOnHomog<>unimod}, but we do not need this fact.)

\end{exercises}

\section{\texorpdfstring{$\Gamma$ is finitely presented}{Γ is finitely presented}}
\label{finpresSect}

\begin{defns}
Let $\Lambda$ be a group.
\noprelistbreak
	\begin{enumerate}
	\item $\Lambda$ is \defit[finitely!generated]{finitely generated} it has a finite generating set. That is, there is a finite subset of~$\Lambda$ that is not contained in any proper subgroup of~$\Lambda$.
	\item $\Lambda$ is \defit[finitely!presented]{finitely presented} it has a presentation with only finitely many generators and finitely many relations. In other words, there exist:
		\begin{itemize}
		\item a finitely generated free group~$F$, 
		\item a surjective homomorphism $\phi \colon F \to \Lambda$,
		and 
		\item a finite subset~$R$ of the kernel of~$\phi$,
		\end{itemize}
	such that $\ker\phi$ is the smallest \emph{normal} subgroup of~$F$ that contains~$R$.
	\end{enumerate}
It is easy to see that every finitely presented group is finitely generated. However, the converse is not true.
\end{defns}

In this section, we describe the proof that $\Gamma$ is
finitely presented. Much like the usual proof that the fundamental group of any compact manifold is finitely presented, it is based on the existence of a nice set that is close to being a fundamental domain for the action of~$\Gamma$ on~$G$.

\begin{defn} \label{CoarseFundDomDefn}
 Suppose $\Gamma$ acts properly discontinuously on a
topological space~$Y$.
 A subset $\fund$ 
 	\nindex{$\fund$ = (coarse) fundamental domain}
  of~$Y$ is a \defit[fundamental!domain!coarse]{coarse fundamental
domain} for~$\Gamma$ if%
\noprelistbreak
 \begin{enumerate}
 \item $\Gamma \fund = Y$, and
 \item \label{CoarseFundDomDefn-finite}
 $\{\, \gamma \in \Gamma \mid \gamma \fund \cap \fund \neq
\emptyset \,\}$ is finite.
 \end{enumerate}
 \end{defn}

\begin{other}
Some authors call~$\fund$ a \defit[fundamental!set]{fundamental set}, rather than a coarse fundamental domain.
\end{other}

The following general principle will be used to show that $\Gamma$ is finitely generated:

\begin{prop} \label{FundDom->fg}
 Suppose a discrete group~$\Lambda$ acts properly
discontinuously on a topological space~$Y$. If $Y$ is
connected, and $\Lambda$ has a coarse fundamental
domain~$\fund$ that is an open subset of~$Y$, then
$\Lambda$ is finitely generated.
 \end{prop}

\begin{proof}
 Let $S = \{\, s \in \Lambda \mid s \fund \cap
\fund \neq \emptyset \,\}$. We know that $S$ is finite
\fullcsee{CoarseFundDomDefn}{finite}, so it suffices to show
that $S$ generates~$\Lambda$. Here is the idea: think of
$\{\, \lambda \fund \mid \lambda \in \Lambda \,\}$ as a
tiling of~$Y$. The elements of~$S$ can move~$\fund$ to any
adjacent tile, and $Y$~is connected, so a composition of
elements of~$S$ can move~$\fund$ to any tile. Therefore
$\langle S \rangle$ is transitive on the set of tiles.
Since $S$ also contains the entire stabilizer of the
tile~$\fund$, we conclude that $\langle S \rangle = \fund$.

Now, here is the formal proof. 
Consider some $\lambda \in \Lambda\mk$, such that $\lambda \fund \cap \langle S \rangle \fund \neq \emptyset$. This means there exists $s \in \langle S \rangle$, such that $\lambda \fund \cap s \fund \neq \emptyset$, so $s^{-1} \lambda \fund \cap \fund \neq \emptyset$. Therefore, by the definition of~$S$, we have $s^{-1} \lambda \in S$, so $\lambda \in s S \subseteq \langle S\rangle$. 
Thus, we have shown that $\bigl( \Lambda \smallsetminus \langle S \rangle \bigr) \fund$ is disjoint from $\langle S \rangle \fund$. 

However, both of these sets are open (since $\fund$ is open), and their union is all of~$Y$ (since $\Lambda \fund = Y$). Therefore, since $Y$ is connected, the two sets cannot both be nonempty. Since $\langle S \rangle$ is obviously nonempty, we conclude that $\Lambda \smallsetminus \langle S \rangle = \emptyset$, so $\langle S \rangle = \Lambda$.
 \end{proof}

\begin{cor} \label{cocpct->fg}
 If\/ $\Gamma \backslash G$ is compact, then\/ $\Gamma$ is
finitely generated.
 \end{cor}

\begin{proof}
 Since $\Gamma \backslash G$ is compact, there is a
compact subset~$C$ of~$G$, such that $\Gamma C = G$
\csee{CpctInHomog}. Let $\fund$ be a precompact, open
subset of~$G$, such that $C \subseteq \fund$. Because $C
\subseteq \fund$, we have $\Gamma \fund = G$. Also, because
$\fund$ is precompact, and $\Gamma$ acts properly
discontinuously on~$G$, we know that
Condition~\fullref{CoarseFundDomDefn}{finite} holds. Therefore, $\fund$ is
a coarse fundamental domain for~$\Gamma$.
By passing to a subgroup of finite index, we may assume $G$~is connected \csee{FinInd->FG}, so \cref{FundDom->fg} applies.
 \end{proof}

%See \cref{CpctGen->fg} for an alternate proof of
%\cref{cocpct->fg} that does not require
%$G$~to be connected.

\begin{eg}
 Let $\fund$ be the closed unit square in~$\real^2$, so
$\fund$ is a coarse fundamental domain for the usual action
of~$\integer^2$ on~$\real^2$ by translations. Define~$S$
as in the proof of \cref{FundDom->fg}, so
 $$ S = \bigl\{\, (m,n) \in \integer^2 \mid m,n \in
\{-1,0,1\} \, \bigr\} 
 = \{0, \pm a_1, \pm a_2, \pm a_3 \} , $$
 where $a_1 = (1,0)$, $a_2 = (0,1)$, and $a_3 = (1,1)$.
Then $S$ generates $\integer^2$; in fact, the subset
$\{a_1,a_2\}$ is already a generating set.

\Cref{FundDom->fg} does not apply to this
situation, because $\fund$~is not open. We could enlarge
$\fund$ slightly, without changing~$S$. Alternatively, the
proposition can be proved under the weaker hypothesis that
$\fund$~is in the interior of $\bigcup_{s \in S} \fund$
\csee{FundDomNotOpen->fg}.

Note that $\integer^2$ has the presentation
 $$ \integer^2 = \langle\, x_1,x_2,x_3 \mid x_1 x_2 = x_3, \ 
x_2 x_1 = x_3 \,\rangle .$$
 (More precisely, if $F_3$ is
the free group on $3$~generators $x_1,x_2,x_3$, then there
is a surjective homomorphism $\phi \colon F_3 \to
\integer^2$, defined by
 $$ \phi(x_1) = a_1, \qquad \phi(x_2) = a_2, \qquad
\phi(x_3) = a_3 ,$$
 and the kernel of~$\phi$ is the smallest normal subgroup
of~$F_3$ that contains both $x_1 x_2 x_3^{-1}$ and $x_2
x_1 x_3^{-1}$.) Each of the relations in this presentation is of a very simple form, merely stating that the product of two elements of~$S$ is equal to
another element of~$S$. The proof of the following
\lcnamecref{FundDom->FinPres} shows that relations of this type suffice to
define~$\Lambda$ in a very general situation.
 \end{eg}

\begin{prop} \label{FundDom->FinPres}
 Suppose\/ $\Lambda$ acts properly discontinuously on a
topological space~$Y$. If 
	\begin{itemize}
	\item $Y$ is both connected and simply connected, 
	and 
	\item there is a coarse fundamental domain~$\fund$ for\/~$\Lambda$ that is a
connected, open subset of~$Y$, 
	\end{itemize}
then\/ $\Lambda$ is finitely presented.
 \end{prop}

\begin{proof}
 This is somewhat similar to the proof of \cref{FundDom->fg}, but is more elaborate.
 As before, let $S = \{\, \lambda
\in \Lambda \mid \lambda \fund \cap \fund \neq \emptyset
\,\}$. For each $s \in S$, define a formal symbol~$x_s$,
and let $F$ be the free group on $\{x_s\}$. Finally, let
	$$R = \{\, x_s x_t x_{st}^{-1} \mid s,t,st \in S \,\} ,$$
so $R$ is a finite subset of~$F$.

We have a homomorphism $\phi \colon F \to \Lambda$
determined by $\phi(x_s) = s$. From the proof of 
\cref{FundDom->fg}, we know that $\phi$ is
surjective, and it is clear that $R \subseteq \ker\phi$.
The main part of the proof is to show that $\ker \phi$ is
the smallest normal subgroup of~$F$ that contains~$R$.
(Since $R$ is finite, and $F/\ker \phi \iso \Lambda$, this
implies that $\Lambda$ is finitely presented, as desired.)

Let $N$ be the smallest normal subgroup of~$F$ that
contains~$R$. (It is clear that $N \subseteq \ker(\phi)$; we
wish to show $\ker(\phi) \subseteq N$.)
 \begin{itemize}
 
 \item Define an equivalence relation $\sim$ on $(F/N)
\times \fund$, by stipulating that
 $(f N, y) \sim (f' N , y')$
 if and only if there exists $s \in S$, such that $x_s fN =
f'N$ and $sy = y'$ \csee{EquivRelnOnFxY}. 

 \item Let $\widetilde{Y}$ be the quotient space $\bigl( (F/N) \times \fund \bigr) /
{\sim}$.

 \item Define a map $\psi \colon (F/N)
\times \fund \to Y$ by $\psi(fN,y) = \phi(f^{-1}) y$. (Note
that, because $N \subseteq \ker(\phi)$, the map~$\psi$ is
well defined.)
 \end{itemize}
 Because 
 $$\psi(x_s fN, sy) 
 = \bigl( \phi(f^{-1}) s^{-1} \bigr) (sy)
 = \psi(fN,y) ,$$
 we see that $\psi$ factors through to a
well-defined map $\widetilde \psi \colon \widetilde{Y} \to
Y$.

Let $\widetilde \fund$ be the image of $\bigl(\ker(\phi)/N
\bigr) \times \fund$ in~$\widetilde{Y}$. Then it is obvious,
from the definition of~$\psi$, that $\widetilde \psi(
\widetilde \fund) = \fund$. In fact, it is not difficult
to see that $\widetilde \psi^{-1}(\fund) = \widetilde
\fund$ \csee{FinPresInvImg}.

For each $f \in F$, the image~$\fund_f$ of $(f N) \times
\fund$ in~$\widetilde{Y}$ is open \csee{FinPresOpenSet}, and,
for $f_1,f_2 \in \ker(\phi)$, one can show that
$\fund_{f_1} \cap \fund_{f_2} = \emptyset$ if $f_1
\not\equiv f_2 \pmod{N}$ \ccf{FinPresPancakes}. Therefore,
from the preceding paragraph, we see that $\widetilde \psi$
is a covering map over~$\fund$. Since $Y$ is covered by
translates of~$\fund$ (and $\fund$ is open) it follows that $\widetilde \psi$
is a covering map. 
%Furthermore, the degree of this covering map is $|\ker(\phi)/N|$.

Because $\fund$ is connected, it is not difficult to see
that $\widetilde{Y}$ is connected \csee{FinPresYConn}.
Since $Y$ is simply connected, and $\widetilde \psi$ is a
covering map, this implies that $\widetilde \psi$ is a
homeomorphism. 
%Hence, the degree of the cover is~$1$
Hence, $\widetilde \psi$~is injective, and it is easy to see that this implies
$\ker(\phi) = N$, as desired.
 \end{proof}

\begin{rem} \label{LocConn->FinPres}
 The assumption that $\fund$ is connected can be replaced
with the assumption that $Y$ is \term[connected!locally]{locally connected}.
However, the proof is somewhat more complicated in this
setting.
 \end{rem}

\begin{cor}
 If\/ $\Gamma \backslash G$ is compact, then\/ $\Gamma$ is
finitely presented.
 \end{cor}

\begin{proof}
 Let $K$ be a maximal compact subgroup of~$G$,
 so $\Gamma$ acts properly discontinuously on $G/K$.
 Arguing as in the proof of
\cref{cocpct->fg}, we see that $\Gamma$ has a
coarse fundamental domain~$\fund$ that is an open subset of~$G/K$.
From the ``Iwasawa decomposition'' $G = KAN$ \csee{IwasawaDecomp}, we see that $G/K$ is connected and simply connected \fullcsee{G=KxRn}{G/K}.
So \cref{FundDom->FinPres} implies that $\Gamma$ is finitely presented.
 \end{proof}

If $\Gamma \backslash G$ is not compact, then it is more
difficult to prove that $\Gamma$ is finitely presented (or
even finitely generated). 

\begin{thm} \label{GammaFinPres} \label{GammaFinGen}
 $\Gamma$ is finitely presented.
 \end{thm}

\begin{proof}[Idea of proof]
 It suffices to find a coarse fundamental domain for~$\Gamma$
that is a connected, open subset of~$G/K$. Assume, without
loss of generality, that $\Gamma$ is irreducible. 

 In each of the following two cases, a coarse fundamental
domain~$\fund$ can be constructed as the union of
finitely many translates of ``\term{Siegel set}s\zz.'' (This will
be discussed in \cref{SLnZLattChap,ReductionChap}.)
 \begin{enumerate}
 \item \label{GammaFPrank2}
 $\Gamma$ is ``arithmetic\zz,'' as defined
in~\pref{ArithDefn}, or
 \item \label{GammaFinPresPf-Rank1}
 $G$ has a simple factor of real rank one, or,
more generally, we have $\Qrank \Gamma \le 1$ \csee{RrankDefn,Qrank-arithmetic}.
 \end{enumerate}
 The (amazing!) \thmindex{Margulis!Arithmeticity}{Margulis Arithmeticity Theorem}~\pref{MargulisArith} implies that these two cases
are exhaustive, which completes the proof.
 \end{proof}

\begin{rem}
 It is not necessary to appeal to the Margulis
Arithmeticity Theorem in order to prove only that
$\Gamma$ is finitely generated (and not that it is finitely presented).
Namely, if \pref{GammaFinPresPf-Rank1} does not apply,
then the real rank of every simple factor of~$G$ is at
least two, so \term[Kazhdan!property T@property~$(T)$]{Kazhdan's Property~$(T)$} implies that $\Gamma$ is finitely
generated  \fullcsee{KazhdanEasy}{fg}.
 \end{rem}

\begin{rem} \label{typeFn}
For $n \ge 2$, $\Gamma$~is said to be of \defit[type Fn@type~$F_n$]{type~$F_n$} if there is a compact CW complex~$X$, such that:
	\begin{itemize}
	\item the fundamental group $\pi_1(X)$ is isomorphic to~$\Gamma$,
	and
	\item the homotopy group $\pi_k(X)$ is trivial for $2 \le k < n$.
	\end{itemize}
Since $\Gamma$ is finitely presented \csee{GammaFinPres}, it is easy to see that $\Gamma$ is of type~$F_2$.
(In fact, a group is of type~$F_2$ if and only if it is finitely presented.)
Borel and Serre proved the much stronger result that $\Gamma$ is of type~$F_n$ for every~$n$. 
%In fact, there is a compact CW complex~$X$ as above, such that $\pi_k(X)$ is trivial for all $k \ge 2$.
(If $G/\Gamma$ is compact, and $\Gamma$~is torsion free, then one may let $X = \Gamma \backslash G/K$, where $K$~is a maximal compact subgroup of~$G$. (In other words, $X$ is the \index{locally!symmetric}{locally symmetric space} associated to~$\Gamma$.) When $G/\Gamma$ is not compact, but $\Gamma$~is torsion free, then $X$ is a certain space called the ``\term{Borel-Serre compactification}'' of $\Gamma \backslash G/K$.)
\end{rem}

\begin{exercises}

\item \label{FinInd->FG}
Show that if some finite-index subgroup of~$\Lambda$ is finitely generated, then $\Lambda$ is finitely generated.
\hint{If $F$ is a finite set of coset representatives for $\langle S \rangle$ in~$\Lambda$, then $S \cup F$ generates~$\Lambda$.}

\item \label{FGforFinInd}
Assume $\Lambda$ is abstractly commensurable to~$\Lambda'$.  Show $\Lambda$ is finitely generated if and only if $\Lambda'$ is finitely generated.
\hint{\Cref{FinInd->FG} is half of the proof. For the other half, suppose $F$ is a finite set of coset representatives for the subgroup~$\Lambda'$ of~$\Lambda\,$, and $S$~is a finite generating set for~$\Lambda\,$. For each $f \in F$ and $s \in S$, we have $fs \in \Lambda' f'$, for some $f'_{f,s} \in F$. Then $\{\,f s (f'_{f,s})^{-1} \mid f \in F, s \in S\,\}$ generates~$\Lambda'$. Alternatively, it is easy to prove this topologically: If $\Lambda$ is the fundamental group of a CW-complex~$\Sigma$ with only finitely many $1$-cells, then $\Lambda'$ is the fundamental group of a finite cover of~$\Sigma\,$, which must also have only finitely many $1$-cells.}

\item Assume $\Lambda$ is abstractly commensurable to~$\Lambda'$.  Show $\Lambda$ is finitely presented if and only if $\Lambda'$ is finitely presented. 
\hint{Suppose $\Lambda'$ has finite index in~$\Lambda$, and $F$ is a set of coset representatives. Let $S$ and $S'$ be finite generating sets of $\Lambda$ and~$\Lambda'$, respectively \csee{FGforFinInd}.
\par  
($\Leftarrow$) For each $s \in S \cup F$ and $f \in F$, there exist $g \in \Lambda'$ and $f' \in F$, such that $f s = g f'$. Adding these relations to a presentation of~$\Lambda'$ yields a presentation of~$\Lambda$.
\par
($\Rightarrow$) Proving this direction algebraically is somewhat more complicated, but there is an easy topological proof: If $\Lambda$ is the fundamental group of a CW-complex~$\Sigma$ whose $2$-skeleton is finite, then $\Lambda'$ is the fundamental group of a finite cover of~$\Sigma\,$, which must also have only finitely many $1$-cells and $2$-cells.}

\item \label{CpctGen->fg}
 Suppose $\Lambda$ is a discrete subgroup of a locally
compact group~$H$. Show that if $H/\Lambda$ is compact, and $H$~is
\defit{compactly generated} (that is, there is a compact
subset~$C$ of~$H$, such that $\langle C \rangle = H$),
then $\Lambda$ is finitely generated.
(This provides an alternate proof of \cref{FundDom->fg} that does not require $G$ to be connected.)
 \hint{Assume $e \in C$. Choose a compact subset~$\fund$ of~$H$, such
that $\fund \Lambda = H$ (and $e \in \fund$), and let $S =
\Lambda \cap (\fund^{-1} C^{\pm 1} \fund)$.
If $\lambda = c_1 c_2 \cdots c_n$ with $c_i \in C^{\pm 1}$, then $\lambda = \lambda_1 \cdots \lambda_n$, with $\lambda_i \in S$.}

\item \label{FundDomNotOpen->fg}
 Prove \cref{FundDom->fg}, replacing the assumption
that $\fund$ is open with the weaker assumption that
$\fund$ is in the interior of $\bigcup_{s \in S} s \fund$
(where $S$ is as defined in the proof of
\cref{FundDom->fg}).

\item \label{EquivRelnOnFxY}
 Show that the relation~$\sim$ defined in the proof of
\cref{FundDom->FinPres} is an equivalence relation.

\item \label{FinPresInvImg}
 In the notation of the proof of
\cref{FundDom->FinPres}, show that if $\psi(fN,y) \in
\fund$, then $(fN,y) \sim (f'N,y')$, for some $f' \in
\ker(\phi)$ and some $y' \in \fund$.
 \hint{We have $\phi(f) \in S$, because $\phi(f^{-1}) y \in
\fund$.}

\item \label{FinPresOpenSet}
 In the notation of the proof of
\cref{FundDom->FinPres},  show that the inverse image
of~$\fund_f$ in $(F/N) \times \fund$ is
 $$ \bigcup_{s \in S} \bigl( (x_s fN/N ) \times (\fund
\cap s\fund) \bigr) ,$$
 which is open.

\item \label{FinPresPancakes}
 In the notation of the proof of
\cref{FundDom->FinPres}, show that if we have $\fund_f \cap
\fund_e \neq \emptyset$ and $f \in \ker(\phi)$, then $f
\in N$.
 \hint{If $(f N, y_1) \sim (N, y_2)$, then there is some
$s \in S$ with $x_s N = fN$. Since $f \in \ker(\phi)$, we
have $s = \phi(x_s) = \phi(f) = e$.}

\item \label{FinPresYConn}
 Show that the set~$\widetilde{Y}$ defined in the proof of
\cref{FundDom->FinPres} is connected.
 \hint{For $s_1,\ldots,s_r \in S$, define $\fund_j =
\{ x_{s_j} \cdots x_{s_1} N \} \times \fund$. Show
there exist $a \in \fund_j$ and $b \in \fund_{j+1}$, such
that $a \sim b$.}

\item \label{CoarseFundIff}
Assume $\Lambda$ acts properly discontinuously on a
topological space~$Y$.
Show that a Borel subset $\fund$ of~$Y$ is a \term[fundamental!domain!coarse]{coarse fundamental domain} for~$\Lambda$ if and only if
 \begin{enumerate}
 \item $\fund$ contains a strict fundamental domain~$\fund_0$ for~$\Lambda$,
 and
 \item there is a finite subset~$F$ of~$\Lambda$, such that $\fund \subseteq F \fund_0$.
 \end{enumerate}

\end{exercises}

\section{\texorpdfstring
	{$\Gamma$ has a torsion-free subgroup of finite index}%
	{Γ has a torsion-free subgroup of finite index}}
	 \label{TorsionFreeSect}

\begin{defn}
 A group is \defit{torsion free} if it has no
nontrivial finite subgroups. Equivalently, the identity
element~$e$ is the only element of finite
order.
 \end{defn}

\begin{thm}[{}{(\thmindex{Selberg's Lemma}Selberg's Lemma)}] \label{torsionfree}
$\Gamma$ has a torsion-free
subgroup of finite index.
 \end{thm}

\begin{proof}
 From \cref{standassump}, we know $\Gamma \subseteq \SL(\ell,\real)$, for some~$\ell$. Let us start with an illustrative special case.

\setcounter{case}{0}

\begin{case} \label{torsionfree-SLnZ}
 Assume $\Gamma = \SL(\ell,\integer)$.
 \end{case}
 For any positive integer~$n$, the natural ring
homomorphism $\integer \to \integer/n\integer$ induces a
group homomorphism
 $\Gamma \to \SL(\ell,\integer/n\integer)$ \csee{SLfunctor};
 let $\Gamma_n$ be the kernel of this homomorphism. (This
is called the \label{PrincCongSubgrp}\defit{principal congruence subgroup} of
$\SL(\ell,\integer)$ of level~$n$.) Since it is the kernel of a 
group homomorphism, we know that $\Gamma_n$ is a normal 
subgroup of~$\Gamma$. It is also not difficult to see that
$\Gamma_n$ has finite index in~$\Gamma$ \csee{CSGfinite}. It therefore
suffices to show that $\Gamma_n$ is torsion free, for
some~$n$.
In fact, $\Gamma_n$ is torsion free whenever $n \ge
3$ \csee{TorsFreeAlln}, but, for simplicity, we will 
assume $n = p$ is an odd prime.

Given $\gamma \in \Gamma_p \smallsetminus \{\Id\}$ and $k
\in \natural \smallsetminus \{0\}$, we wish to show that
$\gamma^k \neq \Id$. We may write
 $$\gamma = \Id + p^d T ,$$
 where
 \begin{itemize}
 \item $d \ge 1$,
 \item $T \in
\Mat_{\ell \times \ell}(\integer)$, and
 \item $p \nmid T$ (that is, not every matrix entry of~$T$ is divisible by~$p$).
 \end{itemize}
 Also, we may assume $k$~is prime \csee{torsion->prime}.
Therefore, either $p \nmid k$ or $p = k$.

\begin{subcase}
 Assume $p \nmid k$.
 \end{subcase}
 Noting that
 $$ (p^d T)^2 = p^{2d} T^2 \equiv 0 \pmod{p^{d+1}} ,$$
 and using the Binomial Theorem, we see that
 $$ \gamma^k = (\Id + p^d T)^k
 \equiv \Id + k(p^d T) 
 \not\equiv \Id
 \pmod{p^{d+1}} ,$$
 as desired.

\begin{subcase}
 Assume $p = k$.
 \end{subcase}
 Using the Binomial Theorem (and noting that $\binom{p}{i} p^{di}$ is divisible by $p^{d + 2}$ for $i > 1$ \csee{TorsionErrorTerm}), we have
 $$ \gamma^k = \gamma^p = (\Id + p^d T)^p
 \equiv \Id + p (p^d T) 
 = \Id + p^{d+1} T
 \not\equiv \Id
 \pmod{p^{d+2}} .$$
% (Note that if $p = 2$, then the congruence requires $d
%\ge 2$ \csee{Gamma2HasTorsion}.)

\begin{case}
 Assume $\Gamma \subseteq \SL(\ell,\integer)$.
 \end{case}
  From \cref{torsionfree-SLnZ}, we know there is a
torsion-free, finite-index subgroup~$\Gamma_n$ of
$\SL(\ell, \integer)$. Then $\Gamma \cap \Gamma_n$ is a
torsion-free subgroup of finite index in~$\Gamma$.

\begin{case}
 The general case.
 \end{case}
 The proof is very similar to \cref{torsionfree-SLnZ},
with the addition of some commutative algebra (or
algebraic number theory) to account for the more general
setting.

 We know that $\Gamma$ is finitely generated \csee{GammaFinPres},
so there exist $a_1,\ldots,a_r \in \complex$, such that every
matrix entry of every element of~$\Gamma$ is contained in
the ring $Z = \integer[a_1,\ldots,a_r]$ generated by
$\{a_1,\ldots,a_r\}$ \csee{GammaFG->fgRing}. Therefore, letting
 $\Lambda = \SL(\ell, Z)$,
 we have $\Gamma \subseteq \Lambda$.

Now let $\mathfrak{p}$ be a maximal ideal in~$Z$. Then
$Z/\mathfrak{p}$ is a field, so, because $Z/\mathfrak{p}$
is also known to be a finitely generated ring, it must be a
finite field. Therefore, the kernel of the natural homomorphism
$\Lambda \to \SL(\ell, Z/\mathfrak{p})$ has finite index
in~$\Lambda$. Basic facts of Algebraic Number Theory allow
us to work with the prime ideal~$\mathfrak{p}$ in very
much the same way as we used the prime number~$p$ in
\cref{torsionfree-SLnZ}.
 \end{proof}

\begin{warn} \label{NoTorsionFreeSubrgpWarn}
 Our standing assumption that $G \subseteq \SL(\ell,\real)$ is needed for \cref{torsionfree}.
For example, the group $\Sp(4,\real)$ has an $8$-fold cover, which we call~$H$.
The inverse image of $\Sp(4,\integer)$ in~$H$ is a lattice~$\Lambda$ in~$H$.
It can be shown that every finite-index subgroup of~$\Lambda$ contains an element
of order~$2$, so no subgroup of finite index is torsion free.
%$\SL(4,\real)$ has a double cover, which we call~$H$.
%The inverse image of $\SL(4,\integer)$ in~$H$ is a lattice~$\Lambda$ in~$H$.
%It can be shown that every finite-index subgroup of~$\Lambda$ contains an element
%of order~$2$, so no subgroup of finite index is torsion free.
This does not contradict \cref{torsionfree}, because $H$~is not linear: it has no faithful embedding in any $\SL(\ell,\real)$.
 \end{warn}

If $\gamma^k = \Id$, then every eigenvalue of~$\gamma$
must be a $k^{\text{th}}$~root of unity. If, in addition,
$\gamma \neq \Id$, then at least one of these roots of
unity must be nontrivial. Therefore, the following is a
strengthening of \cref{torsionfree}. 

\begin{thm} \label{weaklynet}
 There is a finite-index subgroup\/~$\Gamma'$ of\/~$\Gamma$,
such that no eigenvalue of any element of\/~$\Gamma'$ is a
nontrivial root of unity.
 \end{thm}

\begin{proof}
 Assume $\Gamma = \SL(\ell,\integer)$. Let
 \begin{itemize}
 \item $n$~be some (large) natural number, 
 \item $\Gamma_n$ be the principal congruence subgroup
of~$\Gamma$ of level~$n$,
 \item $\omega$ be a nontrivial $k^{\text{th}}$~root of
unity, for some~$k$,
 \item $\gamma$ be an element of~$\Gamma_n$, such that
$\omega$ is an eigenvalue of~$\gamma$,
 \item $T = \gamma - \Id$, 
 \item $Q(x)$ be the characteristic polynomial of~$T$, and
 \item $\lambda = \omega - 1$, so $\lambda$ is a nonzero
eigenvalue of~$T$.
 \end{itemize}

Since $\gamma \in \Gamma_n$, we know that $n|T$, so $Q(x)
= x^\ell + n R(x)$, for some integral polynomial $R(x)$.
Since $Q(\lambda) = 0$, we conclude that $\lambda^\ell = n
\zeta$, for some $\zeta \in \integer[\lambda]$. Therefore,
$\lambda^\ell$ is divisible by~$n$, in the ring of
algebraic integers.

The proof can be completed by noting that any particular
nonzero algebraic integer is divisible by only finitely
many natural numbers, and there are only finitely many
roots of unity that satisfy a monic integral polynomial of
degree~$\ell$. See \cref{NoTorsOmega/p} for a slightly
different argument.
 \end{proof}

 \begin{rems} \ \label{SelbergRems}
 \noprelistbreak
	 \begin{enumerate}
	
	\item The proof of \cref{torsionfree} shows that $\Gamma$
	has nontrivial, proper, normal subgroups, so $\Gamma$ is
	not simple. However, the normal subgroups
	constructed there all have finite index. In fact, it is often the case that every nontrivial, normal
	subgroup of~$\Gamma$ has finite index
	\csee{MargNormalSubgrpsThm}. 
	
	Moreover, although it will not be proved in this book,
	it is often the case that all of the normal
	subgroups of finite index are close to being of the type
	constructed in the course of the proof. 
	More precisely, the
	``\thmindex{Congruence Subgroup Property}Congruence Subgroup Property''
	asserts there is a constant~$C$, such that if $N$ is any finite-index, normal subgroup of~$\Gamma$, then there is a principal congruence subgroup~$\Gamma'$ of~$\Gamma$, such that $|\Gamma' : N'| < C$. This is not always true, but it has been proved for the lattices in many groups.
	
	 \item \label{SelbergRems-net}
	 Arguing more carefully, one can obtain a
	finite-index subgroup~$\Gamma'$ with the stronger property
	that, for every $\gamma \in \Gamma'$, the multiplicative
	group generated by the (complex) eigenvalues of~$\gamma$
	does not contain any nontrivial roots of unity. Such a
	subgroup is sometimes called \defit[net subgroup]{net}.
	%(the French word for ``nice'').
	
	 \item \label{SelbergRems-F}
	 If $F$~is any field of characteristic zero, then \cref{torsionfree} remains valid (with the same proof) when $\Gamma$~is replaced with any finitely generated subgroup~$\Lambda$ of $\SL(\ell,F)$.
	  \end{enumerate}
 \end{rems}

Let us now present an alternate approach to the general
case of \cref{torsionfree}. It requires only the
Nullstellensatz, not Algebraic Number Theory.

\begin{proof}[Another proof of \cref{torsionfree} \optional]
 Let 
 \begin{itemize}
 \item $Z$ be the subring of~$\complex$ generated by the
matrix entries of the elements of~$\Gamma$, and
 \item $F$ be the quotient field of~$Z$.
 \end{itemize}
 Because $\Gamma$ is a finitely generated group
\csee{GammaFinPres}, we know that $Z$~is a finitely
generated ring \csee{GammaFG->fgRing}, so $F$ is a
finitely generated extension of~$\rational$.

\setcounter{step}{0}

\begin{step}
 We may assume that $F = \rational(x_1,\ldots,x_r)$ is a
purely transcendental extension of\/~$\rational$.
 \end{step}
 Choose a subfield $L = \rational(x_1,\ldots,x_r)$ of~$F$,
such that 
 \begin{itemize}
 \item $L$~is a purely transcendental extension
of~$\rational$, and 
 \item $F$ is an algebraic extension of~$L$.
 \end{itemize}
 Let $d$~be the degree of~$F$ over~$L$. Because $F$ is
finitely generated (and algebraic over~$L$), we know that
$d < \infty$. Therefore, we may identify $F^\ell$
with~$L^{d\ell}$, so there is an embedding
 $$ \Gamma \subseteq \SL(\ell,F) \hookrightarrow \SL(d\ell,
L) .$$
 Hence, by replacing~$F$ with~$L$ (and replacing~$\ell$
with~$d\ell$), we may assume that $F$ is purely
transcendental.
(Identifying $F^\ell$ with~$L^{d\ell}$ is the foundation of an important technique called ``\term{Restriction of Scalars}'' that will be introduced in \cref{RestrictScalarsSect}.)

\begin{step} \label{TorsionAltPf-traces}
 If $\gamma$ is any element of finite order in
$\SL(\ell,F)$, then $\trace(\gamma) \in \integer$, and
$|\trace(\gamma)| \le \ell$.
 \end{step}
 There is a positive integer~$k$ with $\gamma^k = \Id$, so
every eigenvalue of~$\gamma$ is a $k^{\text{th}}$~root of
unity. The trace of~$\gamma$ is the sum of these
eigenvalues, and any root of unity is an algebraic
integer, so we conclude that the trace of~$\gamma$ is an
algebraic integer. 

Since $\trace(\gamma)$ is the sum of the diagonal entries
of~$\gamma$, we know $\trace(\gamma) \in F$.  Since
$\trace(\gamma)$ is algebraic, but~$F$ is a purely
transcendental extension of~$\rational$, this implies
$\trace(\gamma) \in \rational$. Since $\trace(\gamma)$ is
an algebraic integer, this implies $\trace(\gamma) \in
\integer$.

Since $\trace(\gamma)$ is the sum of~$\ell$ roots of
unity, and every root of unity is on the unit circle, we
see, from the triangle inequality, that $|\trace(\gamma)|
\le \ell$.

\begin{step}
 There is a prime number $p > 2 \ell$, such that $1/p
\notin Z$.
 \end{step}
 From the \thmindex{Nullstellensatz}Nullstellensatz \pref{Nullstellensatz-ringhomo},
we know that there is a nontrivial homomorphism $\phi
\colon Z \to \overline{\rational}$, where
$\overline{\rational}$ is the algebraic closure
of~$\rational$ in~$\complex$. Replacing $Z$ with
$\phi(Z)$, let us assume that $Z \subset
\overline{\rational}$. Thus, for each $z \in Z$, there is
some nonzero integer~$n$, such that $nz$~is an algebraic
integer. More precisely, because $Z$~is finitely
generated, there is an integer~$n$, such that, for each $z
\in Z$, there is some positive integer~$k$, such that $n^k
z$~is an algebraic integer. It suffices to choose~$p$ so
that it is not a divisor of~$n$.

\begin{step}
 There is a finite field~$E$ of characteristic~$p$, and a
nontrivial homomorphism $\phi_p \colon Z \to E$.
 \end{step}
 Because $1/p \notin Z$, there is a maximal
ideal~$\mathfrak{p}$ of~$Z$, such that $p \in
\mathfrak{p}$. Then $E = Z/\mathfrak{p}$ is a field of
characteristic~$p$. Because it is a finitely generated
ring, $E$~must be a finite extension of the prime
field~$\integer/p \integer$ \csee{Nullstellensatz-fg}, so $E$~is
finite.

\begin{step}
 Let $\Lambda$ be the kernel of the homomorphism
$\hat\phi_p \colon \SL(\ell,Z) \to \SL(\ell, E)$ that is induced by~$\phi_p$. Then
$\Lambda$ is torsion free.
 \end{step}
 Let $\gamma$ be an element of finite order in~$\Lambda$.
Then 
 $$ \trace \bigl( \hat\phi_p(\gamma) \bigr) = \trace(\Id) =
\ell \pmod{p} ,$$
 so $p \mid \bigl( \ell - \trace(\gamma) \bigr)$ (since \cref{TorsionAltPf-traces} tells us that $\trace(\gamma) \in \integer$). Since we also know that $|\trace(\gamma)| \le \ell$ and $p > 2 \ell$, we conclude that $\trace(\gamma) = \ell$. 
 Since the $\ell$~eigenvalues of~$\gamma$ are roots of
unity, and $\trace(\gamma)$ is the sum of these
eigenvalues, we conclude that $1$~is the only eigenvalue
of~$\gamma$. Since $\gamma^k = \Id$, we know that
$\gamma$~is elliptic (hence, diagonalizable
over~$\complex$), so this implies $\gamma = \Id$, as
desired.
 \end{proof}

\begin{exercises}

\item \label{TorsionErrorTerm}
Show that if $p$ is an odd prime, $d \ge 1$, and $2 \le i \le p$, then $\binom{p}{i} p^{di}$ is divisible by $p^{d + 2}$.
\hint{If either $d > 1$ or $i > 2$, then $di \ge d + 2$.}

\item \label{SLfunctor}
 Show that $\SL(\ell,\cdot)$ is a (covariant) functor
from the category of rings with identity to the category of
groups. That is, show:
 \begin{enumerate}
 \item if $A$ is any ring with identity, then
$\SL(\ell,A)$ is a group,
 \item for every ring homomorphism $\phi \colon A \to B$
(with $\phi(1) = 1$), there is a group
homomorphism $\phi_* \colon \SL(\ell,A) \to \SL(\ell,B)$,
and
 \item if $\phi \colon A \to B$ and $\psi \colon B \to C$
are ring homomorphisms (with $\phi(1) = 1$ and $\psi(1) =
1$), then $(\psi \circ \phi)_* = \psi_* \circ \phi_*$.
 \end{enumerate}

\item \label{CSGfinite}
 Show that if $B$ is a finite ring with identity, then
$\SL(\ell,B)$ is finite. Use this fact to show, for every
positive integer~$n$, that if $\Gamma_n$ denotes the
principal congruence subgroup of $\SL(\ell,\integer)$ of
level~$n$ (cf.\ \cref{torsionfree-SLnZ} of the proof of 
\cref{torsionfree}), then $\Gamma_n$ has finite index
in $\SL(\ell,\integer)$.

\item \label{torsion->prime}
 Show that if $\Gamma$~has a nontrivial element of finite
order, then $\Gamma$~has an element of prime order.

\item \label{TorsFreeAlln}
Show the principal congruence subgroup~$\Gamma_n$ is torsion free if $n \ge 3$.
\hint{Since $\Gamma_m \subseteq \Gamma_n$ whenever $n \mid m$, you may assume, without loss of generality, that $n$~is either $4$ or an odd prime.}

\item \label{Gamma2HasTorsion}
 In the notation of \cref{torsionfree-SLnZ} of
the proof of  \cref{torsionfree}, show that $\Gamma_2$
is not torsion free. Where does your solution of
\cref{TorsFreeAlln} fail?

\item\label{GammaFG->fgRing}
 Show that if $\Lambda$ is a finitely generated subgroup of
$\SL(\ell,\complex)$, then there is a finitely generated
subring~$B$ of~$\complex$, such that $\Lambda \subset
\SL(\ell,B)$.
 \hint{Let $B$ be the subring of~$\complex$ generated by
the matrix entries of the generators of~$\Lambda$.}

\item \label{NoTorsOmega/p}
 Suppose $\omega$ is a nontrivial root of unity, and $(\omega-1)^\ell = n \zeta$, for some $n,\ell
\in \integer^+$ and
$\zeta \in \integer + \integer
\omega + \cdots + \integer \omega^{\ell-1}$. Show $n < 2^{(\ell+1)!}$.
 \hint{Let $F$ be the Galois closure of the field
extension $\rational(\omega)$ of~$\rational$ generated
by~$\omega$, and define $N \colon F
\to \rational$ by $N(x) = \prod_{\sigma \in
\operatorname{Gal}(F/\rational)} \sigma(x)$. Then
$N(\omega - 1)^\ell = n^d N(\zeta)$, and $|N(\omega - 1)|
\le 2^d \le 2^{\ell!}$, where $d$~is the degree of~$F$
over~$\rational$.}

\item \label{GammaResidFinite}
 Show that $\Gamma$ is \term{residually finite}. That is,
for every $\gamma \in \Gamma \smallsetminus \{e\}$, show
that there is a finite-index, normal subgroup~$\Gamma'$
of~$\Gamma$, such that $\gamma \notin \Gamma'$. (In
particular, if $\Gamma$ is infinite, then $\Gamma$ is \emph{not} a simple group.)

\item \label{GammaNormSubgrpsToE}
 Show there is a sequence $N_1, N_2, \ldots$ of
subgroups of~$\Gamma$, such that
 \begin{enumerate}
 \item $N_1 \supset N_2 \supset \cdots$,
 \item each $N_k$ is a finite-index, normal subgroup
of~$\Gamma$,
 and
 \item $N_1 \cap N_2 \cap \cdots = \{e\}$.
 \end{enumerate}
 \hint{Use \cref{GammaResidFinite}.}

\item \label{SLQCommSLZ}
 Show that $\SL(n,\rational)$ commensurates
$\SL(n,\integer)$. 
\hint{For each $g \in
\SL(n,\rational)$, there is a principal congruence
subgroup~$\Gamma_m$ of $\SL(n,\integer)$, such that $g^{-1}
\Gamma_m g \subseteq \SL(n,\integer)$.}

\item \label{ImgNotTorsion}
Show that if $\varphi \colon \Gamma \to \GL(n,\real)$ is a homomorphism, such that every element of $\varphi(\Gamma)$ has finite order, then $\varphi(\Gamma)$ is finite.  
\hint{\fullcref{SelbergRems}{F}.}

\end{exercises}

\section{\texorpdfstring{$\Gamma$ has a nonabelian free subgroup}%
	{Γ has a nonabelian free subgroup}}
\label{TitsAlternativeSect}

In this section, we describe the main ideas in the proof
of the following important result.

\begin{thm}[(\thmindex{Tits Alternative}Tits Alternative)]
\label{TitsAlternative}
 If $\Lambda$ is a subgroup of\/ $\SL(\ell,\real)$, then
either%
 \begin{enumerate}
 \item $\Lambda$ contains a nonabelian free group, or
 \item $\Lambda$ has a solvable subgroup of finite index.
 \end{enumerate}
 \end{thm}

Since $\Gamma$ is not solvable when $G$ is not compact \csee{GammaNotSolvable},
the following is an immediate corollary.

\begin{cor} \label{FreeInGamma}
 If $G$ is not compact, then\/ $\Gamma$ contains a nonabelian
free group.
 \end{cor}

\begin{defn}
 Let us say that a homeomorphism $\phi$ of a topological
space $M$ is \defit[contracting homeomorphism]{$(A^-,B,A^+)$-contracting} if $A^-$,
$B$ and~$A^+$ are nonempty, disjoint, open subsets of~$M$,
such that $\phi(B \cup A^+) \subseteq A^+$ and $\phi^{-1}(B
\cup A^-) \subseteq A^-$.
 \end{defn}

In a typical example, $A^-$ and~$A^+$ are small
neighborhoods of points~$p^-$ and~$p^+$, such that $\phi$
collapses a large open subset of~$M$ into~$A^+$, and
$\phi^{-1}$ collapses a large open subset of~$M$
into~$A^-$ \csee{Acontracting}.

\begin{figure}[ht]
 \includegraphics{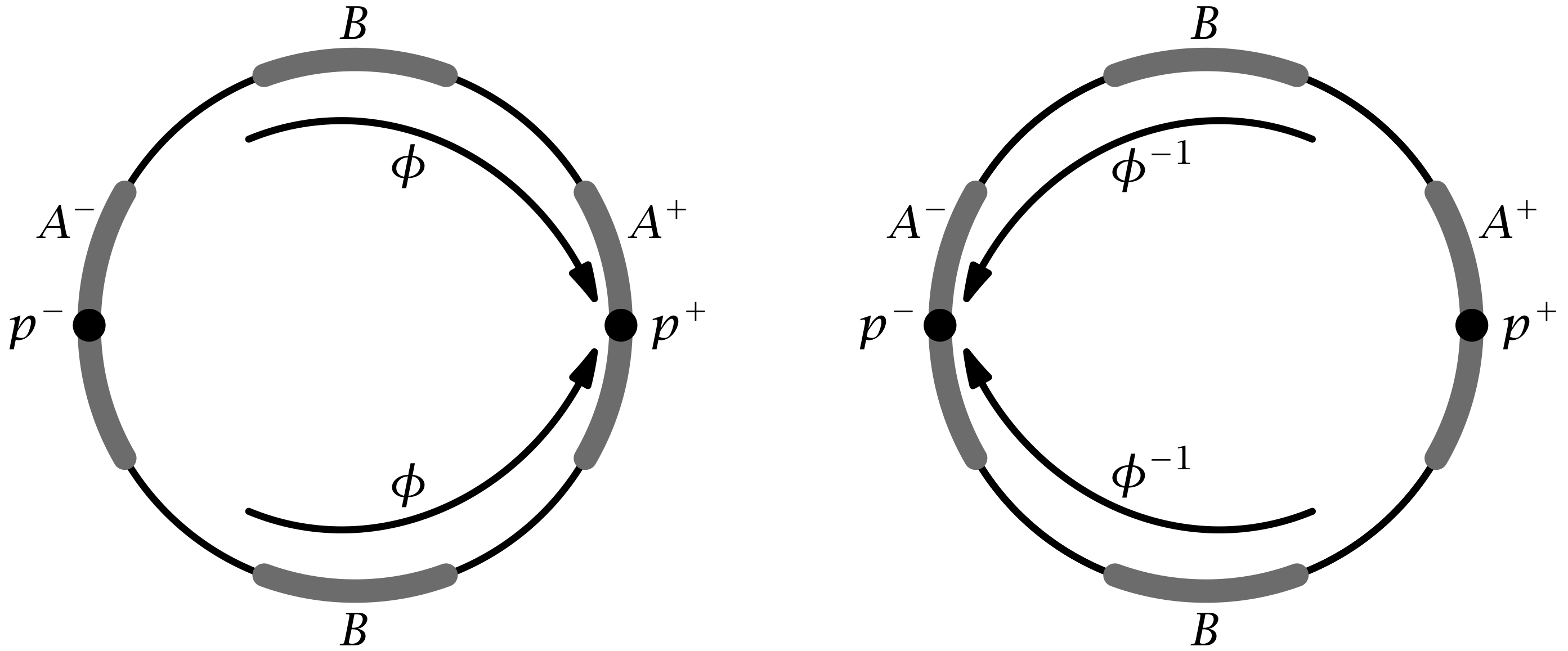}
 \caption{A typical $(A^-,B,A^+)$-contracting
homeomorphism of the circle.} \label{Acontracting}
 \end{figure}
%texpreamble
%(" \usepackage[LY1]{fontenc}
% \usepackage[expert, LY1, mylucidascale]{mylucidabr} % I adjusted the scaling
% \usepackage{amsmath}
% \everymath{\displaystyle}
% ");
%defaultpen(  fontcommand("\normalfont") + fontsize(10) ); 
%
%from graph access *;
%
%unitsize(2cm);
%
%real r = 1, b = 3.2;
%
%currentpen=linewidth(1.5);
%
%draw( circle( (b,0) , r) );
%
%for (int i = 0; i <=1; ++i){
%	draw( circle( (i*b,0) , r) );
%
%	draw( arc( (i*b,0), r, 150, 210) , linewidth(5)+gray);
%	draw( (i*b-r,0) , linewidth(7)); 
%	label("$p^-$", (i*b-r,0) , 1.5*W);
%	label( "$A^-$", (i*b-r-0.08, 0.4*r) );
%
%	draw( arc( (i*b,0), r, -30, 30) , linewidth(5)+gray);
%	draw( (i*b+r,0) , linewidth(7)); 
%	label("$p^+$", (i*b+r,0) , 2*E);
%	label( "$A^+$", (i*b+r, 0.4*r), 0.5*E );
%
%	draw( arc( (i*b,0), r, 70, 110) , linewidth(5)+gray);
%	label( "$B$", (i*b, r), 1.5*N );
%
%	draw( arc( (i*b,0), r, -70, -110) , linewidth(5)+gray);
%	label( "$B$", (i*b, -r), 1.5*S );
%	}
%
%draw( (-0.4*r, -0.7*r){ESE}..{NNE}(0.9*r, -0.1*r), EndArrow(7) );
%draw( (-0.4*r, 0.7*r){ENE}..{SSE}(0.9*r, 0.1*r), EndArrow(7) );
%label( "$\phi$", (0.2*r, -0.6*r) ); label( "$\phi$", (0.2*r, 0.6*r) );
%
%draw( (b + 0.4*r, -0.7*r){WSW}..{NNW}(b -0.9*r, -0.1*r), EndArrow(7) );
%draw( (b + 0.4*r, 0.7*r){WNW}..{SSW}(b -0.9*r, 0.1*r), EndArrow(7) );
%label( "$\phi^{-1}$", (b - 0.2*r, -0.55*r) ); label( "$\phi^{-1}$", (b - 0.2*r, 0.6*r) );

\begin{eg} \label{22contracts}
 Let 
 \noprelistbreak
 \begin{itemize}
 \item $M$ be the real projective line
$\projective(\real^2)$, 
 \item $\gamma = 
 \begin{bmatrix}
 2 & 0 \\
 0 & 1/2
 \end{bmatrix}
 \in \SL(2,\real)$,
 \item $A^-$ be any (small) neighborhood of $p^- = [0:1]$
in $\projective(\real^2)$,
 \item $A^+$ be any (small) neighborhood of $p^+ = [1:0]$
in $\projective(\real^2)$, and
 \item $B$ be any precompact, open subset of
$\projective(\real^2) \smallsetminus \{p^-,p^+\}$.
 \end{itemize}
 For any $(x,y) \in \real^2$ with $x \neq 0$, we have
  $$\text{$\gamma^n[x:y] = [2^nx:2^{-n}y] = [1:2^{-2n}y/x] \to
[1:0] = p^+$ \  as \ $n \to \infty$} , $$
 and the convergence is uniform on
compact subsets. Similarly, we have $\gamma^{-n}[x:y] \to p^-$ as
$n \to \infty$. Hence, for sufficiently large~$n$, the
homeomorphism $\gamma^n$ is $(A^-,B,A^+)$-contracting
on $\projective(\real^2)$.
 \end{eg}

More generally, if $\gamma$~is any nontrivial, hyperbolic
element of $\SL(2,\real)$, then $\gamma^n$ is
$(A^-,B,A^+)$-contracting on $\projective(\real^2)$, for
some appropriate choice of $A^-$, $B$, and~$A^+$
\csee{hypercontracts}.

The following is easy to prove by induction on~$n$.

\begin{lem} \label{TitsPhi(A)}
 If $\phi$ is $(A^-,B,A^+)$-contracting, then
 \begin{enumerate}
 \item \label{TitsPhi(A)-A+}
 $\phi^n(B) \subseteq A^+$ for all $n > 0$,
 \item $\phi^n(B) \subseteq A^-$ for all $n < 0$,
 \item \label{TitsPhi(A)-A}
 $\phi^n(B) \subseteq A^- \cup A^+$ for all $n \neq 0$.
 \end{enumerate}
 \end{lem}

The following lemma is the key to the proof of
\cref{TitsAlternative}.

\begin{lem}[(\thmindex{Ping-Pong Lemma}Ping-Pong Lemma)] \label{PingPong}
 Suppose 
 \begin{itemize}
 \item $\phi$ and $\psi$ are homeomorphisms of a
topological space~$M$,
 \item $A^-$, $A^+$, $B^-$, and~$B^+$ are nonempty,
pairwise-disjoint, open subsets of~$M$,
 \item $\phi$ is $(A^-,B,A^+)$-contracting, where $B = B^-
\cup B^+$, and
 \item $\psi$ is $(B^-,A,B^+)$-contracting, where $A = A^-
\cup A^+$.
 \end{itemize}
 Then $\phi$ and~$\psi$ have no nontrivial relations; so
$\langle \phi, \psi \rangle$ is free.
 \end{lem}

\begin{proof}
 Consider a word of the form $w = \phi^{m_1}
\psi^{n_1} \ldots \phi^{m_k} \psi^{n_k}$, with each $m_j$
and~$n_j$ nonzero. We wish to show $w \neq e$.

From \fullcref{TitsPhi(A)}{A}, we have
 $$ \mbox{$\phi^{m_j}(B) \subseteq A$\qquad
 and\qquad
 $\psi^{n_j}(A) \subseteq B$,\qquad
 for $j = 1,2,\ldots,k$.} $$

Therefore
 \begin{align*}
 \psi^{n_k}(A) &\subseteq B, \\
 \phi^{m_k} \psi^{n_k}(A) &\subseteq A, \\
 \psi^{n_{k-1}} \phi^{m_k} \psi^{n_k}(A) &\subseteq B, \\
 \phi^{m_{k-1}} \psi^{n_{k-1}} \phi^{m_k} \psi^{n_k}(A)
&\subseteq A, 
 \end{align*}
 and so on: points bounce back and forth between~$A$
and~$B$. (Hence, the name of the \lcnamecref{PingPong}.) In the
end, we see that $w(A) \subseteq A$.

 Assume, for definiteness, that $m_1 > 0$. Then, by
applying \fullref{TitsPhi(A)}{A+} in the last step, instead
of~\fullref{TitsPhi(A)}{A}, we obtain the
more precise conclusion that $w(A) \subseteq A^+$. Since $A
\not\subseteq A^+$ (recall that $A^-$~is disjoint
from~$A^+$), we conclude that $w \neq e$, as desired.
 \end{proof}

\begin{cor} \label{hyperSL2Rfree}
 If $\gamma_1$ and~$\gamma_2$ are two nontrivial hyperbolic
elements of $\SL(2,\real)$ that have no common
eigenvector, then, for sufficiently large $n \in
\integer^+$, the group $\langle (\gamma_1)^n, (\gamma_2)^n
\rangle$ is free.
 \end{cor}

\begin{proof}
 Let 
 \noprelistbreak
 \begin{itemize}
 \item $v_j$ and~$w_j$ be linearly independent eigenvectors
of~$\gamma_j$, with eigenvalues $\lambda_j$ and
$1/\lambda_j$, such that $\lambda_j > 1$,
 \item $A^+$ and~$A^-$ be small neighborhoods of $[v_1]$
and~$[w_1]$ in $\projective(\real^2)$, and
 \item $B^+$ and~$B^-$ be small neighborhoods of $[v_2]$
and~$[w_2]$ in $\projective(\real^2)$.
 \end{itemize}
 By the same argument as in \cref{22contracts}, we see
that if $n$~is sufficiently large, then
 \begin{itemize}
 \item $(\gamma_1)^n$ is $(A^-, B^- \cup B^+,
A^+)$-contracting, and
 \item $(\gamma_2)^n$ is $(B^-, A^- \cup A^+,
B^+)$-contracting
 \end{itemize}
 \csee{hypercontracts}. Therefore, the Ping-Pong Lemma
\pref{PingPong} implies that $\langle (\gamma_1)^n,
(\gamma_2)^n \rangle$ is free.
 \end{proof}

We can now give a direct proof of \cref{FreeInGamma},
in the special case where $G = \SL(2,\real)$.

\begin{cor}
 If $G = \SL(2,\real)$, then\/ $\Gamma$ contains a
nonabelian, free group.
 \end{cor}

\begin{proof}
 By passing to a subgroup of finite index, we may assume
that $\Gamma$ is torsion free \csee{torsionfree}. Hence,
$\Gamma$ has no elliptic elements. Not every element
of~$\Gamma$ is unipotent \csee{GammaHasNonUnip}, so we
conclude that some nontrivial element~$\gamma_1$
of~$\Gamma$ is hyperbolic.

Let $v$ and~$w$ be linearly independent eigenvectors
of~$\gamma_1$. The Borel Density Theorem
\pref{BDT(Zardense)} implies that there is some $\gamma
\in \Gamma$, such that $\{\gamma v , \gamma w \} \cap
(\real v \cup \real w) = \emptyset$
\csee{gammavnotinW}. Let $\gamma_2 = \gamma \gamma_1
\gamma^{-1}$, so $\gamma_2$~is a hyperbolic element
of~$\Gamma$ with eigenvectors $\gamma v$ and~$\gamma w$.

From \cref{hyperSL2Rfree}, we conclude that
$\langle (\gamma_1)^n, (\gamma_2)^n \rangle$ is a
nonabelian, free subgroup of~$\Gamma$, for some $n \in
\integer^+$.
 \end{proof}

The same ideas work in general:

\begin{proof}[Idea of direct proof of
\cref{FreeInGamma}]
 Assume $G \subseteq \SL(\ell,\real)$. Choose some
nontrivial, hyperbolic element~$\gamma_1$ of~$\Gamma$, and let $\lambda_1 \ge \lambda_2 \ge \cdots \ge
\lambda_\ell$ be its eigenvalues. We may assume, without loss of generality,
that $\lambda_1 > \lambda_2$. (If the
eigenvalue~$\lambda_1$ has multiplicity~$d$, then we may
pass to the $d^{\text{th}}$ exterior power
$\wedge^d(\real^\ell)$, to obtain a representation in which
the largest eigenvalue of~$\gamma_1$ is simple.)

Let us assume that the smallest eigenvalue~$\lambda_\ell$
is also simple; that is, $\lambda_\ell <
\lambda_{\ell-1}$. (One can show that this is a generic
condition in~$G$, so it can be achieved by
replacing~$\gamma_1$ with some other element of~$\Gamma$.)

Let $v$ be an eigenvector corresponding to the
eigenvalue~$\lambda_1$ of~$\gamma_1$, and let $w$ be an
eigenvector for the eigenvalue~$\lambda_\ell$. Assume, to
simplify the notation, that all of the eigenspaces
of~$\gamma_1$ are orthogonal to each other. Then, for any
$x \in \real^\ell \smallsetminus v^\perp$, we have
$(\gamma_1)^n[x] \to [v]$ in $\projective(\real^\ell)$, as
$n \to \infty$ \csee{SimpleEig->Contract}. Similarly, if
$x \notin w^\perp$, then $(\gamma_1)^{-n}[x] \to [w]$.

We may assume, by replacing $\real^\ell$ with a minimal
$G$-invariant subspace, that $\real^\ell$ has no nontrivial, 
proper, $G$-invariant subspaces. Then the Borel Density Theorem implies
that there exists $\gamma \in \Gamma$, such that we have
 $\{\gamma v, \gamma w\} \cap (\real v \cup \real w) =
\emptyset$.

Then, for any small neighborhoods $A^-$, $A^+$, $B^-$,
and~$B^+$ of $[v]$, $[w]$, $[\gamma v]$, and~$[\gamma w]$,
and any sufficiently large~$n$, the Ping-Pong
Lemma implies that the subgroup $\langle (\gamma_1)^n, (\gamma \gamma_1
\gamma^{-1})^n \rangle$ is free. 
 \end{proof}

\begin{rem}
 The proof of \cref{TitsAlternative} is similar, but
involves additional complications.
 \begin{enumerate}
 \item In order to replace $\real^\ell$ with an
irreducible subspace~$W$, it is necessary to have $\dim W >
1$ (otherwise, there do not exist two linearly independent
eigenvectors $v$ and~$w$). Unfortunately, the minimal
$\Lambda$-invariant subspaces may be $1$-dimensional. 
After modding these out, the minimal subspaces in the
quotient may also be $1$-dimensional, and so on. In this
case, the group~$\Lambda$ consists entirely of
upper-triangular matrices (after a change of basis), so
$\Lambda$ is solvable.
 \item The subgroup~$\Lambda$ may not have any hyperbolic
elements. Even worse, it may be the case that $1$~is the
absolute value of every eigenvalue of every element
of~$\Lambda$. (For example, $\Lambda$ may be a subgroup of
the compact group $\SO(n)$, so that every element
of~$\Lambda$ is elliptic.) In this case, the proof replaces the
usual absolute value with an appropriate $p$-adic norm.
Not all eigenvalues are roots of unity \ccf{weaklynet}, so
Algebraic Number Theory tells us that some element
of~$\Lambda$ has an eigenvalue whose $p$-adic norm is
greater than~$1$. The proof is completed by using this eigenvalue, and the
corresponding eigenvector, just as we used~$\lambda_1$ and
the corresponding eigenvector~$v$.
 \end{enumerate}
 \end{rem}

\begin{exercises}
\item \label{hypercontracts}
 In the notation of the proof of \cref{hyperSL2Rfree},
show that if $A^-$, $A^+$, $B^-$, and~$B^+$ are disjoint,
then, for all large~$n$, the homeomorphism $(\gamma_1)^n$
is $(A^-, B^- \cup B^+, A^+)$-contracting on
$\projective(\real^2)$.

\item \label{gammavnotinW}
 Assume that $G$ is irreducible in $\SL(\ell,\real)$
\csee{irredrepDefn}, and that $\Gamma$~projects densely
into the maximal compact factor of~$G$. If $F$~is a finite
subset of $\real^\ell \smallsetminus \{0\}$, and
$\mathcal{W}$~is a finite set of proper subspaces
of~$\real^\ell$, show that there exists $\gamma \in
\Gamma$, such that 
 $$ \gamma F \cap \bigcup_{W \in \mathcal{W}} W =
\emptyset .$$
 \hint{For $v \in F$ and $W \in \mathcal{W}$, the set
 $ A_{v,W} = \{\, g \in G \mid gv \in W \,\} $
 is Zariski closed, so $\bigcup_{v,W} A_{v,W}$ is Zariski
closed. Apply the Borel Density Theorem and
\cref{Conn->ZarConn}.}

\goodbreak % @@@

\item \label{SimpleEig->Contract}
 Let 
 \noprelistbreak
 \begin{itemize}
 \item $\gamma$ be a hyperbolic element of
$\SL(\ell,\real)$, 
\item  $\lambda_1 > \lambda_2
\ge \cdots \ge \lambda_\ell$ be the eigenvalues of~$\gamma$,
 \item $v$~be an eigenvector of~$\gamma$ corresponding to
the eigenvalue~$\lambda_1$, and 
 \item $W$~be the sum of the other eigenspaces. 
 \end{itemize}
 Show that if $x \in \projective(\real^\ell)
\smallsetminus [W]$, then $\gamma^n x \to [v]$ as $n \to
\infty$. Furthermore, the convergence is uniform on
compact subsets.

\item \label{PingPongABonly}
(another version of the Ping-Pong Lemma)
Suppose $A$ and~$B$ are disjoint, nonempty subsets of~$M$, such that $\phi^n(B) \subseteq A$ and $\psi^n(A) \subseteq B$, for every \emph{nonzero} integer~$n$.
Show $\langle \phi, \psi \rangle$ is free.
\hint{If every $m_j$ and~$n_j$ is nonzero, then $\phi^{m_1} \psi^{n_1} \cdots \phi^{m_k} \psi^{n_k} \phi^{m_{k+1}}(B) \subseteq A$.}

\item \label{SanovIsFree}
Let 
	$\Gamma' = \left\langle 
	\begin{Smallbmatrix} 1 & 2 \\ 0 & 1 \end{Smallbmatrix} ,
	\begin{Smallbmatrix} 1 & 0 \\ 2 & 1 \end{Smallbmatrix} \right\rangle $
be the \defit{Sanov subgroup} of $\SL(2,\integer)$.
Show that $\Gamma'$ is a free subgroup of finite index in $\SL(2,\integer)$.
\hint{\Cref{PingPongABonly} implies $\Gamma'$ is free. 
Matrices of the form $\begin{Smallbmatrix} 4k+1 & 2\ell \\ 2m & 4n+1 \end{Smallbmatrix}$ are in~$\Gamma'$.}

\item \label{SL2ZLotsFree}
Generalizing \cref{SanovIsFree}, show that every torsion-free subgroup of $\SL(2,\integer)$ is a free group.
\hint{%Many different proofs are possible, but here is one suggestion. 
Let $\bdry\fund$ be the boundary of the usual fundamental domain for the action of $\SL(2,\integer)$ on the upper half plane~$\hyperbolic^2$ \csee{FundDomSL2R}. Then $\bigcup_{\gamma \in \Gamma} \gamma \cdot \bdry\fund$ is a contractible $1$-dimensional simplicial complex; in other words, it is a tree. $\Gamma$~acts properly on this tree, so any subgroup of~$\Gamma$ that acts freely must be a free group.}

\item \label{irredinSL2xSO3}
 Show there is an irreducible lattice~$\Gamma$ in
$\SL(2,\real) \times \SO(3)$, such that $\Gamma \cap
\SL(2,\real)$ is infinite.
 \hint{There is a free group~$F$ and a homomorphism $\phi
\colon F \to \SO(3)$, such that $\phi(F)$ is dense in
$\SO(3)$.}

\end{exercises}

\section{Moore Ergodicity Theorem} \label{MooreErgBasicSect}

All mathematicians encounter situations in which they would like to prove that some function~$\varphi$ on some space~$X$ is constant. If $X = G/\Gamma$, this means that they would like to prove $\varphi $ is $G$-invariant. 

\begin{defn}
Suppose $f$ is a function on $G/\Gamma$, and $H$~is a subgroup of~$G$. We say that $f$~is \defit[invariant function@invariant function]{$H$-invariant} if $f(hx) = f(x)$ for all $h \in H$ and $x \in G/\Gamma$.
\end{defn}

The following fundamental result shows that it suffices to prove $\varphi $ is invariant under a much smaller subgroup of~$G$ (if we make the very weak assumption that $\varphi$ is measurable). 
It does not suffice to prove that $\varphi $ is invariant under a compact subgroup, because it is easy to find a non-constant, continuous function on $G/\Gamma$ that is invariant under any given compact subgroup of~$G$ (unless $G$ itself is compact) \csee{ExistKInvtFuncs}, so the following result is optimal.

\begin{thm}[(Moore Ergodicity Theorem)] \label{MooreErgThmGsimple}
Suppose
\noprelistbreak
	\begin{itemize}
	\item $G$ is connected and simple,
	\item $H$ is a closed, noncompact subgroup of~$G$,
	and
	\item $\varphi \colon G/\Gamma \to \complex$ is $H$-invariant and measurable.
	\end{itemize}
Then $\varphi$ is constant \textup(a.e.\textup).
\end{thm}

\Cref{MooreLattFromDiscreteLp} shows how to derive this theorem from the following more general result that replaces $\Gamma$ with a discrete subgroup that need not be a lattice. (This generalization will be a crucial ingredient in \cref{SLNZISLATTSlickSect}'s proof of the important fact that $\SL(n,\integer)$ is a lattice in $\SL(n,\real)$.) In this more general situation, we impose an $\LL{p}$-integrability hypothesis on~$\varphi$, in order to compensate for the fact that $G/\Lambda$ is not assumed to have finite measure \ccf{ExistInvtNotLp}.

\begin{thm} \label{MooreErgBasicThm}
Suppose
\noprelistbreak
	\begin{itemize}
	\item $G$ is connected and simple,
	\item $H$ is a closed, noncompact subgroup of~$G$,
	\item $\Lambda$ is a discrete subgroup of~$G$,
	and
	\item $\varphi$ is an $H$-invariant $\LL{p}$-function on~$G/\Lambda$ \textup(with $1 \le p < \infty$\textup).
	\end{itemize}
Then $\varphi$ is constant \textup(a.e.\textup).
\end{thm}

\begin{proof}[Idea of proof]
To illustrate the key ingredient in the proof, let us consider only the special case where 
	$G = \SL(2,\real)$
	and 
	 $H$~is the group of diagonal matrices. 
(A proof of the general case will be given in \cref{MooreErgPfSect}.) Let
	$$ a^t = \begin{Smallbmatrix} e^t & 0 \\ 0 & e^{-t} \end{Smallbmatrix} \in H
	\qquad \text{and} \qquad
	u \in \begin{Smallbmatrix} 1& 0 \\ * & 1 \end{Smallbmatrix} .$$
Note that straightforward matrix multiplication \csee{MooreErg-commutationEx} shows
	\begin{align} \label{MooreErg-commutation}
	\lim_{t \to \infty} a^t u a^{-t} = e
	. \end{align}

For $g \in G$, define $g \varphi \colon G/\Lambda \to \complex$ by $g\varphi(x) = \varphi(g^{-1} x)$.
We are assuming that $\varphi$~is $a^t$-invariant (which means $a^t \varphi = \varphi$), and the crux of the proof is the observation that we can use \pref{MooreErg-commutation} to show that $\varphi$~must also be $u$-invariant:
we have
	\begin{align*}
	 \|u \varphi - \varphi \|_p
	&= \| a^t u \varphi - a^t \varphi \|_p 
	&& \text{(\cref{TranslateHasSameIntegral})}
	\\&= \| (a^t u a^{-t}) a^t \varphi - a^t \varphi\|_p 
	&& \text{(inserting $a^{-t} a^t$)}
	\\&= \| (a^t u a^{-t}) \varphi - \varphi \|_p 
	&& \begin{pmatrix} \text{$a^t \varphi = \varphi$ because} \\  \text{$\varphi$ is $H$-invariant} \end{pmatrix}
	\\&\to \| e \varphi - \varphi \|_p \quad \text{as $t \to \infty$} 
	&& \begin{pmatrix} \text{\pref{MooreErg-commutation} and} \\ \text{\cref{GContOnLp}} \end{pmatrix}
	\\&= 0
	, \end{align*}
so $u\varphi = \varphi$ (a.e.). 

Thus, from the fact that $\varphi$~is $H$-invariant, we have shown that 
	$$ \text{$\varphi$~must also be $\left[\begin{smallmatrix} 1& 0 \\ * & 1 \end{smallmatrix}\right]$-invariant (a.e.).} $$
The same calculation, but with $t \to -\infty$, shows that 
	$$ \text{$\varphi$ must also be $\left[\begin{smallmatrix} 1& * \\ 0 & 1 \end{smallmatrix}\right]$-invariant (a.e.).} $$
Since $\begin{Smallbmatrix} 1& 0 \\ * & 1 \end{Smallbmatrix}$ and $\begin{Smallbmatrix} 1& * \\ 0 & 1 \end{Smallbmatrix}$ generate $\SL(2,\real) = G$ \csee{ElemsGenSL2R}, we conclude that $\varphi$~is $G$-invariant (a.e.). Since $G$ is transitive on $G/\Lambda$, this implies that $\varphi$ is constant (a.e.) \csee{InvtOnTrans}.
\end{proof}

%See \cref{ErgodicChap} for more discussion of results of this type.

\begin{exercises}

\item \label{ExistKInvtFuncs}
Show that if $K$ is any compact subgroup of~$G$, and $G$ is not compact, then there is a continuous, $K$-invariant function on $G/\Gamma$ that is not constant.
%\hint{One possibility is to fix a point $p \in G/\Gamma$, and let $f(x)$ be the distance from $x$ to the $K$-orbit of~$p$.}

\item \label{ExistInvtNotLp}
Show there is a counterexample to \cref{MooreErgBasicThm} if we remove the assumption that the measurable function~$\varphi$ is~$\LL{p}$.
\hint{It is easy to construct an counterexample by taking $\Lambda$ to be trivial (or finite).}

\item \label{MooreLattFromDiscreteLp}
Derive \cref{MooreErgThmGsimple} from \cref{MooreErgBasicThm}.
\hint{If there is a nonconstant $H$-invariant function on $G/\Gamma$, then there is one that is bounded.}

\item \label{MooreErg-commutationEx}
Verify \cref{MooreErg-commutation}.

\item \label{MooreErgBasicInvtSetEx}
Suppose
	\begin{itemize} 
	\item $G$ and $H$ are as in the Moore Ergodicity Theorem \pref{MooreErgThmGsimple},
	and
	\item $X$ is an $H$-invariant, measurable subset of $G/\Gamma$.
	\end{itemize}
Show that either $X$ has measure~$0$, or the complement of~$X$ has measure~$0$.
%% \hint{The characteristic function of~$X$ is $H$-invariant.}
%\item Where does your proof use the assumption that $\Gamma$ is a lattice, rather than an arbitrary discrete subgroup?
%%\item Show there is a counterexample if we omit the assumption  that $\Lambda$ is a lattice.

\item \label{ElemsGenSL2R}
Show $\SL(2,\real)$ is generated by the subgroups
$\left\{\begin{Smallbmatrix} 1& 0 \\ * & 1 \end{Smallbmatrix}\right\}$ and $\left\{\begin{Smallbmatrix} 1& * \\ 0 & 1 \end{Smallbmatrix}\right\}$.
\hint{If a matrix can be reduced to the identity matrix by a sequence of elementary row operations, then it is a product of elementary matrices.}

\item \label{InvtOnTrans}
Suppose $\varphi$ is a measurable function on $G/\Lambda$, and for each $g \in G$, we have $\varphi(gx) = \varphi(x)$ for a.e.\ $x \in G/\Lambda$. Show $\varphi$ is constant (a.e.).
\hint{Use Fubini's Theorem to reverse the quantifiers.}

\item Suppose $G$ and $H$ are as in the Moore Ergodicity Theorem \pref{MooreErgThmGsimple}.
Show that $Hx$ is dense in $G/\Gamma$, for a.e.\ $x \in G/\Gamma$.
\hint{Use \cref{MooreErgBasicInvtSetEx}. For any open subset~$\open$ of $G/\Gamma$, the set $H\open$ is measurable (why?) and $H$-invariant.}

\item Assume $G$ is simple, and let $H$ be a subgroup of~$G$ (not necessarily closed). Show that every (real-valued) $H$-invariant measurable function on $G/\Gamma$ is constant (a.e.) if and only if the closure of $H$ is not compact.

\item \label{TranslateHasSameIntegral}
Show that if $\Lambda$ is a discrete subgroup of~$G$, and $\mu$ is a $G$-invariant measure on $G/\Lambda$, then $\int_X g\varphi \, d\mu = \int_X \varphi \, d\mu$, for every $g \in G$ and measurable $\varphi \colon G/\Lambda \to \complex$.
\hint{Apply a change of variables, and use the fact that $g_*\mu = \mu$ (since $\mu$ is $G$-invariant).}

\item \label{GContOnLp}
Show that $G$ acts continuously on $\LL{p}(G/\Lambda)$ if $1 \le p < \infty$. More precisely, show that if $\Lambda$ is a discrete subgroup of~$G$, and we define $\alpha \colon G \times \LL{p}(G/\Lambda) \to \LL{p}(G/\Lambda)$ by $\alpha(g,\varphi) = g\varphi$, then $\alpha$ is continuous.
\hint{To show $\alpha$ is continuous in~$g$, use \thmindex{Lusin's}{Lusin's Theorem} \pref{LusinsThm} to approximate~$\varphi$ by a uniformly continuous function. Then use \cref{TranslateHasSameIntegral} (and the Triangle Inequality) to complete the proof.}

\end{exercises}

\begin{notes}

Raghunathan's book \cite{RaghunathanBook} is the standard
reference for the basic properties of lattices. It contains
almost all of the material in this chapter, except the Tits
Alternative (\cref{TitsAlternativeSect}) and the Moore Ergodicity Theorem (\cref{MooreErgBasicSect}).

%Isogenies are discussed in texts, such as
%\cite{Borel-AlgicGrps}, on algebraic groups and in the book
%of Platonov and Rapinchuk \cite{PlatonovRapinchukBook}.
 
\Cref{GodementConverse} (the existence of unipotent elements in noncompact lattices) was proved by Kazhdan and Margulis \cite{KazhdanMargulis-PfSelberg}. Expositions can be found in \cite{Borel-KazhdanMargulisBourbaki} and \cite[Cor.~11.13, p.~180]{RaghunathanBook}. 

The Borel Density Theorem \pref{BDT} was proved by
Borel~\cite{Borel-BDT}. It appears in
\cite[Thm.~2.4.4, p.~93]{MargulisBook}, \cite[Thm.~5.5, p.~79]{RaghunathanBook}, and
\cite[Thm.~3.2.5, pp.~41--42]{ZimmerBook}. Several authors have published
generalizations or alternative proofs
(for example, \cite{Dani-BDT, Furstenberg-BDT,
Wigner-BDT}).
%\cite{CowlingSteger-BDT, IozziNevo-BDT, Shalom-BDT, Stuck-BDT}

Our presentation of \cref{FundDom->fg,FundDom->FinPres} is based on
\cite[pp.~195--199]{PlatonovRapinchukBook}. A proof of
\cref{LocConn->FinPres} can also be found there.
A proof of \cref{GammaFinPres} for the case where
$\Gamma$ is arithmetic can be found in
\cite{Borel-IntroGrpArith} or \cite[Thm.~4.2,
p.~195]{PlatonovRapinchukBook}. For the case where $\Qrank
\Gamma = 1$, see \cite{GarlandRaghunathan} or
\cite[Cor.~13.20, p.~210]{RaghunathanBook}.

Borel and Serre \cite[\S11.1]{BorelSerre-corners} proved $\Gamma$ is of type~$F_n$ \csee{typeFn}.
%This also follows from the theorem of Raghunathan \cite{Raghunathan-quotients} that $\Gamma \backslash G / K$ is homeomorphic to the interior of a compact manifold with boundary.
(We remark that there is no harm in assuming $\Gamma$ is torsion free, since being of type~$F_n$ is invariant under passage to finite-index subgroups \cite[Cor.~7.2.4, p.~170]{Geoghegan-TopMethGrpThy}.)

\Cref{torsionfree}
is proved in \cite[Thm.~6.11, p.~93]{RaghunathanBook} and
\cite[Cor.~17.7, p.~119]{Borel-IntroGrpArith}, in stronger forms that establish \cref{SelbergRems}(\ref{SelbergRems-net},\ref{SelbergRems-F}).
Our alternate proof of
\cref{torsionfree} is excerpted from the elementary
proof in \cite{Alperin-SelbergLemma}.

\Cref{NoTorsionFreeSubrgpWarn} is due to P.\,Deligne \cite{Deligne-torsion}.
See also \cite{Raghunathan-torsion}.

For an introduction to the Congruence Subgroup Property, see \cite[Chap.~6]{Humphreys-ArithmeticGroups} or \cite{Sury-CSPBook}.

The Tits Alternative \pref{TitsAlternative} was proved by
Tits \cite{Tits-TitsAlt}. A nice introduction (and a proof
of some special cases) can be found in
\cite{delaHarpe-TitsAlt}.

See \cref{ErgodicitySect} for more on the Moore Ergodicity Theorem \pref{MooreErgThmGsimple} and related results.

\end{notes}

  	\standassumptrue

 %!TEX root = IntroArithGrps.tex

\mychapter{\texorpdfstring{What is an\\Arithmetic Group?}%
	{What is an Arithmetic Group?}}
\label{ArithGrpsChap}

\prereqs{\Cref{LatticeDefn,CommensDefn} (lattice subgroups and commensurability).}

$\SL(n,\integer)$ is the most basic example of an ``arithmetic group\zz.'' We will see that, by definition, the other arithmetic groups are obtained by intersecting $\SL(n,\integer)$ with some semisimple subgroup~$G$ of $\SL(n,\real)$. More precisely, if $G$ is a subgroup of $\SL(n,\real)$ that satisfies certain technical conditions (to be explained in \cref{ArithLattDefnSect}), then $G \cap \SL(n,\integer)$ (the group of ``integer points'' of~$G$) is said to be an \defit[arithmetic!subgroup]{arithmetic subgroup} of~$G$. However, the official definition \pref{ArithDefn} also allows certain modifications of this subgroup to be called \defit[arithmetic!subgroup]{arithmetic}.

Different embeddings of~$G$ into $\SL(n,\real)$ can yield different intersections with $\SL(n,\integer)$, so $G$ has many different arithmetic subgroups. 
(Examples can be found in \cref{EgArithGrpsChap}.) 
\Cref{arith->latt} tells us that all of them are lattices in~$G$. In particular, $\SL(n,\integer)$ is a lattice in $\SL(n,\real)$.

\section{Definition of arithmetic subgroups} \label{ArithLattDefnSect}

We are assuming that $G$ is a subgroup of $\SL(\ell,\real)$ \csee{standassump}, and we are interested in $\Gamma = G \cap \SL(\ell,\integer)$, the set of 
``integer points'' of~$G$. However, in order for the integer
points to form a lattice, $G$~needs to be well-placed with
respect to $\SL(\ell,\integer)$. (If we replace $G$ by a
conjugate under some terrible irrational matrix, perhaps $G
\cap \SL(\ell,\integer)$ would become trivial
\csee{conjnoZ}.) The following \lcnamecref{RnDefdQ<>Latt} is an elementary
illustration of this idea.

\begin{prop} \label{RnDefdQ<>Latt}
 The following are
equivalent, for every subspace~$W$ of\/~$\real^\ell$:
 \begin{enumerate}
 \item \label{RnDefdQ<>Latt-latt}
 $W \cap \integer^\ell$ is a cocompact lattice in~$W$.
 \item \label{RnDefdQ<>Latt-span}
 $W$ is spanned by $W \cap \integer^\ell$.
 \item \label{RnDefdQ<>Latt-dense}
 $W \cap \rational^\ell$ is dense in~$W$.
 \item \label{RnDefdQ<>Latt-eqs}
 $W$ can be defined by a set of linear equations with coefficients in~$\rational$.
\end{enumerate}
 \end{prop}

\begin{proof}
Let $k = \dim W$.

($\ref{RnDefdQ<>Latt-latt} \Rightarrow
\ref{RnDefdQ<>Latt-span}$) Let $V$ be the $\real$-span of $W
\cap \integer^\ell$. Then $W/V$, being a vector space
over~$\real$, is homeomorphic to~$\real^d$, for some~$d$. On
the other hand, we know that $W \cap \integer^\ell \subset
V$, and that $W/(W \cap \integer^\ell)$ is compact, so $W/V$
is compact. Hence $d = 0$, so $V = W$.

($\ref{RnDefdQ<>Latt-span} \Rightarrow \ref{RnDefdQ<>Latt-latt}$) 
Let $\{\varepsilon_1,\ldots,\varepsilon_k\}$ be the standard
basis of~$\real^k$. 
Because $W \cap \integer^\ell$ contains a basis of~$W$, there is a
linear isomorphism $T \colon \real^k \to W$, such that
$T\bigl( \{\varepsilon_1,\ldots,\varepsilon_k\} \bigr) \subseteq W \cap \integer^\ell$. This implies that $T(\integer^k) \subseteq W \cap
\integer^\ell$. Since $\real^k/\integer^k$ is compact, and $T$~is continuous, we
conclude that $W/(W \cap \integer^\ell)$ is compact.

 ($\ref{RnDefdQ<>Latt-span} \Rightarrow \ref{RnDefdQ<>Latt-dense}$) 
 As in the proof of ($\ref{RnDefdQ<>Latt-span} \Rightarrow \ref{RnDefdQ<>Latt-latt}$), 
there is a linear isomorphism $T \colon \real^k \to W$, such that
$T(\integer^k) \subseteq W \cap \integer^\ell$. 
Then $T(\rational^k) \subseteq W \cap
\rational^\ell$. Since $\rational^k$ is dense in~$\real^k$,
and $T$~is continuous, we conclude that $T(\rational^k)$ is
dense in~$W$.

 ($\ref{RnDefdQ<>Latt-eqs} \Rightarrow
\ref{RnDefdQ<>Latt-span}$) By assumption, $W$~is the
solution space of a system of linear equations whose
coefficients belong to~$\rational$. (Since $\real^\ell$ is 
finite dimensional, only finitely many of the equations are 
necessary.) Therefore, by elementary linear algebra (row 
reductions), we may find a basis for~$W$ that
consists entirely of vectors in~$\rational^\ell$.
Multiplying by a scalar to clear the denominators, we may
assume that the basis consists entirely of vectors
in~$\integer^\ell$.

 ($\ref{RnDefdQ<>Latt-dense} \Rightarrow
\ref{RnDefdQ<>Latt-eqs}$) 
Since $W \cap \rational^\ell$ is
dense in~$W$, we know that the orthogonal complement $W^\perp$ is defined by a set
of linear equations with rational coefficients. (For each 
$w \in W \cap \rational^\ell$, we write the equation 
$w \cdot x = 0$.) Thus, from 
($\ref{RnDefdQ<>Latt-eqs} \Rightarrow
\ref{RnDefdQ<>Latt-span}$), we conclude that there is a
basis $v_1,\ldots,v_m$ of~$W^\perp$, such that each $v_j \in
\rational^\ell$. Then $W = (W^\perp)^\perp$ is defined by the
system of equations $v_1 \cdot x= 0$, \dots, $v_m \cdot x =
0$.
 \end{proof}

With the above \lcnamecref{RnDefdQ<>Latt} in mind, we make the following
definition.

% cref won't combine these references, because one is defn and the other is defns !!!
\begin{defn}[(cf.\ Definitions \ref{AlgicGrpDefn} and~\ref{AlmZarDefn})] \label{DefdQDefn}
 Let $H$ be a closed subgroup of $\SL(\ell,\real)$. We say
 $H$~is \index{algebraic!group!over Q@over~$\rational$}\defit[defined!over Q@over~$\rational$]{defined over\/~$\rational$} (or that $H$ is a
\defit[Q-@$\rational$-!subgroup]{$\rational$-subgroup\/}) if there is a
subset~$\mathcal{Q}$ of $\rational[x_{1,1}, \ldots,
x_{\ell,\ell}]$, such that 
 \begin{itemize}
 \item \nindex{$\Var(\mathcal{Q})$ = $\{\, g \in \SL(\ell,\real) \mid
Q(g) = 0, \ \forall Q \in \mathcal{Q} \,\}$}
 $\Var(\mathcal{Q}) = \{\, g \in \SL(\ell,\real) \mid
Q(g) = 0, \ \forall Q \in \mathcal{Q} \,\}$ is a subgroup of
$\SL(\ell,\real)$,
 \item $H^\circ = \Var(\mathcal{Q})^\circ$, and
 \item $H$ has only finitely many components.
 \end{itemize}
 In other words, $H$~is commensurable to the variety
$\Var(\mathcal{Q})$, for some set~$\mathcal{Q}$ of
$\rational$-polynomials.
 \end{defn}

\begin{egs} \ 
\noprelistbreak
 \begin{enumerate}
 \item $\SL(\ell,\real)$ is defined over~$\rational$: let
$\mathcal{Q} = \emptyset$.
 \item If $n < \ell$, we may embed $\SL(n,\real)$ in the top
left corner of $\SL(\ell,\real)$. This copy of $\SL(n,\real)$ is defined over~$\rational$: let
 $$ \mathcal{Q} = 
 \{\, x_{i,j} - \delta_i^j \mid \max\{i,j\} > n \,\} .$$
 \item For $A \in \SL(\ell,\rational)$, 
 the group 
 $\SO_\ell(A; \real) = \{ g \in \SL(\ell,\real) \mid g A g^T = A \}$ % @@@ should have \, 
 is defined over~$\rational$: let 
 $$ \mathcal{Q} =
 \bigset{
 \sum_{1\le p,q \le m+n}
 x_{i,p} A_{p,q} x_{j,q}
 - A_{i,j}
 }{1 \le i,j \le m+n}
 .$$
 In particular, $\SO(m,n)$, under its usual embedding in
$\SL(m+n,\real)$, is defined over~$\rational$.
 \item $\SL(n,\complex)$, under its usual embedding in
$\SL(2n,\real)$, is defined over~$\rational$
\fullccf{EgZarClosed}{SLnC}.
 \end{enumerate}
 \end{egs}

\begin{rems} \label{DefQRems}\ 
\noprelistbreak
 \begin{enumerate}
 \item
 There is always a subset~$\mathcal{Q}$ of $\real[x_{1,1},
\ldots, x_{\ell,\ell}]$, such that $G$ is commensurable to
$\Var(\mathcal{Q})$ \csee{GisAlgic}; that is, $G$~is
\defit[defined!over R@over~$\real$]{defined
over\/~$\real$}. However, it may not be possible to find a
set~$\mathcal{Q}$ that consists entirely of polynomials
whose coefficients are rational, so $G$ may not be defined
over~$\rational$.
 \item If $G$ is defined over~$\rational$, then the
set~$\mathcal{Q}$ of \cref{DefdQDefn} can be chosen to
be finite (because the ring $\rational[x_{1,1}, \ldots,
x_{\ell,\ell}]$ is \term[Noetherian ring]{Noetherian}). 
 \end{enumerate}
 \end{rems}

\begin{prop} \label{hasQform}
 $G$ is isogenous to a group that is defined over\/~$\rational$.
 \end{prop}

\begin{proof}
 It is easy to handle direct products, so the crucial case
is when $G$ is simple. This is easy if $G$ is classical.
Indeed, the groups in \cref{classical-fulllinear,classical-orthogonal} are defined over~$\rational$
(after identifying $\SL(\ell,\complex)$ and
$\SL(\ell,\quaternion)$ with appropriate subgroups of
$\SL(2\ell,\real)$ and $\SL(4\ell,\real)$, in a natural
way). 

The general case is not difficult for someone familiar with
exceptional groups. Namely, since $\Ad G$ is a
finite-index subgroup of $\Aut(\Lie G)$,
it suffices to find a basis of~$\Lie G$, for which the
structure constants of the Lie algebra are rational. We omit the details.
 \end{proof}

\begin{notation} 
For each subring~$\ints$
of~$\real$ (containing~$1$), we construct
\nindex{$G_A$ = $G \cap \SL(n,A)$}$G_{\ints} = G \cap
\SL(n,\ints)$. That is, $G_\ints$ is the subgroup consisting of
the elements of~$G$ whose matrix entries all belong to~$\ints$. 
%(The same notation can be applied when $G \hookrightarrow \SL(\ell,\complex)$, and $A$~is any subring of~$\complex$.)
 \end{notation}

%\begin{notation} \nindex{$G_A$ = $G \cap \SL(n,A)$}
%If we have an embedding of~$G$ in some $\SL(n,\complex)$, and
% $A$~is any subring
%of~$\complex$ (containing~$1$), we let $G_{A} = G \cap
%\SL(n,A)$. That is, $G_A$ is the subgroup consisting of
%the elements of~$G$ whose matrix entries all belong to~$A$. 
% \end{notation}

\begin{eg} \label{SLnCQ}
 Let $\phi \colon \SL(n,\complex) \to \SL(2n,\real)$ be the
natural embedding. Then
 $$ \phi \bigl( \SL(n,\complex) \bigr)_{\rational} 
 = \phi \bigl( \SL(n,\rational[i]) \bigr) .$$
 Therefore, if we think of $\SL(n,\complex)$ as a Lie group
over~$\real$, then $\SL(n,\rational[i])$ represents the
``$\rational$-points'' of $\SL(n,\complex)$.
 \end{eg}

The following result provides an alternate point of view on
being defined over~$\rational$.
It is the nonabelian version of  ($\ref{RnDefdQ<>Latt-dense} \Leftrightarrow \ref{RnDefdQ<>Latt-eqs}$) of \cref{RnDefdQ<>Latt}.

\begin{prop} \label{QptsDense}
 Let $H$ be a connected subgroup of\/ $\SL(\ell,\real)$ that
is almost Zariski closed. The group $H$ is defined
over\/~$\rational$ if and only if $H_{\rational}$ is dense
in~$H$.
 \end{prop}

\begin{proof}
($\Leftarrow$) Let $\mathcal{Q}_\complex = 
	\{\, Q \in \complex[x_{1,1},\ldots,x_{\ell,\ell}] \mid Q(h) = 0, \ \forall h \in H \,\}$. Also, for $d \in \natural$, let $\mathcal{Q}_\complex^d = \{\, Q \in \mathcal{Q}_\complex \mid \deg Q \le d \,\}$. Since $H_\rational$ is dense in~$H$ (and polynomials are continuous), it is clear that $\mathcal{Q}_\complex^d$ is invariant under the Galois group $\Gal(\complex/\rational)$, so it is not difficult to see that $\mathcal{Q}_\complex^d$ is spanned (as a vector space over~$\complex$) by a collection $\mathcal{Q}^d$ of polynomials with rational coefficients \csee{InvtIffDefdQ}. Since $H$ is almost Zariski closed, and polynomial rings are Noetherian, we have $H^\circ = \Var(\mathcal{Q}^d)\!^\circ$ for $d$ sufficiently large. The polynomials in $\mathcal{Q}^d$ all have rational coefficients, so this implies that $H$ is defined over~$\rational$.

($\Rightarrow$)  See \cref{QDenseGSimple} for a proof when $G$ is simple and $G/G_\integer$ is not compact. The general case utilizes a fact from the theory of algebraic groups that will not be proved in this book \csee{GuniratSoQDense}. 
\end{proof}

\begin{warn} \label{HQnotdense}
 \Cref{QptsDense} requires the assumption
that $H$ is connected; there are subgroups~$H$ of
$\SL(\ell,\real)$, such that $H$ is defined
over~$\rational$, but $H_{\rational}$ is not dense in~$H$.
For example, let 
 $$ H = \{\, h \in \SO(2) \mid h^8 = \Id \,\} .$$
 \end{warn}
 
 \begin{rem} \label{JacobsonMorosovOverF}
The \thmindex{Jacobson-Morosov}{Jacobson-Morosov Lemma}~\pref{JacobsonMorosov}
has a relative version: if $G$~is defined over~$\rational$, and $u$ is a nontrivial, unipotent element of $G_\rational$, then there is a (polynomial)
homomorphism $\phi \colon \SL(2,\real) \to G$, such that
 $\phi 
\left( \left[ \begin{smallmatrix}
 1& 1 \\
 0 & 1
 \end{smallmatrix}
 \right] \right)
 = u$
 and
%$\phi$ is defined over~$F$
$\phi \bigl( \SL(2,\rational) \bigr) \subseteq G_\rational$.
\end{rem}

We now state a theorem of fundamental importance in the theory of lattices and
arithmetic groups. It is a nonabelian analogue of the obvious fact that $\integer^\ell$ is a lattice in~$\real^\ell$, and of ($\ref{RnDefdQ<>Latt-eqs} \Rightarrow \ref{RnDefdQ<>Latt-latt}$) of \cref{RnDefdQ<>Latt}.

\begin{major} \label{arith->latt}
	\thmindex{arithmetic subgroups are lattices}%
 If $G$ is defined over~$\rational$, then $G_{\integer}$ is
a lattice in~$G$.
 \end{major}

\begin{proof}
The statement of this theorem is more important than its proof, so, for
most purposes, the reader could accept this fact as an axiom, without 
learning the proof.\label{arith->lattNotPf}
For those who do not want to take this on faith, a discussion of two different proofs 
can be found in \cref{SLnZLattChap} (with some additional details in \cref{ReductionChap}).
\end{proof}

\begin{eg} \label{ArithLattEg}
 Here are some standard cases of
\cref{arith->latt}.%
\noprelistbreak
 \begin{enumerate}
 \item $\SL(2,\integer)$ is a lattice in $\SL(2,\real)$. 
 (We proved this in \cref{SL2Zlatt}.)
 \item \label{ArithLattEg-SLnZ}
$\SL(n,\integer)$ is a lattice in $\SL(n,\real)$.
  (We will prove this in \cref{SLnZLattChap}.) 

 \item $\SO(m,n)_{\integer}$ is a lattice in $\SO(m,n)$.
 \item $\SL(n, \integer[i])$ is a lattice in
$\SL(n,\complex)$ \ccf{SLnCQ}.
 \end{enumerate}
 \end{eg}

\begin{eg}
 As an additional example, let 
 	$$G = \SO(7 x_1^2 - x_2^2 - x_3^2; \real) \iso \SO(1,2). $$
Then \cref{arith->latt}
implies that $G_{\integer}$ is a lattice in~$G$. This
illustrates that the \lcnamecref{arith->latt} is a highly nontrivial result.
For example, in this case, it may not even be obvious to the
reader that $G_{\integer}$ is infinite. 
 \end{eg}

\begin{warn}
 \Cref{arith->latt} requires our standing assumption
that $G$ is semisimple; there are subgroups~$H$ of
$\SL(\ell,\real)$, such that $H$ is defined
over~$\rational$, but $H_{\integer}$ is not a lattice in~$H$.
For example, if $H$ is the group of diagonal matrices in
$\SL(2,\real)$, then $H_{\integer}$ is finite, not a lattice in~$H$.
 \end{warn}

 \begin{rem} The converse of \cref{arith->latt} holds
when $G$ has no compact factors \csee{arithlatt->defdQ}.
 \end{rem}

Combining \cref{hasQform} with \cref{arith->latt}
yields the following important conclusion:

\begin{cor} \label{GHasLatt}
 $G$ has a lattice.
 \end{cor}

 In fact, a more careful look at the proof shows that if $G$ is not compact, then the lattice we constructed is not cocompact:

\begin{cor} \label{GHasNoncpctLatt}
If $G$ is not compact, then $G$ has a noncocompact lattice.
\end{cor}

\begin{proof}
Assume that $G$ is classical, which means it is one of the groups listed in \cref{classical-fulllinear,classical-orthogonal}. As was mentioned in the proof of \cref{hasQform}, each of these groups has an obvious $\rational$-form $G_\rational$, obtained by replacing $\real$ with~$\rational$ (or replacing $\complex$ with $\rational[i]$), in a natural way.
Whenever $G$ is noncompact, it is not difficult to see that $G_\rational$ has a nontrivial unipotent element \csee{UnipInClassicalQ}, so \cref{GammaUnip->notcpct} tells us that $G/G_\integer$ is not compact. 
\end{proof}

\begin{rem}
We will show in \cref{GHasCpctLatt} that $G$ also has a cocompact lattice, and a special case that illustrates the main idea of the proof will be seen much earlier, in \cref{SO(12;Z[sqrt2])}.
\end{rem}

A lattice of the form~$G_{\integer}$ is said to be
\emph{arithmetic}.  However, for the following reasons, a
somewhat more general class of lattices is also said to be
arithmetic. The idea is that there are some obvious
modifications of~$G_{\integer}$ that are also lattices, and
any subgroup that is obviously a lattice should be called
arithmetic.
 \begin{itemize}
 \item If $\phi \colon G_1 \to G_2$ is an isomorphism, and
$\Gamma_1$ is an arithmetic subgroup of~$G_1$, then we wish
to be able to say that $\phi(\Gamma_1)$ is an arithmetic
subgroup in~$G_2$.
 \item We wish to ignore compact
groups; that is, modding out a compact subgroup should not
affect arithmeticity. So we wish to be able to say that if
$K$ is a compact normal subgroup of~$G$, and $\Gamma$ is a
lattice in~$G$, then $\Gamma$ is arithmetic if and only if
$\Gamma K/K$ is an arithmetic subgroup of~$G/K$
 \item Arithmeticity should be independent of
commensurability. 
 \end{itemize}
 The following formal definition implements these
considerations.

\begin{defn} \label{ArithDefn}
 $\Gamma$ is an \defit[arithmetic!subgroup]{arithmetic} subgroup of~$G$ if and
only if there exist
 \begin{itemize}
 \item a closed, connected, semisimple subgroup~$G'$ of
some $\SL(n,\real)$, such that $G'$ is defined
over~$\rational$,
 \item compact normal subgroups $K$ and~$K'$ of $G^\circ$
and~$G'$, respectively, and
 \item an isomorphism $\phi \colon G^\circ/K \to G'/K'$, 
 \end{itemize}
 such that $\phi(\overline{\Gamma})$ is commensurable to
$\overline{G'_{\integer}}$, where $\overline{\Gamma}$
and~$\overline{G'_{\integer}}$ are the images of $\Gamma \cap G^\circ$
and~$G'_{\integer}$ in $G^\circ/K$ and $G'/K'$, respectively.
 \end{defn}

\begin{rems} \label{ArithDefnRem} \ 
\noprelistbreak
	\begin{enumerate}
	\item \label{ArithDefnRem-noK}
	If $G$ has no compact factors, then it is obvious that the subgroup~$K$ in \cref{ArithDefn} must be finite. 

	\item \Cref{IrredNoncpct->NoROS} will show that if $G/\Gamma$ is not compact (and $\Gamma$ is irreducible), then the annoying compact subgroups are not needed in \cref{ArithDefn}.

	\item On the other hand, if $\Gamma$ is cocompact, then a nontrivial
(connected) compact group~$K'$ may be required (even if $G$
has no compact factors). We will see many examples of this
phenomenon, starting with \cref{SO(12;Z[sqrt2])}.

	\item  Up to conjugacy, there are only countably many arithmetic
lattices in~$G$, because there are only countably many
finite subsets of the polynomial ring $\rational[x_{1,1},\ldots,x_{\ell,\ell}]$.

	\end{enumerate}
 \end{rems}

\begin{terminology}
 Our definition of \defit[arithmetic!subgroup]{arithmetic subgroup} assumes the perspective of Lie theory, where $\Gamma$ is assumed
to be embedded in some Lie group~$G$. The theory of algebraic groups
has a more strict definition, which requires
$\Gamma$ to be commensurable to~$G_{\integer}$: arbitrary
isomorphisms are not allowed, and compact subgroups cannot
be ignored. At the other extreme, abstract group theory 
has a much looser definition, which completely ignores~$G$:
if an abstract group~$\Lambda$ is abstractly commensurable to a
group that is arithmetic in our sense, then $\Lambda$~is considered to be arithmetic.
 \end{terminology}

\begin{exercises}

\item \label{conjnoZ}
 Show that if
 $G$ is connected, 
 $G \subseteq \SL(\ell,\real)$,
 and
 $-\Id \notin G$,
 then there exists $h \in \SL(\ell,\real)$, such that
$(h^{-1} G h) \cap \SL(\ell,\integer)$ is trivial.
 \hint{For each nontrivial $\gamma \in \SL(\ell,\integer)$,
let 
 \ $ X_\gamma = \{\, h \in \SL(\ell,\real) \mid h \gamma
h^{-1} \in G \,\} $. \ 
 Then each $X_\gamma$ is nowhere dense in $\SL(\ell,\real)$
\fullcsee{VarClosed}{nodense}.}

\item \label{InvtIffDefdQ}
Let $W$ be a vector subspace of~$\complex^n$, for some~$n$. Show that $W$ is invariant under $\Gal(\complex/\rational)$ if and only if $W$~is spanned by a set of vectors with rational coordinates. 
\hint{($\Rightarrow$) Choose $w \in W \smallsetminus \{0\}$ with a minimal number of nonzero coordinates, and multiply by a scalar to assume at least one coordinate is a nonzero rational. Since $\sigma(w) - w \in W$ for all $\sigma \in \Gal(\complex/\rational)$, the minimality implies $w \in \rational^n$. Mod out $w$ and induct on the dimension.}

\item \label{GuniratSoQDense} 
It can be shown that $G^\circ$ is \defit{unirational}. This means there exists an open subset~$U$ of some~$\real^n$, and a function $f \colon U \to G^\circ$, such that 
	\begin{itemize}
	\item $f(U)$ contains an open subset of~$G$, 
	and
	\item each matrix entry of $f(x)$ is a rational function of~$x$ (that is, a quotient of two polynomials).
	\end{itemize}
Furthermore, if $G$ is defined over~$\rational$, then $f$ can be chosen to be defined over~$\rational$ (that is, all of the coefficients of~$f$ are in~$\rational$).

Assuming the above, show that $G_\rational$ is dense in~$G$ if $G$~is connected and $G$~is defined over~$\rational$.
\hint{Unirationality implies that $\closure{G_\rational}$ contains an open subset of~$G$.}

\item \label{HQnotdenseEx}
 For $H$ as in \cref{HQnotdense}, show that
$H_{\rational}$ is not dense in~$H$.
\hint{$H$ is finite, and $H_\rational \neq H$.}

\item \label{arithlatt->defdQ}
 Show that if $G \subseteq \SL(\ell,\real)$, $G$ has no
compact factors, and $G_{\integer}$ is a lattice in~$G$,
then $G$ is defined over~$\rational$. 
\hint{See the proof of \cref{QptsDense}($\Leftarrow$). % !!!
Since $G_{\integer}$ is a lattice in~$G$, the Borel Density Theorem \pref{BDT(Zardense)} implies that $\mathcal{Q}_\complex^d$ is invariant under the Galois group.} 

\item \label{GZzardense->latt}
 Show that if 
 \begin{itemize}
 \item $G \subseteq \SL(\ell,\real)$, and
 \item $G_{\integer}$ is Zariski dense in~$G$,
 \end{itemize}
 then $G_{\integer}$ is a lattice in~$G$.
 \hint{It suffices to show that $G$~is defined
over~$\rational$.}

\item \label{IntZarDense->EqualEx}
 Show that if 
 \begin{itemize}
 \item $G$ has no compact factors,
 \item $\Gamma_1$ and~$\Gamma_2$ are arithmetic subgroups
of~$G$, and
 \item $\Gamma_1 \cap \Gamma_2$ is Zariski dense in~$G$,
 \end{itemize}
 then $\Gamma_1$ is commensurable to~$\Gamma_2$.
 \hint{Suppose $\phi_j \colon G \to H_j$ is an isomorphism, such
that $\phi_j(\Gamma_j) = (H_j)_{\integer}$. Define $\phi
\colon G \to H_1 \times H_2$ by $\phi(g) = \bigl( \phi_1(g),
\phi_2(g) \bigr)$. 
 Then $\phi(G)_{\integer}
 = \phi(\Gamma_1 \cap \Gamma_2)$
 is Zariski dense in $\phi(G)$, so $\Gamma_1 \cap \Gamma_2$
is a lattice in~$G$ \csee{GZzardense->latt}.
A similar (but slightly more complicated) argument applies if $\phi_j \colon G \to H_j/K_j$, where $K_j$ is compact.}

\item \label{UnipInClassicalQ}
For each classical simple group~$G$ in \cref{classical-fulllinear,classical-orthogonal}, let $G_\rational$ be the subgroup obtained by replacing $\real$ with~$\rational$, $\complex$ with $\rational[i]$, or $\quaternion$ with $\quaternion_\rational = \rational + \rational i+ \rational j+ \rational k$, as appropriate. Show that if $G$ is not compact, then $G_\rational$ contains a nontrivial unipotent element.
\hint{Show that $G_\rational$ contains a copy of either $\SL(2,\rational)$, $\SO(1,2)_\rational$, or $\SU(1,1)_\rational$ \ccf{SL2RinG}.}

\end{exercises}

\section{Margulis Arithmeticity Theorem}

The following astonishing theorem shows that taking integer
points is usually the only way to make a lattice. (See \cref{MargArithPf}
for a sketch of the proof.)

\begin{thm}[(\thmindex{Margulis!Arithmeticity}{Margulis Arithmeticity Theorem})]
\label{MargulisArith}
 If
 \begin{itemize}
 \item $G$ is not isogenous to\/ $\SO(1,n) \times K$ or\/
$\SU(1,n) \times K$, for any compact group~$K$, and
 \item $\Gamma$ is irreducible,
 \end{itemize}
 then $\Gamma$ is arithmetic.
 \end{thm}

\begin{warn}
 Unfortunately,
 \begin{itemize}
 \item $\SL(2,\real)$ is isogenous to $\SO(1,2)$, and 
 \item $\SL(2,\complex)$ is isogenous to $\SO(1,3)$,
 \end{itemize}
 so the arithmeticity theorem says nothing about the lattices in these two
important groups.
 \end{warn}

\begin{rem} \label{MargArithThmCpctFacts}
 The conclusion of \cref{MargulisArith} can be
strengthened: the subgroup~$K$ of \cref{ArithDefn} can
be taken to be finite. More precisely, if $G$ and~$\Gamma$
are as in \cref{MargulisArith}, and $G$~is noncompact
and has trivial center, then there exist
 \begin{itemize}
 \item a closed, connected, semisimple subgroup~$G'$ of
some $\SL(\ell,\real)$, such that $G'$ is defined
over~$\rational$, and
 \item a surjective (continuous) homomorphism $\phi \colon
G' \to G$, 
 \end{itemize}
 such that
 \begin{enumerate}
 \item $\phi( G'_{\integer} )$ is commensurable
to~$\Gamma$; and
 \item the kernel of~$\phi$ is compact.
 \end{enumerate}
 \end{rem}

\begin{rems} \ 
\noprelistbreak
	\begin{enumerate}
	\item For any~$G$, it is possible to give a reasonably complete
description of the arithmetic subgroups of~$G$ (up to
conjugacy and commensurability).
Some examples are worked out in fair detail in \cref{EgArithGrpsChap}.
More generally, \cref{ArithLattsAreClassical} (or the table on \cpageref{IrredInG})
essentially provides a list of all the
irreducible arithmetic subgroups of almost all of the
classical groups. Thus, for most
groups, the Margulis Arithmeticity Theorem provides a
list of all the lattices in~$G$.

	\item Furthermore, knowing that $\Gamma$ is arithmetic provides a
foothold to use algebraic and number-theoretic techniques to
explore the detailed structure of~$\Gamma$. For example, we saw that it is easy to show $\Gamma$ is torsion free if $\Gamma$ is
arithmetic \csee{torsionfree}. A more important example is that (apparently) the
only known proof that every lattice is finitely presented
\csee{GammaFinPres} relies on the Margulis Arithmeticity Theorem.

	\item It is known that there are nonarithmetic lattices in
$\SO(1,n)$ for every~$n$ \csee{NonarithInSO1n}, but we do
not yet have a theory that describes them all when $n \ge 3$.
 Also, nonarithmetic lattices have been constructed in
$\SU(1,n)$ for $n \in \{1,2,3\}$, but (apparently) it is
still not known whether they exist when $n \ge 4$.
	\end{enumerate}
 \end{rems}
 
 \begin{rem} \label{GQComm}
The subgroup \nindex{$\Comm_G(\Gamma)$ = commensurator of~$\Gamma$ in~$G$}
	$$\Comm_G(\Gamma) = \{\, g \in G \mid \text{$g \Gamma g^{-1}$ is commensurable to~$\Gamma$} \,\} $$
is called the \defit{commensurator} of~$\Gamma$ in~$G$. It is easy to see that if $G$ is defined over~$\rational$, 
then $G_{\rational} \subseteq \Comm_G(G_\integer)$
\ccf{SLQCommSLZ}. 
\begin{enumerate}

\item \label{GQComm-criterion}
This implies that if $\Gamma$ is arithmetic (and $G$ is connected, with no compact factors), then $\Comm_G(G_{\integer})$ is dense in~$G$ \csee{QptsDense}.
Margulis proved a converse. Namely, if $G$ is connected and has no compact factors, then
{\bf\mathversion{bold} $\Gamma$ is arithmetic iff\/ $\Comm_G(\Gamma)$ is dense in~$G$}
\csee{CommCriterion}.
This is known as the \thmindex{Commensurability Criterion for Arithmeticity}{Commensurability Criterion for Arithmeticity}.

\item \label{GQComm-bigger}
In some cases, the commensurator of $G_{\integer}$ is much larger than $G_{\rational}$ \csee{Comm(SL2Z)}. However, it was observed by Borel that this never happens when the ``\term{complexification}'' of~$G^\circ$ has trivial center (and other minor conditions are satisfied) \csee{Comm=GC}.
%(Intuitively, the \defit{complexification}~$G_\complex$ of~$G$ is the
%complex Lie group that is obtained from~$G$ by replacing
%real numbers with complex numbers. 
% For example,
%the complexification of $\SL(n,\real)$ is $\SL(n,\complex)$, and the complexification of $\SO(2,3)$ is (isomorphic to) $\SO(5)$.
(See \cref{RFormsOfCGrps} for an explanation of the complexification.)

\end{enumerate}
 \end{rem}

\begin{exercises}
\noprelistbreak
\item \label{Comm(SL2Z)}
Let $G = \SL(2,\real)$ and $G_\integer = \SL(2,\integer)$.
Show  $\Comm_G(G_{\integer})$ is \emph{not}
commensurable to~$G_{\rational}$.
\hint{The diagonal matrix $\diag \bigl( \sqrt{p} , \, 1/ \! \sqrt{p} \bigr)$
 commensurates $G_\integer$, for all $p \in \integer^+$.}

\item \label{Comm(SL3Z)}
Show $\Comm_{\SL(3,\real)} \bigl( \SL(3,\integer) \bigr)$
is \emph{not} commensurable to $\SL(3,\rational)$. 
% \hint{See \cref{Comm(SL2Z)}.}

\item \label{Comm(Gamma)discrete}
 Show that if $G$ is simple and $\Gamma$ is not arithmetic,
then $\Gamma$, $\nzer_G(\Gamma)$, and $\Comm_G(\Gamma)$ are
commensurable to each other.

\item \label{Comm=GC}
(\emph{requires some knowledge of algebraic groups})
Assume $G$ is connected and $G_{\integer}$ is Zariski dense in~$G$ \ccf{BDT(Zardense)}. 
The complexification $G \otimes \complex$ is defined in \cref{GxCNot}.

Show that if $Z(G \otimes \complex) = \{e\}$, then $\Comm_G(G_{\integer}) = G_{\rational}$.
\hint{For $g \in \Comm_G(G_{\integer})$, we know that $\Ad g$ is an automorphism of the Lie algebra~$\Lie G$ that is defined over~$\rational$, so $\Ad g \in (\Ad G)_{\rational}$. However, the assumptions imply that the adjoint representation is an isomorphism (and it is defined over~$\rational$).}

\item
Show that the assumption $Z(G \otimes \complex) = \{e\}$ cannot be replaced with the weaker assumption $Z(G) = \{e\}$ in \cref{Comm=GC}.
\hint{%Let $G = \SL(3,\real)$. 
Any matrix in $\GL(3,\rational)$ has a scalar multiple that is in $\SL(3,\real)$, but $\SL(3,\rational)$ has infinite index in $\GL(3,\rational)$.}

\end{exercises}

\section{Unipotent elements of noncocompact lattices} \label{GodementSect}

The following result answers one of the most basic topological questions about the manifold $G/G_\integer\mk$: is it compact?

\begin{prop}[(\thmindex{Godement Compactness Criterion}{Godement Compactness Criterion})]
\label{GodementCriterion}
 Assume that $G$ is defined over\/~$\rational$. The homogeneous
space $G/G_{\integer}$ is compact if and only if\/
$G_{\integer}$~has no nontrivial unipotent elements.
 \end{prop}

\begin{proof}
 ($\Rightarrow$) This is the easy direction
\csee{GammaUnip->notcpct}.

($\Leftarrow$) We prove the contrapositive: suppose
$G/G_{\integer}$ is not compact. (We wish to show
that $G_{\integer}$ has a nontrivial unipotent element.)
 From \cref{diverge<>contract} (and the fact that
$G_{\integer}$ is a lattice in~$G$ \csee{arith->latt}), we
know that there exist nontrivial $\gamma \in G_{\integer}$
and $g \in G$, such that ${}^g \gamma \approx \Id$. Because the
characteristic polynomial of a matrix is a continuous
function of the matrix entries of the matrix, we conclude
that the characteristic polynomial of ${}^g \gamma$ is
approximately $(x-1)^\ell$ (the characteristic polynomial
of~$\Id$). On the other hand, similar matrices have the same
characteristic polynomial, so this means that the
characteristic polynomial of~$\gamma$ is approximately
$(x-1)^\ell$. Now all the coefficients of the characteristic
polynomial of~$\gamma$ are integers (because $\gamma$ is an
integer matrix), so the only way this polynomial can be
close to $(x-1)^\ell$ is by being exactly equal to
$(x-1)^\ell$. Therefore, the characteristic polynomial
of~$\gamma$ is $(x-1)^\ell$, so $\gamma$~is unipotent.
 \end{proof}

The following important consequence of the Godement Criterion tells us that there is often no need for compact subgroups in \cref{ArithDefn}, the definition of an arithmetic group:

\begin{cor} \label{IrredNoncpct->NoROS}
 Assume 
 	\begin{itemize}
	\item $\Gamma$ is an irreducible, arithmetic subgroup of~$G$, 
	\item $G/\Gamma$ is not compact,
	and
	\item $G$ is connected and has no compact factors.
	\end{itemize}
Then, perhaps after replacing $G$ by an isogenous group, there is an
embedding of~$G$ in some\/ $\SL(\ell,\real)$, such that 
 	\begin{enumerate}
 	\item $G$ is defined over\/~$\rational$, and
	 \item $\Gamma$ is commensurable to $G_{\integer}$.
 	\end{enumerate}
\end{cor}

\begin{proof}
 From \cref{ArithDefn} (and \fullcref{ArithDefnRem}{noK}) we know that (up to isogeny and
commensurability) there is a compact group~$K'$, such that we
may embed $G' = G \times K'$ in some $\SL(\ell,\real)$, such
that $G'$~is defined over~$\rational$, and $\Gamma K' =
G'_{\integer} K'$.

Let $N$ be the almost-Zariski closure of the subgroup
of~$G'$ generated by all of the unipotent elements of
$G'_{\integer}$. 
Since $G/\Gamma$ is not compact, the \lcnamecref{GodementCriterion} % @@@
implies $N$ is infinite.
However, $K'$ has no unipotent elements \fullcsee{SSeltRem}{cpct}, so $N \subseteq G$. 
Also, the definition of $N$ implies that it is normalized by the Zariski closure of~$G'_{\integer}$.
Therefore, the Borel Density Theorem \pref{BDTNotProj} implies that $N$
is a normal subgroup of~$G$. 

Assume, for simplicity, that $G$ is simple \csee{NoROSPf-NotSimpleEx}. Then the conclusion of the preceding paragraph tells us that $N = G$.  
Therefore, $G$ is the almost-Zariski closure of a subset
of~$G'_{\integer}$, which
implies that $G$ is defined over~$\rational$
\ccf{arithlatt->defdQ}. Hence, $G_{\integer}$ is a lattice
in~$G$, and it is easy to see that it is commensurable to~$\Gamma$ \csee{NoROSPf-GZ=GammaEx}.
 \end{proof}

In the special case where $\Gamma$ is arithmetic, the following
result is an easy consequence of \cref{GodementCriterion}, but we will not prove the general case (which is more difficult). The
assumption that $G$~has no compact factors cannot be
eliminated \csee{NonCocpctNoUnip}.

\begin{thm} \label{GodementNoCpctFactor}
 Assume $G$~has no
compact factors. The homogeneous space $G/\Gamma$
is compact if and only if\/ $\Gamma$~has no nontrivial
unipotent elements.
 \end{thm}

The above proof of \cref{GodementCriterion} relies on
the fact that $G_{\integer}$ is a lattice in~$G$, which
will not be proved until \cref{SLnZLattChap}. The following result illustrates that
the cocompactness of~$G_\integer$ can sometimes be proved quite
easily from the Mahler Compactness Criterion
\pref{MahlerCpct}, without assuming that it is a lattice.

\begin{prop} \label{anis->cocpct}
 If 
 \noprelistbreak
 	\begin{itemize}
	\item $B(x,y)$ is a symmetric, bilinear form on\/~$\rational^\ell$, 
	such that
	\item  $B(x,x) \neq 0$ for all nonzero $x \in \rational^\ell$,
	\end{itemize}
then\/ $\SO(B)_{\integer}$
is cocompact in $\SO(B)_{\real}$.
 \end{prop}

\begin{proof}
 Let $G = \SO(B)$ and $\Gamma = \SO(B)_{\integer} = G_\integer$.
(Our proof will not use the fact that $\Gamma$~is a
lattice in~$G$.)
Replacing $B$ by an integer multiple to clear the
denominators, we may assume $B(\integer^\ell,\integer^\ell)
\subseteq \integer$.

\goodbreak % @@@

\setcounter{step}{0}

\begin{step} \label{anis->cocpct-precpct}
 The image of~$G$ in\/ $\SL(\ell,\real)
/ \SL(\ell,\integer)$ is precompact.
 \end{step}
 Let 
 \begin{itemize}
 \item $\{g_n\}$ be a sequence of elements of~$G$ and
 \item $\{v_n\}$ be a sequence of elements of $\integer^\ell
\smallsetminus \{0\}$. 
 \end{itemize}
 Suppose that $g_n v_n \to 0$. (This will lead to a
contradiction, so the desired conclusion follows from the
Mahler Compactness Criterion \pref{MahlerCpct}.)

Since $B(v,v) \neq 0$ for all
nonzero $v \in \integer^\ell$, and $B(\integer^\ell,\integer^\ell)
\subseteq \integer$, we have $|B(v_n,v_n)| \ge 1$
for all~$n$. Therefore
 $$ 1 \le |B(v_n,v_n)| = |B(g_n v_n, g_n v_n)| \to |B(0,0)|
= 0 .$$
 This is a contradiction.

\begin{step} \label{anis->cocpct-closed}
 The image of~$G$ in\/ $\SL(\ell,\real) /
\SL(\ell,\integer)$ is closed.
 \end{step}
 Suppose 
 	$$ \text{$g_n \gamma_n \to h \in \SL(\ell,\real)$, \ with $g_n
\in G$ and $\gamma_n \in \SL(\ell,\integer)$} . $$
We wish to show $h \in G \, \SL(\ell,\integer)$.

Let $\{\varepsilon_1,\cdots,\varepsilon_\ell\}$ be the standard basis
of~$\real^\ell$ (so each $\varepsilon_j \in \integer^\ell$).  Then
 $$ B(\gamma_n \varepsilon_j, \gamma_n \varepsilon_k)
 \in B(\integer^\ell,\integer^\ell)
 \subseteq \integer .$$
 We also have
 $$ B(\gamma_n \varepsilon_j, \gamma_n \varepsilon_k)
 =  B(g_n \gamma_n \varepsilon_j, g_n \gamma_n \varepsilon_k)
 \to B(h \varepsilon_j, h \varepsilon_k) .$$
 Since $\integer$ is discrete, we conclude that
 $ B(\gamma_n \varepsilon_j, \gamma_n \varepsilon_k)
 = B(h \varepsilon_j, h \varepsilon_k) $
 for any sufficiently large~$n$.
 Therefore $h \gamma_n^{-1} \in \SO(B)$ \csee{hgammainG},
 so we have $h \in G \gamma_n \subseteq G \, \SL(\ell,\integer)$.

\begin{step}
 Completion of the proof.
 \end{step}
 Define $\phi \colon G/\Gamma \to \SL(\ell,\real) /
\SL(\ell,\integer)$ by $\phi(g\Gamma) = g
\SL(\ell,\integer)$. By combining
\cref{anis->cocpct-precpct,anis->cocpct-closed}, we see that the image
of~$\phi$ is compact.  Therefore, it suffices to show that $\phi$
is a homeomorphism onto its image.

Given a sequence $\{g_n\}$ in~$G$, such that $\{\phi(g_n
\Gamma)\}$ converges, we wish to show that $\{g_n \Gamma\}$
converges. There is a sequence $\{\gamma_n\}$ in
$\SL(\ell,\integer)$, and some $h \in G$, such that $g_n
\gamma_n \to h$. The proof of \cref{anis->cocpct-closed}
shows, for all large~$n$, that $h \in G \gamma_n$. Then
$\gamma_n \in Gh = G$ (and we know $\gamma_n \in \SL(\ell,\integer)$), so $\gamma_n \in G_{\integer} =
\Gamma$. Therefore, $\{g_n \Gamma\}$ converges (to
$h\Gamma$), as desired.
 \end{proof}

\begin{exercises}

\item \label{NoROSPf-NotSimpleEx}
 The proof of \cref{IrredNoncpct->NoROS} assumes that $G$ is simple. Eliminate this hypothesis.
 \hint{The proof shows that $N \cap \Gamma$ is a lattice in~$N$. Since $\Gamma$ is irreducible, this implies $N = G$.}
 
\item \label{NoROSPf-GZ=GammaEx}
 At the end of the proof of \cref{IrredNoncpct->NoROS}, show that $G_{\integer}$ is commensurable to~$\Gamma$.
 \hint{We know $G'_{\integer} K' = \Gamma K'$, and $G_{\integer}$ has finite index in $G'_{\integer}$ \csee{finext->latt}. Mod out~$K'$.}

\item \label{NonCocpctNoUnip}
 Show there is a noncocompact lattice~$\Gamma$ in
$\SL(2,\real) \times \SO(3)$, such that no nontrivial
element of~$\Gamma$ is unipotent.
\noprelistbreak % @@@
\hint{$\SL(2,\real)$ has a lattice~$\Gamma'$ that is free. Let $\Gamma$ be the graph of a  homomorphism from $\Gamma'$ to $\SO(3)$.}

\item Suppose $G \subseteq \SL(\ell,\real)$ is defined
over~$\rational$.
 \begin{enumerate}
 \item Show that if $N$ is a closed, normal subgroup of~$G$,
and $N$~is defined over~$\rational$, then $G_{\integer} N$
is closed in~$G$.
 \item Show that $G_{\integer}$ is irreducible if and only
if no proper, closed, connected, normal subgroup of~$G$ is
defined over~$\rational$. (That is, if and only if $G$~is
\defit[Q-@$\rational$-!simple group]{$\rational$-simple}.)
 \item Let $H$ be the Zariski closure of the subgroup
generated by the unipotent elements of~$G_{\integer}$. Show
that $H$~is defined over~$\rational$.
% \item Provide a direct proof of \cref{noncocpct->nocpct}, without using
%\cref{GZirred->scalars}, or any other results on
%Restriction of Scalars.
 \end{enumerate}

\item \label{Cocpct->AllSS}
Show that if every element of~$\Gamma$ is semisimple, then $G/\Gamma$ is compact. 
\hint{There is no harm in assuming that $G$ has no compact factors (why?), so \cref{GodementNoCpctFactor} applies.}

\item \label{CocpctAllSS}
(\emph{assumes some familiarity with reductive groups})
Prove the converse of \cref{Cocpct->AllSS}.
\hint{Let $kau$ be the real Jordan decomposition of an element~$g$ of~$\Gamma$. Since $C_G(ka)$ is reductive \csee{C(T)reductive}, the Jacobson-Morosov Lemma provides a subgroup~$L$ of $C_G(ka)$ that contains~$u$ and is isogenous to $\SL(2,\real)$. So $ka$ is in the closure of ${}^G \!g$. However, ${}^G \! g$ is closed, since $\Gamma$ is discrete and cocompact. Therefore $g = ka$ is semisimple.}

 \item \label{GZCocpctIff}
Assuming $\Gamma = G_\integer$ is arithmetic (and $G$ is defined over~$\rational$), prove the following are equivalent:
 	\begin{enumerate}

	\item \label{GZCocpctIff-cpct}
	$G/G_\integer$ is compact.

	 \item \label{GZCocpctIff-unip}
	 $G_{\rational}$ has no nontrivial unipotent elements.

	 \item \label{GZCocpctIff-GQss}
	 Every element of~$G_\rational$ is semisimple.

	 \item \label{GZCocpctIff-ss}
	 Every element of\/~$\Gamma$ is semisimple.

	 \item \label{GZCocpctIff-SL2Q}
	 $G_\rational$ does not contain a subgroup isogenous to $\SL(2,\rational)$. 
	 (More precisely, there does not exist a continuous homomorphism $\rho \colon \SL(2,\real) \to G$, such that $\rho \bigl( \SL(2,\rational) \bigr) \subseteq G_\rational$.)
	
	 \end{enumerate}
 \hint{($\ref{GZCocpctIff-unip} \Rightarrow \ref{GZCocpctIff-GQss}$)~Jordan decomposition.
 ($\ref{GZCocpctIff-SL2Q} \Rightarrow \ref{GZCocpctIff-unip}$)~Jacobson-Morosov Lemma \pref{JacobsonMorosovOverF}.}

\item \label{QDenseGSimple}
Show that $G_\rational$~is dense in~$G$ if $G$ is defined over~$\rational$, $G$~is simple, and $G/G_\integer$ is not compact.
\hint{The Godement Criterion implies that $G_\rational$ has a nontrivial unipotent element~$u$. Write $u = \exp T = \sum_{k=0}^\ell T^k/k!$ (where $T \in \Mat_{\ell\times\ell}(\rational)$ and $T^{\ell+1} = 0$). Then $\exp(rT) \in G_\rational$ for all $r \in \rational$, so the identity component of $\closure{G_\rational}$ is nontrivial. Combining \cref{arith->latt} with the Borel Density Theorem \pref{BDT-normalize} implies that $\closure{G_\rational} = G$.}
 
 \item \label{hgammainG}
 Let $B(x,y)$ be a symmetric, bilinear form on~$\real^\ell$, 
 let $\{v_1,\cdots,v_\ell\}$ be a  basis of~$\real^\ell$,
and let $\gamma,h \in \SL(\ell,\real)$.
 If $B(\gamma v_j, \gamma v_k) = B(h v_j, h v_k)$ for all $j$ and~$k$, 
 show that $h \gamma^{-1} \in \SO(B)$.
\hint{$\{\gamma v_1,\ldots,\gamma v_\ell\}$ is a basis of~$\real^\ell$.}

\end{exercises}

\section{How to make an arithmetic subgroup}
\label{MakeArithLattSect}

The definition that (modulo commensurability, isogenies, and
compact factors) an arithmetic subgroup must be the
$\integer$-points of~$G$ has the virtue of being concrete.
However, this concreteness imposes a certain lack of
flexibility. (Essentially, we have limited ourselves to the
standard basis of the vector space~$\real^n$, ignoring the
possibility that some other basis might be more convenient
in some situations.) We now describe a more abstract
viewpoint that makes the construction of general arithmetic
lattices more transparent. (In particular, this approach
will be used in \S\ref{RestrictScalarsSect}.) The key point is that there are analogues of $\integer^\ell$ and~$\rational^\ell$ in any real
vector space, not just~$\real^\ell$ \fullcsee{VQ=Qd}{VQ}.

\begin{defns} \label{DefdQAbstract}
 Let $V$ be a real vector space.
 \begin{enumerate}
 \item
 A $\rational$-subspace $V_{\rational}$ of~$V$ is a
\defit[Q-@$\rational$-!form]{$\rational$-form} of~$V$ if the natural $\real$-linear
map $V_{\rational} \otimes_{\rational} \real \to V$ is an isomorphism
\csee{Qform<>basis}. (The map is defined
by $v \otimes t \mapsto tv$.)
 \item 
 A polynomial~$f$ on~$V$ is \defit[defined!over Q@over~$\rational$]{defined
over $\rational$} (with respect to the $\rational$-form
$V_{\rational}$) if $f(V_{\rational}) \subseteq \rational$
\csee{Qdefd<>Qcoeffs}.
 \item
 A subgroup~\nindex{$\Zlatt$ = $\integer$-lattice in a vector space}$\Zlatt$ of the additive
group of~$V_{\rational}$ is a \defit[Z-lattice@$\integer$-lattice]{$\integer$-lattice}
in~$V_{\rational}$ if it is finitely generated and the natural $\rational$-linear map $\Zlatt \otimes_{\integer} \rational  \to V_{\rational}$ is an isomorphism
\csee{vsLatt<>basis}.  (The map is defined
by $v \otimes t \mapsto tv$.)

 \item Each $\rational$-form~$V_{\rational}$ of~$V$
yields a corresponding $\rational$-form of the real
vector space $\End(V)$ by
 $ \End(V)_{\rational} = \{\, A \in \End(V) \mid
A(V_{\rational}) \subseteq V_{\rational} \,\} $
 \csee{VQ->End(V)Q}.
 \item A function~$Q$ on a real vector space~$W$ is a
\defit[polynomial!on a vector space]{polynomial} if for
some (hence, every) $\real$-linear isomorphism $\phi \colon
\real^\ell \iso W$, the composition $f \circ \phi$ is a
polynomial function on~$\real^\ell$.
 \item A subgroup~$H$ of $\SL(V)$ is
\defit[defined!over Q@over~$\rational$]{defined over~$\rational$} (with
respect to the $\rational$-form $V_{\rational}$) if there
exists a set~$\mathcal{Q}$ of polynomials on $\End(V)$, such
that 
 \begin{itemize}
 \item every $Q \in \mathcal{Q}$ is defined
over~$\rational$ (with respect to the $\rational$-form
$V_{\rational}$),
 \item 
 $\Var(\mathcal{Q}) = \{\, g \in \SL(V) \mid
\mbox{$Q(g) = 0$ for all $Q \in \mathcal{Q}$} \,\}$
 is a subgroup of $\SL(V)$, and
 \item $\Var(\mathcal{Q})^\circ$ is a finite-index subgroup of~$H$.
 \end{itemize}
 \end{enumerate}
 \end{defns}

\begin{rems} \ 
\noprelistbreak
 \begin{enumerate}
 \item Suppose $G \subseteq \SL(\ell,\real)$, as usual. For the standard
$\rational$-form $\rational^\ell$ of~$\real^\ell$, it is easy to see
that $G$ is defined over~$\rational$ in terms of
\cref{DefdQAbstract} if and only if it is defined
over~$\rational$ in terms of \cref{DefdQDefn}. 
 \item Some authors simply call $\Zlatt$ a 
 \defit[lattice!in~$V_{\rational}$|indsee{$\integer$-lattice}]{lattice in~$V_{\rational}$}, 
 but this could cause confusion, because
$\Zlatt$ is \emph{not} a lattice in~$V_{\rational}$, in the
sense of \cref{LatticeDefn} (although it \emph{is} a
lattice in~$V$).
 \end{enumerate}
 \end{rems}

A $\rational$-form~$V_{\rational}$
and $\integer$-lattice~$\Zlatt$ simply represent
$\rational^\ell$ and~$\integer^\ell$, under some
identification of~$V$ with~$\real^\ell$:

\begin{lem} \label{VQ=Qd}
 Let\/ $V$ be an $\ell$-dimensional real vector space.
 \begin{enumerate}
 \item \label{VQ=Qd-VQ}
 If\/ $V_{\rational}$ is a\/ $\rational$-form of\/~$V$,
then there exists an\/ $\real$-linear isomorphism $\phi \colon V
\to \real^\ell$, such that $\phi(V_{\rational}) =
\rational^\ell$. Furthermore, if $\Zlatt$ is any
$\integer$-lattice in\/~$V_{\rational}$, then $\phi$ may be
chosen so that $\phi(\Zlatt) = \integer^\ell$.
 \item A polynomial~$f$ on\/~$\real^\ell$ is defined
over\/~$\rational$ \textup(with respect
to the standard\/ $\rational$-form~$\rational^\ell$\textup) if and
only if every coefficient of~$f$ is in\/~$\rational$
\csee{Qdefd<>Qcoeffs}.
 \end{enumerate}
 \end{lem}

Also note that any two $\integer$-lattices
in~$V_{\rational}$ are commensurable:

\begin{lem}[\csee{VQ->vslatt}] \label{pLambda1inLambda2}
 If $\Zlatt_1$ and~$\Zlatt_2$ are two $\integer$-lattices
in $V_{\rational}$, then there is some nonzero $p \in
\integer$, such that $p \Zlatt_1 \subseteq \Zlatt_2$ and $p
\Zlatt_2 \subseteq \Zlatt_1$.
 \end{lem}

It is now easy to prove the following more abstract
characterization of arithmetic subgroups \csee{VQ->LattInG,GLambdaUnique}.

\begin{prop} \label{AbstractArith}
 Suppose $G \subseteq \GL(V)$, and $G$ is defined
over\/~$\rational$, with respect to the\/
$\rational$-form~$V_{\rational}$.
 \begin{enumerate}
 \item \label{AbstractArith-arith}
 If $\Zlatt$ is any $\integer$-lattice in~$V_{\rational}$,
then
 $$ G\!_{\Zlatt} = \{\, g \in G \mid g \Zlatt = \Zlatt \,\}
$$
 is an arithmetic subgroup of~$G$.
 \item \label{AbstractArith-comm}
 If $\Zlatt_1$ and~$\Zlatt_2$ are
$\integer$-lattices in~$V_{\rational}$, then
$G\!_{\Zlatt_1}$ is commensurable to~$G\!_{\Zlatt_1}$.
 \end{enumerate}
 \end{prop}

From \fullcref{AbstractArith}{comm}, we see that the
arithmetic subgroup~$G\!_{\Zlatt}$ is almost entirely determined
by the $\rational$-form~$V_{\rational}$; choosing a different $\integer$-lattice 
in $V_\rational$ will yield a commensurable arithmetic subgroup.

\begin{exercises}

\item \label{Qform<>basis}
 Show that a $\rational$-subspace $V_{\rational}$ of~$V$ is
a $\rational$-form if an only if there is a subset~$\basis$
of~$V_{\rational}$, such that $\basis$ is both a
$\rational$-basis of~$V_{\rational}$ and an $\real$-basis
of~$V$.

\item \label{Qdefd<>Qcoeffs}
 For the standard $\rational$-form $\rational^\ell$
of~$\real^\ell$, show that a polynomial is defined
over~$\rational$ if and only if all of its coefficients are
rational.

\item \label{vsLatt<>basis}
 Show that a subgroup~$\Zlatt$ of~$V_{\rational}$ is a
$\integer$-lattice in~$V_{\rational}$ if and only if there
is a $\rational$-basis~$\basis$ of~$V_{\rational}$, such
that $\Zlatt$ is the additive abelian subgroup
of~$V_{\rational}$ generated by~$\basis$.

\item \label{VSLatt<>Rank}
 Let $V$ be a real vector space of dimension~$\ell$,
and let $\Zlatt$~be a discrete subgroup of the additive
group of~$V$. Recall that the \defit[rank!of an abelian group]{rank} of an abelian group is the largest~$r$, such that the group contains a copy of~$\integer^r$.
 \begin{enumerate}
 \item Show that $\Zlatt$ is a finitely generated, abelian
group of rank $\le \ell$, with equality if and only if the
$\real$-span of~$\Zlatt$ is~$V$.
 \item Show that if the rank of~$\Zlatt$ is~$\ell$, then
the $\rational$-span of~$\Zlatt$ is a $\rational$-form
of~$V$, and $\Zlatt$~is a $\integer$-lattice
in~$V_{\rational}$.
 \end{enumerate}
 \hint{Induction on~$\ell$. For $\lambda \in \Zlatt$, show
that the image of~$\Zlatt$ in $V/\real \lambda$ is
discrete.}

\item \label{VQ->End(V)Q}
 Verify: if $V_{\rational}$ is a $\rational$-form
of~$V$, then $\End(V)_{\rational}$ is a
$\rational$-form of $\End(V)$.

\item \label{VQ->vslatt}
 Prove \cref{pLambda1inLambda2}. Conclude that
$\Lambda_1$ and~$\Lambda_2$ are commensurable.

\item \label{VQ->LattInG}
 Prove \fullcref{AbstractArith}{arith}. [{\it Hint:}
Use \cref{VQ=Qd}.]

\item \label{GLambdaUnique}
 Prove \fullcref{AbstractArith}{comm}. [{\it Hint:}
Use \cref{pLambda1inLambda2}.]

%\item \label{ChevalleyStabOverQ} % cf. \ref{ChevalleyStabDiscreteZ}
%Assume $G$ is defined over~$\rational$ (and connected).
%Show there exist
%	\begin{itemize}
%	\item a $\rational$-form $V_{\rational}$ of some finite-dimensional real vector space~$V$, 
%	\item a $\integer$-lattice~$\Zlatt$ in~$V_{\rational}$,
%	\item a vector~$v$ in~$\Zlatt$,
%	and
%	\item a homomorphism $\rho \colon \SL(\ell,\real) \to \SL(V)$,
%	\end{itemize}
%such that 
%	\begin{enumerate}
%	\item $G = \Stab_{\SL(\ell,\real)}(v)^\circ$,
%	\item $\rho \bigl( \SL(\ell,\rational) \bigr) v \subseteq V_{\rational}$,
%	and 
%	\item $\rho \bigl( \SL(\ell,\integer) \bigr) v \subseteq \Zlatt$.
%	\end{enumerate}
%\hint{See the hint to \cref{ChevalleyStabEx}.}

\end{exercises}

\section{Restriction of scalars} \label{RestrictScalarsSect}

We know that $\SL(2,\integer)$ is an arithmetic subgroup of
$\SL(2,\real)$. In this section, we explain that $\SL \bigl(
2, \integer[\sqrt{2}] \bigr)$ is an arithmetic subgroup of the group
$\SL(2,\real) \times \SL(2,\real)$ (see
\cref{SL(2Z[sqrt2])}). More generally, recall that any finite extension of~$\rational$
is called an \defit[algebraic!number field]{algebraic number field}. We will see that if
$\ints$ is the ring of algebraic integers in any
algebraic number field~$F$, and $G$~is defined over~$F$,
then $G_{\ints}$ is an arithmetic subgroup of a certain
group~$G'$ that is related to~$G$. 

\begin{rem}
 In practice, we do not require $\ints$ to be the
entire ring of algebraic integers in~$F$: it suffices for the ring~$\ints$ to have finite index in the ring of integers
(as an additive group); equivalently, the $\rational$-span
of~$\ints$ should be all of~$F$, or, in other words,
the ring $\ints$ should be a $\integer$-lattice
in~$F$. (A $\integer$-lattice in~$F$ that is also a
subring is called an \defit[order (in an algebraic number
field)]{order} in~$F$.)
 \end{rem}

Any complex vector space can be thought of as a real vector
space (of twice the dimension). Similarly, any complex Lie
group can be thought of as a real group (of twice the
dimension). Restriction of scalars is the generalization of
this idea to any field extension $F/L$, not just
$\complex/\real$. This yields a general method to construct
arithmetic subgroups. 

\begin{eg}
 Let 
 \begin{itemize}
 \item $F = \rational[\!\sqrt{2}]$,
 \item $\ints = \integer[\!\sqrt{2}]$, and
 \item $\sigma$ be the nontrivial Galois automorphism of~$F$,
 \end{itemize}
 and define a ring homomorphism $\Delta \colon F \to
\real^2$ by 
 $\Delta(x) = \bigl( x, \sigma(x) \bigr) $.

 It is easy to show that $\Delta(\ints)$ is discrete
in~$\real^2$. Namely, for $x \in \ints$, the product
of the coordinates of~$\Delta(x)$ is the product $x \cdot
\sigma(x)$ of all the Galois conjugates of~$x$. This is the
\defit[norm!of an algebraic number]{norm} of the algebraic
number~$x$. Because $x$ is an algebraic integer, its norm is
an ordinary integer; hence, its norm is bounded away
from~$0$. So it is impossible for both coordinates of
$\Delta(x)$ to be small simultaneously.

More generally, if $\ints$ is the ring of integers of
any algebraic number field~$F$, this same argument shows that if we let
$\{ \sigma_1, \ldots, \sigma_r \}$ be the set of all
embeddings of~$\ints$ in~$\complex$, and define
 $\Delta \colon \ints \to \complex^r$
 by 
 $$ \Delta(x) = \bigl( \sigma_1(x), \ldots, \sigma_r(x) \bigr) , $$
 then $\Delta( \ints )$ is a \emph{discrete} subring
of~$\complex^r$. 

Now $\Delta$ induces a homomorphism $\Delta_* \colon \SL(\ell, \ints) \to \SL(\ell, \complex^r)$ (because $\SL(\ell, \,{\cdot}\,)$ is a functor from the category of commutative rings to the category of groups). Furthermore, the group $\SL(\ell, \complex^r)$ is naturally isomorphic to $\SL(\ell, \complex)^r$. Therefore, we have a homomorphism (again called~$\Delta$) from $\SL(\ell, \ints)$ to  $\SL(\ell, \complex)^r$. Namely, for $\gamma \in \SL(\ell, \ints)$, we let $\sigma_i(\gamma) \in \SL(\ell,\complex)$ be obtained by applying $\sigma_i$ to each entry of~$\gamma$, and then
 $$ \Delta(\gamma) = \bigl( \sigma_1(\gamma), \ldots, \sigma_r(\gamma) \bigr) .$$
Since $\Delta(\ints)$ is discrete in~$\complex^r$, it is obvious that the image of~$\Delta_*$ is discrete in $\SL(\ell, \complex^r)$, so $\Delta(\Gamma)$ is a discrete subgroup
of $\SL(\ell,\complex)^r$, for any subgroup~$\Gamma$ of
$\SL(\ell,\ints)$. 
 \end{eg}

The main goal of this \lcnamecref{RestrictScalarsSect} is to show that if $\Gamma = G_{\ints}$, and $G$~is defined
over~$F$, then the discrete group~$\Delta(\Gamma)$ is an arithmetic subgroup of a certain subgroup of
$\SL(\ell,\complex)^r$.

To illustrate, let us show that
$\SL \bigl( 2,\integer[\!\sqrt{2}] \bigr)$ is isomorphic to an
arithmetic subgroup of $\SL(2,\real) \times
\SL(2,\real)$.

\begin{eg} \label{SL(2Z[sqrt2])}
 Let
 \noprelistbreak
 \begin{itemize}
 \item $\Gamma = \SL \bigl( 2,\integer[\!\sqrt{2}] \bigr)$,
 \item $G = \SL(2, \real) \times \SL(2, \real)$, and
 \item $\sigma$ be the conjugation on $\rational[\!\sqrt{2}]$
\textup(so $\sigma \bigl( a + b \sqrt{2} \bigr) = a - b
\sqrt{2}$, for $a,b \in \rational$\textup),
 \end{itemize}
 and define $\Delta \colon \Gamma \to G$ by
$\Delta(\gamma) = \bigl(\gamma, \sigma(\gamma) \bigr)$.

Then $\Delta(\Gamma)$ is an irreducible, arithmetic subgroup
of~$G$.
 \end{eg}

\begin{proof}
 Let $F = \rational[\!\sqrt{2}]$ and $\ints =
\integer[\!\sqrt{2}]$. Then $F$ is a 2-dimensional vector
space over~$\rational$, and $\ints$ is a
$\integer$-lattice in~$F$.

Since $\bigl\{ (1,1) , (\sqrt{2},-\sqrt{2}) \bigr\}$ is both
a $\rational$-basis of~$\Delta(F)$ and an $\real$-basis
of~$\real^2$, we see that $\Delta(F)$ is a $\rational$-form
of~$\real^2$. Therefore,
 $$\Delta(F^2) = \bigl\{\, \bigl( u, \sigma(u) \bigr) \in F^4
\mid u \in F^2 \,\bigr\} $$
 is a $\rational$-form of~$\real^4$, and
$\Delta(\ints^2)$ is a $\integer$-lattice in
$\Delta(F^2)$.

Now $G$ is defined over~$\rational$ \csee{SL2xSL2Defd/Q}, so
$G_{\Delta(\ints^2)}$ is an arithmetic subgroup of~$G$.
It is not difficult to see that $G_{\Delta(\ints^2)} =
\Delta(\Gamma)$ \csee{SL2xSL2-Delta(Gamma)=}. Furthermore,
because $\Delta(\Gamma) \cap \bigl( \SL(2,\real) \times e
\bigr)$ is trivial, we see that the lattice $\Delta(\Gamma)$
must be irreducible in~$G$ \csee{prodirredlatt}.
 \end{proof}

More generally, the proof of \cref{SL(2Z[sqrt2])}
shows that if $G$ is defined over~$\rational$, then
$G_{\integer[\!\sqrt{2}]}$ is isomorphic to an (irreducible)
arithmetic subgroup of $G \times G$. 

Here is another sample application of the method.

\begin{eg} \label{SO(12;Z[sqrt2])}
 Let $G = \SO(x^2 + y^2 - \sqrt{2} z^2 ; \real) \iso
\SO(1,2)$. Then $G_{\integer[\!\sqrt{2}]}$ is a
cocompact, arithmetic subgroup of~$G$.
 \end{eg}

\begin{proof}
 As above, let $\sigma$ be the conjugation on
$\rational[\!\sqrt{2}]$. Let $\Gamma = G_{\integer[\!\sqrt{2}]}$.

Let $K' = \SO(x^2 + y^2 + \sqrt{2} z^2) \iso \SO(3)$, so
$\sigma(\Gamma) \subseteq K'$. (However, $\sigma(\Gamma) \not\subseteq G$.) Then, we may 
	$$ \text{define \ $\Delta
\colon \Gamma \to G \times K'$ \ by \ $\Delta(\gamma) = \bigl( \mkern2mu
 \gamma, \sigma(\gamma) \bigr)$} .$$
Arguing as in the proof of \cref{SL(2Z[sqrt2])} establishes that
$\Delta(\Gamma)$ is an arithmetic subgroup of $G \times K'$.
(See \cref{GxG-Defd/Q} for the technical point of
verifying that $G \times K'$ is defined over~$\rational$.)
Since $K'$ is compact, we see, by modding out~$K'$, that
$\Gamma$ is an arithmetic subgroup of~$G$. (This type of
example is the reason for including the compact normal
subgroup~$K'$ in \cref{ArithDefn}.)

 Let $\gamma$ be any nontrivial element of~$\Gamma$. Since
$\sigma(\gamma) \in K'$, and compact groups have no
nontrivial unipotent elements \fullcsee{SSeltRem}{cpct}, we know that
$\sigma(\gamma)$ is not unipotent. Therefore, $\sigma(\gamma)$
has some eigenvalue $\lambda \neq 1$. Hence, $\gamma$ has
the eigenvalue $\sigma^{-1}(\lambda) \neq 1$, so $\gamma$ is
not unipotent. Therefore, Godement's Criterion
\pref{GodementCriterion} implies that $\Gamma$ is cocompact.
Alternatively, this conclusion can easily be obtained
directly from the Mahler Compactness Criterion
\pref{MahlerCpct} \csee{SO(B)Z2cocpt(Mahler)}.
 \end{proof}

Let us consider one more example before stating the general
result.

\begin{eg}
 Let
 \begin{itemize}
 \item $F = \rational[\!\! \root 4\!\! \of 2]$,
 \item $\ints = \integer[\!\!\root 4\!\!\of 2]$,
 \item $\Gamma = \SL(2,\ints)$, and
 \item $G = \SL(2,\real) \times \SL(2,\real) \times
\SL(2,\complex)$.
 \end{itemize}
 Then $\Gamma$ is isomorphic to an irreducible, arithmetic
subgroup of~$G$.
 \end{eg}

\begin{proof}
 For convenience, let $\alpha = \!\root 4\!\! \of 2$. There are
exactly 4 distinct embeddings $\sigma_0$, $\sigma_1$,
$\sigma_2$, $\sigma_3$ of~$F$ in~$\complex$ (corresponding
to the 4~roots of $x^4 - 2 = 0$); they are determined by:
 $$ \text{$\sigma_0(\alpha) = \alpha$ \  (so $\sigma_0 = \Id$),
 \  $\sigma_1(\alpha) = -\alpha$,
 \  $\sigma_2(\alpha) = i\alpha$, 
 \ and
 \  $\sigma_3(\alpha) = -i\alpha$}
 . $$
 Define $\Delta \colon F \to \real \oplus \real \oplus
\complex$ by $\Delta(x) = \bigl( x, \sigma_1(x), \sigma_2(x)
\bigr)$. Then, arguing much as before, we see that
$\Delta(F^2)$ is a $\rational$-form of~$\real^2 \oplus
\real^2 \oplus \complex^2$, $G$~is defined over~$\rational$,
and $G_{\Delta(\ints^2)} = \Delta(\Gamma)$.
 \end{proof}

These examples illustrate all the ingredients of the general result that will be stated in \cref{ResScal->Latt} after the necessary definitions.

\begin{defn} \label{PlaceDefn}
 Let $F$ be an algebraic number field (or, in other words, let $F$ be a finite extension of~$\rational$).
 \begin{enumerate}
 \item Two distinct embeddings $\sigma_1, \sigma_2 \colon F
\to \complex$ are said to be \defit[equivalent embeddings in
$\complex$]{equivalent} if $\sigma_1(x) =
\overline{\sigma_2(x)}$, for all $x \in F$ (where
$\overline{z}$ denotes the usual complex conjugate of the
complex number~$z$).
 \item A \defit[place!of~$F$]{place} of~$F$ is an equivalence
class of embeddings in~$\complex$. Therefore, each place
consists of either one or two embeddings of~$F$: 
	\begin{itemize}
	\item a \defit[place!real]{real place} consists of only one embedding (with
$\sigma(F) \subset \real$), but
	\item a \defit[place!complex]{complex place} consists
of two embeddings (with $\sigma(F) \not\subset \real$).
	\end{itemize}
 \item We let $S^\infty = \{\, \mbox{places of~$F$} \,\}$,
or, abusing notation, we assume that $S^\infty$ is a set of
embeddings, consisting of exactly one embedding from each
place.
 \item For $\sigma \in S^\infty$, we let
 $$ F_\sigma = 
 \begin{cases}
 \real & \mbox{if $\sigma$ is real}, \\
 \complex & \mbox{if $\sigma$ is complex}.
 \end{cases}
 $$
 Note that $\sigma(F)$ is dense in~$F_\sigma$, so $F_\sigma$
is often called the \defit[completion of~$F$]{completion}
of~$F$ at the place~$\sigma$.
 \item For $\mathcal{Q} \subset
F_\sigma[x_{1,1},\ldots,x_{\ell,\ell}]$, let
 $$ \Var_{F_\sigma}(\mathcal{Q})
 = \{\, g \in \SL(\ell,F_\sigma) \mid Q(g) = 0, \ \forall Q
\in \mathcal{Q} \,\} .$$
 Thus, for $F_\sigma = \real$, we have
$\Var_{\real}(\mathcal{Q}) = \Var(\mathcal{Q})$, and
$\Var_{\complex}(\mathcal{Q})$ is analogous, using the
field~$\complex$ in place of~$\real$.
 \item Suppose $G \subseteq \SL(\ell,\real)$, and $G$ is
defined over~$F$, so there is some subset~$\mathcal{Q}$ of
$F[x_{1,1},\ldots,x_{\ell,\ell}]$, such that $G^\circ =
\Var(\mathcal{Q})^\circ$. For each place~$\sigma$ of~$F$, let
 $$ G^\sigma = \Var_{F_\sigma}\bigl( \sigma(\mathcal{Q})
\bigr)^\circ .$$
 Then $G^\sigma$, the \defit[Galois!conjugate of~$G$]{Galois
conjugate} of~$G$ by~$\sigma$, is defined over~$\sigma(F)$.
 \end{enumerate}
 \end{defn}

\begin{other} \label{onlyinfinite}
Our definition requires places to be \defit[place!infinite]{infinite} (or
\defit[place!archimedean]{archimedean}); that is the
reason for the superscript~$\infty$ on~$S^\infty$. Other authors also
allow places that are \defit[place!finite]{finite} (or
\defit[place!nonarchimedean]{nonarchimedean}, or
\defit[place!p-adic@$p$-adic]{$p$-adic}). These additional places are of
fundamental importance in number theory, and, therefore, in
deeper aspects of the theory of arithmetic groups. For
example, superrigidity at the finite places will play a crucial
role in the proof of the Margulis Arithmeticity Theorem
in \cref{MargArithPf}. 
Finite places are also essential for the definition of the ``$S$-arithmetic'' groups discussed in \cref{SarithChap}.
 \end{other}

\begin{prop} \label{ResScal->Latt}
 If $G$ is defined over an algebraic
number field~$F \subset \real$, and $\ints$ is the
ring of integers of~$F$, then
there is a finite-index subgroup $\dot G_\ints$ of~$G_\ints$, such that
	$$ \text{$\dot G_{\ints}$ embeds as an arithmetic subgroup of
 \ $ \displaystyle \prod_{\sigma \in S^\infty} G^\sigma $} , $$
 via the natural embedding $\Delta \colon \gamma \mapsto
\bigl( \sigma(\gamma) \bigr)_{\sigma \in S^\infty}$ 

Furthermore, if $G$ is simple, then the lattice~$\Delta(G_{\ints})$ is irreducible.
 \end{prop}

\begin{warn}
By our definition, $G^\sigma$ is always connected, since it is the identity component of $\Var_{F_\sigma}\bigl( \sigma(\mathcal{Q})
\bigr)$. If $G$ is assumed to be Zariski closed (so it is equal to $\Var(\mathcal{Q})$, rather than merely being isogenous to it), then it is sometimes more convenient to define $G^\sigma$ to be the entire variety $\Var_{F_\sigma}\bigl( \sigma(\mathcal{Q})
\bigr)$, rather than merely the identity component. In particular, that would eliminate the need to pass to a finite-index subgroup $\dot G_\ints$ in the statement of \cref{ResScal->Latt}. Taking the best of both worlds, we will usually ignore the difference between $G_\ints$ and $\dot G_\ints$, and pretend that the map $\Delta$ of \cref{ResScal->Latt} is defined on all of $G_\ints$. For example, the statements of \cref{scalars->cpct,GZirred->scalars} below omit the dots that should be in $\Delta(\dot G_\ints)$ and $\phi \bigl( \Delta(\dot H_\ints) \bigr)$.
\end{warn}

The argument in the last paragraph of the proof of
\cref{SO(12;Z[sqrt2])} shows the following:

\begin{cor} \label{scalars->cpct}
 If $G^\sigma$ is compact, for some $\sigma \in S^\infty$,
then $\Delta(G_{\ints})$ is cocompact.
 \end{cor}

\begin{rem} \label{ResScal(FnotinR)}
 \Cref{ResScal->Latt} is stated only for real
groups, but the same conclusions hold if
\noprelistbreak
 \begin{itemize}
 \item $G \subseteq \SL(\ell,\complex)$, 
 \item $F$~is an algebraic number field, such that $F \not
\subset \real$, and
 \item $G$ is defined over~$F$, as an algebraic
group over~$\complex$; that is, there is  a
subset~$\mathcal{Q}$ of $F[x_{1,1},\ldots,x_{\ell,\ell}]$,
such that $G^\circ = \Var_{\complex}(\mathcal{Q})^\circ$  \csee{GxCNot}.
 \end{itemize}
 For example, we have the following irreducible arithmetic lattices:
 \begin{enumerate}
 \item $\SO \bigl( n, \integer [ i, \!\sqrt{2} ]
\bigr)$ 
in
$\SO(n,\complex) \times \SO(n,\complex)$, and
 \item
  $\SO \! \left( n, \integer \left[ \! \sqrt{1 - \sqrt{2}} \right]
\right)$ in
$\SO(n,\complex) \times \SO(n,\real) \times \SO(n,\real)$.
 \end{enumerate}
 \end{rem}

The following converse shows that restriction of scalars is
the only way to make a group of $\integer$-points that is irreducible.

\begin{prop} \label{GZirred->scalars}
 If\/ $\Gamma = G_{\integer}$ is an irreducible lattice
in~$G$ \textup(and $G$ is connected\/\textup), then there exist
\noprelistbreak
 \begin{enumerate}
 \item an algebraic number field~$F$, with
completion~$F_\infty$ \textup($= \real$
or\/~$\complex$\textup),
 \item a connected, simple subgroup~$H$ of\/
$\SL(\ell,F_{\infty})$, for some~$\ell$,  such that $H$~is
defined over~$F$ \textup(as an algebraic
group over~$F_{\infty}$\textup), and
 \item an isogeny
 $$\phi \colon \prod_{\sigma \in S^\infty} H^\sigma \to G ,$$
 \end{enumerate}
 such that 
 $\phi \bigl( \Delta(H_{\ints}) \bigr)$
 is commensurable to~$\Gamma$.
 \end{prop}

\begin{proof} 
 It is easier to work with the algebraically closed
field~$\complex$, instead of~$\real$, so, to avoid minor
complications, let us assume that $G \subseteq \SL(\ell,
\complex)$ is defined over~$\rational[i]$ (as an algebraic
group over~$\complex$), and that $\Gamma = G_{\integer[i]}$.
This assumption results in a loss of generality, but similar
ideas apply in general.

Write $G = G_1 \times \cdots \times G_r$, where each
$G_i$~is simple. Let $H = G_1$. We remark that if $r = 1$,
then the desired conclusion is obvious: let $F =
\rational[i]$, and let $\phi$~be the identity map.

Let $\Sigma$ be the Galois group of~$\complex$
over~$\rational[i]$. Because $G$ is defined
over~$\rational[i]$, we have $\sigma(G) = G$ for every
$\sigma \in \Sigma$. Hence, $\sigma$~must permute the simple
factors $\{G_1,\ldots,G_r\}$. 

We claim that $\Sigma$ acts transitively on
$\{G_1,\ldots,G_r\}$. To see this, suppose, for example,
that $r = 5$, and that $\{G_1,G_2\}$ is invariant
under~$\Sigma$. Then $A = G_1 \times G_2$ is invariant
under~$\Sigma$, so $A$ is defined
over~$\rational[i]$. Similarly, $A' = G_3 \times G_4 \times
G_5$ is also defined over~$\rational[i]$. Then
$A_{\integer[i]}$ and $A'_{\integer[i]}$ are lattices in $A$
and~$A'$, respectively, so $\Gamma = G_{\integer[i]} \approx
A_{\integer[i]} \times A'_{\integer[i]}$ is reducible. This
is a contradiction.

Let
 $$ \Sigma_1 = \{\, \sigma \in \Sigma \mid \sigma(G_1) = G_1
\,\}$$
 be the stabilizer of~$G_1$, and let
 $$ F = \{\, z \in \complex \mid \sigma(z) = z, \ \forall
\sigma \in \Sigma_1 \,\} $$
 be the fixed field of~$\Sigma_1$. Because $\Sigma$ is
transitive on a set of $r$~elements, we know that $\Sigma_1$
is a subgroup of index~$r$ in~$\Sigma$, so Galois Theory
tells us that $F$~is an extension of $\rational[i]$ of
degree~$r$.

Since $\Sigma_1$ is the Galois group of~$\complex$ over~$F$,
and $\sigma(G_1) = G_1$ for all $\sigma \in \Sigma_1$, we see
that $G_1$ is defined over~$F$. 

Let $\sigma_1,\ldots,\sigma_r$ be coset representatives
of~$\Sigma_1$ in~$\Sigma$. Then $\sigma_1|_F, \ldots,
\sigma_r|_F$ are the $r$~places of~$F$ and, after
renumbering, we have $G_j = \sigma_j(G_1)$.
So (with $H = G_1)$, we have
 $$\prod_{\sigma \in S^\infty} H^\sigma
 = H^{\sigma_1|_F} \times \cdots \times H^{\sigma_r|_F}
 = \sigma_1(G_1) \times \cdots \times \sigma_r(G_1)
 = G_1 \times \cdots \times G_r
 = G .$$
 Let $\phi$ be the identity map.

For $h \in H_F$, let $\Delta'(h) = \prod_{j=1}^r
\sigma_j(h)$. Then $\sigma \bigl( \Delta'(h) \bigr) =
\Delta'(h)$ for all $\sigma \in \Sigma$, so $\Delta'(h) \in
G_{\rational[i]}$. In fact, it is not difficult
to see that $\Delta'(H_F) = G_{\rational[i]}$, and then one
can verify that
 $\Delta'(H_{\ints}) \approx G_{\integer[i]} =
\Gamma$, so $\phi \bigl( \Delta(H_{\ints})$ is
commensurable to~$\Gamma$.
 \end{proof}

\begin{rem} \label{ROSAbsSimple}
Although it may not be clear from our proof, the group $G'$ in \cref{irred->scalars} can be chosen to be ``\term[simple!absolutely]{absolutely simple}\zz.'' This means that if $F \subset \real$, then the following three equivalent conditions must be true: $G'$ remains simple over~$\complex$, $\Lie G' \otimes_\real \complex$ is simple, and $G'$ is not isogenous to any ``complexification'' $(G'')_\complex$.
\end{rem}

Combining \cref{GZirred->scalars} with
\cref{scalars->cpct} yields the following result.

\begin{cor} \label{noncocpct->nocpct}
 If $G_{\integer}$ is an irreducible lattice in~$G$, and
$G/G_{\integer}$ is not cocompact, then $G$ has no compact
factors.
 \end{cor}

By combining \cref{GZirred->scalars} with \cref{ArithDefn}, we see that every irreducible arithmetic subgroup can be constructed by using restriction of scalars, and then modding out a compact subgroup:

\begin{cor} \label{irred->scalars}
 If\/ $\Gamma$ is an irreducible, arithmetic lattice in~$G$ \textup(and $G$ is connected\/\textup), then there exist
 \begin{enumerate}
 \item an algebraic number field~$F$, with
completion~$F_\infty$ \textup($= \real$
or~$\complex$\textup),
 \item a connected, simple subgroup~$G'$ of\/
$\SL(\ell,F_{\infty})$, for some~$\ell$,  such that $G'$~is
defined over~$F$ \textup(as an algebraic group
over~$F_{\infty}$\textup), and
 \item a continuous surjection 
 $$\phi \colon \prod_{\sigma \in S^\infty} (G')^\sigma \to G ,$$
 with compact kernel,
 \end{enumerate}
 such that 
 $\phi \bigl( \Delta(G'_{\ints}) \bigr)$
 is commensurable to~$\Gamma$.
 \end{cor}
 
 When $G$ is simple, the restriction of~$\phi$ to some simple factor of $\prod_{\sigma \in S^\infty} (G')^\sigma$ must be an isogeny, so the conclusion can be stated in the following much simpler form:

\begin{cor} \label{simple->Arith=Ints}
 If\/ $\Gamma$ is an arithmetic subgroup of~$G$, and $G$ is simple, then there exist
 \begin{enumerate}
 \item an algebraic number field~$F$, with
completion~$F_\infty$ \textup($= \real$
or~$\complex$\textup),
 \item a connected, simple subgroup~$G'$ of\/
$\SL(\ell,F_{\infty})$, for some~$\ell$,  such that $G'$~is
defined over~$F$ \textup(as an algebraic group
over~$F_{\infty}$\textup), and
 \item an isogeny
 $\phi \colon G' \to G$,
 \end{enumerate}
 such that 
 $\phi( G'_{\ints})$
 is commensurable to~$\Gamma$.
 \end{cor}

However, we should point out that this result is of interest only when $\Gamma$ is cocompact (or is reducible with at least one cocompact factor). This is because there is no need for restriction of scalars when the irreducible lattice~$\Gamma$ is not cocompact \csee{IrredNoncpct->NoROS}.

\begin{exercises}

\item \label{End(R4)Q}
 In the notation of the proof of \cref{SL(2Z[sqrt2])},
show, for the $\rational$-form $\Delta(F^2)$ of~$\real^4$,
that
 $$ \End(\real^4)_{\rational}
 = \bigset{
 \begin{bmatrix}
 A & B \\
 \sigma(B) & \sigma(A)
 \end{bmatrix}
 }{
 A,B \in \Mat_{2 \times 2}(F)
 } .$$
\hint{Since the $F$-span of $\Delta(F^2)$ is~$F^4$,
we have $\End(\real^4)_{\rational} \subseteq \Mat_{4 \times
4}(F)$. Thus, for any $T \in \End(\real^4)_{\rational}$, we
may write
 $ T = \begin{Smallbmatrix}
 A & B \\
 C & D
 \end{Smallbmatrix}$,
 with $A,B,C,D \in \Mat_{2 \times 2}(F)$. Now use the fact
that, for all $u \in F^2$, we have $T(u) = \bigl( v,
\sigma(v) \bigr)$, for some $v \in F^2$.}

\item \label{SL2xSL2Defd/Q}
 In the notation of the proof of \cref{SL(2Z[sqrt2])},
let 
 \begin{align*}
  \mathcal{Q} = &\bigset{ x_{i,j+2}  + x_{i+2,j}, ~ x_{i,j+2} x_{i+2,j} }{  1 \le i, j \le 2 }
	\\& \quad \cup \left\{ \frac{1}{\sqrt{2}} \Bigl( (x_{1,1} x_{2,2} - x_{1,2} x_{2,1}) - (x_{3,3} x_{4,4} - x_{3,4} x_{4,3}) \Bigr) \right\} 
	. \end{align*}
 \begin{enumerate}
 \item Use the conclusion of \cref{End(R4)Q} to show
that each $Q \in \mathcal{Q}$ is defined over~$\rational$.
 \item Show that $\Var(\mathcal{Q})^\circ = \SL(2,\real)
\times \SL(2,\real)$.
 \end{enumerate}

\item \label{SL2xSL2-Delta(Gamma)=}
 In the notation of the proof of \cref{SL(2Z[sqrt2])},
use \cref{End(R4)Q} to show that
$G_{\Delta(\ints^2)} = \Delta(\Gamma)$.

\item \label{GxG-Defd/Q}
 Let $F$, $\ints$, $\sigma$, $\Delta$ be as in the
proof of \cref{SL(2Z[sqrt2])}.
If $G \subseteq \SL(\ell,\real)$, and $G$ is defined over~$F$,
show $G \times G$ is defined over~$\rational$ (with
respect to the $\rational$-form on $\End(\real^{2\ell})$
induced by the $\rational$-form $\Delta(F^\ell)$
on~$\real^{2\ell}$).
\hint{For each $Q \in
\rational[x_{1,1},\ldots,x_{\ell,\ell}]$, let us define a
corresponding polynomial $Q^+ \in
\rational[x_{\ell+1,\ell+1},\ldots,x_{2\ell,2\ell}]$ by
replacing every occurrence of each variable $x_{i,j}$ with
$x_{\ell+i,\ell+j}$. For example, if $\ell = 2$, then
 $$ (x_{1,1}^2 + x_{1,2} x_{2,1} - 3 x_{1,1} x_{2,2})^+
 = x_{3,3}^2 + x_{3,4} x_{4,3} - 3 x_{3,3} x_{4,4} .$$
 Choose $\mathcal{Q}_0 \subset
\rational[x_{1,1},\ldots,x_{\ell,\ell}]$ that defines $G$ as
a subgroup of $\SL(\ell,\real)$, and let
 $$ \mathcal{Q}_1 = 
 \{\, Q + \sigma(Q^+), ~ Q \, \sigma(Q^+) \mid Q \in
\mathcal{Q}_0 \,\}.$$
 A natural generalization of \cref{SL2xSL2Defd/Q} shows
that $\SL(\ell,\real) \times \SL(\ell,\real)$ is defined
over~$\rational$: let $\mathcal{Q}_2$ be the corresponding
set of $\rational$-polynomials. Now define $\mathcal{Q} =
\mathcal{Q}_1 \cup \mathcal{Q}_2$.}

\item \label{Delta(O)=VSLatt}
 Suppose $\ints$ is the ring of integers of an
algebraic number field~$F$.
 \begin{enumerate}
 \item \label{Delta(O)=VSLatt-discrete}
 Show $\Delta(\ints)$ is discrete in
$\bigoplus_{\sigma \in S^\infty} F_{\sigma}$.
 \item \label{Delta(O)=VSLatt-F}
 Show $\Delta(F)$ is a $\rational$-form of
$\bigoplus_{\sigma \in S^\infty} F_{\sigma}$.
 \item \label{Delta(O)=VSLatt-O}
 Show $\Delta(\ints)$ is a $\integer$-lattice
in $\Delta(F)$.
 \end{enumerate}

\item \label{SO(B)Z2cocpt(Mahler)}
 Let 
 \begin{itemize}
 \item $B(v,w) = v_1 w_1 + v_2 w_2 - \sqrt{2} v_3 w_3$, for
$v,w \in \real^3$,
 \item $G = \SO(B)^\circ$,
 \item $G^* = G \times G^\sigma$,
 \item $\Gamma = G_{\integer[\sqrt{2}]}$,
 and 
 \item $\Gamma^* = \Delta(\Gamma)$.
 \end{itemize}
 Show:
 \begin{enumerate}
 \item The image of~$G^*$ in $\SL(6,\real) /
\SL(6,\real)_{\Delta(\ints^3)}$ is precompact (by using
the Mahler Compactness Criterion).
 \item The image of~$G^*$ in $\SL(6,\real) /
\SL(6,\real)_{\Delta(\ints^3)}$ is closed.
 \item $G^*/\Gamma^*$~is compact.
 \item $G/\Gamma$~is compact (without using the fact that 
$\Gamma$ is a lattice in~$G$).
 \end{enumerate}
 \hint{This is similar to \cref{anis->cocpct}.}

\item \label{SLnFsigma/Q}
 For any algebraic number field~$F$,
the $\rational$-form $\Delta(F^\ell)$ on $\bigoplus_{\sigma \in
S^\infty} (F_{\sigma})^\ell$ induces a natural $\rational$-form on $\End_{\real} \bigl( \bigoplus_{\sigma
\in S^\infty} (F_{\sigma})^\ell \bigr)$.
 Show the group
 $\prod_{\sigma \in S^\infty} \SL(\ell,F_\sigma)$ is
defined over~$\rational$, with respect to this $\rational$-form.
 \hint{This is a generalization of
\cref{SL2xSL2Defd/Q}. That proof is based on the
elementary symmetric functions of two variables:
$P_1(a_1,a_2) = a_1 + a_2$ and $P_2(a_1,a_2) = a_1 a_2$. For
the general case, use symmetric functions of $d$~variables,
where $d$~is the degree of~$F$ over~$\rational$.}

\item Suppose $G \subseteq \SL(\ell,\real)$, and $G$ is
defined over an algebraic number field $F \subset \real$. Show
$\prod_{\sigma \in S^\infty} G^\sigma$ is defined
over~$\rational$, with respect to the
$\rational$-form on $\End_{\real} \bigl( \bigoplus_{\sigma
\in S^\infty} (F_{\sigma})^\ell \bigr)$ induced by the
$\rational$-form $\Delta(F^\ell)$ on~$\bigoplus_{\sigma \in
S^\infty} (F_{\sigma})^\ell$.
 \hint{This is a generalization of \cref{GxG-Defd/Q}.
See the hint to \cref{SLnFsigma/Q}.}

\item \label{DenseProjSO(n)}
 Show, for all $m,n \ge 1$, with $m + n \ge 3$, that there exist 
a lattice~$\Gamma$ in $\SO(m,n)$,
 and
 a homomorphism $\rho \colon \Gamma \to \SO(m+n)$,
 such that $\rho(\Gamma)$ is dense in $\SO(m+n)$.

\end{exercises}

\section{Only isotypic groups have irreducible lattices}

Intuitively, the \defit{complexification}~$G_\complex$ of~$G$ is the
complex Lie group that is obtained from~$G$ by replacing
real numbers with complex numbers. 
For example, $\SL(n,\real)_\complex = \SL(n,\complex)$, and $\SO(n)_\complex = \SO(n,\complex)$.
(See \cref{RFormsOfCGrps} for more discussion of this.)

\begin{defn} \label{IsotypicDefn}
 $G$ is \defit[isotypic semisimple Lie group]{isotypic} if
all of the simple factors of~$G_\complex$ are isogenous to each other.
 \end{defn}

For example, $\SL(2,\real) \times \SL(3,\real)$ is not
isotypic, because $\SL(2,\complex)$ is not isogenous to $\SL(3,\complex)$.
Similarly, $\SL(5,\real)
\times \SO(2,3)$ is not isotypic, because the complexification of 
 $\SL(5,\real)$ is $\SL(5,\complex)$, but the complexification of $\SO(2,3)$ is (isomorphic to) $\SO(5, \complex)$. Therefore, the
following consequence of the arithmeticity theorem implies
that neither $\SL(2,\real) \times \SL(3,\real)$ nor
$\SL(5,\real) \times \SO(2,3)$ has an irreducible lattice.

\begin{thm}[(Margulis)] \label{irred->isotypic}
 Assume that $G$ has no compact factors. 
If $G$ has an \index{irreducible!lattice}{irreducible lattice}, then $G$ is \term[isotypic semisimple Lie group]{isotypic}.
 \end{thm}

\begin{proof} 
Suppose $\Gamma$ is an irreducible lattice
in~$G$. We may assume that $G$ is not simple (otherwise,
the desired conclusion is trivially true), so $G$~is
neither $\SO(1,n)$ nor $\SU(1,n)$. Therefore, from the Margulis
Arithmeticity Theorem \pref{MargulisArith}, we know that
$\Gamma$ is arithmetic. Then, since $\Gamma$ is irreducible,
\cref{irred->scalars} implies there is a simple subgroup~$G'$ of some $\SL(\ell,\real)$, and a compact group~$K$,
such that
 \begin{itemize}
 \item $G'$ is defined over a number field~$F$, 
 and
 \item $G  \times K$ is isogenous to $\prod_{\sigma \in S^\infty}
(G')^\sigma$.
 \end{itemize}
 So the simple factors of $G \times K$ are all in $\{\, (G')^\sigma \mid \sigma \in S^\infty\,\}$ (up to isogeny). 
 It then follows from \cref{Gsigma=G} below that $G$ is isotypic.
 \end{proof}

\begin{rems} \label{Irred->IsotypicRem} \ 
\noprelistbreak
	\begin{enumerate}
	\item We will prove the converse of \cref{irred->isotypic} in \cref{Irred->Isotypic} (without the assumption that $G$ has no compact factors).
	
	\item By arguing just a bit more carefully, it can be shown
that \cref{irred->isotypic} remains valid when the assumption that $G$ has no compact factors is replaced with the weaker hypothesis that $G$ is not isogenous to $\SO(1,n)
\times K$ or $\SU(1,n) \times K$, for any nontrivial,
connected compact group~$K$ \csee{irred<>isotypic(cpct)}. 

	\end{enumerate}
 \end{rems}

The following example shows that a nonisotypic group can have irreducible lattices, so some restriction on~$G$ is
necessary in \cref{irred->isotypic}.

\begin{eg}
 $\SL(2,\real) \times K$ has an irreducible lattice, for any
connected, compact Lie group~$K$ \ccf{irredinSL2xSO3}.
 \end{eg}

We now complete the proof of \cref{irred->isotypic}:

\begin{lem} \label{Gsigma=G}
 Assume $G$~is defined over an algebraic number field~$F$.
 If $\sigma$ is a place of~$F$,
 and
  $G$ is simple,
  then the complexification of~$G$ is isogenous to the complexification of~$G^\sigma$.
 \end{lem}

%\begin{lem} \label{Gsigma=G(real)}
% If
% \begin{itemize}
% \item $G$~is defined over an algebraic number field $F
%\subset \real$, and
% \item $\sigma$ is an embedding of~$F$ in~$\real$,
% \end{itemize}
% then $G \otimes \complex$ is isogenous to $G^\sigma \otimes
%\complex$.
% \end{lem}

\begin{proof}
 Extend~$\sigma$ to an automorphism~$\widehat\sigma$
of~$\complex$. Then $\widehat\sigma(G_\complex) =
(G^\sigma)_\complex$, so it is clear that $G_\complex$ is isomorphic to $(G^\sigma)_\complex$.
Unfortunately, however, the automorphism~$\widehat\sigma$ is not
continuous (not even measurable) unless it happens to be the
usual complex conjugation, so we have only an isomorphism of
abstract groups, not an isomorphism of Lie groups. Hence, this observation is not
a proof,
although it is suggestive.
To give a rigorous proof, it is easier to work at the Lie algebra
level.

First, let us make an observation that will also be pointed out in \cref{GxC=GxGiff}. If $G = \SL(n,\complex)$, or, more generally, if $G$ is isogenous to a complex group $G'_\complex$, then $G_\complex = G \times G$ (because $\complex \otimes_\real \complex \iso \complex \oplus \complex$). So $G_\complex$ is not simple. However, it can be shown that this is the only situation in which the complexification of a simple group fails to be simple: if $G$ is simple, but $G_\complex$ is not simple, then $G$ is isogenous to a complex simple group $G'_\complex$. 
Therefore, although the complexification of a simple group is not always simple, it is always isotypic.

Now assume, for definiteness, that $F \subset \real$ \csee{CPfOfGsigma=G}.
Since $G$ is defined over~$F$, its Lie algebra~$\Lie G$ is also defined over~$F$. This means there is a basis $\{v_1,\ldots,v_n\}$ of~$\Lie G$, such that the corresponding structure constants $\{c_{j,k}^\ell\}_{j,k,\ell=1}^n$ all belong to~$F$; recall that the structure constants are defined by the formula
 	$$ \textstyle [v_j,v_k] = \sum_{\ell=1}^n c_{j,k}^\ell v_\ell .$$
 
 Because $G$ is isogenous to a group that is defined
over~$\rational$ \csee{hasQform}, there is also a basis
$\{u_1,\ldots,u_n\}$ of~$\Lie G$ whose structure
constants are in~$\rational$. Write
 $ v_k = \sum_{\ell=1}^n \alpha_k^\ell u_\ell $
 with each $\alpha_k^\ell \in \real$, and define
 $$ \textstyle v_k^\sigma = \sum_{\ell=1}^n \widehat\sigma(\alpha_k^\ell)
u_\ell . $$
 Then $v_1^\sigma, \ldots, v_n^\sigma$ is a basis of $\Lie G
\otimes_{\real} \complex$ whose structure constants are
$\bigl\{ \sigma(c_{j,k}^\ell) \bigr\}_{j,k,\ell=1}^n$. These
are obviously the structure constants of the Lie algebra $\Lie G^\sigma$ of~$G^\sigma$.

If $\sigma(F) \subset \real$, then the $\real$-span of $\{v_1^\sigma, \ldots, v_n^\sigma\}$ is (isomorphic to)~$\Lie G^\sigma$, so its $\complex$-span is $\Lie G^\sigma
\otimes_{\real} \complex$. Since $v_1^\sigma, \ldots, v_n^\sigma$ is also a basis of $\Lie G
\otimes_{\real} \complex$, we conclude that $(G^\sigma)_\complex$ is isogenous to~$G_\complex$.

Finally, if $\sigma(F) \not\subset \real$, then the $\complex$-span of $\{v_1^\sigma, \ldots, v_n^\sigma\}$ is (isomorphic to)~$\Lie G^\sigma$, so $\Lie G
\otimes_{\real} \complex = \Lie G^\sigma$. This implies that $(G^\sigma)_\complex$ is isogenous to~$G_\complex$.
 \end{proof}

\begin{rem} \label{Gsigma=GRems} 
The proof of \cref{Gsigma=G} used our standing assumption that $G$~is
semisimple only to show that $G$~is isogenous to a group
that is defined over~$\rational$. 
See \cref{H<>Hsigma} for an example of a Lie
group~$H$, defined over an algebraic number field $F \subset
\real$, and an embedding~$\sigma$ of~$F$ in~$\real$, such
that $H \times H^\sigma$ is not isotypic.
 \end{rem}

\begin{exercises}

\item Show, for $m,n \ge 2$, that $\SL(m,\real) \times
\SL(n,\real)$ has an irreducible lattice if and only if $m =
n$.

\item \label{irred<>isotypic(cpct)}
 Suppose $G$ is not isogenous to $\SO(1,n)
\times K$ or $\SU(1,n) \times K$, for any nontrivial,
connected compact group~$K$. Show that if $G$~has an
irreducible lattice, then $G$~is isotypic.
 \hint{Use \cref{MargArithThmCpctFacts} to modify the proof of \cref{irred->isotypic}.}

\item \label{H<>Hsigma} \optional\ 
 For $\alpha \in \complex \smallsetminus \{0,-1\}$, let
$\Lie H_\alpha$ be the 7-dimensional, nilpotent Lie algebra
over~$\complex$, generated by $\{x_1,x_2,x_3\}$, such that 
 \begin{itemize}
 \item $[\Lie H_\alpha, x_1,x_1] = [\Lie H_\alpha, x_2,x_2]
= [\Lie H_\alpha,x_3,x_3] = 0$, and
 \item $[x_2,x_3,x_1] = \alpha[x_1,x_2,x_3]$.
 \end{itemize}
 \begin{enumerate}
 \item Show that $[x_3,x_1,x_2] = -(1+\alpha) [x_1,x_2,x_3]$.
 \item For $h \in \Lie H_\alpha$, show that $[\Lie
H_\alpha,h,h] = 0$ if and only if there exists $x \in
\{x_1,x_2,x_3\}$ and $t \in \complex$, such that $h \in t x
+ [\Lie H_\alpha,\Lie H_\alpha]$.
 \item Show $\Lie H_\alpha \iso \Lie H_\beta$ iff 
 $\beta \in \left\{\alpha, \frac{1}{\alpha}, -(1+\alpha),
-\frac{1}{1+\alpha},
-\frac{\alpha}{1+\alpha}, -\frac{1+\alpha}{\alpha} \right\} $.
 \item Show that if the degree of $\rational(\alpha)$
over~$\rational$ is at least~$7$, then there is a
place~$\sigma$ of $\rational(\alpha)$,
such that $\Lie H_\alpha$ is not
isomorphic to $(\Lie H_\alpha)^\sigma$.  
 \end{enumerate}

\item \optional\ 
In the notation of \cref{H<>Hsigma}, show that if
the degree of $\rational(\alpha)$ over~$\rational$ is at
least~$7$, then  $\Lie H_\alpha$ is not isomorphic to any
Lie algebra that is defined over~$\rational$.

\item \optional\ 
In the notation of \cref{H<>Hsigma}, show, for
$\alpha = \sqrt{2} - (1/2)$, that $\Lie H_\alpha$ is
isomorphic to a Lie algebra that is defined over~$\rational$.
 \hint{Let $y_1 = x_1 + x_2$ and $y_2 = (x_1 -
x_2)/\sqrt{2}$. Show that the $\rational$-subalgebra
of~$\Lie H_\alpha$ generated by $\{y_1,y_2,x_3\}$ is a
$\rational$-form of~$\Lie H_\alpha$.}

\item \label{CPfOfGsigma=G}
Carry out the proof of \cref{Gsigma=G} for the case where $F \not\subset \real$. 
\hint{Write $\Lie G = \Lie G' \otimes_\real \complex$ and let $\{u_1,\ldots,u_n\}$ be a basis of~$\Lie G'$ with rational structure constants. Show that $G$ is isogenous to either $G^\sigma$ or $(G^\sigma)_\complex$.}

\end{exercises}

\begin{notes}

The fact that $G$ is \term{unirational} (used in \cref{GuniratSoQDense}) is proved in \cite[Thm.~18.2, %(ii), 
p.~218]{Borel-LinAlgGrps}.

The Margulis Arithmeticity Theorem \pref{MargulisArith} was
proved by Margulis \cite{MargulisArithProp, MargulisArith}
under the assumption that $\Rrank G \ge 2$. (Proofs also
appear in \cite[Thm.~A, p.~298]{MargulisBook} and
\cite{ZimmerBook}.) Much later, the superrigidity theorems
of Corlette \cite{Corlette} and Gromov-Schoen
\cite{GromovSchoen} extended this to all groups except
$\SO(1,n)$ and $\SU(1,n)$.
% We remark that these proofs rely on the fact that $\Gamma$
%is finitely generated; Venkataramana \cite{Venky-fg} showed
%how to avoid using this assumption.

\Cref{hasQform} is a weak version of a theorem of Borel \cite{Borel-CK}. (A proof also appears in \cite[Chap.~14]{RaghunathanBook}.)

The Commensurability Criterion (\fullref{GQComm}{criterion}) is due to Margulis \cite{Margulis-DiscGrpMot}. We will see it again in \cref{CommCriterion}, and it is proved in
\cite{A'CampoBurger}, \cite{MargulisBook}, and
\cite{ZimmerBook}.

The fact that all noncocompact lattices have unipotent elements (that is, the generalization of \cref{GodementNoCpctFactor} to the nonarithmetic case) is due to D.\,Kazhdan and G.\,A.\,Margulis \cite{KazhdanMargulis} (or see \cite{Borel-KazhdanMargulisBourbaki} or \cite[Cor.~11.13, p.~180]{RaghunathanBook}).

The standard reference on restriction of scalars is \cite[\S1.3, pp.~4--9]{Weil-AdelesAlgGrps}.
(A discussion can also be found in \cite[\S2.1.2, pp.~49--50]{PlatonovRapinchukBook}.)

\Cref{GZirred->scalars} (and \cref{ROSAbsSimple}) is due to A.\,Borel and J.\,Tits \cite[6.21(ii), p.~113]{BorelTits-GrpRed}.

See \cite[Cor.~IX.4.5, p.~315]{MargulisBook} for a proof of \cref{irred->isotypic}.

\end{notes}

 %!TEX root = IntroArithGrps.tex

\mychapter{\texorpdfstring{Examples of\\Arithmetic Groups}%
	{Examples of Arithmetic Groups}}
\label{EgArithGrpsChap}

\prereqs{definition of arithmetic subgroup (\cref{ArithLattDefnSect}), Godement Criterion (\cref{GodementCriterion}), and restriction of scalars (\cref{RestrictScalarsSect}).}

\section{Arithmetic subgroups of \texorpdfstring{$\SL(2,\real)$}{SL(2,R)} via orthogonal groups} 
\label{ArithLattSL2}

$\SL(2,\integer)$ is the obvious example of an arithmetic
subgroup of $\SL(2,\real)$. Later in this section, we will show that (up to
commensurability and conjugates) it is the only one that is
not cocompact \csee{NoncocpctSL2R=SL2Z}. In contrast, there are infinitely many
cocompact, arithmetic subgroups. They can be constructed by
several different methods. Perhaps the easiest way is to
note that $\SL(2,\real)$ is isogenous to the special orthogonal group $\SO(2,1)$. 

\begin{notation} \label{SO3wayNotation}
In this chapter (and others), we will see many different special orthogonal groups over a field~$F$.
They can be specified in (at least) three different, but equivalent ways:
	\begin{enumerate}
	
	\item (\term{Gram matrix}) For a symmetric, invertible matrix $A \in \Mat_{\ell \times \ell}(F)$, we define
	$$ \SO(A;F) = \{\, g \in \SL(n,F) \mid g^\transpose A g = A \,\} .$$
	This is the approach taken to the definition of $\SO(m,n)$ in \cref{classical-fulllinear}.

	\item (\term[bilinear form]{Bilinear form}) A symmetric, bilinear form~$B$ on~$F^\ell$ is \defit[nondegenerate!bilinear form]{nondegenerate} if, for all nonzero $v \in F^\ell$, there exists $w \in F^\ell$, such that $B(v,w) \neq 0$. We define
	$$ \SO(B;F) = \{\, g \in \SL(\ell,F) \mid B(gv,gw) = B(v,w), \ \forall v,w \in F^\ell \,\} .$$

	\item ({Quadratic form}) A \defit{quadratic form} on~$F^\ell$ is a homogeneous polynomial $Q(x_1,\ldots,x_\ell)$ of degree~$2$.  It is \defit[nondegenerate!quadratic form]{nondegenerate} if the corresponding bilinear form $B_Q$ is nondegenerate, where
	$$ B_Q(v,w) = {\textstyle\frac{1}{4}} \bigl( Q(v + w) - Q(v - w) \bigr) .$$
	We define
	 $$ \SO(Q;F) = \{\, g \in \SL(\ell,F) \mid Q(gv) = Q(v), \ \forall v \in F^\ell \,\} .$$

	\end{enumerate}
The three approaches give rise to exactly the same groups \csee{SO3ways}, and it is straightforward to translate between them, so we will use whichever notation is most convenient in a particular context.
\end{notation}	

\begin{egs} \label{ArithSO21Eg} \ 
\noprelistbreak
 \begin{enumerate}
 \item \label{ArithSO21Eg-Q}
 Fix positive integers $a$ and~$b$, and let 
 $$G = \SO(a x^2 + b y^2 - z^2; \real) \iso
\SO(2,1) .$$
 If $(0,0,0)$ is the only integer solution of the
Diophantine equation $a x^2 + b y^2 = z^2$, then
$G_{\integer}$ is a cocompact, arithmetic subgroup of~$G$
\csee{anis->cocpct}. See \cref{px2+py2<>z2} for some
examples of $a$ and~$b$ satisfying the hypotheses.
 \item Restriction of scalars \csee{RestrictScalarsSect}
allows us to use algebraic number fields other
than~$\rational$. Let
 \begin{itemize}
 \item \label{ArithSO21Eg-F}
 $F \neq \rational$ be a \defit[totally!real number
field]{totally real} algebraic number field (that is, an
algebraic number field with no complex places),
 \item $a,b \in F^+$, such that $\sigma(a)$ and~$\sigma(b)$
are negative, for every place $\sigma \neq \Id$,
 \item $\ints$ be the ring of integers of~$F$, and
 \item $G = \SO(a x^2 + b y^2 - z^2; \real) \iso
\SO(2,1)$. 
 \end{itemize}
 Then the group~$G_{\ints}$ is a cocompact, arithmetic subgroup
of~$G$ (cf.~\cref{SO(12;Z[sqrt2])}, or
see~\cref{ResScal->Latt,scalars->cpct}). See
\cref{sigma(ab)<0} for an example of~$F$, $a$,
and~$b$ satisfying the hypotheses.
 \item In both \pref{ArithSO21Eg-Q}
and~\pref{ArithSO21Eg-F}, the group~$G$ is conjugate to
$\SO(2,1)$, via the diagonal matrix 
 $$ g = \diag \bigl( \sqrt{a}, \sqrt{b}, 1 \bigr) .$$
 Therefore, $g^{-1} (G_{\integer}) g$ or $g^{-1}
(G_{\ints}) g$ is a cocompact, arithmetic subgroup of
$\SO(2,1)$.
 \end{enumerate}
 \end{egs}

\begin{rem} 
 For $a$ and~$b$ as in \fullcref{ArithSO21Eg}{F}, $(0,0,0)$
is the only solution in~$\ints^3$ of the equation
$ax^2 + by^2 = z^2$ \csee{sigma(ab)neg->nosoln}.
Therefore, \fullcref{ArithSO21Eg}{Q} and
\fullcref{ArithSO21Eg}{F} could fairly easily be combined into a
single construction, but we separated them to keep them a bit less complicated.
 \end{rem}

\begin{prop} \label{CocpctArithSO21}
 The only cocompact, arithmetic subgroups of\/ $\SO(2,1)$ are the 
arithmetic subgroups constructed in \cref{ArithSO21Eg}
\textup(up to commensurability and conjugates\textup).

More precisely, any cocompact, arithmetic subgroup of\/
$\SO(2,1)$ has a conjugate that is commensurable
to an arithmetic subgroup constructed in \cref{ArithSO21Eg}.
 \end{prop} 

\begin{proof}
Let $\Gamma$ be a cocompact, arithmetic subgroup of
$\SO(2,1)$. Ignoring the minor technical issue that not all
automorphisms are inner \ccf{OutGFinite},
it suffices to
show that there is an automorphism~$\alpha$ of $\SO(2,1)$,
such that $\alpha(\Gamma)$ is commensurable to one of the
arithmetic subgroups constructed in \cref{ArithSO21Eg}.
 
\setcounter{step}{0}

\begin{step} \label{CocpctArithSO21Pf-isog}
 There are 
  \begin{itemize}
 \item an algebraic number field $F \subset \real$, with
ring of integers~$\ints$,
 \item a symmetric, bilinear form $B(x,y)$ on~$F^3$, and
 \item an isomorphism $\phi \colon \SO(B;\real) \to
\SO(2,1)$,
 \end{itemize}
 such that $\phi \bigl( \SO(B;\ints) \bigr)$ is
commensurable to~$\Gamma$. 
 \end{step}
 We give two proofs.

 First, we note that this follows from the classification results that will be proved in \cref{ArithClassicalChap}. Namely, a group of the form $\SO(m,n)$ does not appear in \cref{GFxC}, and it arises as the right-hand side of two different parts of \cref{GFxR}. However, $m+n = 1+2 = 3$ is odd in our situation, so only one of the listings is
relevant: $G_F$ must be $\SO(A;F)$, for some
algebraic number field $F \subset \real$. This means that
$\Gamma$ is commensurable to $\SO(A;\ints)$, where $\ints$ is the 
ring of integers of~$F$.

Second, let us give a direct proof that does not rely on
the results of \cref{ArithClassicalChap}. Because all
(irreducible) arithmetic subgroups are obtained by
restriction of scalars, and $G$ is simple, \cref{simple->Arith=Ints} tells us
there are
 \begin{itemize}
 \item an algebraic number field $F \subset \real$, with
ring of integers~$\ints$,
 \item a simple Lie group $H \subseteq \SL(\ell,\real)$ that
is defined over~$F$, and
 \item an isogeny $\phi \colon H \to \SO(2,1)$,
 \end{itemize}
 such that $\phi(H_{\ints})$ is commensurable
to~$\Gamma$. All that remains is to show that we may
identify $H_F$ with $\SO(B;F)$, for some symmetric bilinear
form~$B$ on~$F^3$.

The \term{Killing form}
 $$ \kappa(u,v) = \trace \bigl( (\ad_{\Lie H} u) (\ad_{\Lie
H} v) \bigr) $$
 is a symmetric, bilinear form on the Lie algebra~$\Lie H$.
It is invariant under $\Ad H$, so $\Ad_H$~is an isogeny
from~$H$ to $\SO(\kappa;\real)$. Pretending that $\Ad_H$ is
an isomorphism, not just an isogeny, we may identify $H$
with $\SO(\kappa;\real)$. Note that $\kappa(\Lie H_F, \Lie
H_F) \subseteq F$, so, by identifying $\Lie H_F$ with~$F^3$,
we may think of~$\kappa$ as a bilinear form on~$F^3$.

\begin{step}
 We may assume that $B(x,x) = a x_1^2 + b x_2^2 - x_3^2$
for some $a,b \in F^+$.
 \end{step}
 By choosing an orthogonal basis that diagonalizes the form, we may assume $B(x,x) = a
x_1^2 + b x_2^2 + c x_3^2$. Since $\SO(B;\real) \approx
\SO(2,1)$, we know that $\pm B(x,x)$ has signature
$(2,1)$. So we may assume $a,b,-c \in F^+$. Dividing by~$c$
(which does not change the orthogonal group) yields the
desired form.

\begin{step} \label{CocpctArithSO21Pf-totreal}
 $F$ is totally real, and both $\sigma(a)$ and~$\sigma(b)$
are negative, for all places $\sigma \neq \Id$.
 \end{step}
 Since $\Delta(G_{\ints})$ is an
irreducible lattice in $\prod_{\sigma \in S^\infty}
G^\sigma$ \csee{ResScal->Latt}, but the projection to the
first factor, namely~$G$, is~$\Gamma$, which is discrete,
we know that $G^\sigma$~is compact, for all $\sigma \neq
\Id$. This implies $G^\sigma \iso \SO(3)$, so $F_\sigma = \real$,
and the three real numbers $\sigma(a)$, $\sigma(b)$, and
$\sigma(-1)$ all have the same sign.

\begin{step} \label{CocpctArithSO21Pf-anis}
 $B$~is anisotropic over~$F$.
 \end{step} 
 Since $G_{\ints}$ is cocompact, it has no nontrivial unipotent elements \csee{GammaUnip->notcpct}. Therefore
$B(x,x) \neq 0$, for every nonzero $x \in F^3$
\csee{isotrop->unip}.
 \end{proof}

\begin{prop} \label{NoncocpctSL2R=SL2Z}
 $\SL(2,\integer)$ is the only noncocompact, arithmetic
subgroup\/~$\Gamma$ of\/ $\SL(2,\real)$ \textup(up to commensurability
and conjugates\textup).
 \end{prop} 

\begin{proof}
 Let us consider the isogenous group $\SO(2,1)$, instead of
$\SL(2,\real)$.

\setcounter{step}{0}

\begin{step} \label{NoncocptSL2R=SL2Z-isog}
 There are 
 \begin{itemize}
 \item a symmetric, bilinear form $B(x,y)$ on~$\rational^3$,
and
 \item an isogeny $\phi \colon \SO(B;\real) \to
\SO(2,1)$,
 \end{itemize}
 such that $\phi \bigl( \SO(B;\integer) \bigr)$ is
commensurable to~$\Gamma$. 
 \end{step}
 Since $\Gamma$ is not cocompact, there is an isogeny $\phi \colon G \to \SO(2,1)$, such that $G$ is defined over~$\rational$ and $\phi(G_\integer)$ is commensurable to~$\Gamma$ \csee{IrredNoncpct->NoROS}. The argument in \cref{CocpctArithSO21Pf-isog,CocpctArithSO21Pf-totreal} of the proof of
\cref{CocpctArithSO21} shows that we may assume $G = \SO(B;\real)$.

\begin{step} \label{CocpctArithSO21Pf-21}
 We may assume $B(x,x) = x_1^2 + x_2^2 - x_3^2$.
 \end{step}
 Because $\Gamma$ is not cocompact, we know that $B$~is
isotropic over~$F$ \csee{anis->cocpct}. So there is some nonzero $u \in F^3$, such that $B(u,u) = 0$. Choose $v \in F^3$, such that $B(u,v) \neq 0$. By adding a scalar multiple of~$u$ to~$v$, we may assume $B(v,v) = 0$. Now choose a nonzero $w \in F^3$ that is orthogonal to both $u$ and~$v$. After multiplying $B$ and~$u$ by appropriate scalars, we may assume $B(w,w) = 2B(u,v) = 1$. Then $B$ has the desired form with respect to the basis $w, u+v, u-v$.
 \end{proof}

\begin{rem} \label{FreeLattINSL2R}
As a source of counterexamples, it is useful to remember that $\SL(2,\integer)$ contains a free subgroup of finite index (see 
\cref{SanovIsFree} or \ref{SL2ZLotsFree}). % are both still exercises !!!
This implies that (finitely generated) nonabelian free groups are lattices in $\SL(2,\real)$.
\end{rem}

\begin{exercises}

\item \label{SO3ways} 
Show:
	\begin{enumerate}
	\item If $Q$ is nondegenerate quadratic form on~$F^\ell$, and $B_Q$ is defined as in \cref{SO3wayNotation}, then $B_Q$ is a bilinear form, and we have $\SO(B_Q; F) = B(Q;F)$.
	\item If $B$ is a nondegenerate bilinear form on $F^\ell$, and we define $Q(x) = B(x,x)$, then $Q(x)$ is a quadratic form, and $B = B_Q$.
	\item If $A$ is a symmetric, invertible matrix in $\Mat_{\ell \times \ell}(F)$, and we define $B(v,w) = v^\transpose A w$, then $B$ is a nondegenerate bilinear form, and $\SO(B; F) = \SO(A;F)$.
	\item If $B$ is a nondegenerate bilinear form on $F^\ell$, and $\{\varepsilon_1,\ldots,\varepsilon_\ell\}$ is the standard basis of~$F^\ell$, then the matrix $A = \bigl( B(\varepsilon_i,\varepsilon_j) \bigr)$ is invertible and symmetric, and we have $\SO(A;F) = \SO(B;F)$.
	\end{enumerate}

\item \label{px2+py2<>z2}
 Suppose $p$ is a prime, such that $x^2 + y^2 \equiv 0
\pmod{p}$ has only the trivial solution $x \equiv y \equiv 0
\pmod{p}$. (For example, $p$~could be~$3$.) Show that
$(0,0,0)$ is the only integer solution of the Diophantine
equation $p x^2 + p y^2 = z^2$.

\item \label{sigma(ab)<0}
 Let $F = \rational \bigl[\! \sqrt{2}, \sqrt{3} \bigr]$, and
$a = b = \sqrt{2} + \sqrt{3} - 3$. Show 
 \begin{enumerate}
 \item $F$~is a totally real extension of~$\rational$, 
 \item $a$~is positive, and 
 \item $\sigma(a)$~is negative, for every place $\sigma
\neq \Id$.
 \end{enumerate}

\item \label{sigma(ab)neg->nosoln}
If $a$ and~$b$ are elements of an algebraic number
field~$F$, and there is a real place~$\sigma$ of~$F$,
such that $\sigma(a)$ and~$\sigma(b)$ are negative, show
$(0,0,0)$ is the only solution in~$F^3$ of the equation $a
x^2 + b y^2 = z^2$.

\item \label{isotrop->unip}
In \cref{CocpctArithSO21Pf-anis} of the proof of \cref{CocpctArithSO21}, verify the assertion that $B(x,x) \neq 0$, for every nonzero $x \in F^3$.
\hint{If $B(x,x) = 0$ for some nonzero~$x$, then, after a change of basis, $B(x,x)$ is a scalar multiple of the form $x_1x_3 + x_2^2$, which is invariant under 
the unipotent transformation $x_1 \mapsto x_1 - 2x_2 - x_3$, $x_2 \mapsto x_2 + x_3$, $x_3 \mapsto x_3$.}
%$\begin{bmatrix} 1 & 1 & 0 \\ 0 & 1 & -1 \\ 0 & 0 & 1 \end{bmatrix}$.

\end{exercises}

\section{Arithmetic subgroups of \texorpdfstring{$\SL(2,\real)$}{SL(2,R)} via quaternion algebras}
\label{QuaternionAlgSect}

In the preceding section, we constructed the cocompact, arithmetic subgroups of $\SL(2,\real)$ from orthogonal groups. As an alternative approach, we will now explain what quaternion algebras are, and how they can be used to construct those same arithmetic subgroups. In later sections (and later chapters), the use of quaternion algebras will sometimes be necessary, not an alternative approach. 
%More general division algebras will also sometimes be required, so we also explain what those are.

\begin{defns} \label{QuaternionDefn} \ 
\noprelistbreak
 \begin{enumerate}
 \item For any field~$F$, and any nonzero $a,b \in
F$, the corresponding \defit{quaternion algebra} over~$F$ is
the ring%
\nindex{$\quaternion^{a,b}_F$ = quaternion algebra}% no page break here !!!
 $$ \quaternion^{a,b}_F
 = \{\, p + qi + rj + sk \mid p,q,r,s \in F \,\} ,$$
 where
 \begin{itemize}
 \item addition is defined in the obvious way, and 
 \item multiplication is determined by the relations
 $$ i^2 = a,
 \ j^2 = b,
 \ ij = k = - ji ,$$ 
 together with the requirement that every element of~$F$ is
in the center of~$D$.
(Note that $k^2 = k \cdot k = (-ji)(ij) = -a j^2 = -ab$.)
 \end{itemize}
 \item The \defit{reduced norm} of $x = p + qi + rj + sk \in
\quaternion^{a,b}_F$ is
 $$ \Nred(x) = x \, \overline x = p^2 - a q^2 - b r^2 + ab s^2 \in F,$$
 where $\overline x = p - qi - rj - sk$ is the \defit[conjugate (of a quaternion)]{conjugate} of~$x$.
 (Note that $\overline{xy} = \overline{y} \, \overline{x}$.) 
 \end{enumerate}
 \end{defns}

\begin{eg} \label{QuatFEgs} \ 
 \begin{enumerate}
 \item We have $\quaternion^{-1,-1}_\real = \quaternion$.
 \item We have $\quaternion^{t^2 a, t^2 b}_F \iso
\quaternion^{a,b}_F$ for any nonzero $a,b,t \in F$
\csee{Db2=D}.
 \item \label{QuatFEgs-squaresplit}
 We have
$\quaternion^{a^2,b}_F \iso \Mat_{2 \times 2}(F)$, for any
nonzero $a,b \in F$ \csee{Quat1=Mat}.
 \item We have $\Nred(gh) = \Nred(g) \cdot \Nred(h)$ for
$g,h \in \quaternion^{a,b}_F$.
 \end{enumerate}
 \end{eg}

\begin{lem} \label{QuatxF=}
  We have $\quaternion^{a,b}_\complex \iso
 \Mat_{2 \times 2}(\complex)$, for all $a,b \in
\complex$, and
 $$ \quaternion^{a,-1}_\real
 \iso \begin{cases}
 \Mat_{2 \times 2}(\real) & \text{if $a > 0$} , \\
 \hfil \quaternion & \text{if $a < 0$} .
 \end{cases}
 $$
 \end{lem}

\begin{proof}
 This follows from the observations in
\cref{QuatFEgs}.
 \end{proof}

\begin{prop} \label{ArithSL2RQuatQ}
 Fix positive integers $a$ and~$b$, and let
 $$G 
 = \SL \bigl( 1,\quaternion^{a,b}_\real \bigr) 
 = \{\, g \in \quaternion^{a,b}_\real \mid \Nred(g) = 1 \,\} .$$
 Then:
 \begin{enumerate}
 \item \label{ArithSL2RQuatQ-=SL2R}
 $G \iso \SL(2,\real)$,
 \item \label{ArithSL2RQuatQ-latt}
 $G_{\integer} = \SL \bigl( 1,\quaternion^{a,b}_\integer \bigr)$ is
an arithmetic subgroup of~$G$, and
 \item the following are equivalent:
 \begin{enumerate}
 \item \label{ArithSL2RQuatQ-cocpct}
 $G_{\integer}$ is cocompact in~$G$.
 \item \label{ArithSL2RQuatQ-nosoln}
 $(0,0,0,0)$ is the only integer solution $(p,q,r,s)$ of the
Diophantine equation 
	$$w^2 - a x^2 - b y^2 + ab z^2 = 0 .$$
 \item \label{ArithSL2RQuatQ-DivAlg}
 Every nonzero element of $\quaternion^{a,b}_\rational$ has a multiplicative inverse \textup(so $\quaternion^{a,b}_\rational$ is a ``\term{division algebra}''\/\textup).
 \end{enumerate}
 \end{enumerate}
 \end{prop}

\begin{proof}
 \pref{ArithSL2RQuatQ-=SL2R} Define an $\real$-linear
bijection
 $\phi \colon \quaternion^{a,b}_\real \to \Mat_{2 \times 2}(\real)$
by $\phi(1) = \Id$,
 $$ 
  \phi(i) = 
 \begin{bmatrix}
 \sqrt{a} & 0 \\
 0 & -\sqrt{a}
 \end{bmatrix}
 , \quad
  \phi(j) = 
 \begin{bmatrix}
 0 & 1 \\
 b & 0
 \end{bmatrix}
 , \quad
  \phi(k) = 
 \begin{bmatrix}
 0 & \sqrt{a} \\
 -b\sqrt{a} & 0
 \end{bmatrix}
 .$$
 It is straightforward to check that $\phi$ preserves
multiplication, so $\phi$~is a ring isomorphism.

For $g = p + q i + r j + s k \in \quaternion^{a,b}_\real$, we have
 \begin{align*}
  \det \bigl( \phi(g) \bigr) 
 &= (p + q \sqrt{a})(p - q \sqrt{a})
 - (r + s\sqrt{a})(br - bs\sqrt{a})
 \\&= p^2 - a q^2 - b r^2 + ab s^2
 \\&= \Nred(g)
  . \end{align*}
 Therefore, $\phi(G) = \SL(2,\real)$.

\pref{ArithSL2RQuatQ-latt} For $g \in G$, define $T_g
\colon \quaternion^{a,b}_\real \to \quaternion^{a,b}_\real$ by $T_g(v) = gv$.
Then $T_g$ is $\real$-linear. 
For $\gamma \in \quaternion^{a,b}_\real$, we have
$T_\gamma \bigl( \quaternion^{a,b}_\integer \bigr) \subset
\quaternion^{a,b}_\integer$ if and only if $\gamma \in
\quaternion^{a,b}_\integer$. So $G_{\integer} = G \cap
\quaternion^{a,b}_\integer$ is an arithmetic subgroup of~$G$.

(\ref{ArithSL2RQuatQ-DivAlg} $\Rightarrow$
\ref{ArithSL2RQuatQ-cocpct})
 We prove the contrapositive. Suppose $G_{\integer}$ is not
cocompact. Then the Godement Criterion \pref{GodementCriterion} tells us that it has a nontrivial unipotent
element~$\gamma$. So $1$~is an
eigenvalue of~$T_\gamma$; that is, there is some nonzero $v
\in \quaternion^{a,b}_\integer$, such that $T_\gamma(v) = v$. By
definition of~$T_\gamma$, this means $\gamma v = v$. Hence
$(\gamma - 1) v = 0$. Since $\gamma \neq 1$ and $v \neq 0$,
this implies $v$~is a zero divisor, so it certainly does not have a multiplicative inverse.

(\ref{ArithSL2RQuatQ-cocpct} $\Rightarrow$
\ref{ArithSL2RQuatQ-DivAlg}) 
 We prove the contrapositive.
 Suppose $\quaternion^{a,b}_\rational$ is not a division algebra.
Then $\quaternion^{a,b}_\rational \iso \Mat_{2 \times 2}(\rational)$
\csee{QuatZeroDiv=Mat}. This implies $\SL \bigl( 1,
\quaternion^{a,b}_\integer \bigr) \approx \SL(2,\integer)$ is not
cocompact. (It has nontrivial unipotent elements.)

(\ref{ArithSL2RQuatQ-nosoln} $\Leftrightarrow$
\ref{ArithSL2RQuatQ-DivAlg}) See \cref{QuatInverses}.
 \end{proof}

The following can be proved similarly \csee{SL1DCocpctEg}.

\begin{prop} \label{ArithSL2RQuatF}
 Let
 \noprelistbreak
 \begin{itemize}
 \item
 $F$ be a totally real algebraic number
field\/ \textup(with $F \neq \rational$\textup),
 \item $\ints$ be the ring of integers of~$F$,
 \item $a,b \in \ints$, such that $a$
and~$b$ are positive, but $\sigma(a)$ and~$\sigma(b)$ are
negative, for every place $\sigma \neq \Id$, and
 \item $G = \SL \bigl( 1,\quaternion^{a,b}_\real \bigr) $.
 \end{itemize}
 Then:
 \begin{enumerate}
 \item \label{ArithSL2RQuatF-=SL2R}
 $G \iso \SL(2,\real)$, and
 \item \label{ArithSL2RQuatF-latt}
 $G_{\ints} = \SL \bigl( 1, \quaternion^{a,b}_{\ints}
\bigr)$ is a cocompact, arithmetic subgroup of~$G$.
 \end{enumerate}
 \end{prop}

\begin{prop} \label{CocpctArithSL2R}
Every cocompact, arithmetic subgroup of\/ $\SL(2,\real)$
appears in either
\cref{ArithSL2RQuatQ} or\/~\ref{ArithSL2RQuatF} % are both still propositions !!!
\textup(up to commensurability and conjugates\/\textup).
 \end{prop}

\begin{proof}
 This can be proved directly, but we will instead derive it as a
corollary of \cref{CocpctArithSO21}. For each
arithmetic subgroup~$\Gamma$ of $\SO(2,1)$, constructed in
\cref{ArithSO21Eg}, we find an isogeny $\phi \colon \SL(2,\real) \to
\SO(2,1)$, such that $\phi(\Gamma')$ is commensurable
to an arithmetic
subgroup constructed in
\cref{ArithSL2RQuatQ} or~\ref{ArithSL2RQuatF}.

\pref{ArithSO21Eg-Q} First, let us show that every arithmetic subgroup
of type~\fullref{ArithSO21Eg}{Q} appears in
\pref{ArithSL2RQuatQ}. Given positive integers $a$ and~$b$,
such that $(0,0,0)$ is the only rational solution of the
equation $a x^2 + b y^2 = z^2$, let
 $$G 
 = \SL \bigl( 1,\quaternion^{a,b}_\real \bigr) \iso \SL(2,\real) .$$
 One can show that $(0,0,0,0)$ is the only rational
solution of the equation $w^2 - a x^2 - b y^2 + a b z^2 = 0$
\csee{ZeroDiv(k=0)}, so $G_{\integer}$ is a cocompact,
arithmetic subgroup of~$G$ \csee{ArithSL2RQuatQ}.

As a subspace of $\quaternion^{a,b}_\real$, the Lie algebra~$\Lie G$
of~$G$ is
 $$ \Lie G = \{\, v \in \quaternion^{a,b}_\real \mid \Re v = 0 \,\}
$$
 \csee{T1(Nred=1)}. For $g \in G$ and $v \in \Lie G$, we
have $(\Ad_G g)(v) = g v g^{-1}$, so $\Nred|_{\Lie G}$ is a
quadratic form on~$\Lie G$ that is invariant under $\Ad
G_F$. For $v = x i + y j + z k \in \Lie G$, we have
 $$ \Nred(v) = -a x^2 - b y^2 + ab z^2 .$$
 After the change of variables $x \mapsto b y$ and $y
\mapsto a x$, this becomes $-ab(a x^2 + b y^2 - z^2)$,
which is a scalar multiple of the quadratic form in
\fullref{ArithSO21Eg}{Q}. Therefore, after identifying $\Lie G$
with~$\real^3$ by an appropriate choice of basis, the
arithmetic subgroup constructed in \fullcref{ArithSO21Eg}{Q} (for the
given values of~$a$ and~$b$) is commensurable to $\Ad_G
G_{\integer}$.

\pref{ArithSO21Eg-F} Similarly, every arithmetic subgroup of
type~\fullref{ArithSO21Eg}{F} appears in
\pref{ArithSL2RQuatF} \csee{ArithSO21->QuatF}.
 \end{proof}

\begin{exercises}

\item \label{Db2=D}
 Show $\quaternion^{u^2a,v^2\gamma}_F \iso
\quaternion^{a,b}_F$, for any nonzero $u,v \in F$.
 \hint{An isomorphism is given by $1 \mapsto 1$, $i \mapsto
ui$, $j \mapsto vj$, $k \mapsto uvk$.}

\item \label{Quat1=Mat}
 Show $\quaternion^{a^2,b}_F \iso \Mat_{2 \times 2}(F)$, for any
field~$F$, and any $a,b \in F$.
 \hint{See the proof of \fullcref{ArithSL2RQuatQ}{=SL2R}.}

\item \label{QuatZeroDiv=Mat}
 Show that if the ring $\quaternion^{a,b}_\rational$ is not a division
algebra, then it is isomorphic to $\Mat_{2 \times 2}(\rational)$.
 \hint{This follows from {Wedderburn's Theorem} \pref{WedderburnThm}, but can also be proved directly: if $x$ is not invertible, then $xy = 0$ for some~$y$, so the left ideal generated by~$x$ is a $2$-dimensional subspace on which $\quaternion^{a,b}_\rational$ acts faithfully.}
 
\item \label{QuatInverses}
 Show that every nonzero element of $\quaternion^{a,b}_F$ has a
multiplicative inverse if and only if the reduced norm of every nonzero element is nonzero.
 \hint{If $\Nred(x) \neq 0$, then multiply the conjugate
of~$x$ by an element of~$F$ to obtain a multiplicative
inverse of~$x$. If $\Nred(x) = 0$, then $x$~is a zero
divisor.}

\item \label{SL1DCocpctEg}
 For $G$, $F$, $\ints$, $a$, and~$b$~as in
\cref{ArithSL2RQuatF}, show:
 \begin{enumerate}
 \item $G \iso \SL(2,\real)$,
 \item $G_{\ints}$ is an arithmetic subgroup of~$G$,
 \item if $g \in \quaternion^{a,b}_F$ with $\Nred(g) = 0$, then $g =
0$, and
 \item $G_{\ints}$ is cocompact in~$G$.
 \end{enumerate}

\item \label{ZeroDiv(k=0)}
 Let $a$ and~$b$ be nonzero elements of a field~$F$.
Show that if there is a nonzero solution of the equation
$w^2 - a x^2 - b y^2 + a b z^2 = 0$, then there is a
nonzero solution of the equation $w^2 - a x^2 - b y^2 = 0$.
 \hint{By assumption, there is a nonzero element~$g$ of
$\quaternion^{a,b}_F$, such that $\Nred(g) = 0$. There is some
nonzero $\alpha \in F + Fi$, such that the $k$-component
of~$\alpha g$ is zero.}

\item \label{T1(Nred=1)}
 For $a,b \in \real$, the set
 $$ G = \{\, g \in \quaternion^{a,b}_\real \mid \Nred(g) = 1 \,\} $$ 
 is a submanifold of $\quaternion^{a,b}_\real$. Show that the
tangent space $T_1 G$ is
 $$ \{\, v \in \quaternion^{a,b}_\real \mid \Re v = 0 \,\} .$$
 \hint{$T_1 G$ is the kernel of the derivative
$d(\Nred)_1$.}

\item \label{ArithSO21->QuatF}
 Carry out Part~\pref{ArithSO21Eg-F} of the proof of \cref{CocpctArithSL2R}.

\end{exercises}

\section{Arithmetic subgroups of \texorpdfstring{$\SL(2,\real)$}{SL(2,R)} via unitary groups}

Unitary groups provide yet another construction of the cocompact, arithmetic subgroups of $\SL(2,\real)$. In later sections (and later chapters), they will join quaternion algebras as another essential tool, not an alternative approach.

In fact, unitary groups can be applied in two different ways.
The simpler of the two approaches is based on the fact that $\SL(2,\real)$ is isomorphic to $\SU(1,1)$ \csee{CocpctSU11}. (This is very similar to the construction in \cref{ArithLattSL2} that is based on the fact that $\SL(2,\real)$ is isogenous to $\SO(2,1)$.) However, the required isogeny has no higher-dimensional analogue, so this method will not provide any lattices in $\SL(n,\real)$ when $n > 2$.

The following method is much more important, because it will be used in later sections to construct arithmetic subgroups of $\SL(n,\real)$ for all~$n$, not just $n = 2$.

\begin{eg} \label{SU2QCocpct}
Let
	\begin{itemize}
	\item $a,b \in \rational^+$,
	\item $L = \rational \bigl[\! \sqrt{a} \bigr] \subset \real$,
	\item $\ints$ be the ring of integers of~$L$ (so $\ints \doteq \integer\bigl[\! \sqrt{a}\bigr]$),
	\item $\tau$ denote the nontrivial element of $\Gal(L/\rational)$,
	\item $A = \diag(b,-1) = \left[\begin{smallmatrix} b & 0 \\ 0 & -1 \end{smallmatrix}\right]$,
	and
	\item $G_\ints = \SU( A, \tau ; \ints ) = \{\, g \in \SL(2,\ints) \mid \tau(g^\transpose) \, A \, g = A \,\} \subset \SL(2,\real)$.
	\end{itemize}
If $ x = (0,0)$ is the only solution in~$L^2$ of the equation $\tau( x{}^\transpose) \, A  x = 0$, then $G_\ints$ is a cocompact, arithmetic subgroup of $\SL(2,\real)$.
\end{eg}

\begin{proof}
It is not at all difficult to verify that $G_\ints$ is commensurable to an arithmetic group constructed from a quaternion algebra in \cref{ArithSL2RQuatQ} \csee{SU2=Quat}, but a direct proof is more instructive.

To see that $G_\ints$ is an arithmetic subgroup, we apply restriction of scalars.
The Galois automorphism $\tau \colon L \to L$ is $\rational$-linear. Therefore, if we think of $L$ as a ($2$-dimensional) vector space over~$\rational$, then $\tau$ is a polynomial with $\rational$-coefficients (with respect to any basis of~$L$ over~$\rational$). Since matrix multiplication and transpose are also defined by polynomial functions, this implies that if we write $g = X + \sqrt{a} \, Y$, where $X,Y \in \Mat_{2 \times 2}(\rational)$, then the equation $\tau(g^\transpose) \, A \, g = A$ is a system of polynomial equations with $\rational$-coefficients, in terms of the matrix entries of~$X$ and~$Y$. Therefore, it determines a group that is defined over~$\rational$. 
More precisely, letting $G = \SL(2,\real)$, define:
	\begin{itemize}
	\item $\Delta \colon L \to L^2$ by $\Delta(s) = \bigl( s, A \, \tau(s) \bigr)$, so $\Zlatt = \Delta(\ints)$ is a $\integer$-lattice in~$\real^2$,
	and
	\item $\phi \colon G \to G \times G$ by $\phi(g) = \bigl( g, (g^\transpose)^{-1} \bigr)$. 
	\end{itemize}
The import of the above argument is that $\phi(G)$ is defined over~$\rational$, with respect to the $\rational$-form $\Delta(L)$ of~$\real^2$. Since it is not difficult to verify that $G_\ints = \rho^{-1} \bigl( \rho(G)_{\Delta(\ints)} \bigr)$, we see that $G_\ints$ is an arithmetic subgroup of~$G$.

If $G_\ints$ is not cocompact, then it has a nontrivial unipotent element~$u$, so there exist nonzero $ x,  y \in L^2$, such that $u  x =  x$ and $u  y =  x + y$. 
Define $B \colon L^2 \times L^2 \to L$ by $B( x_1 ,  x_2) = \tau( x{}_1^\transpose) \, A  x_2$. Since $u \in G_\ints$, the definition of $G_\ints$ implies
	$$ B( x,  y) = B(u x, u  y) = B( x, x + y) = B( x,  x) + B( x,  y) .$$
Therefore $B( x, x ) = 0$. By assumption, this contradicts the fact that $x \neq  0$.
\end{proof}

\begin{eg} \label{SU2TotReal}
The preceding % @@@
\lcnamecref{SU2QCocpct} can be modified, much as in \cref{ArithSL2RQuatF}, to obtain all of the other cocompact lattices in $\SL(2,\real)$. Namely, replace $\rational$ with a totally real number field $F \neq \rational$, and let:
	\begin{itemize}
	\item $a,b \in F^+$, such that $\sigma(a) < 0$ and $\sigma(b) < 0$, for all nonidentity places of~$F$,
	\item $L = F \bigl[\! \sqrt{a} \bigr] \subset \real$,
	and
	\item $\ints$,  $\tau$, $A$, and $G_\ints$ be defined as in \cref{SU2QCocpct}. 
	\end{itemize}
Then $G_\ints$ is a cocompact, arithmetic subgroup of $\SL(2,\real)$.
\end{eg}

\begin{proof}
Let $\ints_F$ be the ring of integers of~$F$. From the second paragraph % @@@
of the proof of \cref{SU2QCocpct} (with $F$ in the place of~$\rational$) we see that $G_\ints$ is the $\ints_F$-points of a certain $F$-form~$G_F$ of $\SL(2,\real)$. Then restriction of scalars \pref{ResScal->Latt} implies that $\Delta(G_\ints)$ is an arithmetic subgroup of $\prod_{\sigma \in S^\infty} G^\sigma$. 

For any nonidentity place~$\sigma$ of~$F$, we have $\sigma(a) < 0$, so
	$$ L_\sigma = F_\sigma \Bigl[ \sqrt{\sigma \vbox to 7pt{\vss\hbox{$($}} a \vbox to 7pt{\vss\hbox{$)$}}} \Bigr] = \complex .$$
Then, since $\sigma(b)$ and~$-1$ are both negative, we have
	$$ \text{$G^\sigma = \SU \bigl( \sigma(A) , \tau_\complex ; \complex \bigr) = \SU \bigl( \diag( \sigma(b), -1 ), \tau_\complex ; \complex \bigr) \iso \SU(2)$ is compact} .$$
Therefore, all factors of $\prod_{\sigma \in S^\infty} G^\sigma$ other than~$G$ are compact, so we can mod them out, to conclude that $G_\ints$ is an arithmetic subgroup of the group $G = \SL(2,\real)$ \ccf{ArithDefn}. Furthermore, the existence of compact factors implies that the arithmetic subgroup  is cocompact \csee{scalars->cpct}.
\end{proof}

\begin{exercises}

\item \label{SU2=Quat}
Let $a,b \in \integer^+$, let $\phi \colon \quaternion_\real^{a,b} \to \Mat_{2 \times 2}(\real)$ be as in the proof of \cref{ArithSL2RQuatQ}, let $\ints = \integer[\!\sqrt{a}]$, and let $G_\ints$ be as in \cref{SU2QCocpct}. 
Show  $\phi \bigl( \SL \bigl(1, \quaternion_\integer^{a,b} \bigr) \bigr) = G_\ints$. 

\item \label{CocpctSU11}
Let
	\begin{itemize}
	\item $a, b \in \rational^+$,
	\item $L = \rational \bigl[ \sqrt{-a} \,\bigr]$,
	\item $\ints$ be the ring of integers of~$L$ (so $\ints \doteq \integer\bigl[\! \sqrt{-a}\bigr]$),
	\item $\tau$ denote complex conjugation (the only nontrivial element of $\Gal(\complex/\real)$, and also of $\Gal(L/\rational)$),
	\item $A = \diag(b,-1) = \left[\begin{smallmatrix} b & 0 \\ 0 & -1 \end{smallmatrix}\right]$,
	and
	\item $G = \SU( A, \tau; \complex ) = \{\, g \in \SL(2,\complex) \mid \tau(g^\transpose) \, A \, g = A \,\} \iso \SU(1,1)$.
	\end{itemize}
Show that if $x = (0,0)$ is the only solution in~$L^2$ of the equation $\tau( x\!{}^\transpose) A  x = 0$, then $G_\ints$ is a cocompact, arithmetic subgroup of~$G$.

%To obtain the other cocompact arithmetic subgroups of $\SU(1,1)$, we replace $\rational$ with a more general totally real number field~$F$. When $F \neq \rational$, other assumptions imply $ x^\transpose \! A  x = 0$ has no nontrivial solution, so we omit this hypothesis.
%
%\begin{prop}
%Let
%	\begin{itemize}
%	\item $F$ be a totally real algebraic number field, such that $F \not\iso \rational$,
%	\item $a, b \in F^+$,
%	\item $L = F \bigl[ \sqrt{-a} \,\bigr]$,
%	\item $\ints$ be the ring of integers of~$L$,
%	\item $\tau$ denote complex conjugation,
%	\item $A = \diag(b,-1)$,
%	and
%	\item $G = \SU( A, \tau ; L ) \iso \SU(1,1)$.
%	\end{itemize}
%Then $G_\ints$ is a cocompact, arithmetic subgroup of~$G$.
%
%Conversely, every cocompact, arithmetic subgroup of\/ $\SU(1,1)$ arises from either this construction or the construction in \cref{SUQLattSU11} \textup(up to commensurability and conjugates\textup). 
%\end{prop}

\end{exercises}

\section{Arithmetic subgroups of \texorpdfstring{$\SO(1,\lowercase{n})$}{SO(1,n)}}
\label{ArithSO1nSect}

\begin{prop} \label{NonCocptArithSOn1}
 Let
 \begin{itemize}
  \item $a_1,\ldots,a_n$ be positive integers,
  and
   \item $G = \SO(x_0^2 - a_1 x_1^2 - \cdots - a_n x_n^2;
\real) \iso \SO(1,n)$. 
 \end{itemize}
 If $n \ge 4$, then $G_{\integer}$ is an arithmetic subgroup of~$G$ that is \textbf{not} cocompact.
 \end{prop}

\begin{proof}
Since $a_1,\ldots,a_n > 0$ it is obvious that $G \iso \SO(1,n)$. Also, since $a_1,\ldots,a_n \in \rational$, it is clear that $G$ is defined over~$\rational$, so $G_\integer$ is an arithmetic subgroup of~$G$.

Since we are assuming $n \ge 4$, a theorem of Number Theory (called \emph{Meyer's Theorem}\thmindex{Meyer's}) tells us that the equation $a_1 x_1^2 + \cdots + a_n x_n^2 = x_0^2$ has a nontrivial integral solution. (This is related to, but more difficult than, the fact that every integer is a sum of four squares.) Therefore, $G_{\integer}$ is noncocompact.
\end{proof}

In most cases, the above construction is exhaustive:

\begin{prop}[\csee{SO1nNotCpctList}] \label{SO1nNotCpctListStated}
 If $n \notin \{3,7\}$, then the arithmetic subgroups constructed in
\cref{NonCocptArithSOn1} are the only noncocompact, arithmetic subgroups of\/
$\SO(1,n)$ \textup(up to commensurability and
conjugates\textup).
 \end{prop}
 
\begin{rems} \ 
\noprelistbreak
\begin{enumerate}
\item The case $n = 7$ is genuinely exceptional: there exist some exotic
arithmetic subgroups of $\SO(1,7)$ \csee{D4weird}. 
\item The groups $\SO(1,2)$ and $\SO(1,3)$ are isogenous to $\SL(2,\real)$ and $\SL(2,\complex)$, respectively. Therefore, \cref{CocpctArithSL2R,NoncocpctSL2R=SL2Z} describe all of the arithmetic subgroups of $\SO(1,2)$. Similar constructions yield the arithmetic subgroups of $\SL(2,\complex) \approx \SO(1,3)$.
%Hence, there is no real harm in assuming $n \ge 4$.
\end{enumerate}
 \end{rems}

Cocompact arithmetic subgroups of $\SO(1,n)$ can be constructed by using an algebraic extension of~$\rational$, much as in \cref{ArithSO21Eg}:

\begin{prop} \label{ArithSOn1F}
\noprelistbreak
 Let
 \begin{itemize}
 \item $F$ be an algebraic number field that is totally
real,
 \item $\ints$ be the ring of integers of~$F$,
 \item $a_1,\ldots,a_n \in \ints$, such that 
 \begin{itemize}
 \item each $a_j$ is positive, and
 \item each $\sigma(a_j)$ is negative, for every place
$\sigma \neq \Id$,
and
 \end{itemize}
 \item $G = \SO(x_0^2 - a_1 x_1^2 - \cdots - a_n x_n^2;
\real) \iso \SO(1,n)$.
 \end{itemize}
 Then $G_{\ints}$ is a cocompact, arithmetic
subgroup of~$G$.
 \end{prop}

This construction is exhaustive when $n$~is even:

\begin{prop}[\csee{SO1nCpctList}] \label{SO1nCpctListStated}
 If $n$~is even, then the arithmetic subgroups constructed in
\cref{ArithSOn1F} are the only cocompact, arithmetic subgroups of\/
$\SO(1,n)$ \textup(up to commensurability and
conjugates\textup).
 \end{prop}

\begin{rem}
 Theoretically, it is easy to tell whether two choices of
$a_1,\ldots,a_n$ give essentially the same arithmetic subgroup
\csee{Gamma=Gamma'<>BEquiv}.
\end{rem}

%In practice, the condition given by
%\cref{Gamma=Gamma'<>BEquiv} is sometimes rather
%subtle, as the following example indicates.
%
%\begin{egprop}
% Define quadratic forms $B_1(x)$ and $B_2(x)$
%on~$\rational^{n+1}$ by
% $$ B_1(x) = x_1^2 + \cdots + x_n^2 - x_{n+1}^2$$
% and 
% $$ B_2(x) = x_1^2 + \cdots + x_n^2 - 2x_{n+1}^2 ,$$
% and let
% \begin{itemize}
% \item $\Gamma_1 = \SO(B_1;\integer)$ and
% \item $\Gamma_2 = h^{-1} \SO(B_2;\integer) h$, where $h =
%\diag(1,\ldots,1,\sqrt{2}) \in \GL(n+1,\real)$.
% \end{itemize}
% There exists $g \in \Ortho(1,n)$, such that $g \Gamma_1
%g^{-1}$ is commensurable to~$\Gamma_2$ if and only if
%$n$~is even.
% \end{egprop}
%
%\begin{proof}
% ($\Leftarrow$) For 
% $$ g' = 
% \begin{bmatrix}
% 1& 1 \\
% 1 & -1 \\
% && 1 & 1 \\
% && 1 & -1 \\
% &&&& \ddots \\
% &&&&& 1 & 1 \\
% &&&&& 1 & -1 \\
% &&&&&&&2
% \end{bmatrix}
% ,$$
% we have
% $$ (g')^\transpose \diag(1,\ldots,1,-1) \, g'
% = \diag(2,2,\ldots,2,-4)
% = 2 \diag(1,1,\ldots,-2) $$
% so the desired conclusion follows from
%\cref{Gamma=Gamma'<>BEquiv}.
%
%($\Rightarrow$) This follows from
%\cref{Gamma=Gamma'(nodd)->disc} below, because $2$~is
%not a square in~$\rational$.
% \end{proof}

When $n$~is odd, we can construct additional arithmetic subgroups of $\SO(1,n)$  by using quaternion algebras. This requires a definition:

\begin{defn} \label{SU(Quat)Defn}
 Suppose $\quaternion^{a,b}_F$ is a quaternion algebra over
a field~$F$.
\noprelistbreak
 \begin{enumerate}
 \item Define $\tau_r \colon \quaternion^{a,b}_F \to
\quaternion^{a,b}_F$ by 
 $$ \tau_r(x_0 + x_1 i + x_2 j + x_3 k)
 =  x_0 + x_1 i - x_2 j + x_3 k .$$
 This is the \defit[reversion anti-involution]{reversion} anti-involution
of~$\quaternion^{a,b}_F$ \ccf{tau=Antiaut}.
 \item For $A \in \GL \bigl( m,\quaternion^{a,b}_F \bigr)$, with $\tau_r(A^\transpose) =
A$, let
 $$ \SU \bigl( A, \tau_r; \quaternion^{a,b}_F \bigr) =  \{\, g \in \SL \bigl( m,\quaternion^{a,b}_F \bigr) \mid
\tau_r(g^\transpose) A \, g = A \,\}.$$
 \end{enumerate}
 \end{defn}
 
 Now, here is the main idea of the construction:
 
 \begin{prop} \label{SU(H_Q)}
  Let 
  \noprelistbreak
  	\begin{itemize}
	\item $a,b \in \rational \smallsetminus \{0\}$, with $a > 0$,
	\item $a_1,\ldots,a_m$ be invertible elements of\/ $\quaternion^{a,b}$, such that $\tau_r(a_\ell) = a_\ell$, for each~$\ell$,
	\item $A = \diag(a_1,\ldots,a_m) \in \GL \bigl( m , \quaternion^{a,b}_\rational \bigr)$,
 \item $G = \SU \bigl( A, \tau_r; \quaternion^{a,b}_\real \bigr)$,
and
 \item $\ints$ be a $\integer$-lattice in
$\quaternion^{a,b}_F$, such that $\ints$ is also a
subring.
 \end{itemize}
 Then:
 \begin{enumerate}
 \item \label{SU(H_Q)-SO}
 $G \iso \SO(p,q)$, for some $p$ and~$q$ with $p + q = 2m$,
 and
 \item  \label{SU(H_Q)-arith}
$\SU ( A, \tau_r; \ints)$ is an arithmetic subgroup of~$G$.
 \end{enumerate}
\end{prop}

\begin{proof}
To make things a bit easier, let us assume $b < 0$ \csee{SU(H_Q)PfbPos}.
\Cref{QuatEigsGammaNeg} provides an isomorphism $\phi \colon \quaternion^{a,b}_\real \to
\Mat_{2 \times 2}(\real)$, such that:
 \begin{itemize}
 \item $\phi \bigl( \tau_r(x) \bigr) = \phi(x)^\transpose$,
for all $x \in \quaternion^{a,b}_\real$,
and
 \item $\phi(x)$ is symmetric, for all $x \in
\quaternion^{a,b}_\real$, such that $\tau_r(x) = x$.
\end{itemize}
Then $\phi(A)$ is symmetric, and $G$ is isomorphic to $\SO_{2m}\bigl( \phi(A) \bigr)$ \csee{SUQuatisSO}. This establishes~\pref{SU(H_Q)-SO}.

As a vector space over~$\real$, $\bigl( \quaternion^{a,b}_\real \bigr)^m$ is isomorphic to~$\real^{4m}$. With this identification, and considering $\bigl(\quaternion^{a,b}_\real\bigr)^m$ as a vector space over $\quaternion^{a,b}_\real$ via scalar multiplication on the right, we have
	$$ \GL \bigl( m,\quaternion^{a,b}_\real \bigr) 
	= \bigset{ g \in \GL(4m,\real) }{ 
	\begin{matrix} \text{$g(\vector x t) = (g \vector x) t$ for all} \\
	\text{$\vector x \in \bigl( \quaternion^{a,b}_\real \bigr)^m$ and $t \in \quaternion^{a,b}_\real$} 
	\end{matrix} 
	} .$$
Since $\quaternion^{a,b}_\rational$ is dense in $\quaternion^{a,b}_\real$, we may restrict $t$ to belong to $\quaternion^{a,b}_\rational$. This implies $G$ is defined over~$\rational$, with respect to the $\rational$-form $\bigl( \quaternion^{a,b}_\rational \bigr)^m$ of $\bigl( \quaternion^{a,b}_\real \bigr)^m$. For this $\rational$-form, we have $G_\integer = \SU ( A, \tau_r; \ints)$. This establishes~\pref{SU(H_Q)-arith}.
\end{proof}

\begin{rem}
Since $p + q = 2m$ must be even, the preceding \lcnamecref{SU(H_Q)} % @@@
cannot yield any arithmetic subgroups of $\SO(1,n)$ unless $1 + n$ is even, which means that $n$~is odd.
\end{rem}

\Cref{SU(H_Q)} yields an arithmetic subgroup of some $\SO(p,q)$, but not necessarily a subgroup of $\SO(1,n)$. Obtaining a particular value of~$p$ requires us to prescribe the number of positive eigenvalues of the symmetric matrix~$\phi(A)$ that appears in the proof. Since $\phi(A)$ is made from $\phi(a_1),\ldots,\phi(a_m)$, this is achieved by calculating the number $\varepsilon_{a,b}(a_\ell)$ of positive eigenvalues of each $\phi(a_\ell)$; the formula is in \cref{QuatEigs} below. % @@@

However, as in \cref{NonCocptArithSOn1}, Meyer's Theorem implies that arithmetic subgroups obtained in this way are never cocompact (unless $G$ is compact or $m \le 2$). To construct cocompact lattices, restriction of scalars is applied, as usual: choose an extension field~$F$ of~$\rational$, and arrange for $G^\sigma$ to be compact at all but one place. The outcome of these considerations is stated in \cref{CocpctSO1nQuat} below. % @@@

\begin{notation}[{(cf.~\cref{CountQuatEigsGammaNeg,CountQuatEigsGammaPos})}] \label{QuatEigs}
 Suppose
 \noprelistbreak
 \begin{itemize}
 \item $a$ and~$b$ are nonzero elements
of~$\real$, such that either $a$ or~$b$ is
positive, and
 \item $x$ is an invertible element of
$\quaternion^{a,b}_\real$, such that $\tau_r(x) = x$.
 \end{itemize}
 Write $x = p + qi + sk$, for some $p,q,s \in \real$.
 For convenience, let 
 	$$N_{a,b}(x) = x \, \overline{x} = p^2 - a q^2 + a b s^2 , $$
and note that $N_{a,b}(x) \neq 0$ (since $x$ is invertible). Define
 $$
 \varepsilon_{a,b}(x) =
 \begin{cases}
 \,1 & \mbox{if $b \, N_{a,b}(x) > 0$} , \\
 \,2 & \mbox{if $b \, N_{a,b}(x) < 0$, and} \\[-2pt] % !!!
 	& \qquad \text{either } \begin{cases}
	\text{$b < 0$ and $p > 0$, or} \\
	\text{$b > 0$ and $(a+1)q + (a - 1) s \sqrt{b} > 0$},
	\end{cases} \\[-4pt] % !!!
% \,2 & \mbox{if $b < 0$, $N_{a,b}(x) > 0$, and $p > 0$} , \\
%% \,0 & \mbox{if $b < 0$, $N_{a,b}(x) > 0$, and $p < 0$} , \\
% \,2 & \mbox{if $b > 0$, $N_{a,b}(x) < 0$, and $(a+1)q + (a - 1) s \sqrt{b} > 0$} , \\
%% \,0 & \mbox{if $b > 0$, $N_{a,b}(x) < 0$, and $(a+1)q + (a - 1) s \sqrt{b} > 0$} .
 \,0 & \mbox{otherwise} .
 \end{cases}
 $$
 \end{notation}

\begin{prop} \label{CocpctSO1nQuat}
 Let
 \noprelistbreak
 \begin{itemize}
 \item
 $F$ be a totally real algebraic number field\/ \textup(such that $F \neq \rational$\textup),
 \item $a$ and~$b$ be nonzero elements of~$F$,
such that, for each place $\sigma$ of~$F$,
either~$\sigma(a)$ or~$\sigma(b)$ is positive,
 \item $a_1,\ldots,a_m \in \quaternion^{a,b}_F$, such that
 \begin{itemize}
 \item $\tau_r(a_\ell) = a_\ell$ for each~$\ell$,
 \item $\sigma(a_\ell)$ is invertible, for each~$\ell$, and each
place~$\sigma$,
 \item $\sum_{\ell=1}^m \varepsilon_{a,b}(a_\ell) = 1$,
and
 \item $\sum_{\ell=1}^m
\varepsilon_{\sigma(a),\sigma(b)} \bigl(
\sigma(a_\ell) \bigr) \in \{0,2m\}$ for each place~$\sigma
\neq \Id$,
 \end{itemize}
 \item $\ints$ be a $\integer$-lattice in
$\quaternion^{a,b}_F$, such that $\ints$ is also a
subring, and
 \item $G = \SU \bigl( \diag(a_1,\ldots,a_m), \tau_r;
\quaternion^{a,b}_\real \bigr)^\circ$.
 \end{itemize}
 Then:
 \noprelistbreak
 \begin{enumerate}
 \item
 $G \iso \SO(1,2m-1)^\circ$, and
 \item
 $G_{\ints}$ is a cocompact, arithmetic subgroup
of~$G$.
 \end{enumerate}
 \end{prop}

\begin{prop}
 If $n \notin \{3,7\}$, then the arithmetic subgroups constructed in
\cref{ArithSOn1F,CocpctSO1nQuat} are the
only cocompact, arithmetic subgroups of\/ $\SO(n,1)$ \textup(up to
commensurability and conjugates\textup).
 \end{prop}

\Cref{D4weird} briefly explains the
need to assume $n \neq 7$.

\begin{exercises}
 
\item \label{CocpctSO}
 Use restriction of scalars \csee{RestrictScalarsSect} to
construct cocompact arithmetic subgroups $\SO(m,n)$ for
all $m$ and~$n$.

\item \label{QuadFormUnique}
 Suppose $G$ is an irreducible subgroup of
$\GL(\ell,\complex)$. (This means there is no nonzero, proper, $G$-invariant subspace of~$\complex^\ell$.) Show that if $B_1$ and~$B_2$ are (nonzero)
$G$-invariant quadratic forms on~$\complex^\ell$, then
there exists $\lambda \in \complex$, such that $B_1 = \lambda
B_2$.
\hint{Let $A_1$ and~$A_2$ be the symmetric matrices that represent~$B_1$ and~$B_2$, and write $A_2 = A_1 L$. For any $g \in G$, we have
	$A_1 L = g^\transpose A_1 L g = A_1 (g^{-1} L g)$.}

\item \label{Gamma=Gamma'<>BEquiv}
 Let 
 	\begin{itemize}
	\item $F$, $\ints$, $a_1,\ldots,a_n$, and~$G$ be as in \cref{ArithSOn1F},
	\item $\Gamma = h^{-1} G_\ints h$, where $n = \diag(1, \sqrt{a_1}, \ldots,\sqrt{a_n})$, 
	and
	\item $F'$, $\ints'$, $a_1',\ldots,a_n'$, $G'$, $\Gamma'$, and $h'$ be defined similarly.
	\end{itemize}
Show $g^{-1} \Gamma g$ is commensurable to~$\Gamma'$, for some $g \in \Ortho(1,n)$,  if and only if there exists $\lambda \in F^\times$ and $g' \in \GL(n+1,F)$, such that
 $$ (g')^\transpose \, \diag(-1,a_1,\ldots,a_n) \, g'
 = \lambda \diag(-1,a'_1,\ldots,a'_n) .$$
\hint{($\Rightarrow$) For $g' = h' g h^{-1}$, we have 
	$(g')^{-1}\SO(B;\ints)g' \subseteq \SO(B';\real)$,
% \begin{align*}
%  (g')^{-1}\SO(B;\ints)g'
% &= h'g \bigl( h^{-1} \SO(B;\ints) h \bigr) g^{-1}
%(h')^{-1} \\
% &= h' (g \Gamma g^{-1}) (h')^{-1} \\
% &\approx  h' \, \Gamma' (h')^{-1} \\
% &\subseteq h' \, \SO(n,1) (h')^{-1} \\
% &= \SO(B';\real)
% . 
% \end{align*}
so the Borel Density Theorem implies $(g')^{-1}\SO(B;\real)g' \subseteq \SO(B';\real)$. 
%we have $(g')^{-1}\SO(B;\real)g' \subseteq \SO(B \circ g'; \real)$. Since $\SO(B;\real) \otimes \complex \iso \SO(n+1;\complex)$ is irreducible on~$\complex^{n+1}$, 
Apply \cref{QuadFormUnique} with $G = (g')^{-1}\SO(B;\real)g'$.} 
%Hence, there is
%some $\lambda \in \real^\times$ with $B \circ g' = \lambda
%B'$. Since both~$B$ and~$B'$ are defined over~$F$, we must
%have $\lambda \in F^\times$.

%\item \label{BEquiv->Gamma=Gamma'Ex}
% Show, for $F$, $B$, $B'$, $\Gamma$, and~$\Gamma'$ as in
%\cref{Gamma=Gamma'<>BEquiv}, that if there is some
%nonzero $\lambda \in F$ and some $g' \in \GL(n+1,F)$, such
%that
% $$ \mbox{$B(g'x) = \lambda B'(x)$, for all $x \in
%F^{n+1}$,}$$
% then there exists $g \in \Ortho(n,1)$, such that $g^{-1}
%\Gamma g$ is commensurable to~$\Gamma'$.

\item \label{Gamma=Gamma'(nodd)->disc}
 Let $F$, $\ints$, $a_1,\ldots,a_n$,
$a'_1,\ldots,a'_n$, $\Gamma$, and~$\Gamma'$ be as in
\cref{Gamma=Gamma'<>BEquiv}.

Show that if $n$~is odd, and there exists $g \in \Ortho(1,n)$, such
that $g \Gamma g^{-1}$ is commensurable to~$\Gamma'$ then
 $$ \frac{a_1 \cdots a_n}{a'_1 \cdots a'_n} \in
(F^\times)^2 .$$
\hint{The \defit[discriminant (of quadratic
form)]{discriminant} of a quadratic form $B(x)$
on~$F^{n+1}$ is defined to be the determinant of the Gram matrix of~$B$, with respect to any
basis~$\basis$ of~$F^{n+1}$. This is not uniquely
determined by~$B$, but show that it is well-defined up to
multiplication by a nonzero square in~$F^\times$.}
%, because
% $$ \det [B]_{g\basis}
% = \det \bigl( g^\transpose [B]_{\basis} \, g \bigr)
% = (\det g)^2 \det [B]_{\basis},$$
% for any $g \in \GL(n+1,F)$.
%
%Note that $\det [\lambda B]_{\basis} = \lambda^{n+1} \det
%[B]_{\basis}$. Therefore, if $n$~is odd, then $\lambda B$ has
%the same discriminant as~$B$ (up to multiplication by a
%square).

\item \label{QuatEigsGammaNeg}
 Suppose $a$ and~$b$ are real numbers,
such that $a > 0$ and $b < 0$. Show that there is
an isomorphism $\phi \colon \quaternion^{a,b}_\real \to
\Mat_{2 \times 2}(\real)$, such that:
 \begin{enumerate}
 \item $\phi \bigl( \tau_r(x) \bigr) = \phi(x)^\transpose$,
for all $x \in \quaternion^{a,b}_\real$,
and
 \item $\phi(x)$ is symmetric, for all $x \in
\quaternion^{a,b}_\real$, such that $\tau_r(x) = x$.
 \hint{Let $\phi(i) = 
 \begin{bmatrix}
 \sqrt{a} & 0 \\
 0 & - \sqrt{a}
 \end{bmatrix}$
 and
 $\phi(j) = 
 \begin{bmatrix}
 0 & \sqrt{|b|} \\
 -\sqrt{|b|} & 0
 \end{bmatrix}
 $.}
 \end{enumerate}

\item \label{SUQuatisSO}
Assume the notation of the proof of \cref{SU(H_Q)}. Show $G$ is isomorphic to $\SO_{2m}\bigl( \phi(A) \bigr)$.
\hint{Apply the isomorphism $\phi$ to both sides of the equation $\tau_r(g^\transpose) A g = A$.}

\item \label{QuatEigsGammaPos}
 Suppose $a$ and~$b$ are nonzero real numbers,
such that $b > 0$, and let
 $ w = 
 \begin{Smallbmatrix}
 0 & 1 \\ 1 & 0 
 \end{Smallbmatrix}
 $.
 Show there is
an isomorphism $\phi \colon \quaternion^{a,b}_\real \to
\Mat_{2 \times 2}(\real)$, such that:
 \begin{enumerate}
 \item $\phi \bigl( \tau_r(x) \bigr) = w
\phi(x)^\transpose w$, for all $x \in
\quaternion^{a,b}_\real$,
and
 \item $w \phi(x)$ is symmetric, for all $x \in
\quaternion^{a,b}_\real$, such that $\tau_r(x) = x$.
 \end{enumerate}
 \hint{Let $\phi(i) = 
 \begin{bmatrix}
 0 & 1 \\
 a & 0
 \end{bmatrix}$
 and
 $\phi(j) = 
\begin{bmatrix}
 \sqrt{b} & 0 \\
 0 & - \sqrt{b}
 \end{bmatrix}
 $.}

\item \label{SU(H_Q)PfbPos}
Prove \cref{SU(H_Q)} under the additional assumption that $b > 0$.
\hint{Use \cref{QuatEigsGammaPos} and show $G \iso \SO_{2m} \bigl( \diag( w \phi(a_1),\ldots,w \phi(a_m) \bigr)$.}

\item \label{CountQuatEigsGammaNeg}
In the situation of \cref{QuatEigsGammaNeg}, show that $\phi$ can be chosen so that if $x = p + qi + rj +  sk$ is an invertible element of
$\quaternion^{a,b}_\real$, and $\tau_r(x) = x$ (so $r = 0$), then
the number of positive eigenvalues of~$\phi(x)$ is
 $$ \begin{cases}
 \,1 & \mbox{if $N_{a,b}(x) < 0$,} \\
 \,2 & \mbox{if $N_{a,b}(x) > 0$ and $p > 0$,} \\
\, 0 & \mbox{otherwise} 
 . \end{cases}
 $$
\hint{Since both eigenvalues of $\phi(x)$ are real (and nonzero), the number of positive eigenvalues is determined by the determinant and trace.}

\item \label{CountQuatEigsGammaPos}
In the situation of \cref{QuatEigsGammaPos}, show that $\phi$ can be chosen so that if $x = p + qi + rj +  sk$ is an invertible element of
$\quaternion^{a,b}_\real$, and $\tau_r(x) = x$ (so $r = 0$), then
the number of positive eigenvalues of $w \phi(x)$ is
 $$ \begin{cases}
 \,1 & \mbox{if $N_{a,b}(x) > 0$} , \\
 \,2 & \mbox{if $N_{a,b}(x) < 0$ and $(a+1)q + (a - 1) s \sqrt{b} > 0$} , \\
 \,0 & \mbox{otherwise}
 . \end{cases}
 $$
\hint{See the hint to \cref{CountQuatEigsGammaNeg}.}

%\item Suppose $Q$ is a quadratic form on~$\real^\ell$
%that is defined over~$\rational$. Let $G = \SO(Q; \real)$ and
%$\Gamma = G_{\integer}$. Prove directly (without the
%Jacobson-Morosov Lemma) that if $\Gamma$ has a nontrivial
%unipotent element, then $Q$~is isotropic over~$\rational$.
% \hint{If $\gamma v \neq v$, then $\gamma^n v / \bigl(
%\|\gamma^n v \| \bigr)$ converges to an isotropic vector.}

\end{exercises}

\section{Some nonarithmetic lattices in \texorpdfstring{$\SO(1,\lowercase{n})$}{SO(1,n)}}
 \label{NonArithSO1nSect}

\Cref{ArithSO1nSect} describes algebraic methods to
construct all of the arithmetic lattices in $\SO(1,n)$
(when $n \neq 7$). We now present a geometric method that
is sometimes able to produce a new lattice by combining
two known lattices. The result is often nonarithmetic.
We assume some familiarity with hyperbolic geometry.

\subsection{Hyperbolic manifolds}

For geometric purposes, it is more convenient to
consider the locally symmetric space $\Gamma \backslash
\hyperbolic^n$, instead of the lattice~$\Gamma$. 

\begin{defn}
 A connected, Riemannian $n$-manifold~$M$ is
\defit[hyperbolic!manifold]{hyperbolic} if
 \begin{enumerate}
 \item \label{HypMfldDefn-Hn}
 $M$~is locally isometric to~$\hyperbolic^n$ (that
is, each point of~$M$ has a neighborhood that is isometric
to an open set in~$\hyperbolic^n$),
 \item $M$~is complete, and
 \item $M$~is orientable.
 \end{enumerate}
 \end{defn}

\begin{terminology}
 Many authors do not require $M$ to be complete or
orientable. Our requirement~\pref{HypMfldDefn-Hn} is
equivalent to the assertion that $M$~has constant sectional
curvature~$-1$; some authors relax this to require the
sectional curvature to be a negative constant, but do not
require it to be normalized to~$-1$.
 \end{terminology}

\begin{notation}
 Let 
 \nindex{$\PO(1,n)$ = $\Ortho(1,n)/ \{\pm \Id\}$}%
 $\PO(1,n) = \Ortho(1,n)/ \{\pm \Id\}$.

Note that:
 \begin{itemize}
 \item $\PO(1,n)$ is isogenous to $\SO(1,n)$,
 \item $\PO(1,n) \iso \Isom(\hyperbolic^n)$, and
 \item $\PO(1,n)$ has two connected components (one
component consists of orientation-preserving isometries
of~$\hyperbolic^n$, and the other consists of
orientation-reversing isometries).
 \end{itemize}
 \end{notation}

The following observation is easy to prove
\csee{HyperMfld<>latticePfEx}.

\begin{prop} \label{HyperMfld<>lattice}
 A connected Riemannian manifold~$M$ of finite volume is
hyperbolic if and only if there is a torsion-free
lattice~$\Gamma$ in $\PO(1,n)^\circ$, such that $M$~is
isometric to $\Gamma \backslash \hyperbolic^n$.
 \end{prop}

\subsection{Hybrid manifolds and totally geodesic
hypersurfaces}

\begin{figure}[ht]
 \includegraphics{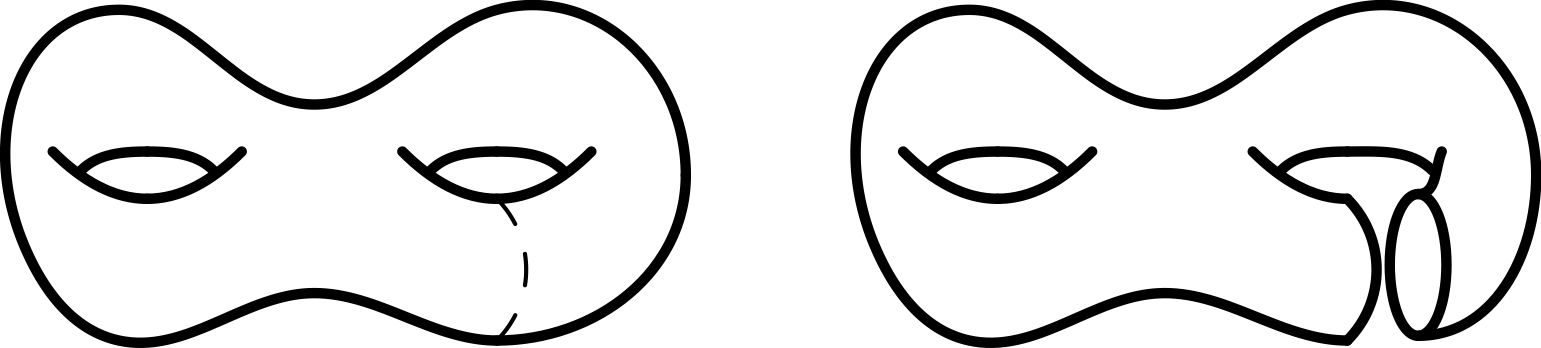}
 \caption{Cutting open a manifold by slicing along a closed
hypersurface (dashed) results in a manifold with boundary.}
 \label{cutopen}
 \end{figure}

We wish to combine two (arithmetic) hyperbolic manifolds
$M_1$ and~$M_2$ into a single hyperbolic manifold. The idea
is that we will choose closed hypersurfaces
$C_1$ and~$C_2$ of $M_1$ and~$M_2$, respectively, such that
$C_1$ is isometric to~$C_2$. Let $M'_j$ be the manifold with
boundary that results from cutting $M_j$ open, by slicing
along~$C_j$ (see \cref{cutopen} and
\cref{NonOrientDblCover}).

%texpreamble
%("  \usepackage{amsmath}
% \usepackage[LY1]{fontenc}
% \usepackage[expert,LY1,mylucidascale]{mylucidabr}
% ");
%defaultpen(  fontcommand("\normalfont") + fontsize(10) ); 
%
%from graph access *;
%unitsize(0.4cm);
%
%pen thickpen = linewidth(1.25);
%
%currentpen = thickpen;
%
%void lefthalf(real x) {
%	draw( (x,0.5){W}..{NW}(x-1, 1) );
%	draw( (x,1){W}..{SW}(x-0.7, 0.8) );
%	}
%void righthalf(real x) {
%	draw( (x,0.5){E}..{NE}(x+1, 1) );
%	draw( (x,1){E}..{SE}(x+0.7, 0.8) );
%	}
%void hole(real x) {lefthalf(x); righthalf(x);}
%
%void M(real m){
%	draw( (m,-1){W}..(m-2,-0.5)..(m-4,-1)..(m-5,0)..(m-4,2.5)..(m-2,1.5)..(m,2.5)..{S}(m+2,0.75) );
%	hole(m-3.7);
%	}
%
%M(0); hole(0); draw( (2,0.75){S}..{W}(0,-1) ); 
%draw( (0,-1){NE}..(0,0.5){NW}, dashed );
%
%M(9); 
%draw( (9,-1){NE}..(9,0.5){NW} );
%draw( ellipse( (9.75,-0.2), 0.3, 0.75 ) );
%lefthalf(9);
%draw( (9,1){E}..{SE}(9.9, 0.8) );
%draw( (9.75,0.55){E}..{NNE}(10, 1) );
%draw( (11,0.75){S}..{W}(9.75, -0.95) );

The boundary of~$M'_1$
(namely, two copies of~$C_1$) is isometric to the boundary
of~$M'_2$ (namely, two copies of~$C_2$)
\csee{NonOrientDblCover}. So we may glue $M'_1$ and $M'_2$
together, by identifying $\bdry M'_1$ with~$\bdry M'_2$
\csee{hybridfig}, as described in the following
well-known proposition.

\begin{figure}[ht]
 \includegraphics{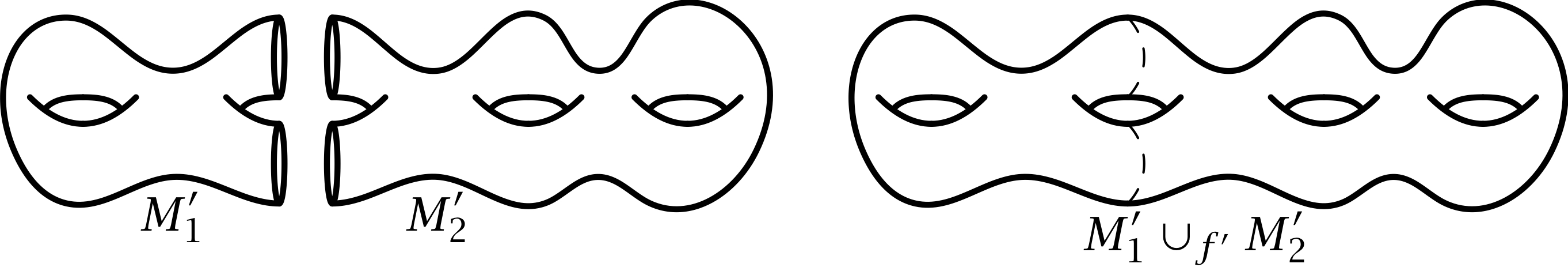}
 \caption{Gluing $M_1'$ to $M_2'$ along their boundaries
results in a manifold without boundary.}
 \label{hybridfig}
 \end{figure}
%texpreamble
%("  \usepackage{amsmath}
% \usepackage[LY1]{fontenc}
% \usepackage[expert,LY1,mylucidascale]{mylucidabr}
% ");
%defaultpen(  fontcommand("\normalfont") + fontsize(10) ); 
%
%from graph access *;
%unitsize(0.39cm);
%
%pen thickpen = linewidth(1.25);
%
%currentpen = thickpen;
%
%void lefthalf(real x) {
%	draw( (x,0.5){W}..{NW}(x-1, 1) );
%	draw( (x,1){W}..{SW}(x-0.7, 0.8) );
%	}
%void righthalf(real x) {
%	draw( (x,0.5){E}..{NE}(x+1, 1) );
%	draw( (x,1){E}..{SE}(x+0.7, 0.8) );
%	}
%void hole(real x) {lefthalf(x); righthalf(x);}
%
%void M(real m){
%	draw( (m,-1){W}..(m-2,-0.5)..(m-4,-1)..(m-5,0)..(m-4,2.5)..(m-2,1.5)..{E}(m,2.5) );
%	lefthalf(m);
%	hole(m-3.7);
%	}
%
%M(0);
%
%draw( ellipse( (0,-0.25), 0.1, 0.75 ) );
%draw( ellipse( (0,1.75), 0.1, 0.75 ) );
%
%void N(real n){
%	draw( (n,-1){E}..(n+2,-0.5)..(n+4,-1)..(n+5,-0.5)..(n+6,-1)..(n+8,0)
%		..(n+6,2.5)..(n+5,1.5)..(n+4,2.5)..(n+2,1.5)..{W}(n,2.5));
%	righthalf(n);
%	hole(n + 3.7);
%	hole(n + 6.7);
%	}
%N(1);
%draw( ellipse( (1,-0.25), 0.1, 0.75 ) );
%draw( ellipse( (1,1.75), 0.1, 0.75 ) );
%
%M(16); N(16);
%draw( (16,-1){NE}..{NW}(16,0.5) , dashed );
%draw( (16,1){NE}..{NW}(16,2.5) , dashed );
%
%label( "$M_1'$", (-2,-1.25) ); 
%label( "$M_2'$", (3,-1.25) ); 
%label( "$M_1' \cup_{f'} M_2'$", (17.3,-1.65) ); 

\begin{prop}
 Suppose 
 \begin{itemize}
 \item $M'_1$ and~$M'_2$ are connected $n$-manifolds with
boundary, and
 \item $f' \colon \bdry M'_1 \to \bdry M'_2$ is any
homeomorphism.
 \end{itemize}
 Define a topological space $M'_1 \cup_{f'} M'_2$, by gluing
$M'_1$ to~$M'_2$ along their boundaries:
 \begin{itemize}
 \item let $M'_1 \disjunion M'_2$ be the disjoint union
of~$M'_1$ and~$M'_2$,
 \item define an equivalence relation on $M'_1 \disjunion
M'_2$ by specifying that we have $m \sim f'(m)$, for every $m \in
\bdry M'_1$, and
 \item let $M'_1 \cup_{f'} M'_2 = (M'_1 \disjunion
M'_2)/\mathord{\sim}$ be the quotient of $M'_1 \disjunion
M'_2$ by this equivalence relation.
 \end{itemize}
 Then $M'_1 \cup_{f'} M'_2$ is an $n$-manifold
\textup(without boundary\textup).
 \end{prop}

\begin{cor}
 Suppose 
 \begin{itemize}
 \item $M_1$ and~$M_2$ are connected, orientable
$n$-manifolds,
 \item $C_j$ is a closed $(n-1)$-submanifold of~$M_j$, and
 \item $f \colon C_1 \to C_2$ is any homeomorphism.
 \end{itemize}
 Define $M_1 \#_f M_2 = M'_1 \cup_{f'} M'_2$, where
 \begin{itemize}
 \item $M'_j$ is the manifold with boundary that is obtained by
slicing $M_j$ open along~$C_j$, and
 \item $f' \colon \bdry M'_1 \to \bdry M'_2$ is
defined by $f'(c,k) = \bigl( f(c), k \bigr)$, under a
natural identification of~$\bdry M'_j$ with $C_j
\times \{1,2\}$.
 \end{itemize}
 Then $M_1 \#_f M_2$ is a \textup(connected\textup)
$n$-manifold \textup(without boundary\textup).

Furthermore, 
 \begin{enumerate}
 \item $M_1 \#_f M_2$ is compact if and only if both
$M_1$ and~$M_2$ are compact, and
 \item $M_1 \#_f M_2$ is connected if and only if either
$M_1 \smallsetminus C_1$ or~$M_2 \smallsetminus C_2$ is
connected.
 \end{enumerate}
 \end{cor}

\begin{terminology}
 Gromov and Piatetski-Shapiro \cite{GromovPS} call the
manifold $M_1 \#_f M_2$ a \defit[hybrid (of hyperbolic
manifolds)]{hybrid} of~$M_1$ and~$M_2$, and they call this
construction \defit[interbreeding hyperbolic
manifolds]{interbreeding}.
 \end{terminology}

Unfortunately, gluing two Riemannian manifolds together does
not always result in a Riemannian manifold (in any natural
way), even if the gluing map~$f$ is an isometry from~$\bdry
M'_1$ to~$\bdry M'_2$.

\begin{eg}
 Let $M'_1 = M'_2$ be the closed unit disk in~$\real^2$, and
let $f \colon \bdry M'_1 \to \bdry M'_2$ be the identity
map. Then $M'_1 \cup_f M'_2$ is homeomorphic to the
$2$-sphere~$S^2$. The Riemannian metrics on~$M'_1$
and~$M'_2$ are flat, so the resulting Riemannian metric
on~$S^2$ would also be flat. However, there is no flat
Riemannian metric on~$S^2$. (This follows, for example,
from the Gauss-Bonnet Theorem.) 
 \end{eg}

We can eliminate this problem by putting a restriction on
the hypersurface~$C_j$.

\begin{defn}
Let $M$ be a hyperbolic $n$-manifold. 
 A \defit[totally!geodesic hypersurface]{totally geodesic hypersurface} in~$M$ is a (closed, nonempty) connected
submanifold~$C$ of~$M$, such that, for each point~$c$
of~$C$, there are
 \begin{itemize}
 \item a neighborhood~$U$ of~$c$ in~$M$,
 \item a point~$x$ in $\hyperbolic^{n-1} = \{\, v \in
\hyperbolic^n \mid v_1 = 0 \,\}$,
 \item a neighborhood~$V$ of~$x$ in~$\hyperbolic^n$, and
 \item a Riemannian isometry $g \colon U \to V$, such that
$g(U \cap C) = V \cap \hyperbolic^{n-1}$.
 \end{itemize}
 \end{defn}

\begin{rem} \label{totgeod=Hn-1}
 If $C$ is a totally geodesic hypersurface in a hyperbolic
$n$-manifold of finite volume, then there are
 \begin{itemize}
 \item a lattice~$\Gamma$ in $\PO(1,n)$, and
 \item an isometry $f \colon M \to \Gamma \backslash
\hyperbolic^n$,
 \end{itemize}
 such that $f(C)$ is the image of $\hyperbolic^{n-1}$ in
$\Gamma \backslash \hyperbolic^n$.
 \end{rem}

\begin{prop} \label{M1cupM2hyper}
 If
 \begin{itemize}
 \item $M_1$ and~$M_2$ are hyperbolic $n$-manifolds,
 \item $C_j$ is a totally geodesic hypersurface in~$M_j$,
 \item $f \colon C_1 \to C_2$ is a Riemannian isometry, 
 and
 \item $M_1$ and~$M_2$ have finite volume,
 \end{itemize}
 then $M_1 \#_f M_2$ is a hyperbolic $n$-manifold of finite
volume.
 \end{prop}

\begin{proof}
 The main issue is to show that each point of $\bdry M'_1$
has a neighborhood~$U$ in $M'_1 \cup_{f'} M'_2$, such that
$U$~is isometric to an open subset of~$\hyperbolic^n$. This
is not difficult \csee{M1cupM2hyperEx}.

We have $\vol(M_1 \#_f M_2) = \vol(M_1) + \vol(M_2) <
\infty$.

If $M_1 \#_f M_2$ is compact, then it is obviously
complete. More generally, since $M'_1$ and~$M'_2$ are
complete, and their union is all of $M'_1 \cup_f M'_2$, it
seems rather obvious that every Cauchy sequence in $M'_1
\cup_f M'_2$ has a convergent subsequence. Hence, it
seems to be more-or-less obvious that $M'_1 \cup_f M'_2$ is
complete.

Unfortunately, if $M_1 \#_f M_2$ is not compact, then there
is a technical difficulty arising from the possibility
that, theoretically, the Riemannian isometry~$f$ may not be
an isometry with respect to the topological metrics that
$C_1$ and~$C_2$ inherit as submanifolds of $M_1$ and~$M_2$,
respectively. 
We will ignore this issue.
 \end{proof}

The following \lcnamecref{Hn-1=totgeod} describes how we will construct the
totally geodesic hypersurface~$C_j$.

\begin{lem} \label{Hn-1=totgeod}
 Suppose
 \begin{itemize}
 \item $\Gamma$ is a torsion-free lattice in\/
$\PO(1,n)^\circ$,
 \item $C$ is the image of\/~$\hyperbolic^{n-1}$ in\/ $\Gamma
\backslash \hyperbolic^n$,
 \item $\tau \colon \hyperbolic^n \to \hyperbolic^n$ is the
reflection across\/ $\hyperbolic^{n-1}$, so
	$$\tau(v_0,v_1,\ldots,v_n) =
(v_0,v_1,\ldots, v_{n-1}, -v_n) ,$$
 \item $\Gamma \cap \PO(1,n-1)$ is a lattice in\/
$\PO(1,n-1)$,
and
 \item $\Gamma$ is contained in a torsion-free lattice\/
$\Gamma'$ of\/ $\PO(1,n)^\circ$, such that\/ $\Gamma'$ is
normalized by~$\tau$.
 \end{itemize}
 Then $C$ is a totally geodesic hypersurface in\/ $\Gamma
\backslash \hyperbolic^n$, and $C$~has finite volume
\textup(as an $(n-1)$-manifold\textup).
 \end{lem}

\begin{proof}
 It is clear, from the definition of~$C$, that we need only
show $C$~is a (closed, embedded) submanifold of $\Gamma
\backslash \hyperbolic^n$.

Let $\Gamma_0 = \{\, \gamma \in \Gamma \mid
\gamma(\hyperbolic^{n-1}) = \hyperbolic^{n-1} \,\}$. (Then
$\Gamma \cap \PO(1,n-1)$ is a subgroup of index at most two
in~$\Gamma_0$.) 
 The natural map
 $$ \phi \colon \Gamma_0 \backslash \hyperbolic^{n-1} \to
\Gamma \backslash \hyperbolic^n$$
 is proper \ccf{LattInH->ProperMap}, so~$C$, being the
image of~$\phi$, is closed.

Because $\phi$ is obviously an immersion (and is a proper
map), all that remains is to show that $\phi$~is injective.
This follows from the assumption on~$\Gamma'$
\csee{Hn-1=totgeodPfinj}.
 \end{proof}

\subsection{Construction of nonarithmetic lattices}

The following theorem is the key to the construction of
nonarithmetic lattices. We postpone the proof until later
in the section \csee{ArithHybrid->Cup=M1Pf,HybridMfromM'}.

\begin{defn} \label{ArithHypMfldDefn}
 A hyperbolic $n$-manifold of finite volume is
\defit[arithmetic!hyperbolic manifold]{arithmetic} if the
corresponding lattice~$\Gamma$ in $\PO(1,n)$
\csee{HyperMfld<>lattice} is arithmetic. (Note that $\Gamma$
is well-defined, up to conjugacy
\csee{HyperMfld->UniqueLatt}, so this definition is
independent of the choice of~$\Gamma$.)
 \end{defn}

\begin{thm} \label{ArithHybrid->hybrid=M1}
 Suppose
 \noprelistbreak
 \begin{itemize}
 \item $M_1$ and~$M_2$ are hyperbolic $n$-manifolds,
 \item $C_j$ is a totally geodesic hypersurface in~$M_j$,
 \item $f \colon C_1 \to C_2$ is a Riemannian isometry,
 \item $M_1$ and~$M_2$ have finite volume \textup(as
$n$-manifolds\textup),
 \item $C_1$ and~$C_2$ have finite volume \textup(as
$(n-1)$-manifolds\textup),
 and
 \item each of $M_1 \smallsetminus C_1$ and $M_2
\smallsetminus C_2$ is connected.
 \end{itemize}
 If the hyperbolic manifold $M_1 \#_f M_2$ is
arithmetic, then $M_1 \#_f M_2$ is commensurable
to~$M_1$; that is, there are
 \begin{enumerate}
 \item a finite cover~$\widetilde{M}$ of $M_1 \cup_f M_2$,
and
 \item a finite cover~$\widetilde{M_1}$ of~$M_1$,
 \end{enumerate}
 such that $\widetilde{M}$ is isometric
to~$\widetilde{M_1}$.
 \end{thm}

\begin{cor} \label{ArithHybrid->M1=M2}
 In the situation of \cref{ArithHybrid->hybrid=M1},
 if the hyperbolic manifold $M_1 \#_f M_2$ is
arithmetic, then $M_1$ is commensurable to~$M_2$.
 \end{cor}

\begin{proof}
 From \cref{ArithHybrid->hybrid=M1}, we know that $M_1
\#_f M_2$ is commensurable to~$M_1$. By interchanging $M_1$
and~$M_2$, we see that $M_1 \#_f M_2$ is also commensurable
to~$M_2$. By transitivity, $M_1$ is commensurable to~$M_2$.
 \end{proof}

\begin{cor} \label{NonarithInSO1n}
 There exist nonarithmetic lattices\/
$\Gamma_\text{\upshape cpct}$ and\/~$\Gamma_\text{\upshape
non}$ in\/ $\SO(1,n)$, such that\/ $\Gamma_\text{\upshape
cpct}$ is cocompact, and\/ $\Gamma_\text{\upshape non}$ is
not cocompact.
 \end{cor}

\begin{proof}
 We construct only $\Gamma_\text{\upshape non}$. (See
\cref{CocpctNonarithInSO1n} for the construction
of~$\Gamma_\text{\upshape cpct}$, which is similar.)

Define quadratic forms $B_1(x)$ and~$B_2(x)$
on~$\rational^{n+1}$ by
 \begin{align*}
 B_1(x) &= x_0^2 - x_1^2 - x_2^2 - \cdots - x_{n-1}^2 - x_n^2 \\
 \intertext{and}
 B_2(x) &= x_0^2 - x_1^2 - x_2^2 - \cdots - x_{n-1}^2 - 2 x_n^2 
 . 
 \end{align*}
 Let
 \begin{itemize}
 \item $\Gamma_1 \approx \SO(B_1;\integer)$,
 \item $\Gamma_2 \approx h^{-1} \SO(B_2;\integer) h$, where
	$$h = \diag(1,1,\ldots,1,\sqrt{2}) \in \GL(n+1,\real) ,$$
 \item $M_j = \Gamma_j \backslash \hyperbolic^n$,
 \item $C_j$ be the image of $\hyperbolic^{n-1}$ in~$M_j$,
 and
 \item $\widehat\Gamma_j = \Gamma_j \cap \SO(1,n-1)$.
 \end{itemize}
 Then \cref{ArithSOn1F} tells us that $\Gamma_1$ and~$\Gamma_2$ are noncocompact
(arithmetic) lattices in $\SO(1,n)$.
By passing to finite-index subgroups, we may assume
$\Gamma_1$ and~$\Gamma_2$ are torsion free
\csee{torsionfree}. Therefore, $M_1$ and~$M_2$ are hyperbolic
$n$-manifolds of finite volume \csee{HyperMfld<>lattice}.

Because $\widehat\Gamma_j \approx \SO(1,n-1;\integer)$ is a
lattice in $\SO(1,n-1)$, and $\SO(B_j;\integer)$ is
normalized by the involution~$\tau$ of
\cref{Hn-1=totgeod}, we know $C_j$~is a totally
geodesic hypersurface in~$M_j$ that has finite volume
\csee{Hn-1=totgeod}.

Let us assume that $M_1 \smallsetminus C_1$ and $M_2
\smallsetminus C_2$ are connected. (See
\cref{NonarithInSO1nPf(notconn)} for a way around this
issue, or note that this hypothesis can be achieved by
passing to finite covers of~$M_1$ and~$M_2$.)

We know that $\widehat \Gamma_1 \approx \widehat \Gamma_2$ (since
both groups are commensurable to $\SO(1,n-1;\integer)$).
By taking a little bit of care in the choice of $\Gamma_1$
and~$\Gamma_2$, we may arrange that $\widehat \Gamma_1 = \widehat
\Gamma_2$ \csee{hybrid:Gamma1=Gamma2Ex}. Then 
 $$ C_1 \iso \widehat \Gamma_1 \backslash \hyperbolic^{n-1}
 =  \widehat \Gamma_2 \backslash \hyperbolic^{n-1}
 \iso C_2 ,$$
 so there is an isometry $f \colon C_1 \to C_2$.

If $n$~is odd, then $M_1$ is not commensurable to~$M_2$
\csee{Gamma=Gamma'(nodd)->disc}, so
\cref{ArithHybrid->M1=M2} implies that $M_1 \#_f M_2$
is not arithmetic; therefore, the corresponding
lattice~$\Gamma_\text{\upshape non}$ is not arithmetic
\csee{ArithHypMfldDefn}. When $n$~is even, an additional
argument is needed; see \cref{hybrid:nEvenNotArith}.
 \end{proof}

\subsection{Proof of \cref{ArithHybrid->hybrid=M1}}
 \label{ArithHybrid->Cup=M1Pf}

Let us recall the following lemma, which was proved in
\cref{IntZarDense->EqualEx}.

\begin{lem} \label{IntZarDense->Equal}
 If 
 \begin{itemize}
 \item $G$ has no compact factors,
 \item $\Gamma_1$ and~$\Gamma_2$ are arithmetic lattices
in~$G$, and
 \item $\Gamma_1 \cap \Gamma_2$ is Zariski dense in~$G$,
 \end{itemize}
 then $\Gamma_1$ is commensurable to~$\Gamma_2$.
 \end{lem}

\begin{defn}
 Let $M'$ be a Riemmanian $n$-manifold with boundary. We say
that $M'$~is a \defit[hyperbolic!manifold!with totally
geodesic boundary]{hyperbolic manifold with totally
geodesic boundary} if 
 \begin{enumerate}
 \item $M'$~is complete,
 \item each point of $M' \smallsetminus \bdry M'$ has a
neighborhood that is isometric to an open set
in~$\hyperbolic^n$, and
 \item for each point~$p$ of~$\bdry M'$, there are
 \begin{itemize}
 \item a neighborhood~$U$ of~$p$ in~$M'$,
 \item a point~$x$ in $\hyperbolic^{n-1} = \{\, v \in
\hyperbolic^n \mid v_1 = 0 \,\}$,
 \item a neighborhood~$V$ of~$x$ in~$\hyperbolic^n$, and
 \item an isometry $g \colon U \to V^+$, where 
 $$ V^+ = \{\, v \in V \mid v_1 \ge 0 \,\} .$$
 \end{itemize}
 (Note that $g(U \cap \bdry M') = V \cap
\hyperbolic^{n-1}$.)
 \end{enumerate}
 \end{defn}

The following is a generalization of
\cref{ArithHybrid->hybrid=M1}  \csee{HybridMfromM'}.

\begin{thm} \label{ArithHybrid->Cup=M1}
 Suppose
 \begin{itemize}
 \item $M_1$ and~$M_2$ are hyperbolic $n$-manifolds,
 \item $M_j'$ is a connected, $n$-dimensional submanifold
of~$M_j$ with totally geodesic boundary,
 \item $f' \colon \bdry M_1' \to \bdry M_2'$ is an isometry,
 \item $M_1$ and~$M_2$ have finite volume \textup(as
$n$-manifolds\textup),
 \item $\bdry M_j'$ has only finitely many components,
 and
 \item $\bdry M_1'$ and~$\bdry M_2'$ have finite volume
\textup(as $(n-1)$-manifolds\textup).
 \end{itemize}
 If the hyperbolic manifold $M_1' \cup_{f'} M_2'$ is
arithmetic, then $M_1' \cup_{f'} M_2'$ is commensurable
to~$M_1$.
 \end{thm}

\begin{proof}
 \ 
 \begin{itemize}
 \item Let $M = M_1' \cup_{f'} M_2'$.
 \item Write $M = \Gamma \backslash \hyperbolic^n$, for
some torsion-free lattice~$\Gamma$ in $\PO(1,n)$.
 \item Let $\phi \colon \hyperbolic^n \to M$ be the
resulting covering map.
 \item Let $B = \phi^{-1}(\bdry M_1')$. Because $M_1'$ has
totally geodesic boundary, we know that $B$~is a union of
disjoint hyperplanes. (That is, each component of~$B$ is of
the form $g(\hyperbolic^{n-1})$, for some $g \in
\Ortho(1,n)$.)
 \item Let $V$ be the closure of some connected component
of $\hyperbolic^n \smallsetminus B$ that contains a point
of $\phi^{-1}(M_1')$.
 \item Let 
 $$\Gamma'
 = \{\, \gamma \in \Gamma \mid \gamma V = V \,\}
 = \{\, \gamma \in \Gamma \mid \text{interior}(\gamma V \cap V) \neq \emptyset \,\} $$
 \csee{VTessellation}, so $M'_1 = \phi(V) \iso \Gamma'
\backslash V$.
 \end{itemize}

By definition, $V$~is an intersection of half-spaces, so it
is (hyperbolically) convex; hence, it is simply connected.
Therefore, $V$~is the universal cover of~$M_1'$, and
$\Gamma'$ can be identified with the fundamental group
of~$M_1'$.

Since $M'_1 \subseteq M_1$, we may define
$\Gamma_1,\phi_1,B_1,V_1,\Gamma_1'$ as above, but with
$M_1$ in the place of~$M$. From the uniqueness of the
universal cover of~$M_1'$, we know that there is an isometry
$\psi \colon V \to V_1$, and an isomorphism $\psi_* \colon
\Gamma' \to \Gamma_1'$, such that $\psi(\gamma v) =
\psi_*(\gamma) \, \psi(v)$, for all $\gamma \in \Gamma'$
and $v \in V$. Since~$\psi$ extends to an isometry
of~$\hyperbolic^n$, we may assume (after replacing
$\Gamma_1$ with $\psi^{-1} \Gamma_1 \psi$) that $V = V_1$
and $\psi_* = \Id$. Hence $\Gamma' = \Gamma_1' \subset
\Gamma \cap \Gamma_1$. It suffices to show (after replacing
$\Gamma$~by a conjugate subgroup) that the Zariski closure
of~$\Gamma'$ contains $\PO(1,n)^\circ$, for then
\cref{IntZarDense->Equal} implies $\Gamma$ is
commensurable to~$\Gamma_1$.

\begin{claim}
 We may assume that the Zariski closure of\/~$\Gamma'$
contains $\PO(1,n)^\circ$.
 \end{claim}
 We may assume $\hyperbolic^{n-1}$ is one of the connected
components of~$\bdry V$. Since $\bdry M_1'$ has finite
volume, this means that
 \begin{equation} \label{ArithHybridPf-BdryLatt}
 \mbox{$\Gamma' \cap \SO(1,n-1)$ is a lattice in
$\PO(1,n-1)$.}
 \end{equation}
 Let $\closure{\Gamma'}$~be the Zariski closure
of~$\Gamma'$. From \pref{ArithHybridPf-BdryLatt} and the
Borel Density Theorem \pref{BDT(Zardense)}, we know that
$\closure{\Gamma'}$~contains $\PO(1,n-1)^\circ$. Then,
since $\PO(1,n-1)^\circ$ is a maximal connected subgroup of
$\PO(1,n)$ \csee{MaxlInPO1n}, we may assume that
$\closure{\Gamma'}^\circ = \PO(1,n-1)^\circ$. (Otherwise,
the claim holds.) Because $\closure{\Gamma'}^\circ$ has
finite index in~$\closure{\Gamma'}$ \see{Zar->AlmConn}, this
implies that $\PO(1,n-1)^\circ$ contains a finite-index
subgroup of~$\Gamma'$. In fact, 
 \begin{equation} \label{ArithHybridPf-FinInd}
 \begin{matrix}
 \text{$\{\, \gamma \in \Gamma' \mid \gamma H = H \,\}$ has
finite index in $\Gamma'$,}
\\ \text{for every connected component~$H$
of~$\bdry V$} .
\end{matrix}
 \end{equation}
 This will lead to a contradiction.

\setcounter{case}{0}

\begin{case}
 Assume $\bdry V$ is connected.
 \end{case}
 We may assume $\bdry V = \hyperbolic^{n-1}$. Then, by
passing to a finite-index subgroup, we may assume that
$\Gamma' \subset \PO(1,n-1)$ \see{ArithHybridPf-FinInd}.
 Define $g \in \Isom(\hyperbolic^n)$ by
 $$g(v_1,v_2,\ldots,v_n) = (-v_1,v_2,\ldots,v_n) .$$
 Then
 \begin{itemize}
 \item $g$~centralizes~$\Gamma'$, and
 \item $\hyperbolic^n = V \cup g(V)$.
 \end{itemize}
  Since $\Gamma' \backslash V \iso M_1'$ has finite
volume, we know that $\Gamma' \backslash g(V)$ also has
finite volume. Therefore
 $$ \Gamma' \backslash \hyperbolic^n 
 = (\Gamma' \backslash V) \cup \bigl( \Gamma' \backslash
g(V) \bigr) $$
 has finite volume, so $\Gamma'$~is a lattice in
$\PO(1,n)$. But this contradicts the Borel Density Theorem
\pref{BDT(Zardense)} (since $\Gamma' \subset \PO(1,n-1)$).

\begin{case} \label{ArithHybrid->Cup=M1PfNotConn}
 Assume $\bdry V$ is not connected.
 \end{case}
 Let $H_1$ and~$H_2$ be two distinct connected components
of~$\bdry V$. Replacing
$\Gamma'$ by a finite-index subgroup, let us assume that
each of~$H_1$ and~$H_2$ is invariant under~$\Gamma'$
\see{ArithHybridPf-FinInd}.

 To simplify the argument, let us assume that $\bdry M_1'$
is compact, rather than merely that it has finite volume.
(See \cref{HybridWithNoncpctBdry} for the general
case.) Therefore, $\Gamma' \backslash H_1$ is compact, so there
is a compact subset~$C$ of~$H_1$, such that $\Gamma' C =
H_1$. Let 
 $$ \delta = \min \{\, \dist(c, H_2) \mid c \in C\,\} > 0
.$$
 Because $\Gamma'$ acts by isometries, we have $\delta =
\dist(H_1,H_2)$. Now, since $\hyperbolic^n$ is negatively
curved, there is a unique point~$p$ in~$H_1$, such that
$\dist(p,H_2) = \delta$. The uniqueness implies that $p$~is
fixed by every element of~$\Gamma'$. Since $\Gamma$ acts
freely on~$\hyperbolic^n$ (recall that it is a group of
deck transformations), we conclude that $\Gamma'$ is
trivial. This contradicts the fact that $\Gamma' \backslash
H_1$ is compact. (Note that $H_1 \iso \hyperbolic^{n-1}$ is
not compact.)
 \end{proof}

\begin{exercises}

\item \label{HyperMfld<>latticePfEx}
 Prove \cref{HyperMfld<>lattice}.

\item \label{HyperMfld->UniqueLatt}
 Show that if $\Gamma_1$ and~$\Gamma_2$ are torsion-free
lattices in $\PO(1,n)$, such that $\Gamma_1 \backslash
\hyperbolic^n$ is isometric to $\Gamma_2 \backslash
\hyperbolic^n$, then $\Gamma_1$ is conjugate to~$\Gamma_2$.
 \hint{Any isometry $\phi \colon \Gamma_1 \backslash
\hyperbolic^n \to \Gamma_2 \backslash \hyperbolic^n$ lifts
to an isometry of~$\hyperbolic^n$.}

 \item \label{NonOrientDblCover}
 Let $C$ be a closed, connected hypersurface in an
orientable Riemannian manifold~$M$, and let $M'$ be the
manifold with boundary that results from cutting~$M$ open,
by slicing along~$C$. Show:
 \begin{enumerate}
 \item If $C$~is orientable, then the boundary of~$M$ is
two copies of~$C$.
 \item If $C$ is not orientable, then the boundary is the
orientable double cover of~$C$.
 \item If $C$ is isometric to a closed, connected
hypersurface~$C_0$ in an orientable Riemannian
manifold~$M_0$, and $M'_0$ is the manifold with boundary
that results from cutting~$M_0$ open, by slicing
along~$C_0$, then the boundary of~$M'$ is isometric to the
boundary of~$M'_0$.
 \end{enumerate}

\item \label{M1cupM2hyperEx}
 For $M_1$, $M_2$, and~$f$ as in
\cref{M1cupM2hyperEx}, show that if $p \in \bdry
M'_1$, then $p$~has a neighborhood~$U$ in $M'_1 \cup_f
M'_2$, such that $U$~is isometric to an open subset
of~$\hyperbolic^n$.
 \hint{Find a ball~$V$ around a point~$x$
in~$\hyperbolic^{n-1}$, and isometries $g_1 \colon U_1 \to
V^+$ and $g_2 \colon U_2 \to V^-$, where $U_j$~is a
neighborhood of~$p$ in~$M'_j$, with $g_1|_{\bdry M'_1}
= (g_2 \circ f)|_{\bdry M'_1}$.}

\item \label{Hn-1=totgeodPfinj}
 For $\phi \colon \Gamma_0 \backslash
\hyperbolic^{n-1} \to \Gamma \backslash \hyperbolic^n$, as defined in the proof of \cref{Hn-1=totgeod}, show that $\phi$ is
injective.
 \hint{Suppose $\gamma x = y$, for some $\gamma \in \Gamma$
and $x,y \in \hyperbolic^{n-1}$. Then $\gamma^{-1} \tau
\gamma \tau$ is an element of~$\Gamma'$ that fixes~$x$, so
it is trivial. Hence, the fixed-point set of~$\tau$ is
$\gamma$-invariant.}

\item \label{CocpctNonarithInSO1n}
 Assume $n$~is odd, and construct a cocompact,
nonarithmetic lattice $\Gamma$ in $\SO(1,n)$.
 \hint{Let $F = \rational[\sqrt{2}]$, define $B_1(x) = \sqrt{2} x_0^2 - x_1^2 - x_2^2 - \cdots -
x_{n-1}^2 - x_n^2$
 and
 $B_2(x) = \sqrt{2} x_0^2 - x_1^2 - x_2^2 - \cdots -
x_{n-1}^2 - 3x_{n+1}^2$,
 and use the proof of \cref{NonarithInSO1n}.}

\item \label{hybrid(notconn)}
 In the notation of the proof of \cref{NonarithInSO1n},
assume that $M_1 \smallsetminus C_1$ and $M_2
\smallsetminus C_2$ are \emph{not} connected; let $M'_j$ be
the closure of a component of $M_j \smallsetminus C_j$.
Show that if $f' \colon C_1 \to C_j$ is any isometry (and
$n$~is odd), then $M'_1 \cup_{f'} M'_2$ is a
\emph{nonarithmetic} hyperbolic $n$-manifold of finite
volume.

\item \label{NonarithInSO1nPf(notconn)}
 Eliminate the assumption that $M_1 \smallsetminus C_1$ and
$M_2 \smallsetminus C_2$ are connected from the proof of
\cref{NonarithInSO1n}.
 \hint{Define
 $B_3(x) = x_0^2 - x_1^2 - x_2^2 - \cdots -
x_{n-1}^2 - 3
x_n^2$. If $M_j \smallsetminus C_j$ has the same number
of components as $M_k \smallsetminus C_k$ (and $j \neq k$),
then either \cref{hybrid(notconn)} or the proof of
\cref{NonarithInSO1n} applies.}

\item \label{hybrid:Gamma1=Gamma2Ex}
 For $B_1(x)$ and $B_2(x)$ as in the proof of
\cref{NonarithInSO1n}, show that there are finite-index
subgroups~$\Gamma_1$ and~$\Gamma_2$ of $\SO(B_1;\integer)$
and $\SO(B_2;\integer)$, respectively, such that
 \begin{enumerate}
 \item $\Gamma_1$ and~$\Gamma_2$ are torsion free, and
 \item $\Gamma_1 \cap \SO(1,n-1) = \Gamma_2 \cap
\SO(1,n-1)$.
 \end{enumerate}
 \hint{Let $\Gamma_j = \Lambda \cap \SO(B_j;\integer)$,
where $\Lambda$ is a torsion-free subgroup of finite index
in $\SL(n+1,\integer)$.}

\item \label{hybrid:nEvenNotArith}
 In the notation of the proof of \cref{NonarithInSO1n},
show that if $n$~is even (and $n \ge 4$), then
$\Gamma_\text{\upshape non}$ is not arithmetic.
 \hint{If $\Gamma_\text{\upshape non}$ is arithmetic, then
its intersection with $\SO(1,n-1)$ is arithmetic in
$\SO(1,n-1)$, and $n-1$ is odd.}

\item \label{HybridMfromM'}
 Derive \cref{ArithHybrid->hybrid=M1} as a corollary of
\cref{ArithHybrid->Cup=M1}.
 \hint{Apply \cref{ArithHybrid->Cup=M1} to
$\widetilde{M_j} = M_j \#_{f_j} M_j$, where $f_j \colon C_j
\to C_j$ is the identity map. Note that $\widetilde{M_j}$
is a double cover of~$M_j$, so $\widetilde{M_j}$ is
commensurable to~$M_j$.}

\item \label{VTessellation}
 For $\Gamma$ and~$V$ as in the proof of
\cref{ArithHybrid->Cup=M1}, let $\interior V$ be the
interior of~$V$, and show, for each $\gamma \in \Gamma$,
that if $\gamma \interior V \cap \interior V \neq \emptyset$,
then $\gamma \interior V = \interior V$.

\item \label{MaxlInPO1n}
 Show that if $H$ is a connected subgroup of $\PO(1,n)$
that contains $\PO(1,n-1)^\circ$, then $H =
\PO(1,n-1)^\circ$.

\item \label{HybridWithNoncpctBdry}
 Eliminate the assumption that $\bdry M_1'$ is compact from
\cref{ArithHybrid->Cup=M1PfNotConn} of the proof of
\cref{ArithHybrid->Cup=M1}.
 \hint{The original proof applies unless $\dist(H_1,H_2) =
0$, which would mean that $H_1$ and~$H_2$ intersect at
infinity. This intersection is a single point, and it is
invariant under~$\Gamma'$, which contradicts the Zariski
density of~$\Gamma'$.}

\end{exercises}

\section{Noncocompact arithmetic subgroups of \texorpdfstring{$\SL(3,\real)$}{SL(3,R)}}
\label{NonLattinSL3Sect}

We saw in \cref{NoncocpctSL2R=SL2Z} that
$\SL(2,\integer)$ is essentially the only noncocompact,
arithmetic subgroup of $\SL(2,\real)$. So it may be
surprising that $\SL(3,\integer)$ is \emph{not} the only
one in $\SL(3,\real)$.

\begin{prop} \label{NoncocpctInSL3Eg}
 Let
 \begin{itemize}
 \item $L$ be a real quadratic extension of\/~$\rational$, so
$L = \rational\bigl[\! \sqrt{r} \bigr]$, for some square-free positive
integer $r \ge 2$,
 \item $\sigma$ be the nontrivial Galois automorphism of~$L$,
 \item $\widetilde\sigma$ be the automorphism of\/ $\Mat_{3
\times 3}(L)$ induced by applying~$\sigma$ to each entry
of a matrix,
 \item $J_3 = 
 \begin{Smallbmatrix}
 0 & 0 & 1 \\
 0 & 1 & 0 \\
 1 & 0 & 0
 \end{Smallbmatrix}$,
 and
 \item $\Gamma = \SU \bigl( J_3, \sigma; \integer \bigl[\! \sqrt{r}
\bigr] \bigr)
 = 
 \bigset{
 g \in \SL \bigl( 3, \integer \bigl[\! \sqrt{r} \bigr] \bigr)
 }{
 \widetilde\sigma(g^{\transpose}) J_3 \, g = J_3
 }$.
 \end{itemize}
 Then:
 \begin{enumerate}
 \item \label{NoncocpctInSL3Eg-latt}
 $\Gamma$ is an arithmetic subgroup of\/ $\SL(3,\real)$,
 \item \label{NoncocpctInSL3Eg-noncocpct}
 $\Gamma$ is not cocompact, and
 \item \label{NoncocpctInSL3Eg-Qrank}
 no conjugate of\/~$\Gamma$ is commensurable to\/ $\SL(3,\integer)$.
% $\Qrank(\Gamma) = 1$.
 \end{enumerate}
 \end{prop}

\begin{proof}
 \pref{NoncocpctInSL3Eg-latt} This is a special case of
\cref{GFxR}\pref{FClassicalDefn-SUSL}, but we
provide a concrete, explicit proof (using the
methods of \cref{MakeArithLattSect,RestrictScalarsSect}).

Define
\noprelistbreak
 \begin{itemize}
 \item $\Delta \colon L^3 \to \real^6$ by $\Delta(v) =
\bigl( v,J_3 \, \sigma(v) \bigr)$,
 \item $V_{\rational} = \Delta(L^3)$,
 \item $\Zlatt = \Delta \bigl( \integer[\sqrt{r}]^3
\bigr)$, and
 \item $\rho \colon \SL(3,\real) \to \SL(6,\real)$ by
	$$ \text{$\rho(A)(v,w) = \bigl( Av, (A^T)^{-1}w \bigr)$ 
	\ for $v,w \in \real^3$} .$$
 \end{itemize}
 Then 
 \begin{itemize}
 \item $V_{\rational}$ is a $\rational$-form of~$\real^6$
\fullccf{Delta(O)=VSLatt}{F},
 \item $\Zlatt$ is a $\integer$-lattice
in~$V_{\rational}$ \fullccf{Delta(O)=VSLatt}{O},  
 \item $\rho$~is a homomorphism, 
 \item $\rho \bigl( \SL(3,\real) \bigr)$ is defined
over~$\rational$ (with respect to the
$\rational$-form~$V_{\rational}$)
(see \pref{NoncocpctInSL3Eg-/Q} below), % !!!
and
 \item $\Gamma = \{\, g \in \SL(3,\real) \mid \rho(g)
\Zlatt = \Zlatt \,\} $ \ccf{End(R4)Q}.
 \end{itemize}
 Hence, \fullcref{AbstractArith}{arith} (together with
\cref{arith->latt}) implies that $\Gamma$ is an
arithmetic subgroup of $\SL(3,\real)$.

Now let us show that 
 \begin{equation} \label{NoncocpctInSL3Eg-/Q}
 \mbox{$\rho \bigl( \SL(3,\real) \bigr)$ is defined
over~$\rational$.}
 \end{equation}
 This can be verified directly, by finding appropriate
$\rational$-polynomials, but let us, instead,
show that $\rho \bigl( \SL(3,\real) \bigr)_{\rational}$ is
dense in $\rho \bigl( \SL(3,\real) \bigr)$.

Define $U_1$ as in \eqref{NoncocpctInSL3EgPf-UGamma} below, % @@@
but allowing $a,b,c$ to range over all of~$\rational$, instead
of only~$2\integer$. Then $\rho(U_1) V_{\rational} \subset
V_{\rational}$ \csee{NoncocpctInSL3Ex-U1/Q}, so we have $\rho(U_1)
\subseteq \rho \bigl( \SL(3,\real) \bigr)_{\rational}$.
Furthermore, $U_1$ is dense in
 $$ U = 
 \begin{Smallbmatrix}
 1 & * & * \\
 0 & 1 & * \\
 0 & 0 & 1 
 \end{Smallbmatrix}
 .$$
 Similarly, there is a dense subgroup~$U_2$
of~$U^\transpose$, with $\rho(U_2) \subseteq \rho \bigl(
\SL(3,\real) \bigr)_{\rational}$
\csee{NoncocpctInSL3Ex-U2}. Since $\langle U, U^\transpose
\rangle = \SL(3,\real)$, we know that $\langle U_1,U_2
\rangle$ is dense in $\SL(3,\real)$, so $\rho \bigl(
\SL(3,\real) \bigr)_{\rational}$ is dense in $\rho \bigl(
\SL(3,\real) \bigr)$. Therefore $\rho \bigl( \SL(3,\real)
\bigr)$ is defined over~$\rational$ \see{QptsDense}.

\pref{NoncocpctInSL3Eg-noncocpct}
 By calculation, one may verify, directly from the
definition of~$\Gamma$, that the subgroup
 \begin{equation} \label{NoncocpctInSL3EgPf-UGamma}
 U_\Gamma = \bigset{
 \begin{bmatrix}
 1& a + b \sqrt{r} & -(a^2 - r b^2)/2 + c \sqrt{r} \\
 0& 1                 &  -a + b \sqrt{r} \\
 0& 0                 & 1 
 \end{bmatrix}
 } { a,b,c \in 2 \integer }
 \end{equation}
 is contained in~$\Gamma$. Then, since every element
of~$U_\Gamma$ is unipotent, it is obvious that $\Gamma$ has
nontrivial unipotent elements. So the Godement Criterion
\pref{GodementCriterion} implies that $G/\Gamma$
is not compact.

\pref{NoncocpctInSL3Eg-Qrank}
We sketch a proof. Choose an element $\omega \in \integer \bigl[ \sqrt{r} \bigr]$, such that $\sigma(\omega) = 1/\omega$. Then $\diag(\omega, 1, \omega^{-1})$ is a hyperbolic element of~$\Gamma$ that normalizes the maximal unipotent subgroup~$U_\Gamma$. On the other hand, it is easy to see that if $U'_\integer$ is any subgroup of $\SL(3,\integer)$ that is commensurable to the maximal unipotent subgroup $U_\integer$, then $U'_\integer$ has finite index in $\nzer_{\SL(3,\integer)}(U'_\integer)$. Since all of the maximal unipotent subgroups of $\SL(3,\integer)$ are conjugate (up to commensurability) under $\SL(3,\rational)$, this implies that  no conjugate of~$\Gamma$ is commensurable to $\SL(3,\integer)$.

Here is a more complete argument that is based on the notion of $\rational$-rank, which will be explained in \cref{QrankChap}.
 Define a nondegenerate $\sigma$-Hermitian form $B(x,y)$
on~$L^3$ by $B(x,y) = \sigma(x^\transpose) \, J_3 \, y$.
Then $v = (1,0,0)$ is an \term{isotropic vector} for~$B$ (i.e., $B(v,v) = 0$). On the other
hand, because $B$~is nondegenerate, the dimension of the orthogonal complement
of any subspace is equal to the codimension of the subspace.
Since $L^3$~is $3$-dimensional, this implies there is no 2-dimensional subspace that consists entirely of isotropic vectors. Therefore, $\Qrank \Gamma = 1$ \fullccf{QrankEg}{SOQ}. However, we have $\Qrank \SL(3,\integer) = 2$ \fullccf{QrankEg}{SL}. Two lattices with different $\rational$-ranks cannot be conjugate. (They cannot even be abstractly commensurable.)
 \end{proof}

\begin{rems} \ 
\noprelistbreak
 \begin{enumerate}
 \item From 
% Recall that $\Qrank \bigl( \SL(3,\integer) \bigr) =
%2$ (see~\ref{QrankIntroEgs} and~\fullref{QrankEg}{SL}).
%Therefore, 
\fullcref{NoncocpctInSL3Eg}{Qrank}, we know that none
of the arithmetic subgroups in \cref{NoncocpctInSL3Eg} are
conjugate to a subgroup that is commensurable to
$\SL(3,\integer)$. 

Indeed, let $X = \SL(3,\real)/\SO(3)$ be the symmetric space
associated to $\SL(3,\real)$. \Cref{HattoriThm}
implies that if $\Gamma$ is one of the arithmetic subgroups constructed
in \cref{NoncocpctInSL3Eg}, then the geometry of the
locally symmetric space $\Gamma \backslash X$ is very
different from that of $\SL(3,\integer) \backslash X$.
Namely, $\Gamma \backslash X$ is only mildly noncompact: it
merely has cusps, which means that its asymptotic cone is a
union of finitely many rays. In contrast, the asymptotic cone
of $\SL(3,\integer) \backslash X$ is a
2-complex, not just a union of rays. Even from a distance,
$\Gamma \backslash X$ and $\SL(3,\integer) \backslash X$
look completely different.

\item Different values of~$r$ always give essentially
different arithmetic subgroups \csee{GammaSL3R(r1<>r2)}, but this is
not so obvious.
 \end{enumerate}
 \end{rems}

The classification results in
\cref{ArithClassicalChap} imply that these are the
only arithmetic subgroups of $\SL(3,\real)$ that are not cocompact:

\begin{prop}[\csee{AllNoncocptSL3R}] \label{AllNoncocptSL3RStated}
 $\SL(3,\integer)$ and the arithmetic subgroups constructed in
\cref{NoncocpctInSL3Eg} are the only noncocompact
arithmetic subgroups of\/ $\SL(3,\real)$ \textup(up to
commensurability and conjugates\textup).
 \end{prop}

\begin{exercises}

\item \label{NoncocpctInSL3Ex-U1/Q}
 For $U_1$, $\rho$, and~$V_{\rational}$ as in the proof
of~\fullcref{NoncocpctInSL3Eg}{latt}, show that $\rho(U_1)
V_{\rational} \subseteq V_{\rational}$.

\item \label{NoncocpctInSL3Ex-U2}
 In the notation of the proof of
\cref{NoncocpctInSL3Eg}, find a dense subgroup~$U_2$
of 
 $$ \begin{Smallbmatrix}
 1 & 0 & 0 \\
 * & 1 & 0 \\
 * & * & 1 
 \end{Smallbmatrix} ,$$
 such that $\rho(U_2) \subseteq \rho \bigl( \SL(3,\real)
\bigr)_{\rational}$.

\item Assume the notation of the proof of
\cref{NoncocpctInSL3Eg}, and let $G = \rho \bigl(
\SL(3,\real) \bigr)$.
 \begin{enumerate}
 \item Show that $G$ is \defit[quasisplit $\rational$-form of
$G$]{quasisplit}. That is, show that some Borel subgroup of
$G$ is defined over~$\rational$.
 \item Show that every proper parabolic $\rational$-subgroup
of~$G$ is a Borel subgroup of~$G$.
 \end{enumerate}
 \hint{Let $B$ be the group of upper-triangular matrices in
$\SL(3,\real)$. Then $B$~is a Borel subgroup of
$\SL(3,\real)$, and $\rho(B)$ is defined over~$\rational$.}

\item \label{GammaSL3R(r1<>r2)}
 Let $\Gamma_1$ and~$\Gamma_2$ be noncocompact
arithmetic subgroups of $\SL(3,\real)$ that correspond to two different
values of~$r$, say~$r_1$ and~$r_2$. Show that $\Gamma_1$ is
not commensurable to any conjugate of~$\Gamma_2$.
 \hint{There is a diagonal matrix in~$\Gamma_1$ whose
trace is not in $\integer \bigl[\! \sqrt{r_2} \bigr]$.}

\end{exercises}

\section{Cocompact arithmetic subgroups of \texorpdfstring{$\SL(3,\real)$}{SL(3,R)}}\label{CocpctLattSL3R}

\Cref{SU2TotReal} used unitary groups over a totally real extension to construct cocompact, arithmetic subgroups of $\SL(2,\real)$. The same technique can be applied to $\SL(3,\real)$:

\begin{prop} \label{CocpctSL3RSU}
 Let
 \begin{itemize}
 \item $F$ be a totally real algebraic number field, such that $F \neq \rational$,
 \item $t,a,b \in F$, such that
 	\begin{itemize}
	\item $t, a, b > 0$,
	but
	\item $\sigma(t), \sigma(a), \sigma(b) < 0$ for every place $\sigma \neq \Id$,
	\end{itemize}
 \item $L = F \bigl[\! \sqrt{t}\bigr]$,
 \item $\tau$~be the Galois automorphism of~$L$ over~$F$,
 \item $\ints$ be the ring of integers of~$L$, and
 \item $\Gamma = \SU \bigl( \diag(a,b,-1), \tau ; \ints)$. 
 \end{itemize}
 Then $\Gamma$ is a cocompact, arithmetic subgroup of\/
$\SL(3,\real)$.
 \end{prop}

Here is a specific example:

\begin{cor} \label{EgCocpctUnitarySL3R}
Let 
	\begin{itemize}
	\item $t = \sqrt{2}$,
	\item $F = \rational \bigl[ t \bigr] = \rational \bigl[\! \sqrt{2} \bigr]$,
	\item $L = F \bigl[ \sqrt{t} \bigr] =  \rational \bigl[ \! \! \root 4 \of {2} \bigr]$,
	\item $\tau$ be the Galois automorphism of~$L$ over~$F$,
	\item $\ints \doteq \integer \bigl[ \! \! \root 4 \of {2}\bigr]$ be the ring of integers of~$L$,
	and
	\item $\Gamma = \SU( \Id_{3 \times 3}, \tau ; \ints )$.
	\end{itemize}
Then\/ $\Gamma$ is a cocompact, arithmetic subgroup of\/ $\SL(3,\real)$.
\end{cor}

 \begin{rem}
 It is necessary to assume $F \neq \rational$ in \cref{CocpctSL3RSU} (in other words, there is no analogue of \cref{SU2QCocpct} for $\SL(3,\real)$), because unitary groups over~$\rational$ yield only noncompact lattices in $\SL(3,\real)$ (as in \cref{NoncocpctInSL3Eg}), not cocompact ones \csee{SU3NotCocpct}.
 \end{rem}

Here is a quite different construction (not using unitary groups) that yields additional examples of cocompact, arithmetic subgroups. See
\cref{GoodLpForSL3Rcocpct} for explicit examples of~$L$
and~$p$ that satisfy the hypotheses.

\begin{prop} \label{CocpctSL3Rbands}
 Let
 \begin{itemize}
 \item $L$~be a cubic, Galois extension of~$\rational$
\textup(that is, a Galois extension of~$\rational$, such
that $|L:\rational| = 3$\textup),
 \item $\sigma$~be a generator of $\Gal(L/\rational)$
\textup(note that $\Gal(L/\rational)$, being of
order~$3$, is cyclic\textup),
 \item $\ints$ be the ring of integers of~$L$,
 \item $p \in \integer^+$,
 \item $\phi \colon L^3 \to \Mat_{3 \times 3}(L)$ be given
by
 \begin{equation} \label{embedDcyclic}
 \phi(x,y,z) =
 \begin{bmatrix}
 x & y & z \\
 p \, \sigma(z) & \sigma(x) & \sigma(y) \\
 p \, \sigma^2(y) & p \, \sigma^2(z) & \sigma^2(x)
 \end{bmatrix}
 ,
 \end{equation}
 and
 \item $\Gamma = \{\, \gamma \in \phi(\ints^3) \mid
\det \gamma = 1 \,\}$.
 \end{itemize}
 Then:
 \begin{enumerate}
 \item  \label{CocpctSL3Rbands-latt}
 $\Gamma$~is an arithmetic subgroup of\/ $\SL(3,\real)$.
 \item \label{CocpctSL3Rbands-cocpct}
 $\Gamma$ is cocompact if and only if $p \neq t \,
\sigma(t) \, \sigma^2(t)$, for all $t \in L$.
 \end{enumerate}
 \end{prop}

\begin{proof}
 \pref{CocpctSL3Rbands-latt} 
 It is easy to see that:
 \begin{itemize}
 \item $L \subset \real$ \csee{L/Qodd->LinR}.
 \item $\phi(L^3)$ and $\phi(\ints^3)$ are subrings
of $\Mat_{3 \times 3}(L)$ (even though $\phi$~is \emph{not}
a ring homomorphism).
 \item $\phi(L^3)$ is a $\rational$-form of $\Mat_{3 \times
3}(\real)$.
 \item $\phi(\ints^3)$ is a $\integer$-lattice in
$\phi(L^3)$.
 \item If we define $\rho \colon \Mat_{3 \times 3}(\real) \to
\End_{\real} \bigl( \Mat_{3 \times 3}(\real) \bigr)$ by
$\rho(g)(v) = gv$, then $\rho \bigl( \SL(3,\real) \bigr)$ is
defined over~$\rational$ (with respect to the
$\rational$-form $\phi(L^3)$ \csee{CocpctSL3RbandsEx/Q}).
 \item $\Gamma = \{\, g \in \SL(3,\real)
\mid g \, \phi(\ints^3) = \phi(\ints^3) \,\}$.
 \end{itemize}
 So $\Gamma$ is an arithmetic subgroup of $\SL(3,\real)$
\fullcsee{AbstractArith}{arith}.

(\ref{CocpctSL3Rbands-cocpct} $\Leftarrow$)
 If $\SL(3,\real)/\Gamma$ is not compact, then the Godement Criterion \pref{GodementCriterion} tells us there is a nontrivial
unipotent element~$u$ in~$\Gamma$. 
 This means $1$~is an eigenvalue of~$u$ (indeed, it is the
only eigenvalue of~$u$), so there is some nonzero $v
\in \real^3$ with $uv = v$. Hence $(u-1)v = 0$. Since $u
\neq \Id$ and $v \neq 0$, we conclude that $\phi(L^3)$ has a
nonzero element that is not invertible.

Hence, letting $D = \phi(L^3)$, it suffices to show that every nonzero element of $D$ is
invertible. (That is, $D$ is a
``division algebra\zz.'') For convenience, define $N \colon L \to
\rational$ by $N(t) = t \, \sigma(t) \, \sigma^2(t)$.
(In Algebraic Number Theory, $N$~is called the ``norm'' from~$L$ to~$\rational$.) We know that $p \neq N(t)$, for all
$t \in L$. It is easy to see that $N(t_1 t_2) = N(t_1) \,
N(t_2)$.
 
Note that if $xyz = 0$, but $(x,y,z) \neq (0,0,0)$, then
$\phi(x,y,z)$ is invertible. For example, if $z = 0$, then
$\det \phi(x,y,z) = N(x) + p\, N(y)$. Since $p \neq N(-x/y)
= -N(x)/N(y)$ (assuming $y \neq 0$), we see that $\det
\phi(x,y,z) \neq 0$, as desired. The other cases are
similar.

For any $x,y,z \in L$, with $z \neq 0$, we have
 $$ \phi \left( 1, - \frac{x}{p\, \sigma(z)}, 0 \right)
 \, \phi(x,y,z)
 = \phi(0,*,*)$$
 is invertible, so $\phi(x,y,z)$ is invertible.

(\ref{CocpctSL3Rbands-cocpct} $\Rightarrow$)
 If $p = t \, \sigma(t) \, \sigma^2(t)$, for some $t \in
L$, then $\phi(L^3) \iso \Mat_{3 \times
3}(\rational)$ \csee{phi(L3)=Mat}. 
From this, it is easy clear that
$\phi(\ints^3)$ contains nonidentity unipotent matrices. 
Since the determinant of any unipotent matrix is~$1$, these unipotents
belong to~$\Gamma$.
Therefore $\Gamma$ is not cocompact.
 \end{proof}

\begin{eg} \label{GoodLpForSL3Rcocpct}
 Let
 \begin{itemize}
 \item $\zeta = 2 \cos(2\pi/7)$,
 \item $L = \rational[\zeta]$, and
 \item $p$~be any prime that is congruent to either~$3$
or~$5$, modulo~$7$.
 \end{itemize}
 Then
 \begin{enumerate}
 \item $L$~is a cubic, Galois extension of~$\rational$, and
 \item $p \neq t \, \sigma(t) \, \sigma^2(t)$, for all $t
\in L$, and any generator~$\sigma$ of $\Gal(L/\rational)$.
 \end{enumerate}
 To see this, let $\omega = e^{2\pi i/7}$ be a primitive
$7^\text{th}$~root of unity, so $\zeta = \omega +
\omega^6$. Now it is well known that the Galois group of
$\rational[\omega]$ is cyclic of order~$6$, generated by
$\tau(\omega) = \omega^3$ \csee{Gal(cyclotomic)}. So the
fixed field~$L$ of $\tau^3$ is a cyclic extension of degree
$6/2 = 3$.

Now suppose $t \, \sigma(t) \, \sigma^2(t) = p$, for some
$t \in L^\times$. Clearing denominators, we have $s \,
\sigma(s) \, \sigma^2(s) = pm$, where
 \begin{itemize}
 \item $m \in \integer^+$,
 \item $s = a + b (\omega+\omega^6) + c
(\omega+\omega^6)^2$, with $a,b,c \in \integer$ and $p
\nmid \gcd(a,b,c)$.
 \end{itemize}
 Replacing~$\omega$ with the variable~$x$, we obtain
integral polynomials $s_1(x)$, $s_2(x)$, and~$s_3(x)$, such that 
 $$ s_1(x) s_2(x) s_3(x) = p m = 0 \mbox{ in }
 \frac{\integer_p[x]}{\langle x^6 + x^5 + \cdots + 1
\rangle}. $$
 This implies that $x^6 + x^5 + \cdots + 1$ is \emph{not}
irreducible in $\integer_p[x]$. This contradicts the choice
of~$p$ \csee{CycloqRedModp}.
 \end{eg}

\begin{rem} \label{AbelExtsAreCyclotomic}
 The famous \thmindex{Kronecker-Weber}{Kronecker-Weber Theorem} tells us that if $L$ is a Galois extension of~$\rational$, with abelian Galois group, then $L$ is contained in an extension obtained by adjoining an $n$th root of unity to~$\rational$ (for some~$n$). (\emph{Warning:} this does not hold for abelian extensions of algebraic number fields other than~$\rational$.) As a very special case, this implies that all of the cubic, Galois extension fields~$L$ of~$\rational$ can be
constructed quite explicitly, in the manner of \cref{GoodLpForSL3Rcocpct}:
 \begin{itemize}
 \item Choose $n \in \integer^+$, such that $\varphi(n)$ is
divisible by~$3$ (where 
 $$\varphi(n) = \# \{\, k \mid 1 \le k \le n, \ \gcd(k,n) =
1 \,\} $$
 is the Euler $\varphi$-function).
 \item Let $\omega = e^{2\pi i/n}$ be a primitive
$n^\text{th}$~root of unity.
 \item Let $H$ be any subgroup of index~$3$ in the
multiplicative group $(\integer_n)^\times$ of units
modulo~$n$.
 \item Let $\zeta = \sum_{k \in H} \omega^k = \sum_{k \in
H} \cos(2 \pi k/n)$.
 \item Let $L = \rational[\zeta]$.
 \end{itemize}
 \end{rem}

We have now seen that cocompact arithmetic subgroups of $\SL(3,\real)$ can be constructed
by two different methods: some are constructed by using unitary groups (as in \cref{CocpctSL3RSU}) and others are constructed by using ``division 
algebras'' $D = \phi(L^3)$ (as in \cref{DCocpctinSL3R}). We will see in the following % @@@
section that these two methods can be
combined: some cocompact arithmetic subgroups are constructed by using \emph{both}
unitary groups \emph{and} division algebras. The classification results in
\cref{QFormsOfSLnSect} show that all
cocompact arithmetic subgroups of $\SL(3,\real)$ can be obtained from these methods, using either unitary groups, division algebras, or 
a combination of the two (perhaps also combined with restriction of scalars). The same is true for the cocompact arithmetic subgroups of any $\SL(n,\real)$, with $n \ge 3$.

%\begin{prop} \label{CocpctSL3R}
% The arithmetic subgroups constructed in \cref{CocpctSL3RSU,CocpctSL3Rbands} are the only cocompact arithmetic subgroups
%in $\SL(3,\real)$ \textup(up to commensurability and
%conjugates\textup).
% \end{prop}

\begin{exercises}

\item \label{SU3NotCocpct}
Assume the situation of \cref{CocpctSL3RSU}, except that $F = \rational$.
More precisely, let 
	\begin{itemize}
	\item $a,b,c \in \rational$ (all nonzero),
	\item $L$ be a real quadratic extension of~$\rational$,
	\item $\tau$ be the Galois automorphism of~$L$ over~$\rational$,
	\item $\ints$ be the ring of integers of~$L$,
	\item $A = \diag( a,b,c )$,
	and
	\item $\Gamma = \SU ( A, \tau ; \ints )$, so $\Gamma$ is an arithmetic subgroup of\/ $\SL(3,\real)$.
	\end{itemize}
Show that $\Gamma$ is \underline{not} cocompact.
\hint{The equation $\tau( x{}^\transpose) A  x = 0$ for $ x \in L^3$ can be considered as an equation in $6$ variables over~$\rational$, so the Number Theory fact mentioned in the proof of \cref{NonCocptArithSOn1} implies it has a nontrivial solution.}

\item \label{L/Qodd->LinR}
 Let $L$~be a Galois extension of~$\rational$\,,
with $|L:\rational|$ odd. Show $L \subset \real$.

\item \label{CocpctSL3RbandsEx/Q}
 Assume the notation of the proof of
\cref{CocpctSL3Rbands}. For $h \in L^3$, define $T_h
\in \End_{\real} \bigl( \Mat_{3 \times 3}(\real) \bigr)$ by
$T_h(v) = \phi(h) \, v$.
 \begin{enumerate}
 \item Show that $\phi(h) \in \End_{\real} \bigl( \Mat_{3
\times 3}(\real) \bigr)_{\rational}$, where the $\rational$-form is  induced by the $\rational$-form
$\phi(L^3)$ of $\Mat_{3 \times 3}(\real)$. 
\hint{Show $\phi(h) \, \phi(L^3) \subseteq \phi(L^3)$.}
 \item Show that $\rho \bigl( \Mat_{3 \times 3}(\real)
\bigr)$ is the centralizer of $\{\, T_h \mid h \in L^3
\,\}$.
 \item Show that $\rho \bigl( \SL(3,\real) \bigr)$ is
defined over~$\rational$.
 \end{enumerate}

\item In the notation of \cref{CocpctSL3Rbands}, show
that if $p = t \, \sigma(t) \, \sigma^2(t)$, then the element 
$\phi\Bigl( 1,1/t, 1/ \bigl(t \, \sigma(t) \bigr) \Bigr)$ of $\phi(L^3)$ is not invertible.

\item \label{phi(L3)=Mat}
In the notation of \cref{CocpctSL3Rbands}, show that if 
$p = t \, \sigma(t) \, \sigma^2(t)$, for some $t \in
L$, then $\phi(L^3) \iso \Mat_{3 \times 3}(\rational)$.
\hint{$\bigl\{ \bigl( a, t \sigma(a), t \sigma(t) \sigma^2(a) \bigr) \bigr\}$
is a $3$-dimensional, $\phi(L)$-invariant $\rational$-subspace of~$L^3$.}

\item \label{CycloqRedModp}
 Let $p$ and~$q$ be distinct primes, and
 $$f(x) = x^{q-1} + \cdots + x + 1 .$$
 Show that $f(x)$ is reducible over~$\integer_p$ if and
only if there exists $r \in \{1,2,\ldots, q-2\}$, such that
$p^r \equiv 1 \pmod{q}$.
 \hint{Let $g(x)$ be an irreducible factor of $f(x)$, and
let $r = \deg g(x) < q-1$. Then $f(x)$ has a root~$\alpha$
in a finite field~$F$ of order~$p^r$. Since $\alpha$ is an
element of order~$q$ in~$F^\times$, we must have $q \mid \#
F^\times$.}

\end{exercises}

\section{Arithmetic subgroups of \texorpdfstring{$\SL(\lowercase{n},\real)$}{SL(n,R)}} \label{LattSlnRSect}

We will briefly explain how the previous results on $\SL(3,\real)$
can be generalized to higher dimensions. (The group $\SL(2,\real)$
is a special case that does not fit into this pattern.)
The proofs are similar to those for $\SL(3,\real)$.

In \cref{QuaternionAlgSect}, and in the proof of \cref{CocpctSL3Rbands}, we have seen that cocompact arithmetic subgroups of $\SL(n,\real)$ can sometimes be constructed by using rings in which every nonzero element has a multiplicative inverse. More such ``division algebras'' will be needed in order to construct all the arithmetic subgroups of $\SL(n,\real)$ with $n > 3$. 

%The list of classical Lie groups includes the groups
%$\SL(n,\quaternion)$ and $\SO(n,\quaternion)$, that are based on the ring of
%quaternions, which is a division algebra over~$\real$.
%Analogously, we have seen in \cref{EgArithGrpsChap} that $\rational$-forms are sometimes based on division algebras over~$\rational$. In fact, because $\rational$-forms can come from Restriction of Scalars \pref{ResScal->Latt}, the study of arithmetic groups involves division algebras not only over~$\rational$, but over extension fields of~$\rational$.
%In this section, we provide some basic information
%on division algebras and the associated unitary
%groups. Usually (but, unfortunately, not for $\SL(n,\real)$), it suffices to have some
%familiarity with the special case of quaternion algebras $\quaternion^{a,b}_F$
%\csee{QuaternionDefn}. 
%The
%interested reader can find a more substantial introduction to
%the theory of central division algebras (including how to
%construct them) in \cref{DivAlgsChap}. 

\subsection{Division algebras}

\begin{defn} \label{DivAlgDefn}
 An associative ring~$D$ is a \defit{division
algebra} over a field~$F$ if%
\noprelistbreak
 \begin{enumerate}
 \item \label{DivAlgDefn-F}
 $D$ contains~$F$ in its center 
 	(that is, $xf = fx$ for all $x \in D$ and $f \in F$),
 \item \label{DivAlgDefn-1}
 the element $1 \in F$ is the identity element of~$D$,
 \item \label{DivAlgDefn-dim}
 $D$ is finite-dimensional as a vector space over~$F$,
 and
 \item \label{DivAlgDefn-inv}
 every nonzero element of~$D$ has a multiplicative
inverse.
 \end{enumerate}
Furthermore, it is \defit[division algebra!central]{central} over~$F$ if the entire center of~$D$ is precisely~$F$.
 \end{defn}

\begin{rems} \ 
\noprelistbreak
 \begin{enumerate}
 \item $D$ is an \emph{algebra} over~$F$ if
\pref{DivAlgDefn-F} and~\pref{DivAlgDefn-1} hold. 
% \item The word \emph{central} requires the center of~$D$ to
%be exactly~$F$, not some larger field.
  \item The word \emph{division} requires
\pref{DivAlgDefn-inv}.
We also assume \pref{DivAlgDefn-dim}, although not all authors require this.
% \item We consider only associative algebras here, but the
%algebra of \term{octonions}, which
%is nonassociative, also arises in the theory of arithmetic
%groups \cf{Rrank1list}.
 \end{enumerate}
 \end{rems}

\begin{terminology}
 Division algebras are also called \defit[skew field]{skew fields}.
 	% or \defit[division ring]{division rings}.
 \end{terminology}

\begin{egs} \label{DivAlgEgs} \ 
\noprelistbreak
 \begin{enumerate}
 \item Any extension field of~$F$ is a division algebra over~$F$ (but is not usually central).
 \item \label{DivAlgEgs-H}
 $\quaternion = \quaternion^{-1,-1}_\real$ is a central division algebra over~$\real$.
 \item \label{DivAlgEgs-quat}
 More generally, a quaternion algebra $\quaternion^{a,b}_F$ is a
central division algebra over~$F$ if and only if $\Nred(x) \neq 0$,
for every nonzero $x \in \quaternion^{a,b}_F$
\csee{QuatInverses}.
 Note that this is consistent with \pref{DivAlgEgs-H}.
 \end{enumerate}
 \end{egs}

The following famous theorem shows that division algebras are the building blocks of simple algebras:

\begin{thm}[(\thmindex{Wedderburn's}Wedderburn's Theorem)] \label{WedderburnThm}
Let $A$ be a finite-dimensional algebra over a field~$F$. If $A$ is \defit[simple!ring]{simple} \textup(that is, if $A$ has no nonzero, proper, two-sided ideals\textup), then $A \iso \Mat_{n \times n}(D)$, for some~$n$ and some division algebra~$D$ over~$F$.
\end{thm}

\begin{proof}
Since $A$ is finite-dimensional, we may let $I$ be a minimal left ideal. 
Then $I$ is a left $A$-module that is \defit[simple!module]{simple} (that is, has no nonzero, proper submodules). So $\End_A(I)$ is a division algebra \csee{SchursLemma}; call it~$D$. 

We have $IA = A$, since $IA$ is a $2$-sided ideal and $A$~is simple. Hence, the minimality of~$I$ implies $A = Ia_1 \oplus \cdots \oplus Ia_n$, for some $a_1,\ldots,a_n \in A$ \csee{DirSumIdeals}. 

For $A$ considered as a left $A$-module, it is easy to see that each element of $\End_A(A)$ is multiplication on the right by an element of~$A$ \csee{End(A)=A}; therefore $\End_A(A) \iso A$.
On the other hand, it is easy to see that $Ia_i$ is isomorphic to~$I$ as a left $A$-module \csee{Ia=I}, so we have
	\begin{align*}
	\End_A(A) 
	&= \End_A(Ia_1 \oplus \cdots \oplus Ia_n)
	\iso \End_A(I^n)
	\\&\iso \Mat_{n \times n} \bigl( \End_A(I) \bigr)
	= \Mat_{n \times n}(D) 
	.  \end{align*}
Therefore $A \iso \Mat_{n \times n}(D) $.
%If we let $n = \dim_D(A)$, then the action of $A$ on~$I$ provides a homomorphism $A \to \Mat_{n \times n}(D)$. The homomorphism is injective (because its kernel must be trivial, since $A$ is simple). 
%The hard part (which we omit) is to show that the homomorphism is surjective; this  assertion is known as the \emph{\thmindex{Jacobson Density}{Jacobson Density Theorem}}.
\end{proof}

\begin{cor}
If $D$ is a central division algebra over~$F$, then we have $\dim_F D = d^2$, for some $d \in \integer^+$. \textup(This integer~$d$ is called the \defit[degree!of a division algebra]{\textit{degree}} of~$D$ over~$F$.\textup)
\end{cor}

\begin{proof}
Let $\overline{F}$ be the algebraic closure of~$F$. Then (from {Wedderburn's Theorem}), we see that $D \otimes_F \overline{F} \iso \Mat_{d \times d}(D')$, for some~$d$ and some central division algebra~$D'$ over~$\overline{F}$. Since $\overline{F}$ is algebraically closed, we must have $D' = \overline{F}$ \csee{DivAlg/AlgClosed}, so
	\begin{align*}
	\dim_F D = \dim_{\overline{F}} (D \otimes_F \overline{F})
	= \dim_{\overline{F}} \Mat_{d \times d}(\overline{F})
	= d^2 
	. & \qedhere \end{align*}
\end{proof}

In order to produce arithmetic groups from division algebras, the following \lcnamecref{OinDivAlg} provides an analogue of the ring of integers in an
algebraic number field~$F$. 

\begin{lem} \label{OinDivAlg}
 If $D$ is a division algebra over an algebraic
number field~$F$, then there is a subring~$\ints_D$
of~$D$, such that $\ints_D$ is a $\integer$-lattice
in~$D$. Any such subring is called an {\upshape\defit[order (in an algebraic number field)]{order}} in~$D$.
 \end{lem}

\begin{proof}
 Let $\{v_0,v_1,\ldots,v_r\}$ be a basis of~$D$
over~$\rational$, with $v_0 = 1$. Let
$\{c_{j,k}^\ell\}_{j,k,\ell=0}^r$ be the structure constants
of~$D$ with respect to this basis. That is, for $j,k \in
\{0,\ldots,r\}$, we have
 $ v_j v_k = \sum_{\ell=0}^r c_{j,k}^\ell v_\ell $.
 There is some nonzero $m \in \integer$, such that $m
c_{j,k}^\ell \in \integer$, for all $j,k,\ell$. Let
$\ints_D$ be the $\integer$-span of
$\{1, mv_1,\ldots,mv_r\}$.
 \end{proof}

In the proof of \cref{CocpctSL3Rbands}, we showed that
$\phi(L^3)$ is a division algebra if
 $p \neq t\, \sigma(t) \, \sigma^2(t)$.
 Conversely, it is known that every division algebra of degree~$3$ arises from the above construction. (In the terminology of ring theory, this means that every central division algebra of degree~$3$ is ``\label{CubicAlgIsCyclic}\term[division algebra!cyclic]{cyclic}\zz.'') Therefore, we can restate the \lcnamecref{CocpctSL3Rbands} in the following more abstract form.  
 
\begin{prop} \label{DCocpctinSL3R}
 Let
 \begin{itemize}
 \item $L$~be a cubic, Galois extension of\/~$\rational$,
 \item $D$ be a central division algebra of degree\/~$3$
over\/~$\rational$, such that $D$~contains~$L$ as a subfield,
and
 \item $\ints_D$ be an order in~$D$ \csee{OinDivAlg}.
 \end{itemize}
 Then there is an embedding $\phi \colon D \to \Mat_{3
\times 3}(\real)$, such that 
 \begin{enumerate}
 \item $\phi \bigl( \SL(1,D) \bigr)$ is a\/ $\rational$-form
of\/ $\Mat_{3 \times 3}(\real)$, and
 \item \label{DCocpctinSL3R-cocpct}
 $\phi \bigl( \SL(1,\ints_D) \bigr)$ is a cocompact,
arithmetic subgroup of\/ $\SL(3,\real)$.
 \end{enumerate}

Furthermore, $\phi \bigl( \SL(1,\ints_D) \bigr)$ is
essentially independent of the choice of~$\ints_D$ or
of the embedding~$\phi$. Namely, if $\ints_D'$
and~$\phi'$ are some other choices, then there is an
automorphism~$\alpha$ of\/ $\SL(3,\real)$, such that
$\alpha \phi' \bigl( \SL(1,\ints_D') \bigr)$ is
commensurable to $\phi \bigl( \SL(1,\ints_D)
\bigr)$.
 \end{prop}

This generalizes in an obvious
way to provide cocompact, arithmetic subgroups of
$\SL(n,\real)$. By replacing $\SL(1,\ints_D)$ with the
more general $\SL(m,\ints_D)$, we can also obtain arithmetic subgroups
that are not cocompact (if $n$ is not prime).

\begin{prop} \label{LattSLnRDivAlg}
 Let
 \noprelistbreak
 \begin{itemize}
 \item $D$ be a central division algebra of degree~$d$
over\/~$\rational$, such that $D$~splits over\/~$\real$,
 \item $m \in \integer^+$, and
 \item $\ints_D$ be $\integer$-lattice in~$D$ that is
also a subring of~$D$.
 \end{itemize}
 Then $\phi \bigl( \SL(m,\ints_D) \bigr)$ is an
arithmetic subgroup of\/ $\SL(dm,\real)$, for any embedding
$\phi \colon D \to \Mat_{d \times d}(\real)$, such that
$\phi(D)$ is a\/ $\rational$-form of\/ $\Mat_{d \times
d}(\real)$.

It is cocompact if and only if $m = 1$.
 \end{prop}

\subsection{Unitary groups over division algebras}
The definition of a unitary group is based on the Galois automorphism of a quadratic extension. This is a field automorphism of order~$2$. The following analogue makes it possible to define unitary groups over division algebras that are not required to be fields.

\begin{defn}
 Let $D$ be a central division algebra. A map $\tau
\colon D \to D$ is an \defit{anti-involution} if $\tau^2 =
\Id$ and $\tau$ is an \textbf{anti}-automorphism; that is,
$\tau(x+y) = \tau(x) + \tau(y)$ and $\tau(xy) =
\tau(y) \, \tau(x)$. (Note that $\tau$ reverses the order of
the factors in a product.)
 \end{defn}

\begin{terminology}
 Some authors call $\tau$ an involution, rather than an
anti-involution, but, to avoid confusion, our terminology emphasizes
the fact that $\tau$ is \emph{not} an automorphism (unless
$D$~is commutative).
 \end{terminology}

\begin{egs} \label{QuatConjRevDefn}
 Let $D$ be a quaternion division algebra. Then:
 \noprelistbreak
 \begin{enumerate}
 \item The map 
 \nindex{$\tau_c$ = conjugation on quaternion algebra}$\tau_c \colon D \to D$ defined by 
 $$\tau_c(a + bi + cj + dk) = a - bi - cj - dk$$
 is an anti-involution. It is called the \defit[anti-involution!standard]{standard
anti-involution} of~$D$, or the \defit[conjugation on
quaternion algebra]{conjugation} on~$D$, so $\tau_c(x)$ can also be denoted~$\overline{x}$.
 \item The map 
  \nindex{$\tau_r$ = reversion anti-involution on quaternion algebra}$\tau_r \colon D \to D$ defined by 
 $$\tau_r(a + bi + cj + dk) = a + bi - cj + dk$$
 is an anti-involution. It is called the
\defit[reversion anti-involution]{reversion} on~$D$.
 \end{enumerate}
 \end{egs}

\begin{defns} \label{SUDDefn}
Let $\tau$ be an anti-involution of a division
algebra~$D$ over~$F$. 
\noprelistbreak
 \begin{enumerate}
 \item A matrix $A \in \Mat_{n \times n}(D)$ is said to be \defit[Hermitian!matrix]{Hermitian} (or,
more precisely, \emph{$\tau$-Hermitian}) if $(A^\tau)^\transpose = A$.
\item Given a Hermitian matrix~$A$, we let
	$$ \SU(A,\tau; D) = \{\, g \in \SL(n,D) \mid (A^\tau)^\transpose A g = A \,\} .$$
 \end{enumerate}
 \end{defns}

This notation makes it possible to state a version of \cref{CocpctSL3RSU} that
replaces the quadratic extension~$L$ with a larger division algebra.

\begin{prop} \label{LattSLnRSU}
 Let
 \begin{itemize}
 \item $L$ be a real quadratic extension of\/~$\rational$,
 \item $D$~be a central simple division algebra of
degree~$d$ over~$L$,
 \item $\tau$~be an anti-involution of~$D$, such that
$\tau|_L$ is the Galois automorphism of~$L$ over~$\rational$,
 \item $b_1,\ldots,b_m \in D^\times$, such that $\tau(b_j)
= b_j$ for each~$j$,
 \item $\ints_D$ be an order in~$D$, 
 and
 \item $\Gamma = \SU \bigl( \diag(b_1,b_2,\ldots,b_m) , \tau ;
\ints_D)$. 
 \end{itemize}
 Then:
 \begin{enumerate}
 \item  $\Gamma$ is an arithmetic subgroup of\/
$\SL(md,\real)$.
 \item $\Gamma$ is cocompact if and only if, for all nonzero $x \in D^m$, we have
 $$\tau \bigl(x^\transpose) \, \diag(b_1,b_2,\ldots,b_m) \,
x \neq 0 .$$
 \end{enumerate}
 \end{prop}

Additional examples of cocompact arithmetic subgroups can be obtained
by generalizing \cref{LattSLnRSU} to allow $L$ to be a
totally real quadratic extension of a totally real
algebraic number field~$F$ (as in
\cref{CocpctSL3RSU}). However, in this situation, one
must require $b_1,\ldots,b_m$ to be chosen in such a
way that $\SU \bigl( \diag(b_1,b_2,\ldots,b_m), \tau ;
\ints_D)^\sigma$ is compact, for every
place~$\sigma$ of~$F$, such that $\sigma \neq \Id$. For $n
\ge 3$, every arithmetic subgroup of $\SL(n,\real)$ is
obtained either from this unitary construction or from
\cref{LattSLnRDivAlg} \csee{QformsOfSLn}.

\begin{exercises}

\item Show $\Nred(xy) = \Nred(x) \, \Nred(y)$ for all
elements $x$ and~$y$ of a quaternion algebra 
$\quaternion^{a,b}_F$.

\item \label{SchursLemma} 
(\thmindex{Schur's Lemma}\emph{Schur's Lemma})
Suppose $A$ is a (finite-dimensional) $F$-algebra, and $M$~is a simple $A$-module (that is finite-dimensional as a vector space over~$F$). 
Show $\End_A(M)$ is a division algebra, where
	$$ \End_A(M) = \{\, \varphi \colon M \to M \mid \varphi(am) = a \, \varphi(m), \ \forall a \in A, \ m \in M \,\} .$$
\hint{If $\varphi$ is not invertible, then it has a nontrivial kernel, which is a nontrivial $A$-submodule of~$M$.}

\item \label{DirSumIdeals}
Show that if $I$ is any minimal left-ideal of a finite-dimensional algebra~$A$, then there exist $a_1,\ldots,a_n \in A$, such that $A = Ia_1 \oplus \cdots \oplus Ia_n$.
\hint{By finite-dimensionality, $A = Ia_1 + \cdots + Ia_n$ for some $a_1,\ldots,a_n \in A$. If $n$ is minimal, then $Ia_n \notsubset a_1 + \cdots + Ia_{n-1}$, so the minimality of~$I$ implies $Ia_n \cap (a_1 + \cdots + Ia_{n-1}) = \{0\}$.}

\item \label{End(A)=A}
Show that if $A$ is a ring with identity, and we consider $A$ to be a left $A$-module, then, for every $\varphi \in \End_A(A)$, there exists $a \in A$, such that $\varphi(x) = xa$ for all $x \in A$.
\hint{Let $a = \varphi(1)$.}

\item \label{Ia=I}
For any minimal left ideal~$I$ of a ring~$A$, and any $a \in A$, such that $Ia \neq \{0\}$, show $I \iso Ia$ as left $A$-modules.
\hint{$i \mapsto ia$ is a homomorphism of modules that is obviously surjective. The minimality of~$I$ implies it is also injective.}

\item \label{DivAlg/AlgClosed}
Show that if $D$ is a division algebra over an algebraically closed field~$\overline{F}$, then $D = \overline{F}$.
\hint{Multiplication on the left by any $x \in D$ is a linear transformation, which must have an eigenvalue $\lambda \in \overline{F}$. Then $x - \lambda$ is not invertible.}

\item Suppose $\quaternion^{a,b}_F$ is a quaternion algebra
over some field~$F$, and let $L = F + Fi \subseteq \quaternion^{a,b}_F$. 
 \begin{enumerate}
 \item Show
that if $a$ is not a square in~$F$, then $L$ is a
subfield of $\quaternion^{a,b}_F$. 
 \item Show that $\quaternion^{a,b}_F$ is a two-dimensional (left) vector
space over~$L$.
 \item For each $x \in \quaternion^{a,b}_F$, define $R_x \colon \quaternion^{a,b}_F \to \quaternion^{a,b}_F$ by
$R_x(v) = vx$, and show that $R_x$ is an $L$-linear
transformation.
 \item For each $x \in \quaternion^{a,b}_F$, show $\det(R_x) = \Nred(x)$.
 \end{enumerate}

\item Let $\tau$ be an anti-involution on a division
algebra~$D$.
 \begin{enumerate}
 \item For any $J \in \Mat_{n \times n}(D)$, define $B_J
\colon D^n \times D^n \to D$ by 
 $$ B_J(x,y) = \tau(x^\transpose) J y $$
 for all $x,y \in D^n = \Mat_{n \times 1}(D)$. Show that
$B_J$ is a Hermitian form if and only if $\tau
(J^\transpose) = J$.
 \item Conversely, show that if $B$ is a Hermitian form
on~$D^n$, then  $B = B_J$, for some $J \in
\Mat_{n \times n}(D)$.
 \end{enumerate}

\item Let $D$ be a finite-dimensional algebra over
a field~$F$. Show that $D$ is a \term{division
algebra} if and only if $D$~has no proper, nonzero
left ideals. (We remark that, by definition, $D$ is
\defit[simple!algebra]{simple} if and only if it has no
proper, nonzero \emph{two-sided} ideals.)

 \end{exercises}

\begin{notes}

Generalizing the examples considered here, see \cref{ArithClassicalChap} for the construction of all arithmetic subgroups of classical groups (except some strange arithmetic subgroups of groups, such as $\SO(1,7)$, whose complexification has $\SO(8,\complex)$ as a simple factor).

The construction of all arithmetic subgroups of
$\SL(2,\real)$ is discussed (from the point of view of
quaternion algebras) in \cite[Chap.~5]{SKatok-Fuchsian}.

See 
%\cite[Thm.~1 of \S1.7 and Thm.~5 of \S1.6, pp.~61 and~51]{BorevichShafarevich} or 
\cite[Cor.~2 of \S4.3.2, p.~43]{Serre-CourseArith} for a proof of the fact (used in \cref{NonCocptArithSOn1,SU3NotCocpct}) that if $a_1,\ldots,a_{n+1} \in \rational$ are not all of the same sign, and $n \ge 4$, then the equation $a_1 x_1^2 + \cdots a_{n+1} x_{n+1}^2 = 0$ has a nontrivial integer solution. It is called \thmindex{Meyer's}{Meyer's Theorem}, and will be used again in \cref{QrankGap}.

The original paper of Gromov and Piatetski-Shapiro
\cite{GromovPS} on the construction of nonarithmetic
lattices in $\SO(1,n)$ (\S\ref{NonArithSO1nSect}) is highly
recommended. The exposition there is very understandable,
especially for a reader with some knowledge of arithmetic
groups and hyperbolic manifolds. A brief treatment also
appears in \cite[App.~C.2, pp.~362--364]{MargulisBook}.

It was known quite classically that there are nonarithmetic
lattices in $\SO(1,2)$ (or, in other words, in
$\SL(2,\real)$). This was extended to $\SO(1,n)$, for $n \le
5$, by Makarov \cite{Makarov-nonarith} and Vinberg
\cite{Vinberg-nonarith}. The nonarithmetic
lattices of Gromov and Piatetski-Shapiro
\cite{GromovPS} came later.
Nonarithmetic lattices in $\SU(1,n)$ were constructed by
Mostow \cite{Mostow-nonarith} for $n = 2$, and by Deligne and
Mostow \cite{DeligneMostow-nonarith} for $n = 3$. These
results on $\SO(1,n)$ and $\SU(1,n)$ are presented briefly
in \cite[App.~C, pp.~353--368]{MargulisBook}.

The Kronecker-Weber Theorem can be found in books on Class Field Theory, such as \cite[Thm.~5.1.10, p.~324]{Neukirch-AlgNumThy} (or see \cite{Greenberg-ElemPfKroneckerWeber}). 
%The special case of cubic extensions, which is all that is needed for  \cref{AbelExtsAreCyclotomic}, can be found in
% \url{http://sbseminar.wordpress.com/2008/05/05/the-cubic-kronecker-weber-theorem/} 

Wedderburn's Theorem \pref{WedderburnThm} is proved in 
\cite[Thm.~3.5, p.~49]{Pierce-AssocAlgs},
 %\cite[pp.~15--16]{Draxl-SkewFields},
 and other introductory texts on noncommutative rings (often in the more general setting of semisimple Artinian rings).

The fact that division algebras of degree~$3$ are cyclic (mentioned on page~\pageref{CubicAlgIsCyclic}) is due to Wedderburn \cite{Wedderburn-OnDivAlgs}, and a proof can be found in \cite[Thm.~2.9.17, p.~69]{Jacobson-FDDivAlgs}.
Much more generally, the famous (and much more difficult) \thmindex{Albert-Brauer-Hasse-Noether}\text{Albert-Brauer-Hasse-Noether Theorem} states that any division algebra (of any degree) over a finite extension of~$\rational$ is cyclic. It was first proved in \cite{AlbertHasee-AlgsOverNumFld,BrauerHasseNoether-Hauptsatzes}.
See \cite[proof of Thm.~32.20, p.~280]{Reiner-MaximalOrders} for references to more modern expositions of Class Field Theory that provide proofs.

\end{notes}

  \immediate\addtocontents{toc}{\protect\toceject} % @@@
 %!TEX root = IntroArithGrps.tex

\standassumpfalse
\mychapter{\texorpdfstring{$\SL(\lowercase{n},\integer)$ is a lattice\\in $\SL(\lowercase{n},\real)$}%
	{SL(n,Z) is a lattice in SL(n,R)}}
 \label{SLnZLattChap}
\standassumptrue

\prereqs{definition of $\integer$-lattices (\cref{DefdQAbstract}) and Moore Ergodicity Theorem (\cref{MooreErgBasicSect}).}

In this chapter, we describe two different proofs of the following crucial fact, which is the basic case of the fundamental fact that if $G$ is defined over~$\rational$, then $G_\integer$ is a lattice in~$G$ \csee{arith->latt}. This special case was specifically mentioned (without proof) in \fullcref{ArithLattEg}{SLnZ}.

\begin{thm} \label{SLNZISLATT}
$\SL(n,\integer)$ is a lattice in\/ $\SL(n,\real)$.
\end{thm}

The case $n = 2$ of \cref{SLNZISLATT} was established in \cref{SL2Zlatt}, by constructing a subset~$\fund$ of $\SL(2,\real)$, such that
	\noprelistbreak
	\begin{enumerate}
	\item $\SL(2,\integer) \cdot \fund = \SL(2,\real)$,
	and
	\item $\fund$ has finite measure. 
	\end{enumerate}
Our first proof of \cref{SLNZISLATT} shows how to generalize this approach to other values of~$n$, by choosing $\fund$ to be an appropriate ``Siegel set'' \csee{SiegelSLnZSect,SLNZISLATTSiegelPfSect}.

\begin{rems} \ 
\noprelistbreak
	\begin{enumerate}
	\item As was mentioned on \cpageref{arith->lattNotPf}, the \emph{statement} that $G_\integer$ is a lattice in~$G$ is more important than the \emph{proof}. The same is true of the special case in \cref{SLNZISLATT}, but it is advisable to understand at least the statements of the three main ingredients of our first proof:
	\begin{enumerate}
	\item the definition of a Siegel set \csee{SiegelSLnZSect},
	\item the fact that every Siegel set has finite measure \csee{SiegelSLnRFinMeas}, 
	and
	\item the fact that some Siegel set is a coarse fundamental domain for $\SL(n,\integer)$ in $\SL(n,\real)$ \ccf{SiegelFundDomSLnZ}.
	\end{enumerate}
	
	\item This subject is often called \defit{Reduction Theory}.
The idea is that, given an element~$g$ of~$G$, we would like to multiply $g$ by an element~$\gamma$ of~$\Gamma$ to make the matrix $\gamma g$ as simple as possible. That is, we would like to ``reduce''~$g$ to a simpler form by multiplying it by an element of~$\Gamma$. This is a generalization of the classical reduction theory of quadratic forms, which goes back to Gauss and others.

	\end{enumerate}
\end{rems}

Unfortunately, serious complications arise when using Siegel sets to establish in general that $G_\integer$ is a lattice in~$G$ \cref{arith->latt} (see the proof in \cref{RedThyPfSect}). Therefore, we will give a second proof with the virtue that it can easily be extended to establish that all arithmetic subgroups are lattices (see \cref{SLNZISLATTSlickSect} for this proof of \cref{SLNZISLATT}, and see \cref{GZLattEx} for the generalization to a proof of \cref{arith->latt}). However, this argument relies on a fact about $\SL(n,\real)/\SL(n,\integer)$ that we will not prove in general \csee{DaniMargulisUnipReturns}.

\begin{warn}
 \textbf{The Standing Assumptions (\ref{standassump} on page~\pageref{standassump}) are \emph{not} in effect in this chapter}, because we are \emph{proving} that $\Gamma = \SL(n,\integer)$ is a lattice, instead of \emph{assuming} that it is a lattice. 
 %On the other hand, we do assume, as always, that $G$ is a semisimple Lie group that is contained in $\SL(\ell,\real)$, for some~$\ell$.
 \end{warn}

\section{Iwasawa decomposition: \texorpdfstring{$\SL(\lowercase{n},\real) = KAN$}{SL(n,R) = KAN}} \label{IwasawaSLnZ}

The definition of a ``Siegel set'' is based on the following fundamental structure theorem:

\begin{thm}[(\thmindex{Iwasawa decomposition!of $\SL(n,\real)$}Iwasawa Decomposition of $\SL(n,\real)$)] \label{IwasawaDecompSLnR}
In $G = \SL(n,\real)$, let
\begin{align*}
 K &= \SO(n), 
&
  N &= \left\{
 \begin{Smallbmatrix}
 1 \\
  & 1 &  \vbox to 0pt{\vss \hbox to 0pt{\Huge $*$\hss}} \\
  &  \vbox to 0pt{\vss\hbox to 0pt{\hss\Huge $0$}\vss} & \rotatebox{-10}{$\ddots$} \\
  & & & 1
  \end{Smallbmatrix}
  \right\}
  , &
  A &= \left\{
 \begin{Smallbmatrix}
 a_1 \\
  & a_2 &  \vbox to 0pt{\vss \hbox to 0pt{\ \Huge $0$\hss}} \\
  &  \vbox to 0pt{\vss\hbox to 0pt{\hss\Huge $0$\ }\vss} & \ddots \\
  & & & a_n
  \end{Smallbmatrix}
  \right\}^{\lower5pt\hbox{\!\larger$\circ$}}
  . \end{align*}
 Then $G = K A N$.
 In fact, every $g \in G$ has a \bemph{unique} representation of the form $g = k a u$ with $k \in K$, $a \in A$, and $u \in N$. 
 \end{thm}

\begin{proof}
It is important to note that, because of the superscript ``$\circ$'' in its definition, $A$~is only the \emph{identity component} of the group of diagonal matrices; the entire group of diagonal matrices has a nontrivial intersection with~$K$. With this in mind, the uniqueness of the decomposition is easy \csee{KANunique}.

We now prove the existence of $k$, $a$, and~$u$. To get started, let $\varepsilon_1,\ldots,\varepsilon_n$ be the standard basis of~$\real^n$. Then, for any $g \in G$, the set $\{g \varepsilon_1, \ldots,g \varepsilon_n\}$ is a basis of~$\real^n$. 

The Gram-Schmidt Orthogonalization process constructs a corresponding orthonormal basis $w_1,\ldots,w_n$. We briefly recall how this is done: 
for $1 \le i \le n$, inductively define 
	\begin{align*}
	w_i^* &= v_i - \sum_{j = 1}^{i-1} \langle v_i \mid w_j \rangle \, w_j 
	\quad \text{and} \quad
	w_i = \frac{1}{\|w_i^*\|} w_i^* 
	, & \text{where $v_i = g \varepsilon_i$}
	. \end{align*}
It is easy to verify that $w_1,\ldots,w_n$ is an orthonormal basis of~$\real^n$ \csee{GSVecsAreOrtho}.

Since $\{w_1, \ldots,w_n\}$ and $\{\varepsilon_1,\ldots,\varepsilon_n\}$  are orthonormal, there is an orthogonal matrix $k \in \Ortho(n)$, such that $kw_i = \varepsilon_i$ for all~$i$. Then
	$$ k w_i^* = k \cdot \|w_i^*\| \, w_i = \|w_i^*\| (k \, w_i) = \|w_i^*\| \, \varepsilon_i ,$$
so there is a diagonal matrix~$a$ (with positive entries on the diagonal), such that 
	$$ \text{$k w_i^* = a \varepsilon_i$ \ for all~$i$} . $$
Also, it is easy to see (by induction) that $w_i \in \langle v_1,\ldots,v_i\rangle$ for every~$i$. With this in mind, we have
	\begin{align*} g^{-1} w_i^* 
	&= g^{-1} \, v_i - g^{-1} \, \sum\nolimits_{j = 1}^{i-1} \langle v_i \mid w_j \rangle \, w_j 
	\\&\in g^{-1} \, v_i + g^{-1} \, \bigl\langle v_1,\ldots,v_{i-1} \bigr\rangle
	\\&= \varepsilon_i + \bigl\langle \varepsilon_1 , \ldots, \varepsilon_{i-1} \bigr\rangle
	, \end{align*}
so there exists $u \in N$, such that 
	$$ \text{$g^{-1} w_i^* = u \varepsilon_i$ \ for all~$i$} . $$
Therefore
	$$ u^{-1} g^{-1} w_i^* = \varepsilon_i = a^{-1} k w_i^* ,$$
so $u^{-1} g^{-1} = a^{-1} k$. Hence, $g = k^{-1} a u^{-1} \in KAN$ \csee{kInKEx}.
\end{proof}

\begin{exercises}

\item \label{KANunique}
Show that if $k_1a_1u_1 = k_2a_2u_2$, with $k_i \in K$, $a_i \in A$, and $u_i \in N$, then $k_1 = k_2$, $a_1 = a_2$, and $u_1 = u_2$.
\hint{Show $k_1^{-1} k_2 = a_1 u_1 u_2^{-1} a_2^{-1} \in K \cap AN = \{e\}$, so $k_1 = k_2$. This implies $a_1^{-1} a_2 = u_1 u_2^{-1} \in A \cap N = \{e\}$, so $a_1 = a_2$ and $u_1 = u_2$.}

\item \label{GSVecsAreOrtho}
In the notation of the proof of \cref{IwasawaDecompSLnR}, show $\{w_1,\ldots,w_n\}$ is an orthonormal basis of~$\real^n$.
\hint{Calculating an inner product shows that $w_i^* \perp w_k$ whenever $i > k$.}

\item \label{IwasawContinuousSLnR}
Show that the components $k$, $a$, and~$u$ in the Iwasawa decomposition $g = k a u$ are real analytic functions of~$g$.
\hint{The matrix entries of $a$ and~$k^{-1}$ can be written explicitly in terms of the vectors $w_i^*$ and~$w_i$, which are real-analytic functions of~$g$. Then $u = a^{-1} k^{-1} g$ is also real analytic.}

\item \label{kInKEx}
In the proof of \cref{IwasawaDecompSLnR}, note that:
	\begin{itemize}
	\item $a$ is a diagonal matrix, but we do not know the determinant of~$a$, so it is not obvious that $a \in \SL(n,\real)$,
	and
	\item $k \in \Ortho(n)$, but $K = \SO(n)$, so it is not obvious that $k \in K$. 
	\end{itemize}
From the fact that $g = k^{-1} a u^{-1}$, show $a \in A$ and $k \in K$.
\hint{We know $\det k \in \{\pm1\}$, $\det a > 0$, $\det u = 1$, and the determinant of a product is the product of determinants.} % $(\det k)^{-1} (\det a)(\det u)^{-1}$

\item \label{PermuteIwasawaSLnR}
Let $G = \SL(n,\real)$.
	\begin{enumerate}
	\item Show $ G = KNA = ANK = NAK$.
		\hint{We have $AN = NA$ and $G = G^{-1}$.}
	\item \label{PermuteIwasawaSLnR-NKA} \optional\ \harder  Show $G \neq NKA$ (if $n \ge 2$). 
	\hint{For $n = 2$, the action of~$G$ by isometries on~$\hyperbolic^2$ yields a simply transitive action on the set of unit tangent vectors. Let $v$~be a vertical tangent vector at the point~$i$, and let $w$ be a horizontal tangent vector at the point $2i$. The $N$-orbit of~$w$ consists of horizontal vectors at points on the line $\real + 2i$, but vectors in the $KA$-orbit of~$v$ are horizontal only on the line $\real + i$.}
	\end{enumerate}

\item \label{ConjToSOn}
Show that every compact subgroup of $\SL(n,\real)$ is conjugate to a subgroup of $\SO(n)$.
\hint{For every compact subgroup~$C$ of $\SL(n,\real)$, there is a $C$-invariant inner product on~$\real^n$, defined by
	$ \langle v \mid w \rangle = \int_C (cv \cdot cw) \, dc$.
Since $\langle v \mid w \rangle = gv \cdot gw$ for some $g \in \SL(n,\real)$, the usual dot product is invariant under some conjugate of~$C$. This conjugate is contained in $\SO(n)$.}

\end{exercises}

\section{Siegel sets for \texorpdfstring{$\SL(\lowercase{n},\integer)$}{SL(n,Z)}}  \label{SiegelSLnZSect}

\begin{eg}
Let $\Gamma = \SL(2,\integer)$ and $G = \SL(2,\real)$.
\Cref{FunDomSL2ZFig}  (on page~\pageref{FunDomSL2ZFig}) 
depicts a well-known fundamental domain for the action of~$\Gamma$ on the upper half plane~$\hyperbolic$. (We have already seen this in \cref{FundDomSL2R}.)
For convenience, let us give it a name, say~$\overline{\fund_0}$. There is a corresponding fundamental domain~$\fund_0$ for $\Gamma$ in~$G$, namely
	$$ \fund_0 = \{\, g \in G \mid g(i) \in \overline{\fund_0} \,\} $$
\ccf{FundInX->FundInG}.
\end{eg}

\setlength{\templength}{0.425\textwidth}
\let\oldthefigure=\thefigure
\begin{figure}[ht]
\refstepcounter{figure} \label{WeakAndFunDomSL2ZFig}
\addtocounter{figure}{-1}
\renewcommand{\thefigure}{\oldthefigure(a)}
\refstepcounter{figure} \label{FunDomSL2ZFig}
\begin{minipage}{\templength}
	$$ \includegraphics{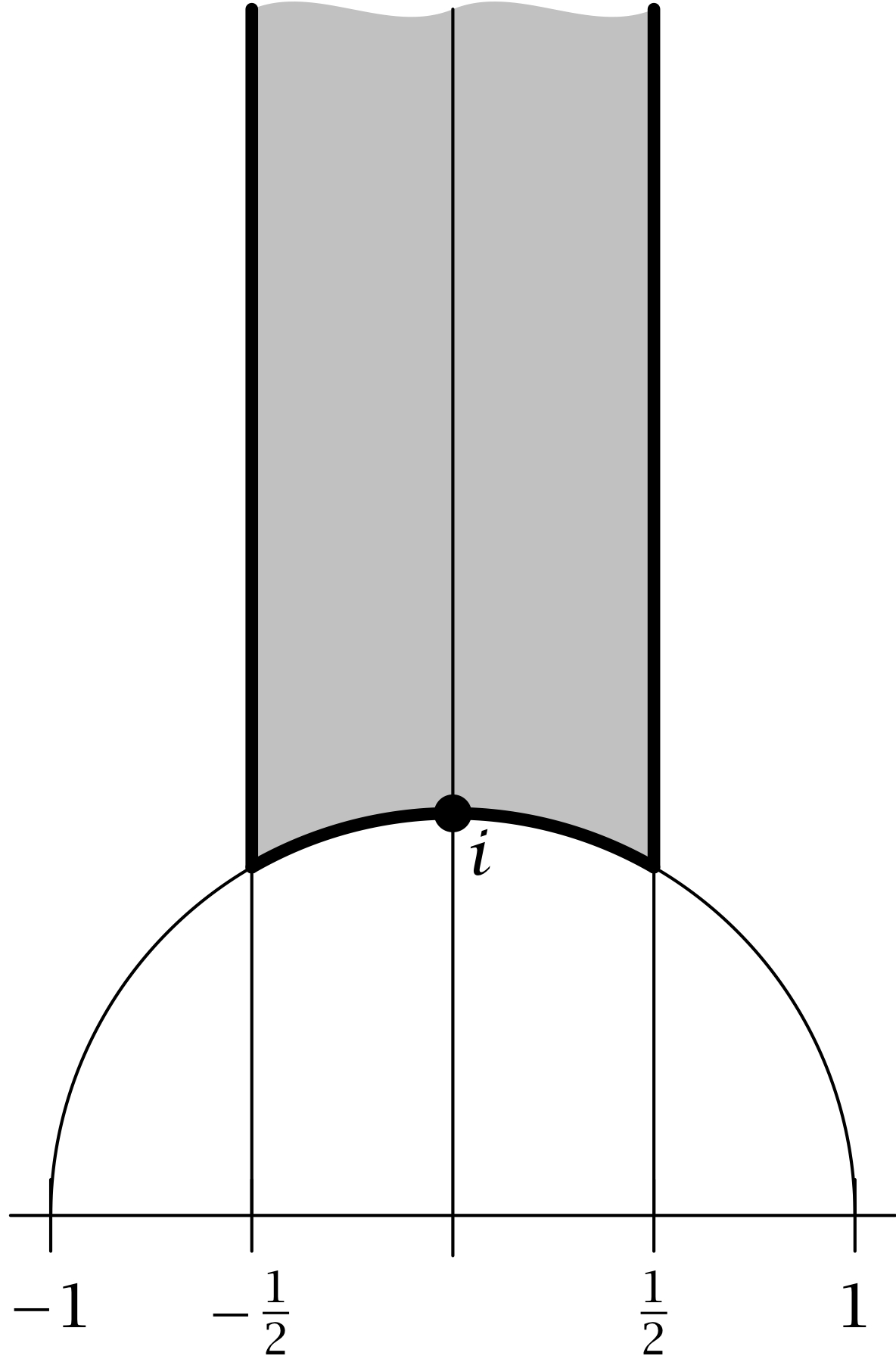} $$
\subfigurename{a}. A fundamental domain $\overline{\fund_0}$.
\end{minipage}
%texpreamble
%(" \usepackage[LY1]{fontenc}
% \usepackage[expert, LY1, mylucidascale]{mylucidabr} % I adjusted the scaling
% \usepackage{amsmath}
% \everymath{\displaystyle}
% ");
%defaultpen(  fontcommand("\normalfont") + fontsize(10) ); 
%
%from graph access *;
%
%size(0,3inch);
%pair i = (0,1);
%real h = sqrt(3)/2;
%pair a = (-1/2, h);
%pair b = (1/2, h);
%real top = 3;
%filldraw( (1/2,h) .. i .. (-1/2,h)--(-1/2,top){ENE}..{ENE}(0,top){ENE}..{ENE}(1/2,top)--cycle , gray(0.8) , invisible  );
%draw( (-1/2,h)--(-1/2,top) , linewidth(2) );
%draw( (1/2,h)--(1/2,top) , linewidth(2) );
%draw( (-1/2,h) .. i .. (1/2,h) , linewidth(2) );
%draw( (-1,0){N}..(-1/2,h) .. i .. (1/2,h) .. {S}(1,0));
%draw( (-1/2,0) -- (-1/2,h) );
%draw( (1/2,0) -- (1/2,h) );
%dotfactor = 12;
%dot( i ); label( "$i$", (0,1), SE );
%real[] xticklist = {-1, -1/2, 1/2, 1};
%string labelfunc(real x){ 
%	if (x < -0.75){return "$-1$";} 
%	else if (x > 0.75) { return "$1$" ;}
%	else if (x < 0) { return "$\textstyle-\frac{1}{2}$" ;}
%	else if (x == 0) { return "$0$" ;}
%	else { return "$\textstyle\frac{1}{2}$" ;}
%	}
%xaxis(-1.1, 1.1, Ticks(xticklist, ticklabel=labelfunc) );
%yaxis(-0.1, top, true);
\qquad
\addtocounter{figure}{-1}
\renewcommand{\thefigure}{\oldthefigure(b)}
\refstepcounter{figure} \label{WeakFunDomSL2ZFig}
\begin{minipage}{\templength}
	$$ \includegraphics{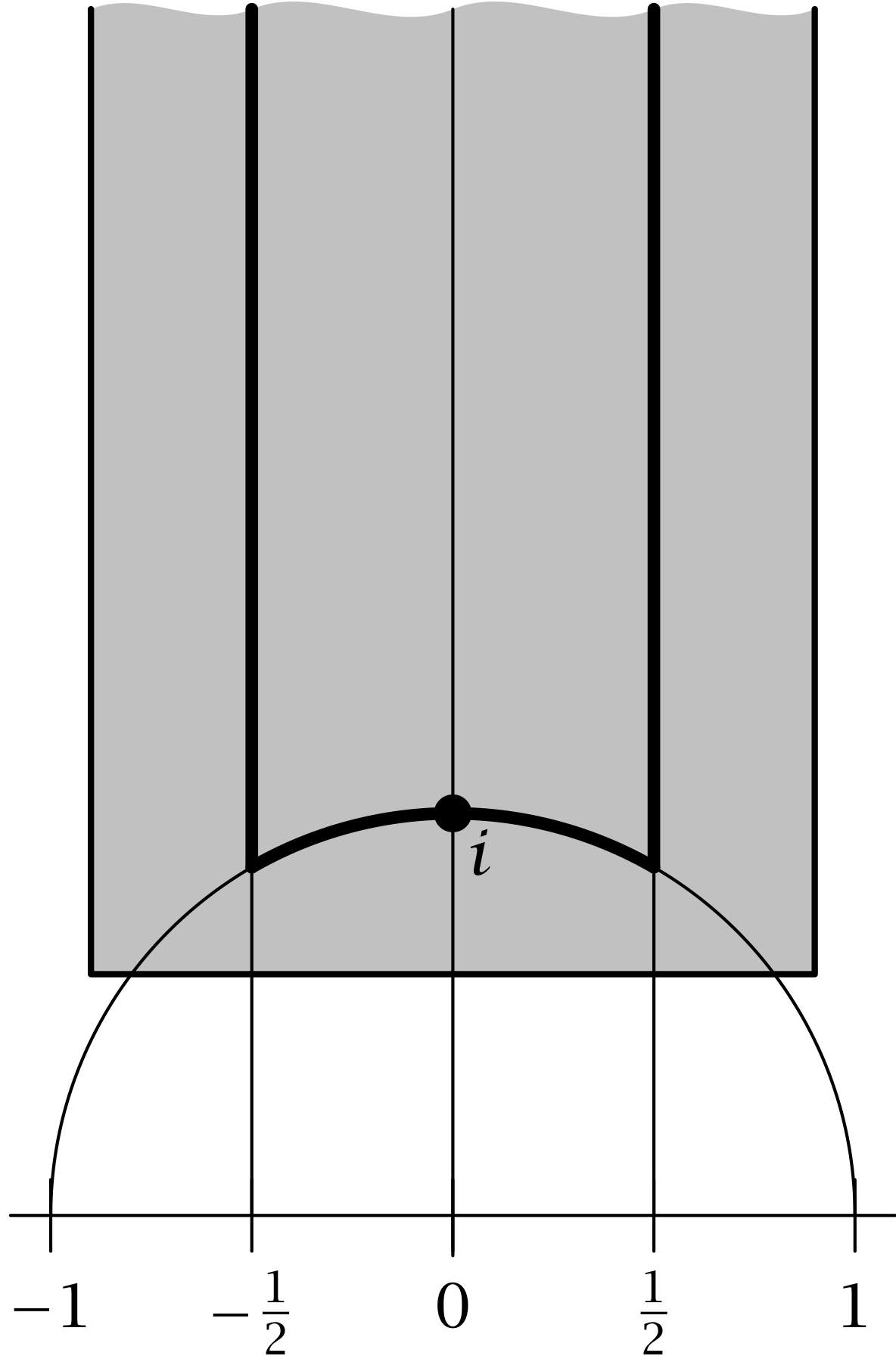} $$
\subfigurename{b}. A coarse fundamental domain~$\overline{\fund}$.
\end{minipage}
\end{figure}
%texpreamble
%(" \usepackage[LY1]{fontenc}
% \usepackage[expert, LY1, mylucidascale]{mylucidabr} % I adjusted the scaling
% \usepackage{amsmath}
% \everymath{\displaystyle}
% ");
%defaultpen(  fontcommand("\normalfont") + fontsize(10) ); 
%
%from graph access *;
%
%size(0,3inch);
%pair i = (0,1);
%real h = sqrt(3)/2;
%pair a = (-1/2, h);
%pair b = (1/2, h);
%real top = 3, x = 0.9, y = 0.6;
%filldraw( (x,y) -- (-x,y)--(-x,top){ENE}.. {ENE}(-1/2,top){ENE}..{ENE}(0,top){ENE}..{ENE}(1/2,top){ENE}..{ENE}(x,top)--cycle , gray(0.8) , invisible  );
%draw( (-x,top)--(-x,y)--(x,y)--(x,top), linewidth(1) );
%draw( (-1/2,h)--(-1/2,top) , linewidth(2) );
%draw( (1/2,h)--(1/2,top) , linewidth(2) );
%draw( (-1/2,h) .. i .. (1/2,h) , linewidth(2) );
%draw( (-1,0){N}..(-1/2,h) .. i .. (1/2,h) .. {S}(1,0));
%draw( (-1/2,0) -- (-1/2,h) );
%draw( (1/2,0) -- (1/2,h) );
%dotfactor = 12;
%dot( i ); label( "$i$", (0,1), SE );
%real[] xticklist = {-1, -1/2, 0, 1/2, 1};
%string labelfunc(real x){ 
%	if (x < -0.75){return "$-1$";} 
%	else if (x > 0.75) { return "$1$" ;}
%	else if (x < 0) { return "$\textstyle-\frac{1}{2}$" ;}
%	else if (x == 0) { return "$0$" ;}
%	else { return "$\textstyle\frac{1}{2}$" ;}
%	}
%xaxis(-1.1, 1.1, Ticks(xticklist, ticklabel=labelfunc) );
%yaxis(-0.1, top, true);

Unfortunately, the shape of~$\overline{\fund_0}$ is not entirely trivial, because the bottom edge is curved. Furthermore, the shape of a fundamental domain gets much more complicated when $G$ is larger than just $\SL(2,\real)$. Therefore, we will content ourselves with finding a set that is easier to describe, and is close to being a fundamental domain. 

\begin{eg}
To construct a region that is simpler than~$\overline{\fund_0}$, we can replace the curved edge with an edge that is straight. Also, because we do not need to find precisely a fundamental domain, we can be a bit sloppy about exactly where to place the edges, so we can enlarge the region slightly by moving the edges out a bit. The result is depicted in 
\cref{WeakFunDomSL2ZFig}. %on page~\pageref{WeakFunDomSL2ZFig}.
This new region~$\overline{\fund}$ is slightly larger than a fundamental domain, but it is within bounded distance of a fundamental domain,  and it suffices for many purposes. In particular, it is a coarse fundamental domain, in the sense of \cref{CoarseFundDomDefn} \csee{FIsWeakFundDom}.

An important virtue of this particular coarse fundamental domain is that it can be specified quite easily:
	$$ \overline{\fund} = \bigset{x + y i }{ \begin{matrix} c_1 \le x \le c_2, \\ y \ge c_3 \end{matrix} } $$
for appropriate $c_1, c_2, c_3 \in \real$.
\end{eg}

By using the Iwasawa decomposition $G = N A \, K$ \csee{{PermuteIwasawaSLnR}}, we can give a fairly simple description of the corresponding coarse fundamental domain~$\fund$ in $\SL(2,\real)$:

\begin{eg}  \label{SiegelSL2ZinSL2R}
Let
\noprelistbreak
	\begin{itemize}
	\item
	$\fund = \{\, g \in G \mid g(i) \in \overline{\fund} \,\}$, 
	\item $N_{c_1,c_2} = \bigset{
	\begin{bmatrix} 1 & t \\ 0 & 1 \end{bmatrix}
	}{
	c_1 \le t \le c_2}$, 
	\item $A_{c_3} = \bigset{
	\begin{bmatrix} e^t \\ &e^{-t} \end{bmatrix}
	}{
	e^{2t} \ge c_3}$,
	and
	\item $K = \SO(2)$.
	\end{itemize}
Then \csee{FIsSiegelSet}
	$$\fund = N_{c_1,c_2} A_{c_3} K  .$$
\end{eg}

Any set of the form $N_{c_1,c_2} A_{c_3} K$ is called a ``Siegel set\zz,''  so we can summarize this discussion by saying that Siegel sets provide good examples of coarse fundamental domains for $\SL(2,\integer)$ in $\SL(2,\real)$. 
%As a matter of notation, we mention that the gothic letter~$\Siegel$ is usually used to denote a Siegel set.

To construct a coarse fundamental domain for $\SL(n,\integer)$ (with $n > 2$), we generalize the notion of Siegel set to $\SL(n,\real)$.

\begin{defn}[(Siegel sets for $\SL(n,\integer)$)] \label{SiegelSLnZDefn}
Let $G = \SL(n,\real)$, and consider the Iwasawa decomposition $G = NAK$ \csee{PermuteIwasawaSLnR}.
To generalize \cref{SiegelSL2ZinSL2R}, we construct a ``Siegel set'' $\Siegel$ by choosing appropriate subsets $\overline{N}$ of~$N$ and $\overline{A}$ of~$A$, and letting $\Siegel = \overline{N} \, \overline{A} \, K$. 

\begin{itemize}
\item The set $\overline{N}$ can be any (nonempty) compact subset of~$N$. For example, we could let 
	$$\overline{N} = N_{c_1,c_2} = \{\, u \in N \mid \text{$c_1 \le u_{i,j} \le c_2$ for $i < j$} \,\} .$$
\item Note that the set $A_{c_3}$ of \cref{SiegelSL2ZinSL2R} has the following alternate description:
	$$ A_{c_3} = \{\, a \in A \mid a_{1,1} \ge c_3 \, a_{2,2} \,\}. $$
%since $e^t/e^{-t} = e^{2t}$. 
Therefore, we can generalize to $\SL(n,\real)$ by defining
	$$ A_c = \{\, a \in A \mid \text{$a_{i,i}   \ge c \, a_{i+i,i+1}$ for $i = 1,\ldots,n-1$} \,\} .$$
\end{itemize}
Thus, for $c_1,c_2 \in \real$ and $c_3 \in \real^+$, we have a Siegel set%
\nindex{$\Siegel_{c_1,c_2,c_3} = N_{c_1,c_2} A_{c_3} K$ is a Siegel set for $\SL(n,\integer)$}% no page break here !!!
	$$\Siegel_{c_1,c_2,c_3} = N_{c_1,c_2} A_{c_3} K .$$
\end{defn}

By calculating an appropriate multiple integral, it is not difficult to see that Siegel sets have finite measure:

\begin{prop}[\csee{SiegelSLnRFinMeasEx}] \label{SiegelSLnRFinMeas}
$\Siegel_{c_1,c_2,c_3}$ has finite measure\/ \textup(with respect to the Haar measure on\/ $\SL(n,\real)$\textup).
\end{prop}

\begin{exercises}

\item \label{FundInX->FundInG}
Suppose $H$ is a closed subgroup of~$G$, and $\overline{\fund}$ is a strict fundamental domain for the action of~$\Gamma$ on~$G/H$. For every $x \in G/H$, show that
	$$ \fund = \{\, g \in G \mid gx \in \overline{\fund}\,\} $$
is a strict fundamental domain for~$\Gamma$ in~$G$.

\item \label{F1inFinF2}
Suppose $\fund_1$ and $\fund_2$ are coarse fundamental domains for~$\Gamma$ in~$G$. Show that if $\fund_1 \subseteq \fund \subseteq \fund_2$, then $\fund$ is also a coarse fundamental domain for~$\Gamma$ in~$G$.

\item \label{FinUnionFundDomsIsWeak}
Suppose 
	\begin{itemize}
	\item $\fund$ is a coarse fundamental domain for the action of~$\Gamma$ on~$X$, 
	and
	\item  $F$ is a nonempty, finite subset of~$\Gamma$.
	\end{itemize}
Show that $F \fund = \bigcup_{f \in F} f \fund$ is also a coarse fundamental domain.

\item \label{CoarseInX->CoarseInG}
Suppose $H$ is a closed subgroup of~$G$, and $\overline{\fund}$ is a coarse fundamental domain for the action of~$\Gamma$ on~$G/H$. For every $x \in G/H$, show that
	$$ \fund = \{\, g \in G \mid gx \in \overline{\fund}\,\} $$
is a coarse fundamental domain for~$\Gamma$ in~$G$.

\item \label{FCoveredByF0}
In the notation of \cref{WeakAndFunDomSL2ZFig}, show that the coarse fundamental domain $\overline{\fund}$ is contained in the union of finitely many $\Gamma$-translates of the fundamental domain~$\overline{\fund_0}$.

\item \label{FIsWeakFundDom}
Show that the set $\overline{\fund}$ depicted in \cref{WeakFunDomSL2ZFig} is indeed a coarse fundamental domain for the action of~$\Gamma$ on~$\hyperbolic$. \hint{\Cref{FinUnionFundDomsIsWeak,FCoveredByF0}.  You may assume (without proof) that $\overline{\fund_0}$ is a fundamental domain.}

\item \label{WeakFundDomFinInd}
Suppose 
	\begin{itemize}
	\item $\fund$ is a coarse fundamental domain for $\Gamma$ in~$G$, 
	and
	\item  $F_1$ is a nonempty, finite subset of~$\Gamma$,
	and
	\item $\Gamma_1$ is a finite-index subgroup of~$\Gamma$. 
	\end{itemize}
Show:
	\begin{enumerate}
	\item that $F_1 \fund$ is also a coarse fundamental domain for $\Gamma$ in~$G$.
	\item \label{WeakFundDomFinInd-IsFund}
	If $\Gamma_1 F_1 \fund = G$, then $F_1 \fund$ is a  coarse fundamental domain for both $\Gamma$ and~$\Gamma_1$ in~$G$. 
	\end{enumerate}

\item Assume $\Gamma$ is infinite (or, equivalently, that $G$ is not compact), and $\Gamma_1$ is a finite-index, proper subgroup of~$\Gamma$. Show there exists a (strict) fundamental domain for $\Gamma_1$ in~$G$ that is \emph{not} contained in any coarse fundamental domain for $\Gamma$ in~$G$.
\hint{Construct a strict fundamental domain for $\Gamma_1$ that contains a strict fundamental domain~$\fund_0$ for~$\Gamma$, but is not covered by finitely many $\Gamma$-translates of $\fund_0$.}

\item \label{ConjFundDomsInGIsWeak}
Suppose 
	\begin{itemize}
	\item $\fund$ is a coarse fundamental domain for $\Gamma$ in~$G$, 
	and
	\item  $g \in \nzer_G(\Gamma)$.
	\end{itemize}
Show that  $\fund^g = g^{-1} \fund g$ is also a coarse fundamental domain.

\item \label{FIsSiegelSet}
Let $\fund$ be the coarse fundamental domain for $\SL(2,\integer)$ in $\SL(2,\real)$ that is  defined in \cref{SiegelSL2ZinSL2R}. Verify that $\fund = N_{c_1,c_2} A_{c_3} K$.
%\hint{We have
%	$K(i) = i$,
%	$A_{c_3} (i) = \{\, y i \mid y \ge c_3 \,\}$,
%	and $N_{c_1,c_2} \bigl( \{\, y i \mid y \ge c_3 \,\} \bigr) = \fund$.}

\item \label{Ac=aA+}
Let $G = \SL(n,\real)$. Given $c > 0$, show there exists $a \in A$, such that $A_c = a A^+$.

\item \label{HaarOnG=KAN}
\emph{This exercise provides a description of the Haar measure on~$G$.}

Let $dg$, $dk$, $da$, and $du$ be the Haar measures on the unimodular groups $G$, $K$, $A$, and~$N$, respectively, where $G = KAN$ is an Iwasawa decomposition. Also, for $a \in A$, let $\rho(a)$ be the modulus (or Jacobian) of the action of~$a$ on~$N$ by conjugation, so
	$$ \text{$\int_N f( a^{-1} u a) \, du = \int_N f(u) \, \rho(a) \, du$ for $f \in C_c(N)$} . $$
Show, for $f \in C_c(G)$, that
	\begin{align*}
	\int_G f \, dg 
	&=  \int_K \, \int_N \, \int_A \, f(kua) \, da \, du \, dk
	\\&= \int_K \, \int_A \, \int_N \, f(kau) \, \rho(a) \, du \, da \, dk
	. \end{align*}
\hint{Since $G$ is unimodular, $dg$ is invariant under left translation by elements of~$K$ and right translation by elements of~$AN$.}

\item \label{ModFuncANinSLnR}
Let $G = \SL(n,\real)$, choose $N$, $A$, and~$K$ as in \cref{SiegelSLnZDefn}, and 
define $\rho$ as in \cref{HaarOnG=KAN}. Show
	$$ \rho \left( \begin{Smallbmatrix} a_{1,1} \\&a_{2,2}
	\\ \BigSymbol{0}{0}{-5}&& \BigSymbol{0}{15}{15}\ddots \\&&& a_{n,n} \end{Smallbmatrix} \right) = \prod_{i<j} \frac{a_{j,j}}{a_{i,i}} .$$
%for all $a_{1,1},\ldots,a_{n,n} \in \real^+$.

\item \label{SiegelSLnRFinMeasEx}
Let $c_1,c_2,c_3 \in \real$, with $c_1 < c_2$ and $c_3 > 0$. Show that the Siegel set $\Siegel_{c_1,c_2,c_3}$ in $\SL(n,\real)$ has finite measure.
\hint{See \cref{HaarOnG=KAN,ModFuncANinSLnR} for a description of the Haar measure on $\SL(n,\real)$.}

\end{exercises}

\section{Constructive proof using Siegel sets} \label{SLNZISLATTSiegelPfSect}

In this section, we prove the following result:

\begin{thm} \label{SiegelFundDomSLnZ}
Let 
	\noprelistbreak
	\begin{itemize}
	\item $G = \SL(n,\real)$, 
	\item $\Gamma = \SL(n,\integer)$, 
	and 
	\item $\Siegel_{0,1,\frac{1}{2}} = N_{0,1} A_{1/2} K$ be the Siegel set defined in \cref{SiegelSLnZDefn}. 
	\end{itemize}
Then $G = \Gamma \, \Siegel_{0,1,\frac{1}{2}}$.
\end{thm}

This establishes \cref{SLNZISLATT}:

\begin{proof}[\bf Proof of \cref{SLNZISLATT}]
Combine the conclusion of \Cref{SiegelFundDomSLnZ} with \cref{SiegelSLnRFinMeas,Latt<>fund}.
\end{proof}

\begin{rems} \label{SiegelRemsLR} \ 
\noprelistbreak
	\begin{enumerate}
	\item \label{SiegelRemsLR-transpose}
	$\Gamma$ is written on the left in the conclusion of \cref{SiegelFundDomSLnZ}, because our definition of Siegel sets is motivated by a fundamental domain for the action of $\SL(2,\integer)$ on~$\hyperbolic^2$, and $\Gamma$~acts on the left there. However, taking the transpose of both sides of the conclusion of \cref{SiegelFundDomSLnZ} yields
		$ G = \Siegel_{0,1,\frac{1}{2}}^\transpose \Gamma$,
	where $\Siegel_{0,1,\frac{1}{2}}^\transpose = K A_{1/2} N_{c_1,c_2}^\transpose$. 
	Thus, $\Gamma$ can be written on the right, if the definition of Siegel set is modified appropriately.
	
	\item Our definition of Siegel sets uses the upper-triangular group~$N$, and \cref{SiegelFundDomSLnZ} puts $\Gamma$ on the left. Then \pref{SiegelRemsLR-transpose} uses the lower-triangular group~$N^\transpose$ (also called $N^-$), and puts $\Gamma$ on the right. Some authors reverse this, using $N^-$ when $\Gamma$ is on the left and using $N$ when the action on the right. However, to accomplish this, the inequality in the definition of~$A_c$ needs to be reversed. (See \cref{FundDomSL2ZinSL2RminEx} and the proof of \cref{SiegelFundDomSLnZ}.) 
	\end{enumerate}
\end{rems}

The following elementary observation is the crux of the proof of \cref{SiegelFundDomSLnZ}:

\begin{lem} \label{RedSLnZProjLem}
If $\Zlatt$ is any $\integer$-lattice in\/ $\real^n$, then there is an ordered basis $v_1,\ldots,v_n$ of\/~$\real^n$, such that
	\begin{enumerate}
	\item \label{RedSLnZProjLem-gen}
	$\{v_1,\ldots,v_n\}$ generates $\Zlatt$ as an abelian group,
	and
	\item \label{RedSLnZProjLem->1/2}
	$\| \proj_i^\perp v_{i+1} \| \ge \frac{1}{2} \| \proj_{i-1}^\perp v_{i} \|$ for\/ $1 \le i < n$, where\/ $\proj_i^\perp \colon \real^n \to V_i^\perp$ is the orthogonal projection onto the orthogonal complement of the subspace~$V_i$ spanned by\/ $\{v_1,v_2,\ldots,v_i\}$.
	\end{enumerate}
\end{lem}

\begin{proof}
Choose $v_1$ to be a nonzero vector of minimal length in~$\Zlatt$. Then define the remaining vectors $v_2,v_3,\ldots,v_n$ by induction, as follows:%
	\begin{itemize}
	\item[] Given $v_1,v_2,\ldots,v_i$, choose $v_{i+1} \in \Zlatt$ to make $\proj_i^\perp v_{i+1}$ as short as possible, subject to the constraint that $v_{i+1}$ is linearly independent from $\{v_1,v_2,\ldots,v_i\}$ \textup(so $\proj_i^\perp v_{i+1}$ is nonzero\/\textup). 
	\end{itemize}
We now verify \pref{RedSLnZProjLem-gen} and \pref{RedSLnZProjLem->1/2}.

\pref{RedSLnZProjLem-gen}
For each~$i$, let $\Zlatt_i$ be the abelian group generated by $v_1,v_2,\ldots,v_i$. If $\Zlatt_n \neq \Zlatt$, we may let $i$ be minimal with $\Zlatt_{i+1} \neq \Zlatt \cap V_{i+1}$. Then we must have $\proj_i^\perp \Zlatt_{i+1} \subsetneq \proj_i^\perp(\Zlatt \cap V_{i+1})$ \csee{RedSLnZProjX=YEx}, so there is some $v \in \Zlatt \cap V_{i+1}$ with $\proj_i^\perp v_{i+1} = k \cdot \proj_i^\perp v$ for some $k \ge 2$ \csee{RedSLnZProjCyclicEx}.  This contradicts the minimality of $\|\proj_i^\perp v_{i+1}\|$.

\pref{RedSLnZProjLem->1/2} 
For simplicity, assume $i = 1$ \csee{RedSLnZProj>1/2Ex}, and let $v_2^* = \proj_1^\perp v_2$, so $v_2 = v_2^* + \alpha v_1$, with $\alpha \in \real$. Obviously, there exists $k \in \integer$, such that $|\alpha - k| \le 1/2$. If 
		$ \| \proj_1^\perp v_2 \| <  \frac{1}{2} \| v_1 \| $,
then
	\begin{align*}
	 \|v_2 - k v_1\| 
	= \|v_2^* + (\alpha - k)v_1 \|
	&\le  \| v_2^*\| + |\alpha - k| \cdot \| v_1 \| 
	\\&< \frac{1}{2} \|v_1\| +  \frac{1}{2} \| v_1 \|
	= \|v_1 \| 
	. \end{align*}
This contradicts the minimality of ~$\|v_1\|$.
\end{proof}

\begin{proof}[\bf Proof of \cref{SiegelFundDomSLnZ}]
We wish to show $G = \Gamma \, \Siegel_{0,1,\frac{1}{2}} = \Gamma \, N_{0,1} \, A_{1/2} \, K$. However, since the proof uses an action of~$G$, and most readers prefer to have this action on the left, we will instead prove an analogous result with $\Gamma$ on the right: $G = \Siegel_{0,1,\frac{1}{2}}^{-} \Gamma$. 
%Roughly speaking, we take the inverse of both sides of the original equation. 
Namely, given $g \in G$, 
	$$ \text{we will show $g \in K \, A^-_{1/2} \, N_{0,1} \, \Gamma$,} $$
where $A^-_c = \{\, a^{-1} \mid a \in A_c \,\} = \{\, a \in A \mid \text{$a_{i,i} \le a_{i+1,i+1}/c$ for all~$i$} \,\}$.

For convenience, let $\Zlatt = g \integer^n$, and let $\{\varepsilon_1,\ldots,\varepsilon_n\}$ be the standard basis of~$\real^n$.  \Cref{RedSLnZProjLem} provides us with a sequence $v_1,\ldots,v_n$ of elements of~$\Zlatt$. From \fullref{RedSLnZProjLem}{gen}, we see that, after multiplying $g$ on the right by an element of~$\Gamma$, we may assume 
	$$ \text{$g \varepsilon_i = v_i$ for $i = 1,\ldots,n$} $$
\csee{eig=pmvi}. 

From the Iwasawa decomposition $G = KAN$ \csee{IwasawaDecompSLnR}, we may write
	$g = k a u$ with $k \in K$, $a \in A$, and $u \in N$.
For simplicity, let us assume $k$ is trivial \csee{SiegelSLnZ-knoteEx}, so
	$$ \text{$g = a u$ with $a \in A$ and $u \in N$.} $$
Since $g \in AN$, we know $g$ is upper triangular (and its diagonal entries are exactly the same as the diagonal entries of~$a$), so
	$$ \text{$\langle \varepsilon_1,\varepsilon_2,\ldots,\varepsilon_i \rangle
	= \langle g \varepsilon_1, g \varepsilon_2, \ldots, g\varepsilon_i \rangle 
	=  \langle v_1,v_2,\ldots,v_i \rangle$ for all~$i$} .$$
This implies that the diagonal entry $a_{i,i}$ of~$a$ is given by
	\begin{align*}
	a_{i,i}
	&=
	g_{i,i}
	=  \| \proj_{i-1}^\perp g \varepsilon_i \| 
	=  \| \proj_{i-1}^\perp v_i \| 
	\\& \le 2 \| \proj_i^\perp v_{i+1} \| 
	=  2\| \proj_i^\perp g \varepsilon_{i+1} \| 
	 = 2 g_{i+1,i+1} 
	 = 2 a_{i+1,i+1}
	 . \end{align*}
Therefore $a \in A_{1/2}^-$.

Also, there exists $ \gamma \in \Gamma \cap N$, such that
	$ u \in N_{0,1} \, \gamma$
\csee{N=NZN01}.
Therefore $g = a u \in  A_{1/2}^- \, N_{0,1} \, \gamma \subseteq  K \, A_{1/2}^- \, N_{0,1} \, \Gamma$, as desired.
\end{proof}

\begin{rem} \label{SiegelIsCoarseFundInSLnR}
	It can be shown that that the Siegel set $\Siegel_{0,1,\frac{1}{2}}$ is a coarse fundamental domain for $\SL(n,\integer)$ in $\SL(n,\real)$ \ccf{SiegelPropertySect}, but this fact is not needed in the proof that $\SL(n,\integer)$ is a lattice in $\SL(n,\real)$.
\end{rem}

\begin{exercises}

\item\label{FundDomSL2ZinSL2RminEx}
Let
	\begin{itemize}
	\item $N^-_{c_1,c_2} = \bigset{
	\begin{Smallbmatrix} 1 & 0 \\ t & 1 \end{Smallbmatrix}
	}{
	c_1 \le t \le c_2}$, 
	\smallskip % @@@
	\item $A^-_{c_3} = \bigset{
	\begin{Smallbmatrix} e^t\\ &e^{-t} \end{Smallbmatrix}
	}{
	e^{2t} \le c_3}$,
	\smallskip % @@@
	\item $K = \SO(2)$,
	and
	\item $\fund' = N^-_{c_1,c_2} A^-_{c_3} K $.
	\end{itemize}
Show that $\fund'$ is a coarse fundamental domain for $\SL(2,\integer)$ in $\SL(2,\real)$ if and only if the set $\fund = N_{c_1,c_2} A_{c_3} K$ of \cref{SiegelSL2ZinSL2R} is a coarse fundamental domain.
\hint{Conjugate by $\begin{Smallbmatrix} 0 & 1 \\ 1 & 0 \end{Smallbmatrix}$.}

\item \label{RedSLnZProjX=YEx}
In the notation of \cref{RedSLnZProjLem}, show that if $X$ and~$Y$ are two subgroups of $V_{i+1}$, such that 
	$$ \text{$X \subseteq Y$, 
	\ $X \cap V_i = Y \cap V_i$,
	\ and
	\ $\proj_i^\perp X = \proj_i^\perp Y$,} $$
then $X = Y$.

\item \label{RedSLnZProjCyclicEx}
In the notation of \cref{RedSLnZProjLem}, show that the group $\proj_i^\perp(\Zlatt \cap V_{i+1})$ is cyclic.
\hint{Since $\dim \proj_i^\perp V_{i+1} = 1$, it suffices to show $\proj_i^\perp(\Zlatt \cap V_{i+1})$ is discrete.}

\item \label{RedSLnZProj>1/2Ex}
Prove \fullcref{RedSLnZProjLem}{>1/2} without assuming $i = 1$.
\hint{Mod out $V_{i-1}$, which is in the kernel of both $\proj_{i-1}^\perp$ and $\proj_i^\perp$.}

\item 
For $g \in \GL(n,\real)$, show $g \in \GL(n,\integer)$ if and only if $g \integer^n \subseteq \integer^n$ and $g^{-1} \integer^n \subseteq \integer^n$.

\item \label{eig=pmvi}
For every $n$-element generating set $\{v_1,\ldots,v_n\}$ of the group~$\integer^n$, show there exists $\gamma \in \SL(n,\integer)$, such that $g \varepsilon_i = \pm v_i$ for every~$i$.
\hint{Show there exists $\gamma \in \GL(n,\integer)$, such that $g \varepsilon_i =  v_i$ for every~$i$.}

\item \label{SiegelSLnZ-knoteEx}
Complete the proof of \cref{SiegelFundDomSLnZ}  (without assuming the element~$k$ is trivial).
\hint{The group~$K$ acts by isometries on~$\real^n$, so replacing $\{v_1,\ldots,v_n\}$ with its image under an element of~$K$ does not affect the validity of \fullref{RedSLnZProjLem}{>1/2}.}

\item \label{N=NZN01}
For all $c \in \real$, show $N =  N_{c, c+1} N_{\integer}$.

\end{exercises}

\section{Elegant proof using nondivergence of unipotent orbits} \label{SLNZISLATTSlickSect}

We now present a very nice proof of \cref{SLNZISLATT} that relies on two key facts: the Moore Ergodicity Theorem \pref{MooreErgBasicThm}, and an important observation about orbits of unipotent elements (\cref{DaniMargulisUnipReturns}).
%that is often called the ``\term{Margulis Lemma}\zz.'' 
The statement of this observation will be more enlightening after some introductory remarks.

\begin{eg} \label{SplitDivergeSL2R}
Let $a = \begin{Smallbmatrix} 2& 0 \\ 0 & 1/2 \end{Smallbmatrix}$, or, more generally, let $a$ be any element of $\SL(2,\real)$ that is diagonalizable over~$\real$ (and is not~$\pm\Id$). Then $a$ has one eigenvalue that is greater than~$1$, and one eigenvalue that is less than~$1$ (in absolute value), so it is obvious that there exist linearly independent vectors $v_+$ and~$v_-$ in~$\real^2$, such that 
	$$ \text{$a^k v_+ \to 0$ \ and \ $a^{-k} v_- \to 0$ \  as $k \to +\infty$} . $$
By the Mahler Compactness Criterion \pref{MahlerCpct}, this implies that some of the orbits of~$a$ on $\SL(2,\real)/\SL(2,\integer)$ are ``divergent'' or ``go off to infinity'' or ``leave compact all sets\zz.'' That is, there exists $x \in \SL(2,\real)/\SL(2,\integer)$, such that, for every compact subset~$C$ of~$\SL(2,\real)/\SL(2,\integer)$,
	$$ \text{$\{\, k \in \integer \mid a^k x \in C \,\}$ is finite} $$
\csee{SplitDivergeSL2REx}. 
\end{eg}

In contrast, if $u = \begin{Smallbmatrix} 1 & 1 \\ 0 & 1 \end{Smallbmatrix}$, then it is clear that there does \emph{not} exist a nonzero vector $v \in \real^2$, such that $u^k v \to 0$ as $k \to \infty$. In fact, if $v$ is not fixed by~$u$ (i.e., if $uv \neq v$), then
	\begin{align} \label{unvtoinfty}
	\text{$\|u^k v \| \to \infty$ as $k \to \pm \infty$}
	\end{align}
\csee{UnipToInftyInR2}.
Therefore, it is not very difficult to show that \emph{none} of the orbits of~$u$ on $\SL(2,\real)/\SL(2,\integer)$ go off to infinity:

\begin{prop} \label{MargulisRecurSL2R}
If $u$ is any unipotent element of\/ $\SL(2,\real)$, then, for all $x \in \SL(2,\real)/\SL(2,\integer)$, there is a compact subset~$C$ of\/ $\SL(2,\real)/\SL(2,\integer)$, such that
	$$ \text{$\{\, k \in \integer^+ \mid u^k x \in C \,\}$ is infinite} .$$
\end{prop}

\begin{proof}
We may assume $u$ is nontrivial. Then, by passing to a conjugate (and perhaps taking the inverse), we may assume $u = \left[\begin{smallmatrix} 1 & 1 \\ 0 & 1 \end{smallmatrix} \right]$.

Choose a small neighborhood~$\open\mk$ of~$0$ in~$\real^2$ so that, for all $g \in \SL(2,\real)$, there do not exist two linearly independent vectors in $\open\mk \cap g \integer^2$ \csee{2SmallNotIndep}. Since $x \integer^2$ is discrete, we may assume $\open\mk$ is small enough that 
	\begin{align} \label{OpenNotXZ2}
	\open\mk \cap x \integer^2 = \{0\}
	. \end{align}
Since $\open\mk$ is open and $0$ is a fixed point of~$u$ (and the action of~$u^{-1}$ is continuous), there exists $r > 0$, such that
	\begin{align} \label{uBrInOpen}
	 B_r(0) \cup u^{-1} B_r(0) \subseteq \open 
	 , \end{align}
where $B_r(0)$ is the open ball of radius~$r$ around~$0$.
Let 
	$$ C = \bigset{ c \in \SL(2,\real)/\SL(2,\integer) }{ c \integer^2 \cap B_r(0) = \{0\} } .$$
The Mahler Compactness Criterion \pref{MahlerCpct} tells us that $C$ is compact.

Given $N \in \integer^+$, it suffices to show there exists $k \ge 0$, such that $u^{N+k} x \in C$. That is,
	$$\text{we wish to show there exists $k \ge 0$, such that $u^{N+k} x \integer^2 \cap B_r(0) = \{0\}$} . $$
Let $v$ be a nonzero vector of smallest length in $u^N x \integer^2$. We may assume $\|v\| < r$ (for otherwise we may let $k = 0$). Hence, \pref{OpenNotXZ2} implies that $v$ is not fixed by~$u$. Then, from \pref{unvtoinfty}, we know there is some $k > 0$, such that $\|u^k v\| \ge r$, and we may assume $k$ is minimal with this property. Therefore $\|u^{k-1} v \|< r$, so $u^{k-1} v \in B_r(0) \subseteq \open$\, by \pref{uBrInOpen}.

From the choice of~$\open$, we know that $\open \cap u^{N+k-1} x \integer$ does not contain any vector that is linearly independent from~$u^{k-1} v$. Therefore $u^{N+k} x \integer^2$ does not contain any nonzero vectors of length less than~$r$ \csee{MargulisSL2RNotLessR}, as desired.
\end{proof}

This result has a natural generalization to $\SL(n,\real)$ (but the proof is more difficult; see \cref{PfUnipOrbitsReturnSect}):

\begin{thm}[(Margulis)] 
\label{MargulisUnipReturns}
Suppose
\noprelistbreak
	\begin{itemize}
	\item $u$ is a unipotent element of\/ $\SL(n,\real)$,
	and
	\item $x \in \SL(n,\real)/\SL(n,\integer)$.
	\end{itemize}
Then there exists a compact subset~$C$ of\/ $\SL(n,\real)/\SL(n,\integer)$, such that
	$$ \text{$\{\, k \in \integer^+ \mid u^k x \in C \,\}$ is infinite} .$$
\end{thm}

In other words, every unipotent orbit visits some compact set infinitely many times. In fact, it can be shown that the orbit visits the compact set quite often --- it spends a nonzero fraction of its life in the set:

\begin{thm}[(Dani-Margulis)]
 \label{DaniMargulisUnipReturns}
Suppose
	\begin{itemize}
	\item $u$ is a unipotent element of\/ $\SL(n,\real)$,
	and
	\item $x \in \SL(n,\real)/\SL(n,\integer)$.
	\end{itemize}
Then there exists a compact subset~$C$ of\/ $\SL(n,\real)/\SL(n,\integer)$, such that
	$$ \liminf_{m \to \infty} \frac{\# \bigset{ k \in \{1,2,\ldots,m\} }{ u^k x \in C} }{m} > 0 .$$
\end{thm}

Before saying anything about the proof of this important fact, let us see how it implies the main result of this chapter:

\begin{proof}[\bf Proof of \cref{SLNZISLATT}]
Let 
	\begin{itemize}
	\item $X = \SL(n,\real)/\SL(n,\integer)$,
	and
	\item $\mu$ be an $\SL(n,\real)$-invariant measure on~$X$ \csee{HaarOnHomog}.
	\end{itemize}
We wish to show $\mu(X) < \infty$.

Fix a nontrivial unipotent element~$u$ of $\SL(n,\real)$. For each $x \in X$ and compact $C \subseteq X$, let
	\begin{align*}
	 \rho_C( x ) 
	 &= \liminf_{m \to \infty} \frac{ \# \bigset{ k \in \{1,2,\ldots,m\} }{ u^k x \in C }  } {m} 
	 . \end{align*}
Since $X$ can be covered by countably many compact sets, \cref{DaniMargulisUnipReturns} implies there is a compact set~$C$, such that
	\begin{align} \label{rhoC>0}
	 \text{$\rho_C > 0$ on a set of positive measure} 
	 \end{align}
\csee{rhoC>0Ex}. 
Letting $\chi_C$ be the characteristic function of~$C$, we have
	\begin{align*}
	\int_X \rho_C  \, d\mu
	&= \int_X \liminf_{m \to \infty} \frac{ \# \! \bigset{ k \in \{1,2,\ldots,m\} }{ u^k x \in C }  } {m} \, d\mu(x)
	\\&\le  \liminf_{m \to \infty} \int_X  \! \frac{ \# \! \bigset{ k \in \{1,2,\ldots,m\} }{ u^k x \in C }  } {m} \, d\mu(x)
	 \vbox to 0pt{\vss\vskip-5pt\hbox{\hskip0.07in$\begin{pmatrix} \text{Fatou's} \\ \text{Lemma} \\ \text{\pref{FatousLemma}} \end{pmatrix}$}\vss} % is the spacing good? @@@
	\\&=  \liminf_{m \to \infty} \frac{1}{m} \int_X  \bigl( \chi_{u^{-1}C} + \chi_{u^{-2}C} + \cdots  + \chi_{u^{-m}C} \bigr) \, d\mu
	\\&= \liminf_{m \to \infty}  \frac{1}{m} \left( \int_X \chi_{u^{-1}C}  \, d\mu + \int_X \chi_{u^{-2}C}  \, d\mu + \cdots  
		+ \int_X \chi_{u^{-m}C}   \, d\mu \right)
	\\&= \liminf_{m \to \infty} \frac{1}{m} \Bigl( \mu( u^{-1} C ) +  \mu( u^{-2} C ) + \cdots  
		+  \mu( u^{-m} C ) \Bigr)
	\\&= \liminf_{m \to \infty} \frac{1}{m} \Bigl( \mu(C) +  \mu(C)  + \cdots  +  \mu(C) \Bigr)
	\\&= \mu(C)
	\\&< \infty
	, \end{align*}
so $\rho_C \in \LL1(X, \mu)$.

It is easy to see that $\rho_C$ is $u$-invariant \csee{rhoCinvtEx}, so the Moore Ergodicity Theorem \pref{MooreErgBasicThm} implies that $\rho_C$ is constant (a.e.).
Also, from \pref{rhoC>0}, we know that the constant is not~$0$. Therefore, we have a nonzero constant function that is in $\LL1(X,\mu)$, which tells us that $\mu(X)$ is finite.
\end{proof}

Now, to begin our discussion of the proof of \cref{DaniMargulisUnipReturns}, we introduce a bit of terminology and notation, and make some simple observations. First of all, let us restate the result by using the Mahler Compactness Criterion \pref{MahlerCpct}, and also replace the discrete times $\{1,2,3,\ldots,m\}$ with a continuous interval $[0,T]$. \Cref{SpendTimeCpctMahlerEx} shows that this new version implies the original.

\begin{defn} \label{UnimodLattDefn}
For any $\integer$-lattice~$\Zlatt$ in~$\real^n$, there is some $g \in \GL(n,\real)$, such that $\Zlatt = g \integer^n$. We say $\Zlatt$ is \defit[Z-lattice@$\integer$-lattice!unimodular]{unimodular} if $\det g = \pm 1$.
\end{defn}

\begin{thm}[(restatement of \cref{DaniMargulisUnipReturns})] \label{SpendTimeCpctMahler}
Suppose
\noprelistbreak
	\begin{itemize}
	\item $\{u^t\}$ is a one-parameter unipotent subgroup of\/ $\SL(n,\real)$,
	\item $\Zlatt$ is a unimodular $\integer$-lattice in\/~$\real^n$,
	and
	\item \nindex{$\leb$ = Lebesgue measure on~$\real$}$\leb$ is the usual Lebesgue measure\/ \textup(i.e., length\/\textup) on\/~$\real$.
	\end{itemize}
Then there exists a neighborhood~$\open$ of\/~$0$ in\/~$\real^n$, such that
	$$ \liminf_{T \to \infty} \frac{\leb \Bigl( \bigset{ t \in [0,T] }{ u^t \Zlatt \cap \open = \{0\} } \Bigr)}{T} > 0 .$$
\end{thm}

\begin{notation} \label{PrimVecDefn}
Suppose $W$ is a discrete subgroup of~$\real^n$. 
\noprelistbreak
	\begin{itemize}
	\item A vector $w \in W$ is \defit[primitive vector]{primitive} in~$W$ if $\lambda w \notin W$, for $0 < \lambda < 1$.
	\item Let \nindex{$\prims{W}$ = $\{\text{primitive vectors in discrete subgroup~$W$ of $\real^n$}\}$}$\prims{W}$ 
	be the set of primitive vectors in~$W$.
	\item Let 
	\nindex{$\primplus{W}$ = set of representatives of $\prims{W}/\{\pm1\}$}$\primplus{W} \subseteq \prims{W}$ be a set of representatives that contains either~$w$ or~$-w$, but not both, for every $w \in \prims{W}$.
	(Note that $\prims{W} = -\prims{W}$; see \cref{PrimVecIffEx}.)
	\end{itemize}
\end{notation}

For simplicity, let us assume now that $n = 2$ (see \cref{PfUnipOrbitsReturnSect} for a discussion of the general case).

\begin{lem} \label{SmallPrimsInUnimodLatt} \ 
\noprelistbreak
	\begin{enumerate}
	\item \label{SmallPrimsInUnimodLatt-1}
There is a neighborhood~$\open_1$ of\/~$0$ in\/~$\real^2$, such that if $W$ is any unimodular $\integer$-lattice in\/~$\real^2$, then\/ $\# \!\left(\primplus{W} \cap \open_1\right) \le 1$.
	\item \label{SmallPrimsInUnimodLatt-compare}
Given any neighborhood~$\open_1$ of\/~$0$ in\/~$\real^2$, and any $\epsilon > 0$, there exists a neighborhood~$\open_2$ of\/~$0$ in\/~$\real^2$, such that if $x \in \real^2$, and\/ $[a,b]$ is an interval in\/~$\real$, such that there exists $t \in [a,b]$ with $u^t x \notin \open_1$, then
		$$ \leb \left( \bigset{ t \in [a,b] }{ u^t x \in \open_2 } \right) 
		\le \epsilon \, \leb \left( \bigset{ t \in [a,b] }{ u^t  x \in \open_1 } \right) . $$
	\end{enumerate}
\end{lem}

\begin{figure}[ht]
\begin{center}
\includegraphics{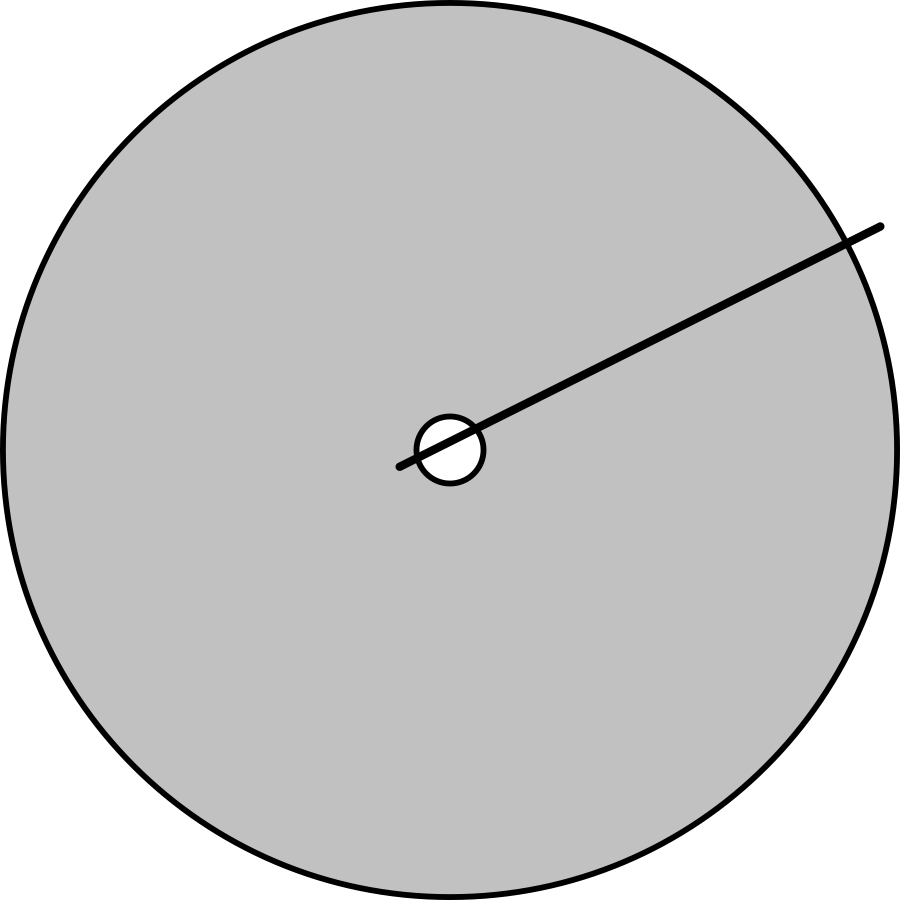}
\end{center}
\caption{Part~\pref{SmallPrimsInUnimodLatt-compare} of \cref{SmallPrimsInUnimodLatt}: Any line segment that reaches the boundary of the large disk has only a small fraction of its length inside the tiny disk.}
\label{LineSegSmallTimeInSmallDisk}
\end{figure}
%texpreamble
%(" \usepackage[T1]{fontenc}
% \usepackage[expert, T1, nofontinfo, mylucidascale]{mylucidabr}
% \usepackage{amsmath}
% ");
%defaultpen(  fontcommand("\normalfont") + fontsize(10) ); 
%
%from graph access *;
%
%size(0,1.5inch);
%
%real R = 2, r = 0.15;
%
%//xaxis(-R-r,R+r);
%//yaxis(-R-r,R+r);
%
%fill( (R,0)..(0,R)..(-R,0)..(0,-R)..cycle, gray(0.8) );
%draw( (R,0)..(0,R)..(-R,0)..(0,-R)..cycle, linewidth(0.7) );
%fill( (r,0)..(0,r)..(-r,0)..(0,-r)..cycle, white );
%draw( (r,0)..(0,r)..(-r,0)..(0,-r)..cycle, linewidth(0.7) );
%
%draw( (-3/2*r,-r/2)--(R-r/2, 0.5*R ) , linewidth(1) );

\begin{proof} 
\pref{SmallPrimsInUnimodLatt-1}~A unimodular $\integer$-lattice in~$\real^2$ cannot contain two linearly independent vectors of norm less than~$1$ \csee{UnimodNot2Short}.

\pref{SmallPrimsInUnimodLatt-compare} Note that $u^t x$ moves at constant velocity along a straight line \csee{UnipOrbIsLine}. So we simply wish to choose~$\open_2$ small enough that every line segment that reaches the boundary of~$\open_1$ has only a small fraction of its length inside~$\open_2$ \ccf{LineSegSmallTimeInSmallDisk}.

By making $\open_1$ smaller, there is no harm in assuming it is a disk centered at~$0$. Let $R$ be the radius of~$\open_1$, and let $\open_2$ be a disk of radius~$r$ centered at~$0$, with $r$ small enough that
	$$ \frac{2r}{R - r} < \epsilon .$$
Then, for any line segment~$L$ that reaches both $\open_2$ and the boundary of~$\open_1$, we have:
	\begin{itemize}
	\item the length of $L \cap \open_2$ is $\le$ the diameter $2r$ of $\open_2$,
	and
	\item the length of $L \cap \open_1$ is $\ge$ the distance $R - r$ from $\bdry \open_1$ to $\bdry \open_2$.
	\end{itemize}
Therefore, the segment of~$L$ that is in~$\open_2$ has length less than $\epsilon$ times the length of the segment that is~$\open_1$  \ccf{LineSegSmallTimeInSmallDisk}.
\end{proof}

\begin{proof}[\mathversion{bold}\bf Proof of \cref{SpendTimeCpctMahler} when $n = 2$]
Let $\open_1$ and $\open_2$ be as described in \cref{SmallPrimsInUnimodLatt}, with $\epsilon = 1/2$.
We may assume $\open_1$ and~$\open_2$ are convex, that they are small enough that they contain no nonzero elements of~$\Zlatt$, and that $\open_2 \subseteq \open_1$.

Fix $T \in \real^+$. For each $x \in \primplus{\Zlatt}$, and $k = 1,2$, let 
	$$ I^k_x = \{\, t \in [0,T] \mid x u^t  \in \open_k \,\} .$$
Since $\open_k$ is convex, and $u^t x$ traces out a line \csee{UnipOrbIsLine}, we know that $I^k_x$ is an interval (possibly empty). Note that:
	\begin{enumerate}
	\item  from \fullcref{SmallPrimsInUnimodLatt}{compare} (and the fact that $\epsilon = 1/2$), we see that 
	$ \leb( I^2_x ) \le \frac{1}{2} \leb( I^1_x )  $,
	and
	\item  from  \fullcref{SmallPrimsInUnimodLatt}{1}, we see that $I^1_{x_1}$ is disjoint from $I^1_{x_2}$ whenever $x_1 \neq x_2$.
	\end{enumerate}
Therefore
	\begin{align*}
	 \leb \bigl( \{\, t \in [0,T] \mid u^t \Zlatt  \cap \open_2 & \neq \{0\} \bigr) 
	= \sum_{x \in \primplus{\Zlatt}} \leb( I^2_x ) 
	\le \sum_{x \in \primplus{\Zlatt}} \frac{\leb( I^1_x )}{2}
	\\&= \frac{1}{2} \leb \left(  \bigcup_{x \in \primplus{\Zlatt}} I^1_x  \right)
	\le \frac{1}{2} \leb \bigl(  [0,T] \bigr)
	=  \frac{T}{2} 
	. \end{align*}
So, passing to the complement, we have
	\begin{align*}
	\leb \bigl( \{\, t \in [0,T] \mid u^t \Zlatt \cap \open_2 = \{0\} \bigr) 
	\ge \frac{T}{2} 
	. & \qedhere \end{align*}
\end{proof}

Unfortunately, \cref{SpendTimeCpctMahler} is not nearly as easy to prove when $n > 2$, because two basic complications arise.
	\begin{enumerate}
	\item The first difficulty is that the $u^t$-orbit of a vector is usually not a straight line (contrary to \cref{UnipOrbIsLine} for $n = 2$). However, the coordinates of $u^t x$ are always polynomials of bounded degree \csee{UnipOrbsArePoly}, so, for any fixed vector~$x$, 
	$$ \text{the function $\| u^t x \|^2$ is a polynomial in~$t$} $$
and the degree of this polynomial is bounded (independent of~$x$).
Therefore, it is easy to prove that the appropriate analogue of \fullcref{SmallPrimsInUnimodLatt}{compare} holds even if $n > 2$ \csee{PolySpendsLittleTimeSmall}, so the nonlinearity is not a major problem.

	\item A much more serious difficulty is the failure of \fullref{SmallPrimsInUnimodLatt}{1}: if $n > 2$, then a unimodular lattice in~$\real^n$ may have two  linearly independent primitive vectors that are very small \csee{UnimodLattCanHaveLotsSmallVectors}. This means that the sets $I^2_x$ in the above proof may not be disjoint, which is a major problem. It is solved by looking at not only single vectors, but at larger sets of linearly independent vectors. More precisely, we look at the subgroups generated by sets of small vectors in~$u^t \Zlatt$. These subgroups can intersect in rather complicated ways, and sorting this out requires a study of chains of these subgroups (ordered by inclusion) and a rather delicate proof by induction. 
Although the proof is completely elementary, using only some observations about polynomial functions, it is very clever and intricate. The main idea is presented in \cref{PfUnipOrbitsReturnSect}.
	\end{enumerate}

\begin{exercises}

\item \label{SplitDivergeSL2REx}
Suppose $a \in \SL(2,\real)$, and there exist linearly independent vectors $v_+$ and~$v_-$ in~$\real^2$, such that 
	$$ \text{$a^n v_+ \to 0$ as $n \to +\infty$ \ and \ $a^n v_- \to 0$ as $n \to -\infty$} . $$
Show $\exists$~$x \in \SL(2,\real)/\SL(2,\integer)$, such that $\{\, n \in \integer \mid a^n x \in C \,\}$ is finite, for every compact subset~$C$ of~$\SL(2,\real)/\SL(2,\integer)$.
\hint{There exists $g \in \SL(2,\real)$ that takes the two standard basis vectors of~$\real^2$ to vectors that are scalar multiples of $v_+$ and~$v_-$.}

\item \label{UnipToInftyInR2}
Let $u = \left[\begin{smallmatrix} 1 & 1 \\ 0 & 1 \end{smallmatrix}\right]$. For every $v \in \real^2$, show that either 
	\begin{itemize}
	\item $u^n v = v$ for all $n \in \integer$,
	or
	\item $\|u^n v\| \to \infty$ as $n \to \infty$.
	\end{itemize}

\item
Generalize the preceding exercise to $\SL(n,\real)$:
\par
Let $u$ be any unipotent element of $\SL(n,\real)$. For every $v \in \real^n$, show that either 
	\begin{itemize}
	\item $u^n v = v$ for all $n \in \integer$,
	or
	\item $\|u^n v\| \to \infty$ as $n \to \infty$.
	\end{itemize}
\hint{Each coordinate of $u^n v$ is a polynomial function of~$n$, and non-constant polynomials cannot be bounded.}

\item \label{2SmallNotIndep}
Suppose $v_1$ and~$v_2$ are linearly independent vectors in~$\integer^2$, and we have $g \in \SL(2,\real)$. Show that if $\|gv_1\| < 1$, then $\|g v_2\| > 1$.
\hint{Since $g \in \SL(2,\real)$, the area of the parallelogram spanned by the vectors $g v_1$ and~$g v_2$ is the same as the area of the parallelogram spanned by $v_1$ and~$v_2$, which is an integer.}

\item \label{MargulisSL2RNotLessR}
Near the end of the proof of \cref{MargulisRecurSL2R}, verify the assertion that $u^{n+N} x \integer^2$ does not contain any nonzero vectors of length less than~$r$.
\hint{If $\|w\| < r$, then $u^{n-1}v$ and $u^{-1} w$ are linearly independent vectors in $\open \cap u^{n+N-1} x \integer$.}

\item \label{rhoC>0Ex}
Prove \pref{rhoC>0}.
\hint{$X$ cannot be the union of countably many sets of measure~$0$.}

\item \label{rhoCinvtEx}
In the proof of \cref{SLNZISLATT}, verify (directly from the definition) that $\rho_C$ is $u$-invariant.

\item \label{UnimodLattDefnWellDef}
Show \cref{UnimodLattDefn} is well-defined. More precisely, given any $g_1 , g_2 \in \GL(n,\real)$, such that $g_1 \integer^n = g_2 \integer^n$, show 
	$$ \det g_1 \in \{ \pm1\}  \iff \det g_2 \in \{ \pm 1\} .$$

\item \label{DiscToContInCpct}
Assume 
	\begin{itemize}
	\item $u^t$ is a one-parameter unipotent subgroup of~$G$, 
	\item $x \in G/\Gamma$,
	and
	\item $C^*$ is a compact subset of $G/\Gamma$.
	\end{itemize}
Show that if 
	$$ \liminf_{T \to \infty} \frac{\leb \Bigl( \bigset{ t \in [0,T] }{  u^t x \in C^* } \Bigr)}{T} > 0 , $$
then there is a compact subset~$C$ of $G/\Gamma$, such that 
	$$ \liminf_{m \to \infty} \frac{\# \bigset{ k \in \{1,2,\ldots,m\} }{ u^k x \in C} }{m} > 0 .$$
\hint{Let $C = \bigcup_{t \in [0,1]} u^t C^*$.}

\item \label{SpendTimeCpctMahlerEx}
Show \cref{SpendTimeCpctMahler} implies \cref{DaniMargulisUnipReturns}.
\hint{Mahler Compactness Criterion \pref{MahlerCpct} and \cref{DiscToContInCpct}.}

\item \label{PrimVecIffEx}
Suppose $w$ is a nonzero element of a discrete subgroup~$W$ of~$\real^n$. Show the following are equivalent:
	\begin{enumerate}
	\item $w$ is primitive in~$W$.
	\item $\real w \cap  W = \{\integer w\}$.
	\item If $kw' = w$, for some $k \in \integer$ and $w' \in W$, then $k \in \{\pm 1\}$.
	\item \label{PrimVecIffEx-minus}
	$-w$ is primitive in~$W$.
	\end{enumerate}

\item \label{UnimodNot2Short}
Suppose $v$ and~$w$ are linearly independent vectors in a unimodular $\integer$-lattice in~$\real^2$. Show $\| v \| \cdot \|w \| \ge 1$.

\item \label{UnipOrbIsLine}
Show that if $x \in \real^2$, and $\{u^t\}$ is any nontrivial one-parameter unipotent subgroup of $\SL(2,\real)$, then $u^t x$ moves at constant velocity along a straight line.
\hint{Calculate the coordinates of $u^t x$ after choosing a basis so that $u^t = \left[\begin{smallmatrix} 1 & t \\ 0 & 1 \end{smallmatrix} \right] $.}

\item \label{UnipOrbsArePoly}
Given $n \in \integer^+$, show there is a constant~$D$, such that
if $x \in \real^n$, and $\{u^t\}$ is any one-parameter unipotent subgroup of $\SL(n,\real)$, then the coordinates of $u^t x$ are polynomial functions of~$t$, and the degrees of these polynomials are $\le D$.
\hint{We have $u^t = \exp( t v)$ for some $v \in \Mat_{n \times n}(\real)$. Furthermore, $v$~is nilpotent, because $u^t$ is unipotent, so the power series $\exp( t v)$ is just a polynomial.}

\item \label{PolySpendsLittleTimeSmall}
Given $R, D, \epsilon > 0$, show there exists $r > 0$, such that if
	\begin{itemize}
	\item $f(x)$ is a (real) polynomial of degree $\le D$,
	and
	\item $[a,b]$ is an interval in~$\real$, with $|f(t)| \ge R$ for some $t \in [a,b]$,
	\end{itemize}
then
	$$ \leb \bigl( \bigset{ t \in [a,b] }{ |f(t)| < r } \bigr) 
		\le \epsilon \, \leb \bigl( \bigset{ t \in [a,b] }{ |f(t)| < R } \bigr) . $$
\hint{If not, then taking a limit yields a polynomial of degree~$D$ that vanishes on a set of positive measure, but is $\ge R$ at some point.}

\item \label{UnimodLattCanHaveLotsSmallVectors}
For every $\epsilon > 0$, find a unimodular $\integer$-lattice~$\Zlatt$ in~$\real^n$ with $n-1$ linearly independent primitive vectors of norm~$\le \epsilon$.

\item \label{ChevalleyStabDiscreteZ}
Assume $G$ is defined over~$\rational$ (and connected).
Show there exist
	\begin{itemize}
	\item a finite-dimensional real vector space~$V$, 
	\item a vector~$v$ in~$V$,
	and
	\item a homomorphism $\rho \colon \SL(\ell,\real) \to \SL(V)$,
	\end{itemize}
such that 
	\begin{enumerate}
	\item $G = \Stab_{\SL(\ell,\real)}(v)^\circ$,
	and 
	\item $\rho \bigl( \SL(\ell,\integer) \bigr) v$ is discrete.
	\end{enumerate}
\hint{See the hint to \cref{ChevalleyStabEx}, and choose $v$ to be the exterior product of polynomials with integer coefficients.}

\item \label{ProperMapToSLnR/SLnZPfEx}
Show that if $G$ is defined over~$\rational$, then the natural embedding $G/G_{\integer} \hookrightarrow \SL(\ell,\real)/\SL(\ell,\integer)$ is a proper map.
\hint{Use \cref{ChevalleyStabDiscreteZ}.}

\item \label{GZLattSimpleEx}
Prove \cref{arith->latt} under the additional assumption that $G$ is simple.
\hint{The natural embedding $G/G_{\integer} \hookrightarrow \SL(\ell,\real)/\SL(\ell,\integer)$ is a proper map \csee{ProperMapToSLnR/SLnZPfEx}, so the $G$-invariant measure on $G/G_{\integer}$ provides a $G$-invariant measure~$\mu$ on $\SL(\ell,\real)/\SL(\ell,\integer)$, such that all compact sets have finite measure. The proof of \cref{SLNZISLATT} (with $u \in G$) implies that $\mu$ is finite.}

\item \label{GZLattEx}
Prove \cref{arith->latt} (without assuming that $G$ is simple).
\hint{You may assume \Cref{MooreErgNonsimple} (without proof). This provides a version of the Moore Ergodicity Theorem for groups that are not simple.}

\end{exercises}

\section{Proof that unipotent orbits return to a compact set} \label{PfUnipOrbitsReturnSect}

The proof of \cref{SpendTimeCpctMahler} is rather complicated. To provide the gist of the argument, while eliminating some of the estimates that obscure the main ideas, we prove only \cref{MargulisUnipReturns}, which is a qualitative version of the result. (The quantitative conclusion in \cref{SpendTimeCpctMahler} makes additional use of observations similar to \fullcref{SmallPrimsInUnimodLatt}{compare} and \cref{PolySpendsLittleTimeSmall}.) 
This section is optional, because none of the material is needed elsewhere in the book.

By the Mahler Compactness Criterion (and an appropriate modification of \cref{DiscToContInCpct}), it suffices to prove the following statement:

\begin{thm}[(restatement of \cref{MargulisUnipReturns})] \label{UnipOrbMissNeigh}
Suppose
	\begin{itemize}
	\item $\{u^t\}$ is a one-parameter unipotent subgroup of\/ $\SL(n,\real)$,
	and
	\item $\Zlatt$ is a unimodular $\integer$-lattice in\/~$\real^n$.
	\end{itemize}
Then there exists a neighborhood~$\open$ of\/~$0$ in\/~$\real^n$, such that
	$$ \text{$\bigset{ t \in \real^+ }{ u^t \Zlatt \cap \open = \{0\} }$ is unbounded.} $$
\end{thm}

\begin{defn} \label{d(W)Defn}
Suppose 
\noprelistbreak
	\begin{itemize}
	\item $W$ is a discrete subgroup of~$\real^n$, 
	and
	\item $k$ is the dimension of the linear span $\langle W \rangle$ of~$W$.
	\end{itemize}
We make the following definitions:
\noprelistbreak
	\begin{enumerate}
	\item We define an inner product on the exterior power $\bigwedge\!\!^k \, \real^n$ by declaring $\{ \varepsilon_{i_1} \wedge \varepsilon_{i_2} \wedge  \cdots \wedge \varepsilon_{i_k} \}$ to be an orthonormal basis, where $\{\varepsilon_1,\ldots,\varepsilon_n\}$ is the standard basis of~$\real^n$.
	\item \label{d(W)Defn-d(W)}
 Since $\bigwedge\!\!^k \, W$ is cyclic \csee{WedgeWCyclicEx}, it has a generator $w_1 \wedge \cdots \wedge w_k$ that is unique up to sign, and we define%
	$$ d(W) = \| w_1 \wedge \cdots \wedge w_k\| 
	\nindex{$d(W) = \| w_1 \wedge \cdots \wedge w_k\|$}. $$
(However, by convention, we let $d \bigl( \{0\} \bigr) = 1$.)
	\end{enumerate}
\end{defn}

\begin{rem} 
If $W$ is the cyclic group generated by a nonzero vector $w \in \real^n$, then it is obvious that $d(W) = \|w\|$. Therefore, \fullcref{d(W)Defn}{d(W)} presents a notion that generalizes the norm of a vector.
\end{rem}

The following generalization of \cref{UnipOrbsArePoly} is straightforward to prove \csee{d(W)isPolyEx}.

\begin{lem} \label{d(W)isPoly}
Suppose
\noprelistbreak
	\begin{itemize}
	\item $\{u^t\}$ is a one-parameter unipotent subgroup of\/ $\SL(n,\real)$,
	and
	\item $W$ is a discrete subgroup of\/~$\real^n$.
	\end{itemize}
Then $d(u^t W)^2$ is a polynomial function of~$t$, and the degree of this polynomial is bounded by a constant~$D$ that depends only on~$n$.
\end{lem}

\Cref{d(W)isPoly} allows us to make good use of  the following two basic properties of polynomials of bounded degree.  (See \cref{ConstsForPolyDiv-alphaEx,ConstsForPolyDiv-betaEx} for the proofs.) The first follows from the observation that polynomials of bounded degree form a finite-dimensional real vector space, so any closed, bounded subset is compact. The second uses the fact that nonzero polynomials of degree~$D$ cannot have more than~$D$ zeroes.

\begin{lem} \label{ConstsForPolyDiv}
Suppose $D \in \integer^+$, $\epsilon > 0$, and $f$ is any real polynomial of degree $\le D$. 
Then there exists $C > 1$, depending only on~$D$ and~$\epsilon$, such that, for all $T,\tau > 0$:
\noprelistbreak
	\begin{enumerate}
	\item \label{ConstsForPolyDiv-alpha}
	If $f(s) \ge \tau$ for some $s \in [0,T]$, and $|f(T)| \le \tau/C$, then there exists $t \in [0, \epsilon T]$, such that $|f(T + t)| = \tau/C$.
	\item \label{ConstsForPolyDiv-beta}
 	If $|f(s)| \le \tau$ for all $s \in [0,T]$, and $f(T) = \tau$,
	then there exists $T_1 \in [T, 4^D T]$, such that 
		$$ \text{$\tau/C \le |f(t)| \le \tau C $ 
		\ for all $t \in [T_1, 2T_1]$} .$$
	\end{enumerate}
\end{lem}

\begin{notation} 
Suppose  $\Zlatt$ is a $\integer$-lattice in~$\real^n$.
	\begin{itemize}
	\item A subgroup~$W$ of~$\Zlatt$ is \defit[full!subgroup of a $\integer$-lattice]{full} if it is the intersection of~$\Zlatt$ with a vector subspace of~$\real^n$. (This is equivalent to requiring $\Zlatt/W$ to be torsion-free.)
	\item Let \nindex{$\subgrps(\Zlatt) = \{\text{full, nontrivial subgroups of $\Zlatt$}\}$}$\subgrps(\Zlatt)$ be the collection of all full, nontrivial subgroups of~$\Zlatt$, partially ordered by inclusion.
%	\item For any totally ordered subset~$S$ of~$\subgrps(\Zlatt)$, we let
%		$$ \subgrps(\Zlatt;S) = \bigl\{\, W \in \subgrps(\Zlatt) \smallsetminus S \mid 
%			\text{$S \cup \{W\}$ is totally ordered} \,\bigr\} .$$
	\item For $W \subseteq \Zlatt$, we let $\langle W \rangle\!_{\Zlatt}$ be the (unique) smallest full subgroup of~$\Zlatt$ that contains~$W$. In other words, $\langle W \rangle\!_{\Zlatt} = \langle W \rangle \cap \Zlatt$.
	\end{itemize}
\end{notation}

The following simple observation uses full subgroups of~$\Zlatt$ to provide a crucial lower bound on the norms of vectors \csee{FullSubgrpLowerBdEx}:

\begin{lem} \label{FullSubgrpLowerBd}
If $W \in \subgrps(\Zlatt)$
	and
	 $v \in \Zlatt \smallsetminus W$,
then\/ 
	$\displaystyle \|v\| \ge \frac{d \bigl( \langle W, v \rangle\!_{\Zlatt}\bigr)}{d \bigl( \langle W \rangle\!_{\Zlatt}\bigr)} $.
\end{lem}

We can now prove \cref{UnipOrbMissNeigh}. However, to avoid the need for a proof by induction, we assume $n = 3$. 

\begin{proof}[\mathversion{bold}\bf Proof of \cref{UnipOrbMissNeigh} when $n = 3$]
It is easy to see that 
	$$ \text{$\{\, W \in \subgrps(\Zlatt) \mid d(W) < 1 \,\}$ is finite} $$
\csee{DisFiniteEx}. Hence, there exists $\tau > 0$, such that
	$$ \text{$d(W) > \tau$, for all $W \in \subgrps(\Zlatt)$.} $$
Let:
\noprelistbreak
	\begin{itemize}
	\item $D$ be the constant provided by \cref{d(W)isPoly},
	\item $\epsilon = 4^{-(D+1)}$,
	and
	\item $C$ be the constant provided by \cref{ConstsForPolyDiv}.
	\end{itemize}
Given $T > 0$, it suffices to find $R \ge 0$, such that $\|u^{T+R} v\| \ge \tau/C^2$ for all nonzero $v \in \Zlatt$.

Let 
	$$ \mathcal{D} = \{\, W  \in \subgrps(\Zlatt) \mid d ( u^T W ) < \tau/C\,\} .$$
We assume $\mathcal{D} \neq \emptyset$ (otherwise, we could let $R = 0$). For each $W \in \mathcal{D}$, \fullcref{ConstsForPolyDiv}{alpha} implies 
	$$ \text{there exists $t_{W} \in [0,  \epsilon T]$, such that $d \bigl( u^{T + t_W} W \bigr) = \tau/C$} .$$
By choosing $t_W$ minimal, we may assume
	$$ \text{$d(u^{T + t}W) < \tau/C$ \ for all $t \in [0,  t_W)$} .$$
Since $\mathcal{D}$ is finite \csee{DisFiniteEx}, we may 
	$$ \text{fix some $W^+ \in \mathcal{D}$ that maximizes $t_{W}$.} $$
 From \fullcref{ConstsForPolyDiv}{beta}, we see that there exists 
 	$$T_1 \in [t_{W^+}, 4^D t_{W^+}] \subseteq \left[t_{W^+}, \frac{T}{2} \right] ,$$
such that 
	$$ \text{$\tau/C^2 \le d ( u^{T+t} W^+ ) \le \tau C^2$
	\ for all $t \in [T_1, 2T_1]$.} $$

Since $\dim \langle \Zlatt \rangle = n = 3$, we know $\dim \langle W^+ \rangle$ is either $1$ or~$2$. To be concrete, let us assume it is~$2$. (See \cref{dimWplus} for the other case.) Then, for any $v \in \Zlatt \smallsetminus W^+$, we have $\langle W^+, v \rangle\!_{\Zlatt} = \Zlatt$, so \cref{FullSubgrpLowerBd} implies $\| u^{T + t} v \| \ge 1/\tau$ for all $t \in [T_1, 2T_1]$. Hence, it is only the vectors in~$W^+$ that can be small anywhere in this interval.

Therefore, we may assume there is some nonzero $v_0 \in W^+$, such that $\| u^{T+T_1} v_0 \| < \tau/C^2$. There is no harm in assuming that $\integer v_0$ is a full subgroup of~$\Zlatt$. Then, since $T_1 \ge t_{W^+}$, the maximality of~$t_{W^+}$ implies $\| u^{T+s} v_0 \| \ge \tau/C$ for some $s \in [0,T_1]$.
Therefore, \fullcref{ConstsForPolyDiv}{alpha} provides some $t \in [T_1, 2T_1]$, such that $\| u^{T+t} v_0 \| = \tau/C^2$.
Now, for any nonzero $v \in \Zlatt$,
	$$ \text{ either \ 
	$\langle v \rangle\!_{\Zlatt} = \langle v_0 \rangle\!_{\Zlatt}$,
	\ or \ 
	$\langle v_0, v \rangle\!_{\Zlatt} = W^+$,
	\ or \ 
	$\langle v, W^+ \rangle\!_{\Zlatt} = \Zlatt$
	}. $$
In each case, we see (by using \cref{FullSubgrpLowerBd} in the latter two cases) that $\|u^{T+t} v\| \ge \tau/C^2$ (if $\tau \le 1$).
\end{proof}

\begin{exercises}

\item 
Show that \cref{UnipOrbMissNeigh} is a corollary of \cref{SpendTimeCpctMahler}.

\item \label{UnipReturnsToCpctSet}
Use \cref{UnipOrbMissNeigh} (and \emph{not} \cref{DaniMargulisUnipReturns} or \cref{SpendTimeCpctMahler}) to show that if
$\{u^t\}$, $X$, and~$x$ are as in \cref{DaniMargulisUnipReturns}, then there exists a compact subset~$K$ of~$X$, such that
	$$ \text{$\bigset{ t \in \real^+ }{ u^t x \in K }$ is unbounded.} $$

\item \label{WedgeWCyclicEx}
In the notation of \cref{d(W)Defn}, show $\bigwedge\!\!^k W$ is cyclic.
\hint{If $\{w_1,\ldots,w_k\}$ generates~$W$, then $\bigwedge\!\!^k W$ is generated by $w_1 \wedge w_2 \wedge \cdots \wedge w_k$.}

\item Suppose 
	\begin{itemize}
	\item $W$ is (nontrivial) discrete subgroup of~$\real^n$, 
	and
	\item $M \in \SO(n)$.
	\end{itemize}
Show $d(MW) = d(W)$.

\item \label{d(W)isPolyEx}
Prove \cref{d(W)isPoly}.

\item \label{ConstsForPolyDiv-alphaEx}
Prove \fullcref{ConstsForPolyDiv}{alpha}.
\hint{Since rescaling does not change the degree of a polynomial, we may assume $T = \tau = 1$. If $C$ does not exist, then taking a limit results in a polynomial of degree~$\le D$ that is~$1$ at some point of $[0,1]$, but vanishes on all of $[1,1+\epsilon]$.}

\item \label{ConstsForPolyDiv-betaEx}
Prove \fullcref{ConstsForPolyDiv}{beta}.
\hint{Assume, without loss of generality, that $T = \tau = 1$. The polynomials of degree $\le D$ that are $\le 1$ on $[0,1]$ form a compact set, so they are uniformly bounded by some constant on $[1, 4^{D+1}]$. For $T_1 \in \{1, 4, \ldots, 4^D\}$, the intervals $[T_1,2T_2]$ are pairwise disjoint. If $f$ is not bounded away from~$0$ on any of these intervals, then taking a limit results in a nonzero polynomial of degree~$\le D$ that vanishes at $D+1$ distinct points.}

\item \label{FullSubgrpLowerBdEx}
Prove \cref{FullSubgrpLowerBd}.
\noprelistbreak % @@@
\hint{This is easy if $W$ is generated by scalar multiples of the standard basis vectors of~$\real^k$, and $v \in \real^{k+1}$.}

\item \label{WedgeDiscrete}
Show that if $\Zlatt$ is a discrete subgroup of~$\real^n$, and $1 \le k \le n$, then $\bigwedge\!\!^k \Zlatt$ is a discrete subset of  $\bigwedge\!\!^k \real^n$.
\hint{By choosing an appropriate basis, you can assume $\Zlatt \subseteq \integer^n$.}

\item \label{DisFiniteEx}
 %In the notation of the proof of \cref{MargLemTechResult}, show $\mathcal{D}$ is finite.
Assume
		\begin{itemize}
		\item $\Zlatt$ is a $\integer$-lattice in~$\real^n$,
		and
		\item $\delta > 0$.
		\end{itemize}
	Show there are only finitely many full subgroups of~$\Zlatt$, such that $d(W) < \delta$.
	\hint{\Cref{WedgeDiscrete}. (If $W_1$ and~$W_2$ are two different $k$-dimensional subspaces of~$\real^n$, then $\bigwedge\!\!^k W_1 \neq \bigwedge\!\!^k W_2$.)}
	
\item \label{dimWplus}
Complete the proof of \cref{UnipOrbMissNeigh} in the special case where $\dim \langle W^+ \rangle = 1$ (and $n = 3$). 
\hint{If there exist $v \in \Zlatt \smallsetminus W^+$ and $t \in [T_1,2T_1]$, such that $\|u^{T+t} v\| < 1/C$, then $d \bigl( u^{T + R}\langle W^+, v \rangle\!_{\Zlatt} \bigr) = \tau/C$ for some $R \in [T_1,2T_1]$.}

\end{exercises}

\begin{notes}

See \cite[\S1]{Borel-IntroGrpArith} or \cite[\S4.2]{PlatonovRapinchukBook} for more information on Siegel sets in $\SL(n,\real)$, and the proof of \cref{SLNZISLATT} that appears in \cref{SiegelSLnZSect,SLNZISLATTSiegelPfSect}.

A brief discussion of the connection with the reduction theory of positive-definite quadratic forms can be found in \cite[\S2, pp.~20--24]{Borel-IntroGrpArith}.

See  \cite[Prop.~3.12, p.~129]{PlatonovRapinchukBook} for a proof of \cref{IwasawaDecompSLnR}. A generalization to other semisimple groups will be stated in \cref{IwasawaDecomp}.

The clever proof in \cref{SLNZISLATTSlickSect} is by G.\,A.\,Margulis \cite[Rem.~3.12(II)]{Margulis-LieGrpsErgThy}.

\Cref{MargulisUnipReturns} is due to G.\,A.\,Margulis \cite{Margulis-OnActionUnip}. 
(\Cref{PfUnipOrbitsReturnSect} is adapted from the nice exposition in \cite[Appendix, pp.~162--173]{DaniMargulis-elementary}, where all details can be found.)
The result had been announced previously (without proof), and
J.\,Tits \cite[p.~59]{Tits-WorkOfMargulis} commented that:
	\par\smallbreak\begin{narrower}\begin{narrower}
	\noindent \sffamily
\llap{``}For a couple of years, Margulis' proof remained unpublished and every attempt
by other specialists to supply it failed. When it  finally appeared \dots, the proof
came as a great surprise, both for being rather short and using no sophisticated
technique: it can be read without any special knowledge and gives a good idea of
the extraordinary inventiveness shown by Margulis throughout his work\zz.''
	\par\end{narrower}\end{narrower}\smallbreak

The quantitative version stated in \cref{DaniMargulisUnipReturns} is due to S.\,G.\,Dani \cite{Dani-LemmaOfMargulis}. 
See \cite{Kleinbock-QuantitativeNondivergence} for a recent generalization, and applications to number theory.

\end{notes}

%!TEX root = IntroArithGrps.tex

\part{Important Concepts} \label{ConceptsPart}

 %!TEX root = IntroArithGrps.tex

\mychapter{Real Rank} \label{RrankChap}

\prereqs{none.}

%Although real rank and parabolic parabolic subgroups are a fundamental notion in the theory of semisimple Lie groups, it will not play a major role in most of this book.
%However, it is important in our discussion of lattices that are Gromov hyperbolic \csee{GromovHyperGrpsSect}, and they provide important background for the notion of $\rational$-rank (and parabolic $\rational$-subgroups), which are absolutely crucial for the construction of coarse fundamental domains in \cref{ReductionChap}. 

%\begin{assump}
% Assume that $G$ is a closed subgroup of $\SL(\ell,\real)$,
%for some~$\ell$. (The definitions and results of this
%chapter are independent of the particular embedding chosen.)
% \end{assump}

\section{\texorpdfstring{$\real$}{R}-split tori and \texorpdfstring{$\real$}{R}-rank}

\begin{defn} \label{RsplitDefn}
  A closed, connected subgroup~$T$ of~$G$ is a \defit{torus} if it is diagonalizable
over~$\complex\mk$; that is, if there exists $g \in
\GL(n,\complex)$, such that $g T g^{-1}$ consists entirely
of diagonal matrices. A torus is 
\defit[torus!$\real$-split]{$\real$-split} if it is diagonalizable
over~$\real\mk$; that is, if $g$~may be chosen to be in $\GL(n,\real)$.
 \end{defn}

%Because commuting diagonalizable matrices can be
%simultaneously diagonalized \csee{SimultDiag}, we have the following:
%
%\begin{prop} \label{torus<>diag}
% A subgroup~$T$ of~$G$ is a torus \index{torus} if and
%only if
% \begin{itemize}
% \item $T$ is closed and connected,
% \item $T$ is abelian, and
% \item every element of~$T$ is semisimple.
% \end{itemize}
% \end{prop}

\begin{egs} \ 
\noprelistbreak 
 \begin{enumerate}
 
 \item Let $A$ be the identity component of the group of
diagonal matrices in $\SL(n,\real)$. Then $A$ is obviously an $\real$-split
torus. 
%However, it is homeomorphic to $\real^{n-1}$, so it
%is not a torus in the usual topological sense.

\item $\SO(1,1)^\circ$ is an $\real$-split torus in $\SL(2,\real)$ \csee{SO(11)torus}.

 \item $\SO(2)$ is a torus in $\SL(2,\real)$ that is \emph{not} $\real$-split. It is diagonalizable over~$\complex$ \csee{SO2diag}, but not over~$\real$ \csee{SO2NotOverR}.
%Since every element of~$T$ is a normal linear
%transformation (that is, commutes with its transpose), we
%know from elementary linear algebra that every element of~$T$
%is diagonalizable.
% \item Generalizing the preceding, we see that $T = \SO(2)^n$
%is a torus in $\SL(2n,\real)$. Note that $T$
%is homeomorphic to the $n$-torus~$\torus^n$. (In fact, any
%compact torus subgroup is homeomorphic to a topological
%torus.) This is the motivation for the term
%\defit[torus subgroup]{torus}.
% \item The hypothesis that $T$ is abelian cannot be omitted
%from \cref{torus<>diag}. For example, every element of
%$\SO(n)$ is semisimple, but $\SO(n)$ is not abelian if $n
%\ge 3$.
 \end{enumerate}
 \end{egs}

\begin{warn}
 An $\real$-split torus is \emph{never} homeomorphic
to the topologist's torus~$\torus^n$ (except in the trivial case $n = 0$).
 \end{warn}
 
 \begin{rems} \ 
 \noprelistbreak
 	\begin{enumerate}

	\item If $T$ is an $\real$-split torus, then every element of~$T$ is hyperbolic \csee{hypelluniDefn}. In particular, no nonidentity element of~$T$ is elliptic or unipotent.

	\item  When $G$ is compact, every torus in~$G$ is isomorphic to $\SO(2)^n$, for some~$n$. This is homeomorphic to~$\torus^n$, which is the reason for the terminology ``torus\zz.''
 
	\end{enumerate}
 \end{rems}

%Recall that an element of~$G$ is
%\defit[hyperbolic!element of~$G$]{hyperbolic} if it is
%diagonalizable over~$\real$ \csee{hypelluniDefn}.
% In this
%terminology, we have the following:
%% analogue of \cref{torus<>diag}.
%
%\begin{prop}[\csee{RSplitTorus<>Ex}] \label{RSplitTorus<>}
% A subgroup~$T$ of~$G$ is an $\real$-split torus if and
%only if
%\noprelistbreak
% \begin{itemize}
% \item $T$ is closed and connected,
% \item $T$ is abelian, and
% \item every element of~$T$ is hyperbolic.
% \end{itemize}
% \end{prop}

It is a key fact in the theory of semisimple Lie groups that \emph{maximal} $\real$-split tori are conjugate:

\begin{thm}
If $A_1$ and~$A_2$ are maximal\/ $\real$-split tori in~$G$, then there exists $g \in G$, such that $A_1 = g A_2 g^{-1}$.
\end{thm}

This implies that all maximal $\real$-split tori have the same dimension, which is called the ``{real rank}'' (or ``$\real$-rank'') of~$G$, and is denoted $\Rrank G$:

 \begin{defn} \label{RrankDefn} \nindex{$\Rrank G$ = real rank of~$G$}%
\index{rank!R-@$\real$- or real}\term[R-@$\real$-!rank]{$\Rrank G$} is the dimension of a maximal $\real$-split torus~$A$ in~$G$. 
This is independent of both the choice of~$A$ and the choice of the embedding of~$G$ in $\SL(\ell,\real)$.
 \end{defn}

\begin{egs} \ \label{RrankEgs}
 \noprelistbreak
 \begin{enumerate}

 \item $\Rrank \bigl( \SL(n,\real) \bigr) = n-1$. (Let $A$
be the identity component of the group of all diagonal
matrices in $\SL(n,\real)$.) 

 \item We have
 $\Rrank \bigl( \SL(n,\complex) \bigr) = \Rrank \bigl(
\SL(n,\quaternion) \bigr) = n-1$. This is because only the
\emph{real} diagonal matrices remain
diagonal when $\SL(n,\complex)$ or $\SL(n,\quaternion)$ is embedded in
$\SL(2n,\real)$ or $\SL(4n,\real)$, respectively.

\item \label{RrankEgs-cpct}
$\Rrank G = 0$ if and only if $G$~is compact \csee{Rrank0Ex}.
 \end{enumerate}
 \end{egs}

\begin{prop} \label{Rrank(SOmn)}
 $\Rrank \SO(m,n)  = \min\{m,n\}$.
 \end{prop}

\begin{proof} 
Since $\SO(m,n)$ contains a copy of $\SO(1,1)^{\min\{m,n\}}$ \csee{SO11inSOmn}, and the identity component of this subgroup is an $\real$-split torus \ccf{SO(11)torus}, we have
	$$ \Rrank \SO(m,n) \ge \dim \bigl( \SO(1,1)^{\min\{m,n\}} \bigr)^\circ = \min\{m,n\}. $$

We now establish the reverse inequality. Let $A$ be a maximal $\real$-split torus. We may assume $A$ is nontrivial. (Otherwise $\Rrank \SO(m,n) = 0$, so the desired inequality is obvious.) Therefore, there is some nontrivial $a \in A$. Since $a$ is diagonalizable over~$\real$, and nontrivial, there is an eigenvector~$v$ of~$a$, such that $av \neq v$; hence, $av = \lambda v$ for some $\lambda \neq 1$. Now, if we let $\langle \cdot \mid \cdot \rangle_{m,n}$ be an $\SO(m,n)$-invariant bilinear form on~$\real^{m,n}$, we have
	$$ \langle v \mid v \rangle_{m,n}
	= \langle av \mid av \rangle_{m,n}
	= \langle \lambda v \mid \lambda v \rangle_{m,n}
	= \lambda^2 \langle v \mid v \rangle_{m,n}
	. $$
By choosing $a$ to be near~$e$, we may assume $\lambda \approx 1$, so $\lambda \neq -1$. Since, by assumption, we know $\lambda \neq 1$, this implies $\lambda^2 \neq 1$. So we must have $ \langle v \mid v \rangle_{m,n} = 0$;
 that is, $v$~is an \defit[isotropic vector]{isotropic} vector.
 Hence, we have shown that if the real rank is $\ge 1$, then there is an isotropic vector in $\real^{m+n}$. 
 
 By arguing more carefully, it is not difficult to see that if the real rank is at least~$k$, then there is a $k$-dimensional subspace of $\real^{m+n}$ that consists entirely of isotropic vectors \csee{SOmnIsotopic}. Such a subspace is said to be \defit[totally!isotropic]{totally isotropic}. The maximum dimension of a totally isotropic subspace is $\min\{m,n\}$ \csee{SOmnLargestIsotropic}, so we conclude that $\min\{m,n\} \ge \Rrank \SO(m,n)$, as desired.
 \end{proof}

\begin{rems} \label{RrankRems} \ 
\noprelistbreak
	\begin{enumerate}
	\item \label{RrankRems-Rrank(SUmn)}
	Other classical groups, not just $\SO(m,n)$, have the property that their real rank is the maximal dimension of a totally isotropic subspace.
More concretely, we have
 $$ %\Rrank \SO(m,n) = 
 \Rrank \SU(m,n) = \Rrank \Sp(m,n) = \min\{m,n\} .$$
 
	\item 
	The Mostow Rigidity Theorem \pref{MostowIso} will tell us that if $\Gamma$ is (isomorphic to) a lattice in both $G$ and~$G_1$, then $G^\circ$ is isomorphic to~$G_1'$, modulo compact groups. Modding out a compact subgroup does not affect the real rank \ccf{Rrank0Ex}, so this implies that the real rank of~$G$ is uniquely determined by the algebraic structure of~$\Gamma$.

\item \label{RrankRems-Rrank=CoverGamma}
Although it is not usually very useful in practice, we now state an explicit relationship between $\Gamma$ and $\Rrank G$. Let $S_r$ be the set of all elements~$\gamma$ of $\Gamma$, such that the centralizer $C_\Gamma(\gamma)$ is commensurable to a subgroup of the free abelian group $\integer^r$ of rank~$r$. Then it can be shown that
	$$ \Rrank G = \min \bigset{ r \ge 0 }{ \begin{matrix} \text{$\Gamma$ is covered by finitely} \\ \text{many translates of~$S_r$}  \end{matrix}} .$$
We omit the proof, which is based on the very useful (and nontrivial) fact that if $T$ is any maximal torus of~$G$, then there exists $g \in G$, such that $g T g^{-1} / (\Gamma \cap g T g^{-1})$ is compact.
\end{enumerate}
\end{rems}

\begin{exercises}

%\item \label{SimultDiag}
%Suppose $\mathcal{T} \subseteq \GL(n,F)$, such that:
%	\begin{itemize}
%	\item the matrices in $\mathcal{T}$ all commute with each other,
%	and, 
%	\item for all $t \in \mathcal{T}$, there exists $g \in \GL(n,F)$, such that $g t g^{-1}$ is a diagonal matrix. 
%	\end{itemize}
%Show the quantifiers can be reversed: there exists $g \in \GL(n,F)$, such that, for all $t \in \mathcal{T}$, $g t g^{-1}$ is a diagonal matrix. 
%\hint{$t$ is diagonalizable iff there is a basis of~$F^n$ that consists of eigenvectors of~$t$. If $s$ commutes with~$t$, then each eigenspace for~$t$ is $s$-invariant.}

\item \label{SO(11)torus}
Show that the identity component of $\SO(1,1)$ is an $\real$-split torus.
\hint{Let $g = \begin{Smallbmatrix} 1 & 1 \\ 1 & -1 \end{Smallbmatrix}$. Alternatively, note that each element of $\SO(1,1)$ is a symmetric matrix (hence, diagonalizable via an orthogonal matrix), and use the fact that any set of commuting diagonalizable matrices is simultaneously diagonalizable.}
%We have $\SO(1,1)^\circ = 
% \begin{bmatrix}
% \cosh t & \sinh t \\
% \sinh t & \cosh t
% \end{bmatrix}
%$,
%where $\cosh t = (e^t + e^{-t})/2$ and $\sinh t = (e^t
%- e^{-t})/2$, since $\cosh^2 t - \sinh^2 t = 1$).

\item \label{SO2diag}
 For $ g =
 \begin{Smallbmatrix}
 1 & -i \\
 1 & i
 \end{Smallbmatrix}$,
show every element of $g \SO(2) g^{-1}$ is diagonal.

\item \label{SO2NotOverR}
Show that $\SO(2)$ is not diagonalizable over~$\real$.
\hint{If $T$ is diagonalizable over~$\real$, then eigenvalues of the elements of~$T$ are real.}

\item \label{TorusAbel}
Show that every $\real$-split torus is abelian.

%\item \label{RSplitTorus<>Ex}
% Prove \cref{RSplitTorus<>}.
% \hint{\Cref{TorusAbel,SimultDiag}.}

\item Suppose
	\begin{itemize}
	\item $T$ is an $\real$-split torus in~$G$,
	and
	\item $A$ is a maximal $\real$-split torus in~$G$.
	\end{itemize}
Show that $T$ is conjugate to a subgroup of~$A$.
\hint{By considering dimension, it is obvious that $T$ is contained in some maximal $\real$-split torus of~$G$.}

\item \label{MaxRsplitZar}
Show that every maximal $\real$-split torus in~$G$ is almost Zariski closed.

\item \label{SO11inSOmn}
Assume $m \ge n$. Then $m + n \ge 2n$, so there is a natural embedding of $\SO(1,1)^n$ in $\SL(m+n,\real)$.
Show that $\SO(m,n)$ contains a conjugate of this copy of $\SO(1,1)^n$.
\hint{Permute the basis vectors.}
%\hint{Let $G = \SO(B)$, where
% $$ B = \sum_{i = 1}^n (x_{2i-1}^2 - x_{2i}^2 ) + \sum_{j=2n+1}^{m+n} x_j^2 ,$$
%so there is an obvious embedding of the $\real$-split torus $\bigl( \SO(1,1)^n \bigr)^\circ$ in~$G$:
% $$ \begin{bmatrix}
% \SO(1,1) \\
% & \SO(1,1) \\
% & & & \ddots \\
% & & & & \SO(1,1) \\
% & & & & & \Id_{m-n,m-n}
% \end{bmatrix}
% \subset G
% ,$$
% and note that $G$ is conjugate to $\SO(m,n)$, just by permuting the basis vectors.}

\item Prove, directly from \cref{RsplitDefn}, that if
$G_1$ is conjugate to~$G_2$ in $\GL(\ell,\real)$, then
$\Rrank(G_1) = \Rrank(G_2)$.

\item \label{Rrank0Ex}
Show $\Rrank G = 0$ if and only if $G$~is compact.
\hint{\Cref{SL2RinG,SSeltRem}\pref{SSeltRem-cpct}.}

\item \label{SOmnIsotopic}
Show that if $\Rrank \SO(m,n) = r$, then there is an $r$-dimensional subspace~$V$ of~$\real^{m+n}$, such that $\langle v \mid w \rangle_{m,n} = 0$ for all $v,w \in V$.
\hint{Because $A$ is diagonalizable over~$\real$, there is a basis $\{v_1,\ldots,v_{m+n}\}$ of~$\real^{m+n}$ whose elements are eigenvectors for every element of~$A$.
Since $\dim A = r$, we may assume, after renumbering, that for all $\lambda_1,\ldots,\lambda_r \in \real^+$, there exists $a \in A$, such that $av_i = \lambda_i v_i$, for $1 \le i \le r$. This implies $\langle v_1,\ldots,v_r \rangle$ is totally isotropic.}

\item \label{SOmnLargestIsotropic}
Show that if $V$ is a subspace of $\real^{m+n}$ that is totally isotropic for 
\text{$\langle \cdot \mid \cdot \rangle_{m,n}$}, % don't allow bad line break
then $\dim V \le \min\{m,n\}$.
\hint{If $v \neq 0$ and the last $n$ coordinates of~$v$ are~$0$, then $\langle v \mid v \rangle_{m,n} > 0$.}

\item Show $ \Rrank(G_1 \times G_2) =  \Rrank G_1 \  + \ \Rrank G_2$.

\item \label{Rrank0<>SL2R}
Show $\Rrank G \ge 1$ if and only if $G$~contains a subgroup that is isogenous to $\SL(2,\real)$.
\hint{\Cref{SL2RinG}.}

\item Show that $\Gamma$ contains a subgroup that is isomorphic to $\integer^r$, where $r = \Rrank G$.
\hint{You may assume the fact stated in the last sentence % !!!
of \fullcref{RrankRems}{Rrank=CoverGamma}.}

\end{exercises}

\section{Groups of higher real rank} \label{HigherRrankSect}

In some situations, there is a certain subset~$S$
of~$G$, such that the centralizer of each
element of~$S$ is well-behaved, and it would be helpful to know
that these centralizers generate~$G$. The results in this section illustrate that an assumption on the real rank of~$G$ may be exactly
what is needed. (However, we will often only prove the special case where $G = \SL(3,\real)$. A reader familiar with the theory of ``\term{real root}s'' should have no difficulty generalizing the arguments.)

\begin{prop} \label{Rrank2<>C(a)gens}
Let $A$ be a maximal $\real$-split torus in~$G$. Then we have
 $\Rrank G \ge 2$
 if and only if there exist nontrivial elements $a_1$ and~$a_2$ of~$A$,
such that $G = \langle \czer_G(a_1), \czer_G(a_2)\rangle$.
\end{prop}

\begin{proof}
($\Rightarrow$) Assume, for simplicity, that $G = \SL(3,\real)$.
(See \fullcref{CentGensG1xG2Ex}{C(a)} for another special case.) 
 Then we may assume $A$ is the group of diagonal matrices (after replacing it by a conjugate). Let
 $$ a_1 =
 \begin{bmatrix}
 2 & 0 & 0 \\
 0 & 2 & 0 \\
 0 & 0 & 1/4 
 \end{bmatrix}
 \mbox{ \qquad and \qquad }
 a_2 =
 \begin{bmatrix}
 1/4 & 0 & 0 \\
 0 & 2 & 0 \\
 0 & 0 & 2 
 \end{bmatrix}
 .$$
 Then
 $$ \czer_G(a_1) =
 \begin{bmatrix}
 * & * & 0 \\
 * & * & 0 \\
 0 & 0 & * 
 \end{bmatrix}
 \mbox{ \qquad and \qquad }
 \czer_G(a_2) =
 \begin{bmatrix}
 * & 0 & 0 \\
 0 & * & * \\
 0 & * & * 
 \end{bmatrix}
 .$$
 These generate~$G$.
 
 ($\Leftarrow$) Suppose $\Rrank G = 1$, so $\dim A = 1$. Then, since $A$ is almost Zariski closed (and contains $\langle a_1 \rangle$), we have 
 	$\czer_G(a_1) = \czer_G(A) = \czer_G(a_2)$,
so 
	$$\langle \czer_G(a_1) , \czer_G(a_2) \rangle = \czer_G(A) .$$
It is obvious that $\czer_G(A) \neq G$ (because the center of~$G$ is finite, and therefore cannot contain the infinite group~$A$).
\end{proof}

The following explicit description of $\czer_G(A)$ will be used in some of the proofs.

\begin{lem} \label{C(A)=AxC}
If $A$ is any maximal $\real$-split torus in~$G$, then $\czer_G(A) = A \times C$, where $C$~is compact.
\end{lem}

\begin{proof}[Proof]  \optional\ 
A subgroup of $\SL(\ell,\real)$ is said to be \defit{reductive} if it is isogenous to $M \times T$, where $M$~is semisimple and $T$~is a torus. It is known that the centralizer of any torus is reductive \csee{C(T)reductive}, so, if we assume, for simplicity, that $\czer_G(A)$ is connected, then we may write $\czer_G(A) = M \times A$, where $M$ is reductive \csee{C=MA}. The maximality of~$A$ implies that $M$ does not contain any $\real$-split tori, so $M$ is compact \csee{Rrank0Ex}.
\end{proof}

 \begin{prop}[\csee{Rrank2<>au=uaEx}] \label{Rrank2<>au=ua}
 $\Rrank G \ge 2$ if and only if there exist a nontrivial hyperbolic element~$a$ and
a nontrivial unipotent element~$u$, such that $au = ua$.
\end{prop}

For use in the proof of the \lcnamecref{Rrank2<>UnipGens} that follows it, we mention a very useful characterization of a somewhat different flavor:

\begin{lem}[\csee{Rrank1UniqMaxUnipEx}] \label{Rrank1UniqMaxUnip}
$\Rrank G \le 1$ if and only if every nontrivial unipotent subgroup of~$G$ is contained in a \textbf{unique} maximal unipotent subgroup.
\end{lem}

\begin{prop} \label{Rrank2<>UnipGens}
$\Rrank G \ge 2$ if and only if there exist nontrivial unipotent subgroups $U_1,\ldots,U_k$, such that 
	\begin{itemize}
	\item $\langle U_1,\ldots,U_k \rangle = G$,
	and
	\item $U_i$ centralizes~$U_{i+1}$ for each~$i$.
	\end{itemize}
\end{prop}

\begin{proof}
($\Rightarrow$) Assume, for simplicity, that $G = \SL(3,\real)$.
Then we take the sequence
 $$ 
 \begin{Smallbmatrix}
 1 & \upast & 0 \\
 0 & 1 & 0 \\
 0 & 0 & 1 
 \end{Smallbmatrix}
 , \ 
 \begin{Smallbmatrix}
 1 & 0 & \upast \\
 0 & 1 & 0 \\
 0 & 0 & 1 
 \end{Smallbmatrix}
 , \ 
 \begin{Smallbmatrix}
 1 & 0 & 0 \\
 0 & 1 & \upast \\
 0 & 0 & 1 
 \end{Smallbmatrix}
 , \ 
 \begin{Smallbmatrix}
 1 & 0 & 0 \\
 \upast & 1 & 0 \\
 0 & 0 & 1 
 \end{Smallbmatrix}
 , \ 
 \begin{Smallbmatrix}
 1 & 0 & 0 \\
 0 & 1 & 0 \\
 \upast & 0 & 1 
 \end{Smallbmatrix}
 , \ 
 \begin{Smallbmatrix}
 1 & 0 & 0 \\
 0 & 1 & 0 \\
 0 & \upast & 1 
 \end{Smallbmatrix}
 .$$

\medbreak

($\Leftarrow$) Since $U_i$ commutes with $U_{i+1}$, we know that $\langle U_i, U_{i+1} \rangle$ is unipotent, so, if $\Rrank G = 1$, then it is contained in a \emph{unique} maximal unipotent subgroup $\overline{U_i}$ of~$G$. Since $\overline{U_i}$ and $\overline{U_{i+1}}$ both contain~$U_{i+1}$, we conclude that $\overline{U_i} = \overline{U_{i+1}}$ for all~$i$. Hence, $\langle U_1,\ldots,U_k \rangle$ is contained in the unipotent group $\overline{U_1}$, and is therefore not all of~$G$.
 \end{proof}

\begin{rem}
See \cref{HighRankGenLi} for yet another result of the same type, which will be used in the proof of the Margulis Superrigidity Theorem in \cref{SuperPfSect}. A quite different characterization, based on the existence of subgroups of the form $\SL(2,\real) \ltimes \real^n$, appears in \cref{SL2RxRnInG}, and is used in proving Kazhdan's Property~$(T)$ in \cref{KazhdanTChap}.
\end{rem}

We know that $\SL(2,\real)$ is the smallest group of real rank one \csee{Rrank0<>SL2R}. However, the smallest group of real rank two is not unique:

\begin{prop} \label{Rrank2<>SL3orSO23}
Assume $G$ is simple. Then\/ $\Rrank G \ge 2$ if and only if $G$~contains a subgroup that is isogenous to either\/ $\SL(3,\real)$ or\/ $\SO(2,3)$.
\end{prop}

\begin{exercises}

\item \label{CentGensG1xG2Ex}
Prove the following results in the special case where $G = G_1 \times G_2$, and $\Rrank G_i \ge 1$ for each~$i$.
	\begin{enumerate}
	\item \label{CentGensG1xG2Ex-C(a)}
	\cref{Rrank2<>C(a)gens}($\Rightarrow$)
	\item \cref{Rrank2<>au=ua}($\Rightarrow$)
	\item \cref{Rrank1UniqMaxUnip}($\Leftarrow$)
	\item \cref{Rrank2<>UnipGens}($\Rightarrow$)
	\end{enumerate}

\item \label{C(T)reductive} \optional\ 
It is known that if $M$ is a subgroup that is almost Zariski closed, and $M^\transpose = M$, then $M$~is reductive \ccf{SelfAdj->SS}. Assuming this, show that if $T$ is a subgroup of the group of diagonal matrices, and $G^\transpose = G$, then $C_G(T)$ is reductive.

\item \label{C=MA} \optional\ 
Suppose $M$ is reductive, and $A$ is an $\real$-split torus in the center of~$M$. Show there exists a reductive subgroup $L$ of~$M^\circ$, such that $M^\circ = L \times A$.
\hint{Up to isogeny, write $M = M_0 \times T$, with $A \subseteq T$. Then it suffices to show $T = E \times A$ for some~$E$. You may assume, without proof, that, since $T$ is a connected, abelian Lie group, it is isomorphic to $\real^m \times \torus^n$ for some $m$ and~$n$.}

\item \label{Rrank2<>au=uaEx}
	\begin{enumerate}
	\item Prove \cref{Rrank2<>au=ua}($\Rightarrow$) under the additional assumption that $G = \SL(3,\real)$.
	\item \label{Rrank2<>au=uaEx-not1}
	Prove \cref{Rrank2<>au=ua}($\Leftarrow$).
%	\hint{\cref{ZarAandCent}.}
	\end{enumerate}

\item \label{Rrank1UniqMaxUnipEx}
Find a nontrivial unipotent subgroup of $\SL(3,\real)$ that is contained in two different maximal unipotent subgroups.

\item (\emph{Assumes the theory of real roots})
Prove the general case of the following results.
	\begin{enumerate}
	\item \cref{C(A)=AxC}
	\item \cref{Rrank2<>C(a)gens}($\Rightarrow$)
	\item \cref{Rrank2<>au=ua}($\Rightarrow$)
	\item \cref{Rrank1UniqMaxUnip}
	\item \cref{Rrank2<>UnipGens}($\Rightarrow$)
	\end{enumerate}

\item Show (without assuming $G$ is simple): $\Rrank G \ge 2$ if and only if $G$~contains a subgroup that is isogenous to either $\SL(3,\real)$ or $\SL(2,\real) \times \SL(2,\real)$.
\hint{\Cref{Rrank2<>SL3orSO23}. You may assume, without proof, that $\SO(2,2)$ is isogenous to $\SL(2,\real) \times \SL(2,\real)$.}

 \end{exercises}

\section{Groups of real rank one}

As a complement to \cref{HigherRrankSect}, here is an explicit list of the simple groups of real rank one.

\begin{thm} \label{rank1simple}
 If $G$ is simple, and\/ $\Rrank G = 1$, then $G$ is isogenous to
either
 \begin{itemize}
 \item $\SO(1,n)$ for some $n \ge 2$,
 \item $\SU(1,n)$ for some $n \ge 2$,
 \item $\Sp(1,n)$ for some $n \ge 2$, or
 \item $F_4^{-20}$ \textup(also known as~$F_{4,1}$\textup), a certain exceptional group.
 \end{itemize}
 \end{thm}

\begin{rem}
 The special linear groups $\SL(2,\real)$,
$\SL(2,\complex)$ and $\SL(2,\quaternion)$ have real rank
one, but they are already on the list under different names, because
 \begin{enumerate}
 \item $\SL(2,\real)$ is isogenous to $\SO(1,2)$ and
$\SU(1,1)$,
 \item $\SL(2,\complex)$ is isogenous to $\SO(1,3)$ and
$\Sp(1,1)$, and
 \item $\SL(2,\quaternion)$ is isogenous to $\SO(1,4)$.
 \end{enumerate}
 \end{rem}

% \begin{cor}
% $\Rrank(G) = 1$ if and only if $G$ is isogenous to a direct
%product $G_0 \times G_1$, where
% \begin{itemize}
% \item $G_0$ is compact and
% \item $G_1$ is of one of the groups listed in
%\cref{rank1simple}.
% \end{itemize}
% \end{cor}

\begin{rem}
 Each of the simple groups of real rank one has a very important geometric realization. Namely, $\SO(1,n)$, $\SU(1,n)$, $\Sp(1,n)$, and $F_{4,1}$ (respectively) are isogenous to the isometry groups of:
	 \begin{enumerate}

	\item (real) \term[hyperbolic!space!real]{hyperbolic $n$-space} $\hyperbolic^n$,
	 \item \defit[hyperbolic!space!complex]{complex hyperbolic $n$-space} $\complex\hyperbolic^n$,
	 \item \defit[hyperbolic!space!complex]{quaternionic hyperbolic $n$-space} $\quaternion\hyperbolic^n$,
	 and
	 \item the \index{hyperbolic!plane, octonionic}\defit[Cayley!plane]{Cayley plane}, which can be thought of as the hyperbolic plane over the (nonassociative) ring~$\octonion$ of Cayley \term{octonions}.
	 \end{enumerate}
%(It is not possible to define octonionic hyperbolic spaces of higher dimension, only the plane, because the ring $\octonion$ is not associative.) 
\end{rem}

\section{Minimal parabolic subgroups} \label{ParabSubgrpSect}

The group of upper-triangular matrices plays a very important role in the study of $\SL(n,\real)$. In this section, we introduce subgroups that play the same role in other semisimple Lie groups:

\begin{defn} \label{MinParabDefn}
Let $A$ be a maximal $\real$-split torus of~$G$, and let $a$ be a \defit[generic element]{generic} element of~$A$, by which we mean that $\czer_G(a) = \czer_G(A)$. Then the corresponding \defit[parabolic!subgroup!minimal]{minimal parabolic subgroup} of~$G$ is
	\begin{align} \label{P=horo}
	P = \bigset{ g \in G }{ \limsup_{n \to \infty} \| a^{-n} g a^n \| < \infty } 
	. \end{align}
This is a Zariski closed subgroup of~$G$.
 \end{defn}

\begin{thm} \label{MinParabConj}
All minimal parabolic subgroups of~$G$ are conjugate.
\end{thm}

\begin{egs} \ \label{ParabEgs}
 \noprelistbreak
 \begin{enumerate}

\item \label{ParabEgs-SLn}
The group of upper triangular matrices is a minimal parabolic subgroup of $\SL(n,\real)$. To see this, let $A$ be the group of diagonal matrices, and choose $a \in A$ with $a_{1,1} > a_{2,2} > \cdots > a_{n,n} > 0$ \csee{MinParabSLnEx}. 

 \item \label{ParabEgs-SO1n}
 It is easier to describe a minimal parabolic subgroup of $\SO(1,n)$ if
we replace $\Id_{m,n}$ with a different symmetric matrix of
the same signature: let 
 $G = \SO(A;\real)$, for
 $$ A =
 \begin{pmatrix}
 0 & 0 & 1 \\
 0 & \Id_{(n-1) \times (n-1)} & 0 \\
 1 & 0 & 0
 \end{pmatrix}
 .$$
 Then $G$ is conjugate to $\SO(1,n)$ \csee{SOA=SOmn},
 the ($1$-dimensional) group of diagonal matrices in~$G$ form a maximal $\real$-split torus, and a minimal parabolic subgroup in~$G$ is
 $$\left\{\begin{pmatrix}
 t&*&* \\
 0& \SO(n-1) & * \\
 0&0&1/t
 \end{pmatrix}
 \right\}
 $$
 \csee{MinParabSO1nEx}.
 \end{enumerate}
 \end{egs}

The following result explains that a minimal parabolic subgroup of a
classical group is simply the stabilizer of a (certain kind
of) flag. 
%(For the general case, including exceptional groups, the theory of real roots gives a good understanding of the parabolic subgroups.)
Recall that
%, if $\langle \cdot \mid \cdot \rangle$ is a bilinear (or Hermitian) form on a vector space~$V$. A 
a subspace $W$
of a vector space~$V$, equipped with a bilinear (or Hermitian) form $\langle \cdot \mid \cdot \rangle$, is said to be \defit[totally!isotropic]{totally isotropic} if $\langle W \mid W \rangle = 0$.

\begin{thm}[\csee{parab=Stab(flag)Ex}] \ \label{parab=Stab(flag)}
\noprelistbreak
 \begin{enumerate}
 \item \label{parab=Stab(flag)-SLn} 
 A subgroup~$P$ of\/ $\SL(n,\real)$ is a minimal parabolic if and
only if there is a chain $V_0 \subsetneq V_1 \subsetneq \cdots
\subsetneq V_n$ of subspaces of\/ $\real^n$ \textup(with $\dim V_i = i$\textup), such that
 $$ P = \{\, g \in \SL(n,\real) \mid
 \forall i, \ g V_i = V_i \,\} .$$
Similarly for\/ $\SL(n,\complex)$ and\/ $\SL(n,\quaternion)$,
taking chains of subspaces in\/ $\complex^n$
or\/~$\quaternion^n$, respectively.

\item \label{parab=Stab(flag)-SOmn} 
A subgroup~$P$ of\/ $\SO(m,n)$ is a minimal parabolic if and only
if there is a chain $V_0 \subsetneq V_1 \subsetneq \cdots
\subsetneq V_r$ of \textbf{totally isotropic} subspaces of $\real^{m+n}$ \textup(with $\dim V_i = i$ and $r = \min\{m,n\}$\textup), such that
 $$ P = \{\, g \in \SO(m,n) \mid
 \forall i, \ g V_i = V_i \,\} .$$
Similarly for\/ $\SO(n,\complex)$, $\SO(n,\quaternion)$,
$\Sp(2m,\real)$, $\Sp(2m,\complex)$, $\SU(m,n)$ and\/
$\Sp(m,n)$.
 \end{enumerate}
 \end{thm}

%Note that, in each of these cases, we may write $P = MAN$, where
%	\begin{itemize}
%	\item $M$ is compact,
%	\item $A$ is a maximal $\real$-split torus of~$G$,
%	and
%	\item $N$ is a maximal unipotent subgroup of~$G$.
%	\end{itemize}
%
%The group of upper triangular matrices in $\SL(n,\real)$ is of the form $AN$, where
%	\begin{itemize}
%	\item $A$ is a maximal $\real$-split torus,
%	and
%	\item $N$ is a maximal unipotent subgroup that is normalized by~$A$.
%	\end{itemize}
%This description can be extended to the minimal parabolic subgroups of any semisimple group, except that we need to replace $A$ with the slightly larger group $\czer_G(A)$:

%\begin{thm}
%$P$ is a minimal parabolic subgroup of~$G$ if and only there exist%
%	\begin{itemize}
%	\item a maximal $\real$-split torus $A$,
%	and
%	\item a maximal unipotent subgroup~$N$ that is normalized by~$A$,
%	\end{itemize}
%such that $P = \czer_G(A) \, N$.
%\end{thm}

Note that any upper triangular matrix in $\SL(n,\real)$ can be written uniquely in the form $mau$, where
	\begin{itemize}
	\item $a$~belongs to the $\real$-split torus~$A$ of diagonal matrices whose nonzero entries are positive,
	\item $m$ is in the finite group~$M$ consisting of diagonal matrices whose nonzero entries are $\pm1$,
	and
	\item $u$ belongs to the unipotent group~$N$ of upper triangular matrices with $1$'s on the diagonal.
	\end{itemize}
The elements of every minimal parabolic subgroup have a decomposition of this form, except that the subgroup~$M$ may need to be compact, instead of only finite:

\begin{thm}[(\thmindex{Langlands decomposition}Langlands decomposition)] \label{LanglandsDecomp}
 If $P$ is a minimal parabolic subgroup of~$G$, then we may write it in the form
 $P = \czer_G(A)\, N = MAN$, where%
	 \begin{itemize}
	 \item $A$~is a maximal $\real$-split torus, 
	 \item $M$ is a compact subgroup of $\czer_G(A)$,
	 and 
	 \item $N$ is the unique maximal unipotent subgroup of~$P$.
	 \end{itemize}
 Furthermore, $N$ is a maximal unipotent subgroup of~$G$, and, for some generic $a \in A$, we have
 	\begin{align} \label{LanglandsDecomp-N}
	N = \bigset{ u \in G }{ \lim_{n \to \infty} a^{-n} u a^n = e } 
	. \end{align}
 \end{thm}
 
 Before discussing the proof (which is not so important for our purposes), let us consider a few examples:

\begin{eg} \label{MinParabEg} \ 
\noprelistbreak
	\begin{enumerate}
	\item If $G = \SL(n,\complex)$, then, for the Langlands decomposition of the group~$P$ of upper-triangular matrices, we may let:
		\begin{itemize}
		\item $A$ be the group of diagonal matrices in~$G$ whose nonzero entries are positive real numbers (just as for $\SL(n,\real)$),
		\item $M$ be the group of diagonal matrices in~$G$ whose nonzero entries have absolute value~$1$,
		and
		\item $N$ be the group of upper triangular matrices with $1$'s on the diagonal.
		\end{itemize}
	The same description applies to $G = \SL(n,\quaternion)$ (and, actually, also to $\SL(n,\real)$).
	
	\item \label{MinParabEg-SOA}
	Assume $m \le m$, and let $G = \SO(A; \real)$, where 
		$$ A = \begin{bmatrix} 
		0 & 0 & J_m \\
		0 & \Id_{(n-m) \times (n-m)} & 0 \\
		J_m & 0 & 0
		\end{bmatrix}
		\quad  \text{and} \quad
		J_m = \left[\begin{smallmatrix}
		&&&&1 \\
		&\vbox to 0pt{\vss\hbox to 0pt{\Large\hskip-5pt0\hss}\vss}&&1 \\
		&&\nedots \\
		& 1 & \vbox to 0pt{\vss\hbox to 0pt{\Large\hskip10pt 0\hss}\vss\vskip3pt}\\
		 1 
		\end{smallmatrix}\right]
		$$
	(and the size of the matrix $J_m$ is $m \times m$). Then $G$ is conjugate to $\SO(m,n)$ \csee{SOA=SOmn}, and a minimal parabolic~$P$ of~$G$ is:
	$$  \bigset{\begin{pmatrix}
	 b&*& * \\
	 0& k & * \\
	 0&0&b^\dagger
	 \end{pmatrix}
	 }{ 
	 \begin{matrix}
	 \text{$b \in \GL(m,\real)$ is upper triangular} , \\
	 k \in \SO(n-m) \\
	 \end{matrix}}
	 $$
	 where $x^\dagger = J_m (x^{-1})^\transpose J_m$,
	 %denotes the reflection of~$x$ across the aniti-diagonal 
	 so, for example, 
	 	$$\diag(b_1,\ldots,b_m)^\dagger = \diag(1/b_m,\ldots,1/b_1) .$$
	Hence, we may let
		\begin{itemize}
		\item $A = \{\, \diag(a_1,\ldots,a_m, 0, \ldots,0, 1/a_m, \ldots,1/a_1) \mid a_i > 0 \,\}$,
		\item $M \iso \SO(n-m) \times \{\pm1\}^m$,
		and
		\item $N$ be the group of upper triangular matrices with $1$'s on the diagonal that are in~$G$.
		\end{itemize}
	\end{enumerate}
\end{eg}
 
 \begin{proof}[Proof of \cref{LanglandsDecomp} \optional]
 Choose a generic element~$a$ of~$A$ satisfying \pref{P=horo}, and define $N$ as in \pref{LanglandsDecomp-N}. Then, since $a$~is diagonalizable over~$\real$, it is not difficult to see that $P = \czer_G(a)\, N$ \csee{P=CN}.
Since $a$~is a generic element of~$A$, this means $P = \czer_G(A)\, N$.

It is easy to verify that $N$ is normal in~$P$ \csee{NnormP}; then, since $P/N \iso \czer_G(A) = A \times (\text{compact})$ \csee{C(A)=AxC}, and therefore has no nontrivial unipotent elements, it is clear that $N$ contains every unipotent element of~$P$. 
Conversely, the definition of~$N$ implies that it is unipotent \csee{NUnip}.
Therefore, $N$~is the unique maximal unipotent subgroup of~$P$.

Suppose $U$ is a unipotent subgroup of~$G$ that properly contains~$N$.
Since unipotent subgroups are nilpotent \csee{UnipNilp}, then $\nzer_U(N)$ properly contains~$N$ \csee{NilpNorm}.
However, it can be shown that $\nzer_G(N) = P$ \csee{P=normalizer}, so this implies $\nzer_U(N)$ is a unipotent subgroup of~$P$ that properly contains~$N$, 
which contradicts the conclusion of the preceding paragraph. % !!!
 \end{proof}

The subgroups $A$ and~$N$ that appear in the Langlands decomposition of~$P$ are two components of the Iwasawa decomposition of~$G$:

\begin{thm}[(\thmindex{Iwasawa decomposition}Iwasawa decomposition)] \label{IwasawaDecomp}
 Let 
 \noprelistbreak
 \begin{itemize}
 \item $K$ be a maximal compact subgroup of~$G$,
 \item $A$ be a maximal $\real$-split torus, 
 and
 \item $N$ be a maximal unipotent subgroup that is normalized by~$A$.
 \end{itemize}
Then $G = K A N$.

 In fact, every $g \in G$ has a \bemph{unique} representation of the form $g = k a u$ with $k \in K$, $a \in A$, and $u \in N$. 
 \end{thm}

\begin{rem} \label{IwasawaDiffeo}
 A stronger statement is true: if we define a function $\varphi \colon K \times A \times N \to G$ by\/
 $\varphi(k,a,u) = kau$, then $\varphi$ is a (real analytic) diffeomorphism. Indeed, \cref{IwasawaDecomp} tells us that $\varphi$ is a bijection, and it is obviously real analytic. It is not so obvious that the inverse of~$\varphi$ is also real analytic, but this is proved in \cref{IwasawContinuousSLnR} when $G = \SL(n,\real)$, and the general case can be obtained by choosing an embedding of~$G$ in $\SL(n,\real)$ for which the subgroups $K$, $A$, and~$N$ of~$G$ are equal to the intersection of~$G$ with the corresponding subgroups of $\SL(n,\real)$.
 \end{rem}

The Iwasawa decomposition implies $KP = G$ (since $AN \subseteq P$), so it has the following important consequence:

\begin{cor} \label{G/Pcpct}
If $P$ is any minimal parabolic subgroup of~$G$, then $G/P$ is compact.
\end{cor}

\begin{rem} \label{ParabRem}
A subgroup of~$G$ is called \defit[parabolic!subgroup]{parabolic} if it contains a minimal parabolic subgroup. 
	\begin{enumerate}

	\item \label{ParabRem-cocpct}
	\Cref{G/Pcpct} implies that if $Q$ is any parabolic subgroup, then $G/Q$ is compact. The converse
does not hold. (For example, if $P = MAN$ is a minimal
parabolic, then $G/(AN)$ is compact, but $AN$~is not
parabolic unless $M$ is trivial.) However, passing to the
``\term{complexification}'' does yield the converse: $Q$ is parabolic if
and only if $G_{\complex}/Q_{\complex}$ is compact.
Furthermore, $Q$ is parabolic if and only if $Q_{\complex}$
contains a maximal solvable subgroup (``\defit[Borel!subgroup]{Borel subgroup}'') of~$G_{\complex}$.

	\item All parabolic subgroups can be described fairly completely (there are only finitely many that contain any given minimal parabolic), but we do not need the more general theory.
	\end{enumerate}
 \end{rem}

\begin{exercises}

\item \label{MinParabSLnEx}
Let $a$ be a diagonal matrix as described in \fullcref{ParabEgs}{SLn}, and show that the corresponding minimal parabolic subgroup is precisely the group of upper triangular matrices.

\item \label{MinParabSO1nEx}
Show that the subgroup at the end of \fullcref{ParabEgs}{SO1n} is indeed a minimal parabolic subgroup of $\SO(A; \real)$.

\item \label{parab=Stab(flag)Ex}
Show the minimal parabolic subgroups of each of the following groups are as described in \cref{parab=Stab(flag)}:
	\begin{enumerate}
	\item $\SL(n,\real)$.
	\item $\SO(m,n)$.
	\end{enumerate}
\hint{It suffices to find one minimal parabolic subgroup in order to understand all of them \csee{MinParabConj}.}

\item \label{SOA=SOmn}
For $A$ as in \fullcref{MinParabEg}{SOA}, show that $\SO(A; \real)$ is conjugate to $\SO(m,n)$.
\hint{Let $\alpha = 1/\sqrt{2}$, and define $v_i$ to be: $\alpha(e_i + e_{n+1-i})$ for $i \le m$, $e_i$~for $m < i \le n$, and $\alpha(e_i - e_{n+1-i})$ for $i > n$. Then $v_i^\transpose A v_i$ is~$1$ for $i \le n$, and is $-1$~for $i > n$.}

\item \label{P=CN} \optional\ 
For $P$, $a$, and~$N$ as in the proof of \cref{LanglandsDecomp}, show $P = \czer_G(a)\, N$.
\hint{Given $g \in P$, show that $a^{-n} g a^n$ converges to some element~$c$ of $C_G(a)$. Also show $c^{-1} g \in N$. You may assume $a$~is diagonal, with $a_{11} \ge a_{22} \ge \cdots \ge a_{\ell\ell}$ (\emph{why?}).}

\item \label{NnormP}
For $P$, $a$, and~$N$ as in the proof of \cref{LanglandsDecomp}, show $N$ is a normal subgroup of~$P$.

\item \label{NUnip}
 Show that a subgroup~$N$ satisfying \pref{LanglandsDecomp-N} must be unipotent.
 \hint{$u$~has the same characteristic polynomial as~$a^{-n} u a^n$.} 

\item \label{P=normalizer}
For $P$ and~$N$ as in \fullcref{parab=Stab(flag)}{SOmn}, show $P = \nzer_G(N)$.
\hint{$P$~is the stabilizer of a certain flag, and the subgroup~$N$ also uniquely determines this same flag.}

\item \label{UnipNilp}
Show that every unipotent subgroup of $\SL(\ell,\real)$ is nilpotent.
(Recall that a group~$N$ is \defit[nilpotent group]{nilpotent} if there is a series 
	$$\{e\} = N_0 \normal \cdots \normal N_r = N$$
of subgroups of~$N$, such that $[N, N_k] \subseteq N_{k-1}$ for each~$k$.)
\hint{Engel's Theorem \pref{EngelUnip}.}

\item \label{NilpNorm}
Show that if $N$ is a proper subgroup of a nilpotent group~$U$, then $\nzer_U(N) \not\subseteq N$.
\hint{If $[N, U_k] \subseteq N$, then $U_k$ normalizes~~$N$.}

\item  \label{G=KxRn}
Assume $K$ is a maximal compact subgroup of~$G$. Show:
	\noprelistbreak
	\begin{enumerate}
	\item $G$ is diffeomorphic to the cartesian product $K \times \real^n$, for
some~$n$,
	\item \label{G=KxRn-G/K}
	$G/K$ is diffeomorphic to $\real^n$, for some~$n$,
	\item $G$ is connected if and only if $K$ is connected,
	and 
	\item $G$ is simply connected if and only if $K$ is simply connected.
	\end{enumerate}
\hint{\cref{IwasawaDiffeo}.}

%\item \label{G/Ksc}
%Show that if $K$ is any maximal compact subgroup of~$G$, then $G/K$ is connected and simply connected.
%\hint{The Iwasawa decomposition implies that $G/K$ is homeomorphic to $AN$. The subgroups $A$ and~$N$ are each homeomorphic to some Euclidean space. So, in fact, $G/K$ is homeomorphic to some~$\real^k$, and is therefore contractible, but we do not need this fact.}

\end{exercises}

\begin{notes}

The comprehensive treatise of Borel and Tits \cite{BorelTits-GrpRed} is the standard reference on rank, parabolic subgroups, and other fundamental properties of reductive groups over any field.
See \cite[\S7.7, pp.~474--487]{Knapp-BeyondIntro} for a discussion of parabolic subgroups of Lie groups (which is the special case in which the field is~$\real$). 

\fullCref{RrankRems}{Rrank=CoverGamma} is due to Prasad-Raghunathan \cite[Thms.~2.8 and 3.9]{PrasadRaghunathan-Cartan+Latts}.

Proofs of the Iwasawa decomposition for both $\SL(n,\real)$ \pref{IwasawaDecompSLnR} and the general case \pref{IwasawaDecomp} can be found in \cite[Prop.~3.12, p.~129, and Thm.~3.9, p.~131]{PlatonovRapinchukBook}. (Iwasawa's original proof is in \cite[\S3]{Iwasawa-SomeTypes}.) The decomposition also appears in many textbooks on Lie groups. In particular, \cref{IwasawaDiffeo} is proved in \cite[Thm.~6.5.1, pp.~270--271]{HelgasonBook}. 
% Since $A$ and~$N$ are connected, the decomposition implies that every maximal compact subgroup of~$G$ intersects all of the components of~$G$.

Regarding \fullcref{ParabRem}{cocpct}, the obvious cocompact subgroups of~$G$ are parabolic subgroups and (cocompact) lattices. See \cite{Witte-cocpct} for a short proof that every cocompact subgroup is a combination of these two types. (A similar result had  been proved previously in \cite[(5.1a)]{GotoWang-NondiscreteUniform}.)

\end{notes}

 %!TEX root = IntroArithGrps.tex

\mychapter{\texorpdfstring{$\rational$}{Q}-Rank} \label{QrankChap}

\prereqs{Real rank and minimal parabolic subgroups (\cref{RrankChap}), and groups defined over~$\rational$ (\cref{DefdQDefn}).}

Algebraically, the definition of real rank extends in a straightforward way to a notion of rank over any field: if $G$ is defined over~$F$, then we can talk about $\rank_F G$. In the study of arithmetic groups, we assume $G$ is defined over~$\rational$, and the corresponding $\rational$-rank is an important invariant of the associated arithmetic group $\Gamma = G_\integer$.

\begingroup \smaller \renewcommand{\baselinestretch}{0.95}
\begin{disclaimer}
The reading of this \lcnamecref{QrankChap} may be postponed without severe consequences (and can even be skipped entirely), because the material here will not arise elsewhere in this book (except marginally) other than in \cref{ReductionChap}, where a coarse fundamental domain for~$\Gamma$ will be constructed. Furthermore, unlike the other chapters in this part of the book, the topic is of importance only for arithmetic groups and closely related subjects, not a broad range of areas of mathematics.
\end{disclaimer}
\par \endgroup

\section{\texorpdfstring{$\rational$}{Q}-rank}

\begin{defn} %[(cf.\ definition of $\real$-split torus in \pref{RsplitDefn})]
Assume $G$ is defined over~$\rational$.
 A closed, connected subgroup~$T$ of~$G$ is a 
 \index{torus!$\rational$-split}%
 \defit[Q-@$\rational$-!split]{$\rational$-split torus} if
 \noprelistbreak
 \begin{itemize}
 \item $T$ is defined over~$\rational$, and
 \item $T$ is diagonalizable over~$\rational$. (That is,
there exists $g \in \GL(\ell,\rational)$, such that $g
T g^{-1}$ consists entirely of diagonal matrices.)
 \end{itemize}
 \end{defn}

\begin{eg} \ 
\noprelistbreak
 \begin{enumerate}
 \item $\SO(1,1)^\circ$ is a $\rational$-split torus, because $g \SO(1,1) g^{-1}$ consists of diagonal matrices if $g =  \begin{Smallbmatrix}
 1 & 1 \\
 1 & -1 
 \end{Smallbmatrix}
 $.
 
\item Although it is obvious that every $\rational$-split torus is an $\real$-split torus, the converse is not true (even if the torus is defined over~$\rational$). For example, let $T = \SO(x^2 - 2y^2; \real)^\circ$. Then $T$ is defined over~$\rational$, and it is $\real$-split (because it is conjugate to $\SO(1,1)^\circ$). However, it is \textbf{not} $\rational$-split. To see this, note that 
 $ \begin{Smallbmatrix}
 3 & 4 \\
 2 & 3 
 \end{Smallbmatrix}
 \in T_{\rational}
 $,
 but the eigenvalues of this matrix are irrational (namely,
$3 \pm 2\sqrt{2}$), so this rational matrix is not
diagonalizable over~$\rational$. 

 \end{enumerate}
\end{eg}

The following key fact implies that the maximal $\rational$-split tori of~$G$ all have the same dimension (which is called the ``$\rational$-rank''):

\begin{thm}
Assume $G$ is defined over~$\rational$.
If $S_1$ and~$S_2$ are maximal\/ $\rational$-split tori in~$G$, then $S_1 = g S_2 g^{-1}$ for some $g \in G_\rational$.
\end{thm}

\begin{defn}[(for arithmetic lattices)] \label{Qrank-arithmetic}
 Assume
 \noprelistbreak
 	\begin{itemize}
	\item $G$ is defined over~$\rational$, 
	and 
	\item $\Gamma$ is commensurable to $G_{\integer}$. 
	\end{itemize}
Then \defit[Q-@$\rational$-!rank]{$\Qrank \Gamma $}%
	\index{rank!Q-@$\rational$-|indsee{$\rational$-rank}}%
	\nindex{$\Qrank \Gamma $ = $\rational$-rank of~$\Gamma$}
 is the dimension of any maximal $\rational$-split torus in~$G$. 
%(Because all
%maximal $\rational$-split tori are conjugate under
%$G_{\rational}$, this definition is independent of the
%particular torus chosen.)

(More generally, if $\phi \colon G/K \stackrel{\iso}{\to} G'/K'$, where $K$
and~$K'$ are compact, and $\phi(\overline{\Gamma})$ is
commensurable to $\overline{G'_{\integer}}$
\csee{ArithDefn}, then $\Qrank \Gamma$ is the
dimension of any maximal $\rational$-split torus in~$G'$.)
 \end{defn}
 
\begin{egs} \ \label{QrankEg} 
 \noprelistbreak
	 \begin{enumerate}
	 
	 \item \label{QrankEg-SL}
	 $\Qrank \bigl( \SL(n,\integer) \bigr) = n-1$. (Let $S$ be the identity component of the group of all diagonal matrices in $\SL(n,\real)$.)
	
	 \item \label{QrankEg-SOQ}
	 Let $G = \SO(Q; \real)$, where $Q(x_1,\ldots,x_\ell)$
	is some quadratic form on~$\real^\ell$, such that $Q$~is
	defined over~$\rational$. (That is, all of the coefficients
	of~$Q$ are rational.) Then $G$ is defined over~$\rational$,
	and the discussion of \cref{RrankEgs}, with $\rational$
	in place of~$\real$, shows that $\Qrank G_{\integer} $ is
	the maximum dimension of a totally isotropic
	$\rational$-subspace of~$\rational^\ell$. 
	
		\begin{enumerate} \itemsep=\smallskipamount

		\item For example, 
		 \ $\Qrank  \SO(m,n)_{\integer} = \min\{m,n\} $. \ 
		 Similarly,
		 $$ \Qrank  \SU(m,n)_{\integer} 
		 = \Qrank  \Sp(m,n)_{\integer}  = \min\{m,n\} .$$
		 So $\Qrank G_\integer = \Rrank G$ for these groups.
		 
		 \item \label{QrankEg-SOQ-anis}
		 Let $G = \SO(x_1^2 + x_2^2 + x_3^2 - 7 x_4^2; \real)$.
		Then, because the congruence $a^2 + b^2 + c^2 + d^2 \equiv 0
		\pmod{8}$ implies that all the variables are even, it is not difficult
		to see that this quadratic form has no nonzero isotropic vectors in~$\rational^4$ \csee{NoSolQ}. This means $\Qrank G_{\integer}  = 0$.

		Note, however, that $G$ is isomorphic to $\SO(3,1)$, so its real rank is~$1$. Therefore, $\Qrank G_\integer \neq \Rrank G$.
%			 Since there are quadratic forms that have nontrivial isotropic subspaces over~$\real$, but not over~$\rational$, we may have $\Rrank \SO(Q; \real) \neq \Qrank \SO(Q; \integer)$ 
%
%		 $$ \begin{bmatrix}
%		 \sqrt{7} & 0 & 0 \\
%		 0 & 1 & 0 & 0 \\
%		 0 & 0 & 1 & 0 \\
%		 0 & 0 & 0 & 1
%		 \end{bmatrix}
%		 .$$
%		 We see from \cref{Qrank0<>} below that $G_{\integer}
%		\backslash G$ is compact; therefore, $G_{\integer} \backslash
%		\hyperbolic^3$ is a \emph{compact} hyperbolic 3-manifold,
%		whereas $\SO(1,3)_{\integer} \backslash \hyperbolic^3$ is
%		not compact. Therefore, the Mostow Rigidity Theorem
%		\pref{MostowIrred} implies that $G_{\integer}$ is not
%		isomorphic to $\SO(1,3)_{\integer}$; this is an illustration
%		of the fact that a different embedding of~$G$ in
%		$\SL(\ell,\real)$ can yield very different integer points.
%		(Technical remark: Actually, $G_{\integer} \backslash
%		\hyperbolic^3$ is an \emph{orbifold}. To obtain a (compact)
%		\emph{manifold,} replace $G_{\integer}$ with a torsion-free
%		subgroup of finite index \see{torsionfree}.)
		\end{enumerate}
	
	\item \label{QrankEg-0}
	$\Qrank \Gamma = 0$ if and only if $G/\Gamma$ is compact \csee{Qrank0Ex}.

	\item \label{QrankEg-SUD}
	$\Qrank \SU(B, \tau; \ints_D)$ is the dimension (over~$D$) of a maximal totally isotropic subspace of~$D^n$, if $B$ is a $\tau$-Hermitian form on~$D^n$, and $D$~is a division algebra over~$F$.
	
	 \end{enumerate}
 \end{egs}

\begin{warn} \label{SL2ZinGamma}
 In analogy with \cref{Rrank0<>SL2R} and
\fullcref{GZCocpctIff}{SL2Q}, one might suppose that
$\Qrank \Gamma \neq 0$ if and only if $\Gamma$
contains a subgroup that is isomorphic to
$\SL(2,\integer)$ (modulo finite groups). However, this
is \textbf{false}: \emph{every} lattice in~$G$ contains a
subgroup that is abstractly commensurable to
$\SL(2,\integer)$ (unless $G$ is compact). Namely, the Tits Alternative tells us that $\Gamma$ contains a nonabelian free subgroup \csee{FreeInGamma},
and it is well known that $\SL(2,\integer)$ has a
finite-index subgroup that is free \csee{SanovIsFree}.
 \end{warn}

\begin{rems} \label{QrankRems} \ 
\noprelistbreak
	\begin{enumerate}

	\item The definition of $\Qrank \Gamma$ is somewhat indirect,
because the $\rational$-split tori of~$G$ are not subgroups
of~$\Gamma$. Therefore, it would be more correct to say that
we have defined $\Qrank G_{\rational} $. 

	\item Although different embeddings
of~$G$ in $\SL(\ell,\real)$ can yield maximal $\rational$-split tori of different dimensions, the theory of algebraic groups shows
that the $\rational$-rank is the same for all of the embeddings in
which $\Gamma$ is commensurable to~$G_{\integer}$
\csee{Qrank(isog)=}; therefore, $\Qrank \Gamma$ is well defined
as a function of~$\Gamma$.

	\item We have $0 \le \Qrank \Gamma \le \Rrank G$, since every $\rational$-split torus is $\real$-split. It can be shown that:
		\begin{enumerate}
	
		\item The extreme values are always realized: there exist lattices $\Gamma_0$ and~$\Gamma_r$ in~$G$, such that $\Qrank \Gamma_0 = 0$ and $\Qrank \Gamma_r = \Rrank G$ \csee{GHasCpctLatt,Qrank=Rrank}.
		
		\item \label{QrankRems-gap} 
		In some cases, there are intermediate values 
			%between $0$ and $\Rrank G$ 
		that are not realized. For example, the $\rational$-rank of every lattice in $\SO(2,5)$ is either $0$ or~$2$ \csee{QrankGap}.
		\end{enumerate}
	
	\item Suppose $\Gamma$ is defined by restriction of scalars \pref{ResScal->Latt}, so $\Gamma$ is commensurable to $G'_{\ints}$, where $G'$~is defined over a finite extension~$F$ of~$\rational$, and $\ints$ is the ring of integers of~$F$. Then $\Qrank \Gamma$ is equal to the ``\index{rank!over a field F@over a field~$F$}$F$-rank'' of~$G'$, or, in other words, the maximal dimension (over~$F_\infty$) of a subgroup of~$G'$ that is diagonalizable over~$F$. For example, the $\rational$-rank of  $\SO(B; \ints)$ is the dimension of a maximal totally isotropic $F$-subspace of~$F^n$.

	\end{enumerate}
 \end{rems}

%Similarly, if $G$ is defined over an algebraic number
%field~$F$, then one can define
%$\operatorname{\mbox{\upshape$F$-rank}}(G)$. The following
%result shows that this can be used to calculate the $\rational$-rank
%of a lattice obtained by Restriction of Scalars.
%
%\begin{lem}
% Suppose 
% \begin{itemize}
% \item $F$ is an algebraic number field, 
% \item $\ints$ is the ring of integers of~$F$,
% \item $G$ is
%defined over~$F$ \textup(as an algebraic group
%over~$F_\infty$\textup), and
% \item $\Delta \colon G_F \to \prod_{\sigma \in S^\infty}
%G^\sigma$ is defined by $\Delta(g) = \bigl( \sigma(g)
%\bigr)_{\sigma \in S^\infty}$, as in
%\cref{ResScal->Latt}.
% \end{itemize}
% Then $\Qrank \bigl( \Delta(G_{\ints}) \bigr) =
%\operatorname{\mbox{\upshape$F$-rank}}(G)$.
% \end{lem}
%
%\begin{proof}
% If $T$ is a torus in~$G$, and $T$~is defined over~$F$, then
%$\prod_{\sigma \in S^\infty} T^\sigma$ is a $\rational$-torus
%in $\prod_{\sigma \in S^\infty} G^\sigma$. Conversely, it is
%not difficult to see that any $\rational$-torus of
%$\prod_{\sigma \in S^\infty} G^\sigma$ is contained in a
%torus of the form $\prod_{\sigma \in S^\infty} T^\sigma$.
%Therefore, the desired conclusion follows from the fact, which
%will be proved in \S\ref{DirichletUnitSect}, that
% $$ \Qrank \left( \prod_{\sigma \in S^\infty} T^\sigma
%\right)
% = \operatorname{\mbox{\upshape$F$-rank}}(T) $$
% \see{Qrank(ROS)=Frank(T)}.
% \end{proof}

\Cref{Qrank-arithmetic} applies only to arithmetic lattices, but the Margulis Arithmeticity Theorem \pref{MargulisArith} allows the definition to be extended to all lattices:

\begin{defn}[\csee{QrankWellDefd}] \label{QrankDefn}
 Up to isogeny, and modulo the maximal compact factor of~$G$,
we may write $G = G_1 \times \cdots \times G_s$, so that
$\Gamma_i = \Gamma \cap G_i$ is an irreducible lattice
in~$G_i$ for $i = 1,\ldots,r$ \csee{prodirredlatt}. We let\term[Q-@$\rational$-!rank]{}
 $$ \Qrank(\Gamma) = \Qrank(\Gamma_1) +
\cdots + \Qrank(\Gamma_s) ,$$
 where:
	 \begin{enumerate}
	 \item If $G/\Gamma_i$ is compact, then $\Qrank \Gamma_i  = 0$.
	 \item If $G/\Gamma_i$ is not compact, and $\Rrank G = 1$, then $\Qrank \Gamma_i = 1$.
	 \item If $G/\Gamma_i$ is not compact, and $\Rrank G \ge 2$, then the Margulis Arithmeticity Theorem \pref{MargulisArith} implies that $\Gamma_i$ is
	arithmetic, so \cref{Qrank-arithmetic} applies.
	 \end{enumerate}
 \end{defn}

\begin{exercises}

\item Show that if $T$ is a $\rational$-split torus, then
$T_{\integer}$ is finite.

\item Give an example of a torus~$T$ (that is defined over~$\rational$), such that
$T_{\integer}$ is infinite.

\item Verify the claim in \fullcref{QrankEg}{SOQ} that $\Qrank \SO(Q ; \integer)$ is the dimension of a maximal totally isotropic subspace of~$\rational^\ell$.

\item \label{NoSolQ}
Verify the claim in \fullcref{QrankEg}{SOQ-anis} that $(0,0,0,0)$ is the only solution in~$\rational^4$ of the equation $x_1^2 + x_2^2 + x_3^2 - 7x_4^2 = 0$.

\item \label{Qrank0Ex}
Prove \fullcref{QrankEg}{0}.
\hint{($\Rightarrow$)~See \cref{GZCocpctIff}.
($\Leftarrow$)~If $a$ is diagonalizable over~$\rational$, then there exists $v \in \integer^\ell$, such that $a^n v \to 0$ as $n \to +\infty$, so the Mahler Compactness Criterion \pref{MahlerCpct} implies $G / G_{\integer}$ is not compact.}

\item \label{QrankWellDefd}
Show that \cref{QrankDefn} is consistent with \cref{Qrank-arithmetic}. More precisely, assume $\Gamma$ is arithmetic, and prove:
	 \begin{enumerate}
	 \item $G/\Gamma$ is compact if and only if $\Qrank \Gamma = 0$. 
	 \item If $G/\Gamma$ is not compact, and $\Rrank G = 1$, then $\Qrank \Gamma = 1$. 
%	$\Qrank(\Gamma) \le \Rrank(G) = 1$ \see{Qrank<Rrank}; so we
%	must have $\Qrank(\Gamma) = 1$.
	 \item If $\Gamma = \Gamma_1 \times \Gamma_2$
	is reducible, then $\Qrank \Gamma  = \Qrank \Gamma_1  + \Qrank \Gamma_2$.
%	 \item If $C$ is a compact, normal subgroup of~$G$, and
%	$\overline{\Gamma}$ is the image of~$\Gamma$ in
%	$\overline{G} = G/C$, then $\Qrank(\Gamma) =
%	\Qrank(\overline{\Gamma})$.
	 \end{enumerate}

\item \label{Qrank=Rrank}
Suppose $G$ is classical. Show that, for the natural embeddings
described in \cref{classical-fulllinear,classical-orthogonal}, we have
$\Qrank G_{\integer}  = \Rrank G$.
\hint{\cref{QrankEg}(\ref{QrankEg-SL},\ref{QrankEg-SOQ})).}
\end{exercises}

%\section{Lattices of $\rational$-rank one}
%
%\begin{thm} \label{Qrank1<>}
% The following are equivalent:
% \begin{enumerate}
% \item $\Qrank(\Gamma) \le 1$.
% \item Any two maximal unipotent subgroups of~$\Gamma$ either
%are equal or are disjoint.
% \end{enumerate}
% \end{thm}
%
%
%\begin{thm} \label{Qrank1<>GQ}
% If $\Gamma$ is commensurable to $G_{\integer}$, then the
%following are equivalent:
% \begin{enumerate}
% \item $\Qrank(\Gamma) \le 1$;
% \item If $U_1$ and $U_2$ are unipotent subgroups
%of~$G_{\rational}$, and $U_1 \cap U_2 \neq e$, then $\langle
%U_1,U_2 \rangle$ is unipotent.
% \item Any two maximal unipotent subgroups of~$G_{\rational}$
%either are equal or are disjoint.
% \item  \label{Qrank1<>GQ-unipP}
% If $P_1$ and $P_2$ are parabolic $\rational$-subgroups
%of~$G$, and $\unip P_1 \cap \unip P_2 \neq e$, then $P_1 = P_2$.
% \item All proper parabolic $\rational$-subgroups of~$G$ are
%conjugate to each other \textup(in fact, they are
%conjugate under $G_{\rational}$\textup).
% \end{enumerate}
% \end{thm}
% 
% \begin{exercises}
% 
%\item Show that if $\Qrank(\Gamma) = 1$, and $U_1$ and
%$U_2$ are unipotent subgroups of~$\Gamma$, such that $U_1
%\cap U_2 \neq e$, then $\langle U_1,U_2 \rangle$ is
%unipotent.
%
%\end{exercises}

\section{Lattices of higher \texorpdfstring{$\rational$}{Q}-rank}

This section closely parallels \cref{HigherRrankSect}, because the results there on semisimple groups of higher real rank can be extended in a natural way to lattices of higher $\rational$-rank.

\begin{assump}
Throughout this section, if the statement of a result mentions $G_\rational$, $G_\integer$, or a $\rational$-split torus in~$G$, then $G$ is assumed to be defined over~$\rational$.
\end{assump}

\begin{prop}[\csee{Qrank2<>C(a)gensEx}] \label{Qrank2<>C(a)gens}
Let $S$ be any maximal\/ $\rational$-split torus in~$G$. Then we have
 $\Qrank G_\integer \ge 2$
 if and only if there exist nontrivial elements $s_1$ and~$s_2$ of~$S_\rational$,
such that $G = \langle \czer_G(s_1), \czer_G(s_2)\rangle$.
\end{prop}

\begin{lem} \label{C(A)=AM}
If $S$~is any maximal $\rational$-split torus in~$G$, then we have $\czer_G(S) = S \times M = S \times CL$, where 
	\begin{itemize}
	\item $M$, $C$, and~$L$ are defined over~$\rational$,
	\item $\Qrank M = 0$,
	\item $L$ is semisimple,
	and
	\item $C$ is a torus that is the identity component of the center of~$M$.
	\end{itemize}
\end{lem}

 \begin{prop}[\csee{Qrank2<>au=uaEx}] \label{Qrank2<>au=ua}
 $\Qrank G_\integer \ge 2$ if and only if there exist nontrivial elements $a$ and~$u$ of~$G_\rational$, such that $a$~belongs to a\/ $\rational$-split torus of~$G$, $u$~is unipotent, and $au = ua$.
\end{prop}

\begin{lem}[\csee{Qrank1UniqMaxUnipEx}] \label{Qrank1UniqMaxUnip}
Assume\/ $\Gamma$ is commensurable to~$G_\integer$. The following are equivalent:
	\begin{enumerate}
	\item $\Qrank \Gamma \le 1$.
	\item Every nontrivial unipotent subgroup of\/~$\Gamma$ is contained in a \textbf{unique} maximal unipotent subgroup of\/~$\Gamma$.
	\item Every nontrivial unipotent\/ $\rational$-subgroup of~$G$ is contained in a \textbf{unique} maximal unipotent\/ $\rational$-subgroup of~$G$.
	\end{enumerate}
\end{lem}

\begin{prop} \label{Qrank2<>UnipGens}
$\Qrank \Gamma \ge 2$ if and only if\/ $\Gamma$ contains nontrivial unipotent subgroups $U_1,\ldots,U_k$, such that 
	\begin{itemize}
	\item $\langle U_1,\ldots,U_k \rangle$ is a finite-index subgroup of\/~$\Gamma$,
	and
	\item $U_i$ centralizes~$U_{i+1}$ for each~$i$.
	\end{itemize}
\end{prop}

% It would be good to give an application of unipotent generators or the following theorem @@@

\begin{prop} \label{Qrank2<>SL3orSO23}
Assume\/ $\Gamma$ is irreducible. Then $\Qrank \Gamma \ge 2$ if and only if\/ $\Gamma$~contains a subgroup that is commensurable to either\/ $\SL(3,\integer)$ or\/ $\SO(2,3)_\integer$.
\end{prop}

%\begin{thm} \label{Qrank2<>SL3orSO23}
%Assume $\Gamma$ is irreducible and $\Rrank G \ge 2$. Then:
%	\begin{enumerate}
%	\item $\Qrank \Gamma \ge 2$ if and only if $\Gamma$~contains a subgroup that is commensurable to either\/ $\SL(3,\integer)$ or\/ $\SO(2,3)_\integer$.
%	\item $\Qrank \Gamma \ge 1$ if and only if $\Gamma$~contains a subgroup that is isomorphic to a noncocompact, irreducible, arithmetic lattice in either\/ $\SL(3,\real)$, $\SL(3,\complex)$, or a group of the form $\SL(2,\real)^r \times \SL(2,\complex)^s$, with $r + s \ge 2$. \textup(These ``minimal'' lattices are either unitary groups or special linear groups.\textup)
%	\item If $\Qrank \Gamma = 0$, then $\Gamma$~contains a subgroup that is commensurable to either a special linear group or a unitary group\/ \textup(perhaps over a division algebra\textup).
%	\end{enumerate}
%\end{thm}

\begin{rems} \ 
\noprelistbreak
 \begin{enumerate}

\item Unfortunately, the list of lattices of $\rational$-rank one is longer and much more complicated than the list of simple groups of real rank one in \cref{rank1simple}. The classical arithmetic groups (of any $\rational$-rank) are described in \cref{ArithClassicalChap} (see the \lcnamecref{IrredInG} on \pageref{IrredInG}), but there are also infinitely many different lattices of $\rational$-rank one in exceptional groups of type $E_6$ and~$F_4$, and the nonarithmetic lattices of $\rational$-rank one in $\SO(1,n)$ and $\SU(1,n)$ have not yet been classified.

 \item Suppose $\Qrank \Gamma \le 1$.
\Cref{Qrank2<>UnipGens} shows that it is impossible to
find a generating set $\{\gamma_1,\ldots,\gamma_r\}$
for~$\Gamma$, such that each $\gamma_i$ is nontrivial and
unipotent, and $\gamma_i$ commutes with~$\gamma_{i+1}$, for
each~$i$. However, it is possible, in some cases, to find a
generating set $\{\gamma_1,\ldots,\gamma_r\}$ that has all
of these properties \emph{except} the requirement that
$\gamma_i$~is unipotent. For example, this is easy (up to
finite index) if $\Gamma$~is reducible \csee{RedCommGens}.

 \end{enumerate}
 \end{rems}

\begin{exercises}

\item \label{Qrank2<>C(a)gensEx}
	\begin{enumerate}
	\item Prove \cref{Qrank2<>C(a)gens}($\Rightarrow$) for the special case where we have $G_\rational = \SL(3,\rational)$.
	\item Prove \cref{Qrank2<>C(a)gens}($\Leftarrow$).
	\end{enumerate}

\item \label{CentGensGamma1xGamma2Ex}
Prove the following results in the special case where $\Gamma = \Gamma_1 \times \Gamma_2$, and $\Qrank \Gamma_i \ge 1$ for each~$i$.
	\begin{enumerate}
	\item \label{CentGensGamma1xGamma2Ex-C(a)}
	\cref{Qrank2<>C(a)gens}($\Rightarrow$)
	\item \cref{Qrank2<>au=ua}($\Rightarrow$)
	\item \cref{Qrank1UniqMaxUnip}($\Leftarrow$)
	\item \cref{Qrank2<>UnipGens}($\Rightarrow$)
	\end{enumerate}

\item \label{Qrank2<>au=uaEx}
	\begin{enumerate}
	\item Prove \cref{Qrank2<>au=ua}($\Rightarrow$) in the special case where we have $G_\rational = \SL(3,\rational)$.
	\item \label{Qrank2<>au=uaEx-not1}
	Prove \cref{Qrank2<>au=ua}($\Leftarrow$).
%	\hint{\cref{ZarAandCent}.}
	\end{enumerate}

\item \label{Qrank1UniqMaxUnipEx}
For each of these groups, find a nontrivial unipotent subgroup that is contained in two different maximal unipotent subgroups.
	\begin{enumerate}
	\item $\SL(3,\rational)$.
	\item $\SL(3,\integer)$.
	\end{enumerate}

\item \label{Qrank2<>UnipGensEx}
Prove \cref{Qrank2<>UnipGens}.

\item (\emph{Assumes the theory of $\rational$-roots})
Prove the general case of the following results.
	\begin{enumerate}
	\item \cref{Qrank2<>C(a)gens}.
	\item \cref{C(A)=AM}.
	\item \cref{Qrank2<>au=ua}($\Rightarrow$).
	\item \cref{Qrank1UniqMaxUnip}.
	\item \cref{Qrank2<>UnipGens}.
	\item \cref{Qrank2<>SL3orSO23}.
	\end{enumerate}

\item \label{RedCommGens}
 Show that if $\Gamma$ is reducible, and $G$ has no compact
factors, then there is a finite subset
$\{\gamma_1,\ldots,\gamma_r\}$ of~$\Gamma$, such that
 \begin{enumerate}
 \item $\{\gamma_1,\ldots,\gamma_r\}$ generates a
finite-index subgroup of~$\Gamma$,
 \item each $\gamma_i$ is nontrivial, and 
 \item $\gamma_i$ commutes with~$\gamma_{i+1}$, for each~$i$.
 \end{enumerate}
 
% \item \label{noCommGens}
%Suppose $\Lambda$ is a group such that:
%	\begin{enumerate}
%	\item if $\gamma_1$ and~$\gamma_2$ are any
%nontrivial elements of~$\Lambda$, such that $\gamma_1$
%commutes with~$\gamma_2$, then $\czer_{\Gamma}(\gamma_1) =
%\czer_{\Gamma}(\gamma_2)$,
%	and
%	\item there is a sequence $\lambda_1,\lambda_2,\ldots,\lambda_n$ of elements of~$\Lambda$, such 
%		\begin{enumerate}
%		\item $\{\lambda_1,\lambda_2,\ldots,\lambda_n\}$ generates~$\Lambda$,
%		and
%		\item $\lambda_i$ commutes with $\lambda_{i+1}$ for $1 \le i < n$.
%		\end{enumerate}
%	\end{enumerate}
%Show that $\Lambda$ is abelian.

\item Let $\Gamma$ be a torsion-free, cocompact lattice in $\SL(3,\real)$, constructed as in \cref{CocpctSL3Rbands}. Show that if $\gamma_1$ and~$\gamma_2$ are any
nontrivial elements of~$\Gamma$, such that $\gamma_1$
commutes with~$\gamma_2$, then $\czer_{\Gamma}(\gamma_1) =
\czer_{\Gamma}(\gamma_2)$. (Hence, it is impossible to find a
sequence of nontrivial generators of~$\Gamma$, such that each generator
commutes with the next.)
\hint{Let $D = \phi(L^3)$, so $D$ is a division ring of degree~$3$ over~$\rational$. Then $\czer_D(\gamma_1)$ is subring of~$D$ that contains the field $\rational[\gamma_1]$ in its center. Because the degree of~$D$ is prime, we conclude that
$\czer_D(\gamma_1) = \rational[\gamma_1] \subseteq \czer_D(\gamma_2)$.}

 \end{exercises}

\section{Minimal parabolic \texorpdfstring{$\rational$}{Q}-subgroups}

%Generally speaking, subgroups of~$G$ are more likely to be useful for the study of arithmetic subgroups if they are defined over~$\rational$. 

Minimal parabolic subgroups of~$G$ play an important role in the study of arithmetic subgroups, even when they are not defined over~$\rational$. However, for some purposes (especially when we construct a coarse fundamental domain in \cref{ReductionChap}), we want a subgroup that is both defined over~$\rational$ and is similar to a minimal parabolic subgroup:

\begin{defn}[(cf.~\cref{MinParabDefn})]
Let $S$ be a maximal $\rational$-split torus of~$G$, and let $a$ be a generic element of~$S$. Then the corresponding \defit[parabolic!Q-subgroup, minimal@$\rational$-subgroup, minimal]{minimal parabolic $\mpmb\rational$-subgroup} % @@@ mypmb
of~$G$ is
	$$P = \bigset{ g \in G }{ \limsup_{n \to \infty} \| a^{-n} g a^n \| < \infty } .$$
This is a Zariski closed subgroup of~$G$ that is defined over~$\rational$.
 \end{defn}

\begin{egs} \ \label{QParabEgs}
 \noprelistbreak
 \begin{enumerate}

	\item \label{QParabEgs-SLn}
	Since the group of upper triangular matrices is a minimal parabolic $\rational$-subgroup of $\SL(n,\real)$, we see that, in this case, the minimal parabolic $\rational$-subgroup is also a minimal parabolic subgroup. 
	%This is because the maximal $\rational$-split torus is also a maximal $\real$-split torus.
	
	\item This is a special case of the fact that if $\Qrank \Gamma = \Rrank G$, then every minimal parabolic $\rational$-subgroup is also a minimal parabolic subgroup \csee{MinParQisMin}. 
	
	 \item \label{QParabEgs-SOQ}
	 \fullCcf{parab=Stab(flag)}{SOmn}
	 Suppose $Q$ is a nondegenerate quadratic form on~$\rational^\ell$ that is defined over~$\rational$. A subgroup~$P$ of $\SO(Q; \real)$ is a minimal parabolic $\rational$-subgroup if and only if there is a chain $V_0 \subsetneq V_1 \subsetneq \cdots \subsetneq V_r$ of \textbf{totally isotropic} subspaces of $\rational^\ell$, such that 
	 	\begin{itemize}
		\item $\dim V_i = i$, for each~$i$,
		\item $V_r$~is a maximal totally isotropic subspace, 
		and
		\item $ P = \{\, g \in \SO(Q; \real) \mid \forall i, \ g V_i = V_i \,\} $.
		\end{itemize}
	 \end{enumerate}
 \end{egs}

We have a Langlands decomposition over~$\rational$. However, unlike in the real case, where the subgroup~$M$ is compact (i.e., $\Rrank M = 0$), we now have a subgroup that may be noncompact (but whose $\rational$-rank is~$0$):

\begin{thm}[(\thmindex{Langlands decomposition}Langlands decomposition)]
\label{LanglandsDecompQ}
 If $P$ is a minimal parabolic\/ $\rational$-subgroup of~$G$, then we may write $P$ in the form $P = MSN = LCSN$, where%
	 \begin{enumerate}
	 \item $M$, $S$, $N$, $L$, and~$C$ are defined over~$\rational$,
	 \item $S$~is a maximal $\rational$-split torus, 
	 \item \label{LanglandsDecompQ-M=0}
	 $\Qrank M = 0$,
	 \item $M S = \czer_G(S)$,
	 \item $M = LC$, where $L$ is semisimple and $C$ is the identity component of the center of~$M$,
	and 
	 \item $N$ is the {\upshape\defit[unipotent!radical]{unipotent radical}} of~$P$; that is, the unique maximal unipotent \textbf{normal} subgroup of~$P$. 
	 \end{enumerate}
Furthermore, 
%the subgroups $S$ and~$N$ are nontrivial if and only if $P \neq G$.
for some $a \in S_{\rational}$, we have
	 \begin{align} \label{aContractsP}
	  P = \bigset{ g \in G }{ \limsup_{n \to \infty} \| a^{-n}
	g a^n \| < \infty } 
	\end{align}
and
	 \begin{align}
	N = \bigset{ g \in G }{ \lim_{n \to \infty} a^{-n} g a^n = e } 
	. \end{align}
 \end{thm}

\begin{proof}
The examples and proof are essentially the same as for the real Langlands decomposition \pref{LanglandsDecomp}, but with $\rational$ in the place of~$\real$, and groups of $\rational$-rank~$0$ in place of compact groups.
\end{proof}

\begin{prop} \label{parab/Q}
 Assume $G$ is defined over~$\rational$, and $P$ is a minimal parabolic $\rational$-subgroup, with Langlands decomposition $P = MSN$. Then:
	 \begin{enumerate}
	 \item \label{parab/Q-conj}
	 Every minimal parabolic $\rational$-subgroup of~$G$ is $G_\rational$-conjugate to~$P$.
	 \item \label{parab/Q-U}
	Every unipotent $\rational$-subgroup of~$G$ is $G_\rational$-conjugate to a subgroup of~$N$.
	\item \label{parab/Q-nzer}
	$P = \nzer_G(N) = \nzer_G(P)$.
	 \end{enumerate}
 \end{prop}

\begin{defn} \label{PosWeylChamberDefn}
For $P$, $M$, $S$, $N$, $L$, and~$C$ as in \cref{LanglandsDecompQ}, the 
\defit[Weyl chamber, positive]{positive Weyl chamber}%
\index{positive!Weyl chamber|indsee{Weyl chamber, positive}}
of~$S$ (with respect to~$P$) is the set~$S^+$ of all elements~$a$ of~$S$, such that $P$ is contained in the right-hand side of \pref{aContractsP}. (Equivalently, it is the closure of the set of elements~$a$ of~$S$ for which equality holds in \pref{aContractsP}.)
\end{defn}

%\begin{cor} \label{MaxUnipQSubgrps}
% Assume $G$ is defined over~$\rational$.
% \begin{enumerate}
% \item \label{MaxUnipQSubgrps-UnipRad}
%The maximal unipotent $\rational$-subgroups of~$G$ are
%precisely the unipotent radicals of the minimal parabolic
%$\rational$-subgroups of~$G$.
% \item \label{MaxUnipQSubgrps-PisNorm}
%The minimal parabolic $\rational$-subgroups of~$G$ are
%precisely the normalizers of the maximal unipotent
%$\rational$-subgroups of~$G$.
% \item All of the maximal unipotent subgroups
%of~$G_{\rational}$ are conjugate \textup(via
%$G_{\rational}$\textup).
% \end{enumerate}
% \end{cor}

\begin{exercises}
\item \label{MinParQisMin}
Show that if we have $\Qrank \Gamma = \Rrank G$, then every minimal parabolic $\rational$-subgroup is also a minimal parabolic subgroup. 
\hint{Choose $A$ to be both a maximal $\rational$-split torus and a maximal $\real$-split torus.}

\item Show that the converse of \cref{MinParQisMin} is not true.
\hint{\cref{NoncocpctInSL3Eg}.}

\item Show that every minimal parabolic $\rational$-subgroup of~$G$ contains a minimal parabolic subgroup.
\hint{Choose a maximal $\rational$-split torus~$S$. Then choose a maximal $\real$-split torus~$A$ that contains~$S$. There is a generic element of~$A$ that is very close to a generic element of~$S$.}

\item \label{G=KP}
If $P$ is any minimal parabolic $\rational$-subgroup of~$G$, and $K$~is any maximal compact subgroup of~$G$, show that $G = KP$.
 \hint{The Iwasawa decomposition \pref{IwasawaDecomp} tells us $G = KAN$, and some conjugate of $AN$ is contained in~$P$.}

\item Assume the notation of \cref{LanglandsDecompQ}. Show that if $\Qrank G = 1$, then there is an isomorphism $\phi \colon S \stackrel{\iso}{\to} \real$, such that $\phi(S^+) = \real^+$.

\item \label{U/UZcpct}
 Show that if $U$ is a unipotent $\rational$-subgroup of $\SL(\ell,\real)$, then $U_{\integer}$ is a cocompact lattice in~$U$.
\hint{Induct on the nilpotence class of~$U$ \csee{UnipNilp}.
Note that the exponential map $\exp \colon \Lie U \to U$ is a polynomial with rational coefficients, as is its inverse, so $U_\integer$ is Zariski dense in~$U$.}
 
\item Show that if $U_1$ and $U_2$ are maximal unipotent
subgroups of~$\Gamma$, and $\Gamma$~is commensurable to
$G_{\integer}$, then there exists $g \in G_{\rational}$, such
that $g^{-1} U_1 g$ is commensurable to~$U_2$.

\end{exercises}

\section{Isogenies over \texorpdfstring{$\rational$}{Q}}

We have seen examples in which $G$ is isogenous (or even isomorphic)
to~$G'$, but the arithmetic subgroup $G_{\integer}$ is very different
from~$G'_{\integer}$. (For example, it may be the case that
$G_{\integer}$ is cocompact, but $G'_{\integer}$ is not.) This
does not happen if the isogeny is defined over~$\rational$, in the following sense:
%Because of our interest
%in restriction of scalars, we describe the theory over any
%algebraic number field~$F$. For the particularly motivated
%reader, we provide some proofs.

\begin{defn} 
\ %To avoid technical issues, assume $G$ is connected.
\noprelistbreak
	 \begin{enumerate}

	 \item A homomorphism $\phi \colon G \to G'$ is 
	 \defit[defined!over Q@over~$\rational$]{defined over~$\mpmb\rational$} %@@@ mypmb
	 if $\phi( G_{\rational} ) \subseteq G'_{\rational}$.
%	 \begin{enumerate}
%	 \item $\phi$ is differentiable,
%	 \item the derivative $d \phi_e$ is $F_\infty$-linear, and
%	 \item $\phi(G_F) \subseteq H_F$.
%	 \end{enumerate}
	 \item $G_1$ is \defit[isogenous!over~$\rational$]{isogenous to~$G_2$
	over~$\mpmb\rational$} %@@@ mypmb
	 (denoted $G_1 \approx_\rational G_2$) if there is a
	group~$G$ that is defined over~$\rational$, and isogenies $\phi_i
	\colon G \to G_i$ that are defined over~$\rational$.
	 \end{enumerate}
 \end{defn}

The following result shows that any isogeny over~$\rational$ can be thought of as a polynomial with rational coefficients.
 
\begin{defn} 
 A function $\phi \colon G \to G'$ is a \defit[polynomial!with rational coefficients]{polynomial with rational coefficients} if 
 	\begin{itemize}
	\item the matrix entries of~$\phi(g)$ can be written as polynomial functions of the coefficients of~$g$,
	and
	\item the polynomials can be chosen so that all of their coefficients are in~$\rational$.
	\end{itemize}
\end{defn}

\begin{prop} \label{isog/Q->poly}
 If $G_1 \approx_\rational G_2$, then there is a group~$G$ that is defined
over\/~$\rational$, and isogenies $\phi_i \colon G
\to G_i$ for $i = 1,2$, that are polynomials with rational coefficients.
 \end{prop}

\begin{proof}
 Given isogenies $\phi_i \colon G \to G_i$ that are defined
over~$\rational$, let 
	 $$G' =
	 \bigset{ \bigl( \phi_1(g), \phi_2(g) \bigr) 
	 }{ 
	 g \in G^\circ }
	 .$$
This is defined over~$\rational$, since $G'_\rational$ is dense \csee{QptsDense}. The projection maps $\phi_i' \colon G' \to G_i$
defined by $\phi_i'(g_1,g_2) = g_i$ are polynomials with rational coefficients.
 \end{proof}

\begin{warn}
 There are examples where $\phi \colon G_1 \to G_2$ is an
isomorphism, and $\phi$ is a polynomial, but
$\phi^{-1}$ is not a polynomial. For example, the natural
homomorphism $\phi \colon \SL(3,\real) \to
\PSL(3,\real)^\circ$ is an isomorphism (because
$\SL(3,\real)$ has no center). However, there is no
isomorphism from $\PSL(3,\complex)$ to $\SL(3,\complex)$
(because one of these groups has a center and the other does
not), so the inverse of~$\phi$ cannot be a polynomial
(because it does not extend to a well-defined map between the
complexifications).
 \end{warn}

The following fundamental result implies that different
embeddings of~$G$ with the same $\rational$-points have
essentially the same $\integer$-points.

\begin{prop} \label{phi(GZ)}
 Suppose $\phi \colon G \to G'$ is a surjective homomorphism
that is defined over~$\rational$. Then $\phi (G_\integer)$ is
commensurable to~$G'_\integer$.
 \end{prop}

\begin{proof}
% Let us assume $F = \rational$ and $\ints = \integer$.
%(The general case follows from this by Restriction of
%Scalars.) 
From the proof of \cref{isog/Q->poly}, we see
that, after replacing~$G$ with an isogenous group, we may
assume that $\phi$ is a polynomial with rational
coefficients. Assume $G \subseteq \SL(\ell,\real)$ and $G'
\subseteq \SL(m,\real)$.

Define $\widetilde\phi \colon G \to \Mat_{m \times m}(\real)$ by
$\widetilde\phi(x) = \phi(x - \Id)$. Then $\widetilde\phi$ is a polynomial, so it is defined on all of $\Mat_{\ell \times
\ell}(\real)$. Since the coefficients are in~$\rational$, there
is some nonzero $n \in \natural$, such that $\widetilde\phi\bigl( n
\Mat_{\ell \times \ell}(\real)  \bigr) \subseteq \Mat_{m \times
m}(\integer)$. Therefore, letting $\Gamma_n$ be the ``\term{principal
congruence subgroup}'' of~$G_{\integer}$ of level~$n$ (see page~\pageref{PrincCongSubgrp}), we have
$\phi(\Gamma_n) \subseteq G'_{\integer}$.

Because $\Gamma_n$ is a lattice in~$G$ (and
$\phi(\Gamma_n)$ is discrete), we know that
$\phi(\Gamma_n)$ is a lattice in~$G'$. Since
$\phi(\Gamma_n) \subseteq G'_{\integer}$, this implies that
$\phi(\Gamma_n)$ is commensurable to $G'_{\integer}$
\csee{finext->latt}.
 \end{proof}

%\begin{cor}
% Suppose $G$ and~$H$ are subgroups of $\SL(\ell,\real)$
%that are defined over~$\rational$. A differentiable
%homomorphism $\phi \colon G \to H$ is defined
%over~$\rational$ if and only if some finite-index subgroup
%of $\phi(G_{\integer})$ is contained in $H_{\rational}$.
% \end{cor}

The following fundamental fact is, unfortunately, not obvious from our definition of ``$\rational$-split\zz.''

\begin{prop} \label{QsplitInvtIsog}
Assume 
	\begin{itemize}
	\item $T$ and~$H$ are connected Lie groups that are defined over\/~$\rational$,
	and
	\item $T \approx_\rational H$.
	\end{itemize}
Then $T$ is a\/ $\rational$-split torus if and only if $H$~is a\/ $\rational$-split torus.
\end{prop}

\begin{cor} \label{Qrank(isog)=}
 If $G \approx_\rational G'$, then $\Qrank G_\integer = \Qrank G'_\integer$.
 \end{cor} 

\begin{proof}
Suppose $G$ is a $\rational$-group, and there is an isogeny $\varphi_i \colon G \to G_i$ that is defined over~$\rational$ for $i = 1,2$.
If $T_1$ is a maximal $\rational$-split torus in~$G_1$, then \cref{QsplitInvtIsog} implies that $\varphi_2 \bigl( \varphi_1^{-1}(T_1)^\circ \bigr)$ is a $\rational$-split torus in~$G_2$. Since isogeny preserves dimension, we conclude that $\Qrank G_1 \le \Qrank G_2$. By symmetry, equality must hold.
\end{proof}

%\begin{exercises}
%
%\item It would be good to have some exercises here @@@
%
%\end{exercises}

\begin{notes}

As was mentioned in the notes of \cref{RrankChap}, the comprehensive treatise of Borel and Tits \cite{BorelTits-GrpRed} is the standard reference on rank, parabolic subgroups, and other fundamental properties of reductive groups over any field (including~$\rational$).
Abbreviated accounts can be found in many texts, including \cite[\S10 and \S11]{Borel-IntroGrpArith} and  \cite[Chap.~5]{Borel-LinAlgGrps}.

See \cite[Rem.~8.11, p.~60]{Borel-IntroGrpArith} for a proof of \cref{phi(GZ)}.

\end{notes}

 %!TEX root = IntroArithGrps.tex

\mychapter{Quasi-Isometries}
\label{QuasiChap}

\prereqs{none.}

\section{Word metric and quasi-isometries} \label{WordMetricSect}

The field of \term{Geometric Group Theory} equips groups with a metric, which allows them to be studied as metric spaces:

\begin{defn} \label{WordMetricDefn}
Fix a finite generating set~$S$ of~$\Gamma$ \csee{GammaFinGen}, and assume, for simplicity, that $S$ is \defit[symmetric!generating set]{symmetric}, which means $s^{-1} \in S$ for every $s \in S$.
\noprelistbreak 
	\begin{enumerate}
	\item For $g \in \Gamma$, the \defit[word!length]{word length} of~$g$ is the length~$\ell$ of the shortest sequence $(s_1,s_2,\ldots,s_\ell)$ of elements of~$S$, such that $s_1s_2\cdots s_\ell = g$. It is denoted $\ell(g)$.
	(By convention, $\ell(e) = 0$.)
	\item For $g,h \in \Gamma$, we let $d(g,h) = \ell(g^{-1} h)$. This defines a metric on~$\Gamma$, called the \defit[word!metric]{word metric} \csee{WordMetricEx}.
	%So $\Gamma$ is a metric space.
	\end{enumerate}
\end{defn}

The word metric has the important property that the action of~$\Gamma$ on itself by left-translations is an action by isometries \csee{LeftTransIsom}.

\begin{rem}
The word metric can be pictured geometrically, by constructing a \defit[Cayley!graph]{Cayley graph}. Namely, $\Cay(\Gamma; S)$ is a certain $1$-dimensional simplicial complex (or ``graph''):
	\begin{itemize} 
	\item its $0$-skeleton is~$\Gamma$,
	and
	\item it has a $1$-simplex (or ``edge'') of length~$1$ joining $v$ to~$vs$, for every $v \in \Gamma$ and $s \in S$.
	\end{itemize}
Define a metric on $\Cay(\Gamma; S)$ by letting $d(x,y)$ be the length of the shortest path from~$x$ to~$y$. Then the restriction of this metric to the $0$-skeleton is precisely the word metric on~$\Gamma$.
\end{rem}

Unfortunately, the word metric on~$\Gamma$ is not canonical, because it depends on the choice of the generating set~$S$ \csee{WordMetricNotWellDefd}. However, it is ``almost'' well-defined, in that changing the generating set can only distort the distances by a bounded factor. This idea is formalized in the following notion:

\begin{defn} \label{QuasiIsomDefn}
 Let $X_1$ and~$X_2$ be metric spaces, with metrics~$d_1$ and~$d_2$, respectively. 
 	\begin{enumerate}
	\item A function $f \colon X_1 \to X_2$ is a \defit[quasi-!isometry]{quasi-isometry} if there is a constant $C > 0$, such that
		 \begin{enumerate}
		 \item
		 for all $x,y \in X_1$ with $d_1(x,y) > C$,
		we have
		 $$ \frac{1}{C}
		 < \frac{ d_2 \bigl( f(x), f(y) \bigr) }
		 { d_1 ( x, y) }
		 < C ,$$ and
		 \item for all $x_2 \in X_2$, there exists $x_1 \in
		X_1$, such that 
			$$d_2 \bigl( f(x_1), x_2 \bigr) < C .$$
		 \end{enumerate}
 	Note that $f$ need not be continuous (and need not be one-to-one or onto, either).
	
	\item We say $X_1$ is \defit[quasi-!isometric]{quasi-isometric} to~$X_2$ (and write $X_1 \QI X_2$) if there is a quasi-isometry from~$X_1$ to~$X_2$. This is an equivalence relation \csee{QIEquivRel}.
	\end{enumerate}
 \end{defn}

\begin{prop}[\csee{WordWellDefEx}] \label{WordWellDef}
Let 
	\begin{itemize}
	\item $S_1$ and~$S_2$ be two finite, symmetric generating sets for~$\Gamma$,
	and
	\item $d_i$ be the word metric on~$\Gamma$ corresponding to the generating set~$S_i$.
	\end{itemize}
Then $(\Gamma,d_1) \QI (\Gamma,d_2)$.
\end{prop}

Therefore, if $\Gamma_1$ and~$\Gamma_2$ are quasi-isometric for some choice of the word metrics on the two groups, then they are quasi-isometric for all choices of the word metrics. So it makes sense to say that two groups are quasi-isometric, without any mention of generating sets (as long as both of the groups are finitely generated).

\begin{rem}
A property is said to be \defit[geometric property]{geometric} if is is invariant under quasi-isometry. For example, we will see (in \cref{AmenQIinvt,TnotQI}, respectively) that amenability is a geometric property, but Kazhdan's property~$(T)$ is not. In other words, if $\Lambda_1 \QI \Lambda_2$, and $\Lambda_1$ is amenable, then $\Lambda_2$~is amenable, but the same cannot be said for Kazhdan's property~$(T)$. 
In general, quasi-isometric groups can be very different from each other, so most of the usual algebraic properties of groups are not geometric.
\end{rem}

Quasi-isometries also arise from cocompact actions. Before stating the result, we introduce some terminology.

\begin{defn}
Let $(X,d)$ be a metric space, and let $C > 0$.
\noprelistbreak
	\begin{enumerate}
	\item $X$ is  \defit[proper metric space]{proper} if the closed ball $B_r(x)$ is compact, for all $r > 0$ and all $x \in X$.
	\item Let $x,y \in X$. A ($C$)-\defit[coarse!geodesic]{coarse geodesic} from $x$ to~$y$ is a finite sequence $x_0,x_1,\ldots,x_n$ in~$X$, such that $x_0 = x$, $x_n = y$, and
			$$ \text{$\bigl| d(x_i, x_j) - |i - j| \bigr| < C$ \ for all $i,j$} .$$
	\item $X$ is ($C$-)\defit{coarsely geodesic} if, for all $x,y \in X$, there is a $C$-coarse geodesic from~$x$ to~$y$.
	\end{enumerate}
\end{defn}

\begin{prop}[\csee{CocpctActIsQIEx}] \label{CocpctActIsQI}
Suppose
	\begin{itemize}
	\item $(X,d)$ is a metric space that is proper and coarsely geodesic,
	\item $\Gamma$ has a properly discontinuous action on~$X$ by isometries, such that $\Gamma \backslash X$ is compact,
	and
	\item $d'$~is a word metric on~$\Gamma$.
	\end{itemize}
Then $(\Gamma,d') \QI (X,d)$.

More precisely, for any basepoint $x_0 \in X$, the map $\gamma \mapsto \gamma x_0$ is a quasi-isometry from\/~$\Gamma$ to~$X$.
\end{prop}

\begin{cor}
 If $G/\Gamma$ is compact, then the inclusion\/
$\Gamma \hookrightarrow G$ is a quasi-isometry, where we
use any word metric on\/~$\Gamma$, and we use any\/
\textup(coarsely geodesic, proper\/\textup) metric on~$G$ that is invariant under left-translations.
 \end{cor}
 
This implies that any two cocompact lattices in the same group are quasi-isometric:

\begin{cor} \label{CpctLattQI}
If\/ $\Gamma_1$ and~$\Gamma_2$ are cocompact lattices in~$G$, then $\Gamma_1 \QI \Gamma_2$.
\end{cor}

\begin{proof}
We have $\Gamma_1 \QI G \QI \Gamma_2$, so $\Gamma_1 \QI \Gamma_2$ by transitivity.
\end{proof}

We will see in \cref{QIRigSect} that the situation is usually very different for lattices that are not cocompact: in most cases, there are infinitely many different (noncocompact) lattices in~$G$ that are not quasi-isometric to each other.

Any (coarsely geodesic, proper) metric on~$G$ provides a metric on~$\Gamma$, by restriction. In most cases, this restriction is the word metric (up to quasi-isometry):

\begin{thm}[(\thmindex{Lubotzky-Mozes-Raghunathan}Lubotzky-Mozes-Raghunathan)] \label{LMRThm}
 If\/ $\Rrank G \ge 2$, and\/ $\Gamma$ is irreducible, then the inclusion\/ $\Gamma \hookrightarrow G$ is a quasi-isometry onto its image.
 \end{thm}

The assumption that $\Rrank G \ge 2$ is essential:

\begin{eg}
Let 
\noprelistbreak
	\begin{itemize}
	\item $G = \SL(2,\real)$,
	\item $\Gamma$ be a free subgroup of finite index in $\SL(2,\integer)$ \csee{SanovIsFree}, 
	\item $u = \begin{bmatrix} 1 & k \\ 0 & 1 \end{bmatrix} \in \Gamma$,
	and
	\item $a^t =  \begin{bmatrix} e^t & 0 \\ 0 & e^{-t} \end{bmatrix} \in G$.
	\end{itemize}
Then:
	\begin{enumerate}
	\item For any word metric $d_\Gamma$ on~$\Gamma$, the function $d_\Gamma(u^n, e)$ grows linearly with~$n$, because $\Gamma$ is free.
	\item For any left-invariant metric $d_G$ on~$G$, the function $d_G(u^n, e)$ grows only logarithmically, because $a^{\log n} u a^{-\log n} = u^{2n}$.
	\end{enumerate}
Therefore, the restriction of~$d_G$ to~$\Gamma$ is not quasi-isometric to~$d_\Gamma$.
\end{eg}

\begin{exercises}

\item \label{WordMetricEx}
Show that the word metric is indeed a metric. More precisely, for $x,y,z \in \Gamma$, show 
%	$d(x,y) \ge 0$,
%	$d(x,y) = 0 \iff x = y$,
%	$d(x,y) = d(y,x)$,
%	and
%	$d(x,y) \le d(x,z) + d(z,y)$.
	$$ \begin{matrix}
	d(x,y) \ge 0, 
	&\quad &
	d(x,y) = 0 \iff x = y,
	\\
	d(x,y) = d(y,x),
	&&
	d(x,y) \le d(x,z) + d(z,y)
	. \end{matrix}$$

\item \label{LeftTransIsom}
Assume $d$ is a word metric on~$\Gamma$ (with respect to a finite, symmetric generating set~$S$). Show that $d(a x, a y) = d(x,y)$ for all $a,x,y \in \Gamma$.
\\ {\smaller{}[\emph{Warning:} $d(xa, ya)$ is usually \emph{not} equal to $d(x,y)$.]}

\item \label{WordMetricNotWellDefd}
Assume $\Gamma$ is infinite. Show there exist two word metrics $d_1$ and~$d_2$ on~$\Gamma$ (corresponding to finite, symmetric generating sets $S_1$ and~$S_2$), such that the metric space $(\Gamma_1,d_1)$ is not isometric to the metric space $(\Gamma_2,d_2)$.
\hint{A ball of radius~$r$ can have different cardinality for the two metrics.}

\item \label{QIEquivRel}
%For convenience, write $\Lambda_1 \sim \Lambda_2$ if $\Lambda_1$ is quasi-isometric to~$\Lambda_2$ (where $\Lambda_1$ and~$\Lambda_2$ are finitely generated groups with generating sets $S_1$ and~$S_2$). 
Show that $\QI$ is an equivalence relation.

\item \label{WordWellDefEx}
Prove \cref{WordWellDef}.
\hint{Show the identity map is a quasi-isometry from $(\Gamma, d_1)$ to $(\Gamma, d_2)$, by choosing $C$ so that $d_1(e, s) \le C$ for each $s \in S_2$.}

%\item \label{CpctLattQIEx}
%Show that if $\Gamma_1$ and~$\Gamma_2$ are two \emph{cocompact} lattices in the same group~$G$, then $\Gamma_1 \QI \Gamma_2$.

\item \label{CocpctActIsQIEx}
Prove \cref{CocpctActIsQI}.
\hint{Assume $S \supseteq \{\, \gamma \mid \exists x \in X, \ d(\gamma x, x) \le 3C\,\}$, and $\Gamma \cdot B_C(x) = X$ for all $x \in X$. Given $x_0,x_1,\ldots,x_n \in X$ with $d(x_i,x_{i+1}) \le C$, there exists $\gamma_i \in \Gamma$, such that $d (\gamma_i x_0, x_i) \le C$, so $\ell(\gamma_n) \le n$.}

\end{exercises}

\section{Hyperbolic groups} \label{GromovHyperGrpsSect}

Manifolds of negative sectional curvature play an important role in differential geometry, both in applications and as a source of examples. 
The triangles in these manifolds have a very special ``thinness'' property that we will now explain. Groups whose triangles have this same property are said to be ``\term[negatively curved group]{negatively curved}\zz,'' or, in our terminology, ``Gromov hyperbolic\zz.''

\begin{defn}[(Gromov)] \label{GromovHyperDefn}
Let $\delta > 0$, and let $X$ be a $C$-coarsely geodesic metric space.
\noprelistbreak
	\begin{enumerate}
	\item A ($C$-coarse) \defit[coarse!triangle]{triangle} $abc$ in~$X$ is a set $\{s_{ab} , s_{bc} , s_{ac}\}$, where $s_{xy}$ is a $C$-coarse geodesic from~$x$ to~$y$ for $x,y \in \{a,b,c\}$.
	\item A triangle $abc$ is 
		\defit[delta-thin@$\delta$-thin]{$\delta$-thin}%
		\index{thin, $\delta$-|indsee{$\delta$-thin}}
	if each of the three sides of the triangle is contained in the (closed) $\delta$-neighborhood of the union of the other two sides. That is, each point in~$s_{ab}$ is at distance no more than~$\delta$ from some point in $s_{ac} \cup s_{bc}$, and similarly for $s_{ac}$ and~$s_{bc}$.
	\item $X$ is \defit[hyperbolic!Gromov]{Gromov hyperbolic} if there exists some $\delta > 0$, such that every ($C$-coarse) triangle in~$X$ is $\delta$-thin. 
	\end{enumerate}
\end{defn}

\begin{thm} \label{NegCurvGromovHyp}
The universal cover of any compact manifold of strictly negative sectional curvature is Gromov hyperbolic.
\end{thm}

\begin{proof}[Idea of proof] As an illustration, let us show that the hyperbolic plane~$\hyperbolic^2$ is Gromov hyperbolic. We use the disk model.

\vskip-0.2\baselineskip % !!!

\vbox{% we can't have a page break in this paragraph !!!
\tolerance=3000
\vbox to 0pt{\vskip 1.5\baselineskip\rightline{\includegraphics{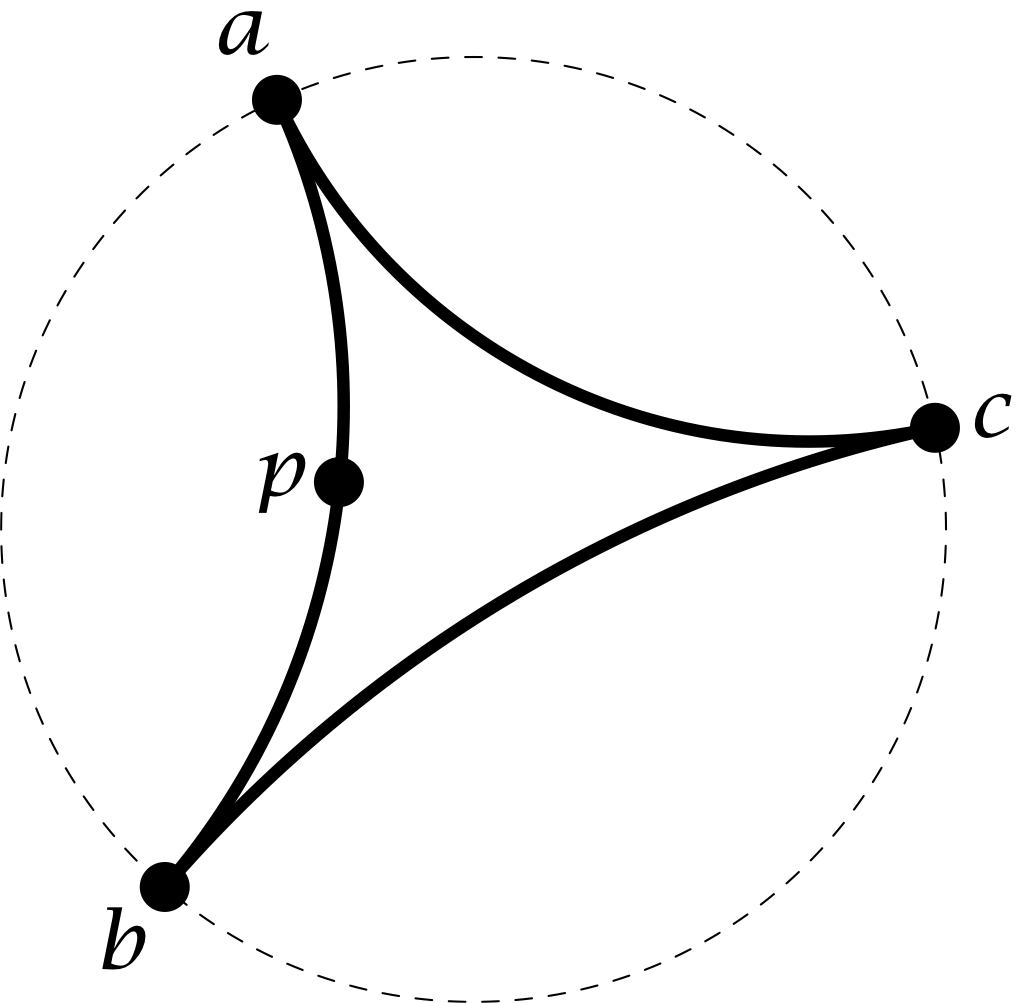}}\vss}%
\hangindent=-1.85in \hangafter=0 
Any three distinct points $a$, $b$, and~$c$ on $\bdry\hyperbolic^2$ are the vertices of an ideal triangle (with geodesic sides).
%texpreamble
%("  \usepackage{amsmath}
% \usepackage[LY1]{fontenc}
% \usepackage[expert,LY1,mylucidascale]{mylucidabr}
% ");
%defaultpen(  fontcommand("\normalfont") + fontsize(10) ); 
%
%from graph access *;
%unitsize(2cm);
%
%real linethick = 1.5;
%real dotthick = 12;
%dotfactor=12;
%
%real atheta = 2;
%pair a = ( cos(atheta), sin(atheta) );
%
%real btheta = 4;
%pair b = ( cos(btheta), sin(btheta) );
%
%real ctheta = 6.5;
%pair c = ( cos(ctheta), sin(ctheta) );
%
%pair p = (-0.285, 0.1);
%dot( p ); label( "$p$", p, W );
%
%draw( circle((0,0), 1), linewidth(0.25)+dashed );
%
%dot( a ); label( "$a$", a, 2*a );
%dot( b );label( "$b$", b, b );
%dot( c );label( "$c$", c, 1.5*c );
%
%draw( a{-a}..{b}b , linewidth(linethick) );
%draw( b{-b}..{c}c , linewidth(linethick));
%draw( c{-c}..{a}a , linewidth(linethick));
%
Choose a point~$p$ on $\overline{ab}$.
Since the geodesic ray $\overrightarrow{pa}$ is asymptotic to $\overrightarrow{ca}$, there is some $\delta > 0$, such that every point of $\overrightarrow{pa}$ is in the $\delta$-neighborhood of $\overrightarrow{ca}$. Similarly, every point of $\overrightarrow{pb}$ is in the $\delta$-neighborhood of $\overrightarrow{cb}$ (after we enlarge~$\delta$). Therefore, all of $\overline{ab}$ is in the $\delta$-neighborhood of the union of the other two sides. By applying the same argument to the sides $\overline{bc}$ and $\overline{ac}$, we see there is some $\delta$, such that the triangle $abc$ is $\delta$-thin.
}

Since the isometry group $\SL(2,\real)$ acts transitively on the (unordered) triples of distinct points on the boundary, we conclude that every ideal triangle is $\delta$-thin for this same value of~$\delta$. Having vertices on the boundary is the worst-case scenario, so this implies that all geodesic triangles in $\hyperbolic^2$ are $\delta$-thin.
\end{proof}

The following important ``shadowing property'' tells us that ``quasi-geodesics'' are always close to geodesics. Unfortunately, its proof is somewhat lengthy, so we omit it. (See \cref{HyperbolicGeodUnique,HyperbolicLargerC} for some weaker results that are easier.) 
% Maybe should also prove logarithmic distance, since that's not hard? @@@

\begin{thm} \label{CoarseShadowing}
Suppose 
	\begin{itemize}
	\item $X$ is a $C$-coarsely geodesic $\delta$-hyperbolic metric space,
	and
	\item $\{x_0,x_1,\ldots,x_n\}$ is finite sequence of points in~$X$, such that, for all $i,j$ we have:
		$$ \frac{|i - j|}{C} - C \le d(x_i, x_j) \le C|i - j| + C .$$
	\end{itemize}
Then there exists $C' > 0$ \textup(depending only on~$C$ and~$\delta$\textup), such that the set $\{x_0,x_1,\ldots,x_n\}$ is contained in the $C'$-neighborhood of every $C$-coarse geodesic from~$x_0$ to~$x_1$.
\end{thm}

This implies that being Gromov hyperbolic is invariant under quasi-isometry:

\begin{cor}[\csee{HyperbolicQIInvtEx}] \label{HyperbolicQIInvt}
Assume
\noprelistbreak
	\begin{itemize}
	\item $X_1$ and~$X_2$ are coarsely geodesic metric spaces,
	and
	\item $X_1 \QI X_2$.
	\end{itemize}
If $X_1$ is Gromov hyperbolic, then $X_2$ is Gromov hyperbolic.
\end{cor}

\begin{cor}
The fundamental group of any compact manifold~$M$ of strictly negative sectional curvature is Gromov hyperbolic.
\end{cor}

\begin{proof}
Since $M$ is compact, the fundamental group $\pi_1(M)$ acts cocompactly on the universal cover~$\cover M$ of~$M$, so $\pi_1(M) \QI \cover M$ \csee{CocpctActIsQI}. Now apply \cref{NegCurvGromovHyp,HyperbolicQIInvt}.
\end{proof}

This observation allows us to determine precisely which lattices are Gromov hyperbolic:

\begin{prop} \label{WhichLattsHyper}
 $\Gamma$ is Gromov hyperbolic if and only if\/
$\Rrank G = 1$ and either $G/\Gamma$ is
compact, or the unique noncompact simple factor of~$G$ is
isogenous to\/ $\SL(2,\real)$.
 \end{prop}

\begin{proof}[Sketch of proof]
With a bit more theory than has been presented here, it is not difficult to show that Gromov hyperbolic groups never contain a subgroup isomorphic to $\integer \times
\integer$, so we may assume $\Gamma$ is irreducible.

\setcounter{case}{0}

\begin{case}
Assume $G/\Gamma$ is compact.
 \end{case}
 From \cref{CocpctActIsQI}, we know that $\Gamma$ is quasi-isometric to the symmetric space $G/K$ associated to~$G$.
 	\begin{itemize}
	\item If $\Rrank G = 1$, then $G/K$ has negative sectional curvature, bounded away from~$0$, so it is Gromov hyperbolic. 
	\item If $\Rrank G \ge 2$, then $G/K$ contains $2$-dimensional
flats, so it is not Gromov hyperbolic.
	\end{itemize}

\begin{case}
Assume $G/\Gamma$ is not compact.
 \end{case}
 $(\Rightarrow)$ We may assume that $G$ has no compact factors. Let $U$ be a maximal unipotent subgroup of~$\Gamma$. We know that $U$ does not
contain a subgroup isomorphic to $\integer \times
\integer$. Also, since $G/\Gamma$ is not compact, we know $U$ is infinite \csee{GodementConverse}. Therefore, since $U$ is nilpotent (and torsion-free), it is easy to see that $U$ must be cyclic. It can be shown that this implies $G$ is isogenous to
$\SL(2,\real)$.

$(\Leftarrow)$ $\Gamma$ is virtually free, so it is Gromov hyperbolic \csee{FreeGrpsAreHyper}.
 \end{proof}

\begin{exercises}

\item \label{HyperbolicGeodUnique}
Show that the coarse geodesic between two points in a Gromov hyperbolic space is coarsely unique. More precisely, given $C$, show there is some $C' > 0$, such that if $\gamma$ and~$\gamma'$ are two $C$-coarse geodesics with the same endpoints, then $\gamma$ is contained in the $C'$-neighborhood of~$\gamma'$.
\hint{If $a$ and~$b$ are the two endpoints, consider the (degenerate) triangle $abb$.}

\item \label{HyperbolicLargerC}
Show that if $X$ is a $C$-coarsely geodesic $\delta$-hyperbolic space, and $C' \ge C$, then there exists~$\delta'$, such that every $C'$-coarse geodesic from~$a$ to~$b$ is in the $\delta'$-neighborhood of every $C$-coarse geodesic from~$a$ to~$b$.
(This is a generalization of \cref{HyperbolicGeodUnique}.)
\hint{For any point~$c$ on a $C'$-coarse geodesic from~$a$ to~$b$, there is a $C$-coarse triangle $abc$. If $c$ is not in the $\delta$-neighborhood of $\overline{ab}$, then there exist $a',b' \in \overline{ab}$ that are distance less than~$\delta$ from points $a''$ and~$b''$ on $\overline{ac}$ and $\overline{bc}$, respectively, such that $d(a',b') < C + 1$. Bound $d(a',c)$ by noting that $c$ is on the $C''$-coarse geodesic $\overline{a''c} \cup \overline{cb''}$.}

\item \label{FreeGrpsAreHyper}
Show that free groups are Gromov hyperbolic.
\hint{The word metric corresponding to a set of free generators is $0$-hyperbolic.}

\item \label{HyperbolicQIInvtEx}
Prove \cref{HyperbolicQIInvt}.
\hint{Use \cref{CoarseShadowing} to show that coarse triangles in~$X_2$ can be approximated by coarse triangles in~$X_1$.}

\end{exercises}

\begin{notes}

Almost all of the material in this chapter can be found in any treatment of geometric group theory, such as \cite{GhysDelaHarpe-InfGrpsGeomObjs}, \cite{delaHarpe-TopicsGeomGrpThy}, or (more elementary)~\cite{Bowditch-CourseGGT}. A detailed treatment of this and much more is in \cite{BridsonHaefliger}. 

The Lubotzky-Mozes-Raghunathan Theorem \pref{LMRThm} is proved in \cite{LubotzkyMozesRaghunathan} (or see \cite{LubotzkyMozesRaghunathan-CR} for an exposition of the special case where $\Gamma = \SL(n,\integer)$). 

See \cite{Ghys-GrpHyp} for an introduction to the theory of Gromov hyperbolic groups, or \cite{GhysDelaHarpe-GrpsHyp,Gromov-HypGrps} for much more information.
The notion of $\delta$-hyperbolic group is credited to E.\,Rips, who also proved some of the basic properties (such as that they are finitely presented), but much of the foundational work in the subject was done by M.\,Gromov \cite{Gromov-HypGrps}.

\Cref{WhichLattsHyper} is well known.

\end{notes}

 %!TEX root = IntroArithGrps.tex

\mychapter{Unitary Representations} \label{UnitaryRepChap}

\prereqs{none.}

Unitary representations are of the utmost importance in the study of Lie groups. For our purposes, one of the main applications is the proof of the Moore Ergodicity Theorem \pref{MooreErgBasicThm} in \cref{MooreErgPfSect}, but they are also the foundation of the definition (and study) of Kazhdan's Property~$(T)$ in \cref{KazhdanTChap}.

\section{Definitions} \label{UnitaryRepSect}

\begin{defn}
Assume $\Hilbert$ is a Hilbert space, with inner product $\langle~\mid~\rangle$,
\nindex{$\Hilbert$ = Hilbert space with inner product $\langle~\mid~\rangle$}
and $H$~is a Lie group.
	\begin{enumerate}
	
	\item $\unitary(\Hilbert)$ is the group of {unitary operators} on~$\Hilbert$. 
	\nindex{$\unitary(\Hilbert)$ = group of unitary operators on~$\Hilbert$}
%	Therefore, 
%	$$ \unitary(\Hilbert) = \bigset{
%		T \colon \Hilbert \to \Hilbert
%		}{
%		\begin{matrix}
%		\text{$T$ is linear, and} \\
%		\langle Tv \mid Tw \rangle = \langle v \mid w \rangle, \ 
%		\forall v,w \in \Hilbert
%		\end{matrix}
%		} .$$
	
	\item A \defit[unitary!representation]{unitary representation} of the Lie group~$H$ on the Hilbert space~$\Hilbert$ is a homomorphism 
	$\pi \colon H \to \unitary(\Hilbert)$, 
	\nindex{$\pi$ = unitary representation}
	such that the map $h \mapsto \pi(h) \varphi$ is continuous, for each $\varphi \in \Hilbert$. (If we wish to spell out that a unitary representation is on a particular Hilbert space~$\Hilbert$, we may refer it as $(\pi, \Hilbert)$, rather than merely~$\pi$.) 
	
	\item The \defit[dimension of a representation]{dimension} of a unitary representation $(\pi,\Hilbert)$ is the dimension of the Hilbert space~$\Hilbert$.

	\item Suppose $(\pi_1, \Hilbert_1)$ and $(\pi_2, \Hilbert_2)$ are unitary representations  of~$H$.
\noprelistbreak
		\begin{enumerate}
		
		\item The \defit[direct sum!of representations]{direct sum} of the representations $\pi_1$ and~$\pi_2$ is the unitary representation $\pi_1 \oplus \pi_2$ of~$H$ on $\Hilbert_1 \oplus \Hilbert_2$ that is defined by 
			$$(\pi_1 \oplus \pi_2)(h)(\varphi_1,\varphi_2) = \bigl( \pi_1(h) \varphi_1, \pi_2(h) \varphi_2 \bigr), $$
		for $h \in H$ and $\varphi_i \in \Hilbert_i$.
		
		\item $\pi_1$ and~$\pi_2$ are \defit[isomorphic unitary representations]{isomorphic} if there is a Hilbert-space isomorphism $T \colon \Hilbert_1 \stackrel{\iso}{\longrightarrow} \Hilbert_2$ that \defit[intertwining operator]{intertwines} the two representations. This means $T \bigl( \pi_1(h) \varphi \bigr) = \pi_2(h) \, T(\varphi)$, for all $h \in H$ and $\varphi \in \Hilbert_1$.
		\end{enumerate}

	\end{enumerate}
\end{defn}

\begin{eg}
Every group~$H$ has a \defit[representation!trivial]{trivial representation}, denoted by~$\trivrep$
	\nindex{$\trivrep$ = trivial representation}%
(or $\trivrep_H$, if it will avoid confusion). It is a unitary representation on the $1$-dimensional Hilbert space~$\complex$, and is defined by $\trivrep(h) \varphi = \varphi$ for all $h \in H$ and $\varphi \in \complex$.
\end{eg}

Here is a more interesting example:

\begin{eg} \label{RepL2}
Suppose 
	\begin{itemize}
	\item $H$ is a Lie group,
	\item $H$ acts continuously on a locally compact, metrizable space~$X$, 
	and 
	\item $\mu$ is an $H$-invariant Radon measure on~$X$.
	\end{itemize}
Then there is a unitary representation of~$H$ on $\LL2(X,\mu)$,
	% be translations yields a representation $\pi \colon G \to \unitary \bigl( \LL2(X) \bigr)$, 
defined by
	$$ \bigl( \pi(h) \varphi \bigr)(x) = \varphi(h^{-1} x) $$
\ccf{GContOnLp}.
For the action of~$H$ on itself by translations (on the left), the resulting representation~$\regrep$\nindex{$\regrep$ = regular representation}
of~$H$ on $\LL2(H)$ is called the (left) \defit[regular!representation]{regular representation}\index{representation!regular|indsee{regular representation}} of~$H$.
\end{eg}

\begin{defn} \label{OrthoCompDefn}
Suppose $\pi$ is a unitary representation of~$H$ on~$\Hilbert$, $H'$~is a subgroup of~$H$, and $\mathcal{K}$ is a closed subspace of~$\Hilbert$.
	\begin{enumerate}
	\item $\mathcal{K}$ is \defit[invariant!subspace]{$H'$-invariant} if $\pi(h') \mathcal{K} = \mathcal{K}$, for all $h' \in H'$. (If the representation~$\pi$ is not clear from the context, we may also say that $\mathcal{K}$ is $\pi(H')$-invariant.)
	
	\item For the special case where $H' = H$, an $H$-invariant subspace is simply said to be \defit[invariant!subspace]{invariant}, and the representation of~$H$ on any such subspace is called a \defit{subrepresentation} of~$\pi$. More precisely, if $\mathcal{K}$ is $H$-invariant (and closed), then the corresponding subrepresentation is the unitary representation $\pi_{\mathcal{K}}$ of~$H$ on~$\mathcal{K}$, defined by $\pi_{\mathcal{K}}(h)\varphi = \pi(h) \varphi$, for all $h \in H$ and $\varphi \in \mathcal{K}$.
	\end{enumerate}
\end{defn}

\begin{lem}[\csee{subspace->sumEx}] \label{subspace->sum}
If $(\pi,\Hilbert)$ is a unitary representation of~$H$, and $\mathcal{K}$ is a closed, $H$-invariant subspace of~$\Hilbert$, then $\pi \iso \pi_{\mathcal{K}} \oplus \pi_{\mathcal{K}^\perp}$.
\end{lem}

The above \lcnamecref{subspace->sum} shows that any invariant subspace leads to a decomposition of the representation into a direct sum of subrepresentations. This suggests that the fundamental building blocks are the representations that do not have any (interesting) subrepresentations. Such representations are called ``irreducible:''

\begin{defns}
Let $H$ be a Lie group.
\noprelistbreak
	\begin{enumerate}
	\item A unitary representation $(\pi,\Hilbert)$ of~$H$ is \defit[irreducible!representation]{irreducible} if it has no nontrivial, proper, closed, invariant subspaces. That is, the only closed, $H$-invariant subspaces of~$\Hilbert$ are $\{0\}$ and~$\Hilbert$.
	\item The set of all irreducible representations of~$H$ (up to isomorphism) is called the \defit[unitary!dual]{unitary dual} of~$H$, and is denoted $\widehat H$.
	\end{enumerate}
\end{defns}

\begin{warns} \label{UnitaryWarns} \ 
\noprelistbreak
	\begin{enumerate}
	\item Unfortunately, it is usually not the case that every unitary representation of~$H$ is a direct sum of irreducible representations. (This is a generalization of the fact that if $U$ is a unitary operator on~$\Hilbert$, then $\Hilbert$~may not be a direct sum of eigenspaces of~$U$.) However, it will be explained in \cref{DirectIntegralSect} that every unitary representation is a ``direct integral'' of irreducible representations. 
In the special case where $H = \integer$, this is a restatement of the Spectral Theorem for unitary operators \ccf{SpectralThm}.

	\item \label{UnitaryWarns-tame}
	Although the unitary dual $\widehat H$ has a fairly natural topology, it can be quite bad. In particular, the topology may not be Hausdorff. Indeed, in some cases, the topology is so bad that there does not exist an injective, Borel measurable function from~$\widehat H$ to~$\real$. Fortunately, though, the worst problems do not arise for semisimple Lie groups: the unitary dual is always ``tame'' (measurably, at least) in this case.
	
	\end{enumerate}
\end{warns}

\begin{exercises}

%\item \label{RepL2Ex}
%Show that the map~$\pi$ of \cref{RepL2} is indeed a unitary representation of~$H$.
%\hint{It is obvious that $\pi(h)$ is linear, and the invariance of~$\mu$ easily implies that it preserves the inner product on $\LL2(X)$. Also, it is easy to check that $\pi$ is a homomorphism. Continuity follows from the fact that continuous functions are dense in~$\LL2(X)$ (by {Lusin's Theorem} \pref{LusinsThm}).}

\item Suppose $\pi$ is a unitary representation of~$H$ on~$\Hilbert$, and define a map $\xi \colon H \times \Hilbert \to \Hilbert$ by $\xi(h,v) = \pi(h) v$. Show that $\xi$ is continuous.
\hint{Use the fact that $\pi(H)$ consists of unitary operators.}

\item \label{subspace->sumEx}
Prove \cref{subspace->sum}.
\hint{If $\mathcal{K}$ is invariant, then $\mathcal{K}^\perp$ is also invariant. %(Since $(\mathcal{K}^\perp)^\perp = \mathcal{K}$, this also implies the converse.)
Define $T \colon \mathcal{K} \oplus \mathcal{K}^\perp \to \Hilbert$ by $T(\varphi,\psi) = \varphi + \psi$.}

\item \label{SchursLemmaUnitary}
(\thmindex{Schur's Lemma}Schur's Lemma)
Suppose $(\pi,\Hilbert)$ is an irreducible unitary representation of~$H$, and $T$~is a bounded operator on~$\Hilbert$ that commutes with every element of $\pi(\Hilbert)$.
Show there exists $\lambda \in \complex$, such that $T \varphi = \lambda \, \varphi$, for every $\varphi \in \Hilbert$.
\hint{Assume $T$ is normal, by considering $T+ T^*$ and $T - T^*$, and apply the Spectral Theorem \pref{SpectralThm}.}

\item Suppose $(\pi,\Hilbert)$ is an irreducible unitary representation of~$H$, and $\langle~\mid~\rangle'$ is another $H$-invariant inner product on~$\Hilbert$ that defines the same topology on~$\Hilbert$ as the original inner product $\langle~\mid~\rangle$. Show there exists $c \in \real^+$, such that $\langle~\mid~\rangle' = c \langle~\mid~\rangle$.
\hint{Each inner product provides an isomorphism of~$\Hilbert$ with~$\Hilbert^*$. \Cref{SchursLemmaUnitary} implies they are the same, up to a scalar multiple.}

\item \label{RegRepNoInvtMeas}
\optional\ 
In the situation of \cref{RepL2}, a weaker assumption on~$\mu$ suffices to define a unitary representation on $\LL2(X,\mu)$. Namely, instead of assuming that $\mu$ is invariant, it suffices to assume only that the \emph{class} of~$\mu$ is invariant. This means, for every measurable subset~$A$, and all $h \in H$, we have $\mu(A) = 0 \Leftrightarrow \mu(hA) = 0$. Then, for each $h \in H$, the Radon-Nikodym Theorem \pref{RadonNikodym} provides a function $D_h \colon X \to \real^+$, such that $h_* \mu = D_h \mu$. Show that the formula
	$$\bigl( \pi(h) \varphi \bigr) \  = \  \sqrt{D_h(x)} \  \varphi( h^{-1} x ) $$
defines a unitary representation of~$H$ on $\LL2(X,\mu)$.

\end{exercises}

\section{Proof of the Moore Ergodicity Theorem} \label{MooreErgPfSect}

Recall the following result that was proved only in a special case:

\begin{thm}[{(\thmindex{Moore Ergodicity}{Moore Ergodicity Theorem} \pref{MooreErgBasicThm})}] \label{MooreErgBasicThmReprise}
Suppose
	\begin{itemize}
	\item $G$ is connected and simple,
	\item $H$ is a closed, noncompact subgroup of~$G$,
	\item $\Lambda$ is a discrete subgroup of~$G$,
	and
	\item $\phi$ is an $H$-invariant $\LL{p}$-function on~$G/\Lambda$ \textup(with $1 \le p < \infty$\textup).
	\end{itemize}
Then $\phi$ is constant\/ \textup(a.e.\textup).
\end{thm}

This is an easy consequence of the following result in representation theory \csee{MooreErgodicitySimplePf}. 

\begin{thm}[(\thmindex{Decay of matrix coefficients}{Decay of matrix coefficients})] \label{DecayMatCoeffSimple}
 If
 \begin{itemize}
 \item $G$ is simple,
 \item $\pi$ is a unitary representation of~$G$ on a Hilbert
space~$\Hilbert$, such that no nonzero vector is fixed by $\pi(G)$, and
 \item $\{g_j\}$ is a sequence of elements of~$G$, such that $\lVert g_j
\rVert \to \infty$,
 \end{itemize}
 then $\langle \pi(g_j) \phi \mid \psi \rangle \to 0$, for every $\phi,\psi \in \Hilbert$.
 \end{thm}

\begin{proof}
Assume, for simplicity, that 
	$$G = \SL(2,\real) .$$
(A reader familiar with the theory of real roots and Weyl chambers should have little difficulty in extending this proof to the general case; cf.~\cref{DecayMatCoeffSL3}.) Let
	$$ A = \begin{Smallbmatrix} \upast&\\ &\upast \end{Smallbmatrix} \subset G .$$
Further assume, for simplicity, that 
	$$ \{g_j\} \subseteq A .$$
(It is not difficult to eliminate this hypothesis; see \cref{DecayMatCoeffSimpleNotAEx}.)
 By passing to a subsequence, we may assume $\pi(g_j)$ converges weakly, to
some operator~$E$; that is,
 $$ \langle \pi(g_j) \phi \mid \psi \rangle
 \to \langle E \phi \mid \psi \rangle
 \mbox{ \ for every $\phi,\psi \in \Hilbert$} $$
 \csee{SubseqConvWeakly}.
 Let 
 \begin{align} \label{horodefn}
 U &= \{\, v \in G \mid g_j^{-1} v g_j \to e \,\} \\
 \intertext{and} 
 U^- &= \{\, u \in G \mid g_j u g_j^{-1} \to e \,\} 
. \end{align}
 For $u \in U^-$ and $\phi,\psi \in \Hilbert$, we have
 \begin{align*}
 \langle E\pi(u) \phi \mid \psi \rangle
 &= \lim \langle \pi(g_j u) \phi \mid \psi \rangle
 \\&= \lim \langle \pi(g_j u g_j^{-1}) \, \pi(g_j) \phi \mid \psi \rangle
 \\&= \lim \langle \pi(g_j) \phi \mid \psi \rangle
 && \text{\csee{MatCoeffConv}}
 \\&= \langle E \phi \mid \psi \rangle 
 , \end{align*}
 so $E \pi(u) = E$. Therefore, letting 
 	$$\Hilbert^{U^-} = \{\, \phi \in \Hilbert \mid \text{$\pi(u) \phi = \phi$ for all $u \in U^-$} \,\} $$
be the space of $U^-$-invariant vectors in~$\Hilbert$, we have
 	\begin{align} \label{UperpInKer}
	(\Hilbert^{U^-})^\perp \subseteq \ker E
	\end{align}
\csee{UperpInKerEx}.
Similarly, since
 $$ \langle E^* \phi \mid \psi \rangle
 =  \langle \phi \mid E \psi \rangle
 = \lim \langle \phi \mid \pi(g_j) \psi \rangle
 = \lim \langle \pi(g_j^{-1}) \phi \mid  \psi \rangle ,$$
the same argument, with $E^*$ in the place of~$E$ and $g_j^{-1}$ in the
place of~$g_j$, shows that 
	$$(\Hilbert^U)^\perp \subseteq \ker E^* .$$

We also have
\begin{align*}
 \langle \pi(g_j) \phi \mid \pi(g_k) \phi \rangle
 &= \langle \pi(g_k^{-1} g_j) \phi \mid \phi \rangle
 	&& (\text{$\pi(g_k^{-1})$ is unitary})
 \\&= \langle \pi(g_j g_k^{-1}) \phi \mid \phi \rangle
 	&& (\text{$A$ is abelian}) 
 \\&= \langle \pi(g_k^{-1}) \phi \mid \pi(g_j^{-1}) \phi \rangle
.\end{align*}
Applying $\lim_{j \to \infty} \lim_{k \to \infty}$ to both sides yields
$ \lVert E\phi\rVert^2 = \lVert E^* \phi\rVert^2 $,
 and this implies $\ker E = \ker E^* $. 
Hence, 
 \begin{align*}
 \ker E
 &= \ker E + \ker E^*
 \supset (\Hilbert^{U^-})^\perp + (\Hilbert^U)^\perp
 \\&= (\Hilbert^{U^-} \cap \Hilbert^U)^\perp
 = (\Hilbert^{\langle U, U^- \rangle})^\perp
 . \end{align*}
 Now, by passing to a subsequence of $\{g_j\}$, we may assume $\langle U, U^-
\rangle = G$ \csee{horosubseq}. Then $\Hilbert^{\langle U, U^- \rangle} =
\Hilbert^G = 0$,
 so $\ker E
 \supset 0^\perp
 = \Hilbert$.
 This implies that, for all $\phi,\psi \in \Hilbert$, we have
 \begin{align*}
  \lim \langle \pi(g_j) \phi \mid \psi \rangle
 = \langle E \phi \mid \psi \rangle
 = \langle 0 \mid \psi \rangle
 = 0 . & \qedhere \end{align*}
 \end{proof}

\begin{rem}
If $A$ is a bounded operator on a Hilbert space~$\Hilbert$, and $\phi,\psi \in \Hilbert$, then the inner product $\langle A \phi \mid \psi \rangle$ is called a \defit{matrix coefficient} of~$A$. The motivation for this terminology is that if $A \in \Mat_{n \times n}(\real)$, and $\varepsilon_1,\ldots,\varepsilon_n$ is the standard basis of $\Hilbert = \real^n$, then $\langle A \varepsilon_j \mid \varepsilon_i \rangle$ is the $(i,j)$ matrix entry of~$A$. 
\end{rem}

The above argument yields the following more general result.

\begin{cor}[(of proof)] \label{DecayMatCoeffNoCpct}
 Assume
 \begin{itemize}
 \item $G$ has no compact factors,
 \item $\pi$ is a unitary representation of~$G$
on a Hilbert space~$\Hilbert$, and
 \item $\{g_n\} \to \infty$ in $G/N$, for every proper,
normal subgroup~$N$ of~$G$.
 \end{itemize}
 Then $\langle g_n \phi \mid \psi \rangle \to 0$, for
every $\phi,\psi \in (\Hilbert^G)^\perp$.
 \end{cor}

This has the following consequence, which implies the {Moore Ergodicity Theorem} \csee{MautnerPhenom->MooreErgodicityEx}.

\begin{cor}[({\thmindex{Mautner phenomenon}Mautner phenomenon \csee{MautnerPhenomPfEx}})] \label{MautnerPhenom}
 Assume
 \begin{itemize}
 \item $G$ has no compact factors,
 \item $\pi$ is a unitary representation of~$G$ on a Hilbert space~$\Hilbert$, and
 \item $H$ is a closed subgroup of~$G$.
 \end{itemize}
 Then there is a closed, normal subgroup~$N$ of~$G$, containing a cocompact
subgroup of~$H$, such that every $\pi(H)$-invariant vector in~$\Hilbert$ is
$\pi(N)$-invariant.
 \end{cor}
 
 \begin{rem} \label{QuantitativeDecay}
 \Cref{DecayMatCoeffSimple} does not give any information about the \emph{rate} at which the function $\langle \phi(g) \phi \mid \psi \rangle$ tends to~$0$ as $\|g\| \to \infty$. For some applications, it is helpful to know that, for many choices of the vectors $\phi$ and~$\psi$, the inner product decays exponentially fast:

 	If $G$, $\pi$, and~$\Hilbert$ are as in \cref{DecayMatCoeffSimple}, then there is a dense, linear subspace~$\Hilbert_\infty$ of~$\Hilbert$, such that, for all $\phi , \psi \in \Hilbert_\infty$, there exist $a,b > 1$, such that 
		$$ \text{$\displaystyle |\langle \phi(g) \phi \mid \psi \rangle| < \frac{b}{a^{\|g\|}}$ \quad for all $g \in G$} . $$
Specifically, if $K$ is a maximal compact subgroup of~$G$, then we may take
	$$\Hilbert_\infty = \bigset{ \phi \in \Hilbert }{  \text{the linear span of $K \phi$ is finite-dimensional}} \! .$$
(So $\Hilbert_\infty$ is the space of ``\defit[K-finite vector@$K$-finite vector]{\mathversion{bold}$K$-finite}'' vectors.)
 \end{rem}

 \begin{exercises}
 
 \item \label{MooreErgodicitySimplePf}
Show that \cref{MooreErgBasicThmReprise} is a corollary of \cref{DecayMatCoeffSimple}.
 \hint{If $\phi$ is an $H$-invariant function in $\LL{p}(G/\Gamma)$, let $\phi' = |\phi|^{p/2} \in \LL2(G/\Gamma)$. Then $\langle \phi', g_j \phi' \rangle = \langle \phi', \phi' \rangle$ for every $g_j \in H$.}
  
\item \label{SubseqConvWeakly}
 Let $\{T_j\}$ be a sequence of unitary operators on a Hilbert space~$\Hilbert$.
 Show there is a subsequence $\{T_{j_i}\}$ of $\{T_j\}$ and a bounded operator~$E$ on~$\Hilbert$, such that $\langle T_{j_i} v \mid w \rangle \stackrel{i \to \infty}{\longrightarrow} \langle E v \mid w \rangle$ for all $v,w \in \Hilbert$.
\hint{Choose an orthonormal basis $\{e_p\}$. For each $p,q$, the sequence $\{\langle T_j e_p \mid e_q \rangle$ is bounded, and therefore has a subsequence that converges to some~$\alpha_{p,q}$. Cantor diagonalization implies that we may assume, after passing to a subsequence, that $\langle T_j e_p \mid e_q \rangle \to \alpha_{p,q}$ for all $p$ and~$q$.}

\item \label{MatCoeffConv}
Suppose 
	\begin{itemize}
	\item $\pi$ is a unitary representation of~$G$ on~$\Hilbert$,
	\item $\{\phi_j\}$ is a sequence of unit vectors in~$\Hilbert$, 
	and 
	\item $u_j \to e$ in~$G$.
	\end{itemize}
Show $\lim \langle \pi(u_j) \phi_j \mid \psi \rangle = \lim \langle \phi_j \mid \psi \rangle$, for all $\psi \in \Hilbert$.
\hint{Move $\pi(u_j)$ to the other side of the inner product. Then use the continuity of~$\pi$ and the boundedness of $\{\phi_j\}$.}

\item \label{UperpInKerEx}
Prove \pref{UperpInKer}.
\hint{Let $\mathcal{K}$ be the closure of $\{\, \pi(u) \phi - \phi \mid u \in U^-, \phi \in (\Hilbert^{U^-})^\perp \,\}$, and note that $\mathcal{K} \subseteq \ker E$. If $\psi \in \mathcal{K}^\perp$, then $\pi(u) \psi - \psi = 0$ for all $u \in U^-$ (why?), so $\psi \in  \Hilbert^{U^-}$.}

 \item \label{DecayMatCoeffSimpleNotAEx}
 Eliminate the assumption that $\{g_j\} \subseteq A$ from the proof of \cref{DecayMatCoeffSimple}.
 \hint{You may assume the \index{Cartan!decomposition}{Cartan decomposition}, which states that $G = KAK$, where $K$ is compact. Hence, $g_j = c_j a_j c'_j$,
with $c_j,c'_j \in K$ and $a_j \in A$. Assume,
by passing to a subsequence, that $\{c_j\}$ and $\{c'_j\}$ converge. Then
	$$ \lim \langle \pi(g_j) \phi \mid \psi \rangle
	= \lim \bigl\langle \pi(a_j) \bigl( \pi(c') \phi \bigr) 
	 \mathrel{\big|} \pi(c)^{-1}
	\psi \bigr\rangle= 0$$
if $c_j \to c$ and $c'_j \to c'$.}

\item \label{DecayMatCoeffSL3}
Prove \cref{DecayMatCoeffSimple} for the special case where $G = \SL(n,\real)$.

\item \label{horosubseq}
For $G$, $A$, $\{g_j\}$, $U$, and~$U^-$ as in the proof of \cref{DecayMatCoeffSimple} (with $\{g_j\} \subseteq A)$, show that if $\{g_j\}$ is replaced by an appropriate subsequence, then $\langle U , U^- \rangle = G$.
\hint{Arrange that $U$ is $\begin{Smallbmatrix} 1&\upast \\ &1 \end{Smallbmatrix}$ and $U^-$ is $\begin{Smallbmatrix} 1&\\ \upast&1 \end{Smallbmatrix}$, or vice versa.}

\item \label{MautnerPhenomPfEx}
Derive \cref{MautnerPhenom} from \cref{DecayMatCoeffNoCpct}.

\item \label{MautnerPhenom->MooreErgodicityEx}
Derive \cref{MooreErgodicity} from \cref{MautnerPhenom}.
\hint{If $f$~is $H$-invariant, then $\langle \pi(h)f \mid f \rangle = \langle f \mid f \rangle$ for all $h \in H$.}

\item \label{MooreErgNonsimple}
Suppose
	\begin{itemize}
	\item $G$ is connected, with no compact factors,
	\item $\Lambda$ is a discrete subgroup of~$G$,
	\item $H$ is a subgroup of $G$ whose projection to every simple factor of~$G$ is not precompact, 
	and
	\item $\phi$ is an $H$-invariant $\LL{p}$ function on $G/\Lambda$, with $1 \le p < \infty$.
	\end{itemize}
Show that $\phi$ is constant (a.e.).

\item \label{MooreErgLatticeEx}
Suppose
	\begin{itemize}
	\item $H$ is a noncompact, closed subgroup of $G$,
	\item $\Gamma$ is irreducible,
	and
	\item $\phi \colon G/\Gamma \to \real$ is any $H$-invariant, measurable function.
	\end{itemize}
Show that $\phi$ is constant (a.e.).
\hint{There is no harm in assuming that $\phi$ is bounded (why?), so it is in $\LL2(G/\Gamma)$ (why?). Apply \cref{MautnerPhenom}.}

 \end{exercises}

\section{Induced representations} \label{InducedRepSect}

It is obvious that if $H$ is a subgroup of~$G$, then any unitary representation of~$G$ can be restricted to a unitary representation of~$H$. (That is, if we define $\pi|_H$ by $\pi|_H(h) = \pi(h)$ for $h \in H$, then $\pi|_H$ is a unitary representation of~$H$.) What is not so obvious is that, conversely, every unitary representation of~$H$ can be ``induced'' to a unitary representation of~$G$. We will need only the special case where $H = \Gamma$ is a lattice in~$G$ (but see \cref{GenInducedRepEx} for the definition in general). This construction will be a key ingredient of the proof in \cref{LattInTSect} that $\Gamma$ often has Kazhdan's property~$(T)$.

\begin{defn}[{(\term[representation!induced]{Induced representation})}] \label{InducedDefn}
Suppose $\pi$ is a unitary representation of~$\Gamma$ on~$\Hilbert$. 
\begin{enumerate}

\item A measurable function $\varphi \colon G \to \Hilbert$ is said to be (essentially) right \defit[equivariant]{$\Gamma$-equivariant} if, for each $\gamma \in \Gamma$, we have 
	$$ \text{$\varphi(g\gamma^{-1}) = \pi(\gamma) \, \varphi(g)$ for a.e.\ $g \in G$.} $$

\item We use 
	\nindex{$\LGamma{}(G; \Hilbert)$ = $\{$~$\Gamma$-equivariant functions from~$G$ to~$\Hilbert$~$\}$}%
	$\LGamma{}(G; \Hilbert)$ 
	to denote the space of right $\Gamma$-equivariant measurable functions from~$G$ to~$\Hilbert$, where two functions are identified if they agree almost everywhere.

\item For $\varphi \in \LGamma{}(G; \Hilbert)$, we have $\| \varphi(g \gamma) \|_\Hilbert = \| \varphi(g) \|_\Hilbert$ for every $\gamma \in \Gamma$ and a.e.\ $g \in G$ \csee{Norm(gGamma)WD}. Hence, $\|\varphi(g)\|_\Hilbert$ is a well-defined function on $G/\Gamma$ (a.e.), so we may define the $\LL2\,$-norm of $\varphi$ by
	$$ \| \varphi \|_2 = \left( \int_{G/\Gamma} \| \varphi(g) \|_\Hilbert^2 \, d g \right) ^{1/2}.$$

\item We use 
	\nindex{$\LGamma2(G; \Hilbert)$ = $\{$~square-integrable functions in $\LGamma{}(G; \Hilbert)$~$\}$}%
	$\LGamma2(G; \Hilbert)$ to denote the subspace of $\LGamma{}(G; \Hilbert)$ consisting of the functions with finite $\LL2\,$-norm. It is a Hilbert space \csee{L2(G;V)Hilbert}.

\item \label{InducedDefn-formula}
Note that $G$ acts by unitary operators on $\LGamma2(G; \Hilbert)$, via
	\begin{align*}
	\text{$(g \cdot \varphi)(x) = \varphi( g^{-1} x)$ \quad for $g \in G$, $\varphi \in \LGamma2(G; \Hilbert)$, and $x \in G$} 
	\end{align*}
 \csee{GActsL2(G;V)}. This unitary representation of~$G$ is called the representation \defit[representation!induced]{induced} from~$\pi$, and it is denoted 
 	\nindex{$\Ind_\Gamma^G(\pi)$ = induced representation}%
 	$\Ind_\Gamma^G(\pi)$.

 \end{enumerate}
\end{defn}

\begin{exercises}

\item \label{GenInducedRepEx}
\optional\ 
Suppose $(\pi,\Hilbert)$ is a unitary representation of a closed subgroup~$H$ of~$G$. 
Define $\Ind_H^G(\pi)$, without assuming that there is a $G$-invariant measure on $G/H$.
\hint{Since $G/H$ is a $C^\infty$ manifold \csee{HomosAreSmooth}, we may use a nowhere-vanishing differential form to choose a measure~$\mu$ on~$G/H$, such that $f_*\mu$ is in the class of~$\mu$, for every diffeomorphism~$f$ of $G/H$. A unitary representation of~$G$ on $\LL2(G/H, \mu)$ can be defined by using Radon-Nikodym derivatives, as in \cref{RegRepNoInvtMeas}, and the same idea yields a unitary representation on a space of $H$-equivariant functions.}

\item \label{Norm(gGamma)WD}
Let $\varphi \in \LGamma{}(G; \Hilbert)$ and $\gamma \in \Gamma$, where $\pi$ is a unitary representation of~$\Gamma$ on~$\Hilbert$. Show $\| \varphi(g \gamma) \|_\Hilbert = \| \varphi(g) \|_\Hilbert$, for a.e.\ $g \in G$.

\item \label{L2(G;V)Hilbert}
Show $\LGamma2(G;\Hilbert)$ is a Hilbert space (with the given norm, and assuming that two functions represent the same element of the space if and only if they are equal a.e.).

\item \label{GActsL2(G;V)}
Show that the formula in \fullcref{InducedDefn}{formula} defines a unitary representation of~$G$ on $\LGamma2(G;\Hilbert)$.

\item \label{Ind1=L2(G/Gamma)}
Show that $\Ind_\Gamma ^G(\trivrep)$ is (isomorphic to) the usual representation of~$G$ on $\LL2(G/\Gamma)$ (by left translation).

\item \label{IndIrred->Irred}
Show that if $\Ind_\Gamma^G(\pi)$ is irreducible, then $\pi$ is irreducible.

\item Show that the converse of \cref{IndIrred->Irred} is false.
\hint{Is the representation of~$G$ on $\LL2(G/\Gamma)$ irreducible?}

\end{exercises}

\section{Representations of compact groups}

\begin{eg} \label{FourierSeriesEg}
Consider the circle $\real/\integer$. For each $n \in \integer$, define
	$$ \text{$e_n \colon \real/\integer \to \complex$ by $e_n(t) = e^{2\pi i n t}$} .$$
The theory of Fourier Series tells us that $\{e_n\}$ is an orthonormal basis of $\LL2(\real/\integer)$, which means we have the direct-sum decomposition
	$$ \LL2(\real/\integer) = \bigoplus_{n \in \integer} \complex e_n . $$
Furthermore, it is easy to verify that each subspace $\complex e_n$ is an invariant subspace for the regular representation \csee{TorusEigenvector}, and, being $1$-dimensional, is obviously irreducible. Hence, we have a decomposition of the regular representation into a direct sum of irreducible representations. In addition, it is not difficult to see that every irreducible representation of~$\torus$ occurs exactly once in this representation.

More generally, it is not difficult to show that every unitary representation of~$\torus$ is a direct sum of $1$-dimensional representations.
\end{eg}

The following theorem generalizes this to any compact group. However, for nonabelian groups, the irreducible representations cannot all be $1$-dimensional \csee{1D->Abel}.

\begin{thm}[(\thmindex{Peter-Weyl}Peter-Weyl Theorem)] \label{PeterWeyl}
Assume $H$ is \textbf{compact}. Then:
	\begin{enumerate}
	\item \label{PeterWeyl-sum}
	Every unitary representation of~$H$ is\/ \textup(isomorphic to\/\textup) a direct sum of irreducible representations.
	\item \label{PeterWeyl-fd}
	Every irreducible representation of~$H$ is finite-dimensional. 
	\item \label{PeterWeyl-countable}
	$\widehat H$ is countable.
	\item \label{PeterWeyl-regrep}
	For the particular case of the regular representation $\bigl(\regrep, \LL2(H) \bigr)$, we have
		$$\regrep \iso \bigoplus_{(\pi,\Hilbert) \in \widehat H} (\dim \Hilbert) \cdot \pi , $$
	where $k \cdot \pi$ denotes the direct sum $\pi \oplus \cdots \oplus \pi$ of $k$~copies of\/~$\pi$. That is, the ``multiplicity'' of each irreducible representation is equal to its dimension.
	\end{enumerate}
\end{thm}

\begin{proof}
In order to establish both \pref{PeterWeyl-sum} and~\pref{PeterWeyl-fd} simultaneously, it suffices to show that if $(\pi,\Hilbert)$ is any unitary representation of~$H$, then $\Hilbert$ is a direct sum of finite-dimensional, invariant subspaces.
Zorn's Lemma \pref{ZornsLemma} provides a subspace~$\mathcal{M}$ of~$\Hilbert$ that is maximal among those that are a direct sum of finite-dimensional, invariant subspaces. By passing to the orthogonal complement of~$\mathcal{M}$, we may assume that $\Hilbert$ has no nonzero, finite-dimensional, invariant subspaces.

Let 
\noprelistbreak
	\begin{itemize}
	\item $P$ be the orthogonal projection onto some nonzero subspace of~$\Hilbert$ that is finite-dimensional, 
	\item $\mu$ be the Haar measure on~$H$,
	and
	\item $ \overline{P} = \int_H \pi(h) \, P \, \pi(h^{-1}) \, d\mu(h)$.
	\end{itemize}
Note that, since $\overline{P}$ commutes with $\pi(H)$ \csee{AvgOpCommutes}, every eigenspace of~$\overline{P}$ is $H$-invariant \csee{Commutes->EigInvt}. 

Since $P$ is self-adjoint and each $\pi(h)$ is unitary (so $\pi(h^{-1}) = \pi(h)^*$), it is not difficult to see that $\overline{P}$ is self-adjoint. It is also compact \csee{PbarCpct} and nonzero \csee{PbarNonzero}. Therefore, the Spectral Theorem \pref{SpectralThmCpct} implies that $\overline{P}$ has at least one eigenspace~$E$ that is finite-dimensional. By contradicting the fact that $\Hilbert$ has no nonzero, finite-dimensional, invariant subspaces, this completes the proof of \pref{PeterWeyl-sum} and~\pref{PeterWeyl-fd}.

\medbreak

Note that \pref{PeterWeyl-countable} is an immediate consequence of \pref{PeterWeyl-regrep}, since Hilbert spaces are assumed to be separable \csee{HilbertSpaceSeparable}, and therefore cannot have uncountably many terms in a direct sum.

\medbreak

We now give the main idea in the proof of \pref{PeterWeyl-regrep}.
Given an irreducible representation $(\pi, \complex^k)$, we will not calculate the exact multiplicity of~$\pi$, but only indicate how to obtain the correct lower bound by using properties of matrix coefficients. 
Write $\pi(x) = \bigl[ f_{i,j}(x) \bigr]$. Then
	\begin{align}  \label{PeterWeylMatCoeff}
	\begin{split}
	\bigl[ \bigl( \regrep(h)f_{i,j} \bigr)(x) \bigr]
	&= \bigl[ f_{i,j}(h^{-1} x) \bigr]
	= \pi(h^{-1} x)
	\\&= \pi(h^{-1}) \, \pi(x)
	= \pi(h^{-1}) \bigl[ f_{i,j}(x) \bigr] 
	. \end{split} \end{align}
Now, for $1 \le j \le k$, define $ T_j \colon \complex^k \to \LL2(H)$, by
	$$ T_j(a_1,\ldots,a_k) = a_1  f_{1,j} + a_2 f_{2,j} + \cdots + a_k f_{k,j} .$$
Equating the $j$th columns of the two ends of \pref{PeterWeylMatCoeff} tells us that
	$$ T_j \bigl( \pi(h) v \bigr) = \regrep(h) \, T_j(v) ,$$
so $T_j(\complex^k)$ is an invariant subspace, and the corresponding subrepresentation is isomorphic to~$\pi$. Therefore, there are (at least) $k$~different copies of~$\pi$ in~$\regrep$ (one for each value of~$j$). Since, by definition, $k = \dim \pi$, this establishes the correct lower bound for the multiplicity of~$\pi$.
\end{proof}

As an illustrative, simple case of the main results in \cref{RepRnSect,DirectIntegralSect}, we present two different reformulations of the Peter-Weyl Theorem for the special case of abelian groups, after some preliminaries.

\begin{defn}
Let $A$ be an abelian Lie group.
	\begin{enumerate}
	\item A \defit[character (of an abelian group)]{character} of~$A$ is a continuous homomorphism $\chi \colon A \to \torus$, where $\torus = \{\, z \in \complex \mid |z| = 1 \,\}$.
	\item The set of all characters of~$A$ is called the \defit{Pontryagin dual} of~$A$, and is denoted~%
	\nindex{$A^*$ = $\{$characters of abelian group~$A$$\}$}%
	$A^*$.
	 It is an abelian group under the operation of pointwise multiplication. (That is, the product $\chi_1 \chi_2$ is defined by $(\chi_1 \chi_2)(a) = \chi_1(a) \, \chi_2(a)$, for $\chi_1, \chi_2 \in A^*$ and $a \in A$.) Furthermore, if $A/A^\circ$ is finitely generated, then $A^*$ is a Lie group (with the topology of uniform convergence on compact sets).
	\end{enumerate}
\end{defn}

\begin{obs}
If $A$ is any abelian Lie group\/ \textup(compact or not\/\textup), then every irreducible representation $(\pi, \Hilbert)$ of~$A$ is $1$-dimensional \csee{AbelIrred1D}. Therefore, the unitary dual~$\widehat A$ can be identified with the Pontryagin dual~$A^*$ \csee{1Drep<>Char}.
\end{obs}

Hence, for the special case where $H = A$ is abelian, we have the following reformulation of the Peter-Weyl Theorem:

\begin{cor}[\csee{PeterWeylAbelEx}] \label{PeterWeylAbel}
Assume $(\pi,\Hilbert)$ is a unitary representation of a compact, abelian Lie group~$A$.  For each $\chi \in A^*$, let
	\begin{itemize}
	\item $ \Hilbert_\chi = \{\, \varphi \in \Hilbert \mid 
		\text{$\phi(a) \varphi = \chi(a) \, \varphi$, for all $a \in A$} \,\} $,
	and
	\item $P_\chi \colon \Hilbert \to \Hilbert_\chi$ be the orthogonal projection.
	\end{itemize}
Then $\Hilbert = \bigoplus_{\chi \in A^*} \Hilbert_\chi$, so, for all $a \in A$, we have
	$$ \pi(a) = \sum_{\chi \in A^*} \chi(a) \, P_\chi. $$
\end{cor}

Here is another way of saying the same thing:

\begin{cor}[\csee{PeterWeylAbelL2Ex}] \label{PeterWeylAbelL2}
Assume $(\pi,\Hilbert)$ is a unitary representation of a compact, abelian Lie group~$A$.  Then there exist
	\begin{itemize}
	\item a Radon measure $\mu$ on a locally compact metric space~$Y$,
	and
	\item a Borel measurable function $\chi \colon Y \to A^* \colon y \mapsto \chi_y$ \textup(where the countable set~$A^*$ is given the discrete topology\/\textup)
	\end{itemize}
such that $\pi$ is isomorphic to the the unitary representation $\rho_\chi$ of~$A$ on $\LL2(Y, \mu)$ that is defined by
	$$ \bigl( \rho_\chi(a) \varphi  \bigr)(y) = \chi_y(a) \, \varphi(y)
	 \quad \text{for $a \in A$, $\varphi \in \LL2(Y, \mu )$, and $y \in Y$} .$$
%where, for convenience, we write $\chi_y$ for $\chi(y) \in A^*$.
\end{cor}

An analogue of this result for semisimple groups will be stated in \cref{DirectIntegralSect}, after we define the ``direct integral'' of a family of representations.

\begin{exercises}

\item \label{TorusEigenvector} 
In the notation of \cref{FourierSeriesEg}, show $\regrep(h) e_n = e^{-2\pi i h} e_n$, for all $h \in \real/\integer$.

\item \label{MatCoeffOfCpct}
Suppose $(\pi,\Hilbert)$ is a unitary representation of a compact group~$H$, let $\varphi,\psi \in \Hilbert$, and define $f \colon H \to \complex$ by $f(h) = \langle \pi(h) \varphi \mid \psi \rangle$. Show $f \in \LL2(H)$.
\hint{It is a bounded function on a compact set.}

\item \label{AvgOpCommutes}
Suppose $(\pi,\Hilbert)$ is a unitary representation of a compact group~$H$. Show that if $T$ is any bounded operator on~$\Hilbert$, then 
	$$\overline{T} = \int_H \pi(h) \, T \, \pi(h^{-1}) \, d\mu(h)$$
is an operator that commutes with every element of $\pi(H)$.
\hint{The invariance of Haar measure implies $\pi(g) \overline{T} \pi(g^{-1}) = \overline{T}$.}

\item \label{PbarCpct}
Show that the operator $\overline{P}$ in the proof of \cref{PeterWeyl} is compact.
\hint{Apply \cref{CpctOpBasics}, by noting that any integral can be approximated by a finite sum, and the finite sum is an operator whose range is finite-dimensional.}

\item \label{PbarNonzero}
Show that the operator $\overline{P}$ in the proof of \cref{PeterWeyl} is nonzero.
\hint{Choose some nonzero $\varphi \in \Hilbert$, such that $P \varphi = \varphi$.
Then $\langle \overline{P} \varphi \mid \varphi \rangle > 0$, since $\langle P \psi \mid \psi \rangle \ge 0$ for all $\psi \in \Hilbert$.}

\item \label{Commutes->EigInvt}
Suppose $(\pi,\Hilbert)$ is a unitary representation of a Lie group~$H$, $T$~is a bounded operator on~$\Hilbert$, $\lambda \in \complex$, and $\varphi \in \Hilbert$. Show that if $T$ commutes with every element of $\pi(H)$, and $T(\varphi) = \lambda \varphi$, then $T \bigl( \pi(h) \varphi \bigr) = \lambda \, \pi(h) \varphi$, for every $h \in H$.

\item Assume $H$ is compact. Show that $H$ is finite if and only if it has only finitely many different irreducible unitary representations (up to isomorphism).
\hint{You may assume \cref{PeterWeyl}.}

\item \label{AbelIrred1D}
Show that every irreducible representation of any abelian Lie group is $1$-dimensional.
\hint{If $\pi(a)$ is not a scalar, for some $a$, then the Spectral Theorem \pref{SpectralThm} yields an invariant subspace.}

\item \label{1Drep<>Char}
Let $H$ be a Lie group. Show there is a natural bijection between the set of $1$-dimensional unitary representations (modulo isomorphism) and the set of continuous homomorphisms from~$H$ to~$\torus$.
\hint{Any $1$-dimensional unitary representation is isomorphic to a representation on~$\complex$.}

\item \label{PeterWeylAbelEx}
Derive \cref{PeterWeylAbel} from \cref{PeterWeyl}.

\item \label{PeterWeylAbelL2Ex}
Prove \cref{PeterWeylAbelL2}.
\hint{If $\Hilbert_\chi \neq \{0\}$, then $\Hilbert_\chi$ is isomorphic to some $\LL2(Y_\chi, \mu_\chi)$. Let $Y = \bigcup_\chi Y_\chi$.}

\end{exercises}

\section{Unitary representations of \texorpdfstring{$\real^{\lowercase{n}}$}{Rn}} \label{RepRnSect}

Any character~$\chi$ of $\real^n$ is of the form 
	$$ \chi(a) = e^{2 \pi i \, (a \cdot y)} \qquad
	\text{for some (unique) $y \in \real^n$} $$
\csee{RnChars}.  Therefore, the Pontryagin dual~$(\real^n)^*$ (or, equivalently, the unitary dual~$\widehat{\real^n}$) can be identified with~$\real^n$ (by matching $\chi$ with the corresponding vector~$y$). In particular, unlike in \cref{PeterWeyl}, the unitary dual is uncountable.

Unfortunately, however, not every representation of~$\real^n$ is a direct sum of irreducibles.
For example, the regular representation~$\regrep$ of $\real^n$ on $\LL2(\real^n)$ has no $1$-dimensional, invariant subspaces \csee{L2RnNo1D}, so it does not even contain a single irreducible representation and is therefore not a sum of them. Indeed, \term{Fourier Analysis} tells us that a function in $\LL2(\real^n)$ is not a \emph{sum} of exponentials, but an \emph{integral}:
	$$ \varphi(a) = \int_{\real^n} \widehat\varphi(y) \, e^{2\pi i \, (a \cdot y)} \, dy ,$$
where $\widehat\varphi$ is the Fourier transform of~$\varphi$.
Now, for each Borel subset~$E$ of~$\real^n$, let 
	\begin{align} \label{HEFourier}
	\Hilbert_E = \{\, f \in \LL2(\real^n) \mid \text{$\widehat f(y) = 0$ for a.e.\ $y \notin E$} \,\} 
	. \end{align}
Then it is not difficult to show that $\Hilbert_E$ is a closed, $\regrep$-invariant subspace \csee{HEinvt}.

Now, let $P(E)\colon \Hilbert \to \Hilbert_E$ be the orthogonal projection. Then we can think of $P$ as a projection-valued measure on $\real^n$ (or on $(\real^n)^*$), and, for all $a \in \real^n$,  we have
	$$\regrep(a) = \int_{\real^n} e^{i (a \cdot y)} \, dP(y)
		= \int_{(\real^n)^*} \chi(a) \, dE(\chi) .$$
If we let $\pi = \regrep$, this is a perfect analogue of the conclusion of \cref{PeterWeylAbel}, but with the sum replaced by an integral.

A version of the \thmindex{Spectral}{Spectral Theorem} tells us that this generalizes in a natural way to all unitary representations of~$\real^n$, or, in fact, of any abelian Lie group:

\begin{prop} \label{ProjValMeas}
Suppose $\pi$ is a unitary representation of an abelian Lie group~$A$ on~$\Hilbert$. Then there is a \textup(unique\/\textup) projection-valued measure~$P$ on\/~$A^*$, such that
	$$ \text{$\pi(a) = \int_{A^*} \chi(a) \, dP(\chi)$ \ for all $a \in A$} . $$
\end{prop}

This can be reformulated as a generalization of \cref{PeterWeylAbelL2}:

\begin{cor} \label{IntIrredAbel}
Let $(\pi,\Hilbert)$ be a unitary representation of an abelian Lie group~$A$.  Then there exist
	\begin{itemize}
	\item a probability measure $\mu$ on a locally compact metric space~$Y$,
	and
	\item a Borel measurable function $\chi \colon Y \to A^* \colon y \mapsto \chi_y$,
	\end{itemize}
such that $\pi$ is isomorphic to the unitary representation $\rho_\chi$ of~$A$ on $\LL2(Y, \mu)$ that is defined by
	$$ \bigl( \rho_\chi(a) \varphi  \bigr)(y) = \chi_y(a) \, \varphi(y)
	 \quad \text{for $a \in A$, $\varphi \in \LL2(Y, \mu )$, and $y \in Y$} .$$
\end{cor}

\begin{proof}
Let $P$ be the projection-valued measure given by \cref{ProjValMeas}.
A closed subspace~$\Hilbert'$ of~$\Hilbert$ is said to be \defit[cyclic, for a projection-valued measure]{cyclic} for~$P$ if there exists $\psi \in \Hilbert'$, such that the span of $\{\, P(E)\,\psi \mid E \subset A^* \,\}$ is a dense subspace of~$\Hilbert'$. It is not difficult to see that $\Hilbert$ is an orthogonal direct sum of countably many cyclic subspaces \csee{IntIrredCyclicEx}. Therefore, we may assume $\Hilbert$ is cyclic \csee{AssumeHCyclicEx} (and nonzero). So we may fix some unit vector $\psi$ that generates a dense subspace of~$\Hilbert$.

Define a probability measure $\mu$ on $A^*$ by 
	$$ \mu(E) = \langle P(E)\psi \mid \psi \rangle = \langle P(E)\psi \mid P(E)\psi \rangle ,$$
and let $\Id$ be the identity map on~$A^*$.

For the characteristic function~$f_E$ of each Borel subset~$E$ of~$A^*$, define $\Phi(f_E) = P(E) \psi$. Then $\langle \Phi(f_{E_1}) \mid  \Phi(f_{E_2}) \rangle =  \langle f_{E_1} \mid  f_{E_2} \rangle$, by the definition of~$\mu$, so $\Phi$ extends to a norm-preserving linear map~$\Phi'$ from $\LL2(A^*,\mu)$ to~$\Hilbert$. Since $\psi$ is a cyclic vector for~$\Hilbert$, we see that $\Phi'$ is surjective, so it is an isomorphism of Hilbert spaces. Indeed, $\Phi'$ is an isomorphism from $\rho_{\Id}$ to~$\pi$ \csee{IsoRhoPi}.
\end{proof}

\begin{exercises}

\item \label{RnChars}
Show that every character~$\chi$ of $\real^n$ is of the form 
	$ \chi(a) = e^{2 \pi i \, (a \cdot t)}$ for some $t \in \real^n$.
\hint{Since $\real^n$ is simply connected, any continuous homomorphism from~$\real^n$ to~$\torus$ can be lifted to a homomorphism into the universal cover, which is~$\real$. Apply \cref{RnHomIsLinear}.}

\item \label{L2RnNo1D}
Let $H$ be a noncompact Lie group. Show that the regular representation of~$H$ has no $1$-dimensional, invariant subspaces.
\hint{If $\varphi$ is in a $1$-dimensional, invariant subspace of $\LL2(H)$, then $|\varphi|$~is constant (a.e.).}

\item \label{HEinvt}
For every measurable subset~$E$ of~$\real^n$, show that the subspace $\Hilbert_E$ defined in \pref{HEFourier} is closed, and is invariant under $\regrep(\real^n)$.
\hint{It is clear that $\{\, \widehat f \mid f \in \Hilbert_E \,\}$ is closed. For invariance, note that the Fourier transform of $\regrep(a)f$ is $e^{-2\pi i (a \cdot y)} \widehat f(y)$.}

\item \label{IntIrredCyclicEx}
Given a projection-valued measure~$P$ on a Hilbert space~$\Hilbert$, show that $\Hilbert$ is the orthogonal direct sum of countably many cyclic subspaces.
\hint{Every vector in~$\Hilbert$ is contained in a cyclic subspace, the orthogonal complement of a $P(E)$-invariant subspace is $P(E)$-invariant, and all Hilbert spaces are assumed to be separable.} 

\item \label{AssumeHCyclicEx}
In the notation of \cref{IntIrredAbel}, suppose $\rho_{\chi_i}$ is the representation on $\LL2(Y_i, \mu_i)$ corresponding to some $\chi_i \colon Y_i \to A^*$. Show $\bigoplus_{i=1}^\infty \rho_{\chi_i} \iso \rho_\chi$, for some $Y$, $\mu$, and~$\chi$.
\hint{Let $(Y,\mu)$ be the disjoint union of (copies of) $(Y_i,\mu_i)$.}

\item \label{IsoRhoPi}
In the notation of the proof of \cref{IntIrredAbel}, show that $\Phi'$ is an isomorphism from $\rho_{\Id}$ to~$\pi$.
\hint{Given $a \in A$ and $E \subset A^*$, write $E$ as the disjoint union of small sets $E_1,\ldots,E_n$ (so $\chi \mapsto \chi(a)$ is almost constant on each~$E_i$). 
Then
	$$ \textstyle \Phi' \bigl( \rho_{\Id}(a) f_E \bigr)
	\approx \Phi' \bigl( \sum_i \chi_i(a) \, f_{E_i} \bigr)
	= \sum_i \chi_i(a) \, P(E_i) \, \psi
	\approx \pi \bigl( P(E) \psi \bigr)
	= \pi \bigl( \Phi'(f_E) \bigr)
	, $$
for any $\chi_i \in E_i$.}

\end{exercises}

\section{Direct integrals of representations} \label{DirectIntegralSect}

Before we define the direct integral of a collection of unitary representations, we first discuss the simpler case of a direct sum of a sequence $\{(\pi_n, \Hilbert_n)\}_{n= 1}^\infty$ of unitary representations.

\begin{defn} \label{HilbertDirSumInfty}
If $\{\Hilbert_n \}_{n=1}^\infty$ is a sequence of Hilbert spaces, then the \defit[direct sum!of Hilbert spaces]{direct sum} $\bigoplus_{n= 1}^\infty \Hilbert_n$ consists of all sequences $\{\varphi_n\}_{n= 1}^\infty$, such that 
	\begin{itemize}
	\item $\varphi_n \in \Hilbert_n$ for each~$n$, 
	and
	\item $\sum_{n= 1}^\infty \| \varphi_n\|^2 < \infty$.
	\end{itemize}
This is a Hilbert space, under the inner product 
	$$ \bigl\langle \{\varphi_n\} \mid \{\psi_n\} \bigr\rangle = \sum_{n= 1}^\infty \langle \varphi_n \mid \psi_n \rangle .$$
It contains a copy of $\Hilbert_n$, for each~$n$, such that $\Hilbert_m \perp \Hilbert_n$, for $m \neq n$.
\end{defn}

Suppose, now, that all of the Hilbert spaces in the sequence are the same; say, $\Hilbert_n = \Hilbert$, for all~$n$. Then $\bigoplus_{n=1}^\infty \Hilbert_n$ is equal to the set of square-integrable functions from $\integer^+$ to $\Hilbert$, which can be denoted $\LL2(\integer^+ ; \Hilbert)$. In this notation, the direct sum of unitary representations is easy to describe:

\begin{defn}
If $\{\, \pi_n\}_{n= 1}^\infty$ is a sequence of unitary representations of~$H$ on a fixed Hilbert space~$\Hilbert$, then $\bigoplus_{n= 1}^\infty \pi_n$ is the unitary representation $\pi$ on $\LL2(\integer^+ ; \Hilbert)$ that is defined by
	$$ \bigl( \pi(h) \varphi  \bigr)(n) = \pi_n(h) \, \varphi(n)
	 \quad \text{for $h \in H$, $\varphi \in \LL2(\integer^+ ; \Hilbert)$, and $n \in \integer^+$} .$$
\end{defn}

This description of the direct sum naturally generalizes to a definition of the direct integral of representations:

\begin{defn} \label{DirectIntegDefn}
Suppose 
	\begin{itemize}
	\item $\Hilbert$ is a Hilbert space,
	\item $\{\pi_x\}_{x \in X}$ is a measurable family of unitary representations of~$H$ on~$\Hilbert$, which means:
		\begin{itemize}
		\item $X$~is a locally compact metric space,
		\item $\pi_x$ is a unitary representation of~$H$ on~$\Hilbert$, for each $x \in X$,
		and
		\item for each fixed $\varphi,\pi \in \Hilbert$, the expression $\langle \pi_x(h) \varphi \mid \psi \rangle$ is a Borel measurable function on $X \times H$,
		\end{itemize}
	and
	\item $\mu$ is a Radon measure on~$X$.
	\end{itemize}
Then $\int_X \pi_x \, d\mu(x)$ is the unitary representation~$\pi$ of~$H$ on $\LL2(X, \mu ; \Hilbert)$ that is defined by
	$$ \bigl( \pi(h) \varphi  \bigr)(x) = \pi_x(h) \, \varphi(x)
	 \quad \text{for $h \in H$, $\varphi \in \LL2(X, \mu ; \Hilbert)$, and $x \in X$} .$$
This is called the \defit{direct integral} of the family of representations $\{\pi_x\}$.
\end{defn}

The above definition is limited by requiring all of the representations to be on the same Hilbert space. The construction can be generalized to eliminate this assumption \csee{GeneralDirectIntegDefn}, but there is often no need: 

\begin{thm} \label{UnitaryGDirectInt}
Assume 	
	\begin{itemize}
	\item $\pi$ is a unitary representation of~$G$, 
	\item $G$ is connected, and has no compact factors,
	and
	\item no nonzero vector is fixed by every element of $\pi(G)$.
	\end{itemize}
Then there exist $\Hilbert$, $\{\pi_x\}_{x \in X}$, and~$\mu$, as in \cref{DirectIntegDefn}, such that 
	\begin{enumerate}
	\item $\pi \iso \int_X \pi_x \, d\mu(x)$, 
	and
	\item $\pi_x$ is irreducible for every $x \in X$.
	\end{enumerate}
\end{thm}

\begin{rem} \label{GeneralDirectIntegDefn}
Up to isomorphism, there are only countably many different Hilbert spaces (since any two Hilbert spaces of the same dimension are isomorphic). It is therefore not difficult to generalize \cref{DirectIntegDefn} to deal with a family of representations in which the Hilbert space varies with~$x$. Such a generalization allows every unitary representation of any Lie group to be written as a direct integral of representations that are irreducible.

Here is one way.
Let us say that
$\{(\pi_x,\Hilbert_x\}_{x \in X}$ is a measurable family of unitary representations of~$H$ if:
	\begin{itemize}
	\item $X = \bigcup_{n=1}^\infty$ is the union of countably many locally compact metric spaces~$X_n$,
	\item for each~$n$, there is Hilbert space~$\Hilbert_n$, such that $\Hilbert_x = \Hilbert_n$ for all $x \in X_n$
	\item $\pi_x$ is a unitary representation of~$H$ on~$\Hilbert_x$, for each $x \in X$,
	\item for each~$n$ and each $\varphi,\pi \in \Hilbert_n$, the expression $\langle \pi_x(h) \varphi \mid \psi \rangle$ is a Borel measurable function on $X_n \times H$,
	and
	\item $\mu$ is a Radon measure on~$X$.
	\end{itemize}
Given such a family of representations, we define
	$$ \int_X \pi_x \, d\mu(x) = \bigoplus_{n=1}^\infty \int_{X_n} \pi_x \, d\mu(x) .$$
With this, more general, notion of direct integral, it can be proved that every unitary representation of any Lie group is isomorphic to a direct integral of a measurable family of irreducible unitary representations.
\end{rem}

\begin{exercises}

\item \label{1D->Abel}
Let $H$ be a Lie group. Show that if the regular representation of~$H$ is a direct integral of $1$-dimensional representations, then $H$ is abelian.

\end{exercises}

\begin{notes}

There are many books on the theory of unitary representations, including the classics of Mackey \cite{Mackey-ThyUnitaryGrpReps,Mackey-UnitaryRepsPhysProbNT}. Several books, such as \cite{Knapp-RepThySSGrps}, specifically focus on the representations of semisimple Lie groups.

The Moore Ergodicity Theorem \pref{MooreErgBasicThmReprise} is due to C.\,C.\,Moore \cite{Moore-ergodicity}. 

\Cref{DecayMatCoeffNoCpct} is due to R.\,Howe and C.\,C.\,Moore
\cite[Thm.~5.1]{HoweMoore} and (independently) R.\,J.\,Zimmer
\cite[Thm.~5.2]{Zimmer-orbitspace}. 
%It appears in \cite[Thm.~2.2.20, p.~23]{ZimmerBook}.
The elementary proof we give here was found by
R.\,Ellis and M.\,Nerurkar \cite{EllisNerurkar}.
%(See \cite{Pittet} for a more complete exposition.) -- preprint hasn't been published? @@@
Other proofs are in \cite[\S2.3, pp.~85--92]{MargulisBook} and \cite[\S2.4, pp.~28--31]{ZimmerBook}.

A more precise form of the quantitative estimate in \cref{QuantitativeDecay} can be found in \cite[Cor.~7.2]{Howe-rank}.
(As stated there, the result requires the matrix coefficient $\langle \pi(g)\phi \mid \psi \rangle$ to be an $\LL{p}$~function of~$g$, for $\phi,\psi$ in a dense subspace of~$\Hilbert$, and for some $p < \infty$, but it was proved in \cite[Thm.~2.4.2]{Cowling-coeffs} that this integrability hypothesis always holds.)

\Cref{PeterWeyl} is proved in \cite[Chap.~3]{Sepanski-CpctLieGrps}

See \cite[Chap.~2]{Helson-SpectralThm} for a nice proof of \cref{ProjValMeas}. (Although most of the proof is written for $n = 1$, it is mentioned on p.~31 that the argument works in general.) 

See \cite[Thm.~2.9, p.~108]{Mackey-ThyUnitaryGrpReps} for a proof of \cref{GeneralDirectIntegDefn}'s statement that every unitary representation is a direct integral of irreducibles. (This is a generalization of \cref{UnitaryGDirectInt}.)

Regarding \fullcref{UnitaryWarns}{tame}, groups for which the set of irreducible unitary representations admits an injective Borel map to $[0,1]$ are called ``Type~I'' (and the others are ``Type~II''). See \cite[\S2.3, pp.~77--85]{Mackey-ThyUnitaryGrpReps} for some discussion of this.

\end{notes}

  \immediate\addtocontents{toc}{\protect\toceject} % @@@
 %!TEX root = IntroArithGrps.tex

\mychapter{Amenable Groups}
\label{AmenableChap}

\prereqs{none.}
	%, but unitary representations (\cref{UnitaryRepSect}) and quasi-isometries (\cref{QuasiChap}) are mentioned.}

The classical {Kakutani-Markov Fixed Point Theorem} \pref{AbelAmen} implies that any abelian group of continuous linear operators has a fixed point in any compact, convex, invariant set. This theorem can be extended to some non-abelian groups; the groups that satisfy such a fixed-point property are said to be ``amenable\zz,'' and they have quite a number of interesting features. 
Many important subgroups of~$G$ are amenable, so the theory is directly relevant to the study of arithmetic groups, even though we will see that $G$ and~$\Gamma$ are usually not amenable. In particular, the theory yields an important equivariant map that will be constructed in \cref{AmenEquiMapSect}.

\section{Definition of amenability}

\begin{assump}
Throughout this chapter, $H$ denotes a Lie group.
The ideas here are important even in the special case where $H$ is discrete.
\end{assump}

\begin{defn} \label{CpctCnvxDefn}
Suppose $H$ acts continuously (by linear maps) on a locally convex topological vector space~$\LocConvex$.
Every $H$-invariant, compact, convex subset of~$\LocConvex$ is called a \defit[compact!convex $H$-space]{compact, convex $H$-space}.
\end{defn}

\begin{defn} \label{AmenDefn}
$H$ is \defit[amenable group]{amenable} if and only if $H$ has a fixed point in every nonempty, compact, convex $H$-space.
\end{defn}

This is just one of many different equivalent definitions of amenability. (A few others are discussed in \cref{OtherAmen}.) The equivalence of these diverse definitions is a testament to the fact that this notion is very fundamental. 

\begin{rems} \ 
\noprelistbreak
	\begin{enumerate}

	\item All locally convex topological vector spaces are assumed to be Hausdorff.

	\item In most applications, the locally convex space~$\LocConvex$ is the dual of a separable Banach space, with the weak$^*$ topology \csee{WeakStarDefn}.  In this situation, every compact, convex subset~$C$ is second countable, and is therefore metrizable \csee{Metrizable<>2ndCount}.
	With these thoughts in mind, we feel free to assume metrizability when it eliminates technical difficulties in our proofs.
In fact, we could restrict to these spaces in our definition of amenability, because it turns out that this modified definition results in exactly the same class of groups (if we only consider groups that are second countable) \csee{MetrizableSuffices}.

\item The choice of the term ``\term[amenable group]{amenable}'' seems to have been motivated by two considerations:%
\noprelistbreak
		\begin{enumerate}
		\item The word ``amenable'' can be pronounced ``a-MEAN-able\zz,'' and we will see in \cref{OtherAmen} that a group is amenable if and only if it admits certain types of means.
		\item One definition of ``amenable'' from the \emph{Oxford American Dictionary} is ``capable of being acted on a particular way\zz.'' In other words, in colloquial English, something is ``amenable'' if it is easy to work with. Classical analysis has averaging theorems and other techniques that were developed for the study of functions on the group~$\real^n$. Many of these methods can be generalized to all amenable groups, so amenable groups are easy to work with.
		\end{enumerate}
	\end{enumerate}
\end{rems}

\begin{exercises}

\item Show that every finite group is amenable.
\hint{For some $c_0 \in C$, let 
	$c = \frac{1}{\#H} \sum_{h\in H} hc_0 $.
Then $c \in C$ and $c$ is fixed by~$H$.}

\item \label{QuotAmen}
Show that quotients of amenable groups are amenable. That is, if $H$ is amenable, and $N$ is any closed, normal subgroup of~$H$, then $H/N$ is amenable.

\item Suppose $H_1$ is amenable, and there is a continuous homomorphism $\varphi \colon H_1 \to H$ with dense image. Show $H$ is amenable.

\end{exercises}

\section{Examples of amenable groups}

In this section, we will see that:
\noprelistbreak
\begin{itemize}
\item abelian groups are amenable \see{AbelAmen},
\item compact groups are amenable \see{CpctAmen},
\item solvable groups are amenable \see{AmenExtCor}, because the class of amenable groups is closed under extensions \see{AmenExtension},  
and
\item closed subgroups of amenable groups are amenable \see{SubgrpAmen}.
\end{itemize}
On the other hand, however, it is important to realize that not all groups are amenable. In particular, we will see in \cref{NonamenSect} that:
	\noprelistbreak
	\begin{itemize}
	\item nonabelian free groups are not amenable, 
	and
	\item $\SL(2,\real)$ is not amenable.
	\end{itemize}

\bigbreak

We begin by showing that $\integer$ is amenable:

\begin{prop} \label{CyclicAmen}
Cyclic groups are amenable.
\end{prop}

\begin{proof}
Assume $H = \langle T \rangle$ is cyclic. Given a nonempty, compact, convex $H$-space~$C$, choose some $c_0 \in C$.  For $n \in \natural$, let
	\begin{equation} \label{AmenInvt-c_n}
	c_n  =  \frac{1}{n+1}  \sum_{k = 0}^{n} T^k(c) 
	. \end{equation}
Since $C$ is compact, the sequence $\{c_n\}$ must have an accumulation point $c \in C$. It is not difficult to see that $c$ is fixed by~$T$ \csee{lim(c_n)invt}. Since $T$ generates~$H$, this means that $c$ is a fixed point for~$H$.
\end{proof}

\begin{cor}[(\thmindex{Kakutani-Markov Fixed Point}Kakutani-Markov Fixed Point Theorem)] \label{AbelAmen}
Every abelian group is amenable.
\end{cor}

\begin{proof}
Let us assume $H = \langle g,h \rangle$ is a $2$-generated abelian group. (See \cref{AbelAmenEx} for the general case.) If $C$ is any nonempty, compact, convex $H$-space, then \cref{CyclicAmen} implies that the set $C^g$ of fixed points of~$g$ is nonempty. It is easy to see that $C^g$ is compact and convex \csee{FixedPtsCpctConv}, and, because $H$ is abelian, that $C^g$ is invariant under~$h$ \csee{FixedPtsInvCent}. Hence, $C^g$ is a nonempty, compact, convex $\langle h \rangle$-space. Therefore, \cref{CyclicAmen} implies that $h$ has a fixed point~$c$ in $C^g$. Now $c$ is fixed by~$g$ (because it belongs to~$C^g$), and $c$ is fixed by~$h$ (by definition), so $c$ is fixed by $\langle g,h \rangle = H$.
\end{proof}

Compact groups are also easy to work with:

\begin{prop} \label{CpctAmen}
Compact groups are amenable.
\end{prop}

\begin{proof}
Assume $H$ is compact, and let $\mu$ be a Haar measure on~$H$. 
Given a nonempty, compact, convex $H$-space~$C$, choose some $c_0 \in C$.
Since $\mu$ is a probability measure, we may let 
	\begin{equation} \label{CenterOfHOrbit}
	c = \int_H h(c_0) \, d\mu(h) \in C
	. \end{equation}
(In other words, $c$ is the center of mass of the $H$-orbit of~$c_0$.)
The $H$-invariance of~$\mu$ implies that $c$ is a fixed point for~$H$ \csee{InvtMeas->FP}.
\end{proof}

It is easy to show that amenable extensions of amenable groups are amenable \csee{AmenExtensionEx}:

\begin{prop} \label{AmenExtension}
 If $H$ has a closed, normal subgroup~$N$, such that $N$ and $H/N$ are amenable, then $H$ is amenable.
\end{prop}

Combining the above results has the following consequences:

\begin{cor} \label{AmenExtCor} \ 
\noprelistbreak
\begin{enumerate}
\item \label{AmenExtCor-Solvable} 
 Every solvable group is amenable.

\item \label{AmenExtCor-CocpctSolv}
If $H$ has a solvable, normal subgroup~$N$, such that $H/N$ is compact, then $H$ is amenable.

\end{enumerate}
\end{cor}

\begin{proof}
\Cref{SolvAmen,CocpctSolv->Amen}.
\end{proof}

The converse of \fullcref{AmenExtCor}{CocpctSolv} is true for connected groups \csee{ConnAmen}.

\begin{prop} \label{SubgrpAmen}
Every closed subgroup of any amenable group is amenable.
\end{prop}

\begin{proof}
This proof employs a bit of machinery, so we postpone it to \cref{AmenSubgrpSect}. 
(For discrete groups, the result follows easily from some other characterizations of amenability; see \cref{FinAddRem,FolnerSubgrp} below.)
\end{proof}

\begin{exercises}

\item \label{lim(c_n)invt}
Suppose $T$ is a continuous linear map on a locally convex space~$\LocConvex$. Show that if $c$ is any accumulation point of the sequence $\{c_n\}$ defined by \pref{AmenInvt-c_n}, then $c$ is $T$-invariant.
\hint{If $\| c_n - c \|$ is small, then $\|T(c_n) - T(c)\|$ is small. Show that $\|T(c_n) - c_n\|$ is small whenever $n$~is large. Conclude that $\|T(c) - c \|$ is smaller than every~$\epsilon$.}

\item \label{FixedPtsCpctConv}
Suppose $C$ is a compact, convex $H$-space. Show that the set $C^H$ of fixed points of~$H$ is compact and convex.
\hint{Closed subsets of~$C$ are compact.}

\item \label{FixedPtsInvCent}
Suppose $H$ acts on a space~$C$, $A$~is a subgroup of~$H$, and $h$~is an element of the centralizer of~$A$. Show that the set $C^A$ of fixed points of~$A$ is invariant under~$h$.

\item \label{FixedPtsInvNorm}
Establish \cref{FixedPtsInvCent} under the weaker assumption that $h$ is an element of the \emph{normalizer} of~$A$, not the centralizer.

\item \label{AbelAmenEx}
Prove \cref{AbelAmen}.
\hint{For each $h \in H$, let $C^h$ be the set of fixed points of~$h$. The given argument implies (by induction) that $\{\, C^h \mid h \in H \,\}$ has the finite intersection property, so the intersection of these fixed-point sets is nonempty.}

\item \label{InvtMeas->FP}
Show that if $\mu$ is the Haar measure on~$H$, and $H$ is compact, then the point~$c$ defined in \pref{CenterOfHOrbit} is fixed by~$H$.

\item \label{AmenExtensionEx}
Prove \cref{AmenExtension}.
\hint{\Cref{FixedPtsCpctConv,FixedPtsInvNorm}.}

\item Show that $H_1 \times H_2$ is amenable if and only if $H_1$ and~$H_2$ are both amenable.

\item \label{SolvAmen}
Prove \fullcref{AmenExtCor}{Solvable}.
\hint{\Cref{AmenExtension}.}

\item \label{CocpctSolv->Amen}
Prove \fullcref{AmenExtCor}{CocpctSolv}.
\hint{\Cref{AmenExtension}.}

\item Suppose $H$ is discrete, and $H_1$ is a finite-index subgroup. Show $H$ is amenable if and only if $H_1$ is amenable.

\item Show that if $\Lambda$ is a lattice in~$H$, and $\Lambda$ is amenable, then $H$ is amenable.
\hint{Let $\mu = \int_{H/\Lambda} hv \, d h$, where $v$ is a fixed point for~$\Lambda$.}

\item Assume $H$ is discrete. Show that if every finitely generated subgroup of~$H$ is amenable, then $H$ is amenable.
\hint{For each $h \in H$, let $C^h$ be the set of fixed points of~$h$. Then $\{\, C^h \mid h \in H \,\}$ has the finite intersection property, so $\bigcap_h C^h \neq \emptyset$.}

%\item Is every finite group amenable?

%\item Is $\real$ amenable.

\item \label{BorelSL3Amen}
Let \ 
	$P = \begin{Smallbmatrix} \upast&& \\ \upast&\upast& \\ \upast&\upast&\upast \end{Smallbmatrix} \subset \SL(3,\real) $. \ 
Show that $P$ is amenable.
\hint{$P$ is solvable.}

\item Assume there exists a discrete group that is not amenable. Show the free group~$\free_2$ on $2$ generators is not amenable.
\hint{$\free_n$ is a subgroup of~$\free_2$.}

\item Assume there exists a Lie group that is not amenable. 
\noprelistbreak
	\begin{enumerate}
	\item Show the free group~$\free_2$ on $2$ generators is not amenable.
	\item Show $\SL(2,\real)$ is not amenable.
	\end{enumerate}
%\hint{You may use \cref{SubgrpAmen}.}

%\item \label{FPNonSep}
%Suppose $H$ acts continuously (in norm) on a Banach space~$\Banach$, and $C$ is a weak$^*$ compact, compact, convex, $H$-invariant subset of the dual space~$\Banach^*$. 
%Show that if $H$ is amenable, then $H$ has a fixed point in~$C$, even though $\Banach^*$ may not be Fréchet (because $\Banach$ may not be separable).
%\hint{For each $H$-invariant, separable subspace~$\LocConvex$ of~$\Banach$, let $C_{\LocConvex}^H \subseteq \Banach^*$ consist of the elements of~$C$ whose restriction to~$\LocConvex$ is $H$-invariant. Amenability implies that $\{ C_{\LocConvex}^H \}$ has the finite-intersection property, so some element of~$C$ belongs to every $C_{\LocConvex}^H$.}

\end{exercises}

\section{Other characterizations of amenability}
\label{OtherAmen}

Here are a few of the many conditions that are equivalent to amenability. The necessary definitions are provided in the discussions that follow.

\begin{thm} \label{AmenEquiv}
The following are equivalent:
\noprelistbreak
\begin{enumerate}
\item \label{AmenEquiv-amen}
 $H$ is amenable.
\item \label{AmenEquiv-FP}
 $H$ has a fixed point in every nonempty, compact, convex $H$-space.
\item \label{AmenEquiv-InvMeas}
 If $H$ acts continuously on a compact, metrizable topological space~$X$, then there is an $H$-invariant probability measure on~$X$.
\item \label{AmenEquiv-Mean}
 There is a left-invariant mean on the space $\Cbdd(H)$ of all real-valued, continuous, bounded functions on~$H$.
 \item \label{AmenEquiv-FinAdd}
 There is a left-invariant finitely additive probability measure~$\rho$ defined on the collection of all Lebesgue measurable subsets of~$H$, such that $\rho(E) = 0$ for every set~$E$ of Haar measure~$0$.
\item \label{AmenEquiv-Vectors}
 The left regular representation of~$H$ on $\LL2(H)$ has almost-invariant vectors.
\item \label{AmenEquiv-Folner}
 There exists a F\o lner sequence in~$H$.
\end{enumerate}
\end{thm}

The equivalence ($\ref{AmenEquiv-amen} \Leftrightarrow \ref{AmenEquiv-FP}$) is the definition of amenability \pref{AmenDefn}. Equivalence of the other characterizations will be proved in the remainder of this section.

%\begin{rem}
%  Also equivalent to amenability:
% If $H$ acts continuously\/ {\rm(}by isometries{\rm)} on a Banach space~$\Banach$, and $C$ is a nonempty, convex, $H$-invariant subset of the dual space~$\Banach^*$ that is compact in the weak$^*$ topology, then $H$ has a fixed point in~$C$.
%\end{rem}

\subsection{Invariant probability measures}

\begin{defns} 
Let $X$ be a complete metric space.
	\begin{enumerate}
	\item A measure~$\mu$ on~$X$ is a \defit[measure!probability]{probability measure} if $\mu(X) = 1$.
	\item \nindex{$\Prob(X)$ = $\{$probability measures on~$X$$\}$}%
	$\Prob(X)$ denotes the space of all probability measures on~$X$.
	\end{enumerate}
Any measure on~$X$ is also a measure on the one-point compactification~$X^+$ of~$X$, so, if $X$ is locally compact, then the {Riesz Representation Theorem} \pref{RieszRepThm} tells us that every finite measure on~$X$ can be thought of as a linear functional on the Banach space $C(X^+)$ of continuous functions on~$X^+$. This implies that $\Prob(X)$ is a subset of the closed unit ball in the dual space $C(X^+)^*$, and therefore has a weak$^*$ topology. If $X$ is compact (so there is no need to pass to~$X^+$), then the Banach-Alaoglu Theorem \pref{BanachAlaogluThm} tells us that $\Prob(X)$ is compact \csee{Prob(X)cpctEx}.
\end{defns}

\begin{eg} \label{Prob(X)CpctConvexEg}
If a group~$H$ acts continuously on a compact, metrizable space~$X$,
%The weak$^*$ topology provides $C(X)^*$ with the structure of a locally convex vector space, such that the set $\Prob(X)$ of probability measures on~$X$ is a compact, convex set.
then $\Prob(X)$ is a compact, convex $H$-space \csee{HonProb(X)cont}.
\end{eg}

\begin{rem}[(\thmindex{Urysohn's Metrization}{Urysohn's Metrization Theorem})] \label{Metrizable<>2ndCount}
Recall that a compact, Hausdorff space is metrizable if and only if it is \term{second countable}, so requiring a compact, separable, Hausdorff space to be metrizable is not a strong restriction.
\end{rem}

\begin{prop}[($\ref{AmenEquiv-amen} \Leftrightarrow \ref{AmenEquiv-InvMeas}$)]
$H$ is amenable if and only if for every continuous action of~$H$ on a compact, metrizable space~$X$, there is an $H$-invariant probability measure~$\mu$ on~$X$.
\end{prop}

\begin{proof}
($\Rightarrow$) If $H$ acts on~$X$, and $X$ is compact, then $\Prob(X)$ is a nonempty, compact, convex $H$-space \csee{Prob(X)CpctConvexEg}. So $H$ has a fixed point in $\Prob(X)$; this fixed point is the desired $H$-invariant measure.

\medbreak

($\Leftarrow$)
Suppose $C$ is a nonempty, compact, convex $H$-space.
By replacing $C$ with the closure of the convex hull of a single $H$-orbit, we may assume $C$ is separable; then $C$ is metrizable  \csee{CpctFrechet->Metrizable}. Since $H$ is amenable, this implies there is an $H$-invariant probability measure~$\mu$ on~$C$.
Since $C$ is convex and compact, the center of mass
	$$ p = \int_C c \, d\mu(c) $$
belongs to~$C$ \csee{CenterInC}. Since $\mu$ is $H$-invariant (and the $H$-action is by linear maps), a simple calculation shows that $p$ is $H$-invariant \csee{CenterInvt}.
\end{proof}

\subsection{Invariant means}

\begin{defn}
Suppose $\LocConvex$ is some linear subspace of $\LL\infty(H)$, and assume $\LocConvex$ contains the constant function~$1_H$ that takes the value~$1$ at every point of~$H$.
A \defit{mean} on~$\LocConvex$ is
a linear functional
$\lambda$ on $\LocConvex$, such that
\noprelistbreak
 \begin{itemize}
 \item $\lambda(1_H) = 1$,
 and
 \item $\lambda$ is positive, i.e., $\lambda(f) \ge 0$
whenever $f \ge 0$.
 \end{itemize}
 \end{defn}

\begin{rem} \label{MeanNorm1}
Any mean is a continuous linear functional; indeed, $\|\lambda\| = 1$ \csee{MeanNorm1Ex}.
\end{rem}
 
 It is easy to construct means:
 
 \begin{eg} \label{EasyMeans}
 If $\phi$ is any unit vector in $\LL1(H)$, and $\mu$ is the left Haar measure on~$H$, then defining
 	$$ \lambda(f) = \int_H f\, |\phi| \, d\mu$$
produces a mean (on any subspace of $\LL\infty(H)$ that contains~$1_H$). Means constructed in this way are (weakly) dense in the set of all means \csee{EasyMeansDense}.
 \end{eg}
 
Compact groups are the only ones with invariant probability measures, but invariant means exist more generally:

\begin{prop}[($\ref{AmenEquiv-amen} \Rightarrow \ref{AmenEquiv-Mean}$)] \label{Amen->Mean}
If $H$ is amenable, then there exists a left-invariant mean on the space $\Cbdd(H)$ of bounded, continuous functions on~$H$.
\end{prop}

\begin{proof}
%To avoid technical problems, let us assume $H$ is discrete.
The set of means on $\Cbdd(H)$ is obviously nonempty, convex and invariant under left translation \csee{SetOfMeans}. Furthermore, it is a weak$^*$ closed subset of the unit ball in $\Cbdd(H)^*$ \csee{MeansCpct}, so it is compact by the Banach-Alaoglu Theorem (\cref{BanachAlaogluThm}).
Therefore, the amenability of~$H$ implies that some mean is left-invariant.
(Actually, there is a slight technical problem here if $H$ is not discrete: the action of~$H$ on $\Cbdd(H)$ may not be continuous in the sup-norm topology, because continuous functions do not need to be uniformly continuous.)
 \end{proof}

\begin{rem} \label{Amen<>MeanLinfty}
With a bit more work, it can be shown that if $H$ is amenable, then there is a left-invariant mean on $\LL\infty(H)$, not just on $\Cbdd(H)$ \csee{MeanLinfty}.
Therefore, $\Cbdd(H)$ can be replaced with $\LL\infty(H)$ in \fullcref{AmenEquiv}{Mean}.
Furthermore, there exists a mean on $\LL\infty(H)$ that is \defit[mean!bi-invariant]{bi-invariant} (both left-invariant \emph{and} right-invariant) \ccf{DiscreteBiInv}.
\end{rem}

\begin{prop}[($\ref{AmenEquiv-Mean} \Rightarrow \ref{AmenEquiv-InvMeas}$)]
Suppose $H$ acts continuously on a compact, metrizable space~$X$. If there is a left-invariant mean on $\Cbdd(H)$, then there is an $H$-invariant probability measure on~$X$.
\end{prop}

\begin{proof}
 Fix some $x \in X$. Then we have a continuous, $H$-equivariant linear map from $C(X)$ to $\Cbdd(H)$, defined by
 	$$ \overline{f}(h) = f(hx) .$$
Therefore, any left-invariant mean on $\Cbdd(H)$ induces an $H$-invariant mean~$\lambda$ on $C(X)$ \csee{MapMeans}. Since $X$ is compact, the {Riesz Representation Theorem} \pref{RieszRepThm} tells us that any continuous, positive linear functional on $C(X)$ is a measure; thus, this $H$-invariant mean~$\lambda$ can be represented by an $H$-invariant measure~$\mu$ on~$X$. Since $\lambda$ is a mean, we have $\lambda(1) = 1$, so $\mu(X) = 1$, which means that $\mu$ is a probability measure.
 \end{proof}

\subsection{Invariant finitely additive probability measures}

The following proposition is based on the observation that, just as probability measures on~$X$ correspond to elements of the dual of $C(X)$, finitely additive probability measures correspond to elements of the dual of $\LL\infty(X)$.

\begin{prop}[($\ref{AmenEquiv-Mean} \Leftrightarrow \ref{AmenEquiv-FinAdd}$)] \label{Amen<>FinAdd}
There is a left-invariant mean on $\LL\infty(X)$ if and only if there is a left-invariant finitely additive probability measure~$\rho$ defined on the collection of all Lebesgue measurable subsets of~$H$, such that $\rho(E) = 0$ for every set~$E$ of Haar measure~$0$.
\end{prop}

\begin{proof}
($\Rightarrow$) Because $H$ is amenable, there exists a left-invariant mean~$\lambda$ on $\LL\infty(H)$ \csee{Amen<>MeanLinfty}. For a measurable subset~$E$ of~$H$, let
$\rho(E) = \lambda(\chi_E)$, where $\chi_E$ is the characteristic function of~$E$.
It is easy to verify that $\rho$ has the desired properties \csee{RhoIsFinAdd}.

($\Leftarrow$) We define a mean~$\lambda$ via an approximation by step functions: for $f \in \LL\infty(H)$, let
	$$ \lambda(f) = \inf \bigset{ \sum_{i=1}^n a_i \rho(E_i) }{ 
	f \le  \sum_{i=1}^n a_i \chi_{E_i} \text{ a.e.}
	} .$$
Since $\rho$ is finitely additive, it is straightforward to verify that $\lambda$ is a mean on $\LL\infty(H)$ \csee{LambdaIsMean}. Since $\rho$ is bi-invariant, we know that $\lambda$ is also bi-invariant.
\end{proof}

\begin{rem} \label{FinAddRem} \ 
\noprelistbreak
\begin{enumerate}
\item \label{FinAddRem-subgrp}
 \cref{Amen<>FinAdd} easily implies that every subgroup of a discrete amenable group is amenable \csee{SubgrpDiscAmenFinAdd}, establishing \cref{SubgrpAmen} for the case of discrete groups. In fact, it is not very difficult to prove the general case of \cref{SubgrpAmen} similarly \csee{SubgrpAmenFinAdd}.
\item Because any amenable group~$H$ has a bi-invariant mean on $\LL\infty(H)$ \csee{Amen<>MeanLinfty}, the proof of \cref{Amen<>FinAdd}($\Rightarrow$) shows that the finitely additive probability measure~$\rho$ can be taken to be bi-invariant.
\end{enumerate}
\end{rem}

\subsection{Almost-invariant vectors}

\begin{defn} \label{AlmInvtVecDefn}
 An action of~$H$ on a normed vector space~$\Banach$ has \emph{almost-invariant} vectors if, for every
compact subset~$C$ of~$H$ and every $\epsilon > 0$, there
is a unit vector
 $v \in \Banach$, such that 
 \begin{equation} \label{epsCInvt}
 \|c v - v \| < \epsilon \text{\quad for all $c \in C$}
.\end{equation}
(A unit vector satisfying \pref{epsCInvt} is said to be $(\epsilon,C)$-invariant.)
 \end{defn}

\begin{eg}
Consider the regular
representation of~$H$ on $\LL2(H)$.
\begin{enumerate}
\item  If $H$ is a compact Lie group, then the constant
function~$1_H$ belongs to $\LL2(H)$, so $\LL2(H)$
has an $H$-invariant unit vector.

\item If $H = \real$, then $\LL2(H)$ does not have any (nonzero)
$H$-invariant vectors \csee{NonCpctNoInvtVect}, but it does have almost-invariant vectors: Given $C$ and~$\epsilon$, choose $n \in
\natural$ so large that $C \subseteq [-n,n]$ and $2/\sqrt{n}
< \epsilon$. Let $\phi = \frac{1}{n} \chi_{n^2}$, where
$\chi_{n^2}$ is the characteristic function of $[0,n^2]$.
 Then $\phi$ is a unit vector and, for $c \in C$, we
have 
  $$ ||c \phi - \phi ||^2 \le \int_{-n}^n \frac{1}{n^2}
\, dx
 + \int_{n^2-n}^{n^2 + n} \frac{1}{n^2} \, dx
 = \frac{4}{n} 
 < \epsilon^2.$$
 \end{enumerate}
 \end{eg}
 
 \begin{rem} \label{L2invtIffL1inv}
 $\LL2(H)$ has almost-invariant vectors if and only if $\LL1(H)$ has almost-invariant vectors \csee{L2invtIffL1invEx}. Therefore,  $\LL2(H)$ may be replaced with $\LL1(H)$ in \fullcref{AmenEquiv}{Vectors}. (In fact, $\LL2(H)$ may be replaced with $\LL{p}(H)$, for any $p \in [1,\infty)$ \csee{L1invtIffLpInvEx}.)
 \end{rem}

 \begin{prop}[($\ref{AmenEquiv-Mean} \Leftrightarrow \ref{AmenEquiv-Vectors}$)]
There is a left-invariant mean on $\LL\infty(H)$ if and only if $\LL2(H)$ has almost-invariant vectors.
 \end{prop}
 
\begin{proof}
Because of \cref{L2invtIffL1inv}, we may replace $\LL2(H)$ with $\LL1(H)$.

($\Leftarrow$) By applying the construction of means in \cref{EasyMeans} to almost-invariant vectors in $\LL1(H)$, we obtain almost-invariant means on $\LL\infty(H)$. A limit of almost-invariant means is invariant \csee{LimAlmInvMean}.

($\Rightarrow$) Because the means constructed in \cref{EasyMeans} are dense in the space of all means, we can approximate a left-invariant mean by an $\LL1$~function. Vectors close to an invariant vector are almost-invariant, so $\LL1(H)$ has almost-invariant vectors. However, there are technical issues here; one problem is that the approximation is in the weak$^*$ topology, but we are looking for vectors that are almost-invariant in the norm topology. See \cref{Amen->VectorsEx} for a correct proof in the case of discrete groups (using the fact that a convex set has the same closure in both the norm topology and the weak$^*$ topology).
\end{proof}

 \subsection{F\o lner sequences}
 
 \begin{defn}
 Let  $\{F_n\}$ be a sequence of measurable sets in~$H$, such that $0 < \mu(F_n) < \infty$ for every~$n$. We say $\{F_n\}$ is a \defit[Folner@F\o lner!sequence]{F\o lner sequence} if, for every compact subset~$C$ of~$H$, we have 
 	\begin{equation} \label{FolnerLim}
	 \lim_{n \to \infty} \max_{c \in C} \frac{ \mu( F_n \symmdiff c F_n) }{ \mu(F_n) } = 0 ,
	 \end{equation}
where $\mu$ is the Haar measure on~$H$.
\end{defn}

\begin{eg} \label{FolnerEg} \ 
	\begin{enumerate}
	\item  \label{FolnerEg-Rl}
	If $F_n = B_n(0)$ is the ball of radius~$n$ in~$\real^\ell$, then $\{F_n\}$ is a F\o lner sequence in~$\real^\ell$ \csee{FolnerRlEx}.
	\item The free group $\free_2$ on 2 generators does not have F\o lner sequences \csee{FreeNoFolnerEx}.
	\end{enumerate}
\end{eg}

The reason that $\real^\ell$ has a F\o lner sequence, but the free group~$\free_2$ does not, is that $\real^\ell$ is amenable, but $\free_2$ is not:

\begin{prop}[($\ref{AmenEquiv-Vectors} \Leftrightarrow \ref{AmenEquiv-Folner}$)] \label{Amen<>Folner}
There is an invariant mean on $\LL2(H)$ if and only if $H$ has a F\o lner sequence.
\end{prop}

\begin{proof}
($\Leftarrow$) Normalized characteristic functions of F\o lner sets are almost invariant vectors in $\LL1(H)$ \csee{Amen<-Folner}.

($\Rightarrow$) Let us assume $H$ is discrete.
Given $\epsilon > 0$, and a finite subset~$C$ of~$H$, we wish to find a finite subset~$F$ of~$H$, such that
	$$ \frac{ \# \bigl( F \symmdiff c(F) \bigr) }{ \#(F_n) } < \epsilon 
	\quad \text{for all $c \in C$} .$$
Since $H$ is amenable, we know $\LL1(H)$ has almost-invariant vectors \csee{L2invtIffL1inv}; hence, there exists $f \in \LL1(H)$, such that 
\noprelistbreak
	\begin{enumerate}
	\item $f \ge 0$,
	\item $\| f \|_1 = 1$,
	and
	\item $\|c f - f \|_1 < \epsilon/\#C$, for every $c \in C$.
	\end{enumerate}
Note that if $f$ were the normalized characteristic function of a set~$F$, then this set~$F$ would be what we want; for the general case, we will approximate $f$ by a sum of such characteristic functions.

Approximating $f$ by a step function, we may assume $f$ takes only finitely many values. Hence, there exist:
\noprelistbreak
	\begin{itemize}
	\item finite subsets $A_1 \subseteq A_2 \subseteq \cdots \subseteq A_n$ of~$H$,
	and
	\item real numbers $\alpha_1, \ldots,\alpha_n > 0$,
	\end{itemize}
such that
\noprelistbreak
	\begin{enumerate}
	\item $\alpha_1 + \alpha_2 + \cdots  + \alpha_n = 1$
	and
	\item $f = \alpha_1 f_1 + \alpha_2 f_2 + \cdots \alpha_n f_n$,
	\end{enumerate}
where $f_i$ is the normalized characteristic function of~$A_i$ \csee{StepFuncIsSumEx}. For all $i$ and~$j$, and any $c \in H$, we have
	\begin{equation} \label{cAiDisjAj}
	\text{$A_i \smallsetminus cA_i$ is disjoint from $cA_j \smallsetminus A_j$} 
	\end{equation}
\csee{cAiDisjAjEx}, so, for any $x \in H$, we have
	$$f_i(x) > (c f_i)(x) \implies f_j(x) \ge (cf_j)(x)$$
and
	$$f_i(x) < (cf_i)(x) \implies f_j(x) \le (cf_j)(x) .$$
Therefore
	$$ \text{$\displaystyle | (cf - f)(x)| = \sum_i  \alpha_i | (cf_i - f_i)(x)|$ \quad for all $x \in H$} .$$
Summing over~$H$ yields
	$$ \sum_i \alpha_i \|cf_i - f_i \|_1
	= \|cf - f \|_1
	< \frac{\epsilon}{\#C} .$$
Summing over~$C$, we conclude that
	$$ \sum_i \alpha_i \sum_{c \in C} \|cf_i - f_i \|_1 < \epsilon .$$
Since $\sum_i \alpha_i = 1$ (and all terms are positive), this implies there is some~$i$, such that
	$$ \sum_{c \in C} \|cf_i - f_i \|_1 < \epsilon .$$
Hence, $\|cf_i - f_i \|_1 < \epsilon$, for every~$c \in C$, so we may let $F = A_i$.
\end{proof}

\begin{rem} \label{FolnerSubgrp}
F\o lner sets provide an easy proof that subgroups of discrete amenable groups are amenable.
\end{rem}

\begin{proof}
Let
\noprelistbreak
	\begin{itemize}
	\item $A$ be a closed subgroup of a discrete, amenable group~$H$,
	\item $C$ be a finite subset of~$A$,
	and
	\item $\epsilon > 0$.
	\end{itemize}
Since $H$ is amenable, there is a corresponding F\o lner set~$F$ in~$H$. 

It suffices to show there is some $h \in H$, such that $Fh \cap A$ is a F\o lner set in~$A$.
We have
	\begin{align*}
	\#F
	= \sum_{Ah \in A \backslash H} \# (F \cap Ah)
	\end{align*}
and, letting $\epsilon' = \epsilon \#C$, we have
	\begin{align*}
	(1 + \epsilon') \#F
	\ge \#(C F) 
	= \sum_{Ah \in A \backslash H} \# \bigl( C (F \cap Ah) \bigr)
	, \end{align*}
so there must be some $Ah  \in A \backslash H$, such that 
	$$ \# \bigl( C (F \cap Ah) \bigr) \le (1 + \epsilon') \# (F \cap Ah)  $$
(and $F \cap Ah \neq \emptyset$). Then, letting $F' = Fh^{-1} \cap A$, we have
	$$ \# (CF') = \# \bigl( C (F \cap Ah) \bigr) \le (1 + \epsilon') \# (F \cap Ah) = (1 + \epsilon') \#F' ,$$
so $F'$ is a F\o lner set in~$A$.
\end{proof}

\begin{exercises}

\item[]
\exersubsection{Invariant probability measures}

\item \label{Prob(X)cpctEx}
In the setting of \cref{Prob(X)CpctConvexEg}, show that $\Prob(X)$ is a compact, convex subset of $C(X)^*$.
\hint{You may assume the Banach-Alaoglu Theorem (\cref{BanachAlaogluThm}).}

\item \label{HonProb(X)cont}
Suppose $H$ acts continuously on a compact, metrizable space~$X$. There is an induced action of~$H$ on $\Prob(X)$ defined by
	$$ \text{$(h_*\mu)(A) = \mu (h^{-1} A)$ \quad for $h \in H$, $\mu \in \Prob(X)$, and $A \subseteq X$} .$$
Show that this induced action of~$H$ on $\Prob(X)$ is continuous (with respect to the weak$^*$ topology on $\Prob(X)$).

\item \label{CpctFrechet->Metrizable}
Let $A$ be a separable subset of a Fréchet space~$\LocConvex$. Show
	\begin{enumerate}
	\item $A$ is second countable.
	\item \label{CpctFrechet->Metrizable-metrizable}
 If $A$ is compact, then $A$ is metrizable.
	\end{enumerate}
\hint{\pref{CpctFrechet->Metrizable-metrizable}~\cref{Metrizable<>2ndCount}.}

\item \label{CenterInC}
Let $\mu$ be a probability measure on a compact, convex subset~$C$ of a Fréchet space~$\LocConvex$. The \defit[center!of mass]{center of mass} of~$C$ is a point $c \in \LocConvex$, such that, for every continuous linear functional~$\lambda$ on~$\LocConvex$, we have
	$$ \lambda (c) = \int_C \lambda(x) \, d\mu(x) .$$
Show the center of mass of~$\mu$ exists and is unique, and is an element of~$C$.

\item
Give an example of a probability measure~$\mu$ on a Fréchet space, such that the center of mass of~$\mu$ does not exist.
\hint{There are probability measures on~$\real^+$, such that the center of mass is infinite.}

\item \label{CenterInvt}
Show that if $p$ is the center of mass of a probability measure~$\mu$ on a Fréchet space~$\LocConvex$, then $p$~is invariant under every continuous, linear transformation of~$\LocConvex$ that preserves~$\mu$.

\item Suppose $H$ acts continuously on a compact, metrizable space~$X$. Show that the map
	$$ H \times \Prob(X) \colon (h,\mu) \mapsto h_* \mu $$
defines a continuous action of~$H$ on $\Prob(X)$.

\exersubsection{Left-invariant means}

\item \label{MeanNorm1Ex}
Verify \cref{MeanNorm1}.
\hint{$\lambda(1_H) = 1$ implies $\|\lambda\| \ge 1$. For the other direction, note that if $\|f\|_\infty \le 1$, then $1_H-f \ge 0$ a.e., so $\lambda(1_H - f) \ge 0$; similarly, $\lambda(f + 1_H) \ge 0$.}

\item \label{MeanLtoC}
Show that the restriction of a mean is a mean. More precisely, let $\LocConvex_1$ and $\LocConvex_2$ be linear subspaces of $\LL\infty(H)$, with $1_H \in \LocConvex_1 \subseteq \LocConvex_2$. Show that if $\lambda$ is a mean on~$\LocConvex_2$, then the restriction of~$\lambda$ to~$\LocConvex_1$ is a mean on~$\LocConvex_1$.

\item Suppose $\lambda$ is a mean on $\Cbdd(H)$, the space of bounded, continuous functions on~$H$. For $f \in \Cbdd(H)$, show 
	$$ \min f \le \lambda(f) \le \max  f .$$

\item For $h \in H$, define $\delta_h \colon \Cbdd(H) \to \real$ by $\delta_h(f) = f(h)$. Show $\delta_h$ is a mean on $\Cbdd(H)$.

\item \label{EasyMeansDense}
Let $\Banach$ be any linear subspace of $\LL\infty(H)$, such that $\Banach$ contains $1_H$ and is closed in the $\LL\infty$-norm. Show that the means constructed in \cref{EasyMeans} are weak$^*$ dense in the set of all means on~$\Banach$.
\hint{If not, then the Hahn-Banach Theorem implies there exist $\epsilon > 0$, a mean~$\lambda$, and some $f \in (\Banach^*)^* = \Banach$, such that 
	$$  \lambda(f) > \epsilon + \int_H f \, |\phi| \, d\mu ,$$
for every unit vector~$\phi$ in $\LL1(H)$. This contradicts the fact that $\lambda(f) \le \mathop{\text{ess.\ sup}} f$.}

\item \label{SetOfMeans}
Let $\mathcal{M}$ be the set of means on $\Cbdd(H)$. Show:
	\begin{enumerate}
	\item \label{SetOfMeans-notempty}
	$\mathcal{M} \neq \emptyset$.
	\item $\mathcal{M}$ is convex.
	\item $\mathcal{M}$ is $H$-invariant.
	\end{enumerate}
\hint{\pref{SetOfMeans-notempty}~Evaluation at any point is a mean.}

\item \label{MeansCpct}
Let $\mathcal{M}$ be the set of means on $\Cbdd(H)$. Show:
	\begin{enumerate}
	\item $\mathcal{M}$ is contained in the closed unit ball of $\Cbdd(H)^*$. (That is, we have $|\lambda(f) \le \|f\|_\infty$ for every $f \in \Cbdd(H)$.)
	\item $\mathcal(M)$ is weak$^*$ closed.
	\item $\mathcal{M}$ is compact in the weak$^*$ topology.
	\end{enumerate}
\hint{You may assume the Banach-Alaoglu Theorem (\cref{BanachAlaogluThm}).}

\label{MeanLinfty}
Show that if $H$ is amenable, then there is a left-invariant mean on $\LL\infty(H)$.
\hint{Define $\lambda(f) = \mu_0(f * \eta)$, where $\lambda_0$ is a left-invariant mean on $\Cbdd(H)$, and $\eta$ is a nonnegative function of integral~$1$.}

\item \label{MapMeans}
Suppose $\psi \colon Y \to X$ is continuous, and $\lambda$ is a mean on $\Cbdd(Y)$. Show that $\psi_*\lambda$ (defined by $(\psi_*\lambda)(f) = \lambda ( f \circ \psi )$) is a mean on $\Cbdd(X)$.

%\item \label{MeanOnCpctIsMeas}
%Suppose $\lambda$ is a mean on $\Cbdd(X)$, and $X$ is compact. Show that $\lambda$ is a probability measure.
%\hint{You may assume the Riesz Representation Theorem, which states that any continuous, positive linear functional on $C_c(X)$ is a measure. (Here, $C_c(X)$ is the space of continuous functions with compact support, under the supremum norm.}

\item \label{DiscreteBiInv}
Assume $H$ is amenable and discrete. Show there is a bi-invariant mean on $\LL\infty(H)$.
\hint{Since $\LL\infty(H) = \Cbdd(H)$, amenability implies there is a left-invariant mean on $\LL\infty(H)$ \fullcsee{AmenEquiv}{Mean}. Now $H$ acts by right translations on the set of all such means, so amenability implies that some left-invariant mean is right-invariant.}

\item \label{MetrizableSuffices}
\harder
Assume $H$ has a fixed point in every \emph{metrizable}, nonempty, compact, convex $H$-space (and $H$ is second countable). Show $H$ is amenable.
\hint{To find a fixed point in~$C$, choose some $c_0 \in C$. For each mean~$\lambda$ on $\Cbdd(H)$ and each $\rho \in \LocConvex^*$, define $\phi_\lambda(\rho) = \lambda \bigl(h \mapsto \rho(hc_0) \bigr)$, so $\phi_\lambda \in (\LocConvex^*)'$, the algebraic dual of~$\LocConvex^*$. 
If $\lambda$ is a convex combination of evaluations at points of~$H$, 
 it is obvious there exists $c_\lambda \in C$, such that $\phi_\lambda(\rho) = \rho(c_\lambda)$. Since the map $\lambda \mapsto \phi_\lambda$ is continuous (with respect to appropriate weak topologies), this implies $c_\lambda$ exists for every~$\lambda$. The proof of \cref{Amen->Mean} shows that $\lambda$ may be chosen to be left-invariant, and then $c_\lambda$ is $H$-invariant.}
%\hint{The proof of \cref{Amen->Mean} shows there is a left-invariant mean~$\lambda$ on $\Cbdd(H)$.
%Now, to find a fixed point in~$C$, fix some $c_0 \in C$. For each $\rho \in \LocConvex^*$, define %$f_\rho \in \Cbdd(H)$ by $f_\rho(h) = \rho(hc_0)$, and then define $\phi(\rho) = \lambda(f_\rho)$, 
%$\phi(\rho) = \lambda \bigl(h \mapsto \rho(hc_0) \bigr)$, 
%so $\phi \in \LocConvex^{*{}*}$. Since $\phi$ is continuous for the weak$^*$ topology, we have $\phi \in \LocConvex$. In fact, $\phi$ is an $H$-invariant element of~$C$.}

\exersubsection{Invariant finitely additive probability measures}

\item \label{RhoIsFinAdd}
Verify that~$\rho$, as defined in the proof of \cref{Amen<>FinAdd}($\Rightarrow$), has the properties specified in the statement of the \lcnamecref{Amen<>FinAdd}.

\item \label{LambdaIsMean}
Let $\rho$ and $\lambda$ be as in the proof of \cref{Amen<>FinAdd}($\Leftarrow$). 
	\begin{enumerate}
	\item \label{LambdaIsMean-well}
 If $\displaystyle \sum_{i=1}^m a_i \, \chi_{E_i} = \sum_{j=1}^n b_j \, \chi_{F_j}$ a.e., show $\displaystyle \sum_{i=1}^m a_i \, \rho(E_i) = \sum_{j=1}^n \, b_j \rho(F_j)$.
	\item \label{LambdaIsMean-le}
 If $\displaystyle \sum_{i=1}^m a_i \, \chi_{E_i} \le \sum_{j=1}^n b_j \, \chi_{F_j}$ a.e., show $\displaystyle \sum_{i=1}^m a_i \, \rho(E_i) \le \sum_{j=1}^n b_j \, \rho(F_j)$.
 	\item Show that $\lambda(1_H) = 1$.
	 \item Show that if $f \ge 0$, then $\lambda(f) \ge 0$.
	 \item Show that 
	 $$ \lambda(f) = \sup \bigset{ \sum_{i=1}^n a_i \rho(E_i) }{ 
		f \ge  \sum_{i=1}^n a_i \chi_{E_i} \text{ a.e.}
		} .$$
	\item Show that $\lambda$ is a mean on $\LL\infty(H)$.
	\end{enumerate}
\hint{(\ref{LambdaIsMean-well},\ref{LambdaIsMean-le})~By passing to a refinement, arrange that $\{E_i\}$ are pairwise disjoint, $\{F_j\}$ are pairwise disjoint, and each $E_i$ is contained in some~$F_j$.
}

\item \label{SubgrpDiscAmenFinAdd}
Use \cref{Amen<>FinAdd} to prove that every subgroup~$A$ of a discrete amenable group~$H$ is amenable.
\hint{Let $X$ be a set of representatives of the right cosets of~$A$ in~$H$,
and let $\lambda$ be a left-invariant finitely additive probability measure on~$H$. For $E \subseteq A$, define
	$ \overline\lambda(E) = \lambda( EX )$.}

\item \label{SubgrpAmenFinAdd}
Use \cref{Amen<>FinAdd} to prove that every closed subgroup~$A$ of an amenable group~$H$ is amenable.
\hint{Let $X$ be a Borel set of representatives of the right cosets of~$A$ in~$H$, and define $\overline\lambda$ as in the solution of \cref{SubgrpDiscAmenFinAdd}. Fubini's Theorem implies that if $E$ has measure~$0$ in~$A$, then $XA$ has measure~$0$ in~$H$. You may assume (without proof) the fact that if $f \colon M \to N$ is a continuous function between manifolds $M$ and~$N$, and $E$ is a Borel subset of~$M$, such that the restriction of~$f$ to~$E$ is one-to-one, then $f(E)$ is a Borel set in~$N$.
}

\exersubsection{Almost-invariant vectors}

\item \label{NonCpctNoInvtVect} 
\begin{enumerate}
\item For $v \in \LL2(H)$, show that $v$ is invariant under translations if and only if $v$~is constant (a.e.).
\item Show that $H$ is compact if and only if  $\LL2(H)$ has a nonzero vector that is invariant under translation.
\end{enumerate}

\item \label{L2invtIffL1invEx}
Show that $\LL2(H)$ has almost-invariant vectors if and only if $\LL1(H)$ has almost-invariant vectors.
\hint{Note that $f^2 - g^2 = (f - g)(f + g)$, so $\|f^2 - g^2 \|_1 \le \| f - g\|_2 \, \| f + g \|_2$. 
Conversely, for $f,g \ge 0$, we have $(f-g)^2 \le |f^2-g^2|$, so $\| f-g \|_2^2 \le \|f^2 - g^2\|_1$.}

\item \label{L1invtIffLpInvEx}
For $p \in [1,\infty)$, show that $\LL1(H)$ has almost-invariant vectors if and only if $\LL{p}(H)$ has almost-invariant vectors. 
\hint{If $p < q$, then almost-invariant vectors in $\LL{p}(H)$ yield almost-invariant vectors in $\LL{q}(H)$, because $|(f - g)|^{q/p} \le |f^{q/p} - g^{q/p}|$. And almost-invariant vectors in $\LL{p}(H)$ yield almost-invariant vectors in $\LL{p/2}(H)$, by the argument of the first hint in \cref{L2invtIffL1invEx}.}

\item \label{LimAlmInvMean}
Let
\noprelistbreak
	\begin{itemize}
	\item $\{C_n\}$ be an increasing sequence of compact subsets of~$H$, such that $\bigcup_n C_n = H$,
	\item $\epsilon_n = 1/n$,
	\item $\phi_n$ be an $(\epsilon_n,C_n)$-invariant unit vector in $\LL1(H)$, 
	\item $\lambda_n$ be the mean on $\LL\infty$ obtained from $\phi_n$ by the construction in \cref{EasyMeans},
	and
	\item $\lambda$ be an accumulation point of $\{\lambda_n\}$.
	\end{itemize}
Show that $\lambda$ is invariant.

\item \label{Amen->VectorsEx}
Assume $H$ is discrete. Let 
	$$\mathcal{P} = \{\, \phi \in \LL1(H) \mid \phi \ge 0, \|\phi\|_1 = 1 \,\} .$$
Suppose $\{\phi_i\}$ is a net in $\mathcal{P}$, such that the corresponding means $\lambda_i$ converge weak$^*$ to an invariant mean~$\lambda$ on $\LL\infty(H)$.
	\begin{enumerate}
	\item For each $h \in H$, show that the net $\{h^*\phi_i - \phi_i\}$ converges weakly to~$0$.
	\item Take a copy $\LL1(H)_h$ of $\LL1(H)$ for each $h \in H$, and let 		$$\LocConvex = \bigtimes_{h \in H} \LL1(H)_h $$
	with the product of the norm topologies. Show that $\LocConvex$ is a Fréchet space.
	\item Show that the weak topology on~$\LocConvex$ is the product of the weak topologies on the factors.
	\item Define a linear map $T \colon \LL1(H) \to \LocConvex$ by $T(f)_h = h^* f - f$.
	\item Show that the net $\{T(\phi_i\})$ converges to $0$ weakly.
	\item Show that $0$ is in the strong closure of $T(\mathcal{P})$.
	\item Show that $\LL1(H)$ has almost-invariant vectors.
	\end{enumerate}

\item Show that if $H$ is amenable, then $H$ has the \defit{Haagerup property}. By definition, this means there is a unitary representation of~$H$ on a Hilbert space~$\Hilbert$, such that there are almost-invariant vectors, and all matrix coefficients decay to~$0$ at~$\infty$ as in the conclusion of \cref{DecayMatCoeffSimple}. (A group with the Haagerup property is also said to be \defit[a-T-menable|indsee{Haagerup property}]{a-T-menable}.)

\exersubsection{F\o lner sequences}

\item Show that $\{F_n\}$ is a F\o lner sequence if and only if, for every compact subset~$C$ of~$H$, we have
	$$ \lim_{n \to \infty} \frac{ \mu( F_n \cup c F_n) }{ \mu(F_n) } = 1 .$$

\item \label{FolnerRlEx}
Justify \fullcref{FolnerEg}{Rl}.
\hint{$C \subseteq B_r(0)$, for some~$r$. We have $\mu\bigl(B_{r+\ell}(0) \bigr)/\mu\bigl(B_{\ell}(0) \bigr) \to 1$.}

\item \label{Amen<-Folner}
Prove \Cref{Amen<>Folner}($\Leftarrow$).
\hint{Normalizing the characteristic function of $F_n$ yields an almost-invariant unit vector.}

\item \label{Folner<>ratio}
Show \pref{FolnerLim} is equivalent to
	 $$ \lim_{n \to \infty} \max_{c \in C} \frac{ \mu( F_n \cup c F_n) }{ \mu(F_n) } = 1 .$$

\item \label{Folner<>CF}
Assume $H$ is discrete. Show that a sequence $\{F_n\}$ of finite subsets of~$H$ is a F\o lner sequence if and only if, for every finite subset~$C$ of~$H$, we have
	 $$ \lim_{n \to \infty} \frac{ \#( C F_n ) }{ \#(F_n) } = 1 .$$

\item \label{StepFuncIsSumEx}
Given a step function~$f$, as in the proof of \cref{Amen<>Folner}($\Rightarrow$), let 
	\begin{itemize}
	\item $a_1 > a_2 > \cdots > a_n$ be the finitely many positive values taken by~$f$,
	\item $A_i = \{\, h \in H \mid f(h) \ge a_i \,\}$,
	and
	\item $f_i$ be the normalized characteristic function of~$A_i$.
	\end{itemize}
Show 
	\begin{enumerate}
	\item $A_1 \subseteq A_2 \subseteq \cdots \subseteq A_n$,
	\item there exist real numbers $\alpha_1,\ldots,\alpha_n > 0$, such that 
	$$ f = \alpha_1 f_1 + \cdots + \alpha_n f_n ,$$
	and
	\item $\alpha_1 + \cdots + \alpha_n = 1$.
	\end{enumerate}

\item \label{cAiDisjAjEx}
Prove \pref{cAiDisjAjEx}.
\hint{Note that either $A_i \subseteq A_j$ or $A_j \subseteq A_i$.}

\item 
\harder
Use F\o lner sets to prove \cref{FolnerSubgrp} (without assuming $H$ is discrete).
\hint{Adapt the proof of the discrete case. There are technical difficulties, but begin by replacing the sum over $A \backslash H$ with an integral over $A \backslash H$.}

\item A finitely generated (discrete) group~$\Lambda$ is said to have \defit{subexponential growth} if there exists a generating set~$S$ for $\Lambda$, such that, for every $\epsilon > 0$, 
	$$ \text{$\# (S \cup S^{-1})^n \le e^{\epsilon n}$ for all large~$n$.} $$
Show that every group of subexponential growth is amenable.

\item Give an example of an finitely generated, amenable group that does \emph{not} have subexponential growth.

\end{exercises}

\section{Some nonamenable groups} \label{NonamenSect}

Other proofs of the following proposition appear in \cref{FreeNotAmenMeas,FreeNoFolnerEx}.

\begin{prop} \label{FreeNotAmen}
Nonabelian free groups are not amenable.
\end{prop}

\begin{proof}
For convenience, we consider only the free group~$\free_2$ on two generators $a$ and~$b$. Suppose $\free_2$ has a left-invariant finitely additive probability measure~$\rho$. (This will lead to a contradiction.)

We may write $\free_2 = A \cup A^- \cup B \cup B^- \cup \{e\}$, where $A$, $A^-$, $B$, and~$B^-$ consist of the reduced words whose first letter is~$a$, $a^{-1}$, $b$, or~$b^{-1}$, respectively. Assume, without loss of generality, that $\rho(A \cup A^-) \le \rho(B \cup B^-)$ and $\rho(A) \le \rho(A^-)$. Then
	$$ \text{$\displaystyle \rho \bigl( B \cup B^- \cup \{e\} \bigr) \ge \frac{1}{2}$
	\quad
	and
	\quad
	$\displaystyle \rho(A) \le \frac{1}{4}$}
	.$$
Then, by left-invariance, we have
	$$ \rho \Bigl( a \bigl( B \cup B^- \cup \{e\} \bigr) \Bigr)
	= \rho \bigl( B \cup B^- \cup \{e\} \bigr)
	\ge \frac{1}{2}
	>  \rho(A) .$$
This contradicts the fact that $a \bigl( B \cup B^- \cup \{e\} \bigr) \subseteq A$.
\end{proof}

Combining this with the fact that subgroups of discrete amenable groups are amenable \csee{SubgrpAmen}, we have the following consequence:

\begin{cor} \label{FreeSubgrp->Nonamen}
Suppose $H$ is a discrete group. 
If $H$ contains a nonabelian, free subgroup, then $H$ is not amenable.
\end{cor}

\begin{rems} \label{FreeSubgrpRem} \ 
\noprelistbreak
\begin{enumerate}
\item  \label{FreeSubgrpRem-Olshanskii}
The converse of \cref{FreeSubgrp->Nonamen} is known as ``\thmindex{von\,Neumann's Conjecture}von\,Neumann's Conjecture\zz,'' but it is false: a nonamenable group with no nonabelian free subgroups was found by  Ol'shanskii in 1980. 
(The name is misleading: apparently, the conjecture is due to M.\,Day, and was never stated by Von Neumann.)
\item The assumption that $H$ is discrete cannot be deleted from the statement of \cref{FreeSubgrp->Nonamen}.
For example, the orthogonal group $\SO(3)$ is amenable (because it is compact), but the Tits Alternative \pref{TitsAlternative} implies that it contains nonabelian free subgroups. 
\item \label{FreeSubgrpRem-BanachTarski}
 The nonamenability of nonabelian free subgroups of $\SO(3)$ is the basis of the famous \thmindex{Banach-Tarski Paradox}Banach-Tarski Paradox: A $3$-dimensional ball~$B$ can be decomposed into finitely many subsets $X_1,\ldots,X_n$, such that these subsets can be reassembled to form the union of two disjoint balls of the same radius as~$B$. (More precisely, the union $B_1 \cup B_2$ of two disjoint balls of the same radius as~$B$ can be decomposed into subsets $Y_1,\ldots,Y_n$, such that $Y_i$ is congruent to~$X_i$, for each~$i$.)
\item If $H$ contains a \emph{closed}, nonabelian, free subgroup, then $H$ is not amenable.
\end{enumerate}
\end{rems}

Here is an example of a nonamenable connected group:

\begin{prop} \label{SL2RNotAmen}
$\SL(2,\real)$ is not amenable.
\end{prop}

\begin{proof}
Let $G = \SL(2,\real)$. The action of $G$ on $\real \cup \{\infty\} \iso \circle$ by linear-fractional transformations is transitive, and the stabilizer of the point~$0$ is the subgroup
	$ P = \begin{Smallbmatrix} *&* \\ &* \end{Smallbmatrix} $,
so $G/P$ is compact. However, the Borel Density Theorem implies there is no $G$-invariant probability measure on $G/P$ \csee{BDT-G/H}. (See \cref{NoProbMeasOnRP1} for a direct proof that there is no $G$-invariant probability measure.) So $G$ is not amenable.
\end{proof}

More generally:

\begin{prop} \label{GNotAmen}
If a connected, semisimple Lie group $G$ is not compact, then $G$ is not amenable.
\end{prop}

\begin{proof}
The {Jacobson-Morosov Lemma} \pref{JacobsonMorosov} tells us that $G$ contains a closed subgroup isogenous to $\SL(2,\real)$. Alternatively, recall that any lattice $\Gamma$ in~$G$ must contain a nonabelian free subgroup \csee{FreeInGamma}, and, being discrete, this is a closed subgroup of~$G$.
\end{proof}

\begin{rem}
Readers familiar with the structure of semisimple Lie groups will see that the proof of \cref{SL2RNotAmen} generalizes to the situation of \cref{GNotAmen}: Since $G$ is not compact, it has a proper parabolic subgroup~$P$. Then $G/P$ is compact, but the Borel Density Theorem implies that $G/P$ has no $G$-invariant probability measure.
\end{rem}

Combining this result with the structure theory of connected Lie groups yields the following classification of connected, amenable Lie groups:

\begin{prop} \label{ConnAmen}
A connected Lie group~$H$ is amenable if and only if $H$ contains a connected, closed, solvable normal subgroup~$N$, such that $H/N$ is compact.
\end{prop}

\begin{proof}
($\Leftarrow$) \fullcref{AmenExtCor}{CocpctSolv}.

($\Rightarrow$) The structure theory of Lie groups tells us that there is a connected, closed, solvable, normal subgroup~$R$ of~$H$, such that $H/R$ is semisimple. (The subgroup~$R$ is called the \defit[radical!of a Lie group]{radical} of~$H$.) Since quotients of amenable groups are amenable \csee{QuotAmen}, we know that $H/R$ is amenable. So $H/R$ is compact \csee{GNotAmen}.
\end{proof}

\begin{exercises}

\item \label{FreeNotAmenMeas}
\noprelistbreak
	\begin{enumerate}
	\item \label{FreeNotAmenMeas-ptmass}
Find a homeomorphism $\phi$ of the circle~$\circle$, such that the only $\phi$-invariant probability measure is the delta mass at a single point~$p$.
	\item \label{FreeNotAmenMeas-diff}
 Find two homeomorphisms $\phi_1$ and~$\phi_2$ of~$\circle$, such that the subgroup $\langle \phi_1,\phi_2 \rangle$ they generate has no invariant probability measure.
	\item Deduce that the free group~$\free_2$ on 2 generators is not amenable.
	\end{enumerate}
\hint{\pref{FreeNotAmenMeas-ptmass}~Identifying $S^1$ with $[0,1]$, let $\phi(x) = x^2$. For any $x \in (0,1)$, we have $\phi \bigl( (0,x) \bigr) = (0,x^2)$, so $\mu \bigl( (x^2,x) \bigr) = 0$. Since $(0,1)$ is the union of countably many such intervals, this implies that $\mu \bigl( (0,1) \bigr) = 0$.}

\item \label{FreeNoFolnerEx}
Show explicitly that free groups do not have F\o lner sequences. More precisely, let $\free_2$ be the free group on two generators~$a$ and~$b$, and show that if $F$ is any nonempty, finite subset of~$\free_2$, then there exists $c \in \{a,b,a^{-1},b^{-1}\}$, such that $\#(F \smallsetminus cF) \ge (1/4) \#F$.
This shows that $\free_2$ free groups is not amenable.
\hint{Suppose $F = A \cup B \cup A^- \cup B^-$, where words in $A,B,A^-,B^-$ start with $a,b,a^{-1},b^{-1}$, respectively. If $\#A \le \# A^-$ and $\# (A \cup A^-) \le \#(B \cup B^-)$, then $\#(aF \smallsetminus F) \ge \#(B \cup B^-) - \#A$.}

\item Assume that $H$ is discrete, and that $H$ is isomorphic to a (not necessarily discrete) subgroup of $\SL(\ell,\real)$. Show: 
\noprelistbreak
	\begin{enumerate}
	\item $H$ is amenable if and only if $H$ has no nonabelian, free subgroups.
	\item $H$ is amenable if and only if $H$ has a solvable subgroup of finite index.
	\end{enumerate}
\hint{Tits Alternative \pref{TitsAlternative}.}

\item \label{NoProbMeasOnRP1}
Let $G = \SL(2,\real)$ act on $\real \cup \{\infty\}$ by linear-fractional transformations, as usual.
\noprelistbreak
	\begin{enumerate}
	\item \label{NoProbMeasOnRP1-u}
 For 
		$ u = \left[\begin{smallmatrix} 1 & 1 \\ 0 & 1 \end{smallmatrix} \right] \in G $,
	show that the only $u$-invariant probability measure on $\real \cup \{\infty\}$ is concentrated on the fixed point of~$u$.
	\item Since the fixed point of~$u$ is not fixed by all of~$G$, conclude that there is no $G$-invariant probability measure on $\real \cup \{\infty\}$.
	\end{enumerate}
\hint{\pref{NoProbMeasOnRP1-u}~The action of~$u$ is conjugate to the homeomorphism~$\phi$ in the hint to \fullcref{FreeNotAmenMeas}{ptmass}, % !!!
so a similar argument applies.}

\item \label{GammaNotAmen}
Show that if a semisimple Lie group~$G$ is not compact, then every lattice~$\Gamma$ in~$G$ is not amenable.

%\item Show that if $G$ is amenable, then $G$ is compact.
%\hint{If $G$ is not compact, then $G$ has a proper parabolic subgroup~$P$. Since $G/P$ is compact, amenability implies there is a $G$-invariant probability measure on $G/P$. The Borel Density Theorem \pref{BDT-G/H} provides a contradiction.}

\item Give an example of a nonamenable Lie group that has a closed, cocompact, amenable subgroup. (By \cref{ConnAmen}, the subgroup cannot be normal.)

\end{exercises}

\section{Closed subgroups of amenable groups}
\label{AmenSubgrpSect}

Before proving that closed subgroups  of amenable groups  are amenable (\cref{SubgrpAmen}), we introduce some notation and establish a \lcnamecref{Linfty(H;C)good}. (Proofs for the case of discrete groups have already been given in \cref{FinAddRem,FolnerSubgrp}.)

\begin{notation} \ 
\noprelistbreak
	\begin{enumerate}
	\item We use \nindex{$\LL\infty(H;C) = \{\text{bounded functions from~$H$ to~$C$} \}$}$\LL\infty(H;C)$ to denote the space of all measurable functions from the Lie group~$H$ to the compact, convex set~$C$, where two functions are identified if they are equal a.e.\ (with respect to the Haar measure on~$H$). 
	\item If $\Lambda$ is a closed subgroup of~$H$, and $C$ is a $\Lambda$-space, then%
	\nindex{$\LLequi\Lambda(H;C) = \{\text{essentially $\Lambda$-equivariant maps in $\LL\infty(H;C)$}\}$} % no page break here !!!
		$$ \LLequi\Lambda(H;C) 
		= \bigset{ \psi \in \LL\infty(H;C) }{ 
		\begin{matrix} \text{$\psi$ is essentially} \\ \text{$\Lambda$-equivariant} \end{matrix}
		} .$$
	(To say $\psi$ is \defit[essentially!$\Lambda$-equivariant]{essentially $\Lambda$-equivariant} means, for each $\lambda \in \Lambda$, that $\psi(\lambda h) = \lambda \cdot \psi(h)$ for a.e.\ $h \in H$.)
	\end{enumerate}
\end{notation}

\begin{egs} \ 
\noprelistbreak
\begin{enumerate}
\item Suppose $H$ is discrete. Then every function on~$H$ is measurable, so $\LL\infty(H;C) = C^H$ is the cartesian product of countably many copies of~$C$. Therefore, in this case, {Tychonoff's Theorem}~\pref{TychonoffThm} implies that $\LL\infty(H;C)$ is compact.
\item If $C$ is the closed unit disk in the complex plane (and $H$ is arbitrary), then $\LL\infty(H;C)$ is the closed unit ball in the Banach space $\LL\infty(H)$, so the Banach-Alaoglu Theorem (\cref{BanachAlaogluThm}) states that it is compact in the weak$^*$ topology. 
\end{enumerate}
\end{egs}

More generally, if we put a technical restriction on~$C$, then there is a weak topology on $\LL\infty(H;C)$ that makes it into a compact, convex $H$-space:

\begin{lem} \label{Linfty(H;C)good}
Assume 
\noprelistbreak
	\begin{itemize}
	\item $\Lambda$ is a closed subgroup of~$H$,
	\item $C$ is a nonempty, compact, convex $H$-space,
	and
	\item $C$ is contained in the dual of some separable Banach space~$\Banach$.
	\end{itemize}
Then $\LL\infty(H; C)$ and $\LLequi\Lambda(H; C)$ are nonempty, compact, convex $H$-spaces.
\end{lem}

\begin{proof}
Let $\LL\infty(H;\Banach^*)$ be the space of all bounded measurable functions from~$H$ to~$\Banach^*$ (where two functions are identified if they are equal a.e.). This is the dual of the (separable) Banach space $\LL1(H;\Banach)$, so it has a natural weak$^*$ topology. Since $\LL\infty(H; C)$ is a closed, bounded, convex subset of $\LL\infty(H;\Banach^*)$, the Banach-Alaoglu Theorem \pref{BanachAlaogluThm} tells us that it is weak$^*$ compact.
In addition, the action of~$H$ by right-translation on $\LL\infty(H; C)$ is continuous \csee{HContOnL(H;C)}.

It is not difficult to see that $\LLequi\Lambda\!(H;C)$ is a nonempty, closed, convex, $H$-invariant subset \csee{Equi(H;C)ClosedConvex}.
\end{proof}

\begin{proof}[\bf Proof of \cref{SubgrpAmen}]
Let $\Lambda$ be a closed subgroup of an amenable Lie group~$H$.
Given any continuous action of~$\Lambda$ on a compact, metrizable space~$X$, it suffices to show there is a $\Lambda$-invariant probability measure on~$X$ \fullcsee{AmenEquiv}{InvMeas}.
From \cref{Linfty(H;C)good}, we know that $\LLequi\Lambda \bigl( H; \Prob(X) \bigr)$ is a nonempty, compact, convex $H$-space. Therefore, the amenability of~$H$ implies that $H$ has a fixed point~$\psi$ in~$\LLequi\Lambda (H;C)$. So $\psi$ is essentially $H$-invariant. If we fix any $\lambda \in \Lambda$, then, for a.e.\ $h \in H$, we have
	\begin{align*}
	\lambda \cdot \psi(h) 
	&= \psi(\lambda h) 
	&& \text{($\psi$ is essentially $\Lambda$-equivariant)} 
	\\&= \psi(h) 
	&& \text{($\psi$ is essentially $H$-invariant)} 
	. \end{align*}
If we assume, for simplicity, that $\Lambda$ is countable \csee{SubgrpIsAmen(NotDiscrete)}, then the quantifiers can be reversed (because the union of countably many null sets is a null set), so we conclude that the probability measure $\psi(h)$ is $\Lambda$-invariant.
\end{proof}

\begin{exercises}

\item \label{HContOnL(H;C)}
Show that the action of~$H$ on $\LL\infty(H;C)$ by right translations is continuous in the weak$^*$-topology.

\item \label{GammaEquiEx}
  Suppose $\Lambda$ is a closed subgroup of~$H$, and that $\Lambda$ acts measurably on a measure
space~$\Omega$. Show there is a $\Lambda$-equivariant,
measurable map $\psi \colon H \to \Omega$.
 \hint{$\psi$ can be defined arbitrarily on a strict
fundamental domain for~$\Lambda$ in~$H$.}

\item \label{Equi(H;C)ClosedConvex}
Show that $\LLequi\Lambda(H;C)$ is a nonempty, closed, convex, $H$-invariant subset of $\LL\infty(H;C)$.

\item \label{SubgrpIsAmen(NotDiscrete)}
Prove \cref{SubgrpAmen} without assuming $\Lambda$ is countable.
\hint{Consider $\lambda$ in a countable dense subset of~$\Lambda$.}

\end{exercises}

 \section{Equivariant maps from \texorpdfstring{$G/P$}{G/P} to \texorpdfstring{$\Prob(X)$}{Prob(X)}}\label{AmenEquiMapSect}

We now use amenability to prove a basic result that has important consequences for the theory of arithmetic groups. In particular, it is an ingredient in two fundamental results of G.\,A.\,Margulis: his Superrigidity Theorem \pref{MargSuperG'} and his Normal Subgroups Theorem \pref{MargNormalSubgrpsThm}.

\begin{prop}[(\thmindex{Furstenberg's Lemma}Furstenberg's Lemma)] \label{G/amen->Meas(X)}
 If 
\noprelistbreak
 \begin{itemize}
% \item $\Gamma$ is a lattice in~$G$,
 \item $P$ is a closed, amenable subgroup of~$G$,
 and
 \item $\Gamma$ acts continuously on a compact metric space~$X$,
 \end{itemize}
 then there is a Borel measurable map $\psi \colon G/P \to
\Prob(X)$, such that $\psi$ is essentially $\Gamma$-equivariant.
 \end{prop}

\begin{proof}
%The argument is similar to the proof of \cref{SubgrpAmen}.
%
\Cref{Linfty(H;C)good} tells us that $\LLequi\Gamma \bigl( G;\Prob(X) \bigr)$
is a nonempty, compact, convex $G$-space. By restriction, it is also a nonempty, compact, convex $P$-space, so $P$ has a fixed point~$\psi_0$ (under the action by right-translation). 
Then $\psi_0$ factors through
to an (essentially) well-defined map $\psi \colon G/P \to
\Prob(X)$. Because $\psi_0$ is $\Gamma$-equivariant, it is
immediate that $\psi$~is $\Gamma$-equivariant.
 \end{proof}
 
In applications of \cref{G/amen->Meas(X)}, the subgroup~$P$ is usually taken to be a \index{parabolic!subgroup!minimal}minimal parabolic subgroup. % \csee{MinParabDefn,MinParabIsAmen}. 
Here is an example of this:

\begin{cor} \label{SL3R/P->Meas(X)}
If
\noprelistbreak
	\begin{itemize}
	\item $G = \SL(3,\real)$,
	\item $ P
	 = \begin{Smallbmatrix} \upast&& \\ \upast&\upast& \\ \upast&\upast&\upast \end{Smallbmatrix}
	 \subset G $,
	 and
	 \item $\Gamma$ acts continuously on a compact metric
space~$X$,
 \end{itemize}
 then there is a Borel measurable map $\psi \colon G/P \to
\Prob(X)$, such that $\psi$ is essentially
$\Gamma$-equivariant.
 \end{cor}

\begin{proof}
$P$ is amenable \csee{BorelSL3Amen}.
\end{proof}

\begin{rem} The function~$\psi$ that is provided by Furstenberg's Lemma \pref{G/amen->Meas(X)} (or \cref{SL3R/P->Meas(X)}) can be thought of as being a ``random'' map from~$G/P$ to~$X$; for each $z \in G/P$, the value of~$\psi(z)$ is a probability distribution that defines a random value for the function at the point~$z$. However, we will see in \cref{QuickProximalitySect} that the theory of proximality makes it possible to show, in certain cases, that $\psi(z)$ is actually a single well-defined point of~$X$, not a random value that varies over some range.
\end{rem}

\begin{exercises}

\item \label{MinParabIsAmen}
Show that every minimal parabolic subgroup of~$G$ is amenable.
\hint{Langlands decomposition \pref{LanglandsDecomp}.}

\end{exercises}

\section{More properties of amenable groups \texorpdfstring{\optional}{(optional)}}

In this section, we mention (without proof, and without even defining all of the terminology) a variety of very interesting properties of amenable groups. For simplicity, 
	$$ \text{we assume $\Lambda$ is a discrete group.} $$

\subsection{Bounded harmonic functions}

\begin{defn} 
Fix a probability measure~$\mu$ on~$\Lambda$.	
	\begin{enumerate}
	\item A function $f \colon \Lambda \to \real$ is \defit[harmonic!function]{$\mu$-harmonic} if $f = \mu *f$. This means, for every $\lambda \in \Lambda$, 
	$$ f(\lambda) = \sum_{x \in \Lambda} \mu(x) \, f( x \lambda )  .$$
	\item $\mu$ is \defit[symmetric!measure]{symmetric} if $\mu(A^{-1}) = \mu(A)$ for every $A \subseteq \Lambda$.
	\end{enumerate}
\end{defn}

\begin{thm} \label{AmenNoBddHarm}
$\Lambda$ is amenable if and only if there exists a symmetric probability measure~$\mu$ on~$\Lambda$, such that 
\noprelistbreak
	\begin{enumerate}
	\item the support of~$\mu$ generates~$\Lambda$,
	and
	\item every bounded, $\mu$-harmonic function on~$\Lambda$ is constant.
	\end{enumerate}
\end{thm}

Because any harmonic function is the Poisson integral of a function on the Poisson boundary (and vice-versa), this result can be restated in the following equivalent form:

\begin{cor}
$\Lambda$ is amenable if and only if there exists a symmetric probability measure~$\mu$ on~$\Lambda$, such that
\noprelistbreak
	\begin{enumerate}
	\item the support of~$\mu$ generates~$\Lambda$,
	and
	\item the Poisson boundary of~$\Lambda$ {\rm(}with respect to~$\mu${\rm)} consists of a single point.
	\end{enumerate}
\end{cor}

\subsection{Norm of a convolution operator}

\begin{defn}
For any probability measure~$\mu$ on~$\Lambda$, there is a corresponding convolution operator $C_\mu$ on $\LL2(\Lambda)$, defined by
	$$ (C_\mu f )(\lambda) = \sum_{x \in \Lambda} \mu(x) \, f(x^{-1} \lambda) .$$
\end{defn}

\begin{thm} \label{Amen<>ConvOp} 
Let $\mu$ be any probability measure on~$\Lambda$, such that the support of~$\mu$ generates~$\Lambda$. Then $\| C_\mu \| = 1$ if and only if $\Lambda$ is amenable.
\end{thm}

\subsection{Spectral radius}

In geometric terms, the following famous result characterizes amenability in terms of the spectral radius of random walks on Cayley graphs.

\begin{thm}[(Kesten)]
Let $\mu$ be a finitely supported, symmetric probability measure on~$\Lambda$, such that the support of~$\mu$ generates~$\Lambda$. Then $\mu$ is amenable if and only if
	$$ \lim_{n \to \infty} \left( 
	\raise13pt\hbox{$\displaystyle % @@@ raise to reduce blank space in the formula
	 \sum_{\begin{matrix}
	g_1,\ldots,g_{2n} \in \mathop{\mathrm{supp}} \mu
	\\
	g_1g_2 \cdots g_{2n} = e
	\end{matrix}}
	\mu(g_1) \, \mu(g_2) \cdots \mu(g_{2n}) 
	$}
	\right)^{1/2n} = 1
	. $$
\end{thm}

\subsection{Positive-definite functions}

\begin{defn}[\ccf{PosDefTerm}]
A $\complex$-valued function~$\varphi$ on~$\Lambda$ is \defit[positive!definite]{positive-definite} if, for all $a_1,\ldots,a_n \in \Lambda$, the matrix
	$$ [a_{i,j}]_{i,j=1}^n = \bigl[ \varphi( a_i^{-1} a_j ) \bigr] $$
is Hermitian and has no negative eigenvalues.
\end{defn}

\begin{thm}
$\Lambda$ is amenable if and only if\/ $\sum_{g \in \Lambda} \varphi(g) \ge 0$ for every\/ \textup(finitely supported\/\textup) positive-definite function~$\varphi$ on~$\Lambda$.
\end{thm}

%Wikipedia: Godement condition. Every bounded positive-definite measure ? on G satisfies ?(1) ³ 0. Valette (1998) improved this criterion by showing that it is sufficient to ask that, for every continuous positive-definite compactly supported function f on G, the function ?Ð?f has non-negative integral with respect to Haar measure, where ? denotes the modular function.

\subsection{Growth}

\begin{defn}
Assume $\Lambda$ is finitely generated, and fix a symmetric generating set~$S$ for~$\Lambda$. 
%Assume, for simplicity, that $e \in S$.
\noprelistbreak
	\begin{enumerate}
	\item For each $r \in \integer^+$, let $B_r(\Lambda)$ be the ball of radius~$r$ centered at~$e$, More precisely,
		$$ B_r(\Lambda; S) = \{\, \lambda \in \Lambda \mid \exists s_1,s_2,\ldots,s_r \in S \cup \{e\}, \lambda = s_1 s_2 \cdots s_r \,\} .$$
%	\item We say $\Lambda$ has \emph{polynomial growth} if there is some polynomial function $f \colon \real \to \real$, such that $\# B_r < f(r)$, for all $r \in \integer^+$.
	\item We say $\Lambda$ has \emph{subexponential growth} if for every $\epsilon > 0$, we have $\# B_r(\Lambda; S) < e^{\epsilon r}$, for all sufficiently large $r \in \integer^+$.
	\end{enumerate}
\end{defn}

\begin{prop}[\csee{SubexpIsAmenEx}] \label{SubexpIsAmen}
If $\Lambda$ has subexponential growth, then $\Lambda$ is amenable.
\end{prop}

\begin{warn}
The implication in \cref{SubexpIsAmen} goes only one direction: there are many groups (including many solvable groups) that are amenable, but do not have subexponential growth \csee{AmenNotSubExpEx}.
\end{warn}

\subsection{Cogrowth}

\begin{defn} \label{CoGrowthDefn}
Assume $\Lambda$ is finitely generated. Let:
	\begin{enumerate}
	\item  $S = \{s_1,s_2,\ldots,s_k\}$ be a finite generating set of~$\Lambda$.  
	\item  $F_k$ be the free group on~$k$ generators $x_1,\ldots,x_k$.
	\item  $\phi_S \colon F_k \to \Lambda$ be the homomorphism defined by $\phi(x_i) = s_i$.
	\end{enumerate}
The \defit[cogrowth of a discrete group]{cogrowth} of~$\Lambda$ (with respect to~$S$) is
	$$ \lim_{r \to \infty} \frac{1}{r} \log_{2k-1}  \# \bigl( (\ker \phi_S) \cap B_r(F_k; x_1^{\pm1},\ldots,x_k^{\pm1} \bigr) .$$
\end{defn}

Note that $\#B_r (F_k; x_1^{\pm1},\ldots,x_k^{\pm1})$ is equal to the number of reduced words of length~$r$ in the symbols $x_1^{\pm1},\ldots,x_k^{\pm1}$, which is approximately $(2k-1)^r$. Therefore, it is easy to see that that the cogrowth of~$\Lambda$ is between $0$ and $1$ \csee{CoGrowthBtwn01Ex}. The maximum value is obtained if and only if $\Lambda$~is amenable:

\begin{thm} \label{AmenCoGrowth}
$\Lambda$ is amenable if and only if the cogrowth of~$\Lambda$ is~$1$, with respect to some\/ \textup(or, equivalently, every\/\textup) finite generating set~$S$.
\end{thm}

\subsection{Unitarizable representations}

\begin{defn}
Let $\rho \colon \Lambda \to \bddop(\Hilbert)$ be a (not necessarily unitary) representation of~$\Lambda$ on a Hilbert space~$\Hilbert$.
	\begin{enumerate}
	\item $\rho$ is \defit[representation!uniformly bounded]{uniformly bounded} if there exists $C > 0$, such that $\| \rho(\lambda) \| < C$, for all $\lambda \in \Lambda$.
	\item $\rho$ is \defit[representation!unitarizable]{unitarizable} if it is conjugate to a unitary representation. This means there is an invertible operator~$T$ on~$\Hilbert$, such that the representation $\lambda \mapsto T^{-1} \, \rho(\lambda) \, T$ is unitary.
	\end{enumerate}
\end{defn}

It is fairly obvious that every unitarizable representation is uniformly bounded \csee{UnitIsUnifBddEx}. The converse is not true, although it holds for amenable groups:

\begin{thm} \label{AmenUnifBddIsUnit}
If $\Lambda$ is amenable, then every uniformly bounded representation of~$\Lambda$ is unitarizable.
\end{thm}

\begin{rem}
The converse of \cref{AmenUnifBddIsUnit} is an open question.
\end{rem}

\subsection{Almost representations are near representations}

\begin{defn}
Fix $\epsilon > 0$, and let $\varphi$ be a function from~$\Lambda$ to the group $\unitop(\Hilbert)$ of unitary operators on a Hilbert space~$\Hilbert$.
	\begin{enumerate}
	\item $\varphi$ is \defit[representation!almost]{$\epsilon$-almost} a unitary representation if 
		$$ \text{$\| \varphi(\lambda_1 \lambda_2) -  \varphi(\lambda_1) \, \varphi( \lambda_2) \| < \epsilon$ for all $\lambda_1,\lambda_2 \in \Lambda$.} $$
	\item $\varphi$ is \defit[representation!near]{$\epsilon$-near} a unitary representation if there exists a unitary representation $\rho \colon \Lambda \colon \unitop(\Hilbert)$, such that 
		$$ \text{$\| \varphi(\lambda) - \rho(\lambda) \| < \epsilon$ for every $\lambda \in \Lambda$.} $$
	\end{enumerate}
\end{defn}

For amenable groups, every almost representation is near a representation:

\begin{thm} \label{AlmRepsAreNear}
Assume $\Lambda$ is amenable. Given $\epsilon > 0$, there exists $\delta > 0$, such that if $\varphi$ is $\delta$-almost a unitary representation, then $\varphi$~is $\epsilon$-near a unitary representation.
\end{thm}

\subsection{Bounded cohomology}

The bounded cohomology groups of~$\Lambda$ are defined just like the ordinary cohomology groups, except that all cochains are assumed to be bounded functions.

\begin{defn}
Assume $\Banach$ is a Banach space.
\noprelistbreak
	\begin{enumerate}
	\item $\Banach $ is a \defit[Banach!$\Lambda$-module]{Banach $\Lambda$-module} if $\Lambda$ acts continuously on~$\Banach$, by linear isometries.
	\item $\cocyc1\bdd(\Lambda;\Banach) = \cocyc1(\Lambda;\Banach) \cap \LL\infty(\Lambda;\Banach)$.
	\item $\coho1\bdd(\Lambda;\Banach) = \cocyc1\bdd(\Lambda;\Banach) / \cobdry1(\Lambda;\Banach)$.
	\end{enumerate}
\end{defn} 

\begin{thm} \label{AmenBddCoho0}
$\Lambda$ is amenable if and only if\/ $\coho1\bdd(\Lambda;\Banach) = 0$ for every Banach $\Lambda$-module~$\Banach$.
\end{thm}

\begin{rem} \label{AmenBddCoho0Alln}
In fact, if $\Lambda$ is amenable, then $\coho {\raise 1pt \hbox{$\scriptstyle n$}}\bdd(\Lambda;\Banach) = 0$ for all~$n$.
\end{rem}

\subsection{Invariance under quasi-isometry}

\begin{prop}[\csee{AmenQIinvtEx}] \label{AmenQIinvt}
Assume $\Lambda_1$ and~$\Lambda_2$ are finitely generated groups, such that $\Lambda_1$~is quasi-isometric to~$\Lambda_2$ \csee{QuasiIsomDefn}.
Then $\Lambda_1$~is amenable if and only if $\Lambda_2$~is amenable.
\end{prop}

\subsection{Ponzi schemes}

Assume $\Lambda$ is finitely generated, and let $d$ be the word metric on~$\Lambda$, with respect to some finite, symmetric generating set~$S$ \csee{WordMetricDefn}.

\begin{defn}
A function $f \colon \Lambda \to \Lambda$ is a \defit{Ponzi scheme on~$\Lambda$} if there is some $C > 0$, such that, for all $\lambda \in \Lambda$, we have:
		\begin{enumerate}
		\item $\# \, f^{-1}(\lambda) \ge 2$,
		and
		\item $d \bigl( f(\lambda), \lambda \bigr) < C$.
		\end{enumerate}
\end{defn}

\begin{thm} \label{AmenNoPonzi}
%Assume $S$ is a finite, symmetric generating set of~$\Lambda$.
 $\Lambda$ is amenable if and only if there does not exist a Ponzi scheme on~$\Lambda$.
\end{thm}

\begin{exercises}

\item \label{SubexpIsAmenEx}
Prove \cref{SubexpIsAmen}.
\hint{If no balls are F\o lner sets, then the group has exponential growth.}

\item \label{AmenNotSubExpEx}
Choose a prime number~$p$, and let
	$$ \Lambda = \bigset{
	\begin{bmatrix} p^k & mp^n \\ 0 & p^{-k} \end{bmatrix} 
	}{ k,m,n \in \integer}
	\subset \SL(2,\rational) ,$$
with the discrete topology. Show $\Lambda$ is an amenable group that does not have subexponential growth.

\item \label{CoGrowthBtwn01Ex}
In the notation of \cref{CoGrowthDefn}, show:
	\begin{enumerate}
	\item $\# B_r(F_k; x_1^{\pm1},\ldots,x_k^{\pm1}) = 2k(2k-1)^{r-1}$.
	\item If $\mathop{\text{\rm cog}} \Lambda$ is the cogrowth of~$\Lambda$, then $0 \le \mathop{\text{\rm cog}} \Lambda \le 1$.
	\end{enumerate}

\item \label{UnitIsUnifBddEx}
Show that every unitarizable representation is uniformly bounded.

\item For every $\epsilon > 0$, show there exists $\delta > 0$, such that if $\varphi$ is $\delta$-near a unitary representation, then $\varphi$ is $\epsilon$-almost a unitary representation.

\item \label{AmenQIinvtEx}
Prove \cref{AmenQIinvt}.
\hint{Show that if $\Lambda$ is not amenable, then, for every $k$, it has a finite subset~$S$, such that $\#(S F) \ge k \cdot \# F$ for every finite subset~$F$ of~$\Lambda$.}

\item Explicitly construct a Ponzi scheme on the free group with two generators.

\item Show (without using \cref{AmenNoPonzi}) that if $\Lambda$ is amenable, then there does not exist a Ponzi scheme on~$\Lambda$.
\hint{F\o lner sets.}

\end{exercises}

\begin{notes}

The notion of amenability is attributed to J.\,von\,Neumann \cite{vonNeumann-AllgemeineMasses}, but he used the German word ``messbar'' (which can be translated as ``measurable''). The term ``amenable'' was apparently introduced into the literature by M.\,Day \cite[\#507, p.~1054]{Day-MeansAbstract} in the announcement of a talk.

The monographs \cite{PatersonBook,PierBook} are standard references on amenability.
Briefer treatments are in \cite[App.~G]{BekkaHarpeValette-T}, \cite{GreenleafBook}, and \cite[\S4.1]{ZimmerBook}.
Quite a different approach to amenability appears in \cite[Chaps.~10--12]{WagonBook} (for discrete groups only).

The fact that closed subgroups of amenable groups are amenable (\cref{SubgrpAmen}) is proved in \cite[Thm.~2.3.2, pp.~30--32]{GreenleafBook}, \cite[Prop.~13.3, p.~118]{PierBook}, and \cite[Prop.~4.2.20, p.~74]{ZimmerBook}.

See \cite[p.~67]{GreenleafBook} for a proof of \cref{Amen<>Folner}($\Rightarrow$) that does not require $H$ to be discrete.

\Cref{L2invtIffL1inv} is proved in \cite[pp.~46--47]{GreenleafBook}.

The solution of \cref{MetrizableSuffices} can be found in \cite[Thm.~5.4, p.~45]{PierBook}.

For a proof of the fact (mentioned in the hint to \cref{SubgrpAmenFinAdd}) that the one-to-one continuous image of a Borel set is Borel, see \cite[Thm.~3.3.2, p.~70]{Arveson-Cstar}.

Our proof of \cref{Amen<>Folner}($\Rightarrow$) is taken from \cite[pp.~66-67]{GreenleafBook}.

\fullCref{FreeSubgrpRem}{Olshanskii}, the existence of a nonamenable group with no nonabelian free subgroup, is due to Olshanskii \cite{Olshanskii-ExInvMean}. (In this example, called an ``Olshanskii Monster'' or ``Tarski Monster\zz,'' every proper subgroup of the group is a cyclic group of prime order, so there is obviously no free subgroup.) A much more elementary example has recently been constructed by N.\,Monod \cite{Monod-PiecewiseProj}.

The book of S.\,Wagon \cite{WagonBook} is one of the many places to read about the \thmindex{Banach-Tarski Paradox}Banach-Tarski Paradox (\fullcref{FreeSubgrpRem}{BanachTarski}).

Furstenberg's Lemma \pref{G/amen->Meas(X)} appears in \cite[Thm.~15.1]{Furstenberg-BdThyStochProc}. Another proof can be found in \cite[Prop.~4.3.9, p.~81]{ZimmerBook}. 

\Cref{AmenNoBddHarm} is due to Kaimanovich and Vershik \cite[Thms.~4.2 and~4.4]{KaimanovichVershik-RandWalksBdryEntropy} and (independently) Rosenblatt \cite[Props.~1.2 and 1.9 and Thm.~1.10]{Rosenblatt-ErgMixRandWalks}.

\Cref{Amen<>ConvOp} is due to H.\,Kesten (if $\mu$ is symmetric). See \cite[G.4.4]{BekkaHarpeValette-T} for a proof.

A proof of \cref{SubexpIsAmen} can be found in \cite[Props.~12.5 and~12.5]{PierBook}.

\Cref{AmenCoGrowth} was proved by R.\,I.\,Grigorchuk \cite{Grigorchuk-SymmRandWalks} and J.\,M.\,Cohen \cite{Cohen-cogrowth} (independently). 

\Cref{AmenUnifBddIsUnit} was proved by J.\,Dixmier and M.\,Day in 1950 (independently). 
See \cite{Pisier-UnitAmen} for historical remarks and progress on the converse.
(Another result on the converse is proved in \cite{MonodOzawa-Unitarizable}.)

\Cref{AlmRepsAreNear} is due to D.\,Kazhdan \cite{Kazhdan-EpsilonRep}.

\Cref{AmenBddCoho0,AmenBddCoho0Alln} are due to B.\,E.\,Johnson \cite{Johnson-CohoBanachAlg}. See \cite{Monod-InvBddCoho} (and its many references) for an introduction to bounded cohomology.

\Cref{AmenNoPonzi} appears in \cite[6.17 and 6.17\raise1pt\hbox{$\scriptstyle\frac{1}{2}$}, p.~328]{Gromov-MetricStructures}.

\end{notes}

 %!TEX root = IntroArithGrps.tex

\mychapter{Kazhdan's Property \texorpdfstring{$(T)$}{(T)}}
\label{KazhdanTChap}

\prereqs{Unitary representations (\cref{UnitaryRepSect,InducedRepSect,RepRnSect}).}
	%Quasi-isometries (\cref{QuasiChap}) are mentioned.

Recall that if a Lie group~$H$ is not amenable, then $\LL2(H)$ does not
have almost-invariant vectors \fullcsee{AmenEquiv}{Vectors}. Kazhdan's property~$(T)$ is the much stronger condition that \textbf{no} unitary representation
of~$H$ has almost-invariant vectors (unless it has a vector that is fixed by~$H$).
Thus, in a sense, Kazhdan's property is the antithesis of amenability.

We already know that $\Gamma$ is not amenable (unless it is finite) \csee{GammaNotAmen}. In this chapter, we will see that $\Gamma$ usually has Kazhdan's Property~$(T)$, and we will look at some of the consequences of this.

\section{Definition and basic properties}

Part~\pref{KazhdanTDefn-alminvt} of the following definition is repeated from \cref{AlmInvtVecDefn}, but the second half is new.

\begin{defn} \label{KazhdanTDefn}
Let $H$ be a Lie group. 
\noprelistbreak
\begin{enumerate}
\item \label{KazhdanTDefn-alminvt}
An action of~$H$ on a normed vector space~$\Banach$ has \defit[almost-invariant vector]{almost-invariant} vectors if, for every
compact subset~$C$ of~$H$ and every $\epsilon > 0$, there
is a unit vector
 $v \in \Banach$, such that 
 \begin{equation} \label{KazhdanTDefn-epsCInvt}
 \|c v - v \| < \epsilon \text{\quad for all $c \in C$}
.\end{equation}
(A unit vector satisfying \pref{KazhdanTDefn-epsCInvt} is said to be $(\epsilon,C)$-invariant.)

\item $H$ has \defit[Kazhdan!property T@property~$(T)$]{Kazhdan's property~$(T)$} if every unitary representation of~$H$ that has almost-invariant vectors also has (nonzero) invariant vectors.
\end{enumerate}
We often abbreviate ``Kazhdan's property~$(T)$'' to ``Kazhdan's property\zz.'' Also, a group that has Kazhdan's property is often said to be a \defit[Kazhdan!group]{Kazhdan group}.
 \end{defn}
 
\begin{warn} \label{TnotBanach}
By definition, unitary representations are actions on Hilbert spaces, so Kazhdan's property says nothing at all about actions on other types of topological vector spaces. In particular, there are actions of Kazhdan groups by norm-preserving linear transformations on some Banach spaces that have almost-invariant vectors, without having invariant vectors \csee{BanachNoInvt}. On the other hand, it can be shown that there are no such examples on $\LL{p}$~spaces (with $1 \le p < \infty$).
\end{warn}

\begin{prop} \label{Kazhdan+amenable}
A Lie group is compact if and only if it is amenable and has Kazhdan's
property.
 \end{prop}

\begin{proof} 
\Cref{Cpct->Kazhdan+amenable,Cpct<-Kazhdan+amenable}.
 \end{proof}
 
 \begin{cor} \label{T+amen->finite}
 A discrete group~$\Lambda$ is finite if and only if it is amenable and has Kazhdan's
property.
\end{cor}

\begin{eg}
$\integer^n$ does not have Kazhdan's property, because it is a discrete, amenable group that is not finite.
\end{eg}

\begin{prop} \label{KazhdanEasy} 
If $\Lambda$ is a discrete group with Kazhdan's property, then:
\begin{enumerate}
\item \label{KazhdanEasy-quotient}
every quotient $\Lambda/N$ of~$\Lambda$ has Kazhdan's property,
\item \label{KazhdanEasy-Abel}
 the abelianization $\Lambda/[\Lambda,\Lambda]$ of~$\Lambda$ is finite,
 and
 \item \label{KazhdanEasy-fg}
 $\Lambda$ is finitely generated.
\end{enumerate}
 \end{prop}

\begin{proof}
 For \pref{KazhdanEasy-quotient} and \pref{KazhdanEasy-Abel}, see \Cref{KazhdanEasy-quotientPf,KazhdanEasy-AbelPf}.
 
 \pref{KazhdanEasy-fg}
  Let $\{ \Lambda_n \}$ be the collection of all finitely
generated subgroups of~$\Lambda$.
  We have a unitary representation of~$\Lambda$ on each $\LL2(\Lambda/\Lambda_n)$, 
given by $(\gamma f)(x \Lambda_n) = f(\gamma^{-1} x \Lambda_n)$.
The direct sum of these is a unitary representation on
 $$\Hilbert = \LL2(\Lambda/\Lambda_1) \oplus \LL2(\Lambda/\Lambda_2) \oplus
\cdots. $$

Any compact set $C \subseteq \Lambda$ is finite, so we have
$C \subseteq \Lambda_n$, for some~$n$.
 Then $C$ fixes the base point $p = \Lambda_n/\Lambda_n$ in $\Lambda
/\Lambda_n$,
 so, letting $f = \delta_p$ be a nonzero function in
$\LL2(\Lambda/\Lambda_n)$ that is supported on $\{p\}$,
 we have $\gamma f = f$ for all $\gamma \in C$. Therefore,
$\Hilbert$ has almost-invariant vectors, so there must be
an $H$-invariant vector in~$\Hilbert$.

 So some $\LL2(\Lambda/\Lambda_n)$  has an invariant vector.
 Since $\Lambda$ is transitive on $\Lambda/\Lambda_n$, 
  an invariant function must be constant.
 So a (nonzero) constant function is in $\LL2(\Lambda/\Lambda_n)$,
 which means $\Lambda /\Lambda_n$ is finite.
 Because $\Lambda_n$ is finitely generated, this implies that 
$\Lambda $ is finitely generated.
 \end{proof}

Since the abelianization of any (nontrivial) free group is infinite, we have the following example:

\begin{cor} \label{FreeNotT}
 Free groups do not have Kazhdan's property.
 \end{cor}
 
\begin{rem}
\Cref{KazhdanEasy} can be generalized to groups that are not required to be discrete, if we replace the word ``finite'' with ``compact'' (see \cref{KazhdanQuotient,KazhdanAbelianization,Kazhdan->CpctGen}). This leads to the following definition:
\end{rem}

\begin{defn} \label{CpctGenDefn}
A Lie group $H$ is \defit{compactly generated} if there exists a compact subset that generates~$H$. 
\end{defn}

 \begin{warn} \label{TIsQuotOfFP}
 Although discrete Kazhdan groups are always finitely generated \fullcsee{KazhdanEasy}{fg}, they need not be finitely presented. (In fact, there are uncountably many non-isomorphic discrete groups with Kazhdan's property~$(T)$, and only countably many of them can be finitely presented.) However, it can be shown that every discrete Kazhdan group is a quotient of a finitely presented Kazhdan group.
 \end{warn}

\begin{exercises}

\item \label{BanachNoInvt}
Let $C_0(H)$ be the Banach space of continuous functions on~$H$ that tend to~$0$ at infinity (with the supremum norm). Show:
	\begin{enumerate}
	\item $C_0(H)$ has almost-invariant vectors of norm~$1$, 
	but
	\item $C_0(H)$ does not have $H$-invariant vectors other than~$0$, unless $H$ is compact.
	\end{enumerate}
\hint{Choose a uniformly continuous function $f(h)$ that tends to~$+\infty$ as $h$~leaves compact sets. For large~$n$, the function $h \mapsto n/ \bigl( n + f(h) \bigr)$ is almost invariant.}

\item \label{Cpct->Kazhdan+amenable}
Prove \Cref{Kazhdan+amenable}($\Rightarrow$).
\hint{If $H$ is compact, then almost-invariant vectors are invariant.}
 
\item \label{Cpct<-Kazhdan+amenable}
Prove \Cref{Kazhdan+amenable}($\Leftarrow$).
\hint{Amenability plus Kazhdan's property implies $\LL2(H)$ has an invariant vector.}

\item \label{KazhdanEasy-quotientPf}
Prove \fullcref{KazhdanEasy}{quotient}
\hint{Any representation of $\Lambda/N$ is also a representation of $\Lambda$.}

\item \label{KazhdanQuotient}
Show that if $H$ has Kazhdan's property, and $N$ is a closed, normal subgroup of~$H$, then $H/N$ has Kazhdan's property.

\item \label{KazhdanEasy-AbelPf}
Prove \fullcref{KazhdanEasy}{Abel}.
\hint{$\Lambda/[\Lambda,\Lambda]$ is amenable and has Kazhdan's property.}

\item \label{KazhdanAbelianization}
Show that if $H$ has Kazhdan's property, then $H/[H,H]$ is compact.

\item Show that if $N$ is a closed, normal subgroup of~$H$, such that $N$ and $H/N$ have Kazhdan's property, then $H$ has Kazhdan's property.
\hint{The space of $N$-invariant vectors is $H$-invariant (why?).}

\begin{warn*}
The converse is not true: there are examples in which a normal subgroup of a Kazhdan group is not Kazhdan \csee{SL3xR3HasT}.
\end{warn*}

\item Show that $H_1 \times H_2$ has Kazhdan's property if and only if $H_1$ and~$H_2$ both have Kazhdan's property.

\item Let $(\pi,V)$ be a unitary representation of a Kazhdan group~$H$. Show that almost-invariant vectors in~$V$ are near invariant vectors. More precisely, given $\epsilon > 0$, find a compact subset~$C$ of~$H$ and $\delta > 0$, such that if $v$ is any $(\delta,C)$-invariant vector in~$V$, then there is an invariant vector~$v_0$ in~$V$, such that $\| v - v_0 \| < \epsilon$.
\hint{There are no almost-invariant unit vectors in~$(V^H)^\perp$, the orthogonal complement of the space of invariant vectors.}

\item \label{KazhdanConstEx}
Suppose $S$ is a generating set of a discrete group~$\Lambda$, and $\Lambda$~has Kazhdan's property. Show there exists $\epsilon > 0$, such that if $\pi$ is any unitary representation of~$\Lambda$ that has an $(\epsilon,S)$-invariant vector, then $\pi$ has an invariant vector. 
(The point here is to reverse the quantifiers: the same $\epsilon$ works for every~$\pi$.)
Such an~$\epsilon$ is called a \defit[Kazhdan!constant]{Kazhdan constant} for~$\Lambda$.

\item \label{T->NotHaagerup}
Recall that we say $H$ has the \defit{Haagerup property} if it has a unitary representation, such that there are almost-invariant vectors, and all matrix coefficients decay to~$0$ at~$\infty$.
Show that if $H$ is a noncompact group with Kazhdan's property, then $H$ does not have the Haagerup property.

\item Assume:
	\begin{itemize}
	\item $\varphi \colon H_1 \to H_2$ is a homomorphism with dense image,
	and
	\item $H_1$ has Kazhdan's property.
	\end{itemize}
Show $H_2$ has Kazhdan's property.

\item \label{CpctGen<>fg}
Show that a Lie group~$H$ is compactly generated if and only if $H/H^\circ$ is finitely generated.
\hint{($\Leftarrow$) Since $H^\circ$ is connected, it is generated by any subset with nonempty interior.}

\item \label{Kazhdan->CpctGen}
Show that every Lie group with Kazhdan's property is compactly generated.
\hint{Either adapt the proof of \fullcref{KazhdanEasy}{fg}, or use \fullcref{KazhdanEasy}{fg} together with \cref{KazhdanQuotient,CpctGen<>fg}.}

\item \label{Kazhdan->Expander}
Assume $\Gamma$ has Kazhdan's property~$(T)$, and $S$~is a finite generating set for~$\Gamma$. Show there exists $\epsilon > 0$, such that if $N$~is any finite-index normal subgroup of~$\Gamma$, and $A$~is any subset of $\Gamma/N$, then
	$$ \#(SA \cup A) \ge \min \left\{\, (1+\epsilon) \cdot \#A, \ {\textstyle\frac{1}{2}} |\Gamma / N| \,\right\} .$$
{\smaller (In graph-theoretic terminology, this means the Cayley graphs $\mathrm{Cay}(\Gamma/N_k; S)$ form a family of \defit{expander graphs} if $N_1,N_2,\ldots$ are finite-index normal subgroups, such that $|\Gamma/N_k| \to \infty$.)\par}

\end{exercises}

 \section{Semisimple groups with Kazhdan's property} \label{SSKazhdanSect}
 
 \begin{thm}[(Kazhdan)] \label{SL3RHasT}
 $\SL(3,\real)$ has Kazhdan's property.
 \end{thm}
 
 This theorem is an easy consequence of the following \lcnamecref{RelTForSL2RxR2}, which will be proved in \cref{RelTForSL2RxR2Sect}.
 
 \begin{lem} \label{RelTForSL2RxR2}
 Assume
 \noprelistbreak
 	\begin{itemize}
	\item $\pi$ is a unitary representation of the natural semidirect product 
		$$\SL(2,\real) \ltimes \real^2
		 =	\begin{Smallbmatrix}
			 \upast & \upast & \upast \\
			 \upast & \upast & \upast \\
			0 & 0 & 1 
			\end{Smallbmatrix}
			\subset \SL(3,\real) ,$$
	and 
	\item $\pi$ has almost-invariant vectors.
	\end{itemize}
Then $\pi$ has a nonzero vector that is invariant under the subgroup\/~$\real^2$.
\end{lem}

\begin{terminology*}
Suppose $R$ is a subgroup of a topological group~$H$. The pair $(H,R)$ is said to have \defit[Kazhdan!property T@property~$(T)$!relative]{relative property~$(T)$} if every unitary representation of~$H$ that has almost-invariant vectors must also have an $R$-invariant vector. In this terminology, \cref{RelTForSL2RxR2} states that the pair $\bigl( \SL(2,\real) \ltimes \real^2, \real^2 \bigr)$ has relative property~$(T)$.
\end{terminology*}

\begin{proof}[Proof of \cref{SL3RHasT}]
Let 
	$$ \text{$G = \SL(3,\real) $,
	\ 
	$R = 
	\begin{Smallbmatrix}
	1 & 0 & \upast \\
	0 & 1 & \upast \\
	0 & 0 & 1 
	\end{Smallbmatrix}
	\iso \real^2$,
	\ 
	and
	\ 
	$H = 
	\SL(2,\real) \ltimes R$}
	, $$
and suppose $\pi$ is a unitary representation of~$G$ that has almost-invariant vectors. Then it is obvious that the restriction of~$\pi$ to~$H$ also has almost-invariant vectors \csee{RestrictAlmInvt}, so \cref{RelTForSL2RxR2} implies there is a nonzero vector~$v$ that is fixed by~$R$.
Then the {Moore Ergodicity Theorem} \pref{MautnerPhenom} implies that $v$ is fixed by all of~$G$.
So $\pi$ has a fixed vector (namely,~$v$).
\end{proof}

If $G$ is simple, and $\Rrank G \ge 2$, then $G$ contains a subgroup isogenous to $\SL(2,\real) \ltimes \real^n$, for some~$n$ \ccf{SL2RxRnInG}, so a modification of the above argument shows that $G$ has Kazhdan's property. On the other hand, it is important to know that not all simple Lie groups have the property:

\begin{eg}
 $\SL(2,\real)$ does \emph{not} have Kazhdan's property.
 \end{eg}
 
 \begin{proof}
 Choose a torsion-free lattice~$\Gamma$ in
$\SL(2,\real)$. Then $\Gamma$ is either a surface group or
a nonabelian free group. In either case,
$\Gamma/[\Gamma,\Gamma]$ is infinite, so $\Gamma$ does not
have Kazhdan's property. Therefore, we conclude from
\cref{Kazhdan:G->Gamma} below that $\SL(2,\real)$ does not
have Kazhdan's property.
\end{proof}

\begin{proof}[Alternate proof]
A reader familiar with the unitary representation theory of $\SL(2,\real)$ can easily construct a sequence of representations in the principal series  whose limit is the trivial representation. The 
direct sum of this sequence of representations has almost-invariant vectors.
 \end{proof}

We omit the proof of the following precise characterization of the semisimple groups that have Kazhdan's property:

\begin{thm} \label{WhichGKazhdan}
 $G$ has Kazhdan's property if and only if no simple
factor of~$G$ is isogenous to\/ $\SO(1,n)$ or\/ $\SU(1,n)$.
\end{thm}

\begin{exercises}

\item \label{RestrictAlmInvt}
Assume $\pi$ is a unitary representation of~$H$ that has almost-invariant vectors, and $L$~is a subgroup of~$H$. Show that the restriction of~$\pi$ to~$L$ has almost-invariant vectors. 

\item \label{SL2RxRnInG}
(\emph{Assumes familiarity with real roots})
Assume $G$ is simple. Show $\Rrank G \ge 2$ if and only if some connected subgroup~$L$ of~$G$ is isogenous to $\SL(2,\real)$ and normalizes (but does not centralize) a nontrivial, unipotent subgroup~$U$ of~$G$.
\hint{($\Rightarrow$)~An entire maximal parabolic subgroup of~$G$ normalizes a nontrivial unipotent subgroup.
($\Leftarrow$)~Construct two unipotent subgroups of~$G$ that both contain~$U$, but generate a subgroup that is not unipotent.}

\item
Suppose $H$ is a closed, noncompact subgroup of~$G$, and $G$~is simple.
Show that the pair $(G,H)$ has relative property $(T)$ if and only if $G$ has Kazhdan's property.

\item Suppose $G$ has Kazhdan's property. Show there is a compact subset~$C$ of~$G$ and some $\epsilon > 0$, such that every unitary representation of~$G$ with $(\epsilon,C)$-invariant vectors has invariant vectors.

\end{exercises}

\section{Proof of relative property \texorpdfstring{$(T)$}{(T)}} \label{RelTForSL2RxR2Sect}

In this \lcnamecref{RelTForSL2RxR2Sect}, we prove \cref{RelTForSL2RxR2}, thereby completing the proof that $\SL(3,\real)$ has Kazhdan's property~$(T)$. The argument relies on a decomposition theorem for representations of~$\real^n$.
%, or as a consequence of elementary multivariate Fourier analysis that we  (or, in our terminology,) that will not be needed elsewhere in the book. 

%\begin{notation}
%For convenience, in this \lcnamecref{RelTForSL2RxR2Sect}, we let $R = \real^2$.
%\end{notation}

\begin{proof}[Proof of \cref{RelTForSL2RxR2}]
For convenience, let $H = \SL(2,\real) \ltimes \real^2$. % and $R = \real^2 \normal H$. 
Given a unitary representation $(\pi,\Hilbert)$ of~$H$ that has almost-invariant vectors, we wish to show that some nonzero vector in~$\Hilbert$ is fixed by the subgroup~$\real^2$ of~$H$. In other words, if we let $E$ be the projection-valued measure provided by \cref{ProjValMeas} (for the restriction of~$\pi$ to~$\real^2$), then we wish to show $E \bigl( \{0\} \bigr)$ is nontrivial.

Letting 
	\nindex{$\bddop(\Hilbert)$ = $\{\text{bounded operators on~$\Hilbert$} \}$}%
	$\bddop(\Hilbert)$ 
	be the algebra of bounded linear operators on~$\Hilbert$, and using the fact that $\pi$ has almost-invariant vectors, \cref{ConstructPhiEx} provides a continuous, linear functional $\lambda \colon \bddop(\Hilbert) \to \complex$, such that 
	\begin{itemize}
	\item $\lambda(\Id) = 1$,
	\item $\lambda( E ) \ge 0$ for every orthogonal projection~$E$,
	and
	\item $\lambda$ is bi-invariant under $H$. (More precisely, for all $h_1,h_2 \in H$ and $T \in \bddop(\Hilbert)$ we have $\lambda \bigl( \pi(h_1) \, T  \, \pi(h_2) \bigr) = \lambda(T)$.)
	\end{itemize}

Now, let $\mu$ be the composition of $\lambda$ with~$E$ (that is, let $\mu(A) = \lambda \bigl(E(A) \bigr)$ for $A \subseteq \real^n$), so $\mu$ is a finitely additive probability measure on~$\real^n$ \csee{MuIsFinAdd}. Since $\real^2 \normal H$, there is an action of~$H$ on~$\real^2$ by conjugation. One can show that the probability measure~$\mu$ is invariant under this action \csee{muIsHInvtEx}.

On the other hand, the only $\SL(2,\real)$-invariant, finitely additive probability measure on~$\real^2$ is the point-mass supported at the origin \csee{InvtMeanOnR2}. 
Therefore, we must have $\mu \bigl( \{0\} \bigr) = 1 \neq 0$. Hence, $E \bigl( \{0\} \bigr)$ is nonzero, as desired. 
\end{proof}

\begin{exercises}

\item \label{ConstructPhiEx}
Prove the existence of the linear functional $\lambda \colon \bddop(\Hilbert) \to \complex$ in the proof of \cref{RelTForSL2RxR2}.
\hint{For $T \in \bddop(\Hilbert)$, define $\lambda_n(T) = \langle Tv_n \mid v_n \rangle$, where $\{v_n\}$ is a sequence of unit vectors in~$\Hilbert$, such that $\| \pi(h) v_n - v_n\| \to 0$ for every $h \in H$. Let $\lambda$ be an accumulation point of $\{\lambda_n\}$ in an appropriate weak topology.}

\item \label{MuIsFinAdd}
Let $\mu$ be as defined near the end of the proof of \cref{RelTForSL2RxR2}. Show:
	\begin{enumerate}
	\item $\mu(\real^2) = 1$.
	\item If $A_1$ and~$A_2$ are disjoint Borel subsets of~$\real^2$, then we have $\mu(A_1 \cup A_2) = \mu(A_1) + \mu(A_2)$.
	\item $\mu(A) \ge 0$ for every Borel set $A \subseteq \real^2$.
	\end{enumerate}

\item \label{muIsHInvtEx}
Show the finitely additive measure~$\mu$ in the proof of \cref{RelTForSL2RxR2} is invariant under the action of~$H$ on~$\dual R$.
\hint{Since 
	\begin{align*}
	\int_{\dual R} \tau(r) \, dE(\tau^h)
	&= \int_{\dual R} \tau^{h^{-1}}(r) \, dE(\tau)
	= \int_{\dual R} \tau(h^{-1} r h) \, dE(\tau)
	= \pi(h^{-1} r h)
	\\&= \pi(h^{-1}) \, \pi(r) \, \pi(h)
	= \int_{\dual R} \tau(r) \, \bigl( \pi(h^{-1}) \,  dE(\tau) \, \pi(h) \bigr) ,
	\end{align*}
we have $E(A^h) = \pi(h^{-1}) \,  E(A)  \, \pi(h)$ for $A \subseteq \dual R$.}

\item \label{InvtMeanOnR2}
Show that any $\SL(2,\real)$-invariant, finitely additive probability measure~$\mu$ on~$\real^2$ is supported on $\{(0,0)\}$.
\hint{Let $V = \{\, (x,y) \mid y > |x| \,\}$ and $h = \text{\smaller[2]$\begin{bmatrix} 1 & 2 \\ 0 & 1 \end{bmatrix}$}$. Then $h^iV$ is disjoint from $h^j V$ for $i \neq j \in \integer^+$, so $\mu(V) = 0$. All of $\real^2 \smallsetminus \{(0,0
)\}$ is covered by finitely many sets of the form $hV$ with $h \in \SL(2,\real)$.}

\item \label{SL3xR3HasT}
Show that the natural semidirect product $\SL(3,\real) \ltimes \real^3$ has Kazhdan's property.
\hint{We know $\SL(3,\real)$ has Kazhdan's property, and the proof of \cref{RelTForSL2RxR2} shows that the pair $\bigl( \SL(3,\real) \ltimes \real^3 , \real^3 \bigr)$ has relative property~$(T)$.}

\item Show that the direct product $\SL(3,\real) \times \real^3$ does \emph{not} have Kazhdan's property. (Comparing this with \cref{SL3xR3HasT} shows that, for group extensions, Kazhdan's property may depend not only the groups involved, but also on the details of the particular extension.)

\end{exercises}

\section{Lattices in groups with Kazhdan's property} \label{LattInTSect}

In this \lcnamecref{LattInTSect}, we will use basic properties of induced representations to prove the following important result:

\begin{prop} \label{Kazhdan:G->Gamma}
 If $G$ has Kazhdan's property, then  $\Gamma$ has
Kazhdan's property.
 \end{prop}
  
 Combining this with \cref{WhichGKazhdan}, we obtain:
 
\begin{cor} \label{GammaHasT}
\label{Kazhdanlattice}
 If no simple factor of~$G$ is isogenous to\/ $\SO(1,n)$ or\/
$\SU(1,n)$, then\/  $\Gamma$ has Kazhdan's property.
\end{cor}

By \cref{KazhdanEasy}, this has two important consequences:

\begin{cor} \label{KazhdanlatticeCor}
 If no simple factor of~$G$ is isogenous to\/ $\SO(1,n)$ or\/
$\SU(1,n)$, then 
 \begin{enumerate}
\item \label{KazhdanlatticeCor-fg}
 $\Gamma $ is finitely generated,
 and
 \item \label{KazhdanlatticeCor-noabel}
 $\Gamma/[\Gamma ,\Gamma ]$ is finite.
   \end{enumerate}
\end{cor}

\begin{rems} \label{KazhdanLattRem} \ 
\noprelistbreak
\begin{enumerate}
\item It was pointed out in \cref{GammaFinGen} that \pref{KazhdanlatticeCor-fg} remains true without any assumption on the simple factors of~$G$. In fact, $\Gamma$ is always finitely presented, not merely finitely generated.
\item \label{KazhdanLattRem-abel}
On the other hand, \pref{KazhdanlatticeCor-noabel} is not always true, because lattices in $\SO(1,n)$ and $\SU(1,n)$ can have infinite abelian quotients. (In fact, it is conjectured that every lattice in $\SO(1,n)$ has a finite-index subgroup with an infinite abelian quotient, and this is known to be true when $n = 3$.) The good news is that the Margulis Normal Subgroup Theorem implies these are the only examples (modulo multiplying $G$ by a compact factor) if we make the additional assumption that $\Gamma$ is irreducible (see \cref{GNoAbelianization} or \cref{GammaHasFiniteAbelianization}).
\end{enumerate}
\end{rems}

The proof of \cref{Kazhdan:G->Gamma} uses some machinery from the theory of unitary representations.

\begin{notation}
 Let $(\pi, \Hilbert)$ and $(\sigma,\mathcal{K})$ be unitary representations of a Lie group~$H$. (In our applications, $H$~will be either $G$ or~$\Gamma$.)
 \begin{enumerate}
 
 \item We write 
 	\nindex{$\sigma \le \pi$ means $\sigma$ is a subrepresentation of~$\pi$}%
	$\sigma \le \pi$ if $\sigma$ is (isomorphic to) a \defit{subrepresentation} of~$\pi$. 
This means there exist
	\begin{itemize}
	\item a closed, $H$-invariant subspace $\Hilbert'$ of~$\Hilbert$,
	and
	\item a bijective, linear isometry $T \colon \mathcal{K} \stackrel{\iso}{\longrightarrow} \Hilbert'$,
	\end{itemize}
such that $T \bigl( \sigma(h) \phi \bigr) = \pi(h) \, T(\phi)$, for all $h \in H$ and $\phi \in \mathcal{K}$.

 \item We write 
  	\nindex{$\sigma \wl \pi$ means $\sigma$ is weakly contained in~$\pi$}%
	$\sigma \wl \pi$ if $\sigma$ is \defit{weakly contained} in~$\pi$. This means that, for every compact set $C$ in~$H$, every $\epsilon > 0$,
 and all unit vectors $\phi _1, \ldots,\phi_n  \in \mathcal{K}$,
 there exist unit vectors
$\psi_1,\ldots,\psi_n \in \Hilbert$, such that, for all $h \in C$ and
all $1 \le i,j \le n$, we have
 $$  \bigl| \langle \sigma(h) \phi_i \mid \phi_j \rangle - \langle \pi(h) \psi_i 
\mid  \psi_j  \rangle \bigr| < \epsilon . $$
 \end{enumerate}
\end{notation}

\begin{rems} \label{WeakContainRem}\ 
\noprelistbreak
	\begin{enumerate}
	\item \label{WeakContainRem-weaker}
	It is obvious that if $\sigma \le \pi$, then $\sigma \wl \pi$. 
	
	\item We have:
		\begin{itemize}
		\item $\pi$ has invariant vectors if and only if $\trivrep\le \pi$,
		and
		\item $\pi$ has almost-invariant vectors if and only if $\trivrep \wl \pi$. 
		\end{itemize}
	Therefore, Kazhdan's property asserts the converse to \pref{WeakContainRem-weaker} in the special case where $\sigma = \trivrep$: for all~$\pi$, if $\trivrep \wl \pi$, then $\trivrep \le \pi$.
	\end{enumerate}
\end{rems}

It is not difficult to show that induction preserves weak containment \csee{IndWeakContEx}:

\begin{lem} \label{IndWeakCont}
If $\sigma \wl \pi$, then $\Ind_\Gamma ^G(\sigma ) \wl  \Ind_\Gamma ^G(\pi )$.
\end{lem}

This (easily) implies the main result of this \lcnamecref{LattInTSect}:

\begin{proof}[Proof of \cref{Kazhdan:G->Gamma}]
Suppose a representation $\pi$ of~$\Gamma $ has
almost-invariant vectors.
 Then $\pi \wg \trivrep$, so 
  $$ \Ind_\Gamma ^G(\pi ) \wg \Ind_\Gamma ^G(\trivrep) =
\LL2(G/\Gamma ) \ge \trivrep $$
\csee{Ind1=L2(G/Gamma),1<L2(G/Gamma)}.
 Because $G$ has Kazhdan's property, we conclude that
$\Ind_\Gamma ^G(\pi ) \ge \trivrep$.
 This implies $\pi \ge \trivrep$ \csee{1<Ind->1<rho}, as desired.
 \end{proof}
 
\begin{rem} \label{ExplicitKazhdanConstants}
If $\Gamma$ has Kazhdan's property, and $S$ is any generating set of~$\Gamma$, then there is some $\epsilon > 0$, such that every unitary representation of~$\Gamma$ with an $(\epsilon,S)$-invariant unit vector must have invariant vectors \csee{KazhdanConstEx}. Our proof does not provide any estimate on~$\epsilon$, but, in many cases, including $\Gamma = \SL(n,\integer)$, an explicit value of~$\epsilon$ can be obtained by working directly with the algebraic structure of~$\Gamma$ (rather than using the fact that $\Gamma$ is a lattice).
\end{rem}

\begin{rem} \label{OtherT}
For many years, lattices (and some minor modifications of them) were the only discrete groups known to have Kazhdan's property~$(T)$, but other constructions are now known. In particular:
\noprelistbreak
	\begin{enumerate}
	\item \label{OtherT-random}
	Groups can be defined by generators and relations. It can be shown that if the relations are selected at random (with respect to a certain probability distribution), then the resulting group has Kazhdan's property~$(T)$ with high probability.
	\item \label{OtherT-Alg}
	An algebraic approach that directly proves Kazhdan's property for $\SL(n,\integer)$, without using the fact that it is a lattice, has been generalized to allow some other rings, such as polynomial rings, in the place of~$\integer$. In particular, $\SL \bigl( n, \integer[X_1,\ldots,X_k] \bigr)$ has Kazhdan's property~$(T)$ if $n \ge k + 3$.
	\end{enumerate}
\end{rem}

\begin{rem} \label{TnotQI}
We saw in \cref{AmenQIinvt} that amenability is invariant under quasi-isometry \csee{QuasiIsomDefn}. In contrast, this is not true for Kazhdan's property~$(T)$.
To see this, let 
\noprelistbreak
	\begin{itemize}
	\item $G$ be a simple group with Kazhdan's property~$(T)$, 
	\item $\cover G$ be the universal cover of~$G$,
	\item $\Gamma$ be a cocompact lattice in~$G$,
	and
	\item $\cover\Gamma$ be the inverse image of~$\Gamma$ in~$\cover G$, so $\cover\Gamma$ is a lattice in~$\cover G$.
	\end{itemize}
Then $\cover\Gamma$ has Kazhdan's property~$(T)$ (because $\cover G$ has the property). However, if $G = \Sp(4,\real)$ (or, more generally, if the fundamental group of~$G$ is an infinite cyclic group), then $\cover\Gamma$ is quasi-isometric to $\Gamma \times \integer$, which obviously does not have Kazhdan's property (because its abelianization is infinite).

Here is a brief explanation of why $\cover\Gamma$ is quasi-isometric to $\Gamma \times \integer$. Note that $\cover\Gamma/\integer \iso \Gamma$ yields a $2$-cocycle $\alpha \colon \Gamma \times \Gamma \to \integer$ of group cohomology. 
Since $G/\Gamma$ is compact, it turns out that $\alpha$ can be chosen to be uniformly bounded, as a function on $\Gamma \times \Gamma$. This implies that the extension~$\cover\Gamma$ is quasi-isometric to the extension corresponding to the trivial cocycle. This extension is $\Gamma \times \integer$.
\end{rem}

\begin{exercises}

\item \label{IndWeakContEx}
Prove \cref{IndWeakCont}.

\item \label{1<L2(G/Gamma)}
Show $\trivrep \le \LL2(G/\Gamma)$.

\item \label{1<Ind->1<rho}
Show that if $\pi$ is a unitary representation of~$\Gamma$, and $\trivrep \le \Ind_\Gamma^G(\pi)$, then $\trivrep \le \pi$.

\item \label{GammaT->GT}
Prove the converse of \cref{Kazhdan:G->Gamma}: Show that if $\Gamma$ has Kazhdan's property, then $G$ has Kazhdan's property.
\hint{Any $\Gamma$-invariant vector~$v$ can be averaged over $G/\Gamma$ to obtain a $G$-invariant vector. If $v$ is $\epsilon$-invariant for a compact set whose projection to $G/\Gamma$ has measure $> 1 - \epsilon$, then the average is nonzero.}

\end{exercises}

 \section{Fixed points in Hilbert spaces}

We now describe an important geometric interpretation of Kazhdan's
property. 
%For simplicity, we assume all Hilbert spaces are \emph{real} (that is, the scalars come from~$\real$, rather than~$\complex$).

 \begin{defn} \label{AffIsomDefn}
 Let $\Hilbert$ be a Hilbert space. A bijection $T
\colon \Hilbert \to \Hilbert$ is an \defit{affine isometry}
of~$\Hilbert$ if there exist a unitary operator~$U$ on~$\Hilbert$, and $b \in \Hilbert$,
such that 
	$$ \text{$T(v) = Uv + b$ \  for all $v \in \Hilbert$} .$$
 \end{defn}

\begin{eg} \label{NoFPEg}
Let $w_0$ be a nonzero vector in a Hilbert space~$\Hilbert$. For $t \in \real$, define an affine isometry $\phi^t$ of~$\Hilbert$ by $\phi^t (v) = v + t w_0\mk$; this yields an action of~$\real$ on~$\Hilbert$ by affine isometries. Since $\phi^1(v) = v + w_0 \neq v$, we know that the action has no fixed point. 
\end{eg}

The main theorem of this section shows that  the groups that do not have Kazhdan's property are characterized by the existence of a fixed-point-free action as in \cref{NoFPEg}. However, before stating the result, let us introduce some notation, so that we can also state it in cohomological terms. 

 \begin{defn} \label{HilbertCohoDefn}
 Suppose $(\pi, \Hilbert)$ is a unitary representation of a Lie group~$H$. Define
 \begin{enumerate}
 \item $C(H;\Hilbert) = \{\, \text{continuous functions $f \colon H \to \Hilbert$} \,\}$,
 
 \item $\cocyc1(H; \pi) \,{=}\, 
 \{ f \in C(H;\Hilbert) \mid  \forall g,h \in H, \, f(gh) = f(g) + \pi(g)  f(h) \}$,

 \item $\cobdry1(H; \pi) = 
 \{ f \in C(H;\Hilbert) \mid  \exists v \in \Hilbert, \, \forall h \in H, \, f(h) = v - \pi(h) v  \}$,
 
 \item $\coho1(H; \pi) = \cocyc1(H; \pi) / \cobdry1(H; \pi)$
\csee{B1inZ1}.
 \end{enumerate}
If the representation~$\pi$ on~$\Hilbert$ is clear from the context, we may write $\cocyc1(H; \Hilbert)$, $\cobdry1(H; \Hilbert)$, and $\coho1(H; \Hilbert)$,
instead of $\cocyc1(H; \pi)$, $\cobdry1(H; \pi)$, and $\coho1(H; \pi)$.
 \end{defn}

\begin{thm} \label{T<>FH}
 For a Lie group~$H$, the following are equivalent:
 \noprelistbreak
 \begin{enumerate}
 \item  \label{T<>FH-T}
 $H$ has Kazhdan's property.
 \item \label{T<>FH-FP}
  For every Hilbert space~$\Hilbert$, every continuous action of~$H$ by affine isometries
on~$\Hilbert$  has a fixed point.
\item   \label{T<>FH-H1=0}
$\coho1(H; \pi) = 0$, for every unitary representation~$\pi$ of~$H$.
 \end{enumerate}
 \end{thm}

\begin{proof}[\bf \mathversion{bold} Proof of \pref{T<>FH-FP}$\implies$\pref{T<>FH-H1=0}]
Given $f \in \cocyc1(H;\pi)$, define an action of~$f$ on~$\Hilbert$ via affine isometries by defining 
	$$ \text{$hv = \pi(h) v + f(h)$ for $h \in H$ and $v \in \Hilbert$} $$
\csee{Z1<>act}. By assumption, this action must have a fixed point~$v_0$. For all $h \in H$, we have
	$ v_0 = h v_0 = \pi(h) v_0 + f(h) $, 
	so $ f(h) = v_0 - \pi(h) v_0$.
Therefore $f \in \cobdry1(H;\pi)$. Since $f$ is an arbitrary element of $\cocyc1(H;\pi)$, this implies $\coho1(H; \pi) = 0$.
\end{proof}

\begin{proof}[\bf \mathversion{bold} Proof of \pref{T<>FH-H1=0}$\implies$\pref{T<>FH-T}]
We prove the contrapositive: assume $H$ does not have Kazhdan's property.
This means a unitary representation of~$H$ on some Hilbert space~$\Hilbert$ has almost-invariant vectors, but does not have invariant vectors. We claim $\coho1(H; \pi^\infty) \neq 0$, where $\pi^\infty$ is the obvious diagonal action of~$H$ on the Hilbert space $\Hilbert^\infty = \Hilbert \oplus \Hilbert \oplus \cdots$.

Choose an increasing chain $C_1 \subseteq C_2 \subseteq \cdots$ of compact subsets of~$G$, such that $G = \bigcup_n C_n$. For each~$n$, since $\Hilbert$ has almost-invariant vectors, there exists a unit vector $v_n \in \Hilbert$, such that
	$$ \text{$\displaystyle \| c v_n - v_n \| < \frac{1}{2^n}$ for all $c \in C_n$} . $$
Now, define $f \colon H \to \Hilbert^\infty$ by 
	$$ f(h)_n = n\bigl( hv_n - v_n \bigr) $$
\csee{HinftyFunc}, so $f \in \cocyc1(H; \pi^\infty)$ \csee{HinftyCocyc}.
However, it is easy to see that $f$ is an unbounded function on~$H$ \csee{HinftyUnbdd}, so $f \notin \cobdry1(H; \Hilbert^\infty)$ \csee{CobdryIsBdd}. Therefore $f$ represents a nonzero cohomology class in $\coho1(H; \pi^\infty)$. 
\end{proof}

\begin{proof}[\bf \mathversion{bold} Alternate proof of \pref{T<>FH-H1=0}$\implies$\pref{T<>FH-T}]
Assume the unitary representation~$\pi$ has no invariant vectors. (We wish to show this implies there are no almost-invariant vectors.)
Define a linear map 
	$$ \text{$F \colon \Hilbert \to \cocyc1(H;\pi)$ by $ F_v(h) = \pi(h)v - v$} .$$
Assume, for simplicity, that $H$ is compactly generated \csee{AltFH->T-NotCpctGen}, so some compact, symmetric set~$C$ generates~$H$. By enlarging~$C$, we may assume $C$~has nonempty interior. 
Then the supremum norm on~$C$ turns $\cocyc1(H;\pi)$ into a Banach space \csee{CocycIsBanach}, and the map~$F$ is continuous in this topology \fullcsee{AltFH->T-FEx}{cont}. 

Since there are no invariant vectors in~$\Hilbert$, we know that $F$ is injective \fullcsee{AltFH->T-FEx}{inj}.
Also, the image of~$F$ is obviously $\cobdry1(H;\pi)$. Since $\coho1(H;\pi) = 0$, this means that $F$ is surjective. Therefore, $F$~is a bijection.
So the Open Mapping Theorem (\fullref{OpenMappingThm}{bij}) provides a constant $\epsilon > 0$, such that $\|F_v\| > \epsilon$ for every unit vector~$v$. This means there is some $h \in C$, such that $\| \pi(h)v - v \| > \epsilon$, so $v$ is not $(C,\epsilon)$-invariant. Therefore, there are no almost-invariant vectors.
\end{proof}

\begin{proof}[\bf \mathversion{bold} Sketch of proof of \pref{T<>FH-T}$\implies$\pref{T<>FH-FP}]
We postpone this proof to \cref{PosDefFuncSect}, where functions of positive type are introduced. They yield an embedding of~$\Hilbert$ in the unit sphere of a (larger) Hilbert space~$\widehat\Hilbert$. This embedding is nonlinear and non-isometric, but there is a unitary representation~$\widehat\pi$ on~$\widehat\Hilbert$ for which the embedding is equivariant. Kazhdan's property provides an invariant vector in~$\widehat\Hilbert$, and this pulls back to a fixed point in~$\Hilbert$.
See \cref{PosDefFuncSect} for more details.
\end{proof}

\begin{rem}
If $H$ satisfies \pref{T<>FH-FP} of \cref{T<>FH}, it is said to have ``\term[property FH@property $(FH)$]{property $(FH)$}'' (because it has \underline{\textbf{F}}ixed points on \underline{\textbf{H}}ilbert spaces).
\end{rem}

In \cref{HilbertCohoDefn}, the subspace $\cobdry1(H; \pi)$ may fail to be closed \csee{AlmInvt->B1NotClosed}. In this case, the quotient space $\coho1(H; \pi)$ does not have a good topology. Fortunately, it can be shown that \cref{T<>FH} remains valid even if we replace $\cobdry1(H;\pi)$ with its closure:

\begin{defn}
In the notation of \cref{HilbertCohoDefn}, let:
	\begin{enumerate}
	\item $\overline{\cobdry1(H; \pi)}$ be the closure of $\cobdry1(H; \pi)$ in $\cocyc1(H;\pi)$,
	and
	\item $\redcoho1(H; \pi) = \cocyc1(H; \pi) / \overline{\cobdry1(H; \pi)}$. This is called the \emph{reduced} $1$st~cohomology.
	\end{enumerate}
\end{defn}

The following result requires the technical condition that $H$ is compactly generated \csee{CpctGenDefn,RedCohoCpctGenEx}.

\begin{thm} \label{T<>H1bar=0}
A compactly generated Lie group~$H$ has Kazhdan's property if and only if\/
 $\redcoho1(H; \pi) = 0$, for every unitary
representation~$\pi$ of~$H$.
\end{thm}

Because reduced cohomology behaves well with respect to the direct integral decomposition of a unitary representation
(although the unreduced cohomology does not), this theorem implies that it suffices to consider only the irreducible representations of~$H$:

\begin{cor} \label{T<>H1irred=0}
A compactly generated Lie group~$H$ has Kazhdan's property if and only if\/
 $\redcoho1(H; \pi) = 0$, for every \textbf{irreducible} unitary representation~$\pi$ of~$H$.
\end{cor}

\begin{rem} \label{T->FPonOther}
We have seen that a group with Kazhdan's property has bounded orbits whenever it acts isometrically on a Hilbert space. The same conclusion has been proved for isometric actions on some other spaces, including real hyperbolic $n$-space~$\hyperbolic^n$, complex hyperbolic $n$-space~$\hyperbolic_{\complex}^n$, and all ``median spaces'' (including all $\real$-trees). (In many cases, the existence of a bounded orbit implies the existence of a fixed point.)
See \cref{WatataniThm} for an example.
\end{rem}

\begin{exercises}

\item \label{AffineIsomAxioms}
Let $T \colon \Hilbert \to \Hilbert$. Show that if $T$ is an affine isometry, then
 \begin{enumerate}
 \item $T(v - w) = T(v) - T(w) + T(0)$,
 and
 \item $\|T(v)-T(w)\| = \|v-w\|$,
 \end{enumerate}
 for all $v,w \in \Hilbert$. 

\item \label{AffIsom<>}
Prove the converse of \cref{AffineIsomAxioms}.

\item \label{B1inZ1}
 In the notation of \cref{HilbertCohoDefn}, show that
$\cobdry1(H;\pi) \subseteq \cocyc1(H;\pi)$ (so the quotient $\cocyc1(H; \pi) / \cobdry1(H; \pi)$ is defined).

\item \label{CobdryIsBdd}
Suppose $f \in \cobdry1(H;\pi)$ so $f \colon H \to \Hilbert$. Show $f$ is bounded.

\item \label{Z1<>act}
Suppose 
	\begin{itemize}
	\item $(\pi,\Hilbert)$ is a unitary representation of~$H$,
	and
	\item $\tau \colon H \to \Hilbert$.
	\end{itemize}
For $h \in H$ and $v \in \Hilbert$, let $\alpha(h)v = \pi(h)\, v + \tau(h)$, so $\alpha(h)$ is an affine isometry of~$\Hilbert$. Show that $\alpha$ defines a continuous action of~$H$ on~$\Hilbert$ if and only if $\tau \in \cocyc1(H;\pi)$ and $\tau$ is continuous.

\item \label{AffineMustFromRep}
Suppose $H$ acts continuously by affine isometries on the Hilbert space~$\Hilbert$. Show there is a unitary representation~$\pi$ of~$H$ on~$\Hilbert$, and some $\tau \in \cocyc1(H;\pi)$, such that $hv =  \pi(h)\, v + \tau(h)$ for every $h \in H$ and $v \in \Hilbert$. 

 \item \label{FP<>SomeOrbitBdd<>EveryOrbitBddEx}
  Suppose $H$ acts continuously by affine isometries
on the Hilbert space~$\Hilbert$. Show the following are equivalent:
	\begin{enumerate}
	\item $H$ has a fixed point in~$\Hilbert$.
	\item The orbit $Hv$ of each vector~$v$ in~$\Hilbert$ is a bounded subset of~$\Hilbert$.
	\item The orbit $Hv$ of some vector~$v$ in~$\Hilbert$ is a bounded subset of~$\Hilbert$.
	\end{enumerate}
\hint{You may use (without proof) the fact that every nonempty, bounded subset~$X$ of a Hilbert space has a unique circumcenter. By definition, the \defit{circumcenter} is a point~$c$, such that, for some $r > 0$, the set $X$ is contained in the closed ball of radius~$r$ centered at~$c$, but $X$~is not contained in any ball of radius $< r$ (centered at any point).}

\item \label{FP->H1=0Ex}
Prove directly that \fullref{T<>FH}{H1=0}$\implies$\fullref{T<>FH}{FP}, without using Kazhdan's property.
\hint{For each $h \in H$, there is a unique unitary operator $\pi(h)$, such that we have $hv = \pi(h)v + h(0)$ for all $v \in \Hilbert$.
Fix $v \in \Hilbert$ and define $f \in \cocyc1(H;\pi)$ by $f(h) = hv - v$. If $f \in \cobdry1(H;\pi)$, then $H$ has a fixed point.}

\item \label{HinftyFunc}
In the notation of the proof of \ref{T<>FH}($\ref{T<>FH-H1=0}\Rightarrow\ref{T<>FH-T}$), show $f \colon H \to \Hilbert^\infty$.
\hint{For each $h \in H$, show the sequence $\bigl\{ \|f(h)_n\| \bigr\}$ is square-summable.}

\item \label{HinftyCocyc}
In the notation of the proof of  \ref{T<>FH}($\ref{T<>FH-H1=0}\Rightarrow\ref{T<>FH-T}$), show $f \in \coho1(H; \pi^\infty)$.

\item \label{HinftyUnbdd}
In the notation of the proof of \ref{T<>FH}($\ref{T<>FH-H1=0}\Rightarrow\ref{T<>FH-T}$), show $f$ is unbounded.
\hint{You may use (without proof) the fact that every nonempty, bounded subset of a Hilbert space has a unique circumcenter, as in \cref{FP<>SomeOrbitBdd<>EveryOrbitBddEx}.}

\item \label{CocycIsBanach}
Assume $C$ is a compact, symmetric set that generates~$H$, and has nonempty interior. For each $f \in \cocyc1(H;\Hilbert)$, let $\xi(f)$ be the restriction of~$f$ to~$C$. Show that $\xi$ is a bijection from $\cocyc1(H;\Hilbert)$ onto a closed subspace of the Banach space of continuous functions from~$C$ to~$\Hilbert$.

\item \label{AltFH->T-FEx}
In the notation of the alternate proof of \ref{T<>FH}($\ref{T<>FH-H1=0}\Rightarrow\ref{T<>FH-T}$), show:
	\begin{enumerate}
	\item  \label{AltFH->T-FEx-cont}
	$F$ is continuous.
	\item  \label{AltFH->T-FEx-inj}
	$F$ is injective.
	\end{enumerate}
\hint{\pref{AltFH->T-FEx-inj}~If $F_v = F_w$, then what is $\pi(h)(v - w)$?}

\item \label{AltFH->T-NotCpctGen}
Remove the assumption that $H$ is compactly generated from the alternate proof of \pref{T<>FH-H1=0}$\implies$\pref{T<>FH-T}.
\hint{The topology of uniform convergence on compact sets makes $\cocyc1(H;\pi)$ into a Fréchet space.}

\item Assume
 \begin{itemize}
 \item $\Gamma$ has Kazhdan's property~$T$,
 \item $V$ is a vector space,
 \item $\Hilbert$ is a Hilbert space that is contained in~$V$,
 \item $v \in V$,
 and
 \item $\sigma \colon \Gamma \to \GL(V)$ is any homomorphism, such
that 
 \begin{itemize}
 \item the restriction $\sigma(\gamma)|_{\Hilbert}$ is unitary, for
every $\gamma \in \Gamma$, 
 and
 \item $\Hilbert + v$ is $\sigma(\Gamma)$-invariant.
 \end{itemize}
 \end{itemize}
 Show $\sigma(\Gamma)$ has a fixed point in $\Hilbert + v$.
 \hint{\fullCref{T<>FH}{FP}.}
 
 \item \label{AlmInvt->B1NotClosed}
 Show that if $\pi$ has almost-invariant vectors, then $\cobdry1(H;\pi)$ is not closed in $\cocyc1(H;\pi)$.
 \hint{See the alternate proof of \cref{T<>FH}(\ref{T<>FH-H1=0}$\implies$\ref{T<>FH-T}).}
 
 \item \label{RedCohoCpctGenEx}
 Show the assumption that $H$ is compactly generated cannot be removed from the statement of \cref{T<>H1bar=0}.
 \hint{Let $H$ be an infinite, discrete group, such that every finitely generated subgroup of~$H$ is finite.}
 
 \item \label{TiffH1(irred)}
 Show $H$ has Kazhdan's property if and only if 
 $\coho1(H; \pi) = 0$, for every \textbf{irreducible} unitary
representation~$\pi$ of~$H$.
\hint{You may assume \cref{T<>H1bar=0}.}

\begin{defn*}
A \defit{tree} is a contractible, $1$-dimensional simplicial complex.
\end{defn*}

\item  \label{WatataniThm}
(Watatani)
Suppose 
	\begin{itemize}
	\item $\Lambda$ is a discrete group that has Kazhdan's property, 
	and 
	\item acts by isometries on a tree~$T$. 
	\end{itemize}
Show $\Lambda$ has a fixed point in~$T$ (without assuming \cref{T->FPonOther}). 
\hint{Fix an orientation of~$T$, and fix a vertex $v$ in~$T$. For each $\lambda \in \Lambda$, the geodesic path in~$T$ from $v$ to $\lambda(v)$ can be represented by a $\{0, \pm1\}$-valued function~$P_\lambda$ on the set~$E$ of edges of~$T$. Verify that $\lambda \mapsto P_\lambda$ is in $\cocyc1\bigl(\Lambda; \LL2(E) \bigr)$, and conclude that the orbit of~$v$ is bounded.}

\item \label{SO1nNotKazhdan}
Show $\SO(1,n)$ and $\SU(1,n)$ do not have Kazhdan's property.
\hint{You may assume the facts stated in \cref{T->FPonOther}.}

\item It is straightforward to verify that all of the results in this \lcnamecref{KazhdanTChap} remain valid if we require $\Hilbert$ to be a \textbf{real} Hilbert space (instead of a Hilbert space over~$\complex$), as in \cref{RealHilbertAssump} below. In this setting, there is no need to restrict attention to \emph{affine} isometries in the statement of \fullcref{T<>FH}{FP}, because \emph{all} isometries are affine:

Let $\Hilbert$ be a real Hilbert space, and let $\varphi \colon \Hilbert \to \Hilbert$ be any distance-preserving bijection (so $\| \varphi(v) - \varphi(w) \| = \| v - w \|$ for all $v,w \in \Hilbert$). Show that $\varphi$ is an affine isometry.
\hint{The main problem is to show that if $\varphi(0) = 0$, then $\varphi$ is $\real$-linear. This is well known (and easy to prove) when $\Hilbert = \real^2$. The general case follows from this.}

\end{exercises}

\section{Functions on\texorpdfstring{~$H$}{ H} that are of positive type} \label{PosDefFuncSect}

This section completes the proof of \cref{T<>FH}, by showing that affine isometric actions of Kazhdan groups on Hilbert spaces always have fixed points. For this purpose, we develop some of the basic theory of functions of positive type. 

\begin{assump} \label{RealHilbertAssump}
To simplify some details, Hilbert spaces in this section are assumed to be real, rather than complex. (That is, the field of scalars is~$\real$, rather than~$\complex$.)
\end{assump}

\begin{defn} \label{PosTypeDefn} \  
\noprelistbreak
	\begin{enumerate}
	\item Let $A$ be an $n \times n$ real symmetric matrix.
		\noprelistbreak
		\begin{enumerate}
		\item \label{PosTypeDefn-Pos}
		$A$ is of \defit[positive!type]{positive type} if $\langle Av \mid v \rangle \ge 0$ for all $v \in \real^n$. Equivalently, this means all of the eigenvalues of~$A$ are $\ge 0$ \csee{PosTypeIffInnerProdEx}.
		\item  \label{PosTypeDefn-Cond}
		$A$ is \defit[positive!type!conditionally]{conditionally of positive type} if 
			\begin{enumerate}
			\item $\langle Av \mid v \rangle \ge 0$ for all $v = (v_1,\ldots,v_n) \in \real^n$, such that we have $v_1 + \cdots + v_n = 0$,
			and
			\item  \label{PosTypeDefn-Cond-0}
			all the diagonal entries of~$A$ are~$0$.
			\end{enumerate}
		(The word ``conditionally'' refers to the fact that the inequality on $\langle Av \mid v \rangle$ is only required to be satisfied when a particular condition is satisfied, namely, when the sum of the coordinates of~$v$ is~$0$.)
		\end{enumerate}
	\item A continuous, real-valued function~$\varphi$ on a topological group~$H$ is said to be of \defit[positive!type]{positive type} (or \defit[positive!type!conditionally]{conditionally of positive type}, respectively) if, for all~$n$ and all $h_1,\ldots,h_n \in H$, the matrix $\bigl( \varphi(h_i^{-1} h_j ) \bigr)$ is a symmetric matrix of the said type.
	\end{enumerate}
\end{defn}

\begin{warn}
A function that is of positive type is almost never conditionally of positive type. This is because a  matrix satisfying \pref{PosTypeDefn-Pos} of \cref{PosTypeDefn}  will almost never satisfy \pref{PosTypeDefn-Cond-0} \csee{PosAndCond->0}.
\end{warn}

\begin{terminology} \label{PosDefTerm}
Functions of positive type are often called \defit[positive!definite]{positive definite} or \defit[positive!semi-definite]{positive semi-definite}.
\end{terminology}

Such functions arise naturally from actions of~$H$ on Hilbert spaces:

\begin{lem} \label{PosTypeFromAction}
Suppose
\noprelistbreak
	\begin{itemize}
	\item $H$ is a topological group,
	\item $H$ acts continuously by affine isometries on a Hilbert space~$\Hilbert$,
	\item $v \in \Hilbert$,
	\item $\varphi \colon H \to \real$ is defined by $\varphi(h) = -\| hv - v \|^2$ for $h \in H$,
	and
	\item $\psi \colon H \to \real$ is defined by $\psi(h) = \langle hv \mid v \rangle$ for $h \in H$.
	\end{itemize}
Then:
	\begin{enumerate}
	\item  \label{PosTypeFromAction-Cond}
	$\varphi$ is conditionally of positive type,
	and
	\item  \label{PosTypeFromAction-Pos}
	$\psi$ is of positive type if $h(0) = 0$ for all $h \in H$.
	\end{enumerate}
\end{lem} 

\begin{proof}
\Cref{CondPosTypeFromActionEx,PosTypeFromActionEx}.
\end{proof}

Conversely, the following result shows that all functions of positive type arise from this construction. 
(The ``GNS'' in its name stands for Gelfand, Naimark, and Segal.) 

\begin{prop}[(``GNS construction'')] \label{GNS}
If $f \colon H \to \real$ is of positive type, then there exist
\noprelistbreak
	\begin{itemize}
	\item a continuous action of~$H$ by linear isometries on a Hilbert space~$\Hilbert$ {\rm(}so $h(0) = 0$ for all~$h${\rm)},
	and
	\item $v \in \Hilbert$,
	\end{itemize}
such that $f(h) = \langle hv \mid v \rangle$ for all $h \in H$.
\end{prop}

\begin{proof}
Let $\real[H]$ be the vector space of functions on~$H$ with finite support. Since the set of delta functions $\{\,\delta_h \mid h \in H \,\}$ is a basis, there is a unique bilinear form on $\real[H]$, such that 
	$$ \text{$\langle \delta_{h_1} \mid  \delta_{h_2} \rangle = f(h_1^{-1} h_2)$ \ for all $h_1,h_2$} .$$
Since $f$ is of positive type, this form is symmetric and satisfies the inequality $\langle w \mid w \rangle \ge 0$ for all~$w$. Let $Z$ be the \defit[radical!of a bilinear form]{radical} of the form, which means
	$$ Z = \{\, z \in \real[H] \mid \langle z \mid z \rangle = 0 \,\} ,$$
so $\langle \, \mid \, \rangle$ factors through to a well-defined positive-definite, symmetric bilinear form on the quotient $\real[H]/Z$. This makes the quotient into a pre-Hilbert space; let $\Hilbert$ be its completion, which is a Hilbert space, and let $v$ be the image of~$\delta_e$ in~$\Hilbert$. 

The group~$H$ acts by translation on $\real[H]$, and it is easy to verify that the action is continuous, and preserves the bilinear form \csee{GNSPf-ActCont&IsomEx}. Therefore, the action extends to a unitary representation of~$H$ on~$\Hilbert$. Furthermore, for any $h \in H$, we have
	$$ 
	f(h)
	= f(e \cdot h)
	= \langle  \delta_e \mid \delta_h \rangle
	= \langle  \delta_e \mid h\delta_e \rangle
	=  \langle v \mid hv \rangle
	= \langle hv \mid v \rangle
	, $$
as desired.
\end{proof}

We will also use the following important relationship between the two concepts:

\begin{lem}[(Schoenberg's Lemma)] \label{Schoenberg}
If $\varphi$ is conditionally of positive type, then $e^\varphi$ is of positive type.
\end{lem}

\begin{proof}
A function $\kappa \colon H \times H \to \real$ is said to be a \defit{kernel of positive type} if the matrix $\bigl( \kappa(h_i,h_j) \bigr)$ is a symmetric matrix of positive type, for all~$n$ and all $h_1,\ldots,h_n \in H$.

Define $\kappa \colon H \times H \to \real$ by
	$$ \kappa(g,h) = \varphi(g^{-1} h) - \varphi(g) - \varphi(h) .$$
Then:
	\begin{itemize}
	\item $\kappa$ is a kernel of positive type \csee{SchoenbergEx-K},
	\item so $e^\kappa$ is a kernel of positive type \csee{SchoenbergEx-exp},
	\item and $e^{\varphi(g)} e^{\varphi(h)}$ is a kernel of positive type \csee{SchoenbergEx-func},
	\item so the product $e^{\kappa(g,h)} \bigl( e^{\varphi(g)} e^{\varphi(h)} \bigr)$ is a kernel of positive type \csee{SchoenbergEx-prod}.
	\end{itemize}
This product is $e^{\varphi(g^{-1} h)}$, so $e^\varphi$ is a function of positive type.
\end{proof}

With these tools, it is not difficult to show that affine isometric actions of Kazhdan groups on Hilbert spaces always have fixed points:

\begin{proof}[\bf \mathversion{bold} Proof of \cref{T<>FH} ($\ref{T<>FH-T} \Rightarrow\ref{T<>FH-FP}$)]
Let $\alpha$ be the given action of~$H$ on~$\Hilbert$ by affine isometries, and let $\pi$ be the corresponding unitary representation \csee{AffineMustFromRep}. Therefore, we have
	$$ \text{$\alpha(h) v = \pi(h)v + \tau(h)$ for $h \in H$ and $v \in \Hilbert$,} $$
where $\tau \in \cocyc1(H;\pi)$.

Let $\widehat H = \Hilbert \rtimes H$ be the semidirect product of (the additive group of)~$\Hilbert$ with~$H$, where $H$ acts on~$\Hilbert$ via~$\pi$. This means the elements of~$\widehat H$ are the ordered pairs $(v,h)$, and, for $v_1,v_2 \in \Hilbert$ and $h_1,h_2 \in H$, we have
	$$ (v_1,h_1) \cdot  (v_2,h_2) = (v_1 + \pi(h_1) v_2, h_1 h_2 ) .$$
This semidirect product is a topological group, so we can apply the above theory of functions of positive type to it.
Define a continuous action~$\widehat\alpha$ of $\widehat H$ on~$\Hilbert$ by
	\begin{align} \label{T<>FH-alphahatdefn}
	\widehat\alpha(v,h)w = \alpha(h)w + v
	\end{align}
\csee{alphahatIsAction}, and define
	$$ \text{$\widehat\varphi \colon \widehat H \to \real$ by $\widehat\varphi(v,h) = -\| \widehat\alpha(v,h)(0) \|^2$.} $$
Since $\widehat\alpha(v,h)$ is an affine isometry for every $v$ and~$h$, we know $\widehat\varphi$ is conditionally of positive type \fullcsee{PosTypeFromAction}{Cond}. Therefore $e^{\widehat\varphi}$ is of positive type \csee{Schoenberg}. Hence, the GNS construction \pref{GNS} provides a unitary representation~$\widehat\pi$ of~$\widehat H$ on a Hilbert space~$\widehat\Hilbert$ and some $\hat v \in \widehat \Hilbert$, such that 
	\begin{align} \label{T->FH-GNSEquality}
	\text{$\langle \widehat\pi(v,h) \hat v \mid \hat v \rangle = e^{\widehat\varphi(v,h)}$ for all $v \in \Hilbert$ and $h \in H$.} 
	\end{align}
We now define
	$$ \text{$\Phi \colon \Hilbert \to \widehat\Hilbert$ by
	$ \Phi(v) = \widehat\pi(v,e) \hat v $} .$$
We have
	\begin{align} \label{T<>FH-PhiIsEqui}
	\text{$\Phi \bigl( \alpha(h)v \bigr) = \widehat\pi( 0,h ) \, \Phi(v)$ for $h \in H$ and $v \in \Hilbert$} 
	\end{align}
\csee{T<>FH-PhiIsEquiEx}, so $\Phi$ converts the affine action of~$H$ on~$\Hilbert$ to a linear action on~$\widehat\Hilbert$.
Since the linear span of $\Phi(\Hilbert)$ contains~$\hat v$ and is invariant under $\widehat\pi \bigl( \widehat H \bigr)$  \csee{SpanIsHInvt}, there is no harm in assuming that its closure is all of~$\widehat\Hilbert$.

It is clear from the definition of $\widehat\varphi$ that $\widehat\varphi(0,e) = 0$, so we know that $\hat v$ is a unit vector. Therefore
	\begin{align}
	\| \widehat\pi(v,h) \hat v - \hat v\|^2 
	&= \bigl\langle \, \widehat\pi(v,h) \hat v - \hat v, \ \widehat\pi(v,h) \hat v - \hat v \,\bigr\rangle 
	\notag
	\\&= 2 \bigl( 1 - \langle \widehat\pi(v,h) \hat v \mid \hat v \rangle \bigr)
	\label{T<>FHPf-AlmInvt}
	\\& = 2 \bigl( 1 - e^{\widehat\varphi(v,h)} \bigr) 
	\notag
	 .\end{align}

Since $H$ has Kazhdan's property, there is a compact subset~$C$ of~$H$ and some $\epsilon > 0$, such that every unitary representation of~$H$ that has a $(C,\epsilon)$-invariant vector must have an invariant vector. There is no harm in multiplying the norm on~$\Hilbert$ by a small positive scalar, so we may assume $\widehat\varphi(0,h)$ is as close to~$0$ as we like, for all $h \in C$. Then \pref{T<>FHPf-AlmInvt} tells us that $\hat v$ is $(C,\epsilon)$-invariant, so $\widehat\Hilbert$ must have a nonzero $H$-invariant vector~$\hat v_0$.

Suppose the affine action~$\alpha$ does not have any fixed points. (This will lead to a contradiction.) This implies that every $H$-orbit on~$\Hilbert$ is unbounded \csee{FP<>SomeOrbitBdd<>EveryOrbitBddEx}. Hence, for any fixed $v \in \Hilbert$, and all $h \in H$, we have
	\begin{align*}
	\langle \Phi(v) \mid \hat v_0 \rangle
	&= \langle \Phi(v) \mid \widehat\pi(0,h^{-1}) \hat v_0 \rangle
	&& \text{($\hat v_0$ is $H$-invariant)}
	\\&= \langle  \widehat\pi(0,h) \, \Phi(v)  \mid \hat v_0 \rangle
	&& \text{($\pi$ is unitary)}
%	\\&= \langle  \widehat\pi \bigl( (0,h) \cdot(v,e) \bigr) \hat v \mid \hat v_0 \rangle
%	\\&= \langle  \widehat\pi \bigl( \pi(h) v,h \bigr) \hat v \mid \hat v_0 \rangle
	\\&= \langle  \Phi( \alpha(h)v \bigr) \mid \hat v_0 \rangle
	&& \text{(\ref{T<>FH-PhiIsEqui})}
	\\&\to 0 \text{\quad as $\|\alpha(h) v\| \to \infty$}
	&& \text{(\cref{T->FH-WeaklyToZeroEx})}
	.\end{align*}
So $\hat v_0$ is orthogonal to the linear span of $\Phi(\Hilbert)$, which is dense in~$\widehat\Hilbert$. Therefore $\hat v_0 = 0$. This is a contradiction.
\end{proof}

\begin{exercises}

\item \label{PosTypeIffInnerProdEx}
Let $A$ be a real symmetric matrix.
Show $A$ is of positive type if and only if all of the eigenvalues of~$A$ are $\ge 0$.
\hint{$A$ is diagonalizable by an orthogonal matrix.}

\item \label{PosAndCond->0}
Suppose $A$ is an $n \times n$ real symmetric matrix that is of positive type and is also conditionally of positive type. Show $A = 0$.
\hint{What does \fullref{PosTypeDefn}{Cond-0} say about the trace of~$A$?}

\item \label{alphahatIsAction}
In the notation of the proof of \cref{T<>FH}, show
	$$ \widehat\alpha(v_1,h_1) \cdot \widehat\alpha(v_2,h_2) 
	= \widehat\alpha \bigl( (v_1,h_1) \cdot (v_2,h_2) \bigr) $$
for all $v_1,v_2 \in \Hilbert$ and $h_1,h_2 \in H$.

%\item 
%In the notation of the proof of \cref{T<>FH}, show
%	$$ \widehat\pi(e,h)\hat v  = \widehat\pi \bigl( \tau(h), e \bigr) \hat v .$$

\item \label{SpanIsHInvt}
In the notation of the proof of \cref{T<>FH}, show that the linear span of $\Phi(\Hilbert)$ is invariant under $\widehat\pi \bigl( \widehat H \bigr)$.
\hint{See~\pref{T<>FH-PhiIsEqui}.}

\item \label{T->FH-WeaklyToZeroEx}
In the notation of the proof of \cref{T<>FH}, show that if $v \in \Hilbert$, and $\{w_n\}$ is a sequence in~$\Hilbert$, such that $\| w_n \| \to \infty$, then
	$\Phi( w_n ) \to 0$ weakly.
\hint{If $\widehat w \in \Phi(\Hilbert)$, then \pref{T->FH-GNSEquality} implies $\bigl\langle \Phi( w_n ) \mid \widehat w \bigr\rangle \to 0$.% We may assume the linear span of $\Phi(\Hilbert)$ is dense in~$\widetilde\Hilbert$.
}

\item \label{CondPosTypeFromActionEx}
Prove \fullcref{PosTypeFromAction}{Cond}.
\hint{If $\sum_i t_i = 0$, then
$\sum_{i,j} t_i \, t_j \, \varphi( h_i^{-1} h_j v ) = 2 \, \bigl\| \sum_i t_i h_i v \bigr\|^2$.}

\item \label{PosTypeFromActionEx}
Prove \fullcref{PosTypeFromAction}{Pos}.
\hint{$\sum_{i,j} t_i \, t_j \, \psi( h_i^{-1} h_j v ) = \bigl\| \sum_i t_i h_i v \bigr\|^2$.}

\item \label{GNSPf-ActCont&IsomEx}
Prove that the action of~$H$ acts on $\real[H]$ by translation is continuous, and preserves the bilinear form defined in the proof of \cref{GNS}. 

\item \label{SchoenbergEx-K} 
Show that the kernel~$\kappa$ in the proof of \cref{Schoenberg} is of positive type.
\hint{For $v_1,\ldots,v_n \in \real$ and $h_1,\ldots,h_n \in H$, let $v_0 = -\sum v_i$ and $h_0 = e$. Then $\sum_{i,j} v_i v_j \kappa(h_i,h_j) \ge 0$ since $\varphi$ is conditionally of positive type and $\sum v_i = 0$. However, the terms with either $i = 0$ or $j = 0$ have no net contribution, since $\varphi(e) = 0$.}

\item \label{SchoenbergEx-prod} 
Show that if $\kappa$ and~$\lambda$ are kernels of positive type, then $\kappa \lambda$ is a kernel of positive type.
\hint{Given $h_1,\ldots,h_n \in H$, show there is a matrix~$L$, such that $L^2 = \bigl( \lambda(h_i,h_j) \bigr)$. Note that
	$\sum\nolimits_{i,j} v_i \, v_j \, \kappa(h_i,h_j) \, \lambda(h_i,h_j) = \sum\nolimits_k \sum\nolimits_{i,j} ( L_{i,k} v_i)  \,( L_{k,j} v_j) \, \kappa(h_i, h_j) \ge 0$.}

\item \label{SchoenbergEx-exp} 
If $\kappa$ is a kernel of positive type, show $e^\kappa$ is a kernel of positive type.
\hint{Since $\kappa^n$ is a kernel of positive type for all~$n$, the same is true of $\sum \kappa^n/n!$.}

\item \label{SchoenbergEx-func}
For every $\varphi \colon H \to \real$, show $\varphi(g)\, \varphi(h)$ is a kernel of positive type.

\item Suppose $\varphi \colon H \to \real$.
Prove the following converse of \cref{Schoenberg}:
If $e^{t\varphi}$ is of positive type for all $t > 0$, then $\varphi$ is conditionally of positive type.
\hint{Verify that $e^{t\varphi} - 1$ is conditionally of positive type.
Then $\lim_{t \to 0^+} (e^{t\varphi} - 1)/t$ has the same type.}

\item \label{T<>FH-PhiIsEquiEx}
Verify \pref{T<>FH-PhiIsEqui}.
\hint{Since $\widehat\alpha\bigl( -\tau(h), h \bigr)(0) = 0$, we have $\bigl\langle \widehat\pi(e,h) \hat v \mid \widehat\pi \bigl( \tau(h), e \bigr) \hat v \bigr\rangle = 1$. Since they are of norm~$1$, the two vectors must be equal.}

\item \label{Haagerup<>proper}
Recall that a Lie group~$H$ has the \defit{Haagerup property} if it has a unitary representation, such that there are almost-invariant vectors, and all matrix coefficients decay to~$0$ at~$\infty$.
It is known that $H$ has the Haagerup property if and only if it has a continuous, \emph{proper} action by affine isometries on some Hilbert space. Prove the implication ($\Leftarrow$) of this equivalence.

\end{exercises}

\begin{notes}

The monograph \cite{BekkaHarpeValette-T} is the standard reference on Kazhdan's property~$(T)$. 
%It includes much material that is not covered here.
%Brief treatments can also be found in \cite[Chap.~3]{MargulisBook} and \cite[Chap.~7]{ZimmerBook}.
The property was defined by D.\,Kazhdan in \cite{KazhdanT}, where \cref{KazhdanEasy,SL3RHasT,Kazhdan:G->Gamma} were proved.

See \cite{BaderFurmanGelanderMonod-BanachT} for a discussion of the generalization of Kazhdan's property to actions on Banach spaces, including \cref{TnotBanach,BanachNoInvt}.
See \cite{CherixEtAl-Haagerup} for a discussion of the \term{Haagerup property} that is mentioned in \cref{T->NotHaagerup,Haagerup<>proper}.
See \cite{Lubotzky-ExpandingGrps} for much more information about expander graphs and their connection with Kazhdan's property, mentioned in \cref{Kazhdan->Expander}.

Regarding \cref{TIsQuotOfFP}:
\noprelistbreak
	\begin{itemize} %\itemsep=\smallskipamount
	\item By showing that $\SL\bigl( 3, \mathbf{F}_q[t] \bigr)$ is not finitely presentable, H.\,Behr \cite{Behr-SL3NotFP} provided the first example of a group with Kazhdan's property that is not finitely presentable.
	\item The existence of uncountably many Kazhdan groups was proved by M.\,Gromov \cite[Cor.~5.5.E, p.~150]{Gromov-HypGrp}. More precisely, any cocompact lattice in $\Sp(1,n)$ has uncountably many different quotients (because it is a ``hyperbolic'' group), and all of these quotient groups have Kazhdan's property.
	\item Y.\,Shalom \cite[p.~5]{Shalom-RigidCommens} proved that every discrete Kazhdan group is a quotient of a finitely presented Kazhdan group. The proof can also be found in 
\cite[Thm.~3.4.5, p.~187]{BekkaHarpeValette-T}.
	\end{itemize}

Our proof of \cref{SL3RHasT} is taken from \cite[\S1.4]{BekkaHarpeValette-T}. 

\Cref{WhichGKazhdan} appears in \cite[Thm.~3.5.4, p.~177]{BekkaHarpeValette-T}. It combines work of several people, including Kazhdan \cite{KazhdanT} %Wang \cite{Wang-DualSpace}, 
and Kostant \cite{Kostant-Tannounce,Kostant-Tpaper}. See \cite[pp.~5--7]{BekkaHarpeValette-T} for an overview of the various contributions to this theorem.

A detailed solution of \cref{ConstructPhiEx} can be found in \cite[Lem.~1.4.1]{BekkaHarpeValette-T}.

See \cite[Thm.~1.7.1, p.~60]{BekkaHarpeValette-T} for a proof of \cref{Kazhdan:G->Gamma} and its converse (\cref{GammaT->GT}).

Regarding \fullcref{KazhdanLattRem}{abel}, see \cite{Lubotzky-BettiNum} (and its references) for a discussion of W.\,Thurston's conjecture that lattices in $\SO(1,n)$ have finite-index subgroups with infinite abelian quotients. (For $n = 3$, the conjecture was proved by Agol \cite{Agol-VirtHakenConj}.)
Lattices in $\SU(1,n)$ with an infinite abelian quotient were found by D.\,Kazhdan \cite{Kazhdan-WeilRep}.

Explicit Kazhdan constants for $\SL(n,\integer)$ \ccf{ExplicitKazhdanConstants} were first found by M.\,Burger \cite[Appendix]{delaHarpeValette-T} (or see \cite[\S4.2]{BekkaHarpeValette-T}). 
An approach developed by Y.\,Shalom (see \cite{Shalom-AlgT}) applies to more general groups (such as $\SL \bigl( n, \integer[x] \bigr)$) that are not assumed to be lattices.

\fullCref{OtherT}{random} is a theorem of Zuk \cite[Thm.~4]{Zuk-TandConsts}.
(Or see \cite{KotowskiKotowski-RandGrps} for a more detailed proof.)
\fullCref{OtherT}{Alg} is explained in \cite{Shalom-AlgT}.

\Cref{TnotQI} is due to S.\,Gersten (unpublished). A proof (based on the same example, but rather different from our sketch) is in \cite[Thm.~3.6.5, p.~182]{BekkaHarpeValette-T}.

\Cref{T<>FH} is due to Delorme \cite[Thm.~V.1]{Delorme-T<>FH} (for ($\ref{T<>FH-T}\Rightarrow\ref{T<>FH-FP}$)) and Guichardet \cite[Thm.~1]{Guichardet-T<>FH} (for ($\ref{T<>FH-H1=0} \Rightarrow \ref{T<>FH-T}$)). 

\Cref{T<>H1bar=0} was proved for discrete groups by Korevaar and Schoen \cite[Cor.~4.1.3]{KorevaarSchoen-Global}. The general case is due to Shalom  \cite[Thm.~6.1]{Shalom-RigidCommens}. 
\Cref{T<>H1irred=0,TiffH1(irred)} are also due to Shalom \cite[proof of Thm.~0.2]{Shalom-RigidCommens}. 

The part of \cref{T->FPonOther} dealing with real or complex hyperbolic spaces is in \cite[Cor.~2.7.3]{BekkaHarpeValette-T}. See \cite{ChatterjiEtAl-Median} for median spaces.

The existence and uniqueness of the circumcenter (mentioned in the hints to \cref{FP<>SomeOrbitBdd<>EveryOrbitBddEx,HinftyUnbdd}) is proved in \cite[Lem.~2.2.7]{BekkaHarpeValette-T}.

\Cref{WatataniThm} is due to Watatani \cite{Watatani-T->FA}, and can also be found in \cite[\S2.3]{BekkaHarpeValette-T}.
Serre's book \cite{Serre-Trees} is a very nice exposition of the theory of group actions on trees, but, unfortunately, does not include this theorem.

See \cite[\S2.10--\S2.12 and \S C.4]{BekkaHarpeValette-T} for the material of \cref{PosDefFuncSect}.
\end{notes}

 %!TEX root = IntroArithGrps.tex

\mychapter{Ergodic Theory} \label{ErgodicChap}

\prereqs{none.}

Ergodic Theory is the study of measure-theoretic aspects of group actions. 
%(Because ergodic theory is essentially synonymous with probability theory, one can also say that it is the study of probabilistic aspects of group actions.) 
%(Classical Ergodic Theory considers only actions of $\integer$ or~$\real$.)
Topologists and geometers may be more comfortable in the category of continuous functions, but important results in \cref{MargulisSuperChap,NormalSubgroupChap} will be proved by using measurable properties of actions of~$\Gamma$, so we will introduce some of the basic ideas. 

\section{Terminology}

The reader is invited to skim through this section, and refer back as necessary.

\begin{assump} \ 
\noprelistbreak
\begin{enumerate}
\item All measures are assumed to be \defit[sigma-finite measure@$\sigma$-finite measure]{$\sigma$-finite}. That is, if $\mu$ is a measure on a measure space~$X$, then we always assume that $X$ is the union of countably many subsets of finite measure.
\item We have no need for abstract measure spaces, so all measures are assumed to be \defit[measure!Borel]{Borel}. That is, when we say $\mu$ is a measure on a measure space~$X$, we are assuming that $X$ is a Borel subset of a complete, separable, metrizable space, and the implied $\sigma$-algebra on~$X$ consists of the subsets of~$X$ that are equal to a Borel set, modulo a set of measure~$0$.
\end{enumerate}
\end{assump}

\begin{defns}
Let $\mu$ be a measure on a measure space~$X$.
\noprelistbreak
	\begin{enumerate}
	\item We say $\mu$ is a \defit[measure!finite]{finite measure} if $\mu(X) < \infty$.
	\item A subset $A$ of~$X$ is:
		\begin{itemize}
		\item \defit[null set]{null} if $\mu(A) = 0$,
		\item \defit[conull set]{conull} if the complement of~$A$ is null.
		\end{itemize}
	\item We often abbreviate ``almost everywhere'' to ``a.e\zz.''
	\item \defit[essentially]{Essentially} is a synonym for ``almost everywhere\zz.'' For example, a function~$f$ is \emph{essentially constant} iff $f$~is constant (a.e.).
%	\item For any measurable function $\varphi \colon X \to Y$, the measure $\varphi_*\mu$ on~$Y$ is defined by
%		\nindex{$\varphi_*\mu$ = push-forward of measure~$\mu$}
%		$$ (\varphi_*\mu)(A) = \mu \bigl( \varphi^{-1}(A) \bigr) .$$
	\item Two measures $\mu$ and~$\nu$ on~$X$ are in the same \defit[measure!class]{measure class}
		%(or $\nu$ is \defit[equivalent!measure]{equivalent} to~$\mu$) 
if they have exactly the same null sets:
		$$ \mu(A) = 0 \iff \nu(A) = 0 .$$
	(This defines an equivalence relation.) Note that if $\nu = f \mu$, for some real-valued, measurable function~$f$, such that $f(x) \neq 0$ for a.e.\ $x \in X$, then $\mu$ and~$\nu$ are in the same measure class \csee{fmuClassOfMuEx}. 
The Radon-Nikodym Theorem \pref{RadonNikodym} implies that the converse is true.
	\end{enumerate}
\end{defns}

\begin{defns}
Suppose $H$ is a Lie group~$H$ that acts continuously on a metrizable space~$X$, $\mu$ is a measure on~$X$, and $A$~is a subset of~$X$.
	\begin{enumerate}
	\item The set $A$ is \defit[invariant!set]{invariant} (or, more precisely, \defit[invariant!set]{$H$-invariant}) if $hA = A$ for all $h \in H$.
	\item The measure~$\mu$ is \defit[invariant!measure]{invariant} (or, more precisely, \defit[invariant!measure]{$H$-invariant}) if $h_*\mu = \mu$ for all $h \in H$. (Recall that the push-forward $h_*\mu$ is defined in \pref{PushForwardDefn}.)
	\item The measure~$\mu$ is \defit[quasi-!invariant measure]{quasi-invariant} if $h_*\mu$ is in the same measure class as~$\mu$, for all $h \in H$.
	\item A (measurable) function $f$ on~$X$ is \defit[essentially!$H$-invariant]{essentially $H$-invariant} if, for every $h \in H$, we have
	$$ \text{$f(hx) = f(x)$ for a.e.\ $x \in X$.} $$
	\end{enumerate}
\end{defns}

The Lebesgue measure on a manifold is not unique, but it determines a well-defined measure class, which is invariant under any smooth action:

\begin{lem}[\csee{LebesgueMeasClassEx}]
If $X$ is a manifold, and $H$ acts on~$X$ by diffeomorphisms, then Lebesgue measure provides a measure on~$X$ that is quasi-invariant for~$H$.
\end{lem}

\begin{exercises}

\item \label{fmuClassOfMuEx}
Suppose $\mu$ is a measure on a measure space~$X$, and $f$~is a real-valued, measurable function on~$X$, such that $f \ge 0$ for a.e.~$x$. Show that $f \mu$ is in the measure class of~$\mu$ iff $f(x) \neq 0$ for a.e.\ $x \in X$.

\item Suppose a Lie group~$H$ acts continuously on a metrizable space~$X$, and $\mu$ is a measure on~$X$. Show that $\mu$ is quasi-invariant iff the collection of null sets is $H$-invariant. (This means that if $A$~is a null set, and $h \in H$, then $h(A)$ is a null set.)

\item \label{Diffble(null)=nullEx}
Suppose 
	\begin{itemize}
	\item $A$ is a null set in~$\real^n$ (with respect to Lebesgue measure),
	and
	\item $f$ is a diffeomorphism of some open subset~$\open$\, of~$\real^n$.
	\end{itemize}
Show that $f(A \cap \open\,)$ is a null set.
\hint{Change of variables.}

\item \label{LebesgueMeasClassEx}
Suppose $X$ is a (second countable) smooth, $n$-dimensional manifold. This means that $X$ can be covered by coordinate patches $(X_i,\varphi_i)$ (where $\varphi_i \colon X_i \to \real^n$, and the overlap maps are smooth).
	\begin{enumerate}
	\item Show there exists a partition $X = \bigcup_{i=1}^\infty \hat X_i$ into measurable subsets, such that $\hat X_i \subseteq X_i$ for each~$i$.
	\item Define a measure $\mu$ on~$X$ by $\mu(X) = \lambda \bigl( \varphi_i(A \cap \hat X_i) \bigr)$, where $\lambda$ is the Lebesgue measure on~$\real^n$. This measure may depend on the choice of $X_i$, $\varphi_i$, and~$\hat X_i$, but show that the measure class of~$\mu$ is independent of these choices. 
	\end{enumerate}
\hint{\Cref{Diffble(null)=nullEx}.}

\end{exercises}

\section{Ergodicity} \label{ErgodicitySect}

Suppose $H$ acts on a topological space~$X$. If $H$ has a dense orbit on~$X$, then it is easy to see that every continuous, $H$-invariant function is constant \csee{DenseOrb->FuncConstEx}. Ergodicity is the much stronger condition that every \emph{measurable} $H$-invariant function is constant (a.e.):

\begin{defn} \label{ErgodicDefn}
Suppose $H$ acts on~$X$ with a quasi-invariant measure~$\mu$.
We say the action of~$H$ is \defit{ergodic} (or that $\mu$~is an \defit{ergodic} measure for~$H$) if every $H$-invariant, real valued, measurable function on~$X$ is essentially constant.
\end{defn}

It is easy to see that transitive actions are ergodic \csee{transitive->erg}. But non-transitive actions can also be ergodic:

\begin{eg}[(Irrational rotation of the circle)] \label{IrratRotErg}
For any $\alpha \in \real$, we may define a homeomorphism~$T_\alpha$ of the circle $\torus = \real/\integer$ by
	$$ T_\alpha(x) = x + \alpha \pmod{\integer} .$$
By considering Fourier series, it is not difficult to show that if $\alpha$~is irrational, then every $T_\alpha$-invariant function in $\LL2(\torus)$ is essentially constant \csee{IrratRotErgEx}. This implies that the $\integer$-action generated by~$T_\alpha$ is ergodic \csee{ErgCheckLp}.
\end{eg}

\Cref{IrratRotErg}  is a special case of the following general result:

\begin{prop} \label{Dense->Erg}
If $H$ is any dense subgroup of a Lie group~$L$, then the natural action of~$H$ on~$L$ by left translation is ergodic\/ \textup(with respect to the Haar measure on~$L$\textup).
\end{prop}

\begin{proof}
For any measurable $f \colon L \to \real$, its \defit[stabilizer!essential]{essential stabilizer} in~$L$ is defined to be:
	$$ \Stab_L(f) = \{\, g \in L \mid \text{$f(gx) = f(x)$ \ for a.e.\ $x \in L$} \,\} .$$
It is not difficult to show that $\Stab_L(f)$ is closed \csee{EssStabClosed}.
(It is also a subgroup of~$L$, but we do not need this fact.)  Hence, if $\Stab_L(f)$ contains a dense subgroup~$H$, then it must be all of~$L$. This implies that $f$ is constant (a.e.) \csee{EssInvt->EssConst}.
\end{proof}

It was mentioned above that transitive actions are ergodic; therefore, $G$ is ergodic on $G/\Gamma$.
What is not obvious, and leads to important applications for arithmetic groups, is that most subgroups of~$G$ are also ergodic on $G/\Gamma$:

\begin{thm}[(\thmindex{Moore Ergodicity}Moore Ergodicity Theorem, see \cref{MooreErgLatticeEx})] \label{MooreErgodicity}
If 
	\begin{itemize}
	\item $H$ is any noncompact, closed subgroup of~$G$, 
	and
	\item $\Gamma$~is irreducible,
	\end{itemize}
then $H$ is ergodic on $G/\Gamma$.
\end{thm}

If $H$ is ergodic on $G/\Gamma$, then $\Gamma$ is ergodic on $G/H$ \csee{MackeySwitchErgodic}. Hence:

\begin{cor} \label{GammaErgOnG/H}
If $H$ and\/~$\Gamma$ are as in \cref{MooreErgodicity},
then\/ $\Gamma$ is ergodic on $G/H$.
\end{cor}

\begin{exercises}

\item \label{EssStabClosed}
Show $\Stab_L(f)$ is closed, for every Lie group~$L$ and measurable $f \colon L \to \real$.
\hint{If $f$~is bounded, then, for any $\varphi \in C_c(L)$ and $\{g_n\} \subseteq \Stab_L(f)$, we have
	$ \int_L {}^g \! f \cdot \varphi \,d \mu
	= \int_L f \cdot {}^{g^{-1}} \! \varphi \,d \mu
	= \lim \int_L f \cdot {}^{g_n^{-1}} \! \varphi \,d \mu
	= \lim \int_L {}^{g_n} \! f \cdot  \varphi \,d \mu
%	= \lim_{n \to \infty} \int_L f \cdot  \varphi \,d \mu
	= \int_L f \cdot  \varphi \,d \mu
	$.} %since $\varphi$~is uniformly continuous.
%There is a weak topology on $\LL\infty(L, \mu)$, defined by 
%	$$ f_n \to f \quad \Leftrightarrow \quad
%	\int_L f_n \varphi \,d\mu \to \int_L f\varphi \,d\mu
%	\ \text{for all $\varphi \in \LL1(L,\mu)$} .$$
%Since functions in $\LL1$ can be approximated by uniformly continuous functions, it is not difficult to see that the action of~$L$ on  $\LL\infty(L, \mu)$ is continuous (with respect to this weak topology). This implies that stabilizers are closed.}

\item \label{DenseOrb->FuncConstEx}
Suppose $H$ acts on a topological space~$X$, and has a dense orbit. Show that every real-valued, continuous, $H$-invariant function on~$X$ is constant.

\item \label{transitive->erg}
Show that $H$ is ergodic on $H/H_1$, for every closed subgroup~$H_1$ of~$H$.
\hint{Every $H$-invariant function is constant, not merely essentially constant.}

\item \label{EssInvt->EssConst}
Suppose $H$ is ergodic on~$X$, and $f \colon X \to \real$ is measurable and essentially $H$-invariant. Show that $f$ is essentially constant.

\item Our definition of ergodicity is not the usual one, but it is equivalent:
show that $H$ is ergodic on~$X$ iff every $H$-invariant measurable subset of~$X$ is either null or conull. 
\hint{The characteristic function of an invariant set is an invariant function. Conversely, the sub-level sets of an invariant function are invariant sets.}

\item \label{IrratRotErgEx}
In the notation of \cref{IrratRotErg}, show (without using \cref{Dense->Erg}):
	\begin{enumerate}
	\item If $\alpha$~is irrational, then every $T_\alpha$-invariant function in $\LL2(\torus)$ is essentially constant.
	\item If $\alpha$~is rational, then there exist $T_\alpha$-invariant functions in $\LL2(\torus)$ that are not essentially constant.
	\end{enumerate}
\hint{Any $f \in \LL2(\torus)$ can be written as a unique Fourier series:
	$ f = \sum_{n=-\infty}^\infty a_n e^{in\theta}$.
If $f$ is invariant and $\alpha$~is irrational, then uniqueness implies $a_n = 0$ for $n \neq 0$.}

\item \label{ErgCheckLp}
Suppose $\mu$ is an $H$-invariant, \emph{finite} measure on~$X$. For all $p \in [1,\infty]$, show that $H$ is ergodic iff every $H$-invariant element of $\LL{p}(X,\mu)$ is essentially constant.

\item Let $H = \integer$ act on $X = \real$ by translation, and let $\mu$ be Lebesgue measure. Show:
	\begin{enumerate}
	\item $H$ is not ergodic on~$X$,
	and
	\item for every $p \in [1,\infty)$, every $H$-invariant element of $\LL{p}(X,\mu)$ is essentially constant.
	\end{enumerate}
Why is this not a counterexample to \cref{ErgCheckLp}?

\item Let $H$ be a dense subgroup of~$L$. Show that if $L$ is ergodic on~$X$, then $H$ is also ergodic on~$X$.
\hint{\cref{EssStabClosed}.}

\item Show that if $H$ acts continuously on~$X$, and $\mu$~is a quasi-invariant measure on~$X$, then the support of~$\mu$ is an $H$-invariant subset of~$X$.

\item \label{AEOrbitDense}
Ergodicity implies that a.e.\ orbit is dense in the support of~$\mu$. More precisely, show that if $H$ is ergodic on~$X$, and the support of~$\mu$ is all of~$X$ (in other words, no open subset of~$X$ has measure~$0$), then a.e.\ $H$-orbit in~$X$ is dense. (That is, for a.e.\ $x \in X$, the orbit $Hx$ of~$x$ is dense in~$X$.
\hint{The characteristic function of the closure of any orbit is invariant.}

\item \label{MackeySwitchErgodic}
Suppose $H$ is a closed subgroup of~$G$. Show that $H$ is ergodic on $G/\Gamma$ iff $\Gamma$ is ergodic on $G/H$.
\hint{$H \times \Gamma$ acts on~$G$ (by letting $H$ act on the left and $\Gamma$~act on the right). Show $H$ is ergodic on $G/\Gamma$ iff $H \times \Gamma$ is ergodic on~$G$ iff $\Gamma$ is ergodic on $G/H$.}

\item The Moore Ergodicity Theorem has a converse: Assume $G$ is not compact, and show that if $H$ is any compact subgroup of~$G$, then $H$~is not ergodic on $G/\Gamma$.
\hint{$\Gamma$ acts properly discontinuously on $G/H$, so the orbits are not dense.}

\item Show that if $n \ge 2$, then
	\begin{enumerate}
	\item the natural action of $\SL(n,\integer)$ on~$\real^n$ is ergodic,
	and
	\item the $\SL(n,\integer)$-orbit of a.e.\ vector in~$\real^n$ is dense in $\real^n$.
	\end{enumerate}
\hint{Identify $\real^n$ with a homogeneous space of $G = \SL(n,\real)$ (a.e.), by noting that $G$ is transitive on the nonzero vectors of~$\real^n$.}

\item
Let 	
	\begin{itemize}
	\item $G = \SL(3,\real)$,
	\item $\Gamma$ be a lattice in~$G$,
	and
	\item $P = \begin{Smallbmatrix} \upast&& \\ \upast&\upast& \\ \upast&\upast&\upast\end{Smallbmatrix} \subset G$.
	\end{itemize}
Show:
\begin{enumerate}
\item The natural action of~$\Gamma$ on the homogeneous space $G/P$ is ergodic.
\item The diagonal action of~$\Gamma$ on $(G/P)^2 = (G/P) \times (G/P)$ is ergodic.
\item The diagonal action of~$\Gamma$ on $(G/P)^3 = (G/P) \times (G/P) \times (G/P)$ is \emph{not} ergodic.
\end{enumerate}
\hint{$G$ is transitive on a conull subset of $(G/P)^k$, for $k \le 3$. What is the stabilizer of a generic point in each of these spaces?}

 \item \label{AEOrbitDenseInG/Gamma}
 Assume $\Gamma$ is irreducible, and let $H$ be a closed, noncompact subgroup of~$G$. Show, for a.e.\ $x \in G/\Gamma$, that $Hx$ is dense in $G/\Gamma$.

\item \label{RErgodic->ZErgodic}
Suppose $H$ acts ergodically on~$X$, with invariant measure~$\mu$. 
Show that if $\mu(X) < \infty$ and $H \iso \real$, then some cyclic subgroup of~$H$ is ergodic on~$X$.
\hint{For each $t \in \real$, choose a nonzero, $h^t$-invariant function $f_t \in L^2(X)$, such that $f_t \perp 1$. The projection of $f_r$ to the space of $h^s$-invariant functions is invariant under both $h^r$ and~$h^s$. 
Therefore, if $r$ and~$s$ are linearly independent over~$\rational$, then $f_r \perp f_s$. This is impossible, because $\LL2(X,\mu)$ is separable.}

\end{exercises}

\section{Consequences of a finite invariant measure}

Measure-theoretic techniques are especially powerful when the action has an invariant  measure that is finite. One example of this is the Poincar\'e Recurrence Theorem \pref{PoincareRecurThm}. Here is another.

We know that almost every orbit of an ergodic action is dense \csee{AEOrbitDense}.
For the case of a $\integer$-action with a finite, invariant measure, the orbits are not only dense, but uniformly distributed:

\begin{defn}
Let 
	\begin{itemize}
	\item $\mu$ be a finite measure on a topological space~$X$,
		% such that $\mu(X) = 1$,
	and
	\item $T$ be a homeomorphism of~$X$. 
	\end{itemize}
The $\langle T \rangle$-orbit of a point $x$ in~$X$ is \defit{uniformly distributed} with respect to~$\mu$ if
	$$ \lim_{n \to \infty} \frac{1}{n} \sum_{k=1}^n f \bigl( T^k(x) \bigr) 
	= \frac{1}{\mu(X)} \int_X f \, d\mu ,$$
for every bounded, continuous function~$f$ on~$X$.
\end{defn}

\begin{thm}[{(\thmindex{Ergodic!Pointwise}{Pointwise Ergodic Theorem})}] 
	%\index{Pointwise Ergodic Theorem|indsee{Ergodic Theorem, Pointwise}}
\label{PointwiseErgThm}
Suppose 
	\begin{itemize}
	\item $\mu$ is a finite measure on a second countable, metrizable space~$X$,
	and
	\item $T$ is an ergodic, measure-preserving homeomorphism of~$X$.
	\end{itemize}
Then a.e.\ $\langle T \rangle$-orbit in~$X$ is uniformly distributed\/ \textup(with respect to~$\mu$\textup).
\end{thm}

It is tricky to show that $\lim_{n \to \infty} \frac{1}{n} \sum_{k=1}^n f \bigl( T^k(x) \bigr)$ converges %to $\int_X f \, d\mu$ 
pointwise \csee{PtwiseErgZ}.
Convergence in norm is much easier \csee{MeanErgThmEx}. 

\begin{rem} \label{PtwiseErgAmenRem}
Although the Pointwise Ergodic Theorem was stated only for actions of a cyclic group, it generalizes very nicely to the ergodic actions of any \term{amenable group}. (The values of~$f$ are averaged over an appropriate \term[Folner@F\o lner!set]{F\o lner set} in the amenable group.) See \cref{PointwiseErgThmForREx} for actions of~$\real$.
\end{rem}

\begin{exercises}

\item Suppose the $\langle T \rangle$-orbit of~$x$ is uniformly distributed with respect to a finite measure~$\mu$ on~$X$. Show that if the support of~$\mu$ is all of~$X$, then the $\langle T \rangle$-orbit of~$x$ is dense in~$X$.

\item \label{IterateUnitaryEx}
Suppose 
	\begin{itemize}
	\item $U$ is a unitary operator on a \term{Hilbert space}~$\Hilbert$, 
	\item $v \in \Hilbert$,
	and
	\item $\langle v \mid w \rangle = 0$, for every vector~$w$ that is fixed by~$U$.
	\end{itemize}
Show $\frac{1}{n} \sum_{k=1}^n U^k v \to 0$ as $n \to \infty$.
\hint{Apply the {Spectral Theorem} to diagonalize the unitary operator~$U$.}

\item \label{MeanErgThmEx}
(\thmindex{Ergodic!Mean}{Mean Ergodic Theorem})
Assume the setting of the Pointwise Ergodic Theorem \pref{PointwiseErgThm}.
Show that if $f \in \LL2(X,\mu)$, then
	$$\frac{1}{n} \sum_{k=1}^n f \bigl( T^k(x) \bigr) \to \frac{1}{\mu(X)} \int_X f \, d\mu \text{\quad in $\LL2$} .$$
That is, show
	$$ \lim_{n \to \infty} \left\| \ \frac{1}{n} \sum_{k=1}^n f \bigl( T^k(x) \bigr) 
	\  - \  \frac{1}{\mu(X)} \int_X f \, d\mu \ \right\|_2 = 0 .$$
\emph{Do not assume \cref{PointwiseErgThm}.}
\hint{\Cref{IterateUnitaryEx}.}

\item \label{PtwiseErgZMaxl}
 Assume $X$, $\mu$, and~$T$ are as in \cref{PointwiseErgThm}, and that $\mu(X) = 1$. For $f \in L^1(X,\mu)$, define
	$$ S_n(x) =  f(x) + f \bigl( T(x) \bigr) + \cdots + f \bigl( T^{n-1}(x) \bigr) .$$
  Prove the \thmindex{Ergodic!Maximal}\textit{\textbf{Maximal Ergodic
Theorem}}: for every $\alpha \in \real$, if we let
 $$ E_\alpha = \bigset{ x \in X }{ \sup_n \frac{S_n(x)}{n} > \alpha} ,$$
 then $\int_E f \, d\mu \ge \alpha \, \mu(E)$.
 \hint{Assume $\alpha = 0$. Let $S^+_n(x) = \max_{0\le k\le n} S_k(x)$,
and $E_n = \{\, x \mid S^+_n > 0\,\}$, so $E = \cup_n E_n$. For $x \in
E_n$, we have $f(x) \ge S^+_n(x) - S^+_n \bigl( T(x) \bigr)$, so
$\int_{E_n} f\, d\mu \ge 0$.}

\item \label{PtwiseErgZ}
 Prove the Pointwise Ergodic Theorem \pref{PointwiseErgThm}.
 \hint{Either $\{\, x \mid \limsup S_n(x)/n >
\alpha \,\}$ or its complement must have measure~$0$. If it is the complement, then \cref{PtwiseErgZMaxl} implies $\int_X f \, d\mu \ge \alpha$.}

\item \label{PointwiseErgBddL1Ex}
Assume the setting of the Pointwise Ergodic Theorem \pref{PointwiseErgThm}. 
For every bounded $\varphi \in \LL1(X,\mu)$, show, for a.e.\ $x \in X$, that 
	$$ \lim_{n \to \infty} \frac{1}{n} \sum_{k=1}^n \varphi \bigl( T^k(x) \bigr) 
	= \frac{1}{\mu(X)} \int_X \varphi \, d\mu .$$
\hint{You may assume the Pointwise Ergodic Theorem. Use {Lusin's Theorem}~\pref{LusinsThm} to approximate~$\varphi$ by a continuous function.}

\item \harder 
Remove the assumption that $\varphi$ is bounded in \cref{PointwiseErgBddL1Ex}.

\item \label{PointwiseErgThmForREx}
Suppose 
	\begin{itemize}
	\item $\{a^t\}$ is a (continuous) $1$-parameter group of homeomorphisms of a topological space~$X$,
and
	\item $\mu$ is an ergodic, $a^t$-invariant, finite measure on~$X$.
	\end{itemize}
For every bounded, continuous function~$f$ on~$X$, show that
	$$ \text{$\displaystyle \lim_{T \to \infty} \frac{1}{T} \int_0^t f(a^t x) \, dt = \frac{1}{\mu(X)} \int_X f\, d\mu$ \quad  for a.e.~$x \in X$} .$$
\hint{Apply \cref{PointwiseErgThm} to $\overline{f}(x) = \int_0^1 f(a^t x)\, dt$ if $a^1$ is ergodic \ccf{RErgodic->ZErgodic}.}

\end{exercises}

\section{Ergodic decomposition} \label{ErgDecompSect}

In this \lcnamecref{ErgDecompSect}, we briefly explain that every group action (with a quasi-invariant measure) can be decomposed into ergodic actions.

\begin{eg}[(Irrational rotation of the plane)]
For any irrational $\alpha \in \real$, define a homeomorphism~$T_\alpha$ of~$\complex$ by
	$ T_\alpha(z) = e^{2\pi i \alpha} x $.
Then $|T_\alpha(z)| = |z|$, so each circle centered at the origin is invariant under~$T_\alpha\mk$. The restriction of $T_\alpha$ to any such circle is an irrational rotation of the circle, so it is ergodic \csee{IrratRotErg}. Thus, the entire action can be decomposed as a union of ergodic actions. 
\end{eg}

A similar decomposition is always possible, as long as we work with nice spaces:

\begin{defn}
A topological space~$X$ is \defit[Polish topological space]{Polish} if it is homeomorphic to a complete, separable metric space.
\end{defn}

\begin{thm}[(\thmindex{Ergodic!decomposition}Ergodic decomposition)] \label{ErgodicDecomp}
Suppose a Lie group~$H$ acts continuously on a Polish space~$X$, and $\mu$ is a quasi-invariant, finite measure on~$X$. Then there exist
	\begin{itemize}
	\item a measurable function $\zeta \colon X \to [0,1]$, 
	and
	\item a finite measure~$\mu_z$ on $\zeta^{-1}(z)$, for each $z \in [0,1]$, 
	\end{itemize}
such that $\mu = \int_{[0,1]} \mu_z \, d \nu(z)$, where $\nu = \zeta_* \mu$. For $f \in C_c(X)$, this means
	$$ \int_X f \, d\mu = \int_Z \, \int_{\zeta^{-1}(z)} f \, d\mu_z \, d\nu(z) .$$
Furthermore, for a.e.\ $z \in [0,1]$, 
	\begin{enumerate}
	\item $\zeta^{-1}(z)$ is $H$-invariant,
	and
	\item $\mu_z$ is quasi-invariant and ergodic for the action of~$H$.
	\end{enumerate}
\end{thm}

\begin{rem}
The ergodic decomposition is unique (a.e.). More precisely, if $\zeta'$ and~$\mu_z'$ also satisfy the conclusions, then there is a measurable bijection $\pi \colon [0,1] \to [0,1]$, such that 
	\begin{enumerate}
	\item $\zeta' = \pi \circ \zeta$ a.e.,
	and
	\item $\mu_{\pi(z)}' = \mu_z$ for a.e.~$z$.
	\end{enumerate}
\end{rem}

\begin{defn}
In the notation of \cref{ErgodicDecomp}, each set $\zeta^{-1}(z)$ is called an \defit[ergodic!component]{ergodic component} of the action.
\end{defn}

The remainder of this \lcnamecref{ErgDecompSect} sketches two different proofs of \cref{ErgodicDecomp}.

\subsection{First proof}
The main problem is to find the function~$\zeta$, because the following general Fubini-like theorem will then provide the required decomposition of~$\mu$ into an integral of measures~$\mu_z$ on the fibers of~$\zeta$. (In Probability Theory, each $\mu_z$ is called a \defit{conditional measure} of~$\mu$.)

\begin{prop} %[(\thmindex{Rokhlin's}Rokhlin)]
\label{RokhlinDecompMeas}
Suppose
	\begin{itemize}
	\item $X$ and~$Z$ are Polish spaces,
	\item $\zeta \colon X \to Z$ is a Borel measurable function,
	and
	\item $\mu$ is a probability measure on~$X$.
	\end{itemize}
Then there is a Borel map $\lambda \colon Z \to \Prob(X)$, such that 
 \begin{enumerate}
 \item $\mu = \int_Z \lambda_z \, d\nu(z)$, where $\nu = \zeta_*\mu$,
 and
 \item $\lambda_z \bigl( \zeta^{-1}(z) \bigr) = 1$, for all $z \in Z$.
 \end{enumerate}
 Furthermore, $\lambda$ is unique\/ \textup(a.e.\textup).
 \end{prop}

The map $\zeta$ is a bit difficult to pin down, since it is not completely well-defined --- it can be changed on an arbitrary set of measure zero. We circumvent this difficulty by looking not at the value of~$\zeta$ on individual points (which is not entirely well-defined), but at its effect on an algebra of functions (which is completely well defined). 

\begin{defns} \label{BoolDefn}
Assume $\mu$ is a finite measure on a Polish space~$X$.
	\begin{enumerate}
	\item Let 
	\nindex{$\Bool(X)$ = $\{ \text{Borel subsets} \}$, modulo sets of measure~$0$}%
	$\Bool(X)$ be the collection of all Borel subsets of~$X$, where two sets are identified if they differ by a set of measure~$0$. This is a $\sigma$-algebra.
%	Alternatively, by identifying a set with its characteristic function, we can think of $\Bool(X)$ as the set of $\{0,1\}$-valued measurable functions, considered as a subset of $\LL\infty(X, \mu)$ (which means that two functions are identified if they are equal almost everywhere). 
	\item $\Bool(X)$ is a complete, separable metric space, with respect to the metric $d(A,B) = \mu(A \symmdiff B)$, where 
	\nindex{$A \symmdiff B$ = symmetric difference of $A$ and~$B$}%
	$A \symmdiff B = (A \smallsetminus B) \cup (B \smallsetminus A)$ is the \defit[symmetric!difference]{symmetric difference} of $A$ and~$B$.
%	Being a subset of the unit ball in $\LL\infty(X,\mu)$, $\Bool(X)$ inherits a weak$^*$ topology.  The Banach-Alaoglu Theorem \pref{BanachAlaogluThm} implies that $\Bool(X)$ is compact.
 	\item If a Lie group~$H$ acts continuously on~$X$, we let $\Bool(X)^H$ be the set of $H$-invariant elements of~$\Bool(X)$. This is a sub-$\sigma$-algebra of $\Bool(X)$.
	\end{enumerate}
\end{defns}

%Note that $\Bool(X)$ is a \defit{Boolean algebra}. (This means that it is closed under the set-theoretic operations of union and complement.). Also, $\Bool(X)^H$ is a topologically closed Boolean subalgebra \csee{B(X)HClosedBoolean}. 
The map~$\zeta$ is constructed by the following result:
 
 \begin{lem}
 Suppose 
 	\begin{itemize}
	\item $\mu$ is a finite measure on a Polish space~$X$,
	and
	\item $\Bool$ is a sub-$\sigma$-algebra of $\Bool(X)$.
	\end{itemize}
Then there is a Borel map $\zeta \colon X \to [0,1]$, such that 
	$$ \Bool = \{\, \zeta^{-1}(E) \mid \text{$E$ is a Borel subset of\/ $[0,1]$} \,\} .$$
\end{lem}

\begin{proof}[Idea of proof]
Let 
	\begin{itemize}
	\item $\{E_n\}$ be a countable, dense subset of~$\Bool$,
	\item $\chi_n$ be the characteristic function of~$E_n$, for each~$n$,
	and
	\item $\displaystyle \zeta(x) = \sum_{n = 1}^\infty \frac{\chi_n(x)}{3^n}$.
	\end{itemize}

It is clear from the definition of~$\zeta$ that if $I$ is any open interval in $[0,1]$, then $\zeta^{-1}(I)$ is a Boolean combination of elements of~$\{E_n\}$; therefore, it belongs to~$\Bool$. Since $\Bool$ is a $\sigma$-algebra, this implies that $\zeta^{-1}(E) \in \Bool$ for every Borel subset~$E$ of $[0,1]$.

Conversely, it is clear from the definition of~$\zeta$ that each $E_n$ is the inverse image of a (one-point) Borel subset of $[0,1]$. Since $\{E_n\}$ generates $\Bool$ as a $\sigma$-algebra \csee{DenseGensSigmaAlg}, this implies that every element of~$\Bool$ is the inverse image of a Borel subset of $[0,1]$.
\end{proof}

We will use the following very useful fact:

\begin{thm}[(\thmindex{von\,Neumann Selection}von\,Neumann Selection Theorem)] \label{vonNeumannSelectionThm}
  Let
  \begin{itemize}
  \item $X$ and~$Y$ be Polish spaces,
  \item $\mu$ be a finite measure on~$X$,
  \item $\mathcal{F}$ be a Borel subset of $X \times Y$,
  and
  \item $X_{\mathcal{F}}$ be the projection of~$\mathcal{F}$ to~$X$.
  \end{itemize}
  Then there is a Borel function $\Phi \colon X \to Y$,
  such that $\bigl( x, \Phi(x) \bigr) \in \mathcal{F}$, for a.e.\ $x \in X_{\mathcal{F}}$.
  \end{thm}

\begin{proof}[The first proof of \cref{ErgodicDecomp}]
We wish to show that $\mu_z$ is ergodic (a.e.). If not, then there is a set~$E$ of positive measure in $[0,1]$, such that, for each $z \in E$, the action of~$H$ on $\bigl( \zeta^{-1}(z), \mu_z \bigr)$ is not ergodic. This means there exists an $H$-invariant, measurable, $\{0,1\}$-valued function $\varphi_z \in \LL\infty \bigl( \zeta^{-1}(z), \mu_z \bigr)$ that is \textbf{not} constant (a.e.).
There are technical problems that we will ignore, but, roughly speaking, the  von\,Neumann Selection Theorem \pref{vonNeumannSelectionThm} implies that the selection of $\varphi_z$ can be done measurably, so we have a Borel subset $A$ of~$X$, defined by
	$$ A = \{\, x \in X \mid \text{$\zeta(x) \in E$ and $\varphi_z(x) = 1 $} \,\} .$$

Since $\varphi_z$ is not constant on the fiber $\zeta^{-1}(z)$, we know that $A$ is not of the form $\zeta^{-1}(E)$, for any Borel subset~$E$ of $[0,1]$. On the other hand, we have $A \in \Bool(X)^H$ (since each $\varphi_z$ is $H$-invariant). This contradicts the choice of the function~$\zeta$.
\end{proof}

\subsection{Second proof} \label{ErgDecompPfChoquet}
We now describe a different approach. Instead of obtaining the decomposition of~$\mu$ from the map~$\zeta$, we reverse the argument, and obtain the map~$\zeta$ from a direct-integral decomposition of~$\mu$. 
For simplicity, however, we assume that the space we are acting on is compact. We consider only invariant measures, instead of quasi-invariant measures, so we do not have to keep track of Radon-Nikodym derivatives.

\begin{defns}
Suppose $C$ is a convex subset of a vector space~$V$.
	\begin{enumerate}
	\item A point $c \in C$ is an \defit{extreme point} of~$C$ if there do not exist $c_0,c_1 \in C \smallsetminus \{c\}$ and $t \in (0,1)$, such that $c = t c_0 + (1-t)c_1$.
	\item Let \nindex{$\ext C$ = $\{\text{extreme points of $C$}\}$}$\ext C$ be the set of extreme points of~$C$.
	\end{enumerate}
\end{defns}

\begin{eg}
Suppose $H$ acts continuously on a compact, separable metric space~$X$, and let 
	$$\Prob(X)^H = \{\, \mu \in \Prob(X) \mid \text{$\mu$ is $H$-invariant} \,\} . $$ 
This is a closed, convex subset of $\Prob(X)$, so $\Prob(X)^H$ is a compact, convex subset of $C(X)^*$ (with the weak$^*$ topology). The extreme points of this set are precisely the $H$-invariant probability measures that are ergodic \csee{Erg<>Extreme}.
\end{eg}

The well-known \thmindex{Krein-Milman}Krein-Milman Theorem states that every compact, convex set~$C$ is the closure of the convex hull of the set of extreme points of~$C$. (So, in particular, if $C$ is nonempty, then there exists an extreme point.) We will use the following strengthening of this fact:

\begin{thm}[{(\thmindex{Choquet's}Choquet's Theorem)}] \label{ChoquetsThm}
Suppose 
\noprelistbreak
	\begin{itemize}
	\item $\LocConvex$ is a locally convex topological vector space over~$\real$,
	\item $C$ is a metrizable, compact, convex subset of~$\LocConvex$, 
	and 
	\item $c_0 \in C$.
	\end{itemize}
Then there is a probability measure $\nu$ on $\ext C$, such that
	$$ c_0 = \int_{\ext C} x \, d\nu(x) .$$
\end{thm}

We will also need a corresponding uniqueness result.

\begin{defns}[(Choquet)]
Suppose $\LocConvex$ and $C$ are as in the statement of \cref{ChoquetsThm}. 
\noprelistbreak
	\begin{enumerate}
	\item Let \nindex{$\Sigma C = \{\, t c \mid t \in [0,\infty), \ c \in C \,\}$}%
	$\Sigma C = \{\, t c \mid t \in [0,\infty), \ c \in C \,\} \subseteq \LocConvex$.
	\item Define a partial order $\le$ on~$\Sigma C$ by $a \le b$ iff $b - a \in \Sigma C$.
	\item Two elements $a_1,a_2 \in \Sigma C$ have a \defit{least upper bound} if there exists $b \in \Sigma C$, such that
		\begin{itemize}
		\item $a_i \le b$ for $i = 1,2$,
		and
		\item for all $c \in \Sigma C$, such that $a_i \le c$ for $i = 1,2$, we have $b \le c$.
		\end{itemize}
	\end{enumerate}
\end{defns}

\begin{eg} \label{ProbUB}
Any two elements of $\Sigma \Prob(X)$ have a least upper bound \csee{ProbLUB}, so the same is true for $\Prob(X)^H$.
\end{eg}

\begin{thm}[(Choquet)]
Suppose $\LocConvex$, $C$, and $c_0$ are as in the statement of \cref{ChoquetsThm}. If every two elements of~$\Sigma C$ have a least upper bound, then the measure~$\nu$ provided by \cref{ChoquetsThm} is unique.
\end{thm}

\begin{proof}[The second proof of \cref{ErgodicDecomp}]
Assume, for simplicity, that $\mu$ is $H$-invariant, and that $X$ is compact. By normalizing, we may assume $\mu(X) = 1$, so $\mu \in \Prob(X)$. Then Choquet's Theorem \pref{ChoquetsThm} provides a probability measure~$\nu$, such that 
	$$ \mu = \int_{\ext \Prob(X)^H} \sigma \, d\nu(\sigma) .$$
By identifying $\ext \Prob(X)$ with a Borel subset of $[0,1]$, we may rewrite this as:
	$$ \mu = \int_{[0,1]} \mu_z \, d\nu(z) ,$$
where $\nu$ is a probability measure on $[0,1]$. Furthermore, \cref{Erg<>Extreme} tells us that each $\sigma \in \ext \Prob(X)^H$ is ergodic, so $\mu_z$ is an ergodic $H$-invariant measure for a.e.~$z$.

All that remains is to define a map $\zeta \colon X \to [0,1]$, such that $\mu_z$ is concentrated on $\zeta^{-1}(z)$.
For each Borel subset~$E$ of $\ext \Prob(X)^H$, let $\mu_E = \int_E \sigma \, d\nu(\sigma)$. Then $\mu_E$ is absolutely continuous with respect to~$\mu$, so we may write $\mu_E = f_E \mu$, for some measurable $f_E \colon X \to [0,\infty)$. Then 
	$$\psi(E) = \{\, x \in X \mid f_E(x) \neq 0 \,\}$$
is a well-defined element of $\Bool(X)$. Therefore, we have defined a map $\psi \colon \Bool \bigl( [0,1] \bigr) \to \Bool(X)$, and it can be verified that this is a homomorphism of $\sigma$-algebras. Hence, there is a measurable function $\zeta \colon X \to [0,1]$, such that $\psi(E) = \zeta^{-1}(E)$, for all~$E$ \csee{BoolPtwise}.
By using the uniqueness of~$\nu$, it can be shown that $\mu_z \bigl( \zeta^{-1}(z) \bigr) = 1$ for a.e.~$z$.
\end{proof}

\begin{exercises}

\item \label{sigma<>closed}
Let $\Bool$ be a (nonempty) subset of $\Bool(X)$ that is closed under complements and finite unions. Show that $\Bool$ is closed under countable unions (so $\Bool$ is a sub-$\sigma$-algebra of $\Bool(X)$) if and only if $\Bool$ is a closed set with respect to the topology determined by the metric on $\Bool(X)$.
\hint{($\Leftarrow$)~If $E = \bigcup_{i=1}^\infty E_i$, then $\bigcup_{i=1}^n E_i \to E$ in the topology on $\Bool(X)$.
\par
($\Rightarrow$)~If $d(E_i,E) < 2^{-i}$, then $E = \bigcap_{n = 1}^\infty \bigcup_{i=n}^\infty E_i$ (up to a set of measure~$0$.}

\item \label{DenseGensSigmaAlg}
Show that if $\mathcal{E}$ is dense in a sub-$\sigma$-algebra~$\Bool$ of $\Bool(X)$, then $\mathcal{E}$ is not contained in any proper sub-$\sigma$-algebra of~$\Bool$.
\hint{If $E_n \to E$, then $\bigcup_{n=1}^\infty (E_n \cap E) = E$ (up to a set of measure~$0$).}

\item \label{Erg<>Extreme}
Prove that a measure $\mu \in \Prob(X)^H$ is ergodic iff it is an extreme point.
\hint{If $E$ is an $H$-invariant set, then $\mu$ is a convex combination of the restrictions to $E$ and its complement.
Conversely, if $\mu = t\, \mu_1 + (1-t)\,\mu_2$, then the Radon-Nikodym Theorem implies $\mu_1 = f \, \mu$ for some ($H$-invariant) function~$f$.}

\item \label{ProbLUB}
Show that any two elements of $\Sigma \Prob(X)$ or $\Sigma \Prob(X)^*$ have a least upper bound.
\hint{Write $\nu = f \, \mu + \nu_s$ (uniquely), where $\nu_s$ is concentrated on a set of measure~$0$.}

\item \label{BoolPtwise}
Suppose $\psi \colon \Bool(Z) \to \Bool(X)$ is a function that respects complements and countable unions (and $\psi(\emptyset) = \emptyset$). Show there is a Borel function $\zeta \colon X \to Z$, such that $\psi(E) = \zeta^{-1}(E)$, for every Borel subset~$E$ of~$Z$.
\hint{Assume, for simplicity, that $Z = \bigl\{\, \sum a_k 3^{-k} \mid a_k \in \{0,1\} \, \bigr\} \subset [0,1]$. Then $\zeta = \sum \chi_{E_k} 3^{-k}$ for an appropriate collection $\{E_k\}$ of Borel subsets of~$X$.}

\end{exercises}

\section{Mixing} \label{MixingSect}

It is sometimes important to know that a product of group actions is ergodic. To discuss this issue (and related matters), let us fix some notation.

\begin{notation}
Throughout this \lcnamecref{MixingSect}, we assume:
	\begin{enumerate}
	\item $X_i$ is a Polish space, for every~$i$,
	\item $H$ is a Lie group that acts continuously on~$X_i$, for each~$i$,
	and
	\item $\mu_i$ is an $H$-invariant probability measure on~$X_i$, for each~$i$.
	\end{enumerate}
Furthermore, we use $X$ and~$\mu$ as abbreviations for $X_0$ and~$\mu_0$, respectively.
\end{notation}

\begin{defns} \ 
	\begin{enumerate}
	\item The \defit[product!action]{product action} of~$H$ on $X_1 \times X_2$ is the $H$-action defined by $h(x_1,x_2) = (hx_1, hx_2)$. The product measure $\mu_1 \times \mu_2$ is an invariant measure for this action.
	\item The action on~$X$ is said to be \defit[mixing!weak or weakly]{weak mixing} (or \defit[mixing!weak or weakly]{weakly mixing}) if the product action on $X \times X$ is ergodic.
	\end{enumerate}
\end{defns}

We have the following very useful characterizations of weakly mixing actions:

\begin{thm} \label{WMIff}
The action of~$H$ on~$X$ is weak mixing if and only if the \textup(one-dimensional\/\textup) space of constant functions is the only nontrivial, finite-dimensional, $H$-invariant subspace of $\LL2(X,\mu)$.
\end{thm}

\begin{proof}
We prove the contrapositive of each direction.

($\Rightarrow$)
Suppose $V$ is a nontrivial, finite-dimensional, $H$-invariant subspace of $\LL2(X)$. If the functions in~$V$ are not all constant (a.e.), then we may assume (by passing to a subspace) that $V \perp 1$. Choose an orthonormal basis $\{\varphi_1,\ldots,\varphi_r\}$ of~$V$, and define $\varphi \colon X \times X \to \complex$ by
	\begin{align} \label{WMIffPf-DefinePhi}
	 \varphi(x,y) = \sum_{i=1}^r \varphi_i(x) \, \overline{\varphi_i(y)} 
	 . \end{align}
Then $\varphi$ is an $H$-invariant function that is not constant \csee{WMIffPfDefinePhiEx}, so $H$ is not ergodic on $X \times X$.

($\Leftarrow$) Suppose $\varphi$ is a nonconstant, $H$-invariant, bounded function on $X \times X$. We may assume $\varphi(x,y) = \overline{\varphi(y,x)}$ by replacing $\varphi$ with either $\varphi(x,y) + \overline{\varphi(y,x)}$ or $\sqrt{-1} \, \bigl( \varphi(x,y) - \overline{\varphi(y,x)} \bigr)$. Therefore, we have a compact, self-adjoint operator on $\LL2(X,\mu)$, defined by
	$$ (T \psi)(x) = \int_X \varphi(x,y) \, \psi(y) \, d\mu(y) .$$
The Spectral Theorem \pref{SpectralThmCpct} implies $T$ has an eigenspace~$V$ that is finite-dimensional (and contains a nonconstant function). This eigenspace is $H$-invariant, since $T$~commutes with~$H$ (because $\varphi$ is $H$-invariant).
\end{proof}

We often have the following stronger condition:

\begin{defn} \label{StrongMixingDefn}
The action of $H$ on~$X$ is said to be \defit{mixing} (or, alternatively, \defit[mixing!strongly]{strongly mixing}) if $H$~is noncompact, and, for all $\varphi,\psi \in \LL2(X,\mu)$, such that $\varphi \perp 1$, we have
	$$ \lim_{h \to \infty} \langle h \varphi \mid \psi \rangle = 0 .$$
(We are using $1$ to denote the constant function of value~$1$, and, as usual, $(h \varphi)(x) = \varphi(h^{-1} x)$ \csee{RepL2}.)
\end{defn}

\begin{prop}[\csee{MixingIffEx}] \label{MixingIff}
If $H$ is not compact, then the following are equivalent:
	\begin{enumerate}
	\item \label{MixingIff-mixing}
	The action of $H$ on~$X$ is mixing.
	\item \label{MixingIff-funcs}
	For all $\varphi,\psi \in \LL2(X,\mu)$, we have
	$$ \lim_{h \to \infty} \langle h \varphi \mid \psi \rangle = \langle \varphi \mid 1 \rangle \, \langle 1 \mid \psi \rangle .$$
	\item \label{MixingIff-sets}
	For all Borel subsets $A$ and~$B$ of~$X$, we have
		$$ \lim_{h \to \infty} \mu(hA \cap B) = \mu(A) \, \mu(B) .$$ 
	\end{enumerate}
\end{prop}

\begin{rem}
Condition~\pref{MixingIff-sets} is the motivation for the choice of the term ``mixing:'' as $h \to \infty$, the space~$X$ is getting so stirred up (or well-mixed) that $hA$ is becoming uniformly distributed throughout the entire space.
\end{rem} 

When $G$ is simple, decay of matrix coefficients \pref{DecayMatCoeffSimple} implies that every action of~$G$ (with finite invariant measure) is mixing. In fact, we can say more.

\begin{defn}
Generalizing \cref{StrongMixingDefn}, we say that the action of~$H$ on~$X$ is \defit[mixing!of order~$r$]{mixing of order~$r$} if $H$ is not compact, and, for all Borel subsets $A_1,\ldots,A_r$ of~$X$, we have
	$$ \lim_{h_i^{-1} h_j \to \infty} \mu \left( \bigcap_{i=1}^r h_i A_i \right)
		= \mu(A_1) \, \mu(A_2) \, \cdots \, \mu(A_r) .$$
In particular, 
	\begin{itemize}
	\item every action of~$H$ is mixing of order~$1$ (if $H$ is not compact),
	and
	\item ``mixing'' is the same as ``mixing of order~$2$\zz.''
	\end{itemize}
\end{defn}

\begin{warn}
Some authors use a different numbering, for which this is ``mixing of order $r-1$\zz,'' instead of ``mixing of order $r$\zz.'' 
\end{warn}

Ledrappier constructed an action of $\integer^2$ that is mixing of order~$2$, but not of order~$3$. However, there are no such examples for semisimple groups:

\begin{thm} \label{GMixing}
Every mixing action of~$G$ \textup(with finite invariant measure\/\textup) is mixing of all orders.
\end{thm}

In the special case where $H = \integer$, we mention the following additional characterizations, some of which are weaker versions of \fullcref{MixingIff}{sets}:

\begin{thm} \label{WeakMixingZ}
If $H = \langle T \rangle$ is an infinite cyclic group, then the following are equivalent:
	\begin{enumerate}
	
	\item \label{WeakMixingZ-wm}
	$H$~is weak mixing on~$X$.

	\item \label{WeakMixingZ-eig}
	Every eigenfunction of~$T$ in $\LL2(X,\mu)$ is constant\/ \textup(a.e.\textup). That is, if $f \in \LL2(X,\mu)$, and there is some\/ $\lambda \in \complex$, such that $f(Tx) = \lambda \, f(x)$ a.e., then $f$~is constant\/ \textup(a.e.\textup).
			
	\item \label{WeakMixingZ-spectrum}
	The spectral measure of~$T$ has no atoms other than~$1$, and the eigenvalue~$1$ is simple\/ \textup(that is, the corresponding eigenspace is\/ $1$-dimensional\/\textup). 
%That is, if we let $\nu$ be the spectral measure of the unitary operator~$T$ on $\LL2(X,\mu)$ \csee{SpectralMeasDefn}, then $\nu(\{z\}) = 0$, for every singleton set~$\{z\}$ in~$\torus$.

	\item \label{WeakMixingZ-XxAny}
	The action of~$H$ on $X \times X_1$ is ergodic, whenever the action of~$H$ on~$X_1$ is ergodic.

%	\item \label{WeakMixingZ-IndepSets}
%The weakly independent sets are dense in $\Bool(X)$, where a set~$A$ is said to be \defit[weakly independent set]{\emph{weakly independent}} if there exist $n_1 < n_2 < \cdots$, such that 
%		$$ \text{$\mu \bigl( T^{n_i} A \cap T^{n_j} A \bigr) = \mu(A)^2$ whenever $i \neq j$} . $$
		
%	\item \label{WeakMixingZ-fulldensity}
%For any Borel subsets $A$ and~$B$ of~$X$, there is a subset $\mathcal{N}$ of full density in~$\integer^+$, such that 
%	$$ \mu(T^k A \cap B) \longrightarrow \mu(A) \, \mu(B)  
%	\qquad \text{as $k \to \infty$ with $k \in \mathcal{N}$}
%	.$$
%To say $\mathcal{N}$ has \defit[full density]{\emph{full density}} means 
%	$$ \lim_{n \to \infty} \frac{\# \bigl( \mathcal{N} \cap \{1,2,3,\ldots,n\} \bigr)}{n} = 1 .$$

	\item \label{WeakMixingZ-AbsVal}
For all Borel subsets $A$ and~$B$ of~$X$, 
		$$ \lim_{n \to \infty} \frac{1}{n} \sum_{k=1}^n 
			\bigl| \mu( T^k A \cap B) - \mu(A) \, \mu(B) \bigr| = 0 .$$

	\item \label{WeakMixingZ-fulldensitySameset}
%The quantifiers in \pref{WeakMixingZ-fulldensity} can be reversed:
There is a subset $\mathcal{N}$ of full density in~$\integer^+$, such that, for all Borel subsets $A$ and~$B$ of~$X$, we have
	$$ \mu(T^k A \cap B) \longrightarrow \mu(A) \, \mu(B)  
	\qquad \text{as $k \to \infty$ with $k \in \mathcal{N}$}
	.$$
To say $\mathcal{N}$ has \defit[full!density]{\emph{full density}} means 
	$$ \lim_{n \to \infty} \frac{\# \bigl( \mathcal{N} \cap \{1,2,3,\ldots,n\} \bigr)}{n} = 1 .$$
		
	\item \label{WeakMixingZ-multi}
If $A_0,A_1,\dots,A_r$ are any Borel subsets of~$X$, then there is a subset $\mathcal{N}$ of full density in~$\integer^+$, such that
		$$ \lim_{\begin{smallmatrix}k \to \infty \\ k \in \mathcal{N}\end{smallmatrix}} 
			\mu \bigl(A_0 \cap T^kA_1 \cap T^{2k} A_2 \cap \cdots \cap T^{rk} A_r \bigr)
		= \mu(A_0) \, \mu(A_1) \cdots \mu(A_r)
		%\quad \text{as $k \to \infty$ with $k \in \mathcal{N}$}
		.$$

	\end{enumerate}
\end{thm}

\begin{proof}[Sketch of proof]
($\ref{WeakMixingZ-wm} \Leftrightarrow \ref{WeakMixingZ-eig}$) By \cref{WMIff}, it suffices to observe that every finite-dimensional representation contains an irreducible subrepresentation, and that the irreducible representations of~$\integer$ (or, more generally, of any abelian group) are one-dimensional.

($\ref{WeakMixingZ-eig} \Leftrightarrow \ref{WeakMixingZ-spectrum}$)
These are two different ways of saying the same thing.

($\ref{WeakMixingZ-spectrum} \Rightarrow \ref{WeakMixingZ-XxAny}$)
Since $\LL2(X \times X_1) \iso \LL2(X) \otimes \LL2(X_1)$, the spectral measure~$\nu$ of $\LL2(X \times X_1)$ is the product $\nu_1 \times \nu_2$ of the spectral measures of $\LL2(X)$ and $\LL2(X_1)$. Therefore, any point mass in $\nu$ is obtained by pairing a point mass in~$\nu_1$ with a point mass in~$\nu_2$.

($\ref{WeakMixingZ-XxAny} \Rightarrow \ref{WeakMixingZ-wm}$)
Take $X_1 = X$.

($\ref{WeakMixingZ-wm} \Rightarrow \ref{WeakMixingZ-AbsVal}$)
For simplicity, let $a = \mu(A)$ and $b = \mu(B)$.
By \cref{MeanErgThmEx} (the Mean Ergodic Theorem), % @@@
ergodicity on~$X$ implies
	$$  \frac{1}{n} \sum_{k=1}^n \mu(T^k A \cap B) \stackrel{n \to \infty}{\longrightarrow} ab .$$
For the same reason, ergodicity on $X \times X$ implies
	\begin{align*}
	\sum_{k=1}^n (\mu \times \mu) \bigl (T^k A \times T^k A) \cap (B \times B) \bigr)
	\stackrel{n \to \infty}{\longrightarrow} 
	(\mu \times \mu)(A \times A) \cdot (\mu \times \mu) ( B \times B) 
	. \end{align*}
By simplifying both sides, we see that
	$$ \frac{1}{n} \sum_{k=1}^n \mu(T^k A \cap B)^2 
	\stackrel{n \to \infty}{\longrightarrow} 
	a^2 b^2 .$$
Therefore, simple algebra yields
	$$  \frac{1}{n} \sum_{k=1}^n \bigl| \mu(T^k A \cap B) - \mu(A) \, \mu(B) \bigr|^2
	\ \stackrel{n \to \infty}{\longrightarrow}  \ 
	a^2 \,b^2 - 2 (ab)(ab) + (ab)^2
	= 0
	. $$
\Cref{CesaroIff} implies that we have the same limit without squaring the absolute value.

($\ref{WeakMixingZ-AbsVal} \Rightarrow \ref{WeakMixingZ-eig}$)
Approximating by linear combinations of characteristic functions implies
	$$ \lim_{n \to \infty} \frac{1}{n} \sum_{k=1}^n \bigl| \langle T^k \varphi \mid \varphi \rangle  \bigr|
	= 0 
	\qquad \text{for all $\varphi \perp 1$} 
	. $$
However, if $\varphi$ is an eigenfunction for an eigenvalue $\neq 1$, then it is easy to see that the limit is nonzero (either directly, or by applying \cref{CesaroIff}).

($\ref{WeakMixingZ-AbsVal} \Leftrightarrow \ref{WeakMixingZ-fulldensitySameset}$)
\Cref{CesaroIff} implies that the two assertions are equivalent, up to reversing the order of the quantifiers in~\pref{WeakMixingZ-fulldensitySameset}. To reverse the quantifiers, note that a variant of Cantor diagonalization provides a set~$\mathcal{N}$ of full density that works for all $A$ and~$B$ in a countable dense subset of $\Bool(X)$.

 ($\ref{WeakMixingZ-multi} \Rightarrow \ref{WeakMixingZ-fulldensitySameset}$)
Take $r = 1$.

 ($\ref{WeakMixingZ-fulldensitySameset} \Rightarrow \ref{WeakMixingZ-multi}$)
The proof proceeds by induction on~$r$ (with \pref{WeakMixingZ-fulldensitySameset} as the starting point), and is nontrivial. We have no need for this result, so we omit the proof. % @@@
\end{proof}

\begin{exercises}

\item Show (directly from the definitions) that if the action of~$H$ on~$X$ is weak mixing, then it is ergodic.

\item \label{WMIffPfDefinePhiEx}
Let $\varphi \colon X \times X \to \complex$ be as in \pref{WMIffPf-DefinePhi} of the proof of \cref{WMIff}.
	\begin{enumerate}
	\item \label{WMIffPfDefinePhiEx-invt}
	Show $\varphi$ is $H$-invariant (a.e.).
	\item \label{WMIffPfDefinePhiEx-NotConst}
	Show $\varphi$ is not constant (a.e.). 
	\end{enumerate}
\hint{\pref{WMIffPfDefinePhiEx-invt}~Write $h \, \varphi_i = \sum_{i,j} h_{i,j} \varphi_j$, and observe that $[h_{i,j}]$ is a unitary matrix.
\\ 
\pref{WMIffPfDefinePhiEx-NotConst}~$\varphi(x,x) > 0$, but $\int \varphi(x,y) \, d\mu(y) = 0$.}

\item \label{MixingIffEx}
Prove \cref{MixingIff}.
\hint{($\ref{MixingIff-mixing} \Rightarrow \ref{MixingIff-funcs}$)
For $c = \langle \varphi \mid 1 \rangle$, we have $\langle (\varphi - c) \mid 1 \rangle = 0$. Now calculate $\lim_{h \to \infty} \langle h (\varphi - c) \mid \psi \rangle$ in two ways.
%		= \lim_{h \to \infty} \langle h \varphi \mid \psi \rangle 
%			- c \langle 1 \mid \psi \rangle .$$
\\ ($\ref{MixingIff-funcs} \Rightarrow \ref{MixingIff-sets}$)
Let $\varphi$ and~$\psi$ be the characteristic functions of $A$ and~$B$.
\\ ($\ref{MixingIff-sets} \Rightarrow \ref{MixingIff-funcs}$)
Approximate $\varphi$ and~$\psi$ by linear combinations of characteristic functions.}

\item \label{CesaroIff}
For every bounded sequence $\{a_k\} \subset [0,\infty)$, show
	$$ \lim_{n \to \infty} \frac{1}{n} \sum_{k=1}^n a_k = 0
	\ \Leftrightarrow \ 
	\text{$a_k \to 0$ as $k \to \infty$ in some set of full density.} $$
\hint{($\Rightarrow$)~For each $m > 0$, the set $A_m = \{\, k \mid a_k > 1/m\,\}$ has density~$0$, so there exists $N_m > N_{m-1}$, such that, for all $n \ge N_m$, we have
	$$N_{m-1} + \#\bigl( A_m \cap \{1,2,\ldots,n\} \bigr) < n/m .$$
Let $\mathcal{N}$ be the complement of $\bigcup_m \bigl( A_m \cap [N_m, N_{m+1}) \bigr)$.}

\end{exercises}

\begin{notes}

The focus of classical Ergodic Theory is on actions of $\integer$ and~$\real$ (or other abelian groups). A few of the many introductory books on this subject are \cite{EinsiedlerWard-ErgThyViewNumThy, Halmos, KatokHasselblatt-intro, Walters}. They include proofs of the Poincar\'e Recurrence Theorem \pref{PoincareRecurThm} and the Pointwise Ergodic Theorem \pref{PointwiseErgThm}.

Some basic results on the Ergodic Theory of noncommutative groups can be found in \cite[\S2.1]{ZimmerBook}.

The Moore Ergodicity Theorem \pref{MooreErgodicity} is due to C.\,C.\,Moore \cite{Moore-ergodicity}. 

See \cite{Lindenstrauss-Ptwise} for a very nice version of the Pointwise Ergodic Theorem that applies to all amenable groups (\cref{{PtwiseErgAmenRem}}).

See \cite[Thm.~1.1 (and Thm.~5.2)]{GreschonigSchmidt} for a proof of the ergodic decomposition \pref{ErgodicDecomp}, using Choquet's Theorem as in \cref{ErgDecompPfChoquet}.
\Cref{RokhlinDecompMeas} can be found in \cite[\S3]{Rohlin-FundIdeas}.
See \cite[\S5.5]{Srivastava-BorelSets} for a proof of \cref{vonNeumannSelectionThm}.
%A version of \cref{vonNeumannSelectionThm} appears in \cite[Thm.~3.4.3, p.~77]{Arveson-InvitationCstar}.

See \cite[Prop.~1, pp.~157--158]{Moore-ergodicity} for a proof of \cref{MixingIff}.

The standard texts on ergodic theory only prove \cref{MixingIff} for the special case $H = \integer$, but the same arguments apply in general.

\Cref{GMixing} is due to S.\,Mozes \cite{Mozes-MixingAllOrders}. Ledrappier's counterexample for $\integer^2$ is in 
\cite{Ledrappier-ChampMarkovien}.

\Cref{WeakMixingZ} is in the standard texts on ergodic theory, except for Part~\pref{WeakMixingZ-multi}, which is a ``\thmindex{multiple recurrence}multiple recurrence theorem'' that plays a key role in Furstenberg's proof of \thmindex{Szemeredi's}Szemeredi's theorem that there are arbitrarily long arithmetic progressions in every set of positive density in~$\integer^+$. For a proof of ($\ref{WeakMixingZ-fulldensitySameset} \Rightarrow \ref{WeakMixingZ-multi}$), see 
\cite[Prop.~7.13, p.~191]{EinsiedlerWard-ErgThyViewNumThy} 
or 
\cite[Thm.~4.10]{Furstenberg-RecErgThyCombNumThy}.

A proof of \cref{CesaroIff} is in \cite[Lem.~2.41, p.~54]{EinsiedlerWard-ErgThyViewNumThy}.
\end{notes}

%!TEX root = IntroArithGrps.tex

\part{Major Results} \label{MajorPart}

 %!TEX root = IntroArithGrps.tex

\mychapter{Mostow Rigidity Theorem}
\label{MostowChap}

\prereqs{Quasi-isometries (\cref{QuasiChap}).}

%\begin{narrower} \hyphenpenalty=3000
%\em The Mostow Rigidity Theorem plays a fundamental role in nearly every paper on the geometry of Lie groups, and the techniques and ideas [Mostow] % he @@@
% introduced have influenced a lot of further work in geometry, group theory and dynamics.
%\\ \hbox{ } \hfil --- Yair Minsky, Yale University\break
%% \hbox{ } \hfil on the awarding of the Wolf Prize to G.\,D.\,Mostow in January 2013\break
%\par \end{narrower}

\section{Statement of the theorem}

In its simplest form, the Mostow Rigidity Theorem says that a single group~$\Gamma$ cannot be a lattice in two different semisimple groups~$G_1$ and~$G_2$ (except for  minor modifications involving compact factors, the center, and passing to a finite-index subgroup):

\begin{thm}[(Weak version of the Mostow Rigidity Theorem)] \label{MostowIso}
	\thmindex{Mostow Rigidity!weak version}
Assume
\noprelistbreak
	\begin{itemize}
	\item $G_1$ and~$G_2$ are connected, with trivial center and no compact factors,
	and
	\item $\Gamma_i$ is a lattice in~$G_i$, for $i = 1,2$.
	\end{itemize}
If\/ $\Gamma_1 \iso \Gamma_2$, then $G_1 \iso G_2$.
\end{thm}

In other words, if there is an isomorphism from~$\Gamma_1$ to~$\Gamma_2$, then there is also an isomorphism from~$G_1$ to~$G_2$. In fact, it is usually the case that something much stronger is true: 
%if no simple factor of~$G_1$ is isogenous to $\SL(2,\real)$, then every isomorphism between~$\Gamma_1$ and~$\Gamma_2$ actually comes from an isomorphism between $G_1$ and~$G_2$:

\begin{namedthm}[Mostow Rigidity Theorem]
\label{MostowRigidity}
Assume
\noprelistbreak
	\begin{itemize}
	\item $G_1$ and~$G_2$ are connected, with trivial center and no compact factors,
	\item $\Gamma_i$ is a lattice in~$G_i$, for $i = 1,2$,
	and
	\item there does not exist a simple factor~$N$ of~$G_1$, such that $N \iso \PSL(2,\real)$ and $N \cap \Gamma_1$ is a lattice in~$N$.
	\end{itemize}
Then any isomorphism from~$\Gamma_1$ to~$\Gamma_2$ extends
to a continuous isomorphism from~$G_1$ to~$G_2$.
 \end{namedthm}

\begin{rems}  \ \label{MostowRems}
\noprelistbreak
	\begin{enumerate}

	\item \label{MostowRems-deformationrigidity}
	Assume $G$ is connected, and has no simple factors that are either compact or isogenous to $\SL(2,\real)$. Then the Mostow Rigidity Theorem implies that lattices in~$G$ have no nontrivial deformations. More precisely, if $\Gamma_t$ is a continuous family of lattices in~$G$, then $\Gamma_t$ is conjugate to~$\Gamma_0$, for every~$t$ \csee{deformationrigidityEx}. 

This is not always true when $G$ is isogenous to $\SL(2,\real)$ \csee{ModuliSpaceSL2}, which explains why the statement of the Mostow Rigidity Theorem must forbid factors that are isogenous to $\SL(2,\real)$.

	\item In geometric terms, the Mostow Rigidity Theorem tells us that the topological structure of any irreducible finite-volume locally symmetric space of noncompact type completely determines its geometric structure as a Riemannian manifold (up to multiplying the metric by a scalar on each irreducible factor of the universal cover), if the manifold is not $2$-dimensional \ccf{MostowGeomEx}.

	\item \label{MostowRems-FurmanRigidity}
%\Cref{MostowIso} tells us that if we ignore some obvious minor modifications coming from compact groups, then $\Gamma$ does not embed as a lattice in any connected, semisimple Lie group other than~$G$. However, $\Gamma$~is also a lattice in a quite different Lie group; namely, $\Gamma$~is a lattice in itself.
%
	Assume $\Gamma$ is cocompact in~$G$. Then, as a strengthening of \cref{MostowIso}, it can be shown that if $\Gamma$ is a cocompact lattice in some Lie group~$H$ (not assumed to be semisimple), then $H$ must be either $G$ or~$\Gamma$ (modulo the usual minor modifications involving compact groups). In fact, this remains true even if we allow~$H$ to be any locally compact group, not necessarily a Lie group.

	\end{enumerate}
\end{rems}

\begin{exercises}

\item Suppose 
	\begin{itemize}
	\item $G$ has trivial center and no compact factors,
	\item  $G \not\iso \PSL(2,\real)$,
	and
	\item $\Gamma$ is irreducible. 
	\end{itemize}
Show that every automorphism of~$\Gamma$ extends to a continuous automorphism of~$G$.

\item Show that \cref{MostowIso} is a corollary of \cref{MostowRigidity}.
\hint{This is obvious when no simple factor of~$G_1$ is isomorphic to $\PSL(2,\real)$. The problem can be reduced to the case where $G_1$ and~$G_2$ are irreducible.}

\item \label{NotIsoForCpctFactors}
Assume $G_1$ and~$G_2$ are isogenous. Show that if $K$ is any compact group, then some lattice in~$G_1$ is isomorphic to a lattice in $G_2 \times K$ (even though $G_1$ may not be isomorphic to $G_2 \times K$). This is why \cref{MostowIso} assumes $G_1$ and~$G_2$ are connected, with trivial center and no compact factors
\hint{Any torsion-free lattice in~$G_1^\circ$ is isomorphic to a lattice in $G_2 \times K$.}

\item For $i = 1,2$, suppose 
	\begin{itemize}
	\item $\Gamma_i$ is a lattice in~$G_i$, 
	and 
	\item $G_i$ is connected and has no compact factors. 
	\end{itemize}
Show that if $\Gamma_1$ is isomorphic to~$\Gamma_2$, then $G_1$ is isogenous to~$G_2$.

\item Let $\Gamma_i$ be a lattice in~$G_i$ for $i = 1,2$. Show that if $\Gamma_1 \iso \Gamma_2$, then there is a compact, normal subgroup~$K_i$ of~$G_i^\circ$, for $i = 1,2$, such that $G_1^\circ/K_1 \iso G_2^\circ/K_2$.

%\item Let $G_1 = \SL(4,\real)$ and $G_2 = \PSL(4,\real)$. 
%\begin{enumerate}
%\item Show $G_1 \not\iso G_2$.
%\item Show there exist lattices $\Gamma_1$ and~$\Gamma_2$ in $G_1$ and~$G_2$, respectively, such that $\Gamma_1 \iso \Gamma_2$.
%\item Why is this not a counterexample to \cref{MostowRigidity}?
%\end{enumerate}

\item Let $G = \PSL(2,\real)$. 
Find an automorphism~$\varphi$ of some lattice~$\Gamma$ in~$G$, such that $\varphi$ does not extend to an automorphism of~$G$.
\par
Why is this not a counterexample to \cref{MostowRigidity}?

%\item Let $G_1 = \PSO(2,3)$ and $G_2 = \PSO(2,3) \times \PSO(5)$. 
%\begin{enumerate}
%\item Show $G_1 \not\iso G_2$.
%\item Show there exist irreducible lattices $\Gamma_1$ and~$\Gamma_2$ in $G_1$ and~$G_2$, respectively, such that $\Gamma_1 \iso \Gamma_2$.
%\item Why is this not a counterexample to \cref{MostowRigidity}?
%\end{enumerate}

%\item In the statement of \cref{MostowRigidity}:
%\begin{enumerate}
%\item Show the assumption that $\Gamma_1$ is irreducible can be replaced with the assumption that $\Gamma_2$ is irreducible.
%\item Show the assumption that $\Gamma_1$ is irreducible can be replaced with assumption that no simple factor of~$G_1$ is isomorphic to $\PSL(2,\real)$.
%\end{enumerate}

%\item \label{TrivialDeformEx}
%Verify that the construction of \cref{TrivialDeformLatt} defines a continuous deformation of~$\Gamma$.

\item \label{deformationrigidityEx}
Assume there is a continuous function $\rho \colon \Gamma \times [0,1] \to G$, such that, if we let $\rho_t(\gamma) = \rho(\gamma,t)$, then:
 	\begin{itemize}
	\item $\rho_t$ is a homomorphism, for all~$t$,
	\item $\rho_t(\Gamma)$ is a lattice in~$G$, for all~$t$,
	and
	\item $\rho_0(\gamma) = \gamma$, for all~$\gamma \in \Gamma$.
	\end{itemize}
Show that if $G$ is as in the first sentence of \fullcref{MostowRems}{deformationrigidity}, then $\Gamma_t$ is conjugate to~$\Gamma$, for every~$t$.
\hint{Reduce to the case where $G$ has trivial center. You may use, without proof, the fact that the identity component of the automorphism group of~$G$ consists of inner automorphisms \csee{OutGFinite}.}

\item \label{MostowGeomEx}
(\emph{Requires some familiarity with locally symmetric spaces})
Assume, for $i = 1,2$:
	\begin{itemize}
	\item $G_i$ is connected and simple, with trivial center,
	\item $K_i$ is a maximal compact subgroup of~$G_i$,
%	\item $X_i = K_i \backslash G_i$ is the symmetric space associated to~$G_i$,
	\item $\Gamma_i$ is a torsion-free, irreducible lattice in~$G_i$,
	\item $X_i = K_i \backslash G_i / \Gamma_i$ is the corresponding locally symmetric space of finite volume,
	\item the metric on~$X_i$ is normalized so that $\vol(X_1) = 1$,
	and
	\item $\dim X_1 \ge 3$.
	\end{itemize}
Show that any homotopy equivalence from~$X_1$ to~$X_2$ is homotopic to an isometry.
\hint{Since the universal cover $K_i \backslash G_i$ is contractible, a homotopy equivalence is determined, up to homotopy, by its effect on the fundamental group.}

\end{exercises}

\section{Sketch of the proof for \texorpdfstring{$\SO(1,\lowercase{n})$}{SO(1,n)} 
\texorpdfstring{\optional}{(optional)}} \label{MostowPfSect}

In most cases, the conclusion of the Mostow Rigidity Theorem \pref{MostowRigidity} is an easy consequence of the \thmindex{Margulis!Superrigidity}{Margulis Superrigidity Theorem}, which will be discussed in \cref{MargulisSuperChap}.
More precisely, if we assume, for simplicity, that the lattices $\Gamma_1$ and~$\Gamma_2$ are irreducible \csee{MostowRigIrredEnough}, then the Margulis Superrigidity Theorem applies unless $G_1$ and~$G_2$ are isogenous to either $\SO(1,n)$ or $\SU(1,n)$ \csee{ProveMostMostowEx}.
To illustrate the main ideas involved in completing the proof, we discuss a special case:

\begin{proof}[\mathversion{bold}Proof  of Mostow Rigidity Theorem for cocompact lattices in $\SO(1,n)$]
Assume, for $i = 1,2$:
\noprelistbreak
	\begin{itemize}
	\item $G_i \iso \PSO(1,n_i)$, for some $n_i \ge 3$,
	\item $\Gamma_i$ is cocompact in~$G_i$,
	and
	\item $\rho \colon \Gamma_1 \to \Gamma_2$ is an isomorphism.
	\end{itemize}
In order to show that $\rho$ extends to a continuous homomorphism from $G_1$ to~$G_2$, we take a geometric approach that uses the action of $G_i$ on its associated symmetric space, which is the hyperbolic space~$\hyperbolic^{n_i}$.
This assumes some understanding of hyperbolic space (and other matters) that is not required elsewhere in this book.

\begin{claim}
We have $n_1 = n_2$.
\end{claim}
By passing to subgroups of finite index, we may assume $\Gamma_1$ and~$\Gamma_2$ are torsion free, so $\Gamma_i$ acts freely on~$\hyperbolic^{n_i}$. The action is also properly discontinuous, so, since hyperbolic space is contractible, this implies that $X_i = \Gamma_i \backslash \hyperbolic^{n_i}$ is a $K(\Gamma_i, 1)$-space. Since $X_i$ is a compact manifold of dimension~$n_i$, we conclude that the cohomological dimension of~$\Gamma_i$ is~$n_i$. However, the groups $\Gamma_1$ and~$\Gamma_2$ are isomorphic, so they must have the same cohomological dimension. This completes the proof of the claim.
\qed

\medbreak

Therefore, $\Gamma_1$ and~$\Gamma_2$ are two lattices in the same group $G = \PSO(1,n)$. To simplify matters, let us assume $n = 3$.

Since $\Gamma_i$ is cocompact in~$G_i$ (so $\Gamma_i$ is quasi-isometric to~$\hyperbolic^3$), it is not difficult to construct a quasi-isometry $\varphi \colon \hyperbolic^3 \to \hyperbolic^3$, such that
	\begin{align} \label{MostowPfSO13Equi}
	\text{$\varphi(\gamma x) = \rho(\gamma) \cdot \varphi(x)$ 
	\ for all $\gamma \in \Gamma_1$ and $x \in \hyperbolic^3$} 
	\end{align}
\csee{MostowPfSO13EquiEx}.

Consider the ball model of~$\hyperbolic^3$, whose boundary $\bdry \hyperbolic^3$ is the round $2$-sphere~$S^2$. It is easy to see that any isometry~$\phi$ of~$\hyperbolic^3$ induces a well-defined homeomorphism $\overline\phi$ of $\bdry \hyperbolic^3$. 
Furthermore, it is well known that this boundary map is \defit{conformal} (i.e., it is angle-preserving). This implies that if $C$ is a very small circle in $\bdry \hyperbolic^3$, then $\overline\phi(C)$ is very close to being a circle; more precisely, if we let $S_r(p)$ be the sphere of radius~$r$ around the point~$p$, then
	$$ \limsup_{r \to 0^+} \frac{\sup_{x \in S_r(p)} d \bigl( \overline\phi(x), \overline\phi(p) \bigr)}{\inf_{y \in S_r(p)} d \bigl( \overline\phi(y), \overline\phi(p) \bigr)}
	= 1 .$$
With this background in mind, it should not be difficult to believe (and it is not terribly difficult to prove) that the quasi-isometry $\varphi$ induces a well-defined homeomorphism $\overline\varphi$ of $\bdry \hyperbolic^3$, such that 
	\begin{align} \label{MostowPfSO13OverlineVarphi}
	\text{$\overline\varphi(\gamma p) = \overline{\rho(\gamma)} \cdot \overline\varphi(p)$ 
	\ for all $\gamma \in \Gamma_1$ and $p \in \bdry\hyperbolic^3$} 
	. \end{align}
Furthermore, this boundary map is \defit[quasi-!conformal]{quasi-conformal}, which means that if $C$ is a very small circle in $\bdry \hyperbolic^3$, then $\overline\varphi(C)$ is approximated by an ellipse of bounded eccentricity; more precisely, there is some constant $\kappa > 0$, such that, for all $p \in \bdry \hyperbolic^3$, we have
	$$ \limsup_{r \to 0^+} \frac{\sup_{x \in S_r(p)} d \bigl( \overline\varphi(x), \overline\varphi(p) \bigr)}{\inf_{y \in S_r(p)} d \bigl( \overline\varphi(y), \overline\varphi(p) \bigr)}
	< \kappa .$$ 

It is a fundamental, but highly nontrivial, fact that quasi-conformal maps are differentiable almost everywhere. Therefore, if $C_p$ is a circle (centered at the origin) in the tangent plane at almost any point $p \in \bdry\hyperbolic^3$, then $\overline\varphi(C_p)$ is an ellipse in the tangent plane at $\overline\varphi(p)$. Furthermore, since multiplying all the vectors in a tangent plane by a scalar does not change the eccentricity of any ellipse in the plane, we see that the eccentricity $e_p$ of this ellipse is independent of the choice of the circle~$C_p$. So~$e_p$ is a well-defined, measurable function on (almost all of)~$\bdry\hyperbolic^3$.

\setcounter{case}{0}

\begin{case}
Assume $e_p = 1$ for almost all $p \in \bdry\hyperbolic^3$.
\end{case}
This implies that the quasi-conformal map~$\overline\varphi$ is actually conformal. So there is an isometry~$\alpha$ of~$\hyperbolic^3$, such that $\overline\alpha = \overline\phi$. Then, for $\gamma \in \Gamma_1$ and $p \in \bdry\hyperbolic^3$, we have
	$$ \overline\alpha( \gamma p)
	= \overline\varphi( \gamma p)
	= \overline{\rho(\gamma)} \cdot \overline\varphi(p)
	= \overline{\rho(\gamma)} \cdot \overline\alpha(p)
	. $$
On the other hand, since $G$ is the identity component of $\Isom(\hyperbolic^3)$, it is normalized by~$\alpha$, so there is an automorphism $\widehat\alpha$ of~$G$, such that, for all $g \in G$ and $p \in \bdry\hyperbolic^3$, we have
	$$ \overline\alpha(gp) = \overline{\widehat\alpha(g)} \cdot \overline\alpha(p) .$$
By comparing the displayed equations (and letting $g = \gamma$), we see that $\widehat\alpha(\gamma) = \rho(\gamma)$ for all $\gamma \in \Gamma_1$. Therefore, $\widehat\alpha$ is the desired extension of~$\rho$ to an isomorphism defined on all of~$G$.

\begin{case} \label{MostowPfSO13Not1}
Assume $e_p$ is not almost always equal to~$1$.
\end{case}
Since $\overline{\rho(\gamma)}$ is conformal (because $\rho(\gamma)$ is an isometry), we see from \pref{MostowPfSO13OverlineVarphi} that $e_{\gamma p} = e_p$ for all $\gamma \in \Gamma_1$ and (almost) all $p \in \bdry\hyperbolic^3$. However, since $\overline{G}$ is transitive on $\bdry\hyperbolic^3$ with noncompact point-stabilizers, the Moore Ergodicity Theorem \pref{GammaErgOnG/H} implies that $\overline{\Gamma_1}$ is ergodic on $\bdry\hyperbolic^3$, so we conclude that the function $e_p$ is  constant (a.e). 

Thus, the assumption of this \lcnamecref{MostowPfSO13Not1} implies that, for almost every~$p$, the ellipse $\overline\varphi(C_p)$ is not a circle, and therefore has a well-defined major axis~$\ell_p$, which is a line through the origin in the tangent plane at~$\overline\varphi(p)$. Hence, $\{\ell_p\}$ is a (measurable) section of a certain bundle $\mathbf{P}\hyperbolic^3$, namely, the $\real \mathbf{P}^1$-bundle over $\hyperbolic^3$ whose fiber at each point is the projectivization of the tangent space. Furthermore, since $\overline{\Gamma_2} = \overline{\rho(\Gamma_1)}$ is conformal, we see from \pref{MostowPfSO13OverlineVarphi} that this section is $\overline{\Gamma_2}$-invariant. In fact, if we rotate $\{\ell_p\}$ by any angle~$\theta$, then the resulting section $\{\ell_p^\theta\}$  is also invariant (since $\overline{\Gamma_2}$ is conformal). Then $\bigcup_{0 \le \theta \le \pi/2} \{\ell_p^\theta\}$ is a measurable, $\overline{\Gamma_2}$-invariant subset of $\mathbf{P}\hyperbolic^3$. However, the Moore Ergodicity Theorem \pref{GammaErgOnG/H} implies there are no such (nontrivial) subsets, since $\overline{G}$ is transitive on $\mathbf{P}\hyperbolic^3$ with noncompact point-stabilizers. This contradiction completes the proof, by showing that this \lcnamecref{MostowPfSO13Not1} does not occur.
\end{proof}

\begin{rems} \label{MostowPfRems} \ 
\noprelistbreak
	\begin{enumerate}
	\item The above proof makes the simplifying assumption that $n = 3$. Essentially the same proof works for larger values of~$n$, but $\overline\varphi(C_p)$ will be an ellipsoid, rather than an ellipse, so the space $\ell_p$ of major directions may be a higher-dimensional subspace of the tangent space, instead of just being a line.

	\item \label{MostowPfRems-QuasiC}
	Mostow was able to modify the proof to deal with $\SU(1,n)$, instead of $\SO(1,n)$, by developing a theory of  maps that are quasiconformal over~$\complex$. (In fact, replacing $\complex$ with the quaternions and octonions yields proofs for lattices in  the other simple groups of real rank one, namely, $\Sp(1,n)$ and $F_{4,1}$. However, this is not necessary, because the Margulis Superrigidity Theorem applies to these groups.)

	\item \label{MostowPfRems-Qrank1}
	For lattices in $\SO(1,n)$ that are \textbf{not} cocompact, it is not at all obvious that an isomorphism $\Gamma_1 \iso \Gamma_2$ should yield a quasi-isometry $\hyperbolic^n \to \hyperbolic^n$. This was proved by G.\,Prasad, by using the ``\term{Siegel set}'' description of a coarse fundamental domain for the action of~$\Gamma_i$ on $\hyperbolic^n$ \ccf{ReductionChap}. 
	The same method also works for noncocompact lattices in $\SU(1,n)$, or, more generally, whenever $\Qrank \Gamma_1 = 1$.
	\end{enumerate}
\end{rems}

\begin{exercises}

\item \label{MostowRigIrredEnough}
Show that the proof of \cref{MostowRigidity} can be reduced to the special case where the lattices $\Gamma_1$ and~$\Gamma_2$ are irreducible in $G_1$ and~$G_2$, respectively.
In other words, assume the conclusion of \cref{MostowRigidity} holds whenever $\Gamma_i$ is irreducible in~$G_i$, for $i = 1,2$, and show that this additional hypothesis can be eliminated.

\item \label{MostowPfSO13EquiEx}
Construct a quasi-isometry $\varphi \colon \hyperbolic^3 \to \hyperbolic^3$ that satisfies \pref{MostowPfSO13Equi}.
\hint{Let $\fund$ be a precompact strict fundamental domain for the action of~$\Gamma_1$ on~$\hyperbolic^3$, and choose some $x_0 \in \hyperbolic^3$. For $x \in \hyperbolic^3$, let $\varphi(x) = \rho(\gamma)x_0$, where $x \in \gamma \cdot \fund$.}

\end{exercises}

\section{Moduli space of lattices in \texorpdfstring{$\SL(2,\real)$}{SL(2,R)}} \label{ModuliSpaceSL2}

Suppose $\Gamma$ is a torsion-free, cocompact lattice in $\SL(2,\real)$. In contrast to the Mostow Rigidity Theorem \pref{MostowRigidity}, we will see that an isomorphism $\Gamma \iso \Gamma'$ need not extend to an automorphism of $\SL(2,\real)$. In fact, there are uncountably many different embeddings of~$\Gamma$ in $\SL(2,\real)$ that are not conjugate to each other.

To see this, we take a geometric approach.
Since $\SL(2,\real)$ acts transitively on the hyperbolic plane~$\hyperbolic^2$ (by isometries), the quotient $\Gamma \backslash \hyperbolic^2$ is a compact surface~$M$. We will show there are uncountably many different possibilities for~$M$ (up to isometry). 
The proof is based on the fact that the hyperbolic plane has uncountably many different right-angled hexagons.

\begin{lem} \label{ManyHexagons}
The hyperbolic plane\/~$\hyperbolic^2$ has uncountably many different right-angled hexagons \textup(no two congruent to each other\/\textup).
\end{lem}

\begin{proof}
Fix a basepoint $p \in \hyperbolic^2$, and a starting direction~$\vector v$. For $\vector s \in (\real^+)^6$, construct the piecewise-linear path $L = L(\vector s)$ in~$\hyperbolic^2$ determined by:

\begin{minipage}{2.2in} %{% no page break in this paragraph !!!
	\raggedright
	\vbox to -3pt{\vskip 0.075in \rightline{\hbox to 0pt{\hskip8pt 
		\includegraphics{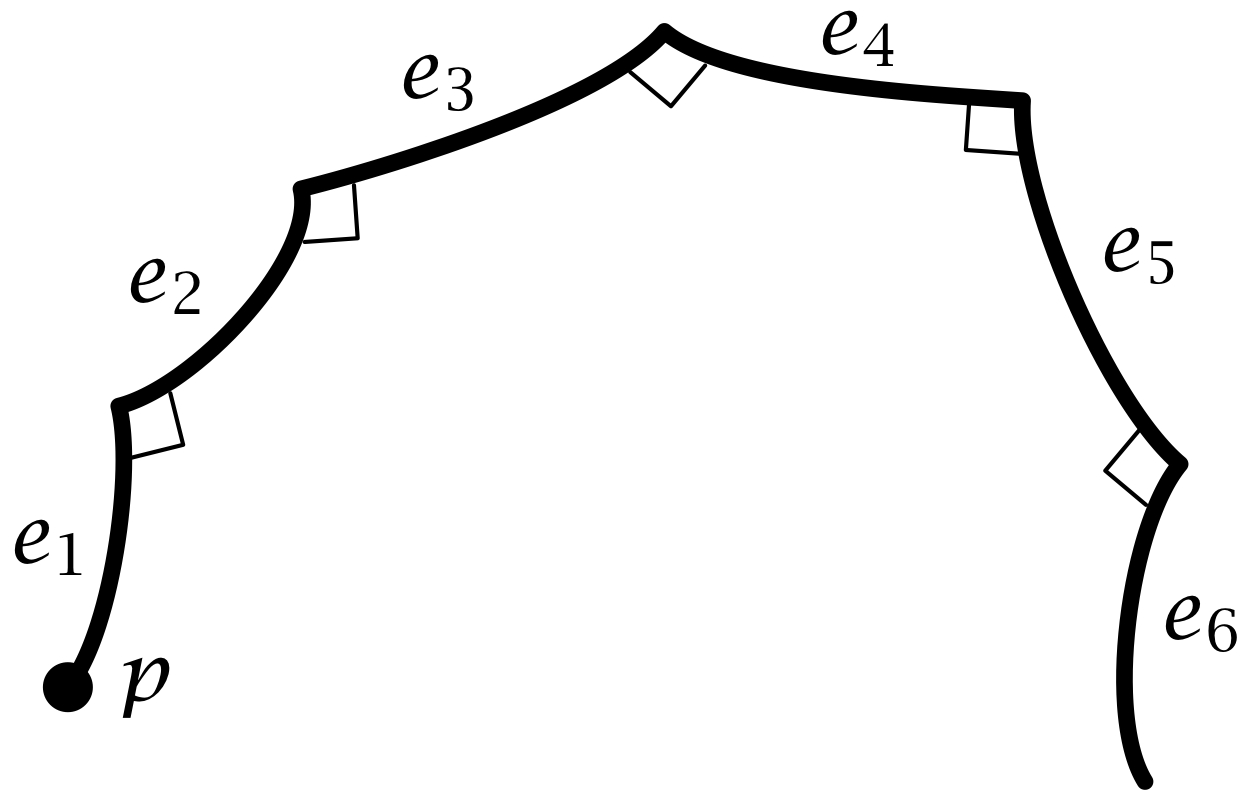}\hss}}\vss}
	\begin{itemize}
	\hangindent=-2.2in
	\hangafter=0
	\item The path starts at the point~$p$.
	\item \hangindent=-2in \hangafter=0
	 The path consists of $6$ geodesic segments (or ``edges'') $e_1,\ldots,e_6$, of lengths $s_1,\ldots,s_6$, respectively.\par
	\item The first edge~$e_1$ starts at~$p$ and heads in the direction~$\vector v$.
	\item For $i \ge 2$, edge~$e_i$ %is of length~$s_i$, and 
	makes a (clockwise) right angle with~$e_{i-1}$.
	\end{itemize}
\end{minipage}
\par\smallskip\noindent 
The variables $s_1,\ldots,s_6$ provide $6$~degrees of freedom in the construction of~$L$. Requiring that $L$ be a closed path (i.e., that the terminal endpoint of~$L$ is equal to~$p$) takes away two degrees of freedom (because $\hyperbolic^2$ is $2$-dimensional). Then, requiring the angle between $e_6$ and~$e_1$ to be right angle takes away one more degree of freedom. Hence (from the 
\thmindex{Implicit Function}{Implicit Function Theorem}), % should this be in the Background appendix? @@@
we see that there are $3$~degrees of freedom in the construction of a right-angled hexagon in~$\hyperbolic^2$. 
\end{proof}

\begin{rem} \label{AltEdgesInRtHex}
In fact, calculations using the trigonometry of triangles in~$\hyperbolic^2$ yields the much more precise fact that, for any $s_2,s_4,s_6 \in \real^+$, there exists a unique right-angled hexagon, with edges $e_1,\ldots,e_6$, such that the length of edge $e_{2i}$ is exactly~$s_{2i}$, for $i = 1,2,3$. (That is, the lengths of the three edges $e_2$, $e_4$, and~$e_6$ can be chosen completely arbitrarily, and they uniquely determine the lengths of the other three edges in the right-angled hexagon.)
\end{rem}

\begin{defn}
A \defit[hyperbolic!surface]{hyperbolic surface} is a compact Riemannian manifold (without boundary) whose universal cover is the hyperbolic plane~$\hyperbolic^2$.
\end{defn}

\begin{cor} \label{ManyHypSurfs}
There are uncountably many non-isometric hyperbolic surfaces of any given genus $g \ge 2$.
\end{cor}

\begin{proof}
Choose a right-angled hexagon~$P$ in~$\hyperbolic^2$. Call its edges $e_1,\ldots,e_6$, and let $s_i$ be the length of~$e_i$. Make a copy~$P'$ of~$P$, and form a surface~$\mathcal{P}$ by gluing $e_{2i}$ to the corresponding edge~$e_{2i}'$ of~$P'$, for $i = 1,2,3$. 
(Topologically, this surface~ $\mathcal{P}$ is a disk with two holes, and is usually called a ``pair of pants\zz.'')
$$ \includegraphics{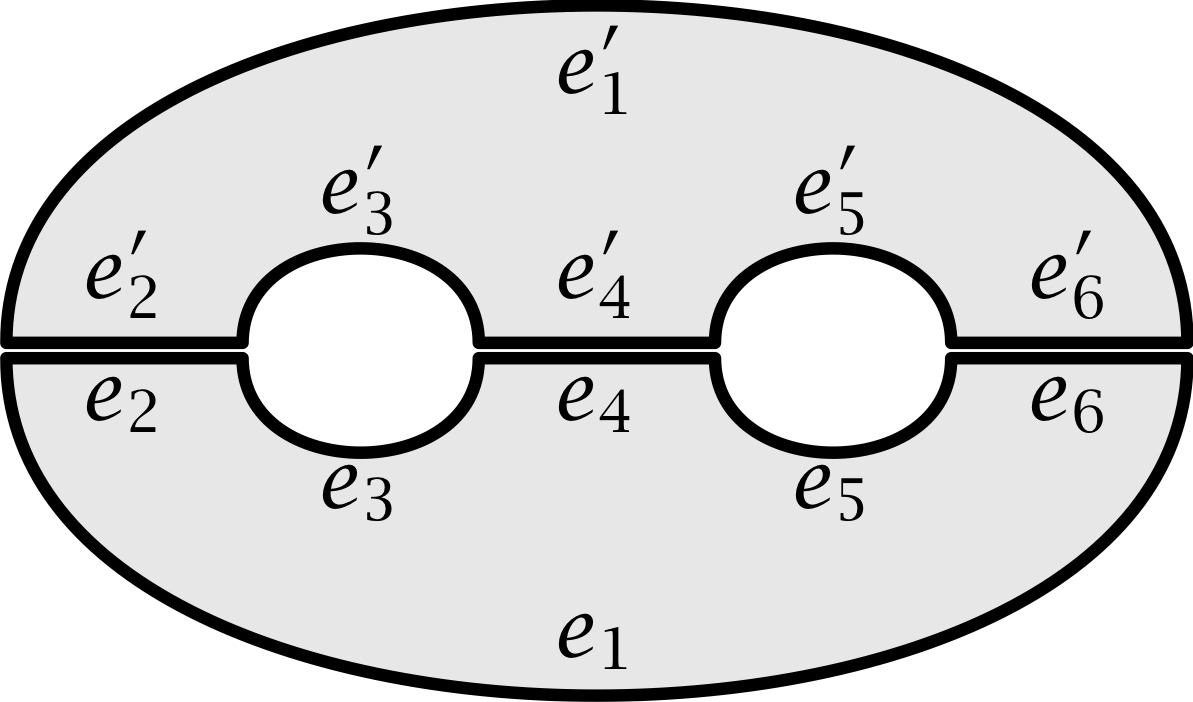}
%texpreamble
%("  \usepackage{amsmath}
% \usepackage[LY1]{fontenc}
% \usepackage[expert,LY1,mylucidascale]{mylucidabr}
% ");
%defaultpen(  fontcommand("\normalfont") + fontsize(10) ); 
%
%from graph access *;
%unitsize(1cm);
%
%real linethick = 1.5;
%real dotthick = 12;
%dotfactor=12;
%
%
%// draw( rotate(60)*ellipse( (0,0), 1,2 ) ) ;
%
%real eps = -0.0325;
%
%filldraw(
%	(
%	(0,0)
%	{N}..tension(1.75)..{S}
%	(5,0)--(4,0)
%	{N}..tension(1.25)..{S}
%	(3,0)--(2,0)
%	{N}..tension(1.25)..{S}
%	(1,0)--cycle
%	),
%	gray(0.925),
%	linewidth(1.5)
%	);
%
%filldraw(
%	reflect((0,eps), (1,eps))*(
%	(0,0)
%	{N}..tension(1.75)..{S}
%	(5,0)--(4,0)
%	{N}..tension(1.25)..{S}
%	(3,0)--(2,0)
%	{N}..tension(1.25)..{S}
%	(1,0)--cycle
%	),
%	gray(0.925),
%	linewidth(1.5)
%	);
%
%label( "$e_6'$", (4.5,0), N );
%label( "$e_5'$", (3.5,0.65) );
%label( "$e_4'$", (2.5,0), N );
%label( "$e_3'$", (1.5,0.65) );
%label( "$e_2'$", (0.5,0), N );
%label( "$e_1'$", (2.5,1.45), S);
%
%label( "$e_6$", (4.5,-0 + eps), S );
%label( "$e_5$", (3.5,-0.6 + eps) );
%label( "$e_4$", (2.5,-0 + eps), S );
%label( "$e_3$", (1.5,-0.6 + eps) );
%label( "$e_2$", (0.5,-0 + eps), S );
%label( "$e_1$", (2.5,-1.45 + eps) , N );
\qquad
 \includegraphics{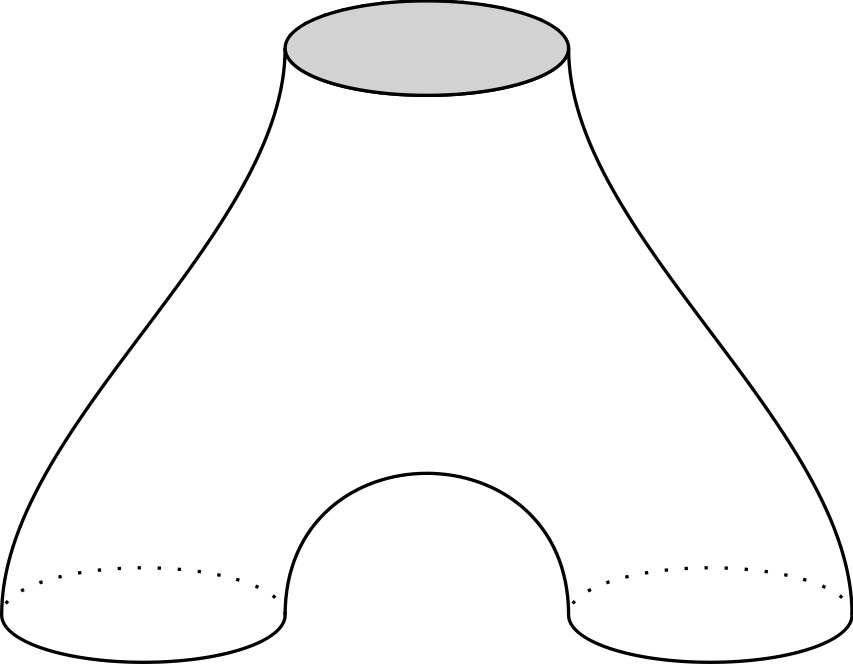} $$
 % http://commons.wikimedia.org/wiki/File:Pair_of_pants_cobordism_%28pantslike%29.svg
 % @@@ I do not have asymptote code for this figure, since it was converted from svg to PDF
Since $P$ is right-angled, the three boundary curves of~$\mathcal{P}$ are geodesics. Their lengths are $2s_1$, $2s_3$, and~$2s_5$. 

\vbox{ % can't allow a page break in this paragraph !!!
\vbox to 0pt{\vskip 0.3in \rightline{$\begin{matrix}
\includegraphics{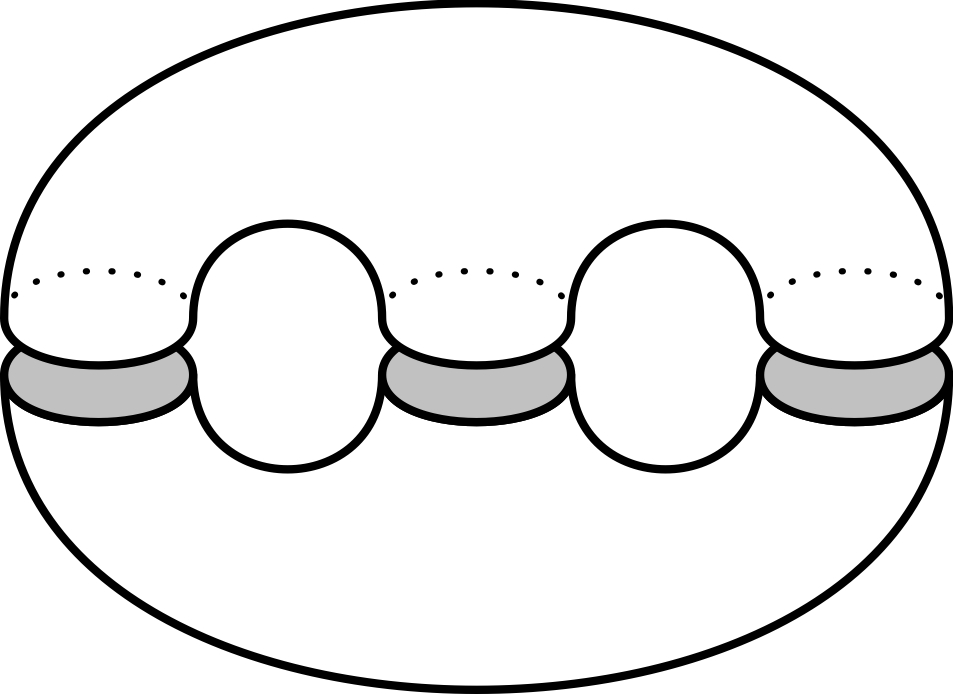}
\\
\text{A surface of genus 2 can be} 
\\ \text{made from two pairs of pants.}
\end{matrix}
$}\vss}
%from graph access *;
%
%unitsize(0.4cm);
%
%currentpen=linewidth(1);
%
%real a = 0, b = 2, c = 4, d = 6, e = 8, f = 10;
%real z = -0.6;
%real s = 1;
%real t = 2;
%real gl = 0.8;
%real lw = 0.75;
%string lt = "3 3";
%
%
%draw(  (a,z){S}.. tension t .. {N}(b,z){S}..tension s ..{N}(c,z){S}..tension t ..{N}(d,z){S}..tension s ..{N}(e,z){S}..tension t ..{N}(f,z){S}.. tension 1.5 ..{N}cycle );
%
%
%
%filldraw( (a,z){S}.. tension t .. {N}(b,z){N}.. tension t .. {S}cycle , gray(gl) );
%filldraw( (c,z){S}.. tension t .. {N}(d,z){N}.. tension t .. {S}cycle , gray(gl) );
%filldraw( (e,z){S}.. tension t .. {N}(f,z){N}.. tension t .. {S}cycle , gray(gl) );
%
%
%filldraw(  (a,0){S}..tension t ..{N}(b,0){N}..tension s ..{S}(c,0){S}..tension t ..{N}(d,0){N}..tension s ..{S}(e,0){S}..tension t ..{N}(f,0){N}.. tension 1.5 ..{S}cycle, white );
%
%draw( (a,0){N}.. tension t ..{S}(b,0), linewidth(lw)+dotted );
%draw( (c,0){N}.. tension t ..{S}(d,0), linewidth(lw)+dotted );
%draw( (e,0){N}.. tension t ..{S}(f,0), linewidth(lw)+dotted );

\hangindent=-2.2in
\hangafter=0
Construct a closed surface~$M$ of genus~$2$ from two copies of~$\mathcal{P}$, by gluing corresponding boundary components to each other. 
(For a discussion of higher genus, see \cref{HigherGenusPants}.) Since the only curves that have been glued together are geodesics, it is easy to see that each point in~$M$ has a neighborhood that is isometric to an open subset of~$\hyperbolic^2$. Therefore, the universal cover of~$M$ is~$\hyperbolic^2$ (since $M$ is complete). So $M$ is a hyperbolic surface.
\par}

Furthermore, from the construction, we see that $M$ has a closed geodesic of length~$2s_1$. (In fact, there is a geodesic of length~$2s_i$, for $1 \le i \le 6$.) Since a single closed surface has closed geodesics of only countably many different lengths, but \cref{ManyHexagons} implies that there are uncountably many possible values of~$s_1$, this implies there must be uncountably many different isometry classes of surfaces.
\end{proof}

\begin{rem} \label{HigherGenusPants}
A hyperbolic surface~$M$ of any genus $g \ge 2$ can be constructed by gluing together $2g - 2$ pairs of pants:
$$ \includegraphics{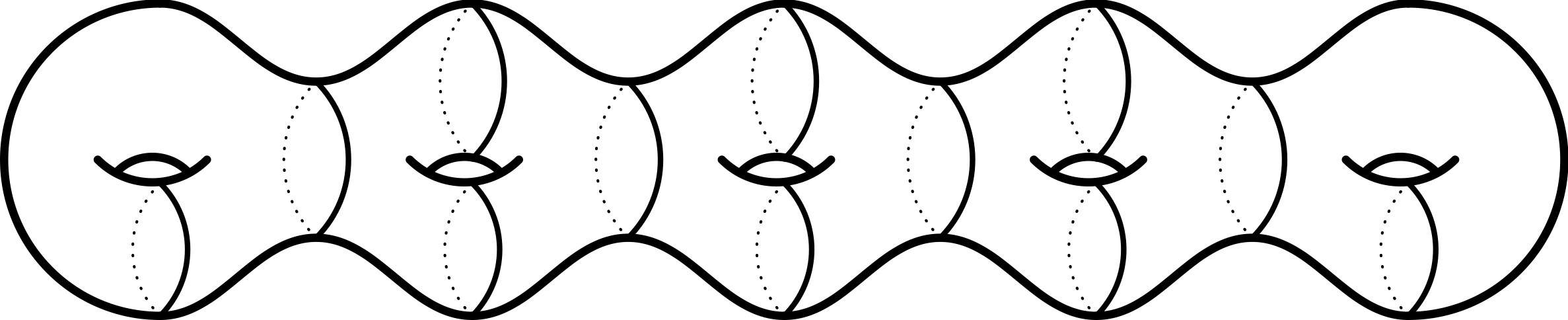} $$
The lengths of the three boundary curves of each pair of pants can be varied independently \ccf{AltEdgesInRtHex}, except that curves that will be glued together need to have the same length. The surface can also be modified by rotating any boundary curve through an arbitrary angle~$\theta$ before it is glued to its mate. This yields $6g-6$ degrees of freedom in the construction of~$M$. It can be shown that this is precisely the dimension of the space of hyperbolic surfaces of genus~$g$. In other words, $6g-6$ is the dimension of the \defit{moduli space} of hyperbolic surfaces of genus~$g$. 
\end{rem}
%texpreamble
%(" \usepackage[LY1]{fontenc}
% \usepackage[expert, LY1, mylucidascale]{mylucidabr} % I adjusted the scaling
% \usepackage{amsmath}
% \everymath{\displaystyle}
% ");
%defaultpen(  fontcommand("\normalfont") + fontsize(10) ); 
%
%from graph access *;
%
%unitsize(1cm);
%
%
%real a = 1, b = 0.5;
%
%currentpen=linewidth(1.5);
%draw( (1,-a){W}..(0,0)..{E}(1,a) );
%draw( (9,-a){E}..(10,0)..{W}(9,a) );
%for (int i = 1; i <= 7; i = i + 2){
%	draw( (i,a){E}..(i+1,b)..{E}(i+2,a) );
%	draw( (i,-a){E}..(i+1,-b)..{E}(i+2,-a) );
%	}
%
%for (int i = 0; i <= 8; i = i + 2 ){
%	currentpen=linewidth(1.5);
%	draw( (i + 0.6,0){SE}..{NE}(i + 1.3,0) );
%	draw( (i + 0.725,-0.075){NE}..{SE}(i + 1.175,-0.075) );
%
%	currentpen=linewidth(1);
%	draw( (i + 1, -0.15){(1,-1)}..{(-1,-1)}(i + 1,-a) );
%	draw( (i + 1, -0.15){(-1,-1)}..{(1,-1)}(i + 1,-a), dotted );
%	}
%
%currentpen=linewidth(1);
%
%for (int i = 3; i <= 7; i = i + 2 ){
%	draw( (i, 0.005){(1,1)}..{(-1,1)}(i,a) );
%	draw( (i, 0.005){(-1,1)}..{(1,1)}(i,a), dotted );
%	}
%
%for (int i = 2; i <= 8; i = i + 2 ){
%	draw( (i, -b){(1,1)}..{(-1,1)}(i,b) );
%	draw( (i, -b){(-1,1)}..{(1,1)}(i,b), dotted );
%	}

\begin{cor} \label{LotsOfLattsSL2}
If\/ $\Gamma$ is any lattice in\/ $\SL(2,\real)$, then there are uncountably many nonconjugate embeddings of\/~$\Gamma$ as a lattice in\/ $\SL(2,\real)$.
\end{cor}

\begin{proof}
Let us assume $\SL(2,\real)/\Gamma$ is compact. (Otherwise, $\Gamma$ is a free group, so it is easy to find embeddings.) Let us also assume, for simplicity, that $\Gamma$~is torsion free, so it is the fundamental group of the hyperbolic surface $\Gamma \backslash \hyperbolic^2$, which has some genus~$g$. Then $\Gamma$ is isomorphic to the fundamental group~$\Gamma'$ of any hyperbolic surface $\Gamma' \backslash \hyperbolic^2$ of genus~$g$. However, if $\Gamma \backslash \hyperbolic^2$ is not isometric to $\Gamma' \backslash \hyperbolic^2$, then $\Gamma$ cannot be conjugate to~$\Gamma'$ in the isometry group of~$\hyperbolic^2$. Therefore, \cref{ManyHypSurfs} implies that there must be uncountably many different conjugacy classes of subgroups~$\Gamma'$ that are isomorphic to~$\Gamma$.
\end{proof}

\section{Quasi-isometric rigidity} \label{QIRigSect}

The Mostow Rigidity Theorem's weak form \pref{MostowIso} tells us (under mild hypotheses) that lattices in two different semisimple Lie groups cannot be isomorphic. In fact, they cannot even be quasi-isometric \csee{QuasiIsomDefn}:

\begin{thm} \label{QI->GIso}
Assume
\noprelistbreak
	\begin{itemize}
	\item $G_1$ and~$G_2$ are connected, with trivial center and no compact factors,
	and
	\item $\Gamma_i$ is an irreducible lattice in~$G_i$, for $i = 1,2$.
	\end{itemize}
If\/ $\Gamma_1 \QI \Gamma_2$, then $G_1 \iso G_2$.
\end{thm}

Although nothing more than \cref{QI->GIso} can be said about quasi-isometric lattices that are cocompact \csee{CpctLattQI}, there is a much stronger conclusion for noncocompact lattices. Namely, not only are the Lie groups $G_1$ and~$G_2$ isomorphic, but the isomorphism can be chosen to make the lattices commensurable (unless $G_1 = G_2 = \PSL(2,\real)$): 

\begin{thm} \label{QuasiMostowRig}
Assume
\noprelistbreak
	\begin{itemize}
	\item $G_1$ and~$G_2$ are connected, with trivial center and no compact factors,
	and
	\item $\Gamma_i$ is an irreducible lattice in~$G_i$, for $i = 1,2$.
	\end{itemize}
Then\/ $\Gamma_1 \QI \Gamma_2$ if and only if
$G_1 \iso G_2$ and either
 \begin{enumerate}
 \item both $G_1/\Gamma_1$ and $G_2/\Gamma_2$ are compact, or
 \item there is an isomorphism $\phi \colon G_1 \to G_2$,
such that $\phi(\Gamma_1)$ is commensurable
to\/~$\Gamma_2$, or
 \item $G_1$ and $G_2$ are isomorphic to\/ $\PSL(2,\real)$,
and neither $G_1/\Gamma_1$ nor $ G_2/\Gamma_2$ is compact.
 \end{enumerate}
 \end{thm}
 
%We will say a little about its proof later in the section. 
One of the key ingredients in the proof of this \lcnamecref{QuasiMostowRig} is the fact that any group quasi-isometric to a lattice is isomorphic to a lattice, modulo finite groups:

\begin{thm} \label{QItoLatt}
 If $\Lambda$ is a finitely generated group that is quasi-isometric to an irreducible lattice\/~$\Gamma$ in~$G$, then there are
 	\begin{itemize}
	\item a finite-index subgroup~$\Lambda'$ of~$\Lambda$,
	and
	\item a finite, normal subgroup~$N$ of~$\Lambda'$, 
	\end{itemize}
such that $\Lambda'/N$ is isomorphic to a lattice in~$G$.
 \end{thm}
 
 \begin{proof}[\mathversion{bold}Sketch of proof of a special case]
 Assume $G = \SO(1,3)$ and $\Gamma$~is cocompact. 
The group~$\Lambda$ acts on itself by translation. Since $\Lambda \QI \Gamma \QI \hyperbolic^3$, this provides an action of~$\Lambda$ by quasi-isometries on~$\hyperbolic^3$. Let $\overline{\Lambda}$ be the corresponding group of quasiconformal maps on $\bdry\hyperbolic^3$. Note that the quasiconformality constant~$\kappa$ is uniformly bounded on~$\overline{\Lambda}$.

Let $(\bdry\hyperbolic^3)^3_\circ$ be the space of ordered triples of distinct points in $\bdry\hyperbolic^3$, and define $p \colon (\bdry\hyperbolic^3)^3_\circ \to \hyperbolic^3$ by letting $p(a,b,c)$ be the point on the geodesic $\overline{ab}$ that is closest to~$c$. It is not difficult to see that $p$ is compact-to-one. Since the action of~$\Lambda$ is cocompact on $\hyperbolic^3$, this implies that the action of~$\overline{\Lambda}$ is cocompact on $(\bdry\hyperbolic^3)^3_\circ$. 

The above information allows us to apply a theorem of Tukia to conclude that $\overline{\Lambda}$ is quasiconformally conjugate to a subgroup of $\overline{\SO(1,n)}$. Hence, after conjugating the action of~$\Lambda$ by a quasi-isometry, 
%modding out the (finite) kernel of the map $\lambda \mapsto \overline{\lambda}$, and choosing a basepoint $x_0 \in \hyperbolic^3$, 
we may assume $\Lambda \subseteq \SO(1,n)$. Furthermore, if we fix a basepoint~$x_0$, then the map $\lambda \mapsto \lambda x_0$ is a quasi-isometry from~$\Lambda$ to~$\hyperbolic^3$. This implies that $\Lambda$ is a cocompact lattice in $\SO(1,n)$.
 \end{proof}

Many additional ideas are needed to prove \cref{QuasiMostowRig}, but this suffices for the weaker version:

\begin{proof}[Proof of \cref{QI->GIso}]
\Cref{QItoLatt} tells us that $\Gamma_2$ is isomorphic to a lattice in~$G_1$ (if we ignore some finite groups). So $\Gamma_2$ is a lattice in both~$G_1$ and~$G_2$. Therefore, the Mostow Rigidity Theorem \pref{MostowIso} implies $G_1 \iso G_2$, as desired. 
\end{proof}

From the Mostow Rigidity Theorem, we know that every
automorphism of~$\Gamma$ extends to an automorphism of~$G$
(if $G$ is not isogenous to $\SL(2,\real)$ and we ignore compact factors and the center). The following
analogue of this result for quasi-isometries is another key ingredient in the proof of
Theorem~\ref{QuasiMostowRig}.

\begin{thm} \label{QIClose}
 Assume
 \noprelistbreak
 \begin{itemize}
 \item $\Gamma$ is irreducible, and \textbf{not} cocompact,
 \item $G$ has trivial center and no compact factors,
 \item $G$ is not isogenous to $\SL(2,\real)$,
 and
 \item $f \colon \Gamma \to \Gamma$.
 \end{itemize}
 Then $f$ is a quasi-isometry if and only if it is at bounded distance from some automorphism of~$G$.
%  is an
%automorphism~$\phi$ of~$G$, such that $f \approx_{QI}
%\phi|_\Gamma$.
%
%Therefore, the group $\QIsom(\Gamma)$ is naturally isomorphic to
%$\Comm_{\Aut(G)}(\Gamma)$.
 \end{thm}

\begin{exercises}

\item \label{QI->Mostow}
Assume the hypotheses of  \cref{MostowRigidity}, and also assume $G_1/\Gamma_1$ is not compact. Show that the conclusion of \cref{MostowRigidity} can be obtained by combining \cref{QuasiMostowRig} with  \cref{QIClose}.

%\item \label{QItoLattEx}
%Prove \cref{QItoLatt} under the additional assumption that the lattice~$\Gamma$ is cocompact.
%\hint{$\Lambda$ acts on itself by translations. Since $\Lambda \QI \Gamma$, this provides a homomorphism from $\Lambda$ to $\QIsom(\Gamma)$. After passing to finite index, this is a homomorphism into~$G$.}

\end{exercises}

\begin{notes}

The Mostow Rigidity Theorem \pref{MostowRigidity} is a combination of (overlapping) special cases proved by three authors:
	\begin{itemize}
	\item Mostow: $\Gamma$ is cocompact
		(case where $G = \SO(1,n)$ \cite{Mostow-QIhyperbolic},
		 general case \cite{MostowRigidity}),
	\item Prasad: $\Qrank \Gamma = 1$ (which includes the case where $\Gamma$ is \emph{not} cocompact and $\Rrank G = 1$) \cite{PrasadMostowRig},
	\item Margulis: $\Rrank G \ge 2$ (and $\Gamma$ is irreducible) \csee{StateMostowSubsect}.
	\end{itemize}
In recognition of this, many authors call \cref{MostowRigidity} the 
	\thmindex{Mostow-Prasad Rigidity}Mostow-Prasad Rigidity Theorem, 
or the Mostow-Prasad-Margulis Rigidity Theorem. 
%Furthermore, some authors, including Mostow himself, insert the adjective "Strong" before "Rigidity\zz.''

The results described in \fullcref{MostowRems}{FurmanRigidity} are due to A.\,Furman \cite{Furman-MostowMargulisRigidity}.

See the exposition in \cite[\S5.9, pp.~106--112]{Thurston-GeomTop3Mflds} for more details of the proof in \cref{MostowPfSect}.  A different proof of this special case was found by Gromov, and is described in \cite[\S6.3, pp.~129--130]{Thurston-GeomTop3Mflds}. Yet another nice proof (which applies to cocompact lattices in any groups of real rank one) appears in \cite[\S5.2]{BessonCourtoisGallot-MinEntropy}.

Regarding \fullcref{MostowPfRems}{QuasiC}, see \cite[\S21, esp.\ (21,18)]{MostowRigidity} for a discussion of the notion of maps that are quasiconformal over~$\complex$ (or $\quaternion$, or~$\octonion$).

Regarding \fullcref{MostowPfRems}{Qrank1}, see \cite{PrasadMostowRig} for Prasad's proof of Mostow rigidity for lattices of $\rational$-rank one.

Regarding \cref{ModuliSpaceSL2}, see \cite[Thm.~5.3.5]{Thurston-GeomTop3Mflds} for a more complete discussion of ``pairs of pants'' and the dimension of the space of hyperbolic metrics on a surface of genus~$g$. (The calculations to justify \cref{AltEdgesInRtHex} can be found in \cite[\S2.6]{Thurston-GeomTop3Mflds}.)

The results on quasi-isometric rigidity in \cref{QIRigSect} include work of A.\,Eskin, B.\,Farb, B.\,Kleiner, B.\,Leeb, P.\,Pansu, R.\,Schwarz, and others. See the survey \cite{FarbQISurvey} for references, discussion of the proofs, and other information.

\end{notes}

  \immediate\addtocontents{toc}{\protect\toceject} % @@@
 %!TEX root = IntroArithGrps.tex

\mychapter{Margulis Superrigidity\texorpdfstring{\\}{ }Theorem}
\label{MargulisSuperChap}

\prereqs{none, except that the proof of the main theorem requires real rank (\cref{RrankChap}), amenability (Furstenberg's Lemma \pref{G/amen->Meas(X)}), and the Moore Ergodicity Theorem (\cref{MooreErgPfSect}).}

Roughly speaking, the Margulis Superrigidity Theorem tells us that homomorphisms defined on~$\Gamma$ can be extended to be defined on all of~$G$ (unless $G$ is either $\SO(1,n)$ or $\SU(1,n)$). 
In cases where it applies, this fundamental theorem is much stronger than the Mostow Rigidity Theorem \pref{MostowRigidity}.
It also implies the Margulis Arithmeticity Theorem (\ref{MargulisArith} or~\ref{MargArithFromSuper}).

\section{Statement of the theorem} \label{MargSuperStatementSect}

It is not difficult to see that every group homomorphism from~$\integer^k$ to~$\real^n$ can be extended to a continuous homomorphism from~$\real^k$ to~$\real^n$ \csee{ZkSuperrigEx}. Noting that $\integer^k$ is a lattice in~$\real^k$, it is natural to hope that, analogously, homomorphisms defined on~$\Gamma$ can be extended to be defined on all of~$G$. The Margulis Superrigidity Theorem shows this is true if $G$ has no simple factors isomorphic to $\SO(1,m)$ or $\SU(1,m)$, except that the conclusion may only be true modulo finite groups and up to a bounded error. Here is an illustrative special case that is easy to state, because the bounded error does not arise. 

\begin{thm}[(Margulis)] \label{MargSuperSL3RNonCpct}
Assume
\noprelistbreak
	\begin{itemize}
	\item $G = \SL(k,\real)$, with $k \ge 3$,
	\item $G/\Gamma$ is \textbf{not} compact,
	and
	\item $\varphi \colon \Gamma \to \GL(n,\real)$ is any homomorphism.
	\end{itemize}
Then there exist:
\noprelistbreak
	\begin{itemize}
	\item a continuous homomorphism $\widehat\varphi \colon G \to \GL(n,\real)$, 
	and
	\item a finite-index subgroup~$\Gamma'$ of~$\Gamma$,
	\end{itemize}
such that  $\widehat\varphi(\gamma) = \varphi(\gamma)$ for all $\gamma \in \Gamma'$.
\end{thm}

\begin{proof}
See \cref{MargSuperNonCpctSCEx}.
\end{proof}

Here is a much more general version of the theorem that has a slightly weaker conclusion. To simplify the statement, we preface it with a definition. 

\begin{defn} \label{AlgSCDefn}
$G$ is \defit[simply connected, algebraically]{algebraically simply connected} if (for every~$\ell$\,) every Lie algebra homomorphism $\Lie G \to \LieSL(\ell,\real)$ is the derivative of a well-defined Lie group homomorphism $G \to \SL(\ell,\real)$.
\end{defn}

\begin{rem}
Every simply connected Lie group is algebraically simply connected, but the converse is not true \csee{SLnASC}. In general, if $G$ is connected, then some finite cover of~$G$ is algebraically simply connected. Therefore, assuming that $G$ is algebraically simply connected is just a minor technical assumption that avoids the need to pass to a finite cover. 
\end{rem}

\begin{thm}[(Margulis Superrigidity Theorem)] \label{MargSuperC}
Assume
\noprelistbreak
	\begin{enumerate} \renewcommand{\theenumi}{\roman{enumi}}
	\item $G$ is connected, and algebraically simply connected,
	\item \label{MargSuperC-notSOSU}
	$G$ is not isogenous to any group that is of the form\/ $\SO(1,m) \times K$ or\/ $\SU(1,m) \times K$, where $K$~is compact,
	\item $\Gamma$ is irreducible,
	and
	\item $\varphi \colon \Gamma \to \GL(n,\real)$ is a homomorphism.
	\end{enumerate}
Then there exist:
\noprelistbreak
	\begin{enumerate}
	\item a continuous homomorphism $\widehat\varphi \colon G \to \GL(n,\real)$,
	\item a compact subgroup~$C$ of\/ $\GL(n,\real)$ that centralizes $\widehat\varphi(G)$,
	and
	\item a finite-index subgroup\/~$\Gamma'$ of\/~$\Gamma$,
	\end{enumerate}
such that $\varphi(\gamma) \in \widehat\varphi(\gamma) \, C$, for all $\gamma \in \Gamma'$.
\end{thm}

\begin{proof}
See \cref{SuperPfSect}.
\end{proof}

\begin{rems} \label{SuperRem} \ 
\noprelistbreak
	\begin{enumerate}
	\item Since $\varphi(\gamma) \in \widehat\varphi(\gamma) \, C$, we have $\widehat\varphi(\gamma)^{-1} \, \varphi(\gamma) \in C$ for all~$\gamma$. Therefore, although $\widehat\varphi(\gamma)$ might not be exactly equal to $\varphi(\gamma)$, the error is an element of~$C$, which is a bounded set (because $C$ is compact). Hence, the size of the error is uniformly bounded on all of~$\Gamma'$.
 
	\item \label{SuperRem-NotSOSU}
	Assumption~\pref{MargSuperC-notSOSU} cannot be removed. For example, if $G = \PSL(2,\real)$, then the lattice $\Gamma$ can be a free group \csee{FreeLattINSL2R}. In this case, there exist many, many homomorphisms from~$\Gamma$ into any group~$G'$, and many of them will not extend to~$G$ \csee{FreeNotSuperrig}.

	\end{enumerate}
 \end{rems}

If we make an appropriate assumption on the range of~$\varphi$, then there is no need for the compact error term~$C$ or the finite-index subgroup~$\Gamma'$:

\begin{cor} \label{MargSuperG'}
Assume
\noprelistbreak
	\begin{enumerate} \renewcommand{\theenumi}{\roman{enumi}}
	\item \label{MargSuperG'-notSOSU}
	$G$ is not isogenous to any group that is of the form\/ $\SO(1,m) \times K$ or\/ $\SU(1,m) \times K$, where $K$~is compact,
	\item  \label{MargSuperG'-irred}
 $\Gamma$ is irreducible, 
 and
	\item  \label{MargSuperG'-G}
 $G$ and $G'$ are connected, with trivial center, and no compact factors.
	\end{enumerate}
If $\varphi \colon \Gamma \to G'$ is any homomorphism, such that $\varphi(\Gamma)$ is Zariski dense in~$G'$, then $\varphi$ extends to a continuous homomorphism $\widehat\varphi \colon G \to G'$.
\end{cor}

\begin{proof}
See \cref{MargSuperG'PfEx}.
\end{proof}
 
Because of our standing assumption \pref{standassump} that $G'$ is semisimple, \cref{MargSuperG'} implicitly assumes that the Zariski closure $\Zar{\varphi(\Gamma)} = G'$ is semisimple. In fact, that is automatically the case:
  
 \begin{cor} \label{MargImgSS}
 Assume
 \noprelistbreak
 	\begin{enumerate} \renewcommand{\theenumi}{\roman{enumi}}
	\item  \label{MargImgSS-notSOSU}
	$G$ is not isogenous to any group that is of the form\/ $\SO(1,m) \times K$ or\/ $\SU(1,m) \times K$, where $K$~is compact,
	\item  \label{MargImgSS-irred}
	 $\Gamma$ is irreducible,
	 and
 	\item \label{MargImgSS-phi}
	$\varphi \colon \Gamma \to \GL(n,\real)$ is a homomorphism.
	\end{enumerate}
Then\/ $\Zar{\varphi(\Gamma)}$ is semisimple.
\end{cor}

\begin{proof}
See \cref{MargImgSSPfEx}.
\end{proof}

\begin{exercises}

\item \label{ZkSuperrigEx}
Suppose $\varphi$ is a homomorphism from~$\integer^k$ to~$\real^n$. Show that $\varphi$ extends to a continuous homomorphism from~$\real^k$ to~$\real^n$.
\hint{Let $\widehat\varphi \colon \real^k \to \real^n$ be a linear transformation, such that $\widehat\varphi(\varepsilon_i) = \varphi(\varepsilon_i)$, where $\{\varepsilon_1,\ldots,\varepsilon_k\}$ is the standard basis of~$\real^k$.}

\item \label{SLnASC}
Show that $\SL(n,\real)$ is algebraically simply connected. (On the other hand, $\SL(n,\real)$ is not simply connected, because its fundamental group is nontrivial.)
\hint{By tensoring with~$\complex$, any homomorphism $\LieSL(n,\real) \to \LieSL(\ell,\real)$ extends to a homomorphism $\LieSL(n,\complex) \to \LieSL(\ell,\complex)$, and $\SL(n,\complex)$ is simply connected.}

\item \label{GNoAbelianization}
Assume 
	\begin{itemize}
	\item $G$ is not isogenous to $\SO(1,n)$ or $\SU(1,n)$, for any~$n$,
	\item $\Gamma$ is irreducible,
	and
	\item $G$ has no compact factors.
	\end{itemize}
Use the Margulis Superrigidity Theorem to show that the abelianization $\Gamma / [\Gamma,\Gamma]$ of~$\Gamma$ is finite.
(When $G$ is simple, this was already proved from Kazhdan's property~$(T)$ in \fullcref{KazhdanlatticeCor}{noabel}. We will see yet another proof in \cref{GammaHasFiniteAbelianization}.)

\item Assume $G$, $\Gamma$, $\varphi$, $\widehat\varphi$, $C$, and~$\Gamma'$ are as in \cref{MargSuperC}. Show there is a homomorphism $\epsilon \colon \Gamma' \to C$, such that $\varphi(\gamma) = \widehat\varphi(\gamma) \cdot \epsilon(\gamma)$, for all $\gamma \in \Gamma'$.

\item \label{FreeNotSuperrig}
Suppose $G = \PSL(2,\real)$ and $\Gamma$ is a free group.
Construct a homomorphism $\varphi \colon \Gamma \to \GL(n,\real)$ (for some~$n$), such that, for every continuous homomorphism $\widehat\varphi \colon G \to \GL(n,\real)$, and every finite-index subgroup~$\Gamma'$ of~$\Gamma$, the set $\{\, \widehat\varphi(\gamma)^{-1} \, \varphi(\gamma) \mid \gamma \in \Gamma' \,\}$ is not precompact.
\hint{$\varphi$ may have an infinite kernel.}

\item \label{MargSuperG'PfEx}
Prove \cref{MargSuperG'} from \cref{MargSuperC}.

\item %Assume $G$ has no compact factors. 
Show that the extension~$\widehat\varphi$ in \cref{MargSuperG'} is unique.
\hint{Borel Density Theorem.}

\item \label{MargSuperG'NeedEx}
In each case, find
	\begin{itemize}
	\item a lattice~$\Gamma$ in~$G$
	and
	\item a homomorphism $\varphi \colon \Gamma \to G'$,
	\end{itemize}
such that 
	\begin{itemize}
	\item $\varphi(\Gamma)$ is Zariski dense in~$G'$,
	and
	\item $\varphi$ does not extend to a continuous homomorphism $\widehat\varphi \colon G \to G'$.
	\end{itemize}
Also explain why they are not counterexamples to \cref{MargSuperG'}.
	\begin{enumerate}
	\item $G = G' =  \PSL(2,\real) \times \PSL(2,\real)$.
	\item $G = \PSL(4,\real)$ and $G' = \SL(4,\real)$.
	\item $G = \SO(2,3)$ and $G' = \SO(2,3) \times \SO(5)$.
	\end{enumerate}

\item \label{MargImgSSPfEx}
Prove \cref{MargImgSS} from \cref{MargSuperC}.

\item \label{MargSuperFromG'PfEx}
Derive \cref{MargSuperC} from the combination of \cref{MargSuperG'} and \cref{MargImgSS}. (This is a converse to \cref{MargSuperG'PfEx,MargImgSSPfEx}.)

%\item Suppose the homomorphism $\varphi \colon \Gamma \to \SL(n,\real)$ extends to a continuous homomorphism $\widehat\varphi \colon G \to \SL(n,\real)$. Show that the almost-Zariski closure $\Zar{\varphi(\Gamma)}$ is connected and semisimple.

%\item \label{MargSuperNoncpctEx}
%Show that if 
%\noprelistbreak
%	\begin{itemize}
%	\item $\Rrank G \ge 2$,
%	\item $G/\Gamma$ is not compact,
%	\item $\varphi$ is any homomorphism from~$\Gamma$ to $\SL(n,\real)$,
%	and
%	\item $G'$ is the identity component of $\Zar{\varphi(\Gamma)}$,
%	\end{itemize}
%then there exist:
%\noprelistbreak
%	\begin{enumerate}
%	\item a finite-index subgroup~$\Gamma'$ of~$\Gamma$,
%	\item a finite subgroup~$Z$ of the center of~$G'$,
%	and
%	\item a homomorphism $\widehat \varphi \colon G \to G'/Z$,
%	\end{enumerate}
%such that $\varphi(\gamma) Z = \widehat\varphi(\gamma)$, for all $\gamma \in \Gamma'$. 

\end{exercises}

\section{Applications}

We briefly describe a few important consequences of the Margulis Superrigidity Theorem.

\subsection{Mostow Rigidity Theorem}
\label{StateMostowSubsect}
\thmindex{Mostow Rigidity}

The special case of the Margulis Superrigidity Theorem \pref{MargSuperG'} in which the homomorphism~$\varphi$ is assumed to be an isomorphism onto a lattice~$\Gamma'$ in~$G'$ is very important:

\begin{thm}[(\thmindex{Mostow Rigidity Theorem}Mostow Rigidity Theorem, cf.\  \pref{MostowRigidity})] \label{MostowRigidityIrred}
Assume
\noprelistbreak
	\begin{itemize}
	\item $G_1$ and~$G_2$ are connected, with trivial center and no compact factors,
	\item $G_1 \not\iso \PSL(2,\real)$,
	\item $\Gamma_i$ is an irreducible lattice in~$G_i$, for $i = 1,2$,
	and 
	\item $\varphi \colon \Gamma_1 \to \Gamma_2$ is a group isomorphism.
	\end{itemize}
Then $\varphi$ extends to a continuous isomorphism from~$G_1$ to~$G_2$.
\end{thm}

This theorem has already been discussed in \cref{MostowChap}.
In most cases, it follows easily from the Margulis Superrigidity Theorem \csee{ProveMostMostowEx}. However, since the superrigidity theorem does not apply when $G_1$ is either $\SO(1,m)$ or $\SU(1,m)$, a different argument is needed for those cases; see \cref{MostowPfSect} for a sketch of the proof.

\subsection{Triviality of flat vector bundles over $G/\Gamma$}

\begin{defn} \label{FlatVecBundleDefn}
For any homomorphism $\varphi \colon \Gamma \to \GL(n,\real)$, there is a diagonal action of~$\Gamma$ on $G \times \real^n$, defined by 
	$$ (x,v) \cdot \gamma = \bigl( x \gamma, \varphi(\gamma^{-1}) v \bigr) . $$
Let $\bundle_\varphi = (G \times \real^n) / \Gamma$ 
	\nindex{$\bundle_\varphi$ = flat vector bundle over $G/\Gamma$}
be the space of orbits of this action. Then there is a well-defined map 
	$$ \text{$\pi \colon \bundle_\varphi \to G/ \Gamma$, defined by $\pi \bigl( [x,v] \bigr) = x \Gamma$,} $$
and this makes $\bundle_\varphi$ into a vector bundle over $G/\Gamma$ (with fiber~$\real^n$) \csee{FlatVecBundleDefnEx}.
A vector bundle defined from a homomorphism in this way is said to be a \defit[flat!vector bundle]{flat vector bundle}.
\end{defn}

The Margulis Superrigidity Theorem implies (in some cases) that every flat vector bundle over $G/\Gamma$ is nearly trivial. Here is an example:

\begin{prop} \label{Super->VecBdlTrivial}
Let $G = \SL(n,\real)$ and\/ $\Gamma = \SL(n,\integer)$. If $\bundle_\varphi$ is any flat vector bundle over\/ $G/\Gamma$, then there is a finite-index subgroup\/~$\Gamma'$ of\/~$\Gamma$, such that the lift of $\bundle_\varphi$ to the finite cover\/ $G/\Gamma'$ is trivial.

In other words, if we let $\varphi'$ be the restriction of~$\varphi$ to\/~$\Gamma'$, then the vector bundle $\bundle_{\varphi'}$ is isomorphic to the trivial vector bundle\/ $(G/\Gamma') \times \real^n$.
\end{prop}

\begin{proof}
From \cref{MargSuperSL3RNonCpct}, we may choose $\Gamma'$ so that the restriction~$\varphi'$ extends to a homomorphism $\widehat\varphi \colon G \to \GL(n,\real)$. Define a continuous function $T \colon G \times \real^n \to G \times \real^n$ by
	$$ T(g,v) = \bigl( g, \widehat\varphi(g) v \bigr) .$$
Then, for any $\gamma \in \Gamma'$, a straightforward calculation shows
	\begin{align} \label{FlatBundleTequi}
	 T \bigl(\, (g,v) \cdot \gamma \, \bigr) 
	 &= T(g,v) * \gamma
	 , \text{ where $(g,v) * \gamma = (g \gamma, v)$}
	 \end{align}
%	\begin{align*}
%	 T \bigl(\, (g,v) \cdot \gamma \, \bigr) 
%	 &= T \bigl(\, g \gamma, \varphi(\gamma^{-1}) v \, \bigr) 
%	 =  \bigl(\, g \gamma, \hat\varphi( g \gamma)   \hat\varphi(\gamma^{-1}) v \, \bigr) 
%	 \\&= \bigl(\,g \gamma, \hat\varphi( g \gamma\gamma^{-1}) v \, \bigr) 
%	 = \bigl(\,g \gamma, \hat\varphi( g)  v \, \bigr) 
%	 \\&=  \bigl(\, g, \hat\varphi( g)  v \, \bigr) * \gamma
%	 = T(g,v) * \gamma
%	\end{align*}
\csee{FlatBundleTequiEx}.
Therefore $T$ factors through to a well-defined bundle isomorphism $\bundle_{\varphi'} \stackrel{\cong}{\to} (G/\Gamma') \times \real^n$.
\end{proof}

\subsection{Embeddings of locally symmetric spaces from embeddings of lattices} \label{TotGeodSect}

Let $M = \Gamma \backslash G / K$ and $M' = \Gamma ' \backslash G' \! / K'$.
Roughly speaking, the Mostow Rigidity Theorem \pref{MostowRigidity} tells us that if $\Gamma$ is isomorphic to~$\Gamma'$, then $M$ is isometric to~$M'$.  More generally, superrigidity implies that if $\Gamma$ is isomorphic to a subgroup of~$\Gamma'$, then $M$ is isometric to a submanifold of~$M'$ (modulo finite covers).

\begin{prop} \label{TotGeodProp}
Suppose 
	\begin{itemize}
	\item $M = \Gamma \backslash G / K$ and $M' = \Gamma ' \backslash G' \! / K'$ are irreducible locally symmetric spaces with no compact factors,
	\item $\Gamma$ is isomorphic to a subgroup of\/~$\Gamma'$,
	and
	\item the universal cover of~$M$ is neither the real hyperbolic space~$\hyperbolic^n$ nor the complex hyperbolic space~$\complex\hyperbolic^n$.
	\end{itemize}
Then some finite cover of\/ $\Gamma \backslash G / K$ embeds as a totally geodesic submanifold of a finite cover of\/ $\Gamma' \backslash G' \! / K'$.
\end{prop}

\begin{proof}[Idea of proof]
There is no harm in assuming that $G$ and~$G'$ have trivial center and no compact factors.
After passing to a finite-index subgroup of~$\Gamma$ 
	(and ignoring a compact group~$C$),
the Margulis Superrigidity Theorem tells us that the embedding $\Gamma \hookrightarrow \Gamma'$ extends to a continuous embedding $\varphi \colon G \hookrightarrow G'$. Conjugate $\varphi$ by an element of~$G'$, so that $\varphi(K) \subseteq K'$, and $\varphi(G)$ is invariant under the Cartan involution of~$G'$ corresponding to the maximal compact subgroup~$K'$. Then $\varphi$ induces an embedding $\Gamma \backslash G / K \to \Gamma ' \backslash G' \! / K'$ whose image is a totally geodesic submanifold.
\end{proof}

\begin{exercises}

\item \label{ProveMostMostowEx}
Prove the Mostow Rigidity Theorem \pref{MostowRigidityIrred} under the additional assumption that $G_1$ is neither $\PSO(1,n)$ nor $\PSU(1,n)$.

\item The statement of the Mostow Rigidity Theorem in \cref{MostowRigidity} is slightly different from \cref{MostowRigidityIrred}. Show that these two theorems are corollaries of each other.
\hint{\Cref{MostowRigIrredEnough}.}

\item \label{FlatVecBundleDefnEx}
In the notation of \cref{FlatVecBundleDefn}:
\noprelistbreak
	\begin{enumerate}
	\item Show the map~$\pi$ is well defined.
	\item Show $\bundle_\varphi$ is a vector bundle over $G/\Gamma$ with fiber $\real^n$.
	\end{enumerate}

\item \label{FlatBundleTequiEx}
Verify \pref{FlatBundleTequi} for all $\gamma \in \Gamma'$.

\end{exercises}

\section{Why superrigidity implies arithmeticity}
 \label{MargArithPf}

Recall the following major theorem that was stated without proof in \cref{MargulisArith}:

\begin{namedthm}[\thmindex{Margulis!Arithmeticity}{Margulis Arithmeticity Theorem}]
\label{MargArithFromSuper}
Every irreducible lattice in~$G$ is arithmetic, except, perhaps, when $G$ is isogenous to\/ $\SO(1,m) \times K$ or\/ $\SU(1,m) \times K$, for some compact group~$K$.
\end{namedthm}

This important fact is a consequence of the Margulis Superrigidity Theorem, but the implication is not at all obvious.  In this section, we will explain the main ideas that are involved.

In addition to our usual assumption that $G \subseteq \SL(\ell,\real)$, let us also assume, for simplicity:
\noprelistbreak
 \begin{itemize}
% \item $G \subseteq \SL(\ell,\real)$, for some~$\ell$, 
 \item $G \iso \SL(3,\real)$ (or, more generally, $G$ is algebraically simply connected; see \cref{AlgSCDefn}),
 and 
 \item $G / \Gamma$ is not compact. 
 \end{itemize}
 We wish to show that $\Gamma$ is arithmetic. It suffices to
show $\Gamma  \subseteq G_{\integer}$, that is, that every
matrix entry of every element of~$\Gamma$ is an integer, for
then $\Gamma$ is commensurable to $G_{\integer}$
\csee{finext->latt}.

Here is a loose description of the 4 steps of the proof:
	\begin{enumerate}
	\item The Margulis Superrigidity Theorem \pref{MargSuperSL3RNonCpct} implies that every matrix entry of every element of~$\Gamma$ is an algebraic number.
	\item By \term[Restriction of Scalars]{restriction of scalars}, we may assume that these algebraic numbers are rational; that is, $\Gamma \subseteq G_{\rational}$.
	\item For every prime~$p$, a ``$p$-adic'' version of the 
	\thmindex{Margulis!Superrigidity!p-adic@$p$-adic}Margulis Superrigidity Theorem 
	provides a natural number~$N_p$, such that no element of~$\Gamma$ has a matrix entry whose denominator is divisible by~$p^{N_p}$. 
	\item This implies that some finite-index subgroup~$\Gamma'$ of~$\Gamma$ is contained in~$G_{\integer}$.
	\end{enumerate}

\setcounter{step}{0}

\begin{step} \label{ArithThmPf-algic}
 Every matrix entry of every element of\/~$\Gamma$ is an
algebraic number.
 \end{step}
 Suppose some $\gamma_{i,j}$ is transcendental.
 Then, for any transcendental number~$\alpha$, there is a
field automorphism~$\phi$ of~$\complex$ with
$\phi(\gamma_{i,j}) = \alpha$. Applying~$\phi$ to all the
entries of a matrix induces an automorphism~$\widetilde\phi$
of $\SL(\ell,\complex)$. Let
	$$ \text{$\varphi$ be the restriction of~$\widetilde\phi$ to~$\Gamma$,} $$
so $\varphi$ is a homomorphism from~$\Gamma$ to $\SL(\ell,\complex)$.
The Margulis Superrigidity Theorem implies there is a
continuous homomorphism $\widehat\varphi \colon G \to
\SL(\ell,\complex)$, such that $\widehat\varphi = \varphi$ on
a finite-index subgroup of~$\Gamma$ \csee{MargSuperNonCpctSCEx}. By passing to this finite-index subgroup,
we may assume $\widehat\varphi = \varphi$ on all of~$\Gamma$.

Since there are uncountably many transcendental
numbers~$\alpha$, there are uncountably many different
choices of~$\phi$, so there must be uncountably many
different $n$-dimensional representations~$\widehat\varphi$
of~$G$. However, it is well known from the the theory of
``roots and weights'' that $G$ (or, more generally, any
connected, simple Lie group) has
only finitely many non-isomorphic representations
of any given dimension, so this is a contradiction.%
\footnote{Actually, this is not quite a
contradiction, because it is possible that two different
choices of~$\varphi$ yield the same representation of~$\Gamma$,
up to isomorphism; that is, after a change of basis. The
trace of a matrix is independent of the basis, so the
preceding argument really shows that the trace
of~$\varphi(\gamma)$ must be algebraic, for every
$\gamma \in \Gamma$. Then one can use some algebraic methods
to construct some other matrix representation~$\varphi'$
of~$\Gamma$, such that the matrix entries of~$\varphi'(\gamma)$
are algebraic, for every $\gamma \in \Gamma$.}

\begin{step}
 We have\/ $\Gamma \subseteq \SL(\ell,\rational)$.
 \end{step}
 Let $F$ be the subfield of~$\complex$ generated by the 
matrix entries of the elements of~$\Gamma$, so $\Gamma 
\subseteq \SL(\ell,F)$. From \cref{ArithThmPf-algic}, we know
that this is an algebraic extension of~$\rational$.
Furthermore, because $\Gamma$ is finitely generated \csee{GammaFinGen}, 
we see that this field extension is finitely
generated. Therefore, $F$ is finite-degree field extension
of~$\rational$ (in other words, $F$ is an algebraic number
field). This means that $F$ is almost the same
as~$\rational$, so it is only a slight exaggeration to say
that we have proved $\Gamma  \subseteq \SL(\ell,\rational)$.

Indeed, restriction of scalars \pref{ResScal->Latt} provides a way to change $F$
into~$\rational$: there is a representation $\rho \colon G
\to \SL(r,\complex)$, for some~$r$, such that 
$\rho \bigl( G \cap \SL(\ell,F) \bigr) \subseteq \SL(r,\rational)$ 
\csee{ROSPutsGFinGQ}. Therefore, after
replacing~$G$ with $\rho(G)$, we have the
desired conclusion (without any exaggeration).

\begin{step} \label{MargArithPf-BddPowerP}
For every prime~$p$, there is a natural number~$N_p$, such 
that no element of\/~$\Gamma$ has a matrix entry whose 
denominator is divisible by~$p^{N_p}$. 
\end{step}
The fields $\real$ and~$\complex$ are complete (that is,
every Cauchy sequence converges), and they obviously 
contain~$\rational$. For any prime~$p$, the $p$-adic 
numbers~$\rational_p$ are another field that has these
same properties.

As we have stated it, the Margulis Superrigidity Theorem
deals with homomorphisms into $\SL(\ell,\F)$, where $\F = \real$,
but Margulis also proved a version of the theorem
that applies when $\F$ is a $p$-adic field \csee{padicSuper}. 
Now $G$ is connected,
but $p$-adic fields are totally disconnected, so every continuous 
homomorphism from~$G$ to $\SL(\ell,\rational_p)$ is trivial.
Therefore, superrigidity tells us that $\varphi$ is trivial, after we mod
out a compact group \ccf{MargSuperC}.
In other words, the closure of 
$\varphi(\Gamma)$ is compact in $\SL(\ell,\rational_p)$.

This conclusion can be rephrased in more elementary terms,
without any mention of $p$-adic
numbers. Namely, it says that there is a bound on the 
highest power of~$p$ that divides the denominator of any matrix entry of
any element of~$\Gamma$. This is what we wanted.

\begin{step}
Some finite-index subgroup\/~$\Gamma'$ of\/~$\Gamma$ is 
contained in\/ $\SL(\ell,\integer)$.
\end{step}
 Let $D \subseteq \natural$ be the set consisting of the
denominators of the matrix entries of the elements of
$\varphi(\Gamma)$. 

We claim there exists $N \in \natural$, such that every
element of~$D$ is less than~$N$.
Since $\Gamma$ is known to be finitely generated,
some finite set of primes $\{p_1,\ldots,p_r\}$ contains all
the prime factors of every element of~$D$. (If $p$~is in the
denominator of some matrix entry of $\gamma_1
\gamma_2$, then it must
appear in a denominator somewhere in either $\gamma_1$
or~$\gamma_2$.) Therefore, every element of~$D$ is of the
form $p_1^{m_1} \cdots p_r^{m_r}$, for some $m_1,\ldots,m_r
\in \natural$. From \cref{MargArithPf-BddPowerP}, 
we know $m_i < N_{p_i}$,
for every~$i$. Thus, every element of~$D$ is less than 
$p_1^{N_{p_1}} \cdots p_r^{N_{p_r}}$. This establishes
the claim.

 From the preceding paragraph, we see that 
 $\Gamma \subseteq \frac{1}{N!} \Mat_{\ell
\times \ell}(\integer)$. Note that if $N = 1$, then $\Gamma
\subseteq \SL(\ell,\integer)$. In general, $N$ is a finite 
distance from~$1$, so it should not be hard to believe
(and it can indeed be shown) that some 
finite-index subgroup of~$\Gamma$ must be contained
in $\SL(\ell,\integer)$ \csee{GammaBddDenomsEx}. Therefore, a finite-index subgroup
of~$\Gamma$ is contained in~$G_{\integer}$, as desired.
\qed

\medbreak

For ease of reference, we officially record the key fact used in \cref{MargArithPf-BddPowerP}:

\begin{thm}[(\thmindex{Margulis!Superrigidity!p-adic@$p$-adic}%
	Margulis superrigidity over $p$-adic fields)] \label{padicSuper}
Assume
\noprelistbreak
	\begin{enumerate} \renewcommand{\theenumi}{\roman{enumi}}
	\item $G$ is not isogenous to any group that is of the form\/ $\SO(1,m) \times K$ or\/ $\SU(1,m) \times K$, where $K$~is compact,
	\item $\Gamma$ is irreducible,
	\item $\rational_p$ is the field of $p$-adic numbers, for some prime~$p$,
	and
	\item $\varphi \colon \Gamma \to \GL(n,\rational_p)$ is a homomorphism.
	\end{enumerate}
Then $\closure{\varphi(\Gamma)}$ is compact.

In other words, there is some $N \in \integer$, such that every matrix entry of every element of $\varphi(\Gamma)$ is in $p^N \ \integer_p$, where $\integer_p$ is the ring of $p$-adic integers.
\end{thm}

The Margulis Arithmeticity Theorem \pref{MargArithFromSuper} does not apply to lattices in $\SO(1,n)$ or $\SU(1,n)$, but, for those groups, Margulis proved the following characterization of the lattices that are arithmetic:

\begin{namedthm}[\thmindex{Commensurability Criterion for Arithmeticity}{Commensurability Criterion for Arithmeticity} {\normalfont (Margulis)}] \label{CommCriterion}
Assume
	\begin{itemize}
	\item $G$ is connected, with no compact factors,
	and
	\item $\Gamma$ is irreducible.
	\end{itemize}
Then\/ $\Gamma$ is arithmetic if and only if the commensurator\/ $\Comm_G(\Gamma)$ of\/~$\Gamma$ is dense in~$G$.
\end{namedthm}

As was already mentioned in \fullcref{GQComm}{criterion}, the direction ($\Rightarrow$) follows from the simple observation that $\Comm_G(G_\integer)$ contains $G_\rational$. 

The proof of ($\Leftarrow$) is more difficult. It is the same as the proof of the Margulis Arithmeticity Theorem, but replacing the Margulis Superrigidity Theorem \pref{MargSuperG'} with the following superrigidity theorem (and also replacing the $p$-adic superrigidity theorem with a suitable commensurator analogue):

\begin{thm}[(\thmindex{Commensurator Superrigidity}{Commensurator Superrigidity})] \label{CommSuper}
Assume
\noprelistbreak
	\begin{enumerate} \renewcommand{\theenumi}{\roman{enumi}}
	\item $\Gamma$ is irreducible, 
	\item $\Comm_G(\Gamma)$ is dense in~$G$,
 	and
	\item  $G$ and $G'$ are connected, with trivial center, and no compact factors.
	\end{enumerate}
If $\varphi \colon \Comm_G(\Gamma) \to G'$ is any homomorphism whose image is Zariski dense in~$G'$, then $\varphi$ extends to a continuous homomorphism $\widehat\varphi \colon G \to G'$.
\end{thm}

\begin{exercises}

\item \label{ROSPutsGFinGQ}
Suppose 
	\begin{itemize}
	\item $G \subseteq \SL(\ell,\complex)$, 
	\item $\Gamma \subseteq \SL(\ell,F)$, for some algebraic number field~$F$,
	and
	\item $G$ has no compact factors.
	\end{itemize}
Show there is a continuous homomorphism $\rho \colon G
\to \SL(r,\complex)$, for some~$r$, such that $\rho \bigl( G \cap
\SL(\ell,F) \bigr) \subseteq \SL(r,\rational)$.
\hint{Apply restriction of scalars (\S\ref{RestrictScalarsSect}) after noting that the Borel Density Theorem (\S\ref{BDTSect}) implies $G$ is defined over~$F$.}

\item \label{GammaBddDenomsEx}
Show that if $\Lambda$ is a subgroup of $\SL(\ell,\rational)$, and $\Lambda \subseteq \frac{1}{N} \Mat_{\ell\times\ell}(\integer)$, for some $N \in \natural$, then $\SL(\ell,\integer)$ contains a finite-index subgroup of~$\Lambda$.
\hint{The additive group of~$\rational^\ell$ contains a $\Lambda$-invariant subgroup~$V$, such that we have $\integer^\ell \subseteq V \subseteq \frac{1}{N} \integer^\ell$. Choose $g \in \GL(\ell,\rational)$, such that $g(V) = \integer^\ell$. Then $g$ commensurates $\SL(\ell,\integer)$ and we have $g \Lambda g^{-1} \subseteq \SL(\ell,\integer)$.}

%\item \label{BddDenom->finind}
% Suppose $\Lambda$ is a subgroup of $\SL(n,\rational)$, and
%$k$~is a positive integer, such that $k \lambda \in \Mat_{n
%\times n}(\integer)$ for every $\lambda \in \Lambda$. Show
%that $\Lambda \cap \SL(n,\integer)$ has finite index
%in~$\Lambda$.
% \hint{If $k \gamma \equiv k \lambda \pmod{k}$,
%then $\gamma \lambda^{-1} \in \Mat_{n \times
%n}(\integer)$.}

\item \label{EigenValsAreInts}
%(\emph{requires some Algebraic Number Theory})
Assume, as usual, that 
	\begin{itemize}
	\item $G$ is not isogenous to any group that is of the form\/ $\SO(1,m) \times K$ or\/ $\SU(1,m) \times K$, where $K$~is compact,
	and
	\item $\Gamma$ is irreducible.
	\end{itemize}
Use the proof of the Margulis Arithmeticity Theorem to show that if $\varphi \colon \Gamma \to \SL(n,\complex)$ is any homomorphism, then every eigenvalue of every element of $\varphi(\Gamma)$ is an algebraic integer.
%\hint{After a change of basis, and passing to a subgroup of finite index, every entry of every matrix in $\varphi(\Gamma)$ is an algebraic integer.}

\end{exercises}

\section{Homomorphisms into compact groups} \label{SuperCpctSect}

The Margulis Superrigidity Theorem \pref{MargSuperC} does not say anything about homomorphisms whose image is contained in a compact subgroup of $\GL(n,\real)$. (This is because all of $\varphi(\Gamma)$ can be put into the error term~$C$, so the homomorphism $\widehat\alpha$ can be taken to be trivial.) Fortunately, there is a different version that completely eliminates the error term (and applies very generally). Namely, from the Margulis Arithmeticity Theorem \pref{MargulisArith}, we know that the lattice $\Gamma$ must be arithmetic (if no simple factors of~$G$ are $\SO(1,m)$ or $\SU(1,m)$). This means that if we add some compact factors to~$G$, then we can assume that $\Gamma$ is commensurable to~$G_{\integer}$. In this situation, there is no need for the error term~$C$:

\begin{cor} \label{GZSuper}
Assume
\noprelistbreak
	\begin{itemize}
	\item $G$ is connected, algebraically simply connected, and defined over~$\rational$,
	\item $G$ is not isogenous to any group that is of the form\/ $\SO(1,m) \times K$ or\/ $\SU(1,m) \times K$, where $K$~is compact,
	\item $\Gamma$ is irreducible,
	and
	\item $\varphi \colon \Gamma \to \GL(n,\real)$ is a homomorphism.
	\end{itemize}
If\/ $\Gamma$ is commensurable to~$G_{\integer}$, then there exist:
\noprelistbreak
	\begin{enumerate}
	\item a continuous homomorphism $\widehat\varphi \colon G \to \GL(n,\real)$,
	and
	\item a finite-index subgroup\/~$\Gamma'$ of\/~$\Gamma$,
	\end{enumerate}
such that $\varphi(\gamma) = \widehat\varphi(\gamma)$, for all $\gamma \in \Gamma'$.
\end{cor}

Here is a less precise version of \cref{GZSuper} that may be easier to apply in situations where the lattice~$\Gamma$ is not explicitly given as the $\integer$-points of~$G$.  
%It states that homomorphisms of~$\Gamma$ into compact groups are Galois twists of homomorphisms that extend to~$G$ (if we consider each simple factor of $\closure{\alpha(\Gamma)}$ individually).
However, it only applies to the homomorphism into each simple component of $\closure{\alpha(\Gamma)}$, not to the entire homomorphism all at once.

\begin{cor} \label{MargSuperCpct}
Assume
\noprelistbreak
	\begin{itemize}
	\item $G$ is connected, and algebraically simply connected, 
	\item $G$ is not isogenous to any group that is of the form\/ $\SO(1,m) \times K$ or\/ $\SU(1,m) \times K$, where $K$~is compact,
	\item $\Gamma$ is irreducible,
	\item $\varphi \colon \Gamma \to \GL(n,\complex)$ is a homomorphism,
	and
	\item $\closure{\varphi(\Gamma)}$ is simple.
	\end{itemize}
Then there exist:
\noprelistbreak
	\begin{enumerate}
	\item a continuous homomorphism $\widehat\varphi \colon G \to \GL(n,\complex)$,
	\item a finite-index subgroup\/~$\Gamma'$ of\/~$\Gamma$,
	and
	\item a Galois automorphism~$\sigma$ of\/~$\complex$,
	\end{enumerate}
such that $\varphi(\gamma) = \sigma \bigl( \widehat\varphi(\gamma) \bigr)$, for all $\gamma \in \Gamma'$.
\end{cor}

\begin{proof}
We may assume $\closure{\varphi(\Gamma)}$ is compact, for otherwise \cref{MargSuperG'} applies (after modding out the centers of $G$ and $\closure{\alpha(\Gamma)}$). Then every element of $\closure{\varphi(\Gamma)}$ is semisimple. 

Choose some $h \in \varphi(\Gamma)$, such that $h$ has infinite order \csee{ImgNotTorsion}. Then the conclusion of the preceding paragraph implies that some eigenvalue~$\lambda$ of~$h$ is not a root of unity. On the other hand, if $\lambda$~is algebraic, then $p$-adic superrigidity \pref{padicSuper} implies that $\lambda$~is an algebraic integer \csee{EigenValsAreInts}. So there is a Galois automorphism~$\sigma$ of~$\complex$, such that $|\sigma(\lambda)| \neq 1$ \csee{KroneckerNot1}. Then $\{\, \sigma(\lambda)^k \mid k \in \integer \,\}$ is an unbounded subset of~$\complex$, so $\bigl\langle \sigma(h) \bigr\rangle$ is not contained in any compact subgroup of $\GL(n,\complex)$.

Now, let 
\noprelistbreak
	\begin{itemize}
	\item $\varphi'$ be the composition $\sigma \circ \varphi$,
	and
	\item $G'$ be the Zariski closure of $\varphi'(\Gamma)$.
	\end{itemize}
Then $G'$ is simple, and the conclusion of the preceding paragraph implies that $G'$ is not compact (since $\sigma(h) \in G'$). After passing to a finite-index subgroup (so $G'$ is connected), \cref{MargSuperG'} provides a continuous homomorphism $\widetilde\varphi \colon G \to G'$, such that $\varphi'(\gamma) = \widehat\varphi(\gamma)$, for all $\gamma$ in some finite-index subgroup of~$\Gamma$.
\end{proof}

\begin{warn}
Assume $\Gamma$ is irreducible, and $G$ is not isogenous to any group of the form $\SO(1,m) \times K$ or $\SU(1,m) \times K$.
\Cref{MargSuperCpct} implies that if there exists a homomorphism~$\varphi$ from~$\Gamma$ to a compact Lie group (and $\varphi(\Gamma)$ is infinite), then $G/\Gamma$ must be compact \csee{MargNoncpct}. However, \textbf{the converse is not true}.  Namely, \cref{GZSuper} tells us that if $\Gamma$ is commensurable to $G_{\integer}$, where $G$~is defined over~$\rational$, and $G_{\real}$ has no compact factors, then $\Gamma$~does not have any homomorphisms to compact groups (with infinite image). It does not matter whether $G/\Gamma$ is compact or not.
\end{warn}

\begin{exercises}

\item \label{MargNoncpct}
Assume, as usual, that the lattice $\Gamma$~is irreducible, that $G$~is not isogenous to any group of the form $\SO(1,m) \times K$ or $\SU(1,m) \times K$, and that $\varphi \colon \Gamma \to \GL(n,\real)$ is a homomorphism. If $G/\Gamma$ is \textbf{not} compact, show the semisimple group $\Zar{\varphi(\Gamma)}$ has no compact factors.
\hint{Godement's Criterion \pref{GodementCriterion}.}

\item \label{MargSuperNonCpctSCEx}
Assume
	\begin{itemize}
	\item $G$ is algebraically simply connected,
	\item $G$ is not isogenous to any group that is of the form\/ $\SO(1,m) \times K$ or\/ $\SU(1,m) \times K$, where $K$~is compact,
	\item $\Gamma$ is irreducible,
	\item $G/\Gamma$ is not compact,
	and
	\item $\varphi \colon \Gamma \to \SL(n,\real)$ is a homomorphism.
	\end{itemize}
Show there is a continuous homomorphism $\widehat\varphi \colon G \to \SL(n, \real)$, such that $\varphi(\gamma) = \widehat\varphi(\gamma)$ for all $\gamma$ in some finite-index subgroup of~$\Gamma$.
\hint{\cref{MargSuperC,MargImgSS,MargNoncpct}.}

%\item \label{SuperSLnEx}
%Prove \cref{MargSuperSL3RNonCpct}.
%\hint{\cref{MargSuperC,MargImgSS,SLnASC,MargNoncpct}.}

\item 
Assume $\Gamma$ is irreducible, and $G$ has no factors isogenous to $\SO(1,m)$ or $\SU(1,m)$.
Show that if $N$ is an infinite normal subgroup of~$\Gamma$, such that $\Gamma/N$ is linear (i.e., isomorphic to a subgroup of $\GL(\ell,\complex)$, for some~$\ell$), then $\Gamma/N$ is finite.

\item \label{KroneckerNot1}
(\emph{\thmindex{Kronecker's}{Kronecker's Theorem}})
Assume $\lambda$ is an algebraic integer. Show that if $|\sigma(\lambda)| = 1$ for every Galois automorphism~$\sigma$ of~$\complex$, then $\lambda$ is a root of unity.
\hint{The powers of $\lambda$ form a set that (by restriction of scalars) is discrete in $\bigtimes_{\sigma \in S^\infty} F_\sigma^\times$. Alternate proof: there are only finitely many polynomials of degree $n$ with integer coefficients that are all $\le C$ in absolute value.}

\end{exercises}

\section{Proof of the Margulis Superrigidity Theorem} \label{SuperPfSect}

In order to establish \cref{MargSuperG'}, it suffices to prove the following special case \csee{MargSuperHSimpleEx}:

\begin{thm} \label{MargSuperHSimple}
Suppose
\noprelistbreak
\begin{itemize}
\item $G$ is connected, and it is not isogenous to any group that is of the form\/ $\SO(1,m) \times K$ or\/ $\SU(1,m) \times K$, where $K$~is compact, 
\item the lattice $\Gamma$ is irreducible in~$G$,
\item $H$ is a connected, noncompact, simple subgroup of\/ $\SL(n,\real)$, for some~$n$ {\rm(}and $H$ has trivial center{\rm)},
\item $\varphi \colon \Gamma \to H$ is a homomorphism,
and
\item $\varphi(\Gamma)$ is Zariski dense in~$H$.
\end{itemize}
Then $\varphi$ extends to a continuous homomorphism $\widehat\varphi \colon G \to H$.
\end{thm}

Although it does result in some loss of generality, we assume:

\begin{assump} \label{MargSuperRrankAssump}
$\Rrank G \ge 2$.
\end{assump}

The case where $\Rrank G = 1$ requires quite different methods --- see \cref{SuperRank1Sect} for a very brief discussion.

\subsection{Geometric reformulation}
To set up the proof of \cref{MargSuperHSimple}, let us translate the problem into a geometric setting, by replacing the homomorphism~$\varphi$ with the corresponding flat vector bundle~$\bundle_\varphi$ over~$G/\Gamma$ \csee{FlatVecBundleDefn}.

\begin{rem} \label{pfsuper-sect<>equi}
The sections of the vector bundle~$\bundle_\varphi$ are in natural one-to-one correspondence with the right $\Gamma$-equivariant maps from~$G$ to~$\real^n$ \csee{Equi<>SectionEx}.
\end{rem}

\begin{lem} \label{Extend<>SectionBij}
$\varphi$ extends to a homomorphism $\cover\varphi \colon G \to \GL(n,\real)$ if and only if there exists a $G$-invariant subspace $V \subseteq \Sect(\bundle_\varphi)$, such that the evaluation map $V \to V_{[e]}$ is bijective.
\end{lem}

\begin{proof}
($\Leftarrow$) Since $V$ is $G$-invariant, we have a representation of~$G$ on~$V$; let us say $\pi \colon G \to \GL(V)$. Therefore, the isomorphism $V \to V_{[e]} = \real^n$ yields a representation $\widehat\pi$ of~$G$ on~$\real^n$. It is not difficult to verify that $\widehat\pi$ extends~$\varphi$ \csee{PiHatExtendsPhi}.

($\Rightarrow$) For $v \in \real^n$ and $g \in G$, let 
	$$\xi_v(g) = \cover\varphi(g^{-1}) v .$$
It is easy to verify that $\xi_v \colon G \to \real^n$ is right $\Gamma$-equivariant \csee{XivEqui}, so we may think of~$\xi_v$ as a section of~$\bundle_\varphi$ \csee{pfsuper-sect<>equi}. Let 
	$$V = \{\, \xi_v \mid v \in \real^n\,\} \subseteq \Sect(\bundle_\varphi) .$$
Now the map $v \mapsto \xi_v$ is linear and $G$-equivariant \csee{XiGEqui}, so $V$ is a $G$-invariant subspace of $\Sect(\bundle_\varphi)$. Since 
	$$ \xi_v \bigl( [e] \bigr) = \cover\varphi(e) v = v ,$$
it is obvious that the evaluation map is bijective.
\end{proof}

In fact, if we assume the representation $\varphi$ is irreducible, then it is not necessary to have the evaluation map $V \to V_{[e]}$ be bijective. Namely, in order to show that $\varphi$ extends, it suffices to have $V$ be finite dimensional (and nonzero):

\begin{lem}[\csee{Extend<>FDEx}] \label{Extend<>FD}
Assume that the representation~$\varphi$ is irreducible. If there exists a\/ \textup(nontrivial\/\textup) $G$-invariant subspace $V$ of\/ $\Sect(\bundle_\varphi)$ that is finite dimensional,
then $\varphi$ extends to a continuous homomorphism $\cover\varphi \colon G \to \GL(n,\real)$.
\end{lem}

\subsection{The need for higher real rank}

We now explain how \cref{MargSuperRrankAssump} comes into play. 

\begin{notation}
Let $A$ be a maximal $\real$-split torus of~$G$. For example, if $G = \SL(3,\real)$, we let
	$$ A = \begin{Smallbmatrix} \upast&0&0 \\ 0&\upast&0 \\ 0&0&\upast \end{Smallbmatrix} .$$
By definition, the assumption that $\Rrank G \ge 2$ means $\dim A \ge 2$.
\end{notation}

It is the following result that relies on our assumption $\Rrank G \ge 2$. It is easy to prove if $G$ has more than one noncompact simple factor \csee{G1xG2GenLiEx}, and is not difficult to verify for the case $G = \SL(\ell,\real)$ \ccf{SL3RGenLiEx}.  Readers familiar with the structure of semisimple groups (including the theory of real roots) should have little difficulty in generalizing to any semisimple group of real rank $\ge 2$ \csee{HighRankGenLiEx}.

\begin{lem} \label{HighRankGenLi}
If\/ $\Rrank G \ge 2$, then, for some $r \in \natural$, there exist closed subgroups $L_1,L_2,\ldots,L_r$ of~$G$, such that
	\begin{enumerate}
%	\item $L_i$ is a closed subgroup of~$G$ that is isogenous to $\SL(2,\real)$,
	\item $ G = L_r L_{r-1} \cdots L_1$,
	and
	\item both $H_i$ and $H_i^\perp$ are noncompact, where
		\begin{itemize}
		\item $H_i = L_i \cap A$,
		and
		\item $H_i^\perp = \czer_A(L_i)$ {\rm(}so $L_i$ centralizes~$H_i^\perp${\rm)}.
		\end{itemize}
	\end{enumerate}
\end{lem}

\subsection{Outline of the proof}

The idea for proving \cref{MargSuperHSimple} is quite simple. We begin by finding a (nonzero) $A$-invariant section of~$\bundle_\varphi$; this section spans a ($1$-dimensional) subspace~$V_0$ of~$\bundle_\varphi$ that is invariant under~$A$. Since (by definition) the subgroup~$H_1$ of \cref{HighRankGenLi} is contained in~$A$, we know that $V_0$ is invariant under~$H_1$, so \cref{pfsuper-C(H)V} below provides a subspace of $\Sect(\bundle_\varphi)$ that is invariant under a larger subgroup of~$G$, but is still finite dimensional. Applying the lemma repeatedly yields finite-dimensional subspaces that are invariant under more and more of~$G$. Eventually, the lemma yields a finite-dimensional subspace that is invariant under all of~$G$. Then \cref{Extend<>FD} implies that $\varphi$ extends to a homomorphism that is defined on~$G$, as desired.

\begin{lem} \label{pfsuper-C(H)V}
If
\noprelistbreak
	\begin{itemize}
	\item $H$ is a closed, noncompact subgroup of~$A$,
	and
	\item $V$ is an $H$-invariant subspace of\/ $\Sect(\bundle_\varphi)$ that is finite dimensional,
	\end{itemize}
then $\langle \czer_G(H) \cdot V \rangle$ is finite dimensional.
\end{lem}

\begin{proof}[Idea of proof]
To illustrate the idea of the proof, let us assume that $V = \real \sigma$ is the span of an $H$-invariant section \csee{pfsuper-C(H)VEx}.
Since $H$ is noncompact, the Moore Ergodicity Theorem \pref{MooreErgodicity} tells us that $H$ has a dense orbit on $G/\Gamma$ \csee{AEOrbitDenseInG/Gamma}. (In fact, almost every orbit is dense.) This implies that any continuous $H$-invariant section of~$\bundle_\varphi$ is determined by its value at a single point \csee{InvtSectDetPtEx}, so the space of $H$-invariant sections is finite-dimensional \csee{InvtSectFDEx}. Since this space contains $\langle \czer_G(H) \cdot V \rangle$ \csee{CG(H).HInvtEx}, the desired conclusion is immediate.
\end{proof}

Here is a more detailed outline:

\begin{proof}[Idea of the proof of \cref{MargSuperHSimple}]
Assume there exists a nonzero $A$-invariant section~$\sigma$ of~$\bundle_\varphi$. Let 
	$$ \text{$H_0 = A$ and $V_0 = \langle \sigma \rangle$.} $$
Thus, $V_0$ is a $1$-dimensional subspace of $\Sect(\bundle_\varphi)$ that is $H_0$-invariant.

Now, for $i = 1,\ldots,r$, let
	$$ V_i = \langle L_i \cdot A \cdot L_{i-1} \cdot A \cdots L_1 \cdot A \cdot V_0 \rangle .$$
Since $L_r L_{r-1} \cdots L_1 = G$, it is clear that $V_r$ is $G$-invariant. Therefore, it will suffice to show (by induction on~$i$) that each $V_i$ is finite dimensional.

Since $H_{i-1} \subseteq L_{i-1}$, it is clear that $V_{i-1}$ is $H_{i-1}$-invariant. Therefore, since $A$~centralizes~$H_{i-1}$, \Cref{pfsuper-C(H)V} implies that $\langle A \cdot V_{i-1} \rangle$ is finite dimensional. Now, since $H_i^\perp \subseteq A$, we know that $\langle A \cdot V_{i-1} \rangle$ is $H_i^\perp$-invariant. Then, since $L_i$ centralizes~$H_i^\perp$, \Cref{pfsuper-C(H)V} implies that the subspace $V_i = \langle  L_i \cdot A \cdot V_{i-1} \rangle$ is finite dimensional.
\end{proof}

Therefore, the key to proving \cref{MargSuperHSimple} is finding a nonzero $A$-invariant section~$\sigma$ of~$\bundle_\varphi$. Unfortunately, the situation is a bit more complicated than the above would indicate, because  we will not find a \emph{continuous} $A$-invariant section, but only a \emph{measurable} one \csee{pfsuper-keyfact-meas}. Then the proof appeals to \cref{pfsuper-C(H)V-meas} below, instead of \cref{pfsuper-C(H)V}. We leave the details to the reader \csee{ProveMargSuperHSimpleEx}.

\begin{defn}
Let $\Sectm(\bundle_\varphi)$ be the vector space of measurable sections of~$\bundle_\varphi$, where two sections are identified if they agree almost everywhere.
\end{defn}

\begin{lem}[\csee{pfsuper-C(H)V-measPfEx,pfsuper-C(H)V-measFullPfEx}] \label{pfsuper-C(H)V-meas}
If
\noprelistbreak
	\begin{itemize}
	\item $H$ is a closed, noncompact subgroup of~$A$,
	and
	\item $V$ is a finite-dimensional, $H$-invariant subspace of\/ $\Sectm(\bundle_\varphi)$,
	\end{itemize}
then $\langle \czer_G(H) \cdot V \rangle$ is finite dimensional.
\end{lem}

\begin{exercises}

\item \label{MargSuperHSimpleEx}
Derive \cref{MargSuperG'} as a corollary of \cref{MargSuperHSimple}.

\item \label{Equi<>SectionEx}
Suppose $\xi \colon G \to \real^n$. Show that $\overline\xi \colon G/\Gamma \to \bundle_\varphi$, defined by
	$$ \overline\xi(g \Gamma) = \bigl[ \bigl( g,{\xi}(g) \bigr) \bigr] , $$
is a well-defined section of~$\bundle_\varphi$ if and only if $\xi$ is right $\Gamma$-equivariant; i.e., $\xi(g \gamma) = \varphi(\gamma^{-1}) \,\xi(g)$.

\item \label{PiHatExtendsPhi}
In the notation of the proof of \cref{Extend<>SectionBij}($\Leftarrow$), show $\widehat\pi(\gamma) = \varphi(\gamma)$ for every $\gamma \in \Gamma$.

\item \label{XivEqui}
In the notation of the proof of \cref{Extend<>SectionBij}($\Rightarrow$), show that we have $\xi_v(gh) = \cover\varphi(h^{-1}) \, \xi_v(g)$. Since $\cover\varphi(\gamma^{-1}) = \varphi(\gamma^{-1})$ for all $\gamma \in \Gamma$, this implies that $\xi_v$ is right $\Gamma$-equivariant.

\item \label{XiGEqui}
In the notation of the proof of \cref{Extend<>SectionBij}($\Rightarrow$), show that we have
$\xi_{\cover\varphi(g) v} = g \cdot \xi_v$, where the action of~$G$ on $\Sect(\bundle_\varphi)$ is defined by $(g \cdot \xi_v)(x) = \xi_v(g^{-1} x)$, as usual.

\item \label{Extend<>FDEx}
Prove \cref{Extend<>FD}.
\hint{By choosing $V$ of minimal dimension, we may assume it is an irreducible $G$-module, so the evaluation map is either $0$ or injective. It cannot be~$0$, and then it must also be surjective, since $\varphi$ is irreducible.}

\item \label{G1xG2GenLiEx}
Prove \cref{HighRankGenLi} under that additional assumption that we have $G = G_1 \times G_2$, where $G_1$ and~$G_2$ are noncompact (and semisimple).
\hint{Let $L_i = G_i$ for $i = 1,2$.}

\item \label{SL3RGenLiEx}
Prove the conclusion of \cref{HighRankGenLi} for $G = \SL(3,\real)$.
\hint{A \defit[unipotent!elementary matrix]{unipotent elementary matrix} is a matrix with $1$'s on the diagonal and only one nonzero off-diagonal entry. Every element of $\SL(3,\real)$ is a product of $\le 10$ unipotent elementary matrices, and any such matrix is contained in a subgroup isogenous to $\SL(2,\real)$ that has a $1$-dimensional intersection with~$A$.}

\item \label{HighRankGenLiEx}
Prove \cref{HighRankGenLi}.

\item \label{InvtSectDetPtEx}
Let $H$ be a subgroup of~$G$. Show that if $\sigma_1$ and~$\sigma_2$ are $H$-invariant, continuous sections of~$\bundle_\varphi$, and there is some $x \in G/\Gamma$, such that
	\begin{itemize}
	\item $Hx$ is dense in $G/\Gamma$
	and
	\item $\sigma_1(x) = \sigma_2(x)$,
	\end{itemize}
then $\sigma_1 = \sigma_2$.

\item \label{InvtSectFDEx}
Let $H$ be a subgroup of~$G$, and assume $H$ has a dense orbit in $G/\Gamma$. Show the space of $H$-invariant, continuous sections of $\bundle_\varphi$ has finite dimension.

\item \label{CG(H).HInvtEx}
Let $H$ be a subgroup of~$G$. Show that if $\sigma$ is an $H$-invariant section of~$\bundle_\varphi$, and $c$~is an element of~$G$ that centralizes~$H$, then $\sigma \cdot c$ is also $H$-invariant.

\item \label{pfsuper-C(H)VEx}
Prove \cref{pfsuper-C(H)V} without assuming that $V$ is $1$-dimensional.
\hint{Fix $x \in G$.
For $c \in C_G(H)$ and $\sigma \in V$, define $T \colon V \to \real^n$ by $T(\xi) = \xi (x)$, and note that $(c \sigma)(hx) = T(h^{-1} \cdot \sigma)$ for all $h \in H$. If $H x \Gamma$ is dense, this implies that $c\sigma$ is determined by $\sigma$ and~$T$. So $\dim \bigl( C_G(H) \cdot V) \le ( \dim V) \cdot \bigl( \dim \Hom(V,\real^n) \bigr)$.}
%For $h \in H$, we have $(c\sigma)[h] = \bigl( h^{-1}(c\sigma)\bigr) [e] = \bigl(c (h^{-1}\sigma) \bigr)[e] = T(h^{-1}\sigma)$. 

\item \label{pfsuper-C(H)V-measPfEx}
Prove \cref{pfsuper-C(H)V-meas} in the special case where $V = \real \sigma$ is the span of an $H$-invariant measurable section.
\hint{This is similar to \cref{pfsuper-C(H)V}, but use the fact that $H$ is ergodic on $G/\Gamma$.}

\item \label{pfsuper-C(H)V-measFullPfEx}
Prove \cref{pfsuper-C(H)V-meas} (without assuming $\dim V = 1$).
\hint{This is similar to \cref{pfsuper-C(H)VEx}.}

\item \label{ProveMargSuperHSimpleEx}
Prove \cref{MargSuperHSimple}.

\end{exercises}

\section{An \texorpdfstring{$A$}{A}-invariant section}

This section sketches the proof of the following result, which completes the proof of \cref{MargSuperHSimple} (under the assumption that $\Rrank G \ge 2$).

\begin{keyfact} \label{pfsuper-keyfact-meas}
For some~$n$, there is an embedding of~$H$ in $\SL(n,\real)$, such that%
\noprelistbreak
	\begin{enumerate}
	\item the associated representation $\varphi \colon \Gamma \to H \subseteq \SL(n,\real)$ is irreducible,
	and
	\item there exists a nonzero $A$-invariant $\sigma \in \Sectm(\bundle_\varphi)$.
	\end{enumerate}
\end{keyfact}

\Cref{pfsuper-sect<>equi} allows us to restate this as follows:

\begin{thmref}{pfsuper-keyfact-meas}
\begin{keyfact} \label{pfsuper-keyfact-equi}
For some embedding of~$H$ in $\SL(n,\real)$, 
\begin{enumerate}
\item $H$ acts irreducibly on~$\real^n$, 
and
\item there exists a\/ $\Gamma$-equivariant, measurable function\/ $\xi \colon G/A \to \real^n$ \textup(and $\xi$~is nonzero\/\textup).
\end{enumerate}
\end{keyfact}
\end{thmref}

In this form, the result is closely related to the following consequence of amenability (from \cref{AmenableChap}). For simplicity, it is stated only for the case $G = \SL(3,\real)$.

\begin{thmref}{SL3R/P->Meas(X)}
\begin{prop}[(Furstenberg)] \label{SL3R/P->Meas(X)'}
If
\noprelistbreak
	\begin{itemize}
	\item $G = \SL(3,\real)$,
	\item $ P
	 =  \begin{Smallbmatrix} \upast&\upast&\upast \\ &\upast&\upast \\ &&\upast \end{Smallbmatrix}
	 \subset G $,
	 and
	 \item $\Gamma$ acts continuously on a compact metric space~$X$,
 \end{itemize}
 then there is a Borel measurable map $\psi \colon G/P \to
\Prob(X)$, such that $\psi$ is essentially
$\Gamma$-equivariant.
 \end{prop}
 \end{thmref}

For convenience, let $W = \real^n$. There are 3 steps in the proof of \cref{pfsuper-keyfact-equi}:
	\begin{enumerate}
	\item (amenability) Letting $X$ be the projective space $\projective(W)$, which is compact, \cref{SL3R/P->Meas(X)'} provides a $\Gamma$-equivariant, measurable map $ \widehat \xi \colon G/P \to \Prob \bigl( \projective(W) \bigr)$.
	\item (proximality) The representation of~$\Gamma$ on~$W$ induces a representation of~$\Gamma$ on any exterior power $\bigwedge^k W$. By replacing $W$ with an appropriate subspace of such an exterior power, we may assume there is some $\gamma \in \Gamma$, such that $\gamma$ has a unique eigenvalue of maximal absolute value \csee{UniqMaxEigEx}. Therefore, the action of $\gamma$ on $\projective(W)$ is ``proximal'' \csee{ProxLem}. The theory of proximality (discussed in \cref{QuickProximalitySect}) now tells us that the $\Gamma$-equivariant random map~$\widehat\xi$ must actually be a well-defined map into $\projective(W)$ \csee{pfsuper-xi(x)ptmass}. 

\item (algebra trick) We have a $\Gamma$-equivariant map $\widehat\xi \colon G/P \to \projective(W)$. By the same argument, there is a $\Gamma$-equivariant map $\widehat\xi^* \colon G/P \to \projective(W^*)$, where $W^*$ is the dual of~$W$. Combining these yields a $\Gamma$-equivariant map 
	$$\overline\xi \colon G/P \times G/P \to \projective(W \otimes W^*) \iso \projective \bigl( \End(W) \bigr) .$$
We can lift $\overline\xi$ to a well-defined map 
	$$\xi \colon G/P \times G/P \to \End(W), $$
by specifying that $\trace \bigl( \xi(x) \bigr) = 1$ \csee{trace(xi)=1Ex}.
Since the action of $\Gamma$ on $\End(W)$ is by conjugation \csee{GammaOnEnd(W)ConjEx}
and the trace of conjugate matrices are equal, we see that $\xi$ is $\Gamma$-equivariant \csee{XiGammaEquiEx}.

Finally, note that there is a $G$-orbit in $G/P \times G/P$ whose complement is a set of measure~$0$, and the stabilizer of a point is (conjugate to) the group~$A$ of diagonal matrices \csee{GTransOnFlags,GTransOnP/PxG/P}. Therefore, after discarding a set of measure~$0$, we may identify $G/P \times G/P$ with $G/A$, so $\xi \colon G/A \to \End(W)$.
	\end{enumerate}

\begin{exercises}

\item \label{UniqMaxEigEx}
Let 
	\begin{itemize}
	\item $\gamma$ be a semisimple element of~$\Gamma$, such that some eigenvalue of~$\gamma$ is not of absolute value~$1$. 
	\item $\lambda_1,\ldots,\lambda_k$ be the eigenvalues of~$\gamma$ (with multiplicity) that have maximal absolute value.
	\item $W' = \bigwedge^k W$.
	\end{itemize}
Show that, in the representation of~$\Gamma$ on~$W'$, the element~$\gamma$ has a unique eigenvalue of maximal absolute value.

\item \label{trace(xi)=1Ex}
Let $\cover\xi \colon G/P \to W$ and $\cover\xi^* \colon G/P \to W^*$ be well-defined, measurable lifts of $\widehat\xi$ and $\widehat\xi^*$. 
	\begin{enumerate}
	\item Show, for a.e.\ $x,y \in G/P$, that $\cover\xi(x)$ is not in the kernel of the linear functional $\cover\xi^*(y)$.
	\item Show, for a.e.\ $x,y \in G/P$, that, under the natural identification of $W \otimes W^*$ with $\End(W)$, we have 
		$$\trace \bigl( \cover\xi(x) \otimes \cover\xi^*(y) \bigr) \neq 0 .$$
	\item Show $\overline\xi$ can be lifted to a well-defined measurable function $\xi \colon G/P \times G/P \to \End(W)$, such that $\trace \bigl( \xi(x,y) \bigr) = 1$, for a.e.\ $x,y \in G/P$.
	\end{enumerate}
\hint{$\Gamma$ acts irreducibly on~$W$, and ergodically on $G/P \times G/P$.}

\item \label{GammaOnEnd(W)ConjEx}
Show that the action of $\Gamma$ on $\End(W) \iso W \otimes W^*$ is given by conjugation:
	$ \overline\varphi(\gamma) T = \varphi(\gamma) \, T \, \varphi(\gamma)^{-1} $.

\item \label{XiGammaEquiEx}
Show that $\xi$ is $\Gamma$-equivariant.

\item \label{GTransOnFlags}
Recall that a \defit[flag in $\real^3$]{flag} in~$\real^3$ is a pair $(\ell, \Pi)$, where 
	\begin{itemize}
	\item $\ell$\, is a line through the origin (in other words, a $1$-dimensional linear subspace),
	and
	\item $\Pi$ is a plane through the origin (in other words, a $2$-dimensional linear subspace),
	such that
	\item $\ell \subset \Pi$.
	\end{itemize}
Show:
	\begin{enumerate}
	\item $\SL(3,\real)$ acts transitively on the set of all flags in~$\real^3$,
	and
	\item the stabilizer of any flag is conjugate to the subgroup~$P$ of
 \cref{SL3R/P->Meas(X)'}.
 	\end{enumerate}
Therefore, the set of flags can be identified with $G/P$.

\item \label{GTransOnP/PxG/P}
Two flags $(\ell_1,\Pi_1)$ and $(\ell_2,\Pi_2)$ are in \defit{general position} if 
	$$ \text{
	$\ell_1 \notin \Pi_2$, \  
	and \ 
	$\ell_2 \notin \Pi_1$
	} .$$
Letting $\mathcal{G}$ be the subset of $G/P \times G/P$ corresponding to the pairs of flags that are in general position, show:
	\begin{enumerate}
	\item \label{GTransOnP/PxG/P-trans}
	$\SL(3,\real)$ is transitive on~$\mathcal{G}$,
	\item \label{GTransOnP/PxG/P-stab}
	the stabilizer of any point in~$\mathcal{G}$ is conjugate to the group of diagonal matrices,
	and
	\item the complement of~$\mathcal{G}$ has measure zero in $G/P \times G/P$.
	\end{enumerate}
	\hint{For \pref{GTransOnP/PxG/P-trans} and \pref{GTransOnP/PxG/P-stab}, identify $\mathcal{G}$ with the set of triples $(\ell_1,\ell_2,\ell_3)$ of lines  that are in general position, by letting $\ell_3 = \Pi_1 \cap \Pi_2$.}
\end{exercises}

\section{A quick look at proximality} \label{QuickProximalitySect}

\begin{assump} \label{ProxAssump}
Assume 
	\begin{enumerate}
	\item $\Gamma \subset \SL(\ell,\real)$, 
	\item every finite-index subgroup of~$\Gamma$ is irreducible on~$\real^\ell$,
	and 
	\item there exists a semisimple element $\overline\gamma \in \Gamma$, such that $\overline\gamma$ has a unique eigenvalue~$\overline\lambda$ of maximal absolute value (and the eigenvalue is simple, which means the corresponding eigenspace is $1$-dimensional).
	\end{enumerate}
\end{assump}

\begin{notation} \ 
\begin{enumerate}
\item Let $\overline{v}$ be an eigenvector associated to the eigenvalue~$\overline\lambda$.
\item For convenience, let $W = \real^\ell$.
\end{enumerate}
\end{notation}

\begin{lem}[(Proximality)] \label{ProxLem}
The action of\/~$\Gamma$ on\/ $\projective(W)$ is {\upshape\defit{proximal}}. This means that, for every $[w_1],[w_2] \in \projective(W)$, there exists a sequence $\{\gamma_n\}$ in~$\Gamma$, such that $d\bigl( [\gamma_n(w_1)],[\gamma_n(w_2)] \bigr) \to 0$ as $n \to \infty$.
\end{lem}

\begin{proof}
Assume, to simplify the notation, that all of the eigenspaces of~$\overline\gamma$ are orthogonal to each other. 
Then, for any $w \in W \smallsetminus \overline{v}^\perp$, we have $\overline\gamma^n[w] \to [\overline{v}]$, as $n \to \infty$ \csee{PowerToEigEx}. Since the finite-index subgroups of~$\Gamma$ act irreducibly, there is some $\gamma \in \Gamma$, such that $\gamma(w_1), \gamma(w_2) \notin \overline{v}^\perp$ \csee{MoveBothOutOfOrthCompEx}. Therefore, 
	$$d\bigl( \overline\gamma^n\gamma ([w_1]), \overline\gamma^n\gamma ([w_2]) \bigr) 
	\to d( [\overline{v}],[\overline{v}]) 
	= 0 ,$$
as desired.
\end{proof}

In the above proof, it is easy to see that the convergence $\overline\gamma^n[w] \to [\overline{v}]$ is uniform on compact subsets of $W \smallsetminus \overline{v}^\perp$ \csee{PowerToEigUnifEx}. This leads to the following stronger assertion \csee{MeasureProxEx}:

\begin{prop}[(Measure proximality)] \label{MeasureProx}
Let $\mu$ be any probability measure on\/ $\projective(W)$. Then there is a sequence\/ $\{\gamma_n\}$ in\/~$\Gamma$, such that\/ $(\gamma_n)_*\mu$ converges to a delta-mass supported at a single point of\/ $\projective(W)$.
\end{prop}

It is obvious from \cref{MeasureProx} that there is no $\Gamma$-invariant probability measure on $\projective(W)$. However, it is easy to see that there does exist a probability measure that is invariant ``on average\zz,'' in the following sense \csee{ExistStatMeasEx}:

\begin{defn} \ 
\noprelistbreak
\begin{enumerate}
\item Fix a finite generating set~$S$ of~$\Gamma$, such that $S^{-1} = S$. A probability measure~$\mu$ on $\projective(W)$ is \defit{stationary} for~$S$ if 
	$$ \frac{1}{\#S} \sum_{\gamma \in S} \gamma_* \mu = \mu .$$
\item More generally, let $\nu$ be a probability measure on~$\Gamma$.
 A probability measure~$\mu$ on $\projective(W)$ is \defit[stationary]{$\nu$-stationary} if $\nu * \mu = \mu$. More concretely, this means
	$$  \sum_{\gamma \in \Gamma} \nu(\gamma) \, \gamma_* \mu = \mu .$$
\end{enumerate}
(Some authors call $\mu$ ``\index{harmonic!measure}{harmonic}\zz,'' rather than ``stationary\zz.'')
\end{defn}

\begin{rem}
A random walk on $\projective(W)$ can be defined as follows: Choose a sequence $\gamma_1,\gamma_2, \ldots$ of elements of~$\Gamma$, independently and with distribution~$\nu$. Also choose a random $x_0 \in \projective(W)$, with respect to some probability distribution~$\mu$ on $\projective(W)$. Then $x_n \in \projective(W)$ is defined  by 
	$$x_n = \gamma_1 \gamma_2 \cdots \gamma_n(x_0) ,$$
so $\{x_n\}$ is a random walk on $\projective(W)$.
A stationary measure represents a ``stationary state'' (or equilibrium distribution) for this random walk. Hence the terminology.  
\end{rem}

If the initial distribution~$\mu$ is stationary, then a basic result of probability (the ``\thmindex{Martingale Convergence}Martingale Convergence Theorem'') implies, for almost every sequence $\{\gamma_n\}$, that the resulting random walk $\{x_n\}$ has a limiting distribution; that is, 
	$$ \text{for a.e.\ $\{\gamma_n\}$, \ $(\gamma_1 \gamma_2 \cdots \gamma_n)_*\mu$ converges in 
$\Prob \bigl( \projective(W) \bigr)$}. $$
This theorem applies to stationary measures on any space, with no need for \cref{ProxAssump}. By using measure proximality, we will now show that the limiting distribution is almost always a point mass.

\begin{defn}
A closed, nonempty, $\Gamma$-invariant subset of $\projective(W)$ is \defit[minimal!closed, invariant set]{minimal} if it does not have any nonempty, proper, closed, $\Gamma$-invariant subsets. (Since $\projective(W)$ is compact, the finite-intersection property implies that every nonempty, closed, $\Gamma$-invariant subset of $\projective(W)$ contains a minimal set.)
\end{defn}

\begin{thm}[(\thmindex{mean proximality}Mean proximality)] \label{meanprox}
Assume 
	\begin{itemize}
	\item $\nu$ is a probability measure on~$\Gamma$, such that $\nu(\gamma) > 0$ for all $\gamma \in \Gamma$,
	\item $C$ is a minimal closed, $\Gamma$-invariant subset of $\projective(W)$,
	and
	\item $\mu$ is a $\nu$-stationary probability measure on~$C$. 
	\end{itemize}
Then, for a.e.\ $\{\gamma_n\} \in \Gamma^\infty$, there exists $c \in \projective(W)$, such that 
	$$ \text{$(\gamma_1 \gamma_2 \cdots \gamma_n)_*(\mu) \to \delta_c$ as $n \to \infty$.} $$
\end{thm}

\begin{proof}
It was mentioned above that the Martingale Convergence Theorem implies $(\gamma_1 \gamma_2 \cdots \gamma_n)_*(\mu)$ has a limit (almost surely), so it suffices to show there is (almost surely) a subsequence that converges to a measure of the form~$\delta_c$.

\Cref{MeasureProx} provides a sequence $\{g_k\}$ of elements of~$\Gamma$, such that $(g_k)_*\mu \to \delta_{c_0}$, for some $c_0 \in C$. To extend this conclusion to a.e.\ sequence $\{\gamma_n\}$, we use equicontinuity: we may write $\Gamma$ is the union of finitely many sets $E_1,\ldots,E_r$, such that each $E_i$ is equicontinuous on some nonempty open subset~$U_i$ of~$C$ \csee{EquiOnProj}.

The minimality of~$C$ implies $\Gamma U_i = C$ for every~$i$. Then, by compactness, there is a finite subset $F = \{f_1,\ldots,f_s\}$ of~$\Gamma$, such that $F U_i = C$ for each~$i$. Since $\nu(\gamma) > 0$ for every $\gamma \in \Gamma$, there is (almost surely) a subsequence $\{\gamma_{n_k}\}$ of~$\{\gamma_n\}$, such that, for every~$k$, we have
	$$ \text{$\gamma_{n_k+1} \gamma_{n_k+2}\cdots \gamma_{n_k+j} = f_j^{-1} \, g_k$ \ for $1 \le j \le s$} .$$
By passing to a subsequence, we may assume there is some~$i$, such that 
	$$ \text{$\gamma_1 \gamma_2 \cdots \gamma_{n_k} \in E_i$, for all~$k$} .$$
To simplify the notation, let us assume $i = 1$.

Since $F U_1 = C$, we may write $c_0 = f_j u$, for some $f_j \in F$ and $u \in U_1$.
Then
	$$ (\gamma_{n_k+1} \gamma_{n_k+2}\cdots \gamma_{n_k+j})_* \nu = (f_j^{-1} \, g_k)_* \nu \to  (f_j^{-1})_* \delta_{c_0} = \delta_{f_j^{-1} c_0} = \delta_u.$$
By passing to a subsequence, we may assume $(\gamma_1 \gamma_2 \cdots \gamma_{n_k})u$ converges to some $c \in C$. Then, since $\gamma_1 \gamma_2 \cdots \gamma_{n_k} \in E_1$, 
and $E_1$ is equicontinuous on $U_1$, this implies 
	\begin{align*}
	(\gamma_1 \gamma_2 \cdots \gamma_{n_k+j})_* \nu 
	= (\gamma_1 \gamma_2 \cdots \gamma_{n_k})_*  
	\bigl( (\gamma_{n_k+1} \ \cdots \gamma_{n_k+j})_* \nu \bigr)
	\to \delta_c 
	. & \qedhere \end{align*}
\end{proof}

In order to apply this theorem, we need a technical result, whose proof we omit:

\begin{lem} \label{pfsuper-LebesgueStat}
There exist:
	\begin{itemize}
	\item a probability measure~$\nu$ on\/~$\Gamma$,
	and
	\item a $\nu$-stationary probability measure~$\mu$ on $G/P$,
	\end{itemize}
such that 
	\begin{enumerate}
	\item  the support of~$\nu$ generates\/~$\Gamma$,
	and
	\item $\mu$ is in the class of Lebesgue measure. \textup(That is, $\mu$ has exactly the same sets of measure\/~$0$ as Lebesgue measure does.\textup)
	\end{enumerate}
\end{lem}

Also note that if $C$ is any nonempty, closed, $\Gamma$-invariant subset of $\projective(W)$, then $\Prob(C)$ is a nonempty, compact, convex $\Gamma$-space, so Furstenberg's Lemma \pref{G/amen->Meas(X)} provides a $\Gamma$-equivariant map $\overline{\xi} \colon G/P \to \Prob(C)$. This observation allows us to replace $\projective(W)$ with a minimal subset.

We can now fill in the missing part of the proof of \cref{pfsuper-keyfact-equi}:

\begin{cor} \label{pfsuper-xi(x)ptmass}
Suppose 
\noprelistbreak
	\begin{itemize}
	\item $C$ is a minimal closed, $\Gamma$-invariant subset of\/ $\projective(W)$,
	and
	\item $\overline{\xi} \colon G/P \to \Prob(C)$ is $\Gamma$-equivariant.
	\end{itemize}
Then $\overline{\xi}(x)$ is a point mass, for a.e.\ $x \in G/P$.

Hence, there exists $\hat\xi \colon G/P \to \projective(W)$, such that $\overline{\xi}(x) = \delta_{\hat\xi(x)}$, for a.e.\ $x \in G/P$.
\end{cor}

\begin{proof}
Let 
	\begin{itemize}
	\item $\delta_{\projective(W)} = \{\, \delta_{x} \mid x \in \projective(W)\,\}$
		be the set of all point masses in the space $\Prob \bigl( \projective(W) \bigr)$,
	and
	\item $\mu$ be a $\nu$-stationary probability measure on $G/P$ that is in the class of Lebesgue measure \csee{pfsuper-LebesgueStat}.
	\end{itemize}
We wish to show $\overline{\xi}(x) \in \delta_{\projective(W)}$, for a.e.\ $x \in G/P$. In other words, we wish to show that $\overline{\xi}_*(\mu)$ is supported on $\delta_{\projective(W)}$.

Note that:
	\begin{itemize}
	\item $\delta_{\projective(W)}$ is a closed, $\Gamma$-invariant subset of $\Prob \bigl( \projective(W) \bigr)$,
	and
	\item because $\overline{\xi}$ is $\Gamma$-equivariant, we know that $\overline{\xi}_*(\mu)$ is a $\nu$-stationary probability measure on $\Prob \bigl( \projective(W) \bigr)$. 
	\end{itemize}

Roughly speaking, the idea of the proof is that almost every trajectory of the random walk on $\Prob \bigl( \projective(W) \bigr)$ converges to a point in $\delta_{\projective(W)}$ \see{meanprox}. On the other hand, being stationary, $\overline{\xi}_*(\mu)$ is invariant under the random walk. Therefore, we conclude that $\overline{\xi}_*(\mu)$ is supported on $\delta_{\projective(W)}$, as desired. 

We now make this rigorous. Let
	$$ \mu_{\projective(W)} = \int_{G/P} \overline{\xi}(x) \, d\mu(x) , $$
so $\mu_{\projective(W)}$ is a stationary probability measure on $\projective(W)$. By mean proximality \pref{meanprox}, we know, for a.e.\ $(\gamma_1,\gamma_2,\ldots) \in \Gamma^\infty$, that
	$$ d \bigl( \, (\gamma_1 \gamma_2 \cdots \gamma_n)_*(\mu_{\projective(W)}), \ \delta_{\projective(W)} \, \bigr) \ \stackrel{n \to \infty}{\longrightarrow} \ 0 .$$
For any $\epsilon > 0$, this implies, by using the definition of~$\mu_{\projective(W)}$, that
	$$ \mu \Bigl( \bigset{ x \in G/P }{ d \bigl( \gamma_1 \gamma_2 \cdots \gamma_n \bigl( \overline{\xi}(x) \bigr) 
	, \delta_{\projective(W)} \bigr)  > \epsilon } \Bigr) \ \stackrel{n \to \infty}{\longrightarrow} \ 0 .$$
Since $\overline{\xi}$ is $\Gamma$-equivariant, we may
	$$ \text{replace \ $\gamma_1 \gamma_2 \cdots \gamma_n \bigl( \overline{\xi}(x) \bigr)$
	\ with \ $\overline{\xi}( \gamma_1 \gamma_2 \cdots \gamma_n  x)$} .$$
Then, since the measure $\mu$ on $G/P$ is stationary, we can delete $\gamma_1 \gamma_2 \cdots \gamma_n$, and conclude that
	\begin{align} \label{DeleteGamma}
	 \mu
	\bigset{ x \in G/P }{ 
	d \bigl( \overline{\xi}(x) 
	, \delta_{\projective(W)} \bigr)  > \epsilon } \ \stackrel{n \to \infty}{\longrightarrow} \ 0 
	\end{align}
\csee{DeleteGammaEx}.
Since the left-hand side does not depend on~$n$, but tends to~$0$ as $n \to \infty$, it must be~$0$. Since $\epsilon > 0$ is arbitrary, we conclude that $\overline{\xi}(x)  \in \delta_{\projective(W)}$ for a.e.~$x$, as desired.
\end{proof}

\begin{exercises}

\item \label{PowerToEigEx}
In the notation of \cref{ProxLem}, show, for every $w \in W \smallsetminus \overline{v}^\perp$, that $\overline\gamma^n[w] \to [\overline{v}]$, as $n \to \infty$.

\item \label{MoveBothOutOfOrthCompEx}
Show, for any nonzero $w_1,w_2 \in W$, that there exists $\gamma \in \Gamma$, such that neither $\gamma w_1$ nor $\gamma w_2$ is orthogonal to~$\overline{v}$.
\hint{Let $H$ be the Zariski closure of~$\Gamma$ in $\SL(\ell,\real)$, and assume, by passing to a finite-index subgroup, that $H$ is connected. Then $W_i = \{\, h \in H \mid h w_i \in \overline{v}^\perp \,\}$ is a proper, Zariski-closed subset. Since $\Gamma$ is Zariski dense in~$H$, it must intersection the complement of $W_1 \cup W_2$.}

\item \label{PowerToEigUnifEx}
Show that the convergence in \cref{PowerToEigEx} is uniform on compact subsets of $W \smallsetminus \overline{v}^\perp$.

\item \label{MeasureProxEx}
Prove \cref{MeasureProx}.
\hint{Show $\max_{w \in \projective(W) \\ \nu \in \overline{\Gamma \mu}} \nu(w) = 1$.}

\item \label{ExistStatMeasEx}
Show there exists a stationary probability measure on $\projective(W)$.
\hint{Kakutani-Markov Fixed-Point Theorem \cf{CyclicAmen}.}

\item \label{EquiOnProj}
Let $C$ be a subset of $\projective(\real^n)$, and assume that $C$ is not contained in any $(n-1)$-dimensional hyperplane.
Prove that $\GL(n,\real)$ is the union of finitely many sets $E_1,\ldots,E_r$, such that each $E_i$ is equicontinuous on some nonempty open subset~$U_i$ of~$C$.
\hint{Each matrix $T \in \Mat_{n\times n}(\real)$ induces a well-defined, continuous function $\overline{T} \colon \bigl( \projective(\real^n) \smallsetminus \projective(\ker T) \bigr) \to \projective(\real^n)$. 
If $B_T$ is a small ball around  $\overline{T}$ in $\projective \bigl( \Mat_{n\times n}(\real) \bigr)$, then $B_T$ is equicontinuous on an open set. A compact set can be covered by finitely many balls.}

\item \label{DeleteGammaEx}
Establish \pref{DeleteGamma}.
\hint{Since $\mu$ is stationary, the map
	$$ \Gamma^\infty \times G/P \to G/P \colon \bigl( (\gamma_1,\gamma_2,\ldots), x \bigr)
	\mapsto \gamma_1 \gamma_2 \cdots \gamma_n x $$
is measure preserving.}

\item
Show that if $\projective(W)$ is minimal, then the $\Gamma$-equivariant measurable map $\xi \colon G/P \to \projective(W)$ is unique (a.e.).
\hint{If $\psi$ is another $\Gamma$-equivariant map, then define $\overline\xi \colon G/P \to \Prob \bigl( \projective(W) \bigr)$ by $\overline\xi(x) = \frac{1}{2} ( \delta_{\xi(x)} + \delta_{\psi(x)})$.}

\end{exercises}

\section{Groups of real rank one} \label{SuperRank1Sect}

The Margulis Superrigidity Theorem \pref{MargSuperC} was proved for groups of real rank at least two in \cref{SuperPfSect}. Suppose, now, that $\Rrank G = 1$ (and $G$ has no compact factors). The classification of simple Lie groups tells us that $G$ is isogenous to the isometry group of either:
	\begin{itemize}
	\item real hyperbolic space $\hyperbolic^n$,
	\item complex hyperbolic space $\complex\hyperbolic^n$,
	\item quaternionic hyperbolic space $\quaternion\hyperbolic^n$,
	or
	\item the Cayley hyperbolic plane $\mathbb{O}\hyperbolic^2$ (where $\mathbb{O}$ is the ring of ``Cayley numbers'' or ``octonions'')
	\end{itemize}
 \ccf{rank1simple}.
Assumption~\fullref{MargSuperC}{notSOSU} rules out $\hyperbolic^n$ and $\complex\hyperbolic^n$, so, from the connection of superrigidity with totally geodesic embeddings \ccf{TotGeodSect}, the following result completes the proof:

\begin{thm}
Assume
	\begin{itemize}
	\item $X = \quaternion\hyperbolic^n$ or\/ $\mathbb{O}\hyperbolic^2$,
	\item $\Gamma$ is a torsion-free, discrete group of isometries of~$X$, such that\/ $\Gamma \backslash X$ has finite volume,
	\item $X'$ is an irreducible symmetric space of noncompact type,
	and
	\item $\varphi \colon \Gamma \to \Isom(X')^\circ$ is  a homomorphism whose image is Zariski dense.
	\end{itemize}
Then there is a map $f \colon X \to X'$, such that 
	\begin{enumerate}
	\item $f(X)$ is totally geodesic,
	and
	\item $f$ is $\varphi$-equivariant, which means $f(\gamma x) = \varphi(\gamma) \cdot f(x)$.
	\end{enumerate}
\end{thm}

\begin{proof}[Brief outline of proof]
Choose a (nice) fundamental domain~$\fund$ for the action of~$\Gamma$ on~$X$. For any $\varphi$-equivariant map $f \colon X \to X'$, define the \defit{energy} of~$f$ to be the $L^2$-norm of the derivative of~$f$ over~$\fund$. Since $f$~is $\varphi$-equivariant, and the groups $\Gamma$ and $\varphi(\Gamma)$ act by isometries, this is independent of the choice of the fundamental domain~$\fund$.

It can be shown that this energy functional attains its minimum at some function~$f$. The minimality implies that $f$ is harmonic. Then, by using the geometry of~$X$ and the negative curvature of~$X'$, it can be shown that $f$ must be totally geodesic.
\end{proof}

\begin{notes}

This chapter is largely based on \cite[Chaps.~6 and~7]{MargulisBook}. (However, we usually replace the assumption that $\Rrank G \ge 2$ with the weaker assumption that $G$ is not $\SO(1,m) \times K$ or $\SU(1,m) \times K$. (See \cite[Thm.~5.1.2, p.~86]{ZimmerBook} for a different exposition that proves version \pref{MargSuperG'} for $\Rrank G \ge 2$.)
In particular:
\begin{itemize}
\item For $\Rrank G \ge 2$, our statement of the Margulis Superrigidity Theorem \pref{MargSuperG'} is a special case of \cite[Thm.~7.5.6, p.~228]{MargulisBook}.

\item For $\Rrank G \ge 2$, \cref{MargImgSS} is stated in \cite[Thm.~9.6.15(i)(a), p.~332]{MargulisBook}.

\item For $\Rrank G \ge 2$, \cref{MostowRigidity} is stated in \cite[Thm.~7.7.5, p.~254]{MargulisBook}. (See \cite[Thm.~B]{PrasadMostowRig} for the general case, which does not follow from superrigidity.)

\item \Cref{Extend<>FD} is a version of \cite[Prop.~4.6, p.~222]{MargulisBook}
\item \Cref{HighRankGenLi} is a version of \cite[Lem.~7.5.5, p.~227]{MargulisBook}.
\item \Cref{pfsuper-C(H)V-meas} is \cite[Prop.~7.3.6, p.~219]{MargulisBook}.

\item \Cref{pfsuper-keyfact-equi} is adapted from \cite[Thm.~6.4.3(b)2, p.~209]{MargulisBook}.

\item \Cref{meanprox} is based on \cite[Prop.~6.2.13, pp.~202--203]{MargulisBook}.
\item \Cref{pfsuper-LebesgueStat} is taken from \cite[Prop.~6.4.2, p.~209]{MargulisBook}.
\item \Cref{pfsuper-xi(x)ptmass} is based on \cite[Prop.~6.2.9, p.~200]{MargulisBook}.
\item \Cref{EquiOnProj} is \cite[Lem.~6.3.2, p.~203]{MargulisBook}.
\end{itemize}

Long before the general theorem of Margulis for groups of real rank $\ge 2$, it was proved by Bass, Milnor, and Serre \cite[Thm.~61.2]{BassMilnorSerre-CSP} that the \thmindex{Congruence Subgroup Property}Congruence Subgroup Property implies $\SL(n,\integer)$ is superrigid in $\SL(n,\real)$. 

\defit[Geometric superrigidity]{\normalfont ``Geometric superrigidity''} is the study of differential geometric versions of the Margulis Superrigidity Theorem, such as \cref{TotGeodProp}. (See, for example, \cite{MokSiuYeung-GeomSuper}.)

Details of the derivation of arithmeticity from superrigidity (\cref{MargArithPf}) appear in \cite[Chap.~9]{MargulisBook} and \cite[\S6.1]{ZimmerBook}. 

Proofs of the Commensurability Criterion \pref{CommCriterion} and Commensurator Superrigidity \pref{CommSuper} can be found in \cite{ACampoBurger-reseaux}, \cite[\S9.2.11, pp.~305\emph{ff}, and Thm.~7.5.4, pp.~226--227]{MargulisBook}, and \cite[\S6.2]{ZimmerBook}.

Much of the material in \cref{QuickProximalitySect} is due to Furstenberg \cite{Furstenberg-BdThyStochProc}.

The superrigidity of lattices in the isometry groups of $\quaternion\hyperbolic^n$ and $\mathbb{O}\hyperbolic^2$ \csee{SuperRank1Sect} was proved by Corlette \cite{Corlette-superrig}. The $p$-adic version \pref{padicSuper} for these groups was proved by Gromov and Schoen \cite{GromovSchoen-padicsuper}.

\end{notes}

 %!TEX root = IntroArithGrps.tex

\mychapter{Normal Subgroups of\texorpdfstring{~$\Gamma$}{ Γ}}
\label{NormalSubgroupChap}

\prereqs{amenability (Furstenberg's Lemma \pref{G/amen->Meas(X)}) and Kazhdan's Property~$(T)$ (\cref{KazhdanTChap}). 
\emph{Also used:} the $\sigma$-algebra of Borel sets modulo sets of measure~$0$ (\cref{ErgDecompSect}) and manifolds of negative curvature.}

This chapter presents a contrast between the lattices in groups of real rank~$1$ and those of higher real rank:
\begin{itemize}
\item If $\Rrank G = 1$, then $\Gamma$ has many, many normal subgroups, so $\Gamma$ is very far from being simple.
\item If $\Rrank G > 1$ (and $\Gamma$ is irreducible), then $\Gamma$ is simple modulo finite groups. More precisely, if $N$ is any normal subgroup of~$\Gamma$, then either $N$ is finite, or $\Gamma/N$ is finite.
\end{itemize}

\section{\texorpdfstring{Normal subgroups in lattices of real rank $\ge 2$}{Normal subgroups in lattices of real rank at least two}}

\begin{thm}[(\thmindex{Margulis!Normal Subgroups}Margulis Normal Subgroups Theorem)] \label{MargNormalSubgrpsThm}
	\thmindex{Margulis!Normal Subgroups}
Assume 
\noprelistbreak
	\begin{itemize}
	\item $\Rrank G \ge 2$,
	\item $\Gamma$ is an irreducible lattice in~$G$,
	and
	\item $N$ is a normal subgroup of\/~$\Gamma$.
	\end{itemize}
Then either $N$ is finite, or\/ $\Gamma/N$ is finite.
\end{thm}

\begin{eg}
Every lattice in $\SL(3,\real)$ is simple, modulo finite groups.
In particular, this is true of $\SL(3,\integer)$.
\end{eg}

\begin{rems} \label{MargNormSubgrpRems} \ 
\noprelistbreak
\begin{enumerate}
\item The hypotheses on $G$ and~$\Gamma$ are essential:
	\noprelistbreak
	\begin{enumerate}
	\item If $\Rrank G = 1$, then every lattice in~$G$ has an infinite normal subgroup of infinite index \csee{Rrank1->GammaNotAlmSimple}. 
	\item \label{MargNormSubgrpRems-hyp-red}
	If $\Gamma$ is reducible (and $G$ has no compact factors), then $\Gamma$ has an infinite normal subgroup of infinite index \csee{MargNormSubgrpRems-hyp-redEx}.
	\end{enumerate}
\item The finite normal subgroups of~$\Gamma$ are easy to understand (if $\Gamma$ is irreducible): the Borel Density Theorem implies that they are the subgroups of the finite abelian group $\Gamma \cap Z(G)$ \csee{NinZ(G)}.
\item \label{MargNormSubgrpRems-finind}
If $\Gamma$ is infinite, then $\Gamma$ has infinitely many normal subgroups of finite index \csee{MargNormSubgrpRems-finindEx}, so $\Gamma$ is \emph{not} simple.

\item \label{MargNormSubgrpRems-CSP}
In most cases, the subgroups of finite index are described by the ``\thmindex{Congruence Subgroup Property}Congruence Subgroup Property\zz.'' For example, if $\Gamma = \SL(3,\integer)$, then the principal congruence subgroups are obvious subgroups of finite index \csee{CSGfinite}. More generally, any subgroup of~$\Gamma$ that contains a principal congruence subgroup obviously has finite index. The Congruence Subgroup Property is the assertion that every finite-index subgroup is one of these obvious ones. It is true for $\SL(n,\integer)$, whenever $n \ge 3$, and a similar (but slightly weaker) statement is conjectured to be true whenever $\Rrank G \ge 2$ and $\Gamma$~is irreducible.
\end{enumerate}
\end{rems}

The remainder of this section presents the main ideas in the proof of \cref{MargNormalSubgrpsThm}. In a nutshell, we will show that if $N$ is an infinite, normal subgroup of~$\Gamma$, then
	\begin{enumerate}
	\item $\Gamma/N$ has Kazhdan's property $(T)$,
	and
	\item $\Gamma/N$ is amenable.
	\end{enumerate}
This implies that $\Gamma/N$ is finite \csee{T+amen->finite}.

In most cases, it is easy to see that $\Gamma/N$ has Kazhdan's property (because $\Gamma$ has the property), so the main problem is to show that $\Gamma/N$ is amenable. This amenability follows easily from an ergodic-theoretic result that we will now describe. 

\begin{assump}
To minimize the amount of Lie theory needed, let us assume 
	$$G = \SL(3,\real) .$$
\end{assump}

\begin{notation} \label{MargNormalP=upper}
Let 
	$$P = \begin{bmatrix}  *&& \\ *&*& \\ *&*&* \end{bmatrix} \subset \SL(3,\real) = G .$$
Hence, $P$ is a (minimal) parabolic subgroup of~$G$.
\end{notation}

Note that if $Q$ is any closed subgroup of~$G$ that contains~$P$, then the natural map $G/P \to G/Q$ is $G$-equivariant, so we may say that $G/Q$ is a $G$-equivariant quotient of $G/P$. Conversely, it is easy to see that spaces of the form $G/Q$ are the only $G$-equivariant quotients of $G/P$. In fact, these are the only quotients even if we only assume that quotient map is equivariant \emph{almost} everywhere \csee{GQuot(G/H)}.

Furthermore, since $\Gamma$ is a subgroup of~$G$, it is obvious that every $G$-equivariant map is $\Gamma$-equivariant. Conversely, the following surprising result shows that every $\Gamma$-equivariant quotient of $G/P$ is $G$-equivariant (up to a set of measure~$0$):

\begin{thm}[(Margulis)] \label{blackbox} Suppose
\noprelistbreak
\begin{itemize}
\item $\Rrank G \ge 2$,
\item $P$ is a minimal parabolic subgroup of~$G$,
\item $\Gamma$ is irreducible,
\item $\Gamma$ acts by homeomorphisms on a compact, metrizable space~$Z$,
and
\item  $\psi \colon G/P \to Z$ is essentially\/ $\Gamma$-equivariant\/ \textup(and measurable\/\textup).
\end{itemize}
Then the action of\/~$\Gamma$ on~$Z$ is measurably isomorphic to the natural action of\/~$\Gamma$ on~$G/Q$\/ \textup(a.e.\textup), for some closed subgroup~$Q$ of~$G$ that contains~$P$.
\end{thm}

\begin{rem} \label{BlackboxRem} \ 
\noprelistbreak
\begin{enumerate}
\item \label{BlackboxRem-meas}
Perhaps we should clarify the choice of measures in the statement of \cref{blackbox}.
(A measure class on $G/P$ is implicit in the assumption that $\psi$ is \emph{essentially} $\Gamma$-equivariant. Measure classes on~$Z$ and~$G/Q$ are implicit in the ``(a.e.)'' in the conclusion of the theorem.)
	\begin{enumerate}
	\item Because $G/P$ and $G/Q$ are $C^\infty$ manifolds, Lebesgue measure supplies a measure class on each of these spaces. The Lebesgue class is invariant under all diffeomorphisms, so, in particular, it  is $G$-invariant.
	\item There is a unique measure class on~$Z$ for which $\psi$ is measure-class preserving \csee{UniqMeasClassPresEx}.
	\end{enumerate}
\item The proof of \cref{blackbox} will be presented in \cref{QuotG/PPfSect}. It may be skipped on a first reading.
%\item Although \cref{blackbox} is stated only for $G = \SL(3,\real)$, it is valid (and has essentially the same proof) whenever $\Rrank G \ge 2$ and $\Gamma$~is irreducible.
\end{enumerate}
\end{rem}

\begin{proof}[\textbf{\upshape Proof of \cref{MargNormalSubgrpsThm}}]
Let $N$ be a normal subgroup of~$\Gamma$, and assume $N$ is infinite. We wish to show $\Gamma/N$ is finite.  
Let us assume, for simplicity, that $\Gamma$ has Kazhdan's Property $(T)$. (For example, this is true if $G = \SL(3,\real)$, or, more generally, if $G$ is simple \csee{GammaHasT}.)
Then $\Gamma/N$ also has Kazhdan's Property $(T)$ \csee{KazhdanEasy}, so it suffices to show that $\Gamma/N$ is amenable \csee{T+amen->finite}.

Suppose $\Gamma/N$ acts by homeomorphisms on a compact, metrizable space~$X$.
In order to show that $\Gamma/N$ is amenable, it suffices to find an invariant probability measure on~$X$ \fullcsee{AmenEquiv}{InvMeas}. In other words, we wish to show that $\Gamma$ has a fixed point in $\Prob(X)$.
\begin{itemize}
\item Because $P$ is amenable, there is an (essentially) $\Gamma$-equivariant measurable map $\psi \colon G/P \to \Prob(X)$ \csee{SL3R/P->Meas(X)}.
\item From \cref{blackbox}, we know there is a closed subgroup~$Q$ of~$G$, such that the action of~$\Gamma$ on $\Prob(X)$ is measurably isomorphic (a.e.)\  to the natural action of~$\Gamma$ on~$G/Q$.
 \end{itemize}
 Since $N$ acts trivially on~$X$, we know it acts trivially on $\Prob(X) \iso G/Q$.
Hence, the kernel of the $G$-action on $G/Q$ is infinite \csee{NormalPf-NTrivEx}.
However, $G$ is simple (modulo its finite center), 
 so this implies that the action of~$G$ on $G/Q$ is trivial \csee{NormalPf-GTrivEx}. 
 (It follows that $G/Q$ is a single point, so $Q = G$, but we do not need quite such a strong conclusion.) Since $\Gamma \subseteq G$, then the action of $\Gamma$ on~$G/Q$ is trivial. In other words, every point in $G/Q$ is fixed by~$\Gamma$.
 Since $G/Q \iso \Prob(X)$ (a.e.), we conclude that almost every point in $\Prob(X)$ is fixed by~$\Gamma$; therefore, $\Gamma$ has a fixed point in $\Prob(X)$, as desired.
 \end{proof}

\begin{rem} \label{AE1Pt}
The proof of \cref{MargNormalSubgrpsThm} concludes that ``almost every point in $\Prob(X)$ is fixed by~$\Gamma$\zz,'' so it may seem that the proof provides not just a single $\Gamma$-invariant measure, but many of them. This is not the case: The proof implies that $\psi$ is essentially constant \csee{AE1PtEx}. This means that the $\Gamma$-invariant measure class $[\psi_* \mu]$ is supported on a single point of $\Prob(X)$, so ``a.e\zz.'' means only one point.
\end{rem}

\begin{exercises}

\item \label{GammaHasFiniteAbelianization}
Assume 
	\begin{itemize}
	\item $G$ is not isogenous to $\SO(1,n)$ or $\SU(1,n)$, for any~$n$,
	\item $\Gamma$ is irreducible,
	and
	\item $G$ has no compact factors.
	\end{itemize}
In many cases, Kazhdan's property $(T)$ implies that the abelianization $\Gamma / [\Gamma,\Gamma]$ of~$\Gamma$ is finite \fullcsee{KazhdanlatticeCor}{noabel}. Use \cref{MargNormalSubgrpsThm} to prove this in the remaining cases. (We saw a different proof of this in \cref{GNoAbelianization}.)

\item \label{MargNormSubgrpRems-hyp-redEx}
Verify \fullcref{MargNormSubgrpRems}{hyp-red}.
\hint{\Cref{prodirredlatt}.}

\item Suppose $\Gamma$ is a lattice in $\SL(3,\real)$. Show that $\Gamma$ has no nontrivial, finite, normal subgroups.

\item Suppose $\Gamma$ is an irreducible lattice in~$G$. Show that $\Gamma$ has only finitely many finite, normal subgroups.

\item \label{MargNormSubgrpRems-finindEx}
Show that if $\Gamma$ is infinite, then it has infinitely many normal subgroups of finite index.
\hint{\Cref{GammaResidFinite}.}

\item \label{GQuot(G/H)}
Suppose 
\noprelistbreak
	\begin{itemize}
	\item $H$ is a closed subgroup of~$G$, 
	\item $G$ acts continuously on a metrizable space~$Z$,
	and
	\item $\psi \colon G/H \to Z$ is essentially $G$-equivariant (and measurable).
	\end{itemize}
Show the action of~$G$ on~$Z$ is measurably isomorphic to the action of~$G$ on $G/Q$ (a.e.), for some closed subgroup~$Q$ of~$G$ that contains~$H$. More precisely, show there is a measurable $\phi \colon Z \to G/Q$, such that:%
\noprelistbreak
	\begin{enumerate}
	\item $\phi$ is measure-class preserving (i.e., a subset $A$ of~$G/Q$ has measure~$0$ if and only if its inverse image $\phi^{-1}(A)$ has measure~$0$),
	\item $\phi$ is one-to-one (a.e.) (i.e., $\phi$ is one-to-one on a conull subset of~$Z$),
	and
	\item $\phi$ is essentially $G$-equivariant.
	\end{enumerate}
\hint{See \fullcref{BlackboxRem}{meas} for an explanation of the measure classes to be used on $G/H$, $G/Q$, and~$Z$. For each $g \in G$, the set $\{\, x \in G/H \mid \psi(gx) = g \cdot \psi(x) \,\}$ is conull. By Fubini's Theorem, there is some $x_0 \in G/H$, such that $\psi(gx_0) = g \cdot \psi(x_0)$ for a.e.~$g$. Show the $G$-orbit of $\psi(x_0)$ is conull in~$Z$, and let $Q = \Stab_G \bigl( \psi(x_0) \bigr)$.}

\item \label{UniqMeasClassPresEx}
Suppose 
\noprelistbreak
	\begin{itemize}
	\item $\psi \colon Y \to Z$ is measurable,
	and
	\item $\mu_1$ and~$\mu_2$ are measures on~$Y$ that are in the same measure class.
	\end{itemize}
Show:
\noprelistbreak
	\begin{enumerate}
	\item The measures $\psi_*(\mu_1)$ and $\psi_*(\mu_2)$ on~$Z$ are in the same measure class.
	\item For any measure class on~$Y$, there is a unique measure class on~$Z$ for which $\psi$ is measure-class preserving.
	\end{enumerate}

\item In the setting of \cref{blackbox}, show that $\psi$ is essentially onto. That is, the image $\psi(G/P)$ is a conull subset of~$Z$.
\hint{By choice of the measure class on~$Z$, we know that $\psi$ is measure-class preserving.}

\item Let $G = \SL(3,\real)$ and $\Gamma = \SL(3,\integer)$.
Show that the natural action of~$\Gamma$ on $\real^3 / \integer^3 = \torus^3$ is a $\Gamma$-equivariant quotient of the action on~$\real^3$, but is not a $G$-equivariant quotient.

\item \label{NormalPf-NTrivEx}
In the proof of \cref{MargNormalSubgrpsThm}, we know that $\Prob(X) \iso G/Q$ (a.e.), so each element of~$N$ fixes a.e.\ point in $G/Q$. Show that $N$ acts trivially on $G/Q$ (everywhere, not only a.e.).
\hint{The action of $N$ is continuous.}

\item \label{NormalPf-GTrivEx}
In the notation of the proof of \cref{MargNormalSubgrpsThm}, show that the action of~$G$ on $G/Q$ is trivial.
\hint{Show that the kernel of the action of~$G$ on $G/Q$ is closed. You may assume, without proof, that $G$ is an almost simple Lie group. This means that every proper, closed, normal subgroup of~$G$ is finite.}

\item  \label{AE1PtEx}
In the setting of the proof of \cref{MargNormalSubgrpsThm}, show that $\psi$ is constant (a.e.).
\hint{The proof shows that a.e.\ point in the image of~$\psi$ is fixed by~$G$. Because $\psi$ is $G$-equivariant, and $G$ is transitive on $G/P$, this implies that $\psi$ is constant (a.e.).}
\end{exercises}

\section{Normal subgroups in lattices of rank one}

\Cref{MargNormalSubgrpsThm} assumes $\Rrank G \ge 2$. The following result shows that this condition is necessary:

\begin{thm} \label{Rrank1->GammaNotAlmSimple}
If\/ $\Rrank G = 1$, then\/ $\Gamma$ has a normal subgroup~$N$, such that neither $N$ nor\/ $\Gamma/N$ is finite.
\end{thm}

\begin{proof}[Proof \normalfont (assumes familiarity with manifolds of negative curvature)]
For simplicity, assume:
\noprelistbreak
	\begin{itemize}
	\item $\Gamma$ is torsion free, so it is the fundamental group of the locally symmetric space $M = \Gamma \backslash G / K$ (where $K$ is a maximal compact subgroup of~$G$).
	\item $M$ is compact.
	\item The locally symmetric metric on~$M$ has been normalized to have sectional curvature $\le -1$.
	\item The injectivity radius of~$M$ is $\ge 2$.
	\item There are closed geodesics $\gamma$ and $\lambda$ in~$M$, such that
		 $\length(\lambda) > 2 \pi$
		and
		 $\dist(\gamma,\lambda) > 2$.
	\end{itemize}
The geodesics $\gamma$ and~$\lambda$ represent (conjugacy classes of) nontrivial elements $\widehat\gamma$ and $\widehat\lambda$ of the fundamental group~$\Gamma$ of~$M$. Let $N$ be the smallest normal subgroup of~$\Gamma$ that contains~$\widehat\lambda$.

It suffices to show that $\widehat\gamma^n$ is nontrivial in $\Gamma/N$, for every $n \in \integer^+$ \csee{EnoughGammaNontrivEx}.
Construct a CW complex $\overline{M}$ by gluing the boundary of a $2$-disk~$D_\lambda$ to~$M$ along the curve~$\lambda$, so the fundamental group of~$\overline{M}$ is $\Gamma/N$. 

We wish to show that $\gamma^n$ is not null-homotopic in~$\overline{M}$. 
Suppose there is a continuous map $f \colon D^2 \to \overline{M}$, such that the restriction of~$f$ to the boundary of~$D^2$ is~$\gamma^n$. Let 
	$$D^2_0 = f^{-1}(M) ,$$
so $D^2_0$ is a surface of genus~$0$ with some number~$k$ of boundary curves. 
We may assume $f$ is \index{minimal!surface}{minimal} (i.e., the area of~$D^2$ under the pull-back metric is minimal). 
Then $D^2_0$ is a surface of curvature $\kappa(x) \le -1$ whose boundary curves are geodesics. 
Note that $f$ maps
	\begin{itemize}
	\item one boundary geodesic onto~$\gamma^n$,
	and
	\item the other $k-1$ boundary geodesics onto multiples of~$\lambda$.
	\end{itemize}
This yields a contradiction:
\begin{align*}
2 \pi (k-2)
&= -2 \pi \, \chi(D^2_0)
&& \text{\csee{Chi(PuncturedDiskEx)}}
\\&= - \int_{D^2_0} \kappa(x) \, dx 
&& \text{(Gauss-Bonnet Theorem)}
\\&\ge  \int_{D^2_0} 1 \, dx 
\\&\ge (k-1) \length(\lambda) 
&& \text{\csee{IntBdryCollarsEx}}
\\&> 2 \pi (k-1)
.  &&\qedhere \end{align*}
\end{proof}

%\begin{thm}
%If\/ $\Rrank G = 1$, then\/ $\Gamma$ has a normal subgroup~$N$, such that\/ $\Gamma/N$ is infinite, but every element of\/ $\Gamma/N$ has finite order. 
% Or maybe this needs $\Gamma$ to be cocompact?
%\end{thm}

\begin{rem}
Perhaps the simplest example of \cref{Rrank1->GammaNotAlmSimple} is when $G = \SL(2,\real)$ and $\Gamma$~is a free group \csee{FreeLattINSL2R}. In this case, it is easy to find a normal subgroup~$N$, such that $N$ and~$\Gamma/N$ are both infinite. (For example, we could take $N = [\Gamma,\Gamma]$.)
\end{rem}

There are numerous strengthenings of \cref{Rrank1->GammaNotAlmSimple} that provide infinite quotients of~$\Gamma$ with various interesting properties (if $\Rrank G = 1$). 
We will conclude this section by briefly describing just one such example.

A classical theorem of 
	\label{HNNThm}\thmindex{Higman-Neumann-Neumann}%
	Higman, Neumann, and Neumann states that every countable group can be embedded in a $2$-generated group. Since $2$-generated groups are precisely the quotients of the free group~$\free_2$ on $2$~generators, this means that $\free_2$ is ``SQ-universal'' in the following sense:

\begin{defn}
$\Gamma$ is \emph{SQ-universal} if every countable group is isomorphic to a subgroup of a quotient of~$\Gamma$. 
 (The letters ``SQ'' stand for ``subgroup-quotient\zz.'')

More precisely, the SQ-universality of~$\Gamma$ means that if $\Lambda$ is any countable group, then there exists a normal subgroup~$N$ of~$\Gamma$, such that $\Lambda$ is isomorphic to a subgroup of $\Gamma/N$.
\end{defn}

\begin{eg} \label{FnSQUnivEg}
$\free_n$ is SQ-universal, for any $n \ge 2$ \csee{FnSQUnivEx}.
\end{eg}

SQ-universality holds not only for free groups, which are lattices in $\SL(2,\real)$ \csee{FreeLattINSL2R}, but for any other lattice of real rank one:

\begin{thm} \label{Rank1SQUniv}
If\/ $\Rrank G = 1$, then\/ $\Gamma$ is SQ-universal.
\end{thm}

\begin{rem}
Although the results in this section have been stated only for~$\Gamma$, which is a lattice, the theorems are valid for a much more general class of groups. 
This is because normal subgroups can be obtained from an assumption of negative curvature (as is illustrated by the proof of \cref{Rrank1->GammaNotAlmSimple}). Indeed, \cref{Rrank1->GammaNotAlmSimple,Rank1SQUniv} remain valid when $\Gamma$ is replaced with 
%the fundamental group of any compact manifold of strictly negative curvature. More generally, the theorems are valid for 
any group that is \term[hyperbolic!Gromov]{Gromov hyperbolic} \csee{GromovHyperDefn}, or even \index{hyperbolic!group, relatively}``relatively'' hyperbolic (and not commensurable to a cyclic group).
\end{rem}

\begin{exercises}

\item \label{EnoughGammaNontrivEx}
Suppose 
	\begin{itemize}
	\item $\gamma$ and~$\lambda$ are nontrivial elements of~$\Gamma$,
	\item $\Gamma$ is torsion free,
	\item $N$ is a normal subgroup of~$\Gamma$,
	\item $\lambda \in N$,
	and
	\item $\gamma^n \notin N$, for every positive integer~$n$.
	\end{itemize}
Show that neither $N$ nor $\Gamma/N$ is infinite.

\item \label{Chi(PuncturedDiskEx)}
Show that the {Euler characteristic} of a $2$-disk with $k-1$ punctures is $2 - k$.

\item \label{IntBdryCollarsEx}
In the notation of the proof of \cref{Rrank1->GammaNotAlmSimple}, show
	$$  \int_{D^2_0} 1 \, dx \ge  (k-1) \length(\lambda) .$$
\hint{All but one of the boundary components are at least as long as~$\lambda$, and a boundary collar of width~$1$ is disjoint from the collar around any other boundary component.}

\item \label{FnSQUnivEx}
Justify \cref{FnSQUnivEg}.
\hint{You may assume the theorem of Higman, Neumann, and Neumann on embedding countable groups in $2$-generated groups.}

\end{exercises}

\section{\texorpdfstring{$\Gamma$}{Γ}-equivariant quotients of \texorpdfstring{$G/P$}{G/P}
	\texorpdfstring{\optional}{(optional)}} \label{QuotG/PPfSect}

In this section, we explain how to prove \cref{blackbox}. However, we will assume $G = \SL(2,\real) \times \SL(2,\real)$, for simplicity.

The space~$Z$ is not known explicitly, so it is difficult to study directly. Instead, as in the proof of the ergodic decomposition in \cref{ErgDecompSect}, we will look at the $\sigma$-algebra $\Bool(Z)$ of Borel sets, modulo the sets of measure~$0$. (We will think of this as the set of $\{0,1\}$-valued functions in $\LL\infty(Z)$, by identifying each set with its characteristic function.)
Note that $\psi$ induces a $\Gamma$-equivariant inclusion
	 $$\psi^* \colon \Bool(Z) \hookrightarrow \Bool(G/P) $$
\csee{B(Z)InjectsEx}. 
Via the inclusion $\psi^*$, we can identify $\Bool(Z)$ with a sub-$\sigma$-algebra of $\Bool(G/P)$:
	 $$\Bool(Z) \subseteq \Bool(G/P) .$$
%Now $\Bool(Z)$ is a $\Gamma$-invariant sub-$\sigma$-algebra of $\Bool(G/P)$ \cref{B(Z)ClosedBoolean}. 
In order to establish that $Z$ is a $G$-equivariant quotient of $G/P$, we wish to show that $\Bool(Z)$ is $G$-invariant \csee{ZGInvt->GQuot}. Therefore, \cref{blackbox}  can be reformulated as follows:

\begin{thmrefer}{blackbox}
\begin{thm} \label{blackbox'}
If $\Bool$ is any\/ $\Gamma$-invariant sub-$\sigma$-algebra of $\Bool(G/P)$, then $\Bool$ is $G$-invariant.
\end{thm}
\end{thmrefer}

To make things easier, let us settle for a lesser goal temporarily:

\begin{defn} 
The \emph{trivial} Boolean sub-$\sigma$-algebra of $\Bool(G/P)$ is $\{0,1\}$ (the set of constant functions).
\end{defn}

\begin{prop} \label{lessergoal}
If $\Bool$ is any nontrivial, $\Gamma$-invariant sub-$\sigma$-algebra of $\Bool(G/P)$, then $\Bool$ contains a nontrivial $G$-invariant Boolean algebra.
\end{prop}

\begin{rem} \ 
\noprelistbreak
\begin{enumerate}
\item To establish \cref{lessergoal}, we will find a characteristic function $\overline f \in \Bool(G/P) \setminus \{0,1\}$, such that $G \, \overline f \subseteq \Bool$.
\item The proof of \cref{blackbox'} is similar: let $\Bool_G$ be the (unique) maximal $G$-invariant Boolean subalgebra of~$\Bool$. If $\Bool_G \neq \Bool$, we will find some $\overline f \in \Bool(G/P) \setminus \Bool_G$, such that $G \, \overline f \subseteq \Bool$. (This is a contradiction.)
\end{enumerate}
\end{rem}

\begin{assump} 
To  simplify the algebra in the proof of \cref{lessergoal}, let us assume $G = \SL(2,\real) \times \SL(2,\real)$.
\end{assump}

\begin{notation} \label{BlackBoxPfNotn} \ 
\noprelistbreak
\begin{itemize}
\item $G = G_1 \times G_2$, where $G_1 = G_2 = \SL(2,\real)$,
\item $P = P_1 \times P_2$, where $P_i = \begin{Smallbmatrix} \upast& \\ \upast&\upast \end{Smallbmatrix}
\subset G_i$,
\item $U = U_1 \times U_2$, where $U_i = \begin{Smallbmatrix} 1& \\ \upast&1 \end{Smallbmatrix}
\subset P_i$,
\item $V = V_1 \times V_2$, where $V_i = \begin{Smallbmatrix} 1&\upast \\ &1 \end{Smallbmatrix}
\subset G_i$,
\item $\Gamma = $ some irreducible lattice in~$G$,
and
\item $\Bool$ = some $\Gamma$-invariant sub-$\sigma$-algebra of $\Bool(G/P)$ .
\end{itemize}
\end{notation}

\begin{rem} \label{G/PisR2}
We have $G/P = (G_1/P_1) \times (G_2/ P_2)$. Here are two useful, concrete descriptions of this space:
\begin{itemize}
\item $G/P = \RP1 \times \RP1 \iso \real^2$ (a.e.),
and
\item $G/P \iso V_1 \times V_2$ (a.e.) \csee{Gi/Pi=ViEx}.
\end{itemize}
Note that, if we identify $G/P$ with $\real^2$ (a.e.), then, for the action of~$G_1$ on~$G/P$, we have
\begin{itemize}
\item $\begin{bmatrix} k& \\ &k^{-1} \end{bmatrix} (x,y) = (k^2 x, y)$,
and
\item $\begin{bmatrix} 1& t \\ &1 \end{bmatrix} (x,y) = (x + t, y)$
\end{itemize}
\csee{GiActsByLinFracEx}.
\end{rem}

The proof of \cref{lessergoal} employs two preliminary results.
The first is based on a standard fact from first-year analysis:

\begin{lem}[(\thmindex{Lebesgue Differentiation}Lebesgue Differentiation Theorem)] \label{LebDiffThm}
Let 
	\begin{itemize}
	\item $f \in \LL1(\real^n)$, 
	\item $\lambda$ be the Lebesgue measure on~$\real^n$,
	and
	\item $B_r(p)$ be the ball of radius~$r$ centered at~$p$.
	\end{itemize}
For a.e.\ $p \in \real^n$, we have
	\begin{equation} \label{LebDiffThmEq}
	\lim_{r \to 0} \frac{1}{\lambda \bigl( B_r(p) \bigr)}\int_{B_r(p)} f \, d\lambda = f(p) 
	. \end{equation}
%where $\lambda$ is Lebesgue measure, and $B_r(x)$ is the ball of radius~$r$ centered at~$x$.
\end{lem}

Letting $n = 1$ and applying Fubini's Theorem yields:

\begin{cor} \label{LebConvMeasToProj}
Let
\noprelistbreak
\begin{itemize}
\item $f \in \LL\infty(\real^2)$,
\item $a = \begin{bmatrix} k & \\ &k^{-1} \end{bmatrix} \in G_1$, for some $k > 1$,
and
\item $\pi_2 \colon \real^2 \to \{0\}\times \real$ be the projection onto the $y$-axis.
\end{itemize}
Then, for a.e.\ $v \in V_1$,
	$$ \text{$a^n v f$ converges in measure to $(vf) \circ \pi_2$ as $n \to \infty$} . $$
%\textup(Here, $(a^n v f)(p) = f \bigl( a^{-n}(p-v) \bigr)$.\textup)
\end{cor}

\begin{proof}
\Cref{LebConvMeasToProjPfEx}.
\end{proof}

The other result to be used in the proof of \cref{lessergoal} is a consequence of the Moore Ergodicity Theorem:

\begin{prop} \label{GammaVADense}
For a.e.\ $v \in V_1$, \ $\Gamma v^{-1} a^{-\natural}$ is dense in~$G$.
\end{prop}

\begin{proof}
Taking inverses, we wish to show $\closure{a^{\natural} v \Gamma} = G$; 
i.e., the (forward) $a$-orbit of $v \Gamma$ is dense in $G/\Gamma$, for a.e.\ $v \in V_1$.
We will show that 
	$$ \text{$\closure{a^{\natural} g \Gamma} = G$, for a.e.\ $g \in G$,} $$
and leave the remainder of the proof to the reader \csee{a^nvGammaDenseEx}.

Given a nonempty open subset~$\open$\, of~$G/\Gamma$,
let 
	$$E = \bigcup_{n>0} a^{-n} \open .$$
Clearly, $a^{-1} E \subseteq E$. Since $\mu(a^{-1} E) = \mu(E)$ (because the measure on $G/\Gamma$ is $G$-invariant), this implies
	 $E$ is $a$-invariant (a.e.).
Since the Moore Ergodicity Theorem \pref{MooreErgodicity} tells us that $a$ is ergodic on $G/\Gamma$, 
we conclude that $E = G/\Gamma$ (a.e.). This means that, for a.e.\ $g \in G$, the forward $a$-orbit of~$g$ intersects~$\open$\,. 

Since $\open$\, is an arbitrary open subset, and $G/\Gamma$ is second countable, we conclude that the forward $a$-orbit of a.e.~$g$ is dense.
\end{proof}

\begin{proof}[\upshape \textbf{\mathversion{bold}Proof of \cref{lessergoal} for $G = \SL(2,\real) \times \SL(2,\real)$}] 
Identify $G/P$ with~$\real^2$, as in \cref{G/PisR2}.
Since $\Bool$ is nontrivial, it contains some nonconstant~$f$.
Now $f$ cannot be essentially constant both on almost every vertical line and on almost every horizontal line \csee{ConstVert+Horiz->Const}, so we may assume there is a non-null set of vertical lines on which it is not constant. This means that
	 $$ \text{$\bigset{v \in V_1 }{ \begin{matrix} \text{$(vf) \circ \pi_2$ is not} \\ \text{essentially constant} \end{matrix} }$ 
	has positive measure.} $$
\Cref{LebConvMeasToProj,GammaVADense} tell us we may choose~$v$ in this set, with the additional properties that 
\begin{itemize}
\item $a^n v f \to (vf) \circ \pi_2$,
and
\item $\Gamma v^{-1} a^{-\natural}$ is dense in~$G$.
\end{itemize}
Let $\overline{f} = (vf) \circ \pi_2$, so
	$$ a^n v f \to \overline{f} .$$
Now, for any $g \in G$, there exist $\gamma_i \in \Gamma$ and $n_i \to \infty$, such that 
	 $$g_i \defd \gamma_i v^{-1} a^{-n_i} \to g .$$
Then we have 
	$$g_i a^{n_i} v = \gamma_i \in \Gamma ,$$
so the $\Gamma$-invariance of~$\Bool$ implies
	$$\Bool \ni \gamma_i f = g_i \, a^{n_i} v f\to g \, \overline {f} $$
\csee{ContOnB(G/P)}.
Since $\Bool$ is closed \csee{WeakClosed->ConvMeasEx}, we conclude that $g \, \overline {f} \in \Bool$. Since $g$ is an arbitrary element of~$G$, this means $G \, \overline {f} \subseteq \Bool$.
Also, from the choice of~$v$, we know that $\overline {f} = (vf) \circ \pi_2$ is not essentially constant.
\end{proof}

Combining the above argument with a list of the $G$-invariant Boolean subalgebras of $\Bool(G/P)$ yields \cref{blackbox'}:

\begin{proof}[\upshape \textbf{\mathversion{bold}Proof of \cref{blackbox'} for $G = \SL(2,\real) \times \SL(2,\real)$}] \ 
Let $\Bool_G$ be the largest $G$-invariant subalgebra of~$\Bool$, and suppose $\Bool \neq \Bool_G$. (This will lead to a contradiction.)

It is shown in \cref{EquivQuotsSL2xSL2/PEx} that the only $G$-invariant subalgebras of $\Bool(G/P) = \Bool(\real^2)$ are  
\begin{itemize}
\item $\Bool(\real^2)$,
\item $\{\, \text{functions constant on horizontal lines (a.e.)} \,\}$,
\item $\{\, \text{functions constant on vertical lines (a.e.)} \,\}$,
and
\item $\{\, 0,1 \,\}$.
\end{itemize}
So $\Bool_G$ must be one of these $4$ subalgebras.

We know $\Bool_G \neq \Bool(\real^2)$ (otherwise $\Bool = \Bool_G$). 
Also, we know $\Bool$ is nontrivial (otherwise $\Bool = \{0,1\} = \Bool_G$), so \cref{lessergoal} tells us that $\Bool_G \neq \{0,1\}$.
Hence, we may assume, by symmetry, that 
	 \begin{equation} \label{BoolG=ConstVertEq}
	 \Bool_G = \{\, \text{functions constant on vertical lines (a.e.)} \,\}  
	 . \end{equation}
Since $\Bool \neq \Bool_G$, there is some $f \in \Bool$, such that $f$ is \emph{not} essentially constant on vertical lines.
Applying the proof of \cref{lessergoal} yields $\overline{f}$, such that
	\begin{itemize}
	\item $G \, \overline{f} \subseteq \Bool$, so $\overline{f} \in \Bool_G$,
	and
	\item $\overline{f}$ is \emph{not} essentially constant on vertical lines.
	\end{itemize}
This contradicts \pref{BoolG=ConstVertEq}.
\end{proof}

Very similar ideas yield the general case of \cref{blackbox}, if one is familiar with real roots and parabolic subgroups. To illustrate this, without using extensive Lie-theoretic language, let us explicitly describe the setup for $G = \SL(3,\real)$.

\subsection*{\hskip-\parindent Modifications for $\SL(3,\real)$.}
\noprelistbreak
\begin{itemize}
\smallskip
\item $P = \begin{Smallbmatrix} \upast&& \\ \upast&\upast& \\ \upast&\upast&\upast \end{Smallbmatrix} $,
	\quad
	$V = \begin{Smallbmatrix} 1&\upast&\upast \\ &1&\upast \\ &&1 \end{Smallbmatrix}$,
	\quad
	$V_1 = \begin{Smallbmatrix} 1&\upast& \\ &1& \\ &&1 \end{Smallbmatrix}$,
	\quad
	$V_2 = \begin{Smallbmatrix} 1&& \\ &1&\upast \\ &&1 \end{Smallbmatrix}$.
\\[3pt] Note that $V = \langle V_1, V_2 \rangle$.

\item There are exactly four subgroups containing $P$, namely,
$$ \text{$P$, \quad $G$, \quad $P_1 =
 \begin{Smallbmatrix} \vphantom{1}\upast&\upast& \\ \vphantom{1}\upast&\upast& \\ \vphantom{1}\upast&\upast&\upast \end{Smallbmatrix} = \langle V_1,  P \rangle$, \quad 
%\quad and \quad
 $P_2 = \begin{Smallbmatrix} \vphantom{1}\upast&& \\ \vphantom{1}\upast&\upast&\upast \\ \vphantom{1}\upast&\upast&\upast \end{Smallbmatrix} = \langle V_2, P \rangle$.} $$
Hence, there are precisely four $G$-invariant subalgebras of $\Bool(G/P)$. Namely, if we identify $\Bool(G/P)$ with $\Bool(V)$, then the $G$-invariant subalgebras of $\Bool(V)$ are 
\begin{itemize}
\item $\Bool(V)$, 
\item $\{0,1\}$, 
\item right $V_1$-invariant functions,
\item right $V_2$-invariant functions.
\end{itemize}

\begin{rem} \label{GeomInterpG/Pi}
The homogeneous spaces $G/P_1$ and $G/P_2$ are $\RP2$ and the Grassmannian $\Grass23$ of $2$-planes in~$\real^3$ \csee{GeomInterpG/PiEx}. Hence, in geometric terms, the $G$-invariant Boolean subalgebras of $\Bool(G/P)$ are 
	 $\Bool(G/P)$, $\{0,1\}$, $\Bool(\RP2)$, and $\Bool(\Grass23)$.
\end{rem}

\item Let $\pi_2$ be the projection onto~$V_2$ in the natural semidirect product 
 $V = V_2 \ltimes V_2^\perp$, where $V_2^\perp = \begin{Smallbmatrix} 1&&\upast \\ &1&\upast \\ &&1 \end{Smallbmatrix}$. 

\item For $a = \begin{Smallbmatrix} k && \\ & k & \\ &&1/k^2 \end{Smallbmatrix} \in G$, \cref{ADilateInSL3Ex} tells us
	 \begin{equation} \label{ADilateInSL3}
	 a \begin{bmatrix} 1&x&z \\  &1&y \\ &&1 \end{bmatrix} P
	= \begin{bmatrix} 1&x&k^3 z \\  &1&k^3 y \\ &&1 \end{bmatrix} P 
	. \end{equation}
	\item A generalization of the Lebesgue Differentiation Theorem tells us, for $f \in \Bool(G/P) = \Bool(V)$ and a.e.\ $v \in V_2^\perp$, that
	 $$ \text{$a^n v f$ converges in measure to $(vf) \circ \pi_2$.} $$

\end{itemize}
With these facts in hand, it is not difficult to prove \cref{blackbox'} under the assumption that $G = \SL(3,\real)$ \csee{BlackBoxSL3PfEx}.

\begin{exercises}

\item \label{B(Z)InjectsEx}
In the setting of \cref{blackbox}, define $\psi^* \colon \Bool(Z) \to \Bool(G/P)$ by $\psi^*(f) = f \circ \psi$.
Show that $\psi^*$ is injective and $\Gamma$-equivariant.
\hint{Injectivity relies on the fact that $\psi$ is measure-class preserving.}

%\item \label{B(Z)ClosedBoolean}
%In the setting of \cref{blackbox}, show that $\psi^* \bigl( \Bool(Z) \bigr)$ is a $\Gamma$-invariant sub-$\sigma$-algebra of $\Bool(G/P)$.

\item \label{ZGInvt->GQuot}
In the setting of \cref{blackbox}, show that if the sub-$\sigma$-algebra $\psi^* \bigl( \Bool(Z) \bigr)$ of $\Bool(G/P)$ is $G$-invariant, then $Z$ is a $G$-equivariant quotient of $G/P$ (a.e.).
\hint{To reduce problems of measurability, you may pretend that $G$ is countable. More precisely, use \cref{BoolPtwise} to show that if $H$ is any countable subgroup of~$G$ that contains~$\Gamma$, then the $\Gamma$-action can be extended to an action of~$H$ on~$Z$ by Borel maps, such that, for each $h \in H$, we have $\psi(hx) = h \, \psi(x)$ for a.e.\ $x \in G/P$.}

\item \label{GiActsByLinFracEx}
Let $G_i$ and~$P_i$ be as in \cref{BlackBoxPfNotn}.  Show that choosing appropriate coordinates on $\RP1 = \real \cup \{\infty\}$ identifies the action of $G_i$ on $G_i/P_i$ with the action of $G_i = \SL(2,\real)$ on $\real \cup \{\infty\}$ by linear-fractional transformations:
	$$ \begin{bmatrix} a & b \\ c & d \end{bmatrix} (x) = \frac{ax + b}{cx + d} .$$
In particular, 
	$$\begin{bmatrix} k & \\ & k^{-1} \end{bmatrix} (x) = k^2 x 
	\text{\qquad and\qquad}
	 \begin{bmatrix} 1 & t\\ & 1 \end{bmatrix} (x) = x + t . $$
\hint{Map a nonzero vector $(x_1,x_2) \in \real^2$ to its reciprocal slope $x_1/x_2 \in \real \cup \{\infty\}$.}

\item \label{Gi/Pi=ViEx}
Let $G_i$, $P_i$, and~$V_i$ be as in \cref{BlackBoxPfNotn}. 
Show that the map $V_i \to G_i/P_i \colon v \mapsto v P_i$ injective and measure-class preserving.
\hint{\Cref{GiActsByLinFracEx}.} 

\item \label{LebDiffThmEquivEx}
Show that \cref{LebDiffThmEq} is equivalent to
	$$ \lim_{k \to \infty} \frac{1}{\lambda \bigl( B_1(0) \bigr)} \int_{B_1(0)} f\left(p + \frac{x}{k}\right) \, d\lambda(x)= f(p) .$$
\hint{A change of variables maps $B_1(0)$ onto $B_r(p)$ with $r = 1/k$.}

\item \label{LebConvMeasToProjPfEx}
Prove \cref{LebConvMeasToProj}.
\hint{\Cref{LebDiffThmEquivEx}.}
 
 \item \label{a^nvGammaDenseEx}
 Complete the proof of \cref{GammaVADense}: assume, for a.e.\ $g \in G$, that $a^{\natural} g \Gamma$ is dense in~$G$, and show, for a.e.\ $v \in V_1$, that $a^{\natural} v \Gamma$ is dense in~$G$.
 \hint{If $a^{\natural} g \Gamma$ is dense, then the same is true when $g$ is replaced by any element of $\czer_G(a) \, U_1 \, g$.}

\item \label{ConstVert+Horiz->Const}
Let $f \in \Bool(\real^2)$. Show that if $f$ is essentially constant on a.e.\ vertical line and on a.e.\ horizontal line, then $f$~is constant (a.e.).
 
 \item \label{ParabolicsInSL2xSL2Ex}
Assume \cref{BlackBoxPfNotn}. Show that the only subgroups of~$G$ containing~$P$ are  $P$, $G_1\times P_2$, $P_1\times G_2$, and~$G$. 
\hint{$P$ is the stabilizer of a point in $\RP1 \times \RP1$, and has only $4$ orbits.} 

\item \label{EquivQuotsSL2xSL2/PEx}
Assume \cref{BlackBoxPfNotn}. 
Show that the only $G$-equivariant quotients of $G/P$ are
	$G/P$, $G_2/P_2$, $G_1/P_1$, and $G/G$.
\hint{\cref{ParabolicsInSL2xSL2Ex}.}

\item \label{WeakClosed->ConvMeasEx}
Suppose $\Bool$ is a sub-$\sigma$-algebra of $\Bool(G/P)$. Show that $\Bool$ is closed under \term{convergence in measure}. 
\par
More precisely, fix a probability measure $\mu$ in the Lebesgue measure class on $G/P$, and show that $\Bool$ is a closed in the topology corresponding to the metric on $\Bool(G/P)$ that is defined by $d(A_1,A_2) = \mu( A_1 \symmdiff A_2)$. 

\item \label{ContOnB(G/P)}
Show that the action of~$G$ on $\Bool(G/P)$ is continuous.
\hint{Suppose $g_n \to e$ and $\mu(A_n \symmdiff A) \to 0$. The Radon-Nikodym derivative $d (g_n)_*\mu/d\mu$ tends uniformly to~$1$, so $\mu(g_n A_n \symmdiff g_nA) \to 0$. To bound $\mu(g_n A \symmdiff A)$, note that $\int_{g_n A} \varphi \, d\mu \to \int_A \varphi d\mu$, for every $\varphi \in C_c(G/P)$.}

\item \label{GeomInterpG/PiEx}
In the notation of \cref{GeomInterpG/Pi}, show that $G/P_1$ and $G/P_2$ are $G$-equivariantly diffeomorphic to  $\RP2$ and~$\Grass23$, respectively.
\hint{Verify that the stabilizer of a point in $\RP2$ is~$P_1$, and the stabilizer of a point in $\Grass23$ is~$P_2$.}

\item \label{ADilateInSL3Ex}
Verify \cref{ADilateInSL3}.
\hint{Since $a \in P$, we have $agP = (a g a^{-1}) P$, for any $g \in G$.}

\item \label{BlackBoxSL3PfEx}
Prove \cref{blackbox'} under the assumption that $G = \SL(3,\real)$.
\hint{You may assume (without proof) the facts stated in the ``Modifications for $\SL(3,\real)$\zz.''}
 
 \end{exercises}

\begin{notes}

The Normal Subgroups Theorem \pref{MargNormalSubgrpsThm} is due to G.\,A.\,Margulis \cite{MargulisFactorGroupsDoklady,MargulisFactorGroups,MargulisFactorGroupsRank1}.
Expositions of the proof appear in \cite[Chap.~4]{MargulisBook} and \cite[Chap.~8]{ZimmerBook}.
(However, the proof in \cite{ZimmerBook} assumes that $G$ has Kazhdan's property~$(T)$.) 

When $\Gamma$ is not cocompact, the Normal Subgroups Theorem can be proved by algebraic methods derived from the proof of the Congruence Subgroup Problem (see \cite[Thms.~A and~B, p.~109]{RaghunathanCSP} and \cite[Cor.~1, p.~75]{RaghunathanCSP2}). 
On the other hand, it seems that the ergodic-theoretic approach of Margulis provides the only known proof in the cocompact case.

Regarding \fullcref{MargNormSubgrpRems}{CSP}, see 
	%\cite[Chap.~6]{Humphreys-ArithmeticGroups} or 
\cite{Sury-CSPBook} for an introduction to the Congruence Subgroup Property. 

\Cref{blackbox} is stated for general $G$ of real rank $\ge 2$ in \cite[Cor.~2.13]{MargulisBook} and \cite[Thm.~8.1.4]{ZimmerBook}.
\Cref{blackbox'} is in \cite[Thm.~4.2.11]{MargulisBook} and \cite[Thm.~8.1.3]{ZimmerBook}. See \cite[\S8.2 and \S8.3]{ZimmerBook} and \cite[\S4.2]{MargulisBook} for expositions of the proof.

The proof of \cref{Rrank1->GammaNotAlmSimple} is adapted from \cite[5.5.F, pp.~150--152]{Gromov-HypGrp}.

The Higman-Neumann-Neumann Theorem on SQ-universality of~$\free_2$ (see p.~\pageref{HNNThm}) was proved in \cite{HNN-F2IsSQUniv}. A very general version of \cref{Rank1SQUniv} that applies to all relatively hyperbolic groups was proved in \cite{HyperIsSQUniv}.
(The notion of a \index{hyperbolic!group, relatively}relatively hyperbolic group was introduced in~\cite{Gromov-HypGrp}, and generalized in \cite{Farb-RelHypGrp}.)

\end{notes}

%!TEX root = IntroArithGrps.tex

\mychapter{\texorpdfstring{Arithmetic Subgroups\\of Classical Groups}%
	{Arithmetic Subgroups of Classical Groups}}
\label{ArithClassicalChap}

\prereqs{Restriction of Scalars (\cref{RestrictScalarsSect}) and examples of arithmetic subgroups (\cref{EgArithGrpsChap}).}

This \lcnamecref{ArithClassicalChap} will give a quite explicit description (up to commensurability) of all the arithmetic subgroups of almost every classical Lie group~$G$ \csee{ArithLattsAreClassical}. (Recall that a simple Lie group~$G$ is ``classical'' if it is either a special linear group, an orthogonal group, a unitary group, or a symplectic group \csee{ClassicalDefn}.) The key point is that all the $\rational$-forms of~$G$ are also classical, not exceptional, so they are fairly easy to understand. 
However, there is an exception to this rule: some $8$-dimensional orthogonal groups have $\rational$-forms of so-called ``\index{triality}triality type'' that are not classical and will not be discussed in any detail here \csee{D4weird}.
%From the correspondence between $\rational$-forms and arithmetic subgroups, this can be thought of as a list of all the arithmetic subgroups of~$G$ (up to commensurability).

Given~$G$, which is a Lie group over~$\real$, we would like to know all of its $\rational$-forms (because, by definition, arithmetic groups are made from $\rational$-forms). However, we will start with the somewhat simpler problem that replaces the fields $\rational$ and~$\real$ with the fields $\real$ and~$\complex$: finding the $\real$-forms of the classical Lie groups over~$\complex$.

\section{\texorpdfstring{$\real$}{R}-forms of classical simple groups over%
\texorpdfstring{~$\complex$}{ C}} \label{RFormsOfCGrps}

To set the stage, let us recall the classical result that almost all complex simple groups are classical:

\begin{thm}[(Cartan, Killing)] \label{AlmostAllOverC}
All but finitely many of the simple Lie groups over\/~$\complex$ are isogenous to either\/ $\SL(n,\complex)$, $\SO(n,\complex)$, or\/ $\Sp(2n,\complex)$, for some~$n$.
\end{thm}

\begin{rem} \label{EFGOverCRem}
Up to isogeny, there are exactly five simple Lie groups over~$\complex$ that are not classical. They are the ``\index{exceptional Lie group}exceptional'' simple groups, and are called $E_6$, $E_7$, $E_8$, $F_4$, and~$G_2$.
%Although we do not need this terminology, each simple group has a \defit[type (of a simple Lie group)]{type}:
%	\begin{itemize}
%	\item $\SL(n,\complex)$ is of type~$A$, 
%	\item $\SO(n,\complex)$ is of type~$B$ or~$D$, depending on whether $n$ is odd or even, respectively,
%	and
%	\item $\Sp(2n,\complex)$ is of type~$C$.
%	\end{itemize}
%The remaining simple groups over~$\complex$ are of type~$E$, $F$, or~$G$. Up to isogeny, there are exactly three of type~$E$, one of type~$F$, and one of type~$G$.
\end{rem}

Now, we would like to describe the $\real$-forms of each of the classical groups. For example, finding all the $\real$-forms of $\SL(n,\complex)$ would mean making a list of the (simple) Lie groups~$G$, such that the ``complexification'' of~$G$ is $\SL(n,\complex)$. This is not difficult, but we should perhaps begin by explaining more clearly what it means.

It has already been mentioned that, intuitively, the complexification of~$G$ is the
complex Lie group that is obtained from~$G$ by replacing
real numbers with complex numbers. 
 For example,
the complexification of $\SL(n,\real)$ is $\SL(n,\complex)$. In general, $G$ is
(isogenous to) the set of real solutions of a certain set of
equations, and we let $G_{\complex}$ be the set of
complex solutions of the same set of equations:

\begin{notation} \label{GxCNot}
 Assume $G \subseteq \SL(\ell,\real)$, for some~$\ell$. Since $G$ is almost Zariski closed \csee{GisAlgic}, there
is a certain subset~$\mathcal{Q}$ of $\real[x_{1,1}, \ldots,
x_{\ell,\ell}]$, such that $G^\circ = \Var(\mathcal{Q})^\circ$. Let
 $$ G_{\complex} = \Var_{\complex}(\mathcal{Q})
 = \{\, g \in \SL(\ell,\complex) \mid \mbox{$Q(g) = 0$, for
all $Q \in \mathcal{Q}$} \,\} .$$
 Then $G_{\complex}$ is a (complex, semisimple) Lie
group.
 \end{notation}

\begin{eg} \label{GxCEg}
 \ 
 \begin{enumerate}
 \item \label{GxCEg-SL}
 $\SL(n,\real)_{\complex} = \SL(n,\complex)$. 
 \item \label{GxCEg-SO}
 $\SO(n)_{\complex} = \SO(n,\complex)$.
 \item \label{GxCEg-SOpq}
 $\SO(m,n)_{\complex} \iso \SO(m+n,\complex)$
\csee{SOpqxC}.
 \end{enumerate}
 \end{eg}

\begin{defn}
 If $G_{\complex}$ is isomorphic to~$H$, then we say
that 
 \begin{itemize}
 \item $H$ is the \defit{complexification} of~$G$, and that
 \item $G$ is an \defit[R-@$\real$-!form]{$\real$-form} of~$H$.
 \end{itemize}
 \end{defn}

The following result lists the complexification of each classical group.
It is not difficult to
memorize the correspondence. For example, it is obvious
from the notation that the complexification of $\Sp(m,n)$ should be
symplectic. Indeed, the only case that really requires
memorization is the complexification of $\SU(m,n)$
\fullcsee{GxC=}{SUmn}.

\begin{prop} \label{GxC=}
 Here is the complexification of each classical Lie group.
 \noprelistbreak
	 \begin{enumerate} \renewcommand{\theenumi}{\Alph{enumi}}
			\renewcommand{\theenumii}{\roman{enumii}}
	 \item \label{GxC=-SL}
	 Real forms of special linear groups:
	 \noprelistbreak
	 	\begin{enumerate} 
		 \item $\SL(n,\real)_{\complex} = \SL(n,\complex)$,
		 \item $\SL(n,\complex)_{\complex} \iso \SL(n,\complex)
		\times \SL(n,\complex)$,
		 \item $\SL(n,\quaternion)_{\complex} \iso
		\SL(2n,\complex)$,
		 \item \label{GxC=-SUmn}
		 $\SU(m,n)_{\complex} \iso \SL(m+n,\complex)$.
		 \end{enumerate}

	 \item Real forms of orthogonal groups:
	 \noprelistbreak
	 	\begin{enumerate}
		 \item $\SO(m,n)_{\complex} \iso \SO(m+n,\complex)$,
		 \item $\SO(n,\complex)_{\complex} \iso \SO(n,\complex)
		\times \SO(n,\complex)$,
		 \item $\SO(n,\quaternion)_{\complex} \iso
		\SO(2n,\complex)$.
		 \end{enumerate}

	 \item Real forms of symplectic groups:
	 \noprelistbreak
	 	\begin{enumerate}
		 \item $\Sp(n,\real)_{\complex} = \Sp(n,\complex)$,
		 \item $\Sp(n,\complex)_{\complex} \iso \Sp(n,\complex)
		\times \Sp(n,\complex)$,
		 \item $\Sp(m,n)_{\complex} \iso
		\Sp \bigl( 2(m+n),\complex \bigr)$.
		 \end{enumerate}
	 \end{enumerate}
 \end{prop}

Some parts of this \lcnamecref{GxC=} are more-or-less obvious (such as $\SL(n,\real)_{\complex} = \SL(n,\complex)$). A few
other examples appear in \cref{CalcGxCSect} below, and the
methods used there can be applied to all of the cases. In
fact, all of the calculations are straightforward
adaptations of the examples, except perhaps the determination of $\SO(n,\quaternion)_{\complex}$ \csee{SOnHxC}.

Nothing in \cref{GxC=} is very surprising. What is not at all obvious is that this list of real forms is complete:

\begin{thm}[(\'E.\,Cartan)] \label{RformsComplete}
Every real form of\/ $\SL(n,\complex)$, $\SO(n,\complex)$, or\/ $\Sp(n,\complex)$ appears in \cref{GxC=}\/ \textup(up to isogeny\/\textup).
\end{thm}

We will discuss a proof of this \lcnamecref{RformsComplete} in \cref{GaloisCohoRealFormsSect}.

\begin{rems} \label{GxC=tensorprod} \ 
\noprelistbreak
 \begin{enumerate}
 
 \item From \cref{GxC=}, we see that a single complex
group may have several different real forms. However, there are always only finitely many (even for exceptional groups).

 \item \label{GxC=tensorprod-LieAlg}
 The Lie algebra of $G_{\complex}$ is the tensor
product $\Lie G \otimes \complex$ \csee{Lie(GxC)}. This is
independent of the embedding of~$G$ in $\SL(\ell,\complex)$,
so, up to isogeny, $G_{\complex}$ is
independent of the embedding of~$G$ in $\SL(\ell,\complex)$.
%Thus, we will often speak of the complexification of~$G$,
%even if $G$~is not a group of matrices, with the
%understanding that the complexification is not entirely well
%defined \see{GCnotWellDefd}.

 \item \label{GxC=tensorprod-Q}
 We ignored a technical issue in \cref{GxCNot}:
there may be many different choices of~$\mathcal{Q}$ (having
the same set of real solutions), and it may be the case that
different choices yield different sets of complex
solutions. (In fact, a bad choice of~$\mathcal{Q}$ can yield
a set of complex solutions that is not a group.) To
eliminate this problem, we should insist that $\mathcal{Q}$
be maximal; that is,
 $$ \mathcal{Q} = \{\, Q \in \real[x_{1,1}, \ldots,
x_{\ell,\ell}] \mid \mbox{$Q(g) = 0$, for all $g \in G$}
\,\} .$$
 Then $G_{\complex}$ is the Zariski closure of~$G$
(over the field~$\complex$), from which it follows that $G_\complex$, like~$G$, is a semisimple Lie group.
 \end{enumerate}
 \end{rems}

\begin{eg} \label{GCnotWellDefd}
 Because the center of $\SL(3,\real)$ is trivial, we see
that $\SL(3,\real)$ is the same Lie group as $\PSL(3,\real)$.
On the other hand, we have
 $$ \SL(3,\real)_{\complex} = \SL(3,\complex) \not\iso
\PSL(3,\complex) = \PSL(3,\real)_{\complex} .$$
 This is a concrete illustration of the fact that different
embeddings of $G$ can yield different complexifications.
Note, however, that $\SL(3,\complex)$ is isogenous to
$\PSL(3,\complex)$, so the difference between the
complexifications is negligible
\fullccf{GxC=tensorprod}{LieAlg}.
 \end{eg}

\begin{exercises}

\item \label{SOpqxC}
 Show that $\SO(m,n)_{\complex} \iso
\SO(m+n,\complex)$. 
 \hint{$\SO(m,n)_{\complex}$ is conjugate to
$\SO(m+n,\complex)$ in $\SL(m+n,\complex)$, because $-1$~is a
square in~$\complex$.}

\item \label{Lie(GxC)}
 Show that the Lie algebra of $G_{\complex}$ is $\Lie G \otimes \complex$.

\end{exercises}

\section{Calculating the complexification of%
	\texorpdfstring{~$G$}{ G}}
\label{CalcGxCSect}

This section justifies \cref{GxC=}, by calculating the complexification of each classical group.

Let us start with $\SL(n,\complex)$. This is already a complex Lie group, but
we can think of it as a real Lie group of twice the
dimension. As such, it has a complexification:

\begin{prop} \label{SLnCxC}
 $\SL(n,\complex)_{\complex} \iso \SL(n,\complex)
\times \SL(n,\complex)$.
 \end{prop}

\begin{proof}
 We should embed $\SL(n,\complex)$ as a subgroup of
$\SL(2n,\real)$, find the corresponding set~$\mathcal{Q}$ of
defining polynomials, and determine the complex solutions.
However, it is more convenient to sidestep some of these
calculations by using restriction of scalars, the method
described in \S\ref{RestrictScalarsSect}.

Define $\Delta \colon \complex \to \complex \oplus \complex$
by $\Delta(z) = (z, \cjg{z})$. Then the vectors $\Delta(1) =
(1,1)$ and $\Delta(i) = (i,-i)$ are linearly independent
(over~$\complex$), so they form a basis of $\complex \oplus
\complex$. Thus, $\Delta(\complex)$ is the $\real$-span
of a basis, so it is a $\real$-form of $\complex \oplus
\complex$. Therefore, letting $V = \complex^{2n}$, we see
that
 $$ V_{\real} = \Delta(\complex^n)
 = \bigl\{\, ( v, \overline{v} )
 \mid v \in \complex^n \, \bigr\} $$
 is a real form of~$V$. Let
 $$ \bigl( \SL(n,\complex) \times 
\SL(n,\complex) \bigr)_{\real}
 = \bigset{
 g \in \SL(n,\complex) \times \SL(n,\complex)
 }{
 g (V_{\real}) = V_{\real}
 } .$$
 Then we have an isomorphism 
 $$\tilde \Delta \colon
\SL(n,\complex) \towith{\iso} \bigl(
\SL(n,\complex) \times \SL(n,\complex) \bigr)_{\real} ,$$
 defined by $\tilde
\Delta(g) = (g,\cjg{g})$, so
 \begin{align*}
  \SL(n,\complex)_{\complex}
 \iso \bigl( [\SL(n,\complex) \times \SL(n,\complex)]_{\real}
\bigr)_{\complex}
 = \SL(n,\complex) \times \SL(n,\complex) 
 . & \qedhere \end{align*}
 \end{proof}

\begin{rems} \ \label{GxC=GxGiff}
 \noprelistbreak
 \begin{enumerate}
 \item Generalizing \cref{SLnCxC}, one can show that if
$G$~is isogenous to a complex Lie group, then $G_\complex$ is isogenous to $G \times G$. 
%For example,
% $$ \SO(1,3)_{\complex}
% \iso \SO(4,\complex)
% \sim \SL(2,\complex) \times \SL(2,\complex)
% \sim \SO(1,3) \times \SO(1,3) $$
% (see \pref{isogTypes-D2-SO4C} and~\pref{isogTypes-D2-SO13}
%of \cref{isogTypes}).
 \item From \cref{SLnCxC}, we see that $G_\complex$ need not be simple, even if $G$~is simple.
However, this only happens when $G$~is complex: if $G$~is
simple, and $G$~is not isogenous to a complex Lie group,
then $G_{\complex}$ is simple.
 \end{enumerate}
 \end{rems}

Although not stated explicitly there, the proof of
\cref{SLnCxC} is based on the fact that $\complex
\otimes_{\real} \complex \iso \complex \oplus \complex$.
Namely, the map
 $$ \mbox{$\complex \otimes_{\real} \complex \to \complex
\oplus \complex$
 defined by $v \otimes \lambda \mapsto \Delta(v) \, \lambda$} $$
 is an isomorphism of $\complex$-algebras. Analogously, understanding the complexification of a group defined from the algebra~$\quaternion$ of quaternions will be
based on a calculation of $\quaternion \otimes_{\real}
\complex$.

\begin{lem} \label{HxC}
 The tensor product $\quaternion \otimes_{\real}
\complex$ is isomorphic to $\Mat_{2 \times 2}(\complex)$.
 \end{lem}

\begin{proof}
  Define an $\real$-linear map $\phi \colon \quaternion \to
\Mat_{2 \times 2}(\complex)$ by
 $$ \phi(1) = \Id, 
 \quad \phi(i)
 = \begin{bmatrix}
 i & 0 \\
 0 & -i
 \end{bmatrix}
 ,
 \quad \phi(j)
 = \begin{bmatrix}
 0 & 1 \\
 -1 & 0
 \end{bmatrix}
 ,
 \quad \phi(k)
 = \begin{bmatrix}
 0 & i \\
 i & 0
 \end{bmatrix}
 .
 $$
 It is straightforward to verify that $\phi$ is an injective
ring homomorphism. Furthermore, $\phi \bigl( \{1,i,j,k\}
\bigr)$ is a $\complex$-basis of $\Mat_{2 \times
2}(\complex)$. Therefore, the map $\hat\phi \colon
\quaternion \otimes \complex \to \Mat_{2 \times
2}(\complex)$ defined by $\hat{\phi}(v \otimes \lambda) = \phi(v)
\, \lambda$ is a ring isomorphism \csee{HxCExer}.
 \end{proof}

\begin{prop} \label{SLnHxC}
 $\SL(n,\quaternion)_{\complex} \iso \SL(2n,\complex)$.
 \end{prop}

\begin{proof}
 From \cref{HxC}, we have
 $$ \SL(n,\quaternion)_{\complex}
 \iso \SL \bigl( n, \Mat_{2 \times 2}(\complex) \bigr)
 \iso \SL(2n,\complex) $$
 \csee{SLnHxCExer,Matd(Matn)}.
 \end{proof}

As additional examples, let us look at the
complexifications of the classical simple Lie
groups that are compact, namely, $\SO(n)$, $\SU(n)$, and $\Sp(n)$. 
As observed in \fullcref{GxCEg}{SO}, we
have $\SO(n)_{\complex} = \SO(n,\complex)$. The other
cases are not as obvious.

\begin{prop} \label{SUnxC}
 $\SU(n)_{\complex} = \SL(n,\complex)$.
 \end{prop}

\begin{proof}
 Let 
 \begin{itemize}
 \item $\sigma \colon \complex \to \complex$,
 \item $\vector{\sigma} \colon \complex^n \to \complex^n$,
and
 \item $\widetilde{\sigma} \colon \SL(n,\complex) \to
\SL(n,\complex)$
 \end{itemize}
 be the
usual complex conjugations $\sigma(z) = \cjg{z}$,
$\vector{\sigma}(v) = \cjg{v}$, and
$\widetilde{\sigma}(g) = \cjg{g}$. We have
 \begin{align*}
 \SU(n) &= \{\, g \in \SL(n,\complex) \mid g^* g = \Id \,\}
\\
 &= \{\, g \in \SL(n,\complex) \mid \widetilde{\sigma}(g^\transpose)
g = \Id \,\} ,
 \end{align*}
 so, in order to calculate $\SU(n)_{\complex}$, we
should determine the map $\widetilde{\eta}$ on $\SL(n,\complex)
\times \SL(n,\complex)$ that corresponds
to~$\widetilde{\sigma}$ when we identify $\complex^n$ with
$(\complex^n \oplus \complex^n)_{\real}$ under the
map~$\vector{\Delta}$. 

First, let us determine $\vector{\eta}$. That
is, we wish to identify $\complex^n$ with~$\real^{2n}$, and
extend $\vector{\sigma}$ to a $\complex$-linear map
on~$\complex^{2n}$. However, as usual, we use the
$\real$-form $\vector{\Delta}(\complex^n)$, in place of
$\real^{2n}$. It is obvious that if we 
	$$ \text{define $\vector{\eta}
\colon \complex^n \oplus \complex^n \to \complex^n \oplus
\complex^n$ by $\vector{\eta}(x,y) = (y,x)$,} $$
then
$\vector{\eta}$ is $\complex$-linear, and the following
diagram commutes:
 $$ 
 \begin{matrix}
 \complex^n & \towith{\vector{\Delta}} & \complex^n
\oplus \complex^n \\
 \mapdown{\vector{\sigma}} & & \mapdown{\vector{\eta}} \\
 \complex^n & \towith{\vector{\Delta}} & \complex^n \oplus
\complex^n 
 . \end{matrix}
 $$
 Thus, it is fairly clear that $\widetilde{\eta}(g,h) = (h,g)$. 
 Hence
 \begin{align*}
 \SU(n)_{\complex}
 &= \{\, (g,h) \in \SL(n,\complex)
\times \SL(n,\complex) \mid \widetilde{\eta} \bigl( g^\transpose,
h^\transpose \bigr) (g,h) = (\Id,\Id) \,\} \\
 &= \{\, (g,h) \in \SL(n,\complex)
\times \SL(n,\complex) \mid \bigl( h^\transpose,
g^\transpose \bigr) (g,h) = (\Id,\Id) \,\} \\
 &= \{\, \bigl( g,(g^\transpose)^{-1} \bigr) 
\mid g  \in
\SL(n,\complex) \,\} \\
 &\iso \SL(n,\complex) .
&& \qedhere \end{align*}
 \end{proof}

\begin{prop} \label{SpnxC}
 $\Sp(n)_{\complex} = \Sp(2n,\complex)$.
 \end{prop}

\begin{proof}
 Let 
 \begin{itemize}
 \item $\phi \colon \quaternion \hookrightarrow \Mat_{2
\times 2}(\complex)$ be the embedding that is described in the proof
of \cref{HxC},
 \item $\tau$ be the usual conjugation on~$\quaternion$, 
 \item $J = 
 \begin{bmatrix}
 0 & 1 \\
 -1 & 0 
 \end{bmatrix}
 $, and
 \item $\eta \colon \Mat_{2 \times 2}(\complex)  \to \Mat_{2
\times 2}(\complex)$ be defined by 
 $\eta(x) = J^{-1} x^{\transpose} J$.
 \end{itemize}
 Then $\eta$ is $\complex$-linear, and the following diagram
commutes:
 $$ 
 \begin{matrix}
 \quaternion & \towith{\phi} & \Mat_{2 \times 2}(\complex) \\
 \mapdown{\tau} & & \mapdown{\eta} \\
 \quaternion & \towith{\phi} & \Mat_{2 \times 2}(\complex)
 . \end{matrix}
 $$

Hence, because
 \begin{align*}
 \Sp(2)
 &= \bigl\{\, g \in \SL(2,\quaternion) \mid g^* g = \Id \,
\bigr\} \\
 &= 
 \bigset{
 \begin{bmatrix}
 a & b \\
 c & d
 \end{bmatrix}
 \in \SL(2,\quaternion)
 }{
 \begin{bmatrix}
 \tau(a) & \tau(c) \\
 \tau(b) & \tau(d)
 \end{bmatrix}
 \begin{bmatrix}
 a & b \\
 c & d
 \end{bmatrix}
 = \Id
 }
 ,
 \end{align*}
 we see that
 \begin{align*}
 \Sp(2)_{\complex}
 &=
 \bigset{\!
 \begin{bmatrix}
 a & b \\
 c & d
 \end{bmatrix}
 \in \SL\bigl( 2, \Mat_{2 \times 2}(\complex) \bigr)
 }{
 \begin{bmatrix}
 \eta(a) & \eta(c) \\
 \eta(b) & \eta(d)
 \end{bmatrix}
 \begin{bmatrix}
 a & b \\
 c & d
 \end{bmatrix}
 = \Id
 \! } \\
 &=
 \bigset{\!
 \begin{bmatrix}
 a & b \\
 c & d
 \end{bmatrix}
 \in \SL\bigl( 2, \Mat_{2 \times 2}(\complex) \bigr)
 }{
 J^{-1}
 \begin{bmatrix}
 a^\transpose & c^\transpose \\
 b^\transpose & d^\transpose
 \end{bmatrix}
 J
 \begin{bmatrix}
 a & b \\
 c & d
 \end{bmatrix}
 = \Id
 \! } \\
 &=
 \bigset{
 g
 \in \SL(4,\complex)
 }{
 J^{-1}
 g^\transpose
 J
 g
 = \Id
 } \\
 &=
 \bigset{
 g
 \in \SL(4,\complex)
 }{
 g^\transpose
 J
 g
 = J
 } \\
 &= \Sp(4,\complex)
 .
 \end{align*}
 Similarly, letting
 $$ \hat{J}_n =
 \begin{bmatrix}
 J \\
 & J \\
 & & \ddots \\
 & & & J 
 \end{bmatrix}
 \in \SL(2n,\complex)
 ,$$
 the same calculations show that
	 \begin{align*} \Sp(n)_{\complex}
	 = \{\, g \in \SL(2n,\complex) \mid g^\transpose \hat{J}_n g
	= \hat{J}_n \,\}
	 \iso \Sp(2n,\complex)
	 . & \qedhere \end{align*}
 \end{proof}

\begin{exercises}

\item \label{HxCExer}
 In the proof of \cref{HxC}, verify:
 \begin{enumerate}
 \item $\phi$ is an injective ring homomorphism,
 \item $\phi \bigl( \{1,i,j,k\} \bigr)$ is a
$\complex$-basis of $\Mat_{2 \times 2}(\complex)$, and
 \item $\hat\phi$ is an isomorphism of $\complex$-algebras.
 \end{enumerate}

\item \label{SLnHxCExer}
 Show $\SL(n,\quaternion)_{\complex}
 \iso \SL \bigl( n, \Mat_{2 \times 2}(\complex) \bigr)$.
 \hint{Define $\phi$ as in the proof of \cref{HxC}. Use
the proof of \cref{SLnCxC}, with $\phi$ in the place
of~$\Delta$.}

\item \label{Matd(Matn)}
 Show $\SL \bigl( n, \Mat_{d \times d}(\complex) \bigr)
 \iso \SL(dn,\complex)$.

\item \label{SOnHxC}
 Show that $\SO(n,\quaternion)_{\complex} \iso
\SO(2n,\complex)$. 
 \hint{Similar to \pref{SpnxC}. To calculate~$\tau_r \otimes
\complex$, note that $\tau_r(x) = j^{-1} \, \tau(x) \, j$,
for $x \in \quaternion$.}

\end{exercises}

\section{How to find the real forms of complex groups} \label{GaloisCohoRealFormsSect}

In this \lcnamecref{GaloisCohoRealFormsSect}, we will explain how to find all of the possible $\real$-forms of $\SL(n,\complex)$. (Similar techniques can be used to justify the other cases of \cref{RformsComplete}, but additional calculations are needed, and we omit the details.)
We take an algebraic approach, based on \index{Galois!Theory}{Galois theory}, and we first review the most basic terminology from the theory of (nonabelian) group cohomology.

\subsection{Definition of the first cohomology of a group}
 \label{GroupCohoSect}

\begin{defns}
Suppose a group~$X$ acts (on the left) by automorphisms on a group~$M$. (For $x \in X$ and $m \in M$, we write ${}^x m$ for the image of~$m$ under~$x$.)
	\begin{enumerate}
	\item A function $\alpha \colon X \to M$ is a \defit[cocycle]{$1$-cocycle} (or ``\term{crossed homomorphism}'') if 
		$$ \text{$\alpha(xy) = \alpha(x) \cdot {}^x \! \alpha(y)$ for all $x,y \in X$} .$$
	\item Two $1$-cocycles $\alpha$ and~$\beta$ are equivalent (or ``\term[cocycle!cohomologous]{cohomologous}'') if there is some $m \in M$, such that 
		$$ \text{$\alpha(x) = m^{-1} \cdot \beta(x) \cdot {}^x m$ \ for all $x \in X$} .$$
	\item \nindex{$\coho1(X;M)$ = 1st cohomology of~$X$}%
	$\coho1(X;M)$ is the set  of equivalence classes of all $1$-cocycles. It is called the \defit[cohomology (1st)]{first cohomology} of~$X$ with coefficients in~$M$.
	\item A $1$-cocycle is a \defit{coboundary} if it is cohomologous to the trivial $1$-cocycle defined by $\tau(x) = e$ for all $x \in X$.
	\end{enumerate}
\end{defns}

\begin{warn}
In our applications, the coefficient group~$M$ is sometimes \textbf{non}abelian. In this case, $\coho1(X;M)$ is a set with no obvious algebraic structure. However, if $M$ is an abelian group (as is often assumed in textbooks on group cohomology), then $\coho1(X;M)$ is an abelian group. 
\end{warn}

\subsection{How Galois cohomology comes into the picture}
For convenience, let $G_{\complex} = \SL(n,\complex)$. Suppose $\rho \colon G_{\complex} \to \SL(N,\complex)$ is an embedding, such that $\rho \bigl( G_{\complex} \bigr)$ is defined over~$\real$. We wish to find all the possibilities for the group 
	$ \rho(G_{\complex})_\real = \rho (G_{\complex}) \cap \SL(N,\real) $
that can be obtained by considering all the possible choices of~$\rho$. 
%(The $\real$-form corresponding to this embedding is $G_{\real} = \rho^{-1} \bigl[ \rho(G_{\complex})_\real \bigr)$.)

Let $\sigma$ denote complex conjugation, the nontrivial Galois automorphism of~$\complex$ over~$\real$. Since
	$ \real = \bigset{ z \in \complex }{ \sigma(z) = z } $,
we have
	$$ \SL(N,\real) = \bigset{ g \in \SL(N,\complex) }{ \sigma(g) = g } ,$$
where we apply $\sigma$ to a matrix by applying it to each of the matrix entries. Therefore
	$$ \rho(G_{\complex})_\real 
	= \rho(G_{\complex}) \cap \SL(N,\real)
	= \bigset{ g \in \rho(G_{\complex}) }{ \sigma(g) = g} .$$
Since $\rho(G_{\complex})$ is defined over~$\real$, we know that it is invariant under~$\sigma$, so we have
	$$ %\SL(3,\complex) \stackrel{\sigma}{\longrightarrow}
	 G_{\complex} \stackrel{\rho}{\longrightarrow}
	  \rho(G_\complex) \stackrel{\sigma}{\longrightarrow}
	   \rho(G_\complex) \stackrel{\rho^{-1}}{\longrightarrow}
	    G_{\complex} .$$
Let $\widetilde\sigma = \rho^{-1} \sigma \rho \colon G_{\complex} \to G_{\complex}$ be the composition. 
Then the real form corresponding to~$\rho$ is
	\begin{align*}
	G_{\real} &= \rho^{-1} \bigl( \rho (G_{\complex}) \cap \SL(N,\real) \bigr)  
	= \{\, g \in G_{\complex} \mid \widetilde\sigma(g) = g  \,\} 
	. \end{align*}
To summarize, the obvious $\real$-form of $G_{\complex}$ is the set of fixed points of the usual complex conjugation, and any other $\real$-form is the set of fixed points of some other automorphism of~$G_\complex$.

Now let
	\begin{align} \label{CocAlphaRFormDefn}
	 \alpha(\sigma) = \widetilde\sigma \, \sigma^{-1} \colon G_{\complex} \to G_{\complex} 
	 . \end{align}
It is not difficult to see that 
	\begin{itemize}
	\item $\alpha(\sigma)$ is an automorphism of~$G_{\complex}$ (as an abstract group),
	and
	\item $\alpha(\sigma)$ is holomorphic (since $\rho^{-1}$ and $\sigma \rho \sigma^{-1}$ are holomorphic --- in fact, they can be represented by polynomials in local coordinates).
	\end{itemize}
So $\alpha(\sigma) \in \Aut (G_{\complex})$. Thus, by defining $\alpha(1)$ to be the trivial automorphism, we obtain a function $\alpha \colon \Gal(\complex/\real) \to \Aut ( G_{\complex} )$.

Let $\Gal(\complex/\real)$ act on $\Aut ( G_{\complex} )$, by defining 
	$$ \text{${}^\sigma \! \varphi = \sigma \varphi \sigma^{-1}$ \ for $\varphi \in \Aut ( G_{\complex} )$} .$$
Then $\alpha(\sigma) = \varphi^{-1} \, {}^\sigma \!\varphi$, so 
$\alpha(\sigma) \cdot {}^\sigma \! \alpha(\sigma) = \alpha(1)$ (since $\sigma^2 = 1$). This means that $\alpha$ is $1$-cocycle of group cohomology, and therefore defines an element of the cohomology set $\coho1\bigl( \Gal(\complex/\real), \Aut ( G_{\complex} ) \bigr)$. 
%In general, the cohomology of a Galois group is called a \emph{Galois cohomology group}, so each $\real$-form of~$G$ defines a Galois cohomology class. 
In fact:
	\begin{align} \label{RFormsG=Coho}
	 \hbox{\em$\begin{matrix}
	\text{This construction provides a one-to-one correspondence}
	\\\text{between $\coho1\bigl( \Gal(\complex/\real), \Aut ( G_{\complex} ) \bigr)$ and the set of\/ $\real$-forms of~$G_{\complex}$} \end{matrix}$ } 
	\end{align}
\csee{RFormsG=CohoEx}.
Thus, finding all of the $\real$-forms of~$G_{\complex}$ amounts to calculating the cohomology of a Galois group, or, in other words, ``Galois cohomology\zz.''

\begin{obs} \label{H1=Rforms}
\normalfont The above discussion is an example of a fairly general principle: if $X$ is an algebraic object that is defined over~$\real$, then $\coho1\bigl( \Gal(\complex/\real), \Aut(X_\complex) \bigr)$ is in one-to-one correspondence with the set of $\real$-isomorphism classes of $\real$-defined objects whose $\complex$-points are isomorphic to~$X_\complex$.
\end{obs}

\begin{eg} \label{H1GL=0} \label{H1SL=0}
Suppose $V_1$ and~$V_2$ are two vector spaces over~$\real$, and they are isomorphic over~$\complex$.  (I.e., $V_1 \otimes \complex \iso V_2 \otimes \complex$.) Then the two vector spaces have the same dimension, so elementary linear algebra tells us that they are isomorphic over~$\real$. This means that the $\real$-form of any complex vector space~$V_\complex$ is unique (up to isomorphism), so the general principle \pref{H1=Rforms} tells us
	$$\coho1\bigl( \Gal(\complex/\real), \Aut(V_\complex) \bigr) = 0 .$$
In other words, we have 
	$$\coho1\bigl( \Gal(\complex/\real), \GL(n,\complex) \bigr) = 0 .$$
A similar argument shows $\coho1\bigl( \Gal(\complex/\real),  \SL(n,\complex) \bigr) = 0$ \csee{H1SL=0Ex}.
\end{eg}

\begin{warn} \label{NotGeneralWarn}
The ``fairly general principle'' \pref{H1=Rforms} is not completely general. Although almost nothing needs to be assumed in order to construct a well-defined, injective map from the set of $\rational$-forms to the cohomology set \ccf{RFormsG=CohoEx}, this map might not be surjective. That is, there might be cohomology classes that do not come from $\rational$-forms, unless some (fairly mild) hypotheses are imposed on the class of algebraic objects.
\end{warn}

\begin{warn}
Up to now, we have usually ignored finite groups in this book: an answer up to isogeny or commensurability was good enough.  However, such sloppiness is unacceptable when calculating Galois cohomology groups. For example, even though $\SL(n,\complex)$ is isogenous to $\PSL(n,\complex)$, the two groups have completely different cohomology. Namely:
	\begin{itemize}
	\item we saw in \cref{H1SL=0} that $\coho1\bigl( \Gal(\complex/\real),  \SL(n,\complex) \bigr)$ is trivial,
	but
	\item \cref{RFormFromCoho} will show $\coho1\bigl( \Gal(\complex/\real),  \PSL(n,\complex) \bigr)$ is infinite.
	\end{itemize}
\end{warn}

\subsection{Constructing explicit $\real$-forms from cohomology classes} \label{RFormFromCoho}
Given $\alpha \in \coho1\bigl( \Gal(\complex/\real), \Aut ( G_{\complex} ) \bigr)$, we will now see how to find the corresponding $\real$-form~$G_\real$.

It is known that the \textbf{outer} automorphism group of $G_\complex = \SL(n,\complex)$ has only one nontrivial element, namely, the ``transpose-inverse'' automorphism, defined by $\omega(g) = (g^\transpose)^{-1}$. So
	$$ \Aut  \bigl( G_{\complex} \bigr) = \PSL(n,\complex) \rtimes \langle \omega \rangle .$$
We consider two cases.

\setcounter{case}{0}

\begin{case} \label{RFormFromCoho-innerCase}
Assume $\alpha \in \coho1\bigl( \Gal(\complex/\real), \PSL(n,\complex) \bigr)$. 
\end{case}
%We have a short exact sequence 
%	$$0 \to \mu_n \to G_{\complex} \to \PSL(n,\complex) \to 0$$
%where $\mu_n = Z \bigl( G_{\complex} \bigr)$ is the group of $n$th roots of unity in~$\complex$. This yields a long exact sequence of cohomology \csee{LongExactCoho}:
%	$$ \coho1\bigl( \Gal(\complex/\real), G_{\complex} \bigr)
%		\to \coho1\bigl( \Gal(\complex/\real), \PSL(n,\complex) \bigr)
%		\stackrel{\delta}{\to} \coho2\bigl( \Gal(\complex/\real), \mu_n \bigr) .$$
%
%We know that the first term in this sequence is~$0$ \csee{H1SL=0}, so we will focus on the last term. 
It is a fundamental fact in the theory of finite-dimensional algebras that every $\complex$-linear automorphism of the matrix algebra $\Mat_{n \times n}(\complex)$ is inner \csee{AutMatInnerEx}. Since the center acts trivially, this means $\Aut \bigl( \Mat_{n \times n}(\complex) \bigr) = \PSL(n,\complex)$. Therefore,  
	$$\coho1\bigl( \Gal(\complex/\real), \PSL(n,\complex) \bigr) = \coho1\bigl( \Gal(\complex/\real), \Aut \bigl( \Mat_{n \times n}(\complex) \bigr) \bigr) ,$$
	so, by the general principle \pref{H1=Rforms}, we can identify this cohomology set with the set of $\real$-forms of $\Mat_{n \times n}(\complex)$.  More precisely, it is the set of algebras $A$ over~$\real$, such that $A \otimes \complex \iso \Mat_{n \times n}(\complex)$. Such an algebra must be simple (since $\Mat_{n \times n}(\complex)$ is simple), so, by {Wedderburn's Theorem} \pref{WedderburnThm}, it is a matrix algebra over a division algebra: $A \iso \Mat_k (D)$, where $D$~is a division algebra over~$\real$. 
The corresponding $\real$-form~$G_\real$ is $\SL(k, D)$. It is well known that the only division algebras over~$\real$ are $\real$, $\complex$, and~$\quaternion$ \csee{DivAlg/R}, so the real form must be either $\SL(k,\real)$, $\SL(k,\complex)$, or $\SL(k,\quaternion)$, all of which are on the list in \fullcref{GxC=}{SL}.

\begin{case} \label{RFormSLOuterCase}
Assume the image of~$\alpha$ is \textbf{not} contained in $\PSL(n,\complex)$.
\end{case}
%Since the outer automorphism group $\Out ( G_{\complex} )$ is of order~$2$, it has no nontrivial automorphisms. Therefore, the action of the Galois group $\Gal(\complex/\real)$ on $\Out ( G_{\complex} )$ must be trivial, so if we let $\overline\alpha \colon \Gal(\complex/\real) \to \Out (G_{\complex} )$ be the  $1$-cocycle obtained from~$\alpha$ by modding out $\PSL(n,\complex)$, then $\overline\alpha$ is an actual homomorphism (not merely a ``crossed homomorphism''). 
%%The assumption of this \lcnamecref{RFormSLOuterCase} implies that the domain and range of~$\overline\alpha$ both have order~$2$, so the kernel of~$\overline\alpha$ is trivial.
%
In this case, we have $\alpha(\sigma) = (\text{conjugation~by~$A$}) \, \omega$ for some $A \in \GL(n,\real)$. Hence, for every $g \in G_\real$, 
	$$g = \widetilde\sigma(g) 
	=  \bigl( \alpha(\sigma) \, \sigma\bigr) (g) 
	= A \, \omega \bigl( \sigma(g) \bigr) \, A^{-1}
	= A \, \bigl( ({}^{\sigma} \! g)^\transpose \bigr) ^{-1} \, A^{-1} ,$$
which means $g \, A \, ({}^{\sigma} \! g)^\transpose = A$. In other words, $g$ is in the unitary group  $\SU( A, \sigma)$ corresponding to the Hermitian form on~$\complex^n$ that is defined by the matrix~$A$. Since every Hermitian form on~$\complex^n$ is determined (up to isometry) by the number of positive and negative eigenvalues of~$A$, we conclude that $G_{\real} \iso \SU(m,n)$ for some $m$ and~$n$. So $G_\real$ is listed in \fullcref{GxC=}{SL}.

\begin{exercises}

%\item \label{LongExactCohoEx}
%Prove \cref{LongExactCoho}.
%\hint{For a $1$-cocycle $\alpha \colon X \to M/Z$, let $\delta \alpha$ be the cohomology class represented by $c_\alpha(x,y) = \alpha(xy) - \alpha(x) - {}^x \! \alpha(y)$.}

\item \label{RFormsG=CohoEx}
Suppose $\rho_1(G_\complex)_\real$ and $\rho_2(G_\complex)_\real$ are two $\real$-forms of $\SL(n,\complex)$, with corresponding $1$-cocycles $\alpha_1$ and~$\alpha_2$.
	\begin{enumerate}
	\item \label{RFormsG=CohoEx-well}
	Show that if $\rho_1(G_\complex)_\real \iso \rho_2(G_\complex)_\real$, then $\alpha_1$ and~$\alpha_2$ are cohomologous.
	(So the correspondence in \pref{RFormsG=Coho} is well-defined.) 
	\item Conversely, show that if $\alpha_1$ is cohomologous to~$\alpha_2$, then we have $\rho_1(G_\complex)_\real \iso \rho_2(G_\complex)_\real$.
	(So the correspondence in \pref{RFormsG=Coho} is one-to-one.) 
	\end{enumerate}
In \cref{RFormFromCoho}, a real form of $\SL(n,\complex)$ is constructed for each cohomology class~$\alpha$. This shows that the correspondence is onto, and therefore completes the proof of \pref{RFormsG=Coho}.
\hint{In \pref{RFormsG=CohoEx-well}, you may assume, without proof, that every isomorphism from $\rho_1(G_\complex)_\real$ to $\rho_2(G_\complex)_\real$ extends to an isomorphism from $\rho_1(G_\complex)$ to $\rho_2(G_\complex)$.}

\item \label{H1SL=0Ex}
Show $\coho1\bigl( \Gal(\complex/\real),  \SL(n,\complex) \bigr) = 0$, by identifying $\SL(n,\complex)$ with the automorphism group of a pair $(V,\xi)$, where $V$ is an $n$-dimensional vector space and $\xi$~is a nonzero element of the exterior power ${\bigwedge}\!^n \, V$.

\item The short exact sequence
	$$ 1 \to \SL(n,\complex) \hookrightarrow \GL(n,\complex) \stackrel{\det}{\longrightarrow} \complex^\times \to 1 $$
gives rise to the following long exact sequence of cohomology:
	\begin{align*}
	&\coho0\bigl( \Gal(\complex/\real), \GL(n,\complex) \bigr)
	\to \coho0\bigl( \Gal(\complex/\real),  \complex^\times \bigr)
	\\&\qquad \to \coho1\bigl( \Gal(\complex/\real),  \SL(n,\complex) \bigr)
	\to \coho1\bigl( \Gal(\complex/\real),  \GL(n,\complex) \bigr) 
	. \end{align*}
Show that the first map in this sequence is surjective, and combine this with the vanishing of the last term to provide another proof that $\coho1\bigl( \Gal(\complex/\real),  \SL(n,\complex) \bigr) = 0$.
\hint{The $0$th cohomology group is the set of fixed points of the action.}

\item Show that if $n$ is odd, then every $\real$-form of $\SO(n,\complex)$ is isogenous to $\SO(p,q)$, for some $p$ and~$q$.
\hint{You may assume, without proof, that every automorphism of $\SO(n,\complex)$ is inner. Also note that $\SO(n,\complex) = \PSO(n,\complex)$ (why?). Both of these observations require the assumption that $n$ is odd.}

%\item Show that every $\real$-form of $\SO(2n,\complex)$ is isogenous either to $\SO(n,\quaternion)$ or to $\SO(p,q)$, for some $p$ and~$q$.
%\hint{You may assume, without proof, that every outer automorphism is conjugation by an orthogonal matrix of determinant~$-1$.}

\item \label{AutMatInnerEx}
Show that if $\alpha$ is any $\complex$-linear automorphism of the ring $\Mat_{n\times n}(\complex)$, then there exists $T \in \GL(n,\complex)$, such that $\varphi(X) = T X T^{-1}$ for all $X \in \Mat_{n\times n}(\complex)$.
\hint{For $A = \Mat_{n\times n}(\complex)$, make $\complex^n$ into a simple $A$-module via $a * v = \alpha(a) v$. However, the usual action on $\complex^n$ is the unique simple $A$-module (up to isomorphism), because $A$~is a direct sum of submodules that are isomorphic to~$\complex^n$.}

\item \label{DivAlg/R}
Show:
	\begin{enumerate}
	\item \label{DivAlg/R-C}
	$\complex$ is the only finite field extension of~$\real$ (other than $\real$ itself).
	\item \label{DivAlg/R-H}
	$\quaternion$ is the only division algebra over~$\real$ that is not commutative.
	\end{enumerate}
\hint{\pref{DivAlg/R-C}~You may assume, without proof, that $\complex$ is algebraically closed. This implies that every irreducible real polynomial is either linear or quadratic.
	\pref{DivAlg/R-H}~If $x \in D \smallsetminus \real$, then $\real[x]$ is a field extension of~$\real$; identify it with~$\complex$. Then conjugation by~$i$ is a $\complex$-linear map on~$D$. Choose $j$ to be in the $-1$-eigenspace, and let $b = j^2$. Show $b \in \real$ and $D \iso \quaternion_\real^{-1,b}$.}

\end{exercises}

\section{The \texorpdfstring{$\rational$}{Q}-forms of 
	\texorpdfstring{$\SL(\lowercase{n},\real)$}{SL(n,R)}} 
\label{QFormsOfSLnSect}

To illustrate how the method of the preceding section is used to find $\rational$-forms, instead of $\real$-forms, we prove the following result that justifies the claims made in \cref{EgArithGrpsChap}
%\cref{ArithLattSL2,NonLattinSL3Sect,CocpctLattSL3R,LattSlnRSect}
 about arithmetic subgroups of $\SL(n,\real)$:

\begin{thm}[\ccf{LattSlnRSect}] \label{QformsOfSLn}
Every\/ $\rational$-form $G_\rational$ of\/ $\SL(n,\real)$ is either a special linear group or a unitary group\/ \textup(perhaps over a division algebra\/\textup).
\end{thm}

\begin{rem}
More precisely, $G_\rational$ is isomorphic to either:
	\begin{enumerate}
	\item $\SL(m,D)$, for some $m$ and some division algebra~$D$ over~$\rational$,
	or
	\item $\SU(A, \tau ; D) = 
		 \bigset{ g \in \SL ( k, D ) }{ g A ({}^\tau \!\! g)^\transpose = A }$, where
		\begin{itemize} 
		\item $D$ is a division algebra over~$\rational$,
		\item $\tau$ is an anti-involution of~$D$ that acts nontrivially on the center of~$D$,
		and
		\item $A$ is a matrix in $\Mat_{k \times k}(D)$ that is Hermitian (i.e., $({}^\tau \!\! A)^\transpose = A$).
		\end{itemize}
	\end{enumerate}
\end{rem}

The proof is based on the following connection with Galois cohomology. 
We will work with $G_\complex$, instead of~$G$, because algebraically closed fields are much more amenable to Galois Theory. (That is, we are replacing $\SL(n,\real)$ with $\SL(n,\complex)$ to avoid technical issues.)

\begin{prop} \label{GQ<>H1}
There is a one-to-one correspondence between the\/ $\rational$-forms of $G_\complex$ and the Galois cohomology set\/ $\coho1\bigl( \Gal(\complex/\rational), \Aut(G_\complex) \bigr)$.
\end{prop}

\begin{proof}
We assume familiarity with the proof in \cref{GaloisCohoRealFormsSect}, and highlight the changes that need to be made.

Suppose we have an embedding $\rho \colon G_\complex \to \SL(N,\complex)$, such that $\rho(G_\complex)$ is defined over~$\rational$. 
%
%The first difference is a technical issue: Let us assume (without proof) that $\rho$ extends to $\rho_\complex \colon \SL(n,\complex) \stackrel{\iso}{\longrightarrow} \rho(G)_\complex$. (At the Lie algebra level, this is easy: any isomorphism $\Lie G \to \Lie G'$ extends to an isomorphism from $\Lie G \otimes \complex$ to $\Lie G' \otimes \complex$.)
%
The main difference from \cref{GaloisCohoRealFormsSect} is that, unlike $\Gal(\complex/\real)$, the Galois group $\Gal(\complex/\rational)$ has infinitely many nontrivial elements, and we need to consider all of them: since
	$$ \rational = \bigset{ z \in \complex }{ \sigma(z) = z, \ \forall \sigma \in \Gal(\complex/\rational) } \! , $$
we have
%	$$ \SL(N,\rational) = \bigset{ g \in \SL(N,\complex) }{ \text{$\sigma(g) = g$ \ for all $\sigma \in \Gal(\complex/\rational)$} } ,$$
%so
	$$ \rho(G_\complex)_\rational 
%	= \rho_\complex (G_\complex) \cap \SL(N,\rational)
	= \bigset{ g \in \rho_\complex(G_\complex) }{ \sigma(g) = g, \ \forall \sigma \in \Gal(\complex/\rational) } \! .$$
For each $\sigma \in \Gal(\complex/\rational)$, let 
	$$ \text{$\widetilde\sigma = \rho^{-1} \sigma \rho \colon G_{\complex} \to G_{\complex}$
\ and \ $ \alpha(\sigma) = \widetilde\sigma \, \sigma^{-1} \colon G_{\complex} \to G_{\complex}$} .$$
Then 
	$$ G_{\rational} = \{\, g \in G_{\complex} \mid \widetilde\sigma(g) = g, \ \forall\sigma \in \Gal(\complex/\rational) \,\} ,$$
and $\alpha(\sigma) \in \Aut(G_{\complex})$. 
Furthermore, since 
	$\alpha(\sigma) = \rho^{-1} \sigma \rho \sigma^{-1} = \rho^{-1} \, {}^\sigma \! \!\rho$
is formally a $1$-coboundary,  
%	\begin{align*}
%	 \alpha({\sigma \tau})
%	&= \rho^{-1} \, (\sigma\tau) \, \rho\, (\sigma\tau)^{-1}
%	\\&= \rho^{-1} \sigma\tau  \rho  \tau^{-1}  \sigma^{-1}
%	\\&= (\rho^{-1} \sigma  \rho \sigma^{-1} ) \cdot \sigma ( \rho^{-1} \tau  \rho  \tau^{-1})  \sigma^{-1}
%	\\&= \alpha(\sigma) \, {}^\sigma \! \alpha(\tau)
%	. \end{align*}
it is easily seen to be a $1$-cocycle, and therefore represents a cohomology class in $\coho1\bigl( \Gal(\complex/\rational), \Aut(G_{\complex}) \bigr)$. 

This defines the desired map from the set of $\rational$-forms to the Galois cohomology set. It can be proved to be well-defined and injective by replacing $\real$ with~$\rational$ in \cref{RFormsG=CohoEx}. That the map is surjective will be established in the proof of \cref{QformsOfSLn} below, % is it still below !!!
where we explicitly describe the $\rational$-form corresponding to each cohomology class.
\end{proof}

More generally, we have the following natural analogue of \cref{H1=Rforms}:

\begin{obs} \label{H1=Qforms}
\normalfont If $X$ is an algebraic object that is defined over~$\rational$ (and satisfies mild hypotheses; cf.\ \cref{NotGeneralWarn}), then the Galois cohomology set $\coho1\bigl( \Gal(\complex/\rational), \Aut(X_\complex) \bigr)$ is in one-to-one correspondence with the set of $\rational$-isomorphism classes of $\rational$-defined objects whose $\complex$-points are isomorphic to~$X_\complex$.
\end{obs}

\begin{cor}[\ccf{H1SL=0}]
$$ \text{$\coho1\bigl( \Gal(\complex/\rational), \GL(n,\complex) \bigr) = 0$
\ and \  $\coho1\bigl( \Gal(\complex/\rational), \SL(n,\complex) \bigr) = 0$.} $$
\end{cor}

\begin{proof}[Proof of \cref{QformsOfSLn}]
Let $G_\complex = \SL(n,\complex)$. As in \cref{RFormFromCoho}, we have 
	$$ \Aut (G_{\complex}) = \PSL(n,\complex) \rtimes \langle \omega \rangle ,$$
where $\omega(g) = (g^\transpose)^{-1}$. 
Given $\alpha \in \coho1\bigl( \Gal(\complex/\rational), \Aut(G_{\complex}) \bigr)$, corresponding to a $\rational$-form~$G_\rational$, we consider two cases.

\setcounter{case}{0}

\begin{case} \label{SLQforms-InnerCase}
Assume $\alpha \in \coho1\bigl( \Gal(\complex/\rational), \PSL(n,\complex) \bigr)$. 
\end{case}
By arguing exactly as in \cref{RFormFromCoho-innerCase} of \cref{RFormFromCoho} (but with $\rational$ in the place of~$\real$), we see that $G_\rational \iso \SL(k, D)$, for some~$k$ and some division algebra~$D$ over~$\rational$. 

\begin{case} \label{SLQforms-OuterCase}
Assume the image of~$\alpha$ is \textbf{not} contained in\/ $\PSL(n,\complex)$.
\end{case}
Since the outer automorphism group $\Out ( G_{\complex} )$ is of order~$2$, it has no nontrivial automorphisms. Therefore, the action of the Galois group $\Gal(\complex/\rational)$ on $\Out ( G_{\complex} )$ must be trivial. Hence, if we let $\overline\alpha \colon \Gal(\complex/\rational) \to \Out (G_{\complex} )$ be the  $1$-cocycle obtained from~$\alpha$ by modding out $\PSL(n,\complex)$, then $\overline\alpha$ is an actual homomorphism (not merely a ``crossed homomorphism''). 
%The assumption of this \lcnamecref{RFormSLOuterCase} implies that the domain and range of~$\overline\alpha$ both have order~$2$, so the kernel of~$\overline\alpha$ is trivial.

By the assumption of this \lcnamecref{SLQforms-OuterCase} (and the fact that $| {\Out(G_{\complex})} | = 2$), the kernel of~$\overline\alpha$ is a subgroup of index~$2$ in $\Gal(\complex/\rational)$. This means that the fixed field of $\ker \overline\alpha$ is a quadratic extension $L = \rational \bigl[ \!\sqrt{r}\bigr]$ of~$\rational$. Then, by construction, we have
	$ \Gal ( \complex / L ) = \ker \overline\alpha $.

For any $\sigma \in \Gal(\complex/L)$, the conclusion of the preceding paragraph tells us that $\overline\alpha(\sigma)$ is trivial. For simplicity, let us assume that the bar can be removed, so $\alpha(\sigma)$ is trivial \fullcsee{QFormCorrections}{mixed}. Since, by definition, we have $ \alpha(\sigma) = \widetilde\sigma \, \sigma^{-1}$ (see \pref{CocAlphaRFormDefn}), this implies $\sigma = \widetilde\sigma$. Therefore, for any $g \in G_\rational$, we have
	$ g^\sigma = g^{\widetilde\sigma} = g $.
Since this holds for all $\sigma \in \Gal(\complex/L)$, we conclude that $g \in \SL( n, L)$. 

Now, for the unique nontrivial $\tau \in \Gal(L/\rational)$, we have $\tau \notin \ker \overline\alpha$, so $\alpha(\tau) = (\text{conj~by~$A$}) \, \omega$ for some $A \in \GL(n,\real)$. Hence,  for any $g \in G_\rational$, we have
	$$g = \widetilde\tau(g) 
	=  \bigl( \alpha(\tau) \, \tau \bigr) (g) 
	= A \, \omega \bigl( \tau(g) \bigr) \, A^{-1}
	= A \, \bigl( ({}^{\tau} \! g)^\transpose \bigr) ^{-1} \, A^{-1} ,$$
so $g \, A \, ({}^{\tau} \! g)^\transpose = A$, which means $g \in \SU( A, \tau;  L)$. Furthermore, the equation $\widetilde\tau^2 = 1$ provides an equation that can be used to show $A$ is Hermitian (or, more precisely, can be chosen to be Hermitian) \csee{SLQforms-OuterCase-AHermEx}.
\end{proof}

\begin{corrections} \label{QFormCorrections} \ 
\noprelistbreak
\begin{enumerate}
\item \label{QFormCorrections-mixed} 
\emph{Mixed case.}
We seem to have shown that all $\rational$-forms of $\SL(n,\real)$ can be constructed from either division algebras (\cref{SLQforms-InnerCase}) or unitary groups (\cref{SLQforms-OuterCase}). However, the discussion in \cref{SLQforms-OuterCase} assumes that $\alpha(\sigma)$ is trivial, when all we actually know is that $\overline\alpha(\sigma)$ is trivial. Removing this assumption means that $\alpha$ can map a part of the Galois group into $\PSL(n,\complex)$. In other words, in addition to the homomorphism~$\overline\alpha$, there is a nontrivial cocycle from $\Gal(\complex/L)$  to $\PSL(n,\complex)$. By the argument of \cref{SLQforms-InnerCase}, this cocycle yields a division algebra~$D$ over~$L$. The resulting $\rational$-form	
	$ G_\rational = \SU( A, \tau ; D)$
is obtained by combining division algebras with unitary groups.

%All non-cocompact (arithmetic) lattices in $\SL(n,\real)$ are one of these, constructed from either unitary groups, or division algebras, or a combination of the two.

\item \label{QFormCorrections-Qbar} 
\emph{$\complex$ vs.\ $\overline{\rational}$.} We should really be using the algebraic closure~$\overline{\rational}$ of~$\rational$, instead of~$\complex$. The Galois cohomology set $\coho1 \bigl( \Gal ( \overline\rational /\rational ), \Aut(G_{\overline\rational}) \bigr)$ is defined to be the natural limit of the sets $\coho1 \bigl( \Gal ( F /\rational ), \Aut(G_{\overline\rational}) \bigr)$, where $F$~ranges over all finite Galois extensions of~$\rational$.

%\item \emph{$\SL(n,\real)$ vs.\ $\SL(n,\complex)$.}
%To discuss Galois cohomology, we replaced~$\real$ with the algebraically closed field~$\complex$. Thus, some of the groups we found might not be $\rational$-forms of $\SL(n,\real)$ (although we know that their complexification is $\SL(3,\complex)$). For example, if $G_\rational = \SU \Bigl( J, \sigma; \rational\bigl[ \sqrt{r} \, \bigr] \Bigr)$, and $r < 0$, then $G_\real$ is $\SU(2,1)$, not $\SL(3,\real)$. 
%	\begin{itemize} \itemsep=\smallskipamount
%	\item In practice, one can determine which of the groups we constructed are $\rational$-forms of $\SL(3,\real)$. 
%	\item Abstractly, $\SL(3,\real)$ is a $\real$-form of $\SL(3,\complex)$, so, by the general principle, it is represented by a cohomology class $\beta \in \coho1\bigl( \Gal(\complex/\real), \Aut( G_{\complex} ) \bigr)$. There is a natural restriction homomorphism 
%		$$r \colon \coho1\bigl( \Gal(\complex/\rational), \Aut(G_{\complex}) \bigr) \to 
%		\coho1\bigl( \Gal(\complex/\real), \Aut(G_{\complex}) \bigr) ,$$ 
%and the $\rational$-forms of $\SL(3,\real)$ are represented by the elements of $r^{-1}(\beta)$.
%	\end{itemize}
\end{enumerate}
\end{corrections}

\begin{exercises}

\item \label{SLQforms-OuterCase-AHermEx}
In \cref{SLQforms-OuterCase} of the proof of \cref{QformsOfSLn}, show that the matrix~$A$ can be chosen to be Hermitian.
\hint{$A$ must be a scalar multiple~$\lambda$ of a Hermitian matrix (since $\widetilde\tau^2 = 1$). Use the fact that $\coho1 \bigl( \Gal(\complex/L) ; \complex^\times \bigr)$ is trivial (why?) to replace $A$ with a scalar multiple of itself that makes $\lambda = 1$.}

\end{exercises}

\section{\texorpdfstring{$\rational$}{Q}-forms of classical groups}
\label{QFormClassicalSect}

By arguments similar to the ones applied to $\SL(n,\real)$ in \cref{QFormsOfSLnSect}, it can be shown that the $\rational$-forms of almost any classical group come from special linear groups, unitary groups, orthogonal groups, or symplectic groups. However, the special linear groups and unitary groups may involve division algebras, and restriction of scalars \pref{ResScal->Latt} implies that the groups may be over an extension~$F$ of~$\rational$. 
(Recall that unitary groups over division algebras were defined in \cref{SUDDefn}, and the involutions $\tau_c$ and~$\tau_r$ on the quaternion algebra $\quaternion_F^{a,b}$ were defined in \cref{QuatConjRevDefn}.)
Here is a list of the groups that arise:

\begin{defn} \label{FClassicalDefn}
For any algebraic number field~$F$, and any~$n$, the following groups are said to be of \defit[classical!type]{classical type}:
	 \begin{enumerate}

	 \item \label{FClassicalDefn-SL}
	 $\SL(n,D)$, where $D$ is a division algebra whose center is~$F$.

	 \item \label{FClassicalDefn-Sp}
	 $\Sp(2n,F)$.

	 \item \label{FClassicalDefn-SO}
	 $\SO(A;F)$, where $A$ is an invertible, symmetric matrix 
	 in $\Mat_{n \times n}(F)$.

	 \item \label{FClassicalDefn-SUSp}
	 $\SU(A, \tau_c; \quaternion_F^{a,b})$, where $\quaternion_F^{a,b}$ is a quaternion division algebra over~$F$, and $A$ is an invertible, $\tau_c$-Hermitian matrix in  $\Mat_{n \times n}(\quaternion_F^{a,b})$.

	 \item \label{FClassicalDefn-SUSO}
	 $\SU(A, \tau_r; \quaternion_F^{a,b})$, where $\quaternion_F^{a,b}$ is a quaternion division algebra, and $A$ is an invertible, $\tau_r$-Hermitian matrix in  $\Mat_{n \times n}(\quaternion_F^{a,b})$.

	 \item \label{FClassicalDefn-SUSL}
	 $\SU(A,\tau;D)$, where 
	 	\begin{itemize}
		\item $D$ is a division algebra whose center is a quadratic extension~$L$ of~$F$, 
		\item $\tau$ is an anti-involution whose restriction to~$L$ is the Galois automorphism of~$L$ over~$F$, 
		and 
		\item $A$ is an invertible, $\tau$-Hermitian matrix in  $\Mat_{n \times n}(D)$.
		\end{itemize}
	 \end{enumerate}
\end{defn}

\begin{rem}
 \Cref{FClassicalDefn} is directly analogous to
the list of classical simple Lie groups
\csee{classical-fulllinear,classical-orthogonal}.
Specifically:
 \begin{enumerate}
 \item[\ref{FClassicalDefn-SL})]
 $\SL(n,D)$ is the analogue of $\SL(n,\real)$,
$\SL(n,\complex)$, and $\SL(n,\quaternion)$.
 \item[\ref{FClassicalDefn-Sp})]
 $\Sp(2n,F)$ is the analogue of $\Sp(2n,\real)$ and
$\Sp(2n,\complex)$.
 \item[\ref{FClassicalDefn-SO})]
 $\SO(A;F)$ is the analogue of $\SO(m,n)$ and
$\SO(n,\complex)$.
 \item[\ref{FClassicalDefn-SUSp})]
 $\SU(A, \tau_c; \quaternion_F^{a,b})$ is the analogue of $\Sp(m,n)$.
 \item[\ref{FClassicalDefn-SUSO})]
 $\SU(A, \tau_r; \quaternion_F^{a,b})$ is the analogue of $\SO(n,\quaternion)$.
 \item[\ref{FClassicalDefn-SUSL})]
 $\SU(A,\tau;D)$ (with $\tau$ nontrivial on the center) is
the analogue of $\SU(m,n)$.
 \end{enumerate}
 \end{rem}

\begin{thm} \label{ArithLattsAreClassical}
Suppose 
\noprelistbreak
	\begin{itemize}
	\item $G$ is classical, 
	and 
	\item no simple factor of $G_\complex$ is isogenous to\/ $\SO(8,\complex)$.
	\end{itemize}
Then every irreducible, arithmetic lattice in~$G$ is commensurable to the integer points of some group\/ \textup(of classical type\/\textup) listed in \cref{FClassicalDefn}.
\end{thm}

\begin{rem}
To state the conclusion of \cref{ArithLattsAreClassical} more explicitly, let us assume, for simplicity, that the center of~$G$ is trivial. Then \cref{ArithLattsAreClassical} states that there exist:
	\begin{itemize}
	\item algebraic number field~$F$, with places~$S^\infty$ and ring of integers~$\ints$,
	\item a group $\widehat G_F$ listed in \cref{FClassicalDefn}, with corresponding semisimple Lie group~$\widehat G$ that is defined over~$F$,  
	and
	\item a homomorphism $\varphi \colon \prod_{\sigma \in S^\infty} \widehat G^\sigma \to G$, with compact kernel,
	\end{itemize}
such that $\varphi \bigl( \Delta(G_{\ints}) \bigr)$ is commensurable to~$\Gamma$ \ccf{ResScal->Latt}.
\end{rem}

\begin{warn}
Although $\varphi \bigl( \Delta(G_{\ints}) \bigr)$ is commensurable to~$\Gamma$, this does \textbf{not} imply that $\varphi \bigl( \Delta(G_F) \bigr)$ is
commensurable to $G_\rational$. For example, the image of
$\SL(2,\rational)$ in $\PSL(2,\rational)$ has infinite index
\ccf{Comm(SL2Z)}.
 \end{warn}

Each of the groups in \cref{FClassicalDefn} has a corresponding semisimple Lie group~$G$ that is defined over~$F$. Before determining which Lie group corresponds to each $F$-group, we first find the complexification of~$G$. This is
similar to calculations that we have already seen,
so we omit the details.

\begin{prop}[\ccf{CalcGxCSect}] \label{GFxC}
 The notation of each part of this \lcnamecref{GFxC} is taken from the corresponding part of \cref{FClassicalDefn}. 
 We use $d$~to denote the
degree of the central division algebra~$D$, and the matrix~$A$ is assumed to be $n \times n$.
 \begin{enumerate}
 \item[\ref{FClassicalDefn-SL})]
 $\SL(n,D \otimes_F \complex) \iso \SL(dn,\complex)$.
 \item[\ref{FClassicalDefn-Sp})]
 $\Sp(2n, \complex) = \Sp(2n,\complex)$ \textup(obviously!\/\textup).
 \item[\ref{FClassicalDefn-SO})]
 $\SO(A;\complex) \iso \SO(n,\complex)$.
 \item[\ref{FClassicalDefn-SUSp})]
 $\SU(A,\tau_c; \quaternion^{a,b}_F \otimes_F \complex) \iso \Sp(2n,\complex)$.
 \item[\ref{FClassicalDefn-SUSO})]
 $\SU(A,\tau_r; \quaternion^{a,b}_F \otimes_F \complex) \iso \SO(2n,\complex)$.
 \item[\ref{FClassicalDefn-SUSL})]
 $\SU(A,\tau; D \otimes_F \complex) \iso \SL(dn,\complex)$.
 \end{enumerate}
 \end{prop}

If $F \notsubset \real$, then the semisimple Lie group~$G$ corresponding to $G_F$ is the complex Lie group in the corresponding line of the above % !!!
\lcnamecref{GFxC}. However, if $F \subset \real$, then $G$ is some $\real$-form of that complex group. The following result lists the correct $\real$-form for each of the groups of classical type.

\begin{prop} \label{GFxR}
 The notation of each part of this
proposition is taken from the corresponding part of
\cref{FClassicalDefn}. We use $d$~to denote the
degree of the central division algebra~$D$, and the matrix~$A$ is assumed to be $n \times n$.

 Assume $F$ is an algebraic number field, and that $F \subset \real$. Then:

 \begin{enumerate} \itemsep=\medskipamount
 
 \item[\ref{FClassicalDefn-SL})]
 $\SL(n,D \otimes_F \real) \iso
 \begin{cases}
 \SL(dn,\real) & \text{if $D  \otimes_F \real \iso \Mat_{d \times d}(\real)$}, \\
 \SL(dn/2,\quaternion) & \text{otherwise} .
 \end{cases}$
 
 \item[\ref{FClassicalDefn-Sp})]
 $\Sp(2n, \real) = \Sp(2n,\real)$ \textup(obviously!\/\textup).

 \item[\ref{FClassicalDefn-SO})]
 $\SO(A,\real) \iso \SO(p,n-p)$.
 
 \item[\ref{FClassicalDefn-SUSp})] \label{GFxR-SUSp}
 $\SU(A,\tau_c; \quaternion^{a,b}_F \otimes_F \real) \iso 
 \begin{cases}
 \Sp(2n,\real) & \text{if\/ $\quaternion^{a,b}_\real \iso \Mat_{2 \times 2}(\real)$}, \\
 %\SU(A,\tau_c;\quaternion) \iso 
 \Sp(p,n-p) & \text{if\/ $\quaternion^{a,b}_\real \iso \quaternion$} .
 \end{cases}$
 
 \item[\ref{FClassicalDefn-SUSO})] 
 $\SU(A,\tau_r; \quaternion^{a,b}_F \otimes_F \real) \iso 
 \begin{cases}
 \SO(p,2n-p) & \text{if\/ $\quaternion^{a,b}_\real \iso \Mat_{2 \times 2}(\real)$}, \\
 \SO(n,\quaternion) & \text{if\/ $\quaternion^{a,b}_\real \iso \quaternion$} .
 \end{cases}$
 
 \item[\ref{FClassicalDefn-SUSL})]
 $\SU(A,\tau;D \otimes_F \real) \iso 
 \begin{cases}
 \SU(p, dn-p)&  \text{if $L \not\subset \real$} , \\
 \SL(dn,\real) & 
 \begin{matrix} \text{if $L \subset \real$ and} \hfill
 	\\[-4pt] \quad \text{$D \otimes_F \real \iso \Mat_{d \times d}(\real)$} , 
	\end{matrix}\\
 \SL(dn/2,\quaternion) & \text{otherwise} .
 \end{cases}$
 
 \end{enumerate}
 \end{prop}

\begin{rem}
\Cref{GFxR} does not specify the value of~$p$, where it appears. However, it can be calculated for any particular matrix~$A$.
For example, to calculate $p$ in
\pref{FClassicalDefn-SUSL}, note that, because $L
\not\subset \real$, we have 
	$$D \otimes_F \real 
	\iso D \otimes_L \complex 
	\iso \Mat_{d \times d}(\complex) ,$$
so we may think of $A \in \Mat_{n
\times n}(D)$ as a $(dn) \times (dn)$ Hermitian matrix. Then
$p$~is the number of positive eigenvalues of this Hermitian
matrix (and $dn - p$ is the number of negative eigenvalues).
%
%In each of the other cases, it is not difficult to give a
%fairly uniform calculation of~$p$ and~$q$ (cf.\ 
%\cref{Calcpq(SOAF),Calcpq(Sp)}, and
%Propositions~\ref{CocpctSO1nQuat}, \fullref{FormAisDiag}{diag},
%and~\fullref{DFormAisDiag}{FormAisDiag-diag}). 
We have already seen this type of consideration in \cref{QuatEigs,CocpctSO1nQuat}.
\end{rem}

 \begin{figure}
\begin{center}  \tabcolsep=4pt
\medskip
 \begin{tabular}{|c|c|c|c|c|c|c|c|c|}
 \noalign{\hrule} 
$\begin{matrix} \text{Lie group $G$} \end{matrix}$ 
& $\begin{matrix} \text{$F$-form $G_F$} \end{matrix} $
& reference 
& $m$ or $p + q$ 
& $\Qrank \Gamma$ \\

\betweengrpseparator

$\SL(m,\real)$
 & $\SL(n, D)$
 &
 \twoline
 {\pref{FClassicalDefn-SL}, $F \subset \real$,}
 {$D$ split/$\real$}
 & $m = dn$
 & $n-1$ \\

\withingrpseparator

 & $\SU(B, \tau; D)$
 &
 \twoline
 {\pref{FClassicalDefn-SUSL}, $F \subset L \subset
\real$,}
 {$D$ split/$\real$}
 & $m = dn$
 & $D$-subspace \\

\betweengrpseparator

\tstrut $\SL(m,\complex)$
 & $\SL(n, D)$
 & \pref{FClassicalDefn-SL}, $F \not\subset \real$
 & $m = dn$
 & $n-1$ \\

\withingrpseparator

& $\SU(B, \tau; D)$
 & 
 \twoline{\pref{FClassicalDefn-SUSL}, $F \not\subset
\real$}
 {(so $L \not\subset \real$)}
 & $m = dn$
 & $D$-subspace \\

\betweengrpseparator

$\SL(m,\quaternion)$
 & $\SL(n, D)$
 &
 \twoline{\pref{FClassicalDefn-SL}, $F \subset \real$,}
 {$D$ not split/$\real$}
 & 
 \twoline
 {$m = dn/2$,}
 {$d$~even}
 & $n-1$
 \\

\withingrpseparator

 & $\SU(B, \tau; D)$
 &
 \twoline{\pref{FClassicalDefn-SUSL}, $F \subset L
\subset \real$,}
 {$D$ not split/$\real$}
 &
 \twoline
 {$m = dn/2$,}
 {$d$~even}
 & $D$-subspace
 \\

\betweengrpseparator

$\SU(p,q)$
 & $\SU(B, \tau; D)$
 &
 \twoline
 {\pref{FClassicalDefn-SUSL}, $F \subset \real$,}
 {$L \not\subset \real$}
 & $p + q = dn$
 & $D$-subspace \\

\betweengrpseparator

\tstrut $\SO(p,q)$
 & $\SO(B; F)$
 &
\pref{FClassicalDefn-SO}, $F \subset \real$
 & $p+q = n$
 & $F$-subspace
\\

\withingrpseparator

 & $\SU(B, \tau_r; D)$
 &
 \twoline
 {\pref{FClassicalDefn-SUSO}, $F \subset \real$,}
 {$D$ split/$\real$}
 & 
 \twoline{$p+q = 2n$,}
 {$d = 2$}
 & $D$-subspace
  \\

\withingrpseparator

 & ?
 &
 \cref{D4weird}
 & 
 $p + q = 8$
 & ?
  \\

\betweengrpseparator

\tstrut $\SO(m,\complex)$
 & $\SO(B; F)$
 & \pref{FClassicalDefn-SO}, $F \not\subset \real$
 & $m = n$
 & $F$-subspace \\

\withingrpseparator

& $\SU(B, \tau_r; D)$
 &
\pref{FClassicalDefn-SUSO}, $F \not\subset \real$
 & 
 \twoline
 {$m = 2n$,}
 {$d = 2$}
 & $D$-subspace
 \\

\withingrpseparator

 & ?
 &
 \cref{D4weird}
 & 
 $m = 8$
 & ?
  \\

\betweengrpseparator

$\SO(m,\quaternion)$
 & $\SU(B, \tau_r; D)$
 & 
 \twoline{\pref{FClassicalDefn-SUSO}, $F \subset \real$,}
 {$D$ not split/$\real$}
 & 
 \twoline
 {$m = n$,}
 {$d = 2$}
 & $D$-subspace \\

\withingrpseparator

 & ?
 &
 \cref{D4weird}
 & 
 $m = 4$
 & ?
  \\

\betweengrpseparator

\tstrut $\Sp(2m,\real)$
 & $\Sp(2n, F)$
 &
\pref{FClassicalDefn-Sp}, $F \subset \real$
 & $m = n$
 & $n$ \\

\withingrpseparator

& $\SU(B, \tau_c; D)$
 &
 \twoline{\pref{FClassicalDefn-SUSp}, $F \subset \real$,}
 {$D$ split/$\real$}
 &
 \twoline
 {$m = n$,}
 {$d = 2$}
 & $D$-subspace
 \\

\betweengrpseparator

\tstrut $\Sp(2m,\complex)$
 & $\Sp(2n, F)$
 & \pref{FClassicalDefn-Sp}, $F \not\subset \real$
 & $m = n$
 & $n$ \\

\withingrpseparator

 & $\SU(B, \tau_c; D)$
 &
\pref{FClassicalDefn-SUSp}, $F \not\subset \real$
 & 
 \twoline
 {$m = n$,} 
 {$d = 2$}
 & $D$-subspace
 \\

\betweengrpseparator

$\Sp(p,q)$
 & $\SU(B, \tau_c; D)$
 &
 \twoline{\pref{FClassicalDefn-SUSp}, $F \subset \real$,}
 {$D$ not split/$\real$}
 &
 \twoline
 {$p+q = n$,}
 {$d = 2$}
 & $D$-subspace
 \\

 \noalign{\hrule}
\end{tabular}
%\caption{The irreducible arithmetic lattices in~$G$, except
%that (as indicated by ``?'') the list is not complete for
%groups of type~$D_4$. The reference is to
%\ref{kFormClassicalThm} and to \ref{GFxC} (for $F_\sigma =
%\complex$) and \ref{GFxR} (for $F_\sigma = \real$).}
\refstepcounter{table}
\label{IrredInG}

{See \cref{IrredInGRem} on \cpageref{IrredInGRem} for an explanation of this table.}
 \end{center}
\vskip-2.01pt % eliminate overfull vbox @@@
 \end{figure}

\begin{rem} \label{IrredInGRem}
The \lcnamecref{IrredInG} % @@@
 on \cpageref{IrredInG} summarizes the above results in a format that makes it easy to find the arithmetic subgroups of any given simple Lie group~$G$ (or, by restriction of scalars, to find the irreducible arithmetic subgroups of any semisimple Lie group that has $G$ as a simple factor), except
that (as indicated by ``?'') the list is not complete for groups whose complexification is isogenous to $\SO(8,\complex)$.

The arithmetic group~$\Gamma$ that corresponds to a given $F$-form~$G_F$ is obtained by:
	\begin{itemize}
	\item replacing $F$ with its ring of integers~$\ints$,
	or
	\item replacing $D$ with an order~$\ints_D$ \csee{OinDivAlg}.
	\end{itemize}
By restriction of scalars \pref{ResScal->Latt}, $\Gamma$~is an arithmetic subgroup of $\prod_{\sigma \in S^\infty} G^\sigma$.

A parenthetical reference indicates the corresponding part of \cref{FClassicalDefn}, and also of  \cref{GFxC} (for $F_\sigma = \complex$) and \cref{GFxR} (for $F_\sigma = \real$).
The reference column (combined with the ``$m$ or $p + q$'' column) also lists additional conditions that determine $G_F \otimes_F F_\infty$ is the desired simple Lie group~$G$ (except that the parameters $p$ and~$q$ will need to be calculated, if they arise).
%As usual, $d$~is the degree of~$D$ over its center (which is either $F$ or~$L$).

The $\rational$-rank of the corresponding arithmetic group~$\Gamma$ is either given explicitly (as a function of~$n$), or  is the dimension of a maximal isotropic subspace (of the associated vector space over either the field~$F$ or the division algebra~$D$, as indicated).
\end{rem}

\begin{rem}[(``triality'')] \label{D4weird} 
Perhaps we should explain why the statement of \cref{ArithLattsAreClassical} assumes no simple factor of $G_\complex$ is isogenous to $\SO(8,\complex)$. Fundamentally, the reason  $\PSO(8,\complex)$ is special is that, unlike all the other simple Lie groups over~$\complex$, it has an
outer automorphism~$\phi$ of order~$3$, called ``\term{triality}\zz.''
For all of the other simple groups, the outer automorphism group is either trivial or has order~$2$.

Here is how the triality automorphism~$\phi$ can be used to construct $\rational$-forms that are not listed in \cref{ArithLattsAreClassical}.
We first choose any homomorphism $\alpha \colon \Gal(\complex/\rational) \stackrel{\text{onto}}{\longrightarrow} \langle \phi \rangle$ (so the kernel of~$\alpha$ is a cubic, Galois extension of~$\rational$).  The triality automorphism can be chosen so that it commutes with the action of the Galois group (in other words, $\phi$ is ``defined over~$\rational$''), so the homomorphism~$\alpha$ is a $1$-cocycle into $\Aut \bigl( \PSO(8,\complex) \bigr)$.
Therefore, by the correspondence between cohomology and $\rational$-forms \pref{GQ<>H1}, there is a corresponding $\rational$-form~$G_\rational$. This $\rational$-form is not any of the groups listed in \cref{ArithLattsAreClassical}, because, for all those groups, the image of the induced homomorphism $\overline\alpha \colon \Gal(\complex/\rational) \to \mathrm{Out}(G_\complex)$ has order $1$ or~$2$, not~$3$.

Mathematicians who understand the triality automorphism can construct the corresponding $\rational$-form explicitly, by reversing the steps in the proof of \cref{GQ<>H1}. Namely, for each $\sigma \in \Gal(\complex/\rational)$, let 
	$$\widetilde\sigma = \alpha(\sigma) \cdot \sigma \in \Aut\bigl( \PSO(8,\complex) \!\bigr) .$$
Then 
	$$G_\rational =  \{\, g \in \PSO(8,\complex) \mid \widetilde\sigma(g) = g, \ \forall\sigma \in \Gal(\complex/\rational) \,\} . $$
 \end{rem}

\section{Applications of the classification of arithmetic groups}

Several results that were stated without proof in previous chapters are easy consequences of the above classification of $F$-forms.

\begin{cor}[\ccf{SO1nCpctListStated}] \label{SO1nCpctList}
  Suppose\/ $\Gamma$ is an arithmetic subgroup of\/ $\SO(m,n)$, and $m+n \ge 5$ is odd. Then there is a finite extension~$F$ of\/~$\rational$, with ring of integers~$\ints$, such that\/ $\Gamma$ is commensurable to\/ $\SO(A; \ints)$, for some invertible, symmetric matrix~$A$ in\/ $\Mat_{n \times n}(F)$.
 \end{cor}
 
 \begin{proof}
 Let $G = \SO(m,n)$.
 Restriction of scalars \pref{simple->Arith=Ints} implies there is a group~$\widehat G$ that is defined over an algebraic number field~$F$ and has a simple factor that is isogenous to~$G$, such that $\Gamma$ is commensurable to~$\widehat G_\ints$.  
 By inspection, we see that a group of the form $\SO(m,n)$ never appears in \cref{GFxC}, and appears at two places in \cref{GFxR}. However, in our situation, we know that $m + n$~is odd, so the only possibility for $\widehat G_F$ is $\SO(A; F)$. Therefore, $\Gamma$ is commensurable to\/ $\SO(A; \ints)$. 
 \end{proof}

\begin{cor}[\fullccf{QrankRems}{gap}] \label{QrankGap}
 If $G = \SO(2,n)$, with $n \ge 5$, and
$n$~is odd, then\/ $\Rrank G = 2$, but there is no
lattice\/~$\Gamma$ in~$G$, such that\/ $\Qrank \Gamma = 1$.
 \end{cor}

\begin{proof}
 We have $\Rrank \SO(2,n) = \min\{2,n\} = 2$ \csee{Rrank(SOmn)}.

From the Margulis Arithmeticity Theorem \pref{MargulisArith}, we know that $\Gamma$ is
arithmetic, so \cref{SO1nCpctList} tells us that $\Gamma$ is of the form $\SO(B; \ints)$, where 
	\begin{itemize}
	\item $\ints$ is the ring of integers of some algebraic number field~$F$,
	and
	\item $B$ is a symmetric bilinear form on~$F^{n+2}$.
	\end{itemize}
If $\Qrank \Gamma = 1$, then $G/\Gamma$  is not compact, so \cref{IrredNoncpct->NoROS} tells us that we may take $F = \rational$; therefore $\ints = \integer$.
We see that:
	 \begin{enumerate}
	 \item $B$ has signature $(2,n)$ on $\real^{n+2}$ (because $G = \SO(2,n)$), 
	 and
	 \item no $2$-dimensional $\rational$-subspace of~$\rational^{n+2}$ is
	totally isotropic (because we have $\Qrank \Gamma < 2$).
	 \end{enumerate}
Recall the following important fact that was used in the proof of \cref{NonCocptArithSOn1}:
	 \begin{itemize} \renewcommand{\labelitemi}{}
	  \it
	 \item \label{Isot/R->Isot/Q}
	 {\bf Meyer's Theorem\thmindex{Meyer's}.}
	 If $B_0(x,y)$ is any nondegenerate, symmetric bilinear form
	on\/ $\real^d$, such that
	 \begin{itemize}
	 \item $B$ is defined over\/~$\rational$, 
	 \item $d \ge 5$, and
	 \item $B_0$ is isotropic over\/~$\real$ \textup(that is, $B(v,v) =
	0$ for some nonzero $v \in \real^d$\textup), 
	 \end{itemize}
	 then $B_0$~is also isotropic over\/~$\rational$ \textup(that is,
	$B(v,v) = 0$ for some nonzero $v \in \rational^d$\textup).
	 \end{itemize}
This theorem, % @@@
tells us there is a
nontrivial isotropic vector $v \in \rational^{n+2}$. Then,
because $B$ is nondegenerate, there is a vector $w \in
\rational^{n+2}$, such that $B(v,w) = 1$ and $B(w,w) = 0$.
Let $V$ be the $\real$-span of $\{ v,w \}$. Because the restriction
of~$B$ to~$V$ is nondegenerate, we have $\real^{n+2} = V
\oplus V^\perp$. This direct sum is obviously orthogonal
(with respect to~$B$), and the restriction of~$B$ to~$V$ has
signature $(1,1)$, so we conclude that the restriction
of~$B$ to~$V^\perp$ has signature $(1,n-1)$. Hence, there is
an isotropic vector in~$V^\perp$.
By applying Meyer's Theorem again, we conclude that there is an
isotropic vector~$z$ in $(V^\perp)_{\rational}$. Then
$\langle v,z \rangle$ is a $2$-dimensional totally isotropic
subspace of~$\rational^{n+2}$. This is a contradiction.
 \end{proof}

\begin{cor}[\ccf{SO1nNotCpctListStated}] \label{SO1nNotCpctList}
 If\/ $n \notin \{3,7\}$, then every noncocompact, arithmetic subgroup of\/ $\SO(1,n)$ is commensurable to a conjugate of\/ $\SO(A; \integer)$, for some invertible, symmetric matrix $A \in \Mat_{n \times n}(\rational)$.
 \end{cor}

\begin{proof}
Assume, for simplicity, that $n = 5$, and let $\Gamma$ be a noncocompact, arithmetic subgroup of $G = \SO(1,5)$.
Since $\Gamma$ is not cocompact, there is no need for restriction of scalars: $\Gamma$ corresponds to a $\rational$-form~$G_\rational$ on~$G$ itself \csee{IrredNoncpct->NoROS}. We may assume the $\rational$-form is not $\SO(A; \rational)$; otherwise, $\Gamma$ is as described. Therefore, by inspection of \cref{GFxC} and \cref{GFxR}, we see that $G_\rational$ must be of the form $\SU(A, \tau_r ; \quaternion^{a,b}_\rational)$, where $A \in \Mat_{3 \times 3}( \quaternion^{a,b}_\rational)$, and 
	$\quaternion^{a,b}_\rational \otimes_\rational
\real \iso \Mat_{2 \times 2}(\real)$.

Because $G/\Gamma$ is not compact, there is a vector $v \in ( \quaternion^{a,b}_\rational )^3$, such that $\tau_r(v)^\transpose A v = 0$ \csee{NotCpct->UnipD}. Hence, it is not difficult to see that, by making a change of basis, we may assume
	 $$ A = \begin{bmatrix}
	 1 & 0 & 0 \\
	 0 & -1 & 0 \\
	 0 & 0 & p
	 \end{bmatrix}, 
	 \ \text{ for some $p \in \quaternion^{a,b}_\rational$} .$$
Since the identity matrix $\Id_{2 \times 2}$ is the image of $1 \in \quaternion^{a,b}_\rational$ under any
isomorphism $\quaternion^{a,b}_\rational \otimes_\rational \real \to \Mat_{2 \times
2}(\real)$, this means
 $$ G
 = \SU \bigl( A ; \quaternion^{a,b}_{\rational}
\otimes_\rational \real \bigr)
 \iso \SO ( A_\real ; \real )
 , \text{ \ where }
 A_\real = 
 \begin{Smallbmatrix}
 1 &   &   &   &   &   \\
   & 1 &   &   &   &   \\
   &   & -1 &   &   &   \\
   &   &   & -1 &   &   \\
   &   &   &   & \upast & \upast \\
   &   &   &   & \upast & \upast \\
 \end{Smallbmatrix}
 .$$
 Therefore, $G$ is isomorphic to either $\SO(2,4)$ or
$\SO(3,3)$; it is not isogenous to $\SO(1,5)$. This is a
contradiction.
\end{proof}

\begin{prop}[\csee{AllNoncocptSL3RStated}] \label{AllNoncocptSL3R}
Every noncocompact, arithmetic subgroup of\/ $\SL(3,\real)$ is commensurable to a conjugate of either\/
 $\SL(3,\integer)$ or a subgroup of the form\/ $\SU(J_3, \sigma; \ints)$, where
 	\begin{itemize}
	\item $J_3 = 
		 \begin{bmatrix}
		 0 & 0 & 1 \\
		 0 & 1 & 0 \\
		 1 & 0 & 0
		 \end{bmatrix}$,
 	\item $L$ is a real quadratic extension of\/~$\rational$, with Galois automorphism~$\sigma$,
	and
	\item $\ints$ is the ring of integers of~$L$.
	\end{itemize}
 \end{prop}

 \begin{proof}
 Let $\Gamma$ be an arithmetic subgroup of $G = \SL(3,\real)$, such that
$G/\Gamma$ is not compact. We
know, from the Margulis Arithmeticity Theorem
\pref{MargulisArith}, that $\Gamma$ is arithmetic. Since
$G/\Gamma$ is not compact, there is no need for restriction of scalars \csee{IrredNoncpct->NoROS}, so there is a $\rational$-form $G_\rational$ of~$G$, such that $\Gamma$ is commensurable to~$G_\integer$.
By inspection of \cref{GFxC,GFxR}, we see that there are only two
possibilities for~$G_{\rational}$. We consider them individually, as separate cases.

\setcounter{case}{0}

\begin{case}
 Assume $G_{\rational} = \SL(n,D)$, for some central
division algebra~$D$ of degree~$d$ over~$\rational$, with
$dn = 3$.
 \end{case}
 Because $3$~is prime, there are only two possibilities
for~$n$ and~$d$.

\begin{subcase}
 Assume $n = 3$ and $d = 1$.
 \end{subcase}
 Since $d = 1$, we have $\dim_{\rational} D = 1$, so $D =
\rational$. Therefore, $G_{\rational} = \SL(3,\rational)$.
So $\Gamma \approx \SL(3,\integer)$.

\goodbreak

\begin{subcase} \label{AllNoncocptSL3RPf-SL1D}
 Assume $n = 1$ and $d = 3$.
 \end{subcase}
 We have $G_{\rational} = \SL(1,D)$. Therefore $\Gamma \approx
\SL(1,\ints_D)$ is cocompact
\fullcsee{DCocpctinSL3R}{cocpct}. This is a contradiction.

\begin{case}
 Assume $G_{\rational} = \SU(A,\tau;D)$, for $A,D,\sigma$
as in \fullref{FClassicalDefn}{SUSL}, with $F =
\rational$, $L \subset \real$, and $dn = 3$.
 \end{case}
 If $n = 1$, then $G_{\rational} \subseteq \SL(1,D)$, so it has no unipotent elements, which contradicts the fact that $G/\Gamma$ is not compact. Thus, we may
assume that $n = 3$ and $d = 1$. 

Since $d = 1$, we have $D = L$, so $G_{\rational} =
\SU(A,\sigma;L)$, where $\sigma$~is the (unique) Galois
automorphism of~$L$ over~$\rational$, and $B$~is a
$\sigma$-Hermitian form on~$L^3$.

Since $\Gamma$ is not cocompact, we know $\Qrank \Gamma \ge 1$, so there is some nonzero $v \in L^3$ with $v^\transpose A v = 0$ \fullccf{QrankEg}{SUD}. From this, it is not difficult to construct a basis of~$L^3$ in which $A$ is a scalar multiple of~$J_3$ \csee{A=J}. 
 \end{proof}

\begin{exercises}

\item \label{Calcpq(SOAF)}
 Assume that $F \subset \real$, and that $A$ is an invertible, symmetric matrix in $\Mat_{n \times n}(F)$. Show that
if exactly~$p$ of the eigenvalues of $A$ are positive,
then $\SU(A;\real) \iso \SO(p,n-p)$.

\item \label{Calcpq(Sp)}
 Suppose 
 	\begin{itemize}
	\item $F \subset \real$, 
	$\quaternion^{a,b}_F$ is a quaternion division algebra over~$F$, 
	and 
	\item $A$ is an invertible $\tau_c$-Hermitian matrix in $\Mat_{n \times n}(D)$ that is diagonal. 
	\end{itemize}
Show:
	 \begin{enumerate}
	 \item every entry of the matrix $A$ belongs to~$F$ (and, hence, to~$\real$),
	and
	 \item if exactly~$p$ of the diagonal entries of $A$ are positive, then we have $\SU(A, \tau_c; \quaternion^{a,b}_F \otimes_F \real) \iso \Sp(p,n-p)$.
	 \end{enumerate}

\item \label{A=J}
Complete the proof of \cref{AllNoncocptSL3R}, by showing that we may assume $A = J$.
\hint{Assume $v_1$ and~$v_3$ are isotropic, and $v_2$ is orthogonal to both $v_1$ and~$v_3$. Multiply $A$ by a scalar in~$\rational$, so $v_2^* A v_2 = 1$. Then normalize~$v_3$, so $v_1^* A v_3 = 1$.}

\item (B.\,Farb) 
For each $n \ge 2$, find a \emph{cocompact} lattice~$\Gamma_n$ in $\SL(n,\real)$, such that $\Gamma_2 \subseteq \Gamma_3 \subseteq \Gamma_4 \subseteq \cdots$.  (If we did not require $\Gamma_n$ to be cocompact, we could let $\Gamma_n = \SL(n, \integer)$.)

\item \label{G=SU(skewHerm)}
 Show that if $G = \SU(A,\tau_r;\quaternion^{a,b}_F)$, as in
\fullcref{FClassicalDefn}{SUSO}, then there exists $A' \in \Mat_{n \times n}(\quaternion^{a,b}_F)$, such that $A'$ is  skew-Hermitian with respect to the
standard anti-involution~$\tau_c$, and $G = \SU(A', \tau_c;\quaternion^{a,b}_F)$.
 \hint{Use the fact that $\tau_r(x) = j^{-1} \tau_c(x) j$.}

 \end{exercises}

\section{\texorpdfstring{$G$}{G} has a cocompact lattice}
 \label{CocptSect}

We have already seen that if $G$ is not compact, then it has a noncocompact lattice \csee{GHasNoncpctLatt}. In this \lcnamecref{CocptSect}, we will show there is also a lattice that is cocompact:

\begin{thm} \label{GHasCpctLatt}
 $G$ has a cocompact, arithmetic lattice.
\end{thm}

To illustrate the main idea, we briefly recall the prototypical case, which is a generalization of \cref{SO(12;Z[sqrt2])}.

\begin{prop} \label{SOmnHasCocpct}
 $\SO(m,n)$ has a cocompact, arithmetic lattice.
 \end{prop}

\begin{proof}
 Let 
 \noprelistbreak
 \begin{itemize}
 \item $F = \rational(\sqrt{2})$,
 \item $\sigma$ be the Galois automorphism of~$F$
over~$\rational$,
 \item $\ints = \integer[\!\sqrt{2}]$, 
 \item $B(x,y) = \sum_{j=1}^p x_j y_j - \sqrt{2}
\sum_{j=1}^q x_{p+j} y_{p+j}$, for $x,y \in F^{p+q}$,
 \item $G = \SO(B)^\circ$,
 \item $\Gamma = G_{\ints}$, and
 \item $\Delta \colon G_F \to G \times G^\sigma$ defined by
$\Delta(g) = \bigl( g, \sigma(g) \bigr)$.
 \end{itemize}
 We know (from restriction of scalars) that $\Delta(\Gamma)$
is an irreducible, arithmetic lattice in $G \times G^\sigma$
\csee{ResScal->Latt}. Since $G^\sigma \iso \SO(p+q)$ is
compact, we may mod it out, to conclude that $\Gamma$ is an
arithmetic lattice in~$G \iso \SO(m,n)^\circ$. Also, since
$G^\sigma$ is compact, we know that $\Gamma$ is cocompact
\csee{scalars->cpct}.
 \end{proof}

More generally, the same idea can be used to prove that any simple group~$G$ has a
cocompact, arithmetic lattice. Namely, start by letting $K$ be a compact group that has the same complexification as~$G$. (For classical groups, the correct choice of~$K$ can be found by looking at \cref{GxC=}.) Then show that $G \times K$
has an irreducible, arithmetic lattice. More precisely, construct 
	\begin{itemize}
	\item an extension~$F$ of~$\rational$ that has exactly two places $\sigma$ and~$\tau$, 
	and 
	\item a group~$\widehat G$ that is defined over~$F$, 
	\end{itemize}
such that $\widehat G^\sigma$ and $\widehat G^\tau$ are isogenous to $G$ and~$K$, respectively.

Although we could do this explicitly for the classical groups, we will save ourselves a lot of work (and also be able to handle the exceptional groups at the same time) by quoting the following powerful theorem. 

\begin{thm}[(Borel-Harder)] \label{BorelHarderLocGlob}
Suppose 
\noprelistbreak
	\begin{itemize}
	\item $F$ is an algebraic number field,
	\item $S^\infty$ is the set of places of~$F$,
	\item $G$ is defined over~$F$,
	\item $G_\complex$ is connected and simple, and has trivial center,
	and
	\item for each $\sigma \in S^\infty$, we are given some $F_\sigma$-form~$G_\sigma$ of~$G_\complex$.
	\end{itemize}
Then there exists a group~$\widehat G$ that is defined over~$F$, such that $\widehat G^\sigma \iso G_\sigma$, for each $\sigma \in S^\infty$.
\end{thm}

\begin{proof}[Remark on the proof]
For any place~$\sigma$ of~$F$, there is a natural map
	$$ \sigma^* \colon \coho1\bigl( \Gal(\complex/F), \Aut ( G_{\complex} ) \bigr) \to \coho1\bigl( \Gal(\complex/F_\sigma), \Aut ( G_{\complex} ) \bigr) .$$
Namely, any element of the domain corresponds to an $F$-form~$\widehat G$ of~$G_\complex$. The twisted group $\widehat G^\sigma$ is defined over~$\sigma(F)$, and is therefore also defined over~$F_\sigma$. Hence, it determines an element of the range. (The map can also be defined directly, in terms of $1$-cocycles, by restricting a cocycle $\alpha \colon \Gal(\complex/F) \to \Aut ( G_{\complex} )$ to the subgroup $\Gal(\complex/F_\sigma)$.)

However, we should replace $\complex$ with~$\overline{\rational}$ in the domain \fullcsee{QFormCorrections}{Qbar}. Making this correction, and putting together the maps for the various choices of~$\sigma$, we obtain a map
	$$ \coho1\bigl( \Gal(\overline{\rational}/F), \Aut ( G_{\overline{\rational}} ) \bigr) \to \bigtimes_{\sigma \in S^\infty} \coho1\bigl( \Gal(\complex/F_\sigma), \Aut ( G_\complex ) \bigr) .$$
The theorem is proved by showing that this map is surjective.
\end{proof}

%\begin{defn}
%$G$ is \defit[isotypic semisimple Lie group]{isotypic} if the simple factors of $G_\complex$ are all isogenous to each other.
%\end{defn}

\begin{cor} \label{Isotypic->irred}
If $G$ is isotypic \csee{IsotypicDefn}, then $G$ has a cocompact, irreducible, arithmetic lattice.
 \end{cor}
 
 \begin{proof}
 Assume $G$ has trivial center (by passing to an isogenous group), and write $G = G^1 \times \cdots \times G^n$, where each $G^i$ is simple.  
Let $\widehat G = G^1$,
 %We may assume \widehat G is defined over~$\rational$ \csee{hasQform}. Also 
 and assume, for simplicity, that $G^i_{\complex}$ is simple, for every~$i$ \csee{Isotypic->irredNotAbsSimpleEx}. Then, by assumption, $G^i_{\complex} \iso \widehat G_\complex$ for every~$i$, which means that $G^i$ is an $\real$-form of $\widehat G_\complex$. 

Let $F$ be an extension of~$\rational$, such that $F$ has exactly $n$ places $v_1,\ldots,v_n$, and all of them are real \csee{rReal+sImag}. The Borel-Harder Theorem \pref{BorelHarderLocGlob} provides some group~$\widehat G$ that is defined over~$F$, such that $\widehat G_{F_{v_i}} \iso G_i$, for each~$i$. Then restriction of scalars \pref{ResScal->Latt} tells us that $\widehat G_{\ints}$ is (isomorphic to) an irreducible, arithmetic lattice in $\prod_{i=1}^n G_{F_{v_i}} \iso G$.

We may assume that $G^n$ is compact (by replacing $G$ with $G \times K$ for a compact group~$K$, such that $G \times K$ is isotypic). Then \cref{scalars->cpct} tells us that every irreducible, arithmetic lattice in~$G$ is cocompact.
 \end{proof}

\begin{proof}[Proof of \cref{GHasCpctLatt}]
 We may assume $G$ is simple. (If $\Gamma_1$ and~$\Gamma_2$
are cocompact, arithmetic lattices in $G_1$ and~$G_2$, then
$\Gamma_1 \times \Gamma_2$ is a cocompact, arithmetic
lattice in $G_1 \times G_2$.) Then $G$ is isotypic, so
\cref{Isotypic->irred} applies.
 \end{proof}

The converse of \cref{Isotypic->irred} is true (even without assuming cocompactness):

\begin{prop}[\ccf{GZirred->scalars}] \label{Irred->Isotypic}
If $G$ has an irreducible, arithmetic lattice, then $G$ is isotypic.
\end{prop}

By the Margulis Arithmeticity Theorem \pref{MargulisArith}, there is usually no need to assume that the lattice is arithmetic:

\begin{cor}
$G$ is isotypic if 
	\begin{itemize}
	\item it has an irreducible lattice, 
	and 
	\item it is not isogenous to a group of the form $\SO(1,n) \times K$ or $\SU(1,n) \times K$, where $K$ is compact.
	\end{itemize}
\end{cor}

We know that if $G$ has an irreducible, arithmetic lattice that is not cocompact, then $G$ is isotypic \csee{Irred->Isotypic} and has no compact factors \csee{scalars->cpct}. However, the converse is not true:

\begin{eg} \label{SO15xSO3H->Cocpct}
 Every irreducible lattice in $\SO(3,\quaternion) \times \SO(1,5)$ is cocompact.
 \end{eg}

\begin{proof}
 Suppose $\Gamma$ is an irreducible lattice in $\SO(3,\quaternion) \times \SO(1,5)$, such that $G/\Gamma$ is \textbf{not} compact. This will lead to a contradiction.

The Margulis Arithmeticity Theorem \pref{MargulisArith}
implies that $\Gamma$ is arithmetic, so
\cref{irred->scalars} implies that $\Gamma$ can be
obtained by restriction of scalars. Hence, there exist:
	\begin{itemize}
	\item an algebraic number field $F$ with exactly two places $1$ and~$\sigma$, 
	and
	\item a (connected) group~$\widehat G$ that is defined over~$F$,
 	\end{itemize}
such that 
	\begin{itemize}
	\item $\widehat G$ is isogenous to $\SO(3,\quaternion)$, 
	\item $\widehat G^\sigma$ is isogenous to $\SO(1,5)$,
	and
	\item $\Gamma$ is commensurable to $\Delta(G_{\ints})$ in $\widehat G \times \widehat G^\sigma$.
	\end{itemize}
Since $\SO(n,\quaternion)$ occurs only once in \cref{GFxC,GFxR} combined, $\widehat G_F$ must be of the form $\SU(A, \tau_r ; \quaternion^{a,b}_F)$, for some $A \in \Mat_{3 \times 3}( \quaternion^{a,b}_F )$.

However, since $G/\Gamma$ is not compact, the proof of \cref{SO1nNotCpctList} implies that $\widehat G^\sigma$ is isomorphic to either $\SO(2,4)$ or
$\SO(3,3)$; it is not isogenous to $\SO(1,5)$. This is a
contradiction.
 \end{proof}

We close with a stronger version of a fact that was used in the proof of \cref{Isotypic->irred}:

\begin{lem} \label{rReal+sImag}
 For any natural numbers~$r$ and~$s$, not both\/~$0$, there is
an algebraic number field~$F$ with exactly $r$~real places
and $s$~complex places.
 \end{lem}

\begin{proof}
 Let $n = r + 2s$. It suffices to find an irreducible
polynomial $f(x) \in \rational[x]$ of degree~$n$, such
that $f(x)$ has exactly~$r$ real roots. (Then we may let $F
= \rational(\alpha)$, where $\alpha$~is any root of $f(x)$.)

Choose a monic polynomial $g(x) \in \integer[x]$, such that
\noprelistbreak
	 \begin{itemize}
	 \item $g(x)$ has degree $n$,
	 \item $g(x)$ has exactly~$r$ real roots, and
	 \item all of the real roots of $g(x)$ are simple.
	 \end{itemize}
 For example, choose distinct integers
$a_1,\ldots,a_r$, and let 
	$$g(x) = (x-a_1)\cdots(x-a_r)(x^{2s}+ 1) .$$

Fix a prime~$p$. \Cref{rescalepoly} allows us to assume
 \begin{enumerate}
 \item \label{rReal+sImagPf-modp2}
 $g(x) \equiv x^n \pmod{p^2}$, and
 \item \label{rReal+sImagPf-critpt}
 $\min \{\, g(t) \mid g'(t) = 0 \,\} > p$,
 \end{enumerate}
by replacing $g(x)$ with $k^n g(x/k)$, for an
appropriate integer~$k$.

 Let $f(x) = g(x) - p$. From  \pref{rReal+sImagPf-modp2}, we
know that $f(x) \equiv x^n - p \pmod{p^2}$, so the Eisenstein
Criterion \pref{Eisenstein} implies that $f$~is irreducible.
From~\pref{rReal+sImagPf-critpt}, we know that $f(x)$ has the same number of real roots as $g(x)$ (see the figure below). % @@@
Therefore $f(x)$ has exactly~$r$ real roots.
 \end{proof}
 
\centerline{\includegraphics{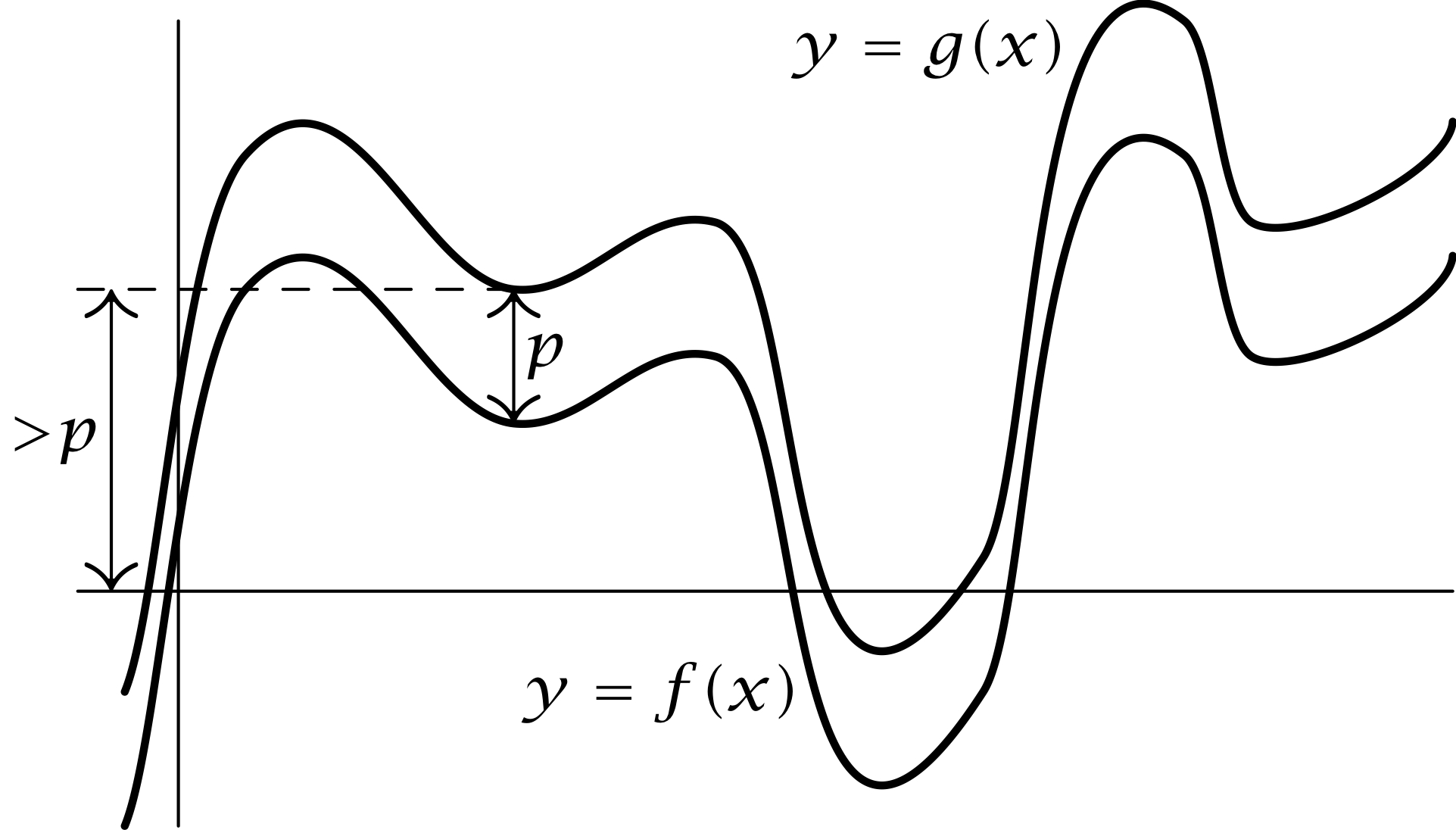}}
%texpreamble
%("  \usepackage{amsmath}
% \usepackage[LY1]{fontenc}
% \usepackage[expert,LY1,mylucidascale]{mylucidabr}
% ");
%defaultpen(  fontcommand("\normalfont") + fontsize(10) ); 
%
%from graph access *;
%unitsize(0.75cm);
%
%pen thickpen = linewidth(1.25);
%
%
%path f = (0.1,-2){NNE}.. tension 2 ..(1,2)..(3,1)..(4.5,1.5) ..(6.5,-1).. tension 2 ..(8,3).. tension 2 ..(8.5,1.5)..tension 2 ..(10,2.25);
%
%draw(f, thickpen);
%draw( shift(0,1)*f, thickpen);
%
%label("$p$", (3,1.5) , 0.5*E);
%draw( (3,1)--(3,2), Arrows(TeXHead,2));
%
%draw( (-0.25,-0.25)--(10,-0.25) );
%draw( (0.5,-2)--(0.5,4) );
%
%draw( (-0.25,2)--(3,2), dashed);
%draw( (0,-0.25)--(0,2), Arrows(TeXHead,2));
%label("${>} p$", (0,.875), 0.5*W);
%
%label("$y = f(x)$", (4.1,-1) );
%label("$y = g(x)$", (6.1,3.8) );

\begin{exercises}

\item \label{CocpctSUSp}
 Use restriction of scalars to
construct cocompact arithmetic subgroups of  $\SU(m,n)$ and
$\Sp(m,n)$ for all $m$ and~$n$.
\hint{See the proof of \cref{SOmnHasCocpct}.}

\item \label{Isotypic->irredNotAbsSimpleEx}
The proof of \cref{Isotypic->irred} assumes that $G^i_\complex$ is simple for every~$i$. Remove this assumption.
\hint{You may use \cref{GxC=GxGiff} (without proof), and you will need the full strength of \cref{rReal+sImag}.}

\item \label{CpctFactorNotCocpct}
 Construct a noncocompact, irreducible lattice in $\SL(2,\real)
\times \SO(3)$.
 \hint{The free group $F_2$ is a noncocompact lattice in
$\SL(2,\real)$. Let $\Gamma$ be the graph of a homomorphism
$F_2 \to \SO(3)$ that has dense image.}

\item \label{NotCpct->UnipD}
In the proof of \cref{SO15xSO3H->Cocpct}, show there exists $v \in ( \quaternion^{a,b}_F )^3$, such that $\tau_r(v)^\transpose A v = 0$.
\hint{If $g$ is a nontrivial unipotent element of~$G_F$, then there is some $w$, such that $g$ fixes the nonzero vector $v = g w - w$.}

\item \label{rescalepoly}
 \begin{enumerate}
 \item Suppose $g(x)$ is a monic polynomial of degree~$n$,
and assume $k \in \integer \smallsetminus \{0\}$, such that $k \, g(x) \in
\integer[x]$. Show $k^n g(x/k) \in \integer[x]$.
 \item Suppose $g(x)$ is a monic, integral polynomial of
degree~$n$, and $p$~is a prime. Show that
 $p^{2n}g(x/p^2) \equiv x^n \pmod{p^2}$.
 \item Suppose $g(x)$ and $h(x)$ are monic polynomials, and
$k$~and~$n$ are nonzero integers, such that $h(x) =
k^ng(x/k)$. Show that
 $$\min \{\, |h(t)| \mid h'(t) = 0 \,\}
 = k^n \min \{\, |g(t)| \mid g'(t) = 0 \,\} .$$
 \end{enumerate}

 \end{exercises}

\begin{notes}

A proof of the classification of complex semisimple Lie algebras (\ref{AlmostAllOverC},\ref{EFGOverCRem}) can be found in standard texts, such as \cite[Thm.~11.4, pp.~57--58, and Thm.~18.4, p.~101]{Humphreys-LieAlg}.

The classification of real simple Lie algebras
(\cref{RformsComplete}) was
obtained by \'E.\,Cartan \cite{Cartan-RSimple}. (The
intervening decades have led to enormous simplifications in
the proof.) 

In \cref{GaloisCohoRealFormsSect,QFormsOfSLnSect}, our cohomological approach to the classification of $F$-forms of the classical groups is based on \cite[\S2.3]{PlatonovRapinchukBook}, where full details can be found. See \cite{Tits-Classification} for a list of all $F$-forms (mostly without proof), including exceptional groups (intended for readers familiar with root systems). The special case of real forms (including exceptional groups) is proved, by a somewhat different approach, in \cite[Chap.~10]{HelgasonBook}. 

\Cref{ArithLattsAreClassical} is due to Weil \cite{Weil-Classical}. 
A proof (together with
\cref{GxC=,GFxC}) is in \cite[\S2.3,
pp.~78--92]{PlatonovRapinchukBook}. We copied \pref{GxC=},
\pref{FClassicalDefn}, and~\pref{GFxC} from
\cite[p.~92]{PlatonovRapinchukBook}, except that
\cite{PlatonovRapinchukBook} uses a different description of
the groups in~\fullref{FClassicalDefn}{SUSO}
\csee{G=SU(skewHerm)}.

\Cref{AutMatInnerEx} is an easy special case of the \thmindex{Skolem-Noether}{Skolem-Noether Theorem}, which can be found in texts on ring theory, such as \cite[\S12.6, p.~230]{Pierce-AssocAlgs}.

\Cref{DivAlg/R}, the classification of division algebras over~$\real$, is due to Frobenius (1878), and a proof can be found in \cite[pp.~452--453]{Jacobson-BasicAlgebra1}. 

\Cref{BorelHarderLocGlob} is due to A.\,Borel and G.\,Harder \cite{BorelHarder-exist}.
See \cite{Johnson-isotypic} for an explicit construction of~$\widehat G$ in the special case where the simple factors of~$G$ are classical.

G.\,Prasad (personal communication) supplied \cref{SO15xSO3H->Cocpct}. It is a counterexample to the noncocompact part of \cite[Thm.~C]{Johnson-isotypic}, which erroneously states that isotypic groups with no compact factors have both cocompact and noncocompact irreducible lattices.

\Cref{GxC=GxGiff} is a consequence of \cite[Prop.~1 of App.~2, p.~385]{Bourbaki-LieGrpsLieAlgs2}, since a connected Lie group is simple if and only if its adjoint representation has no nonzero, proper, invariant subspaces. 

\thmindex{Meyer's}{Meyer's Theorem} (used in
the proof of \cref{QrankGap}) can be found in
\cite[Thm.~1 of \S1.7 and Thm.~5 of \S1.6, pp.~61
and~51]{BorevichShafarevich} or \cite[Cor.~2 of \S4.3.2,
p.~43]{Serre-CourseArith}.

\end{notes}

 %!TEX root = IntroArithGrps.tex

\mychapter{Construction of a Coarse\texorpdfstring{\\}{ }Fundamental Domain} 
\label{ReductionChap}

\prereqs{$\rational$-rank (\cref{QrankChap}). \emph{Recommended:} Siegel sets for $\SL(n,\integer)$ (\cref{IwasawaSLnZ,SiegelSLnZSect,SLNZISLATTSiegelPfSect}).}

The ordinary $2$-torus is often depicted as a square with opposite sides identified, and it would be useful to have a similar representation of $\Gamma  \backslash G$, so we would like to construct a fundamental domain for $\Gamma$ in~$G$. 
Unfortunately, it is usually not feasible to do this explicitly, so, as in \cref{SLnZLattChap}, where we showed that $\SL(n,\integer)$ is a lattice in $\SL(n,\real)$, we will make do with a nice set that is close to being a fundamental domain:

\begin{defn}[\ccf{CoarseFundDomDefn}] \label{CoarseFundDomDefnRedux}
 A subset $\fund$ of~$G$ is called a \defit[fundamental!domain!coarse]{coarse fundamental
domain} for~$\Gamma$ in~$G$ if 
 \begin{enumerate}
 \item $\Gamma \fund = G$, and
 \item \label{CoarseFundDomDefnRedux-finite}
 $\{\, \gamma \in \Gamma \mid \fund \cap \gamma \fund \neq \emptyset \,\}$ is finite.
 \end{enumerate}
\end{defn}

The main result is \cref{ReductThyArithGrps}, which states that the desired set~$\fund$ can be constructed as a finite union of (translates of) ``Siegel sets'' in~$G$. 
Applications of the construction are described in \cref{ReductionAppsSect}.

%We assume familiarity with the special case where $G = \SL(n,\real)$ and $\Gamma = \SL(n,\integer)$, which was discussed in \Cref{SLnZLattChap}.

\section{What is a Siegel set?}

Before defining Siegel sets in every semisimple group, we recall the following special case:

\begin{defn}[\ccf{SiegelSLnZDefn}] \label{WhatIsSiegelSLnROnly}
A \defit{Siegel set} for $\SL(n,\integer)$ is a set of the form $\Siegel_{\overline{N},c} = \overline{N} \, A_c \, K \subseteq \SL(n,\real)$, where 
	\begin{itemize}
	\item $\overline{N}$ is a compact subset of the group~$N$ of upper-triangular unipotent matrices,
%		$$ N = \begin{bmatrix} 1\\ 
%		\BigSymbol{0}{0}{-12}& \BigSymbol{*}{15}{9}\ddots& \\ && 1 \end{bmatrix} $$
	\item $ A_c = \{\, a \in A \mid \text{$a_{i-1,i-1}   \ge c \, a_{i,i}$ for $i = 1,\ldots,n-1$} \,\} $, where $A$ is the group of positive-definite diagonal matrices (and $c > 0$),
%		$$ A = \begin{bmatrix} *\\ 
%		\BigSymbol{0}{0}{-12}& \BigSymbol{0}{15}{7}\ddots& \\ && * \end{bmatrix} 
%		\qquad \begin{pmatrix} \text{with all diagonal} \\ \text{entries positive} \end{pmatrix} ,$$
	and
	\item $K = \SO(n)$.
	\end{itemize}
\end{defn}

In this section, we generalize this notion by replacing $\SL(n,\integer)$ with any arithmetic subgroup (or, more generally, any lattice) in any semisimple Lie group~$G$.
To this end, note that the subgroups $N$, $A$, and~$K$ above are the components of the Iwasawa decomposition $G = KAN$ (or $G = NAK$), which can be defined for any semisimple group \csee{IwasawaDecomp}:
	\noprelistbreak
	\begin{itemize}
	\item %A subgroup of $\SL(n,\real)$ is \defit[unipotent!subgroup]{unipotent} iff it is conjugate to a subgroup of~$N$ \fullcsee{unipEx}{SubgrpInN}.	Therefore, 
	$N$ is a {maximal} unipotent subgroup of~$G$.
	\item %A connected subgroup of $\SL(n,\real)$ is an \defit[torus!R-split@$\real$-split]{$\real$-split torus} iff it is conjugate to a subgroup of~$A$ \csee{RsplitDefn}. Therefore, 
	$A$ is a {maximal} $\real$-split torus of~$G$ that normalizes~$N$,
	and
	\item %Every compact subgroup of $\SL(n,\real)$ is conjugate to a subgroup of $\SO(n)$ \csee{ConjToSOn}. Therefore, 
	$K$ is a {maximal} compact subgroup of~$G$.
	\end{itemize}
%Hence, we have natural characterizations of the subgroups $N$, $A$, and~$K$ in any semisimple group, and the following result generalizes the Iwasawa decomposition of $\SL(n,\real)$ that was stated in \cref{IwasawaDecompSLnR}:
%
%\begin{thm}[(Iwasawa Decomposition)] %% moved to Rrank
% Let 
% \noprelistbreak
% \begin{itemize}
% \item $K$ be a maximal compact subgroup of~$G$,
% \item $A$ be a maximal $\real$-split torus of~$G$, 
% and
% \item $N$ be a maximal unipotent subgroup of~$G$,
% \end{itemize}
% such that $AN$ is a subgroup of~$G$. 
%
%Then $G = K A N$.
% In fact, every $g \in G$ has a \bemph{unique} representation of the form $g = k a u$ with $k \in K$, $a \in A$, and $u \in N$. 
% \end{thm}

Now, to construct Siegel sets in the general case, we will do two things. First, we rephrase \cref{WhatIsSiegelSLnROnly} in a way that does not refer to any specific realization of~$G$ as a matrix group. To this end, recall that, for $G = \SL(n,\real)$, the \emph{\term[Weyl chamber, positive]{positive Weyl chamber}} is
	$$ A^+ = \{\, a \in A \mid \text{$a_{i,i} > a_{i+1,i+1}$ for $i = 1,\ldots,n-1$} \,\} . $$
Therefore, in the notation of \cref{WhatIsSiegelSLnROnly}, we have
	$ A^+ = A_1 $, and, for any $c > 0$, it is not difficult to see that there exists some $a \in A$, such that $A_c =  a A^+$ \csee{Ac=aA+}.
Therefore, letting $C = \overline{N} a$, we see that 
	$$ \text{$\Siegel_{\overline{N},c} = C A^+ K$, and $C$ is a compact subset of~$NA$.} $$
This description of Siegel sets can be generalized in a natural way to any semisimple group~$G$. 

However, all of the above is based entirely on the structure of~$G$, with no mention of~$\Gamma$, but a coarse fundamental domain~$\fund$ needs to be constructed with a particular arithmetic subgroup~$\Gamma$ in mind.
For example, if $\Gamma \backslash G$ is compact (or, in other words, if $\Qrank \Gamma = 0$), then our coarse fundamental domain needs to be compact, so none of the factors in the definition of a Siegel set can be unbounded. 
Therefore, we need to replace the maximal $\real$-split torus~$A$ with a smaller torus~$S$ that reflects the choice of a particular subgroup~$\Gamma$. In fact, $S$ will be the trivial torus when $G/\Gamma$ is compact. In general, $S$ is a maximal $\rational$-split torus of~$G$ (hence, $S$~is compact if and only if $\Gamma \backslash G$ is compact).

Now, if $S$ is properly contained in~$A$, then $NSK$ is not all of~$G$. Hence, $NS$ will usually not be the appropriate replacement for the subgroup $NA$. Instead, if we note that $NA$ is the identity component of a minimal parabolic subgroup of $\SL(n,\real)$ \fullcsee{ParabEgs}{SLn}, and that $NA$ is obviously defined over~$\rational$, then it is natural to replace $NA$ with a minimal parabolic $\rational$-subgroup~$P$ of~$G$.

The following definition implements these considerations. 

\begin{defn} \label{ArithSiegelDefn} 
Assume
\noprelistbreak
	\begin{itemize}
	\item $G$ is defined over~$\rational$,
	\item $\Gamma$ is commensurable to~$G_{\integer}$,
	\item $P$ is a minimal parabolic $\rational$-subgroup of~$G$,
	\item $S$ is a maximal $\rational$-split torus that is contained in~$P$,
	\item $S^+$ is the positive Weyl chamber in~$S$ (with respect to~$P$), 
	\item $K$ is a maximal compact subgroup of~$G$,
	and 
	\item $C$ is any nonempty, compact subset of~$P$.
	\end{itemize}
%We fix the choice of $S$, $S^+$, $P$, and~$K$ for the remainder of this chapter.
Then 
	\ $\Siegel = \Siegel_C = C \, S^+ K $ \ 
is a \defit{Siegel set} for~$\Gamma$ in~$G$.
\end{defn}

\begin{warn}
Our definition of a Siegel set is slightly more general than what is usually found in the literature, because other authors place some restrictions on the compact set~$C$.
For example, it is often assumed that $C$ has nonempty interior.
\end{warn}

\begin{exercises}

\item \label{SiegelinS^+C}
Show that if $\Siegel = C \, S^+ K$ is a Siegel set, then there is a compact subset~$C'$ of~$G$, such that $\Siegel \subseteq S^+ C'$.
\hint{Conjugation by any element of~$S^+$ centralizes $MS$ and contracts~$N$ (where $P = MSN$ is the Langlands decomposition).}

\item \label{CpctInSiegel}
For every compact subset~$C$ of~$G$, show there is a Siegel set that contains~$C$.
\hint{\Cref{G=KP}.}

\end{exercises}

\section{Coarse fundamental domains made from Siegel sets} \label{FundFromSiegelSect}

%It is well-known that every Riemann surface of finite volume has only finitely many cusps. In particular, since the Siegel set pictured in \cref{WeakFunDomSL2ZFig} is a coarse fundamental domain for the action of $\SL(2,\integer)$ on~$\hyperbolic$, it is clear that $\SL(2,\integer) \backslash \hyperbolic $ has only one cusp. Indeed, if $\Gamma$ is any lattice in $G = \SL(2,\real)$, and some Siegel set is a coarse fundamental domain for~$\Gamma$ in~$G$, then $\Gamma \backslash\hyperbolic$ must have only one cusp.
%Since there are examples of (noncompact) hyperbolic surfaces with more than one cusp, this implies that 
%not every lattice has a coarse fundamental domain that consists of a Siegel set.

\begin{eg} \label{SiegelNotFundEg}
Let
\noprelistbreak
	\begin{itemize}
	\item $G = \SL(2,\real)$,
	\item $\Siegel$ be a Siegel set that is a coarse fundamental domain for $\SL(2,\integer)$ in~$G$ \csee{WeakFunDomSL2ZFig},
	and
	\item $\Gamma$ be a subgroup of finite index in~$\SL(2,\integer)$.
	\end{itemize}
Then $\Siegel $ may not be a coarse fundamental domain for~$\Gamma$, because $\Gamma \Siegel$ may not be all of~$G$. In fact, if the hyperbolic surface $\Gamma \backslash \hyperbolic$ has more than one cusp, then no Siegel set is a coarse fundamental domain for~$\Gamma$.

However, if we let $F$ be a set of coset representatives for~$\Gamma$ in $\SL(2,\integer)$, then
	$F \Siegel $
is a coarse fundamental domain for~$\Gamma$ \fullcsee{WeakFundDomFinInd}{IsFund}.
\end{eg}

From the above example, we see that a coarse fundamental domain can sometimes be the union of several translates of a Siegel set, even in cases where it cannot be a single Siegel set.
In fact, this construction always works (if $\Gamma$ is arithmetic):

\begin{namedthm}[\thmindex{Reduction Theory for Arithmetic Groups}Reduction Theory for Arithmetic Groups]
\label{ReductThyArithGrps}
If\/ $\Gamma$ is commensurable to~$G_{\integer}$,
then there exist a Siegel set\/~$\Siegel$ and a finite subset~$F$ of~$G_{\rational}$, such that $\fund = F \, \Siegel$ is a coarse fundamental domain for\/ $\Gamma$ in~$G$.
\end{namedthm}

The proof will be given in \cref{RedThyPfSect}.

\medbreak

Although the statement of this result only applies to arithmetic lattices, it can be generalized to the non-arithmetic case. However, this extension requires a notion of Siegel sets in groups that are not defined over~$\rational$. The following definition reduces this problem to the case where $\Gamma$ is irreducible.

\begin{defn}
If $\Siegel_i$ is a Siegel set for~$\Gamma_i$ in~$G_i$, for $i = 1,2,\ldots,n$, then 
	$$ \Siegel_1 \times \Siegel_2 \times \cdots \times \Siegel_n $$ 
is a \defit{Siegel set} for the lattice $\Gamma_1 \times \cdots \times \Gamma_n$ in $G_1 \times \cdots \times G_n$.
\end{defn}

Then, by the Margulis Arithmeticity Theorem \pref{MargulisArith}, all that remains is to define Siegel sets for lattices in $\SO(1,n)$ and $\SU(1,n)$, but we can use the same definition for all simple groups of real rank one:

\begin{defn} \label{rank1SiegelDefn}
Assume $G$ is simple, $\Rrank G = 1$, and $K$ is a maximal compact subgroup of~$G$.
\noprelistbreak
	\begin{enumerate} \setcounter{enumi}{-1}
	\item If $\Qrank \Gamma = 0$, and $C$ is any compact subset of~$G$, then 
	$ \Siegel = C K $ is a \defit{Siegel set} in~$G$.
	\item Assume now that $\Qrank \Gamma = 1$. Let
		$P$ be a minimal parabolic subgroup of~$G$, with Langlands decomposition $P = MAN$, such that
			\begin{align} \tag{\ref{rank1SiegelDefn}$N$} \label{rank1SiegelDefnN}
			\text{$\Gamma \cap N$ is a maximal unipotent subgroup of~$\Gamma$.}
			\end{align}
If 
	\begin{itemize}
	\item $C$ is any compact subset of~$P$, 
	and 
	\item $A^+$ is the positive Weyl chamber of~$A$ \textup(with respect to~$P$\textup), 
	\end{itemize}
then
	\ $ \Siegel = C A^+ K $ \ 
	is a \defit[Siegel set!generalized]{generalized Siegel set} in~$G$.
	\end{enumerate}
\end{defn}

\begin{rem} \label{Rank1NLattIffQ}
If $\Gamma$ is commensurable to~$G_\integer$ (and $G$ is defined over~$\rational$), then \pref{rank1SiegelDefnN} holds if and only if $P$~is defined over~$\rational$ (and is therefore a minimal parabolic $\rational$-subgroup).
\end{rem}

We can now state a suitable generalization of \cref{ReductThyArithGrps}:

\begin{thm} \label{ReductThyNonarith}
If $G$ has no compact factors, then there exist a generalized Siegel set\/~$\Siegel$ and a finite subset~$F$ of~$G$, such that $\fund = F \, \Siegel$ is a coarse fundamental domain for\/~$\Gamma$ in~$G$.
\end{thm}

The proof is essentially the same as for \cref{ReductThyArithGrps}.

\begin{exercises}

\item \label{GZCpctSiegel=G}
Without using any of the results in this chapter (other than the definitions of ``Siegel set'' and ``coarse fundamental domain''), show that if $\Qrank \Gamma = 0$, then some Siegel set is a coarse fundamental domain for~$\Gamma$ in~$G$.
%\hint{There is a compact subset~$C_0$ of~$G$ with $\Gamma C_0 = G$, and $C_0 K$ is a Siegel set.}

\item Suppose $\fund_1$ and~$\fund_2$ are coarse fundamental domains for $\Gamma_1$ and~$\Gamma_2$ in $G_1$ and~$G_2$, respectively. Show that $\fund_1 \times \fund_2$ is a coarse fundamental domain for $\Gamma_1 \times \Gamma_2$ in $G_1 \times G_2$.

\item Suppose $N$ is a compact, normal subgroup of~$G$, and let $\overline{\Gamma}$ be the image of~$\Gamma$ in $\overline{G} = G/N$. Show that if $\overline{\fund}$ is a coarse fundamental domain for~$\overline{\Gamma}$ in~$\overline{G}$, then 
	$$ \fund = \{\, g \in G \mid gN \in \overline{\fund} \,\} $$
is a coarse fundamental domain for~$\Gamma$ in~$G$.

\item If $G$ is simple, $\Rrank G = 1$, and $G$ is defined over~$\rational$, then \cref{ArithSiegelDefn,rank1SiegelDefn} give two different definitions of the Siegel sets for~$G_\integer$. Show that \cref{rank1SiegelDefn} is more general: any Siegel set according to \cref{ArithSiegelDefn} is also a Siegel set by the other definition.
\hint{\cref{Rank1NLattIffQ}.}

\end{exercises}

\section{Applications of reduction theory} \label{ReductionAppsSect}

Having a coarse fundamental domain is very helpful for understanding the geometry and topology of $\Gamma \backslash G$. Here are a few examples of this (with only sketches of the proofs).

\subsection{$\Gamma$ is finitely presented} \label{FinPresSect}
\Cref{FundDom->FinPres} tells us that if $\Gamma$ has a coarse fundamental domain that is a connected, open subset of~$G$, then $\Gamma$ is finitely presented. The coarse fundamental domains constructed in \cref{ReductThyArithGrps,ReductThyNonarith} are closed, rather than open, but it is easy to deal with this minor technical issue: 

\begin{defn}
A subset $\openSiegel$ of~$G$ is an \defit[Siegel set!open]{open Siegel set} if
$\openSiegel = \open \, S^+ K$, where $\open$\, is a nonempty, precompact, open subset of~$P$.
\end{defn}

Choose a maximal compact subgroup~$K$ of~$G$ that contains a maximal compact subgroup of $\czer_G(S)$. Then we may let:
	\begin{itemize}
	\item $\fund = F \, \Siegel$ be a coarse fundamental domain, with $\Siegel = C \, S^+ K$, such that $C \subseteq P^\circ$ and $\fund$ is connected \csee{CanHaveCinP}, 
	\item $\open$\, be a connected, open, precompact subset of~$P^\circ$ that contains~$C$,
	\item $\openSiegel = \open \,S^+ K$ be the corresponding open Siegel set,
	and
	\item $\openfund = F \, \openSiegel$.
	\end{itemize}
Then $\openfund$ is a coarse fundamental domain for~$\Gamma$ \csee{openfundIsFund}, and $\openfund$ is both connected and open.

This establishes \cref{GammaFinPres}, which stated (without proof) that $\Gamma$~is finitely presented.

\subsection{Mostow Rigidity Theorem}
When $\Qrank \Gamma_1 = 1$, G.\,Prasad constructed a quasi-isometry $\varphi \colon G_1/K_1 \to G_2/K_2$ from an isomorphism $\rho \colon \Gamma_1 \to \Gamma_2$, by using the Siegel-set description of the coarse fundamental domain for $\Gamma_i \backslash G_i$. % would be good to give more explanation @@@
This completed the proof of the Mostow Rigidity Theorem \pref{MostowRigidity}.

\subsection{Divergent torus orbits} \label{DivTorusSect}

\begin{defn}
Let $T$ be an $\real$-split torus in~$G$, and let $x \in G/\Gamma$. We say the $T$-orbit of~$x$ is \defit[divergent torus orbit]{divergent} if the natural map $T \to Tx$ is proper.
\end{defn}

\begin{thm} \label{Qrank=DivTorus}
$\Qrank \Gamma$ is the maximal dimension of an $\real$-split torus that has a divergent orbit on $G/\Gamma$.
\end{thm}

We start with the easy half of the proof:

\begin{lem}[\ccf{SplitDiv}]
Assume $G$ is defined over~$\rational$, and let $S$ be a maximal $\rational$-split torus in~$G$ \textup(so $\dim S = \Qrank G_\integer$\textup). Then the $S$-orbit of $e G_\integer$ is divergent.
\end{lem}

Now, the other half:

\begin{thm}
If $T$ is an\/ $\real$-split torus, and\/ $\dim T > \Qrank \Gamma$, then no $T$-orbit in\/ $G/\Gamma$ is divergent.
\end{thm}

\begin{proof}[Proof \normalfont(assuming $\Qrank \Gamma = 1$)]
Let $T$ be a $2$-dimensional, $\real$-split torus~$T$ of~$G$, and define $\pi \colon T \to G/\Gamma$, by $\pi(t) = t \Gamma$. Suppose $\pi$ is proper. (This will lead to a contradiction.)

Let $P$ be a minimal parabolic $\rational$-subgroup of~$G$, and let $S$ be a maximal $\rational$-split torus in~$P$. For simplicity, let us assume that $\Gamma = G_{\integer}$, and also that a single open Siegel set $\openSiegel = K S^- \open$\, provides a coarse fundamental domain for~$\Gamma$ in~$G$. (Note that, since we are considering $G/\Gamma$, instead of $\Gamma \backslash G$, we have reversed the order of the factors in the definition of the Siegel set, and we use the opposite Weyl chamber.)

Choose a large, compact subset~$C$ of $G/\Gamma$, and let $T_R$ be a large circle in~$T$ that is centered at~$e$. Since $\pi$ is proper, we may assume $T_R$ is so large that $\pi(T_R)$ is disjoint from~$C$. Since $T_R$ is connected, this implies $\pi(T_R)$ is contained in a connected component of the complement of~$C$. So there exists $\gamma \in \Gamma$, such that $T_R \subseteq \Siegel P_{\integer}  \gamma$ (cf.\ \cref{CuspGroup} below). 

Let $t \in T_R$, and assume, for simplicity, that $\gamma = e$. Then $t \in \Siegel P_{\integer}$, and, since $T_R$ is closed under inverses, we see that $\Siegel P_{\integer}$ also contains~$t^{-1}$. However, it is not difficult to see that conjugation by any large element of $\Siegel P_{\integer}$ expands the volume form on~$P$ \csee{SiegelContracts}. Since the inverse of an expanding element is a contracting element, not an expanding element, this is a contradiction.
\end{proof}

\Cref{Qrank=DivTorus} can be restated in the following geometric terms:

\begin{thm}[\csee{QrankFlats}]
$\Qrank \Gamma$ is the largest natural number~$r$, such
that some finite cover of the locally symmetric space\/ $\Gamma \backslash G/K$ contains a closed, simply connected, $r$-dimensional flat.
 \end{thm}

\subsection{The large-scale geometry of locally symmetric spaces} \label{LargeScaleSect}

If we let $\pi \colon G \to \Gamma \backslash G / K$ be the natural map,
then it is not difficult to see that the restriction of~$\pi$ to any Siegel set is proper \csee{ProperOnSiegel}. 
In fact, with much more work (which we omit), it can be shown that the restriction of~$\pi$ is very close to being an isometry:

\begin{thm} \label{AlmIsomOnSiegel}
If\/ $\Siegel = C \, S^+ K$ is any Siegel set, and
	$$ \text{$\pi \colon G \to \Gamma \backslash G / K$ is the natural map,} $$
then there exists $c \in \real^+$, such that, for all $x,y \in \Siegel$, we have
	$$ d \bigl( \pi(x), \pi(y) \bigr) \le d( x,y) \le d \bigl( \pi(x), \pi(y) \bigr) + c .$$
\end{thm}

This allows us to describe the precise shape of the the locally symmetric space associated to~$\Gamma$, up to quasi-isometry:

\begin{thm}
Let 
	\begin{itemize}
	\item $X = \Gamma \backslash G/K$ be the locally symmetric space associated to~$\Gamma$, 
	and 
	\item $r = \Qrank \Gamma$.
	\end{itemize}
Then $X$ is quasi-isometric to the cone on a certain $(r-1)$-dimensional simplicial complex at~$\infty$.
\end{thm}

\begin{proof}[Idea of proof]
Modulo quasi-isometry, any features of bounded size in~$X$ can be completely ignored. Note that:
	\begin{itemize}
	\item \cref{AlmIsomOnSiegel} tells us that, up to a bounded error, $\Siegel$ looks the same as its image in~$X$. 	
	\item There is a compact subset~$C'$ of~$G$, such that $\Siegel \subseteq S^+C'$ \csee{SiegelinS^+C}, so every element of~$\Siegel$ is within a bounded distance of~$S^+$. 
Therefore, $\Siegel$ and~$S^+$ are indistinguishable, up to quasi-isometry.
	\end{itemize}
Then, since $\fund = F \, \Siegel$ covers all of~$X$, we conclude that $X$ is quasi-isometric to $ \bigcup_{f \in F} f S^+ $.

The Weyl chamber $S^+$ is a cone; more precisely, it is the cone on an $(r-1)$-simplex at~$\infty$. Therefore, up to quasi-isometry, $X$~is the union of these finitely many cones, so it is the cone on some $(r-1)$-dimensional simplicial complex at~$\infty$.
\end{proof}
 
\begin{rems} \label{LargeScaleRems} \ 
\noprelistbreak
	\begin{enumerate}

	\item The same argument shows that we get the same picture if, instead of looking at~$X$ modulo quasi-isometry, we look at it from farther and farther away, as in the definition of the asymptotic cone of~$X$ in \cref{AsympConeDefn}. Therefore, the asymptotic cone of~$X$ is the cone on a certain $(r-1)$-dimensional simplicial complex at~$\infty$. This establishes \cref{HattoriThm}.

	\item \label{LargeScaleRems-Tits}
	For a reader familiar with ``Tits buildings\zz,'' the proof (and the construction of~$\fund$) shows that this simplicial complex at~$\infty$ can be constructed by taking the \term{Tits building} of parabolic $\rational$-subgroups of~$G$, and modding out by the action of~$\Gamma$.

	\end{enumerate}
\end{rems}

\begin{exercises}

\item  \label{FinManyFinSubgrps}
Show that
$\Gamma$ has only finitely many conjugacy classes of finite subgroups.
\hint{If $H$ is a finite subgroup of~$\Gamma$, then $H^g \subseteq K$, for some $g \in G$. Write $g = \gamma x$, with $\gamma \in \Gamma$ and $x \in \fund$. Then
	$ H^\gamma x 
%	=  \gamma^{-1} H \gamma x
%	= \gamma^{-1} H g
%	= (\gamma^{-1} g) (g^{-1} H g) 
	= x \cdot H^g
%	\subseteq \fund \cdot K
%	= \fund
	\subseteq \fund$,
so $H$ is conjugate to a subset of
	$ \{\, \gamma \in \Gamma \mid \fund \cap \gamma\fund \neq \emptyset \,\} $.}

\item \label{SplitDiv}
Let $G = \SL(n,\real)$, $\Gamma = \SL(n,\integer)$, and $S$~be the group of positive-definite diagonal matrices. Show the $S$-orbit of $\Gamma e$ is proper.
\hint{If $s_{j,j}/s_{i,i}$ is large, then conjugation by~$s$ contracts a unipotent matrix~$\gamma$ whose only off-diagonal entry is~$\gamma_{i,j}$.}

%\item Show that if $\fund = F \Siegel$, where $\Siegel$ is a Siegel set and $F$ is a finite set, then $\fund = \fund K$.

\item Show that every open Siegel set is an open subset of~$G$ (so the terminology is consistent).

\item \label{CanHaveCinP}
Assume 
	\begin{itemize}
	\item $K$ contains a maximal compact subgroup of $\czer_G(S)$,
	\item $C$ is a compact subset of~$P$, 
	and 
	\item $F$ is a finite subset of~$G$.
	\end{itemize}
Show there is a compact subset $C_\circ$ of~$P^\circ$, such that $C S^+ K \subseteq C_\circ S^+ K$ and $F C_\circ S^+ K$ is connected.
\hint{Show $P^\circ \bigl( K \cap \czer_P(S) \bigr) = P$.}

\item \label{openfundIsFund}
Show the set~$\openfund$ constructed in \cref{FinPresSect} is indeed a coarse fundamental domain for~$\Gamma$ in~$G$.
\hint{\Cref{F1inFinF2}.}

\item \label{ProperOnSiegel}
Let $\pi \colon G \to \Gamma \backslash G$ be the natural map.
Show that if $\Siegel = C S^+ K$ is a Siegel set for~$\Gamma$, then the restriction of~$\pi$ to~$\Siegel$ is proper.
\hint{Let $v$ be a nontrivial element of $N \cap \Gamma$. If $g$ is a large element of $\Siegel$, then $v^g \approx e$.}

\item \label{SiegelContracts}
Show that if $\open$\, is contained in a compact subset of~$P$, then conjugation by any large element of $K S^- \open\, P_{\integer}$ expands the Haar measure on~$P$.
\hint{Conjugation by any element of $M \cup N$ preserves the measure, conjugation by an element of~$\open$\, is bounded, and $S^-$~centralizes $SM$ and expands~$N$. Also note that $P_\integer \doteq M_\integer N_\integer$ \csee{MNZinPZ}.}

\end{exercises}

\section{Outline of the proof of reduction theory} \label{RedThyPfSect}

\begin{notation}
Throughout this \lcnamecref{RedThyPfSect}, we assume 
	\begin{itemize}
	\item $G$~is defined over~$\rational$, 
	\item $\Gamma$~is commensurable to~$G_\integer$,
	and 
	\item $P$ is a minimal parabolic $\rational$-subgroup of~$G$, with Langlands decomposition $P = MSN$.
	\end{itemize}
\end{notation}

In order to use Siegel sets to construct a coarse fundamental domain, a bit of care needs to be taken when choosing a maximal compact subgroup~$K$. Before stating the precise condition, we recall that the \defit[Cartan!involution]{Cartan involution} corresponding to~$K$ is an automorphism~$\tau$ of~$G$, such that $\tau^2$ is the identity, and $K$~is the set of fixed points of~$\tau$. (For example, if $G = \SL(n,\real)$ and $K = \SO(n)$, then $\tau(g) = (g^\transpose)^{-1}$ is the transpose-inverse of~$g$.)

\begin{defn}
A Siegel set $\Siegel = C \, S^+ K$ is \defit[Siegel set!normal]{normal} if $S$ is invariant under the Cartan involution corresponding to~$K$.
\end{defn}

Fix a normal Siegel set $\Siegel = C \, S^+ K$, and some finite $F \subseteq G_\rational$. Then, letting $\fund = F \, \Siegel$, the proof of \cref{ReductThyArithGrps} has two parts, corresponding to the two conditions in the definition of coarse fundamental domain \pref{CoarseFundDomDefnRedux}:
	\begin{enumerate} \renewcommand{\theenumi}{\roman{enumi}}
	\item $\Siegel$ and~$F$ can be chosen so that $\Gamma \fund = G$ \csee{GZFSiegel=G,GammaG/P}, 
	and
	\item for all choices of $\Siegel$ and~$F$, the set 
		$\{\, \gamma \in \Gamma \mid \fund \cap \gamma \fund \neq \emptyset \,\}$ is finite \csee{SiegelProperty}.
	\end{enumerate}
We will sketch proofs of both parts (assuming $\Qrank \Gamma = 1$).
However, as a practical matter, the methods of proof are not as important as understanding the construction of the coarse fundamental domain as a union of Siegel sets \csee{FundFromSiegelSect}, and being able to use this in applications (as in \cref{ReductionAppsSect}). 
%So this material should be considered optional.

\subsection{Proof that $\Gamma \fund = G$}
Here is the rough idea: Fix a base point in $\Gamma \backslash G$.  A Siegel set can easily cover all of the nearby points \csee{CpctInSiegel}, so consider a point $\Gamma g$ that is far away. Godement's Criterion \pref{GodementCriterion} implies there is some nontrivial unipotent $v \in \Gamma$, such that $v^g \approx e$. Replacing $g$ with a different representative of the coset replaces $v$ with a conjugate element. If we assume all the maximal unipotent $\rational$-subgroups of~$G$ are conjugate under~$\Gamma$, this implies that we may assume $v \in N$. If we furthermore assume, for simplicity, that the maximal $\rational$-split torus~$S$ is actually a maximal $\real$-split torus, then the Iwasawa decomposition \pref{IwasawaDecomp} tells us $G = NSK$. The compact group~$K$ is contained in our Siegel set~$\Siegel$, and the subgroup~$N$ is contained in $\Gamma \Siegel$ if $\Siegel$ is sufficiently large, so let us assume $g \in S$. Since $g$~contracts the element~$v$ of~$N$, and, by definition, $S^+$~consists of the elements of~$S$ that contract~$N$, we conclude that $g \in S^+ \subseteq \Siegel$.

\medbreak

We now explain how to turn this outline into a proof.

Recall that all minimal parabolic $\rational$-subgroups of~$G$ are conjugate under~$G_{\rational}$ \fullcsee{parab/Q}{conj}. 
The following technical result from the algebraic theory of arithmetic groups asserts that there are only finitely many conjugacy classes under the much smaller group~$G_{\integer}$. In geometric terms, it is a generalization of the fact that hyperbolic manifolds of finite volume have only finitely many cusps.

\begin{thm} \label{GammaG/P}
There is a finite subset~$F$ of~$G_{\rational}$, such that\/ $\Gamma \, F \, P_{\rational} = G_{\rational}$.
\end{thm}

The finite subset~$F$ provided by the theorem can be used to construct the coarse fundamental domain~$\fund$:

\begin{thm} \label{GZFSiegel=G}
If $F$ is a finite subset of~$G_{\rational}$, such that\/ $\Gamma F P_{\rational} = G_{\rational}$, then there is a\/ \textup(normal\/\textup) Siegel set $\Siegel = C \, S^+ K$, such that\/ $\Gamma F \Siegel = G$.
\end{thm}

\begin{proof}[Idea of proof \normalfont (assuming $\Qrank \Gamma \le 1$)]
For simplicity, assume $\Gamma = G_{\integer}$, and that $F = \{e\}$ has only one element \csee{GZFSiegel=GBigFEx}, so
	\begin{align} \label{AssumeGzPQ=GQ}
	\Gamma \, P_{\rational} = G_{\rational}
	. \end{align}
The theorem is trivial if $\Gamma$ is cocompact \csee{GZCpctSiegel=G}, so let us assume $\Qrank \Gamma = 1$. 

From the proof of the Godement Compactness Criterion \pref{GodementCriterion}, we have a compact subset~$C_0$ of~$G$, such that, for each $g \in G$, either $g \in \Gamma C_0$, or there is a nontrivial unipotent element~$v$ of~$\Gamma$, such that $v^g \approx e$. By choosing $C$ large enough, we may assume $C_0 \subseteq \Siegel$ \csee{CpctInSiegel}.

Now suppose some element $g$ of~$G$ is not in $\Gamma \Siegel$.
Then $g \notin \Gamma C_0$, so there is a nontrivial unipotent element~$v$ of~$\Gamma$, such that 
	\begin{align} \label{Siegelvsmall}
	v^g \approx e 
	. \end{align}
From \pref{AssumeGzPQ=GQ}, we see that we may assume $v \in N$, after multiplying~$g$ on the left by an element of~$\Gamma$ \csee{SiegelvInN}.

We have $G = PK$ \ccf{G=KP}. Furthermore, $P = MSN$, and $\Gamma$ intersects both~$M$ and~$N$ in a cocompact lattice (see \fullcref{QrankEg}{0}, 
\fullcref{LanglandsDecompQ}{M=0}, and \cref{U/UZcpct}).
Therefore, if we multiply $g$ on the left by an element of $\Gamma \cap P$, and ignore a bounded error, we may assume $g \in S$ \csee{SiegelgInS}. Then, since $\Qrank \Gamma = 1$, we have either $g \in S^+$ or $g^{-1} \in S^+$ \csee{S+orS-}. From \pref{Siegelvsmall}, we conclude it is $g$ that is in~$S^+$. So $g \in \Siegel$, which contradicts the fact that $g \notin \Gamma \Siegel$.
\end{proof}

\begin{rem} \label{SiegelMReductive}
The above proof overlooks a technical issue: in the Langlands decomposition $P = MSN$, the subgroup~$M$ may be reductive, rather than semisimple. However, the maximality of~$S$ implies that the central torus~$T$ of~$M$ has no $\real$-split subtori \fullccf{LanglandsDecompQ}{M=0}, so it can be shown that this implies $T/T_{\integer}$ is compact.
Therefore $M/M_\integer$ is compact, even if $M$ is not semisimple.
\end{rem}

\subsection{Proof that $\Siegel$ intersects only finitely many $\Gamma$-translates} \label{SiegelPropertySect}
We know that $\Gamma$ is commensurable to~$G_\integer$. Therefore, if we make the minor assumption that $G_\complex$ has trivial center, then $\Gamma \subseteq G_\rational$ \csee{Comm=GC}. Hence, the following result establishes Condition~\fullref{CoarseFundDomDefnRedux}{finite} for $\fund = F \Siegel$:

\begin{thm}[(``Siegel property'')] \label{SiegelProperty}
If
\noprelistbreak
	\begin{itemize}
%	\item $\Gamma$ is commensurable to $G_{\integer}$,
%	\item $S$ is a maximal $\rational$-split torus in~$G$,
	\item $\Siegel = C S^+ K$ is a normal Siegel set,
	and
	\item $q \in G_{\rational}$,
	\end{itemize}
then 
	$\{\, \gamma \in G_\integer \mid q \Siegel \cap  \gamma \Siegel \neq \emptyset \,\}$ 
is finite.
\end{thm}

\begin{proof}[Proof \normalfont (assuming $\Qrank \Gamma \le 1$)]
The desired conclusion is obvious if $\Siegel$ is compact, so we may assume $\Qrank \Gamma = 1$. To simplify matters, let $\Gamma = G_{\integer}$, and
	$$ \text{assume $q = e$ is trivial.} $$

The proof is by contradiction: assume 
	$$\sigma = \gamma \sigma'  , $$
for some large element $\gamma$ of~$\Gamma$, and some $\sigma, \sigma' \in \Siegel$.
Since $\gamma$ is large, we may assume $\sigma$~is large (by interchanging $\sigma$ with $\sigma'$ and  replacing $\gamma$ with~$\gamma^{-1}$, if necessary).
Let 
	$$ \text{$u$ be an element of $N_{\integer}$ of bounded size.} $$
Since $\sigma \in \Siegel = C S^+ K$, we may write 
	$$ \text{$\sigma = c s k$ with $c \in C$, $s \in S^+$, and $k \in K$. } $$
Then $s$ must be large (since $K$ and~$C$ are compact), so conjugation by~$s$ performs a large contraction on~$N$. 
Since $u^c$ is an element of~$N$ of bounded size, and $K$~is compact, this implies that 
	$u^\sigma \approx e $.
In other words, 
	$$u^{\gamma \sigma'} \approx e .$$
In addition, we know that $u^{\gamma} \in G_{\integer}$. Since $\sigma' \in \Siegel$, we conclude that $u^{\gamma} \in N$ \csee{OnlyUnipPShrinks}.

Now we use the assumption that $\Qrank \Gamma = 1$: since 
	$$u^\gamma \in N \cap N^\gamma ,$$
\cref{Qrank1UniqMaxUnip} tells us that $N = N^\gamma$, so \fullcref{parab/Q}{nzer} implies
	$$ \gamma \in \nzer_G(N) = P .$$
Then, since $(MN)_{\integer}$ has finite index in~$P_{\integer}$ \csee{MNZinPZ}, we may assume $\gamma \in (MN)_{\integer}$.

This implies that we may work inside of~$MN$: if we choose a compact subset $\overline{C} \subseteq MN$, such that 
	$C S^+ \subseteq \overline{C} S$, 
then we have
	\begin{align} \label{SiegelInP}
	\overline{C} (K \cap P) \cap \gamma \overline{C} (K \cap P) \neq \emptyset 
	\end{align}
\csee{SiegelInPEx}. 
Since $\overline{C} (K \cap P)$ is compact (and $\Gamma$ is discrete), we conclude that there are only finitely many possibilities for~$\gamma$.
\end{proof}

\begin{rem} \label{CuspGroup}
When $\Qrank \Gamma = 1$, the first part of the proof establishes the useful fact that there is a compact subset~$C_0$ of~$G$, such that if $\gamma \in \Gamma$, and $\Siegel \cap \gamma \Siegel \not\subseteq C_0$, then $\gamma \in P$.
\end{rem}

\begin{exercises}

\item \label{CanChooseNormalSiegel}
Show that every $\rational$-split torus of~$G$ is invariant under some Cartan involution of~$G$. (Therefore, for any maximal $\rational$-split torus~$S$, there exists a maximal compact subgroup~$K$, such that the resulting Siegel sets $C S^+ K$ are normal.)
\hint{If $\tau$ is any Cartan involution, then there is a maximal $\real$-split torus~$A$, such that $\tau(a) = a^{-1}$ for all $a \in A$. Any $\real$-split torus is contained in some conjugate of~$A$.}

\item \label{SiegelvInN}
In the proof of \cref{GZFSiegel=G}, explain why it may be assumed that $v \in N$.
\hint{Being unipotent, $v$ is contained in the unipotent radical of some minimal parabolic $\rational$-subgroup \fullcsee{parab/Q}{U}. Since $\Gamma \, P_{\rational} = G_{\rational}$, we know that all minimal parabolic $\rational$-subgroups are conjugate under~$\Gamma$.}
%Let $P_0$ be a minimal parabolic $\rational$-subgroup of~$G$, such that $v \in \unip P_0$ \fullsee{parab/Q}{U}. From \pref{AssumeGzPQ=GQ} and the fact that all minimal parabolic $\rational$-subgroups of~$G$ are conjugate under~$G_\rational$ \fullcsee{parab/Q}{conj}, we know that all minimal parabolic $\rational$-subgroups of~$G$ are conjugate under~$\Gamma$. Hence, there is some $\gamma \in \Gamma$ with $P^\gamma = P_0$. Then, by replacing $g$ with $\gamma g$ and $v$ with $v^{\gamma^{-1}}$, we may assume $P_0 = P$. 

\item \label{SiegelgInS}
In the proof of \cref{GZFSiegel=G}, complete the proof without assuming that $g \in S$.
\hint{Write $g = pk \in PK$. If $C$ is large enough that $MN \subseteq (\Gamma \cap P) C$, then we have $g \in \Gamma cs k \subseteq \Gamma CSK$, so $v^s \approx e$.}

\item \label{MNZinPZ}
Show that if $P = MAN$ is a Langlands decomposition of a parabolic $\rational$-subgroup of~$G$, then $(MN)_{\integer}$ contains a finite-index subgroup of~$P_{\integer}$.
\hint{A $\rational$-split torus can have only finitely many integer points.}

\item \label{S+orS-}
Show that if $\Qrank \Gamma = 1$, then $S = S^+ \cup (S^+)^{-1}$.
\hint{If $s \in S^+$, then $s^t \in S^+$ for all $t \in \rational^+$ (and, hence, for all $t \in \real^+$).}

\item \label{GZFSiegel=GBigFEx}
Prove \cref{GZFSiegel=G} in the case where $\Qrank \Gamma = 1$. 
\hint{Replace the imprecise arguments of the text with rigorous statements, and do not assume $F$ is a singleton. (In order to assume $v \in N$, multiply $g$ on the left by an element~$x$ of $F\Gamma$. Then $x \Gamma x^{-1}$ contains a finite-index subgroup of $\Gamma \cap P$.)}

\item \label{OnlyUnipPShrinks}
Let $\Siegel = C S^+ K$ be a Siegel set, and let $P$ be the minimal parabolic $\rational$-subgroup corresponding to~$S^+$. Show there is a neighborhood~$W$ of~$e$ in~$G$, such that if $\gamma \in G_{\integer}$ and $\gamma^\sigma \in W$, for some $\sigma \in \Siegel$, then $\gamma \in \unip P$.

\item \label{NormalSiegel->KcapP}
Suppose 
	\begin{itemize}
%	\item $P$ is a minimal parabolic $\rational$-subgroup with Langlands decomposition $P = MSN$,
	\item $\tau$ is the Cartan involution corresponding to the maximal compact subgroup~$K$, 
	and
	\item $S$ is a $\tau$-invariant.
	\end{itemize}
Show $K \cap P \subseteq M \subseteq \czer_G(S)$.
\hint{Since $\czer_G(S)$ is $\tau$-invariant, the restriction of~$\tau$ to the semisimple part of~$M$ is a Cartan involution. 
Therefore $K \cap M$ contains a maximal compact subgroup of~$M$, which is a maximal compact subgroup of~$P$. The second inclusion is immediate from the definition of the Langlands decomposition.}

\item  \label{SiegelInPEx}
Establish \pref{SiegelInP}.
\hint{$\Siegel \cap P \cap \gamma(\Siegel \cap P) \neq \emptyset$ (since $\gamma \in P$)
and $\Siegel \cap P = C S^+ (K \cap P) = C (K \cap P) S^+$ \csee{NormalSiegel->KcapP}.}

\item \label{SiegelPropertyFullPf}
Give a complete proof of \cref{SiegelProperty}.

\item Show that $\Gamma$ has only finitely many conjugacy classes of
maximal unipotent subgroups.
\hint{You may assume, for simplicity, that $\Gamma = G_{\integer}$ is arithmetic. 
Use \fullCref{parab/Q}{U} and \cref{GammaG/P}.}
%Suppose $U$ is a maximal unipotent subgroup
%of~$G_{\integer}$. Then the Zariski closure~$\Zar{U}$
%of~$U$ is a unipotent $\rational$-subgroup of~$G$. Hence,
%there is a parabolic $\rational$-subgroup $P = MAN$, such
%that $U \subseteq N$. Then $U \subseteq \Zar{U}_{\integer}
%\subseteq N_{\integer}$, so the maximality of~$U$ implies that
%$U = N_{\integer}$. 
%
%The converse is similar (and uses \cref{U/UZcpct}).

\end{exercises}

\begin{notes}

The main results of this chapter were obtained for many classical groups by L.\,Siegel (see, for example, \cite{Siegel-DiscGrps}), and the general results are due to A.\,Borel and Harish-Chandra \cite{BorelHarishChandra-ArithSubgrps}.

The book of A.\,Borel \cite{Borel-IntroArithGrps} is the standard reference for this material; see \cite[Thm.~13.1, p.~90]{Borel-IntroArithGrps} for the construction of a fundamental domain for $G_\integer$ as a union of Siegel sets. The proof there does not assume $G_\integer$ is a lattice, so this provides a proof of the fundamental fact that every arithmetic 
	\thmindex{arithmetic subgroups are lattices}
subgroup of~$G$ is a lattice \csee{arith->latt}. See \cite[\S4.6]{PlatonovRapinchukBook} for an exposition of Borel and Harish-Chandra's original proof of this fact (using the Siegel set $\Siegel _{c_1,c_2,c_3}$ for $\SL(n,\integer)$ from \cref{SiegelSLnZDefn}).

\Cref{AlmIsomOnSiegel} was conjectured by C.\,L.\,Siegel in 1959, and was proved by L.\,Ji \cite[Thm.~7.6]{Ji-MetricCompact}. Another proof is in E.\,Leuzinger \cite[Thm.~B]{Leuzinger-TitsGeometry}.

The construction of coarse fundamental domains for non-arithmetic lattices in groups of real rank one is due to Garland and Raghunathan \cite{GarlandRaghunathan-Rrank1}.
(An exposition appears in \cite[Chap.~13]{RaghunathanBook}.) Combining this result with \cref{ReductThyArithGrps} yields \cref{ReductThyNonarith}.

\Cref{FinManyFinSubgrps} can be found in \cite[Thm.~4.3, p.~203]{PlatonovRapinchukBook}.

%\Cref{HattoriThm,HattoriTits} are due to T.\,Hattori \cite{Hattori-AsympGeom} (and they are also proved in \cite{Leuzinger-TitsGeometry}).

Regarding \cref{SiegelMReductive}, see \cite[Prop.~8.5, p.~55]{Borel-IntroArithGrps} or \cite[Thm.~4.11, p.~208]{PlatonovRapinchukBook} for a proof that if $T$ is a $\rational$-torus that has no $\rational$-split subtori, then $T/T_\integer$ is compact.
\end{notes}

  \immediate\addtocontents{toc}{\protect\toceject} % @@@
 %!TEX root = IntroArithGrps.tex

\mychapter{Ratner's Theorems\texorpdfstring{\\}{ }on Unipotent Flows}
\label{RatnerChap}

\prereqs{none.}

This chapter presents three theorems that strengthen and vastly generalize the  well-known and useful observation that if $V$ is any straight line in the Euclidean plane~$\real^2$, then the closure of the image of~$L$ in~$\torus^2$ is a very nice submanifold \csee{LineDenseInT2}. The plane can be replaced with any Lie group~$H$, and $V$ can be any subgroup of~$H$ that is generated by unipotent elements.

\section{Statement of Ratner's Orbit-Closure Theorem}

\begin{eg} \label{LineDenseInT2}
 Let 
 	\begin{itemize}
	\item $V$ be any $1$-dimensional subspace of the vector space~$\real^2$,
	\item $\torus^2   = \real^2 / \integer^2$ be the ordinary $2$-torus,
	\item $x \in \torus^2$,
	and
	\item $\pi \colon \real^2 \to \torus^2$ be the natural covering map.
	\end{itemize}
Geometrically, $V + x$~is a straight line in the plane, and it is classical (and not difficult to prove) that:
	\begin{enumerate}
	\item \label{LineDenseInT2-notdense}
If the slope of the line $V + x$ is rational, then $\pi(V + x)$ is closed, and it is homeomorphic to the circle~$\torus^1$ \csee{LineDenseInT2-notdenseExer}. 
	\item \label{LineDenseInT2-dense}
If the slope of the line $V + x$ is irrational, then the closure of $\pi(V + x)$ is the entire torus~$\torus^2$ \csee{LineDenseInT2-denseExer}. 
	\end{enumerate}
\end{eg}

An analogous result holds in higher dimensions: 
if we take any vector subspace of~$\real^\ell$, the closure of its image in~$\torus^\ell$ will always be a nice submanifold of~$\torus^\ell$. Indeed, the closure will be a subtorus of~$\torus^\ell$: 

\begin{eg}[\csee{SubspaceInTellExer}] \label{SubspaceInTell}
Let 
	\begin{itemize}
	\item $V$ be a vector subspace of~$\real^\ell$,
	\item $\torus^\ell   = \real^\ell / \integer^\ell$ be the ordinary $\ell$-torus,
	\item $x \in \torus^\ell$,
	and
	\item $\pi \colon \real^\ell \to \torus^\ell$ be the natural covering map.
	\end{itemize}
Then the closure of $\pi(V + x)$ in~$\torus^\ell$ is homeomorphic to a torus~$\torus^k$ (with $0\le k \le \ell$\,).  

More precisely, there is a vector subspace~$L$ of~$\real^\ell$, such that
	\begin{itemize}
	\item the closure of $\pi(V + x)$ is $\pi(L + x)$,
	and
	\item $L$ is defined over~$\rational$ (or, in other words, $L \cap \integer^\ell$ is a $\integer$-lattice in~$L$), so $\pi(L + x) \iso L/ (L \cap \integer^\ell)$ is a torus.
	\end{itemize}
 \end{eg}

The above observation about tori generalizes in a natural way to much more general homogeneous spaces, by replacing:
	\begin{itemize}
	\item $\real^3$ with any connected Lie group~$H$,
	\item $\integer^3$ with a lattice~$\Lambda$ in~$H$,
	\item the vector subspace~$V$ of~$\real^3$ with any subgroup of~$H$ that is generated by unipotent elements,
	\item $x + V$ with the coset $x V$,
	\item the map $\pi \colon \real^\ell \to \torus^\ell$ with the natural covering map $\pi \colon H \to H/\Lambda$,
	and
	\item the vector subspace~$L$ of~$\real^\ell$ with a closed subgroup~$L$ of~$H$.
	\end{itemize}
Because it suffices for our purposes, we state only the case where $H$ is semisimple (so we call the group~$G$, instead of~$H$):

\begin{namedthm}[Ratner's Orbit-Closure Theorem] \label{Ratner-OrbitClosure}
	\thmindex{Ratner's!Orbit-Closure}
Suppose
	\begin{itemize}
	\item $V$ is a subgroup of~$G$ that is generated by unipotent elements, 
	and
	\item $x \in G/\Gamma$.
	\end{itemize}
Then there is a closed subgroup~$L$ of~$G$, such that  the closure of\/ $V x$ is~$L x$. 

Furthermore, $L$~can be chosen so that:
	\begin{enumerate}
	\item $L$ contains the identity component of\/~$V$,
	\item \label{Ratner-OrbitClosure-AlmConn}
	$L$ has only finitely many connected components,
	and
	\item \label{Ratner-OrbitClosure-probmeas}
there is an $L$-invariant, finite measure on $L x$.
	\end{enumerate}
\end{namedthm}

(Also note that $L x$ is closed in $G/\Gamma$, since it is the closure of~$Vx$.)

\begin{rem} 
Write $x = g \Gamma$, for some $g \in G$, and let $\Lambda = (g\Gamma g^{-1}) \cap L$.
	\begin{enumerate}
	\item The theorem tells us that the closure of $V x$ is a very nice submanifold of $G/\Gamma$. Indeed, the closure is homeomorphic to the homogeneous space $L/\Lambda$\,.
	\item Conclusion~\pref{Ratner-OrbitClosure-probmeas} of the theorem is equivalent to the assertion that $\Lambda$ is a lattice in~$L$.
	\end{enumerate}
\end{rem}

\begin{warn} \label{MustBeUnip}
The assumption that $V$ is generated by unipotent elements cannot be eliminated. For example, it is known that if $V$ is the group of diagonal matrices in $G = \SL(2,\real)$, then there are points $x \in G/\SL(2,\integer)$, such that the closure of $V x$ is a \term{fractal}. This means that the closure of $V x$ can be a very bad set that is not anywhere close to being a submanifold.
\end{warn}

Unfortunately, the known proofs of Ratner's Orbit-Closure Theorem are rather long. 
One of the paramount ideas in the proof will be described in \cref{RatnerShearingSect}, but, first, we will present a few of the theorem's applications (in \cref{RatnerApplsSect}) and state two other variants of the theorem (in \cref{RatnerVariantsSect}).

\begin{exercises}

\item \label{LineDenseInT2-notdenseExer}
Verify \fullcref{LineDenseInT2}{notdense}.

\item \label{LineDenseInT2-denseExer}
Verify \fullcref{LineDenseInT2}{dense}.

\item \label{SubspaceInTellExer}
Verify \cref{SubspaceInTell}.

\item Show that if $V$ is connected, then the subgroup~$L$ in the conclusion of \cref{Ratner-OrbitClosure} can also be taken to be connected.

\item \label{NonDivergeEx}
(Non-divergence of unipotent flows)
Suppose 
	\begin{itemize}
	\item $\{u^t\}$ is a unipotent one-parameter subgroup of~$G$,
	and
	\item $x \in G/\Gamma$.
	\end{itemize}
Use \cref{Ratner-OrbitClosure} to show there is a compact subset~$K$ of~$G$, such that 
	$$ \text{$\{\, t \in \real \mid u^t x \in K \,\}$ is unbounded.} $$
{\smaller (Hence, \cref{MargulisUnipReturns} is logically a corollary of \cref{Ratner-OrbitClosure}. However, in practice, \cref{MargulisUnipReturns} is used in the proof of \cref{Ratner-OrbitClosure}.\par}
\hint{Conclusion~\pref{Ratner-OrbitClosure-probmeas} of \cref{Ratner-OrbitClosure} is crucial.}

\item Show, by providing an explicit counterexample, that the assumption that $\{u^t\}$ is unipotent cannot be eliminated in \cref{NonDivergeEx}.
\hint{Consider the one-parameter group of diagonal matrices in $\SL(2,\real)$, and let $\Gamma = \SL(2,\integer)$.}

\end{exercises}

\section{Applications} \label{RatnerApplsSect}

We will briefly describe just a few of the many diverse applications of Ratner's Orbit-Closure Theorem \pref{Ratner-OrbitClosure}.

\subsection{Closures of totally geodesic subspaces} \ 

\noindent \begin{minipage}{2.6in}
\hbox to 0pt{\hskip 2.65in\vbox to 0pt{
\begin{minipage}{2in}
\vskip 0.1in
\includegraphics{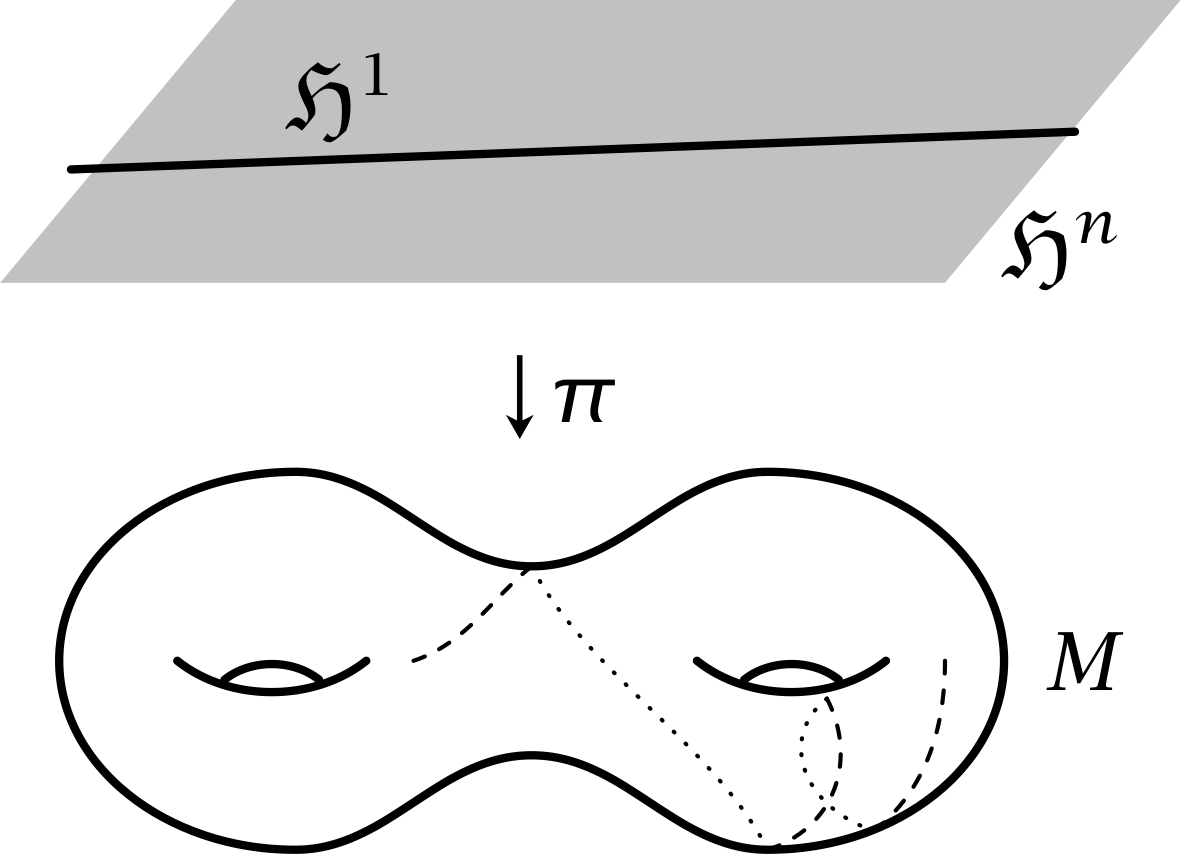}\hss
\end{minipage}
 \vss}\hss}
 \raggedright
\begin{eg} \label{FractalInHyper3Mfld}
 Let $M$ be a  compact, hyperbolic $n$-manifold (with $n \ge 2$).
 \begin{itemize} \itemindent=0pt \leftskip=\parindent
\item There is a covering map $\pi \colon \hyperbolic^n \to M$ that is a local isometry.
\item There is a natural embedding $\iota \colon \hyperbolic^1 \hookrightarrow
\hyperbolic^n$.
\item Let $f_1$ be the composition $\pi \circ \iota$,
 \newline so $f_1 \colon \hyperbolic^1  \to M$.
 \end{itemize}
 \end{eg}
 \end{minipage}

\smallskip

\noindent Then the image of~$f_1$ is a curve in~$M$. It is well known (though not at all obvious) that the closure of this curve can be a very bad set; in fact, even though $\hyperbolic^1$ and~$M$ are nice, smooth manifolds, this closure can be a \term{fractal}. (This is a higher-dimensional analogue of the example in \cref{MustBeUnip}. In the literature, it is the fact that the closure of a geodesic in a compact manifold of negative curvature can be a \term{fractal}.)

It is a consequence of Ratner's Theorem that this pathology never occurs if we replace $\hyperbolic^1$ with a higher-dimensional hyperbolic space:

\begin{cor} \label{GoodClosure(hyper)}
 Let:
 \begin{itemize}
 \item $m,n \in \natural$, with $m \le n$,
 \item $M$ be a  compact, hyperbolic $n$-manifold,
\item $\pi \colon \hyperbolic^n \to M$ be a covering map that is a local isometry,
\item $\iota \colon \hyperbolic^m \hookrightarrow
\hyperbolic^n$ be a totally geodesic embedding,
and
\item $f_m = \pi \circ \iota$, so $f_m \colon \hyperbolic^m  \to M$.
 \end{itemize}
If $m \ge 2$, then the closure $\closure{f_m(\hyperbolic^m)}$ of the image of~$f_m$ is a \textup(totally geodesic\textup) immersed submanifold of~$M$.
\end{cor}

\begin{proof}
We prove only that the closure is a submanifold, not that it is totally geodesic.
Let 
	$$ \text{$V = \SO(1,m)$, \  $G =  \SO(1,n)$, 
	\ and \  $x \in \iota(\hyperbolic^m)$} ,$$
so 
	\begin{itemize}
	\item$G$ acts by isometries on~$\hyperbolic^n$,
	\item $M = \Gamma \backslash \hyperbolic^n$, for some lattice~$\Gamma$ in~$G$,
	and
	\item $\iota(\hyperbolic^m) = V x$.
	\end{itemize}
From Ratner's Orbit-Closure Theorem \pref{Ratner-OrbitClosure}, we know there is a subgroup~$L$ of~$G$, such that 
	$ \closure{\Gamma V} = \Gamma L $.
So 
	$$ \closure{f_m(\hyperbolic^m)}  = \closure{ \pi \bigl( \iota( \hyperbolic^m ) \bigr)}= \closure{\pi(V x)} = \pi(Lx) $$
is an immersed submanifold of~$M$.
\end{proof}

\begin{rems} \label{GoodClosure(LocSymm)} \ 
\begin{enumerate}
\item \label{GoodClosure(LocSymm)-closure}
The same conclusion holds (with the same proof) when $\hyperbolic^m$ and $\hyperbolic^n$ are replaced with much more general symmetric spaces $\cover{X}$ and~$\cover{Y}$ that have no compact factors, except that the closure may not be totally geodesic if $\rank\cover{X} < \rank\cover{Y}$.
\item When $\rank\cover{X} = \rank\cover{Y}$, one proves that the submanifold is totally geodesic by showing that the subgroup~$L$ in the above proof is invariant under the appropriate Cartan involution of~$G$.
\end{enumerate}
\end{rems}

\subsection{Values of quadratic forms}

Many of the most impressive applications of Ratner's Orbit-Closure Theorem (and the related results that will be described in \cref{RatnerVariantsSect}) are in Number Theory. As an example, we present a famous result on values of quadratic forms. It is now an easy corollary of Ratner's Orbit-Closure Theorem, but, historically, it was proved by Margulis before this major theorem was available (by proving the relevant special case of the general theorem).

Let 
	$$ \text{$Q(\vector{x}) = Q(x_1,x_2,\ldots,x_n)$ be a quadratic form in $n$ variables} $$
(in other words, $Q(\vector{x})$ is a homogeneous polynomial of degree~$2$). 

Classical number theorists were interested in determining the values of~$c$ for which the equation $Q( \vector{x}) = c$ has an integer solution; that is, a solution with $\vector{x} \in \integer^n$. For example:
	\begin{enumerate}
	\item Lagrange's 4-Squares Theorem tells us that if
		$$Q(x_1,x_2,x_3,x_4) = x_1^2 + x_2^2 + x_3^2 + x_4^2 ,$$
	then $Q(\vector{x}) = c$ has a solution iff $c \in \integer^{\ge 0}$.
	\item Fermat's 2-Squares Theorem tells us that if $Q(x_1,x_2) = x_1^2 + x_2^2$, and $p$~is an odd prime, then $Q(\vector{x}) = p$ has a solution iff $p \equiv 1 \pmod{4}$.
	\end{enumerate}
These very classical results consider only forms whose coefficients are integers, but we can also look at forms with \emph{irrational} coefficients, such as
	$$ Q(x_1,x_2,x_3,x_4) = 3 x_1^2 - \sqrt{2} x_2 x_3 + \pi x_4^2 .$$
For a given quadratic form $Q(\vector{x})$, it is clear that
	$$ \begin{matrix}
	\text{the equation $Q(\vector{x}) = c$ does not have an integral solution,} \\
	\text{for most real values of~$c$,} 
	\end{matrix} $$
for the simple reason that there are only countably many possible integer values of the variables $x_1,x_2,\ldots,x_n$, but there are uncountably many possible choices of~$c$. Therefore, instead of trying to solve the equation \emph{exactly}, we must be content with solving the equation \emph{approximately}. That is, we will be satisfied with knowing that we can find a value of $Q(\vector{x})$ that is within~$\epsilon$ of~$c$, for every~$\epsilon$ (and every~$c$). In other words, we would like to know that $Q(\integer^n)$ is \emph{dense} in~$\real$.

\begin{egs}
There are some simple reasons that $Q(\integer^n)$ may fail to be dense in~$\real$:
	\begin{enumerate}
	\item Suppose all of the coefficients of $Q(\vector{x})$ are integers. Then we have $Q(\integer^n) \subseteq \integer$, so $Q(\integer^n)$ is obviously not dense in~$\real$. More generally, if $Q(\vector{x})$ is a scalar multiple of a form with integer coefficients, then $Q(\integer^n)$ is not dense in~$\real$.
	\item Suppose all values of $Q(\vector{x})$ are $\ge 0$ (or all are $\le 0$). (In this case, we say that $Q(\vector{x})$ is \defit[positive!definite]{positive-definite} (or \defit{negative-definite}, respectively). For example, this is the case if 
		$$Q(\vector{x}) = a_1 x_1^2 + a_2 x_2^2 + a_3 x_3^2 + \cdots + a_n x_n^2 ,$$
with all coefficients~$a_i$ of the same sign. Then it is clear that $Q(\integer^n)$ is not dense in all of~$\real$.
	\item Let $Q(x_1,x_2) = x_1^2 - \alpha x_2^2$, where $\alpha = 3 + 2 \sqrt{2}$. Then, although it is not obvious, one can show that $Q(\integer^2)$ is not dense in~$\real$ \csee{NotDenseEgEx}. Certain other choices of~$\alpha$ also provide examples where $Q(\integer^2)$ is not dense \csee{BadlyApproxNotDenseEx}, so having only $2$~variables in the quadratic form can cause difficulties.
	\item Even if a form has many variables, there may be a linear change of coordinates that turns it into a form with fewer variables. (For example, letting $z = x + \sqrt{2} y$ transforms $x^2 + 2 \sqrt{2} x y + 2 y^2$ into $z^2$.) A form that admits such a change of coordinates is said to be \defit[quadratic form!degenerate]{degenerate}. Therefore, a degenerate form with more than~$2$ variables could merely be a disguised version of a form with $2$~variables whose image is not dense in~$\real$. 
	\end{enumerate}
\end{egs}

The following result shows that any quadratic form avoiding these simple obstructions does have values that are dense in~$\real$.
It is often called the ``\thmindex{Oppenheim Conjecture}Oppenheim Conjecture\zz,'' because it was an open problem under that name for more than 50 years, but that terminology is no longer appropriate, since it is now a theorem.

\begin{cor}[(Margulis' Theorem on Values of Quadratic Forms)] \label{OppenheimConj}
	\thmindex{Margulis!Values of Quadratic Forms}
Let $Q(\vector{x})$ be a quadratic form  in $n \ge
3$~variables, and assume $Q(\vector{x})$ is:
	\begin{itemize}
	\item not a scalar multiple of a form with integer coefficients,
	\item neither positive-definite nor negative-definite,
	and
	\item nondegenerate.
	\end{itemize}
Then  $Q (\integer^n)$  is dense  in~$\real$.
\end{cor}

\begin{proof}
For simplicity, assume $n = 3$. Let
	\begin{itemize}
	\item $G = \SL(3,\real)$,
	and
	\item $H = \SO(Q)^\circ = \{\, h \in G \mid \text{$Q( h \vector{x} ) = Q(\vector{x})$ for all $\vector{x} \in \real^3$} \, \}^\circ$.
	\end{itemize}
Since $Q(\vector{x})$ is nondegenerate, and neither positive-definite nor negative-definite, we have $H \iso \SO(1,2)^\circ \iso \PSL(2,\real)$, so $H$~is generated by unipotent elements. Furthermore, calculations in Lie theory (which we omit) show that the only connected subgroups of~$G$ containing~$H$ are the obvious ones: $H$ and~$G$. Therefore, Ratner's Orbit-Closure Theorem \pref{Ratner-OrbitClosure} tells us that either:
	\begin{itemize}
	\item $H G_{\integer}$ is closed, and $G_{\integer} \cap H$ is a lattice in~$H$,
	or
	\item the closure of $H G_{\integer}$ is all of~$G$.
	\end{itemize}
However, if $H_{\integer} = G_{\integer} \cap H$ is a lattice in~$H$, then the Borel Density Theorem \pref{BDT(Zardense)} implies that $H$~is defined over~$\rational$ \csee{arithlatt->defdQ}.
Then, since $H = \SO(Q)^\circ$, a bit of algebra shows that $Q(\vector{x})$ is a scalar multiple of a form with integer coefficients \csee{DefdQ->ScalarMult}. This is a contradiction.

Therefore, we conclude that  the closure of $H G_{\integer}$ is all of~$G$. In other words,
	$$ \text{$H G_{\integer}$ is dense in~$G$,} $$
so
	$$ \text{$H G_{\integer} (1,0,0)$ is dense in~$G(1,0,0)$.} $$
Since $G_{\integer} (1,0,0)  \subseteq \integer^3$, and $G (1,0,0) = \real^3 \smallsetminus \{0\}$, this tells us that
	$$ \text{$H \integer^3$ is dense in $\real^3$.} $$
Then, since $Q(\vector{x})$ is continuous, we conclude that
	$$ \text{$Q( H \integer^3 )$ is dense in $Q(\real^3)$.} $$
We also know:
	\begin{itemize}
	\item $Q( H \integer^3 ) = Q(\integer^3)$, by the definition of~$H$,
	and
	\item $Q(\real^3) = \real$, because $Q(\vector{x})$ is neither positive-definite nor negative-definite \csee{Indef->Q(R)=R}.
	\end{itemize}
Therefore $Q(\integer^3)$ is dense in~$\real$.
\end{proof}

\subsection{Products of lattices}

\begin{cor}[\csee{ProdLattsDenseEx}] \label{ProdLattsDense}
If\/ $\Gamma_1$ and\/~$\Gamma_2$ are any two lattices in~$G$, and $G$~is simple, then either
	\begin{enumerate}
	\item $\Gamma_1$ and\/ $\Gamma_2$ are commensurable, so the product\/ $\Gamma_1 \, \Gamma_2$ is discrete,
	or
	\item $\Gamma_1 \, \Gamma_2$ is dense in~$G$.
	\end{enumerate}
\end{cor}

\begin{exercises}

\item \label{BadlyApproxDefnEx}
Suppose $\beta$ is a quadratic irrational. (This means that $\alpha$ is irrational, and that $\alpha$ is a root of a quadratic polynomial with integer coefficients.) Show that $\beta $ is \defit{badly approximable}: i.e., there exists $\epsilon > 0$, such that if $p/q$ is any rational number, then
	$$ \left| \frac{p}{q} - \beta \right| > \frac{\epsilon}{q^2} .$$

\item \label{BadlyApproxNotDenseEx}
Let $Q(x_1,x_2) = x_1^2 - \beta ^2 x_2^2$, where $\beta $ is any badly approximable number \ccf{BadlyApproxDefnEx}. Show that $Q(\integer^2)$ is not dense in~$\real$.
\hint{There exists $\delta > 0$, such that $|Q(p,q)| \ge \delta $ for $p,q \in \integer \smallsetminus \{0\}$.}

\item \label{NotDenseEgEx}
Let $Q(x_1,x_2) = x_1^2 - \alpha x_2^2$, where $\alpha = 3 + 2 \sqrt{2}$. Show $Q(\integer^2)$ is not dense in~$\real$.
\hint{Use previous exercises, and note that $3 + 2 \sqrt{2} = \bigl( 1 + \sqrt{2} \bigr)^2$.}

\item \label{DefdQ->ScalarMult}
Suppose $Q(\vector{x})$ is a nondegenerate quadratic form in $n$~variables.
Show that if $\SO(Q)^\circ$ is defined over~$\rational$, then $Q(\vector{x})$ is a scalar multiple of a form with integer coefficients.
\hint{Up to scalar multiples, there is a unique quadratic form that is invariant under $\SO(Q)^\circ$, and the uniqueness implies that it is invariant under the Galois group $\Gal(\complex/\rational)$.}

\item \label{Indef->Q(R)=R}
Suppose $Q(\vector{x})$ is a quadratic form in $n$~variables that is neither positive-definite nor negative-definite. Show $Q(\real^n) = \real$.
\hint{$Q(\lambda \vector{x}) = \lambda^2 Q(\vector{x})$.}

\item \label{ProdLattsDenseEx}
Prove \cref{ProdLattsDense}.
\hint{Let $\Gamma = \Gamma_1 \times \Gamma_2 \subset G \times G$, and $H = \{\, (g,g) \mid g \in G\,\}$. Show the only connected subgroups of $G \times G$ that contain~$H$ are the two obvious ones: $H$ and $G \times G$. 
Therefore, Ratner's Theorem implies that either $H \cap \Gamma$ is a lattice in~$H$, %(so $\Gamma_1$ and~$\Gamma_2$ are commensurable) 
or $\Gamma H$ is dense in $G \times G$.%(so $\Gamma_1 \, \Gamma_2$ is dense in~$G$
}

\end{exercises}

\section{Two measure-theoretic variants of the theorem} \label{RatnerVariantsSect}

Ratner's Orbit-Closure Theorem \pref{Ratner-OrbitClosure} is purely topological, or qualitative. In some situations, it is important to have quantitative information. 

\begin{eg}
We mentioned earlier that if $V$ is a line with irrational slope in~$\real^2$, then the image $\pi(V)$ of~$V$ in~$\torus^2$ is dense \csee{LineDenseInT2}. 
For applications in analysis, it is often necessary to know more, namely, that $\pi(V)$ is \defit{uniformly distributed} in~$\torus^2$. Roughly speaking, this means that a long segment of~$\pi(V)$ visits all parts of the torus equally often \csee{UnifDistVisitsOpenSets}. 
\end{eg}

\begin{defn} \label{UnifDistDefn}
Let
	\begin{itemize}
	\item $\mu$ be a probability measure on a topological space~$X$,
	and
	\item $c \colon [0,\infty) \to X$ be a continuous curve in~$X$.
	\end{itemize}
We say that $c$ is \defit{uniformly distributed} in~$X$ (with respect to~$\mu$) if, for every continuous function $f \colon X \to \real$ with compact support, we have
	$$ \lim_{T \to \infty} \frac{1}{T} \int_0^T f \bigl( c(t) \bigr) \, dt = \int_X f \, d\mu . $$
\end{defn}

Ratner's Orbit-Closure Theorem tells us that if 
	\begin{itemize}
	\item $U$ is any one-parameter subgroup of~$G$, 
	and 
	\item $x$~is any point in $G/\Gamma$, 
	\end{itemize}
then the closure of the $U$-orbit $U x$ is a nice submanifold of $G/\Gamma$. The following theorem tells us that the $U$-orbit is uniformly distributed in this submanifold.

\begin{namedthm}[Ratner's Equidistribution Theorem] \label{Ratner-Equidistribution}
	\thmindex{Ratner's!Equidistribution}
Let 
	\begin{itemize}
	\item $\{u^t\}$ be any one-parameter unipotent subgroup of~$G$,
	\item $x \in G/\Gamma$,
	and
	\item $c(t) = u^t x$, for $t \in [0,\infty)$.
	\end{itemize}
Then there is a connected, closed subgroup~$L$ of~$G$, such that
	\begin{enumerate}
	\item there is a\/ \textup(unique\textup) $L$-invariant probability measure~$\mu$ on $L x$,
	\item the curve $c$ is uniformly distributed in $L x$, with respect to~$\mu$,
	\item the closure of\/ $\{\, c(t) \mid t \in [0,\infty) \,\}$ is~$L x$ \textup(so $L x$ is closed in $G/\Gamma$\textup),
	and
	\item $\{u^t\} \subseteq L$.
	\end{enumerate}
\end{namedthm}

In the special case where $V$ is a one-parameter unipotent subgroup of~$G$, the following theorem is a consequence of the above Equidistribution Theorem \csee{EquiDist->MeasuresEx}.

\begin{namedthm}[Ratner's Classification of Invariant Measures] \label{Ratner-MeasClass}
	\thmindex{Ratner's!Classification of Invariant Measures}
Suppose
\noprelistbreak
	\begin{itemize}
	\item $V$ is a subgroup of~$G$ that is generated by unipotent elements, 
	and
	\item $\mu$ is any ergodic $V$-invariant probability measure on\/ $G/\Gamma$.
	\end{itemize}
Then there is a closed subgroup~$L$ of~$G$, and some $x \in G/\Gamma$, such that $\mu$ is the unique $L$-invariant probability measure on $L x$. 

Furthermore, $L$ can be chosen so that:
	\begin{enumerate}
	\item $L$ has only finitely many connected components,
	\item $L$ contains the identity component of~$V$,
	and
	\item $L x$ is closed in $G/\Gamma$.
	\end{enumerate}
\end{namedthm}

Here is a sample consequence of the Measure-Classification Theorem:

\begin{cor} \label{UnipotentIsomorphismRigidity}
Suppose 
	\begin{itemize}
	\item $u_i$ is a nontrivial unipotent element of~$G_i$, for $i = 1,2$,
	\item $f \colon G_1/\Gamma_1 \to G_2/\Gamma_2$ is a measurable map that intertwines the translation by~$u_1$ with the translation by~$u_2$; that is,
		$$ \text{$f( u_1 x) = u_2 \, f(x)$, \ for a.e.\ $x \in G_1/\Gamma_1$,} $$
	and
	\item $G_1$ is connected and almost simple.
	\end{itemize}
Then $f$ is a continuous function\/ \textup(a.e.\textup).
\end{cor}

\begin{proof}
Let
	\begin{itemize}
	\item $G = G_1 \times G_2$,
	\item $\Gamma = \Gamma_1 \times \Gamma_2$,
	\item $u = (u_1,u_2) \in G$,
	and
	\item $\graph(f) = \bigset{ \bigl( x, f(x) \bigr) }{ x \in G_1/\Gamma_1 } \subset G/\Gamma$.
	\end{itemize}
The projection from $G/\Gamma$ to $G_1/\Gamma_1$ defines a natural one-to-one correspondence between $\graph(f)$ and $G_1/\Gamma_1$. In fact, since $f$~intertwines $u_1$ with~$u_2$, it is easy to see that the projection provides an isomorphism between the action of~$u_1$ on $G_1/\Gamma_1$ and the action of $u$ on $\graph(f)$.
In particular, the $u_1$-invariant probability measure~$\mu_1$ on $G_1\Gamma_1$ naturally corresponds to a $u$-invariant probability measure~$\mu$ on $\graph(f)$.

The 
	\thmindex{Moore Ergodicity}Moore Ergodicity Theorem \pref{MooreErgodicity} 
tells us that $\mu_1$ is ergodic for~$u_1$, so $\mu$~is ergodic for~$u$. Hence, Ratner's Measure-Classification Theorem \pref{Ratner-MeasClass} provides
	\begin{itemize}
	\item a closed subgroup $L$ of~$G$,
	and
	\item $x \in G/\Gamma$,
	\end{itemize}
such that 
	\begin{itemize}
	\item $\mu$ is the $L$-invariant measure on $L x$,
	and
	\item $L x$ is closed.
	\end{itemize}
Since the definition of~$\mu$ implies that  $\mu \bigl( \graph(f) \bigr) = 1$, and the choice of~$L$ implies that the complement of $L x$ has measure~$0$, we may assume, by changing $f$ on a set of measure~$0$, that 
	$$\graph(f) \subseteq Lx .$$

Assume, for simplicity, that $L$ is connected and $G_1$ is simply connected. Then the natural projection from~$L$ to~$G_1$ is an isomorphism \csee{UnipotentIsomorphismRigidity-coveringmapEx}, so there is a (continuous) homomorphism $\rho \colon G_1 \to G_2$, such that 
	$$L = \graph(\rho) .$$
Assuming, for simplicity, that $x = (e,e)$, this implies that $f(g \Gamma) = \rho(g) \Gamma$ for all $g \in G$. So $f$, like~$\rho$, is continuous.
\end{proof}

As an example of the many important consequences of Ratner's Measure-Classification Theorem \pref{Ratner-MeasClass}, we point out that it implies the Equidistribution Theorem \pref{Ratner-Equidistribution}. The proof is not at all obvious, and we will not attempt to explain it here, but the following simple example illustrates the important precept that knowing all of the invariant measures can lead to an equidistribution theorem.

\begin{prop} \label{UniqErg->UnifDist}
Let  
\noprelistbreak
	\begin{itemize}
	\item $\{g^t\}$ be a one-parameter subgroup of~$G$, 
	\item $\mu$ be the $G$-invariant probability measure on $G/\Gamma$,
	\item $x \in G/\Gamma$,
	and
	\item $c(t) = g^t x$.
	\end{itemize}
If 
	\begin{itemize}
	\item $\mu$ is the only $g^t$-invariant probability measure on $G/\Gamma$, 
	and
	\item $G/\Gamma$ is compact,
	\end{itemize}
then the curve $c$ is uniformly distributed in $G/\Gamma$, with respect to~$\mu$.
\end{prop}

\begin{proof}
Suppose $c$ is not uniformly distributed. Then there is a sequence $T_k \to \infty$, and some continuous function $f_0 \in C(G/\Gamma)$, such that
	\begin{align} \label{UniqErg->UnifDistPf-notmu}
	\lim_{k \to \infty} \frac{1}{T_k} \int_0^{T_k} f_0 \bigl( c(t) \bigr) \, dt \neq \mu(f_0) 
	\end{align}
By passing to a subsequence, we may assume that
	$$ \lambda(f) = \lim_{k \to \infty} \frac{1}{T_k} \int_0^{T_k} f \bigl( c(t) \bigr) \, dt $$
exists for every $f \in C(G/\Gamma)$ \csee{UniqErg->UnifDistPf-LimExistsEx}. Then:
	\begin{enumerate}
	\item $\lambda$ is a continuous linear functional on the space $C(G/\Gamma)$ of continuous functions on $G/\Gamma$, so the Riesz Representation Theorem \pref{RieszRepThm} tells us that $\lambda$ is a measure on $G/\Gamma$. 
	\item $\lambda(1) = 1$, so $\lambda$ is a probability measure. 
	\item From the definition of~$\lambda$, it is not difficult to see that $\lambda$ is $g^t$-invariant \csee{UniqErg->UnifDistPf-gInvt}. 
	\end{enumerate}
Since $\mu$ is the only $g^t$-invariant probability measure, we must have $\lambda = \mu$. However, \pref{UniqErg->UnifDistPf-notmu} says $\lambda(f_0) \neq \mu(f_0)$, so this is a contradiction.
\end{proof}

\begin{rem}
Here is a rough outline of how the three theorems are proved:
	\begin{enumerate}
	\item Measure-Classification is proved in the case where $V$ is unipotent.
		\begin{itemize}
		\item The general case of Measure-Classification follows from this.
		\end{itemize}
	\item Equidistribution is a consequence of Measure-Classification.
	\item Equidistribution easily implies Orbit-Closure in the special case where $V$ is a one-parameter unipotent subgroup.
		\begin{itemize}
		\item The general case of Orbit-Closure can be deduced from this.
		\end{itemize}
	\end{enumerate}
\end{rem}

\begin{rems} \label{MeasClassRems} \ 
\noprelistbreak
	\begin{enumerate}
	\item Ratner's Measure-Classification Theorem \pref{Ratner-MeasClass} remains valid if the lattice~$\Gamma$ is replaced with any closed subgroup of~$G$. However, the other two theorems do not remain valid in this generality.

	\item \label{MeasClassRems-BenoistQuint}
Ratner's Theorems assume that $V$ is generated by unipotent elements, but
N.\,Shah suggested that they might remain valid under the much weaker assumption that the Zariski closure of~$V$ is generated by unipotent elements. For the important case where the Zariski closure of~$V$ is semisimple (with no compact factors), this was recently proved by Y.\,Benoist and J.--F.\,Quint.
	\end{enumerate}
\end{rems}

\begin{exercises}

\item \label{UnifDistVisitsOpenSets}
Suppose $\mu$ is a probability measure on a compact metric space~$X$. Show that a curve $c \colon [0,\infty) \to X$ is uniformly distributed with respect to~$\mu$ if and only if, for every open subset~$\open$ of~$X$, such that $\mu( \bdry \open\,) = 0$, we have
	$$ \lim_{T \to \infty} \frac{1}{T} \leb \{\, t \in [0,T] \mid c(t) \in \open \,\} = \mu(\open) .$$
\hint{Bound the characteristic function of~$\open$\, above and below by continuous functions and apply \cref{UnifDistDefn}.}

\item Show that if $v$ is a nonzero vector of irrational slope in~$\real^2$, then the curve
	$ c(t) = \pi(tv)$
is uniformly distributed in~$\torus^2$ (with respect to the usual Lebesgue measure on the torus).
\hint{Any continuous function on~$\torus^2$ can be approximated by a trigonometric polynomial $\sum a_{m,n} e^{2\pi i m x + 2\pi i n y}$.}

\item \label{EquiDist->MeasuresEx}
Show that if $V = \{u^t\}$ is a one-parameter unipotent subgroup of~$G$, then the conclusions of \cref{Ratner-MeasClass} follow from \cref{Ratner-Equidistribution}.
\hint{Pointwise Ergodic Theorem \csee{PointwiseErgThmForREx}.}

\item \label{UnipotentIsomorphismRigidity-coveringmapEx}
In the setting of the proof of \cref{UnipotentIsomorphismRigidity}, show that the projection $L \to G_1$ is a (surjective) covering map.
\hint{Show $L \cap \bigl(\{e\} \cap G_2 \bigr)$ is discrete, by using Fubini's Theorem and the fact that $\graph(f)$ has nonzero measure.}

\item \label{UniqErg->UnifDistPf-LimExistsEx}
Suppose 
	\begin{itemize}
	\item $T_k \to \infty$,
	\item $G/\Gamma$ is compact,
	and
	\item $c \colon [0,\infty) \to G/\Gamma$ is a continuous curve.
	\end{itemize}
Show there is a subsequence $T_{k_i} \to \infty$, such that
	$$ \lambda(f) = \lim_{i \to \infty} \frac{1}{T_{k_i}} \int_0^{T_{k_i}} f \bigl( c(t) \bigr) \, dt $$
exists for all $f \in C(G/\Gamma)$.
\hint{It suffices to consider a countable subset of $C(G/\Gamma)$ that is dense in the topology of uniform convergence.}

\item \label{UniqErg->UnifDistPf-gInvt}
In the notation used in the proof of \cref{UniqErg->UnifDist}, show, for every $f \in C(G/\Gamma)$ and every $s \in \real$, that
	$$ \lambda(f) = \lim_{k \to \infty} \frac{1}{T_k} \int_0^{T_k} f \bigl( g^s \, c(t) \bigr) \, dt .$$

\end{exercises}

\section{Shearing --- a key ingredient in the proof} \label{RatnerShearingSect}

The known proofs of any of the three variants of Ratner's Theorem are quite lengthy, so we will just illustrate one of the main ideas that is involved. To keep things simple, we will assume $G = \SL(2,\real)$.

\begin{notation} \label{ShearingNotation}
 Throughout this section,
 \begin{itemize}
 \item $G = \SL(2,\real)$,
 \item $u^t = \begin{bmatrix} 1 & 0 \\ t & 1 \end{bmatrix}$ is a one-parameter
\term[one-parameter subgroup!unipotent]{unipotent subgroup} of~$G$,
 \item $a^t= \begin{bmatrix} e^t & 0 \\ 0 & e^{-t} \end{bmatrix}$ is a  one-parameter
\term[one-parameter subgroup!diagonal]{diagonal subgroup} of~$G$,
and
 \item $X = \SL(2,\real)/\Gamma$.
 \end{itemize}
 \end{notation}

The proofs of Ratner's Theorems depend on an understanding of what happens to two \term[nearby points|(]{nearby points} of~$X$ as they are moved by the one-parameter subgroup~$u^t$.

\begin{defn}
 If $x$ and~$y$ are any two points of~$G/\Gamma$, then there
exists $q \in G$, such that $y = qx$. If $x$ is close to~$y$ (which we
denote 
	\nindex{$x \approx y$ means $x$ is close to~$y$}$x \approx y$), 
then $q$~may be chosen close to the identity. Therefore,
we may define a {metric}~$d$ on $G/\Gamma$ by
 $$ d(x,y) = \min \bigset{ \|q - \Id \| }{
 \begin{matrix} q \in G, \\ qx = y \end{matrix}
  } ,$$
 where
 \begin{itemize}
 \item $\Id$ is the identity matrix,
 and
 \item  $\| \cdot \|$ is any (fixed) \index{norm!matrix}matrix norm on
 $\Mat_{2 \times 2}(\real)$. For example, one may take
 $$ \left\| \begin{bmatrix} \mathsf{a} & \mathsf{b} \\ \mathsf{c} &
\mathsf{d} \end{bmatrix} \right \|
 = \max \bigl\{ |\mathsf{a}|, |\mathsf{b}|, |\mathsf{c}|, |\mathsf{d}|
\bigr\} .$$
(Actually, this definition does not guarantee $d(x,y) = d(y,x)$, so it may not define a metric, but let us ignore this minor issue.)
 \end{itemize}
 \end{defn}

Now, we consider two points $x$ and~$qx$, with $q
\approx \Id$, and we wish to calculate
 $ d( u^t x, u^t q x ) $
 \csee{movapart}. 
 \begin{figure}[ht]
 \begin{center}
 \includegraphics{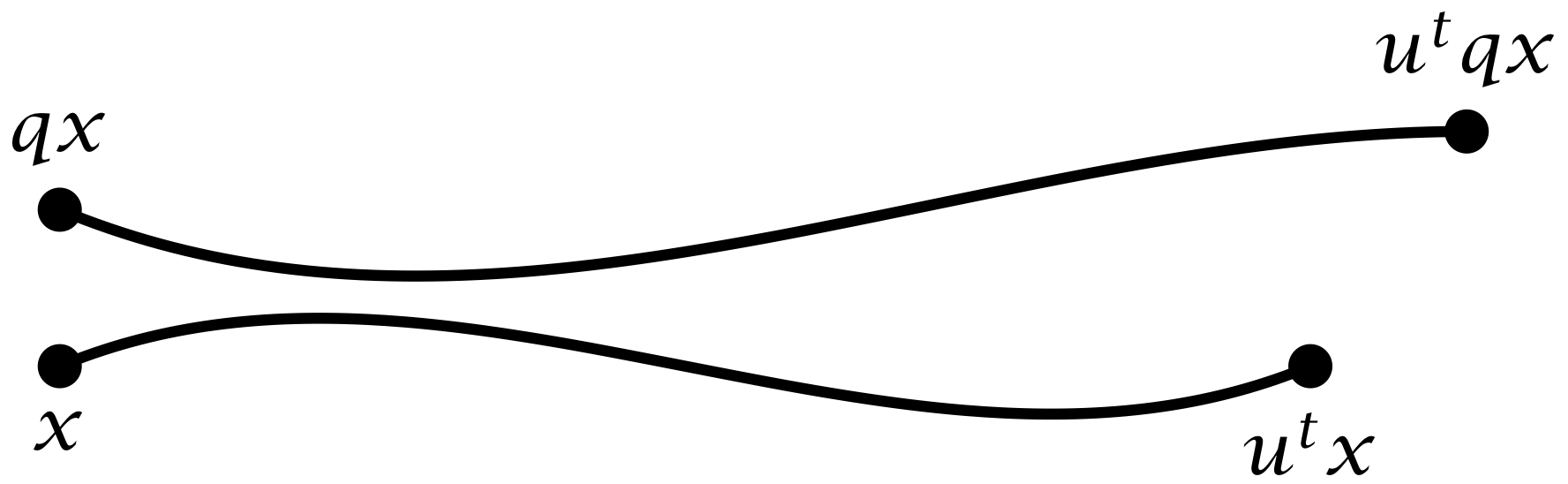}
 \caption{The $u^t$-orbits of two nearby orbits.}
 \label{movapart}
 \end{center}
 \end{figure}
%texpreamble
%("  \usepackage{amsmath}
% \usepackage[LY1]{fontenc}
% \usepackage[expert,LY1,mylucidascale]{mylucidabr}
% ");
%defaultpen(  fontcommand("\normalfont") + fontsize(10) ); 
%
%from graph access *;
%unitsize(1.5cm);
%
%real linethick = 1.5;
%real dotthick = 12;
%dotfactor=12;
%
%pair x = (0,0), qx = (0,0.5);
%dot( x ); dot( qx); label( "$x$", x, 2*S); label( "$qx$", qx, 2*N); 
%
%pair ux = (4,0), uqx = (4.5,0.75);
%dot( ux ); dot( uqx); label( "$u^tx$", ux, 2*S); label( "$u^tqx$", uqx, 2*N); 
%
%draw( x{ENE}..{ENE}ux , linewidth(linethick) );
%draw( qx{ESE}..{E}uqx , linewidth(linethick) );
 \begin{itemize}
 \item To get from~$x$ to~$qx$, one multiplies by~$q$; therefore 
 	$$d(x,qx) = \|q - \Id \| .$$
 \item To get from $u^t x$ to~$u^t q x$, one multiplies by $u^t q u^{-t}$;
therefore 
 $$d( u^t x, u^t q x) = \|u^t q u^{-t} - \Id \| .$$
 (Actually, this equation only holds when the right-hand side is small --- there are infinitely many elements~$g$ of~$G$ with $g u^t x = u^t q x$, and the distance is obtained by choosing
the smallest one, which may not be $u^t q u^{-t}$ if $t$ is large.)
 \end{itemize}
 Letting
 $$ q - \Id = \begin{bmatrix}\mathsf{a} &\mathsf{b} \\ \mathsf{c} & \mathsf{d} \end{bmatrix}
,$$
 a simple matrix calculation shows that
 \begin{align} \label{ConjByU}
 u^t q u^{-t} - \Id = \begin{bmatrix}
 \mathsf{a} -\mathsf{b} t & \mathsf{b} \\
 \mathsf{c}  + (\mathsf{a} -\mathsf{d})t - \mathsf{b} t^2 & \mathsf{d} + \mathsf{b} t
 \end{bmatrix}
 .
  \end{align}

\begin{notation}
 For convenience, let $x_t = u^t x$ and $y_t = u^t y$.
 \end{notation}

Consider the right-hand side of \cref{ConjByU}, with $\mathsf{a}$, $\mathsf{b}$,
$\mathsf{c}$, and~$\mathsf{d}$ very small. Indeed, let us say they are
\index{infinitesimal}infinitesimal (too small to see). As $t$~grows, it is
the quadratic term in the bottom left corner that will be the first matrix entry to attain
macroscopic size. Comparing with the definition
of~$u^t$ \csee{ShearingNotation}, we see that this is exactly the direction
of the $u^t$-orbit. Therefore:

\begin{prop}[{(\index{shearing property}{Shearing Property})}]
\label{ShearingSL2R}
 The fastest \term{relative motion} between two
\term[nearby points|)]{nearby points} is parallel to the orbits of the flow.
 \end{prop}

 \begin{figure}[ht]
 \begin{center}
 \includegraphics{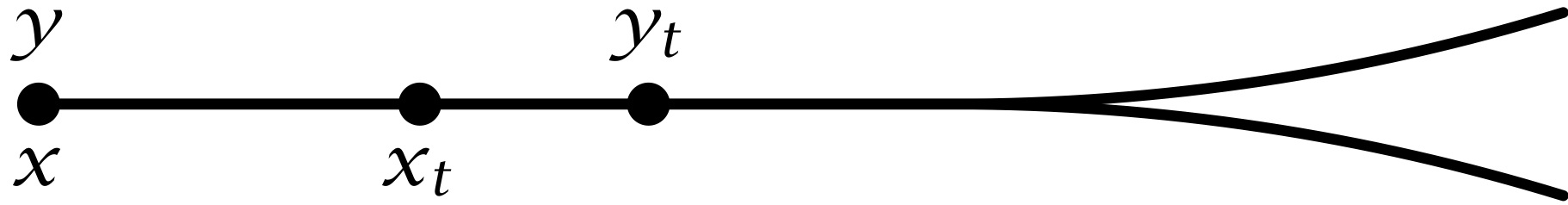}
 \caption{Shearing: If two points start out so close together that we
cannot tell them apart, then the first difference we see will be that one
gets ahead of the other, but (apparently) following the same path. It is
only much later that we will be able to detect any difference between their paths.}
 \label{shearing}
 \end{center}
 \end{figure}
%texpreamble
%("  \usepackage{amsmath}
% \usepackage[LY1]{fontenc}
% \usepackage[expert,LY1,mylucidascale]{mylucidabr}
% ");
%defaultpen(  fontcommand("\normalfont") + fontsize(10) ); 
%
%from graph access *;
%unitsize(1.5cm);
%
%real linethick = 1.5;
%real dotthick = 12;
%dotfactor=12;
%
%real b = 0.3;
%pair x = (0,0), xt = (1.25,0), yt = (2,0);
%pair split = (3,0), xend = (5,-b), yend = (5,b);
%
%dot( x ); dot( xt); dot(yt);
%
%draw( x{E}..{E}split{E}..yend , linewidth(linethick) );
%draw( x{E}..{E}split{E}..xend , linewidth(linethick) );
%
%label( "$x$", x, 2*S); label( "$y$", x, 2*N); 
%label( "$x_t$", xt, 2*S); label( "$y_t$", yt, 2*N); 

\begin{rems} \ 
\noprelistbreak
\begin{enumerate}
\item The only exception to \cref{ShearingSL2R} is that if $q$ is in the centralizer $\czer_G(u^t)$, then $u^t q u^{-t} = q$ for
all~$t$; in this case, the points $x_t$ and~$y_t$ simply move along
together at exactly the same speed, with no \term{relative motion}.
\item  In contrast to the above discussion of~$u^t$, 
 \begin{itemize}
 \item the matrix $a^t$ is diagonal, 
 but
 \item the \term[largest!term]{largest entry} in 
 	$$ a^t q a^{-t} = \begin{bmatrix}
 \mathsf{a} & e^{2t} \mathsf{b} \\
e^{-2t} \mathsf{c}  & \mathsf{d} 
 \end{bmatrix} $$
is an off-diagonal entry,
 \end{itemize}
 so, under the action of the diagonal group, points move apart (at
exponential speed) in a direction transverse to the orbits
\csee{expdivfig}. 
\end{enumerate}
 \end{rems}

 \begin{figure}[ht]
 \begin{center}
 \includegraphics{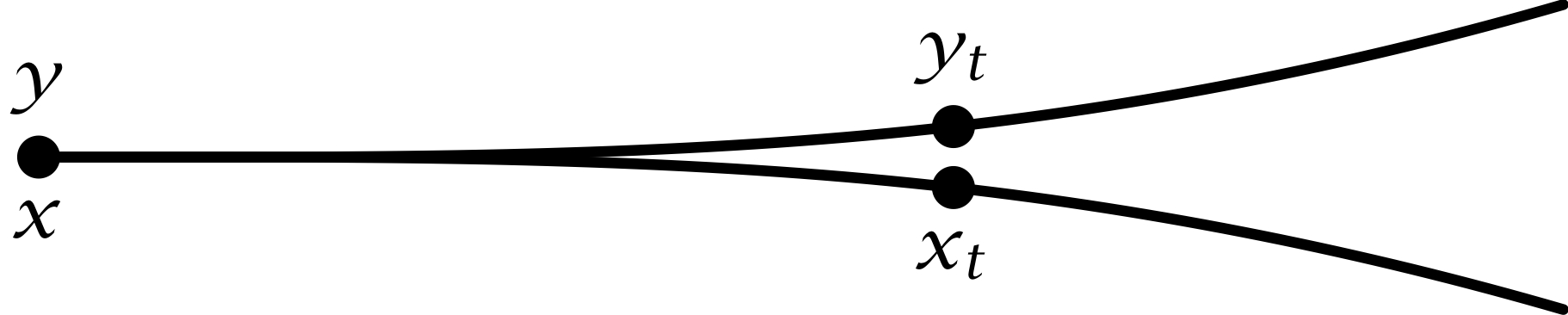}
 \caption{Divergence under a diagonal subgroup: when two points start out so close
together that we cannot tell them apart, the first difference we see
will be in a direction transverse to the orbits.}
 \label{expdivfig}
 \end{center}
 \end{figure}
%texpreamble
%("  \usepackage{amsmath}
% \usepackage[LY1]{fontenc}
% \usepackage[expert,LY1,mylucidascale]{mylucidabr}
% ");
%defaultpen(  fontcommand("\normalfont") + fontsize(10) ); 
%
%from graph access *;
%unitsize(1.5cm);
%
%real linethick = 1.5;
%real dotthick = 12;
%dotfactor=12;
%
%real a = 0.1, b = 0.5;
%pair x = (0,0), xt = (3,-a), yt = (3,a);
%pair xend = (5,-b), yend = (5,b);
%
%dot( x ); dot( xt); dot(yt);
%
%draw( x{E}..yt..yend , linewidth(linethick) );
%draw( x{E}..xt..xend , linewidth(linethick) );
%
%label( "$x$", x, 2*S); label( "$y$", x, 2*N); 
%label( "$x_t$", xt, 2*S); label( "$y_t$", yt, 2*N); 

The Shearing Property \pref{ShearingSL2R} shows that the direction of \term[relative
motion!fastest]{fastest relative motion} is along~$u^t$. However, in the
proof of Ratner's Theorems, it turns out that we wish to ignore
motion \emph{along} the orbits, and consider, instead, only the
\index{component, transverse}component of the relative motion that is
\defit[transverse divergence]{transverse} (or perpendicular) to the
orbits of~$u^t$. This direction, by definition, does not belong to
$\{u^t\}$. 

\begin{defn}
Suppose, as before, that $x$ and~$y$ are two points in
$X$ with $x \approx y$. Then, by continuity, $x_t \approx
y_t$ for a long time. Eventually, we will be able to see a difference
between $x_t$ and~$y_t$. The \term[shearing property!in
$\SL(2,\real)$]{Shearing Property} \pref{ShearingSL2R} tells us that, when
this first happens, $y_t$ will be indistinguishable from some point on the
orbit of~$x$; that is, $y_t \approx x_{t'}$ for some~$t'$. This will
continue for another long time (with $t'$ some function of~$t$), but we can
expect that $y_t$ will eventually diverge from the orbit of~$x$ --- this
is \defit[transverse divergence]{transverse divergence}. (Note that this
transverse divergence is a \term[second-order effect]{second-order
effect}; it is only apparent after we mod out the relative motion along
the orbit.) Letting $x_{t'}$ be the point on the orbit of~$x$ that is
closest to~$y_t$, we write $y_t = g x_{t'}$ for some $g \in G$. Then $g -
\Id$ represents the transverse divergence. When this transverse divergence
first becomes macroscopic, we wish to understand which of the matrix
entries of $g - \Id$ are macroscopic.
 \end{defn}

 In the matrix on the right-hand side of \cref{ConjByU}, we have already observed
that the \term[largest!entry]{largest entry} is in the bottom left corner,
the direction of~$\{u^t\}$. If we ignore that entry, then the two diagonal
entries are the largest of the remaining entries. The diagonal corresponds to the
subgroup~$\{a^t\}$. Therefore, the fastest
\defit[transverse divergence]{transverse} divergence is in the
direction of $\{a^t\}$. Notice that $\{a^t\}$
\term[normalizer]{normalizes}~$\{u^t\}$.

\begin{prop} \label{High1DTransN}
 The fastest \defit[transverse divergence]{transverse} motion is
along some direction in the \term{normalizer} of~$u^t$.

More precisely, if $x,y \in X$, with $x \approx y$, and $r > 0$ is much smaller than the injectivity radius of~$X$, then either:
	\begin{enumerate}
	\item there exist large $t,t' \in \real$ and $g \in \nzer_G\bigl( \{u^t\} \bigr)$ such that
	$$ \text{$u^t y \approx g u^{t'} x$ and\/ $\| g \| = d \bigl(g, \{u^t\} \bigr) = r$,} $$
	or
	\item for all $t \in \real$, there exists $t' \in \real$, such that $u^t y \approx u^{t'} x$ \textup(i.e., there is no transverse motion, only shearing\textup).
	\end{enumerate}
 \end{prop}

To illustrate how understanding the transverse motion can be useful, let us prove a very special case of Ratner's Orbit-Closure Theorem \pref{Ratner-OrbitClosure}.

\begin{prop} \label{MinSetIsUOrbit}
Let $C = \closure{\{u^t\}x}$, for some $x \in X$, and assume
	\begin{itemize}
	\item $C$ is a \defit[minimal!closed, invariant set]{minimal} $u^t$-invariant closed subset of~$X$ 
	\textup(this means that no nonempty, proper, closed subset of~$C$ is\/ $\{u^t\}$-invariant\/\textup),
	and
	\item $\{\, g \in G \mid gC = C \,\} =  \{u^t\}$.
	\end{itemize}
Then $C = \{u^t\}x$, so $C$~is a submanifold of~$X$.
\end{prop}

\begin{proof}
We wish to show $C \subseteq \{u^t\}x$, but \cref{MinSetIsUOrbit-CcapCg=0Ex} implies that it suffices to prove only the weaker statement that $C \subseteq \nzer_G\bigl( \{u^t\} \bigr) \, x$.

Suppose $C \not\subseteq \nzer_G\bigl( \{u^t\} \bigr) \, x$. Then, since $C$ is connected, there exists $y \in C$, with $y \approx x$, but $y \notin \nzer_G\bigl( \{u^t\} \bigr) \, x$.
From \cref{High1DTransN}, we see that there exist $t,t' \in \real$ and 
	$$ \text{$g \in \nzer_G\bigl( \{u^t\} \bigr)$, with $g \notin \{u^t\}$, such that
	$ u^t y \approx g u^{t'} x $.} $$
For simplicity, let us pretend that 
	$$\text{$u^t y$ is \emph{equal} to $g u^{t'} x$} , $$
rather than merely being approximately equal \csee{MinSetIsUOrbit-xu=yuEx}. Then we have $C \cap gC \neq \emptyset$ (because $u^t y \in C$ and $g u^{t'} x \in g C$).
This contradicts \cref{MinSetIsUOrbit-CcapCg=0Ex}.
\end{proof}

\begin{exercises}

\item \label{MinSetIsUOrbit-CcapCg=0Ex}
Under the assumptions of \cref{MinSetIsUOrbit}, show: 
	$$ \text{if $g \in \nzer_G\bigl( \{u^t\} \bigr)$, but $g \notin \{u^t\}$, then $C \cap gC = \emptyset$.} $$ (In particular, $g x \notin C$.)
\hint{$gC$ is $u^t$-invariant (because $g \in \nzer_G\bigl( \{u^t\} \bigr)$), so $C \cap gC$ is a $u^t$-invariant subset of the minimal set~$C$.}

\item \label{MinSetIsUOrbit-xu=yuEx}
Complete the proof of \cref{MinSetIsUOrbit} by eliminating the pretense that $u^t y$ is \emph{equal} to $g u^{t'} x$.
\hint{The compact sets $C$ and $\{\, g \in \nzer_G\bigl( \{u^t\} \bigr) \mid \|g\| = d \bigl( g,\{u^t\} \bigr) = r \,\} \cdot C$ are disjoint, so it is impossible for a point in one set to be arbitrarily close to a point in the other set.}

\end{exercises}

\begin{notes}

\Cref{Ratner-OrbitClosure} is due to M.\,Ratner \cite{Ratner-OrbitClosure}, under the assumption that $V$ is either unipotent or connected. 
(A shorter proof can be found in \cite{MargulisTomanov-InvtMeas}.)
This additional hypothesis was removed by N.\,Shah \cite{Shah-GenByUnip} (except for a technical problem involving Conclusions \pref{Ratner-OrbitClosure-AlmConn} and~\pref{Ratner-OrbitClosure-probmeas} that was resolved in \cite[Cor.~3.5.4]{KleinbockShahStarkov}).

See \cite{Morris-RatnersThms} for a more thorough introduction to Ratner's Theorems, their proofs, and some applications.

See \cite[Lem.~2]{Starkov-StructOrbAndRagConj} for the construction of orbits whose closure is not a submanifold, demonstrating the pathology in \cref{MustBeUnip,FractalInHyper3Mfld}.
	%(originally Morse?)
However, it was conjectured by G.\,A.\,Margulis \cite[\S1.1]{Margulis-problems} that certain analogues of Ratner's Theorems are valid in some situations where the subgroup~$V$ is a split torus of dimension~$> 1$; see \cite[\S4.4c]{KleinbockShahStarkov} and~\cite{Maucourant-Nonhomog} for references on this open problem and its applications.

\Cref{GoodClosure(hyper)} was proved by N.\,Shah \cite{Shah-TotallyGeodesic}. The generalization in \fullcref{GoodClosure(LocSymm)}{closure} is due to T.\,Payne \cite{Payne-Closures}.

\Cref{OppenheimConj} is due to G.\,A.\,Margulis \cite{Margulis-Oppenheim}.
See \cite{Margulis-OppenheimSurvey} for a survey of its history and later related developments.

\Cref{ProdLattsDense} was discovered by N.\,Shah \cite[Cor.~1.5]{Shah-GenByUnip}. This consequence of Ratner's Theorem played an important role in \cite{Vatsal-Heegner}.

\Cref{Ratner-Equidistribution} is due to M.\,Ratner \cite{Ratner-OrbitClosure}.

\Cref{Ratner-MeasClass} was proved by M.\,Ratner \cite{Ratner-MeasClass} in the case where $V$ is either unipotent or connected. 
(See \cite{Einsiedler-RatnerSL2R} for a shorter and more self-contained proof in the case where $V \iso \SL(2,\real)$.)
The general case is due to N.\,Shah \cite{Shah-GenByUnip}.

\Cref{UnipotentIsomorphismRigidity} was proved by M.\,Ratner \cite{Ratner-UnipRigidity} if $G_1 \iso G_2 \iso \SL(2,\real)$. The general case is due to D.\,Witte \cite{Witte-UnipRigidity}.

\Cref{UniqErg->UnifDist} is a special case of a classical result in Ergodic Theory that can be found in textbooks such as \cite[Thm.~6.19]{Walters}.

See \cite{BenoistQuint-StatMeas3} for the work of Y.\,Benoist and J.--F.\,Quint mentioned in \fullcref{MeasClassRems}{BenoistQuint}. Shah's suggestion about Zariski closures appears in \cite[end of \S1, p.~232]{Shah-GenByUnip}.

The discussion of shearing in \cref{RatnerShearingSect} is excerpted from \cite[\S1.5]{Morris-RatnersThms}, except that \cref{MinSetIsUOrbit} is a variant of \cite[Prop.~1.6.10]{Morris-RatnersThms}.

\end{notes}

\begingroup 
%!TEX root = IntroArithGrps.tex

\part*{Appendices} % \label{AppendixesPart}

% don't want \S on sections in the appendix
\immediate\addtocontents{toc}{\protect\notappendixtocfalse}

\startappendix  

 %!TEX root = IntroArithGrps.tex

\mychapter{Basic Facts about Semisimple Lie Groups}
\label{SSGrpsChap}

\section{Definitions} \label{SSLieGrpDefnSect}

We are interested in groups of matrices that are (topologically) closed:

\begin{defns} \label{LieDefn} \ 
\noprelistbreak
	\begin{enumerate}
	\item Let $\Mat_{\ell \times \ell}(\real)$ 
	\nindex{$\Mat_{\ell \times \ell}(\real)$ = $\{ \text{$\ell \times \ell$ matrices with real entries} \}$}%
	be the set of all $\ell \times \ell$ matrices with real entries. This has a natural topology, obtained by identifying it with the Euclidean space $\real^{\ell^2}$. 
	%This is a vector space over~$\real$.
	
	\item Let 
	\nindex{$\SL(\ell,\real)$ = $\{ \text{$\ell \times \ell$ matrices of determinant~$1$} \}$}%
	$\SL(\ell,\real) = \{\, g \in \Mat_{\ell \times \ell}(\real) \mid \det g = 1 \,\}$. This is a group under matrix multiplication \csee{SLisGrp}, and it is a closed subset of $\Mat_{\ell \times \ell}(\real)$ \csee{SLisClosed}.

	\item \label{LieDefn-LieGrp}
	A \defit[Lie!group]{Lie group} is any (topologically) closed subgroup of some $\SL(\ell,\real)$. 
	\end{enumerate}
\end{defns}

Recall that an abstract group is \defit[simple!abstract
group]{simple} if it has no nontrivial, proper, normal
subgroups. For Lie groups, we relax this to allow normal
subgroups that are discrete (except that the
one-dimensional abelian groups $\real$ and~$\torus$ are not
considered to be simple).

\begin{defn} \label{simpleDefn}
A Lie group~$G$ is \defit[simple!Lie group]{simple} if it has no
nontrivial, connected, closed, proper, normal subgroups,
and $G$~is not abelian.
 \end{defn}

\begin{eg}
It can be shown that $G = \SL(\ell,\real)$ is a simple Lie group (when $\ell > 1$). If $\ell$ is even, then $\{\pm \Id\}$ is a subgroup of~$G$, and it is normal,
but, because this subgroup is not connected, it does not
disqualify~$G$ from being simple as a Lie group.
 \end{eg}

\begin{rem}
Although \cref{simpleDefn} only refers to
\emph{closed} normal subgroups, it turns out that, except
for the center, there are no normal subgroups at all:
if $G$ is simple, then every proper, normal subgroup of~$G$ is
contained in the center of~$G$.
\end{rem}

\begin{terminology}
 Some authors say that $\SL(n,\real)$ is
\defit[simple!almost]{almost simple}, and reserve
the term ``simple'' for groups that have no (closed) normal
subgroups at all, not even finite ones.
 \end{terminology}

A Lie group is said to be \emph{semisimple} if it is a direct product of simple groups, modulo passing to a finite-index subgroup and/or modding out a finite group: 

\begin{defns} \label{SSDefn} \ 
\noprelistbreak
	\begin{enumerate}
	\item $G_1$ is \defit{isogenous} to~$G_2$ if there is a finite, normal subgroup~$N_i$ of a finite-index subgroup $G_i'$ of~$G_i$, for $i = 1,2$, such that $G_1'/N_1$ is isomorphic to $G_2'/N_2$.
	\item \label{SSDefn-SS}
	$G$ is \defit[semisimple!Lie group]{semisimple} if it is
isogenous to a direct product of simple Lie groups. That
is, $G$ is isogenous to $G_1 \times \cdots \times G_r$,
where each $G_i$ is simple.
	\end{enumerate}
 \end{defns}
 
 \begin{eg}
 $\SL(2,\real) \times \SL(3,\real)$ is a semisimple Lie group that is not simple (because $\SL(2,\real)$ and $\SL(3,\real)$ are normal subgroups).
 \end{eg}
 
 \begin{rem}[\csee{SSHasNoCenter}]
 If $G$ is semisimple, then the center of~$G$ is finite.
 \end{rem}

\begin{assump}[\ccf{standassump}] \label{FinCompsAssump}
Now that we have the definition of a semisimple group, we will henceforth assume in this chapter that the symbol~$G$ always denotes a semisimple Lie group with only finitely many connected components (but the symbol~$\Gamma$ will never appear).
\end{assump}

\begin{warn}
A \emph{Lie group} is usually defined to be any group that is also a smooth manifold, such that the group operations are~$C^\infty$ functions. \fullCref{HomosAreSmooth}{subgrp} below shows that every closed subgroup of $\SL(\ell,\real)$ is a Lie group in the usual sense. However, the converse is false: not every Lie group (in the usual sense) can be realized as a subgroup of some $\SL(\ell,\real)$. (In other words, not every Lie group is \defit[Lie!group!linear]{linear}.) Therefore, our \fullcref{LieDefn}{LieGrp} is more restrictive than the usual definition.
% (That is, the maps $H \times H \to H \colon (g,h) \mapsto gh$ and $H \to H \colon h \mapsto h^{-1}$ are $C^\infty$.) 
%(Manifolds are assumed to have only countably many connected components, so they are second countable.)
(However, every connected Lie group is ``locally isomorphic'' to a linear Lie group.) 
\end{warn}

\begin{exercises}

\item \label{SLisGrp}
Show that $\SL(\ell,\real)$ is a group under matrix multiplication.
\hint{You may assume (without proof) basic facts of linear algebra, such as the fact that a square matrix is invertible if and only if its determinant is not~$0$.}

\item \label{SLisClosed}
Show that $\SL(\ell,\real)$ is a closed subset of $\Mat_{\ell\times\ell}(\real)$.
\hint{For a continuous function, the inverse image of a closed set is closed.}

\item Recall that
	\nindex{$\GL(\ell,\real)$ = $\{ \text{invertible $\ell \times \ell$ matrices} \}$}%
	$\GL(\ell,\real) = \{\, g \in \Mat_{\ell\times\ell}(\real) \mid \det g \neq 0 \,\}$, and that this is a group under matrix multiplication. Show that it is (isomorphic to) a Lie group, by showing it is isomorphic to a closed subgroup of $\SL(\ell+1,\real)$.

\item \label{Normal=ProdGi}
Suppose $G = G_1 \times \cdots \times G_r$, where each $G_i$~is
simple, and $N$~is a connected, closed, normal subgroup
of~$G$. Show there is a subset~$S$ of
$\{1,\ldots,r\}$, such that $N = \prod_{i \in S} G_i$.
\hint{If the projection of~$N$ to~$G_i$ is all of~$G_i$, then $G_i = [G_i,G_i] = [N, G_i] \subseteq N$.}

\item Show that if $G$ is semisimple, and $N$ is any closed,
normal subgroup of~$G$, then $G/N$ is semisimple.
\hint{\Cref{Normal=ProdGi}.}

\item \label{G=NxH}
 Suppose $N$ is a connected, closed, normal subgroup
of $G$. Show that there is a connected, closed, normal
subgroup~$H$ of $G$, such that $G$~is isogenous to $N \times
H$.
\hint{\Cref{Normal=ProdGi}.}

\end{exercises}

\section{The simple Lie groups}

It is clear from \fullcref{SSDefn}{SS} that the study of semisimple groups requires a good understanding of the simple groups.
Probably the most elementary examples of simple Lie groups are special linear groups and orthogonal groups, but symplectic groups and unitary groups are also fundamental. A group of any of these types is called ``classical\zz.'' (The other simple groups are ``\index{exceptional Lie group}exceptional\zz,'' and are less easy to construct.)

\begin{defn} \label{ClassicalDefn}
 $G$ is a \defit[classical!group]{classical group} if it is isogenous to the
direct product of any collection of the groups constructed in
\cref{classical-fulllinear,classical-orthogonal} below. That is, each simple
factor of~$G$ is either a special linear group or the
isometry group of a bilinear, Hermitian, or skew-Hermitian form, over
$\real$, $\complex$, or~$\quaternion$ (where $\quaternion$
is the algebra of quaternions). 
%See \cref{SOmnIsom,SUmnIsom,SpIsom}. 
 \end{defn}

%\begin{warn}
% Contrary to the usage of some authors, we do not
%require a form to be positive-definite in order to be
%called ``Hermitian\zz.''
% \end{warn}

\begin{notation}
Let
\noprelistbreak
% For natural numbers $m$ and~$n$, let  
 \begin{itemize}
 \item \nindex{$g^\transpose$ = transpose of the matrix~$g$}$g^\transpose$ denote the transpose of the
matrix~$g$, 
 \item $g^*$ denote the adjoint (that
is, the conjugate-transpose) of~$g$,
	\item $G^\circ$ denote the identity component of the Lie group~$G$,
	and

%	 \item $I_{m,n}$ be the $(m
%+ n) \times (m + n)$ diagonal matrix whose diagonal entries
%are $m$~$1$'s followed by $n$~$-1$'s,
	 \item \nindex{$\diag(a_1,a_2\ldots,a_n)$ = diagonal matrix with entries $a_1,\ldots,a_n$}
	 $I_{m,n} = \diag(1,1,\ldots,1,-1,-1,\ldots,-1) \in \Mat_{(m+n)\times(m+n)}(\real)$, 
	 \\ % @@@
	 where the number of~$1$'s is~$m$, and the number of~$-1$'s is~$n$.

 \end{itemize}
 \end{notation}

\begin{eg} \label{classical-fulllinear} \ 
\noprelistbreak \nindex{$\SL(n,\real)$, $\SO(m,n)$, $\SU(m,n)$, $\Sp(2m,\real)$: classical Lie groups}
 \begin{enumerate}
 \item The \defit{special linear group}\/ $\SL(n,\real)$ is a
simple Lie group (if $n \ge 2$). 
%It is connected.

 \item \defit[orthogonal!group, special]{Special orthogonal
group}. Let
 $$ \SO(m,n) = \{\, g \in \SL(m+n,\real) \mid g^\transpose
I_{m,n} \, g = I_{m,n} \,\} .$$
This is always semisimple (if $m + n \ge 3$).
It may not be connected, but the identity component $\SO(m,n)^\circ$ is simple if either $m+n = 3$ or $m+n \ge 5$. 
(Furthermore, the index of $\SO(m,n)^\circ$ in $\SO(m,n)$ is $\le 2$.)
%However, $\SO(m+n)$ is
%not connected unless either $m = 0$ or $n = 0$
%\ccf{SO1nNotConn,SOmnNotConn}; otherwise, it has
%exactly two components. 
%See \cref{SO(m+n=4)} for a discussion of the special case where $m + n = 4$.

We use $\SO(n)$ to denote
$\SO(n,0)$ (or $\SO(0,n)$, which is the same group).

 \item \defit[unitary!group, special]{Special unitary
group}. Let
 $$ \SU(m,n) = \{\, g \in \SL(m+n,\complex) \mid g^*
I_{m,n} \, g = I_{m,n} \,\} .$$
 Then $\SU(m,n)$ is simple if $m+n \ge 2$. 
 	%It is connected. 

We use $\SU(n)$ to denote $\SU(n,0)$ (or $\SU(0,n)$).

 \item \defit[symplectic group]{Symplectic group}. Let 
 $$J_{2m} = \begin{pmatrix}
 0 & \Id_{m\times m} \\
 -\Id_{m\times m} & 0
 \end{pmatrix}
 \in \GL(2m,\real) $$
 (where $\Id_{m\times m}$ denotes the $m \times m$ identity
matrix),
 and let 
 $$ \Sp(2m,\real) = \{\, g \in \SL(2m,\real) \mid
g^\transpose J_{2m} \, g = J_{2m} \,\} .$$
 Then $\Sp(2m,\real)$ is simple if $m \ge 1$. 
 %It is connected.
 \end{enumerate}
 \end{eg}

\begin{eg} \label{classical-orthogonal}
 Additional simple groups can be constructed by replacing
the field~$\real$ with either the field~$\complex$ of complex
numbers or the division ring~$\quaternion$ of quaternions:
 \begin{enumerate}
 	\nindex{$\SL(n,\complex)$, $\SL(n,\quaternion)$, $\SO(n,\complex)$, $\SO(n,\quaternion)$, $\Sp(2m,\complex)$, $\Sp(m,n)$: more classical groups}

 \item \defit[special linear group]{Complex and
quaternionic special linear groups}: $\SL(n,\complex)$ and
$\SL(n,\quaternion)$ are simple Lie groups (if $n \ge 2$).
%Each is connected.

 {\it Note:} The noncommutativity of~$\quaternion$ causes
some difficulty in defining the determinant of a
quaternionic matrix. To avoid this problem, we define the
\defit{reduced norm} of a quaternionic $n \times n$
matrix~$g$ to be the determinant of the $2n \times 2n$
complex matrix obtained by identifying $\quaternion^n$
with~$\complex^{2n}$. Then, by definition, $g$ belongs to
$\SL(n,\quaternion)$ if and only if its reduced norm
is~$1$. It is not difficult to see that the reduced norm of a
quaternionic matrix is always a (nonnegative) real number \csee{RedNormPosEx}.

 \item \label{classical-orthogonal-SOCH}
 \defit[orthogonal!group, special]{Complex and
quaternionic special orthogonal groups}: 
 $$\SO(n,\complex) 
 = \{\, g \in \SL(n,\complex) \mid g^\transpose \Id\, g = \Id
\,\} $$
 and
 $$ \SO(n,\quaternion)
 = \{\, g \in \SL(n,\quaternion) \mid \tau_r(g^\transpose)
\Id\, g = \Id \,\} ,$$
 where $\tau_r$ is the \defit[reversion anti-involution]{reversion} on~$\quaternion$
defined by
 $$ \tau_r(a_0 + a_1 i + a_2 j + a_3 k)
 =  a_0 + a_1 i - a_2 j + a_3 k .$$
 (Note that $\tau_r(ab) = \tau_r(b) \, \tau_r(a)$
\csee{tau=Antiaut}; $\tau_r$ is included in the
definition of $\SO(n,\quaternion)$ in order to compensate
for the noncommutativity of~$\quaternion$
\csee{tau(transpose)}.)

 \item \defit[symplectic group!complex]{Complex symplectic group}:
Let
 $$\Sp(2m,\complex) 
 = \{\, g \in \SL(2m,\complex)
 \mid g^\transpose J_{2m} \, g = J_{2m} \,\} .$$

 \item \defit[symplectic group!unitary]{Symplectic unitary
groups}: Let
 $$\Sp(m,n)
 = \{\, g \in \SL(m+n,\quaternion)
 \mid g^* I_{m,n} \, g = I_{m,n} \,\}
 .$$
 Here, as usual, $g^*$ denotes the conjugate-transpose of~$g$; recall that the \defit[conjugate (of a quaternion)]{conjugate} of a quaternion is defined by%
 \nindex{$\overline{x}$ = conjugate of the quaternion~$x$}% no page break here !!!
 	$$ \overline{a + bi + cj + dk} = a - bi - cj - dk $$
(and that $\overline{xy} = \overline{y} \, \overline{x}$).
  We use $\Sp(n)$ to denote $\Sp(n,0)$ (or $\Sp(0,n)$).
 \end{enumerate}
 \end{eg}

\begin{terminology} \label{SOstarTerminology}
 Some authors use
 \begin{itemize}
 \item  $\SU^*(2n)$ to denote $\SL(n,\quaternion)$,
 \item $\SO^*(2n)$ to denote $\SO(n,\quaternion)$,
 or
 \item $\Sp(n,\real)$ to denote $\Sp(2n,\real)$.
 \end{itemize}
 \end{terminology}

\begin{rem} \label{SL2RinG}
 $\SL(2,\real)$ is the smallest connected, noncompact, simple
Lie group; it is contained (up to isogeny) in any other. For example:
 \begin{enumerate}
 
 \item If $\SL(n,\real)$, $\SL(n,\complex)$, or
$\SL(n,\quaternion)$ is not compact, then $n \ge 2$, so the group contains $\SL(2,\real)$.

 \item If $\SO(m,n)$ is semisimple and  not compact, then
$\min\{m,n\} \ge 1$ and $\max\{m,n\} \ge 2$, so it contains $\SO(1,2)$, which is isogenous to $\SL(2,\real)$.

 \item If $\SU(m,n)$ or $\Sp(m,n)$ is not compact, then $\min\{m,n\} \ge
1$, so the group contains $\SU(1,1)$, which is isogenous to $\SL(2,\real)$.
 
\item $\Sp(2m,\real)$ and $\Sp(2m,\complex)$ both contain $\Sp(2,\real)$, which is equal to $\SL(2,\real)$.

\item If $\SO(n,\complex)$ is semisimple and not compact, then $n \ge 3$, so the group contains $\SO(1,2)$, which is isogenous to $\SL(2,\real)$.

\item If $\SO(n,\quaternion)$ is not compact, then $n \ge 2$, so it contains a subgroup conjugate to $\SU(1,1)$, which is isogenous to $\SL(2,\real)$.

 \end{enumerate}
% Note that $\SO(1,2)$ and $\SU(1,1)$ are isogenous to
%$\SL(2,\real)$ \fullsee{isogTypes}{A1SO12}; and $\Sp(1,1)$
%is isogenous to $\SO(1,4)$ \fullsee{isogTypes}{C2Sp11},
%which contains $\SO(1,2)$.
 \end{rem}

%\begin{rem}
% There is some redundancy in the above lists. (For example,
% $\SL(2,\real)$, $\SU(1,1)$, $\SO(1,2)$, and
%$\Sp(2,\real)$ are isogenous to each other
%\fullsee{isogTypes}{A1SO12}.) A complete list of these
%redundancies is given in \cref{ClassicalIsogChap} below.
% \end{rem}

The classical groups are just examples, so one would expect
there to be many other (more exotic) simple Lie groups.
Amazingly, that is not the case --- there are only finitely many others:

\begin{thm}[(\'E.\,Cartan)] \label{RealSimpleGrps}
 Every simple Lie group is isogenous to either 
 \begin{enumerate}
 \item a classical group, or
 \item one of the finitely many exceptional groups.
 \end{enumerate}
 \end{thm}

 See \cref{RFormsOfCGrps,GaloisCohoRealFormsSect} for an indication of the proof of \cref{RealSimpleGrps}.

\begin{exercises}

\item \label{RedNormPosEx} For all nonzero $g \in \Mat_{n \times n}(\quaternion)$, show that the reduced norm of~$g$ is a nonnegative real number.
\hint{Use row and column operations in $\Mat_{n \times n}(\quaternion)$ to reduce to the case where $g$ is upper triangular. For $n = 1$, the reduced norm of~$g$ is $g \, \overline{g}$.}

\item \label{tau=Antiaut}
 In the notation of \cref{classical-orthogonal},
show that $\tau_r(ab) = \tau_r(b) \, \tau_r(a)$ for all $a,b \in
\quaternion$.
 \hint{Calculate explicitly, or note that $\tau_r(x) =
j \, \overline{x}  j^{-1}$ (and $\overline{xy} = \overline{y} \, \overline{x}$).}

%\item Give an example of two matrices $g,h \in
%\SL(n,\quaternion)$, such that $(gh)^\transpose \neq
%h^\transpose g^{\transpose}$ and $(gh)^\transpose \neq
%g^\transpose h^{\transpose}$.

\item \label{tau(transpose)}
 For $g,h \in \Mat_{n \times n}(\quaternion)$, show that
$\tau_r \bigl( (gh)^\transpose \bigr) = \tau_r(
h^\transpose) \, \tau_r( g^{\transpose})$.

\item Show that $\SO(n,\quaternion)$ is a subgroup of
$\SL(n,\quaternion)$.

%\item \label{conjInvSO}
% For any $g \in \Ortho(\ell)$, the map $\phi_g \colon
%\SO(\ell) \to \SO(\ell)$, defined by $\phi_g(x) = g x
%g^{-1}$, is an automorphism of $\SO(n)$. 
% Show that $\ell$ is odd if and only if, for every $g \in
%\Ortho(\ell)$, there exists $h \in \SO(\ell)$, such that
%$\phi_h = \phi_g$.

\end{exercises}

\section{Haar measure}

Standard texts on real analysis construct a translation-invariant measure on~$\real^n$. this is called \defit[Lebesgue!measure]{Lebesgue measure}, but the analogue for other Lie groups is called ``Haar measure:''

\begin{prop}[(\thmindex{Existence of Haar Measure}{Existence and Uniqueness of Haar Measure})]\index{Haar measure} \label{HaarMeasure}
 If $H$ is any Lie group, then there
is a unique\/ \textup(up to a scalar multiple\textup)
$\sigma$-finite Borel measure~$\mu$ on~$H$, such that
 \begin{enumerate}
 \item \label{HaarMeasure-mu(C)}
 $\mu(C)$ is finite, for every compact subset~$C$
of~$H$, 
%\item \sigma(\open) > 0$, for every nonempty open set~$\open$,
and
 \item $\mu(hA) = \mu(A)$, for every Borel subset~$A$ of~$H$,
and every $h \in H$.
 \end{enumerate}
 \end{prop}
 
 \begin{defns} \ 
 \noprelistbreak
 \begin{enumerate}
 \item The measure~$\mu$ of
\cref{HaarMeasure} is called the \defit[Haar
measure]{left Haar measure} on~$H$.
Analogously, there is a unique \defit[Haar measure!right]{right Haar measure}
with $\mu(Ah) = \mu(A)$ \csee{RightHaarMeasure}.
 \item $H$ is \defit[unimodular group]{unimodular} if the left Haar measure is also a right Haar measure. (This means $\mu(hA) = \mu(Ah) = \mu(A)$.)
 \end{enumerate}
 \end{defns}

\begin{rem}
Haar measure is always \defit[regular!measure, inner]{inner regular}: $\mu(A)$ is the supremum of the measures of the compact subsets of~$A$.
\end{rem}

\begin{prop} \label{modularfuncofG}
 There is a continuous homomorphism $\Delta \colon H \to \real^+$, such that, if $\mu$ is any \textup(left or right\textup) Haar measure on~$H$, then
 	$$ \text{$\mu(h A h^{-1}) = \Delta(h)\, \mu(A)$,
	 for all $h \in H$ and any Borel set $A \subseteq H$}
	 . $$
 \end{prop}

\begin{proof}
Let $\mu$ be a left Haar measure.
 For each $h \in H$, define $\phi_h \colon H \to H$ by $\phi_h(x) = h x h^{-1}$.
Then $\phi_h$ is an automorphism of~$H$, so $(\phi_h)_*\mu$
is a left Haar measure. By uniqueness, we
conclude that there exists $\Delta(h) \in \real^+$, such that
$(\phi_h)_*\mu = \Delta(h) \, \mu$. It is easy to see that
$\Delta$ is a continuous homomorphism. By using the construction of right Haar measure in \cref{RightHaarMeasure}, it is easy to verify that the same formula also applies to it.
 \end{proof} 

\begin{defn}
 The function~$\Delta$ 
 	\nindex{$\Delta$ = modular function of~$H$}
 defined in \cref{modularfuncofG}
is called the \defit{modular function} of~$H$.
 \end{defn}

\begin{terminology}
 Some authors call $1/\Delta$ the modular function, because
they use the conjugation $h^{-1} A h$, instead of $h A
h^{-1}$.
 \end{terminology}

\begin{cor}
 Let $\Delta$ be the modular function of~$H$, and let $A$ be a Borel subset of~$H$. 
 \begin{enumerate}
 \item If $\mu$ is a right Haar measure on~$H$, then
$\mu(hA) = \Delta(h)\, \mu(A)$, for all $h \in H$.
 \item If $\mu$ is a left Haar measure on~$G$, then
$\mu(Ah) = \Delta(h^{-1})\, \mu(A)$, for all $h \in H$.
 \item $H$ is unimodular if and only if $\Delta(h) = 1$, for
all $h \in H$. 
 \item $\Delta(h) = | \det (\Ad_H h)|$ for all $h \in H$ \csee{AdDefn}.
 \end{enumerate}
 \end{cor}

\begin{rem} \label{SS->unimod}
 $G$ is unimodular, because semisimple groups have no nontrivial (continuous)
homomorphisms to $\real^+$ \csee{SS->perfect}.
 \end{rem}

\begin{prop}
 Let $\mu$ be a left Haar measure on a Lie group~$H$. Then $\mu(H) < \infty$ if and only if $H$ is compact.
 \end{prop}

\begin{proof}
 ($\Leftarrow$) See \fullcref{HaarMeasure}{mu(C)}.
 %Haar measure is finite on compact sets \fullcsee{HaarMeasure}{mu(C)}.

($\Rightarrow$) Since $\mu(H) < \infty$ (and the measure $\mu$ is inner regular),
there is a compact subset~$C$ of~$H$, such that $\mu(C) > \mu(H)/2$. 
Then, for any $h \in H$, we have
	$$ \mu(hC) + \mu(C) = \mu(C) + \mu(C) = 2 \mu(C) > \mu(H) ,$$
so $hC$ cannot be disjoint from~$C$. This implies that $h$~belongs to the set $C \cdot C^{-1}$, which is compact. Since $h$~is an arbitrary element of~$H$, we conclude that $H = C \cdot C^{-1}$ is compact.
%
%We prove the contrapositive. Let $C$ be a
%compact subset of nonzero measure. Because $G \times G \to G
%\colon (g,h) \mapsto gh^{-1}$ is continuous, and the
%continuous image of a compact set is compact, we know $C
%C^{-1}$ is compact. Since $G$ is not compact, then there
%exists $g_1 \notin C C^{-1}$; therefore $g_1 C$ is disjoint
%from~$C$. Continuing, we construct, by induction on~$n$, a
%sequence $\{g_n\}$ of elements of~$G$, such that $\{g_n C\}$
%is a collection of pairwise disjoint sets. They all have the
%same measure (since $\mu$ is $G$-invariant), so we conclude
%that 
% $$ \mu(G) \ge \mu \left( \bigcup_{n=1}^\infty g_n C \right) =
%\infty .$$
 \end{proof}

\begin{exercises}

\item Prove the existence (but not uniqueness) of Haar measure on~$H$, without using \cref{HaarMeasure}, under the additional assumption that the Lie group~$H$ is a $C^\infty$ submanifold of $\SL(\ell,\real)$ \fullccf{HomosAreSmooth}{subgrp}.
\hint{For $k = \dim H$, there is a differential $k$-form on~$H$ that is invariant under left translations.}

\item \label{RightHaarMeasure}
Suppose $\mu$ is a left Haar measure on~$H$, and define $\widetilde\mu(A) = \mu(A^{-1})$. Show $\widetilde\mu$ is a right Haar measure.

\item \label{SS->perfect}
% \begin{enumerate}
% \item 
Assume $G$ is connected. Show that if $\phi \colon G \to A$ is a continuous homomorphism, and $A$~is abelian, then $\phi$~is trivial.
% \item Show that $[G,G] = G$.
% \end{enumerate}
\hint{The kernel of a continuous homomorphism is a closed, normal subgroup.}

\end{exercises}

\section{\texorpdfstring{$G$}{G} is almost Zariski closed} \label{ZariskiSect}

%Every classical group has only
%finitely many connected components. This is a special case
%of the following much more general result.

\begin{defns} \  \label{AlgicGrpDefn}
 \noprelistbreak 
 \begin{enumerate}
 \item We use $\real[x_{1,1}, \ldots, x_{\ell,\ell}]$ to
denote the set of real polynomials in the $\ell^2$ variables
$\{\, x_{i,j} \mid 1 \le i,j \le \ell\, \}$. 
 \item For any $Q \in \real[x_{1,1}, \ldots, x_{\ell,\ell}]$,
and any $g \in \Mat_{\ell\times \ell}(\complex)$, we use
$Q(g)$ to denote the value obtained by substituting the
matrix entries~$g_{i,j}$ into the variables~$x_{i,j}$. For
example, if $Q= x_{1,1} x_{2,2} - x_{1,2} x_{2,1}$, then
$Q(g)$ is the determinant of the first principal $2 \times
2$ minor of~$g$.
 \item For any subset~$\mathcal{Q}$ of
$\real[x_{1,1}, \ldots, x_{\ell,\ell}]$, let
  $$\Var(\mathcal{Q}) = \{\, g \in \SL(\ell,\real) \mid
Q(g) = 0, \ \forall Q \in \mathcal{Q} \,\} .$$
 This is the \defit{variety} associated to~$\mathcal{Q}$.
 \item A subset~$H$ of $\SL(\ell,\real)$ is \defit[Zariski!closed]{Zariski
closed} if there exists a subset~$\mathcal{Q}$ of $\real[x_{1,1}, \ldots,
x_{\ell,\ell}]$, such that $H = \Var(\mathcal{Q})$. (In the
special case where $H$ is a sub\emph{group} of
$\SL(\ell,\real)$, we may also say that $H$ is a \defit[algebraic!group!real]{real
algebraic group} or an 
\defit[algebraic!group!over R@over~$\real$]{algebraic group that is defined over~$\real$}.)
 \item The \defit[Zariski!closure]{Zariski closure} of a subset~$H$ of
$\SL(\ell,\real)$ is the (unique) smallest Zariski closed
subset of $\SL(\ell,\real)$ that contains~$H$. This is sometimes denoted~%
	\nindex{$\Zar{H}$ = Zariski closure of~$H$}$\Zar{H}$.
(It can also be denoted $\overline{H}$, if this will not lead to confusion with the closure of~$H$ in the ordinary topology.)
 \end{enumerate}
 \end{defns}

\begin{eg} \label{EgZarClosed}
 \ 
 \begin{enumerate}
 \item $\SL(\ell,\real)$ is Zariski closed. Let $\mathcal{Q}
= \emptyset$.
 \item \label{EgZarClosed-diag}
 The group of diagonal matrices in $\SL(\ell,\real)$
is Zariski closed. Let 
 $\mathcal{Q} = \{\, x_{i,j} \mid i \neq j \,\}$.
 \item \label{EgZarClosed-central}
 For any $A \in \GL(\ell,\real)$, the centralizer
of~$A$ is Zariski closed. Let
 $$ \mathcal{Q} = \bigset{
 \sum_{k=1}^\ell ( x_{i,k} A_{k,j}
-  A_{i,k} x_{k,j}) 
 }{
 1 \le i,j\le \ell 
 }.$$
 \item \label{EgZarClosed-SLnC}
 If we identify $\SL(n,\complex)$ with a subgroup of
$\SL(2n,\real)$, by identifying $\complex$ with~$\real^2$,
then $\SL(n,\complex)$ is Zariski closed, because it is the
centralizer of~$T_i$, the linear transformation in
$\GL(2n,\real)$ that corresponds to scalar multiplication
by~$i$.
 \item The classical groups of \cref{classical-fulllinear,classical-orthogonal} are Zariski closed (if we
identify $\complex$ with~$\real^2$ and $\quaternion$
with~$\real^4$ where necessary).
 \end{enumerate}
 \end{eg}

\begin{terminology} \ 
\noprelistbreak
 \begin{itemize}
 \item Other authors use $\GL(\ell,\real)$ in the definition
of $\Var(\mathcal{Q})$, instead of $\SL(\ell,\real)$. Our
choice leads to no loss of generality, and simplifies the
theory slightly. (In the $\GL$ theory, one should, for
technical reasons, stipulate that the function $1/\det(g)$ is
considered to be a polynomial. In our setting, $\det g$ is the constant
function~$1$, so this is not an issue.)
 \item What we call $\Var(\mathcal{Q})$ is actually only the
\emph{real} points of the variety. Algebraic geometers
usually consider the solutions in~$\complex$, rather
than~$\real$, but our preoccupation with real Lie groups leads
to our emphasis on real points.
 \end{itemize}
 \end{terminology}

\begin{eg} \label{EgNotZarClosed}
 Let 
 $$ H =
 \bigset{
 \begin{pmatrix}
 e^{t} & 0 & 0 & 0\\
 0 & e^{-t} & 0 & 0\\
 0 & 0 & 1 & t \\
 0 & 0 & 0 & 1 \\
 \end{pmatrix}
 }{
 t \in \real
 }
 \subset \SL(4,\real)
 .$$
 Then $H$ is a 1-dimensional subgroup that is not Zariski
closed. Its Zariski closure is
 $$ \Zar{H} =
 \bigset{
 \begin{pmatrix}
 a & 0 & 0 & 0\\
 0 & 1/a & 0 & 0\\
 0 & 0 & 1 & t \\
 0 & 0 & 0 & 1 \\
 \end{pmatrix}
 }{
 \begin{matrix}
 a \in \real \smallsetminus \{0\}, \\
 t \in \real
 \end{matrix}
 }
 \subset \SL(4,\real)
 .$$
 The point here is that the exponential function is
transcendental, not polynomial, so no polynomial can capture
the relation that ties the diagonal entries to the
off-diagonal entry in~$H$. Therefore, as far as polynomials are
concerned, the diagonal entries in the upper left are independent of the
off-diagonal entry, as we see in the Zariski closure.
 \end{eg}

\begin{rem}
 If $H$ is Zariski closed, then the set~$\mathcal{Q}$ of
\cref{AlgicGrpDefn} can be chosen to be finite
(because the ring $\real[x_{1,1}, \ldots, x_{\ell,\ell}]$ is
Noetherian). 
 \end{rem}

Everyone knows that a (nonzero) polynomial in one variable
has only finitely many roots. The following important fact
generalizes this observation to any collection of
polynomials in any number of variables.

\begin{thm} \label{Zar->AlmConn}
 Every Zariski closed subset of\/ $\SL(\ell,\real)$ has only
finitely many connected components.
 \end{thm}

\begin{defn} \label{AlmZarDefn}
 A closed subgroup~$H$ of $\SL(\ell,\real)$ is \emph{almost
Zariski closed} if it has only finitely many components, and
there is a Zariski closed subgroup~$H_1$ of $\SL(\ell,\real)$,
%(which also has only finitely many components, by \cref{Zar->AlmConn}), 
such that $H^\circ = H_1^\circ$.
In other words, in the terminology of \cref{CommensDefn},
$H$~is \defit{commensurable} to a Zariski closed subgroup.
 \end{defn}

\begin{egs} \ 
\noprelistbreak
 \begin{enumerate}
 \item Let $H$ be the group of diagonal matrices in
$\SL(2,\real)$. Then $H$ is Zariski closed
\fullcsee{EgZarClosed}{diag}, but $H^\circ$ is not: any
polynomial that vanishes on the diagonal matrices with
positive entries will also vanish on the diagonal matrices
with negative entries. So $H^\circ$ is almost Zariski closed,
but it is not Zariski closed.

\item Let $G = \SO(1,2)^\circ$. Then $G$ is almost Zariski
closed (because $\SO(1,2)$ is Zariski closed), but $G$ is
not Zariski closed \csee{SO12notZar}.
 \end{enumerate}
 These examples are typical: a connected Lie group is almost Zariski closed if
and only if it is the identity component of a group that is Zariski closed.
 \end{egs}

The following fact gives the Zariski closure a central role
in the study of semisimple Lie groups.

\begin{thm} \label{GisAlgic}
 If $G \subseteq \SL(\ell,\real)$, then $G$ is almost Zariski closed.
 \end{thm}

\begin{proof}
 Let $\Zar{G}$ be the Zariski closure of~$G$. Then
$\Zar{G}$ is semisimple. (For example, if $G$ is
irreducible in $\SL(\ell,\complex)$, then $\Zar{G}$
is also irreducible, so \cref{irred->SS} below implies that
$\Zar{G}^\circ$ is semisimple.)

Since $G$ has only finitely many connected components \csee{FinCompsAssump}, we may assume, by passing to a subgroup of finite index, that it is connected. This implies that the normalizer $\nzer_{\SL(\ell,\real)}(G)$
is Zariski closed \csee{N(G)ZarClosed}. Therefore
$\Zar{G}$ is contained in the normalizer, which means that $G$ is a normal subgroup
of~$\Zar{G}$. 

Hence (up to isogeny), we have
$\Zar{G} = G \times H$, for some closed, normal
subgroup~$H$ of~$\Zar{G}$ \csee{G=NxH}. 
So $G = \czer_{\Zar{G}}(H)^\circ$ is almost Zariski closed
\fullcsee{EgZarClosed}{central}.
 \end{proof}

\begin{warn}
 \Cref{GisAlgic} relies on our standing
assumption that $G$ is semisimple \csee{EgNotZarClosed}.
(Actually, it suffices to know that, besides being
connected, $G$ is perfect; that is, $G = [G,G]$ is equal to
its commutator subgroup.)
 \end{warn}

\begin{exercises}

\item \label{SO12notZar}
 Show that $\SO(1,2)^\circ$ is not Zariski closed.
 \hint{We have
 $$ \frac{1}{2}
 \begin{pmatrix}
 s + \frac{1}{s} & s - \frac{1}{s} & 0 \\
 s - \frac{1}{s} & s + \frac{1}{s} & 0 \\
 0 & 0 & 2
 \end{pmatrix}
 \in \SO(1,2)^\circ
 \qquad \Leftrightarrow \qquad
 s > 0 .$$
 If a rational function $f \colon \real \smallsetminus \{0\}
\to \real$ vanishes on~$\real^+$, then it also vanishes
on~$\real^-$.}

\item \label{N(G)ZarClosed}
 Show that if $H$~is a connected Lie subgroup of
$\SL(\ell,\real)$, then the normalizer $\nzer_{\SL(\ell,\real)}(H)$ is Zariski
closed.
 \hint{$g \in \nzer(H)$ if and only if $g \Lie H g^{-1} = \Lie
H$, where $\Lie H \subseteq \Mat_{\ell\times\ell}(\real)$ is
the Lie algebra of~$H$.}

\item Show that if $\Zar{H}$ is the Zariski closure of a
subgroup~$H$ of~$G$, then $g \Zar{H} g^{-1}$ is the
Zariski closure of $g H g^{-1}$, for any $g \in G$.

\item \label{VarClosed}
 Suppose $G$ is a connected subgroup of
$\SL(\ell,\real)$ that is almost Zariski closed, and that $\mathcal{Q} \subset
\real[x_{1,1},\ldots,x_{\ell,\ell}]$.
 \begin{enumerate}
 \item \label{VarClosed-closed}
 Show that $G \cap \Var(\mathcal{Q})$ is a closed subset of $G$.
 \item \label{VarClosed-nodense}
 Show that if $G \not\subseteq \Var(\mathcal{Q})$, then $G
\cap \Var(\mathcal{Q})$ does not contain any nonempty open
subset of~$G$.
 \item \label{VarClosed-meas0}
 Show that if $G \not\subseteq \Var(\mathcal{Q})$, then $G
\cap \Var(\mathcal{Q})$ has measure zero, with respect to the
Haar measure on~$G$.
 \end{enumerate}
 \hint{For \pref{VarClosed-nodense} and \pref{VarClosed-meas0}, you may assume, 
 without proof, % @@@
 that, for some $d$, there exist 
 	$$ \emptyset = \Var(\mathcal{Q}_{-1}) \subseteq \Var(\mathcal{Q}_0) \subseteq \Var(\mathcal{Q}_1) \subseteq \cdots \subseteq \Var(\mathcal{Q}_d) = \Var(\mathcal{Q}) ,$$
such that $G \cap \bigl( \Var(\mathcal{Q}_k) \smallsetminus \Var(\mathcal{Q}_{k-1}) \bigr)$ is a (possibly empty) $k$-dimensional $C^\infty$ submanifold of~$G$, for $0 \le k \le d$. ($G \cap \Var(\mathcal{Q}_{k-1})$ is called the \defit{singular set} of the variety $G \cap \Var(\mathcal{Q}_k)$.)} 

\item Show, for any subspace~$V$ of~$\real^\ell$, that 
 $$ \Stab_{\SL(\ell,\real)}(V) = \{\, g \in \SL(\ell,\real)
\mid gV = V \,\}$$
 is Zariski closed.

\item \label{IrredCompsOfVariety}
 A Zariski-closed subset of $\SL(\ell,\real)$ is
\defit[irreducible!Zariski-closed subset]{irreducible} if it
cannot be written as the union of two Zariski-closed, proper
subsets. Show that every Zariski-closed subset~$A$ of
$\SL(\ell,\real)$ has a unique decomposition as an
irredundant, finite union of irreducible, Zariski-closed
subsets. (By irredundant, we mean that no one of the sets is
contained in the union of the others.)
 \hint{The ascending chain condition on ideals of
$\real[x_{1,1},\ldots,x_{\ell,\ell}]$ implies the descending
chain condition on Zariski-closed subsets, so $A$~can be
written as a finite union of irreducibles. To make the union
irredundant, the irreducible subsets must be maximal.}

\item \label{Conn->ZarConn}
 Let $H$ be a connected subgroup of $\SL(\ell,\real)$. Show
that if $H \subseteq A_1 \cup A_2$, where $A_1$ and~$A_2$ are
Zariski-closed subsets of $\SL(\ell,\real)$, then either $H
\subseteq A_1$ or $H \subseteq A_2$.
 \hint{The Zariski closure $\Zar{H} = B_1 \cup \cdots
\cup B_r$ is an irredundant union of irreducible,
Zariski-closed subsets \csee{IrredCompsOfVariety}. For $h
\in H$, we have $\Zar{H} = hB_1 \cup \cdots \cup hB_r$,
so uniqueness implies that $h$~acts as a permutation of
$\{B_j\}$. Because $H$~is connected, conclude that
$\Zar{H} = B_1$ is irreducible.}

\item \label{ChevalleyStabEx}
Assume $G$ is connected, and $G \subseteq \SL(\ell,\real)$.
Show there exist
	\begin{itemize}
	\item a finite-dimensional real vector space~$V$,
	\item a vector~$v$ in~$V$,
	and
	\item a continuous homomorphism $\rho \colon \SL(\ell,\real) \to \SL(V)$,
	\end{itemize}
such that $G = \Stab_{\SL(\ell,\real)}(v)^\circ$.
\hint{Let $V_n$ be the vector space of polynomial functions on $\SL(\ell,\real)$, and let $W_n$ be the subspace consisting of polynomials that vanish on~$G$. Then $\SL(\ell,\real)$ acts on~$V_n$ by translation, and $W_n$ is $G$-invariant. For $n$ sufficiently large, $W_n$ contains generators of the ideal of all polynomials vanishing on~$G$, so $G = \Stab_{\SL(\ell,\real)}(W_n)^\circ$. Now let $V$ be the exterior power $\bigwedge^d V_n$, where $d = \dim W_n$, and let $v$ be a nonzero vector in $\bigwedge^d W_n$.}

\item \label{SSHasNoCenter}
Show that the center of~$G$ is finite.
\hint{The identity component of the Zariski closure of~$Z(G)$ is a connected, normal subgroup of~$G$.}

\end{exercises}

\section{Three useful theorems}

\subsection{Real Jordan decomposition}

\begin{defn} \label{hypelluniDefn}
 Let $g \in\GL(n,\real)$. We say that $g$ is
 \begin{enumerate}
 \item \defit[semisimple!element of $G$]{semisimple} if $g$ is
diagonalizable (over~$\complex$),
 \item  \defit[hyperbolic!element of~$G$]{hyperbolic} if
 \begin{itemize}
 \item $g$ is semisimple, and
 \item every eigenvalue of~$g$ is real and positive,
 \end{itemize}
 \item \defit[elliptic element of~$G$]{elliptic} if 
 \begin{itemize}
 \item $g$ is semisimple, and
 \item every eigenvalue of~$g$ is on the unit circle
in~$\complex$,
 \end{itemize}
 \item \defit[unipotent!element]{unipotent} (or
\defit[parabolic!element of~$G$]{parabolic}) if $1$~is the
only eigenvalue of~$g$ over~$\complex$.
 \end{enumerate}
 \end{defn}

\begin{rem} \label{SSeltRem}
 A matrix $g$~is semisimple if and only
if the minimal polynomial of~$g$ has no repeated factors.
 \begin{enumerate}
 \item Because its eigenvalues are real, any hyperbolic~$g$
element is diagonalizable over~$\real$. That is, there is
some $h \in \GL(\ell,\real)$, such that $h^{-1} g h$ is a
diagonal matrix.
 \item \label{SSeltRem-cpct}
 An element is elliptic if and only if it is
contained in some compact subgroup of $\GL(\ell,\real)$. In
particular, if $g$ has finite order (that is, if $g^n = \Id$
for some $n > 0$), then $g$ is elliptic.
 \item A matrix $g \in \GL(\ell,\real)$ is unipotent if and
only if the characteristic polynomial of~$g$ is $(x-1)^\ell$.
(That is, $1$~is the only root of the characteristic
polynomial, with multiplicity~$\ell$.) Another way of saying
this is that $g$~is unipotent if and only if $g - \Id$ is
nilpotent (that is, if and only if $(g - \Id)^n = 0$ for some
$n \in \natural$).
 \end{enumerate}
 \end{rem}
 
 \begin{rem}
 \Cref{SL2RinG} implies that if $G$ is not compact, then it contains nontrivial hyperbolic elements, nontrivial elliptic elements, and nontrivial unipotent elements.
\end{rem}

\begin{prop}[{(\thmindex{Jordan Decomposition, real}Real Jordan Decomposition)}] \label{RealJordanDecomp} 
 Any element~$g$ of~$G$ can be written uniquely
as the product $g = aku$ of three \textbf{commuting} elements
$a,k,u$ of~$G$, such that $a$~is hyperbolic, $k$~is elliptic,
and $u$~is unipotent.
 \end{prop}

\subsection{Engel's Theorem on unipotent subgroups}

\begin{defn} \label{unipDefn}
A subgroup~$U$ of $\SL(\ell,\real)$ is said to be \defit[unipotent!subgroup]{unipotent} if all of its elements are unipotent.
\end{defn}

\begin{eg} \label{Nunip}
Let~$N$ be the group of upper-triangular matrices with $1$'s on the diagonal; that is,
	$$N = \left\{
 	 \begin{bmatrix}
 	1 \\
 	 & 1 &  \vbox to 0pt{\vss \hbox to 0pt{\ \Huge $*$\hss}} \\
 	 &  \vbox to 0pt{\vss\hbox to 0pt{\hss\Huge $0$\ }\vss} & \ddots \\
  	& & & 1
  	\end{bmatrix}
  	\right\}
	\subseteq \SL(\ell,\real) 
	. $$
It is obvious that $N$ is unipotent.
\end{eg}

Therefore, it is obvious that every subgroup of~$N$ is unipotent. Conversely:

\begin{thm}[(\thmindex{Engel's}Engel's Theorem)] \label{EngelUnip}
Every unipotent subgroup of\/ $\SL(\ell,\real)$ is conjugate to a subgroup of the group~$N$ of \cref{Nunip}.
\end{thm}

\subsection{Jacobson-Morosov Lemma}

%We now mention (without proof) a fundamental result that is often useful in the study of Lie groups.

%\begin{defn}
%A matrix $g \in \GL(\ell,\real)$ is unipotent if and
%only if the characteristic polynomial of~$g$ is $(x-1)^\ell$.
%(That is, $1$~is the only root of the characteristic
%polynomial, with multiplicity~$\ell$.) Another way of saying
%this is that $g$~is unipotent if and only if $g - \Id$ is
%nilpotent (that is, if and only if $(g - \Id)^n = 0$ for some
%$n \in \natural$).
%\end{defn}

\begin{thm}[(\thmindex{Jacobson-Morosov}{Jacobson-Morosov Lemma})]
\label{JacobsonMorosov}
 For every unipotent element~$u$ of~$G$, there is a
subgroup~$H$ of~$G$ isogenous to $\SL(2,\real)$, such that
$u \in H$.
\end{thm}

\begin{exercises}

\item \label{NunipEx}
Show that an element of $\SL(\ell,\real)$ is unipotent if and only if it is conjugate to an element of the subgroup~$N$ of \cref{Nunip}.
	\hint{If $g$ is unipotent, then all of its eigenvalues are real, so it can be triangularized over~$\real$.}

\item Show that the Zariski closure of every unipotent subgroup is unipotent. 
%(See \cref{UnipExpLog} for other useful facts about unipotent groups.)
%
%\item \label{UnipExpLog}
%Let 
%	\begin{itemize}
%	\item $N$ be as in \cref{Nunip},
%	\item $\Lie N$ be the space of strictly upper-triangular matrices in $\SL(\ell,\real)$ (with $0$'s on the diagonal),
%	\item $\exp \colon \Lie N \to N$ be defined by $\exp x = \Id + x + x^2 + \cdots + x^\ell$,
%	and
%	\item $\log \colon N \to \Lie N$ be defined by $\Log u = \hat u - \frac{1}{2} \hat u^2 + \frac{1}{3} \hat u^3 - \cdots \pm \hat u^\ell$, where $\hat u %= u - \Id$.
%	\end{itemize}
%Show:
%	\begin{enumerate}
%	
%	\item $\log$ is the inverse of $\exp$. (Since these functions are polynomials, this implies that $\Lie N$ and~$N$ are isomorphic as varieties.)
%	
%	\item For every nontrivial $u \in N$, there is a unique $C^\infty$ homomorphism $f \colon \real \to N$, such that $f(1) = u$. (In other words, $u$~is %contained in a unique one-parameter subgroup of~$N$.)
%
%	\item \label{UnipExpLog-conn}
%	Every connected, unipotent subgroup of $\SL(\ell,\real)$ is Zariski closed.
%	
%	\item \label{UnipExpLog-closed}
%	Every Zariski-closed, unipotent subgroup of $\SL(\ell,\real)$ is connected.
%	
%	\end{enumerate}
%	\hint{@@@}

\end{exercises}

\section{The Lie algebra of a Lie group}

\begin{defn}
 A map~$\rho$ from one Lie group to another is a \defit[homomorphism!of Lie groups]{homomorphism} if%
		\begin{itemize}
		\item it is a homomorphism of abstract groups (i.e., $\rho(ab) = \rho(a) \, \rho(b)$), 
		and
		\item it is continuous.
		\end{itemize}
	(Hence, an \defit[isomorphism!of Lie groups]{isomorphism} of Lie groups is a continuous isomorphism of abstract groups, whose inverse is also continuous.)
\end{defn}

Although the definition only requires homomorphisms to be continuous, it turns out that they are always infinitely differentiable:

\begin{prop} \label{HomosAreSmooth}
Suppose $H_1$ and~$H_2$ are closed subgroups of\/ $\GL(\ell_i,\real)$, for $i = 1,2$. Then
	\begin{enumerate}
	\item \label{HomosAreSmooth-subgrp}
	$H_i$ is a $C^\infty$ submanifold of\/ $\GL(\ell_i,\real)$,
	\item every\/ \textup(continuous\textup) homomorphism from~$H_1$ to~$H_2$ is~$C^\infty$
	and
	\item if $H_1 \subseteq H_2$, then the coset space $H_2/H_1$ is a $C^\infty$ manifold.
	\end{enumerate}
\end{prop}

\begin{rem}
In fact, the submanifolds and homomorphisms are real analytic, not just~$C^\infty$, but we will have no need for this stronger statement.
\end{rem}

\begin{rem} \label{OutGFinite}
If $H$ is any Lie group, then conjugation by any element~$h$ of~$H$ is an automorphism. That is, if we define a map $\varphi_h \colon H \to H$ by $\varphi_h(x) = h^{-1} x h$, then $\varphi_h$ is a continuous automorphism of~$H$. Any such automorphism is said to be ``\defit[inner!automorphism]{inner}\zz.'' The group of all inner automorphisms is isomorphic to $H/Z(H)$, where $Z(H)$ is the center of~$H$. For some groups, there are many other automorphisms. For example, every inner automorphism of an abelian group is trivial, but the automorphism group of~$\real^n$ is $\GL(n,\real)$, which is quite large. In contrast, it can be shown that the group of inner automorphisms of~$G$ has finite index in $\Aut(G)$ (since $G$ is semisimple).
\end{rem}

\begin{defns} \ 
\noprelistbreak
	\begin{enumerate}

	\item For $A,B \in \Mat_{\ell \times \ell}(\real)$, the \defit{commutator} (or \defit[Lie!bracket]{Lie bracket}) of $A$ and~$B$ is the matrix
		$ [A,B] = AB - BA $.
	\item A vector subspace $\Lie H$ of $\Mat_{\ell \times \ell}(\real)$ is a \defit[Lie!algebra]{Lie algebra} if it is closed under the Lie bracket. That is, for all $A,B \in \Lie H$, we have $[A,B] \in \Lie H$.
	\item A map~$\rho$ from one Lie algebra to another is a 
	\defit[homomorphism!of Lie algebras]{homomorphism}
	if
		\begin{itemize}
		\item it is a linear transformation,
		and
		\item it preserves brackets (that is, $[\rho(A), \rho(B)] = \rho \bigl( [A,B] \bigr)$).
		\end{itemize}
	\item Suppose $H$ is a closed subgroup of $\GL(\ell,\real)$. Then $H$ is a $C^\infty$ manifold, so it has a tangent space at every point; the tangent space at the identity element~$e$ is called the \defit[Lie!algebra]{Lie algebra} of~$H$. Note that, since $H$ is contained in the vector space $\Mat_{\ell \times \ell}(\real)$, its Lie algebra can be identified with a vector subspace of $\Mat_{\ell \times \ell}(\real)$.
	\end{enumerate}
\end{defns}

\begin{notation}
Lie algebras are usually denoted by lowercase German letters: the Lie algebras of $G$ and~$H$ are 
	\nindex{$\Lie G$, $\Lie H$ = Lie algebras of $G$ and~$H$}%
$\Lie G$ and~$\Lie H$, respectively.
\end{notation}

\begin{egs} \label{LieAlgEgs} \ 
\noprelistbreak
	\begin{enumerate}
	
%	\item Since $\GL(\ell,\real)$ is open in $\Mat_{\ell \times \ell}(\real)$, the Lie algebra of $\GL(\ell,\real)$ is all of $\Mat_{\ell \times \ell}(\real)$.
	
	\item \label{LieAlgEgs-SL}
	The Lie algebra $\LieSL(\ell,\real)$ of $\SL(\ell,\real)$ is the set of matrices whose trace is~$0$ \csee{LieAlgSLEx}.
	
	\item \label{LieAlgEgs-SO}
	The Lie algebra $\Lie{SO}(n)$ of $\SO(n)$ is the set of $n \times n$ skew-symmetric matrices of trace~$0$ \csee{LieAlgSOEx}.

%	Let $\Ortho(n) = \{\, A \in \Mat_{n \times n}(\real) \mid A^\transpose A = \Id \,\}$, where $A^\transpose$ is the transpose of~$A$, and $\Id$~is the $n \times n$ identity matrix.
	
	\end{enumerate}
\end{egs}

It is an important fact that the Lie algebra of~$H$ is indeed a Lie algebra:

\begin{prop}
If $H$ is a closed subgroup of\/ $\SL(\ell,\real)$, then the Lie algebra of~$H$ is closed under the Lie bracket.
\end{prop}

Here is a very useful reformulation:

\begin{cor}
Suppose $H_1$ and~$H_2$ are Lie groups, with Lie algebras\/ $\Lie H_1$ and\/~$\Lie H_2$. If $H_1$ is a subgroup of~$H_2$, then\/ $\Lie H_1$ is a Lie subalgebra of\/~$\Lie H_2$.
\end{cor}

Hence, for every closed subgroup of~$H$, there is a corresponding Lie subalgebra of~$\Lie H$. Unfortunately, the converse may not be true: although every Lie subalgebra corresponds to a subgroup, the subgroup might not be closed. 

\begin{eg}
The $2$-torus $\torus^2 = \real^2/\integer^2$ can be identified with the Lie group $\SO(2) \times \SO(2)$. For any line through the origin in~$\real^2$, there is a corresponding $1$-dimensional subgroup of~$\torus^2$. However, if the slope of the line is irrational, then the corresponding subgroup of~$\torus^2$ is dense, not closed.
\end{eg}

Therefore, in order to obtain a subgroup corresponding to each Lie subalgebra, we need to allow subgroups that are not closed:

\begin{defn}
Suppose $H_1$ and~$H_2$ are Lie groups, and $\rho \colon H_1 \to H_2$ is a homomorphism. Then $\rho(H_1)$ is a \defit[Lie!subgroup]{Lie subgroup} of~$H_2$.
\end{defn}

\begin{prop}
If $H$ is a Lie group with Lie algebra\/~$\Lie H$, then there is a one-to-one correspondence between the connected Lie subgroups of~$H$ and the Lie subalgebras of\/~$\Lie H$.
\end{prop}

%\begin{notation}
%\nindex{$H^\circ$ = identity component of~$H$}
%$H^\circ$ denotes the \defit{identity component} of the Lie group~$H$ (that is, the connected component of~$H$ that contains the identity element~$e$).
%\end{notation}

%\begin{prop}
%Suppose\/ $\Lie H_1$ and\/~$\Lie H_2$ are the Lie algebras of closed subgroups $H_1$ and~$H_2$ of\/ $\SL(\ell,\real)$. If\/ $\Lie H_1 = \Lie H_2$, then $H_1^\circ = H_2^\circ$.
%\end{prop}

\begin{defns}
Let $H$ be a Lie group in $\SL(\ell,\real)$.
	\begin{enumerate}
	
	\item If $h \colon \real \to H$ is any (continuous) homomorphism, we call $h$ a \defit{one-parameter subgroup} of~$H$, and we usually write $h^t$, instead of $h(t)$.
	
	\item We define $\exp \colon \Mat_{\ell \times \ell}(\real) \to \GL(\ell,\real)$ by
		$$ \exp X = \sum_{k=0}^\infty \frac{1}{k!} X^k .$$
	This is called the \defit{exponential map}.
	\end{enumerate}
\end{defns}

\begin{prop}
Let\/ $\Lie H$ be the Lie algebra of a Lie group $H \subseteq \SL(\ell,\real)$.
	\begin{enumerate}
	\item For any $X \in \Lie H$, the function $x^t = \exp(tX)$ is a one-parameter subgroup of~$H$.
	
	\item Conversely, every one-parameter subgroup of~$H$ is of this form, for some unique $X \in \Lie H$.
	\end{enumerate}
Furthermore, for $X \in \Lie \Mat_{\ell \times \ell}(\real)$, we have \ 
			$$ X \in \Lie H \ \Leftrightarrow \ \forall t \in \real, \ \exp(tX) \in H .$$
\end{prop}

\begin{defn}
Lie groups $H_1$ and~$H_2$ are \defit[locally!isomorphic]{locally isomorphic} if there is a connected Lie group~$H$ and homomorphisms $\rho_i \colon H \to H_i^\circ$, for $i = 1,2$, such that each $\rho_i$ is a covering map.
\end{defn}

\begin{prop} \label{LocIsoIff}
Two Lie groups are locally isomorphic if and only if their Lie algebras are isomorphic. 
\end{prop}

\begin{notation}[(adjoint representation)] \label{AdDefn}
Suppose $\Lie H$ is the Lie algebra of a closed subgroup $H$ of $\SL(\ell,\real)$.
For $h \in H$ and $x \in \Lie H \subseteq \Mat_{\ell \times \ell}(\real)$, we define%
	\nindex{$\Ad_H$ = adjoint representation of~$H$}
	$$ (\Ad_H h)(x) = h x h^{-1} \in \Lie H.$$
Then $\Ad_H \colon H \to \GL(\Lie H)$ is a (continuous) homomorphism. It is called the \defit[adjoint!representation]{adjoint representation} of~$H$.
\end{notation}

%\begin{rem} % \label{Gislinear}
%In most of this book, the assumption that $G$ is linear could
%be replaced with the weaker condition that the center $Z(G)$ is 
%finite. This is because  the kernel of the adjoint representation is
%$Z(G)$, so $G/Z(G)$ is isomorphic to a group of matrices;
%that is, $G/Z(G)$ is always linear. If $Z(G)$ is finite, this implies that
%$G$ is isogenous to $G/Z(G)$, so $G$ is isogenous to a
%linear group.
% \end{rem}

\begin{exercises}

\item \label{RnHomIsLinear}
Suppose $\rho \colon \real^m \to \real^n$ is a continuous map that preserves addition. (That is, we have $\rho(x + y) = \rho(x) + \rho(y)$.) Show (without using \cref{HomosAreSmooth}) that $\rho$ is a linear transformation (and is therefore~$C^\infty$). This is a very special case of \cref{HomosAreSmooth}.
\hint{By assumption, we have $\rho(k x) = k \rho(x)$ for all $k \in \integer$, so $\rho(tx) = t \rho(x)$ for all $t \in \rational$ (why?). Then continuity implies this is true for all $t \in \real$.}

\item \label{LieAlgSLEx}
Verify \fullcref{LieAlgEgs}{SL}.
\hint{$A \in \LieSL(\ell,\real)$ iff $\frac{d}{dt} \det(\Id + tA)\bigr|_{t = 0} = 0$,
and, letting $\lambda = 1/t$, we have 
$\det(\Id + tA) = t^\ell \det(\lambda I + A) = t^\ell \bigl( \lambda^\ell + (\trace A) \lambda^{\ell-1} + \cdots = 1 + (\trace A) t + \cdots$.}

\item \label{LieAlgSOEx}
Verify \fullcref{LieAlgEgs}{SO}.
\hint{A matrix $A$ of trace~$0$ is in $\Lie{SO}(n)$ iff $\frac{d}{dt} (\Id + tA)^\transpose (\Id + tA)\bigr|_{t = 0} = 0$. Calculate the derivative by using the Product Rule.}

\item In the notation of \cref{AdDefn}, show, for all $h \in H$, that $\Ad_H h$ is an automorphism of the Lie algebra~$\Lie H$. (In particular, it is an invertible linear transformation, so it is in $\GL(\Lie H)$.
\hint{The map $a \mapsto h a h^{-1}$ is a diffeomorphism of~$H$ that fixes~$e$, so its derivative is a linear transformation of the tangent space at~$e$.}

\end{exercises}

\section{How to show a group is semisimple}

A semisimple group $G = G_1 \times \cdots G_r$ will often have
connected, normal subgroups (such as the simple
factors~$G_i$). However, these normal subgroups cannot be
abelian \csee{SS->Nnotabel}. The converse is a major
theorem in the structure theory of Lie groups:

\begin{thm} \label{SS<>noAbelNorm}
 A connected Lie group~$H$ is semisimple if and only if it
has no nontrivial, connected, abelian, normal subgroups.
 \end{thm}

\begin{rem}
 A connected Lie group~$R$ is \defit[solvable Lie
group]{solvable} if every nontrivial quotient of~$R$ has a
nontrivial, connected, abelian, normal subgroup. (For
example, abelian groups are solvable.) It can be shown that
every connected Lie group~$H$ has a unique maximal connected,
closed, solvable, normal subgroup. This subgroup is called the
\defit[radical!of a Lie group]{radical} of~$H$, and is
denoted $\Rad H$. Our statement of
\cref{SS<>noAbelNorm} is equivalent to the more usual
statement that $H$~is semisimple if and only if $\Rad H$ is
trivial \csee{SS<>Nnotsolv}. 
 \end{rem}

%\begin{defn}
%To say that $G$ is \defit[Lie!group!linear]{linear} means $G$ is
%isomorphic to a group of matrices; more precisely, $G$ is
%isomorphic to a closed subgroup of $\SL(\ell,\real)$, for
%some~$\ell$. 
% \end{defn}

The following result makes it easy to see that the
classical groups, such as $\SL(n,\real)$, $\SO(m,n)$, and
$\SU(m,n)$, are semisimple (except a few abelian groups in
small dimensions).

\begin{defn} \label{irredrepDefn}
 A subgroup~$H$ of $\GL(\ell,\real)$ (or
$\GL(\ell,\complex)$) is \defit[irreducible!linear group]{irreducible} if there are no nontrivial,
proper, $H$-invariant subspaces of~$\real^\ell$
(or~$\complex^\ell$, respectively).
 \end{defn}

\begin{eg}
 $\SL(\ell,\real)$ is an irreducible subgroup of
$\SL(\ell,\complex)$ \csee{SLnRIrred}.
 \end{eg}

\begin{warn}
 In a different context, the adjective ``irreducible'' can
have a completely different meaning when it is applied to a
group. For example, saying that a lattice is irreducible
(as in \cref{irreducibleLattice}) has nothing to do
with \cref{irredrepDefn}.
 \end{warn}

\begin{rem}
 If $H$ is a subgroup of $\GL(\ell,\complex)$ that is
\emph{not} irreducible (that is, if $H$ is
\defit[reducible!linear group]{reducible}), then, after a
change of basis, we have
 $$ H \subseteq 
 \begin{pmatrix}
 \GL(k,\complex) & * \\
 0 & \GL(n-k,\complex)
 \end{pmatrix}
 ,$$
 for some~$k$ with $1 \le k \le n-1$.

Similarly for $\GL(\ell,\real)$.
 \end{rem}

\begin{cor} \label{irred->SS}
 If $H$ is a nonabelian, closed, connected, irreducible
subgroup of\/ $\SL(\ell,\complex)$, then $H$ is semisimple.
 \end{cor}

\begin{proof}
 Suppose $A$ is a connected, abelian, normal subgroup
of~$H$. For each function $w \colon A \to \complex^\times$,
let 
 $$ V_w = \{\, v \in \complex^\ell \mid \forall a \in A, \
a(v) = w(a) \, v \,\} .$$
 That is, a nonzero vector~$v$ belongs to~$V_w$ if 
 \begin{itemize}
 \item $v$ is an eigenvector for every element of~$A$, and 
 \item the corresponding eigenvalue for each element of~$a$
is the number that is specified by the function~$w$.
 \end{itemize}
 Of course, $0 \in V_w$ for every function~$w$; let $W =
\{\, w \mid V_w \neq 0\,\}$. (This is called the set of
\defit[weights of a representation]{weights} of~$A$
on~$\complex^\ell$.)

Each element of~$a$ has an eigenvector (because $\complex$
is algebraically closed), and the elements of~$A$ all
commute with each other, so there is a common eigenvector
for the elements of~$A$. Therefore, $W \neq \emptyset$. From
the usual argument that the eigenspaces of any linear
transformation are linearly independent, one can show that
the subspaces $\{\,V_w \mid w \in W \,\}$ are linearly
independent. Hence, $W$ is finite.

For $w \in W$ and $h \in H$, a straightforward calculation
shows that 
 $h V_w = V_{h(w)}$, where  $\bigl( h(w) \bigr) (a) =
w(h^{-1} a h)$. That is, $H$ permutes the subspaces
$\{V_w\}_{w \in W}$. Because $H$ is connected and $W$ is
finite, this implies $h V_w = V_w$ for each~$w$; that is,
$V_w$ is an $H$-invariant subspace of~$\complex^\ell$. Since
$H$ is irreducible, we conclude that $V_w = \complex^\ell$.

Now, for any $a \in A$, the conclusion of the preceding
paragraph implies that $a(v) = w(a)\, v$, for all $v \in
\complex^\ell$. Therefore, $a$ is a scalar matrix. 

Since $\det a = 1$, this scalar is an
$\ell^{\text{th}}$~root of unity. So $A$ is a subgroup of
the group of $\ell^{\text{th}}$~roots of unity, which is
finite. Since $A$ is connected, we conclude that $A =
\{e\}$, as desired.
 \end{proof}

Here is another useful characterization of semisimple
groups.

\begin{cor} \label{SelfAdj->SS}
 Let $H$ be a closed, connected subgroup of\/
$\SL(\ell,\complex)$. If
 \begin{itemize} 
 \item the center $Z(H)$ is finite, and
 \item $H^* = H$ \textup(where $*$ denotes the ``adjoint\zz,''
or conjugate-transpose\textup),
 \end{itemize}
 then $H$ is semisimple.
 \end{cor}

\begin{proof}
 Because $H^* = H$, it is not difficult to show that $H$ is
\defit{completely reducible}: there is a
direct sum decomposition 
 $ \complex^\ell = \bigoplus_{j=1}^r V_j $,
 such that the restriction $H|_{V_j}$ is irreducible, for each~$j$
\csee{H*=H->CompRed}.
 
  Let $A$ be a connected, normal subgroup of~$H$. 
The proof of \cref{irred->SS} (omitting the
final paragraph) shows that $A|_{V_j}$ consists of scalar
multiples of the identity, for each~$j$. Hence $A \subset
Z(H)$. Since $A$ is connected, but (by assumption) $Z(H)$
is finite, we conclude that $A$~is trivial.
 \end{proof}

\begin{rem} \label{SS->SelfAdj}
 There is a converse: if $G$ is semisimple (and connected),
then $G$ is conjugate to a subgroup~$H$, such that $H^* =
H$. However, this is more difficult to prove.
 \end{rem}

\begin{exercises}

\item \label{SS->Nnotabel}
 Prove ($\Rightarrow$) of \cref{SS<>noAbelNorm}.

\item \label{SS<>Nnotsolv}
 Show that a connected Lie group~$H$ is semisimple if and
only if $H$ has no nontrivial, connected, solvable, normal
subgroups.
 \hint{If $R$ is a solvable, normal subgroup of~$H$, then
$[R,R]$ is also normal in~$H$. Repeating this eventually yields an abelian, normal subgroup.}

\item \label{SLnRIrred}
 Show that no nontrivial, proper $\complex$-subspace of
$\complex^\ell$ is invariant under $\SL(\ell,\real)$.
 \hint{Suppose $v,w \in \real^\ell$, not both~$0$.
 If they are linearly independent, then there exists $g \in \SL(\ell,\real)$ with
$g(v + iw) = v - iw$. Otherwise, there exists nonzero $\lambda \in
\complex$ with $\lambda(v + iw) \in \real^\ell$.}

\item Give an example of a nonabelian, closed, connected,
irreducible subgroup~$H$ of $\SL(\ell,\real)$, such that $H$ is
not semisimple. 
 \hint{$\mathrm{U}(2)$ is an irreducible subgroup of $\SO(4)$.}

\item Suppose $H \subseteq \SL(\ell,\complex)$. Show that $H$
is completely reducible if and only if, for every
$H$-invariant subspace~$W$ of $\complex^\ell$, there is an
$H$-invariant subspace~$W'$ of~$\complex^\ell$, such that
$\complex^\ell = W \oplus W'$.
 \hint{($\Rightarrow$) If $W' = V_1 \oplus \cdots \oplus
V_s$, and $W' \cap W = \{0\}$, but $(W' \oplus V_j) \cap W
\neq \{0\}$ for every $j > s$, then $W' + W =
\complex^\ell$. ($\Leftarrow$) Let $W$ be maximal among the 
subspaces that are direct sums of irreducibles, and let $V$
be a minimal $H$-invariant subspace of~$W'$. Then $W \oplus
V$ contradicts the maximality of~$W$.}

\item \label{H*=H->CompRed}
 Suppose $H = H^* \subseteq \SL(\ell,\complex)$. 
 \begin{enumerate}
 \item Show that if $W$ is an $H$-invariant subspace of
$\complex^\ell$, then the orthogonal complement $W^\perp$
is also $H$-invariant.
 \item Show that $H$ is completely reducible.
 \end{enumerate}

\end{exercises}

\begin{notes}

See \cite{Howe-VeryBasic} for a very brief introduction to Lie groups, compatible with \fullcref{LieDefn}{LieGrp}. Similar elementary approaches are taken in the books \cite{Baker-MatGrps} and \cite{Hall-LieGrps}.

Almost all of the material in this \lcnamecref{SSGrpsChap} (other than \S\ref{ZariskiSect}) % @@@
 can be found in
Helgason's book \cite{HelgasonBook}. However, we do not
follow Helgason's notation for some of the classical groups
\csee{SOstarTerminology}.
%The lists of isogenies
%in \cref{isogTypes} and \cref{IsogInType} are
%largely copied from there \cite[\S10.6.4,
%pp.~519--520]{HelgasonBook}.

\Cref{GisAlgic} is proved in \cite[Thm.~8.3.2, p.~112]{Hochschild-AlgicGrps}.

%Helgason's book \cite[\S10.2, pp.~447--451]{HelgasonBook}
%proves that all of the (simply connected) classical simple groups except
%$\SO(m,n)$ are connected. For a geometric proof of
%\fullref{isogTypes}{D2-Sp4SO23}, see \cite[\S10.A.2,
%p.~521]{HelgasonBook}.

\Cref{RealJordanDecomp} can be found in
\cite[Lem.~IX.7.1, p.~430]{HelgasonBook}.

See \cite[Prop.~2 in \S11.2 of Chapter~8, p.~166]{BourbakiLie7-9} or \cite[Thm.~3.17, p.~100]{Jacobson-LieAlgs} for a proof of the Jacobson-Morosov Lemma \pref{JacobsonMorosov}.

See \cite[Thm.~2.7.5, p.~71]{VaradarajanBook} for a proof of \cref{LocIsoIff}.

\Cref{SS->SelfAdj} is due to Mostow
\cite{Mostow-SelfAdjoint}.

%The fact that measurable homomorphisms are
%continuous \see{MeasHomoIsCont} was proved by G.~Mackey.
% %% need reference
% A proof  can be found in \cite[Thm.~B.3,
%p.~198]{ZimmerBook}. 

%A proof that continuous homomorphisms are
%real analytic (hence~$C^\infty$) can be found in
%\cite[Prop.~1 of \S4.8, pp.~128--129]{Chevalley-LieGroups}.

\end{notes}

 %!TEX root = IntroArithGrps.tex

\mychapter{Assumed Background} \label{BackChap}

Since the target audience of this book includes mathematicians from a variety of backgrounds (and because very different theorems sometimes have names that are similar, or even identical), this chapter lists (without proof or discussion) specific notations, definitions, and theorems of graduate-level mathematics that are assumed in the main text. (Undergraduate-level concepts, such as the definitions of groups, metric spaces, and continuous functions, are generally not included.)
All of this material is standard, so proofs can be found in graduate textbooks (and on the internet).

\section{Groups and group actions}

\begin{notation} 
Let $H$ be a group, and let $K$ be a subgroup.
	\begin{enumerate}
	\item \nindex{$e$ = identity element of a group}%
	We usually use $e$ to denote the identity element.
	
	\item \nindex{$Z(H)$ = center of the group~$H$}%
	$Z(H) = \{\, z \in H \mid \text{$hz = zh$ for all $h \in H$} \,\}$ is the \defit[center!of a group]{center} of~$H$.
	
	\item \nindex{$\czer_H(K)$ = centralizer of~$K$ in~$H$}%
	$\czer_H(K) = \{\,h \in H \mid \text{$h k = kh$ for all $k \in K$}\,\}$ is the \defit{centralizer} of~$K$ in~$H$.
	
	\item \nindex{$\nzer_H(K)$ = normalizer of~$K$ in~$H$}%
	$\nzer_H(K) = \{\, h \in H \mid h K h^{-1} = K \,\}$ is the \defit{normalizer} of~$K$ in~$H$.
	\end{enumerate}
\end{notation}

\begin{defn}
An \defit[action of a group]{action} of a Lie group~$H$ on a topological space~$X$ is a
%homomorphism $\phi \colon G \to \Homeo(X)$, where
%$\Homeo(X)$ is the group of all permutations
%of~$X$. Equivalently, an \defit[action of~$G$]{action} is a
continuous function $\alpha \colon H \times X \to X$, such that
 \begin{itemize}
 \item $\alpha(e,x) = x$ for all $x \in X$, and
 \item $\alpha \bigl( g, \alpha(h,x) \bigr) = \alpha(gh,x)$
for $g,h \in H$ and $x \in X$.
\end{itemize}
% \item The action is \emph{proper}\index{proper
%action} if the map~$\alpha$ is proper; that is, if the
%inverse image of every compact set is compact.
% \end{enumerate}
 \end{defn}

%\begin{rem}
%If $\alpha$ is an action of~$G$ on~$X$, then we obtain a homomorphism $\varphi \colon G \to \Homeo(X)$, by $\varphi(g)(x) = \alpha(g,x)$.
%\end{rem}

\begin{defns}
 Let a (discrete) group $\Lambda$ act on a topological space~$M$. 
 \begin{enumerate}
 \item The action is \defit[free action]{free} if no
nonidentity element of~$\Lambda$ has a fixed point.
 \item It is \defit{properly discontinuous} if, for every
compact subset~$C$ of~$M$,
 $$ \text{the set $\{\, \lambda \in \Lambda \mid C \cap (\lambda C) \neq \emptyset
\,\}$ is finite} .$$
 \item For any $p \in M$, we define 
 	\nindex{$\Stab_\Lambda(p)$ = stabilizer of~$p$}
	$\Stab_\Lambda(p)
 = \{\, \lambda \in \Lambda \mid \lambda p = p \,\}$. This is a
subgroup of~$\Lambda$ called the \defit{stabilizer} of~$p$ in~$\Lambda$.

\item $M$ is \defit{connected} if it is \textbf{not} the union of two nonempty, disjoint, proper, open subsets.
\item $M$ is \defit[locally!connected]{locally connected} if every neighborhood of every $p \in M$ contains a connected neighborhood of~$p$.
 \end{enumerate}
 \end{defns}

\begin{prop} \label{PropdiscFree->cover}
 If $\Lambda$ acts freely and properly discontinuously on a
topological space~$M$, then the natural map $\pi \colon M \to \Lambda
\backslash M$ is a {\normalfont\defit{covering map}}.

Under the simplifying assumption that $M$ is locally connected, this means that every $p \in \Lambda \backslash M$ has a connected neighborhood~$U$, such that the restriction of~$\pi$ to each connected component of $\pi^{-1}(U)$ is a homeomorphism onto~$U$.
 \end{prop}

%\begin{prop} \label{Isom->proper}
% If $M$ is any locally compact metric space, then $\Isom(M)$ is a
%locally compact topological group, under the compact-open topology
%\textup(that is, under the topology of uniform convergence
%on compact sets\textup). The action of $\Isom(M)$ on~$M$
%is proper; that is, for every compact subset~$C$ of~$M$, the
%set
% $$ \{\, \phi \in \Isom(M) \mid \phi(C) \cap C \neq \emptyset \,\}$$
% is compact.
%
%If $M$ is a smooth manifold, then the topological group $\Isom(M)$
%can be given the structure of a Lie group, so that the action of
%$\Isom(M)$ on~$M$ is smooth.
% \end{prop}

\section{Galois theory and field extensions}

\begin{thm}[{(\thmindex{Fundamental Theorem of Algebra}{Fundamental Theorem of Algebra})}]
 The field~$\complex$ of complex numbers is algebraically
closed; that is, every nonconstant polynomial $f(x) \in
\complex[x]$ has a root in~$\complex$.
 \end{thm}

%\begin{proof}
% This can be proved algebraically, by combining Galois
%Theory with the elementary fact that every real polynomial of
%odd degree has a real zero \csee{AlgicPfFundThmAlg}, but we
%use a bit of complex analysis.
%
%Suppose $f(x)$ has no root. Then $1/f$ is holomorphic
%on~$\complex$. Furthermore, because $f(z) \to \infty$ as $z
%\to \infty$, it is easy to see that $1/f$ is bounded
%on~$\complex$. Hence, Liouville's Theorem asserts that $1/f$
%is constant. This contradicts the fact that $f$~is not
%constant.
% \end{proof}

%\begin{prop}
% Suppose $F$ is a field, and $\alpha$ is a root of some
%irreducible polynomial $f(x) \in F[x]$. Then the extension
%field $F[\alpha]$ is isomorphic to $F[x]/I$, where $I = f(x)
%\, F[x]$ is the principal ideal of $F[x]$ generated by $f(x)$.
% \end{prop}
%
%\begin{proof}
% The map $\phi \colon F[x] \to F[\alpha]$ defined by
%$\phi \bigl( g(x) \bigr) = g(\alpha)$ is a surjective ring
%homomorphism whose kernel is~$I$.
% \end{proof}
%
%\begin{cor}
%  Suppose $F$ is a field. If $\alpha$ and~$\beta$ are two
%roots of an irreducible polynomial $f(x) \in F[x]$, then
%there is an isomorphism $\sigma \colon F[\alpha] \to
%F[\beta]$ with $\sigma(\alpha) = \beta$.
% \end{cor}

\begin{prop}
 Let $F$ be a subfield of~$\complex$, and let $\sigma \colon
F \to \complex$ be any embedding. Then $\sigma$~extends to an
automorphism~$\widehat\sigma$ of~$\complex$.
 \end{prop}

\begin{notation}
 If $F$ is a subfield of a field~$L$, then 
 	\nindex{$|L:F|$ = degree of field extension}$|L:F|$
denotes
$\dim_F L$, the dimension of~$L$ as a vector space over~$F$. This is called the \defit[degree! of a field extension]{degree} of~$L$ over~$F$.
 \end{notation}

\begin{prop}
  If $F$ and~$L$ are subfields of~$\complex$, such that $F
\subseteq L$, then $|L:F|$ is equal to the number of
embeddings~$\sigma$ of~$L$ in~$\complex$, such that
$\sigma|_F = \Id$.
 \end{prop}

\begin{defn}
 An extension~$L$ of a field~$F$ (of characteristic zero) is \defit[Galois!extension]{Galois} if, for every irreducible polynomial $f(x)
\in F[x]$, such that $f(x)$ has a root in~$L$, there exist
$\alpha_1,\ldots,\alpha_n \in L$, such that 
 $$ f(x) = (x - \alpha_1) \cdots (x - \alpha_n) .$$
 That is, if an irreducible polynomial in $F[x]$ has a root
in~$L$, then all of its roots are in~$L$.
 \end{defn}

\begin{defn}
 Let $L$ be a Galois extension of a field~$F$. Then
 	\nindex{$\Gal(L/F)$ = Galois group}
 $$ \Gal(L/F) = \{\, \sigma \in \Aut(L) \mid \sigma|_F = \Id
\,\} .$$
 This is the \defit[Galois!group]{Galois group} of~$L$ over~$F$.
 \end{defn}

\begin{prop}
  If $L$ is a Galois extension of a field~$F$ of
characteristic~$0$, then $|{\Gal(L/F)}| = |L:F|$.
 \end{prop}

\begin{cor}
 If $L$ is a Galois extension of a field~$F$ of
characteristic\/~$0$, then there is a one-to-one correspondence
between
 \begin{itemize}
 \item the subfields~$K$ of~$L$, such that $F \subseteq K$, and
 \item the subgroups~$H$ of\/ $\Gal(L/F)$.
 \end{itemize}
 Specifically, the subgroup of\/ $\Gal(L/F)$ corresponding to the subfield~$K$ is\/ $\Gal(L/K)$.
 \end{cor}

%\begin{proof}
% Given $K$, with $F \subseteq K \subseteq L$, let $H_K =
%\Gal(L/K)$.
% Conversely, given a subgroup~$H$ of $\Gal(L/F)$, let
% $$ K_H = \{\, x \in L \mid \sigma(x) = x, \ \forall \sigma
%\in H \,\} $$
% be the fixed field of~$H$.
% \end{proof}

%\begin{exercises}
%
%\item \label{AlgicPfFundThmAlg}
% Let $F$ be a field of characteristic zero, such that
% \begin{enumerate} 
% \item if $f \in F[x]$ has odd degree, then $f$ has a root
%in~$F$,
% \item if $a \in F$, then either $a$ or~$-a$ has a square
%root in~$F$, and
% \item $-1$ does not have a square root in~$F$.
% \end{enumerate}
% Show that $F[i]$ is algebraically closed (where $i =
%\sqrt{-1}$).
% \hint{Let $L$ be a finite, Galois extension of $F[i]$. If
%$P$ is a Sylow $2$-subgroup of $\Gal(L/F)$, then the fixed
%field of~$P$ has odd degree over~$F$, so this fixed field must
%be trivial. Therefore, $|L:F|$ is a power of~$2$. Hence, $L$ can be
%obtained by a series of quadratic extensions. Since every
%element of~$F$ has a square root in $F[i]$, the half-angle
%formulas show that every element of $F[i]$ has a square root
%in $F[i]$. Therefore $L \subseteq F[i]$.}
%
%\end{exercises}

\section{Algebraic numbers and transcendental numbers}

\begin{defns} \ 
 \noprelistbreak
 \begin{enumerate}
 \item A complex number~$z$ is \defit[algebraic!number]{algebraic} if there is a nonzero polynomial $f(x) \in
\integer[x]$, such that $f(z) = 0$. 
 \item A complex number is \defit[transcendental number]{transcendental} if it is not algebraic.
 \item A (nonzero) polynomial is \defit[polynomial!monic]{monic} if its leading coefficient is~$1$; that is,
we may write $f(x) = \sum_{k=0}^n a_k x^k$ with $a_n = 1$.
 \item A complex number~$z$ is an \defit[algebraic!integer]{algebraic integer} if
there is a \emph{monic} polynomial $f(x) \in \integer[x]$,
such that $f(z) = 0$. 
 \end{enumerate}
 \end{defns}

%\begin{prop}[{}{(\thmindex{Z@$\integer$ is integrally closed}$\integer$ is
%integrally closed}}] \label{Zintegclosed}
% A rational number $t \in \rational$ is an algebraic integer if and only if
%$t \in \integer$.
% \end{prop}

%\begin{proof}
% ($\Leftarrow$) $t$~is a root of the monic polynomial $x - t$.
%
%($\Rightarrow$) Suppose $f(t) = 0$, where $f(x) = x^n +
%\sum_{k=0}^{n-1} a_k x^k$ with each $a_k \in \integer$.
%Writing $t = p/q$ (in lowest terms) with $p,q \in \integer$,
%we have
% $$ 0 = q^n \cdot 0 = q^n f(t)
% = q^n \left( \frac{p^n}{q^n} + \sum_{k=0}^{n-1} a_k
%\frac{p^k}{q^k} \right)
% = p^n + \sum_{k=0}^{n-1} a_k p^k q^{n-k}
% \equiv p^n \pmod{q} .$$
% Since $p^n$~is relatively prime to~$q$ (recall that $t = p/q$
%is in lowest terms), we conclude that $q = 1$, so $t = p/1
%\in \integer$.
% \end{proof}

\begin{prop}
 If $\alpha$ is an algebraic number, then there is some
nonzero $m \in \integer$, such that $m \alpha$ is an
algebraic integer.
 \end{prop}

%\begin{proof}
% Suppose $g(\alpha) = 0$, where $g(x) = \sum_{k=0}^n b_k
%x^k$, with each $b_k \in \integer$, and $b_n \neq 0$. Let 
% \begin{itemize}
% \item $m = a_n$,
% \item $a_k = m^{n-k-1} b_k$, and
% \item $f(x) = \sum_{k=0}^n a_k x^k$.
% \end{itemize}
% Then $f(x)$ is a monic, integral polynomial, and 
% $$ f(m \alpha) = \sum_{k=0}^n (m^{n-k-1} b_k) (m\alpha)^k
% = m^{n-1}  \sum_{k=0}^n b_k \alpha^k
% = m^{n-1} g(\alpha)
% = m^{n-1} \cdot 0
% = 0 .$$
% \end{proof}

%\begin{lem}
% For $t \in \complex$, the following are equivalent:
% \begin{enumerate}
% \item $t$ is an algebraic integer.
% \item $\integer[t]$ is a finitely-generated
%$\integer$-module.
%  \item $\integer[t]$ is a Noetherian
%$\integer$-module.
% \end{enumerate}
% \end{lem}

\begin{prop}
 The set of algebraic integers is a subring of\/~$\complex$.
 \end{prop}

\begin{prop} \label{Gal(cyclotomic)}
Fix some $n \in \natural^+$.
 Let 
 \begin{itemize}
 \item $\omega$ be a primitive $n^\text{th}$~root of unity,
and
 \item $\integer_n^\times$ be the multiplicative group
of units modulo~$n$.
 \end{itemize}
Then there is an isomorphism 
 $$f \colon \integer_n^\times \to \Gal \bigl(
\rational[\omega]/\rational \bigr) 
\colon k \mapsto f_k,$$
 such that $f_k(\omega) = \omega^k$, for all $k \in
\integer_n^\times$.
 \end{prop}

\section{Polynomial rings}

\begin{defn}
 A commutative ring~$R$ is \defit[Noetherian ring]{Noetherian}
if the following equivalent conditions hold:
 \begin{enumerate}
 \item Every ideal of~$R$ is finitely generated.
% \item Every nonempty collection of ideals of~$R$ has a
%maximal element.
 \item If $I_1 \subseteq I_2 \subseteq \cdots$ is any increasing
chain of ideals of~$R$, then there is some~$m$, such that $I_m = I_{m+1} = I_{m+2} = \cdots$.
 \end{enumerate}
 \end{defn}

%\begin{thm} \label{R[x]Noetherian}
% If $R$ is Noetherian, then the polynomial ring $R[x]$ is
%Noetherian.
% \end{thm}
%
%\begin{proof}
% Suppose $J$ is an idea of $R[x]$. (We wish to show that
%$J$~is finitely generated.) For $d \in \natural$, let
% $$ J_d = \{0\} \cup \{\, \operatorname{lead}(f) \mid f \in
%J, \ \deg f = d \,\} .$$
% Then $J_d$ is an ideal of~$R$, and we have
% $J_1 \subseteq J_2 \subseteq \cdots$, so there is some $d_0$,
%such that $J_d = J_{d_0}$, for all $d \ge d_0$.
%
%For each $d$, let $F_d$ be a finite set of polynomials of
%degree~$d$, such that $\{\, \operatorname{lead}(f) \mid f \in
%F_d \,\}$ generates~$J_d$.
%
%Then $F_0 \cup F_1 \cup \ldots \cup F_{d_0}$ generates~$J$.
%(For any $f \in J$, there exists $f' \in \langle F_0 \cup F_1
%\cup \ldots \cup F_{d_0} \rangle$, such that $\deg f' = \deg
%f$, and $\operatorname{lead}(f') = \operatorname{lead}(f)$.
%Then $\deg(f - f') < \deg f$, so we may assume, by induction,
%that $f - f' \in \langle F_0 \cup F_1 \cup \ldots \cup
%F_{d_0} \rangle$.
% \end{proof}

\begin{prop}[(\thmindex{Hilbert Basis}Hilbert Basis Theorem)]
 For any field~$F$, the polynomial ring $F[x_1,\ldots,x_s]$
\textup(in any number of variables\textup) is
Noetherian.
 \end{prop}

%\begin{proof}
% Note that $F$ has only one proper ideal, namely~$\{0\}$, so
%it is obviously Noetherian. Now use \cref{R[x]Noetherian}
%to induct on~$s$.
% \end{proof}

%There are many equivalent formulations of the following
%important theorem.

\begin{thm} \label{Nullstellensatz-fg}
 Let $F$ be a subfield of a field~$L$. If $L$ is finitely
generated as an $F$-algebra \textup(that is, if there
exist $c_1,\ldots,c_r \in L$, such that $L =
F[c_1,\ldots,c_r]$\textup), then $L$~is algebraic over~$F$.
 \end{thm}

\begin{prop}[(\thmindex{Nullstellensatz}Nullstellensatz)] \label{Nullstellensatz}
 Let 
 \begin{itemize}
 \item $F$ be an algebraically closed field, 
 \item $F[x_1,\ldots,x_r]$ be a polynomial ring over~$F$, and
 \item $I$ be any proper ideal of $F[x_1,\ldots,x_r]$.
 \end{itemize}
 Then there exist $a_1,\ldots,a_r \in F$, such that
 $f(a_1,\ldots,a_r) = 0$ for all $f(x_1,\ldots,x_r) \in I$.
 \end{prop}

%\begin{proof}
% Let $\mathfrak{m}$ be a maximal ideal that contains~$I$, and
%choose $a_1,\ldots,a_r \in F$ as in
%\fullref{Nullstellensatz-maxideal}{roots}. Then
%$f(a_1,\ldots,a_r) = 0$ for all $f(x_1,\ldots,x_r) \in
%\mathfrak{m}$, so, since $I \subseteq \mathfrak{m}$, the desired
%conclusion follows.
% \end{proof}

\begin{cor} \label{Nullstellensatz-ringhomo}
 If $B$ is any finitely generated subring of\/~$\complex$, then
there is a nontrivial homomorphism from~$B$ to the
algebraic closure\/~$\overline{\rational}$ of\/~$\rational$.
 \end{cor}

\begin{lem}[(\thmindex{Eisenstein Criterion}{Eisenstein Criterion})] \label{Eisenstein}
 Let $f(x) \in \integer[x]$. If there
is a prime number~$p$, and some $a \in \integer_p
\smallsetminus \{0\}$, such that 
 \begin{itemize}
 \item $f(x) \equiv a x^n \pmod{p}$, where $n = \deg f(x)$,
and
 \item $f(0) \not\equiv 0 \pmod{p^2}$,
 \end{itemize}
 then $f(x)$ is irreducible over~$\rational$.
 \end{lem}

%\begin{proof}
% Suppose $f(x)$ is reducible over~$\rational$. (This will
%lead to a contradiction.) Then $f(x)$ is also reducible
%over~$\integer$ \see{irredZ->irredQ}, so we may write $f(x)
%= g_1(x) g_2(x)$, with $g_j(x) \in \integer[x]$ and $\deg
%g_j(x) \ge 1$.  Then 
% $$ g_1(x) g_2(x) = f(x) \equiv a x^n \pmod{p} .$$
% From the unique factorization of polynomials in
%$\integer_p[x]$ (recall that $\integer_p[x]$ is a Euclidean
%domain, because $\integer_p$~is a field), we conclude that
%there exist $b_1,b_2 \in \integer_p \smallsetminus \{0\}$
%and $m_1 , m_2 \in \natural$, such that $g_j(x) \equiv b_j
%x^{m_j} \pmod{p}$. Since
% $$ m_1 + m_2 = n = \deg f(x) = \deg \bigl( g_1(x) g_2(x)
%\bigr) = \deg g_1(x) + \deg g_2(x) ,$$
% and $m_j \le \deg g_j(x)$, we conclude that $m_j = \deg
%g_j(x) > 1$. Therefore $g_j(0) \equiv 0 \pmod{p}$, so $f(0)
%= g_1(0) g_2(0) \equiv 0 \pmod{p^2}$. This is a
%contradiction.
% \end{proof}

\section{General topology} \label{TopologySect}

\begin{defns}
Let $X$ be a topological space.
	\begin{enumerate}
	\item A subset~$C$ of~$X$ is \defit{precompact} (or \defit{relatively compact}) if the closure of~$C$ is compact.
	\item $X$ is \defit[locally!compact]{locally compact} if every point of~$X$ is contained in a precompact, open subset.
	\item $X$ is \defit[separable topological space]{separable} if it has a countable, dense subset.
%	\item $X$ is \defit[metrizable topological space]{metrizable} if there is a metric on~$X$, such that the open sets in~$X$ are precisely the sets that are unions of open balls with respect to the metric.
	\item If $I$ is an index set (of any cardinality), and $X_i$ is a topological space, for each $i \in I$, then the Cartesian product $\bigtimes_{i \in I} X_i$ has a natural ``\defit[product!topology]{product topology}\zz,'' in which a set is open if and only if it is a union (possibly infinite) of sets of the form $\bigtimes_{i \in I} U_i$, where each $U_i$ is an open subset of~$X_i$, and we have $U_i = X_i$ for all but finitely many~$i$.
	
	\end{enumerate}
\end{defns}

\begin{thm}[(\thmindex{Tychonoff's}Tychonoff's Theorem)] \label{TychonoffThm}
If $X_i$ is a compact topological space, for each $i \in I$, then the Cartesian product $\bigtimes_{i \in I} X_i$ is also compact\/ \textup(with respect to the product topology\/\textup).
\end{thm}

\begin{prop}[(\thmindex{Zorn's Lemma}Zorn's Lemma)] \label{ZornsLemma}
Suppose $\le$ is a binary relation on a set~$\mathcal{P}$, such that:
	\begin{itemize}
	\item If $a \le b$ and $b \le c$, then $a \le c$.
	\item If $a \le b$ and $b \le a$, then $a = b$.
	\item $a \le a$ for all~$a$.
	\item If $\mathcal{C} \subseteq \mathcal{P}$, such that, for all $c_1,c_2 \in \mathcal{C}$, either $c_1 \le c_2$ or $c_2 \le c_1$, then there exists $b \in \mathcal{P}$, such that $c \le b$, for all $c \in \mathcal{C}$.
	\end{itemize}
Then there exists $a \in \mathcal{P}$, such that $a \not\le b$, for all $b \in \mathcal{P}$.
\end{prop}

\section{Measure theory} \label{MeasThySect}

\begin{assump}
Throughout this section, % @@@ \lcnamecref{MeasThySect}, 
$X$ and~$Y$ are complete, separable metric spaces. (Recall that \defit[complete metric space]{complete} means all Cauchy sequences converge.)
\end{assump}

\begin{defns} \ 
%Suppose $X$ and~$Y$ are topological spaces.
\noprelistbreak
	\begin{enumerate}
	\item The \defit[Borel!sigma-algebra@$\sigma$-algebra]{Borel $\sigma$-algebra} 
	\nindex{$\Borel(X)$ = $\{$Borel subsets of~$X$$\}$}
	$\Borel(X)$ of~$X$ is the smallest collection of subsets of~$X$ that:
		\begin{itemize}
		\item contains every open set,
		\item is closed under countable unions (that is, if $A_1,A_2,\ldots \in \Borel$, then $\bigcup_{i=1}^\infty A_i \in \Borel$),
		and
		\item is closed under complements (that is, if $A \in \Borel$, then $X \smallsetminus A \in \Borel$).
		\end{itemize}
	\item Each element of $\Borel(X)$ is called a \index{Borel!set}\defit[Borel!set]{Borel set}.
	\item A function $f \colon X \to Y$ is \defit[Borel!function]{Borel measurable} if $f^{-1}(A)$ is a Borel set in~$X$, for every Borel set~$A$ in~$Y$.
	\item A function $\mu \colon \Borel(X) \to [0,\infty]$ is called a \defit{measure} if it is \defit{countably additive}. This means that if $A_1,A_2,\ldots$ are pairwise disjoint, then 
		$$ \mu \left( \bigcup_{i=1}^\infty A_i \right) = \sum_{i = 1}^\infty \mu(A_i) .$$
	
	\item A measure~$\mu$ on~$X$ is \defit[measure!Radon]{Radon} if $\mu(C) < \infty$, for every compact subset~$C$ of~$X$.

	\item A measure~$\mu$ on~$X$ is \defit[sigma-finite measure@$\sigma$-finite measure]{$\sigma$-finite} if $X$ is the union of countably many sets of finite measure. This means $X = \bigcup_{i = 1}^\infty A_i$, with $\mu(A_i) < \infty$ for each~$i$.
	
	\end{enumerate}
\end{defns}

\begin{prop}
If $\mu$ is a measure on~$X$, and $f$~is a measurable function on~$X$, such that $f \ge 0$, then the integral $\int_X f \, d\mu$ is a well-defined element of $[0,\infty]$, such that:
	\begin{enumerate}
	\item $\int_X \chi_A \, d\mu = \mu(A)$ if $\chi_A$ is the characteristic function of~$A$.
	\item $\int_X (a_1 f_1 + a_2 f_2) \, d\mu = a_1 \int_X f_1 \, d\mu + a_2 \int_X f_2 \, d\mu$ for $a_1,a_2 \in [0,\infty)$.
	\item if $\{f_n\}$ is a sequence of measurable functions on~$X$, such that we have\/ $0 \le f_1 \le f_2 \le \cdots$, then
		$$ \int_X \ \lim_{n\to \infty} f_n \  d\mu = \lim_{n\to \infty} \int_X f_n \, d\mu . $$
	\end{enumerate}
\end{prop}

\begin{cor}[(Fatou's Lemma)] \label{FatousLemma}
If $\{f_n\}_{n=1}^\infty$ is a sequence of measurable functions on~$X$, with $f_n \ge 0$ for all~$n$, and $\mu$~is a measure on~$X$, then
	$$ \int_X \liminf_{n \to \infty} f_n \, d\mu \le \liminf_{n \to \infty} \int_X f_n \, d\mu .$$
\end{cor}

\begin{prop}
If $X$~is locally compact and separable, then every Radon measure~$\mu$ on~$X$ is \defit[regular!measure, inner]{inner regular}. This means 
	$$ \text{$\mu(E) = \sup\{\, \mu(C) \mid \text{$C$ is a compact subset of~$E$} \,\}$,
	\ for every Borel set~$E$}. $$
\end{prop}

\begin{prop}[(\thmindex{Lusin's}Lusin's Theorem)] \label{LusinsThm}
Assume $\mu$ is a Radon measure on~$X$, and $X$~is locally compact. Then, for every measurable function $f \colon X \to \real$, and every $\epsilon > 0$, there is a continuous function $g \colon X \to \real$, such that 
	$$ \mu \bigl( \{\, x \in X \mid f(x) \neq g(x) \,\} < \epsilon .$$
\end{prop}

\begin{defn} \label{PushForwardDefn}
If $\mu$ is a measure on~$X$, and $f \colon X \to Y$ is measurable, then the \defit[push-forward of measure~$\mu$]{push-forward} of~$\mu$ is the measure 
	\nindex{$f_*\mu$ = push-forward of measure~$\mu$}%
	$f_*\mu$ on~$Y$ that is defined by
		$$ (f_* \mu)(A) = \mu \bigl( f^{-1}(A) \bigr) 
		\quad \text{for $A \subseteq Y$} . $$
\end{defn}

\begin{prop}[{(\thmindex{Fubini's}Fubini's Theorem)}]
Assume
	\begin{itemize}
	\item $X_1$ and $X_2$ are complete, separable metric spaces, 
	and
	\item $\mu_i$ is a $\sigma$-finite measure on~$X_i$, for $i = 1,2$.
	\end{itemize}
Then there is a measure $\nu = \mu_1 \times \mu_2$ on $X_1 \times X_2$, such that:
	\begin{enumerate}
	\item $\nu(E_1 \times E_2) = \nu(E_1) \cdot \nu(E_2)$ when $E_i$ is a Borel subset of~$X_i$ for $i = 1,2$,
	and
	\item $\int_{X_1 \times X_2} f \, d\nu = \int_{X_1} \int_{X_2} f(x_1,x_2) \, d\mu_2(x_2) \, d \mu_1(x_1)$ when the function
	$f \colon X_1 \times X_2 \to [0,\infty]$ is Borel measurable.
	\end{enumerate}
\textup(In particular, $\int_{X_2} f(x_1,x_2) \, d\mu_2(x_2)$ is a measurable function of~$x_1$.\textup)
\end{prop}

\begin{defns} \ 
\noprelistbreak
\begin{enumerate} 
	\item The \defit{support} of a function $f \colon X \to \complex$ is defined to be the closure of $\{\, x \in X \mid f(x) \neq 0 \,\}$.
	\item \nindex{$C_c(X)$ = $\{$continuous functions with compact support$\}$}%
	$C_c(X) = \{\, \text{continuous functions $f \colon X \to \complex$ with compact support} \,\}$.

	\item $\lambda \colon C_c(X) \to \complex$ is a \defit[linear functional! positive]{positive linear functional} on $C_c(X)$ if:
		\begin{itemize}
		\item it is linear (that is, $\lambda( a_1 f_1 + a_2 f_2) = a_1 \lambda(f_1) + a_2 \lambda(f_2)$ for $a_1,a_2 \in \complex$ and $f_1,f_2 \in C(X)$),
		and
		\item it is positive (that is, if $f(x) \ge 0$ for all~$x$, then $\lambda(f) \ge 0$).
		\end{itemize}
\end{enumerate}
\end{defns}

\begin{thm}[(\thmindex{Riesz Representation}Riesz Representation Theorem)] \label{RieszRepThm}
Assume $X$~is locally compact and separable.
If\/ $\lambda$ is any positive linear functional on $C_c(X)$, then there is a Radon measure~$\mu$ on~$X$, such that
	$$ \lambda(f) = \int_X f \, d\mu 
	\quad \text{for all $f \in C_c(X)$} .$$
%Furthermore, $\mu(C) < \infty$ for every compact subset~$C$ of~$X$.
\end{thm}

\begin{defns}
Assume  $\mu$ is a measure on~$X$, and the function $\varphi \colon X \to \complex$ is measurable.

	\begin{enumerate}
	\item For $1 \le p < \infty$, the \defit[Lp-@$\LL{p}$-!norm]{$\pmbLL{p}$-norm} of~$\varphi$ is 
	$$ \|\varphi\|_p = \left( \int_{X} |\varphi(x)|^p \, d\mu(x) \right) ^{1/p} .$$
	
	\item An assertion $P(x)$ is said to be true for \defit[almost!all]{almost all} $x \in X$ (or to be true \defit[almost!everywhere]{almost everywhere}, which is usually abbreviated to \defit{a.e.}),  if $\mu \bigl( \{\, x \mid \text{$P(x)$ is false}\,\} \bigr) = 0$.
	
	\item In particular, two functions $\varphi_1$ and $\varphi_2$ are equal (a.e.) if 
		$$\mu\bigl( \{\, x \mid \varphi_1(x) \neq \varphi_2(x) \,\} \bigr) = 0 .$$ 
	This defines an equivalence relation on the set of (measurable) functions on~$X$.

	\item The \defit[Lq-norm@$\LL{\infty}$-norm]{$\pmbLL{\infty}$-norm} (or \defit{essential supremum}) of~$\varphi$ is 
	$$ \|\varphi\|_\infty = \min \bigl\{\, a \in (-\infty, \infty] \mid \text{$\varphi(x) \le a$ for a.e.~$x$} \,\bigr\} .$$

	\item \nindex{$ \LL{p}(X,\mu)$ = $\{$$\LL{p}$-functions on~$X$$\}$}
	$ \LL{p}(X,\mu) = \{\, \varphi \colon X \to \complex \mid \| \varphi \|_p < \infty \,\} $, for $1 \le p \le \infty$. An element of $\LL{p}(X,\mu)$ is called an \defit[Lp-@$\LL{p}$-!function]{$\pmbLL{p}$-function} on~$X$. Actually, two functions in $\LL{p}(X,\mu)$ are identified if they are equal almost everywhere, so, technically, $\LL{p}(X,\mu)$ should be defined to be a set of equivalence classes, instead of a set of functions.
	
	\end{enumerate}
%Note that $\|\varphi\|_p = 0$ if and only if $\varphi = 0$~a.e.
\end{defns}

\begin{defn}
Two measures $\mu$ and~$\nu$ on~$X$ are in the same \defit[measure!class]{measure class}
		%(or $\nu$ is \defit[equivalent!measure]{equivalent} to~$\mu$) 
if they have exactly the same sets of measure~$0$:
		$$ \mu(A) = 0 \iff \nu(A) = 0 .$$
	(This defines an equivalence relation.) 
%	Note that if $\nu = f \mu$, for some real-valued, measurable function~$f$, such that $f(x) \neq 0$ for a.e.\ $x \in X$, then $\mu$ and~$\nu$ are in the same measure class \csee{fmuClassOfMuEx}. 
%(The \term{Radon-Nikodym Theorem} implies that the converse is true,
%	% if the measures are $\sigma$-finite, 
%but we do not need this fact.)
\end{defn}

\begin{thm}[(\thmindex{Radon-Nikodym}Radon-Nikodym Theorem)] \label{RadonNikodym}
Two $\sigma$-finite measures $\mu$ and~$\nu$ on~$X$ are in the same class if and only if there is a measurable function $D \colon X \to \real^+$, such that $\mu = D \nu$. That is, for every measurable subset~$A$ of~$X$, we have
	$ \mu(A) = \int_A D \, d\nu $.
\end{thm}

The function~$D$ is called the \defit{Radon-Nikodym derivative} $d\mu/d\nu$.

\section{Functional analysis}

\begin{defns}
Let  $\field$ be either $\real$ or~$\complex$, and let $V$ be a vector space over~$\field$.
	\begin{enumerate}
	
	\item A \defit{topological vector space} is a vector space~$V$, with a topology, such that the operations of scalar multiplication and vector addition are continuous (that is, the natural maps $\field \times V \to V$ and $V \times V \to V$ are continuous).
	
	\item A subset~$C$ of~$V$ is \defit[convex set]{convex} if, for all $v,w \in C$ and $0 \le t \le 1$, we have $t v + (1-t) w \in C$.
	
	\item A topological vector space~$V$ is \defit[locally!convex]{locally convex} if every neighborhood of~$0$ contains a convex neighborhood of~$0$.
	
	\item A locally convex topological vector space~$V$ is \defit[Fréchet space]{Fréchet} if its topology can be given by a metric that is \defit[complete metric space]{complete} (that is, such that every Cauchy sequence converges to a limit point).
	
	\item A \defit{norm} on~$V$ is a function $\| \ \| \colon V \to [0,\infty)$, such that:
	\noprelistbreak
		\begin{enumerate}
		\item $\| v + w \| \le \| v \| + \| w \|$ for all $v , w \in V$,
		\item $\| av \| = |a| \, \|v\|$ for $a \in \field$ and $v \in V$,
		and
		\item $\|v\| = 0$ if and only if $v = 0$.
		\end{enumerate}
	Note that any norm $\|\ \|$ on~$V$ provides a metric that is defined by $d(v,w) = \| v - w \|$. Thus, the norm determines a topology on~$V$.
	
	\item A \defit[Banach!space]{Banach space} is a vector space~$\Banach$, together
	\nindex{$\Banach$ = a Banach space}
	with a norm~$\| \ \|$, such that the resulting metric is complete.
	 (Banach spaces are Fréchet.)
	
	\item An \defit[inner!product]{inner product} on~$V$ is a function $\langle \, \mid \, \rangle \colon V \times V \to \field$, such that
		\begin{enumerate}
		\item $\langle a v + b w \mid x \rangle = a \langle v | x \rangle + b \langle w | x \rangle$ for $a,b \in \field$ and $v,w,x \in V$,
		\item $\langle v \mid w \rangle = \overline{\langle w \mid v \rangle}$ for $v,w \in V$, where $\overline{a}$ denotes the complex conjugate of~$a$,
		and
		\item $\langle v \mid v \rangle \ge 0$ for all $v \in V$, with equality iff $v = 0$.
		\end{enumerate}
	Note that if $\langle \, \mid \, \rangle$ is an inner product on~$V$, then a norm on~$V$ is defined by the formula $\|v\| = \sqrt{\langle v \mid v \rangle}$.
	
	\item A \defit{Hilbert space} is a vector space~$\Hilbert$, 
		\nindex{$\Hilbert$ = a Hilbert space}
	together with an inner product $\langle \, \mid \, \rangle$, such that the resulting normed vector space is complete. (Hence, every Hilbert space is a Banach space.)
	
	\item An \defit[isomorphism!of Hilbert spaces]{isomorphism} between Hilbert spaces $\bigl( \Hilbert_1, \langle \, \mid \, \rangle_1 \bigr)$ and $\bigl( \Hilbert_2, \langle \, \mid \, \rangle_2 \bigr)$ is an invertible linear transformation~$T \colon \Hilbert_1 \to \Hilbert_2$, such that 
		$$ \text{$\langle Tv \mid Tw \rangle_2 = \langle v \mid w \rangle_1$ 
		\ for all $v,w \in \Hilbert_1$} .$$
		An isomorphism from~$\Hilbert$ to itself is called a \defit[unitary!operator]{unitary operator} on~$\Hilbert$.

	\end{enumerate}
\end{defns}

\begin{eg}
If $\mu$ is a measure on~$X$, then the $\LL{p}$-norm makes $\LL{p}(X,\mu)$ into a Banach space (for $1 \le p \le \infty$). Furthermore, $\LL2(X,\mu)$ is a Hilbert space, with the inner product
	$$ \langle \varphi \mid \psi \rangle = \int_X \varphi(x) \, \overline{\psi(x)} \, d\mu(x) .$$
\end{eg}

\begin{defns} \label{WeakStarDefn}
Let $\Banach$ be a Banach space (over $\field \in \{\real, \complex\}$).
	\begin{enumerate}
	\item A \defit[linear functional!continuous]{continuous linear functional} on~$\Banach$ is a continuous function $\lambda \colon \Banach \to \field$ that is linear (which means $\lambda (a v + b w) = a \lambda(v) + b \lambda(w)$ for $a,b \in \field$ and $v,w \in \Banach$).
	
	\item \nindex{$\Banach^*$ = dual of~$\Banach$}%
	$\Banach^* = \{\, \text{continuous linear functionals on~$\Banach$} \,\}$ is the \defit[dual of a Banach space]{dual} of~$\Banach$. This is a Banach space: the norm of a linear functional~$\lambda$ is
		$$ \| \lambda \| = \sup \{\, |\lambda(v)| \mid v \in \Banach, \ \|v\| \le 1 \,\} .$$
	
	\item For each $v \in \Banach$, there is a linear function $e_v \colon \Banach^* \to \field$, defined by $e_v(\lambda) = \lambda(v)$. The \defit{weak$^*$ topology} on $\Banach^*$ is the coarsest topology for which every $e_v$ is continuous. 
	
	In other words, the basic open sets in the weak$^*$ topology are of the form 
		$ %\mathcal{O}_v^{U} = 
		\{\, \lambda \in \Banach^* \mid \lambda(v) \in U \,\}$, for some $v \in \Banach$ and some open subset~$U$ of~$\field$. A set in $\Banach^*$ is open if and only if it is a union of sets that are finite intersections of basic open sets.

	\item Any continuous, linear transformation from~$\Banach$ to itself is called a \defit{bounded operator} on~$\Banach$.

	\item The set of bounded operators on~$\Banach$ is itself a Banach space, with the \defit{operator norm}
	\nindex{$\| T \|$ = operator norm}
	$$ \| T \| = \sup \{\, \|T(v)\| \mid \| v \| \le 1 \,\} .$$

	\end{enumerate}
\end{defns}

\begin{prop}[(\thmindex{Banach-Alaoglu}{Banach-Alaoglu Theorem})]  \label{BanachAlaogluThm}
If $\Banach$ is any Banach space, then the closed unit ball in $\Banach ^*$ is compact in the weak$^*$ topology.
%\hint{Let $D = B_1(\complex)$ be the unit disk in~$\complex$. There is a natural embedding of the closed unit ball $B_1(\Banach^*)$ in the infinite Cartesian product $\bigtimes_{v \in B_1(\Banach)} D$. This product is compact (by Tychonoff's Theorem~\pref{TychonoffThm}), and the map is a homeomorphism onto its image, which is closed.}
\end{prop}

\begin{prop}[(\thmindex{Hahn-Banach}Hahn-Banach Theorem)]
Suppose 
	\begin{itemize}
	\item $\Banach$ is a Banach space over~$\field$,
	\item $W$ is a subspace of~$\Banach$  \textup(not necessarily closed\/\textup),
	and
	\item $\lambda \colon W \to \field$ is linear.
	\end{itemize}
If\/ $|\lambda(w)| \le \|w\|$ for all $w \in W$, then $\lambda$~extends to a linear functional $\widehat \lambda \colon \Banach \to \field$, such that $|\widehat \lambda(v)| \le \|v\|$ for all $v \in \Banach$.
\end{prop}

\begin{prop}[(\thmindex{Open Mapping}Open Mapping Theorem)] \label{OpenMappingThm}
Assume $X$ and~$Y$ are Fréchet spaces, and $f \colon X \to Y$ is a continuous, linear map.
	\begin{enumerate}
	\item \label{OpenMappingThm-surj}
	If $f$ is surjective, and $\mathcal{O}$ is any open subset of~$X$, then $f(\mathcal{O})$ is open.
	\item \label{OpenMappingThm-bij}
	If $f$ is bijective, then the inverse $f^{-1} \colon Y \to X$ is continuous.
	\end{enumerate} 
\end{prop}

\begin{assump} \label{HilbertSpaceSeparable}
Hilbert spaces are always assumed to be separable.
\end{assump}

This has the following consequence:

\begin{prop}
There is only one infinite-dimensional Hilbert space\/ \textup(up to isomorphism\/\textup). In other words, every infinite-dimensional Hilbert space is isomorphic to $\LL2(\real, \mu)$, where $\mu$~is Lebesgue measure.
\end{prop}

\begin{defns} \ 
\noprelistbreak
	\begin{enumerate}
	\item If $\Hilbert_1$ and~$\Hilbert_2$ are Hilbert spaces, then the \defit[direct sum!of Hilbert spaces]{direct sum} $\Hilbert_1 \oplus \Hilbert_2$ is a Hilbert space, under the inner product 
	$$ \bigl\langle (\varphi_1,\varphi_2) \mid (\psi_1,\psi_2)  \bigr\rangle = \langle \varphi_1 \mid \psi_1 \rangle + \langle \varphi_2 \mid \psi_2 \rangle .$$
By induction, this determines the direct sum of any finite number of Hilbert spaces; see \cref{HilbertDirSumInfty} for the direct sum of infinitely many.

	\item We use 
	\nindex{$\perp$ = ``is orthogonal to''}%
	``$\perp$'' as an abbreviation for ``is orthogonal to\zz.'' Therefore, if $\varphi,\psi \in \Hilbert$, then $\varphi \perp \psi$ means $\langle \varphi \mid \psi \rangle = 0$.  For subspaces $\mathcal{K}, \mathcal{K}'$ of~$\Hilbert$, we write $\mathcal{K} \perp \mathcal{K}'$ if $\varphi \perp \varphi'$ for all $\varphi \in \mathcal{K}$ and $\varphi' \in \mathcal{K}'$.

	\item The \defit[orthogonal!complement]{orthogonal complement} of a subspace~$\mathcal{K}$ of~$\Hilbert$ is%
		\nindex{$\mathcal{K}^\perp$ = orthogonal complement of the subspace~$\mathcal{K}$}
	$$\mathcal{K}^\perp = \{\, \varphi \in \Hilbert \mid \varphi \perp \mathcal{K} \,\} .$$
	This is a closed subspace of~$\Hilbert$. We have $\Hilbert = \mathcal{K} + \mathcal{K}^\perp$ and $\mathcal{K} \perp \mathcal{K}^\perp$, so $\Hilbert = \mathcal{K} \oplus \mathcal{K}^\perp$.
	
	\item The \defit[orthogonal!projection]{orthogonal projection} onto a closed subspace~$\mathcal{K}$ of~$\Hilbert$ is the (unique) bounded operator $P \colon \Hilbert \to \mathcal{K}$, such that 
		\begin{itemize}
		\item $P(\varphi) = \varphi$ for all $\varphi \in \mathcal{K}$,
		and
		\item $P(\psi) = 0$ for all $\varphi \in \mathcal{K}^\perp$.
		\end{itemize}

	\end{enumerate}
\end{defns}

\begin{defns}
Let $T \colon \Hilbert \to \Hilbert$ be a bounded operator on a Hilbert space~$\Hilbert$.%
\noprelistbreak
	\begin{enumerate}
	\item The \defit[adjoint!of a linear transformation]{adjoint} of~$T$ is the bounded operator~$T^*$ on~$\Hilbert$, such that 
		$$ \text{$\langle T \varphi \mid \psi \rangle = \langle \varphi \mid T^* \psi \rangle$ \ for all $\varphi,\psi \in \Hilbert$} .$$
It does not always exist, but $T^*$ is unique if it does exist.
	
	\item $T$ is \defit[self-adjoint operator]{self-adjoint} (or \defit[Hermitian!operator]{Hermitian}) if $T = T^*$.
	
	\item $T$ is \defit[normal operator]{normal} if $T T^* = T^* T$.
	
	\item $T$ is \defit[compact!linear operator]{compact} if there is a nonempty, open subset~$\open$\, of~$\Hilbert$, such that $T(\open\,)$ is precompact.
	
	\end{enumerate}
\end{defns}

\begin{prop} \label{CpctOpBasics}
Let $T$ be a bounded operator on a Hilbert space~$\Hilbert$.
	\begin{enumerate}
	\item If $T(\Hilbert)$ is finite-dimensional, then $T$ is compact.
	\item The set of compact operators on~$\Hilbert$ is closed\/ \textup(in the topology defined by the operator norm\/\textup).
	\end{enumerate}
\end{prop}

\begin{prop}[(\thmindex{Spectral}Spectral Theorem)] \label{SpectralThm}
If $T$ is any bounded, normal operator on any Hilbert space~$\Hilbert$, then there exist
	\begin{itemize}
	\item a finite measure~$\mu$ on $[0,1]$, 
	\item a bounded, measurable function $f \colon [0,1] \to \complex$,
	and
	\item an isomorphism\/ $U \colon \Hilbert \to \LL2([0,1],\mu)$,
	\end{itemize}
such that $U( T \varphi ) = f \, U(\varphi)$, for all $\varphi \in \Hilbert$ \ \textup(where $f \, U(\varphi)$ denotes the pointwise multiplication of the functions $f$ and~$U(\varphi)$.
 %\in \LL2(X,\mu)$ is defined by $\bigl( f \, U(\varphi) \bigr)(x) = f(x) \, \bigl( \psi(v) \bigr)(x)$\textup).

Furthermore:
	\begin{enumerate} 
	\item $T$ is unitary if and only if $|f(x)| = 1$ for a.e.\ $x \in [0,1]$.
	\item $T$ is self-adjoint if and only if $f(x) \in \real$ for a.e.\ $x \in [0,1]$.
	\end{enumerate}
\end{prop}

\begin{defn} \label{SpectralMeasDefn}
In the situation of \cref{SpectralThm}, the \defit{spectral measure} of~$T$ is $f_*\mu$.
\end{defn}

\begin{cor}[(\thmindex{Spectral}Spectral Theorem for compact, self-adjoint operators)] \label{SpectralThmCpct}
Let $T$ be a bounded operator on any Hilbert space~$\Hilbert$. Then $T$ is both self-adjoint and compact if and only if there exists an orthonormal basis $\{e_n\}$ of~$\Hilbert$, such that
\noprelistbreak
	\begin{enumerate}
	\item each $e_n$ is an eigenvector of~$T$, with eigenvalue~$\lambda_n$,
	\item $\lambda_n \in \real$,
	and
	\item $\lim_{n \to \infty} \lambda_n = 0$.
	\end{enumerate}
\end{cor}

\begin{prop}[(\thmindex{Fréchet-Riesz}Fréchet-Riesz Theorem)]
If $\lambda$ is any continuous linear functional on a Hilbert space~$\Hilbert$, then there exists $\psi \in \Hilbert$, such that $\lambda(\varphi) = \langle \varphi \mid \psi \rangle$ for all $\varphi \in \Hilbert$.
\end{prop}

%\begin{notes}
%
%\end{notes}

 %!TEX root = IntroArithGrps.tex

\standassumpfalse

\mychapter{A Quick Look at \texorpdfstring{$S$}{S}-Arithmetic Groups}
\label{SarithChap}

%Much of the importance of arithmetic groups comes from being lattices in semisimple Lie groups.
Classically, and in the main text of this book, the Lie groups under consideration were manifolds over the field~$\real$ of real numbers. However, in some areas of modern mathematics, especially Number Theory and Geometric Group Theory, it is important to understand the lattices in Lie groups not only over the classical field~$\real$ (or~$\complex$), but also over ``nonarchimedean'' fields of $p$-adic numbers. The natural analogues of arithmetic groups in this setting are called ``$S$-arithmetic groups\zz.'' Roughly speaking, this generalization is obtained by replacing the ring~$\integer$ with a slightly larger ring.

\begin{defn}
For any finite set \nindex{$S$ = a finite set of prime numbers}$S = \{p_1,p_2,\ldots,p_n\}$ of prime numbers, let%
\nindex{$\integer_S = \integer \bigl[ 1/p_1, 1/p_2, \ldots, 1/p_n \bigr]$ is the ring of $S$-integers}
% no page break here !!!
	\begin{align*}
	\integer_S 
	&= \bigset{ \frac{p}{q} \in \rational }{ 
		\begin{matrix}
		\text{every prime factor} \\
		\text{of~$q$ is in~$S$} 
		\end{matrix}} 
	 =  \integer \bigl[ 1/p_1, 1/p_2, \ldots, 1/p_n \bigr]
	. \end{align*}
This is called the ring of \defit[S-@$S$-!integers]{$S$-integers}.
\end{defn}

\begin{eg} \label{SArithPrototype} \ 
\noprelistbreak
	\begin{enumerate}
	\item The prototypical example of an arithmetic group is $\SL(\ell, \integer)$.
	\item The corresponding example of an $S$-arithmetic group is $\SL \bigl( \ell, \integer_S \bigr)$ (where $S$ is a finite set of prime numbers).
	\end{enumerate}
That is, while arithmetic groups do not allow their matrix entries to have denominators, $S$-arithmetic groups allow their matrix entries to have denominators that are products of certain specified primes.
\end{eg}

Most of the results in this book can be generalized in a natural way to $S$-arithmetic groups. (The monographs \cite{MargulisBook} and~\cite{PlatonovRapinchukBook} treat $S$-arithmetic groups alongside arithmetic groups throughout.) We will now give a very brief description of these more general results.

\begin{rem}
The monograph \cite{MargulisBook} of Margulis deals with a more general class of $S$-arithmetic groups that allows $\real$ to be replaced with certain ``local'' fields of characteristic~$p$, but we discuss only the fields of characteristic~$0$. 
\end{rem}

\section{Introduction to \texorpdfstring{$S$}{S}-arithmetic groups} \label{IntroSArithSect}

Most of the theory in this book (and much of the importance of the theory of arithmetic groups) arises from the fundamental fact that $G_{\integer}$ is a lattice in~$G$. Since the ring $\integer_S$ is not discrete (unless $S = \emptyset$), the group $G_{\integer_S}$ is usually not discrete, so it is usually not a lattice in~$G$. Instead, it is a lattice in a group~$G_S$ that will be defined in this section. % \lcnamecref{IntroSArithSect}. @@@

The construction of~$\real$ as the completion of~$\rational$ can be generalized as follows: 

\begin{defn}[($p$-adic numbers)]
Let $p$ be a prime number. 
\noprelistbreak
	\begin{enumerate}
	\item If $x$ is any nonzero rational number, then there is a unique integer $v = v_p(x)$, such that we may write
		$$ x = p^v \, \frac{a}{b} ,$$
	where $a$ and~$b$ are relatively prime to~$p$. (We let $v_p(0) = \infty$.) Then $v_p(x)$ is called the \defit[p-adic@$p$-adic!valuation]{$p$-adic valuation} of~$x$.
%	The \defit[p-adic@$p$-adic!valuation]{$p$-adic valuation} $v_p(k)$ of a nonzero integer~$k$ is defined by the equation $k = p^{v_p(k)} k'$, where $k'$~is relatively prime to~$k$. This extends to the nonzero rational numbers by letting
%		$$ v_p(a/b) = v_p(a) - v_p(b) .$$
	\item Let
		$$ d_p(x,y) = p^{-v_p(x-y)} .$$
	It is easy to verify that $d_p$ is a metric on~$\rational$. It is called the \defit[p-adic@$p$-adic!metric]{$p$-adic metric}.
	\item Let \nindex{$\rational_p$ = field of $p$-adic numbers}$\rational_p$ be the completion of~$\rational$ with respect to this metric. (That is, $\rational_p$ is the set of equivalence classes of convergent Cauchy sequences.) This is a field that naturally contains~$\rational$. It is called the \defit[p-adic@$p$-adic!field]{field of $p$-adic numbers}.
	\item If $\GG$ is an algebraic group over~$\rational$, we can define the group $\GG(\rational_p)$ of $\rational_p$-points of~$\GG$.
	\end{enumerate}
\end{defn}

\begin{notation}
To discuss real numbers and $p$-adic numbers uniformly, it is helpful to let 
\nindex{$\rational_\infty = \real$}$\rational_\infty = \real$.
\end{notation}

The construction of arithmetic subgroups by restriction of scalars \csee{RestrictScalarsSect} is based on the fact that the ring~$\ints$ of integers in a number field~$F$ embeds as a cocompact, discrete subring in $\bigoplus_{v \in S_\infty} F_v$. Using this fact, it was shown that $\GG(\ints)$ is a lattice in $\bigtimes_{v \in S_{\infty}} \GG(F_v)$. 

Similarly, to obtain a lattice in a $p$-adic group $\GG(\integer_p)$, or, more generally, in a product $\bigtimes_{v \in S \cup \{\infty\}} \F_v$ of $p$-adic groups and real  groups, we note that 
	$$ \text{$\integer_S$ embeds as a cocompact, discrete subring in $\bigoplus_{p \in S \cup \{\infty\}} \rational_p$} .$$
Using this fact, it can be shown that 
	$$ \text{$\GG(\ints_S)$ is a lattice in $\GG_S = \bigtimes_{p \in S \cup \{\infty\}} \GG(\rational_p)$} . $$% no page break here !!!
\nindex{$\GG_S = \bigtimes_{p \in S \cup \{\infty\}} \GG(\rational_p)$}%
We call $\GG(\ints_S)$ an \defit[S-@$S$-!arithmetic]{$S$-arithmetic subgroup}.

\begin{eg}
Let $\GG$ be the special linear group $\SL_n$.
\noprelistbreak
	\begin{enumerate}
	\item Letting $S = \emptyset$, we have $\integer_S = \integer$ and $G_S = \SL(n,\real)$. So $\SL(n,\integer)$ is an $S$-arithmetic lattice in $\SL(n,\real)$. This is a special case of the fact that every arithmetic lattice is an $S$-arithmetic lattice (with $S = \emptyset$).
	\item Letting $S = \{p\}$, where $p$ is a prime, we see that $\SL \bigl( n, \integer[1/p] \bigr)$ is a lattice in $\SL(n,\real) \times \SL(n,\rational_p)$.
	\item More generally, letting $S =  \{p_1,p_2,\ldots,p_r\}$, where $p_1,\ldots,p_r$ are primes, we see that $\SL \bigl( n, \integer_S \bigr)$ is a lattice in 
		$$\SL(n,\real) \times \SL(n,\rational_{p_1}) \times \SL(n,\rational_{p_2}) \times \cdots  \times \SL(n,\rational_{p_r}) .$$
	\end{enumerate}
This is an elaboration of our previous comment that $\SL(\ell, \integer_S)$ is the prototypical example of an $S$-arithmetic group \csee{SArithPrototype}.
\end{eg}

\begin{rem}[{\cite[Chap.~7]{Brown-BuildingsBook}}]
In the study of arithmetic subgroups of a Lie group~$G$, the symmetric space $G/K$ is a very important tool.
In the theory of $S$-arithmetic subgroups of $\GG_S$, this role is taken over by a space called the \defit{Bruhat-Tits building} of~$\GG_S$.
It is a Cartesian product
	$$ X_S = (G/K) \times \bigtimes_{p \in S} X_p , $$
where $X_p$ is a contractible simplicial complex on which $\GG(\rational_p)$ acts properly (but not transitively).
\end{rem}

%\begingroup \smaller \baselineskip=10pt
\subsection*{Optional:}
Readers familiar with the basic facts of Algebraic Number Theory will realize that the above discussion has the following natural generalization:

\begin{defn}[({\cite[p.~61]{MargulisBook}, \cite[p.~267]{PlatonovRapinchukBook}})]
\label{SarithFDefn} 
Let
\noprelistbreak
	\begin{itemize}
	\item $\ints$ be the ring of integers of an algebraic number field~$F$, 
	\item $S$~be a finite set of finite places of~$F$, 
	and
	\item $\GG$ be a semisimple algebraic group over~$F$,
	and
	\item $\GG_S = \bigtimes_{v \in S \cup S_\infty} \GG(F_v)$.
	\end{itemize}
Then $\GG(\ints_S)$ is an \defit[S-@$S$-!arithmetic]{$S$-arithmetic subgroup} of~$\GG_S$.
\end{defn}

\begin{rem}
More generally, much as in \cref{ArithDefn}, if 
\noprelistbreak
	\begin{itemize}
	\item $\Gamma'$ is an $S$-arithmetic subgroup of~$\GG'_S$,
	and
	\item $\varphi \colon \GG'_S \to \GG_S$ is a surjective, continuous homomorphism, with compact kernel,
	\end{itemize}
then any subgroup of~$\GG_S$ that is commensurable to $\varphi(\Gamma')$ may be called an $S$-arithmetic subgroup of~$\GG_S$.
\end{rem}

\begin{thm}[{}{\cite[Thm.~5.7, p.~268]{PlatonovRapinchukBook}}]
Every $S$-arithmetic subgroup of\/~$\GG_S$ is a lattice in\/~$\GG_S$.
\end{thm}

%\endgroup

\section{List of results on \texorpdfstring{$S$}{S}-arithmetic groups}

\begin{warn}
The Standing Assumptions \pref{standassump} do \textbf{not} apply in this \lcnamecref{SarithChap}, because $G$ is not assumed to be a \emph{real} Lie group.
\end{warn}

Instead:

\begin{assump}
Throughout the remainder of this \lcnamecref{SarithChap}:
\noprelistbreak
	\begin{itemize}
	\item $\GG$ is a semisimple algebraic group over~$\rational$,
	\item $S$ is a finite set of prime numbers,
	and
	\item $\Gamma$ is an $S$-arithmetic lattice in $\GG_S = \bigtimes_{p \in S \cup \{\infty\}} \GG(\rational_p)$.
	\end{itemize}
To avoid trivialities, we assume $\GG_S$ is not compact, so $\Gamma$ is infinite.
\end{assump}

\begin{defn}
As a substitute for real rank in this setting, let%
\nindex{$\Srank \GG = \sum_{p \in S \cup \{\infty\}} \Qprank \GG$} % no page break here !!!
		$$\Srank \GG = \sum_{p \in S \cup \{\infty\}} \Qprank \GG .$$
\end{defn}

\begin{rem}
All of these results generalize to the setting of \cref{SarithFDefn}, but we restrict our discussion to~$\rational$ for simplicity.
\end{rem}

The following theorems on $S$-arithmetic groups are all stated without proof, but each result is provided with a reference for further reading. The reader should be aware that these references are almost always secondary sources, not the original appearance of the result in the literature.

%\subsection*{Definitions related to \cref{RrankChap} (Real Rank)} 
%\label{RrankForSArith}
%
%	\begin{itemize}
%	\sitem[\spref{TorusDefn}] Let $L$ be any field containing~$F$. A Zariski-closed, Zariski-connected subgroup~$T$ of $\GG$ is a \defit[torus subgroup]{torus} if $T$ is diagonalizable
%over the algebraic closure of~$L$; that is, if there exists $g \in
%\GL(n,\algclosure{L})$, such that $g^{-1} T g$ consists entirely of diagonal matrices.
%	\sitem[\spref{RsplitDefn}] For any $v \in S$, a torus~$T$ in~$\GG$ is \defit[Qp-split torus]{$F_v$-split}
%if $T$ is diagonalizable over~$F_v$.
% 	\sitem[\spref{RrankDefn}] For $v \in S$, 
% $\Fvrank(\GG)$ is the dimension of any maximal $F_v$-split
%torus of~$\GG$. (This does not depend on the choice of the
%maximal torus, because all maximal $F_v$-split tori
%of~$\GG$ are conjugate.)
%	\item As a substitute for $\Rrank$ in the setting of $S$-arithmetic groups, let 
%		$$\Srank G = \sum_{v \in S} \Fvrank(\GG) .$$
%	\end{itemize}

%\begin{notation}
%For convenience, we use $G$ as an abbreviation for $\GG_S$.
%\end{notation}

\subsection*{Results related to \cref{BasicLatticesChap} (Basic Properties of Lattices)} 

	\begin{slist}
	\sitem[\spref{G/GammaCpct<>NoAccPt}]
 $\Gamma \backslash \GG_S$ is compact if and only if the
identity element~$e$ is \textbf{not} an accumulation point
of $\Gamma^{\GG_S}$
\cite[Thm.~1.12, p.~22]{RaghunathanBook}.
	\sitem[\spref{GammaUnip->notcpct}]
 If $\Gamma$ has a nontrivial, unipotent element, then
$\Gamma \backslash \GG_S$ is not compact. In fact, Godement's Criterion \spref{GodementCriterion} tells us that the converse is also true.
	\sitem[\spref{BDT}] The Borel Density Theorem holds, for any continuous homomorphism $\rho \colon \GG_S \to \GL(V)$, where $V$ is a vector space over $\real$, $\complex$, or any $p$-adic field~$\rational_p$
\cite[Thm.~II.2.5 (and Lem.~II.2.3), p.~84]{MargulisBook}.
	\sitem[\spref{GammaFinPres}] $\Gamma$ is finitely presented 
	\cite[Thm.~5.11, p.~272]{PlatonovRapinchukBook}.
	\sitem[\spref{torsionfree}] (Selberg Lemma) $\Gamma$ has a torsion-free
subgroup of finite index
\cite[Thm.~6.11, p.~93]{RaghunathanBook}.
	\sitem[\spref{FreeInGamma}] (Tits Alternative) $\Gamma$ has a nonabelian
free subgroup
\cite[App.~B, pp.~351--353]{MargulisBook}.
	\end{slist}

\subsection*{Results related to \cref{ArithGrpsChap} (What is an Arithmetic Group?).}

	\begin{slist}
	\sitem[\spref{MargulisArith}] (Margulis Arithmeticity Theorem) If $\Srank \GG \ge 2$, then every irreducible lattice in~$\GG_S$ is $S$-arithmetic
\cite[Thm.~IX.1.10, p.~298, and Rem.~(vi) on p.~290]{MargulisBook}. 
(Note that our definition of irreducibility is stronger than the one used in \cite{MargulisBook}.)
	\sitem[\spref{GodementCriterion}] (Godement Criterion) $\Gamma \backslash \GG_S$ is compact if and only if $\Gamma$ has no nontrivial unipotent elements
\cite[Thm.~5.7(2), p.~268]{PlatonovRapinchukBook}.
	\end{slist}

\begin{rem*} \Cref{CpctOpenSubgrp->LattsCocpct} (easily) implies:
	\begin{enumerate}
	\item {\cite[Thm.~1]{Tamagawa-DiscSubgrpsPAdic}} If $v$~is any nonarchimedean place of~$F$, then every lattice in $\GG(F_v)$ is cocompact.
	\item If $\GG(S_\infty)$ is compact, then every lattice in~$G$ is cocompact.
	\end{enumerate}
\end{rem*}

\begin{warn*}
We know that if $\ints$ is the ring of integers of~$F$, then $\GG(\ints)$ embeds as an arithmetic lattice in $\bigtimes_{v \in S_\infty} \GG(F_v)$, but that restriction of scalars allows us to realize this same lattice as the $\integer$-points of an algebraic group defined over~$\rational$ \ccf{ResScal->Latt}. This means that all arithmetic groups can be found by using only algebraic groups that are defined over~$\rational$, not other number fields. It is important to realize that the same cannot be said for $S$-arithmetic groups: most extensions of~$\rational$ provide many $S$-arithmetic groups that cannot be obtained from~$\rational$. 

For example, suppose $p$~is a prime in~$\integer$, but $p$ factors in the integers~$\ints$ of an extension field, and $a$ is a prime factor of~$p$ in~$\ints$. Then the subgroup $\SL \bigl( 2, \ints[1/p] \bigr)$ can be obtained by restriction of scalars, but $\SL \bigl( 2, \ints[1/a] \bigr)$ is an $\{a\}$-arithmetic subgroup that cannot be obtained by this method.
\end{warn*}

\subsection*{A result related to \cref{AmenableChap} (Amenable Groups)}
	\begin{slist}
	\sitem[\spref{GNotAmen}]
For $v \in S$, if $\GG(F_v)$ is not compact, then $\GG(F_v)$ is not amenable
	 \cite[Rem.~8.7.11, p.~260]{ReiterStegeman-HarmAnalLocCpctGrps}.
	\end{slist}

\subsection*{Results related to \cref{KazhdanTChap} (Kazhdan's Property ($T$))}
	\begin{slist}
	\sitem[\spref{WhichGKazhdan}]
If $\Fvrank G \ge 2$, for every simple factor~$G$ of $\GG(F_v)$, and every $v \in S$, then
 $\GG_S$ has Kazhdan's property
 \cite[Cor.~III.5.4, p.~130]{MargulisBook}.
	\sitem[\spref{Kazhdan:G->Gamma}]
If $\GG_S$ has Kazhdan's property, then $\Gamma$ also has Kazhdan's property
\cite[Thm.~III.2.12, p.~117]{MargulisBook}.
	\sitem[\spref{KazhdanlatticeCor}]
 If $\Gamma$ has Kazhdan's property, then $\Gamma/[\Gamma ,\Gamma ]$ is finite
 \cite[Thm.~III.2.5, p.~115]{MargulisBook}.
	\end{slist}

\subsection*{A result related to \cref{MargulisSuperChap} (Margulis Superrigidity Theorem)}

\begin{assump*}
Assume 
	\begin{itemize}
	\item $\Srank G \ge 2$, 
	\item $\Gamma$ is irreducible,
	and
	\item $w$ is a place of some algebraic number field~$F'$.
	\end{itemize}
\end{assump*}

	\begin{slist}
	\sitem[\spref{MargSuperG'}] (Margulis Superrigidity Theorem
	\cite[Prop.~VII.5.3, p.~225]{MargulisBook})
	If
		\begin{itemize} \leftskip=1.25\parindent % !!!
		\item $\GG'$ is a Zariski-connected, noncompact, simple algebraic group over~$F'_w$, with trivial center, 
		and
		\item $\varphi \colon \Gamma \to \GG'(F'_w)$ is a homomorphism, such that $\varphi(\Gamma)$ is:
			\begin{itemize}  \leftskip = 2em % !!!
			\item Zariski dense in~$\GG'$,
			and \par
			\item not contained in any compact subgroup of~$\GG'(F'_w)$,
			\end{itemize}
		\end{itemize}
\leftskip=\sindent % slist lost its hangindent !!!
	then $\varphi$ extends to a continuous homomorphism $\widehat\varphi \colon \GG_S \to \GG'(F'_w)$. \par

\smallskip

Furthermore,  there is some $v \in S$, such that $F_v$ is isomorphic to a subfield of a finite extension of~$F'_w$.
	
%	\sitem[\spref{MargImg}]
%	If $\varphi \colon \Gamma \to \GL(n,F'_w)$ is any homomorphism, then the Zariski closure of $\varphi(\Gamma)$ is semisimple. 
%See \cite[Cor.~VII.6.18, pp.~249--150]{MargulisBook} for cocompact case

	\end{slist}

\begin{warn*}
\Cref{MargNoncpct} does not extend to the setting of $S$-arithmetic groups: for example,
the lattice $\SL(n,\integer)$ is not cocompact, but the image of the natural inclusion $\SL(n,\integer) \hookrightarrow \SL(n,\rational_p)$ is precompact.
\end{warn*}

\subsection*{Results related to \cref{NormalSubgroupChap} (Normal Subgroups of $\Gamma$)}

	\begin{slist}
	\sitem[\spref{MargNormalSubgrpsThm}] (Margulis Normal Subgroups Theorem \cite[Thm.~VIII.2.6, p.~265]{MargulisBook})
	Assume 
		\begin{itemize} \leftskip=1.25\parindent % !!!
		\item $\Srank \GG \ge 2$,
		\item $\Gamma$ is irreducible,
		and
		\item $N$ is a normal subgroup of~$\Gamma$.
		\end{itemize}
\leftskip=\sindent % slist lost its hangindent !!!
	Then either $N$ is finite, or $\Gamma/N$ is finite.

\leftskip=0pt % end of the hangindent problems !!!

	\sitem[\spref{Rrank1->GammaNotAlmSimple}]
	If $\Srank \GG = 1$, then $\Gamma$ has (many) normal subgroups~$N$, such that neither $N$ nor $\Gamma/N$ is finite
		 \cite[Cor.~7.6]{Lubotzky-LattRank1LocalFlds}.
	\end{slist}

\subsection*{Results related to \cref{ArithClassicalChap} (Arithmetic Subgroups of Classical Groups)}

	\begin{slist}
	
	\sitem[\spref{AlmostAllOverC}]
	Let $\overline{\rational_p}$ be the algebraic closure of~$\rational_p$. Then all but finitely many of the simple Lie groups over\/~$\overline{\rational_p}$ are isogenous to either\/ $\SL(n,\overline{\rational_p})$, $\SO(n,\overline{\rational_p})$, or\/ $\Sp(2n,\overline{\rational_p})$, for some~$n$ \cite[Thm.~11.4, pp.~57--58, and Thm.~18.4, p.~101]{Humphreys-LieAlg}.

	\sitem[\spref{RformsComplete}, \spref{ArithLattsAreClassical}] Every $\rational$-form or $\rational_p$-form of $\SL(n,\overline{\rational_p})$, $\SO(n,\overline{\rational_p})$, or $\Sp(n,\overline{\rational_p})$ is of classical type, except for some ``triality'' forms of $\SO(8,\overline{\rational_p})$ \ccf{D4weird} \cite[\S2.3]{PlatonovRapinchukBook}.

	\sitem[\spref{D4weird}] Unlike $\real$, the field $\rational_p$ has extensions of degree~$3$, so some $\rational_p$-forms of $\SO(8,\overline{\rational_p})$ are \term{triality} groups, even though there are no such $\real$-forms of $\SO(8,\complex)$.
	
	\sitem[\spref{GHasCpctLatt}]
	 $\GG_S$ has a cocompact, $S$-arithmetic lattice \cite{BorelHarder-exist}.

	\sitem[\spref{Isotypic->irred}]
	If $\GG$ is isotypic, then $\GG_S$ has a cocompact, irreducible lattice that is $S$-arithmetic \cite{BorelHarder-exist}.

	\sitem[\spref{Irred->Isotypic}]
	If $\GG_S$ has an irreducible, $S$-arithmetic lattice, then $\GG$ is isotypic.

	\end{slist}

\subsection*{A result related to \cref{ReductionChap} (Construction of a Coarse Fundamental Domain)}

 If $\fund$ is any coarse fundamental domain for $\GG(\integer)$ in $\GG(\real)$, then there is a compact subset~$C$ of $\bigtimes_{p \in S} \GG(\rational_p)$, such that $\fund \times C$ is a coarse fundamental domain for $\GG(\integer_S)$ in~$\GG_S$
\cite[Prop.~5.11, p.~267]{PlatonovRapinchukBook}.
% Let $\ints$ be the ring of integers in~$F$. If $\fund$ is any coarse fundamental domain for $\GG(\ints)$ in $\GG(S_\infty)$, then there is a compact subset~$C$ of $\bigtimes_{v \in S \smallsetminus S_\infty} \GG(F_v)$, such that $\fund \times C$ is a coarse fundamental domain for $\GG(\ints_S)$ in~$G$
%\cite[Prop.~5.11, p.~267]{PlatonovRapinchukBook}.

% L.\,Ji (arith grps: what, why, how) gives references, he gives:
% 
%A.Borel, J.P.Serre, Cohomologie dÕimmeubles et groups S-arithm«etiques, Topology 15
%(1976), 211-232.
% 
%A.Borel, Some finiteness properties of adele groups over number fields, Inst. Hautes «Etudes
%Sci. Publ. Math. 16 (1963) 5-30.
% 
%A.Borel, Arithmetic properties of linear algebraic groups, 1963 Proc. Internat. Congr.
%Mathematicians (Stockholm, 1962) pp.~10-22.

This implies that every $S$-arithmetic subgroup of~$\GG_S$ is a lattice \cite[Thm.~5.7, p.~268]{PlatonovRapinchukBook}, but the short proof outlined in \cref{SLNZISLATTSlickSect} does not seem to generalize to this setting.

\subsection*{Results related to \cref{RatnerChap} (Ratner's Theorems on Unipotent Flows)}
Ratner's three main theorems (\ref{Ratner-OrbitClosure}, \ref{Ratner-Equidistribution}, and \ref{Ratner-MeasClass})
%Orbit-Closure Theorem \pref{Ratner-OrbitClosure}, Equidistribution Theorem \pref{Ratner-Equidistribution}, and Classification of Invariant Measures \pref{Ratner-MeasClass} 
have all been generalized to the $S$-arithmetic setting by Ratner \cite{Ratner-CartProd,Ratner-Sarith} and Margulis-Tomanov \cite{MargulisTomanov-LocField,MargulisTomanov-AlmLinear} (independently).

\endgroup

\ifwantindex

% Notation Index
 \NotationIndex

% Index
\printindex

% List of Theorems
\begingroup \makeatletter
\renewcommand{\indexspace}{\par\vskip1.4\medskipamount}
\renewcommand{\indexname}{List of Named Theorems}
\@input@{thmlist.ind}
\endgroup
 
  \fi

\end{document}